\documentclass[lang=en,11pt]{elegantbook}

\title{~~~~~~~~~~~~~~~~~~~~~~~~Intercenter Geometry}
\subtitle{~~~~~~~~~~~~~~~~~~~~~~~~~~~~~~~~~~~~~~~~~~~~~~~~~~~~~~~~~~Daiyuan Zhang}

\date{\textbf{April 30, 2024}}
\version{\textbf{11. Download website of Chinese version: https://www.researchgate.net/profile/Daiyuan-Zhang}}

\bioinfo{\textbf{Features:}}{\textbf{innovation, unique ideas, unified approach, solved some problems that were not solved before!}}

\extrainfo{Welcome to \textbf{Intercenter Geometry!} \\ \textbf{Intercenter Geometry} is created by Daiyuan Zhang.\\ The copyright belongs exclusively to Daiyuan Zhang. \\Plagiarism is strictly prohibited!}

\cover{ZdyGongShi.png}   

\usepackage{array}
\usepackage{halloweenmath}

\allowdisplaybreaks[4] 

\usepackage{longtable} 
\usepackage{appendix}  

\begin{document}
\maketitle

\frontmatter
\chapter*{Abstract}
 \markboth{Abstract}{Abstract}
 \addcontentsline{toc}{chapter}{Abstract}

The author proposes a new geometry in this book. The author named this new geometry Intercenter Geometry. Intercenter Geometry is different from traditional Euclidean geometry and analytic geometry (coordinate geometry). The idea of Intercenter Geometry is that the geometric quantities on a plane will be expressed by the lengths of the three sides of a given triangle. The geometric quantities in space will be expressed by the lengths of the six edges of a given tetrahedron. In Intercenter Geometry, a unified approach is used to deal with the calculation of geometric quantities of triangles and tetrahedrons. Many new theorems and formulas on geometry and geometric inequalities have been obtained. In particular, Intercenter Geometry has solved some problems that have not been solved in Euclidean geometry and analytical geometry (coordinate geometry) for a long time, such as the distance between the centroid and the incenter of a given tetrahedron, etc. Intercenter Geometry enriches geometry, which is not only of theoretical significance to open up the field of mathematical research, to explore new ideas and methods of mathematical research, but also of positive significance to promote the spirit of innovation. The fruitful achievements of Intercenter Geometry also have practical application value.


%
\tableofcontents

\mainmatter
\chapter*{Preface}
\markboth{Preface}{Preface}
\addcontentsline{toc}{chapter}{Preface}

%
%
%

The theories and methods in this book are my original creation, and I name this set of theories and methods Intercenter Geometry.


Intercenter Geometry is divided into Plane Intercenter Geometry and Space Intercenter Geometry. The first 15 chapters of this book introduce Plane Intercenter Geometry, and Space Intercenter Geometry is introduced from Chapter 16.

\section{Update information of this version}
The English version of this book has been published on arXiv. The download website is: https://arxiv.org/abs/2107.08388.

%

The download website for the Chinese version is: https://www.researchgate.net/profile/Daiyuan-Zhang.


This edition has added some content on the basis of the tenth edition, and the main added content can be found in Chapter \ref{Ch20} (Frame components on tetrahedral frame)


%
%
%



Readers are strongly advised to read the latest version.

%
%

\section{The idea of Intercenter Geometry}
%
%

The idea of Intercenter Geometry is that the geometric quantities on a plane can be expressed by the lengths of the three sides of a given triangle. The geometric quantities in space can be expressed by the lengths of the six edges of a given tetrahedron. 

This is a brand new idea. I guess that perhaps no one else on this planet has studied geometry like me so far.

The new idea has yielded new achievements. Intercenter Geometry uses vector method to solve the calculation problem of the distance between two points in plane and space. Specifically, Intercenter Geometry uses vector method to solve the calculation problem of geometric quantities related to two points in triangle and tetrahedron, and has obtained many innovative achievements. I dare to guess that many new theorems and formulas in this book cannot be obtained by traditional methods.

\section{Why did I propose Intercenter Geometry?}

\subsection{Thanks to my daughter}
Why should I be original and create Intercenter Geometry? Why should I spend my time studying Intercenter Geometry? It started with tutoring my daughter. In order to help my daughter learn mathematics, I carefully read and studied the relevant materials of contemporary school mathematics. In the process of reading and research, we found that there were more problems in geometry. The geometry of contemporary school basically introduces Euclidean geometry, but whether it is plane geometry of junior high school stage or solid geometry of senior high school stage, there are always some shortcomings. Euclidean geometry is based on a set of axioms, but the solution of many problems is quite individualized (one problem with one strategy), which requires high skills and lacks a general unified method. Although there is no universal solution to all mathematical problems in our global village, the pursuit of a general solution to a certain kind of problems has always been the goal of mathematicians. Therefore, it is particularly important to find a general solution to a certain kind of problems. For example, when discussing triangles or tetrahedrons (triangular pyramids), the textbook will involve the concepts of centroid, incenter, circumcenter, etc., however, it was difficult to get a satisfactory answer when asking a question like how to solve the vector or the distance between the centroid and the incenter with a unified point of view. I have been asking myself, Is there a unified way to deal with such problems? I checked a lot of data, but did not get satisfactory results, so I calmed down and studied myself. Although it is impossible to find a unified method to deal with all the problems of triangles and tetrahedrons, after a period of thinking and research, I still found a unified theoretical method, which is good to be able to deal with problems such as the vector and the distance between the centroid and the incenter for a given triangle or a given tetrahedron, as well as some other problems in the calculation of geometric quantities. This is the Intercenter Geometry to be introduced in this book. Not only that, I have also obtained some new theorems of geometry and constructed some new inequalities. Therefore, this book is fruitful.

I would like to dedicate this book to my dear daughter. If it wasn't for tutoring my daughter, I probably wouldn't have studied this problem. Thank you, my dear daughter!

\subsection{Advantages and disadvantages of Euclidean geometry}

Euclidean geometry belongs to the traditional synthetic geometry. This kind of geometry establishes some axioms, and then obtains theorems or conclusions by logical reasoning. 

The biggest advantage of Euclidean geometry is that its results are often expressed by “natural” parameters such as length, angle, area and volume, so the results are clear both in form and geometry. 

Unfortunately, very few formulae (such as Helen's formula, etc.) can be derived from Euclidean geometry. Firstly, this is because Euclidean geometry has no connection with algebra. Point is the most basic concept in geometry, but in Euclidean geometry, there is no corresponding relationship between a point in space and a set of real numbers, so it is difficult to use algebraic theories and methods based on real number operations to deal with some problems, such as distance, area and so on. The inability to establish an internal connection with algebra is the “weakness” of Euclidean geometry, or the “natural defect” of Euclidean geometry. In fact, for complex geometric problems, abstract algebraic operations must be introduced. Euclidean geometry does not give the calculation formula of the distance between two points, which makes it very difficult to solve the geometric problems related to the distance between two points. For example, how to calculate the distance from the vertex of a triangle to the incenter or orthocenter? How to calculate the distance between the centroid and the incenter of a tetrahedron? These problems are often helpless and very difficult by using Euclidean geometry. The distance between two points is the most basic quantity of geometry and the basis of many other geometric quantities. This also means that Euclidean geometry would encounter great challenges in solving the geometric quantities related to distance. Secondly, the method of solving problems in Euclidean geometry has great limitations. Essentially, Euclidean geometry solves geometric problems directly from intuitive graphics. Specifically, it is to use the positional relationship between the elements of the intuitive graphics to carry out appropriate translation, rotation, expansion and other transformations, introduce auxiliary lines or auxiliary surfaces, intuitively display or calculate the geometric quantities of some graphics  through congruence and similarity, and then deduce the relationship between geometric quantities. However, for some problems, the relationship between geometric quantities may be hidden deeply and the calculation is very complex. For these complex geometric quantities, it is difficult to find their internal relationships only by virtue of intuitive graphics, so it is difficult to get results. From the point of view of the ideological method, it has a great limitation to solve geometric problems with intuitive graphics. Up to now, human beings do not know which problems can be solved by Euclidean geometry and which problems cannot be solved by Euclidean geometry. In addition, in the face of specific problems, how should the translation, rotation, scaling and other transformations be done? How to introduce auxiliary lines or auxiliary surfaces? Euclidean geometry does not provide an effective unified method. It often solves problems by experience. In many cases, it is one solving method for one problem, which is difficult to form a unified problem-solving strategy. 

\subsection{Advantages and disadvantages of analytic geometry}
Analytic geometry introduces the coordinates of points, establishes the corresponding relationship between the basic elements of space “points” and real numbers, so that we can use algebraic methods to study geometry. As a result of the use of the algebraic method, some problems in geometry that were difficult to calculate in Euclidean geometry can be solved.

Using analytic geometry method, first of all, we need to establish the coordinate system, and the establishment of coordinate system is greatly influenced by human factors. For different coordinate systems, the coordinates of points are different. Even in the same kind of coordinate system, the different position of coordinate origin will lead to the different coordinates of space points. Compared with Euclidean geometry, analytic geometry is more convenient to calculate some geometric quantities, but because the basic parameters of analytic geometry are the coordinates of points, its calculation results often only include the coordinates of points, and do not include some basic geometric quantities such as length, area and so on, which makes the calculation results of analytic geometry not “natural”. For example, when calculating the area of a triangle, the analytic geometry method usually establishes a coordinate system, and then selects the vertex coordinates of the triangle as variables (parameters) for calculation. The result is also a function of the coordinates of three points, which seems not directly related to the lengths of the three sides of the triangle. It is difficult to imagine that the function of the coordinates of these three points will be as graceful and symmetrical as the formula I proposed (see theorem \ref{thm:Thm1.5.1}) and Helen's formula. This is because the independent variables of theorem \ref{thm:Thm1.5.1} and Helen's formula are the “natural” quantities such as the lengths of the three sides of the triangle, rather than the coordinates of space points. Although the results of analytic geometry can be transformed into geometric quantities expressed by side lengths, areas or other parameters by algebraic operations, it usually needs complex and difficult algebraic operations. In particular, if there is no hint of known results (such as the known Helen formula, etc.), algebraic operations will appear to be at a loss and have no goal. In a word, the calculation results of analytic geometry usually can not see the internal relationship between some geometric quantities expressed by “natural” parameters such as side lengths or areas, and it is more difficult to understand the symmetry and beauty among those “natural” quantities. The geometric quantities calculated by analytic geometry seem to be covered with a layer of artificial mystery, covering up the original beautiful and charming “true face” of what she looks like.

In a word, Euclidean geometry obtains the conclusion through the transformation of intuitive graphics and the logical deduction of basic principles. Analytic geometry is based on the establishment of coordinate system and the use of algebraic methods to obtain the results. They have their own advantages and disadvantages.

Then, is there a kind of geometry that not only retains the advantages of Euclidean geometry and analytic geometry, but also overcomes their disadvantages? The Intercenter Geometry I founded is an attempt.

\section{Innovative achievements of Intercenter Geometry}
The prospect of a new theoretical method depends not only on whether the new theoretical method is convenient and effective in solving traditional problems, but also on whether some new results can be obtained. The latter is often more important. If a new theoretical method can get some results that can not be obtained or difficult to obtain by traditional methods, then the new theoretical method will show its value and vitality, and will be very attractive. The Intercenter Geometry I proposed should be full of charm.


The Intercenter Geometry that I created has yielded a lot of innovations.

\subsection{Intercenter Geometry solves the problems difficult to be solved by traditional geometry}
The distance between two points in a plane is expressed by the lengths of three sides of a given triangle. The distance between two points in space is expressed by the lengths of the six edges of a given tetrahedron.

Intercenter Geometry solves many problems that cannot be solved by traditional geometry (such as Euclidean geometry, analytical geometry, differential geometry, etc.).

Here is just one example to illustrate the innovations in Intercenter Geometry. The following example is just the tip of the iceberg in this book's vast collection of innovative achievements, and there are a large number of new theorems and new formulas in this book. If you are interested, please read this book carefully, and I believe you will be surprised.

Now let's take a look at an innovative conclusion in this book. Given a tetrahedron, the square of the distance between the centroid and the incenter of the tetrahedron is (see theorem \ref{thm:Thm26.2.1}):
\[G{{I}^{2}}=-\frac{1}{{{S}^{2}}}\left( \begin{aligned}
	& \left( {{S}^{A}}-\bar{S} \right)\left( {{S}^{B}}-\bar{S} \right)A{{B}^{2}}+\left( {{S}^{A}}-\bar{S} \right)\left( {{S}^{C}}-\bar{S} \right)A{{C}^{2}} \\ 
	& +\left( {{S}^{A}}-\bar{S} \right)\left( {{S}^{D}}-\bar{S} \right)A{{D}^{2}}+\left( {{S}^{B}}-\bar{S} \right)\left( {{S}^{C}}-\bar{S} \right)B{{C}^{2}} \\ 
	& +\left( {{S}^{B}}-\bar{S} \right)\left( {{S}^{D}}-\bar{S} \right)B{{D}^{2}}+\left( {{S}^{C}}-\bar{S} \right)\left( {{S}^{D}}-\bar{S} \right)C{{D}^{2}}  
\end{aligned} \right).\]

Where ${{S}^{A}}$, ${{S}^{B}}$, ${{S}^{C}}$, ${{S}^{D}}$ are the area of the triangle opposite to the four vertices $A$, $B$, $C$, $D$ of the tetrahedron $ABCD$ respectively; $S$ is the surface area of tetrahedron $ABCD$; $AB$, $AC$, $AD$, $BC$, $BD$, $CD$ are the lengths of six edges of tetrahedron $ABCD$; $\Bar {S}$ is the average surface area of tetrahedron $ABCD$, i.e.
\[\bar{S}=\frac{1}{4}S=\frac{1}{4}\left( {{S}^{A}}+{{S}^{B}}+{{S}^{C}}+{{S}^{D}} \right).\]

The expression of $G {{I}^{2}}$ shows that the distance between the centroid and the incenter of a tetrahedron can be expressed by the areas of four triangular faces and the lengths of six edges of the tetrahedron. These quantities are the “natural” parameters of the tetrahedron and do not contain any information related to coordinates. If Euclidean geometry method is used, can you prove the above results? If the analytic geometry method is used to establish a coordinate system, assume the coordinates of four vertices, and then calculate the coordinates of the centroid and the incenter (which are functions of the coordinates of four vertices respectively), and finally use the distance formula to calculate the distance between the centroid and the incenter, can you get the above formula of $G{{I}^{2}}$?

Without the hint of the above formula, it is very difficult to obtain the above result by Euclidean geometry or analytical geometry, or it may be impossible to obtain the above result by by Euclidean geometry or analytical geometry at all.

On the contrary, the results of Intercenter Geometry can help other geometries. For example, given the Cartesian coordinates of the four vertices of a tetrahedron, it is easy to transform the above formula of $G{{I}^{2}}$ into the Cartesian coordinate representation (that is, in the form of analytic geometry).

Is this formula for calculating $G{{I}^{2}}$ symmetrical and graceful? Have you seen this formula before? I am sure that this is a new formula that I put forward, and I publish it here for the first time. If you have seen this formula before, please tell me. I will certainly revise the discussion here, but I guess no one can tell me in this world.

\subsection{Intercenter Geometry can create a new class of geometric inequalities}
Intercenter Geometry is good at calculating the distance between two points. Many new geometric inequalities can be created based on the principle that the distance between two points is non-negative. It can be said that Intercenter Geometry is a “generator” that creates a large class of geometric inequalities.

I'll just illustrate it with an example in the following.

If the following notation is used to represent some geometric quantities of $\triangle ABC$:
\begin{flalign*}
	p=\frac{1}{2}\left( a+b+c \right),
	W=\sin 2A+\sin 2B+\sin 2C.	
\end{flalign*}
then, I proved the following inequality(see corollary \ref{cor:Cor15.2.2}): 
\[\begin{aligned}
	& \left( bW-2p\sin 2B \right)\left( cW-2p\sin 2C \right){{a}^{2}} \\ 
	& +\left( cW-2p\sin 2C \right)\left( aW-2p\sin 2A \right){{b}^{2}} \\ 
	& +\left( aW-2p\sin 2A \right)\left( bW-2p\sin 2B \right){{c}^{2}}\le 0. \\ 
\end{aligned}\]	
and, if and only if $a=b=c$, the equal sign holds.

Have you seen this inequality before? You may be surprised, how are these results obtained? Did the author guess? I want to tell you, I can't guess. I'm not so good at it. These results are direct corollaries of some theorems in Intercenter Geometry. Intercenter Geometry can obtain results unknown before, and can construct new inequalities. You can understand the mystery after reading this book. Of course, you can also try to solve them using Euclidean geometry or analytic geometry to see if you can get the above results.

\section{Research methods of Intercenter Geometry}

In order to develop the advantages of Euclidean geometry and analytic geometry and overcome their disadvantages, it is necessary to establish a coordinate system (I call it frame in this book). This coordinate system can establish a corresponding relationship between a point on a plane (or space) and a set of real numbers. At the same time, it can express the coordinates (I call it frame component in this book) as “natural parameters”, that is to say, the coordinates are expressed by the parameters that people's senses can touch (such as side length, area, etc.).

The research method of Intercenter Geometry is different from Euclidean geometry and analytic geometry. By establishing triangular frame (Plane Intercenter Geometry) or tetrahedral frame (Space Intercenter Geometry), all frame components (coordinates) are independent of the origin of the frame. This book takes the intersecting center, intersecting ratio and frame component as the main clues, takes the vector as the analysis tool, and uses a unified method to deal with the calculation of some geometric quantities of triangles and tetrahedrons, and obtains many new results, some of which are not easily obtained by traditional methods. 

A journey of a thousand miles begins with the first step, and any original research should be based on the most basic elements. As we all know, triangle and tetrahedron are the most basic figures of plane and space respectively. This book is divided into Plane Intercenter Geometry and Space Intercenter Geometry. It studies triangles and tetrahedrons respectively. They begins with the study of triangles.


Given a triangle on the plane, it is denoted as $\triangle ABC$. Select a point $P$ on the plane where $\triangle ABC $ is located, and give any point $O$ (called the frame origin) in space, then the necessary and sufficient condition for the four points $A$, $B$, $C$ and $P$ to be coplanar is that there is a unique set of real numbers $\alpha_{A}^{P}$, $\alpha_{B}^{P}$ and $\alpha_{C}^{P}$ (called the frame components of point $P$) so that:

\[\overrightarrow{OP}=\alpha _{A}^{P}\overrightarrow{OA}+\alpha _{B}^{P}\overrightarrow{OB}+\alpha _{C}^{P}\overrightarrow{OC},\]
\[\alpha _{A}^{P}+\alpha _{B}^{P}+\alpha _{C}^{P}=1.\]


The above is a well-known conclusion. I have made an in-depth study of this well-known conclusion, found the secret and put forward Intercenter Geometry.


The modulus of a vector represents the distance between two points, and the distance between two points is the basis of other quantities in geometry. In order to achieve my goal,  the frame components $\alpha_{A}^{P}$, $\alpha_{B}^{P}$ and $\alpha_{C}^{P}$ have to be expressed by “natural” parameters such as side length, so as to ensure that $\left| \overrightarrow{OP} \right|$ can also be expressed by “natural” parameters such as side length. Since the side length of a triangle is independent of the origin $O$ of the frame, it is also expected that the frame components $\alpha_{A}^{P}$, $\alpha_{B}^{P}$ and $\alpha_{C}^{P}$ are also independent of the origin of the frame $O$. How to find the frame components $\alpha_{A}^{P}$, $\alpha_{B}^{P}$ and $\alpha_{C}^{P}$ that meets these conditions has become the key problem of Plane Intercenter Geometry. My research shows that the frame components $\alpha_{A}^{P}$, $\alpha_{B}^{P}$ and $\alpha_{C}^{P}$ can be expressed by intersecting ratio (see Chapter \ref{Ch4}).

Given a triangle, three straight lines are drawn from each of the three vertices respectively. The lines drawn from each vertex intersect the opposite side of the vertex, and the three lines intersect at one point. This intersection is the “intersecting center” of three straight lines. I named this intersection “intersecting center”, which is the reason why I named this geometry “Intercenter Geometry”. Since the line drawn from each vertex intersects with the opposite side of the vertex, the side of the triangle will be divided into two parts. Taking the ratio of the two parts as a parameter, the three sides of the triangle will have three such parameters, and only two of them are independent. Obviously, this group of parameters (ratios) can carry out algebraic operations, and this group of parameters and the intersecting center on the plane form a corresponding relationship. Therefore, Intercenter Geometry has the “gene” of analytic geometry, which can be studied by algebraic methods. On the other hand, the ratio of the sides of the triangle may be related to the lengths of the three sides of the triangle. Therefore, the results of Intercenter Geometry can be expressed by the “natural” parameter such as the length, angle and area. Therefore, Intercenter Geometry inherits the advantages of Euclidean geometry.

Obviously, the quantities expressed by “natural” parameters such as side length, angle and area can be easily transformed into the quantities expressed by the coordinates of vertices. Therefore, the results of Intercenter Geometry can be easily transformed into the results of analytic geometry. However, the results of analytic geometry are usually difficult to be transformed into the results of Intercenter Geometry.


There is a similar method for the tetrahedron of space.


I have proved the following theorem (see theorem \ref{thm:Thm18.2.1}):


Suppose that given a tetrahedron $ABCD$, the point $O$ is an arbitrary point, and the point $P\in {{\pi }_{ABCD}}$, then the vector $\overrightarrow{OP}$ can be uniquely and linearly expressed in the following form:
\begin{equation*}\label{Eq18.2.1}
	\overrightarrow{OP}=\beta _{A}^{P}\overrightarrow{OA}+\beta _{B}^{P}\overrightarrow{OB}+\beta _{C}^{P}\overrightarrow{OC}+\beta _{D}^{P}\overrightarrow{OD},
\end{equation*}
\begin{equation*}\label{Eq18.2.2}
	\beta _{A}^{P}+\beta _{B}^{P}+\beta _{C}^{P}+\beta _{D}^{P}=1.
\end{equation*}


Through in-depth research, I found that the frame components $\beta_{A}^{P}$, $\beta_{B}^{P}$, $\beta_{C}^{P}$ and $\beta_{D}^{P}$ have nothing to do with the frame origin $O$ and are constant (when the point $P$ is given). In fact, if the point $P$ is selected on some special points of the tetrahedron $ABCD$ (such as centroid, incenter, circumcenter, etc.), its frame components $\beta_{A}^{P}$, $\beta _ {B}^{P}$, $\beta_{C}^{P}$ and $\beta_{D}^{P}$ can be expressed by the lengths of six edges of tetrahedron $ABCD$. In fact, there are infinitely many such points, and their frame components can be expressed by the lengths of six edges of the tetrahedron $ABCD$.

\section{Writing style of this book}
The writing style of this book follows the following points.

1. Follow the principle from simplicity to complexity, starting from the basic concepts, readers are gradually guided to understand Intercenter Geometry.

2.The author tries his best to make the discussion clear and the argument rigorous, which is also the author's consistent style. Explain in a plain and understandable language and combine figures and text as much as possible, so as to reduce the misunderstanding of the readers and make the book easier for them to read.


3.Explanation of theorem or corollary. The author explained all the theorems and corollaries in this book. The explanation consists of two parts, one is the name given to the theorem, and the other is the name of the person who first proposed the theorem. For example, in the description of the theorem \ref{thm:Thm6.2.1}, there are words such as “Theorem of vector of two intersecting centers, Daiyuan Zhang”, where “Theorem of vector of two intersecting centers” is the name given by the author to the theorem \ref{thm:Thm6.2.1}; The "Daiyuan Zhang" after the comma is the name of a person, indicating that the first person who proposed the theorem is Daiyuan Zhang.

To name a theorem is to make it easy for readers to read and search. Obviously, the “Theorem of vector of two intersecting centers” is much easier to associate with the content of the theorem than “theorem \ref{thm:Thm6.2.1}”, and is also easier to remember.

%

Since I am the founder of Intercenter Geometry, and this book is my original, many new theorems and corollaries in the book are proposed and proved by me, so my name (Daiyuan Zhang) is written in the explanation of these theorems and corollaries. Of course, there are also some conclusions omitting names that can be easily deduce. I have pressure to write my name in the explanation of theorem and corollary, which makes me dare not neglect. I need to make every effort to ensure the correctness of these theorems and deductions. However, due to my limited ability and the large number of theorems and corollaries, it is difficult to ensure that there is no mistake at all. I hope those theorems, corollaries or formulas related to my name can stand the test. Practice is the only criterion for testing truth. Readers should verify these theorems, corollaries and formulas by themselves to ensure the reliability of use.

The theorems or corollaries with my name (Daiyuan Zhang) and related formulas are independently proved by me. I (Daiyuan Zhang) have also independently proved the theorems or corollaries with other people's names (that is, the results that have existed in history). For example, many historical documents have given the proof of Ceva's theorem, but the proof method I use in this book is different from the existing historical documents. I use the method of Intercenter Geometry to prove Ceva's theorem (see theorem  \ref{thm:Thm10.1.1}). In short, all the theorems, corollaries and related formulas in this book are independently proved and derived by myself (Daiyuan Zhang).

In order to complete the content of this book, a small amount of basic knowledge is also introduced in the previous chapters. For some relatively basic theorems, if a person's name is not given in the description of the theorem, which means that the author is not sure who originated the theorem.

I believe in myself and respect history. I have repeatedly considered the names of people written in the explanation of theorem and corollary. However, due to my limited knowledge and energy, it is impossible for me to search all the historical documents. I can only tell you very honestly that all the theorems and corollaries with my name are my original creations, and I have never seen the same or similar discussion in any other literature. If you have any doubts about the origin of any theorem or corollary written with my name (Daiyuan Zhang), or you think that someone in history has reached the same conclusion, please do not hesitate to give your advice and present your evidence, thank you!

%
%

\section{Background knowledge and references that must be understood when reading this book}

Any reader with the foundation of Euclidean geometry and vector algebra can understand the main content of this book. Of course, if you have a certain understanding of linear algebra, mathematical analysis, you will better understand the content of this book.

Intercenter Geometry is my first creation. I can't find other similar references except this book. As for some basic knowledge needed, such as Euclidean geometry, vector algebra, linear algebra and mathematical analysis, etc., such references are very numerous and easy to obtain. Readers can choose appropriate references according to their own needs, which are not listed in this book.


%
%

\section{On the New Concept in Intercenter Geometry}
Intercenter Geometry is a new word of English I created, and I also created the corresponding Chinese word.

In order to describe Intercenter Geometry, I introduced many new concepts in the book, such as intersecting center, intersecting ratio, intersecting center of tetrahedron, intersecting center of face, frame component, etc. These new concepts are introduced in detail in the book. 



%
%

\section{Summary and prospect}
New ideas have produced new results. This book is my original work and summarizes my latest research achievements. If there are new results in the future, I will release a new version.

Intercenter Geometry has unique advantages in dealing with the calculation of geometric quantities of triangle and tetrahedron, but how to extend it to the calculation of geometric quantities of other graphics requires further in-depth study. If you have any good suggestions, please don't hesitate to give me advice, I'm here to express my gratitude in advance! My two contacts are:

dyzhang@njupt.edu.cn,  zhangdaiyuan2011@sina.com

%
%
%
%
%

\section{Disclaimer}
All other individuals and organizations other than me (Daiyuan Zhang) are hereinafter referred to as “others”; “This book” refers to the Intercenter Geometry published here.

This book is for academic exchange only.

All the theorems, corollaries and formulas in this book have been independently proved by me (Daiyuan Zhang). Although I have repeatedly verified them, errors inevitably occur due to the large number of them.

I (Daiyuan Zhang) don't make any promises to others.

Others need to verify the applicability of anything in this book (e.g. theorems, formulas, etc.).

I hereby declare that all harmful consequences caused by the use of any content in this book by others are the responsibility of others themselves and that 
I don't bear any responsibility for the result of using this book.

\vspace{\baselineskip}
\begin{flushright}\hfill {\it Daiyuan Zhang}\\
Time and place of the first version: July 18, 2021 in Nanjing, China.

The update time and place of this version: April 30, 2024 in Nanjing, China.
\end{flushright}

\begin{figure}[h]
	\centering
	\includegraphics[totalheight=8cm]{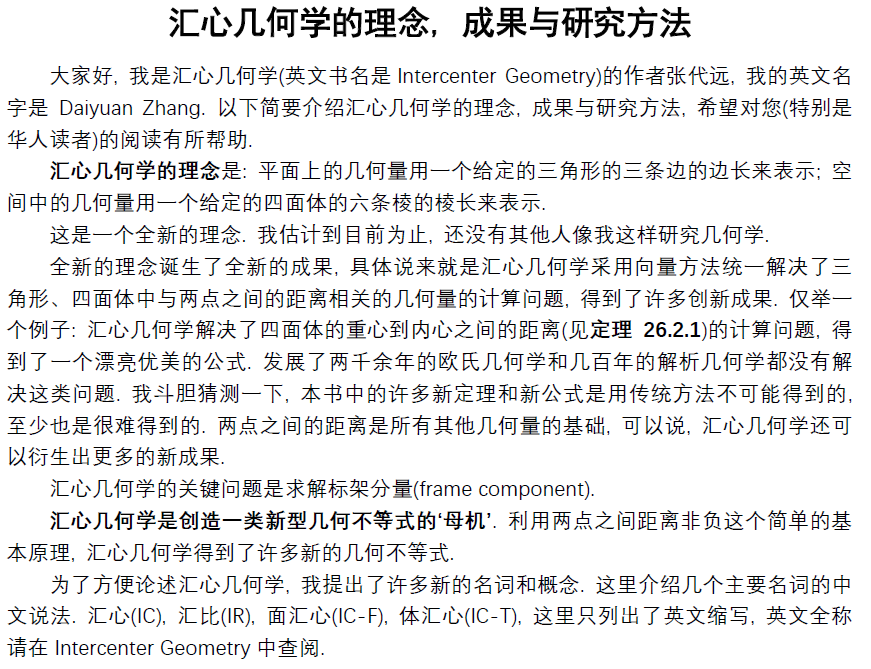}
\end{figure}

\begin{figure}[h]
	\centering
	\includegraphics[totalheight=8cm]{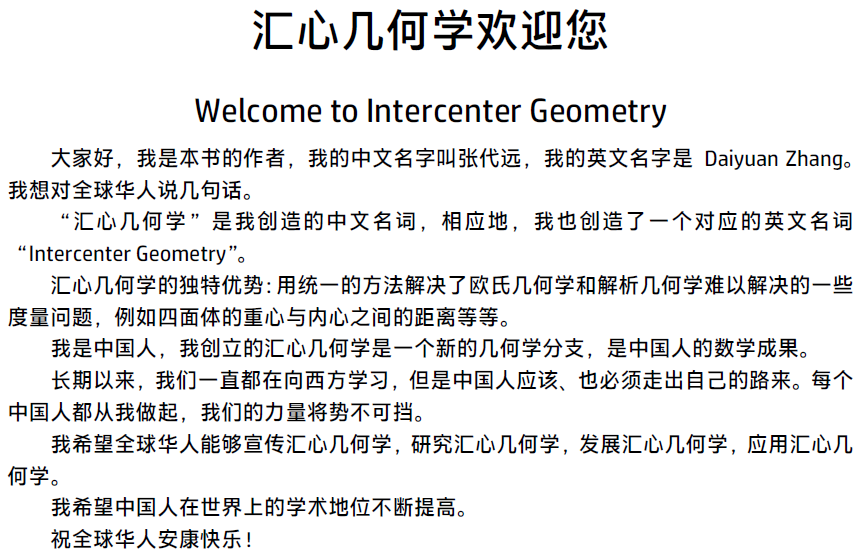}
\end{figure}

\chapter{Ratio of line segments}\label{Ch2} 

In order to study the intersecting ratio, this chapter briefly reviews the concept and basic properties of line segment ratio.

\section{Vector representation of fractional ratio and integral ratio on line segments}\label{Sec2.1}

Generally speaking, vectors in different directions cannot define division operations, but vectors in the same direction can introduce division formally. Suppose three points, $A$, $B$, $P$ are on the same straight line and are called three points collinear. Given points $A$ and $B$, point $A$ is called the starting point, point $B$ is called the ending point (or terminal point), and point $P$ is called the moving point.

This book uses the following symbols to represent fractional ratio and integral ratio of line segments.

\textbf{Fractional ratio}:
\begin{flalign}\label{Eq2.1.1}
	{{\lambda }_{AB}}=\frac{\overrightarrow{AP}}{\overrightarrow{PB}}, \overrightarrow{PB}\ne \overrightarrow{0}.
\end{flalign}

\textbf{Integral ratio}:
\begin{flalign}\label{Eq2.1.2}
	{{\kappa }_{AB}}=\frac{\overrightarrow{AP}}{\overrightarrow{AB}}, \overrightarrow{AB}\ne \overrightarrow{0}.
\end{flalign}

The formula (\ref{Eq2.1.1}) denotes that the vector $\overrightarrow {AP} $ is collinear with the vector $\overrightarrow {PB} $ and $\overrightarrow{AP}={{\lambda }_{AB}}\overrightarrow{PB}$. For example, ${{\lambda}_{ AB}}$ is positive as $\overrightarrow {AP} $ and $\overrightarrow {PB} $ in the same direction; and ${{\lambda}_{ AB}}$ is negative as $\overrightarrow {AP} $ and $\overrightarrow {PB} $ in the opposite direction.

\section{Value range of fractional ratio and integral ratio}\label{Sec2.2}

Obviously, a straight line $\overleftrightarrow{AB}$ can be expressed as the union of the following three line segments, i.e.
\[\overleftrightarrow{AB}={{A}_{-\infty }}{{A}^{*}}\bigcup AB\bigcup {{B}^{\text{*}}}{{B}_{\infty }}.\]
Where
${{A}_{-\infty }}\to -\infty$, ${{A}_{-\infty }}<{{A}^{*}}<A$, ${{A}^{*}}\to A$, $B<{{B}^{*}}<{{B}_{\infty }}$, ${{B}^{*}}\to B$, ${{A}_{-\infty }}{{A}^{*}}\bigcap AB=\varnothing $,
$AB\bigcap {{B}^{\text{*}}}{{B}_{\infty }}=\varnothing $, ${{A}_{-\infty }}{{A}^{*}}\bigcap {{B}^{\text{*}}}{{B}_{\infty }}=\varnothing $.

Obviously, according to the definition of fractional ratio (\ref{Eq2.1.1}), the relationship between the position of point $P$ on the straight line and the numerical range of the fractional ratio ${{\lambda }_{AB}}$ is:
\[\left\{ \begin{aligned}
	& P\in \left( -\infty ,A \right)\Leftrightarrow {{\lambda }_{AB}}\in \left( -1,0 \right) \\ 
	& P\in \left[ A,B \right)\Leftrightarrow {{\lambda }_{AB}}\in \left[ 0,+\infty  \right) \\ 
	& P\in \left( B,+\infty  \right)\Leftrightarrow {{\lambda }_{AB}}\in \left( -\infty ,-1 \right) \\ 
\end{aligned} \right.\]

\section{One-to-one correspondence between the fractional ratio of the segment and the point on a straight line}\label{Sec2.3}
For any point $O$, we have
\[{{\lambda }_{AB}}=\frac{\overrightarrow{AP}}{\overrightarrow{PB}}=\frac{\overrightarrow{OP}-\overrightarrow{OA}}{\overrightarrow{OB}-\overrightarrow{OP}}.\]	

Therefore, 
\[{{\lambda }_{AB}}\left( \overrightarrow{OB}-\overrightarrow{OP} \right)=\overrightarrow{OP}-\overrightarrow{OA},\]	
i.e.
\[\left( 1+{{\lambda }_{AB}} \right)\overrightarrow{OP}=\overrightarrow{OA}+{{\lambda }_{AB}}\overrightarrow{OB}.\]	

Therefore, 
\[\overrightarrow{OP}=\frac{\overrightarrow{OA}+{{\lambda }_{AB}}\overrightarrow{OB}}{1+{{\lambda }_{AB}}}={{\alpha }_{A}}\overrightarrow{OA}+{{\alpha }_{B}}\overrightarrow{OB}.\]	
Where
\begin{flalign*}
	{{\alpha }_{A}}=\frac{1}{1+{{\lambda }_{AB}}},  {{\alpha }_{B}}=\frac{{{\lambda }_{AB}}}{1+{{\lambda }_{AB}}}.
\end{flalign*}


In order to establish the corresponding relationship between the fractional ratio of segment and the one-dimensional Cartesian coordinates, suppose that the points $O$, $A$, $B$, $P$ are collinear, the point $O$ is the origin of the one-dimensional Cartesian coordinates, and the one-dimensional Cartesian coordinates of the starting point $A$, the ending point $B$, and the moving point $P$ are ${{x}_{A}}$, ${{x}_{B}}$, $x$, respectively, then
\begin{flalign*}
	\overrightarrow{OP}=x\mathbf{e}, 	
	\overrightarrow{OA}={{x}_{A}}\mathbf{e}, 	
	\overrightarrow{OB}={{x}_{B}}\mathbf{e}.
\end{flalign*}
Where $\mathbf{e}$ is the unit vector. Therefore,
\[x\mathbf{e}={{\alpha }_{A}}{{x}_{A}}\mathbf{e}+{{\alpha }_{B}}{{x}_{B}}\mathbf{e},\]	
i.e.
\[x={{\alpha }_{A}}{{x}_{A}}+{{\alpha }_{B}}{{x}_{B}}=\frac{1}{1+{{\lambda }_{AB}}}{{x}_{A}}+\frac{{{\lambda }_{AB}}}{1+{{\lambda }_{AB}}}{{x}_{B}},\]	
\[\left( 1+{{\lambda }_{AB}} \right)x={{x}_{A}}+{{\lambda }_{AB}}{{x}_{B}},\]	
\[{{\lambda }_{AB}}\left( x-{{x}_{B}} \right)={{x}_{A}}-x.\]	

Therefore,
\[{{\lambda }_{AB}}=\frac{x-{{x}_{A}}}{{{x}_{B}}-x}.\]

Given the fixed points $A$ and $B$, then ${{x}_{A}}$ and ${{x}_{B}}$ are constants. If $x\ne {{x}_{B}}$, the above formula gives the corresponding relationship between the fractional segment ratio ${{\lambda }_{AB}}$ and the one-dimensional Cartesian coordinate $x$ of the point $P$. The following is proved by mathematical analysis: when $x\in \mathbb{R}-\left\{ {{x}_{B}} \right\}$, the correspondence is one-to-one.


Firstly, ${{\lambda }_{AB}}$ is the elementary function of $x$, and the domain of ${{\lambda }_{AB}}$ is
\[\left( -\infty ,{{x}_{B}} \right)\cup \left( {{x}_{B}},+\infty  \right)=\mathbb{R}-\left\{ {{x}_{B}} \right\}.\]


Since elementary functions are continuous functions in the domain, when $x\in \mathbb{R}-\left\{ {{x}_{B}} \right\}$, ${{\lambda }_{AB}}$ is a continuous function of $x$.
	

Secondly, without losing generality, assuming ${{x}_{B}}>{{x}_{A}}$, then
\[\frac{d{{\lambda }_{AB}}}{dx}=\frac{d}{dx}\left( \frac{x-{{x}_{A}}}{{{x}_{B}}-x} \right)=\frac{\left( {{x}_{B}}-x \right)+\left( x-{{x}_{A}} \right)}{{{\left( {{x}_{B}}-x \right)}^{2}}}=\frac{{{x}_{B}}-{{x}_{A}}}{{{\left( {{x}_{B}}-x \right)}^{2}}}>0.\]


Therefore, when $x\in \mathbb{R}-\left\{ {{x}_{B}} \right\}$, ${{\lambda }_{AB}}$ is a monotone increasing function of $x$, indicating that ${{\lambda }_{AB}}$ is a single valued function of $x$.


To sum up, when $x\in \mathbb{R}-\left\{ {{x}_{B}} \right\}$,  ${{\lambda }_{AB}}$ is a single valued continuous function of $x$, that is, ${{\lambda }_{AB}}$ corresponds to $x$ one by one (one-to-one). This means that if $P\ne B$, the position of the moving point $P$ corresponds to ${{\lambda }_{AB}}$ one by one (one-to-one).

%

\section{Basic properties and relationship about fractional ratio and integral ratio for a given line segment}\label{Sec2.4}

There are some basic properties and some relations between the fraction ratio and the integral ratio for a given line segment. The following theorem gives those relations.

\begin{theorem}{Property theorem of fractional ratio and integral ratio}{Thm2.4.1}\label{Thm2.4.1} 
For the fractional ratio and integral ratio of a given line segment $AB$, they have the following properties:
	\begin{equation}\label{Eq2.4.1}
		{{\lambda }_{AB}}\cdot {{\lambda }_{BA}}=1,		
	\end{equation}
	\begin{equation}\label{Eq2.4.2}
		{{\kappa }_{AB}}+{{\kappa }_{BA}}=1,	
	\end{equation}
	\begin{equation}\label{Eq2.4.3}
		{{\kappa }_{AB}}=\frac{{{\lambda }_{AB}}}{1+{{\lambda }_{AB}}},	
	\end{equation}
	\begin{equation}\label{Eq2.4.4}
		\frac{1}{{{\kappa }_{AB}}}-\frac{1}{{{\lambda }_{AB}}}=1,	
	\end{equation}
	\begin{equation}\label{Eq2.4.5}
		{{\lambda }_{AB}}=\frac{{{\kappa }_{AB}}}{1-{{\kappa }_{AB}}}=\frac{{{\kappa }_{AB}}}{{{\kappa }_{BA}}},	
	\end{equation}		
	\begin{equation}\label{Eq2.4.6}
		{{\lambda }_{AB}}=-{{\kappa }_{PB}},
	\end{equation}		
	\begin{equation}\label{Eq2.4.7}
		{{\kappa }_{AB}}=-{{\lambda }_{PB}}.
	\end{equation}		
\end{theorem}

\begin{proof}
	\[{{\lambda }_{AB}}\cdot {{\lambda }_{BA}}=\frac{\overrightarrow{AP}}{\overrightarrow{PB}}\cdot \frac{\overrightarrow{BP}}{\overrightarrow{PA}}=\frac{\overrightarrow{AP}}{\overrightarrow{PB}}\cdot \left( \frac{-\overrightarrow{PB}}{-\overrightarrow{AP}} \right)=1.\]	
	\[{{\kappa }_{AB}}+{{\kappa }_{BA}}=\frac{\overrightarrow{AP}}{\overrightarrow{AB}}+\frac{\overrightarrow{BP}}{\overrightarrow{BA}}=\frac{\overrightarrow{AP}}{\overrightarrow{AB}}+\left( \frac{-\overrightarrow{PB}}{-\overrightarrow{AB}} \right)=\frac{\overrightarrow{AP}+\overrightarrow{PB}}{\overrightarrow{AB}}=1.\]	
	\[\frac{1}{{{\kappa }_{AB}}}=\frac{\overrightarrow{AB}}{\overrightarrow{AP}}=\frac{\overrightarrow{AP}+\overrightarrow{PB}}{\overrightarrow{AP}}=1+\frac{\overrightarrow{PB}}{\overrightarrow{AP}}=1+\frac{1}{{{\lambda }_{AB}}}=\frac{1+{{\lambda }_{AB}}}{{{\lambda }_{AB}}}.\]	
	\[{{\kappa }_{AB}}=\frac{{{\lambda }_{AB}}}{1+{{\lambda }_{AB}}}.\]	
	\[\frac{1}{{{\kappa }_{AB}}}-\frac{1}{{{\lambda }_{AB}}}=1.\]	
	\[\frac{1}{{{\lambda }_{AB}}}=\frac{1-{{\kappa }_{AB}}}{{{\kappa }_{AB}}}.\]	
	\[{{\lambda }_{AB}}=\frac{{{\kappa }_{AB}}}{1-{{\kappa }_{AB}}}=\frac{{{\kappa }_{AB}}}{{{\kappa }_{BA}}}.\]
	\[{{\lambda }_{AB}}=\frac{\overrightarrow{AP}}{\overrightarrow{PB}}=-\frac{\overrightarrow{PA}}{\overrightarrow{PB}}=-{{\kappa }_{PB}}.\]	
	\[{{\kappa }_{AB}}=\frac{\overrightarrow{AP}}{\overrightarrow{AB}}=-\frac{\overrightarrow{PA}}{\overrightarrow{AB}}=-{{\lambda }_{PB}}.\]	
\end{proof}
\hfill $\square$\par

\chapter{Collinear vectors and coplanar vectors}\label{Ch3}
\thispagestyle{empty}

In many cases, Intercenter Geometry will involve both collinear vectors and coplanar vectors. This chapter introduces several theorems related to collinear vectors and coplanar vectors, which lay a necessary foundation for the following discussion on Intercenter Geometry.

\section{Unique representation theorem of collinear vectors}\label{Sec3.1}

Firstly, the collinearity theorem of three vectors is given.

\begin{theorem}{Unique representation theorem of three collinear vectors}{Thm3.1.1}\label{Thm3.1.1} 
	Assuming two given points $A$, $B$ do not coincide, and point $O$ is any point, then three given points $A$, $B$, $P$ are collinear if and only if there is a unique set of real numbers ${{\alpha }_{A}}$ and ${{\alpha }_{B}}$ makes the following true:
	
	\begin{equation}\label{Eq3.1.1}
		\overrightarrow{OP}={{\alpha }_{A}}\overrightarrow{OA}+{{\alpha }_{B}}\overrightarrow{OB},	
	\end{equation}
	\begin{equation}\label{Eq3.1.2}
		{{\alpha }_{A}}+{{\alpha }_{B}}=1.
	\end{equation}
\end{theorem}

\begin{proof}
	\textbf{Prove of sufficiency}. That is to prove: under the general premise of “assuming that the given two points $A$ and $B $ do not coincide, and the point $O $ is any point”, the collinearity of three terminal points $A $, $B $, $P $ can be deduced from the conditions of (\ref{Eq3.1.1}) and (\ref{Eq3.1.2}).
	
	If there is a set of real numbers ${{\alpha}_{A}}$ and ${{\alpha}_{B}}$ that satisfies conditions (\ref{Eq3.1.1}) and (\ref{Eq3.1.2}), the following formula is obtained
	\[\overrightarrow{OP}={{\alpha }_{A}}\overrightarrow{OA}+{{\alpha }_{B}}\overrightarrow{OB}={{\alpha }_{A}}\overrightarrow{OA}+\left( 1-{{\alpha }_{A}} \right)\overrightarrow{OB},\]	
	i.e.
	\[\overrightarrow{OP}-\overrightarrow{OB}={{\alpha }_{A}}\left( \overrightarrow{OA}-\overrightarrow{OB} \right),\]	
	i.e.
	\begin{equation}\label{Eq3.1.3}
		\overrightarrow{BP}={{\alpha }_{A}}\overrightarrow{BA}.		
	\end{equation}
	
	Because the points $A$ and $B$ do not coincide, so $\overrightarrow{BA}\ne \overrightarrow{0}$, therefore $\overrightarrow {BA}$ is a basis vector of one-dimensional space. The above formula shows that $\overrightarrow {BP}$ can be the linear combination of the basis vector $\overrightarrow {BA}$ with coefficient ${{\alpha}_{A}}$, that is, $\overrightarrow {BP} $ and $\overrightarrow {BA} $ are collinear, so the three points $A $, $B $, $P $ are collinear.
	
	Moreover, this set of real numbers ${{\alpha}_{A}}$ and ${{\alpha }_{B}}$ is unique. In fact, when all three points $A$, $B$, $P$ are given, ${{\alpha}_{A}}$ is unique. This is because if two real numbers ${{\alpha}_{A}}$ and ${{\beta}_{A}}$ make (\ref{Eq3.1.3}) valid, i.e. $\overrightarrow{BP}={{\alpha }_{A}}\overrightarrow{BA}$, $\overrightarrow{BP}={{\beta }_{A}}\overrightarrow{BA}$, then $\left( {{\alpha }_{A}}-{{\beta }_{A}} \right)\overrightarrow{BA}=\overrightarrow{0}$. But $\overrightarrow{BA}\ne\overrightarrow{0}$, we have ${{\alpha }_{A}}={{\beta }_{A}}$, therefore ${{\alpha}_{A}}$ is unique. So ${{\alpha }_{B}}=1-{{\alpha }_{A}}$ is also unique.

	\textbf{Necessity}. That is to prove: under the general premise of “assuming that two given points $A$ and $B$ do not coincide, and point $O$ is any point”, formulas (\ref{Eq3.1.1}) and (\ref{Eq3.1.2}) can be derived according to the collinearity of three points $A$, $B$, $P$.
		
	Assume that the given three points, $A$, $B$, $P$ are collinear, that is, $\overrightarrow{BP}$ and $\overrightarrow{BA}$ are collinear. Because point $A$ and $B$ do not coincide, $\overrightarrow{BA}\ne \overrightarrow{0}$, so $\overrightarrow{BA}$ is a basis vector in one-dimensional space, so there is unique number ${{\alpha}_{A}}$ makes the following true:
	\[\overrightarrow{BP}={{\alpha }_{A}}\overrightarrow{BA}.\]

	For any point $O$, the above expression can be written as
	\[\overrightarrow{OP}-\overrightarrow{OB}={{\alpha }_{A}}\left( \overrightarrow{OA}-\overrightarrow{OB} \right),\]
	i.e.
	\[\overrightarrow{OP}={{\alpha }_{A}}\overrightarrow{OA}+\left( 1-{{\alpha }_{A}} \right)\overrightarrow{OB}={{\alpha }_{A}}\overrightarrow{OA}+{{\alpha }_{B}}\overrightarrow{OB},\]
	where
	\[{{\alpha }_{A}}+{{\alpha }_{B}}=1.\]	
	
	While ${{\alpha }_{B}}=1-{{\alpha }_{A}}$, according to the uniqueness of ${{\alpha }_{A}}$, we known that ${{\alpha }_{B}}$ is also unique, then formulas (\ref{Eq3.1.1}) and (\ref{Eq3.1.2}) can be obtained.
\end{proof}
\hfill $\square$\par

	The real numbers ${{\alpha }_{A}}$ and ${{\alpha }_{B}}$ are called the \textbf{frame components} of the vector $\overrightarrow {OP}$ on the basis vectors $\overrightarrow {OA}$ and $\overrightarrow {OB}$ respectively, also called the coordinates or coefficients. And each of ${{\alpha }_{A}}\overrightarrow{OA}$, ${{\alpha }_{B}}\overrightarrow{OB}$ is called the component vector along the direction of the basis vector $\overrightarrow {OA}$, $\overrightarrow {OB}$ respectively. The point $O$ is called the origin of the frame (referred to as the origin for short), and the vectors $\overrightarrow {OA} $, $\overrightarrow {OB} $ are called the two-dimensional frame. The two-dimensional frame is denoted as $\left (O; A,B \right)$.
	
	Formula (\ref{Eq3.1.1}) states that the vector $\overrightarrow{OP}$ can be a linear combination of the
	basis vectors $\overrightarrow{OA}$ and $\overrightarrow{OB}$ so that three position vectors $\overrightarrow{OP}$, $\overrightarrow{OA}$ and $\overrightarrow{OB}$ are located on the same plane. According to the basic theorem of plane vectors, any vector on the plane spanned by a set of basis vectors $\overrightarrow{OA}$ and $\overrightarrow{OB}$ can be a linear combination of the bisis vectors $\overrightarrow{OA}$ and $\overrightarrow{OB}$, and the representation factor ${{\alpha}_{A}}$ and ${{\alpha}_{B}}$ in formula (\ref{Eq3.1.1}) are also unique, but ${{\alpha}_{A}}$ and ${{\alpha}_{B}}$ are independent of each other and unconstrained. When the  coefficients ${{\alpha}_{ A} } $ and ${{\alpha}_{ B}} $ of formula (\ref{Eq3.1.1}) are constrained by condition (\ref{Eq3.1.2}), then the position of the end point of the vector $\overrightarrow{OP}$ will be constrained. The conclusion of the theorem is that the given three points $A$, $B$, $P$ are collinear.

	In the above theorem, point $O$ is any point, which is useful because we can choose point $O$ to coincide with point $A$ (or point $B$, or point $C$) to make some problems easier.


	Theorem \ref{thm:Thm3.1.1} has been discussed in many materials, and the formulation and proof process of the theorem are also different. However, the formulation and proof process of theorem \ref{thm:Thm3.1.1} in this book are given by the author independently.

	Now that it has been proven that there is a unique set of real numbers ${{\alpha }_{A}}$ and ${{\alpha }_{B}}$, the question is how to calculate those two real numbers ${{\alpha }_{A}}$ and ${{\alpha }_{B}}$? Now let's answer this question.

\section{Representation of frame components for collinear vectors by line segment ratio}\label{Sec3.2}

According to the above definition, the following theorem can be obtained.
\begin{theorem}{Representation of frame components for collinear vectors by line segment ratio}{Thm3.2.1}\label{Thm3.2.1} 
	Assume points $A$, $B$, $P$ are collinear and
	\[{{\lambda }_{AB}}=\frac{\overrightarrow{AP}}{\overrightarrow{PB}},\]	
	then, for any point $O$, the vector $\overrightarrow{OP}$ can be uniquely expressed as:
	\begin{equation}\label{Eq3.2.1}
		\overrightarrow{OP}={{\alpha }_{A}}\overrightarrow{OA}+{{\alpha }_{B}}\overrightarrow{OB},		
	\end{equation}
	where
	\begin{flalign}\label{Eq3.2.2}
		{{\alpha }_{A}}=\frac{1}{1+{{\lambda }_{AB}}}, {{\alpha }_{B}}=\frac{{{\lambda }_{AB}}}{1+{{\lambda }_{AB}}}.
	\end{flalign}
\end{theorem}

\begin{proof}
	
	For any point $O$,
	\[{{\lambda }_{AB}}=\frac{\overrightarrow{AP}}{\overrightarrow{PB}}=\frac{\overrightarrow{OP}-\overrightarrow{OA}}{\overrightarrow{OB}-\overrightarrow{OP}},\]	
	therefore,
	\[{{\lambda }_{AB}}\left( \overrightarrow{OB}-\overrightarrow{OP} \right)=\overrightarrow{OP}-\overrightarrow{OA},\]	
	i.e.
	\[\left( 1+{{\lambda }_{AB}} \right)\overrightarrow{OP}=\overrightarrow{OA}+{{\lambda }_{AB}}\overrightarrow{OB}.\]	
	
	Therefore,
	\[\overrightarrow{OP}=\frac{\overrightarrow{OA}+{{\lambda }_{AB}}\overrightarrow{OB}}{1+{{\lambda }_{AB}}}={{\alpha }_{A}}\overrightarrow{OA}+{{\alpha }_{B}}\overrightarrow{OB},\]	
	where
	\begin{flalign*}
		{{\alpha }_{A}}=\frac{1}{1+{{\lambda }_{AB}}}, {{\alpha }_{B}}=\frac{{{\lambda }_{AB}}}{1+{{\lambda }_{AB}}}.	
	\end{flalign*}
	
	Obviously,
	\[{{\alpha }_{A}}+{{\alpha }_{B}}=1.\]	
	
	According to the above formula and formula (\ref{Eq3.2.1}), using theorem \ref{thm:Thm3.1.1}, we know that the expression (\ref{Eq3.2.1}) is unique.
\end{proof}
\hfill $\square$\par

Obviously, for any point $O$, the above formula holds, that is, formula (\ref{Eq3.2.1}) is independent of the choice of point $O$.

According to formula (\ref{Eq3.2.2}), two frame components (coefficients) ${{\alpha}_{A}}$ and ${{\alpha}_{B}}$ are only related to segment fraction ${{\lambda}_{AB}}$, not to the location of coordinate origin $O$. This means that: if the position of end point is given, then each of the components (coefficients) of the position vector $\overrightarrow{OP}$ in the frame $\overrightarrow{OA}$ and $\overrightarrow{OB}$ is independent of the position of the coordinate origin $O$ (see formula (\ref{Eq3.2.1})), and is only related to the segment ratio ${{\lambda}_{AB}}$. This is an important feature. It is precisely by grasping this characteristic that I put forward the Intercenter Geometry. Chapter 4 of this book introduces the concept of intersecting center. As we all know, in analytic geometry, the components (coefficients) of a position vector determined by the position of the end point in a frame are related to the position of the coordinate origin $O$, so Intercenter Geometry is different from analytic geometry.

In order to deal with the problems on the plane, we need to generalize the above concepts and methods.

%
\section{Unique representation theorem of coplanar vectors}\label{Sec3.3}
In this book, the notation $O\in {{\mathbb{R}}^{3}}$ indicates that point $O$ is an arbitrary point in a three-dimensional finite space; and the notation $P\in {{\mathbb{R}}^{3}}$ indicates that point $P$ is also an arbitrary point in the same three-dimensional finite space as that of point $O$.

The coplanar theorem of space vector is given in the following. 

\begin{theorem}{Unique representation theorem of coplanar vectors}{Thm3.3.1}\label{Thm3.3.1} 
	Suppose that given three points $A$, $B$, $C$ are not coincident with each other and are not collinear, and point $O\in {{\mathbb{R}}^{3}}$, point $P\in {{\mathbb{R}}^{3}}$, then the sufficient and necessary condition that four points $A$, $B$, $C$ and $P$ are coplanar is that there exists a unique set of real numbers ${{\alpha }_{A}}$, ${{\alpha }_{B}}$ and ${{\alpha }_{C}}$ such that:	
	\begin{equation}\label{Eq3.3.1}
		\overrightarrow{OP}={{\alpha }_{A}}\overrightarrow{OA}+{{\alpha }_{B}}\overrightarrow{OB}+{{\alpha }_{C}}\overrightarrow{OC},	
	\end{equation}
	\begin{equation}\label{Eq3.3.2}
		{{\alpha }_{A}}+{{\alpha }_{B}}+{{\alpha }_{C}}=1.
	\end{equation}
\end{theorem}

\begin{proof}
	\textbf{Sufficiency}. That is, to prove that four points $A$, $B$, $C$ and $P$ are coplanar can be derived from the conditions (\ref{Eq3.3.1}) and (\ref{Eq3.3.2}) under the overall precondition that “three points $A$, $B$, $C$ are not coincident and collinear, point $O\in {{\mathbb{R}}^{3}}$, point $P\in {{\mathbb{R}}^{3}}$”.
	
	If there is a set of real numbers ${{\alpha}_{A}}$, ${{\alpha }_{B}}$ and ${{\alpha }_{C}}$ satisfies the conditions (\ref{Eq3.3.1}) and (\ref{Eq3.3.2}), then
	\[\begin{aligned}
		\overrightarrow{OP}& ={{\alpha }_{A}}\overrightarrow{OA}+{{\alpha }_{B}}\overrightarrow{OB}+{{\alpha }_{C}}\overrightarrow{OC} \\
		& =\left( 1-{{\alpha }_{B}}-{{\alpha }_{C}} \right)\overrightarrow{OA}+{{\alpha }_{B}}\overrightarrow{OB}+{{\alpha }_{C}}\overrightarrow{OC},  
	\end{aligned}\]	
	i.e.
	\[\overrightarrow{OP}-\overrightarrow{OA}={{\alpha }_{B}}\left( \overrightarrow{OB}-\overrightarrow{OA} \right)+{{\alpha }_{C}}\left( \overrightarrow{OC}-\overrightarrow{OA} \right),\]	
	i.e.
	\[\overrightarrow{AP}={{\alpha }_{B}}\overrightarrow{AB}+{{\alpha }_{C}}\overrightarrow{AC}.\]	
	
	Because the three points $A$, $B$, $C$ do not coincide with each other and are not collinear, $\overrightarrow{AB}$ and $\overrightarrow{AC}$ can form a set of basis vectors. 
	
	The formula above states that the vector $\overrightarrow{AP}$ can be a linear combination of the
	vectors $\overrightarrow{AB}$ and $\overrightarrow{AC}$, which means that the vector $\overrightarrow{AP}$ is in the plane spanned by $\overrightarrow{AB}$ and $\overrightarrow{AC}$, so the point $P$ and points $A$, $B$, $C$ are on the same plane.
	
	According to the basic theorem of plane vector, ${{\alpha }_{B}}$ and ${{\alpha }_{C}}$ in the formula above is unique. Therefore ${{\alpha }_{A}}=1-{{\alpha }_{B}}-{{\alpha }_{C}}$ is also unique.
	
	\textbf{Necessity}. That is, to prove that the conditions (\ref{Eq3.3.1}) and (\ref{Eq3.3.2}) can be derived from four points $A$, $B$, $C$ and $P$ coplanar under the overall precondition that “three points $A$, $B$, $C$ are not coincident and collinear, point $O\in {{\mathbb{R}}^{3}}$, point $P\in {{\mathbb{R}}^{3}}$”.
	
	Suppose that given four points $A$, $B$, $C$ and $P$ are coplanar, then $\overrightarrow{AP}$ is coplanar with $\overrightarrow{AB}$ and $\overrightarrow{AC}$. Because the three points $A$, $B$, $C$ are not coincident with each other, there must be $\overrightarrow{AB}\ne \overrightarrow{0}$, $\overrightarrow{AC}\ne \overrightarrow{0}$. Because the three points $A$, $B$, $C$ are not collinear, then $\overrightarrow {AB}$ and $\overrightarrow {AC}$ are linearly independent, i.e. $\overrightarrow {AB}$ and $\overrightarrow {AC}$ form a set of basis vectors on the $\triangle ABC $ plane. Therefore $\overrightarrow {AP} $ can be the uniquely linear combination of the basis vectors of $\overrightarrow {AB} $ and $\overrightarrow {AC}$, that is, there exists a unique set of real numbers ${{\alpha }_{B}}$ and ${{\alpha }_{C}}$ such that the following formula holds:
	\[\overrightarrow{AP}={{\alpha }_{B}}\overrightarrow{AB}+{{\alpha }_{C}}\overrightarrow{AC}.\]
	
	For any point $O\in {{\mathbb{R}}^{3}}$, the above formula can be written as
	\[\overrightarrow{OP}-\overrightarrow{OA}={{\alpha }_{B}}\left( \overrightarrow{OB}-\overrightarrow{OA} \right)+{{\alpha }_{C}}\left( \overrightarrow{OC}-\overrightarrow{OA} \right),\]
	i.e.	
	\[\begin{aligned}
		\overrightarrow{OP}& =\overrightarrow{OA}+{{\alpha }_{B}}\left( \overrightarrow{OB}-\overrightarrow{OA} \right)+{{\alpha }_{C}}\left( \overrightarrow{OC}-\overrightarrow{OA} \right) \\
		& =\left( 1-{{\alpha }_{B}}-{{\alpha }_{C}} \right)\overrightarrow{OA}+{{\alpha }_{B}}\overrightarrow{OB}+{{\alpha }_{C}}\overrightarrow{OC} \\
		& ={{\alpha }_{A}}\overrightarrow{OA}+{{\alpha }_{B}}\overrightarrow{OB}+{{\alpha }_{C}}\overrightarrow{OC},  
	\end{aligned}\]
	where
	\[{{\alpha }_{A}}+{{\alpha }_{B}}+{{\alpha }_{C}}=1.\]	
	
	Because ${{\alpha }_{B}}$ and ${{\alpha }_{C}}$ is unique, so ${{\alpha }_{A}}=1-{{\alpha }_{B}}-{{\alpha }_{C}}$ is also unique, this proves that formulas (\ref{Eq3.3.1}) and (\ref{Eq3.3.2}) are unique.
\end{proof}
\hfill $\square$\par

The above theorem holds for any point $O\in {{\mathbb{R}}^{3}}$, which is very useful, because we can choose the point $O$ to coincide with the point $A$ (or point $B$, or point $C$, etc.), which makes some problems simple. You can see this application later.

Theorem \ref{thm:Thm3.3.1} has been discussed in many materials, and the formulation and proof process of the theorem are also different. However, the formulation and proof process of theorem \ref{thm:Thm3.3.1} in this book are given by the author independently.

The above theorem specifies the components (coefficients) ${{\alpha}_{A}}$, ${{\alpha}_{B}}$ and ${{\alpha}_{C}}$ are unique. How do you calculate these three components in detail? Can you find those components (coefficients) that are independent of the origin of the frame $O$ like the segment ratio in theorem \ref{thm:Thm3.2.1}? The answer is yes, which is the question to be studied below. In order to determine the components (coefficients), we need to introduce the concepts of intersecting center and intersecting ratio, and then give the calculation method for components (coefficients) ${{\alpha }_{A}}$, ${{\alpha }_{B}}$ and ${{\alpha }_{C}}$ in subsequent chapters.

\section{Some basic theorems on triangular frame}\label{Sec3.4}
In order to deal with the triangle problem conveniently, a very useful theorem is given below. I name it the collinearity theorem of three vector endpoints.

\begin{theorem}{Theorem of terminal collinearity for three vectors, Daiyuan Zhang}{Thm3.4.1}\label{Thm3.4.1} 
	Assuming that the given three points $A$, $B$, $C$ do not coincide with each other and $O\in {{\mathbb{R}}^{3}}$, the necessary and sufficient conditions for the three points $A$, $B$, $C$ to be collinear are:
			
	\[{{\alpha }_{A}}\overrightarrow{OA}+{{\alpha }_{B}}\overrightarrow{OB}+{{\alpha }_{C}}\overrightarrow{OC}=\vec{0},\]	
	\[{{\alpha }_{A}}+{{\alpha }_{B}}+{{\alpha }_{C}}=0.\]	
	Where the real number ${{\alpha }_{A}}\ne 0$, ${{\alpha }_{B}}\ne 0$, ${{\alpha }_{C}}\ne 0$, and if any one of the three real numbers ${{\alpha }_{A}}$, ${{\alpha }_{B}}$, ${{\alpha }_{C}}$ is 1, the other two are uniquely determined.
\end{theorem}

\begin{proof}
	\textbf {Sufficiency}. According to the conditions of the theorem, we can get:
	\[\overrightarrow{OA}=\left( -\frac{{{\alpha }_{B}}}{{{\alpha }_{A}}} \right)\overrightarrow{OB}+\left( -\frac{{{\alpha }_{C}}}{{{\alpha }_{A}}} \right)\overrightarrow{OC},\]	
	\[\left( -\frac{{{\alpha }_{B}}}{{{\alpha }_{A}}} \right)+\left( -\frac{{{\alpha }_{C}}}{{{\alpha }_{A}}} \right)=1.\]	
	
	Based on the collinearity theorem of three vector endpoints (theorem \ref{thm:Thm3.1.1}), the three endpoints $A$, $B$, $C$ of $\overrightarrow{OA}$, $\overrightarrow{OB}$ and $\overrightarrow{OC}$ are collinear.
	
	
	\textbf{Necessity}. Suppose three vector endpoints $A$, $B$, $C$ are collinear. The three endpoints $A$, $B$, $C$, do not coincide, so $\overrightarrow{CA}\ne \vec{0}$, $\overrightarrow{CB}\ne \vec{0}$. Thus, $\overrightarrow{CA}$ and $\overrightarrow{CB}$ are linearly dependen, according to the basic theorem of vectors, there exists a unique real number $\beta $ such that
	\[\overrightarrow{CA}=\beta \overrightarrow{CB}.\]	
	
	
	The following proves that $\beta \ne 0$ and $\beta \ne 1$. If $\beta = 0 $, the above formula becomes $\overrightarrow{CA}=\vec{0}$, which means that the two points of $A$, $C$ coincide and there is a contradiction. Similarly, if $\beta =1$, the above formula becomes $\overrightarrow{CA}=\overrightarrow{CB}$, which means that the two points of $A$, $B$ coincide, resulting in contradiction.
	
	For the point $O\in {{\mathbb{R}}^{3}}$, the above formula can be written as
	\[\overrightarrow{OA}-\overrightarrow{OC}=\beta \left( \overrightarrow{OB}-\overrightarrow{OC} \right),\]
	i.e. 
	\[\overrightarrow{OA}=\beta \overrightarrow{OB}+\left( 1-\beta  \right)\overrightarrow{OC},\]
	i.e. 
	\begin{equation}\label{Eq3.4.1}
		\overrightarrow{OA}-{{\beta }_{1}}\overrightarrow{OB}-{{\beta }_{2}}\overrightarrow{OC}=\vec{0},
	\end{equation}
	\begin{equation}\label{Eq3.4.2}
		1+\left( -{{\beta }_{1}} \right)+\left( -{{\beta }_{2}} \right)=0,
	\end{equation}
	where 
	\[{{\beta }_{1}}=\beta ,\ {{\beta }_{2}}=1-\beta. \]	
	
	Since ${{\beta }_{1}}=\beta \ne 0$, ${{\beta }_{1}}=\beta \ne 1$, so ${{\beta }_{2}}=1-\beta \ne 1$, ${{\beta }_{2}}=1-\beta \ne 0$.  Formulas (\ref{Eq3.4.1}) and (\ref{Eq3.4.2}) can also be written as
	\[{{\alpha }_{A}}\overrightarrow{OA}+{{\alpha }_{B}}\overrightarrow{OB}+{{\alpha }_{C}}\overrightarrow{OC}=\vec{0},\]	
	\[{{\alpha }_{A}}+{{\alpha }_{B}}+{{\alpha }_{C}}=0.\]	
	
	Where ${{\alpha }_{A}}=\gamma \ne 0$, ${{\alpha }_{B}}=-\gamma {{\beta }_{1}}\ne 0$, ${{\alpha }_{C}}=-\gamma {{\beta }_{2}}\ne 0$. 
	
	
	Obviously, since the arbitrariness of real number $\gamma \ne 0$, ${{\alpha }_{A}}$, ${{\alpha }_{B}}$ and ${{\alpha }_{C}}$ are not unique. However, since non-zero real number ${{\beta }_{1}}=\beta $ is unique, so ${{\beta }_{2}}=1-\beta $ is also unique, so when the coefficient ${{\alpha }_{A}}=\gamma =1$ of $\overrightarrow{OA}$, the other two coefficients ${{\alpha }_{B}}$, ${{\alpha }_{C}}$ are unique.
	
	
	The coefficient ${{\alpha }_{A}}=1$ of $\overrightarrow{OA}$ has been studied above, Similar results can be obtained when ${{\alpha }_{B}}=1$, ${{\alpha }_{C}}=1$, which is left to the reader as an exercise.
\end{proof}
\hfill $\square$\par
According to the above theorem, the following terminal non-collinear theorem can be obtained.


\begin{theorem}{Noncollinearity of terminal points for three vectors, Daiyuan Zhang}{Thm3.4.2}\label{Thm3.4.2} 
	Suppose that the three points $A$, $B$, $C$ do not coincide with each other, $O\in {{\mathbb{R}}^{3}}$, and the following two conditions are satisfied:	
	\begin{equation}\label{Eq3.4.3}
		{{\alpha }_{A}}\overrightarrow{OA}+{{\alpha }_{B}}\overrightarrow{OB}+{{\alpha }_{C}}\overrightarrow{OC}=\vec{0},
	\end{equation}
	\begin{equation}\label{Eq3.4.4}
		{{\alpha }_{A}}+{{\alpha }_{B}}+{{\alpha }_{C}}=0,
	\end{equation}
	then, if the three terminal points $A$, $B$, $C$ of the three position vectors $\overrightarrow{OA}$, $\overrightarrow{OB}$, $\overrightarrow{OC}$ with the same initial point $O$ are not collinear, there must be ${{\alpha }_{A}}={{\alpha }_{B}}={{\alpha }_{C}}=0$.
\end{theorem}

\begin{proof}
	
	When the condition of this theorem is satisfied, from theorem \ref{thm:Thm3.4.1}, if three terminal points $A$, $B$, $C$ are not collinearr, then at least one of the three coefficients ${{\alpha }_{A}}$, ${{\alpha }_{B}}$, ${{\alpha }_{C}}$ is zero. Without losing generality, suppose ${{\alpha }_{A}}=0$. ${{\alpha }_{A}}=0$ is substituted into formula (\ref{Eq3.4.4}) to obtain $-{{\alpha }_{C}}\overrightarrow{OB}+{{\alpha }_{C}}\overrightarrow{OC}=\vec{0}$, i.e. ${{\alpha }_{C}}\overrightarrow{BC}=\vec{0}$.  Since the points $B$ and $C$ do not coincide, so $\overrightarrow{BC}\ne \vec{0}$, so there can only be ${{\alpha }_{C}}=0$.   According to formula (\ref{Eq3.4.3}), we have ${{\alpha }_{B}}=0$. Obviously, The same conclusion is also obtained if more than one of the three coefficients ${{\alpha }_{A}}$, ${{\alpha }_{B}}$, ${{\alpha }_{C}}$ equals zero.
\end{proof}
\hfill $\square$\par
For the triangular frame $\left( O;A,B,C \right)$, its three vertices $A$, $B$, $C$ do not coincide with each other, and the three points $A$, $B$, $C$ are not collinear. Therefore, the following theorem can be obtained directly according to the above theorem.


\begin{theorem}{Theorem of triangular frame, Daiyuan Zhang}{Thm3.4.3}\label{Thm3.4.3} 
	Suppose that given the triangular frame $\left( O;A,B,C \right)$, $O\in {{\mathbb{R}}^{3}}$, and the following two conditions are satisfied:
	\[{{\alpha }_{A}}\overrightarrow{OA}+{{\alpha }_{B}}\overrightarrow{OB}+{{\alpha }_{C}}\overrightarrow{OC}=\vec{0},\]	
	\[{{\alpha }_{A}}+{{\alpha }_{B}}+{{\alpha }_{C}}=0,\]	
	then ${{\alpha }_{A}}={{\alpha }_{B}}={{\alpha }_{C}}=0$.
\end{theorem}

The following theorem makes use of ablve conclusion. 
\begin{theorem}{Unique representation of a vector on triangular frame, Daiyuan Zhang}{Thm3.4.4}\label{Thm3.4.4} 
	Given a triangular frame $\left( O;A,B,C \right)$, $O\in {{\mathbb{R}}^{3}}$, $P\in {{\mathbb{R}}^{3}}$, if the following formulas hold for a set of real numbers ${{\alpha }_{A}}$, ${{\alpha }_{B}}$, ${{\alpha }_{C}}$:
	\[\overrightarrow{OP}={{\alpha }_{A}}\overrightarrow{OA}+{{\alpha }_{B}}\overrightarrow{OB}+{{\alpha }_{C}}\overrightarrow{OC},\]	
	\[{{\alpha }_{A}}+{{\alpha }_{B}}+{{\alpha }_{C}}=1,\]	
	then each of the real numbers ${{\alpha }_{A}}$, ${{\alpha }_{B}}$, ${{\alpha }_{C}}$ is unique.
\end{theorem}

\begin{proof}
	If two sets of real numbers ${{\alpha }_{A}}$, ${{\alpha }_{B}}$, ${{\alpha }_{C}}$ and ${{\beta }_{A}}$, ${{\beta }_{B}}$, ${{\beta }_{C}}$ satisfy the above conditions at the same time, then
	\[\left( {{\alpha }_{A}}-{{\beta }_{A}} \right)\overrightarrow{OA}+\left( {{\alpha }_{B}}-{{\beta }_{B}} \right)\overrightarrow{OB}+\left( {{\alpha }_{C}}-{{\beta }_{C}} \right)\overrightarrow{OC}=\overrightarrow{0},\]
	\[\left( {{\alpha }_{A}}-{{\beta }_{A}} \right)+\left( {{\alpha }_{B}}-{{\beta }_{B}} \right)+\left( {{\alpha }_{C}}-{{\beta }_{C}} \right)=1-1=0.\]
	
	According to theorem \ref{thm:Thm3.4.3}, the following results are obtained: ${{\alpha }_{A}}-{{\beta }_{A}}=0$, ${{\alpha }_{B}}-{{\beta }_{B}}=0$, ${{\alpha }_{C}}-{{\beta }_{C}}=0$, ${{\alpha }_{D}}-{{\beta }_{D}}=0$. 
\end{proof}
\hfill $\square$\par
The above theorem shows that under the condition of the theorem, the vector $\overrightarrow{OP}$ can be expressed by a linear combination of the basis vectors of the triangular frame $\left( O;A,B,C \right)$ uniquely and linearly, where ${{\alpha }_{A}}$, ${{\alpha }_{B}}$, ${{\alpha }_{C}}$ are called the \textbf{affine frame components} of the frame $\overrightarrow{OA}$, $\overrightarrow{OB}$, $\overrightarrow{OC}$ at point $P$, which are called \textbf{components} or \textbf{coefficients} for short.

\chapter{Basic concepts of Plane Intercenter Geometry}\label{Ch4}
\thispagestyle{empty}


For the convenience of discussion, I put forward many new concepts in this chapter, many of which will be frequently used in later chapters of this book. I hope readers will be familiar with those concepts as soon as possible, and the calculation method of quantities related to those new concepts will be studied later.


In order to make readers understand why I put forward these new concepts, I need to briefly explain my original idea of Intercenter Geometry.


The idea of Plane Intercenter Geometry is that the geometric quantities on a plane are expressed by the lengths of the three sides of a given triangle.


Plane Intercenter Geometry uses a unified vector method to solve the calculation problem of geometric quantities related to the distance between two points in the plane.


The key problem of Plane Intercenter Geometry is to establish a frame and calculate the frame components, so that each of the frame components can be expressed as a function of the lengths of the three sides of a given triangle, and the frame components are independent of the origin of the frame.


Firstly, consider the case of straight line (one dimension). Suppose a line segment $AB$ is embedded in a straight line $l$, a point $P$ is selected on the straight line $l$, and a point $O$ is given arbitrarily in space. At this time, the vectors $\overrightarrow{OA}$ and $\overrightarrow{OB}$ constitute the frame. According to theorem \ref{thm:Thm3.2.1}, we have

\[\overrightarrow{OP}={{\alpha }_{A}}\overrightarrow{OA}+{{\alpha }_{B}}\overrightarrow{OB}, \]	
where
\begin{flalign*}
	{{\alpha }_{A}}=\frac{1}{1+{{\lambda }_{AB}}},  {{\alpha }_{B}}=\frac{{{\lambda }_{AB}}}{1+{{\lambda }_{AB}}}.	
\end{flalign*}


It should be noted that the frame components ${{\alpha }_{A}}$, ${{\alpha }_{B}}$ are only related to the fractional ratio of segment ${{\lambda }_{AB}}$, while the fractional ratio of segment ${{\lambda }_{AB}}$ is independent of the position of the frame origin $O$, so the frame components ${{\alpha }_{A}}$, ${{\alpha }_{B}}$ are also independent of the position of the frame origin $O$. On the other hand, as long as the points $A$, $B$ and $P$ are given, the ${{\lambda }_{AB}}$ is a constant. That is, for any point $O$, the vector $\overrightarrow{OP}$ can be the linear combination of the frame vectors $\overrightarrow{OA}$ and $\overrightarrow{OB}$, and the coefficients (frame components) ${{\alpha }_{A}}$ and ${{\alpha }_{B}}$ are constants. This is very important. It can be said that this is the starting point of the germination of the original thought of Intercenter Geometry. It once prompted me to think about a question: can I establish a geometry in which the frame components are all constants? My answer is yes. For triangles and tetrahedrons, I created such geometry and named it Intercenter Geometry.


Now consider the plane case (two-dimensional case).


Suppose a triangle $\triangle ABC$ is embedded in a plane, a point $P$ is selected on the plane of the $\triangle ABC$, and a point $O$ is given at any point in space. We hope to have the following expression:

\[\overrightarrow{OP}={\alpha _{A}^{P}}\overrightarrow{OA}+{\alpha _{B}^{P}}\overrightarrow{OB}+{\alpha _{C}^{P}}\overrightarrow{OC},\]	
and the frame components ${\alpha _{A}^{P}}$, ${\alpha _{B}^{P}}$ and ${\alpha _{C}^{P}}$ are independent of the origin of the frame.


The question now is, can the frame components ${\alpha _{A}^{P}}$, ${\alpha _{B}^{P}}$ and ${\alpha _{C}^{P}}$ be related to the fractional ratio of some line segments similar to the straight line case? The answer I got through my research is yes. The fractional ratio of these line segments is the concept of \textbf{intersecting ratio} to be introduced below. Moreover, if the point $P$ is selected at some special points of $\triangle ABC$ (such as centroid, incenter, orthocenter, etc.), its intersecting ratio can be expressed by the lengths of the three sides of $\triangle ABC$. This means the frame components ${\alpha _{A}^{P}}$, ${\alpha _{B}^{P}}$ and ${\alpha _{C}^{P}}$ of the vector $\overrightarrow{OP}$ are only related to the lengths of the three sides of $\triangle ABC$. in this way, the distance between two points on the $\triangle ABC$ plane can be expressed as a function of the length of the three sides of $\triangle ABC$ in a unified way. Therefore, all geometric quantities that can be expressed by the distance between two points on the plane can also be expressed by the lengths of the three sides of $\triangle ABC$. This is the unique charm of Intercenter Geometry, which is unmatched by Euclidean geometry and analytical geometry.


There are similar methods for space case (three-dimensional case).

\section{Concept of intersecting center of a triangle}\label{Sec4.1}
Intersecting center is a very important concept in this book.

\textbf {Finite plane}: a plane that does not contain infinity.

\textbf {Finite plane formed by a triangle}: a finite plane containing three vertices of the triangle.

The finite plane formed by $\triangle ABC$ is denoted by $P_{\triangle ABC}$.

\textbf{Intersecting center} (abbreviated as IC): If three lines passing through three vertices of a given triangle intersect at one point, the intersection is called the intersecting center (abbreviated as IC) of the triangle.

\begin{figure}[h]
	\begin{minipage}[t]{0.5\linewidth}
		\centering
		\includegraphics[totalheight=1.0in]{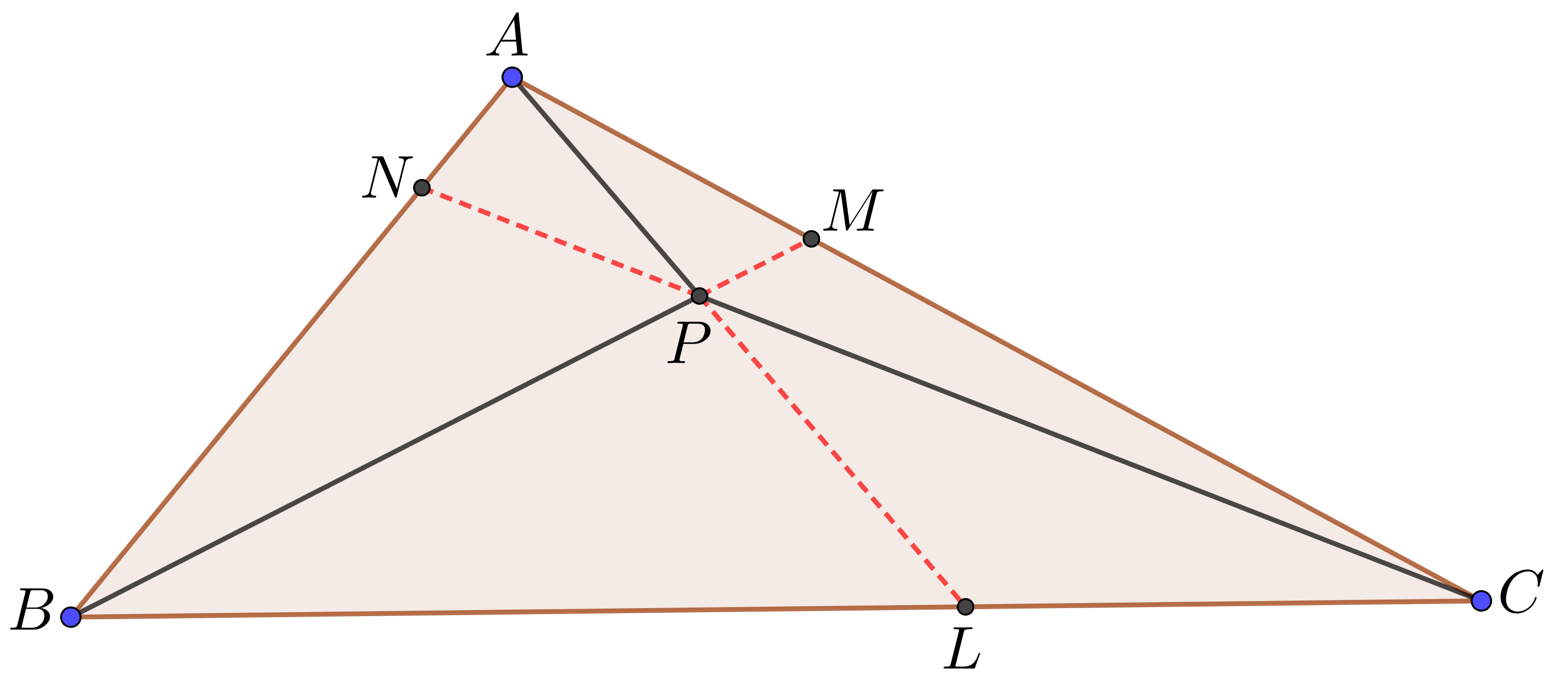}
		\caption{Inner intersecting center $P$} \label{fig:tu4.1.1}
	\end{minipage}
	\begin{minipage}[t]{0.5\linewidth}
		\centering
		\includegraphics[totalheight=1.1in]{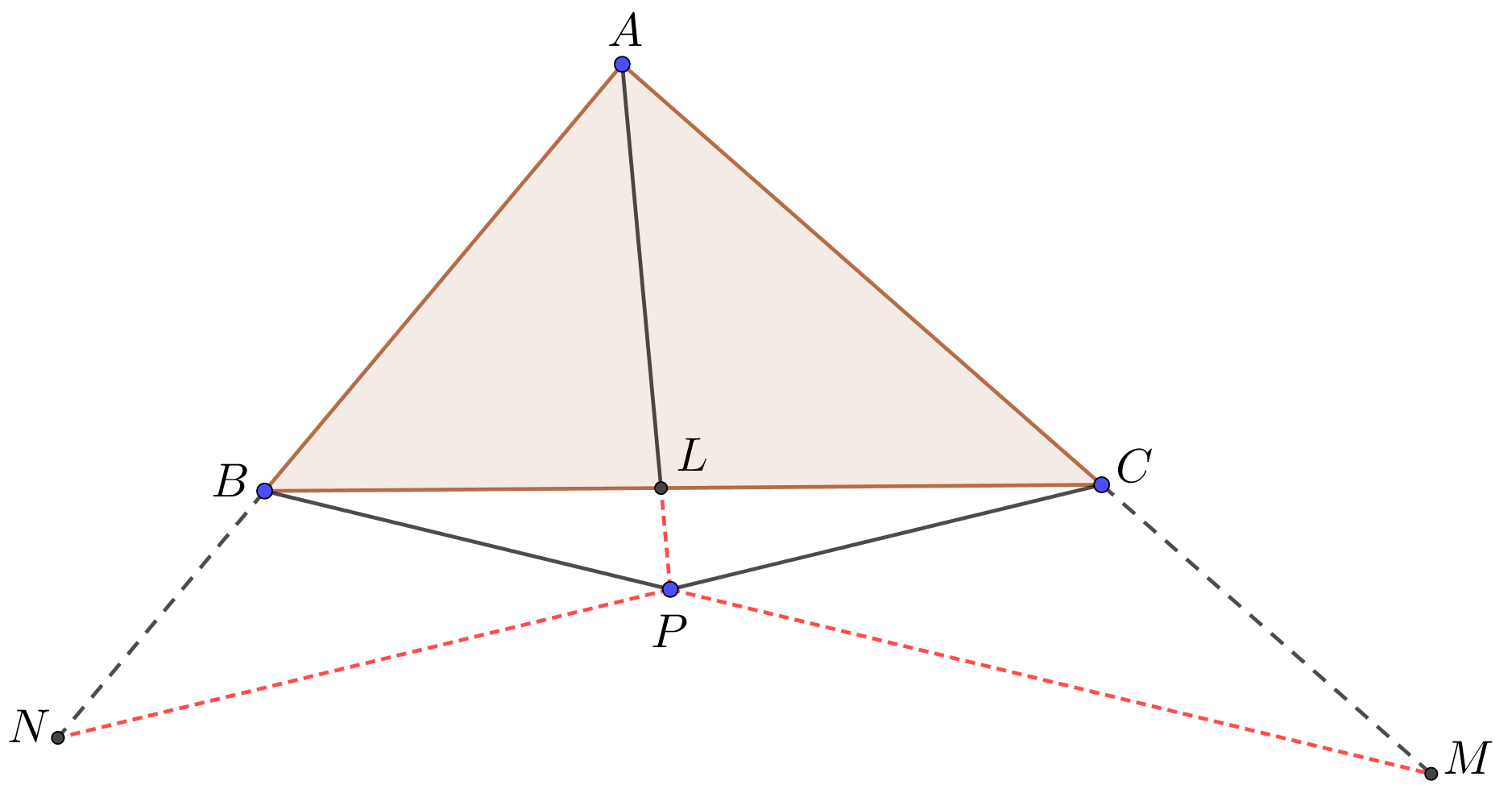}
		\caption{Outer intersecting center $P$} \label{fig:tu4.1.2}
	\end{minipage}
\end{figure}

For example, the point $P$ in Figure \ref{fig:tu4.1.1} and Figure \ref{fig:tu4.1.2} are the ICs of $\triangle ABC$, and
	\[P=\overleftrightarrow{AL}\cap \overleftrightarrow{BM}\cap \overleftrightarrow{CN}.\]
	
Obviously, $P\in {{P}_{\triangle ABC}}$. Not all the three lines passing through the three vertices of a given triangle can intersect at one point. For example, if the three lines passing through the three vertices of a given triangle are parallel to each other, there will be no intersection on the finite plane, or the intersection is at infinity. This book studies the intersecting center on a finite plane, the intersection at infinity needs to be ruled out.

\textbf{Singular intersecting center} (abbreviated as SIC): If at least one intersecting ratio (this concept will be introduced later) of an intersecting center of a triangle does not exist or is 0, the intersecting center is called singular intersecting center. 

The set of all singular intersecting centers of $\triangle ABC$ is denoted as $\pi _{ABC}^{\text{*}}$.


\textbf{Normal intersecting center} (abbreviated as NIC): The IC whose intersecting ratio exists and is not 0.

\textbf{Set of normal intersecting center}: The set of total NICs. 

The set of all NICs of $\triangle ABC$ is denoted as ${{\pi }_{ABC}}$. In the case of no confusion, the NIC is also referred to as IC.

\textbf{Set of intersecting center}: A set of normal intersecting centers.

\textbf{Inner intersecting center}: An intersecting center inside a triangle.

For example, the point $P$ in Figure \ref{fig:tu4.1.1} is the inner intersecting center.

\textbf{Outer intersecting center}: An intersecting center outside a triangle.

For example, the point $P$ in Figure \ref{fig:tu4.1.2} is the outer intersecting center.

\textbf{Line between vertex and intersecting center} (abbreviated as LVIC): The line connecting the vertex and the IC.
	
\textbf{Intersecting foot} (abbreviated as IF): The intersection of the LVIC and the opposite side of the vertex of a triangle.

In Figure \ref{fig:tu4.1.1} and Figure \ref{fig:tu4.1.2}, The points $L$, $M$, $N$ are IFs.

\section{Concept of intersecting ratio of a triangle}\label{Sec4.2}	
	The following concepts can be seen in Figure \ref{fig:tu4.1.1} and Figure \ref{fig:tu4.1.2}.

	\textbf{Intersecting ratio} (abbreviated as IR): Suppose that given a triangle $\triangle ABC$, the point $A$ is the starting point, the point $B$ is the terminal point, and the point $M\in \overleftrightarrow{AB}$, then $\lambda _{AB}^{P}$ is called the intersecting ratio of point $P$ corresponding to edge $AB$ and is defined by the following formula:
\[\lambda _{AB}^{P}=\frac{\overrightarrow{AN}}{\overrightarrow{NB}}.\]	
	
	Similarly, the intersecting ratio of point $P$ corresponding to edge $BC$ is defined by the following:
\[\lambda _{BC}^{P}=\frac{\overrightarrow{BL}}{\overrightarrow{LC}}.\]	
	
	The intersecting ratio of point $P$ corresponding to edge $CA$ is defined by the following:
\[\lambda _{CA}^{P}=\frac{\overrightarrow{CM}}{\overrightarrow{MA}}.\]	
	
	Intersecting ratio is an important concept proposed by the author. The intersecting ratio can be considered as line segment ratio, therefore $\lambda _{AB}^{P}\cdot \lambda _{BA}^{P}=1$. Intersecting ratio is a special ratio of line segment when studying the intersecting center of the triangle.

	In this book, the notation of intersecting ratio has a superscript, which describes the information of the corresponding point. For example, the superscript $P$ of $\lambda $ in the above definition represents the intersecting ratio of point $P$. The subscript of $\lambda$ represents the edge, for example $AB$ represents the $AB$ edge of $\triangle ABC$, and so on.

	In this book, there is no superscript on the ratio of line segment. Readers can compare the notations here with the notations in section \ref{Sec2.1}.

	The intersecting ratio of point $P$ and the fractional ratio of point $P$ are not only different in the form of notation, but also have different meanings. For the fractional ratio of line segment $AB$, the fractional point $P\in \overleftrightarrow{AB}$; For the intersecting ratio of $\triangle ABC$, the point of IC usually satisfies $P\notin \overleftrightarrow{AB}$.
	
	\textbf{Vector of intersecting ratio (abbreviated as vector of IR)}: Given a $\triangle ABC$ and an IC $P$ of the $\triangle ABC$, the vector of IR of the $\triangle ABC$ at point $P$ is defined by the following:	
	\[\boldsymbol{\lambda }_{\triangle ABC}^{P}=\left(\begin{matrix}
		\lambda _{AB}^{P} & \lambda _{BC}^{P} & \lambda _{CA}^{P}  \\
	\end{matrix} \right).\]
	
	The vector of IR is an ordered array of IRs. The vector of IR $\boldsymbol{\lambda }_{\triangle ABC}^{P}$ establishes a corresponding relationship with the point $P$ in space.

\section{Concept of triangular frame}\label{Sec4.3}
	\textbf{Triangular frame} (abbreviated as TF): Given a $\triangle ABC$, point $O$ is any point in space, point $O$ and three vectors $\overrightarrow{OA}$, $\overrightarrow{OB}$ and $\overrightarrow{OC}$ together form a frame, denoted as $\left( O;A,B,C \right)$, and $\left( O;A,B,C \right)$ is also denoted as \textbf{\bm{$\triangle ABC$} frame}. Point $O$ is called the origin of the frame; vectors $\overrightarrow{OA}$, $\overrightarrow{OB}$, $\overrightarrow{OC}$ are called frame vectors. Both the frame vectors and the triangular frame are simply called frame when there is no need to distinguish them strictly.

	Triangular frame $\left( O;A,B,C \right)$ can be regarded as a coordinate system.

	The position of origin of the triangular frame $\left( O;A,B,C \right)$ can be arbitrary, which is also called \textbf{free origin}.

	\textbf{Frame of plane triangle} (abbreviated as FPT): The triangular frame $\left( O;A,B,C \right)$ whose origin $O$ is on the plane of $\triangle ABC$.

	\textbf{Frame of intersecting center} (abbreviated as FIC): The origin $O$ of the frame is the IC of $\triangle ABC$.
	
	\textbf{Frame of circumcenter} (abbreviated as FC): A frame of plane triangle (FPT) whose origin $O$ is the circumcenter of $\triangle ABC$.

	The position of the origin $O$ of the FC cannot be arbitrary, but can only be the circumcenter of $\triangle ABC$. Similar concepts such as frame of centroid can be obtained.

	\textbf{Frame of edge} (abbreviated as FE): Given a $\triangle ABC$, for a vertex $A$ of $\triangle ABC$ and two vectors $\overrightarrow{AB}$, $\overrightarrow{AC}$, they form a frame together, this frame is denoted as $\left(A;B, C \right)$.

	The origin of $\left( A;B,C \right)$ is $A$, two vectors $\overrightarrow{AB}$ and $\overrightarrow{AC}$ of the $\left(A;B,C \right)$ coincide with both sides of $AB$ and $AC$ of the $\triangle ABC$, respectively.

	Similarly, we can define $\left( B;C,A \right)$ and $\left( C;A,B \right)$.

\section{Concept of vector from origin to intersecting center}\label{Sec4.4}	
	\textbf{Vector from origin to intersecting center} (abbreviated as VOIC): Given a $\triangle ABC$, the point $O$ is the origin of a triangular frame $\left( O;A,B,C \right)$ and the point $P$ is an intersecting center of $\triangle ABC$, then the vector $\overrightarrow{OP}$ is called the vector from origin to the intersecting center.

	VOIC $\overrightarrow{OP}$ is the position vector whose origin $O$ points to $P$ (the IC). Since the origin of the frame $O$ can be arbitrary, either not on the plane of $\triangle ABC$ (any point in space) or on the plane of $\triangle ABC$, the origin of the frame $O$ may also be an IC. The key point of the concept of a VOIC is to emphasize that the starting point of the vector $\overrightarrow{OP}$ must be the origin of triangular frame $\left( O;A,B,C \right)$ and the terminal point $P$ must be an intersecting center (IC) of $\triangle ABC$.

	\textbf{Vector of two intersecting centers} (abbreviated as VTICs): Given a $\triangle ABC$, point $O$ is the origin of a triangular frame $\left( O;A,B,C \right)$, point ${{P}_{1}}$ and ${{P}_{2}}$ are the two intersecting centers of $\triangle ABC$, then the vector $\overrightarrow{{{P}_{1}}{{P}_{2}}}$ is called the vector of two intersecting centers (abbreviated as VTICs).

	The VTICs is a vector with two intersecting centers. The position of the origin of the frame is arbitrary, that is, the free origin.
	
	The main point of VTICs concept is: for a vector $\overrightarrow{{{P}_{1}}{{P}_{2}}}$, its starting point ${{P}_{1}}$ and the terminal point ${{P}_{2}}$ are both ICs of $\triangle ABC$, while the position of the origin $O$ is arbitrary and the origin $O$ is usually both different from point ${{P}_{1}}$ and point ${{P}_{2}}$.

	Both the VOIC and the VTICs can be a linear combination of vectors by different frames, which will be discussed later.
	
	\textbf{Vector of vertex to intersecting center} (abbreviated as VVIC): The position vector from the vertex to intersecting center (IC).
	
	In Figure \ref{fig:tu4.1.1} and Figure \ref{fig:tu4.1.2}, $\overrightarrow{AP}$, $\overrightarrow{BP}$ and $\overrightarrow{CP}$ are VVICs.



\chapter{Calculation of intersecting ratio for the special intersecting centers of a triangle}\label{Ch5}
\thispagestyle{empty}

In Intercenter Geometry, the concept of IR is always important. This chapter deals with the calculation of the IR of some special centers of a triangle, which include the centroid, incenter, orthocenter, circumcenter and excenters of the triangle. Those formulas for calculating the IR of the special centers are often used in subsequent studies.

\section{Intersecting ratio of centroid}\label{Sec5.1}

The superscript $G$ of $\lambda$ denotes the centroid to distinguish it from other ICs. The formula for calculating the intersecting ratio of the centroid is as follows:
\[\lambda _{AB}^{G}=\lambda _{BC}^{G}=\lambda _{CA}^{G}=1.\]	

That is, the IR of the centroid is 1.

\section{Intersecting ratio of incenter}\label{Sec5.2}

The superscript $I$ of $\lambda$ denotes the incenter to distinguish it from others. The formula for calculating the intersecting ratio of incenter is as follows:
\begin{flalign}\label{Eq5.2.1}
	\lambda _{AB}^{I}=\frac{b}{a}, \lambda _{BC}^{I}=\frac{c}{b}, \lambda _{CA}^{I}=\frac{a}{c}.	
\end{flalign}

\begin{proof}
	Using $\lambda _{AB}^{I}$ as an example to prove the above formula. The proof processes for $\lambda _{BC}^{I}$ and $\lambda_{ CA}^{I}$ are similar.
	
	As shown in Figure \ref{fig:tu5.2.1}, let $CL$ be the bisector of  $\angle C$ in $\triangle ABC $, $LD\bot CA$ and $LE\bot BC$, and, according to the definition of the intersecting ratio, using the property of the angular bisector, we obtain:
	\[\lambda _{AB}^{I}=\frac{\overrightarrow{AL}}{\overrightarrow{LB}}=\frac{AL}{LB}=\frac{{{S}_{CAL}}}{{{S}_{CLB}}}=\frac{\frac{1}{2}CA\cdot LD}{\frac{1}{2}BC\cdot LE}=\frac{CA}{BC}=\frac{b}{a}.\]	
	
	\begin{figure}[h]
		\centering
		\includegraphics[totalheight=4cm]{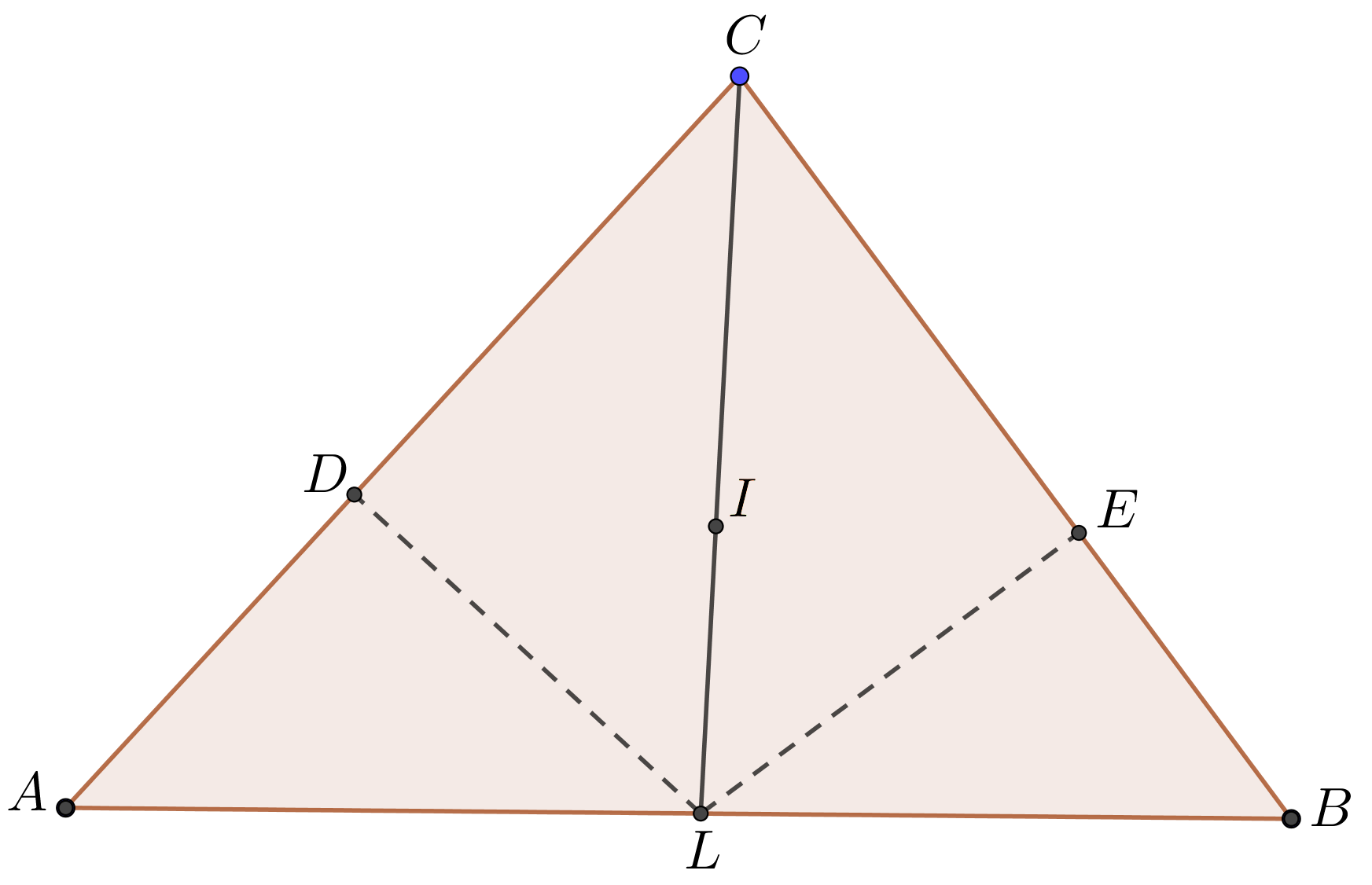}
		\caption{Calculating the intersecting ratio of incenter} \label{fig:tu5.2.1}
	\end{figure}
	
	Similar results are obtained:
	\begin{flalign*}
		\lambda _{BC}^{I}=\frac{c}{b}, \lambda _{CA}^{I}=\frac{a}{c}.	
	\end{flalign*}
\end{proof}
\hfill $\square$\par

\section{Intersecting ratio of orthocenter}\label{Sec5.3}
The formula for calculating the IR of orthocenter $H$ is as follows:

\begin{flalign}\label{Eq5.3.1}
	\lambda _{AB}^{H}=\frac{\tan B}{\tan A}, \lambda _{BC}^{H}=\frac{\tan C}{\tan B}, \lambda _{CA}^{H}=\frac{\tan A}{\tan C}.	
\end{flalign}

\begin{proof}
	As shown in Figure \ref{fig:tu5.3.1}, let $CL$ be the altitude of side $AB$ in $\triangle ABC$ and according to the definition of IR, we use the sine theorem to get the following result:	
	\[\begin{aligned}
	\lambda _{AB}^{H}
	&=\frac{\overrightarrow{AL}}{\overrightarrow{LB}}=\frac{AL}{LB}=\frac{CA\cos A}{BC\cos B}=\frac{b\cos A}{a\cos B} \\ 
		& =\frac{2R\sin B\cos A}{2R\sin A\cos B}=\frac{\tan B}{\tan A}.  
	\end{aligned}\]	
	
	Where $R$ is the radius of the circumscribed circle of $\triangle ABC$. Similar results can be obtained as follows:	
	\[\lambda _{BC}^{H}=\frac{c\cos B}{b\cos C}=\frac{\tan C}{\tan B},\]	
	\[\lambda _{CA}^{H}=\frac{a\cos C}{c\cos A}=\frac{\tan A}{\tan C}.\]	
	
	\begin{figure}[h]
		\centering
		\includegraphics[totalheight=4cm]{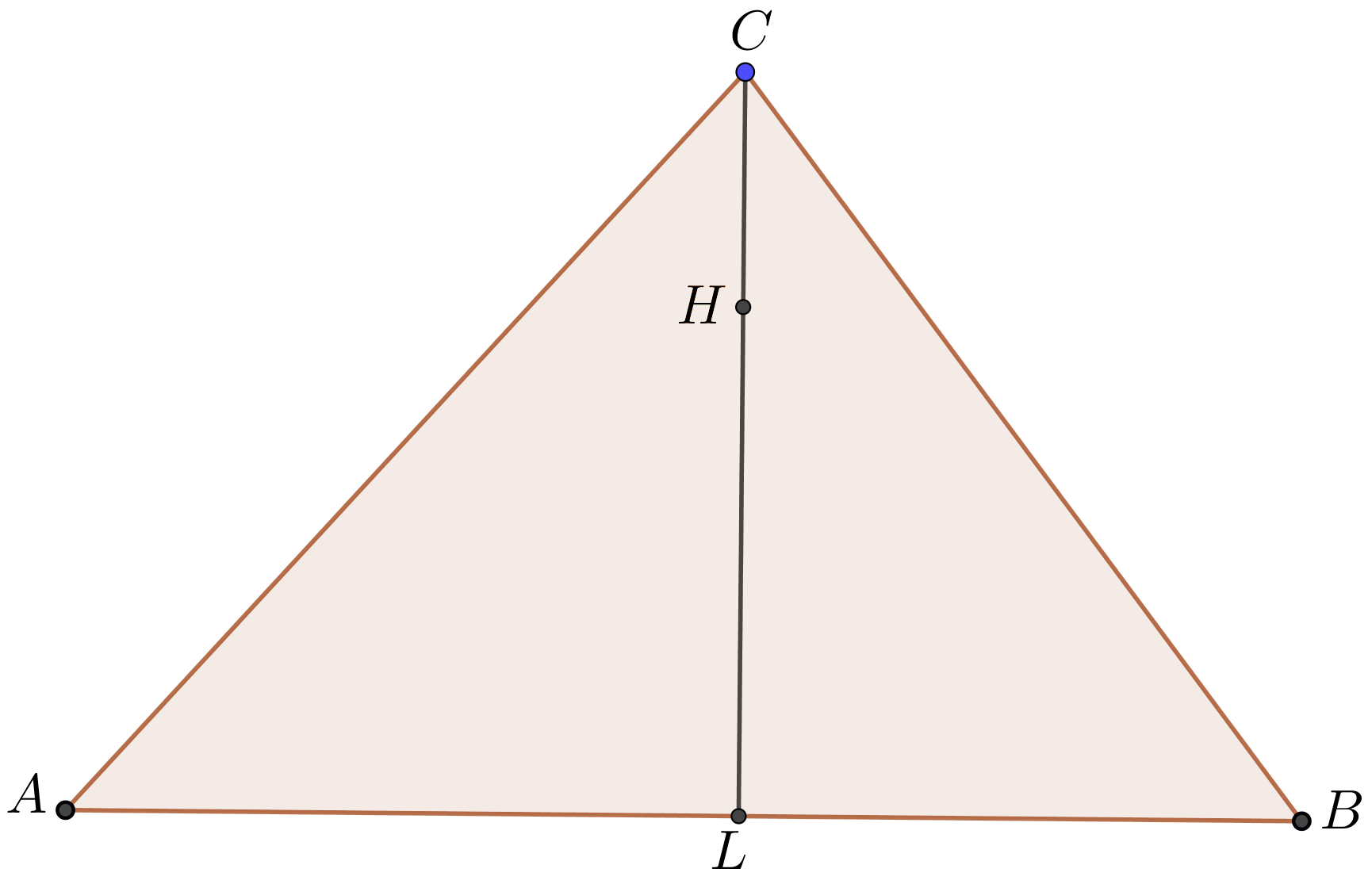}
		\caption{Diagram for calculating the IR of orthocenter}\label{fig:tu5.3.1}
	\end{figure}


%

\end{proof}
\hfill $\square$\par

%
\section{Intersecting ratio of circumcenter}\label{Sec5.4}
The formula for calculating the IR of circumcenter $Q$ is as follows:
\begin{flalign}\label{Eq5.4.1}
	\lambda _{AB}^{Q}=\frac{\sin 2B}{\sin 2A}, \lambda _{BC}^{Q}=\frac{\sin 2C}{\sin 2B}, \lambda _{CA}^{Q}=\frac{\sin 2A}{\sin 2C}.	
\end{flalign}

\begin{figure}[h]
	\centering
	\includegraphics[totalheight=4cm]{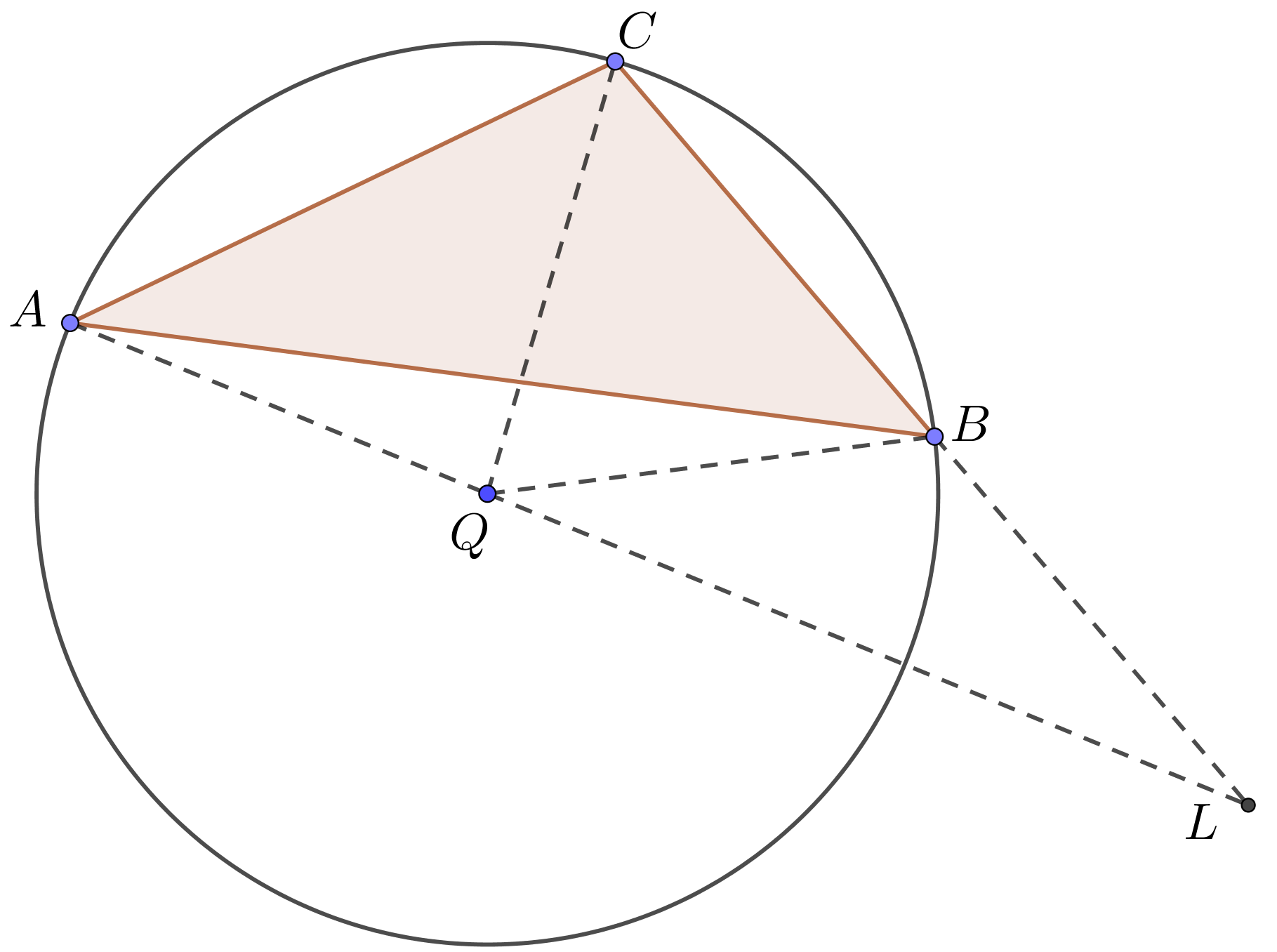}
	\caption{Diagram for calculating the IR of circumcenter} \label{fig:tu5.4.1}
\end{figure}

\begin{proof}
	As shown in Figure \ref{fig:tu5.4.1}, using $\lambda _{BC}^{Q}$ as an example to prove the above formula, $\lambda_{CA}^{Q}$ and $\lambda_{AB}^{Q}$ demonstrate similar procedures.
	\[\lambda _{BC}^{Q}=\frac{\overrightarrow{BL}}{\overrightarrow{LC}}=\frac{\overrightarrow{AB}\times \overrightarrow{AL}}{\overrightarrow{AL}\times \overrightarrow{AC}}=\frac{\overrightarrow{QB}\times \overrightarrow{QL}}{\overrightarrow{QL}\times \overrightarrow{QC}}=\frac{\overrightarrow{AB}\times \overrightarrow{AL}-\overrightarrow{QB}\times \overrightarrow{QL}}{\overrightarrow{AL}\times \overrightarrow{AC}-\overrightarrow{QL}\times \overrightarrow{QC}}.\]
	
	The numerator is
	\[\begin{aligned}
		& \overrightarrow{AB}\times \overrightarrow{AL}-\overrightarrow{QB}\times \overrightarrow{QL}=\left( \overrightarrow{AQ}+\overrightarrow{QB} \right)\times \overrightarrow{AL}-\overrightarrow{QB}\times \overrightarrow{QL} \\ 
		& =\overrightarrow{QB}\times \overrightarrow{AL}-\overrightarrow{QB}\times \overrightarrow{QL}=\overrightarrow{QB}\times \left( \overrightarrow{AL}-\overrightarrow{QL} \right)=\overrightarrow{QB}\times \overrightarrow{AQ}. \\ 
	\end{aligned}\]
	
	The denominator is
	\[\begin{aligned}
		& \overrightarrow{AL}\times \overrightarrow{AC}-\overrightarrow{QL}\times \overrightarrow{QC}=\overrightarrow{AL}\times \left( \overrightarrow{AQ}+\overrightarrow{QC} \right)-\overrightarrow{QL}\times \overrightarrow{QC} \\ 
		& =\overrightarrow{AL}\times \overrightarrow{QC}-\overrightarrow{QL}\times \overrightarrow{QC}=\left( \overrightarrow{AL}-\overrightarrow{QL} \right)\times \overrightarrow{QC}=\overrightarrow{AQ}\times \overrightarrow{QC}. \\ 
	\end{aligned}\]
	
	Therefore,
	\[\lambda _{BC}^{Q}=\frac{\overrightarrow{QB}\times \overrightarrow{AQ}}{\overrightarrow{AQ}\times \overrightarrow{QC}}=\frac{\overrightarrow{QB}\times \overrightarrow{QA}}{\overrightarrow{QA}\times \overrightarrow{QC}},\]	
	i.e.
	\[\lambda _{BC}^{Q}=\frac{\overrightarrow{QB}\times \overrightarrow{QA}}{\overrightarrow{QA}\times \overrightarrow{QC}}=\frac{\frac{1}{2}\left| \overrightarrow{QB} \right|\cdot \left| \overrightarrow{QA} \right|\sin \angle BQA\cdot \mathbf{n}}{\frac{1}{2}\left| \overrightarrow{QA} \right|\cdot \left| \overrightarrow{QC} \right|\sin \angle AQC\cdot \mathbf{n}}=\frac{\sin \angle BQA}{\sin \angle AQC}.\]
	
	
	According to the concepts of directed angle and vector multiplication (see Chapter \ref{SanjiaoxingXianguanZhishi}), the above formula is valid only in the same reference frame (right-handed system). That is, both $\angle BQA$ and $\angle AQC$ are directional angles, and the starting edge of $\angle BQA$ is $\overrightarrow{QB}$, the ending edge is $\overrightarrow{QA}$, and $\angle BQA$ is the angle from the starting edge $\overrightarrow{QB}$ to the ending edge $\overrightarrow{QA}$ in a counterclockwise direction; The starting edge of $\angle AQC$ is $\overrightarrow{QA}$, the ending edge is $\overrightarrow{QC}$, and $\angle AQC$ is the angle from the starting edge $\overrightarrow{QA}$ to the ending edge $\overrightarrow{QC}$ in a counterclockwise direction.
	
	As shown in Figure \ref{fig:tu5.4.1}, it is specified that the initial side of the directed angle $\angle BQA$ is $\overrightarrow{QB}$, the terminal side is $\overrightarrow{QA}$, and the reference positive direction of $\angle BQA$ is counterclockwise, that is, the angle takes a positive value when it rotates counterclockwise.
	
	We specify that the starting point of the arc $\overset\frown{BA}$ is $B$, and the ending point is $A$. The arc $\overset\frown{BA}$ moves from the starting point $B$ to the ending point $A$, and the reference direction of the movement is consistent with the reference direction of the directed angle $\angle BQA$. Such arc $\overset\frown{BA}$ is called the directed arc. In the following discussion about directed arc, it is assumed that the arc does not include the starting point and ending point of the arc.
	
	Since $\odot Q$ is the circumscribed circle of $\triangle ABC$, then there must be $C\in \overset\frown{BA}$ or $C\notin \overset\frown{BA}$, and points $C\ne A$, $C\ne B$. 
	
	When $C\notin \overset\frown{BA}$, First, we consider the directed angle. At this point, the center of the circle $Q$ and the point $C$ are on the opposite side of the segment $AB$. By using the relationship between the center angle and the angle of circumference, select an auxiliary point $D\in \overset\frown{AB}$, we have:	
	\[\angle BQA=2\angle BDA=2\left( \pi -\angle ACB \right)=2\left( \pi -\angle B \right).\]
	
	Secondly, we consider the directed angle $\angle AQC$.
	\[\angle AQC+\angle CQA=2\pi. \]
	
	Since $C\in \overset\frown{BA}$, then $B\notin \overset\frown{CA}$, the center of the circle $Q$ and point $B$ are on the same side of the line segment $AC$. Therefore, we get:
	\[\angle CQA=2\angle CBA,\]
	therefore,
	\[\angle AQC=2\pi -\angle CQA=2\left( \pi -\angle CBA \right)=2\left( \pi -\angle B \right).\]
	
	Therefore,
	\[\frac{\sin \angle BQA}{\sin \angle AQC}=\frac{\sin 2\left( \pi -\angle ACB \right)}{\sin 2\left( \pi -\angle CBA \right)}=\frac{\sin 2\angle ACB}{\sin 2\angle CBA}=\frac{\sin 2C}{\sin 2B}.\]
	
		
	If $C\in \overset\frown{AB}$ (not shown in the Figure \ref{fig:tu5.4.1}), we consider the directed angle $\angle BQA$.
	\[\angle BQA=2\angle BCA=2\angle C.\]
		
	
	Secondly, we consider the directed angle $\angle AQC$, the point $Q$ and the point $C$ are on the same side of the line segment $AB$, and hence:
	\[\angle AQC=2\angle ABC=2\angle B,\]
	therefore,
	\[\frac{\sin \angle BQA}{\sin \angle AQC}=\frac{\sin 2\angle ACB}{\sin 2\angle CBA}=\frac{\sin 2C}{\sin 2B}.\]
	
	That is, in either case, the following formula is obtained:
	\[\lambda _{BC}^{Q}=\frac{\sin \angle BQA}{\sin \angle AQC}=\frac{\sin 2C}{\sin 2B}.\]
		
	
	The above proof is based on the assumption that the edge of the triangle does not pass through the center of the circle. If an edge passes through the center of the circle, the same result will be obtained.
\end{proof}
\hfill $\square$\par

\section{Intersecting ratio of excenters}\label{Sec5.5}

The formulas for calculating the IR of the excenters are as follows:

The IR of the excenter corresponding to $\angle A$:
\begin{flalign}\label{Eq5.5.1}
	\lambda _{AB}^{{{E}_{A}}}=-\frac{b}{a}, \lambda _{BC}^{{{E}_{A}}}=\frac{c}{b}, \lambda _{CA}^{{{E}_{A}}}=-\frac{a}{c}.		
\end{flalign}

The IR of the excenter corresponding to $\angle B$:
\begin{flalign*}
	\lambda _{AB}^{{{E}_{B}}}=-\frac{b}{a}, \lambda _{BC}^{{{E}_{B}}}=-\frac{c}{b}, \lambda _{CA}^{{{E}_{B}}}=\frac{a}{c}.			
\end{flalign*}

The IR of the excenter corresponding to $\angle C$:
\begin{flalign}\label{Eq5.5.2}
	\lambda _{AB}^{{{E}_{C}}}=\frac{b}{a}, \lambda _{BC}^{{{E}_{C}}}=-\frac{c}{b}, \lambda _{CA}^{{{E}_{C}}}=-\frac{a}{c}.			
\end{flalign}

\begin{proof}
	
	We only prove the IR formula of excenter corresponds to $\angle A$, the rest of the proof process are similar.
	
	\begin{figure}[h]
		\centering
		\includegraphics[totalheight=4cm]{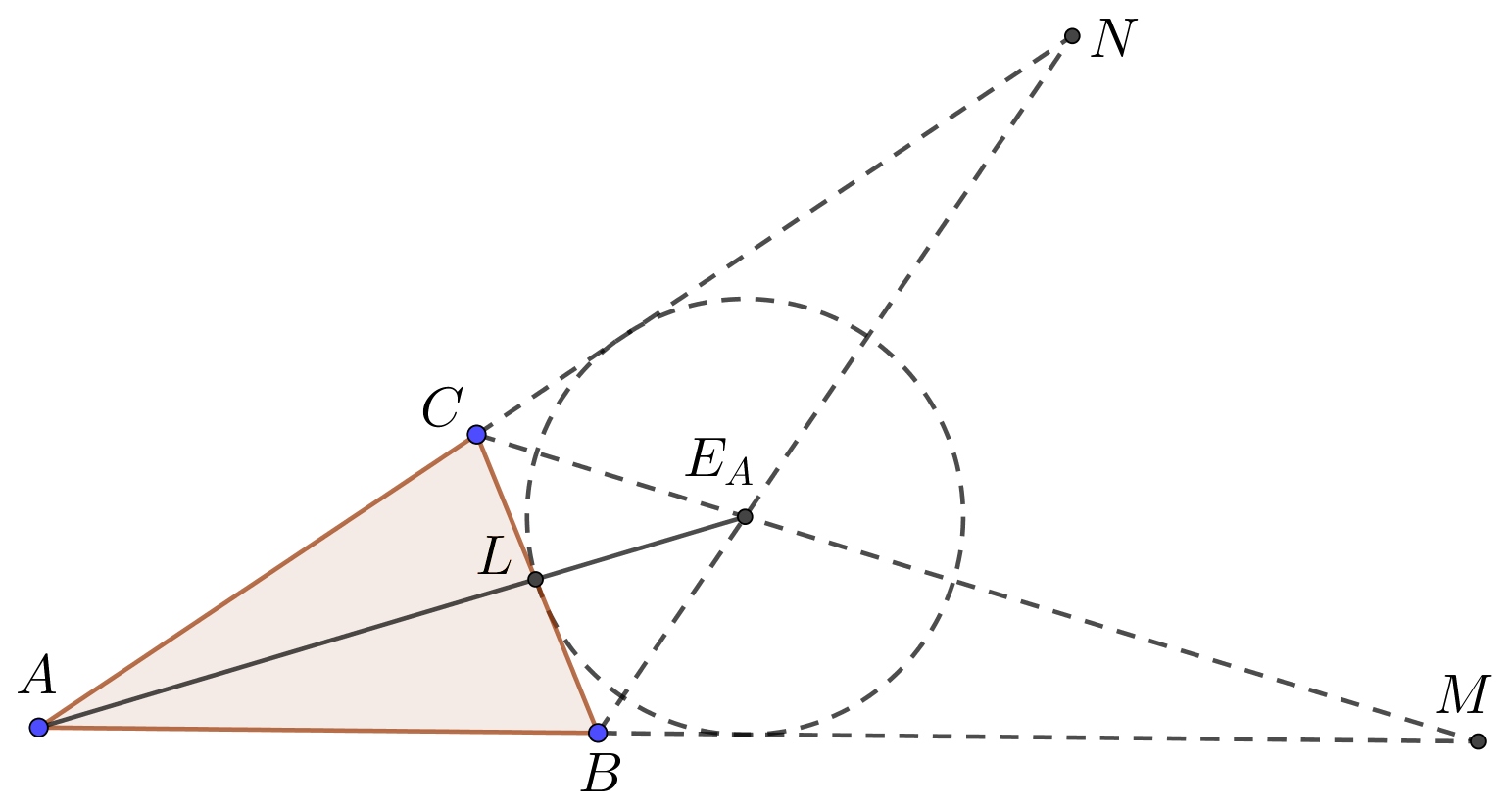}
		\caption{Diagram used to prove the IR of excenter ${{E}_{A}}$ corresponding to $\angle A$}  \label{fig:tu5.5.1}
	\end{figure}	
	
	As shown in Figure \ref{fig:tu5.5.1}, for $\lambda _{AB}^{{{E}_{A}}}$, we have
	\[\begin{aligned}
		\lambda _{AB}^{{{E}_{A}}}& =\frac{\overrightarrow{AM}}{\overrightarrow{MB}}=\frac{\overrightarrow{CA}\times \overrightarrow{CM}}{\overrightarrow{CM}\times \overrightarrow{CB}}=\frac{\frac{1}{2}\left| \overrightarrow{CA} \right|\cdot \left| \overrightarrow{CM} \right|\sin \angle ACM\cdot \mathbf{n}}{\frac{1}{2}\left| \overrightarrow{CM} \right|\cdot \left| \overrightarrow{CB} \right|\sin \angle MCB\cdot \mathbf{n}} \\ 
		& =\frac{\left| \overrightarrow{CA} \right|\sin \angle ACM}{\left| \overrightarrow{CB} \right|\sin \angle MCB}=\frac{CA\sin \angle ACM}{CB\sin \angle MCB},  
	\end{aligned}\]
	while
	\[\angle ACM+\angle MCN=\angle ACN=\pi, \]
	\[\angle ACM=\pi -\angle MCN.\]
	
	Since $M$ is the point on the bisector of the exterior angle of $\angle C$, we have:
	\[\angle MCN=-\angle MCB,\]
	then
	\[\angle ACM=\pi +\angle MCB,\]
	therefore,
	\[\lambda _{AB}^{{{E}_{A}}}=\frac{CA\sin \angle ACM}{CB\sin \angle MCB}=\frac{CA\sin \left( \pi +\angle MCB \right)}{CB\sin \angle MCB}=-\frac{CA}{CB}=-\frac{b}{a}.\]
	
	For $\lambda _{BC}^{{{E}_{A}}}$, we have:
	\[\begin{aligned}
		\lambda _{BC}^{{{E}_{A}}}& =\frac{\overrightarrow{BL}}{\overrightarrow{LC}}=\frac{\overrightarrow{AB}\times \overrightarrow{AL}}{\overrightarrow{AL}\times \overrightarrow{AC}}=\frac{\frac{1}{2}\left| \overrightarrow{AB} \right|\cdot \left| \overrightarrow{AL} \right|\sin \angle BAL\cdot \mathbf{n}}{\frac{1}{2}\left| \overrightarrow{AL} \right|\cdot \left| \overrightarrow{AC} \right|\sin \angle LAC\cdot \mathbf{n}} \\ 
		& =\frac{\left| \overrightarrow{AB} \right|\sin \angle BAL}{\left| \overrightarrow{AC} \right|\sin \angle LAC}=\frac{AB}{AC}=\frac{c}{b}.  
	\end{aligned}\]
	
	For $\lambda _{CA}^{{{E}_{A}}}$, we have:
	\[\begin{aligned}
		\lambda _{CA}^{{{E}_{A}}}& =\frac{\overrightarrow{CN}}{\overrightarrow{NA}}=\frac{\overrightarrow{BC}\times \overrightarrow{BN}}{\overrightarrow{BN}\times \overrightarrow{BA}}=\frac{\frac{1}{2}\left| \overrightarrow{BC} \right|\cdot \left| \overrightarrow{BN} \right|\sin \angle CBN\cdot \mathbf{n}}{\frac{1}{2}\left| \overrightarrow{BN} \right|\cdot \left| \overrightarrow{BA} \right|\sin \angle NBA\cdot \mathbf{n}} \\ 
		& =\frac{\left| \overrightarrow{BC} \right|\sin \angle CBN}{\left| \overrightarrow{BA} \right|\sin \angle NBA}=\frac{BC\sin \angle CBN}{BA\sin \angle NBA},  
	\end{aligned}\]
	and
	\[\angle MBN+\angle NBA=\angle MBA=\pi,\]
	\[\angle NBA=\pi -\angle MBN,\]
	\[\angle MBN=\angle NBC=-\angle CBN,\]
	then
	\[\lambda _{CA}^{{{E}_{A}}}=\frac{BC\sin \angle CBN}{BA\sin \angle NBA}=-\frac{BC}{BA}=-\frac{a}{c}.\]
\end{proof}
\hfill $\square$\par



\chapter{Vector of intersecting center and frame equation}\label{Ch6}
\thispagestyle{empty}

In this chapter, we mainly study the vector of intersecting center (VIC) in triangular frame, which is the basis of Intercenter Geometry and plays a very important role. The vector of intersecting center is mainly divided into the vector from origin to intersecting center (abbreviated as VOIC) and the vector of two intersecting centers (abbreviated as VTICs). In this chapter, I prove two important theorems: theorem \ref{thm:Thm6.1.1} and theorem \ref{thm:Thm6.2.1}. These two theorems are the core theorems of this book and are very important.

The VIC can be expressed by frame system and frame components, which is the basis of other geometric quantities (such as distance, etc.). Many applications in this book are based on these two theorems in this chapter.

\section{Theorem of vector from origin to intersecting center }\label{Sec6.1}

\begin{theorem}{Vector from origin to an intersecting center, Daiyuan Zhang}{Thm6.1.1}\label{Thm6.1.1} 
	Suppose that the three given points $A$, $B$, $C$ do not coincide with each other, the intersecting center $P\in {{\pi }_{ABC}}$, the point $O$ is an arbitrary point in space ($O\in {{\mathbb{R}}^{3}}$), then the vector $\overrightarrow {OP}$ can be the linear and unique combination of the vectors of the triangular frame $\left( O;A,B,C \right)$, i.e.
	\begin{equation}\label{Eq6.1.1}
		\overrightarrow{OP}=\alpha _{A}^{P}\overrightarrow{OA}+\alpha _{B}^{P}\overrightarrow{OB}+\alpha _{C}^{P}\overrightarrow{OC},		
	\end{equation}
	\[\alpha _{A}^{P}+\alpha _{B}^{P}+\alpha _{C}^{P}=1.\]	
	
	Each of $\alpha _{A}^{P}$, $\alpha _{B}^{P}$ and $\alpha _{C}^{P}$ is uniquely determined, where:
	\begin{equation}\label{Eq6.1.2}
		\alpha _{A}^{P}=\frac{1+\lambda _{BA}^{P}+\lambda _{CA}^{P}}{\left( 1+\lambda _{BA}^{P}+\lambda _{CA}^{P} \right)+\left( 1+\lambda _{CB}^{P}+\lambda _{AB}^{P} \right)+\left( 1+\lambda _{AC}^{P}+\lambda _{BC}^{P} \right)},		
	\end{equation}
	\begin{equation}\label{Eq6.1.3}
		\alpha _{B}^{P}=\frac{1+\lambda _{CB}^{P}+\lambda _{AB}^{P}}{\left( 1+\lambda _{BA}^{P}+\lambda _{CA}^{P} \right)+\left( 1+\lambda _{CB}^{P}+\lambda _{AB}^{P} \right)+\left( 1+\lambda _{AC}^{P}+\lambda _{BC}^{P} \right)},	
	\end{equation}
	\begin{equation}\label{Eq6.1.4}
		\alpha _{C}^{P}=\frac{1+\lambda _{AC}^{P}+\lambda _{BC}^{P}}{\left( 1+\lambda _{BA}^{P}+\lambda _{CA}^{P} \right)+\left( 1+\lambda _{CB}^{P}+\lambda _{AB}^{P} \right)+\left( 1+\lambda _{AC}^{P}+\lambda _{BC}^{P} \right)};	
	\end{equation}
	or
	\[\alpha _{A}^{P}=\frac{1}{1+\lambda _{AB}^{P}+\lambda _{AC}^{P}},\]	
	\[\alpha _{B}^{P}=\frac{\lambda _{AB}^{P}}{1+\lambda _{AB}^{P}+\lambda _{AC}^{P}},\]	
	\[\alpha _{C}^{P}=\frac{\lambda _{AC}^{P}}{1+\lambda _{AB}^{P}+\lambda _{AC}^{P}};\]	
	or
	\[\alpha _{B}^{P}=\frac{1}{1+\lambda _{BC}^{P}+\lambda _{BA}^{P}},\]	
	\[\alpha _{C}^{P}=\frac{\lambda _{BC}^{P}}{1+\lambda _{BC}^{P}+\lambda _{BA}^{P}},\]	
	\[\alpha _{A}^{P}=\frac{\lambda _{BA}^{P}}{1+\lambda _{BC}^{P}+\lambda _{BA}^{P}};\]	
	or
	\[\alpha _{C}^{P}=\frac{1}{1+\lambda _{CA}^{P}+\lambda _{CB}^{P}},\]	
	\[\alpha _{A}^{P}=\frac{\lambda _{CA}^{P}}{1+\lambda _{CA}^{P}+\lambda _{CB}^{P}},\]	
	\[\alpha _{B}^{P}=\frac{\lambda _{CB}^{P}}{1+\lambda _{CA}^{P}+\lambda _{CB}^{P}}.\]	
	
	In formulas $(\ref{Eq6.1.2})$, $(\ref{Eq6.1.3})$ and $(\ref{Eq6.1.4})$, $\lambda _{AB}^{P}$ is the IR of point $P$ corresponding to edge $AB$, i.e.
	$\lambda _{AB}^{P}={\overrightarrow{AN}}/{\overrightarrow{NB}}\;$ $($ see Figure $\ref{fig:tu6.1.1}$ $)$. The meanings of the remaining notations are analogous.
\end{theorem}

\begin{proof}
	As shown in Figure \ref{fig:tu6.1.1}, let point $P$ be the IC of $\triangle ABC$, and $AL$, $BM$ and $CN$ be the LVIC of point $P$. For any point $O\in {{\mathbb{R}}^{3}}$, according to the notation of formula (\ref{Eq2.1.1}), the following notation can be naturally introduced:	
	\[\left\{ \begin{aligned}
		& \overrightarrow{AN}=\lambda _{AB}^{P}\overrightarrow{NB} \\ 
		& \overrightarrow{AM}=\lambda _{AC}^{P}\overrightarrow{MC} \\ 
		& \overrightarrow{BP}={{\lambda }_{BM}}\overrightarrow{PM} \\ 
		& \overrightarrow{CP}={{\lambda }_{CN}}\overrightarrow{PN}. \\ 
	\end{aligned} \right.\]	
	
	The $\lambda$ with superscript represents the IR, and the $\lambda$ without superscript represents the ratio of line segment.
	\begin{figure}[h]
		\centering
		\includegraphics[totalheight=6cm]{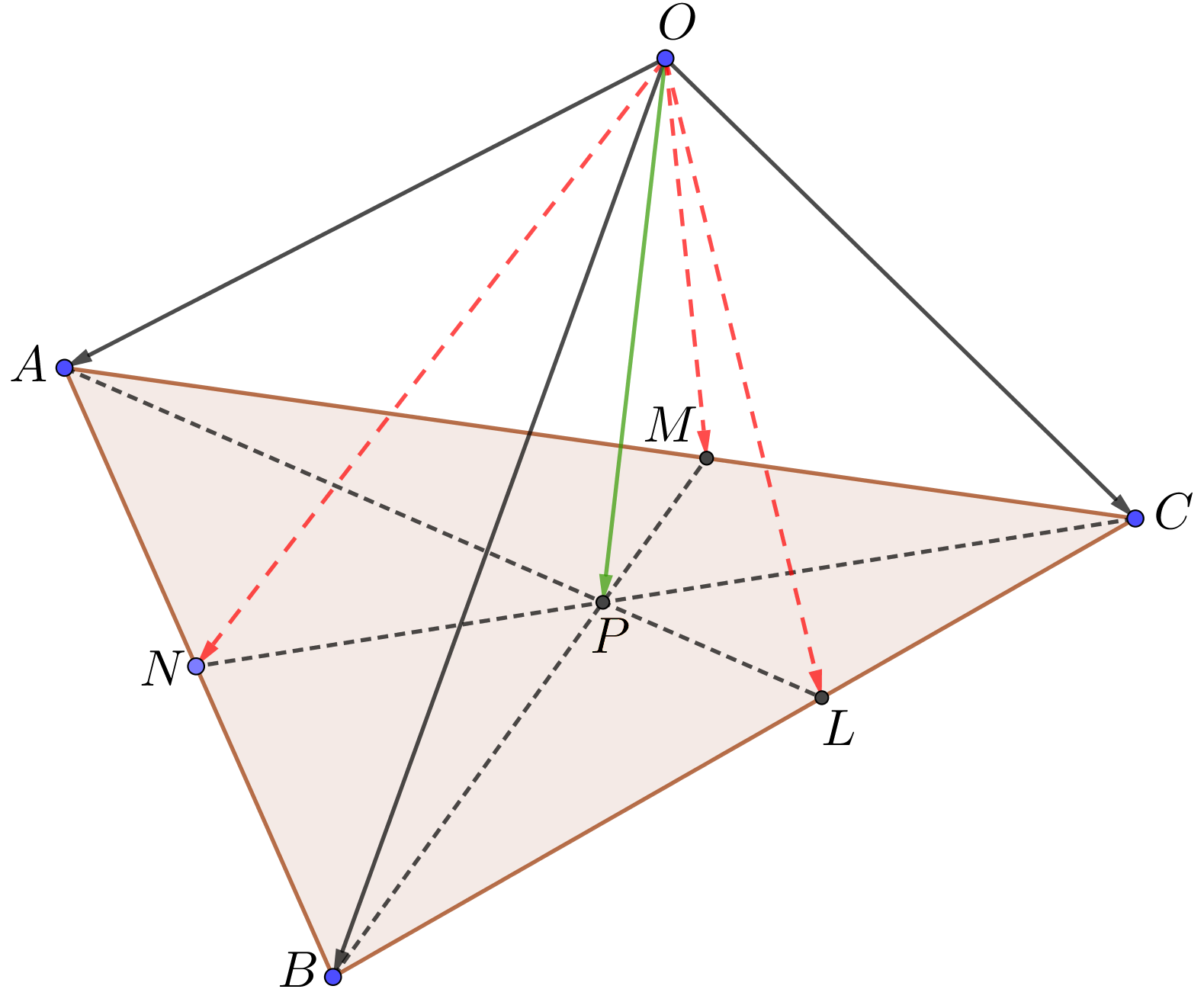}
		\caption{
		Figure for proving theorem \ref{thm:Thm6.1.1}} \label{fig:tu6.1.1}
	\end{figure}

	By using the properties of collinear vector, the following results can be obtained:
	\begin{equation}\label{Eq6.1.5}
		\overrightarrow{OP}=\frac{\overrightarrow{OB}+{{\lambda }_{BM}}\overrightarrow{OM}}{1+{{\lambda }_{BM}}},
	\end{equation}
	\begin{equation}\label{Eq6.1.6}
		\overrightarrow{OP}=\frac{\overrightarrow{OC}+{{\lambda }_{CN}}\overrightarrow{ON}}{1+{{\lambda }_{CN}}}.
	\end{equation}
	
	But
	\begin{equation}\label{Eq6.1.7}
		\overrightarrow{OM}=\frac{\overrightarrow{OA}+\lambda _{AC}^{P}\overrightarrow{OC}}{1+\lambda _{AC}^{P}},	
	\end{equation}
	\begin{equation}\label{Eq6.1.8}
		\overrightarrow{ON}=\frac{\overrightarrow{OA}+\lambda _{AB}^{P}\overrightarrow{OB}}{1+\lambda _{AB}^{P}}.
	\end{equation}
	
	So according to the formula (\ref{Eq6.1.5}) and (\ref{Eq6.1.7}), the following result can be obtained:
	\begin{equation}\label{Eq6.1.9}
		\begin{aligned}
			\overrightarrow{OP}&=\frac{\overrightarrow{OB}+{{\lambda }_{BM}}\overrightarrow{OM}}{1+{{\lambda }_{BM}}}=\frac{\overrightarrow{OB}+{{\lambda }_{BM}}\frac{\overrightarrow{OA}+\lambda _{AC}^{P}\overrightarrow{OC}}{1+\lambda _{AC}^{P}}}{1+{{\lambda }_{BM}}} \\ 
			& =\frac{{{\lambda }_{BM}}\overrightarrow{OA}+\left( 1+\lambda _{AC}^{P} \right)\overrightarrow{OB}+\lambda _{AC}^{P}{{\lambda }_{BM}}\overrightarrow{OC}}{\left( 1+\lambda _{AC}^{P} \right)\left( 1+{{\lambda }_{BM}} \right)} \\ 
			& ={{\gamma }_{1}}\overrightarrow{OA}+{{\gamma }_{2}}\overrightarrow{OB}+{{\gamma }_{3}}\overrightarrow{OC},  
		\end{aligned}
	\end{equation}
	where
	\begin{equation}\label{Eq6.1.10}
		\left\{ \begin{aligned}
			& {{\gamma }_{1}}=\frac{{{\lambda }_{BM}}}{\left( 1+\lambda _{AC}^{P} \right)\left( 1+{{\lambda }_{BM}} \right)} \\ 
			& {{\gamma }_{2}}=\frac{1+\lambda _{AC}^{P}}{\left( 1+\lambda _{AC}^{P} \right)\left( 1+{{\lambda }_{BM}} \right)} \\ 
			& {{\gamma }_{3}}=\frac{\lambda _{AC}^{P}{{\lambda }_{BM}}}{\left( 1+\lambda _{AC}^{P} \right)\left( 1+{{\lambda }_{BM}} \right)}. \\ 
		\end{aligned} \right.
	\end{equation}
	
	Obviously,
	\[{{\gamma }_{1}}+{{\gamma }_{2}}+{{\gamma }_{3}}=1.\]	
	
	According to the formula (\ref{Eq6.1.6}) and (\ref{Eq6.1.8}), the following result can be obtained:
	\begin{equation}\label{Eq6.1.11}
		\begin{aligned}
			\overrightarrow{OP}&=\frac{\overrightarrow{OC}+{{\lambda }_{CN}}\overrightarrow{ON}}{1+{{\lambda }_{CN}}}=\frac{\overrightarrow{OC}+{{\lambda }_{CN}}\frac{\overrightarrow{OA}+\lambda _{AB}^{P}\overrightarrow{OB}}{1+\lambda _{AB}^{P}}}{1+{{\lambda }_{CN}}} \\ 
			& =\frac{{{\lambda }_{CN}}\overrightarrow{OA}+\lambda _{AB}^{P}{{\lambda }_{CN}}\overrightarrow{OB}+\left( 1+\lambda _{AB}^{P} \right)\overrightarrow{OC}}{\left( 1+\lambda _{AB}^{P} \right)\left( 1+{{\lambda }_{CN}} \right)} \\ 
			& ={{\mu }_{1}}\overrightarrow{OA}+{{\mu }_{2}}\overrightarrow{OB}+{{\mu }_{3}}\overrightarrow{OC},  
		\end{aligned}
	\end{equation}
	where
	\begin{equation}\label{Eq6.1.12}
		\left\{ \begin{aligned}
			& {{\mu }_{1}}=\frac{{{\lambda }_{CN}}}{\left( 1+\lambda _{AB}^{P} \right)\left( 1+{{\lambda }_{CN}} \right)} \\ 
			& {{\mu }_{2}}=\frac{\lambda _{AB}^{P}{{\lambda }_{CN}}}{\left( 1+\lambda _{AB}^{P} \right)\left( 1+{{\lambda }_{CN}} \right)} \\ 
			& {{\mu }_{3}}=\frac{1+\lambda _{AB}^{P}}{\left( 1+\lambda _{AB}^{P} \right)\left( 1+{{\lambda }_{CN}} \right)}. \\ 
		\end{aligned} \right. 
	\end{equation}
	
	Obviously,
	\[{{\mu }_{1}}+{{\mu }_{2}}+{{\mu }_{3}}=1.\]	
	
	For the $\triangle ABC$, it is impossible for the three endpoints of the vector $\overrightarrow{OA}$, $\overrightarrow{OB}$ and $\overrightarrow{OC}$ to be collinear. According to the formulas (\ref{Eq6.1.9}) and (\ref{Eq6.1.11}) and the unique representation theorem of coplanar vector (theorem \ref{thm:Thm3.3.1}), it is deduced that:
	${{\gamma }_{1}}={{\mu }_{1}}$, ${{\gamma }_{2}}={{\mu }_{2}}$, ${{\gamma }_{3}}={{\mu }_{3}}$.
	
	Let the components (coefficients) of the corresponding vectors $\overrightarrow{OA}$, $\overrightarrow{OB}$ and $\overrightarrow{OC}$ in formulas (\ref{Eq6.1.10}) and (\ref{Eq6.1.12}) be equal, that is, let ${{\gamma }_{1}}={{\mu }_{1}}$, ${{\gamma }_{2}}={{\mu }_{2}}$, ${{\gamma }_{3}}={{\mu }_{3}}$, the following results can be obtained:	
	\[\frac{{{\lambda }_{BM}}}{\left( 1+\lambda _{AC}^{P} \right)\left( 1+{{\lambda }_{BM}} \right)}=\frac{{{\lambda }_{CN}}}{\left( 1+\lambda _{AB}^{P} \right)\left( 1+{{\lambda }_{CN}} \right)},\]
	\[\frac{1+\lambda _{AC}^{P}}{\left( 1+\lambda _{AC}^{P} \right)\left( 1+{{\lambda }_{BM}} \right)}=\frac{\lambda _{AB}^{P}{{\lambda }_{CN}}}{\left( 1+\lambda _{AB}^{P} \right)\left( 1+{{\lambda }_{CN}} \right)},\]
	\[\frac{\lambda _{AC}^{P}{{\lambda }_{BM}}}{\left( 1+\lambda _{AC}^{P} \right)\left( 1+{{\lambda }_{BM}} \right)}=\frac{1+\lambda _{AB}^{P}}{\left( 1+\lambda _{AB}^{P} \right)\left( 1+{{\lambda }_{CN}} \right)}.\]
	
	From the above three formulas, it is obtained that:	
	\begin{equation}\label{Eq6.1.13}
		\frac{{{\lambda }_{BM}}}{\left( 1+\lambda _{AC}^{P} \right)\left( 1+{{\lambda }_{BM}} \right)}=\frac{{{\lambda }_{CN}}}{\left( 1+\lambda _{AB}^{P} \right)\left( 1+{{\lambda }_{CN}} \right)},	
	\end{equation}
	\begin{equation}\label{Eq6.1.14}
		\frac{1}{1+{{\lambda }_{BM}}}=\frac{\lambda _{AB}^{P}{{\lambda }_{CN}}}{\left( 1+\lambda _{AB}^{P} \right)\left( 1+{{\lambda }_{CN}} \right)},	
	\end{equation}
	\begin{equation}\label{Eq6.1.15}
		\frac{\lambda _{AC}^{P}{{\lambda }_{BM}}}{\left( 1+\lambda _{AC}^{P} \right)\left( 1+{{\lambda }_{BM}} \right)}=\frac{1}{1+{{\lambda }_{CN}}}.
	\end{equation}
	
	Let's use $\lambda _{AB}^{P}$ and $\lambda _{AC}^{P}$ to express ${{\lambda }_{BM}}$ and ${{\lambda }_{CN}}$. According to the formulas (\ref{Eq6.1.13}) and (\ref{Eq6.1.14}), the following result can be obtained:	
	\[\frac{1}{1+{{\lambda }_{BM}}}=\frac{\lambda _{AB}^{P}{{\lambda }_{BM}}}{\left( 1+\lambda _{AC}^{P} \right)\left( 1+{{\lambda }_{BM}} \right)},\]
	i.e.
	\[\frac{\lambda _{AB}^{P}{{\lambda }_{BM}}}{1+\lambda _{AC}^{P}}=1,\]
	therefore,	
	\begin{equation}\label{Eq6.1.16}
		{{\lambda }_{BM}}=\frac{1+\lambda _{AC}^{P}}{\lambda _{AB}^{P}}.
	\end{equation}

	According to the formula (\ref{Eq6.1.13}) and (\ref{Eq6.1.15}), the following result can be obtained:	
	\[\frac{\lambda _{AC}^{P}{{\lambda }_{CN}}}{\left( 1+\lambda _{AB}^{P} \right)\left( 1+{{\lambda }_{CN}} \right)}=\frac{1}{1+{{\lambda }_{CN}}},\]
	therefore,	
	\begin{equation}\label{Eq6.1.17}
		{{\lambda }_{CN}}=\frac{1+\lambda _{AB}^{P}}{\lambda _{AC}^{P}}.
	\end{equation}
	
	By substituting formula (\ref{Eq6.1.16}) into formula (\ref{Eq6.1.9}), the following holds:
	\begin{equation}\label{Eq6.1.18}
		\begin{aligned}
			\overrightarrow{OP}& =\frac{\frac{1+\lambda _{AC}^{P}}{\lambda _{AB}^{P}}\overrightarrow{OA}+\left( 1+\lambda _{AC}^{P} \right)\overrightarrow{OB}+\lambda _{AC}^{P}\frac{1+\lambda _{AC}^{P}}{\lambda _{AB}^{P}}\overrightarrow{OC}}{\left( 1+\lambda _{AC}^{P} \right)\left( 1+\frac{1+\lambda _{AC}^{P}}{\lambda _{AB}^{P}} \right)} \\ 
			& =\frac{\overrightarrow{OA}+\lambda _{AB}^{P}\overrightarrow{OB}+\lambda _{AC}^{P}\overrightarrow{OC}}{1+\lambda _{AB}^{P}+\lambda _{AC}^{P}} \\ 
			& =\alpha _{A}^{P}\overrightarrow{OA}+\alpha _{B}^{P}\overrightarrow{OB}+\alpha _{C}^{P}\overrightarrow{OC},  
		\end{aligned}
	\end{equation}
	where
	\[\alpha _{A}^{P}=\frac{1}{1+\lambda _{AB}^{P}+\lambda _{AC}^{P}},\]	
	\[\alpha _{B}^{P}=\frac{\lambda _{AB}^{P}}{1+\lambda _{AB}^{P}+\lambda _{AC}^{P}},\]	
	\[\alpha _{C}^{P}=\frac{\lambda _{AC}^{P}}{1+\lambda _{AB}^{P}+\lambda _{AC}^{P}}.\]	
	
	We can also substitute formula (\ref{Eq6.1.17}) into formula (\ref{Eq6.1.11}) to get
	\[\begin{aligned}
		\overrightarrow{OP}& =\frac{\frac{1+\lambda _{AB}^{P}}{\lambda _{AC}^{P}}\overrightarrow{OA}+{{\lambda }_{AB}}\frac{1+\lambda _{AB}^{P}}{\lambda _{AC}^{P}}\overrightarrow{OB}+\left( 1+\lambda _{AB}^{P} \right)\overrightarrow{OC}}{\left( 1+\lambda _{AB}^{P} \right)\left( 1+\frac{1+\lambda _{AB}^{P}}{\lambda _{AC}^{P}} \right)} \\ 
		& =\frac{\overrightarrow{OA}+\lambda _{AB}^{P}\overrightarrow{OB}+\lambda _{AC}^{P}\overrightarrow{OC}}{1+\lambda _{AB}^{P}+\lambda _{AC}^{P}}.  
	\end{aligned}\]	
	
	The above formula is consistent with the result of formula (\ref{Eq6.1.18}).
	
	Similarly, the following holds:
	\begin{equation}\label{Eq6.1.19}
		\begin{aligned}
			\overrightarrow{OP}&=\frac{\overrightarrow{OB}+\lambda _{BC}^{P}\overrightarrow{OC}+\lambda _{BA}^{P}\overrightarrow{OA}}{1+\lambda _{BC}^{P}+\lambda _{BA}^{P}} \\ 
			& =\alpha _{A}^{P}\overrightarrow{OA}+\alpha _{B}^{P}\overrightarrow{OB}+\alpha _{C}^{P}\overrightarrow{OC},  
		\end{aligned}
	\end{equation}
	where
	\[\alpha _{B}^{P}=\frac{1}{1+\lambda _{BC}^{P}+\lambda _{BA}^{P}},\]	
	\[\alpha _{C}^{P}=\frac{\lambda _{BC}^{P}}{1+\lambda _{BC}^{P}+\lambda _{BA}^{P}},\]	
	\[\alpha _{A}^{P}=\frac{\lambda _{BA}^{P}}{1+\lambda _{BC}^{P}+\lambda _{BA}^{P}}.\]	
	
	Or
	\begin{equation}\label{Eq6.1.20}
		\begin{aligned}
			\overrightarrow{OP}&=\frac{\overrightarrow{OC}+\lambda _{CA}^{P}\overrightarrow{OA}+\lambda _{CB}^{P}\overrightarrow{OB}}{1+\lambda _{CA}^{P}+\lambda _{CB}^{P}} \\ 
			& =\alpha _{A}^{P}\overrightarrow{OA}+\alpha _{B}^{P}\overrightarrow{OB}+\alpha _{C}^{P}\overrightarrow{OC},  
		\end{aligned}
	\end{equation}
	where
	\[\alpha _{C}^{P}=\frac{1}{1+\lambda _{CA}^{P}+\lambda _{CB}^{P}},\]	
	\[\alpha _{A}^{P}=\frac{\lambda _{CA}^{P}}{1+\lambda _{CA}^{P}+\lambda _{CB}^{P}},\]	
	\[\alpha _{B}^{P}=\frac{\lambda _{CB}^{P}}{1+\lambda _{CA}^{P}+\lambda _{CB}^{P}}.\]	
	
	In order to get the symmetric form, according to the basic properties of proportion, we use the formulas (\ref{Eq6.1.18}), (\ref{Eq6.1.19}) and (\ref{Eq6.1.20}) to get the symmetric form in the following:
	\begin{align*}
		\overrightarrow{OP}&=\frac{\overrightarrow{OA}+\lambda _{AB}^{P}\overrightarrow{OB}+\lambda _{AC}^{P}\overrightarrow{OC}}{1+\lambda _{AB}^{P}+\lambda _{AC}^{P}} \\		&=\frac{\overrightarrow{OB}+\lambda _{BC}^{P}\overrightarrow{OC}+\lambda _{BA}^{P}\overrightarrow{OA}}{1+\lambda _{BC}^{P}+\lambda _{BA}^{P}} \\ &=\frac{\overrightarrow{OC}+\lambda _{CA}^{P}\overrightarrow{OA}+\lambda _{CB}^{P}\overrightarrow{OB}}{1+\lambda _{CA}^{P}+\lambda _{CB}^{P}} \\ 
		& =\frac{\left( \overrightarrow{OA}+\lambda _{AB}^{P}\overrightarrow{OB}+\lambda _{AC}^{P}\overrightarrow{OC} \right)+\left( \overrightarrow{OB}+\lambda _{BC}^{P}\overrightarrow{OC}+\lambda _{BA}^{P}\overrightarrow{OA} \right)+\left( \overrightarrow{OC}+\lambda _{CA}^{P}\overrightarrow{OA}+\lambda _{CB}^{P}\overrightarrow{OB} \right)}{\left( 1+\lambda _{AB}^{P}+\lambda _{AC}^{P} \right)+\left( 1+\lambda _{BC}^{P}+\lambda _{BA}^{P} \right)+\left( 1+\lambda _{CA}^{P}+\lambda _{CB}^{P} \right)} \\ 
		& =\frac{\left( 1+\lambda _{BA}^{P}+\lambda _{CA}^{P} \right)\overrightarrow{OA}+\left( 1+\lambda _{CB}^{P}+\lambda _{AB}^{P} \right)\overrightarrow{OB}+\left( 1+\lambda _{AC}^{P}+\lambda _{BC}^{P} \right)\overrightarrow{OC}}{\left( 1+\lambda _{BA}^{P}+\lambda _{CA}^{P} \right)+\left( 1+\lambda _{CB}^{P}+\lambda _{AB}^{P} \right)+\left( 1+\lambda _{AC}^{P}+\lambda _{BC}^{P} \right)}.
	\end{align*}
	
	The result is as follows:
	\[\overrightarrow{OP}=\frac{\left( 1+\lambda _{BA}^{P}+\lambda _{CA}^{P} \right)\overrightarrow{OA}+\left( 1+\lambda _{CB}^{P}+\lambda _{AB}^{P} \right)\overrightarrow{OB}+\left( 1+\lambda _{AC}^{P}+\lambda _{BC}^{P} \right)\overrightarrow{OC}}{\left( 1+\lambda _{BA}^{P}+\lambda _{CA}^{P} \right)+\left( 1+\lambda _{CB}^{P}+\lambda _{AB}^{P} \right)+\left( 1+\lambda _{AC}^{P}+\lambda _{BC}^{P} \right)}.\]		
	
	The above formula combines the conclusions (\ref{Eq6.1.1}) to (\ref{Eq6.1.4}) of the theorem in one expression.
	
	Obviously, according to (\ref{Eq6.1.2}) to (\ref{Eq6.1.4}), the following result can be obtained
	\[\alpha _{A}^{P}+\alpha _{B}^{P}+\alpha _{C}^{P}=1.\]	
	
	The uniqueness is proved below. According to the unique representation theorem of coplanar vector (theorem \ref{thm:Thm3.3.1}), we immediately deduce that each of $\alpha _{A}^{P}$, $\alpha _{B}^{P}$ and $\alpha _{C}^{P}$ in formulas (\ref{Eq6.1.2}), (\ref{Eq6.1.3}) and (\ref{Eq6.1.4}) is unique, that is to say, the vector $\overrightarrow{OP}$ can be uniquely expressed by the triangular frame $\left( O;A,B,C \right)$.
\end{proof}
\hfill $\square$\par




The $\alpha _{A}^{P}$, $\alpha _{B}^{P}$ and $\alpha _{C}^{P}$ are the frame components (coefficients) of the frame $\overrightarrow{OA}$, $\overrightarrow{OB}$ and $\overrightarrow{OC}$ respectively.

The frame components (coefficients) $\alpha _{A}^{P}$, $\alpha _{B}^{P}$ and $\alpha _{C}^{P}$ of $\overrightarrow{OP}$ contain only the IR of one point $P$. Such a vector $\overrightarrow{OP}$ of the intersecting center $P$ is called the vector from origin to intersecting center (abbreviated as VOIC).

The theorem of unique representation of coplanar vectors (theorem \ref{thm:Thm3.3.1}) points out that there exists a set of unique real numbers $\alpha _{A}^{P}$, $\alpha _{B}^{P}$ and $\alpha _{C}^{P}$, but it does not give the specific expression of this set of real numbers, while the theorem of vector from origin to intersecting centers (theorem \ref{thm:Thm6.1.1}) gives the specific expression of this set of real numbers.

Not only that, theorem \ref{thm:Thm6.1.1} has a deeper meaning. The importance of theorem \ref{thm:Thm6.1.1} lies in the fact that the vector $\overrightarrow{OP}$ can be expressed by the uniquely linear combination of the frame vectors $\overrightarrow{OA}$, $\overrightarrow{OB}$ and $\overrightarrow{OC}$ of the triangular frame $\left( O;A,B,C \right)$. Its expression coefficients $\alpha _{A}^{P}$, $\alpha _{B}^{P}$ and $\alpha _{C}^{P}$ are only the frame components of point $P$. As long as $\triangle ABC$ and point $P$ are given, the IR is a constant. According to the formulas (\ref{Eq6.1.2}), (\ref{Eq6.1.3}) and (\ref{Eq6.1.4}), the frame components (coefficients) $\alpha _{A}^{P}$, $\alpha _{B}^{P}$ and $\alpha _{C}^{P}$ of point $P$ are only determined by the IR, so they are also constant, namely, the set of components including $\alpha _{A}^{P}$, $\alpha _{B}^{P}$ and $\alpha _{C}^{P}$ is only related to the IR of the point $P$, but has nothing to do with the position of the origin $O$ of the frame, and has nothing to do with the position of $\triangle ABC$ in space.

It needs to be explained here that the meaning of given $\triangle ABC$ and point $P$ is that the positional relationship between $\triangle ABC$ and point $P$ is fixed, and you can imagine that point $P$ is pasted together with $\triangle ABC$ such that point $P$ is called the rigid point of $\triangle ABC$.

In fact, in the $\triangle ABC$ plane, if the starting IC is $O$ and the ending IC is $P$, then theorem \ref{thm:Thm6.1.1} calculates the vector $\overrightarrow{OP}$ between the two ICs: $O$ and $P$.

The three vertices and three sides of a triangle have equal status. Therefore, the formula reflecting the triangle should have some symmetry. theorem \ref{thm:Thm6.1.1} has the beauty of mathematical symmetry. It seems complex, but it is easy to remember.

If starting point $O$ and ending point $P$ are two ICs on the plane of $\triangle ABC$, then theorem \ref{thm:Thm6.1.1} states: the vector $\overrightarrow{OP}$ on the $\triangle ABC$ plane can be expressed by the uniquely linear combination of three vectors of $\overrightarrow{OA}$, $\overrightarrow{OB}$ and $\overrightarrow{OC}$, which are the vectors from the starting intersecting center $O$ to the three vertexs, and its representation coefficients are the frame components of the ending IC $P$, but the frame components $\alpha _{A}^{P}$, $\alpha _{B}^{P}$ and $\alpha _{C}^{P}$ of $\overrightarrow{OA}$, $\overrightarrow{OB}$, and $\overrightarrow{OC}$ are not independent of each other, they have a constraint relationship $\alpha _{A}^{P}+\alpha _{B}^{P}+\alpha _{C}^{P}=1$.

%
%
%

In Euclidean geometry, there is no corresponding relationship between a point and a real number, so there is usually no unified method to deal with the vector problem between two points, which can be said to be the “soft rib” of Euclidean geometry.

In analytic geometry, the corresponding relationship between a point and a real number is established, so algebraic method can be used to deal with geometric problems, which overcomes the weakness of Euclidean geometry to a certain extent, and is the advantage of analytic geometry. However, analytic geometry also has some disadvantages.   If the point is $P\left( x,y \right)$, then the frame components (coordinates) $x$, $y$ of the point are usually related to the selection of coordinate system, and the values of $x$, $y$ will be changed as the change of the origin of the coordinate system, which will bring inconvenience to the research and calculation of some problems. Not only that, the expression with point as parameter is often not “natural”. For example, the most “natural” parameter in a triangle is the lengths of three sides. When we want to get the result expressed by the lengths of three sides of a triangle, it is difficult to deal with it by using analytic geometry. One of the most typical examples is the area of a triangle. As we all know, given the lengths of three sides of a triangle, the area can be expressed by Heron's formula, which is a beautiful formula. However, given the coordinates of the three vertices of a triangle, it is difficult to obtain Helen's formula by analytic geometry.

In Intercenter Geometry, the corresponding relationship between point and IR (real number) is established. Therefore, Intercenter Geometry overcomes the weakness of Euclidean geometry to a certain extent and has the advantages of analytic geometry. Moreover, Intercenter Geometry overcomes the shortcomings of analytic geometry to a certain extent. Its frame components are independent of the frame (coordinate system), and the results obtained by Intercenter Geometry can be expressed by “natural” parameters (such as the length of each side of a triangle, etc.). In fact, the vector method is used in the research of Intercenter Geometry. From theorem \ref{thm:Thm6.1.1}, it can be seen that the frame component of vector $\overrightarrow{OP}$ is only determined by the IR, and the distance between two points is obtained after the vector is modeled, which means that the distance between two points is only related to the IR and some “natural” parameters. The keen readers should have realized that many basic geometric quantities can be expressed in terms of IR and “natural” parameters. At this point, readers should have a preliminary understanding of Intercenter Geometry. I believe you will have a deeper understanding after reading this book.

To show the charm of theorem \ref{thm:Thm6.1.1}, let's take a look at an example (see Chapter 9 for more applications).

\begin{example}{(Calculate the frame vector according to the IR)}\label{Exm6.1.1} 
	Find the frame vector of the centroid.	
\end{example}
\begin{solution}
If the IC $P$ is placed at the centroid $G$ of the triangle, then all the IRs are 1 (for example, $\lambda_{BA}^{P}=1$, etc.), then according to the formula (\ref{Eq6.1.2}), (\ref{Eq6.1.3}) and (\ref{Eq6.1.4}), the frame components are obtained as follows:	
	\[\begin{aligned}
		\alpha _{A}^{G}&=\frac{1+\lambda _{BA}^{G}+\lambda _{CA}^{G}}{\left( 1+\lambda _{BA}^{G}+\lambda _{CA}^{G} \right)+\left( 1+\lambda _{CB}^{G}+\lambda _{AB}^{G} \right)+\left( 1+\lambda _{AC}^{G}+\lambda _{BC}^{G} \right)} \\ 
		& =\frac{1+1+1}{\left( 1+1+1 \right)+\left( 1+1+1 \right)+\left( 1+1+1 \right)}=\frac{1}{3},  
	\end{aligned}\]
	\[\begin{aligned}
		\alpha _{B}^{G}& =\frac{1+\lambda _{CB}^{G}+\lambda _{AB}^{G}}{\left( 1+\lambda _{BA}^{G}+\lambda _{CA}^{G} \right)+\left( 1+\lambda _{CB}^{G}+\lambda _{AB}^{G} \right)+\left( 1+\lambda _{AC}^{G}+\lambda _{BC}^{G} \right)} \\ 
		& =\frac{1+1+1}{\left( 1+1+1 \right)+\left( 1+1+1 \right)+\left( 1+1+1 \right)}=\frac{1}{3},  
	\end{aligned}\]
	\[\begin{aligned}
		\alpha _{C}^{G}& =\frac{1+\lambda _{AC}^{G}+\lambda _{BC}^{G}}{\left( 1+\lambda _{BA}^{G}+\lambda _{CA}^{G} \right)+\left( 1+\lambda _{CB}^{G}+\lambda _{AB}^{G} \right)+\left( 1+\lambda _{AC}^{G}+\lambda _{BC}^{G} \right)} \\ 
		& =\frac{1+1+1}{\left( 1+1+1 \right)+\left( 1+1+1 \right)+\left( 1+1+1 \right)}=\frac{1}{3}.  
	\end{aligned}\]
	
	That is, the frame components are all ${1}/{3}\;$, Then according to theorem \ref{thm:Thm6.1.1}, the following holds:
	\[\begin{aligned}
		\overrightarrow{OG}&=\frac{\left( 1+\lambda _{BA}^{G}+\lambda _{CA}^{G} \right)\overrightarrow{OA}+\left( 1+\lambda _{CB}^{G}+\lambda _{AB}^{G} \right)\overrightarrow{OB}+\left( 1+\lambda _{AC}^{G}+\lambda _{BC}^{G} \right)\overrightarrow{OC}}{\left( 1+\lambda _{BA}^{G}+\lambda _{CA}^{G} \right)+\left( 1+\lambda _{CB}^{G}+\lambda _{AB}^{G} \right)+\left( 1+\lambda _{AC}^{G}+\lambda _{BC}^{G} \right)} \\ 
		& =\frac{\overrightarrow{OA}+\overrightarrow{OB}+\overrightarrow{OC}}{3}.  
	\end{aligned}\]		
	
	This formula indicates that the vector $\overrightarrow{OG}$ from any point $O$ to the triangle's centroid $G$ and the sum vector of three frame vectors $\overrightarrow{OA}$, $\overrightarrow{OB}$, $\overrightarrow{OC}$ of $\triangle ABC$ are collinear, and the size is one third of the sum vector, and there is no argument as to where $O$ is located, and above formula always holds regardless of where $\triangle ABC$ is located.
\end{solution}
\hfill $\diamond$\par

For any $\triangle ABC$, three vertices $A$, $B$ and $C$ can't be collinear. Therefore, according to the theorem in this section and the unique representation theorem of coplanar vectors (theorem \ref{thm:Thm3.3.1}), the following very practical theorems can be obtained.

\begin{theorem}{Unique representation theorem of triangular frame, Daiyuan Zhang}{Thm6.1.2}\label{Thm6.1.2} 
	The vector of intersecting center of any triangle can be expressed by the uniquely linear combination of its triangular frame $\left( O;A,B,C \right)$.	
\end{theorem}

\subsection{Representation of vector from origin to an intersecting center by edge frames}\label{Subsec6.1.1}

Theorem \ref{thm:Thm6.1.1} points out that  the origin $O$ of frame $\left( O;A,B,C \right)$ can be anywhere. Sometimes, for the convenience of analysis and calculation, the origin $O$ is selected in some special points. In this section, the origin $O$ is placed at the vertex of $\triangle ABC$ so that the frame vector coincides with the edge of the triangle. This frame is called the edge frame of the triangle, and also called the edge frame for short. For example, when the origin $O$ is placed at the vertex $A$ of $\triangle ABC$, the two edge frames are $\overrightarrow{AB}$ and $\overrightarrow{AC}$.

\begin{theorem}{Representation by edge frames for vector of IC, Daiyuan Zhang}{Thm6.1.3}\label{Thm6.1.3} 
	Given a $\triangle ABC$ and an intersecting center $P\in {{\pi }_{ABC}}$, then the vector between the vertex and the intersecting center of $\triangle ABC$ can be expressed by the uniquely linear combination of the triangle's edge frame, that is
	\[\left\{ \begin{aligned}
		& \overrightarrow{AP}=\alpha _{B}^{P}\overrightarrow{AB}+\alpha _{C}^{P}\overrightarrow{AC} \\ 
		& \alpha _{A}^{P}+\alpha _{B}^{P}+\alpha _{C}^{P}=1, \\ 
	\end{aligned} \right.\]	
	\[\left\{ \begin{aligned}
		& \overrightarrow{BP}=\alpha _{A}^{P}\overrightarrow{BA}+\alpha _{C}^{P}\overrightarrow{BC} \\ 
		& \alpha _{A}^{P}+\alpha _{B}^{P}+\alpha _{C}^{P}=1, \\ 
	\end{aligned} \right.\]	
	\[\left\{ \begin{aligned}
		& \overrightarrow{CP}=\alpha _{A}^{P}\overrightarrow{CA}+\alpha _{B}^{P}\overrightarrow{CB} \\ 
		& \alpha _{A}^{P}+\alpha _{B}^{P}+\alpha _{C}^{P}=1. \\ 
	\end{aligned} \right.\]	
	
	Where $\alpha _{A}^{P}$, $\alpha _{B}^{P}$, $\alpha _{C}^{P}$ are the frame components of $\overrightarrow{OA}$, $\overrightarrow{OB}$, $\overrightarrow{OC}$ corresponding to $P$ in the frame $\left( O;A,B,C \right)$ respectively.
\end{theorem}

\begin{proof}
	This theorem is obtained by replacing the point $O$ in theorem \ref{thm:Thm6.1.1} with point $A$, $B$, $C$, respectively.
\end{proof}
\hfill $\square$\par

The edge frames can reduce the number of vectors, thus reducing the amount of computation, so this is an important special case.

\subsection{Relationship between vector of vertex to intersecting center and vector of vertex to intersecting foot}\label{Subsec6.1.2}

The vector from a vertex to the IC of a triangle is called the vector of vertex to intersecting center (abbreviated as VVIC). In Figure \ref{fig:tu6.1.1}, each of $\overrightarrow{AP}$, $\overrightarrow{BP}$ and $\overrightarrow{CP}$ is the VVIC.

The vector from the vertex of a triangle to its corresponding intersecting foot is called the vector from vertex to intersecting foot (abbreviated as VVIF). In Figure \ref{fig:tu6.1.1}, each of $\overrightarrow{AL}$, $\overrightarrow{BM}$ and $\overrightarrow{CN}$ is the VVIF.

There is a certain relationship between the VVIC and the VVIF, and they can be expressed by the edge frame vectors, and their representation components are functions of the frame components.

\begin{theorem}{Relationship between VVIC and VVIF, Daiyuan Zhang}{Thm6.1.4}\label{Thm6.1.4} 
	For a given $\triangle ABC$, the vector of vertex to intersecting foot (VVIF) is
	\begin{flalign*}
		\overrightarrow{AL}=\frac{\overrightarrow{AP}}{1-\alpha _{A}^{P}},	
		\overrightarrow{BM}=\frac{\overrightarrow{BP}}{1-\alpha _{B}^{P}},	
		\overrightarrow{CN}=\frac{\overrightarrow{CP}}{1-\alpha _{C}^{P}}.	
	\end{flalign*}
	Where $\alpha _{A}^{P}$, $\alpha _{B}^{P}$, $\alpha _{C}^{P}$ are the frame components of $\overrightarrow{OA}$, $\overrightarrow{OB}$, $\overrightarrow{OC}$ corresponding to $P$ in the frame $\left( O;A,B,C \right)$ respectively.
\end{theorem}

\begin{proof}
	According to theorem \ref{thm:Thm2.4.1}, theorem \ref{thm:Thm6.1.1}, and theorem \ref{thm:Thm7.1.1}, the following hold (see Fig.\ref{fig:tu6.1.1}):
	\[\begin{aligned}
		\overrightarrow{AL}&=\frac{\overrightarrow{AP}}{{{\kappa }_{AL}}}=\frac{\left( 1+{{\lambda }_{AL}} \right)\overrightarrow{AP}}{{{\lambda }_{AL}}}=\frac{\left( 1+\lambda _{AB}^{P}+\lambda _{AC}^{P} \right)\overrightarrow{AP}}{\lambda _{AB}^{P}+\lambda _{AC}^{P}} \\ 
		& =\frac{\overrightarrow{AP}}{\frac{\lambda _{AB}^{P}+\lambda _{AC}^{P}}{1+\lambda _{AB}^{P}+\lambda _{AC}^{P}}}=\frac{\overrightarrow{AP}}{\alpha _{B}^{P}+\alpha _{C}^{P}}=\frac{\overrightarrow{AP}}{1-\alpha _{A}^{P}},  
	\end{aligned}\]
	\[\begin{aligned}
		\overrightarrow{BM}&=\frac{\overrightarrow{BP}}{{{\kappa }_{BM}}}=\frac{\left( 1+{{\lambda }_{BM}} \right)\overrightarrow{BP}}{{{\lambda }_{BM}}}=\frac{\left( 1+\lambda _{BC}^{P}+\lambda _{BA}^{P} \right)\overrightarrow{BP}}{\lambda _{BC}^{P}+\lambda _{BA}^{P}} \\ 
		& =\frac{\overrightarrow{BP}}{\frac{\lambda _{BC}^{P}+\lambda _{BA}^{P}}{1+\lambda _{BC}^{P}+\lambda _{BA}^{P}}}=\frac{\overrightarrow{BP}}{\alpha _{C}^{P}+\alpha _{A}^{P}}=\frac{\overrightarrow{BP}}{1-\alpha _{B}^{P}},  
	\end{aligned}\]
	\begin{align*}
		\overrightarrow{CN}&=\frac{\overrightarrow{CP}}{{{\kappa }_{CN}}}=\frac{\left( 1+{{\lambda }_{CN}} \right)\overrightarrow{CP}}{{{\lambda }_{CN}}}=\frac{\left( 1+\lambda _{CA}^{P}+\lambda _{CB}^{P} \right)\overrightarrow{CP}}{\lambda _{CA}^{P}+\lambda _{CB}^{P}} \\ 
		& =\frac{\overrightarrow{CP}}{\frac{\lambda _{CA}^{P}+\lambda _{CB}^{P}}{1+\lambda _{CA}^{P}+\lambda _{CB}^{P}}}=\frac{\overrightarrow{CP}}{\alpha _{A}^{P}+\alpha _{B}^{P}}=\frac{\overrightarrow{CP}}{1-\alpha _{C}^{P}}.  
	\end{align*}
\end{proof}
\hfill $\square$\par

\begin{theorem}{Theorem of VVIF, Daiyuan Zhang}{Thm6.1.5}\label{Thm6.1.5} 
	For a given $\triangle ABC$, its vector of vertex to intersecting foot (VVIF) can be expressed by the edge frame of a triangle.
	\[\overrightarrow{AL}=\frac{\alpha _{B}^{P}}{\alpha _{B}^{P}+\alpha _{C}^{P}}\overrightarrow{AB}+\frac{\alpha _{C}^{P}}{\alpha _{B}^{P}+\alpha _{C}^{P}}\overrightarrow{AC},\]	
	\[\overrightarrow{BM}=\frac{\alpha _{C}^{P}}{\alpha _{C}^{P}+\alpha _{A}^{P}}\overrightarrow{BC}+\frac{\alpha _{A}^{P}}{\alpha _{C}^{P}+\alpha _{A}^{P}}\overrightarrow{BA},\]	
	\[\overrightarrow{CN}=\frac{\alpha _{A}^{P}}{\alpha _{A}^{P}+\alpha _{B}^{P}}\overrightarrow{CA}+\frac{\alpha _{B}^{P}}{\alpha _{A}^{P}+\alpha _{B}^{P}}\overrightarrow{CB}.\]	
	Where $\alpha _{A}^{P}$, $\alpha _{B}^{P}$, $\alpha _{C}^{P}$ are the frame components of $\overrightarrow{OA}$, $\overrightarrow{OB}$, $\overrightarrow{OC}$ corresponding to $P$ in the frame $\left( O;A,B,C \right)$ respectively.
\end{theorem}

\begin{proof}
	According to theorem \ref{thm:Thm6.1.4} and theorem \ref{thm:Thm6.1.3}, the following holds:
	\[\overrightarrow{AL}=\frac{\overrightarrow{AP}}{1-\alpha _{A}^{P}}=\frac{\alpha _{B}^{P}\overrightarrow{AB}+\alpha _{C}^{P}\overrightarrow{AC}}{\alpha _{B}^{P}+\alpha _{C}^{P}}.\]	
	
	Similarly, the other two formulas can be proved.
\end{proof}
\hfill $\square$\par

The above theorem shows that the components of the two edge frames are determined by the components of the triangular frame, and the sum of the components of the two edge frames is 1.

\subsection{Representation of vector from origin to an intersecting center on circumcenter frame}\label{Subsec6.1.3}
In this section, the origin $O$ is placed at the circumcenter $Q$ of $\triangle ABC$ so that the lengths of the three frame vectors are equal, that is, the circumscribed circle radius of $\triangle ABC$, which will bring convenience to the analysis of some problems. Such a frame is called the circumcenter frame of a triangle, abbreviated as the circumcenter frame, and is denoted as $\left( Q;A,B,C \right)$.

\begin{theorem}{Representation by circumcenter frame for vector of IC, Daiyuan Zhang}{Thm6.1.6}\label{Thm6.1.6}  
	Given a $\triangle ABC$, $P\in {{\pi }_{ABC}}$, assuming that the point $Q$ is the circumcenter of $\triangle ABC$, then the vector $\overrightarrow{QP}$ can be expressed by the uniquely linear combination of the circumcenter frame $\left( Q;A,B,C \right)$, i.e.
	\[\overrightarrow{QP}=\alpha _{A}^{P}\overrightarrow{QA}+\alpha _{B}^{P}\overrightarrow{QB}+\alpha _{C}^{P}\overrightarrow{QC},\]	
	\[\alpha _{A}^{P}+\alpha _{B}^{P}+\alpha _{C}^{P}=1.\]	
	Where $\alpha _{A}^{P}$, $\alpha _{B}^{P}$, $\alpha _{C}^{P}$ are the frame components of $\overrightarrow{OA}$, $\overrightarrow{OB}$, $\overrightarrow{OC}$ corresponding to $P$ in the frame $\left( O;A,B,C \right)$ respectively.	
\end{theorem}
\begin{proof}
	This theorem is proved by replacing the point $O$ with the point $Q$ in theorem \ref{thm:Thm6.1.1}.
\end{proof}
\hfill $\square$\par

\section{Theorem of vector of two intersecting centers}\label{Sec6.2}
As known from the preceding discussion, the starting point of the vector $\overrightarrow{OP}$ obtained in theorem \ref{thm:Thm6.1.1} is still the the origin $O$ of the triangular frame $\left( O;A,B,C \right)$. The question now is how to use the three frame vectors $\overrightarrow{OA}$, $\overrightarrow{OB}$ and $\overrightarrow{OC}$ in the triangular frame $\left( O;A,B,C \right)$ to represent the vector $\overrightarrow{{{P}_{1}}{{P}_{2}}}$ uniquely and linearly (where point ${{P}_{1}}$ and ${{P}_{2}}$ is not necessarily coincident with the origin $O$)? This is the theorem to be studied below. The theorem given below has many applications and is very important.

\begin{theorem}{Theorem of vector of two intersecting centers, Daiyuan Zhang}{Thm6.2.1}\label{Thm6.2.1} 
	Suppose that given a $\triangle ABC$, $O$ is an arbitrary point, ${{P}_{1}}\in {{\pi }_{ABC}}$, ${{P}_{2}}\in {{\pi }_{ABC}}$, then vector $\overrightarrow{{{P}_{1}}{{P}_{2}}}$ can be expressed by the uniquely linear combination of the frame $\left( O;A,B,C \right)$, i.e.	
	\[\overrightarrow{{{P}_{1}}{{P}_{2}}}=\alpha _{A}^{{{P}_{1}}{{P}_{2}}}\overrightarrow{OA}+\alpha _{B}^{{{P}_{1}}{{P}_{2}}}\overrightarrow{OB}+\alpha _{C}^{{{P}_{1}}{{P}_{2}}}\overrightarrow{OC},\]	
	\[\alpha _{A}^{{{P}_{1}}{{P}_{2}}}+\alpha _{B}^{{{P}_{1}}{{P}_{2}}}+\alpha _{C}^{{{P}_{1}}{{P}_{2}}}=0,\]  	
	where
	\[\alpha _{A}^{{{P}_{1}}{{P}_{2}}}=\alpha _{A}^{{{P}_{2}}}-\alpha _{A}^{{{P}_{1}}},\]	
	\[\alpha _{B}^{{{P}_{1}}{{P}_{2}}}=\alpha _{B}^{{{P}_{2}}}-\alpha _{B}^{{{P}_{1}}},\]	
	\[\alpha _{C}^{{{P}_{1}}{{P}_{2}}}=\alpha _{C}^{{{P}_{2}}}-\alpha _{C}^{{{P}_{1}}}.\]	
	And, $\alpha _{A}^{{{P}_{1}}}$, $\alpha _{B}^{{{P}_{1}}}$, $\alpha _{C}^{{{P}_{1}}}$ are the frame components of $\overrightarrow{OA}$, $\overrightarrow{OB}$, $\overrightarrow{OC}$ corresponding to ${{P}_{1}}$ respectively in the frame $\left( O;A,B,C \right)$; $\alpha _{A}^{{{P}_{2}}}$, $\alpha _{B}^{{{P}_{2}}}$, $\alpha _{C}^{{{P}_{2}}}$ are the frame components of $\overrightarrow{OA}$, $\overrightarrow{OB}$, $\overrightarrow{OC}$ corresponding to ${{P}_{2}}$  respectively in the frame $\left( O;A,B,C \right)$.	
\end{theorem}

\begin{proof}
	Obviously, for any point $O$, the following result is obtained according to the basic operation of vector.
	\[\overrightarrow{{{P}_{1}}{{P}_{2}}}=\overrightarrow{O{{P}_{2}}}-\overrightarrow{O{{P}_{1}}}.\]
	
	According to theorem \ref{thm:Thm6.1.1}, there exists unique set of real numbers $\alpha_{A}^{{{P}_{1}}}$, $\alpha _{B}^{{{P}_{1}}}$, $\alpha _{C}^{{{P}_{1}}}$ makes the following formula hold:
	\[\overrightarrow{O{{P}_{1}}}=\alpha _{A}^{{{P}_{1}}}\overrightarrow{OA}+\alpha _{B}^{{{P}_{1}}}\overrightarrow{OB}+\alpha _{C}^{{{P}_{1}}}\overrightarrow{OC},\]	
	\[\alpha _{A}^{{{P}_{1}}}+\alpha _{B}^{{{P}_{1}}}+\alpha _{C}^{{{P}_{1}}}=1.\]	
	
	There exists unique set of real numbers $\alpha _{A}^{{{P}_{2}}}$, $\alpha _{B}^{{{P}_{2}}}$, $\alpha _{C}^{{{P}_{2}}}$ makes the following formula hold:
	\[\overrightarrow{O{{P}_{2}}}=\alpha _{A}^{{{P}_{2}}}\overrightarrow{OA}+\alpha _{B}^{{{P}_{2}}}\overrightarrow{OB}+\alpha _{C}^{{{P}_{2}}}\overrightarrow{OC},\]	
	\[\alpha _{A}^{{{P}_{2}}}+\alpha _{B}^{{{P}_{2}}}+\alpha _{C}^{{{P}_{2}}}=1.\]	
	
	Therefore,
	\[\begin{aligned}
		\overrightarrow{{{P}_{1}}{{P}_{2}}}&=\overrightarrow{O{{P}_{2}}}-\overrightarrow{O{{P}_{1}}} \\ 
		& =\alpha _{A}^{{{P}_{2}}}\overrightarrow{OA}+\alpha _{B}^{{{P}_{2}}}\overrightarrow{OB}+\alpha _{C}^{{{P}_{2}}}\overrightarrow{OC}-\left( \alpha _{A}^{{{P}_{1}}}\overrightarrow{OA}+\alpha _{B}^{{{P}_{1}}}\overrightarrow{OB}+\alpha _{C}^{{{P}_{1}}}\overrightarrow{OC} \right) \\ 
		& =\left( \alpha _{A}^{{{P}_{2}}}-\alpha _{A}^{{{P}_{1}}} \right)\overrightarrow{OA}+\left( \alpha _{B}^{{{P}_{2}}}-\alpha _{B}^{{{P}_{1}}} \right)\overrightarrow{OB}+\left( \alpha _{C}^{{{P}_{2}}}-\alpha _{C}^{{{P}_{1}}} \right)\overrightarrow{OC} \\ 
		& =\alpha _{A}^{{{P}_{1}}{{P}_{2}}}\overrightarrow{OA}+\alpha _{B}^{{{P}_{1}}{{P}_{2}}}\overrightarrow{OB}+\alpha _{C}^{{{P}_{1}}{{P}_{2}}}\overrightarrow{OC}.  
	\end{aligned}\]	
	
	And,
	\[\begin{aligned}
		& \alpha _{A}^{{{P}_{1}}{{P}_{2}}}+\alpha _{B}^{{{P}_{1}}{{P}_{2}}}+\alpha _{C}^{{{P}_{1}}{{P}_{2}}}=\left( \alpha _{A}^{{{P}_{2}}}-\alpha _{A}^{{{P}_{1}}} \right)+\left( \alpha _{B}^{{{P}_{2}}}-\alpha _{B}^{{{P}_{1}}} \right)+\left( \alpha _{C}^{{{P}_{2}}}-\alpha _{C}^{{{P}_{1}}} \right) \\ 
		& =\left( \alpha _{A}^{{{P}_{2}}}+\alpha _{B}^{{{P}_{2}}}+\alpha _{C}^{{{P}_{2}}} \right)-\left( \alpha _{A}^{{{P}_{1}}}+\alpha _{B}^{{{P}_{1}}}+\alpha _{C}^{{{P}_{1}}} \right)=1-1=0.  
	\end{aligned}\]	
	
	Since each of $\alpha_{A}^{{{P}_{1}}}$, $\alpha _{B}^{{{P}_{1}}}$, $\alpha _{C}^{{{P}_{1} }}$ is unique, and each of $\alpha_{A}^{{{P}_{2}}}$, $\alpha _{B}^{{{P}_{2}}}$, $\alpha _{C}^{{{P}_{2}}}$ is also unique, thus each of $\alpha_{A}^{{{P}_{1}}{{P}_{2}}}$, $\alpha _{B}^{{{P}_{1}}{{P}_{2}}}$ and $\alpha_{C}^{{{P}_{1}}{{P}_{2}}}$ is unique.
\end{proof}
\hfill $\square$\par


The frame coefficients of $\alpha _{A}^{{{P}_{1}}{{P}_{2}}}$, $\alpha _{B}^{{{P}_{1}}{{P}_{2}}}$ and $\alpha _{C}^{{{P}_{1}}{{P}_{2}}}$ in vector of $\overrightarrow{{{P}_{1}}{{P}_{2}}}$ are related to the two intersecting centers ${{P}_{1}}$ and ${{P}_{2}}$, The IC vector $\overrightarrow{{{P}_{1}}{{P}_{2}}}$ is called the vector of two intersecting centers (abbreviated as VTICs).


The importance of the theorem of VTICs (theorem \ref{thm:Thm6.2.1}) is that the vector $\overrightarrow{{{P}_{1}}{{P}_{2}}}$ can be expressed by the uniquely linear combination of the vectors $\overrightarrow{OA}$, $\overrightarrow{OB}$ and $\overrightarrow{OC}$ in triangular frame $\left( O;A,B,C \right)$, in which their frame components (coefficients) $\alpha _{A}^{{{P}_{1}}{{P}_{2}}}$, $\alpha _{B}^{{{P}_{1}}{{P}_{2}}}$ and $\alpha _{C}^{{{P}_{1}}{{P}_{2}}}$ are only related to two ICs, ${{P}_{1}}$ and ${{P}_{2}}$. If the two ICs are rigid points on the $\triangle ABC $ plane, then each of the IRs of the two ICs is a constant, so their frame components are also constants, which is the same as that of the frame $\left( O;A,B,C \right)$ has nothing to do with the position of the origin $O$, and also has nothing to do with the position of $\triangle ABC $. This conclusion is very important and convenient in many applications. The position of the origin and $\triangle ABC$ can be selected as needed without changing the frame components. For a given triangle, as long as each of the IRs of the two points on the plane of a triangle is known, the vector connecting the two points can be obtained, and thus the distance between the two points can be obtained.

The theorem of VTICs (theorem \ref{thm:Thm6.2.1}) can easily find the vectors between the centers of a triangle (discussed later).

\subsection{Representation of vector of two intersecting centers by edge frame}\label{Subsec6.2.1}
\begin{theorem}{Representation of vector of two intersecting centers by edge frame, Daiyuan Zhang}{Thm6.2.2}\label{Thm6.2.2} 
	Given a $\triangle ABC$, ${{P}_{1}}\in {{\pi }_{ABC}}$, ${{P}_{2}}\in {{\pi }_{ABC}}$, then $\overrightarrow{{{P}_{1}}{{P}_{2}}}$ can be expressed by the uniquely linear combination of the frame of the triangle, i.e.	
	\begin{equation}\label{Eq6.2.1}
		\overrightarrow{{{P}_{1}}{{P}_{2}}}=\alpha _{B}^{{{P}_{1}}{{P}_{2}}}\overrightarrow{AB}+\alpha _{C}^{{{P}_{1}}{{P}_{2}}}\overrightarrow{AC},
	\end{equation}
	\begin{equation}\label{Eq6.2.2}
		\overrightarrow{{{P}_{1}}{{P}_{2}}}=\alpha _{A}^{{{P}_{1}}{{P}_{2}}}\overrightarrow{BA}+\alpha _{C}^{{{P}_{1}}{{P}_{2}}}\overrightarrow{BC},
	\end{equation}	
	\begin{equation}\label{Eq6.2.3}
		\overrightarrow{{{P}_{1}}{{P}_{2}}}=\alpha _{A}^{{{P}_{1}}{{P}_{2}}}\overrightarrow{CA}+\alpha _{B}^{{{P}_{1}}{{P}_{2}}}\overrightarrow{CB},
	\end{equation}	
	\begin{equation}\label{Eq6.2.4}
		\overrightarrow{{{P}_{1}}{{P}_{2}}}=\frac{1}{3}\left( \left( \alpha _{B}^{{{P}_{1}}{{P}_{2}}}-\alpha _{A}^{{{P}_{1}}{{P}_{2}}} \right)\overrightarrow{AB}+\left( \alpha _{C}^{{{P}_{1}}{{P}_{2}}}-\alpha _{B}^{{{P}_{1}}{{P}_{2}}} \right)\overrightarrow{BC}+\left( \alpha _{A}^{{{P}_{1}}{{P}_{2}}}-\alpha _{C}^{{{P}_{1}}{{P}_{2}}} \right)\overrightarrow{CA} \right),
	\end{equation}		
	where
	\begin{equation}\label{Eq6.2.5}
		\alpha _{A}^{{{P}_{1}}{{P}_{2}}}+\alpha _{B}^{{{P}_{1}}{{P}_{2}}}+\alpha _{C}^{{{P}_{1}}{{P}_{2}}}=0,
	\end{equation}	
	\[\alpha _{A}^{{{P}_{1}}{{P}_{2}}}=\alpha _{A}^{{{P}_{2}}}-\alpha _{A}^{{{P}_{1}}},\]	
	\[\alpha _{B}^{{{P}_{1}}{{P}_{2}}}=\alpha _{B}^{{{P}_{2}}}-\alpha _{B}^{{{P}_{1}}},\]	
	\[\alpha _{C}^{{{P}_{1}}{{P}_{2}}}=\alpha _{C}^{{{P}_{2}}}-\alpha _{C}^{{{P}_{1}}}.\]	
	And, $\alpha _{A}^{{{P}_{1}}}$, $\alpha _{B}^{{{P}_{1}}}$, $\alpha _{C}^{{{P}_{1}}}$ are the frame components of $\overrightarrow{OA}$, $\overrightarrow{OB}$, $\overrightarrow{OC}$ corresponding to ${{P}_{1}}$ respectively in the frame $\left( O;A,B,C \right)$; $\alpha _{A}^{{{P}_{2}}}$, $\alpha _{B}^{{{P}_{2}}}$, $\alpha _{C}^{{{P}_{2}}}$ are the frame components of $\overrightarrow{OA}$, $\overrightarrow{OB}$, $\overrightarrow{OC}$ corresponding to ${{P}_{2}}$  respectively in the frame $\left( O;A,B,C \right)$.		
\end{theorem}

\begin{proof}
	According to the theorem of VTICs (theorem \ref{thm:Thm6.2.1}), the formulas (\ref{Eq6.2.1}), (\ref{Eq6.2.2}), (\ref{Eq6.2.3}) can be obtained by coincidence of point $O$ with point $A$, $B$, $C$, respectively. By adding the three formulas together, the formulas (\ref{Eq6.2.4}) can be obtained.
\end{proof}
\hfill $\square$\par

The above theorem illustrates vector $\overrightarrow{{{P}_{1}}{{P}_{2}}}$ can be expressed by the uniquely linear combination of the edge frames whose frame components (coefficients) are constants related only to the IRs of two centers, and are independent of where the triangle is located.
		
The above theorem can be easily used to calculate the length of vector $\overrightarrow{{{P}_{1}}{{P}_{2}}}$. Formulas (\ref{Eq6.2.1}), (\ref{Eq6.2.2}), (\ref{Eq6.2.3}) are simple to calculate, and formula (\ref{Eq6.2.4}) has symmetrical aesthetic feeling and is easy to remember.

\subsection{Representation of vector of two intersecting centers by frames of circumcenter}\label{Subsec6.2.2}

\begin{theorem}{Representation of vector of two ICs by frame of circumcenter, Daiyuan Zhang}{Thm6.2.3}\label{Thm6.2.3} 
	Given a $\triangle ABC$, ${{P}_{1}}\in {{\pi }_{ABC}}$, ${{P}_{2}}\in {{\pi }_{ABC}}$ and point $Q$ is the circumcenter of the triangle, then vector $\overrightarrow{{{P}_{1}}{{P}_{2}}}$ can be expressed by $\left( Q;A,B,C \right)$ uniquely and linearly, i.e.	
	\[\overrightarrow{{{P}_{1}}{{P}_{2}}}=\alpha _{A}^{{{P}_{1}}{{P}_{2}}}\overrightarrow{QA}+\alpha _{B}^{{{P}_{1}}{{P}_{2}}}\overrightarrow{QB}+\alpha _{C}^{{{P}_{1}}{{P}_{2}}}\overrightarrow{QC},\]	
	\[\alpha _{A}^{{{P}_{1}}{{P}_{2}}}+\alpha _{B}^{{{P}_{1}}{{P}_{2}}}+\alpha _{C}^{{{P}_{1}}{{P}_{2}}}=0,\]  	
	where
	\[\alpha _{A}^{{{P}_{1}}{{P}_{2}}}=\alpha _{A}^{{{P}_{2}}}-\alpha _{A}^{{{P}_{1}}},\]	
	\[\alpha _{B}^{{{P}_{1}}{{P}_{2}}}=\alpha _{B}^{{{P}_{2}}}-\alpha _{B}^{{{P}_{1}}},\]	
	\[\alpha _{C}^{{{P}_{1}}{{P}_{2}}}=\alpha _{C}^{{{P}_{2}}}-\alpha _{C}^{{{P}_{1}}}.\]	
	And, $\alpha _{A}^{{{P}_{1}}}$, $\alpha _{B}^{{{P}_{1}}}$, $\alpha _{C}^{{{P}_{1}}}$ are the frame components of $\overrightarrow{OA}$, $\overrightarrow{OB}$, $\overrightarrow{OC}$ corresponding to ${{P}_{1}}$ respectively in the frame $\left( O;A,B,C \right)$; $\alpha _{A}^{{{P}_{2}}}$, $\alpha _{B}^{{{P}_{2}}}$, $\alpha _{C}^{{{P}_{2}}}$ are the frame components of $\overrightarrow{OA}$, $\overrightarrow{OB}$, $\overrightarrow{OC}$ corresponding to ${{P}_{2}}$  respectively in the frame $\left( O;A,B,C \right)$.		
\end{theorem}

\begin{proof}
	According to the theorem of VTICs (theorem \ref{thm:Thm6.2.1}), the result is obtained by overlapping the origin $O$ (of the frame $\left( O;A,B,C \right)$) and the circumcenter $Q$.
\end{proof}
\hfill $\square$\par

\section{Frame equation of intersecting center}\label{Sec6.3}
The frame equation of intersecting center satisfies a certain vector equation.

\begin{theorem}{Frame equation of single intersecting center, Daiyuan Zhang}{Thm6.3.1}\label{Thm6.3.1} 
	Assuming a given $\triangle ABC$, point $P$ is the IC, then
	\[\alpha _{A}^{P}\overrightarrow{PA}+\alpha _{B}^{P}\overrightarrow{PB}+\alpha _{C}^{P}\overrightarrow{PC}=\overrightarrow{0},\]
	\[\alpha _{A}^{P}+\alpha _{B}^{P}+\alpha _{C}^{P}=1.\]	
\end{theorem}

\begin{proof}
	For triangular frame $\left( O;A,B,C \right)$, according to theorem \ref{thm:Thm6.1.1}, let the point $O$ and the point $P$ coincide, and the result can be obtained directly.
\end{proof}
\hfill $\square$\par

According to the frame equation of single intersecting center (abbreviated as FESIC), the frame equation of each IC for a given triangle can be written directly.

The frame equation of two intersecting centers is given below.
\begin{theorem}{Frame equation of two intersecting centers, Daiyuan Zhang}{Thm6.3.2}\label{Thm6.3.2} 
	Assuming a given $\triangle ABC$, points ${{P}_{1}}$ and ${{P}_{2}}$ are two ICs, then
	\[\begin{aligned}
		& \alpha _{A}^{{{P}_{1}}}\overrightarrow{{{P}_{2}}A}+\alpha _{B}^{{{P}_{1}}}\overrightarrow{{{P}_{2}}B}+\alpha _{C}^{{{P}_{1}}}\overrightarrow{{{P}_{2}}C}+\alpha _{A}^{{{P}_{2}}}\overrightarrow{{{P}_{1}}A}+\alpha _{B}^{{{P}_{2}}}\overrightarrow{{{P}_{1}}B}+\alpha _{C}^{{{P}_{2}}}\overrightarrow{{{P}_{1}}C}=\overrightarrow{0}, \\ 
	\end{aligned}\]	
	\[\alpha _{A}^{{{P}_{1}}}+\alpha _{B}^{{{P}_{1}}}+\alpha _{C}^{{{P}_{1}}}=1,\]	
	\[\alpha _{A}^{{{P}_{2}}}+\alpha _{B}^{{{P}_{2}}}+\alpha _{C}^{{{P}_{2}}}=1.\]		
\end{theorem}

\begin{proof}
	Because the ICs have equal status, the following formula can be obtained according to theorem \ref{thm:Thm6.1.1}
	\[\overrightarrow{{{P}_{1}}{{P}_{2}}}=\alpha _{A}^{{{P}_{2}}}\overrightarrow{{{P}_{1}}A}+\alpha _{B}^{{{P}_{2}}}\overrightarrow{{{P}_{1}}B}+\alpha _{C}^{{{P}_{2}}}\overrightarrow{{{P}_{1}}C},\]	
	\[\alpha _{A}^{{{P}_{2}}}+\alpha _{B}^{{{P}_{2}}}+\alpha _{C}^{{{P}_{2}}}=1.\]	
	\[\overrightarrow{{{P}_{2}}{{P}_{1}}}=\alpha _{A}^{{{P}_{1}}}\overrightarrow{{{P}_{2}}A}+\alpha _{B}^{{{P}_{1}}}\overrightarrow{{{P}_{2}}B}+\alpha _{C}^{{{P}_{1}}}\overrightarrow{{{P}_{2}}C},\]	
	\[\alpha _{A}^{{{P}_{1}}}+\alpha _{B}^{{{P}_{1}}}+\alpha _{C}^{{{P}_{1}}}=1.\]	
	
	Since $\overrightarrow{{{P}_{1}}{{P}_{2}}}+\overrightarrow{{{P}_{2}}{{P}_{1}}}=\overrightarrow{0}$, the following result is obtained:
	\[\alpha _{A}^{{{P}_{1}}}\overrightarrow{{{P}_{2}}A}+\alpha _{B}^{{{P}_{1}}}\overrightarrow{{{P}_{2}}B}+\alpha _{C}^{{{P}_{1}}}\overrightarrow{{{P}_{2}}C}+\alpha _{A}^{{{P}_{2}}}\overrightarrow{{{P}_{1}}A}+\alpha _{B}^{{{P}_{2}}}\overrightarrow{{{P}_{1}}B}+\alpha _{C}^{{{P}_{2}}}\overrightarrow{{{P}_{1}}C}=\overrightarrow{0}.\]
\end{proof}
\hfill $\square$\par

The above equation gives the vector equation satisfied by the frame components of two points as well as the corresponding frame.

The frame equation of single intersecting center (abbreviated as FESIC) and the frame equation of two intersecting centers (abbreviated as FETICs) are also called frame equation of intersecting center (abbreviated as FEIC).

\section{Calculate IR by frame components}\label{Sec6.4}

\begin{theorem}{Calculate IR by frame components, Daiyuan Zhang}{Thm6.4.1}\label{Thm6.4.1} 
	Given $\triangle ABC$, point $P$ is an IC, and the three frame components of $P$ are $\alpha _{A}^{P}$, $\alpha _{B}^{P}$, $\alpha _{C}^{P}$ respectively, then its IRs are:
	\[\lambda _{AB}^{P}=\frac{\alpha _{B}^{P}}{\alpha _{A}^{P}},\quad \lambda _{BC}^{P}=\frac{\alpha _{C}^{P}}{\alpha _{B}^{P}},\quad \lambda _{CA}^{P}=\frac{\alpha _{A}^{P}}{\alpha _{C}^{P}}.\]
\end{theorem}

\begin{proof}
	Form theorem \ref{thm:Thm6.1.1} we have
	\[\alpha _{A}^{P}=\frac{1}{1+\lambda _{AB}^{P}+\lambda _{AC}^{P}},\]	
	\[\alpha _{B}^{P}=\frac{\lambda _{AB}^{P}}{1+\lambda _{AB}^{P}+\lambda _{AC}^{P}},\]	
	\[\alpha _{C}^{P}=\frac{\lambda _{AC}^{P}}{1+\lambda _{AB}^{P}+\lambda _{AC}^{P}}.\]	
	
	From the first two formulas we have:
	\[\alpha _{A}^{P}\left( 1+\lambda _{AB}^{P}+\lambda _{AC}^{P} \right)=1,\]	
	\[\alpha _{B}^{P}\left( 1+\lambda _{AB}^{P}+\lambda _{AC}^{P} \right)=\lambda _{AB}^{P}.\]	
	
	The linear equations are obtained as follows:
	\[\left\{ \begin{aligned}
		& \alpha _{A}^{P}\lambda _{AB}^{P}+\alpha _{A}^{P}\lambda _{AC}^{P}=1-\alpha _{A}^{P} \\ 
		& \left( \alpha _{B}^{P}-1 \right)\lambda _{AB}^{P}+\alpha _{B}^{P}\lambda _{AC}^{P}=-\alpha _{B}^{P}. \\ 
	\end{aligned} \right.\]	
	
	Solve the above linear equations to obtain:
	\[\lambda _{AB}^{P}=\frac{\left| \begin{matrix}
			1-\alpha _{A}^{P} & \alpha _{A}^{P}  \\
			-\alpha _{B}^{P} & \alpha _{B}^{P}  \\
		\end{matrix} \right|}{\left| \begin{matrix}
			\alpha _{A}^{P} & \alpha _{A}^{P}  \\
			\alpha _{B}^{P}-1 & \alpha _{B}^{P}  \\
		\end{matrix} \right|}=\frac{\alpha _{B}^{P}}{\alpha _{A}^{P}},\]
	\[\lambda _{AC}^{P}=\frac{\left| \begin{matrix}
			\alpha _{A}^{P} & 1-\alpha _{A}^{P}  \\
			\alpha _{B}^{P}-1 & -\alpha _{B}^{P}  \\
		\end{matrix} \right|}{\left| \begin{matrix}
			\alpha _{A}^{P} & \alpha _{A}^{P}  \\
			\alpha _{B}^{P}-1 & \alpha _{B}^{P}  \\
		\end{matrix} \right|}=\frac{1-\alpha _{A}^{P}-\alpha _{B}^{P}}{\alpha _{A}^{P}}=\frac{\alpha _{C}^{P}}{\alpha _{A}^{P}}.\]
	
	Therefore 
	\[\lambda _{CA}^{P}=\frac{\alpha _{A}^{P}}{\alpha _{C}^{P}}.\]
	
	According to the Ceva's theorem, the following results are obtained:
	\[\lambda _{BC}^{P}=\frac{1}{\lambda _{AB}^{P}\lambda _{CA}^{P}}=\frac{\lambda _{AC}^{P}}{\lambda _{AB}^{P}}=\frac{\alpha _{C}^{P}}{\alpha _{B}^{P}}.\]
\end{proof}
\hfill $\square$\par


\chapter{Some new theorems of plane geometry}\label{Ch7}
\thispagestyle{empty}

In this chapter, two new geometric theorems will be proved by using Intercenter Geometry that I proposed. They are: the theorem of integral ratio from vertex to intersecting foot of a triangle, abbreviated as theorem of IRVIF of a triangle (Theorem \ref{thm:Thm7.3.1}) and the theorem of fractional ratio from vertex to intersecting foot of a triangle, abbreviated as theorem of FRVIF of a triangle (Theorem \ref{thm:Thm7.1.1}). After deriving these two theorems, I have consulted some relevant historical literature on Euclidean geometry and found no similar results. Perhaps these two theorems are new ones I have derived.

\section{Theorem of fractional ratio from vertex to intersecting foot of a triangle}\label{Sec7.1}
This section deduces a new theorem, which is a theorem about triangle. The theorem describes the relationship between the FRVIF of a triangle and the IR of a triangle (See Figure \ref{fig:tu6.1.1}).
\begin{theorem}{Fractional ratio on vertex to intersecting foot of a triangle, Daiyuan Zhang}{Thm7.1.1}\label{Thm7.1.1} 
	For a given $\triangle ABC$,  and $P\in {{\pi }_{ABC}}$, then the relationship between the fractional ratio from vertex to intersecting foot and the intersecting ratio from the triangle is as follows:
	
	\begin{equation}\label{Eq7.1.1}
		{{\lambda }_{AL}}=\lambda _{AB}^{P}+\lambda _{AC}^{P},	
	\end{equation}	
	\begin{equation}\label{Eq7.1.2}
		{{\lambda }_{BM}}=\lambda _{BC}^{P}+\lambda _{BA}^{P},	
	\end{equation}
	\begin{equation}\label{Eq7.1.3}
		{{\lambda }_{CN}}=\lambda _{CA}^{P}+\lambda _{CB}^{P}.	
	\end{equation}
\end{theorem}

\begin{proof}
	As shown in Figure \ref{fig:tu6.1.1}, let the IC be $P\in {{\pi }_{ABC}}$, the point $O$ is any point on the plane $\triangle ABC$. By using theorem \ref{thm:Thm6.1.1}, the following is obtained:	
	\[\overrightarrow{OP}=\frac{\overrightarrow{OA}+\lambda _{AB}^{P}\overrightarrow{OB}+\lambda _{AC}^{P}\overrightarrow{OC}}{1+\lambda _{AB}^{P}+\lambda _{AC}^{P}},\]
	\begin{equation}\label{Eq7.1.4}
		\overrightarrow{OP}-\frac{\overrightarrow{OA}}{1+\lambda _{AB}^{P}+\lambda _{AC}^{P}}=\frac{\lambda _{AB}^{P}\overrightarrow{OB}+\lambda _{AC}^{P}\overrightarrow{OC}}{1+\lambda _{AB}^{P}+\lambda _{AC}^{P}}.
	\end{equation}
	
	First, let the point $O$ coincide with the intersecting foot $L$ of point $A$ (see Figure \ref{fig:tu6.1.1}), that is, $O=L\in \overleftrightarrow{PA}\bigcap \overleftrightarrow{BC}$, i.e., $L$ is both on $\overleftrightarrow{PA}$ and on $\overleftrightarrow{BC}$. Since point $O$ is the intersecting foot $L$ of point $A$, according to the concept of the fractional ratio and integral ratio of segment (see Section 2.1), it follows that:	
	\[\overrightarrow{LP}={{\kappa }_{LA}}\overrightarrow{LA},\]
	\[\overrightarrow{LB}={{\kappa }_{LC}}\overrightarrow{LC},\]
	\[\overrightarrow{BL}=\lambda _{BC}^{P}\overrightarrow{LC},\]
	\[\overrightarrow{LB}=-\lambda _{BC}^{P}\overrightarrow{LC}.\]
	
	Therefore, formula (\ref{Eq7.1.4}) can be written as follows:
	\[\overrightarrow{LP}-\frac{\overrightarrow{LA}}{1+\lambda _{AB}^{P}+\lambda _{AC}^{P}}=\frac{\lambda _{AB}^{P}\overrightarrow{LB}+\lambda _{AC}^{P}\overrightarrow{LC}}{1+\lambda _{AB}^{P}+\lambda _{AC}^{P}},\]	
	i.e.
	\[\left( {{\kappa }_{LA}}-\frac{1}{1+\lambda _{AB}^{P}+\lambda _{AC}^{P}} \right)\overrightarrow{LA}=\frac{-\lambda _{BC}^{P}\lambda _{AB}^{P}+\lambda _{AC}^{P}}{1+\lambda _{AB}^{P}+\lambda _{AC}^{P}}\overrightarrow{LC}.\]
	
	Since the point $L$ is the intersecting foot of the vertex $A$, then $\overrightarrow{LA}$ is linear independent of $\overrightarrow{LC}$, according to the above formula, there must be:	
	\[{{\kappa }_{LA}}-\frac{1}{1+\lambda _{AB}^{P}+\lambda _{AC}^{P}}=0,\]
	\[\frac{-\lambda _{BC}^{P}\lambda _{AB}^{P}+\lambda _{AC}^{P}}{1+\lambda _{AB}^{P}+\lambda _{AC}^{P}}=0,\]
	i.e.
	\[{{\kappa }_{LA}}=\frac{1}{1+\lambda _{AB}^{P}+\lambda _{AC}^{P}},\]	
	\[\frac{-\lambda _{BC}^{P}\lambda _{AB}^{P}+\lambda _{AC}^{P}}{1+\lambda _{AB}^{P}+\lambda _{AC}^{P}}=0.\]
	
	According to formula (\ref{Eq2.4.3}), the following formula holds:
	\[{{\kappa }_{LA}}=\frac{{{\lambda }_{LA}}}{1+{{\lambda }_{LA}}},\]	
	then
	\[\frac{{{\lambda }_{LA}}}{1+{{\lambda }_{LA}}}=\frac{1}{1+\lambda _{AB}^{P}+\lambda _{AC}^{P}},\]
	i.e.
	\[{{\lambda }_{LA}}\left( 1+\lambda _{AB}^{P}+\lambda _{AC}^{P} \right)=1+{{\lambda }_{LA}},\]
	i.e.
	\[{{\lambda }_{LA}}\left( \lambda _{AB}^{P}+\lambda _{AC}^{P} \right)=1,\]
	i.e.
	\[{{\lambda }_{AL}}=\lambda _{AB}^{P}+\lambda _{AC}^{P}.\]	
	
	This proves the formula (\ref{Eq7.1.1}). Similarly, we can prove (\ref{Eq7.1.2}) and (\ref{Eq7.1.3}).
\end{proof}
\hfill $\square$\par

When $L\in \overleftrightarrow{PA}\bigcap \overleftrightarrow{BC}$, the coefficient of $\overrightarrow{OA}$ in theorem \ref{thm:Thm6.1.1} is ${{\kappa }_{LA}}$.

The method of proving (\ref{Eq7.1.2}) and (\ref{Eq7.1.3}) is similar to that of proving formula (\ref{Eq7.1.1}). In fact, because each vertex of the triangle has equal status, we can rotate the subscript of formula (\ref{Eq7.1.1}) to get (\ref{Eq7.1.2}). The specific rotation method is: change the subscript of formula (\ref{Eq7.1.1}) from $A$ to $B$; $B$ to $C$; $C$ to $A$; $L$ to $M$. Similarly, we can rotate the subscript of formula (\ref{Eq7.1.2}) to get (\ref{Eq7.1.3}), and the specific rotation method is: change the subscript of formula (\ref{Eq7.1.2}) from $B$ to $C$; $C$ to $A$; $A$ to $B$; $M$ to $N$. The reason why we can rotate like this is that, after changing $A$ into $B$, $B$ into $C$, $C$ into $A$, $L$ into $M$, it is equivalent to rotating each label of vertex and intersecting foot of the triangle in Figure \ref{fig:tu6.1.1}. So the above process is completely applicable. Of course, we can also prove formulas (\ref{Eq7.1.2}) and (\ref{Eq7.1.3}) just like the previous process of proving formula (\ref{Eq7.1.1}), which is left to the reader as an exercise.

\section{Equation of fractional ratios for triangle}

Using the above theorem, the equation of fractional ratios for a triangle can be obtained.
\begin{theorem}{Equation of fractional ratios for triangle, Daiyuan Zhang}{Thm6.4.2}\label{Thm6.4.2} 
	The fraction ratios of a triangle satisfy the following equation:
	\[\left| \begin{matrix}
		{{\lambda }_{AL}} & -1 & -1  \\
		-1 & {{\lambda }_{BM}} & -1  \\
		-1 & -1 & {{\lambda }_{CN}}  \\
	\end{matrix} \right|=0.\]
	
	i.e.
	\[{{\lambda }_{AL}}{{\lambda }_{BM}}{{\lambda }_{CN}}-\left( {{\lambda }_{AL}}+{{\lambda }_{BM}}+{{\lambda }_{CN}} \right)-2=0.\]
\end{theorem}

\begin{proof}
	From theorem \ref{thm:Thm7.1.1} and theorem \ref{thm:Thm6.4.1}, we have
	\[{{\lambda }_{AL}}=\lambda _{AB}^{P}+\lambda _{AC}^{P}=\frac{\alpha _{B}^{P}}{\alpha _{A}^{P}}+\frac{\alpha _{C}^{P}}{\alpha _{A}^{P}}=\frac{\alpha _{B}^{P}+\alpha _{C}^{P}}{\alpha _{A}^{P}}.\]
	
	Similarly, we have
	\[{{\lambda }_{BM}}=\lambda _{BC}^{P}+\lambda _{BA}^{P}=\frac{\alpha _{C}^{P}}{\alpha _{B}^{P}}+\frac{\alpha _{A}^{P}}{\alpha _{B}^{P}}=\frac{\alpha _{C}^{P}+\alpha _{A}^{P}}{\alpha _{B}^{P}},\]
	\[{{\lambda }_{CN}}=\lambda _{CA}^{P}+\lambda _{CB}^{P}=\frac{\alpha _{A}^{P}}{\alpha _{C}^{P}}+\frac{\alpha _{B}^{P}}{\alpha _{C}^{P}}=\frac{\alpha _{A}^{P}+\alpha _{B}^{P}}{\alpha _{C}^{P}}.\]
	
	i.e.
	\[{{\lambda }_{AL}}\alpha _{A}^{P}-\alpha _{B}^{P}-\alpha _{C}^{P}=0,\]
	\[-\alpha _{A}^{P}+{{\lambda }_{BM}}\alpha _{B}^{P}-\alpha _{C}^{P}=0,\]
	\[-\alpha _{A}^{P}-\alpha _{B}^{P}+{{\lambda }_{CN}}\alpha _{C}^{P}=0.\]
	
	
	Since $\alpha _{A}^{P}$, $\alpha _{B}^{P}$ and $\alpha _{C}^{P}$ has nothing to do with the position of the origin of the frame, if the position of the origin of the frame does not coincide with the point of IC $P$, then $\alpha _{A}^{P}$, $\alpha _{B}^{P}$ and $\alpha _{C}^{P}$ cannot be all 0, that is, there is a non-zero solution for $\alpha _{A}^{P}$, $\alpha _{B}^{P}$ and $\alpha _{C}^{P}$. According to the basic theory of linear equations, the non-zero solution of $\alpha _{A}^{P}$, $\alpha _{B}^{P}$ and $\alpha _{C}^{P}$ means that the determinant of its coefficient matrix is 0, that is, the fractional ratios of the triangle satisfy the following conditions:	
	\[\left| \begin{matrix}
		{{\lambda }_{AL}} & -1 & -1  \\
		-1 & {{\lambda }_{BM}} & -1  \\
		-1 & -1 & {{\lambda }_{CN}}  \\
	\end{matrix} \right|=0.\]
	
	i.e. 
	\[{{\lambda }_{AL}}{{\lambda }_{BM}}{{\lambda }_{CN}}-\left( {{\lambda }_{AL}}+{{\lambda }_{BM}}+{{\lambda }_{CN}} \right)-2=0.\]	
\end{proof}
\hfill $\square$\par
The formula in theorem \ref{thm:Thm6.4.2} is called the equation of fractional ratios for a triangle.

Obviously, in the above formula, if two quantities are known, the other quantity can be obtained. If the three quantities in the above formula are linear with each other, a cubic equation is obtained.

Some interesting formulas can also be obtained by using theorem \ref{thm:Thm6.4.2} For example, the following formula can be proved:
\[\left| \begin{matrix}
	b+c & -a & -a  \\
	-b & c+a & -b  \\
	-c & -c & a+b  \\
\end{matrix} \right|=0.\]

In fact, apply theorem \ref{thm:Thm6.4.2} to the incenter of the triangle, using section \ref{Subsec7.2.2}, we have
\[\left| \begin{matrix}
	\frac{b+c}{a} & -1 & -1  \\
	-1 & \frac{c+a}{b} & -1  \\
	-1 & -1 & \frac{a+b}{c}  \\
\end{matrix} \right|=0,\]

Suppose $a\ne 0$, $b\ne 0$, $c\ne 0$, multiply both sides of the above expression by $abc$, and using the properties of determinant, we have the following result:
\[\left| \begin{matrix}
	b+c & -a & -a  \\
	-b & c+a & -b  \\
	-c & -c & a+b  \\
\end{matrix} \right|=0.\]

Expand the above determinant to obtain:
\[\left( a+b \right)\left( b+c \right)\left( c+a \right)-ab\left( a+b \right)+bc\left( b+c \right)+ca\left( c+a \right)-2abc=0.\]
In fact, the above formula is an identity.


\section{Applications on theorem of fractional ratio from vertex to intersecting foot of a triangle}\label{Sec7.2}
By using theorem \ref{thm:Thm7.1.1}, we can directly get the formula of FRVIF for each IC of a triangle.

\subsection{Fractional ratio from vertex to intersecting foot of a triangle for centroid}\label{Subsec7.2.1}
For the centroid, each IR of the triangle is 1, i.e.
\[\lambda _{AB}^{G}=\lambda _{AC}^{G}=\lambda _{BC}^{G}=\lambda _{BA}^{G}=\lambda _{CA}^{G}=\lambda _{CB}^{G}=1.\]

So the following formula can be obtained from theorem \ref{thm:Thm7.1.1}:
\[{{\lambda }_{AL}}=\frac{\overrightarrow{AG}}{\overrightarrow{GL}}=\lambda _{AB}^{G}+\lambda _{AC}^{G}=2.\]

Similarly, the following formulas can be obtained:
\[{{\lambda }_{BM}}=\frac{\overrightarrow{BG}}{\overrightarrow{GM}}=\lambda _{BC}^{G}+\lambda _{BA}^{G}=2,\]	
\[{{\lambda }_{CN}}=\frac{\overrightarrow{CG}}{\overrightarrow{GN}}=\lambda _{CA}^{G}+\lambda _{CB}^{G}=2.\]	

In this way, we can directly deduce that the length from the vertex to the centroid of the triangle is twice the length of the centroid to the centroid's foot corresponding to the vertex.

\subsection{Fractional ratio from vertex to intersecting foot of a triangle for incenter}\label{Subsec7.2.2}
According to the formula (\ref{Eq5.2.1}), for the incenter $I$ of a triangle, the following formula holds:
\begin{flalign*}
	\lambda _{AB}^{I}=\frac{b}{a}, \lambda _{BC}^{I}=\frac{c}{b}, \lambda _{CA}^{I}=\frac{a}{c}.
\end{flalign*}

Therefore
\begin{flalign*}
	\lambda _{BA}^{I}=\frac{1}{\lambda _{AB}^{I}}=\frac{a}{b}, \lambda _{CB}^{I}=\frac{1}{\lambda _{BC}^{I}}=\frac{b}{c}, \lambda _{AC}^{I}=\frac{1}{\lambda _{CA}^{I}}=\frac{c}{a}.
\end{flalign*}

The following formula is derived directly from theorem \ref{thm:Thm7.1.1}:
\[{{\lambda }_{AL}}=\frac{\overrightarrow{AI}}{\overrightarrow{IL}}=\lambda _{AB}^{I}+\lambda _{AC}^{I}=\frac{b}{a}+\frac{c}{a}=\frac{b+c}{a}.\]

That is, the ratio from the length (between the vertex and the incenter of the triangle) to the length (between the incenter and the incenter's foot corresponding to the vertex) is ${\left( b+c \right)}/{a}\;$, or ${{\lambda }_{AL}}$ is the sum of the length ($b+c$) of the two sides connected to point $A$ divided by the length $a$ of the opposite side of vertex $A$.

Similarly, the following results are obtained from the formula (\ref{Eq5.2.1}):
\[{{\lambda }_{BM}}=\frac{\overrightarrow{BI}}{\overrightarrow{IM}}=\lambda _{BC}^{I}+\lambda _{BA}^{I}=\frac{c}{b}+\frac{a}{b}=\frac{c+a}{b},\]	
\[{{\lambda }_{CN}}=\frac{\overrightarrow{CI}}{\overrightarrow{IN}}=\lambda _{CA}^{I}+\lambda _{CB}^{I}=\frac{a}{c}+\frac{b}{c}=\frac{a+b}{c}.\]

\subsection{Fractional ratio from vertex to intersecting foot of a triangle for orthocenter}\label{Subsec7.2.3}
For the orthocenter $H$ of a triangle, the following results are obtained according to the formula (\ref{Eq5.3.1}):
\[{{\lambda }_{AL}}=\frac{\overrightarrow{AH}}{\overrightarrow{HL}}=\lambda _{AB}^{H}+\lambda _{AC}^{H}=\frac{\tan B}{\tan A}+\frac{\tan C}{\tan A}=\frac{\tan B+\tan C}{\tan A},\]
\[{{\lambda }_{BM}}=\frac{\overrightarrow{BH}}{\overrightarrow{HM}}=\lambda _{BC}^{H}+\lambda _{BA}^{H}=\frac{\tan C}{\tan B}+\frac{\tan A}{\tan B}=\frac{\tan C+\tan A}{\tan B},\]	
\[{{\lambda }_{CN}}=\frac{\overrightarrow{CH}}{\overrightarrow{HN}}=\lambda _{CA}^{H}+\lambda _{CB}^{H}=\frac{\tan A}{\tan C}+\frac{\tan B}{\tan C}=\frac{\tan A+\tan B}{\tan C}.\]	

\subsection{Fractional ratio from vertex to intersecting foot of a triangle for circumcenter}\label{Subsec7.2.4}
For the circumcenter $Q$ of a triangle, the following results are obtained from the formula (\ref{Eq5.4.1}):
\[{{\lambda }_{AL}}=\frac{\overrightarrow{AQ}}{\overrightarrow{QL}}=\lambda _{AB}^{Q}+\lambda _{AC}^{Q}=\frac{\sin 2B}{\sin 2A}+\frac{\sin 2C}{\sin 2A}=\frac{\sin 2B+\sin 2C}{\sin 2A},\]
\[{{\lambda }_{BM}}=\frac{\overrightarrow{BQ}}{\overrightarrow{QM}}=\lambda _{BC}^{Q}+\lambda _{BA}^{Q}=\frac{\sin 2C}{\sin 2B}+\frac{\sin 2A}{\sin 2B}=\frac{\sin 2C+\sin 2A}{\sin 2B},\]	
\[{{\lambda }_{CN}}=\frac{\overrightarrow{CQ}}{\overrightarrow{QN}}=\lambda _{CA}^{Q}+\lambda _{CB}^{Q}=\frac{\sin 2A}{\sin 2C}+\frac{\sin 2B}{\sin 2C}=\frac{\sin 2A+\sin 2B}{\sin 2C}.\]

\subsection{Fractional ratio from vertex to intersecting foot of a triangle for excenter}\label{Subsec7.2.5}
A given triangle has three excenters, and according to the formula (\ref{Eq5.5.1})-(\ref{Eq5.5.2}) the fractional ratio from vertex to intersecting foot of a triangle for excenters can be obtained, which are calculated below.

\subsubsection{Fractional ratio from vertex to intersecting foot of a triangle for excenter ${{E}_{A}}$ corresponding to $\angle A$}

\[{{\lambda }_{AL}}=\frac{\overrightarrow{A{{E}_{A}}}}{\overrightarrow{{{E}_{A}}L}}=\lambda _{AB}^{{{E}_{A}}}+\lambda _{AC}^{{{E}_{A}}}=-\frac{b}{a}-\frac{c}{a}=-\frac{b+c}{a},\]	
\[{{\lambda }_{BM}}=\frac{\overrightarrow{B{{E}_{A}}}}{\overrightarrow{{{E}_{A}}M}}=\lambda _{BC}^{{{E}_{A}}}+\lambda _{BA}^{{{E}_{A}}}=\frac{c}{b}-\frac{a}{b}=\frac{c-a}{b},\]	
\[{{\lambda }_{CN}}=\frac{\overrightarrow{C{{E}_{A}}}}{\overrightarrow{{{E}_{A}}N}}=\lambda _{CA}^{{{E}_{A}}}+\lambda _{CB}^{{{E}_{A}}}=-\frac{a}{c}+\frac{b}{c}=\frac{b-a}{c}.\]

\subsubsection{Fractional ratio from vertex to intersecting foot of a triangle for excenter ${{E}_{B}}$ corresponding to $\angle B$}

\[{{\lambda }_{AL}}=\frac{\overrightarrow{A{{E}_{B}}}}{\overrightarrow{{{E}_{B}}L}}=\lambda _{AB}^{{{E}_{B}}}+\lambda _{AC}^{{{E}_{B}}}=-\frac{b}{a}+\frac{c}{a}=\frac{c-b}{a},\]	
\[{{\lambda }_{BM}}=\frac{\overrightarrow{B{{E}_{B}}}}{\overrightarrow{{{E}_{B}}M}}=\lambda _{BC}^{{{E}_{B}}}+\lambda _{BA}^{{{E}_{B}}}=-\frac{c}{b}-\frac{a}{b}=-\frac{c+a}{b},\]	
\[{{\lambda }_{CN}}=\frac{\overrightarrow{C{{E}_{B}}}}{\overrightarrow{{{E}_{B}}N}}=\lambda _{CA}^{{{E}_{B}}}+\lambda _{CB}^{{{E}_{B}}}=\frac{a}{c}-\frac{b}{c}=\frac{a-b}{c}.\]
\subsubsection{Fractional ratio from vertex to intersecting foot of a triangle for excenter ${{E}_{C}}$ corresponding to $\angle C$}

\[{{\lambda }_{AL}}=\frac{\overrightarrow{A{{E}_{C}}}}{\overrightarrow{{{E}_{C}}L}}=\lambda _{AB}^{{{E}_{C}}}+\lambda _{AC}^{{{E}_{C}}}=\frac{b}{a}-\frac{c}{a}=\frac{b-c}{a},\]	
\[{{\lambda }_{BM}}=\frac{\overrightarrow{B{{E}_{C}}}}{\overrightarrow{{{E}_{C}}M}}=\lambda _{BC}^{{{E}_{C}}}+\lambda _{BA}^{{{E}_{C}}}=-\frac{c}{b}+\frac{a}{b}=\frac{a-c}{b},\]	
\[{{\lambda }_{CN}}=\frac{\overrightarrow{C{{E}_{C}}}}{\overrightarrow{{{E}_{C}}N}}=\lambda _{CA}^{{{E}_{C}}}+\lambda _{CB}^{{{E}_{C}}}=-\frac{a}{c}-\frac{b}{c}=-\frac{a+b}{c}.\]	

\begin{example}\label{Exm7.2.1}
	Proof that the vertex $A$ of $\triangle ABC$, the incenter $I$ and the excenter ${{E}_{A}}$ corresponding to $\angle A$ is on the same straight line; The position of incenter $I$ must lie between vertex $A$ and excenter ${{E}_{A}}$ corresponding to $\angle A$; Finding the ratio from the distance between vertex $A$ and incenter $I$ to the distance between vertex $A$ and the excenter ${{E}_{A}}$.
\end{example}

\begin{solution}
	According to Euclidean geometry, it is easy to prove that the vertex $A$ of $\triangle ABC$, the incenter $I$ and the excenter ${{E}_{A}}$ corresponding to $\angle A$ is on the same straight line. However, if you want to ask a question that how to find the ratio from the distance between vertex $A$ and incenter $I$ to the distance between vertex $A$ and the excenter ${{E}_{A}}$? It is not so easy to use Euclidean geometry to deal with the problem. Let's answer this question by Intercenter Geometry.
		
	First, let's prove by Intercenter Geometry that the vertex $A$ of $\triangle ABC$, the incenter $I$ and the excenter ${{E}_{A}}$ corresponding to $\angle A$ is on the same straight line.
	
	Using theorem \ref{thm:Thm6.1.1}, the following formula is obtained
	\[\overrightarrow{OP}=\frac{\overrightarrow{OA}+\lambda _{AB}^{P}\overrightarrow{OB}+\lambda _{AC}^{P}\overrightarrow{OC}}{1+\lambda _{AB}^{P}+\lambda _{AC}^{P}}.\]
	
	Select the point $O$ at the vertex $A$, and two ICs are the incenter $I$ and the excenter ${{E}_{A}}$. According to the above formulsa and formulas (\ref{Eq5.2.1}), (\ref{Eq5.5.1})-(\ref{Eq5.5.2}), the following results are obtained respectively:
		
	\[\overrightarrow{AI}=\frac{\overrightarrow{AA}+\lambda _{AB}^{I}\overrightarrow{AB}+\lambda _{AC}^{I}\overrightarrow{AC}}{1+\lambda _{AB}^{I}+\lambda _{AC}^{I}}=\frac{\frac{b}{a}\overrightarrow{AB}+\frac{c}{a}\overrightarrow{AC}}{1+\frac{b}{a}+\frac{c}{a}}=\frac{b\overrightarrow{AB}+c\overrightarrow{AC}}{a+b+c},\]
	\[\overrightarrow{A{{E}_{A}}}=\frac{\overrightarrow{AA}+\lambda _{AB}^{{{E}_{A}}}\overrightarrow{AB}+\lambda _{AC}^{{{E}_{A}}}\overrightarrow{AC}}{1+\lambda _{AB}^{{{E}_{A}}}+\lambda _{AC}^{{{E}_{A}}}}=\frac{-\frac{b}{a}\overrightarrow{AB}-\frac{c}{a}\overrightarrow{AC}}{1-\frac{b}{a}-\frac{c}{a}}=-\frac{b\overrightarrow{AB}+c\overrightarrow{AC}}{a-b-c},\]
	therefore
	\[\left( a+b+c \right)\overrightarrow{AI}=\left( b+c-a \right)\overrightarrow{A{{E}_{A}}},\]
	or
	\[\overrightarrow{AI}=\frac{b+c-a}{a+b+c}\overrightarrow{A{{E}_{A}}}.\]
	%
	
	Therefore, $\overrightarrow{AI}\parallel \overrightarrow{A{{E}_{A}}}$. Because vector $\overrightarrow{AI}$ and $\overrightarrow{A{{E}_{A}}}$ has the same starting point $A$, so the vertex $A$ of $\triangle ABC$, the incenter $I$ and the excenter ${{E}_{A}}$ corresponding to $\angle A$ is on the same straight line.
	
	
	Moreover, since all three sides of the triangle are positive real numbers and satisfy $b+c>a$, there must be $a+b+c>b+c-a>0$, then the above formula indicates that the position of incenter $I$ of $\triangle ABC$ must lie between the vertex $A$ and the excenter ${{E}_{A}}$ corresponding to $\angle A$.
		
	
	The ratio from the distance between vertex $A$ and incenter $I$ to the distance between vertex $A$ and excenter ${{E}_{A}}$ is
	\[\rho =\left| \frac{\overrightarrow{AI}}{\overrightarrow{A{{E}_{A}}}} \right|=\frac{b+c-a}{a+b+c}.\]
	
	
	It is not easy to prove the above conclusion by Euclidean geometry. It can be seen that this kind of problem can be easily solved by Intercenter Geometry, which once again highlights the charm of Intercenter Geometry.
	
	
	If it is an equilateral triangle, then, $\rho = {1}/{3}\;$. For equilateral triangles, the distance between vertex $A$ and incenter $I$ is ${1}/{3}\;$ of the distance between vertex $A$ and excenter ${{E}_{A}}$.
\end{solution}
\hfill $\diamond$\par

\section{Theorem of integral ratio from vertex to intersecting foot of a triangle}\label{Sec7.3}

Meaningful conclusion can be obtained by using the uniqueness theorem.
\begin{theorem}{Integral ratio from vertex to intersecting foot of a triangle, Daiyuan Zhang}{Thm7.3.1}\label{Thm7.3.1} 
	Given a $\triangle ABC$, $P\in {{\pi }_{ABC}}$, then the following formula is obtained:
	\[{{\kappa }_{AL}}+{{\kappa }_{BM}}+{{\kappa }_{CN}}=2.\]		
\end{theorem}

\begin{proof}
	According to theorem \ref{thm:Thm7.1.1} and theorem \ref{thm:Thm2.4.1}, there are the following conclusions:
	\[{{\kappa }_{LA}}=\frac{1}{1+\lambda _{AB}^{P}+\lambda _{AC}^{P}},\]	
	\[{{\kappa }_{AL}}=1-{{\kappa }_{LA}}=\frac{\lambda _{AB}^{P}+\lambda _{AC}^{P}}{1+\lambda _{AB}^{P}+\lambda _{AC}^{P}}.\]	
	
	Similarly, we have
	\[{{\kappa }_{BM}}=\frac{\lambda _{BC}^{P}+\lambda _{BA}^{P}}{1+\lambda _{BC}^{P}+\lambda _{BA}^{P}},\]	
	\[{{\kappa }_{CN}}=\frac{\lambda _{CA}^{P}+\lambda _{CB}^{P}}{1+\lambda _{CA}^{P}+\lambda _{CB}^{P}}.\]
	
	Using Ceva's theorem $\lambda _{AB}^{P}\lambda _{BC}^{P}\lambda _{CA}^{P}=1$ can get the following formulas:	
	\begin{align*}
		{{\kappa }_{BM}}&=\frac{\lambda _{BC}^{P}+\lambda _{BA}^{P}}{1+\lambda _{BC}^{P}+\lambda _{BA}^{P}}=\frac{\lambda _{BC}^{P}+\lambda _{BA}^{P}}{1+\lambda _{BC}^{P}+\frac{1}{\lambda _{AB}^{P}}}=\frac{\lambda _{AB}^{P}\left( \lambda _{BC}^{P}+\lambda _{BA}^{P} \right)}{\lambda _{AB}^{P}+\lambda _{AB}^{P}\lambda _{BC}^{P}+1} \\ 
		& =\frac{\lambda _{AB}^{P}\lambda _{BC}^{P}+1}{\lambda _{AB}^{P}+\lambda _{AB}^{P}\lambda _{BC}^{P}+1}=\frac{\frac{1}{\lambda _{CA}^{P}}+1}{\lambda _{AB}^{P}+\frac{1}{\lambda _{CA}^{P}}+1}=\frac{\lambda _{AC}^{P}+1}{\lambda _{AB}^{P}+\lambda _{AC}^{P}+1},  
	\end{align*}	
	\begin{align*}
		{{\kappa }_{CN}}&=\frac{\lambda _{CA}^{P}+\lambda _{CB}^{P}}{1+\lambda _{CA}^{P}+\lambda _{CB}^{P}}=\frac{\lambda _{CA}^{P}+\lambda _{CB}^{P}}{1+\frac{1}{\lambda _{AC}^{P}}+\lambda _{CB}^{P}}=\frac{\lambda _{CA}^{P}+\lambda _{CB}^{P}}{1+\frac{1}{\lambda _{AC}^{P}}+\lambda _{CB}^{P}} \\ 
		& =\frac{\lambda _{AC}^{P}\left( \lambda _{CA}^{P}+\lambda _{CB}^{P} \right)}{\lambda _{AC}^{P}+1+\lambda _{AC}^{P}\lambda _{CB}^{P}}=\frac{1+\lambda _{AC}^{P}\lambda _{CB}^{P}}{\lambda _{AC}^{P}+1+\frac{1}{\lambda _{BA}^{P}}}=\frac{\lambda _{BA}^{P}\left( 1+\lambda _{AC}^{P}\lambda _{CB}^{P} \right)}{\lambda _{BA}^{P}\lambda _{AC}^{P}+\lambda _{BA}^{P}+1} \\ 
		& =\frac{\lambda _{BA}^{P}+1}{\lambda _{BA}^{P}\lambda _{AC}^{P}+\frac{1}{\lambda _{AB}^{P}}+1}=\frac{\lambda _{AB}^{P}\left( \lambda _{BA}^{P}+1 \right)}{\lambda _{AB}^{P}\lambda _{BA}^{P}\lambda _{AC}^{P}+1+\lambda _{AB}^{P}}=\frac{1+\lambda _{AB}^{P}}{\lambda _{AC}^{P}+1+\lambda _{AB}^{P}}. 
	\end{align*}
	
	Then,
	\begin{align*}
		{{\kappa }_{AL}}+{{\kappa }_{BM}}+{{\kappa }_{CN}}&=\frac{\lambda _{AB}^{P}+\lambda _{AC}^{P}}{1+\lambda _{AB}^{P}+\lambda _{AC}^{P}}+\frac{\lambda _{AC}^{P}+1}{\lambda _{AB}^{P}+\lambda _{AC}^{P}+1}+\frac{1+\lambda _{AB}^{P}}{\lambda _{AC}^{P}+1+\lambda _{AB}^{P}}=2. 
	\end{align*}	
\end{proof}
\hfill $\square$\par
Theorem of integral ratio from vertex to intersecting foot of a triangle is also called the theorem of integral ratio from vertex to intersecting foot. Integral ratio from vertex to intersecting foot is abbreviated as IRVIF.

\section{Theorem of intersecting ratio of triangle}\label{Sec7.4}
\begin{theorem}{Intersecting ratio of triangle, Daiyuan Zhang}{Thm7.3.2}\label{Thm7.3.2} 
	Given a $\triangle ABC$, $P\in {{\pi }_{ABC}}$, then the following formula is obtained:	
	\[\frac{1}{1+\lambda _{AB}^{P}+\lambda _{AC}^{P}}+\frac{1}{1+\lambda _{BC}^{P}+\lambda _{BA}^{P}}+\frac{1}{1+\lambda _{CA}^{P}+\lambda _{CB}^{P}}=1.\]			
\end{theorem}

\begin{proof}
	By using the formula (\ref{Eq2.4.3}), the theorem of IRVIF of a triangle (theorem \ref{thm:Thm7.3.1}) can be written in the form of fractional ratio as follows:
	\[\frac{{{\lambda }_{AL}}}{1+{{\lambda }_{AL}}}+\frac{{{\lambda }_{BM}}}{1+{{\lambda }_{BM}}}+\frac{{{\lambda }_{CN}}}{1+{{\lambda }_{CN}}}=2.\]	
	
	Subtracting 1 from both sides, we have:
	\[\frac{1}{1+{{\lambda }_{AL}}}+\frac{1}{1+{{\lambda }_{BM}}}+\frac{1}{1+{{\lambda }_{CN}}}=1.\]	
	
	According to theorem \ref{thm:Thm7.3.1}, the following formula is obtained:
	\[\frac{1}{1+\lambda _{AB}^{P}+\lambda _{AC}^{P}}+\frac{1}{1+\lambda _{BC}^{P}+\lambda _{BA}^{P}}+\frac{1}{1+\lambda _{CA}^{P}+\lambda _{CB}^{P}}=1.\]	
\end{proof}
\hfill $\square$\par
For the IR of the same point in a triangle, the well-known Ceva's theorem (which will be proved by the method of Intercenter Geometry later) gives an identity of IR. I also give an identity of IR here. Ceva's theorem gives the identity in the form of product, and I give the identity in the form of fractional summation. 

Now, two special ICs (centroid and orthocenter) are used to verify the theorem of intersecting ratio of a triangle.

For centroid $G$, we have:
\[\frac{1}{1+\lambda _{AB}^{G}+\lambda _{AC}^{G}}+\frac{1}{1+\lambda _{BC}^{G}+\lambda _{BA}^{G}}+\frac{1}{1+\lambda _{CA}^{G}+\lambda _{CB}^{G}}=\frac{1}{3}+\frac{1}{3}+\frac{1}{3}=1.\]

For orthocenter $H$, we have:
\[\begin{aligned}
	& \frac{1}{1+\lambda _{AB}^{H}+\lambda _{AC}^{H}}+\frac{1}{1+\lambda _{BC}^{H}+\lambda _{BA}^{H}}+\frac{1}{1+\lambda _{CA}^{H}+\lambda _{CB}^{H}} \\ 
	& =\frac{1}{1+\frac{\tan B}{\tan A}+\frac{\tan C}{\tan A}}+\frac{1}{1+\frac{\tan C}{\tan B}+\frac{\tan A}{\tan B}}+\frac{1}{1+\frac{\tan A}{\tan C}+\frac{\tan B}{\tan C}}=1. \\ 
\end{aligned}\]

\chapter{Frame components of triangle and applications}\label{Ch8}
\thispagestyle{empty}

This chapter studies the frame components and their applications. The calculation of frame components is a core problem in Intercenter Geometry. Firstly, the framework components were theoretically studied, and then some applications were presented.

The formulas for calculating the frame components of some intersecting centers (ICs) of a triangle are studied. These special intersecting centers are centroid, incenter, orthocenter, circumcenter and excenter. The frame components of these special centers are represented by the lengths of the three sides of a triangle. In the following chapters, a formula for calculating the distance between two points will be provided, which is represented by the frame component and the length of the three sides of the triangle. This means that if the frame components are represented by the side lengths, then the distance is also represented by the side lengths.

Obviously, obtaining the frame components of the above special points alone is far from enough. In order for the Intercenter Geometry to have broader applications, it is necessary to provide a formula for calculating the frame components at any point.

In order to determine the frame components at any point, theoretical research was first conducted. I have proved that any point on the plane corresponds one-to-one to its frame components (ternary array).

Following the principle of starting from shallow to deep, first study the frame components of the vertices and edges of a triangle. Then solve the following problem: if the frame components of two points are known, the frame components of any point on the line connecting these two points can be obtained.

The frame components of those special points obtained in the previous chapters were all derived from the definition of the IR and obtained using Euclidean geometry methods. However, the use of Euclidean geometry has a serious drawback besides being disconnected from algebra, which is its lack of universality. Usually, it is a strategy of one problem, one solution, and this strategy poses a huge challenge for finding a new point's frame components. In order to determine the frame component at any point on the plane, a new method is needed. The new method I proposed is to use the frame components of known points to calculate the frame components of unknown points, thus opening up a new solution. In fact, the position of any point $P$ on the plane can always be obtained through rotation and stretching transformation based on the known position of point $P_0$. In order to achieve the goal, I proposed an Edge-Axis coordinate system. By utilizing the uniqueness of the frame components and the important characteristic that the frame components are independent of the origin position of the frame, different coordinate systems can be selected in the Edge-Axis coordinate system, and even multiple different coordinate systems can be selected simultaneously to solve the frame components, without the need to calculate in a single coordinate system like in analytic geometry. I think this is a very effective method.

In terms of application, I proposed a theorem that provides a necessary and sufficient condition to solve a class of problems for the concurrent of three straight lines on a plane. Then, the trilinear coordinates were introduced. Unlike traditional trilinear coordinates, I propose the concept of true trilinear coordinates. Then, the relationship between the true trilinear coordinates and the frame components was studied, and the application of the frame components in the true trilinear coordinates was also studied.

I also studied the application of frame components in analytical geometry. In fact, as long as the frame components are obtained, geometric quantities can be easily converted into Cartesian coordinate representations.

The frame component is related to the given triangle. In order to convert between different triangular frames, I studied the transformation of triangular frames.

At the end of this chapter, the center of a triangle was studied. I studied the centers of a triangle from the perspective of fixed points and divided them into several categories. The centers of each category are those fixed points under the corresponding transformation. I have proved that there are infinitely many centers of a given triangle, and that the centers of a triangle cannot fill the entire plane.

This chapter is relatively long, but the content is important.

For convenience, the following notations are used to represent some geometric quantities:

\[p=\frac{1}{2}\left( a+b+c \right),\]
\[W=\sin 2A+\sin 2B+\sin 2C,\]
\[T=\tan A+\tan B+\tan C.\]
\section{Frame components of the centroid of a triangle}\label{Sec8.1}
According to the discussion in the previous chapters, the IR of the centroid is:
\begin{flalign*}
	\lambda _{AB}^{G}=1, \lambda _{BC}^{G}=1, \lambda _{CA}^{G}=1.
\end{flalign*}

Therefore, according to theorem \ref{thm:Thm6.1.1}, it is obtained that:  
\[\alpha _{A}^{G}=\frac{1}{1+\lambda _{AB}^{G}+\lambda _{AC}^{G}}=\frac{1}{1+1+1}=\frac{1}{3},\]	
\[\alpha _{B}^{G}=\frac{\lambda _{AB}^{G}}{1+\lambda _{AB}^{G}+\lambda _{AC}^{G}}=\frac{1}{1+1+1}=\frac{1}{3},\]	
\[\alpha _{C}^{G}=\frac{\lambda _{AC}^{G}}{1+\lambda _{AB}^{G}+\lambda _{AC}^{G}}=\frac{1}{1+1+1}=\frac{1}{3}.\]

\section{Frame components of the incenter of a triangle}\label{Sec8.2}
According to the discussion in the previous chapters, the IR of the incenter is:
\begin{flalign*}
	\lambda _{AB}^{I}=\frac{b}{a}, \lambda _{BC}^{I}=\frac{c}{b}, \lambda _{CA}^{I}=\frac{a}{c}.
\end{flalign*}

Therefore, according to theorem \ref{thm:Thm6.1.1}, it is obtained that: 
\[\alpha _{A}^{I}=\frac{1}{1+\lambda _{AB}^{I}+\lambda _{AC}^{I}}=\frac{1}{1+\frac{b}{a}+\frac{c}{a}}=\frac{a}{a+b+c}=\frac{a}{2p},\]	
\[\alpha _{B}^{I}=\frac{\lambda _{AB}^{I}}{1+\lambda _{AB}^{I}+\lambda _{AC}^{I}}=\frac{\frac{b}{a}}{1+\frac{b}{a}+\frac{c}{a}}=\frac{b}{a+b+c}=\frac{b}{2p},\]	
\[\alpha _{C}^{I}=\frac{\lambda _{AC}^{I}}{1+\lambda _{AB}^{I}+\lambda _{AC}^{I}}=\frac{\frac{c}{a}}{1+\frac{b}{a}+\frac{c}{a}}=\frac{c}{a+b+c}=\frac{c}{2p}.\]

\section{Frame components of the orthocenter of a triangle}\label{Sec8.3}
According to the discussion in the previous chapters, the IR of the orthocenter is:
\[\lambda _{AB}^{H}=\frac{AM}{MB}=\frac{CA\cos A}{BC\cos B}=\frac{b\cos A}{a\cos B}=\frac{\frac{b}{\cos B}}{\frac{a}{\cos A}}=\frac{\tan B}{\tan A}=\frac{\cot A}{\cot B},\]	
\[\lambda _{BC}^{H}=\frac{AQ}{QB}=\frac{AB\cos B}{CA\cos C}=\frac{c\cos B}{b\cos C}=\frac{\frac{c}{\cos C}}{\frac{b}{\cos B}}=\frac{\tan C}{\tan B}=\frac{\cot B}{\cot C},\]	
\[\lambda _{CA}^{H}=\frac{CN}{NA}=\frac{BC\cos C}{AB\cos A}=\frac{a\cos C}{c\cos A}=\frac{\frac{a}{\cos A}}{\frac{c}{\cos C}}=\frac{\tan A}{\tan C}=\frac{\cot C}{\cot A}.\]	

Therefore, according to theorem \ref{thm:Thm6.1.1}, it is obtained that: 
\begin{align*}
	\alpha _{A}^{H}& =\frac{1}{1+\lambda _{AB}^{H}+\lambda _{AC}^{H}}=\frac{1}{1+\frac{\tan B}{\tan A}+\frac{\tan C}{\tan A}} \\ 
	& =\frac{\tan A}{\tan A+\tan B+\tan C}=\frac{\tan A}{T},  
\end{align*}	
\begin{align*}
	\alpha _{B}^{H}& =\frac{\lambda _{AB}^{H}}{1+\lambda _{AB}^{H}+\lambda _{AC}^{H}}=\frac{\frac{\tan B}{\tan A}}{1+\frac{\tan B}{\tan A}+\frac{\tan C}{\tan A}} \\ 
	& =\frac{\tan B}{\tan A+\tan B+\tan C}=\frac{\tan B}{T},  
\end{align*}	
\begin{align*}
	\alpha _{C}^{H}& =\frac{\lambda _{AC}^{H}}{1+\lambda _{AB}^{H}+\lambda _{AC}^{H}}=\frac{\frac{\tan C}{\tan A}}{1+\frac{\tan B}{\tan A}+\frac{\tan C}{\tan A}} \\ 
	& =\frac{\tan C}{\tan A+\tan B+\tan C}=\frac{\tan C}{T}.  
\end{align*}

\section{Frame components of the circumcenter of a triangle}\label{Sec8.4}
According to the discussion in the previous chapters, the IR of the circumcenter is:
\begin{flalign*}
	\lambda _{AB}^{Q}=\frac{\sin 2B}{\sin 2A}, \lambda _{BC}^{Q}=\frac{\sin 2C}{\sin 2B}, \lambda _{CA}^{Q}=\frac{\sin 2A}{\sin 2C}.
\end{flalign*}

Therefore, according to theorem \ref{thm:Thm6.1.1}, it is obtained that: 
\begin{align*}
	\alpha _{A}^{Q}& =\frac{1}{1+\lambda _{AB}^{Q}+\lambda _{AC}^{Q}}=\frac{1}{1+\frac{\sin 2B}{\sin 2A}+\frac{\sin 2C}{\sin 2A}} \\ 
	& =\frac{\sin 2A}{\sin 2A+\sin 2B+\sin 2C}=\frac{\sin 2A}{W},  
\end{align*}	
\begin{align*}
	\alpha _{B}^{Q}& =\frac{\lambda _{AB}^{Q}}{1+\lambda _{AB}^{Q}+\lambda _{AC}^{Q}}=\frac{\frac{\sin 2B}{\sin 2A}}{1+\frac{\sin 2B}{\sin 2A}+\frac{\sin 2C}{\sin 2A}} \\ 
	& =\frac{\sin 2B}{\sin 2A+\sin 2B+\sin 2C}=\frac{\sin 2B}{W},  
\end{align*}	
\begin{align*}
	\alpha _{C}^{Q}& =\frac{\lambda _{AC}^{Q}}{1+\lambda _{AB}^{Q}+\lambda _{AC}^{Q}}=\frac{\frac{\sin 2C}{\sin 2A}}{1+\frac{\sin 2B}{\sin 2A}+\frac{\sin 2C}{\sin 2A}} \\ 
	& =\frac{\sin 2C}{\sin 2A+\sin 2B+\sin 2C}=\frac{\sin 2C}{W}.  
\end{align*}

\section{Frame components of the excenter of a triangle}\label{Sec8.5}
\subsection{Frame components of excenter corresponding to $\angle A$}
According to the discussion in the previous chapters, the IR of the excenter corresponding to $\angle A$ is:
\begin{flalign*}
	\lambda _{AB}^{{{E}_{A}}}=-\frac{b}{a}, \lambda _{BC}^{{{E}_{A}}}=\frac{c}{b}, \lambda _{CA}^{{{E}_{A}}}=-\frac{a}{c}.
\end{flalign*}

Therefore, according to theorem \ref{thm:Thm6.1.1}, it is obtained that: 
\[\alpha _{A}^{{{E}_{A}}}=\frac{1}{1+\lambda _{AB}^{{{E}_{A}}}+\lambda _{AC}^{{{E}_{A}}}}=\frac{1}{1-\frac{b}{a}-\frac{c}{a}}=\frac{a}{a-b-c}=-\frac{a}{2\left( p-a \right)},\]
\[\alpha _{B}^{{{E}_{A}}}=\frac{\lambda _{AB}^{{{E}_{A}}}}{1+\lambda _{AB}^{{{E}_{A}}}+\lambda _{AC}^{{{E}_{A}}}}=\frac{-\frac{b}{a}}{1-\frac{b}{a}-\frac{c}{a}}=-\frac{b}{a-b-c}=\frac{b}{2\left( p-a \right)},\]
\[\alpha _{C}^{{{E}_{A}}}=\frac{\lambda _{AC}^{{{E}_{A}}}}{1+\lambda _{AB}^{{{E}_{A}}}+\lambda _{AC}^{{{E}_{A}}}}=\frac{-\frac{c}{a}}{1-\frac{b}{a}-\frac{c}{a}}=-\frac{c}{a-b-c}=\frac{c}{2\left( p-a \right)}.\]
\subsection{Frame components of excenter corresponding to $\angle B$}
According to the discussion in the previous chapters, the IR of the  excenter corresponding to $\angle B$ is:
\begin{flalign*}
	\lambda _{AB}^{{{E}_{B}}}=-\frac{b}{a}, \lambda _{BC}^{{{E}_{B}}}=-\frac{c}{b}, \lambda _{CA}^{{{E}_{B}}}=\frac{a}{c}.	
\end{flalign*}

Therefore, according to theorem \ref{thm:Thm6.1.1}, it is obtained that: 
\[\alpha _{A}^{{{E}_{B}}}=\frac{1}{1+\lambda _{AB}^{{{E}_{B}}}+\lambda _{AC}^{{{E}_{B}}}}=\frac{1}{1-\frac{b}{a}+\frac{c}{a}}=\frac{a}{a-b+c}=\frac{a}{2\left( p-b \right)},\]
\[\alpha _{B}^{{{E}_{B}}}=\frac{\lambda _{AB}^{{{E}_{B}}}}{1+\lambda _{AB}^{{{E}_{B}}}+\lambda _{AC}^{{{E}_{B}}}}=\frac{-\frac{b}{a}}{1-\frac{b}{a}+\frac{c}{a}}=-\frac{b}{a-b+c}=-\frac{b}{2\left( p-b \right)},\]
\[\alpha _{C}^{{{E}_{B}}}=\frac{\lambda _{AB}^{{{E}_{B}}}}{1+\lambda _{AB}^{{{E}_{B}}}+\lambda _{AC}^{{{E}_{B}}}}=\frac{\frac{c}{a}}{1-\frac{b}{a}+\frac{c}{a}}=\frac{c}{a-b+c}=\frac{c}{2\left( p-b \right)}.\]
\subsection{Frame components of excenter corresponding to $\angle C$}
According to the discussion in the previous chapters, the IR of the  excenter corresponding to $\angle C$ is:
\begin{flalign*}
	\lambda _{AB}^{{{E}_{C}}}=\frac{b}{a}, \lambda _{BC}^{{{E}_{C}}}=-\frac{c}{b}, \lambda _{CA}^{{{E}_{C}}}=-\frac{a}{c}.	
\end{flalign*}

Therefore, according to theorem \ref{thm:Thm6.1.1}, it is obtained that: 
\[\alpha _{A}^{{{E}_{C}}}=\frac{1}{1+\lambda _{AB}^{{{E}_{C}}}+\lambda _{AC}^{{{E}_{C}}}}=\frac{1}{1+\frac{b}{a}-\frac{c}{a}}=\frac{a}{a+b-c}=\frac{a}{2\left( p-c \right)},\]
\[\alpha _{B}^{{{E}_{C}}}=\frac{\lambda _{AB}^{{{E}_{C}}}}{1+\lambda _{AB}^{{{E}_{C}}}+\lambda _{AC}^{{{E}_{C}}}}=\frac{\frac{b}{a}}{1+\frac{b}{a}-\frac{c}{a}}=\frac{b}{a+b-c}=\frac{b}{2\left( p-c \right)},\]
\[\alpha _{C}^{{{E}_{C}}}=\frac{\lambda _{AC}^{{{E}_{C}}}}{1+\lambda _{AB}^{{{E}_{C}}}+\lambda _{AC}^{{{E}_{C}}}}=\frac{-\frac{c}{a}}{1+\frac{b}{a}-\frac{c}{a}}=-\frac{c}{a+b-c}=-\frac{c}{2\left( p-c \right)}.\]


\section{The points on the plane correspond one-to-one with the frame components}\label{PingmianshangDeDianYuBiaojiafenliangYiyiduiying}
Combine the frame components of $\triangle ABC$ into a ternary array, denoted as $\left( \alpha _{A}^{P},\alpha _{B}^{P},\alpha _{C}^{P} \right)$, where $\alpha _{A}^{P}+\alpha _{B}^{P}+\alpha _{C}^{P}=1$. So a correspondence is established between a point $P$ on the plane and this ternary array $\left( \alpha _{A}^{P},\alpha _{B}^{P},\alpha _{C}^{P} \right)$. The current question is whether point $P$ corresponds one-to-one with the frame component $\left( \alpha _{A}^{P},\alpha _{B}^{P},\alpha _{C}^{P} \right)$? The following theorem provides a positive answer.


\begin{theorem}{Points correspond one-to-one with a ternary array of real numbers, Daiyuan Zhang}{PingmianshangDeDianYuSanyuanshishushuzuYiyiduiying}\label{PingmianshangDeDianYuSanyuanshishushuzuYiyiduiying}
	Given a $\triangle ABC$, the triangular frame is $\left( O;A,B,C \right)$, $P$ is a point on the plane of $\triangle ABC$, and $\left( x,y,z \right)$ is a ternary array of real numbers that satisfies the following conditions:
	\[\overrightarrow{OP}=x\overrightarrow{OA}+y\overrightarrow{OB}+z\overrightarrow{OC},\]
	\[x+y+z=1.\]
	Then point $P$ corresponds one-to-one with $\left( x,y,z \right)$.
\end{theorem}

\begin{proof}
	Firstly, let's assume that there is a point $P$ on the plane of $\triangle ABC$, which corresponds to two ternary arrays of real numbers. These two ternary arrays of real numbers are $\left( {{x}_{1}},{{y}_{1}},{{z}_{1}} \right)$ and $\left( {{x}_{2}},{{y}_{2}},{{z}_{2}} \right)$ respectively, and the following conditions are satisfied:
	\[\overrightarrow{OP}={{x}_{1}}\overrightarrow{OA}+{{y}_{1}}\overrightarrow{OB}+{{z}_{1}}\overrightarrow{OC},\]
	\[{{x}_{1}}+{{y}_{1}}+{{z}_{1}}=1.\]
	\[\overrightarrow{OP}={{x}_{2}}\overrightarrow{OA}+{{y}_{2}}\overrightarrow{OB}+{{z}_{2}}\overrightarrow{OC},\]
	\[{{x}_{2}}+{{y}_{2}}+{{z}_{2}}=1.\]
	
	Then
	\[\left( {{x}_{1}}\overrightarrow{OA}+{{y}_{1}}\overrightarrow{OB}+{{z}_{1}}\overrightarrow{OC} \right)-\left( {{x}_{2}}\overrightarrow{OA}+{{y}_{2}}\overrightarrow{OB}+{{z}_{2}}\overrightarrow{OC} \right)=\overrightarrow{OP}-\overrightarrow{OP}=\overrightarrow{0},\]
	
	i.e.
	\[\left( {{x}_{1}}-{{x}_{2}} \right)\overrightarrow{OA}+\left( {{y}_{1}}-{{y}_{2}} \right)\overrightarrow{OB}+\left( {{z}_{1}}-{{z}_{2}} \right)\overrightarrow{OC}=\overrightarrow{0}.\]
	
	And
	\[\left( {{x}_{1}}-{{x}_{2}} \right)+\left( {{y}_{1}}-{{y}_{2}} \right)+\left( {{z}_{1}}-{{z}_{2}} \right)=\left( {{x}_{1}}+{{y}_{1}}+{{z}_{1}} \right)-\left( {{x}_{2}}+{{y}_{2}}+{{z}_{2}} \right)=1-1=0.\]
	
	According to theorem \ref{thm:Thm3.4.2}, it is obtained that ${{x}_{1}}={{x}_{2}}$, ${{y}_{1}}={{y}_{2}}$, ${{z}_{1}}={{z}_{2}}$.
	
	On the other hand, assuming there are two points ${{P}_{1}}$ and ${{P}_{2}}$ on the plane of $\triangle ABC$, they all correspond to the same real  ternary array $\left( x,y,z \right)$, and satisfy the following conditions:
	\[\overrightarrow{O{{P}_{1}}}=x\overrightarrow{OA}+y\overrightarrow{OB}+z\overrightarrow{OC},\]
	\[x+y+z=1.\]
	\[\overrightarrow{O{{P}_{2}}}=x\overrightarrow{OA}+y\overrightarrow{OB}+z\overrightarrow{OC},\]
	\[x+y+z=1.\]
	
	Then
	\[\overrightarrow{{{P}_{1}}{{P}_{2}}}=\overrightarrow{O{{P}_{2}}}-\overrightarrow{O{{P}_{1}}}=\overrightarrow{0}.\]
	
	The above equation means two points ${{P}_{1}}$ and ${{P}_{2}}$ on the plane of $\triangle ABC$ is overlapping together.
	
	Based on the above discussion, it is known that the point $P$ on the plane of $\triangle ABC$ corresponds one-to-one with the real ternary array $\left( x,y,z \right)$.
\end{proof}
\hfill $\square$\par


\begin{theorem}{Points on the plane correspond one-to-one with frame components, Daiyuan Zhang}{PingmianShangDeDianYuBiaojiafenliangYiyiduiying}\label{PingmianShangDeDianYuBiaojiafenliangYiyiduiying}
	Given a $\triangle ABC$, where $P$ is a point on the $\triangle ABC$ plane, the triangular frame is $\left( O;A,B,C \right)$, and the frame components of point $P$ are $\alpha _{A}^{P}$, $\alpha _{B}^{P}$ and $\alpha _{C}^{P}$, then the point $P$ and the frame components $\left( \alpha _{A}^{P},\alpha _{B}^{P},\alpha _{C}^{P} \right)$ of point $P$ corresponds one-to-one. 
\end{theorem}

%

\begin{proof}
	Since $P$ is a point on the plane of $\triangle ABC$, and the frame component of point $P$ is $\left( \alpha _{A}^{P},\alpha _{B}^{P},\alpha _{C}^{P} \right)$, so under the triangular frame of $\left( O;A,B,C \right)$, we have:
	\[\overrightarrow{OP}=\alpha _{A}^{P}\overrightarrow{OA}+\alpha _{B}^{P}\overrightarrow{OB}+\alpha _{C}^{P}\overrightarrow{OC},\]
	\[\alpha _{A}^{P}+\alpha _{B}^{P}+\alpha _{C}^{P}=1.\]
	
	According to theorem \ref{thm:PingmianshangDeDianYuSanyuanshishushuzuYiyiduiying}, it is obtained that the point $P$ and the ordered array $\left( \alpha _{A}^{P},\alpha _{B}^{P},\alpha _{C}^{P} \right)$ of the frame components of point $P$ corresponds one-to-one.
\end{proof}
\hfill $\square$\par


Point $P$ can also be written as $P\left( \alpha _{A}^{P},\alpha _{B}^{P},\alpha _{C}^{P} \right)$ or $\left( \alpha _{A}^{P},\alpha _{B}^{P},\alpha _{C}^{P} \right)$.


The above theorem indicates that for a given $\triangle ABC$, under the triangular frame $\left( O;A,B,C \right)$, the point $P$ on the plane of $\triangle ABC$ can be uniquely represented by the frame component, that is, it can be uniquely represented as $P\left( \alpha _{A}^{P},\alpha _{B}^{P},\alpha _{C}^{P} \right)$, where $\alpha _{A}^{P}+\alpha _{B}^{P}+\alpha _{C}^{P}=1$. This method of representing point is called \textbf{Frame Component Representation of Point}, abbreviated as Frame Component Representation. The representation of frame component is different from the coordinate representation of point in analytical geometry. The frame component is independent of the position of the frame origin, while the coordinate are related to the position of the coordinate origin.


The point $P$ discussed here is the IC $P$. The IC can be seen as a point represented by frame components. This book also refers to IC $P$ as point $P$ or $P$.


\section{The frame components of a vertex and a point on an edge of a triangle}\label{SanjiaoxingDingdianYuBianshangYidianDeBiaojiafenliang} 

This section investigates the frame components of vertices and points on edges of the triangle.
\begin{theorem}{Frame components of the vertices of a triangle, Daiyuan Zhang}{SanjiaoxingDingdianDeBiaojiafenliang}\label{SanjiaoxingDingdianDeBiaojiafenliang}
	Given a $\triangle ABC$, under the triangular frame $\left( O;A,B,C \right)$, the frame components of the three vertices $A$, $B$, and $C$ are:
	\[\alpha _{A}^{A}=1,\quad \alpha _{B}^{A}=0,\quad \alpha _{C}^{A}=0;\]
	\[\alpha _{A}^{B}=0,\quad \alpha _{B}^{B}=1,\quad \alpha _{C}^{B}=0;\]
	\[\alpha _{A}^{C}=0,\quad \alpha _{B}^{C}=0,\quad \alpha _{C}^{C}=1.\]
\end{theorem}

\begin{proof}
	According to the definition of IR, it can be obtained that:
	\[\lambda _{AB}^{P}=\frac{\overrightarrow{AN}}{\overrightarrow{NB}}.\]		
	Select point $P$ at vertex A of $\triangle ABC$, that is, let $P:=A$, then $N=A$  (see Figure \ref{fig:tu6.1.1}), and obtain:
	\[\lambda _{AB}^{A}=\frac{\overrightarrow{AA}}{\overrightarrow{AB}}.\]		
	The above equation is equivalent to $\lambda _{AB}^{A}\overrightarrow{AB}=\overrightarrow{AA}$. Due to $\overrightarrow{AB}\ne \overrightarrow{0}$, $\overrightarrow{AA}=\overrightarrow{0}$, so we can only get $\lambda _{AB}^{A}=0$.

	And
	\[\lambda _{CA}^{P}=\frac{\overrightarrow{CM}}{\overrightarrow{MA}}.\]		
	
	Select point $P$ at vertex $A$ of $\triangle ABC$, that is, make $P:=A$ result in $M=A$  (see Figure \ref{fig:tu6.1.1}), and obtain:
	
	\[\lambda _{CA}^{A}=\frac{\overrightarrow{CM}}{\overrightarrow{MA}}=\frac{\overrightarrow{CA}}{\overrightarrow{AA}}.\]
	
	Therefore
	\[\lambda _{AC}^{A}=\frac{1}{\lambda _{CA}^{A}}=\frac{\overrightarrow{MA}}{\overrightarrow{CM}}=\frac{\overrightarrow{AA}}{\overrightarrow{CA}}.\]
	
	The above equation is equivalent to $\lambda _{AC}^{A}\overrightarrow{CA}=\overrightarrow{AA}$. Due to $\overrightarrow{CA}\ne \overrightarrow{0}$, $\overrightarrow{AA}=\overrightarrow{0}$, so we can only get $\lambda _{AC}^{A}=0$

	According to theorem \ref{thm:Thm6.1.1}, it can be obtained that:
	\[\alpha _{A}^{P}=\frac{1}{1+\lambda _{AB}^{P}+\lambda _{AC}^{P}},\]
	\[\alpha _{B}^{P}=\frac{\lambda _{AB}^{P}}{1+\lambda _{AB}^{P}+\lambda _{AC}^{P}},\]	
	\[\alpha _{C}^{P}=\frac{\lambda _{AC}^{P}}{1+\lambda _{AB}^{P}+\lambda _{AC}^{P}}.\]	
	
	Selecting point $P$ at vertex $A$ of $\triangle ABC$ yields:
	\[\alpha _{A}^{A}=\frac{1}{1+\lambda _{AB}^{A}+\lambda _{AC}^{A}}=\frac{1}{1+0+0}=1,\]
	\[\alpha _{B}^{A}=\frac{\lambda _{AB}^{A}}{1+\lambda _{AB}^{A}+\lambda _{AC}^{A}}=\frac{0}{1+0+0}=0,\]	
	\[\alpha _{B}^{A}=\frac{\lambda _{AB}^{A}}{1+\lambda _{AB}^{A}+\lambda _{AC}^{A}}=\frac{0}{1+0+0}=0.\]	
	
	This proves that $\alpha _{A}^{A}=1$, $\alpha _{B}^{A}=0$, $\alpha _{C}^{A}=0$.
	
	Similar proof can be obtained for the other two cases.
\end{proof}
\hfill $\square$\par


According to the above theorem, the three vertices of the $\triangle ABC$ can be denoted as $A\left( 1,0,0 \right)$, $B\left( 0,1,0 \right)$, $C\left( 0,0,1 \right)$, respectively. They can be regarded as a set of basis vectors in three-dimensional space. The $\triangle ABC$ is called a \textbf{base triangle} or a \textbf{reference triangle}.

%
%
%

\begin{theorem}{Frame components of point on straight line of triangle edges, Daiyuan Zhang}{SanjiaoxingBiandianDeBiaojiafenliang}\label{SanjiaoxingBiandianDeBiaojiafenliang}
	
	Given a $\triangle ABC$, under the triangular frame $\left( O;A,B,C \right)$, the following conclusion holds:
	
	If $P\in \overleftrightarrow{BC}$, then the frame components of point $P$ are:
	\[\alpha _{A}^{P}=0,\quad \alpha _{B}^{P}=\frac{1}{1+\lambda _{BC}^{P}},\quad \alpha _{C}^{P}=\frac{\lambda _{BC}^{P}}{1+\lambda _{BC}^{P}};\]
	
	If $P\in \overleftrightarrow{CA}$, then the frame components of point $P$ are:
	\[\alpha _{B}^{P}=0,\quad \alpha _{C}^{P}=\frac{1}{1+\lambda _{CA}^{P}},\quad \alpha _{A}^{P}=\frac{\lambda _{CA}^{P}}{1+\lambda _{CA}^{P}};\]
	
	If $P\in \overleftrightarrow{AB}$, then the frame components of point $P$ are:
	\[\alpha _{C}^{P}=0,\quad \alpha _{A}^{P}=\frac{1}{1+\lambda _{AB}^{P}},\quad \alpha _{B}^{P}=\frac{\lambda _{AB}^{P}}{1+\lambda _{AB}^{P}}.\]
\end{theorem}

\begin{proof}
	Let's consider the first scenario ($P\in \overleftrightarrow{BC}$).
	
	According to the definition of IR, it can be obtained that:
	\[\lambda _{AB}^{P}=\frac{\overrightarrow{AN}}{\overrightarrow{NB}},\]	
	\[\lambda _{BA}^{P}=\frac{1}{\lambda _{AB}^{P}}=\frac{\overrightarrow{NB}}{\overrightarrow{AN}}.\]	
	
	If $P\in \overleftrightarrow{BC}$, then $M\to C$ and $N\to B$ (see Figure\ref{fig:tu6.1.1}), therefore
	\[\lambda _{BA}^{P}=\frac{1}{\lambda _{AB}^{P}}=\frac{\overrightarrow{NB}}{\overrightarrow{AN}}=\frac{\overrightarrow{BB}}{\overrightarrow{AB}}.\]	
	
	The above equation is equivalent to $\lambda _{BA}^{P}\overrightarrow{AB}=\overrightarrow{BB}$. Due to $\overrightarrow{AB}\ne \overrightarrow{0}$, $\overrightarrow{BB}=\overrightarrow{0}$, so we can only get $\lambda _{BA}^{P}=0$.
	
	According to theorem \ref{thm:Thm6.1.1}, under the triangular frame $\left( O;A,B,C \right)$, if $P\in \overleftrightarrow{BC}$, then the frame components of point $P$ are:
	\[\alpha _{B}^{P}=\frac{1}{1+\lambda _{BC}^{P}+\lambda _{BA}^{P}}=\frac{1}{1+\lambda _{BC}^{P}},\]	
	\[\alpha _{C}^{P}=\frac{\lambda _{BC}^{P}}{1+\lambda _{BC}^{P}+\lambda _{BA}^{P}}=\frac{\lambda _{BC}^{P}}{1+\lambda _{BC}^{P}},\]	
	\[\alpha _{A}^{P}=\frac{\lambda _{BA}^{P}}{1+\lambda _{BC}^{P}+\lambda _{BA}^{P}}=0.\]	
	
	Now, let's consider the second scenario ($P\in \overleftrightarrow{CA}$).
	
	According to the definition of IR, it can be obtained that:
	\[\lambda _{BC}^{P}=\frac{\overrightarrow{BL}}{\overrightarrow{LC}},\]	
	\[\lambda _{CB}^{P}=\frac{1}{\lambda _{BC}^{P}}=\frac{\overrightarrow{LC}}{\overrightarrow{BL}}.\]	
	
	If $P\in \overleftrightarrow{CA}$, then $N\to A$ and $L\to C$ (see Figure\ref{fig:tu6.1.1}), therefore
	\[\lambda _{CB}^{P}=\frac{1}{\lambda _{BC}^{P}}=\frac{\overrightarrow{LC}}{\overrightarrow{BL}}=\frac{\overrightarrow{CC}}{\overrightarrow{BC}}.\]	
	
	
	The above equation is equivalent to $\lambda _{CB}^{P}\overrightarrow{BC}=\overrightarrow{CC}$. Due to $\overrightarrow{BC}\ne \overrightarrow{0}$, $\overrightarrow{CC}=\overrightarrow{0}$, so we can only get $\lambda _{CB}^{P}=0$.
	
	According to theorem \ref{thm:Thm6.1.1}, under the triangular frame $\left( O;A,B,C \right)$, if $P\in \overleftrightarrow{CA}$, then the frame components of point $P$ are:
	\[\alpha _{C}^{P}=\frac{1}{1+\lambda _{CA}^{P}+\lambda _{CB}^{P}}=\frac{1}{1+\lambda _{CA}^{P}},\]	
	\[\alpha _{A}^{P}=\frac{\lambda _{CA}^{P}}{1+\lambda _{CA}^{P}+\lambda _{CB}^{P}}=\frac{\lambda _{CA}^{P}}{1+\lambda _{CA}^{P}},\]	
	\[\alpha _{B}^{P}=\frac{\lambda _{CB}^{P}}{1+\lambda _{CA}^{P}+\lambda _{CB}^{P}}=0.\]	
	
	Let's take a look at the last scenario ($P\in \overleftrightarrow{AB}$).
	
	According to the definition of IR, it can be obtained that:
	\[\lambda _{CA}^{P}=\frac{\overrightarrow{CM}}{\overrightarrow{MA}},\]	
	\[\lambda _{AC}^{P}=\frac{1}{\lambda _{CA}^{P}}=\frac{\overrightarrow{MA}}{\overrightarrow{CM}}.\]	
	
	If $P\in \overleftrightarrow{AB}$, then $L\to B$ and $M\to A$ (see Figure\ref{fig:tu6.1.1}), therefore
	\[\lambda _{AC}^{P}=\frac{1}{\lambda _{CA}^{P}}=\frac{\overrightarrow{MA}}{\overrightarrow{CM}}=\frac{\overrightarrow{AA}}{\overrightarrow{CA}}.\]	
	
	The above equation is equivalent to $\lambda _{AC}^{P}\overrightarrow{CA}=\overrightarrow{AA}$. Due to $\overrightarrow{CA}\ne \overrightarrow{0}$, $\overrightarrow{AA}=\overrightarrow{0}$, so we can only get $\lambda _{AC}^{P}=0$.
		
	According to theorem \ref{thm:Thm6.1.1}, under the triangular frame $\left( O;A,B,C \right)$, if $P\in \overleftrightarrow{AB}$, then the frame components of point $P$ are:
	\[\alpha _{A}^{P}=\frac{1}{1+\lambda _{AB}^{P}+\lambda _{AC}^{P}}=\frac{1}{1+\lambda _{AB}^{P}},\]
	\[\alpha _{B}^{P}=\frac{\lambda _{AB}^{P}}{1+\lambda _{AB}^{P}+\lambda _{AC}^{P}}=\frac{\lambda _{AB}^{P}}{1+\lambda _{AB}^{P}},\]	
	\[\alpha _{C}^{P}=\frac{\lambda _{AC}^{P}}{1+\lambda _{AB}^{P}+\lambda _{AC}^{P}}=0.\]	
\end{proof}
\hfill $\square$\par



The above theorem explains that when a point is located on an edge or its extension of a given triangle, the frame component of the vertex opposite to that edge is zero, and the frame component on that edge is the line segment ratio of that edge.

Question: Under what circumstances is the denominator of the frame component not zero? I will answer this question after introducing the trilinear coordinates later.


\section{The frame component at any point on the line connecting two points}\label{LianggeHuixinShangRenyiYidianDeBiaojiafenliang}
If the frame components of two points are known, how to calculate the frame components at any point on the line connecting these two points? The following theorem answers this question.
\begin{theorem}{Frame component at any point on the line connecting two ICs, Daiyuan Zhang}{HuixinLianxianShangRenyiYidianDeBiaojiafenliang}\label{HuixinLianxianShangRenyiYidianDeBiaojiafenliang}
	Given a $\triangle ABC$, assuming the frame component of IC ${{P}_{1}}$ is $\alpha _{A}^{{{P}_{1}}}$, $\alpha _{B}^{{{P}_{1}}}$, $\alpha _{C}^{{{P}_{1}}}$; the frame component of IC ${{P}_{2}}$ is $\alpha _{A}^{{{P}_{2}}}$, $\alpha _{B}^{{{P}_{2}}}$, $\alpha _{C}^{{{P}_{2}}}$, $P\in \overleftrightarrow{{{P}_{1}}{{P}_{2}}}$ and $\overrightarrow{{{P}_{1}}P}={{\kappa }_{{{P}_{1}}{{P}_{2}}}}\overrightarrow{{{P}_{1}}{{P}_{2}}}$, $\overrightarrow{{{P}_{1}}{{P}_{2}}}\ne \overrightarrow{0}$, then the frame component of the IC $P$ is:
	\[\alpha _{A}^{P}=\alpha _{A}^{{{P}_{1}}}+{{\kappa }_{{{P}_{1}}{{P}_{2}}}}\left( \alpha _{A}^{{{P}_{2}}}-\alpha _{A}^{{{P}_{1}}} \right)=\left( 1-{{\kappa }_{{{P}_{1}}{{P}_{2}}}} \right)\alpha _{A}^{{{P}_{1}}}+{{\kappa }_{{{P}_{1}}{{P}_{2}}}}\alpha _{A}^{{{P}_{2}}},\]
	\[\alpha _{B}^{P}=\alpha _{B}^{{{P}_{1}}}+{{\kappa }_{{{P}_{1}}{{P}_{2}}}}\left( \alpha _{B}^{{{P}_{2}}}-\alpha _{B}^{{{P}_{1}}} \right)=\left( 1-{{\kappa }_{{{P}_{1}}{{P}_{2}}}} \right)\alpha _{B}^{{{P}_{1}}}+{{\kappa }_{{{P}_{1}}{{P}_{2}}}}\alpha _{B}^{{{P}_{2}}},\]
	\[\alpha _{C}^{P}=\alpha _{C}^{{{P}_{1}}}+{{\kappa }_{{{P}_{1}}{{P}_{2}}}}\left( \alpha _{C}^{{{P}_{2}}}-\alpha _{C}^{{{P}_{1}}} \right)=\left( 1-{{\kappa }_{{{P}_{1}}{{P}_{2}}}} \right)\alpha _{C}^{{{P}_{1}}}+{{\kappa }_{{{P}_{1}}{{P}_{2}}}}\alpha _{C}^{{{P}_{2}}}.\]
\end{theorem}

\begin{proof}
	For a given $\triangle ABC$, whose three vertices $A$, $B$, and $C$ do not coincide with each other, let point $O$ be any point in space, and according to theorem \ref{thm:Thm6.1.1}, vector $\overrightarrow{O{{P}_{1}}}$ can be uniquely linearly represented by the triangular frame $\left( O;A,B,C \right)$, i.e
	\[\overrightarrow{O{{P}_{1}}}=\alpha _{A}^{{{P}_{1}}}\overrightarrow{OA}+\alpha _{B}^{{{P}_{1}}}\overrightarrow{OB}+\alpha _{C}^{{{P}_{1}}}\overrightarrow{OC},\]
	\[\alpha _{A}^{{{P}_{1}}}+\alpha _{B}^{{{P}_{1}}}+\alpha _{C}^{{{P}_{1}}}=1.\]
	
	Similarly, the vector $\overrightarrow{O{{P}_{2}}}$ and $\overrightarrow{OP}$ can be uniquely linearly represented by the triangular frame $\left( O;A,B,C \right)$, i.e
	\[\overrightarrow{O{{P}_{2}}}=\alpha _{A}^{{{P}_{2}}}\overrightarrow{OA}+\alpha _{B}^{{{P}_{2}}}\overrightarrow{OB}+\alpha _{C}^{{{P}_{2}}}\overrightarrow{OC},\]
	\[\alpha _{A}^{{{P}_{2}}}+\alpha _{B}^{{{P}_{2}}}+\alpha _{C}^{{{P}_{2}}}=1.\]
	\[\overrightarrow{OP}=\alpha _{A}^{P}\overrightarrow{OA}+\alpha _{B}^{P}\overrightarrow{OB}+\alpha _{C}^{P}\overrightarrow{OC},\]
	\[\alpha _{A}^{P}+\alpha _{B}^{P}+\alpha _{C}^{P}=1.\]
	
	If point $O$ coincides with point $A$, then according to theorem \ref{thm:Thm6.1.3}, it can be obtained that:
	\[\left\{ \begin{aligned}
		& \overrightarrow{A{{P}_{1}}}=\alpha _{B}^{{{P}_{1}}}\overrightarrow{AB}+\alpha _{C}^{{{P}_{1}}}\overrightarrow{AC} \\ 
		& \alpha _{A}^{{{P}_{1}}}+\alpha _{B}^{{{P}_{1}}}+\alpha _{C}^{{{P}_{1}}}=1, \\ 
	\end{aligned} \right.\]
	\[\left\{ \begin{aligned}
		& \overrightarrow{A{{P}_{2}}}=\alpha _{B}^{{{P}_{2}}}\overrightarrow{AB}+\alpha _{C}^{{{P}_{2}}}\overrightarrow{AC} \\ 
		& \alpha _{A}^{{{P}_{2}}}+\alpha _{B}^{{{P}_{2}}}+\alpha _{C}^{{{P}_{2}}}=1. \\ 
	\end{aligned} \right.\]
	
	Therefore
	\[\overrightarrow{{{P}_{1}}{{P}_{2}}}=\overrightarrow{A{{P}_{2}}}-\overrightarrow{A{{P}_{1}}}=\left( \alpha _{B}^{{{P}_{2}}}-\alpha _{B}^{{{P}_{1}}} \right)\overrightarrow{AB}+\left( \alpha _{C}^{{{P}_{2}}}-\alpha _{C}^{{{P}_{1}}} \right)\overrightarrow{AC},\]
	\[{{\kappa }_{{{P}_{1}}{{P}_{2}}}}\overrightarrow{{{P}_{1}}{{P}_{2}}}={{\kappa }_{{{P}_{1}}{{P}_{2}}}}\left( \alpha _{B}^{{{P}_{2}}}-\alpha _{B}^{{{P}_{1}}} \right)\overrightarrow{AB}+{{\kappa }_{{{P}_{1}}{{P}_{2}}}}\left( \alpha _{C}^{{{P}_{2}}}-\alpha _{C}^{{{P}_{1}}} \right)\overrightarrow{AC}.\]
	
	And
	\[\overrightarrow{{{P}_{1}}P}=\overrightarrow{AP}-\overrightarrow{A{{P}_{1}}}=\left( \alpha _{B}^{P}-\alpha _{B}^{{{P}_{1}}} \right)\overrightarrow{AB}+\left( \alpha _{C}^{P}-\alpha _{C}^{{{P}_{1}}} \right)\overrightarrow{AC}.\]
	
	Since \[\overrightarrow{{{P}_{1}}P}={{\kappa }_{{{P}_{1}}{{P}_{2}}}}\overrightarrow{{{P}_{1}}{{P}_{2}}},\] 
	therefore
	\[\left( \left( \alpha _{B}^{P}-\alpha _{B}^{{{P}_{1}}} \right)-{{\kappa }_{{{P}_{1}}{{P}_{2}}}}\left( \alpha _{B}^{{{P}_{2}}}-\alpha _{B}^{{{P}_{1}}} \right) \right)\overrightarrow{AB}+\left( \left( \alpha _{C}^{P}-\alpha _{C}^{{{P}_{1}}} \right)-{{\kappa }_{{{P}_{1}}{{P}_{2}}}}\left( \alpha _{C}^{{{P}_{2}}}-\alpha _{C}^{{{P}_{1}}} \right) \right)\overrightarrow{AC}=\overrightarrow{0}.\]
	
	Since $\overrightarrow{AB}$ is linearly independent of $\overrightarrow{AC}$, the following results are obtained:
	\[\left( \alpha _{B}^{P}-\alpha _{B}^{{{P}_{1}}} \right)-{{\kappa }_{{{P}_{1}}{{P}_{2}}}}\left( \alpha _{B}^{{{P}_{2}}}-\alpha _{B}^{{{P}_{1}}} \right)=0,\]
	\[\left( \alpha _{C}^{P}-\alpha _{C}^{{{P}_{1}}} \right)-{{\kappa }_{{{P}_{1}}{{P}_{2}}}}\left( \alpha _{C}^{{{P}_{2}}}-\alpha _{C}^{{{P}_{1}}} \right)=0.\]
	
	i.e.
	\[\alpha _{B}^{P}=\alpha _{B}^{{{P}_{1}}}+{{\kappa }_{{{P}_{1}}{{P}_{2}}}}\left( \alpha _{B}^{{{P}_{2}}}-\alpha _{B}^{{{P}_{1}}} \right),\]
	\[\alpha _{C}^{P}=\alpha _{C}^{{{P}_{1}}}+{{\kappa }_{{{P}_{1}}{{P}_{2}}}}\left( \alpha _{C}^{{{P}_{2}}}-\alpha _{C}^{{{P}_{1}}} \right).\]
	
	Because
	\[\alpha _{A}^{P}+\alpha _{B}^{P}+\alpha _{C}^{P}=1,\]
	it is obtained that:
	\begin{align*}
		\alpha _{A}^{P}& =1-\alpha _{B}^{P}-\alpha _{C}^{P}=1-\left( \alpha _{B}^{{{P}_{1}}}+{{\kappa }_{{{P}_{1}}{{P}_{2}}}}\left( \alpha _{B}^{{{P}_{2}}}-\alpha _{B}^{{{P}_{1}}} \right) \right)-\left( \alpha _{C}^{{{P}_{1}}}+{{\kappa }_{{{P}_{1}}{{P}_{2}}}}\left( \alpha _{C}^{{{P}_{2}}}-\alpha _{C}^{{{P}_{1}}} \right) \right) \\ 
		& =1-\alpha _{B}^{{{P}_{1}}}-\alpha _{C}^{{{P}_{1}}}-{{\kappa }_{{{P}_{1}}{{P}_{2}}}}\left( \left( \alpha _{B}^{{{P}_{2}}}-\alpha _{B}^{{{P}_{1}}} \right)+\left( \alpha _{C}^{{{P}_{2}}}-\alpha _{C}^{{{P}_{1}}} \right) \right) \\ 
		& =\alpha _{A}^{{{P}_{1}}}-{{\kappa }_{{{P}_{1}}{{P}_{2}}}}\left( 1-\alpha _{A}^{{{P}_{2}}}-\left( 1-\alpha _{A}^{{{P}_{1}}} \right) \right)=\alpha _{A}^{{{P}_{1}}}-{{\kappa }_{{{P}_{1}}{{P}_{2}}}}\left( \alpha _{A}^{{{P}_{1}}}-\alpha _{A}^{{{P}_{2}}} \right) \\ 
		& =\alpha _{A}^{{{P}_{1}}}+{{\kappa }_{{{P}_{1}}{{P}_{2}}}}\left( \alpha _{A}^{{{P}_{2}}}-\alpha _{A}^{{{P}_{1}}} \right)=\left( 1-{{\kappa }_{{{P}_{1}}{{P}_{2}}}} \right)\alpha _{A}^{{{P}_{1}}}+{{\kappa }_{{{P}_{1}}{{P}_{2}}}}\alpha _{A}^{{{P}_{2}}}.  
	\end{align*}
\end{proof}
\hfill $\square$\par


As an application, the following study focuses on the frame components of the midpoint of two ICs. If the midpoint of a line segment ${{P}_{1}}{{P}_{2}}$ is $P$, the following conclusion can be proven.


\begin{corollary}{Frame component of midpoint, Daiyuan Zhang}{ZhongdianDeBiaojiafenliang}\label{ZhongdianDeBiaojiafenliang}
	Given a $\triangle ABC$, assuming the frame component of IC ${{P}_{1}}$ is $\alpha _{A}^{{{P}_{1}}}$, $\alpha _{B}^{{{P}_{1}}}$, $\alpha _{C}^{{{P}_{1}}}$; the frame component of IC ${{P}_{2}}$ is $\alpha _{A}^{{{P}_{2}}}$, $\alpha _{B}^{{{P}_{2}}}$, $\alpha _{C}^{{{P}_{2}}}$, point $P$ is the midpoint of the line segment ${{P}_{1}}{{P}_{2}}$, then the frame components of midpoint $P$ are:
	\[\alpha _{A}^{P}=\frac{\alpha _{A}^{{{P}_{1}}}+\alpha _{A}^{{{P}_{2}}}}{2},\]
	\[\alpha _{B}^{P}=\frac{\alpha _{B}^{{{P}_{1}}}+\alpha _{B}^{{{P}_{2}}}}{2},\]
	\[\alpha _{C}^{P}=\frac{\alpha _{C}^{{{P}_{1}}}+\alpha _{C}^{{{P}_{2}}}}{2}.\]
\end{corollary}


\begin{proof}
	Let ${{\kappa }_{{{P}_{1}}{{P}_{2}}}}={1}/{2}\;$ in theorem \ref{thm:HuixinLianxianShangRenyiYidianDeBiaojiafenliang} to obtain the proof.
\end{proof}
\hfill $\square$\par


The above formula indicates that the frame component of the midpoint is the arithmetic mean of the frame components of the two endpoints.


\begin{example}{}\label{LiyongOulaDingliQiuWaixinDeBiaojiafenliang}
	Let ${{P}_{1}}$ be the centroid, ${{P}_{2}}$ be the orthocenter, $\overrightarrow{{{P}_{1}}P}={{\kappa }_{{{P}_{1}}{{P}_{2}}}}\overrightarrow{{{P}_{1}}{{P}_{2}}}$, ${{\kappa }_{{{P}_{1}}{{P}_{2}}}}=-{1}/{2}\;$, Find the frame component of $P$.
\end{example}

%

\begin{solution}
	According to section \ref {Sec8.1}, section \ref {Sec8.3}, and theorem \ref{thm:HuixinLianxianShangRenyiYidianDeBiaojiafenliang}, it is obtained that
	\begin{align*}
		\alpha _{A}^{P}& =\alpha _{A}^{{{P}_{1}}}+{{\kappa }_{{{P}_{1}}{{P}_{2}}}}\left( \alpha _{A}^{{{P}_{2}}}-\alpha _{A}^{{{P}_{1}}} \right)=\alpha _{A}^{G}+{{\kappa }_{{{P}_{1}}{{P}_{2}}}}\left( \alpha _{A}^{H}-\alpha _{A}^{G} \right) \\ 
		& =\frac{1}{3}-\frac{1}{2}\left( \frac{\tan A}{\tan A+\tan B+\tan C}-\frac{1}{3} \right) \\ 
		& =\frac{1}{3}-\frac{2\tan A-\tan B-\tan C}{6\left( \tan A+\tan B+\tan C \right)} \\ 
		& =\frac{2\left( \tan A+\tan B+\tan C \right)-\left( 2\tan A-\tan B-\tan C \right)}{6\left( \tan A+\tan B+\tan C \right)} \\ 
		& =\frac{\tan B+\tan C}{2\left( \tan A+\tan B+\tan C \right)}.  
	\end{align*}
	
	Similarly, it can be obtained that:
	\begin{align*}
		\alpha _{B}^{P}& =\alpha _{B}^{{{P}_{1}}}+{{\kappa }_{{{P}_{1}}{{P}_{2}}}}\left( \alpha _{B}^{{{P}_{2}}}-\alpha _{B}^{{{P}_{1}}} \right)=\alpha _{B}^{G}+{{\kappa }_{{{P}_{1}}{{P}_{2}}}}\left( \alpha _{B}^{H}-\alpha _{B}^{G} \right) \\ 
		& =\frac{\tan C+\tan A}{2\left( \tan A+\tan B+\tan C \right)},  
	\end{align*}
	\begin{align*}
		\alpha _{C}^{P}& =\alpha _{C}^{{{P}_{1}}}+{{\kappa }_{{{P}_{1}}{{P}_{2}}}}\left( \alpha _{C}^{{{P}_{2}}}-\alpha _{C}^{{{P}_{1}}} \right)=\alpha _{C}^{G}+{{\kappa }_{{{P}_{1}}{{P}_{2}}}}\left( \alpha _{C}^{H}-\alpha _{C}^{G} \right) \\ 
		& =\frac{\tan A+\tan B}{2\left( \tan A+\tan B+\tan C \right)}.  
	\end{align*}
\end{solution}
\hfill $\diamond$\par


According to Euler's theorem (see theorem \ref{thm:Thm10.3.1}), the point $P$ happens to be the circumcenter $Q$, and by the way, the following trigonometric identities are obtained (see section \ref{Sec8.4}):
\[\frac{\sin 2A}{\sin 2A+\sin 2B+\sin 2C}=\frac{\tan B+\tan C}{2\left( \tan A+\tan B+\tan C \right)},\]
\[\frac{\sin 2B}{\sin 2A+\sin 2B+\sin 2C}=\frac{\tan C+\tan A}{2\left( \tan A+\tan B+\tan C \right)},\]
\[\frac{\sin 2C}{\sin 2A+\sin 2B+\sin 2C}=\frac{\tan A+\tan B}{2\left( \tan A+\tan B+\tan C \right)}.\]


\begin{example}{}\label{YuShangyigeLiziJinxingBijiao}
	Find the frame component of the midpoint between the vertex and the orthocenter of a triangle.
\end{example}

\begin{solution}
	The problem can be described as: let ${{P}_{1}}$ be the vertex $A$ of $\triangle ABC$, ${{P}_{2}}$ be the orthocenter $H$ of $\triangle ABC$, and $\overrightarrow{{{P}_{1}}P}={{\kappa }_{{{P}_{1}}{{P}_{2}}}}\overrightarrow{{{P}_{1}}{{P}_{2}}}$ can be written as $\overrightarrow{AP}={{\kappa }_{AH}}\overrightarrow{AH}$, ${{\kappa }_{AH}}={1}/{2}\;$, find the frame component of ${{P}_{m\left( AH \right)}}$. ${{P}_{m\left( AH \right)}}$ refers to the midpoint between vertex $A$ and orthocenter $H$. According to section \ref{Sec8.3} and theorem \ref{thm:HuixinLianxianShangRenyiYidianDeBiaojiafenliang}, it is obtained that:
	\begin{align*}
		\alpha _{A}^{{{P}_{m\left( AH \right)}}}& =\left( 1-{{\kappa }_{{{P}_{1}}{{P}_{2}}}} \right)\alpha _{A}^{{{P}_{1}}}+{{\kappa }_{{{P}_{1}}{{P}_{2}}}}\alpha _{A}^{{{P}_{2}}}=\left( 1-{{\kappa }_{AH}} \right)\alpha _{A}^{A}+{{\kappa }_{AH}}\alpha _{A}^{H} \\ 
		& =\left( 1-\frac{1}{2} \right)+\frac{1}{2}\alpha _{A}^{H}=\frac{1+\alpha _{A}^{H}}{2}=\frac{1+\frac{\tan A}{\tan A+\tan B+\tan C}}{2} \\ 
		& =\frac{2\tan A+\tan B+\tan C}{2\left( \tan A+\tan B+\tan C \right)},  
	\end{align*}
	\begin{align*}
		\alpha _{B}^{{{P}_{m\left( AH \right)}}}& =\left( 1-{{\kappa }_{{{P}_{1}}{{P}_{2}}}} \right)\alpha _{B}^{{{P}_{1}}}+{{\kappa }_{{{P}_{1}}{{P}_{2}}}}\alpha _{B}^{{{P}_{2}}}=\left( 1-{{\kappa }_{AH}} \right)\alpha _{B}^{A}+{{\kappa }_{AH}}\alpha _{B}^{H} \\ 
		& =\frac{1}{2}\alpha _{B}^{H}=\frac{\tan B}{2\left( \tan A+\tan B+\tan C \right)},  
	\end{align*}
	\begin{align*}
		& \alpha _{C}^{{{P}_{m\left( AH \right)}}}=\left( 1-{{\kappa }_{{{P}_{1}}{{P}_{2}}}} \right)\alpha _{C}^{{{P}_{1}}}+{{\kappa }_{{{P}_{1}}{{P}_{2}}}}\alpha _{C}^{{{P}_{2}}}=\left( 1-{{\kappa }_{AH}} \right)\alpha _{C}^{A}+{{\kappa }_{AH}}\alpha _{C}^{H} \\ 
		& =\frac{1}{2}\alpha _{C}^{H}=\frac{\tan C}{2\left( \tan A+\tan B+\tan C \right)}.  
	\end{align*}
	
	Similarly, the frame component of the midpoint between the vertex $B$ and the orthocenter $H$ of $\triangle ABC$ can be obtained as follows:
	\[\alpha _{A}^{{{P}_{m\left( BH \right)}}}=\frac{\tan A}{2\left( \tan A+\tan B+\tan C \right)},\]
	\[\alpha _{B}^{{{P}_{m\left( BH \right)}}}=\frac{\tan A+2\tan B+\tan C}{2\left( \tan A+\tan B+\tan C \right)},\]
	\[\alpha _{C}^{{{P}_{m\left( BH \right)}}}=\frac{\tan C}{2\left( \tan A+\tan B+\tan C \right)}.\]
	
	The frame component of the midpoint between the vertex $C$ and the orthocenter $H$ of $\triangle ABC$ can be obtained as follows:
	\[\alpha _{A}^{{{P}_{m\left( CH \right)}}}=\frac{\tan A}{2\left( \tan A+\tan B+\tan C \right)},\]
	\[\alpha _{B}^{{{P}_{m\left( CH \right)}}}=\frac{\tan B}{2\left( \tan A+\tan B+\tan C \right)},\]
	\[\alpha _{C}^{{{P}_{m\left( CH \right)}}}=\frac{\tan A+\tan B+2\tan C}{2\left( \tan A+\tan B+\tan C \right)}.\]
\end{solution}
\hfill $\diamond$\par

The following example studies the frame components of points on the Line between vertex and intersecting center (LVIC).


\begin{example}{}\label{DinghuixianShangDianDeBiaojiafenliang}
	Given a $\triangle ABC$, denote the IC as $P$, where the frame component of the IC $P$ is $\alpha _{A}^{P}$, $\alpha _{B}^{P}$, $\alpha _{C}^{P}$. Let ${{P}_{AP}}\in \overleftrightarrow{AP}$ and $\overrightarrow{A{{P}_{AP}}}={{\kappa }_{AP}}\overrightarrow{AP}$, $\overrightarrow{AP}\ne \overrightarrow{0}$; ${{P}_{BP}}\in \overleftrightarrow{BP}$ and $\overrightarrow{B{{P}_{BP}}}={{\kappa }_{BP}}\overrightarrow{BP}$, $\overrightarrow{BP}\ne \overrightarrow{0}$; ${{P}_{CP}}\in \overleftrightarrow{CP}$ and $\overrightarrow{C{{P}_{CP}}}={{\kappa }_{CP}}\overrightarrow{CP}$, $\overrightarrow{CP}\ne \overrightarrow{0}$. Find the frame components of ${{P}_{AP}}$, ${{P}_{BP}}$ and ${{P}_{CP}}$
\end{example}

\begin{solution}
	According to theorem \ref{thm:HuixinLianxianShangRenyiYidianDeBiaojiafenliang}, the frame component of point ${{P}_{AP}}$ on the LVIC of $\triangle ABC$ is:
	\[\alpha _{A}^{{{P}_{AP}}}=\left( 1-{{\kappa }_{AP}} \right)\alpha _{A}^{A}+{{\kappa }_{AP}}\alpha _{A}^{P}=1-{{\kappa }_{AP}}+{{\kappa }_{AP}}\alpha _{A}^{P},\]
	\[\alpha _{B}^{{{P}_{AP}}}=\left( 1-{{\kappa }_{AP}} \right)\alpha _{B}^{A}+{{\kappa }_{AP}}\alpha _{B}^{P}={{\kappa }_{AP}}\alpha _{B}^{P},\]
	\[\alpha _{C}^{{{P}_{AP}}}=\left( 1-{{\kappa }_{AP}} \right)\alpha _{C}^{A}+{{\kappa }_{AP}}\alpha _{C}^{P}={{\kappa }_{AP}}\alpha _{C}^{P}.\]
	
	The frame component of point ${{P}_{BP}}$ on the LVIC of $\triangle ABC$ is:
	\[\alpha _{A}^{{{P}_{BP}}}=\left( 1-{{\kappa }_{BP}} \right)\alpha _{A}^{B}+{{\kappa }_{BP}}\alpha _{A}^{P}={{\kappa }_{BP}}\alpha _{A}^{P},\]
	\[\alpha _{B}^{{{P}_{BP}}}=\left( 1-{{\kappa }_{BP}} \right)\alpha _{B}^{B}+{{\kappa }_{BP}}\alpha _{B}^{P}=1-{{\kappa }_{BP}}+{{\kappa }_{BP}}\alpha _{B}^{P},\]
	\[\alpha _{C}^{{{P}_{BP}}}=\left( 1-{{\kappa }_{BP}} \right)\alpha _{C}^{B}+{{\kappa }_{BP}}\alpha _{C}^{P}={{\kappa }_{BP}}\alpha _{C}^{P}.\]
	
	The frame component of point ${{P}_{CP}}$ on the LVIC of $\triangle ABC$ is:
	\[\alpha _{A}^{{{P}_{CP}}}=\left( 1-{{\kappa }_{CP}} \right)\alpha _{A}^{C}+{{\kappa }_{CP}}\alpha _{A}^{P}={{\kappa }_{CP}}\alpha _{A}^{P},\]
	\[\alpha _{B}^{{{P}_{CP}}}=\left( 1-{{\kappa }_{CP}} \right)\alpha _{B}^{C}+{{\kappa }_{CP}}\alpha _{B}^{P}={{\kappa }_{CP}}\alpha _{B}^{P},\]
	\[\alpha _{C}^{{{P}_{CP}}}=\left( 1-{{\kappa }_{CP}} \right)\alpha _{C}^{C}+{{\kappa }_{CP}}\alpha _{C}^{P}=1-{{\kappa }_{CP}}+{{\kappa }_{CP}}\alpha _{C}^{P}.\]
\end{solution}
\hfill $\diamond$\par

Let's study the symmetry of the IC about a point. Specifically, it is necessary to study the IC and the calculation of the frame components of the symmetric IC. 
%
\begin{definition}{Symmetric Center and Symmetric IC}{DuichenzhongxinYuDuichenhuixin}\label{DuichenzhongxinYuDuichenhuixin}
	Given a $\triangle ABC$, assuming ${{P}_{0}}$, $P$ and ${P}'$ are three ICs on the $\triangle ABC$ plane, and ${{P}_{0}}$, $P$ and ${P}'$ are on the same straight line (i.e. ${P}'\in \overleftrightarrow{{{P}_{0}}P}$), if $\overrightarrow{P{P}'}=2\overrightarrow{P{{P}_{0}}}$, then ${{P}_{0}}$ is called the symmetric center, and ${P}'$ is called the symmetric IC of the point $P$ associated with the symmetric center ${{P}_{0}}$.
\end{definition}


\begin{corollary}{Frame component of symmetrical IC, Daiyuan Zhang}{DuichenhuixinDeBiaojiafenliang}\label{DuichenhuixinDeBiaojiafenliang}
	Given a $\triangle ABC$, assuming ${{P}_{0}}$, $P$ and ${P}'$ are three ICs on the $\triangle ABC$ plane, and ${{P}_{0}}$, $P$ and ${P}'$ are on the same straight line (i.e. ${P}'\in \overleftrightarrow{{{P}_{0}}P}$), if $\overrightarrow{P{P}'}=2\overrightarrow{P{{P}_{0}}}$, The frame component of symmetric center ${{P}_{0}}$ is $\alpha _{A}^{{{P}_{0}}}$, $\alpha _{B}^{{{P}_{0}}}$, $\alpha _{C}^{{{P}_{0}}}$; The frame component of IC $P$ is $\alpha _{A}^{P}$, $\alpha _{B}^{P}$, $\alpha _{C}^{P}$; The frame component of symmetric IC ${P}'$ is $\alpha _{A}^{{{P}'}}$, $\alpha _{B}^{{{P}'}}$, $\alpha _{C}^{{{P}'}}$, then the frame component of the symmetric IC ${P}'$ is:
	\[\alpha _{A}^{{{P}'}}=2\alpha _{A}^{{{P}_{0}}}-\alpha _{A}^{P},\]
	\[\alpha _{B}^{{{P}'}}=2\alpha _{B}^{{{P}_{0}}}-\alpha _{B}^{P},\]
	\[\alpha _{C}^{{{P}'}}=2\alpha _{C}^{{{P}_{0}}}-\alpha _{C}^{P}.\]	
\end{corollary}


\begin{proof}
	In theorem \ref{thm:HuixinLianxianShangRenyiYidianDeBiaojiafenliang}, let ${{P}_{1}}=P$, ${{P}_{2}}={{P}_{0}}$, ${{\kappa }_{{{P}_{1}}{{P}_{2}}}}=2$, we obtain this corollary.	
\end{proof}
\hfill $\square$\par


\begin{example}{}\label{ZhongxinGuanyuDingdianDuichendianDeBiaojiafenliang}
	Given a $\triangle ABC$, find the frame component of the symmetric point of centroid $G$ with respect to the vertex $A$ of $\triangle ABC$.
\end{example}


\begin{solution}
	Assuming that the symmetric point of the centroid $G$ with respect to vertex $A$ (symmetric center) is ${G}'$, then based on the above corollary and sections \ref{Sec8.1} and \ref{Sec8.2}, it can be concluded that:
	\[\alpha _{A}^{{{G}'}}=2\alpha _{A}^{A}-\alpha _{A}^{G}=2-\frac{1}{3}=\frac{5}{3},\]
	\[\alpha _{B}^{{{G}'}}=2\alpha _{B}^{A}-\alpha _{B}^{G}=0-\frac{1}{3}=-\frac{1}{3},\]
	\[\alpha _{C}^{{{G}'}}=2\alpha _{C}^{A}-\alpha _{C}^{G}=0-\frac{1}{3}=-\frac{1}{3}.\]
\end{solution}
\hfill $\diamond$\par


\begin{example}{}\label{NeixinGuanyuZhongxinDuichendianDeBiaojiafenliang}
	Given a $\triangle ABC$, find the frame component of the symmetric point of the incenter $I$ with respect to the centroid $G$ of $\triangle ABC$.
\end{example}


\begin{solution}
	The centroid $G$ is the symmetric center. Assuming that the symmetric point of the incenter $I$ with respect to the symmetric center $G$ is ${I}'$, based on the above corollary and sections \ref{Sec8.1} and \ref{Sec8.2}, it can be concluded that:
	\begin{align*}
		\alpha _{A}^{{{I}'}}& =2\alpha _{A}^{G}-\alpha _{A}^{I}=2\times \frac{1}{3}-\frac{a}{a+b+c} \\ 
		& =\frac{2\left( a+b+c \right)-3a}{3\left( a+b+c \right)}=\frac{2\left( b+c \right)-a}{3\left( a+b+c \right)},  
	\end{align*}
	\begin{align*}
		\alpha _{B}^{{{I}'}}& =2\alpha _{B}^{G}-\alpha _{B}^{I}=2\times \frac{1}{3}-\frac{b}{a+b+c} \\ 
		& =\frac{2\left( a+b+c \right)-3b}{3\left( a+b+c \right)}=\frac{2\left( c+a \right)-b}{3\left( a+b+c \right)},  
	\end{align*}
	\begin{align*}
		& \alpha _{C}^{{{I}'}}=2\alpha _{C}^{G}-\alpha _{C}^{I}=2\times \frac{1}{3}-\frac{c}{a+b+c} \\ 
		& =\frac{2\left( a+b+c \right)-3c}{3\left( a+b+c \right)}=\frac{2\left( a+b \right)-c}{3\left( a+b+c \right)}.  
	\end{align*}
\end{solution}
\hfill $\diamond$\par


\section{The frame component of the intersection of two straight lines}

If the frame components of two points on a plane are known, then the frame component of any point on the line connecting these two points can be obtained. If the frame components of two other points on the plane are known, then these two points can also be connected into a straight line. The current question is how to determine the frame component of the intersection point of these two lines? The following theorem I have provided answers to this question

\begin{theorem}{Frame component at the intersection of two straight lines-simplified form, Daiyuan Zhang}{LiangtiaoZhixianJiaodianDeBiaojiafenliang_Jianyuexingshi}\label{LiangtiaoZhixianJiaodianDeBiaojiafenliang_Jianyuexingshi}
	Given a $\triangle ABC$, assuming the frame component of IC ${{P}_{1}}$ is $\alpha _{A}^{{{P}_{1}}}$, $\alpha _{B}^{{{P}_{1}}}$, $\alpha _{C}^{{{P}_{1}}}$; the frame component of IC ${{P}_{2}}$ is $\alpha _{A}^{{{P}_{2}}}$, $\alpha _{B}^{{{P}_{2}}}$, $\alpha _{C}^{{{P}_{2}}}$; the frame component of IC ${{P}_{3}}$ is $\alpha _{A}^{{{P}_{3}}}$, $\alpha _{B}^{{{P}_{3}}}$, $\alpha _{C}^{{{P}_{3}}}$; the frame component of IC ${{P}_{4}}$ is $\alpha _{A}^{{{P}_{4}}}$, $\alpha _{B}^{{{P}_{4}}}$, $\alpha _{C}^{{{P}_{4}}}$. The intersection point of two straight lines $\overleftrightarrow{{{P}_{1}}{{P}_{2}}}$ and $\overleftrightarrow{{{P}_{3}}{{P}_{4}}}$ is $P=\overleftrightarrow{{{P}_{1}}{{P}_{2}}}\cap \overleftrightarrow{{{P}_{3}}{{P}_{4}}}$, and $\overrightarrow{{{P}_{1}}P}={{\kappa }_{{{P}_{1}}{{P}_{2}}}}\overrightarrow{{{P}_{1}}{{P}_{2}}}$, $\overrightarrow{{{P}_{1}}{{P}_{2}}}\ne \overrightarrow{0}$, $\overrightarrow{{{P}_{3}}P}={{\kappa }_{{{P}_{3}}{{P}_{4}}}}\overrightarrow{{{P}_{3}}{{P}_{4}}}$, $\overrightarrow{{{P}_{3}}{{P}_{4}}}\ne \overrightarrow{0}$, then the frame component of IC $P$ is:
	\begin{align*}
		\alpha _{A}^{P}& =\frac{\left( \alpha _{B}^{{{P}_{3}}}-\alpha _{B}^{{{P}_{2}}} \right)\left( \alpha _{C}^{{{P}_{3}}}-\alpha _{C}^{{{P}_{4}}} \right)-\left( \alpha _{B}^{{{P}_{3}}}-\alpha _{B}^{{{P}_{4}}} \right)\left( \alpha _{C}^{{{P}_{3}}}-\alpha _{C}^{{{P}_{2}}} \right)}{\left( \alpha _{B}^{{{P}_{1}}}-\alpha _{B}^{{{P}_{2}}} \right)\left( \alpha _{C}^{{{P}_{3}}}-\alpha _{C}^{{{P}_{4}}} \right)-\left( \alpha _{B}^{{{P}_{3}}}-\alpha _{B}^{{{P}_{4}}} \right)\left( \alpha _{C}^{{{P}_{1}}}-\alpha _{C}^{{{P}_{2}}} \right)}\alpha _{A}^{{{P}_{1}}} \\ 
		& +\frac{\left( \alpha _{B}^{{{P}_{1}}}-\alpha _{B}^{{{P}_{3}}} \right)\left( \alpha _{C}^{{{P}_{3}}}-\alpha _{C}^{{{P}_{4}}} \right)-\left( \alpha _{B}^{{{P}_{3}}}-\alpha _{B}^{{{P}_{4}}} \right)\left( \alpha _{C}^{{{P}_{1}}}-\alpha _{C}^{{{P}_{3}}} \right)}{\left( \alpha _{B}^{{{P}_{1}}}-\alpha _{B}^{{{P}_{2}}} \right)\left( \alpha _{C}^{{{P}_{3}}}-\alpha _{C}^{{{P}_{4}}} \right)-\left( \alpha _{B}^{{{P}_{3}}}-\alpha _{B}^{{{P}_{4}}} \right)\left( \alpha _{C}^{{{P}_{1}}}-\alpha _{C}^{{{P}_{2}}} \right)}\alpha _{A}^{{{P}_{2}}},  
	\end{align*}
	
	\begin{align*}
		\alpha _{B}^{P}& =\frac{\left( \alpha _{C}^{{{P}_{1}}}-\alpha _{C}^{{{P}_{3}}} \right)\left( \alpha _{A}^{{{P}_{3}}}-\alpha _{A}^{{{P}_{4}}} \right)-\left( \alpha _{C}^{{{P}_{3}}}-\alpha _{C}^{{{P}_{4}}} \right)\left( \alpha _{A}^{{{P}_{1}}}-\alpha _{A}^{{{P}_{3}}} \right)}{\left( \alpha _{C}^{{{P}_{1}}}-\alpha _{C}^{{{P}_{2}}} \right)\left( \alpha _{A}^{{{P}_{3}}}-\alpha _{A}^{{{P}_{4}}} \right)-\left( \alpha _{C}^{{{P}_{3}}}-\alpha _{C}^{{{P}_{4}}} \right)\left( \alpha _{A}^{{{P}_{1}}}-\alpha _{A}^{{{P}_{2}}} \right)}\alpha _{B}^{{{P}_{1}}} \\ 
		& +\frac{\left( \alpha _{A}^{{{P}_{1}}}-\alpha _{A}^{{{P}_{3}}} \right)\left( \alpha _{B}^{{{P}_{3}}}-\alpha _{B}^{{{P}_{4}}} \right)-\left( \alpha _{A}^{{{P}_{3}}}-\alpha _{A}^{{{P}_{4}}} \right)\left( \alpha _{B}^{{{P}_{1}}}-\alpha _{B}^{{{P}_{3}}} \right)}{\left( \alpha _{A}^{{{P}_{1}}}-\alpha _{A}^{{{P}_{2}}} \right)\left( \alpha _{B}^{{{P}_{3}}}-\alpha _{B}^{{{P}_{4}}} \right)-\left( \alpha _{A}^{{{P}_{3}}}-\alpha _{A}^{{{P}_{4}}} \right)\left( \alpha _{B}^{{{P}_{1}}}-\alpha _{B}^{{{P}_{2}}} \right)}\alpha _{B}^{{{P}_{2}}},  
	\end{align*}
	
	\begin{align*}
		\alpha _{C}^{P}& =\frac{\left( \alpha _{A}^{{{P}_{3}}}-\alpha _{A}^{{{P}_{2}}} \right)\left( \alpha _{B}^{{{P}_{3}}}-\alpha _{B}^{{{P}_{4}}} \right)-\left( \alpha _{A}^{{{P}_{3}}}-\alpha _{A}^{{{P}_{4}}} \right)\left( \alpha _{B}^{{{P}_{3}}}-\alpha _{B}^{{{P}_{2}}} \right)}{\left( \alpha _{A}^{{{P}_{1}}}-\alpha _{A}^{{{P}_{2}}} \right)\left( \alpha _{B}^{{{P}_{3}}}-\alpha _{B}^{{{P}_{4}}} \right)-\left( \alpha _{A}^{{{P}_{3}}}-\alpha _{A}^{{{P}_{4}}} \right)\left( \alpha _{B}^{{{P}_{1}}}-\alpha _{B}^{{{P}_{2}}} \right)}\alpha _{C}^{{{P}_{1}}} \\ 
		& +\frac{\left( \alpha _{A}^{{{P}_{1}}}-\alpha _{A}^{{{P}_{3}}} \right)\left( \alpha _{B}^{{{P}_{3}}}-\alpha _{B}^{{{P}_{4}}} \right)-\left( \alpha _{A}^{{{P}_{3}}}-\alpha _{A}^{{{P}_{4}}} \right)\left( \alpha _{B}^{{{P}_{1}}}-\alpha _{B}^{{{P}_{3}}} \right)}{\left( \alpha _{A}^{{{P}_{1}}}-\alpha _{A}^{{{P}_{2}}} \right)\left( \alpha _{B}^{{{P}_{3}}}-\alpha _{B}^{{{P}_{4}}} \right)-\left( \alpha _{A}^{{{P}_{3}}}-\alpha _{A}^{{{P}_{4}}} \right)\left( \alpha _{B}^{{{P}_{1}}}-\alpha _{B}^{{{P}_{2}}} \right)}\alpha _{C}^{{{P}_{2}}}.  
	\end{align*}
\end{theorem}

\begin{proof}
	Because the intersection point of two straight lines $\overleftrightarrow{{{P}_{1}}{{P}_{2}}}$ and $\overleftrightarrow{{{P}_{3}}{{P}_{4}}}$ is $P=\overleftrightarrow{{{P}_{1}}{{P}_{2}}}\cap \overleftrightarrow{{{P}_{3}}{{P}_{4}}}$, According to theorem \ref{thm:HuixinLianxianShangRenyiYidianDeBiaojiafenliang}, the following two sets of equations are obtained:
	
	The first set of equations is:
	\[\alpha _{A}^{P}=\left( 1-{{\kappa }_{{{P}_{1}}{{P}_{2}}}} \right)\alpha _{A}^{{{P}_{1}}}+{{\kappa }_{{{P}_{1}}{{P}_{2}}}}\alpha _{A}^{{{P}_{2}}},\]
	\[\alpha _{B}^{P}=\left( 1-{{\kappa }_{{{P}_{1}}{{P}_{2}}}} \right)\alpha _{B}^{{{P}_{1}}}+{{\kappa }_{{{P}_{1}}{{P}_{2}}}}\alpha _{B}^{{{P}_{2}}},\]
	\[\alpha _{C}^{P}=\left( 1-{{\kappa }_{{{P}_{1}}{{P}_{2}}}} \right)\alpha _{C}^{{{P}_{1}}}+{{\kappa }_{{{P}_{1}}{{P}_{2}}}}\alpha _{C}^{{{P}_{2}}}.\]
	
	The second set of equations is:
	\[\alpha _{A}^{P}=\left( 1-{{\kappa }_{{{P}_{3}}{{P}_{4}}}} \right)\alpha _{A}^{{{P}_{3}}}+{{\kappa }_{{{P}_{3}}{{P}_{4}}}}\alpha _{A}^{{{P}_{4}}},\]
	\[\alpha _{B}^{P}=\left( 1-{{\kappa }_{{{P}_{3}}{{P}_{4}}}} \right)\alpha _{B}^{{{P}_{3}}}+{{\kappa }_{{{P}_{3}}{{P}_{4}}}}\alpha _{B}^{{{P}_{4}}},\]
	\[\alpha _{C}^{P}=\left( 1-{{\kappa }_{{{P}_{3}}{{P}_{4}}}} \right)\alpha _{C}^{{{P}_{3}}}+{{\kappa }_{{{P}_{3}}{{P}_{4}}}}\alpha _{C}^{{{P}_{4}}}.\]
	
	So, the first and second equations of the above two sets (containing $\alpha _{A}^{P}$ and $\alpha _{B}^{P}$) are respectively used to obtain the following system of equations (referred to as the $AB$ combination):
	\[\left( 1-{{\kappa }_{{{P}_{1}}{{P}_{2}}}} \right)\alpha _{A}^{{{P}_{1}}}+{{\kappa }_{{{P}_{1}}{{P}_{2}}}}\alpha _{A}^{{{P}_{2}}}=\left( 1-{{\kappa }_{{{P}_{3}}{{P}_{4}}}} \right)\alpha _{A}^{{{P}_{3}}}+{{\kappa }_{{{P}_{3}}{{P}_{4}}}}\alpha _{A}^{{{P}_{4}}},\]
	\[\left( 1-{{\kappa }_{{{P}_{1}}{{P}_{2}}}} \right)\alpha _{B}^{{{P}_{1}}}+{{\kappa }_{{{P}_{1}}{{P}_{2}}}}\alpha _{B}^{{{P}_{2}}}=\left( 1-{{\kappa }_{{{P}_{3}}{{P}_{4}}}} \right)\alpha _{B}^{{{P}_{3}}}+{{\kappa }_{{{P}_{3}}{{P}_{4}}}}\alpha _{B}^{{{P}_{4}}}.\]
	
	The above equation system is a linear equation system with variables ${{\kappa }_{{{P}_{1}}{{P}_{2}}}}$ and ${{\kappa }_{{{P}_{3}}{{P}_{4}}}}$.
	
	From the second and third equations of the above two sets (containing $\alpha _{B}^{P}$ and $\alpha _{C}^{P}$), the following equation system (referred to as the $BC$ combination) is obtained:
	\[\left( 1-{{\kappa }_{{{P}_{1}}{{P}_{2}}}} \right)\alpha _{B}^{{{P}_{1}}}+{{\kappa }_{{{P}_{1}}{{P}_{2}}}}\alpha _{B}^{{{P}_{2}}}=\left( 1-{{\kappa }_{{{P}_{3}}{{P}_{4}}}} \right)\alpha _{B}^{{{P}_{3}}}+{{\kappa }_{{{P}_{3}}{{P}_{4}}}}\alpha _{B}^{{{P}_{4}}}.\]
	\[\left( 1-{{\kappa }_{{{P}_{1}}{{P}_{2}}}} \right)\alpha _{C}^{{{P}_{1}}}+{{\kappa }_{{{P}_{1}}{{P}_{2}}}}\alpha _{C}^{{{P}_{2}}}=\left( 1-{{\kappa }_{{{P}_{3}}{{P}_{4}}}} \right)\alpha _{C}^{{{P}_{3}}}+{{\kappa }_{{{P}_{3}}{{P}_{4}}}}\alpha _{C}^{{{P}_{4}}},\]
	
	From the third and first equations of the above two sets (containing $\alpha _{C}^{P}$ and $\alpha _{A}^{P}$), the following equation system (referred to as the $CA$ combination) is obtained:
	\[\left( 1-{{\kappa }_{{{P}_{1}}{{P}_{2}}}} \right)\alpha _{C}^{{{P}_{1}}}+{{\kappa }_{{{P}_{1}}{{P}_{2}}}}\alpha _{C}^{{{P}_{2}}}=\left( 1-{{\kappa }_{{{P}_{3}}{{P}_{4}}}} \right)\alpha _{C}^{{{P}_{3}}}+{{\kappa }_{{{P}_{3}}{{P}_{4}}}}\alpha _{C}^{{{P}_{4}}},\]
	\[\left( 1-{{\kappa }_{{{P}_{1}}{{P}_{2}}}} \right)\alpha _{A}^{{{P}_{1}}}+{{\kappa }_{{{P}_{1}}{{P}_{2}}}}\alpha _{A}^{{{P}_{2}}}=\left( 1-{{\kappa }_{{{P}_{3}}{{P}_{4}}}} \right)\alpha _{A}^{{{P}_{3}}}+{{\kappa }_{{{P}_{3}}{{P}_{4}}}}\alpha _{A}^{{{P}_{4}}}.\]
	
	It should be noted that due to the uniqueness of the frame components $\alpha _{A}^{P}$, $\alpha _{B}^{P}$ and $\alpha _{C}^{P}$, as well as the uniqueness of ${{\kappa }_{{{P}_{1}}{{P}_{2}}}}$ and ${{\kappa }_{{{P}_{3}}{{P}_{4}}}}$, it can be concluded that regardless of the combination used, the solution ${{\kappa }_{{{P}_{1}}{{P}_{2}}}}$ and ${{\kappa }_{{{P}_{3}}{{P}_{4}}}}$ are all the same.
	
	It is easy to see that the system of equations can be obtained through rotation, that is, through rotation of $A$, $B$, $C$, the combination of $AB$ equations becomes the combination of $BC$ equations, the combination of $BC$ equations becomes the combination of $CA$ equations, and the combination of $CA$ equations becomes the combination of $AB$ equations. Therefore, its solution can also be obtained through rotation of $A$, $B$, and $C$. So all you need to do is to find the solution to one system of equations, and then rotate through $A$, $B$, and $C$ to obtain solutions to other systems of equations. Due to the uniqueness of the frame components and the uniqueness of ${{\kappa }_{{{P}_{1}}{{P}_{2}}}}$ and ${{\kappa }_{{{P}_{3}}{{P}_{4}}}}$, each system of equations has the same solution, but with different expressions.
		
	Below, the first and second equations are combined to obtain a system of equations (referred to as the $AB$ combination):
	\[\left( 1-x \right)\alpha _{A}^{{{P}_{1}}}+x\alpha _{A}^{{{P}_{2}}}=\left( 1-y \right)\alpha _{A}^{{{P}_{3}}}+y\alpha _{A}^{{{P}_{4}}},\]
	\[\left( 1-x \right)\alpha _{B}^{{{P}_{1}}}+x\alpha _{B}^{{{P}_{2}}}=\left( 1-y \right)\alpha _{B}^{{{P}_{3}}}+y\alpha _{B}^{{{P}_{4}}}.\]
	
	i.e.
	\[\left\{ \begin{aligned}
		& \left( \alpha _{A}^{{{P}_{2}}}-\alpha _{A}^{{{P}_{1}}} \right)x-\left( \alpha _{A}^{{{P}_{4}}}-\alpha _{A}^{{{P}_{3}}} \right)y=\alpha _{A}^{{{P}_{3}}}-\alpha _{A}^{{{P}_{1}}} \\ 
		& \left( \alpha _{B}^{{{P}_{2}}}-\alpha _{B}^{{{P}_{1}}} \right)x-\left( \alpha _{B}^{{{P}_{4}}}-\alpha _{B}^{{{P}_{3}}} \right)y=\alpha _{B}^{{{P}_{3}}}-\alpha _{B}^{{{P}_{1}}}. \\ 
	\end{aligned} \right.\]
	
	Solving the above equation system yields:
	\[x=\frac{\left| \begin{matrix}
			\alpha _{A}^{{{P}_{3}}}-\alpha _{A}^{{{P}_{1}}} & -\left( \alpha _{A}^{{{P}_{4}}}-\alpha _{A}^{{{P}_{3}}} \right)  \\
			\alpha _{B}^{{{P}_{3}}}-\alpha _{B}^{{{P}_{1}}} & -\left( \alpha _{B}^{{{P}_{4}}}-\alpha _{B}^{{{P}_{3}}} \right)  \\
		\end{matrix} \right|}{\left| \begin{matrix}
			\alpha _{A}^{{{P}_{2}}}-\alpha _{A}^{{{P}_{1}}} & -\left( \alpha _{A}^{{{P}_{4}}}-\alpha _{A}^{{{P}_{3}}} \right)  \\
			\alpha _{B}^{{{P}_{2}}}-\alpha _{B}^{{{P}_{1}}} & -\left( \alpha _{B}^{{{P}_{4}}}-\alpha _{B}^{{{P}_{3}}} \right)  \\
		\end{matrix} \right|}=\frac{\left( \alpha _{A}^{{{P}_{4}}}-\alpha _{A}^{{{P}_{3}}} \right)\left( \alpha _{B}^{{{P}_{3}}}-\alpha _{B}^{{{P}_{1}}} \right)-\left( \alpha _{A}^{{{P}_{3}}}-\alpha _{A}^{{{P}_{1}}} \right)\left( \alpha _{B}^{{{P}_{4}}}-\alpha _{B}^{{{P}_{3}}} \right)}{\left( \alpha _{A}^{{{P}_{4}}}-\alpha _{A}^{{{P}_{3}}} \right)\left( \alpha _{B}^{{{P}_{2}}}-\alpha _{B}^{{{P}_{1}}} \right)-\left( \alpha _{A}^{{{P}_{2}}}-\alpha _{A}^{{{P}_{1}}} \right)\left( \alpha _{B}^{{{P}_{4}}}-\alpha _{B}^{{{P}_{3}}} \right)},\]
	\[y=\frac{\left| \begin{matrix}
			\alpha _{A}^{{{P}_{2}}}-\alpha _{A}^{{{P}_{1}}} & \alpha _{A}^{{{P}_{3}}}-\alpha _{A}^{{{P}_{1}}}  \\
			\alpha _{B}^{{{P}_{2}}}-\alpha _{B}^{{{P}_{1}}} & \alpha _{B}^{{{P}_{3}}}-\alpha _{B}^{{{P}_{1}}}  \\
		\end{matrix} \right|}{\left| \begin{matrix}
			\alpha _{A}^{{{P}_{2}}}-\alpha _{A}^{{{P}_{1}}} & -\left( \alpha _{A}^{{{P}_{4}}}-\alpha _{A}^{{{P}_{3}}} \right)  \\
			\alpha _{B}^{{{P}_{2}}}-\alpha _{B}^{{{P}_{1}}} & -\left( \alpha _{B}^{{{P}_{4}}}-\alpha _{B}^{{{P}_{3}}} \right)  \\
		\end{matrix} \right|}=\frac{\left( \alpha _{A}^{{{P}_{2}}}-\alpha _{A}^{{{P}_{1}}} \right)\left( \alpha _{B}^{{{P}_{3}}}-\alpha _{B}^{{{P}_{1}}} \right)-\left( \alpha _{A}^{{{P}_{3}}}-\alpha _{A}^{{{P}_{1}}} \right)\left( \alpha _{B}^{{{P}_{2}}}-\alpha _{B}^{{{P}_{1}}} \right)}{\left( \alpha _{A}^{{{P}_{4}}}-\alpha _{A}^{{{P}_{3}}} \right)\left( \alpha _{B}^{{{P}_{2}}}-\alpha _{B}^{{{P}_{1}}} \right)-\left( \alpha _{A}^{{{P}_{2}}}-\alpha _{A}^{{{P}_{1}}} \right)\left( \alpha _{B}^{{{P}_{4}}}-\alpha _{B}^{{{P}_{3}}} \right)}.\]
	
	Therefore
	\[{{\kappa }_{{{P}_{1}}{{P}_{2}}}}=x=\frac{\left( \alpha _{A}^{{{P}_{1}}}-\alpha _{A}^{{{P}_{3}}} \right)\left( \alpha _{B}^{{{P}_{3}}}-\alpha _{B}^{{{P}_{4}}} \right)-\left( \alpha _{A}^{{{P}_{3}}}-\alpha _{A}^{{{P}_{4}}} \right)\left( \alpha _{B}^{{{P}_{1}}}-\alpha _{B}^{{{P}_{3}}} \right)}{\left( \alpha _{A}^{{{P}_{1}}}-\alpha _{A}^{{{P}_{2}}} \right)\left( \alpha _{B}^{{{P}_{3}}}-\alpha _{B}^{{{P}_{4}}} \right)-\left( \alpha _{A}^{{{P}_{3}}}-\alpha _{A}^{{{P}_{4}}} \right)\left( \alpha _{B}^{{{P}_{1}}}-\alpha _{B}^{{{P}_{2}}} \right)},\]
	
	\begin{align*}
		1-{{\kappa }_{{{P}_{1}}{{P}_{2}}}}& =1-\frac{\left( \alpha _{A}^{{{P}_{1}}}-\alpha _{A}^{{{P}_{3}}} \right)\left( \alpha _{B}^{{{P}_{3}}}-\alpha _{B}^{{{P}_{4}}} \right)-\left( \alpha _{A}^{{{P}_{3}}}-\alpha _{A}^{{{P}_{4}}} \right)\left( \alpha _{B}^{{{P}_{1}}}-\alpha _{B}^{{{P}_{3}}} \right)}{\left( \alpha _{A}^{{{P}_{1}}}-\alpha _{A}^{{{P}_{2}}} \right)\left( \alpha _{B}^{{{P}_{3}}}-\alpha _{B}^{{{P}_{4}}} \right)-\left( \alpha _{A}^{{{P}_{3}}}-\alpha _{A}^{{{P}_{4}}} \right)\left( \alpha _{B}^{{{P}_{1}}}-\alpha _{B}^{{{P}_{2}}} \right)} \\ 
		& =\frac{\left( \alpha _{A}^{{{P}_{3}}}-\alpha _{A}^{{{P}_{2}}} \right)\left( \alpha _{B}^{{{P}_{3}}}-\alpha _{B}^{{{P}_{4}}} \right)-\left( \alpha _{A}^{{{P}_{3}}}-\alpha _{A}^{{{P}_{4}}} \right)\left( \alpha _{B}^{{{P}_{3}}}-\alpha _{B}^{{{P}_{2}}} \right)}{\left( \alpha _{A}^{{{P}_{1}}}-\alpha _{A}^{{{P}_{2}}} \right)\left( \alpha _{B}^{{{P}_{3}}}-\alpha _{B}^{{{P}_{4}}} \right)-\left( \alpha _{A}^{{{P}_{3}}}-\alpha _{A}^{{{P}_{4}}} \right)\left( \alpha _{B}^{{{P}_{1}}}-\alpha _{B}^{{{P}_{2}}} \right)},  
	\end{align*}
	and 
	\[{{\kappa }_{{{P}_{3}}{{P}_{4}}}}=y=\frac{\left( \alpha _{A}^{{{P}_{1}}}-\alpha _{A}^{{{P}_{3}}} \right)\left( \alpha _{B}^{{{P}_{1}}}-\alpha _{B}^{{{P}_{2}}} \right)-\left( \alpha _{A}^{{{P}_{1}}}-\alpha _{A}^{{{P}_{2}}} \right)\left( \alpha _{B}^{{{P}_{1}}}-\alpha _{B}^{{{P}_{3}}} \right)}{\left( \alpha _{A}^{{{P}_{1}}}-\alpha _{A}^{{{P}_{2}}} \right)\left( \alpha _{B}^{{{P}_{3}}}-\alpha _{B}^{{{P}_{4}}} \right)-\left( \alpha _{A}^{{{P}_{3}}}-\alpha _{A}^{{{P}_{4}}} \right)\left( \alpha _{B}^{{{P}_{1}}}-\alpha _{B}^{{{P}_{2}}} \right)}.\]
	
	\begin{align*}
		1-{{\kappa }_{{{P}_{3}}{{P}_{4}}}}& =1-\frac{\left( \alpha _{A}^{{{P}_{1}}}-\alpha _{A}^{{{P}_{3}}} \right)\left( \alpha _{B}^{{{P}_{1}}}-\alpha _{B}^{{{P}_{2}}} \right)-\left( \alpha _{A}^{{{P}_{1}}}-\alpha _{A}^{{{P}_{2}}} \right)\left( \alpha _{B}^{{{P}_{1}}}-\alpha _{B}^{{{P}_{3}}} \right)}{\left( \alpha _{A}^{{{P}_{1}}}-\alpha _{A}^{{{P}_{2}}} \right)\left( \alpha _{B}^{{{P}_{3}}}-\alpha _{B}^{{{P}_{4}}} \right)-\left( \alpha _{A}^{{{P}_{3}}}-\alpha _{A}^{{{P}_{4}}} \right)\left( \alpha _{B}^{{{P}_{1}}}-\alpha _{B}^{{{P}_{2}}} \right)} \\ 
		& =\frac{\left( \alpha _{A}^{{{P}_{4}}}-\alpha _{A}^{{{P}_{1}}} \right)\left( \alpha _{B}^{{{P}_{1}}}-\alpha _{B}^{{{P}_{2}}} \right)-\left( \alpha _{A}^{{{P}_{1}}}-\alpha _{A}^{{{P}_{2}}} \right)\left( \alpha _{B}^{{{P}_{4}}}-\alpha _{B}^{{{P}_{1}}} \right)}{\left( \alpha _{A}^{{{P}_{1}}}-\alpha _{A}^{{{P}_{2}}} \right)\left( \alpha _{B}^{{{P}_{3}}}-\alpha _{B}^{{{P}_{4}}} \right)-\left( \alpha _{A}^{{{P}_{3}}}-\alpha _{A}^{{{P}_{4}}} \right)\left( \alpha _{B}^{{{P}_{1}}}-\alpha _{B}^{{{P}_{2}}} \right)}.  
	\end{align*}
	
	After rotation, it is obtained (referred to as the $BC$ combination):
	\[{{\kappa }_{{{P}_{1}}{{P}_{2}}}}=\frac{\left( \alpha _{B}^{{{P}_{1}}}-\alpha _{B}^{{{P}_{3}}} \right)\left( \alpha _{C}^{{{P}_{3}}}-\alpha _{C}^{{{P}_{4}}} \right)-\left( \alpha _{B}^{{{P}_{3}}}-\alpha _{B}^{{{P}_{4}}} \right)\left( \alpha _{C}^{{{P}_{1}}}-\alpha _{C}^{{{P}_{3}}} \right)}{\left( \alpha _{B}^{{{P}_{1}}}-\alpha _{B}^{{{P}_{2}}} \right)\left( \alpha _{C}^{{{P}_{3}}}-\alpha _{C}^{{{P}_{4}}} \right)-\left( \alpha _{B}^{{{P}_{3}}}-\alpha _{B}^{{{P}_{4}}} \right)\left( \alpha _{C}^{{{P}_{1}}}-\alpha _{C}^{{{P}_{2}}} \right)},\]
	\[1-{{\kappa }_{{{P}_{1}}{{P}_{2}}}}=\frac{\left( \alpha _{B}^{{{P}_{3}}}-\alpha _{B}^{{{P}_{2}}} \right)\left( \alpha _{C}^{{{P}_{3}}}-\alpha _{C}^{{{P}_{4}}} \right)-\left( \alpha _{B}^{{{P}_{3}}}-\alpha _{B}^{{{P}_{4}}} \right)\left( \alpha _{C}^{{{P}_{3}}}-\alpha _{C}^{{{P}_{2}}} \right)}{\left( \alpha _{B}^{{{P}_{1}}}-\alpha _{B}^{{{P}_{2}}} \right)\left( \alpha _{C}^{{{P}_{3}}}-\alpha _{C}^{{{P}_{4}}} \right)-\left( \alpha _{B}^{{{P}_{3}}}-\alpha _{B}^{{{P}_{4}}} \right)\left( \alpha _{C}^{{{P}_{1}}}-\alpha _{C}^{{{P}_{2}}} \right)},\]
	\[{{\kappa }_{{{P}_{3}}{{P}_{4}}}}=\frac{\left( \alpha _{B}^{{{P}_{1}}}-\alpha _{B}^{{{P}_{3}}} \right)\left( \alpha _{C}^{{{P}_{1}}}-\alpha _{C}^{{{P}_{2}}} \right)-\left( \alpha _{B}^{{{P}_{1}}}-\alpha _{B}^{{{P}_{2}}} \right)\left( \alpha _{C}^{{{P}_{1}}}-\alpha _{C}^{{{P}_{3}}} \right)}{\left( \alpha _{B}^{{{P}_{1}}}-\alpha _{B}^{{{P}_{2}}} \right)\left( \alpha _{C}^{{{P}_{3}}}-\alpha _{C}^{{{P}_{4}}} \right)-\left( \alpha _{B}^{{{P}_{3}}}-\alpha _{B}^{{{P}_{4}}} \right)\left( \alpha _{C}^{{{P}_{1}}}-\alpha _{C}^{{{P}_{2}}} \right)}.\]
	\[1-{{\kappa }_{{{P}_{3}}{{P}_{4}}}}=\frac{\left( \alpha _{B}^{{{P}_{4}}}-\alpha _{B}^{{{P}_{1}}} \right)\left( \alpha _{C}^{{{P}_{1}}}-\alpha _{C}^{{{P}_{2}}} \right)-\left( \alpha _{B}^{{{P}_{1}}}-\alpha _{B}^{{{P}_{2}}} \right)\left( \alpha _{C}^{{{P}_{4}}}-\alpha _{C}^{{{P}_{1}}} \right)}{\left( \alpha _{B}^{{{P}_{1}}}-\alpha _{B}^{{{P}_{2}}} \right)\left( \alpha _{C}^{{{P}_{3}}}-\alpha _{C}^{{{P}_{4}}} \right)-\left( \alpha _{B}^{{{P}_{3}}}-\alpha _{B}^{{{P}_{4}}} \right)\left( \alpha _{C}^{{{P}_{1}}}-\alpha _{C}^{{{P}_{2}}} \right)}.\]
	
	After rotation, it is obtained (referred to as the $CA$ combination):
	\[{{\kappa }_{{{P}_{1}}{{P}_{2}}}}=\frac{\left( \alpha _{C}^{{{P}_{1}}}-\alpha _{C}^{{{P}_{3}}} \right)\left( \alpha _{A}^{{{P}_{3}}}-\alpha _{A}^{{{P}_{4}}} \right)-\left( \alpha _{C}^{{{P}_{3}}}-\alpha _{C}^{{{P}_{4}}} \right)\left( \alpha _{A}^{{{P}_{1}}}-\alpha _{A}^{{{P}_{3}}} \right)}{\left( \alpha _{C}^{{{P}_{1}}}-\alpha _{C}^{{{P}_{2}}} \right)\left( \alpha _{A}^{{{P}_{3}}}-\alpha _{A}^{{{P}_{4}}} \right)-\left( \alpha _{C}^{{{P}_{3}}}-\alpha _{C}^{{{P}_{4}}} \right)\left( \alpha _{A}^{{{P}_{1}}}-\alpha _{A}^{{{P}_{2}}} \right)},\]
	\[1-{{\kappa }_{{{P}_{1}}{{P}_{2}}}}=\frac{\left( \alpha _{C}^{{{P}_{3}}}-\alpha _{C}^{{{P}_{2}}} \right)\left( \alpha _{A}^{{{P}_{3}}}-\alpha _{A}^{{{P}_{4}}} \right)-\left( \alpha _{C}^{{{P}_{3}}}-\alpha _{C}^{{{P}_{4}}} \right)\left( \alpha _{A}^{{{P}_{3}}}-\alpha _{A}^{{{P}_{2}}} \right)}{\left( \alpha _{C}^{{{P}_{1}}}-\alpha _{C}^{{{P}_{2}}} \right)\left( \alpha _{A}^{{{P}_{3}}}-\alpha _{A}^{{{P}_{4}}} \right)-\left( \alpha _{C}^{{{P}_{3}}}-\alpha _{C}^{{{P}_{4}}} \right)\left( \alpha _{A}^{{{P}_{1}}}-\alpha _{A}^{{{P}_{2}}} \right)},\]
	\[{{\kappa }_{{{P}_{3}}{{P}_{4}}}}=\frac{\left( \alpha _{C}^{{{P}_{1}}}-\alpha _{C}^{{{P}_{3}}} \right)\left( \alpha _{A}^{{{P}_{1}}}-\alpha _{A}^{{{P}_{2}}} \right)-\left( \alpha _{C}^{{{P}_{1}}}-\alpha _{C}^{{{P}_{2}}} \right)\left( \alpha _{A}^{{{P}_{1}}}-\alpha _{A}^{{{P}_{3}}} \right)}{\left( \alpha _{C}^{{{P}_{1}}}-\alpha _{C}^{{{P}_{2}}} \right)\left( \alpha _{A}^{{{P}_{3}}}-\alpha _{A}^{{{P}_{4}}} \right)-\left( \alpha _{C}^{{{P}_{3}}}-\alpha _{C}^{{{P}_{4}}} \right)\left( \alpha _{A}^{{{P}_{1}}}-\alpha _{A}^{{{P}_{2}}} \right)}.\]
	\[1-{{\kappa }_{{{P}_{3}}{{P}_{4}}}}=\frac{\left( \alpha _{C}^{{{P}_{4}}}-\alpha _{C}^{{{P}_{1}}} \right)\left( \alpha _{A}^{{{P}_{1}}}-\alpha _{A}^{{{P}_{2}}} \right)-\left( \alpha _{C}^{{{P}_{1}}}-\alpha _{C}^{{{P}_{2}}} \right)\left( \alpha _{A}^{{{P}_{4}}}-\alpha _{A}^{{{P}_{1}}} \right)}{\left( \alpha _{C}^{{{P}_{1}}}-\alpha _{C}^{{{P}_{2}}} \right)\left( \alpha _{A}^{{{P}_{3}}}-\alpha _{A}^{{{P}_{4}}} \right)-\left( \alpha _{C}^{{{P}_{3}}}-\alpha _{C}^{{{P}_{4}}} \right)\left( \alpha _{A}^{{{P}_{1}}}-\alpha _{A}^{{{P}_{2}}} \right)}.\]
	
	To obtain a symmetrical and aesthetically pleasing form, the following combinations are used to solve for the frame components:
	\begin{align*}
		\alpha _{A}^{P}& =\left( 1-{{\kappa }_{{{P}_{1}}{{P}_{2}}}} \right)\alpha _{A}^{{{P}_{1}}}+{{\kappa }_{{{P}_{1}}{{P}_{2}}}}\alpha _{A}^{{{P}_{2}}} \\ 
		& =\frac{\left( \alpha _{B}^{{{P}_{3}}}-\alpha _{B}^{{{P}_{2}}} \right)\left( \alpha _{C}^{{{P}_{3}}}-\alpha _{C}^{{{P}_{4}}} \right)-\left( \alpha _{B}^{{{P}_{3}}}-\alpha _{B}^{{{P}_{4}}} \right)\left( \alpha _{C}^{{{P}_{3}}}-\alpha _{C}^{{{P}_{2}}} \right)}{\left( \alpha _{B}^{{{P}_{1}}}-\alpha _{B}^{{{P}_{2}}} \right)\left( \alpha _{C}^{{{P}_{3}}}-\alpha _{C}^{{{P}_{4}}} \right)-\left( \alpha _{B}^{{{P}_{3}}}-\alpha _{B}^{{{P}_{4}}} \right)\left( \alpha _{C}^{{{P}_{1}}}-\alpha _{C}^{{{P}_{2}}} \right)}\alpha _{A}^{{{P}_{1}}} \\ 
		& +\frac{\left( \alpha _{B}^{{{P}_{1}}}-\alpha _{B}^{{{P}_{3}}} \right)\left( \alpha _{C}^{{{P}_{3}}}-\alpha _{C}^{{{P}_{4}}} \right)-\left( \alpha _{B}^{{{P}_{3}}}-\alpha _{B}^{{{P}_{4}}} \right)\left( \alpha _{C}^{{{P}_{1}}}-\alpha _{C}^{{{P}_{3}}} \right)}{\left( \alpha _{B}^{{{P}_{1}}}-\alpha _{B}^{{{P}_{2}}} \right)\left( \alpha _{C}^{{{P}_{3}}}-\alpha _{C}^{{{P}_{4}}} \right)-\left( \alpha _{B}^{{{P}_{3}}}-\alpha _{B}^{{{P}_{4}}} \right)\left( \alpha _{C}^{{{P}_{1}}}-\alpha _{C}^{{{P}_{2}}} \right)}\alpha _{A}^{{{P}_{2}}},  
	\end{align*}
	
	\begin{align*}
		\alpha _{B}^{P}& =\left( 1-{{\kappa }_{{{P}_{1}}{{P}_{2}}}} \right)\alpha _{B}^{{{P}_{1}}}+{{\kappa }_{{{P}_{1}}{{P}_{2}}}}\alpha _{B}^{{{P}_{2}}} \\ 
		& =\frac{\left( \alpha _{C}^{{{P}_{1}}}-\alpha _{C}^{{{P}_{3}}} \right)\left( \alpha _{A}^{{{P}_{3}}}-\alpha _{A}^{{{P}_{4}}} \right)-\left( \alpha _{C}^{{{P}_{3}}}-\alpha _{C}^{{{P}_{4}}} \right)\left( \alpha _{A}^{{{P}_{1}}}-\alpha _{A}^{{{P}_{3}}} \right)}{\left( \alpha _{C}^{{{P}_{1}}}-\alpha _{C}^{{{P}_{2}}} \right)\left( \alpha _{A}^{{{P}_{3}}}-\alpha _{A}^{{{P}_{4}}} \right)-\left( \alpha _{C}^{{{P}_{3}}}-\alpha _{C}^{{{P}_{4}}} \right)\left( \alpha _{A}^{{{P}_{1}}}-\alpha _{A}^{{{P}_{2}}} \right)}\alpha _{B}^{{{P}_{1}}} \\ 
		& +\frac{\left( \alpha _{A}^{{{P}_{1}}}-\alpha _{A}^{{{P}_{3}}} \right)\left( \alpha _{B}^{{{P}_{3}}}-\alpha _{B}^{{{P}_{4}}} \right)-\left( \alpha _{A}^{{{P}_{3}}}-\alpha _{A}^{{{P}_{4}}} \right)\left( \alpha _{B}^{{{P}_{1}}}-\alpha _{B}^{{{P}_{3}}} \right)}{\left( \alpha _{A}^{{{P}_{1}}}-\alpha _{A}^{{{P}_{2}}} \right)\left( \alpha _{B}^{{{P}_{3}}}-\alpha _{B}^{{{P}_{4}}} \right)-\left( \alpha _{A}^{{{P}_{3}}}-\alpha _{A}^{{{P}_{4}}} \right)\left( \alpha _{B}^{{{P}_{1}}}-\alpha _{B}^{{{P}_{2}}} \right)}\alpha _{B}^{{{P}_{2}}},  
	\end{align*}
	
	\begin{align*}
		\alpha _{C}^{P}& =\left( 1-{{\kappa }_{{{P}_{1}}{{P}_{2}}}} \right)\alpha _{C}^{{{P}_{1}}}+{{\kappa }_{{{P}_{1}}{{P}_{2}}}}\alpha _{C}^{{{P}_{2}}} \\ 
		& =\frac{\left( \alpha _{A}^{{{P}_{3}}}-\alpha _{A}^{{{P}_{2}}} \right)\left( \alpha _{B}^{{{P}_{3}}}-\alpha _{B}^{{{P}_{4}}} \right)-\left( \alpha _{A}^{{{P}_{3}}}-\alpha _{A}^{{{P}_{4}}} \right)\left( \alpha _{B}^{{{P}_{3}}}-\alpha _{B}^{{{P}_{2}}} \right)}{\left( \alpha _{A}^{{{P}_{1}}}-\alpha _{A}^{{{P}_{2}}} \right)\left( \alpha _{B}^{{{P}_{3}}}-\alpha _{B}^{{{P}_{4}}} \right)-\left( \alpha _{A}^{{{P}_{3}}}-\alpha _{A}^{{{P}_{4}}} \right)\left( \alpha _{B}^{{{P}_{1}}}-\alpha _{B}^{{{P}_{2}}} \right)}\alpha _{C}^{{{P}_{1}}} \\ 
		& +\frac{\left( \alpha _{A}^{{{P}_{1}}}-\alpha _{A}^{{{P}_{3}}} \right)\left( \alpha _{B}^{{{P}_{3}}}-\alpha _{B}^{{{P}_{4}}} \right)-\left( \alpha _{A}^{{{P}_{3}}}-\alpha _{A}^{{{P}_{4}}} \right)\left( \alpha _{B}^{{{P}_{1}}}-\alpha _{B}^{{{P}_{3}}} \right)}{\left( \alpha _{A}^{{{P}_{1}}}-\alpha _{A}^{{{P}_{2}}} \right)\left( \alpha _{B}^{{{P}_{3}}}-\alpha _{B}^{{{P}_{4}}} \right)-\left( \alpha _{A}^{{{P}_{3}}}-\alpha _{A}^{{{P}_{4}}} \right)\left( \alpha _{B}^{{{P}_{1}}}-\alpha _{B}^{{{P}_{2}}} \right)}\alpha _{C}^{{{P}_{2}}}.  
	\end{align*}
\end{proof}
\hfill $\square$\par

In order to obtain a more aesthetically symmetrical form, I have provided the determinant form of the above theorem below.

\begin{theorem}{Frame component of the intersection of two straight lines-determinant form, Daiyuan Zhang}{LiangtiaoZhixianJiaodianDeBiaojiafenliang_Hanglieshixingshi}\label{LiangtiaoZhixianJiaodianDeBiaojiafenliang_Hanglieshixingshi}
	Given a $\triangle ABC$, assuming the frame component of IC ${{P}_{1}}$ is $\alpha _{A}^{{{P}_{1}}}$, $\alpha _{B}^{{{P}_{1}}}$, $\alpha _{C}^{{{P}_{1}}}$; the frame component of IC ${{P}_{2}}$ is $\alpha _{A}^{{{P}_{2}}}$, $\alpha _{B}^{{{P}_{2}}}$, $\alpha _{C}^{{{P}_{2}}}$; the frame component of IC ${{P}_{3}}$ is $\alpha _{A}^{{{P}_{3}}}$, $\alpha _{B}^{{{P}_{3}}}$, $\alpha _{C}^{{{P}_{3}}}$; the frame component of IC ${{P}_{4}}$ is $\alpha _{A}^{{{P}_{4}}}$, $\alpha _{B}^{{{P}_{4}}}$, $\alpha _{C}^{{{P}_{4}}}$. The intersection point of two straight lines $\overleftrightarrow{{{P}_{1}}{{P}_{2}}}$ and $\overleftrightarrow{{{P}_{3}}{{P}_{4}}}$ is $P=\overleftrightarrow{{{P}_{1}}{{P}_{2}}}\cap \overleftrightarrow{{{P}_{3}}{{P}_{4}}}$, and $\overrightarrow{{{P}_{1}}P}={{\kappa }_{{{P}_{1}}{{P}_{2}}}}\overrightarrow{{{P}_{1}}{{P}_{2}}}$, $\overrightarrow{{{P}_{1}}{{P}_{2}}}\ne \overrightarrow{0}$, $\overrightarrow{{{P}_{3}}P}={{\kappa }_{{{P}_{3}}{{P}_{4}}}}\overrightarrow{{{P}_{3}}{{P}_{4}}}$, $\overrightarrow{{{P}_{3}}{{P}_{4}}}\ne \overrightarrow{0}$, then the frame component of IC $P$ is:
	\[\alpha _{A}^{P}=\frac{1}{2\left| \begin{matrix}
			\alpha _{B}^{{{P}_{1}}}-\alpha _{B}^{{{P}_{2}}} & \alpha _{B}^{{{P}_{3}}}-\alpha _{B}^{{{P}_{4}}}  \\
			\alpha _{C}^{{{P}_{1}}}-\alpha _{C}^{{{P}_{2}}} & \alpha _{C}^{{{P}_{3}}}-\alpha _{C}^{{{P}_{4}}}  \\
		\end{matrix} \right|}\left| \begin{matrix}
		1 & 1 & 1 & 1  \\
		\alpha _{A}^{{{P}_{1}}} & \alpha _{A}^{{{P}_{2}}} & -\alpha _{A}^{{{P}_{3}}} & -\alpha _{A}^{{{P}_{4}}}  \\
		\alpha _{B}^{{{P}_{1}}} & \alpha _{B}^{{{P}_{2}}} & \alpha _{B}^{{{P}_{3}}} & \alpha _{B}^{{{P}_{4}}}  \\
		\alpha _{C}^{{{P}_{1}}} & \alpha _{C}^{{{P}_{2}}} & \alpha _{C}^{{{P}_{3}}} & \alpha _{C}^{{{P}_{4}}}  \\
	\end{matrix} \right|,\]
	\[\alpha _{B}^{P}=\frac{1}{2\left| \begin{matrix}
			\alpha _{C}^{{{P}_{1}}}-\alpha _{C}^{{{P}_{2}}} & \alpha _{C}^{{{P}_{3}}}-\alpha _{C}^{{{P}_{4}}}  \\
			\alpha _{A}^{{{P}_{1}}}-\alpha _{A}^{{{P}_{2}}} & \alpha _{A}^{{{P}_{3}}}-\alpha _{A}^{{{P}_{4}}}  \\
		\end{matrix} \right|}\left| \begin{matrix}
		1 & 1 & 1 & 1  \\
		\alpha _{B}^{{{P}_{1}}} & \alpha _{B}^{{{P}_{2}}} & -\alpha _{B}^{{{P}_{3}}} & -\alpha _{B}^{{{P}_{4}}}  \\
		\alpha _{C}^{{{P}_{1}}} & \alpha _{C}^{{{P}_{2}}} & \alpha _{C}^{{{P}_{3}}} & \alpha _{C}^{{{P}_{4}}}  \\
		\alpha _{A}^{{{P}_{1}}} & \alpha _{A}^{{{P}_{2}}} & \alpha _{A}^{{{P}_{3}}} & \alpha _{A}^{{{P}_{4}}}  \\
	\end{matrix} \right|,\]
	\[\alpha _{C}^{P}=\frac{1}{2\left| \begin{matrix}
			\alpha _{A}^{{{P}_{1}}}-\alpha _{A}^{{{P}_{2}}} & \alpha _{A}^{{{P}_{3}}}-\alpha _{A}^{{{P}_{4}}}  \\
			\alpha _{B}^{{{P}_{1}}}-\alpha _{B}^{{{P}_{2}}} & \alpha _{B}^{{{P}_{3}}}-\alpha _{B}^{{{P}_{4}}}  \\
		\end{matrix} \right|}\left| \begin{matrix}
		1 & 1 & 1 & 1  \\
		\alpha _{C}^{{{P}_{1}}} & \alpha _{C}^{{{P}_{2}}} & -\alpha _{C}^{{{P}_{3}}} & -\alpha _{C}^{{{P}_{4}}}  \\
		\alpha _{A}^{{{P}_{1}}} & \alpha _{A}^{{{P}_{2}}} & \alpha _{A}^{{{P}_{3}}} & \alpha _{A}^{{{P}_{4}}}  \\
		\alpha _{B}^{{{P}_{1}}} & \alpha _{B}^{{{P}_{2}}} & \alpha _{B}^{{{P}_{3}}} & \alpha _{B}^{{{P}_{4}}}  \\
	\end{matrix} \right|\text{.}\]
\end{theorem}

\begin{proof}
	Because the intersection point of two straight lines $\overleftrightarrow{{{P}_{1}}{{P}_{2}}}$ and $\overleftrightarrow{{{P}_{3}}{{P}_{4}}}$ is $P=\overleftrightarrow{{{P}_{1}}{{P}_{2}}}\cap \overleftrightarrow{{{P}_{3}}{{P}_{4}}}$, according to theorem \ref{thm:HuixinLianxianShangRenyiYidianDeBiaojiafenliang}, the following two sets of equations are obtained:
	
	\[\alpha _{A}^{P}=\left( 1-{{\kappa }_{{{P}_{1}}{{P}_{2}}}} \right)\alpha _{A}^{{{P}_{1}}}+{{\kappa }_{{{P}_{1}}{{P}_{2}}}}\alpha _{A}^{{{P}_{2}}},\]
	\[\alpha _{B}^{P}=\left( 1-{{\kappa }_{{{P}_{1}}{{P}_{2}}}} \right)\alpha _{B}^{{{P}_{1}}}+{{\kappa }_{{{P}_{1}}{{P}_{2}}}}\alpha _{B}^{{{P}_{2}}},\]
	\[\alpha _{C}^{P}=\left( 1-{{\kappa }_{{{P}_{1}}{{P}_{2}}}} \right)\alpha _{C}^{{{P}_{1}}}+{{\kappa }_{{{P}_{1}}{{P}_{2}}}}\alpha _{C}^{{{P}_{2}}}.\]
	
	And
	\[\alpha _{A}^{P}=\left( 1-{{\kappa }_{{{P}_{3}}{{P}_{4}}}} \right)\alpha _{A}^{{{P}_{3}}}+{{\kappa }_{{{P}_{3}}{{P}_{4}}}}\alpha _{A}^{{{P}_{4}}},\]
	\[\alpha _{B}^{P}=\left( 1-{{\kappa }_{{{P}_{3}}{{P}_{4}}}} \right)\alpha _{B}^{{{P}_{3}}}+{{\kappa }_{{{P}_{3}}{{P}_{4}}}}\alpha _{B}^{{{P}_{4}}},\]
	\[\alpha _{C}^{P}=\left( 1-{{\kappa }_{{{P}_{3}}{{P}_{4}}}} \right)\alpha _{C}^{{{P}_{3}}}+{{\kappa }_{{{P}_{3}}{{P}_{4}}}}\alpha _{C}^{{{P}_{4}}}.\]
	
	So, the system of equations ($AB$ combination) is obtained from the first and second equations:
	\[\left( 1-{{\kappa }_{{{P}_{1}}{{P}_{2}}}} \right)\alpha _{A}^{{{P}_{1}}}+{{\kappa }_{{{P}_{1}}{{P}_{2}}}}\alpha _{A}^{{{P}_{2}}}=\left( 1-{{\kappa }_{{{P}_{3}}{{P}_{4}}}} \right)\alpha _{A}^{{{P}_{3}}}+{{\kappa }_{{{P}_{3}}{{P}_{4}}}}\alpha _{A}^{{{P}_{4}}},\]
	\[\left( 1-{{\kappa }_{{{P}_{1}}{{P}_{2}}}} \right)\alpha _{B}^{{{P}_{1}}}+{{\kappa }_{{{P}_{1}}{{P}_{2}}}}\alpha _{B}^{{{P}_{2}}}=\left( 1-{{\kappa }_{{{P}_{3}}{{P}_{4}}}} \right)\alpha _{B}^{{{P}_{3}}}+{{\kappa }_{{{P}_{3}}{{P}_{4}}}}\alpha _{B}^{{{P}_{4}}}.\]
	
	The system of equations ($BC$ combination) is obtained from the second and third equations:
	\[\left( 1-{{\kappa }_{{{P}_{1}}{{P}_{2}}}} \right)\alpha _{B}^{{{P}_{1}}}+{{\kappa }_{{{P}_{1}}{{P}_{2}}}}\alpha _{B}^{{{P}_{2}}}=\left( 1-{{\kappa }_{{{P}_{3}}{{P}_{4}}}} \right)\alpha _{B}^{{{P}_{3}}}+{{\kappa }_{{{P}_{3}}{{P}_{4}}}}\alpha _{B}^{{{P}_{4}}}.\]
	\[\left( 1-{{\kappa }_{{{P}_{1}}{{P}_{2}}}} \right)\alpha _{C}^{{{P}_{1}}}+{{\kappa }_{{{P}_{1}}{{P}_{2}}}}\alpha _{C}^{{{P}_{2}}}=\left( 1-{{\kappa }_{{{P}_{3}}{{P}_{4}}}} \right)\alpha _{C}^{{{P}_{3}}}+{{\kappa }_{{{P}_{3}}{{P}_{4}}}}\alpha _{C}^{{{P}_{4}}},\]
	
	The system of equations ($CA$ combination) is obtained from the third and first equations:
	\[\left( 1-{{\kappa }_{{{P}_{1}}{{P}_{2}}}} \right)\alpha _{C}^{{{P}_{1}}}+{{\kappa }_{{{P}_{1}}{{P}_{2}}}}\alpha _{C}^{{{P}_{2}}}=\left( 1-{{\kappa }_{{{P}_{3}}{{P}_{4}}}} \right)\alpha _{C}^{{{P}_{3}}}+{{\kappa }_{{{P}_{3}}{{P}_{4}}}}\alpha _{C}^{{{P}_{4}}},\]
	\[\left( 1-{{\kappa }_{{{P}_{1}}{{P}_{2}}}} \right)\alpha _{A}^{{{P}_{1}}}+{{\kappa }_{{{P}_{1}}{{P}_{2}}}}\alpha _{A}^{{{P}_{2}}}=\left( 1-{{\kappa }_{{{P}_{3}}{{P}_{4}}}} \right)\alpha _{A}^{{{P}_{3}}}+{{\kappa }_{{{P}_{3}}{{P}_{4}}}}\alpha _{A}^{{{P}_{4}}}.\]
	
	It can be seen that the system of equations is obtained through rotation, so its solution is also obtained through rotation. Therefore, only one system of equations needs to be solved, and then the form of solutions for other systems of equations can be obtained through rotation.
	
	The following is the system of equations ($AB$ combination) obtained from the first and second equations:
	\[\left( 1-x \right)\alpha _{A}^{{{P}_{1}}}+x\alpha _{A}^{{{P}_{2}}}=\left( 1-y \right)\alpha _{A}^{{{P}_{3}}}+y\alpha _{A}^{{{P}_{4}}},\]
	\[\left( 1-x \right)\alpha _{B}^{{{P}_{1}}}+x\alpha _{B}^{{{P}_{2}}}=\left( 1-y \right)\alpha _{B}^{{{P}_{3}}}+y\alpha _{B}^{{{P}_{4}}}.\]
	
	\[\left\{ \begin{aligned}
		& \left( \alpha _{A}^{{{P}_{2}}}-\alpha _{A}^{{{P}_{1}}} \right)x-\left( \alpha _{A}^{{{P}_{4}}}-\alpha _{A}^{{{P}_{3}}} \right)y=\alpha _{A}^{{{P}_{3}}}-\alpha _{A}^{{{P}_{1}}} \\ 
		& \left( \alpha _{B}^{{{P}_{2}}}-\alpha _{B}^{{{P}_{1}}} \right)x-\left( \alpha _{B}^{{{P}_{4}}}-\alpha _{B}^{{{P}_{3}}} \right)y=\alpha _{B}^{{{P}_{3}}}-\alpha _{B}^{{{P}_{1}}}. \\ 
	\end{aligned} \right.\]
	
	Solving the above equation system yields:
	\[x=\frac{\left| \begin{matrix}
			\alpha _{A}^{{{P}_{3}}}-\alpha _{A}^{{{P}_{1}}} & -\left( \alpha _{A}^{{{P}_{4}}}-\alpha _{A}^{{{P}_{3}}} \right)  \\
			\alpha _{B}^{{{P}_{3}}}-\alpha _{B}^{{{P}_{1}}} & -\left( \alpha _{B}^{{{P}_{4}}}-\alpha _{B}^{{{P}_{3}}} \right)  \\
		\end{matrix} \right|}{\left| \begin{matrix}
			\alpha _{A}^{{{P}_{2}}}-\alpha _{A}^{{{P}_{1}}} & -\left( \alpha _{A}^{{{P}_{4}}}-\alpha _{A}^{{{P}_{3}}} \right)  \\
			\alpha _{B}^{{{P}_{2}}}-\alpha _{B}^{{{P}_{1}}} & -\left( \alpha _{B}^{{{P}_{4}}}-\alpha _{B}^{{{P}_{3}}} \right)  \\
		\end{matrix} \right|}=\frac{\left( \alpha _{A}^{{{P}_{4}}}-\alpha _{A}^{{{P}_{3}}} \right)\left( \alpha _{B}^{{{P}_{3}}}-\alpha _{B}^{{{P}_{1}}} \right)-\left( \alpha _{A}^{{{P}_{3}}}-\alpha _{A}^{{{P}_{1}}} \right)\left( \alpha _{B}^{{{P}_{4}}}-\alpha _{B}^{{{P}_{3}}} \right)}{\left( \alpha _{A}^{{{P}_{4}}}-\alpha _{A}^{{{P}_{3}}} \right)\left( \alpha _{B}^{{{P}_{2}}}-\alpha _{B}^{{{P}_{1}}} \right)-\left( \alpha _{A}^{{{P}_{2}}}-\alpha _{A}^{{{P}_{1}}} \right)\left( \alpha _{B}^{{{P}_{4}}}-\alpha _{B}^{{{P}_{3}}} \right)},\]
	\[y=\frac{\left| \begin{matrix}
			\alpha _{A}^{{{P}_{2}}}-\alpha _{A}^{{{P}_{1}}} & \alpha _{A}^{{{P}_{3}}}-\alpha _{A}^{{{P}_{1}}}  \\
			\alpha _{B}^{{{P}_{2}}}-\alpha _{B}^{{{P}_{1}}} & \alpha _{B}^{{{P}_{3}}}-\alpha _{B}^{{{P}_{1}}}  \\
		\end{matrix} \right|}{\left| \begin{matrix}
			\alpha _{A}^{{{P}_{2}}}-\alpha _{A}^{{{P}_{1}}} & -\left( \alpha _{A}^{{{P}_{4}}}-\alpha _{A}^{{{P}_{3}}} \right)  \\
			\alpha _{B}^{{{P}_{2}}}-\alpha _{B}^{{{P}_{1}}} & -\left( \alpha _{B}^{{{P}_{4}}}-\alpha _{B}^{{{P}_{3}}} \right)  \\
		\end{matrix} \right|}=\frac{\left( \alpha _{A}^{{{P}_{2}}}-\alpha _{A}^{{{P}_{1}}} \right)\left( \alpha _{B}^{{{P}_{3}}}-\alpha _{B}^{{{P}_{1}}} \right)-\left( \alpha _{A}^{{{P}_{3}}}-\alpha _{A}^{{{P}_{1}}} \right)\left( \alpha _{B}^{{{P}_{2}}}-\alpha _{B}^{{{P}_{1}}} \right)}{\left( \alpha _{A}^{{{P}_{4}}}-\alpha _{A}^{{{P}_{3}}} \right)\left( \alpha _{B}^{{{P}_{2}}}-\alpha _{B}^{{{P}_{1}}} \right)-\left( \alpha _{A}^{{{P}_{2}}}-\alpha _{A}^{{{P}_{1}}} \right)\left( \alpha _{B}^{{{P}_{4}}}-\alpha _{B}^{{{P}_{3}}} \right)}.\]
	
	Therefore
	\[{{\kappa }_{{{P}_{1}}{{P}_{2}}}}=x=\frac{\left( \alpha _{A}^{{{P}_{1}}}-\alpha _{A}^{{{P}_{3}}} \right)\left( \alpha _{B}^{{{P}_{3}}}-\alpha _{B}^{{{P}_{4}}} \right)-\left( \alpha _{A}^{{{P}_{3}}}-\alpha _{A}^{{{P}_{4}}} \right)\left( \alpha _{B}^{{{P}_{1}}}-\alpha _{B}^{{{P}_{3}}} \right)}{\left( \alpha _{A}^{{{P}_{1}}}-\alpha _{A}^{{{P}_{2}}} \right)\left( \alpha _{B}^{{{P}_{3}}}-\alpha _{B}^{{{P}_{4}}} \right)-\left( \alpha _{A}^{{{P}_{3}}}-\alpha _{A}^{{{P}_{4}}} \right)\left( \alpha _{B}^{{{P}_{1}}}-\alpha _{B}^{{{P}_{2}}} \right)},\]
	
	\begin{align*}
		1-{{\kappa }_{{{P}_{1}}{{P}_{2}}}}& =1-\frac{\left( \alpha _{A}^{{{P}_{1}}}-\alpha _{A}^{{{P}_{3}}} \right)\left( \alpha _{B}^{{{P}_{3}}}-\alpha _{B}^{{{P}_{4}}} \right)-\left( \alpha _{A}^{{{P}_{3}}}-\alpha _{A}^{{{P}_{4}}} \right)\left( \alpha _{B}^{{{P}_{1}}}-\alpha _{B}^{{{P}_{3}}} \right)}{\left( \alpha _{A}^{{{P}_{1}}}-\alpha _{A}^{{{P}_{2}}} \right)\left( \alpha _{B}^{{{P}_{3}}}-\alpha _{B}^{{{P}_{4}}} \right)-\left( \alpha _{A}^{{{P}_{3}}}-\alpha _{A}^{{{P}_{4}}} \right)\left( \alpha _{B}^{{{P}_{1}}}-\alpha _{B}^{{{P}_{2}}} \right)} \\ 
		& =\frac{\left( \alpha _{A}^{{{P}_{3}}}-\alpha _{A}^{{{P}_{2}}} \right)\left( \alpha _{B}^{{{P}_{3}}}-\alpha _{B}^{{{P}_{4}}} \right)-\left( \alpha _{A}^{{{P}_{3}}}-\alpha _{A}^{{{P}_{4}}} \right)\left( \alpha _{B}^{{{P}_{3}}}-\alpha _{B}^{{{P}_{2}}} \right)}{\left( \alpha _{A}^{{{P}_{1}}}-\alpha _{A}^{{{P}_{2}}} \right)\left( \alpha _{B}^{{{P}_{3}}}-\alpha _{B}^{{{P}_{4}}} \right)-\left( \alpha _{A}^{{{P}_{3}}}-\alpha _{A}^{{{P}_{4}}} \right)\left( \alpha _{B}^{{{P}_{1}}}-\alpha _{B}^{{{P}_{2}}} \right)},  
	\end{align*}
	and 
	\[{{\kappa }_{{{P}_{3}}{{P}_{4}}}}=y=\frac{\left( \alpha _{A}^{{{P}_{1}}}-\alpha _{A}^{{{P}_{3}}} \right)\left( \alpha _{B}^{{{P}_{1}}}-\alpha _{B}^{{{P}_{2}}} \right)-\left( \alpha _{A}^{{{P}_{1}}}-\alpha _{A}^{{{P}_{2}}} \right)\left( \alpha _{B}^{{{P}_{1}}}-\alpha _{B}^{{{P}_{3}}} \right)}{\left( \alpha _{A}^{{{P}_{1}}}-\alpha _{A}^{{{P}_{2}}} \right)\left( \alpha _{B}^{{{P}_{3}}}-\alpha _{B}^{{{P}_{4}}} \right)-\left( \alpha _{A}^{{{P}_{3}}}-\alpha _{A}^{{{P}_{4}}} \right)\left( \alpha _{B}^{{{P}_{1}}}-\alpha _{B}^{{{P}_{2}}} \right)}.\]
	
	\begin{align*}
		1-{{\kappa }_{{{P}_{3}}{{P}_{4}}}}& =1-\frac{\left( \alpha _{A}^{{{P}_{1}}}-\alpha _{A}^{{{P}_{3}}} \right)\left( \alpha _{B}^{{{P}_{1}}}-\alpha _{B}^{{{P}_{2}}} \right)-\left( \alpha _{A}^{{{P}_{1}}}-\alpha _{A}^{{{P}_{2}}} \right)\left( \alpha _{B}^{{{P}_{1}}}-\alpha _{B}^{{{P}_{3}}} \right)}{\left( \alpha _{A}^{{{P}_{1}}}-\alpha _{A}^{{{P}_{2}}} \right)\left( \alpha _{B}^{{{P}_{3}}}-\alpha _{B}^{{{P}_{4}}} \right)-\left( \alpha _{A}^{{{P}_{3}}}-\alpha _{A}^{{{P}_{4}}} \right)\left( \alpha _{B}^{{{P}_{1}}}-\alpha _{B}^{{{P}_{2}}} \right)} \\ 
		& =\frac{\left( \alpha _{A}^{{{P}_{4}}}-\alpha _{A}^{{{P}_{1}}} \right)\left( \alpha _{B}^{{{P}_{1}}}-\alpha _{B}^{{{P}_{2}}} \right)-\left( \alpha _{A}^{{{P}_{1}}}-\alpha _{A}^{{{P}_{2}}} \right)\left( \alpha _{B}^{{{P}_{4}}}-\alpha _{B}^{{{P}_{1}}} \right)}{\left( \alpha _{A}^{{{P}_{1}}}-\alpha _{A}^{{{P}_{2}}} \right)\left( \alpha _{B}^{{{P}_{3}}}-\alpha _{B}^{{{P}_{4}}} \right)-\left( \alpha _{A}^{{{P}_{3}}}-\alpha _{A}^{{{P}_{4}}} \right)\left( \alpha _{B}^{{{P}_{1}}}-\alpha _{B}^{{{P}_{2}}} \right)}.  
	\end{align*}
	
	After rotation, ($BC$ combination) we have:
	\[{{\kappa }_{{{P}_{1}}{{P}_{2}}}}=\frac{\left( \alpha _{B}^{{{P}_{1}}}-\alpha _{B}^{{{P}_{3}}} \right)\left( \alpha _{C}^{{{P}_{3}}}-\alpha _{C}^{{{P}_{4}}} \right)-\left( \alpha _{B}^{{{P}_{3}}}-\alpha _{B}^{{{P}_{4}}} \right)\left( \alpha _{C}^{{{P}_{1}}}-\alpha _{C}^{{{P}_{3}}} \right)}{\left( \alpha _{B}^{{{P}_{1}}}-\alpha _{B}^{{{P}_{2}}} \right)\left( \alpha _{C}^{{{P}_{3}}}-\alpha _{C}^{{{P}_{4}}} \right)-\left( \alpha _{B}^{{{P}_{3}}}-\alpha _{B}^{{{P}_{4}}} \right)\left( \alpha _{C}^{{{P}_{1}}}-\alpha _{C}^{{{P}_{2}}} \right)},\]
	\[1-{{\kappa }_{{{P}_{1}}{{P}_{2}}}}=\frac{\left( \alpha _{B}^{{{P}_{3}}}-\alpha _{B}^{{{P}_{2}}} \right)\left( \alpha _{C}^{{{P}_{3}}}-\alpha _{C}^{{{P}_{4}}} \right)-\left( \alpha _{B}^{{{P}_{3}}}-\alpha _{B}^{{{P}_{4}}} \right)\left( \alpha _{C}^{{{P}_{3}}}-\alpha _{C}^{{{P}_{2}}} \right)}{\left( \alpha _{B}^{{{P}_{1}}}-\alpha _{B}^{{{P}_{2}}} \right)\left( \alpha _{C}^{{{P}_{3}}}-\alpha _{C}^{{{P}_{4}}} \right)-\left( \alpha _{B}^{{{P}_{3}}}-\alpha _{B}^{{{P}_{4}}} \right)\left( \alpha _{C}^{{{P}_{1}}}-\alpha _{C}^{{{P}_{2}}} \right)},\]
	\[{{\kappa }_{{{P}_{3}}{{P}_{4}}}}=\frac{\left( \alpha _{B}^{{{P}_{1}}}-\alpha _{B}^{{{P}_{3}}} \right)\left( \alpha _{C}^{{{P}_{1}}}-\alpha _{C}^{{{P}_{2}}} \right)-\left( \alpha _{B}^{{{P}_{1}}}-\alpha _{B}^{{{P}_{2}}} \right)\left( \alpha _{C}^{{{P}_{1}}}-\alpha _{C}^{{{P}_{3}}} \right)}{\left( \alpha _{B}^{{{P}_{1}}}-\alpha _{B}^{{{P}_{2}}} \right)\left( \alpha _{C}^{{{P}_{3}}}-\alpha _{C}^{{{P}_{4}}} \right)-\left( \alpha _{B}^{{{P}_{3}}}-\alpha _{B}^{{{P}_{4}}} \right)\left( \alpha _{C}^{{{P}_{1}}}-\alpha _{C}^{{{P}_{2}}} \right)}.\]
	\[1-{{\kappa }_{{{P}_{3}}{{P}_{4}}}}=\frac{\left( \alpha _{B}^{{{P}_{4}}}-\alpha _{B}^{{{P}_{1}}} \right)\left( \alpha _{C}^{{{P}_{1}}}-\alpha _{C}^{{{P}_{2}}} \right)-\left( \alpha _{B}^{{{P}_{1}}}-\alpha _{B}^{{{P}_{2}}} \right)\left( \alpha _{C}^{{{P}_{4}}}-\alpha _{C}^{{{P}_{1}}} \right)}{\left( \alpha _{B}^{{{P}_{1}}}-\alpha _{B}^{{{P}_{2}}} \right)\left( \alpha _{C}^{{{P}_{3}}}-\alpha _{C}^{{{P}_{4}}} \right)-\left( \alpha _{B}^{{{P}_{3}}}-\alpha _{B}^{{{P}_{4}}} \right)\left( \alpha _{C}^{{{P}_{1}}}-\alpha _{C}^{{{P}_{2}}} \right)}.\]
	
	After rotation, ($CA$ combination) we have:
	\[{{\kappa }_{{{P}_{1}}{{P}_{2}}}}=\frac{\left( \alpha _{C}^{{{P}_{1}}}-\alpha _{C}^{{{P}_{3}}} \right)\left( \alpha _{A}^{{{P}_{3}}}-\alpha _{A}^{{{P}_{4}}} \right)-\left( \alpha _{C}^{{{P}_{3}}}-\alpha _{C}^{{{P}_{4}}} \right)\left( \alpha _{A}^{{{P}_{1}}}-\alpha _{A}^{{{P}_{3}}} \right)}{\left( \alpha _{C}^{{{P}_{1}}}-\alpha _{C}^{{{P}_{2}}} \right)\left( \alpha _{A}^{{{P}_{3}}}-\alpha _{A}^{{{P}_{4}}} \right)-\left( \alpha _{C}^{{{P}_{3}}}-\alpha _{C}^{{{P}_{4}}} \right)\left( \alpha _{A}^{{{P}_{1}}}-\alpha _{A}^{{{P}_{2}}} \right)},\]
	\[1-{{\kappa }_{{{P}_{1}}{{P}_{2}}}}=\frac{\left( \alpha _{C}^{{{P}_{3}}}-\alpha _{C}^{{{P}_{2}}} \right)\left( \alpha _{A}^{{{P}_{3}}}-\alpha _{A}^{{{P}_{4}}} \right)-\left( \alpha _{C}^{{{P}_{3}}}-\alpha _{C}^{{{P}_{4}}} \right)\left( \alpha _{A}^{{{P}_{3}}}-\alpha _{A}^{{{P}_{2}}} \right)}{\left( \alpha _{C}^{{{P}_{1}}}-\alpha _{C}^{{{P}_{2}}} \right)\left( \alpha _{A}^{{{P}_{3}}}-\alpha _{A}^{{{P}_{4}}} \right)-\left( \alpha _{C}^{{{P}_{3}}}-\alpha _{C}^{{{P}_{4}}} \right)\left( \alpha _{A}^{{{P}_{1}}}-\alpha _{A}^{{{P}_{2}}} \right)},\]
	\[{{\kappa }_{{{P}_{3}}{{P}_{4}}}}=\frac{\left( \alpha _{C}^{{{P}_{1}}}-\alpha _{C}^{{{P}_{3}}} \right)\left( \alpha _{A}^{{{P}_{1}}}-\alpha _{A}^{{{P}_{2}}} \right)-\left( \alpha _{C}^{{{P}_{1}}}-\alpha _{C}^{{{P}_{2}}} \right)\left( \alpha _{A}^{{{P}_{1}}}-\alpha _{A}^{{{P}_{3}}} \right)}{\left( \alpha _{C}^{{{P}_{1}}}-\alpha _{C}^{{{P}_{2}}} \right)\left( \alpha _{A}^{{{P}_{3}}}-\alpha _{A}^{{{P}_{4}}} \right)-\left( \alpha _{C}^{{{P}_{3}}}-\alpha _{C}^{{{P}_{4}}} \right)\left( \alpha _{A}^{{{P}_{1}}}-\alpha _{A}^{{{P}_{2}}} \right)}.\]
	\[1-{{\kappa }_{{{P}_{3}}{{P}_{4}}}}=\frac{\left( \alpha _{C}^{{{P}_{4}}}-\alpha _{C}^{{{P}_{1}}} \right)\left( \alpha _{A}^{{{P}_{1}}}-\alpha _{A}^{{{P}_{2}}} \right)-\left( \alpha _{C}^{{{P}_{1}}}-\alpha _{C}^{{{P}_{2}}} \right)\left( \alpha _{A}^{{{P}_{4}}}-\alpha _{A}^{{{P}_{1}}} \right)}{\left( \alpha _{C}^{{{P}_{1}}}-\alpha _{C}^{{{P}_{2}}} \right)\left( \alpha _{A}^{{{P}_{3}}}-\alpha _{A}^{{{P}_{4}}} \right)-\left( \alpha _{C}^{{{P}_{3}}}-\alpha _{C}^{{{P}_{4}}} \right)\left( \alpha _{A}^{{{P}_{1}}}-\alpha _{A}^{{{P}_{2}}} \right)}.\]
	
	To obtain a symmetrical and aesthetically pleasing form, the following combinations are used to solve for the frame components:
	\begin{align*}
		\alpha _{A}^{P}& =\left( 1-{{\kappa }_{{{P}_{1}}{{P}_{2}}}} \right)\alpha _{A}^{{{P}_{1}}}+{{\kappa }_{{{P}_{1}}{{P}_{2}}}}\alpha _{A}^{{{P}_{2}}} \\ 
		& =\frac{\left( \alpha _{B}^{{{P}_{3}}}-\alpha _{B}^{{{P}_{2}}} \right)\left( \alpha _{C}^{{{P}_{3}}}-\alpha _{C}^{{{P}_{4}}} \right)-\left( \alpha _{B}^{{{P}_{3}}}-\alpha _{B}^{{{P}_{4}}} \right)\left( \alpha _{C}^{{{P}_{3}}}-\alpha _{C}^{{{P}_{2}}} \right)}{\left( \alpha _{B}^{{{P}_{1}}}-\alpha _{B}^{{{P}_{2}}} \right)\left( \alpha _{C}^{{{P}_{3}}}-\alpha _{C}^{{{P}_{4}}} \right)-\left( \alpha _{B}^{{{P}_{3}}}-\alpha _{B}^{{{P}_{4}}} \right)\left( \alpha _{C}^{{{P}_{1}}}-\alpha _{C}^{{{P}_{2}}} \right)}\alpha _{A}^{{{P}_{1}}} \\ 
		& +\frac{\left( \alpha _{B}^{{{P}_{1}}}-\alpha _{B}^{{{P}_{3}}} \right)\left( \alpha _{C}^{{{P}_{3}}}-\alpha _{C}^{{{P}_{4}}} \right)-\left( \alpha _{B}^{{{P}_{3}}}-\alpha _{B}^{{{P}_{4}}} \right)\left( \alpha _{C}^{{{P}_{1}}}-\alpha _{C}^{{{P}_{3}}} \right)}{\left( \alpha _{B}^{{{P}_{1}}}-\alpha _{B}^{{{P}_{2}}} \right)\left( \alpha _{C}^{{{P}_{3}}}-\alpha _{C}^{{{P}_{4}}} \right)-\left( \alpha _{B}^{{{P}_{3}}}-\alpha _{B}^{{{P}_{4}}} \right)\left( \alpha _{C}^{{{P}_{1}}}-\alpha _{C}^{{{P}_{2}}} \right)}\alpha _{A}^{{{P}_{2}}},  
	\end{align*}
	
	\begin{align*}
		\alpha _{B}^{P}& =\left( 1-{{\kappa }_{{{P}_{1}}{{P}_{2}}}} \right)\alpha _{B}^{{{P}_{1}}}+{{\kappa }_{{{P}_{1}}{{P}_{2}}}}\alpha _{B}^{{{P}_{2}}} \\ 
		& \text{=}\frac{\left( \alpha _{C}^{{{P}_{1}}}-\alpha _{C}^{{{P}_{3}}} \right)\left( \alpha _{A}^{{{P}_{3}}}-\alpha _{A}^{{{P}_{4}}} \right)-\left( \alpha _{C}^{{{P}_{3}}}-\alpha _{C}^{{{P}_{4}}} \right)\left( \alpha _{A}^{{{P}_{1}}}-\alpha _{A}^{{{P}_{3}}} \right)}{\left( \alpha _{C}^{{{P}_{1}}}-\alpha _{C}^{{{P}_{2}}} \right)\left( \alpha _{A}^{{{P}_{3}}}-\alpha _{A}^{{{P}_{4}}} \right)-\left( \alpha _{C}^{{{P}_{3}}}-\alpha _{C}^{{{P}_{4}}} \right)\left( \alpha _{A}^{{{P}_{1}}}-\alpha _{A}^{{{P}_{2}}} \right)}\alpha _{B}^{{{P}_{1}}} \\ 
		& +\frac{\left( \alpha _{A}^{{{P}_{1}}}-\alpha _{A}^{{{P}_{3}}} \right)\left( \alpha _{B}^{{{P}_{3}}}-\alpha _{B}^{{{P}_{4}}} \right)-\left( \alpha _{A}^{{{P}_{3}}}-\alpha _{A}^{{{P}_{4}}} \right)\left( \alpha _{B}^{{{P}_{1}}}-\alpha _{B}^{{{P}_{3}}} \right)}{\left( \alpha _{A}^{{{P}_{1}}}-\alpha _{A}^{{{P}_{2}}} \right)\left( \alpha _{B}^{{{P}_{3}}}-\alpha _{B}^{{{P}_{4}}} \right)-\left( \alpha _{A}^{{{P}_{3}}}-\alpha _{A}^{{{P}_{4}}} \right)\left( \alpha _{B}^{{{P}_{1}}}-\alpha _{B}^{{{P}_{2}}} \right)}\alpha _{B}^{{{P}_{2}}},  
	\end{align*}
	
	\begin{align*}
		\alpha _{C}^{P}& =\left( 1-{{\kappa }_{{{P}_{1}}{{P}_{2}}}} \right)\alpha _{C}^{{{P}_{1}}}+{{\kappa }_{{{P}_{1}}{{P}_{2}}}}\alpha _{C}^{{{P}_{2}}} \\ 
		& =\frac{\left( \alpha _{A}^{{{P}_{3}}}-\alpha _{A}^{{{P}_{2}}} \right)\left( \alpha _{B}^{{{P}_{3}}}-\alpha _{B}^{{{P}_{4}}} \right)-\left( \alpha _{A}^{{{P}_{3}}}-\alpha _{A}^{{{P}_{4}}} \right)\left( \alpha _{B}^{{{P}_{3}}}-\alpha _{B}^{{{P}_{2}}} \right)}{\left( \alpha _{A}^{{{P}_{1}}}-\alpha _{A}^{{{P}_{2}}} \right)\left( \alpha _{B}^{{{P}_{3}}}-\alpha _{B}^{{{P}_{4}}} \right)-\left( \alpha _{A}^{{{P}_{3}}}-\alpha _{A}^{{{P}_{4}}} \right)\left( \alpha _{B}^{{{P}_{1}}}-\alpha _{B}^{{{P}_{2}}} \right)}\alpha _{C}^{{{P}_{1}}} \\ 
		& +\frac{\left( \alpha _{A}^{{{P}_{1}}}-\alpha _{A}^{{{P}_{3}}} \right)\left( \alpha _{B}^{{{P}_{3}}}-\alpha _{B}^{{{P}_{4}}} \right)-\left( \alpha _{A}^{{{P}_{3}}}-\alpha _{A}^{{{P}_{4}}} \right)\left( \alpha _{B}^{{{P}_{1}}}-\alpha _{B}^{{{P}_{3}}} \right)}{\left( \alpha _{A}^{{{P}_{1}}}-\alpha _{A}^{{{P}_{2}}} \right)\left( \alpha _{B}^{{{P}_{3}}}-\alpha _{B}^{{{P}_{4}}} \right)-\left( \alpha _{A}^{{{P}_{3}}}-\alpha _{A}^{{{P}_{4}}} \right)\left( \alpha _{B}^{{{P}_{1}}}-\alpha _{B}^{{{P}_{2}}} \right)}\alpha _{C}^{{{P}_{2}}}.  
	\end{align*}
	
	Write $\alpha _{A}^{P}$ in determinant form. First, consider the ${{\kappa }_{{{P}_{1}}{{P}_{2}}}}$ of $BC$ combination. The numerator of ${{\kappa }_{{{P}_{1}}{{P}_{2}}}}$ is
	
	\[\begin{aligned}
		& \left( \alpha _{B}^{{{P}_{1}}}-\alpha _{B}^{{{P}_{3}}} \right)\left( \alpha _{C}^{{{P}_{3}}}-\alpha _{C}^{{{P}_{4}}} \right)-\left( \alpha _{B}^{{{P}_{3}}}-\alpha _{B}^{{{P}_{4}}} \right)\left( \alpha _{C}^{{{P}_{1}}}-\alpha _{C}^{{{P}_{3}}} \right)=\left| \begin{matrix}
			\alpha _{B}^{{{P}_{1}}}-\alpha _{B}^{{{P}_{3}}} & \alpha _{B}^{{{P}_{3}}}-\alpha _{B}^{{{P}_{4}}}  \\
			\alpha _{C}^{{{P}_{1}}}-\alpha _{C}^{{{P}_{3}}} & \alpha _{C}^{{{P}_{3}}}-\alpha _{C}^{{{P}_{4}}}  \\
		\end{matrix} \right| \\ 
		& =\left| \begin{matrix}
			1 & 0 & 0  \\
			\alpha _{B}^{{{P}_{3}}} & \alpha _{B}^{{{P}_{1}}}-\alpha _{B}^{{{P}_{3}}} & \alpha _{B}^{{{P}_{3}}}-\alpha _{B}^{{{P}_{4}}}  \\
			\alpha _{C}^{{{P}_{3}}} & \alpha _{C}^{{{P}_{1}}}-\alpha _{C}^{{{P}_{3}}} & \alpha _{C}^{{{P}_{3}}}-\alpha _{C}^{{{P}_{4}}}  \\
		\end{matrix} \right|=\left| \begin{matrix}
			1 & 1 & -1  \\
			\alpha _{B}^{{{P}_{3}}} & \alpha _{B}^{{{P}_{1}}} & -\alpha _{B}^{{{P}_{4}}}  \\
			\alpha _{C}^{{{P}_{3}}} & \alpha _{C}^{{{P}_{1}}} & -\alpha _{C}^{{{P}_{4}}}  \\
		\end{matrix} \right|=\left| \begin{matrix}
			1 & 1 & 1  \\
			\alpha _{B}^{{{P}_{1}}} & \alpha _{B}^{{{P}_{3}}} & \alpha _{B}^{{{P}_{4}}}  \\
			\alpha _{C}^{{{P}_{1}}} & \alpha _{C}^{{{P}_{3}}} & \alpha _{C}^{{{P}_{4}}}  \\
		\end{matrix} \right|. \\ 
	\end{aligned}\]
	
	The numerator of $1-{{\kappa }_{{{P}_{1}}{{P}_{2}}}}$ is
	
	\[\begin{aligned}
		& \left( \alpha _{B}^{{{P}_{3}}}-\alpha _{B}^{{{P}_{2}}} \right)\left( \alpha _{C}^{{{P}_{3}}}-\alpha _{C}^{{{P}_{4}}} \right)-\left( \alpha _{B}^{{{P}_{3}}}-\alpha _{B}^{{{P}_{4}}} \right)\left( \alpha _{C}^{{{P}_{3}}}-\alpha _{C}^{{{P}_{2}}} \right)=\left| \begin{matrix}
			\alpha _{B}^{{{P}_{3}}}-\alpha _{B}^{{{P}_{2}}} & \alpha _{B}^{{{P}_{3}}}-\alpha _{B}^{{{P}_{4}}}  \\
			\alpha _{C}^{{{P}_{3}}}-\alpha _{C}^{{{P}_{2}}} & \alpha _{C}^{{{P}_{3}}}-\alpha _{C}^{{{P}_{4}}}  \\
		\end{matrix} \right| \\ 
		& =\left| \begin{matrix}
			1 & 0 & 0  \\
			\alpha _{B}^{{{P}_{3}}} & \alpha _{B}^{{{P}_{3}}}-\alpha _{B}^{{{P}_{2}}} & \alpha _{B}^{{{P}_{3}}}-\alpha _{B}^{{{P}_{4}}}  \\
			\alpha _{C}^{{{P}_{3}}} & \alpha _{C}^{{{P}_{3}}}-\alpha _{C}^{{{P}_{2}}} & \alpha _{C}^{{{P}_{3}}}-\alpha _{C}^{{{P}_{4}}}  \\
		\end{matrix} \right|=\left| \begin{matrix}
			1 & -1 & -1  \\
			\alpha _{B}^{{{P}_{3}}} & -\alpha _{B}^{{{P}_{2}}} & -\alpha _{B}^{{{P}_{4}}}  \\
			\alpha _{C}^{{{P}_{3}}} & -\alpha _{C}^{{{P}_{1}}} & -\alpha _{C}^{{{P}_{4}}}  \\
		\end{matrix} \right|=-\left| \begin{matrix}
			1 & 1 & 1  \\
			\alpha _{B}^{{{P}_{2}}} & \alpha _{B}^{{{P}_{3}}} & \alpha _{B}^{{{P}_{4}}}  \\
			\alpha _{C}^{{{P}_{2}}} & \alpha _{C}^{{{P}_{3}}} & \alpha _{C}^{{{P}_{4}}}  \\
		\end{matrix} \right|. \\ 
	\end{aligned}\]
	
	The denominator of the determinant is
	\[\left( \alpha _{B}^{{{P}_{1}}}-\alpha _{B}^{{{P}_{2}}} \right)\left( \alpha _{C}^{{{P}_{3}}}-\alpha _{C}^{{{P}_{4}}} \right)-\left( \alpha _{B}^{{{P}_{3}}}-\alpha _{B}^{{{P}_{4}}} \right)\left( \alpha _{C}^{{{P}_{1}}}-\alpha _{C}^{{{P}_{2}}} \right)=\left| \begin{matrix}
		\alpha _{B}^{{{P}_{1}}}-\alpha _{B}^{{{P}_{2}}} & \alpha _{B}^{{{P}_{3}}}-\alpha _{B}^{{{P}_{4}}}  \\
		\alpha _{C}^{{{P}_{1}}}-\alpha _{C}^{{{P}_{2}}} & \alpha _{C}^{{{P}_{3}}}-\alpha _{C}^{{{P}_{4}}}  \\
	\end{matrix} \right|.\]
	
	So the frame component is obtained as:
	\[\begin{aligned}
		\alpha _{A}^{P}& =\left( 1-{{\kappa }_{{{P}_{1}}{{P}_{2}}}} \right)\alpha _{A}^{{{P}_{1}}}+{{\kappa }_{{{P}_{1}}{{P}_{2}}}}\alpha _{A}^{{{P}_{2}}} \\ 
		& =\frac{1}{\left| \begin{matrix}
				\alpha _{B}^{{{P}_{1}}}-\alpha _{B}^{{{P}_{2}}} & \alpha _{B}^{{{P}_{3}}}-\alpha _{B}^{{{P}_{4}}}  \\
				\alpha _{C}^{{{P}_{1}}}-\alpha _{C}^{{{P}_{2}}} & \alpha _{C}^{{{P}_{3}}}-\alpha _{C}^{{{P}_{4}}}  \\
			\end{matrix} \right|}\left( -\left| \begin{matrix}
			1 & 1 & 1  \\
			\alpha _{B}^{{{P}_{2}}} & \alpha _{B}^{{{P}_{3}}} & \alpha _{B}^{{{P}_{4}}}  \\
			\alpha _{C}^{{{P}_{2}}} & \alpha _{C}^{{{P}_{3}}} & \alpha _{C}^{{{P}_{4}}}  \\
		\end{matrix} \right|\alpha _{A}^{{{P}_{1}}}+\left| \begin{matrix}
			1 & 1 & 1  \\
			\alpha _{B}^{{{P}_{1}}} & \alpha _{B}^{{{P}_{3}}} & \alpha _{B}^{{{P}_{4}}}  \\
			\alpha _{C}^{{{P}_{1}}} & \alpha _{C}^{{{P}_{3}}} & \alpha _{C}^{{{P}_{4}}}  \\
		\end{matrix} \right|\alpha _{A}^{{{P}_{2}}} \right) \\ 
		& =\frac{1}{\left| \begin{matrix}
				\alpha _{B}^{{{P}_{1}}}-\alpha _{B}^{{{P}_{2}}} & \alpha _{B}^{{{P}_{3}}}-\alpha _{B}^{{{P}_{4}}}  \\
				\alpha _{C}^{{{P}_{1}}}-\alpha _{C}^{{{P}_{2}}} & \alpha _{C}^{{{P}_{3}}}-\alpha _{C}^{{{P}_{4}}}  \\
			\end{matrix} \right|}\left| \begin{matrix}
			-\alpha _{A}^{{{P}_{1}}} & -\alpha _{A}^{{{P}_{2}}} & 0 & 0  \\
			1 & 1 & 1 & 1  \\
			\alpha _{B}^{{{P}_{1}}} & \alpha _{B}^{{{P}_{2}}} & \alpha _{B}^{{{P}_{3}}} & \alpha _{B}^{{{P}_{4}}}  \\
			\alpha _{C}^{{{P}_{1}}} & \alpha _{C}^{{{P}_{2}}} & \alpha _{C}^{{{P}_{3}}} & \alpha _{C}^{{{P}_{4}}}  \\
		\end{matrix} \right|.  
	\end{aligned}\]
	
	Consider ($BC$ combination) ${{\kappa }_{{{P}_{3}}{{P}_{4}}}}$. The numerator of ${{\kappa }_{{{P}_{3}}{{P}_{4}}}}$ is
	\[\begin{aligned}
		& \left( \alpha _{B}^{{{P}_{1}}}-\alpha _{B}^{{{P}_{3}}} \right)\left( \alpha _{C}^{{{P}_{1}}}-\alpha _{C}^{{{P}_{2}}} \right)-\left( \alpha _{B}^{{{P}_{1}}}-\alpha _{B}^{{{P}_{2}}} \right)\left( \alpha _{C}^{{{P}_{1}}}-\alpha _{C}^{{{P}_{3}}} \right)=\left| \begin{matrix}
			\alpha _{B}^{{{P}_{1}}}-\alpha _{B}^{{{P}_{3}}} & \alpha _{B}^{{{P}_{1}}}-\alpha _{B}^{{{P}_{2}}}  \\
			\alpha _{C}^{{{P}_{1}}}-\alpha _{C}^{{{P}_{3}}} & \alpha _{C}^{{{P}_{1}}}-\alpha _{C}^{{{P}_{2}}}  \\
		\end{matrix} \right| \\ 
		& =\left| \begin{matrix}
			1 & 0 & 0  \\
			\alpha _{B}^{{{P}_{1}}} & \alpha _{B}^{{{P}_{1}}}-\alpha _{B}^{{{P}_{3}}} & \alpha _{B}^{{{P}_{1}}}-\alpha _{B}^{{{P}_{2}}}  \\
			\alpha _{C}^{{{P}_{1}}} & \alpha _{C}^{{{P}_{1}}}-\alpha _{C}^{{{P}_{3}}} & \alpha _{C}^{{{P}_{1}}}-\alpha _{C}^{{{P}_{2}}}  \\
		\end{matrix} \right|=\left| \begin{matrix}
			1 & -1 & -1  \\
			\alpha _{B}^{{{P}_{1}}} & -\alpha _{B}^{{{P}_{3}}} & -\alpha _{B}^{{{P}_{2}}}  \\
			\alpha _{C}^{{{P}_{1}}} & -\alpha _{C}^{{{P}_{3}}} & -\alpha _{C}^{{{P}_{2}}}  \\
		\end{matrix} \right|=-\left| \begin{matrix}
			1 & 1 & 1  \\
			\alpha _{B}^{{{P}_{1}}} & \alpha _{B}^{{{P}_{2}}} & \alpha _{B}^{{{P}_{3}}}  \\
			\alpha _{C}^{{{P}_{1}}} & \alpha _{C}^{{{P}_{2}}} & \alpha _{C}^{{{P}_{3}}}  \\
		\end{matrix} \right|. \\ 
	\end{aligned}\]
	
	The numerator of $1-{{\kappa }_{{{P}_{3}}{{P}_{4}}}}$ is
	\[\begin{aligned}
		& \left( \alpha _{B}^{{{P}_{4}}}-\alpha _{B}^{{{P}_{1}}} \right)\left( \alpha _{C}^{{{P}_{1}}}-\alpha _{C}^{{{P}_{2}}} \right)-\left( \alpha _{B}^{{{P}_{1}}}-\alpha _{B}^{{{P}_{2}}} \right)\left( \alpha _{C}^{{{P}_{4}}}-\alpha _{C}^{{{P}_{1}}} \right)=\left| \begin{matrix}
			\alpha _{B}^{{{P}_{4}}}-\alpha _{B}^{{{P}_{1}}} & \alpha _{B}^{{{P}_{1}}}-\alpha _{B}^{{{P}_{2}}}  \\
			\alpha _{C}^{{{P}_{4}}}-\alpha _{C}^{{{P}_{1}}} & \alpha _{C}^{{{P}_{1}}}-\alpha _{C}^{{{P}_{2}}}  \\
		\end{matrix} \right| \\ 
		& =\left| \begin{matrix}
			1 & 0 & 0  \\
			\alpha _{B}^{{{P}_{1}}} & \alpha _{B}^{{{P}_{4}}}-\alpha _{B}^{{{P}_{1}}} & \alpha _{B}^{{{P}_{1}}}-\alpha _{B}^{{{P}_{2}}}  \\
			\alpha _{C}^{{{P}_{1}}} & \alpha _{C}^{{{P}_{4}}}-\alpha _{C}^{{{P}_{1}}} & \alpha _{C}^{{{P}_{1}}}-\alpha _{C}^{{{P}_{2}}}  \\
		\end{matrix} \right|=\left| \begin{matrix}
			1 & 1 & -1  \\
			\alpha _{B}^{{{P}_{1}}} & \alpha _{B}^{{{P}_{4}}} & -\alpha _{B}^{{{P}_{2}}}  \\
			\alpha _{C}^{{{P}_{1}}} & \alpha _{C}^{{{P}_{4}}} & -\alpha _{C}^{{{P}_{2}}}  \\
		\end{matrix} \right|=\left| \begin{matrix}
			1 & 1 & 1  \\
			\alpha _{B}^{{{P}_{1}}} & \alpha _{B}^{{{P}_{2}}} & \alpha _{B}^{{{P}_{4}}}  \\
			\alpha _{C}^{{{P}_{1}}} & \alpha _{C}^{{{P}_{2}}} & \alpha _{C}^{{{P}_{4}}}  \\
		\end{matrix} \right|. \\ 
	\end{aligned}\]
	
	So the frame component is obtained as:
	\[\begin{aligned}
		\alpha _{A}^{P}& =\left( 1-{{\kappa }_{{{P}_{3}}{{P}_{4}}}} \right)\alpha _{A}^{{{P}_{3}}}+{{\kappa }_{{{P}_{3}}{{P}_{4}}}}\alpha _{A}^{{{P}_{4}}} \\ 
		& =\frac{1}{\left| \begin{matrix}
				\alpha _{B}^{{{P}_{1}}}-\alpha _{B}^{{{P}_{2}}} & \alpha _{B}^{{{P}_{3}}}-\alpha _{B}^{{{P}_{4}}}  \\
				\alpha _{C}^{{{P}_{1}}}-\alpha _{C}^{{{P}_{2}}} & \alpha _{C}^{{{P}_{3}}}-\alpha _{C}^{{{P}_{4}}}  \\
			\end{matrix} \right|}\left( \left| \begin{matrix}
			1 & 1 & 1  \\
			\alpha _{B}^{{{P}_{1}}} & \alpha _{B}^{{{P}_{2}}} & \alpha _{B}^{{{P}_{4}}}  \\
			\alpha _{C}^{{{P}_{1}}} & \alpha _{C}^{{{P}_{2}}} & \alpha _{C}^{{{P}_{4}}}  \\
		\end{matrix} \right|\alpha _{A}^{{{P}_{3}}}-\left| \begin{matrix}
			1 & 1 & 1  \\
			\alpha _{B}^{{{P}_{1}}} & \alpha _{B}^{{{P}_{2}}} & \alpha _{B}^{{{P}_{3}}}  \\
			\alpha _{C}^{{{P}_{1}}} & \alpha _{C}^{{{P}_{2}}} & \alpha _{C}^{{{P}_{3}}}  \\
		\end{matrix} \right|\alpha _{A}^{{{P}_{4}}} \right) \\ 
		& =\frac{1}{\left| \begin{matrix}
				\alpha _{B}^{{{P}_{1}}}-\alpha _{B}^{{{P}_{2}}} & \alpha _{B}^{{{P}_{3}}}-\alpha _{B}^{{{P}_{4}}}  \\
				\alpha _{C}^{{{P}_{1}}}-\alpha _{C}^{{{P}_{2}}} & \alpha _{C}^{{{P}_{3}}}-\alpha _{C}^{{{P}_{4}}}  \\
			\end{matrix} \right|}\left| \begin{matrix}
			0 & 0 & \alpha _{A}^{{{P}_{3}}} & \alpha _{A}^{{{P}_{4}}}  \\
			1 & 1 & 1 & 1  \\
			\alpha _{B}^{{{P}_{1}}} & \alpha _{B}^{{{P}_{2}}} & \alpha _{B}^{{{P}_{3}}} & \alpha _{B}^{{{P}_{4}}}  \\
			\alpha _{C}^{{{P}_{1}}} & \alpha _{C}^{{{P}_{2}}} & \alpha _{C}^{{{P}_{3}}} & \alpha _{C}^{{{P}_{4}}}  \\
		\end{matrix} \right|.  
	\end{aligned}\]
	
	Therefore
	\[\begin{aligned}
		\alpha _{A}^{P}& =\frac{1}{2}\left( \left( 1-{{\kappa }_{{{P}_{1}}{{P}_{2}}}} \right)\alpha _{A}^{{{P}_{1}}}+{{\kappa }_{{{P}_{1}}{{P}_{2}}}}\alpha _{A}^{{{P}_{2}}}+\left( 1-{{\kappa }_{{{P}_{3}}{{P}_{4}}}} \right)\alpha _{A}^{{{P}_{3}}}+{{\kappa }_{{{P}_{3}}{{P}_{4}}}}\alpha _{A}^{{{P}_{4}}} \right) \\ 
		& =\frac{1}{2\left| \begin{matrix}
				\alpha _{B}^{{{P}_{1}}}-\alpha _{B}^{{{P}_{2}}} & \alpha _{B}^{{{P}_{3}}}-\alpha _{B}^{{{P}_{4}}}  \\
				\alpha _{C}^{{{P}_{1}}}-\alpha _{C}^{{{P}_{2}}} & \alpha _{C}^{{{P}_{3}}}-\alpha _{C}^{{{P}_{4}}}  \\
			\end{matrix} \right|}\left| \begin{matrix}
			-\alpha _{A}^{{{P}_{1}}} & -\alpha _{A}^{{{P}_{2}}} & \alpha _{A}^{{{P}_{3}}} & \alpha _{A}^{{{P}_{4}}}  \\
			1 & 1 & 1 & 1  \\
			\alpha _{B}^{{{P}_{1}}} & \alpha _{B}^{{{P}_{2}}} & \alpha _{B}^{{{P}_{3}}} & \alpha _{B}^{{{P}_{4}}}  \\
			\alpha _{C}^{{{P}_{1}}} & \alpha _{C}^{{{P}_{2}}} & \alpha _{C}^{{{P}_{3}}} & \alpha _{C}^{{{P}_{4}}}  \\
		\end{matrix} \right| \\ 
		& =\frac{1}{2\left| \begin{matrix}
				\alpha _{B}^{{{P}_{1}}}-\alpha _{B}^{{{P}_{2}}} & \alpha _{B}^{{{P}_{3}}}-\alpha _{B}^{{{P}_{4}}}  \\
				\alpha _{C}^{{{P}_{1}}}-\alpha _{C}^{{{P}_{2}}} & \alpha _{C}^{{{P}_{3}}}-\alpha _{C}^{{{P}_{4}}}  \\
			\end{matrix} \right|}\left| \begin{matrix}
			1 & 1 & 1 & 1  \\
			\alpha _{A}^{{{P}_{1}}} & \alpha _{A}^{{{P}_{2}}} & -\alpha _{A}^{{{P}_{3}}} & -\alpha _{A}^{{{P}_{4}}}  \\
			\alpha _{B}^{{{P}_{1}}} & \alpha _{B}^{{{P}_{2}}} & \alpha _{B}^{{{P}_{3}}} & \alpha _{B}^{{{P}_{4}}}  \\
			\alpha _{C}^{{{P}_{1}}} & \alpha _{C}^{{{P}_{2}}} & \alpha _{C}^{{{P}_{3}}} & \alpha _{C}^{{{P}_{4}}}  \\
		\end{matrix} \right|.  
	\end{aligned}\]
	
	The following results were obtained through rotation:
	\[\alpha _{B}^{P}=\frac{1}{2\left| \begin{matrix}
			\alpha _{C}^{{{P}_{1}}}-\alpha _{C}^{{{P}_{2}}} & \alpha _{C}^{{{P}_{3}}}-\alpha _{C}^{{{P}_{4}}}  \\
			\alpha _{A}^{{{P}_{1}}}-\alpha _{A}^{{{P}_{2}}} & \alpha _{A}^{{{P}_{3}}}-\alpha _{A}^{{{P}_{4}}}  \\
		\end{matrix} \right|}\left| \begin{matrix}
		1 & 1 & 1 & 1  \\
		\alpha _{B}^{{{P}_{1}}} & \alpha _{B}^{{{P}_{2}}} & -\alpha _{B}^{{{P}_{3}}} & -\alpha _{B}^{{{P}_{4}}}  \\
		\alpha _{C}^{{{P}_{1}}} & \alpha _{C}^{{{P}_{2}}} & \alpha _{C}^{{{P}_{3}}} & \alpha _{C}^{{{P}_{4}}}  \\
		\alpha _{A}^{{{P}_{1}}} & \alpha _{A}^{{{P}_{2}}} & \alpha _{A}^{{{P}_{3}}} & \alpha _{A}^{{{P}_{4}}}  \\
	\end{matrix} \right|,\]
	
	\[\alpha _{C}^{P}=\frac{1}{2\left| \begin{matrix}
			\alpha _{A}^{{{P}_{1}}}-\alpha _{A}^{{{P}_{2}}} & \alpha _{A}^{{{P}_{3}}}-\alpha _{A}^{{{P}_{4}}}  \\
			\alpha _{B}^{{{P}_{1}}}-\alpha _{B}^{{{P}_{2}}} & \alpha _{B}^{{{P}_{3}}}-\alpha _{B}^{{{P}_{4}}}  \\
		\end{matrix} \right|}\left| \begin{matrix}
		1 & 1 & 1 & 1  \\
		\alpha _{C}^{{{P}_{1}}} & \alpha _{C}^{{{P}_{2}}} & -\alpha _{C}^{{{P}_{3}}} & -\alpha _{C}^{{{P}_{4}}}  \\
		\alpha _{A}^{{{P}_{1}}} & \alpha _{A}^{{{P}_{2}}} & \alpha _{A}^{{{P}_{3}}} & \alpha _{A}^{{{P}_{4}}}  \\
		\alpha _{B}^{{{P}_{1}}} & \alpha _{B}^{{{P}_{2}}} & \alpha _{B}^{{{P}_{3}}} & \alpha _{B}^{{{P}_{4}}}  \\
	\end{matrix} \right|.\]
\end{proof}
\hfill $\square$\par


\begin{example}{}\label{JiaodianBiaojiafenliangJuli}
	Given a $\triangle ABC$, the ccentroid of $\triangle ABC$ is $G$. The intersection point of two straight lines, $\overleftrightarrow{AG}$ and $\overleftrightarrow{BC}$, is $P:=\overleftrightarrow{AG}\cap \overleftrightarrow{BC}$. Find the frame components $\alpha _{A}^{P}$, $\alpha _{B}^{P}$, $\alpha _{C}^{P}$ of point $P$. 
\end{example}

\begin{solution}
	Let ${{P}_{1}}:=A$, ${{P}_{2}}:=G$, ${{P}_{3}}:=B$, ${{P}_{4}}:=C$,  According to theorem \ref{thm:LiangtiaoZhixianJiaodianDeBiaojiafenliang_Hanglieshixingshi}, we have:
	\[\begin{aligned}
		\alpha _{A}^{P}& =\frac{1}{2\left| \begin{matrix}
				\alpha _{B}^{{{P}_{1}}}-\alpha _{B}^{{{P}_{2}}} & \alpha _{B}^{{{P}_{3}}}-\alpha _{B}^{{{P}_{4}}}  \\
				\alpha _{C}^{{{P}_{1}}}-\alpha _{C}^{{{P}_{2}}} & \alpha _{C}^{{{P}_{3}}}-\alpha _{C}^{{{P}_{4}}}  \\
			\end{matrix} \right|}\left| \begin{matrix}
			1 & 1 & 1 & 1  \\
			\alpha _{A}^{{{P}_{1}}} & \alpha _{A}^{{{P}_{2}}} & -\alpha _{A}^{{{P}_{3}}} & -\alpha _{A}^{{{P}_{4}}}  \\
			\alpha _{B}^{{{P}_{1}}} & \alpha _{B}^{{{P}_{2}}} & \alpha _{B}^{{{P}_{3}}} & \alpha _{B}^{{{P}_{4}}}  \\
			\alpha _{C}^{{{P}_{1}}} & \alpha _{C}^{{{P}_{2}}} & \alpha _{C}^{{{P}_{3}}} & \alpha _{C}^{{{P}_{4}}}  \\
		\end{matrix} \right| \\ 
		& =\frac{1}{2\left| \begin{matrix}
				\alpha _{B}^{A}-\alpha _{B}^{G} & \alpha _{B}^{B}-\alpha _{B}^{C}  \\
				\alpha _{C}^{A}-\alpha _{C}^{G} & \alpha _{C}^{B}-\alpha _{C}^{C}  \\
			\end{matrix} \right|}\left| \begin{matrix}
			1 & 1 & 1 & 1  \\
			\alpha _{A}^{A} & \alpha _{A}^{G} & -\alpha _{A}^{B} & -\alpha _{A}^{C}  \\
			\alpha _{B}^{A} & \alpha _{B}^{G} & \alpha _{B}^{B} & \alpha _{B}^{C}  \\
			\alpha _{C}^{A} & \alpha _{C}^{G} & \alpha _{C}^{B} & \alpha _{C}^{C}  \\
		\end{matrix} \right| \\ 
		& =\frac{1}{2\left| \begin{matrix}
				0-\frac{1}{3} & 1-0  \\
				0-\frac{1}{3} & 0-1  \\
			\end{matrix} \right|}\left| \begin{matrix}
			1 & 1 & 1 & 1  \\
			1 & \frac{1}{3} & 0 & 0  \\
			0 & \frac{1}{3} & 1 & 0  \\
			0 & \frac{1}{3} & 0 & 1  \\
		\end{matrix} \right|=\frac{3}{4}\left| \begin{matrix}
			1 & 1 & 1 & 1  \\
			0 & -\frac{2}{3} & -1 & -1  \\
			0 & \frac{1}{3} & 1 & 0  \\
			0 & \frac{1}{3} & 0 & 1  \\
		\end{matrix} \right|=\frac{3}{4}\left| \begin{matrix}
			1 & 1 & 1 & 1  \\
			0 & -\frac{2}{3} & -1 & -1  \\
			0 & \frac{1}{3} & 1 & 0  \\
			0 & \frac{1}{3} & 0 & 1  \\
		\end{matrix} \right|=0,  
	\end{aligned}\]
	\[\begin{aligned}
		\alpha _{B}^{P}& =\frac{1}{2\left| \begin{matrix}
				\alpha _{C}^{{{P}_{1}}}-\alpha _{C}^{{{P}_{2}}} & \alpha _{C}^{{{P}_{3}}}-\alpha _{C}^{{{P}_{4}}}  \\
				\alpha _{A}^{{{P}_{1}}}-\alpha _{A}^{{{P}_{2}}} & \alpha _{A}^{{{P}_{3}}}-\alpha _{A}^{{{P}_{4}}}  \\
			\end{matrix} \right|}\left| \begin{matrix}
			1 & 1 & 1 & 1  \\
			\alpha _{B}^{{{P}_{1}}} & \alpha _{B}^{{{P}_{2}}} & -\alpha _{B}^{{{P}_{3}}} & -\alpha _{B}^{{{P}_{4}}}  \\
			\alpha _{C}^{{{P}_{1}}} & \alpha _{C}^{{{P}_{2}}} & \alpha _{C}^{{{P}_{3}}} & \alpha _{C}^{{{P}_{4}}}  \\
			\alpha _{A}^{{{P}_{1}}} & \alpha _{A}^{{{P}_{2}}} & \alpha _{A}^{{{P}_{3}}} & \alpha _{A}^{{{P}_{4}}}  \\
		\end{matrix} \right| \\ 
		& =\frac{1}{2\left| \begin{matrix}
				\alpha _{C}^{A}-\alpha _{C}^{G} & \alpha _{C}^{B}-\alpha _{C}^{C}  \\
				\alpha _{A}^{A}-\alpha _{A}^{G} & \alpha _{A}^{B}-\alpha _{A}^{C}  \\
			\end{matrix} \right|}\left| \begin{matrix}
			1 & 1 & 1 & 1  \\
			\alpha _{B}^{A} & \alpha _{B}^{G} & -\alpha _{B}^{B} & -\alpha _{B}^{C}  \\
			\alpha _{C}^{A} & \alpha _{C}^{G} & \alpha _{C}^{B} & \alpha _{C}^{C}  \\
			\alpha _{A}^{A} & \alpha _{A}^{G} & \alpha _{A}^{B} & \alpha _{A}^{C}  \\
		\end{matrix} \right| \\ 
		& =\frac{1}{2\left| \begin{matrix}
				0-\frac{1}{3} & 0-1  \\
				1-\frac{1}{3} & 0-0  \\
			\end{matrix} \right|}\left| \begin{matrix}
			1 & 1 & 1 & 1  \\
			0 & \frac{1}{3} & -1 & 0  \\
			0 & \frac{1}{3} & 0 & 1  \\
			1 & \frac{1}{3} & 0 & 0  \\
		\end{matrix} \right|=\frac{3}{4}\left| \begin{matrix}
			1 & 1 & 1 & 1  \\
			0 & \frac{1}{3} & -1 & 0  \\
			0 & \frac{1}{3} & 0 & 1  \\
			0 & -\frac{2}{3} & -1 & -1  \\
		\end{matrix} \right|=\frac{3}{4}\left| \begin{matrix}
			\frac{1}{3} & -1 & 0  \\
			\frac{1}{3} & 0 & 1  \\
			-\frac{2}{3} & -1 & -1  \\
		\end{matrix} \right| \\ 
		& =\frac{3}{4}\left( \frac{2}{3}-\frac{1}{3}+\frac{1}{3} \right)=\frac{1}{2},  
	\end{aligned}\]
	
	Similarly, it can be obtained that:
	\[\begin{aligned}
		\alpha _{C}^{P}& =\frac{1}{2\left| \begin{matrix}
				\alpha _{A}^{{{P}_{1}}}-\alpha _{A}^{{{P}_{2}}} & \alpha _{A}^{{{P}_{3}}}-\alpha _{A}^{{{P}_{4}}}  \\
				\alpha _{B}^{{{P}_{1}}}-\alpha _{B}^{{{P}_{2}}} & \alpha _{B}^{{{P}_{3}}}-\alpha _{B}^{{{P}_{4}}}  \\
			\end{matrix} \right|}\left| \begin{matrix}
			1 & 1 & 1 & 1  \\
			\alpha _{C}^{{{P}_{1}}} & \alpha _{C}^{{{P}_{2}}} & -\alpha _{C}^{{{P}_{3}}} & -\alpha _{C}^{{{P}_{4}}}  \\
			\alpha _{A}^{{{P}_{1}}} & \alpha _{A}^{{{P}_{2}}} & \alpha _{A}^{{{P}_{3}}} & \alpha _{A}^{{{P}_{4}}}  \\
			\alpha _{B}^{{{P}_{1}}} & \alpha _{B}^{{{P}_{2}}} & \alpha _{B}^{{{P}_{3}}} & \alpha _{B}^{{{P}_{4}}}  \\
		\end{matrix} \right| \\ 
		& =\frac{1}{2\left| \begin{matrix}
				\alpha _{A}^{A}-\alpha _{A}^{G} & \alpha _{A}^{B}-\alpha _{A}^{C}  \\
				\alpha _{B}^{A}-\alpha _{B}^{G} & \alpha _{B}^{B}-\alpha _{B}^{C}  \\
			\end{matrix} \right|}\left| \begin{matrix}
			1 & 1 & 1 & 1  \\
			\alpha _{C}^{A} & \alpha _{C}^{G} & -\alpha _{C}^{B} & -\alpha _{C}^{C}  \\
			\alpha _{A}^{A} & \alpha _{A}^{G} & \alpha _{A}^{B} & \alpha _{A}^{C}  \\
			\alpha _{B}^{A} & \alpha _{B}^{G} & \alpha _{B}^{B} & \alpha _{B}^{C}  \\
		\end{matrix} \right| \\ 
		& =\frac{1}{2\left| \begin{matrix}
				1-\frac{1}{3} & 0-0  \\
				0-\frac{1}{3} & 1-0  \\
			\end{matrix} \right|}\left| \begin{matrix}
			1 & 1 & 1 & 1  \\
			0 & \frac{1}{3} & 0 & -1  \\
			1 & \frac{1}{3} & 0 & 0  \\
			0 & \frac{1}{3} & 1 & 0  \\
		\end{matrix} \right|=\frac{3}{4}\left| \begin{matrix}
			1 & 1 & 1 & 1  \\
			0 & \frac{1}{3} & 0 & -1  \\
			0 & -\frac{2}{3} & -1 & -1  \\
			0 & \frac{1}{3} & 1 & 0  \\
		\end{matrix} \right|=\frac{3}{4}\left| \begin{matrix}
			\frac{1}{3} & 0 & -1  \\
			-\frac{2}{3} & -1 & -1  \\
			\frac{1}{3} & 1 & 0  \\
		\end{matrix} \right| \\ 
		& =\frac{3}{4}\left( \frac{2}{3}-\frac{1}{3}+\frac{1}{3} \right)=\frac{1}{2}.  
	\end{aligned}\]
\end{solution}
\hfill $\diamond$\par


\begin{example}{}\label{DinghuixianJiaodianDeBiaojiafenliang}
	Given a $\triangle ABC$, the centroid of $\triangle ABC$ is $G$, and the incenter is $I$. The intersection point of two straight lines $\overleftrightarrow{AG}$ and $\overleftrightarrow{CI}$ is $P$, that is, $P:=\overleftrightarrow{AG}\cap \overleftrightarrow{CI}$, Find the frame component $\alpha _{A}^{P}$, $\alpha _{B}^{P}$, $\alpha _{C}^{P}$ of point $P$.
\end{example}

\begin{solution}
	Let ${{P}_{1}}:=A$, ${{P}_{2}}:=G$, ${{P}_{3}}:=C$, ${{P}_{4}}:=I$, according to theorem \ref{thm:LiangtiaoZhixianJiaodianDeBiaojiafenliang_Jianyuexingshi}, it is obtained that:
	\begin{align*}
		\alpha _{A}^{P}& =\frac{\left( \alpha _{B}^{{{P}_{3}}}-\alpha _{B}^{{{P}_{2}}} \right)\left( \alpha _{C}^{{{P}_{3}}}-\alpha _{C}^{{{P}_{4}}} \right)-\left( \alpha _{B}^{{{P}_{3}}}-\alpha _{B}^{{{P}_{4}}} \right)\left( \alpha _{C}^{{{P}_{3}}}-\alpha _{C}^{{{P}_{2}}} \right)}{\left( \alpha _{B}^{{{P}_{1}}}-\alpha _{B}^{{{P}_{2}}} \right)\left( \alpha _{C}^{{{P}_{3}}}-\alpha _{C}^{{{P}_{4}}} \right)-\left( \alpha _{B}^{{{P}_{3}}}-\alpha _{B}^{{{P}_{4}}} \right)\left( \alpha _{C}^{{{P}_{1}}}-\alpha _{C}^{{{P}_{2}}} \right)}\alpha _{A}^{{{P}_{1}}} \\ 
		& +\frac{\left( \alpha _{B}^{{{P}_{1}}}-\alpha _{B}^{{{P}_{3}}} \right)\left( \alpha _{C}^{{{P}_{3}}}-\alpha _{C}^{{{P}_{4}}} \right)-\left( \alpha _{B}^{{{P}_{3}}}-\alpha _{B}^{{{P}_{4}}} \right)\left( \alpha _{C}^{{{P}_{1}}}-\alpha _{C}^{{{P}_{3}}} \right)}{\left( \alpha _{B}^{{{P}_{1}}}-\alpha _{B}^{{{P}_{2}}} \right)\left( \alpha _{C}^{{{P}_{3}}}-\alpha _{C}^{{{P}_{4}}} \right)-\left( \alpha _{B}^{{{P}_{3}}}-\alpha _{B}^{{{P}_{4}}} \right)\left( \alpha _{C}^{{{P}_{1}}}-\alpha _{C}^{{{P}_{2}}} \right)}\alpha _{A}^{{{P}_{2}}} \\ 
		& =\frac{\left( \alpha _{B}^{C}-\alpha _{B}^{G} \right)\left( \alpha _{C}^{C}-\alpha _{C}^{I} \right)-\left( \alpha _{B}^{C}-\alpha _{B}^{I} \right)\left( \alpha _{C}^{C}-\alpha _{C}^{G} \right)}{\left( \alpha _{B}^{A}-\alpha _{B}^{G} \right)\left( \alpha _{C}^{C}-\alpha _{C}^{I} \right)-\left( \alpha _{B}^{C}-\alpha _{B}^{I} \right)\left( \alpha _{C}^{A}-\alpha _{C}^{G} \right)}\alpha _{A}^{A} \\ 
		& +\frac{\left( \alpha _{B}^{A}-\alpha _{B}^{C} \right)\left( \alpha _{C}^{C}-\alpha _{C}^{I} \right)-\left( \alpha _{B}^{C}-\alpha _{B}^{I} \right)\left( \alpha _{C}^{A}-\alpha _{C}^{C} \right)}{\left( \alpha _{B}^{A}-\alpha _{B}^{G} \right)\left( \alpha _{C}^{C}-\alpha _{C}^{I} \right)-\left( \alpha _{B}^{C}-\alpha _{B}^{I} \right)\left( \alpha _{C}^{A}-\alpha _{C}^{G} \right)}\alpha _{A}^{G} \\ 
		& =\frac{\left( 0-\frac{1}{3} \right)\left( 1-\frac{c}{2p} \right)-\left( 0-\frac{b}{2p} \right)\left( 1-\frac{1}{3} \right)}{\left( 0-\frac{1}{3} \right)\left( 1-\frac{c}{2p} \right)-\left( 0-\frac{b}{2p} \right)\left( 0-\frac{1}{3} \right)} \\ 
		& +\frac{\left( 0-0 \right)\left( 1-\frac{c}{2p} \right)-\left( 0-\frac{b}{2p} \right)\left( 0-1 \right)}{\left( 0-\frac{1}{3} \right)\left( 1-\frac{c}{2p} \right)-\left( 0-\frac{b}{2p} \right)\left( 0-\frac{1}{3} \right)}\cdot \frac{1}{3} \\ 
		& =\frac{-\frac{2p-c}{6p}+\frac{2b}{6p}}{-\frac{2p-c}{6p}-\frac{b}{6p}}+\frac{-\frac{b}{2p}}{-\frac{2p-c}{6p}-\frac{b}{6p}}\cdot \frac{1}{3} \\ 
		& =\frac{-\left( 2p-c \right)+2b}{-\left( 2p-c+b \right)}+\frac{-3b}{-\left( 2p-c+b \right)}\cdot \frac{1}{3} \\ 
		& =\frac{\left( 2p-c \right)-2b+3b}{2p-c+b}=\frac{a}{a+2b},  
	\end{align*}
	\begin{align*}
		\alpha _{B}^{P}& =\frac{\left( \alpha _{C}^{{{P}_{1}}}-\alpha _{C}^{{{P}_{3}}} \right)\left( \alpha _{A}^{{{P}_{3}}}-\alpha _{A}^{{{P}_{4}}} \right)-\left( \alpha _{C}^{{{P}_{3}}}-\alpha _{C}^{{{P}_{4}}} \right)\left( \alpha _{A}^{{{P}_{1}}}-\alpha _{A}^{{{P}_{3}}} \right)}{\left( \alpha _{C}^{{{P}_{1}}}-\alpha _{C}^{{{P}_{2}}} \right)\left( \alpha _{A}^{{{P}_{3}}}-\alpha _{A}^{{{P}_{4}}} \right)-\left( \alpha _{C}^{{{P}_{3}}}-\alpha _{C}^{{{P}_{4}}} \right)\left( \alpha _{A}^{{{P}_{1}}}-\alpha _{A}^{{{P}_{2}}} \right)}\alpha _{B}^{{{P}_{1}}} \\ 
		& +\frac{\left( \alpha _{A}^{{{P}_{1}}}-\alpha _{A}^{{{P}_{3}}} \right)\left( \alpha _{B}^{{{P}_{3}}}-\alpha _{B}^{{{P}_{4}}} \right)-\left( \alpha _{A}^{{{P}_{3}}}-\alpha _{A}^{{{P}_{4}}} \right)\left( \alpha _{B}^{{{P}_{1}}}-\alpha _{B}^{{{P}_{3}}} \right)}{\left( \alpha _{A}^{{{P}_{1}}}-\alpha _{A}^{{{P}_{2}}} \right)\left( \alpha _{B}^{{{P}_{3}}}-\alpha _{B}^{{{P}_{4}}} \right)-\left( \alpha _{A}^{{{P}_{3}}}-\alpha _{A}^{{{P}_{4}}} \right)\left( \alpha _{B}^{{{P}_{1}}}-\alpha _{B}^{{{P}_{2}}} \right)}\alpha _{B}^{{{P}_{2}}} \\ 
		& =\frac{\left( \alpha _{C}^{A}-\alpha _{C}^{C} \right)\left( \alpha _{A}^{C}-\alpha _{A}^{I} \right)-\left( \alpha _{C}^{C}-\alpha _{C}^{I} \right)\left( \alpha _{A}^{A}-\alpha _{A}^{C} \right)}{\left( \alpha _{C}^{A}-\alpha _{C}^{G} \right)\left( \alpha _{A}^{C}-\alpha _{A}^{I} \right)-\left( \alpha _{C}^{C}-\alpha _{C}^{I} \right)\left( \alpha _{A}^{A}-\alpha _{A}^{G} \right)}\alpha _{B}^{A} \\ 
		& +\frac{\left( \alpha _{A}^{A}-\alpha _{A}^{C} \right)\left( \alpha _{B}^{C}-\alpha _{B}^{I} \right)-\left( \alpha _{A}^{C}-\alpha _{A}^{I} \right)\left( \alpha _{B}^{A}-\alpha _{B}^{C} \right)}{\left( \alpha _{A}^{A}-\alpha _{A}^{G} \right)\left( \alpha _{B}^{C}-\alpha _{B}^{I} \right)-\left( \alpha _{A}^{C}-\alpha _{A}^{I} \right)\left( \alpha _{B}^{A}-\alpha _{B}^{G} \right)}\alpha _{B}^{G} \\ 
		& =\frac{\left( 1-0 \right)\left( 0-\frac{b}{2p} \right)-\left( 0-\frac{a}{2p} \right)\left( 0-0 \right)}{\left( 1-\frac{1}{3} \right)\left( 0-\frac{b}{2p} \right)-\left( 0-\frac{a}{2p} \right)\left( 0-\frac{1}{3} \right)}\cdot \frac{1}{3} \\ 
		& =\frac{-\frac{b}{2p}}{-\frac{2b}{6p}-\frac{a}{6p}}\cdot \frac{1}{3}=\frac{3b}{2b+a}\cdot \frac{1}{3}=\frac{b}{a+2b},  
	\end{align*}
	\begin{align*}
		\alpha _{C}^{P}& =\frac{\left( \alpha _{A}^{{{P}_{3}}}-\alpha _{A}^{{{P}_{2}}} \right)\left( \alpha _{B}^{{{P}_{3}}}-\alpha _{B}^{{{P}_{4}}} \right)-\left( \alpha _{A}^{{{P}_{3}}}-\alpha _{A}^{{{P}_{4}}} \right)\left( \alpha _{B}^{{{P}_{3}}}-\alpha _{B}^{{{P}_{2}}} \right)}{\left( \alpha _{A}^{{{P}_{1}}}-\alpha _{A}^{{{P}_{2}}} \right)\left( \alpha _{B}^{{{P}_{3}}}-\alpha _{B}^{{{P}_{4}}} \right)-\left( \alpha _{A}^{{{P}_{3}}}-\alpha _{A}^{{{P}_{4}}} \right)\left( \alpha _{B}^{{{P}_{1}}}-\alpha _{B}^{{{P}_{2}}} \right)}\alpha _{C}^{{{P}_{1}}} \\ 
		& +\frac{\left( \alpha _{A}^{{{P}_{1}}}-\alpha _{A}^{{{P}_{3}}} \right)\left( \alpha _{B}^{{{P}_{3}}}-\alpha _{B}^{{{P}_{4}}} \right)-\left( \alpha _{A}^{{{P}_{3}}}-\alpha _{A}^{{{P}_{4}}} \right)\left( \alpha _{B}^{{{P}_{1}}}-\alpha _{B}^{{{P}_{3}}} \right)}{\left( \alpha _{A}^{{{P}_{1}}}-\alpha _{A}^{{{P}_{2}}} \right)\left( \alpha _{B}^{{{P}_{3}}}-\alpha _{B}^{{{P}_{4}}} \right)-\left( \alpha _{A}^{{{P}_{3}}}-\alpha _{A}^{{{P}_{4}}} \right)\left( \alpha _{B}^{{{P}_{1}}}-\alpha _{B}^{{{P}_{2}}} \right)}\alpha _{C}^{{{P}_{2}}} \\ 
		& =\frac{\left( \alpha _{A}^{A}-\alpha _{A}^{C} \right)\left( \alpha _{B}^{C}-\alpha _{B}^{I} \right)-\left( \alpha _{A}^{C}-\alpha _{A}^{I} \right)\left( \alpha _{B}^{A}-\alpha _{B}^{C} \right)}{\left( \alpha _{A}^{A}-\alpha _{A}^{G} \right)\left( \alpha _{B}^{C}-\alpha _{B}^{I} \right)-\left( \alpha _{A}^{C}-\alpha _{A}^{I} \right)\left( \alpha _{B}^{A}-\alpha _{B}^{G} \right)}\alpha _{C}^{G} \\ 
		& =\frac{\left( 1-0 \right)\left( 0-\frac{b}{2p} \right)-\left( 0-\frac{a}{2p} \right)\left( 0-0 \right)}{\left( 1-\frac{1}{3} \right)\left( 0-\frac{b}{2p} \right)-\left( 0-\frac{a}{2p} \right)\left( 0-\frac{1}{3} \right)}\cdot \frac{1}{3} \\ 
		& =\frac{-\frac{b}{2p}}{-\frac{2b}{6p}-\frac{a}{6p}}\cdot \frac{1}{3}=\frac{3b}{2b+a}\cdot \frac{1}{3}=\frac{b}{a+2b}.  
	\end{align*}
\end{solution}
\hfill $\diamond$\par

\section{The Edge-Axis coordinate system of a triangle}\label{SanjiaoxingBianzhouzuobiaoxi}
In the following chapters, it can be seen that as long as the frame components of two points are given, the distance between these two points can be calculated, and the distance is represented by the lengths of the three sides of the triangle. Therefore, in Intercenter Geometry, solving the frame components is very important. We have already obtained the frame components for some special points (such as centroid, incenter, orthocenter, circumcenter, etc.), but knowing only the frame components for a few special points is far from enough in application.

Looking back at the frame components of those special points that were previously obtained, they were basically obtained using the definition of IR and Euclidean geometry methods. However, using Euclidean geometry methods to solve frame components is not universal, which means that when solving problems in Euclidean geometry, it is usually a special solution method for each problem, lacking a unified approach and solution. This personalized, one problem-one strategy solution is difficult to use the information of the frame components that already calculated to seek the frame components for new points. Or, from the perspective of Euclidean geometry, every new point is a completely new problem that is difficult to solve with the help of historical solutions. Therefore, it is difficult to solve the frame components of more points using Euclidean geometry methods.

As discussed in the previous chapters, if the frame components of two points are known, the frame components of any point on the line connecting these two points can be obtained. This is a great progress, but it is still not enough. For example, on the plane where a triangle is located, rotating the incenter of the triangle counterclockwise around one of its vertices by $30{}^\circ $ yields a new point. How to calculate the frame component of this new point? If we adopt the definition of IR and use Euclidean geometry methods to calculate the frame component, we will definitely encounter huge challenges. Obviously, such points usually do not happen to be exactly on the line connecting two known points (the frame components of these two points are known). So, it is urgent to find a universal method for solving the frame components of any point.

In order to obtain more frame components, a new method is needed. Since the frame components of some special points have been obtained, can we use the frame components of known points to determine the frame components of unknown points? I have come to a positive answer through research.

In order to obtain the frame component of another unknown point from the frame component of a known point, it is necessary to know the positional relationship between the known point and the unknown point. This requires a certain coordinate system, which is the Edge-Axis coordinate system to be discussed below.

I will propose the Edge-Axis coordinate system in this section. It should be pointed out that the frame components are independent of the coordinate system. The advantage of introducing a Edge-Axis coordinate system to solve the frame component is that: 1. By using the Edge-Axis coordinate system, it is convenient to obtain the frame component of another unknown point from the frame component of a known point 2. For the same point, no matter how the coordinate system is selected, the calculated frame components are the same. Therefore, according to the needs of the problem, the most convenient coordinate system can be used, and multiple different coordinate systems can be selected simultaneously for the same problem; 3. The coordinate system is just an intermediate medium, and once the frame components are calculated, the coordinate system can be completely eliminated. The position relationship between known and unknown points will be reflected as parameters in the final formula. For many practical problems, these parameters are constants. In this sense, the final frame component obtained is still a function of the lengths of the three sides of the triangle.

Before introducing the Edge-Axis coordinate system, several conventions need to be introduced.

\textbf {Counterclockwise Arrangement Convention}: When looking at a plane, the three vertices $A$, $B$, and $C$ of $\triangle ABC $ are arranged in counterclockwise order. This arrangement is called $ABC$ arrangement, or in other words, $ABC$ is \textbf{counterclockwise arrangement}.

Obviously, when $ABC$ is arranged counterclockwise, both $BCA$ and $CAB$ are also arranged counterclockwise. Given a $\triangle ABC$, if $ABC$ is arranged counterclockwise, let's assume that a person stands at the vertex $A$ of $\triangle ABC$ and walks along the segment $AB$ towards the vertex $B$, then stands at the vertex $B$ and walks along the segment $BC$ towards the vertex $C$, and then stands at the vertex $C$ and walks along the segment $CA$ towards the vertex $A$, then the interior of the triangle is always on the left side of that person.

\textbf {Directed Angle Convention}: In $\triangle ABC$, $\angle BAC$, $\angle CBA$, $\angle ACB$ are all directed angles, and $A:=\angle BAC$, $B:=\angle CBA$, $C:=\angle ACB$.

The concept of directed angle can be found in the appendix \ref{SanjiaoxingXianguanZhishi}. For the sake of symbol simplicity, $A$($B$, $C$) represents both the vertex $A$($B$, $C$) of $\triangle ABC$ and the angle $A$($B$, $C$) of $\triangle ABC$, which readers can distinguish them based on the context.

\textbf {Vertex-Edge Symbol Convention}: In $\triangle ABC$, the length of the opposite edge of vertex $A$ is always $a$, that is, $a:=BC$; The length of the opposite side of vertex $B$ is always $b$, that is, $b:=CA$; The length of the opposite side of vertex $C$ is always $c$, that is, $c:=AB$.

For the vertices of a triangle, people can freely give labels, such as $A$, $X$, and so on. But if a triangle $\triangle ABC$ is said, it means that the labels of the three vertices of the triangle has been determined, namely $A$, $B$, and $C$, and the length of the opposite side of vertex $A$ can only be $a$, and no other notation can be used. The length of the opposite sides of vertex $B$ and vertex $C$ can only be $b$ and $c$, respectively.



In this book, when using the Edge-Axis coordinate system to analyze problems of a triangle, it is shown that the triangle satisfies the above conventions.

Unless otherwise specified, this book follows the above conventions.

\begin{figure}[h]
	\centering
	\includegraphics[totalheight=6cm]{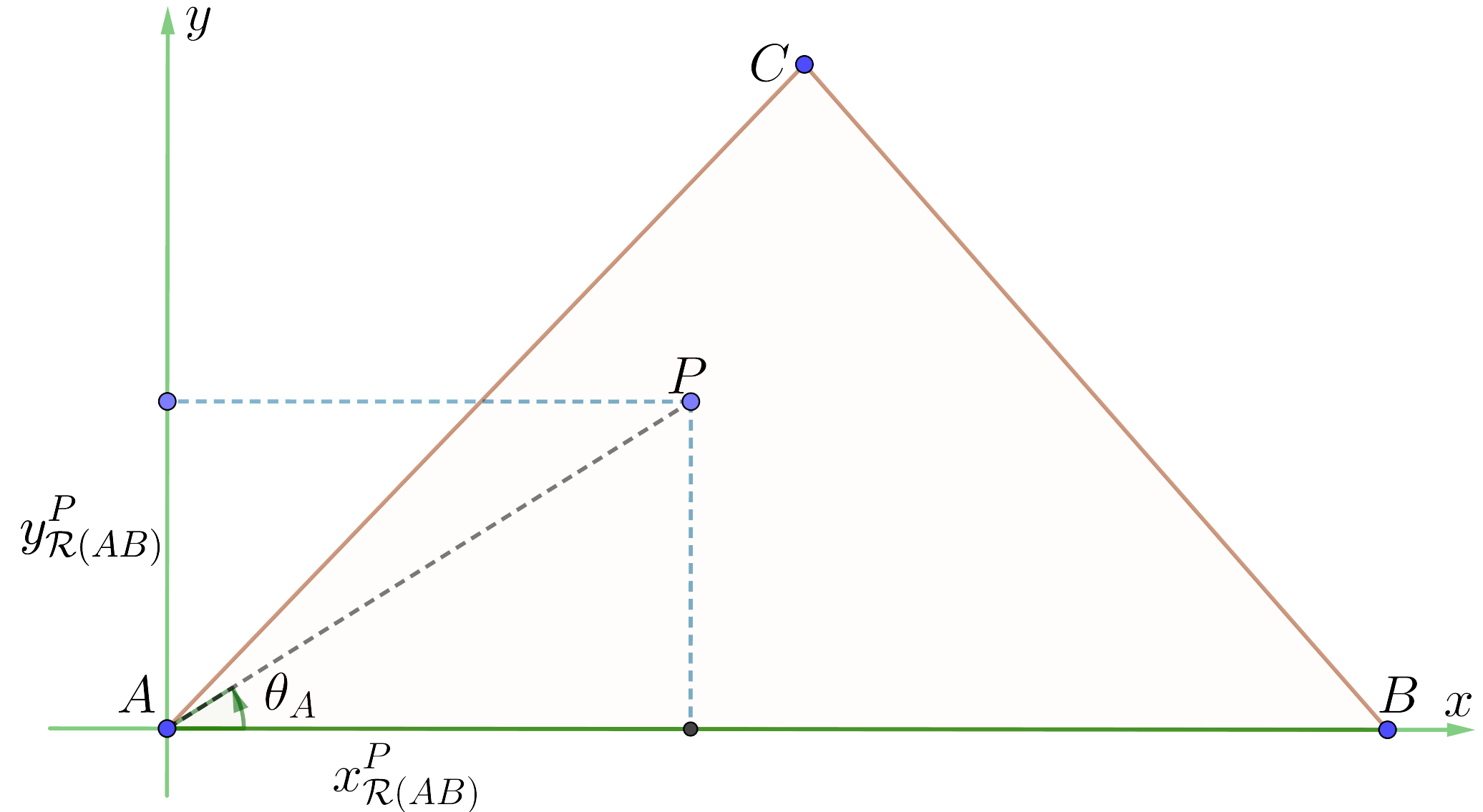}
	\caption{Coordinate system $\mathcal{R}\left( AB \right)$ and $\mathcal{P}\left( AB \right)$ of vertex $A$}
	\label{fig:Bianzhouzuobiaoxi_A}
\end{figure}

\begin{figure}[h]
	\centering
	\includegraphics[totalheight=8cm]{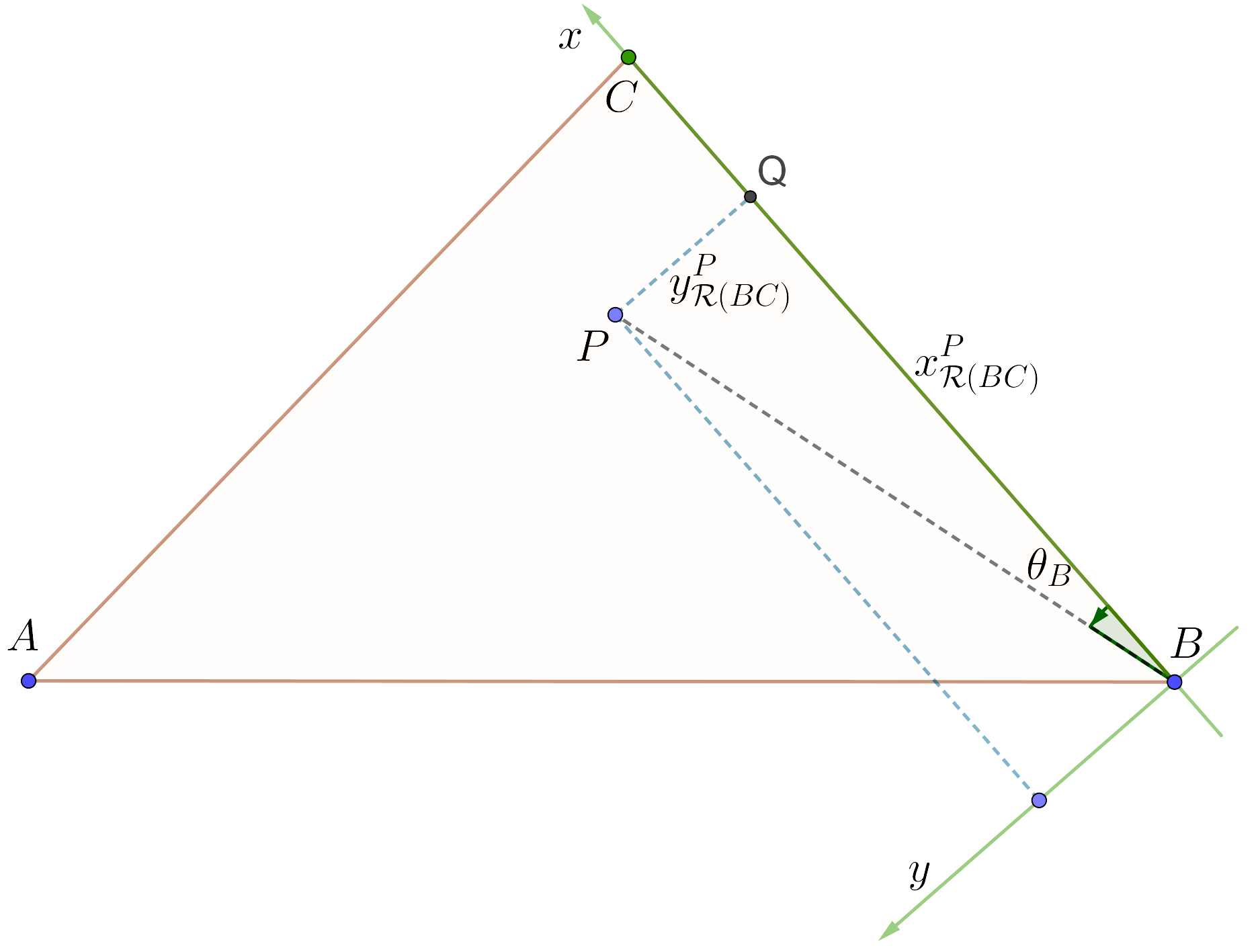}
	\caption{Coordinate system $\mathcal{R}\left( BC \right)$ and $\mathcal{P}\left( BC \right)$ of vertex $B$} \label{fig:Bianzhouzuobiaoxi_B}
\end{figure}

\begin{figure}[h]
	\centering
	\includegraphics[totalheight=6cm]{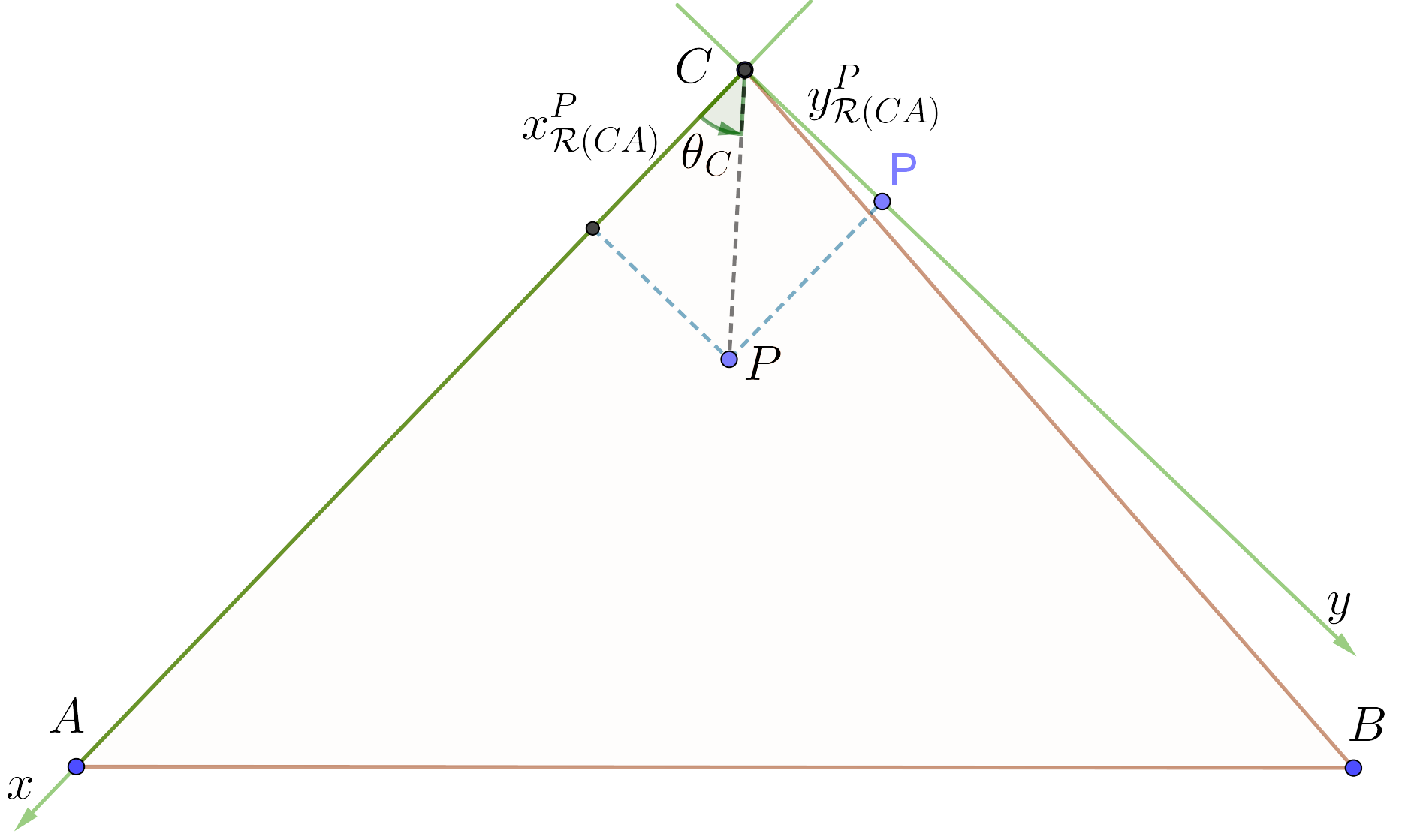}
	\caption{Coordinate system $\mathcal{R}\left( CA \right)$ and $\mathcal{P}\left( CA \right)$ of vertex $C$} \label{fig:Bianzhouzuobiaoxi_C}
\end{figure}


\textbf {Edge-Axis coordinate system} is shown in figure \ref{fig:Bianzhouzuobiaoxi_A}, Figure \ref{fig:Bianzhouzuobiaoxi_B}, and Figure \ref{fig:Bianzhouzuobiaoxi_C}. Taking the example of Figure \ref{fig:Bianzhouzuobiaoxi_A}, we will provide an explanation.


Figure \ref{fig:Bianzhouzuobiaoxi_A} establishes a Cartesian coordinate system (abbreviated as Cartesian coordinate system) and a polar coordinate system. First, let's take a look at the Cartesian coordinate system. The origin of this Cartesian coordinate system is at the vertex $A$ of the $\triangle ABC$, with a horizontal axis of $\overleftrightarrow{AB}$. The positive direction of the horizontal axis of $\overleftrightarrow{AB}$ is the same as the positive direction of the vector $\overrightarrow{AB}$, and such a horizontal axis is also abbreviated as $Ax$. The positive direction of the vertical axis is perpendicular to the positive direction of the vector $\overrightarrow{AB}$, it is in the same direction as the vector $\overrightarrow {AB}$ after rotating $90{}^\circ$ counterclockwise around vertex $A$. This vertical axis is also abbreviated as $Ay$. The convention $\angle xAy$ is a directed angle, and its reference positive direction is counterclockwise. I refer to this coordinate system as the \textbf{Cartesian coordinate system with origin $A$}, abbreviated as the \textbf{Cartesian coordinate system} without causing confusion, denoted as $\mathcal{R}\left(AB\right)$.


In the Cartesian coordinate system $\mathcal{R}\left( AB \right)$, the Cartesian coordinates of point $P$ are denoted as $P\left( x_{\mathcal{R}\left( AB \right)}^{P},y_{\mathcal{R}\left( AB \right)}^{P} \right)$. Similarly, in the Cartesian coordinate system $\mathcal{R}\left( BC \right)$ and $\mathcal{R}\left( CA \right)$, the Cartesian coordinates of point $P$ are denoted as $P\left( x_{\mathcal{R}\left( BC \right)}^{P},y_{\mathcal{R}\left( BC \right)}^{P} \right)$ and $P\left( x_{\mathcal{R}\left( CA \right)}^{P},y_{\mathcal{R}\left( CA \right)}^{P} \right)$, respectively.


The Cartesian coordinate systems of origin $B$ and $C$ are shown in figure \ref{fig:Bianzhouzuobiaoxi_B} and \ref{fig:Bianzhouzuobiaoxi_C}, respectively.

Obviously, any point on the plane corresponds one-to-one to the coordinates in a Cartesian coordinate system, which is a pair of real numbers.


Next, let's take a look at the polar coordinate system. The pole of a polar coordinate system are at vertex $A$, the polar axis is $\overline{AB}$, and the polar angle is defined as $\theta _{\mathcal{R}\left( AB \right)}^{P}:=\angle BAP$. Without causing confusion, the polar angle is also abbreviated as $\theta _A$. The polar angle $\theta _A$ is a directed angle, and its reference positive direction is counterclockwise. I refer to this coordinate system as the \textbf{polar coordinate system with pole of $A$}, abbreviated as the \textbf {polar coordinate system} without causing confusion, denoted as $\mathcal{P}\left(AB\right)$.

In the polar coordinate system $\mathcal{P}\left(AB \right)$, the polar coordinates of point $P$ are denoted as $P\left(AP,\theta _A \right)$. 


Similarly, in the polar coordinate system $\mathcal{P}\left( BC \right)$ and $\mathcal{P}\left( CA \right)$, the polar coordinates of point $P$ are denoted as $P\left( BP,\theta _B \right)$ and $P\left( CP,\theta _C \right)$, respectively. Polar coordinates $P\left( AP,\theta _A \right)$, $P\left( BP,\theta _B \right)$ and $P\left( CP,\theta _C \right)$ are also denoted as $P\left( {{\rho }_{AB}},\theta _A \right)$, $P\left( {{\rho }_{BC}},\theta _B \right)$ and $P\left( {{\rho }_{CA}},\theta _C \right)$, respectively. They are respectively shown in figure \ref{fig:Bianzhouzuobiaoxi_A}, \ref{fig:Bianzhouzuobiaoxi_B}, and Figure \ref{fig:Bianzhouzuobiaoxi_C}.

Obviously, in the polar coordinate system $\mathcal{P}\left( AB \right)$, except for pole $A$, any point on the plane corresponds one-to-one to the polar coordinates in $\mathcal{P}\left( AB \right)$.


For a given $\triangle ABC$, the Cartesian coordinate system with origin $A$ and the polar coordinate system with pole $A$ combined together is referred to as the \textbf{Edge-Axis coordinate system of vertex $A$}, denoted as $\mathcal{T}_{A}$, $\mathcal{T}_{A}:=\mathcal{R}\left( AB \right) \cup{} \mathcal{P}\left( AB \right)$. $\mathcal{R}\left( AB \right)$ is called the Cartesian coordinate system of $\mathcal{T}_{A}$, and $\mathcal{P}\left( AB \right)$ is called the polar coordinate system of $\mathcal{T}_{A}$.


Similarly, the Edge-Axis coordinate systems of vertex $B$ and vertex $C$ can be obtained, denoted as $\mathcal{T}_{B}$ and $\mathcal{T}_{C}$ respectively, where $\mathcal{T}_{B}:=\mathcal{R}\left( BC \right) \cup{} \mathcal{P}\left( BC \right)$, $\mathcal{T}_{C}:=\mathcal{R}\left( CA \right) \cup{} \mathcal{P}\left( CA \right)$. $\mathcal{R}\left( BC \right)$ is called the Cartesian coordinate system of $\mathcal{T}_{B}$, $\mathcal{P}\left( BC \right)$ is called the polar coordinate system of $\mathcal{T}_{B}$; $\mathcal{R}\left( CA \right)$ is called the Cartesian coordinate system of $\mathcal{T}_{C}$, $\mathcal{P}\left( CA \right)$ is called the polar coordinate system of $\mathcal{T}_{C}$. They are shown in figure \ref{fig:Bianzhouzuobiaoxi_B} and \ref{fig:Bianzhouzuobiaoxi_C}.


The Edge-Axis coordinate system $\mathcal{T}_{A}$ of vertex $A$, the Edge-Axis coordinate system $\mathcal{T}_{B}$ of vertex $B$ and the Edge-Axis coordinate system $\mathcal{T}_{C}$ of vertex $C$ combined together, they are collectively referred to as the \textbf {Set of Edge-Axis coordinate system} of the $\triangle ABC$, abbreviated as the Edge-Axis coordinate set, denoted as $\mathcal{T}_{ABC}$, $\mathcal{T}_{ABC}:={{\mathcal{T}_{A}} \cup{} {\mathcal{T}_{B}} \cup{} {\mathcal{T}_{C}}}$.

The important feature of the Edge-Axis coordinate system is that the origin (pole) coincides with a vertex of a triangle. In the same Edge-Axis coordinate system(e.g. $\mathcal{T}_{A}$), the origin and pole are coincident, and the horizontal and polar axes are also coincident, and coincide with one edge of the triangle. This is why I named it the Edge-Axis coordinate system.


For $\mathcal{R}\left( AB \right)$ and $\mathcal{P}\left( AB \right)$ in ${\mathcal{T}_{A}}$, the relationship between the Cartesian and polar coordinates is as follows:
\[x_{\mathcal{R}\left( AB \right)}^{P}={{\rho }_{AP}}\cos {{\theta }_{A}},\ y_{\mathcal{R}\left( AB \right)}^{P}={{\rho }_{AP}}\sin {{\theta }_{A}}, {{\rho }_{AP}}\ge 0.\]


For $\mathcal{R}\left( BC \right)$ and $\mathcal{P}\left( BC \right)$ in ${\mathcal{T}_{B}}$, the relationship between the Cartesian and polar coordinates is as follows:
\[x_{\mathcal{R}\left( BC \right)}^{P}={{\rho }_{BP}}\cos {{\theta }_{B}},\ y_{\mathcal{R}\left( BC \right)}^{P}={{\rho }_{BP}}\sin {{\theta }_{B}}, {{\rho }_{BP}}\ge 0.\]


For $\mathcal{R}\left( CA \right)$ and $\mathcal{P}\left( CA \right)$ in ${\mathcal{T}_{C}}$, the relationship between the Cartesian and polar coordinates is as follows:
\[x_{\mathcal{R}\left( CA \right)}^{P}={{\rho }_{CP}}\cos {{\theta }_{C}},\ y_{\mathcal{R}\left( CA \right)}^{P}={{\rho }_{CP}}\sin {{\theta }_{C}}, {{\rho }_{CP}}\ge 0.\]


Obviously, if $ABC$ is arranged counterclockwise and point $P$ is located inside $\triangle ABC$, then the three vertical coordinates $y_{\mathcal{R}\left( AB \right)}^{P} $, $y_{\mathcal{R}\left( BC \right)}^{P}$ and $y_{\mathcal{R}\left( CA \right)}^{P}$ of point $P$ in the Cartesian coordinate system are all positive real numbers. The three polar angles of point $P$ in polar coordinate system $\theta _A$, $\theta _B$ and $\theta _C$ are also positive real numbers.

It should be noted that in analytic geometry, if a polar coordinate system is used, the poles are usually excluded because at the poles, the pole angle is uncertain, which means there is no definite polar coordinate for the poles, or there is no corresponding relationship between the coordinates of the poles.


In analytic geometry, since the coordinates of each point are closely related to the coordinate system, it is usually necessary to establish a single (global) Cartesian or polar coordinate system before solving a problem, which I refer to as the \textbf{single coordinate system}. A single coordinate system determines the coordinates of all known quantities and then solves them. Once the coordinate system is determined, it will not be easily modified, otherwise the global situation will be affected. On the other hand, different people have different habits, and even for the same problem, they may choose different coordinate systems or coordinate origin. This directly leads to a result where different people have different raw data (coordinates) for the same problem. How embarrassing it is to have different raw data for the same problem. Alternatively, if it is discovered midway that the previously selected single coordinate system is not suitable for solving the problem and a new single coordinate system needs to be established, all previously obtained raw data (coordinates) will be invalidated. Furthermore, if multiple problems need to be solved, the established single coordinate system may not be fully applicable, meaning that a single coordinate system suitable for solving certain problems may not be suitable for solving other problems. For the convenience of calculation, sometimes it is necessary to modify a single coordinate system, which will result in changes in the coordinates of all known quantities, which is not only inconvenient but also requires additional computational costs.


In Intercenter Geometry, this problem does not exist because the goal of Intercenter Geometry is to find the frame components, which are independent of the coordinate system. In other words, no matter what coordinate system you use, no matter where the coordinate origin (pole) is chosen, the obtained frame components are the same.

On the other hand, in Intercenter Geometry, not only can a single coordinate system be established, but more importantly, multiple coordinate systems can be established as needed. That's why I proposed the Edge-Axis coordinate system in this section. The choice of coordinate system depends entirely on the convenience of solving the problem. By utilizing the important feature of Intercenter Geometry that the frame components are independent of coordinate systems, it is possible to collaborate well with multiple coordinate systems to solve the frame components. The characteristic of Intercenter Geometry is very convenient in application. For example, when a polar coordinate system is needed to handle point $A$ in practical problems, $\mathcal{P}\left( BC \right)$ or $\mathcal{P}\left( CA \right)$ can be chosen, because in these two polar coordinate systems, point $A$ is not a pole, which solves the embarrassment of a single polar coordinate system being difficult to handle pole problems. Intercenter Geometry can easily use multiple coordinate systems simultaneously to solve the same problem.

As is well known, trigonometric functions are often used when dealing with triangle problems, and trigonometric functions are built on the basis of Cartesian coordinate systems. The Edge-Axis coordinate system can conveniently use trigonometric functions. For example, in figure \ref{fig:Bianzhouzuobiaoxi_B}, the sine function $\sin {{\theta }_{B}}$ of angle ${{\theta }_{B}}$ can have positive or negative values, depending on the quadrant where ${{\theta }_{B}}$ is located. The quadrant of ${{\theta }_{B}}$ can be easily determined using Figure \ref{fig:Bianzhouzuobiaoxi_B}. The trigonometric function here can also be seen as a \textbf{Correlated trigonometric function} associated with the Cartesian coordinate system $\mathcal{R}\left( BC \right)$.

Additionally, it should be emphasized that vector multiplication can be easily used by Edge-Axis coordinate system.

Readers may ask, since the concept of Plane Intercenter Geometry is to represent various geometric quantities using the lengths of the three sides of a given triangle, does introducing an Edge-Axis coordinate system return to the old path of analytical geometry? The answer is no, Intercenter Geometry will not follow the old path of analytic geometry.

The Edge-Axis coordinate system I introduced is to search for the frame component of another unknown point through the frame component of a known point.

Obviously, the positional relationship between one point and another requires a certain coordinate system to be determined. The Edge-Axis coordinate system I introduced is only an intermediary tool, and the resulting geometric quantity is still represented by the lengths of the three sides of a given triangle. So introducing the coordinate system of Edge-Axis here does not violate the concept of Intercenter Geometry.


\section{Calculate the unknown frame component of a point based on the known frame component of a point}\label{GenjuYizhidianDeBiaojiafenliangQiuWeizhidianDeBiaojiafenliang}

If we know the frame component $\alpha _{A}^{{{P}_{0}}}$, $\alpha _{B}^{{{P}_{0}}}$, $\alpha _{C}^{{{P}_{0}}}$ of an IC, and know the positional relationship of point $P$ and ${{P}_{0}}$, then the frame component of point $P$ can be obtained. The determination of the positional relationship of point $P$ and point ${{P}_{0}}$ requires the assistance of a coordinate system.

\begin{theorem}{Calculating unknown frame component by known frame component, Daiyuan Zhang}{GenjuYizhidianDeBiaojiafenliangQiuWeizhidianDeBiaojiafenliangGongshi}\label{GenjuYizhidianDeBiaojiafenliangQiuWeizhidianDeBiaojiafenliangGongshi}
	Given a $\triangle ABC$, the area of the triangle is $S$, $P\in {{\pi }_{ABC}}$, ${{P}_{0}}\in {{\pi }_{ABC}}$, the frame component of ${{P}_{0}}$ is $\alpha _{A}^{{{P}_{0}}}$, $\alpha _{B}^{{{P}_{0}}}$, $\alpha _{C}^{{{P}_{0}}}$
	
	1. $\angle {{P}_{0}}AP={{\theta }_{A}}$, $A{{P}_{0}}\ne 0$, then the frame component of $P$ is: 
	\[\alpha _{A}^{P}=1-\frac{AP}{A{{P}_{0}}}\frac{\left( \alpha _{B}^{{{P}_{0}}}{{c}^{2}}-\alpha _{C}^{{{P}_{0}}}{{b}^{2}} \right)\sin {{\theta }_{A}}+bc\left( \alpha _{B}^{{{P}_{0}}}\sin \left( A-{{\theta }_{A}} \right)+\alpha _{C}^{{{P}_{0}}}\sin \left( A+{{\theta }_{A}} \right) \right)}{2S},\]	
	\[\alpha _{B}^{P}=\frac{AP}{A{{P}_{0}}}\frac{\alpha _{B}^{{{P}_{0}}}bc\sin \left( A-{{\theta }_{A}} \right)-\alpha _{C}^{{{P}_{0}}}{{b}^{2}}\sin {{\theta }_{A}}}{2S},\]
	\[\alpha _{C}^{P}=\frac{AP}{A{{P}_{0}}}\frac{\alpha _{B}^{{{P}_{0}}}{{c}^{2}}\sin {{\theta }_{A}}+\alpha _{C}^{{{P}_{0}}}bc\sin \left( A+{{\theta }_{A}} \right)}{2S}.\]
	
	2. $\angle {{P}_{0}}BP={{\theta }_{B}}$, $B{{P}_{0}}\ne 0$, then the frame component of $P$ is: 
	\[\alpha _{B}^{P}=1-\frac{BP}{B{{P}_{0}}}\frac{\left( \alpha _{C}^{{{P}_{0}}}{{a}^{2}}-\alpha _{A}^{{{P}_{0}}}{{c}^{2}} \right)\sin {{\theta }_{B}}+ca\left( \alpha _{C}^{{{P}_{0}}}\sin \left( B-{{\theta }_{B}} \right)+\alpha _{A}^{{{P}_{0}}}\sin \left( B+{{\theta }_{B}} \right) \right)}{2S},\]
	\[\alpha _{C}^{P}=\frac{BP}{B{{P}_{0}}}\frac{\alpha _{C}^{{{P}_{0}}}ca\sin \left( B-{{\theta }_{B}} \right)-\alpha _{A}^{{{P}_{0}}}{{c}^{2}}\sin {{\theta }_{B}}}{2S},\]
	\[\alpha _{A}^{P}=\frac{BP}{B{{P}_{0}}}\frac{\alpha _{C}^{{{P}_{0}}}{{a}^{2}}\sin {{\theta }_{B}}+\alpha _{A}^{{{P}_{0}}}ca\sin \left( B+{{\theta }_{B}} \right)}{2S}.\]
	
	3. $\angle {{P}_{0}}CP={{\theta }_{C}}$, $C{{P}_{0}}\ne 0$, then the frame component of $P$ is: 
	\[\alpha _{C}^{P}=1-\frac{CP}{C{{P}_{0}}}\frac{\left( \alpha _{A}^{{{P}_{0}}}{{b}^{2}}-\alpha _{B}^{{{P}_{0}}}{{a}^{2}} \right)\sin {{\theta }_{C}}+ab\left( \alpha _{A}^{{{P}_{0}}}\sin \left( C-{{\theta }_{C}} \right)+\alpha _{B}^{{{P}_{0}}}\sin \left( C+{{\theta }_{C}} \right) \right)}{2S},\]
	\[\alpha _{A}^{P}=\frac{CP}{C{{P}_{0}}}\frac{\alpha _{A}^{{{P}_{0}}}ab\sin \left( C-{{\theta }_{C}} \right)-\alpha _{B}^{{{P}_{0}}}{{a}^{2}}\sin {{\theta }_{C}}}{2S},\]
	\[\alpha _{B}^{P}=\frac{CP}{C{{P}_{0}}}\frac{\alpha _{A}^{{{P}_{0}}}{{b}^{2}}\sin {{\theta }_{C}}+\alpha _{B}^{{{P}_{0}}}ab\sin \left( C+{{\theta }_{C}} \right)}{2S}.\]
\end{theorem}

\begin{proof}
	As shown in figure \ref{fig:Bianzhouzuobiaoxi_A}, only the first case is proven, while the other cases are similar.
		
	If point $O$ coincides with point $A$, then according to theorem \ref{thm:Thm6.1.3}, we have
	\begin{equation}\label{ChajiDianjiHunheFangcheng}
		\left\{ \begin{aligned}
			& \overrightarrow{A{{P}_{0}}}\times \overrightarrow{AP}=\left( \alpha _{B}^{{{P}_{0}}}\overrightarrow{AB}+\alpha _{C}^{{{P}_{0}}}\overrightarrow{AC} \right)\times \left( \alpha _{B}^{P}\overrightarrow{AB}+\alpha _{C}^{P}\overrightarrow{AC} \right) \\ 
			& \overrightarrow{A{{P}_{0}}}\cdot \overrightarrow{AP}=\left( \alpha _{B}^{{{P}_{0}}}\overrightarrow{AB}+\alpha _{C}^{{{P}_{0}}}\overrightarrow{AC} \right)\cdot \left( \alpha _{B}^{P}\overrightarrow{AB}+\alpha _{C}^{P}\overrightarrow{AC} \right) \\ 
			& \alpha _{A}^{P}+\alpha _{B}^{P}+\alpha _{C}^{P}=1. \\ 
		\end{aligned} \right.
	\end{equation}
	
	The above system of equations is called the \textbf {mixed system of equations of cross product and dot product}.
	
	Using the definition of vector product to obtain (see Appendix \ref{SanjiaoxingXianguanZhishi}): 
	\[\overrightarrow{A{{P}_{0}}}\times \overrightarrow{AP}=A{{P}_{0}}\cdot AP\sin {{\theta }_{A}}{{\bm{n}}}.\]
	
	Where ${{\theta }_{A}}:=\angle {{P}_{0}}AP$, which is a directed angle between the starting edge $\overrightarrow{A{{P}_{0}}}$ turning counterclockwise and the ending edge $\overrightarrow{AP}$. The value of the directed angle ${{\theta }_{A}}$ specified in this way can be either a positive real number or a negative real number.
	
	For $\triangle ABC$, it is agreed that the directed angle between the starting edge $\overrightarrow{AB}$ and the ending edge $\overrightarrow{AC}$ is $A:=\angle BAC$.
	
	Therefore
	\[\begin{aligned}
		& \left( \alpha _{B}^{{{P}_{0}}}\overrightarrow{AB}+\alpha _{C}^{{{P}_{0}}}\overrightarrow{AC} \right)\times \left( \alpha _{B}^{P}\overrightarrow{AB}+\alpha _{C}^{P}\overrightarrow{AC} \right) \\ 
		& =\left( \alpha _{B}^{{{P}_{0}}}\alpha _{C}^{P}-\alpha _{C}^{{{P}_{0}}}\alpha _{B}^{P} \right)AB\cdot AC\cdot \sin A{{\bm{n}}}. \\ 
	\end{aligned}\]
	
	Therefore
	\[A{{P}_{0}}\cdot AP\sin {{\theta }_{A}}{{\bm{n}}}=\left( \alpha _{B}^{{{P}_{0}}}\alpha _{C}^{P}-\alpha _{C}^{{{P}_{0}}}\alpha _{B}^{P} \right)cb\sin A\cdot {{\bm{n}}}.\]
	
	i.e.
	\[A{{P}_{0}}\cdot AP\sin {{\theta }_{A}}=\left( \alpha _{B}^{{{P}_{0}}}\alpha _{C}^{P}-\alpha _{C}^{{{P}_{0}}}\alpha _{B}^{P} \right)cb\sin A.\]
	
	Using the dot product definition of vectors to obtain
	\[\overrightarrow{A{{P}_{0}}}\cdot \overrightarrow{AP}=A{{P}_{0}}\cdot AP\cos {{\theta }_{A}}=A{{P}_{0}}\cdot AP\cos {{\theta }_{A}}.\]
	
	And
	\[\begin{aligned}
		& \left( \alpha _{B}^{{{P}_{0}}}\overrightarrow{AB}+\alpha _{C}^{{{P}_{0}}}\overrightarrow{AC} \right)\cdot \left( \alpha _{B}^{P}\overrightarrow{AB}+\alpha _{C}^{P}\overrightarrow{AC} \right) \\ 
		& =\alpha _{B}^{{{P}_{0}}}\alpha _{B}^{P}\overrightarrow{AB}\cdot \overrightarrow{AB}+\alpha _{C}^{{{P}_{0}}}\alpha _{C}^{P}\overrightarrow{AC}\cdot \overrightarrow{AC}+\alpha _{B}^{{{P}_{0}}}\alpha _{C}^{P}\overrightarrow{AB}\cdot \overrightarrow{AC}+\alpha _{C}^{{{P}_{0}}}\alpha _{B}^{P}\overrightarrow{AC}\cdot \overrightarrow{AB} \\ 
		& =\alpha _{B}^{{{P}_{0}}}\alpha _{B}^{P}{{c}^{2}}+\alpha _{C}^{{{P}_{0}}}\alpha _{C}^{P}{{b}^{2}}+\alpha _{B}^{{{P}_{0}}}\alpha _{C}^{P}bc\cos A+\alpha _{C}^{{{P}_{0}}}\alpha _{B}^{P}bc\cos A \\ 
		& =\left( \alpha _{B}^{{{P}_{0}}}{{c}^{2}}+\alpha _{C}^{{{P}_{0}}}bc\cos A \right)\alpha _{B}^{P}+\left( \alpha _{C}^{{{P}_{0}}}{{b}^{2}}+\alpha _{B}^{{{P}_{0}}}bc\cos A \right)\alpha _{C}^{P}. \\ 
	\end{aligned}\]
	
	i.e.
	\[A{{P}_{0}}\cdot AP\cos {{\theta }_{A}}=\left( \alpha _{B}^{{{P}_{0}}}{{c}^{2}}+\alpha _{C}^{{{P}_{0}}}bc\cos A \right)\alpha _{B}^{P}+\left( \alpha _{C}^{{{P}_{0}}}{{b}^{2}}+\alpha _{B}^{{{P}_{0}}}bc\cos A \right)\alpha _{C}^{P}.\]
	
	So the following system of linear equations is obtained:
	\[\left\{ \begin{aligned}
		& \alpha _{A}^{P}+\alpha _{B}^{P}+\alpha _{C}^{P}=1 \\ 
		& \left( \alpha _{B}^{{{P}_{0}}}\alpha _{C}^{P}-\alpha _{C}^{{{P}_{0}}}\alpha _{B}^{P} \right)bc\sin A=A{{P}_{0}}\cdot AP\sin {{\theta }_{A}} \\ 
		& \left( \alpha _{B}^{{{P}_{0}}}{{c}^{2}}+\alpha _{C}^{{{P}_{0}}}bc\cos A \right)\alpha _{B}^{P}+\left( \alpha _{C}^{{{P}_{0}}}{{b}^{2}}+\alpha _{B}^{{{P}_{0}}}bc\cos A \right)\alpha _{C}^{P}=A{{P}_{0}}\cdot AP\cos {{\theta }_{A}}. \\ 
	\end{aligned} \right.\]	
	
	Solving the above linear equation system yields:
	\[\alpha _{A}^{P}=1-\frac{A{{P}_{0}}\cdot AP\left( \begin{aligned}
			& \alpha _{B}^{{{P}_{0}}}\left( \sin \left( {{\theta }_{A}} \right){{c}^{2}}+bc\sin \left( A-{{\theta }_{A}} \right) \right) \\ 
			& -\alpha _{C}^{{{P}_{0}}}{{b}^{2}}\sin \left( {{\theta }_{A}} \right)+\alpha _{C}^{{{P}_{0}}}bc\sin \left( A+{{\theta }_{A}} \right) \\ 
		\end{aligned} \right)}{bc\left( \sin \left( A \right){{\left( \alpha _{B}^{{{P}_{0}}} \right)}^{2}}{{c}^{2}}+\sin \left( 2A \right)\alpha _{B}^{{{P}_{0}}}\alpha _{C}^{{{P}_{0}}}bc+\sin \left( A \right){{\left( \alpha _{C}^{{{P}_{0}}} \right)}^{2}}{{b}^{2}} \right)},\]
	i.e.
	\[\alpha _{A}^{P}=1-\frac{A{{P}_{0}}\cdot AP\left( \left( \alpha _{B}^{{{P}_{0}}}{{c}^{2}}-\alpha _{C}^{{{P}_{0}}}{{b}^{2}} \right)\sin {{\theta }_{A}}+bc\left( \alpha _{B}^{{{P}_{0}}}\sin \left( A-{{\theta }_{A}} \right)+\alpha _{C}^{{{P}_{0}}}\sin \left( A+{{\theta }_{A}} \right) \right) \right)}{\left( {{\left( \alpha _{B}^{{{P}_{0}}} \right)}^{2}}{{c}^{2}}+2\alpha _{B}^{{{P}_{0}}}\alpha _{C}^{{{P}_{0}}}bc\cos A+{{\left( \alpha _{C}^{{{P}_{0}}} \right)}^{2}}{{b}^{2}} \right)bc\sin A},\]
	\[\alpha _{B}^{P}=\frac{A{{P}_{0}}\cdot AP\left( \alpha _{B}^{{{P}_{0}}}c\sin \left( A-{{\theta }_{A}} \right)-\alpha _{C}^{{{P}_{0}}}b\sin {{\theta }_{A}} \right)}{\left( {{\left( \alpha _{B}^{{{P}_{0}}} \right)}^{2}}{{c}^{2}}+2\alpha _{B}^{{{P}_{0}}}\alpha _{C}^{{{P}_{0}}}bc\cos A+{{\left( \alpha _{C}^{{{P}_{0}}} \right)}^{2}}{{b}^{2}} \right)c\sin A},\]
	\[\alpha _{C}^{P}=\frac{A{{P}_{0}}\cdot AP\left( \alpha _{B}^{{{P}_{0}}}c\sin {{\theta }_{A}}+\alpha _{C}^{{{P}_{0}}}b\sin \left( A+{{\theta }_{A}} \right) \right)}{\left( {{\left( \alpha _{B}^{{{P}_{0}}} \right)}^{2}}{{c}^{2}}+2\alpha _{B}^{{{P}_{0}}}\alpha _{C}^{{{P}_{0}}}bc\cos A+{{\left( \alpha _{C}^{{{P}_{0}}} \right)}^{2}}{{b}^{2}} \right)b\sin A}.\]
	
	From theorem \ref{thm:Thm12.1.2} in later chapters, we have
	\[\begin{aligned}
		& {{\left( \alpha _{B}^{{{P}_{0}}} \right)}^{2}}{{c}^{2}}+2\alpha _{B}^{{{P}_{0}}}\alpha _{C}^{{{P}_{0}}}bc\cos A+{{\left( \alpha _{C}^{{{P}_{0}}} \right)}^{2}}{{b}^{2}} \\ 
		& ={{\left( \alpha _{B}^{{{P}_{0}}} \right)}^{2}}{{c}^{2}}+\alpha _{B}^{{{P}_{0}}}\alpha _{C}^{{{P}_{0}}}\left( {{b}^{2}}+{{c}^{2}}-{{a}^{2}} \right)+{{\left( \alpha _{C}^{{{P}_{0}}} \right)}^{2}}{{b}^{2}} \\ 
		& ={{\left( \alpha _{B}^{{{P}_{0}}} \right)}^{2}}{{c}^{2}}+\alpha _{B}^{{{P}_{0}}}\alpha _{C}^{{{P}_{0}}}\left( {{b}^{2}}+{{c}^{2}} \right)+{{\left( \alpha _{C}^{{{P}_{0}}} \right)}^{2}}{{b}^{2}}-\alpha _{B}^{{{P}_{0}}}\alpha _{C}^{{{P}_{0}}}{{a}^{2}} \\ 
		& ={{\left( \alpha _{B}^{{{P}_{0}}} \right)}^{2}}{{c}^{2}}+\alpha _{B}^{{{P}_{0}}}\alpha _{C}^{{{P}_{0}}}{{c}^{2}}+\alpha _{B}^{{{P}_{0}}}\alpha _{C}^{{{P}_{0}}}{{b}^{2}}+{{\left( \alpha _{C}^{{{P}_{0}}} \right)}^{2}}{{b}^{2}}-\alpha _{B}^{{{P}_{0}}}\alpha _{C}^{{{P}_{0}}}{{a}^{2}} \\ 
		& =\left( \alpha _{B}^{{{P}_{0}}}+\alpha _{C}^{{{P}_{0}}} \right)\alpha _{B}^{{{P}_{0}}}{{c}^{2}}+\left( \alpha _{B}^{{{P}_{0}}}+\alpha _{C}^{{{P}_{0}}} \right)\alpha _{C}^{{{P}_{0}}}{{b}^{2}}-\alpha _{B}^{{{P}_{0}}}\alpha _{C}^{{{P}_{0}}}{{a}^{2}} \\ 
		& =\left( \alpha _{B}^{{{P}_{0}}}+\alpha _{C}^{{{P}_{0}}} \right)\left( \alpha _{B}^{{{P}_{0}}}{{c}^{2}}+\alpha _{C}^{{{P}_{0}}}{{b}^{2}} \right)-\alpha _{B}^{{{P}_{0}}}\alpha _{C}^{{{P}_{0}}}{{a}^{2}} \\ 
		& =\left( 1-\alpha _{A}^{{{P}_{0}}} \right)\left( \alpha _{B}^{{{P}_{0}}}{{c}^{2}}+\alpha _{C}^{{{P}_{0}}}{{b}^{2}} \right)-\alpha _{B}^{{{P}_{0}}}\alpha _{C}^{{{P}_{0}}}{{a}^{2}} \\ 
		& =A{{P}_{0}}^{2}.  
	\end{aligned}\]
	
	Therefore
	\begin{align*}
		\alpha _{A}^{P}& =1-\frac{AP}{A{{P}_{0}}bc\sin A}\left( \left( \alpha _{B}^{{{P}_{0}}}{{c}^{2}}-\alpha _{C}^{{{P}_{0}}}{{b}^{2}} \right)\sin {{\theta }_{A}}+bc\left( \alpha _{B}^{{{P}_{0}}}\sin \left( A-{{\theta }_{A}} \right)+\alpha _{C}^{{{P}_{0}}}\sin \left( A+{{\theta }_{A}} \right) \right) \right) \\ 
		& =1-\frac{AP}{A{{P}_{0}}}\frac{\left( \alpha _{B}^{{{P}_{0}}}{{c}^{2}}-\alpha _{C}^{{{P}_{0}}}{{b}^{2}} \right)\sin {{\theta }_{A}}+bc\left( \alpha _{B}^{{{P}_{0}}}\sin \left( A-{{\theta }_{A}} \right)+\alpha _{C}^{{{P}_{0}}}\sin \left( A+{{\theta }_{A}} \right) \right)}{2S},  
	\end{align*}
	\begin{align*}
		\alpha _{B}^{P}& =\frac{AP}{A{{P}_{0}}c\sin A}\left( \alpha _{B}^{{{P}_{0}}}c\sin \left( A-{{\theta }_{A}} \right)-\alpha _{C}^{{{P}_{0}}}b\sin {{\theta }_{A}} \right) \\ 
		& =\frac{AP}{A{{P}_{0}}}\frac{\alpha _{B}^{{{P}_{0}}}bc\sin \left( A-{{\theta }_{A}} \right)-\alpha _{C}^{{{P}_{0}}}{{b}^{2}}\sin {{\theta }_{A}}}{2S},  
	\end{align*}
	\begin{align*}
		\alpha _{C}^{P}& =\frac{AP}{A{{P}_{0}}b\sin A}\left( \alpha _{B}^{{{P}_{0}}}c\sin {{\theta }_{A}}+\alpha _{C}^{{{P}_{0}}}b\sin \left( A+{{\theta }_{A}} \right) \right) \\ 
		& =\frac{AP}{A{{P}_{0}}}\frac{\alpha _{B}^{{{P}_{0}}}{{c}^{2}}\sin {{\theta }_{A}}+\alpha _{C}^{{{P}_{0}}}bc\sin \left( A+{{\theta }_{A}} \right)}{2S}.  
	\end{align*}
	
	Similarly, the other results can be proven.
\end{proof}
\hfill $\square$\par


The above theorem explains (taking the polar coordinate system $\mathcal{P}\left( AB \right)$ as an example) that as long as the \textbf{rotation angle} ${{\theta }_{A}}$ and the \textbf{stretching ratio} ${AP}/{A{{P}_{0}}}\;$ are known, the frame component of any point $P$ can be calculated based on the known frame component of ${{P}_{0}}$.

If specific values of rotation angle ${{\theta }_{A}}$ and stretching ratio ${AP}/{A{{P}_{0}}}\;$ are given, the positional relationship between the point ${{P}_{0}}$ and $P$ is determined. In fact, the position of point $P$ can always be determined through rotation and stretching transformation based on the known point ${{P}_{0}}$, which means that the positional relationship between ${{P}_{0}}$ and $P$ is obtained by the rotation angle ${{\theta }_{A}}$ and the stretching ratio ${AP}/{A{{P}_{0}}}\;$.

The rotation angle ${{\theta }_{A}}$ and the stretching ratio ${AP}/{A{{P}_{0}}}\;$ in the above theorem can be seen as a parameter, which is a known quantity for practical problems and can be regarded as a constant. According to the formula for the frame component of point $P$ obtained from the above theorem, it can be seen that the frame component of point $P$ is derived from a known frame component of point ${{P}_{0}}$ and the three side lengths of $\triangle ABC$, so as long as the frame component of point ${{P}_{0}}$ can be represented by the three side lengths of $\triangle ABC$, the final frame component of point $P$ can also be represented by the three side lengths of $\triangle ABC$.

From the previous discussion, it is known that the frame components of the vertices of a triangle are constant, and the frame components of the centroid are also constant. The frame components of the incenter, orthocenter, and circumcenter are represented by the three side lengths of $\triangle ABC$, so these points can all be used as a known point ${{P}_{0}}$.

The above theorem can be written in the form of trigonometric functions.

\begin{corollary}{Trigonometric form of theorem \ref{thm:GenjuYizhidianDeBiaojiafenliangQiuWeizhidianDeBiaojiafenliangGongshi}, Daiyuan Zhang}{GenjuYizhidianDeBiaojiafenliangQiuWeizhidianDeBiaojiafenliangGongshi_SanjiaohanshuXingshi}\label{GenjuYizhidianDeBiaojiafenliangQiuWeizhidianDeBiaojiafenliangGongshi_SanjiaohanshuXingshi}
	Given a $\triangle ABC$, the frame component of point ${{P}_{0}}$ is $\alpha _{A}^{{{P}_{0}}}$, $\alpha _{B}^{{{P}_{0}}}$, $\alpha _{C}^{{{P}_{0}}}$.
	
	1. $P\in {{\pi }_{ABC}}$, $\angle {{{P}_{0}}}AP={{\theta }_{A}}$, $A{{{P}_{0}}}\ne 0$, then the frame component of point $P$ is: 
	\[\alpha _{A}^{P}=1-\frac{AP}{A{{{P}_{0}}}}\frac{\left( \begin{aligned}
			& \left( \alpha _{B}^{{P}_{0}}{{\sin }^{2}}C-\alpha _{C}^{{P}_{0}}{{\sin }^{2}}B \right)\sin {{\theta }_{A}} \\ 
			& +\sin B\sin C\left( \alpha _{B}^{{P}_{0}}\sin \left( A-{{\theta }_{A}} \right)+\alpha _{C}^{{P}_{0}}\sin \left( A+{{\theta }_{A}} \right) \right) \\ 
		\end{aligned} \right)}{\sin A\sin B\sin C},\]
	\[\alpha _{B}^{P}=\frac{AP}{A{{{P}_{0}}}}\frac{\alpha _{B}^{{P}_{0}}\sin C\sin \left( A-{{\theta }_{A}} \right)-\alpha _{C}^{{P}_{0}}\sin B\sin {{\theta }_{A}}}{\sin C\sin A},\]
	\[\alpha _{C}^{P}=\frac{AP}{A{{{P}_{0}}}}\frac{\alpha _{B}^{{P}_{0}}\sin C\sin {{\theta }_{A}}+\alpha _{C}^{{P}_{0}}\sin B\sin \left( A+{{\theta }_{A}} \right)}{\sin A\sin B}.\]
	
	2. $P\in {{\pi }_{ABC}}$, $\angle {{{P}_{0}}}BP={{\theta }_{B}}$, $B{{{P}_{0}}}\ne 0$, then the frame component of point $P$ is: 
	\[\alpha _{B}^{P}=1-\frac{BP}{B{{{P}_{0}}}}\frac{\left( \begin{aligned}
			& \left( \alpha _{C}^{{P}_{0}}\sin {{}^{2}}A-\alpha _{A}^{{P}_{0}}{{\sin }^{2}}C \right)\sin {{\theta }_{B}} \\ 
			& +\sin C\sin A\left( \alpha _{C}^{{P}_{0}}\sin \left( B-{{\theta }_{B}} \right)+\alpha _{A}^{{P}_{0}}\sin \left( B+{{\theta }_{B}} \right) \right) \\ 
		\end{aligned} \right)}{\sin B\sin C\sin A},\]
	\[\alpha _{C}^{P}=\frac{BP}{B{{{P}_{0}}}}\frac{\alpha _{C}^{{P}_{0}}\sin A\sin \left( B-{{\theta }_{B}} \right)-\alpha _{A}^{{P}_{0}}\sin C\sin {{\theta }_{B}}}{\sin A\sin B},\]
	\[\alpha _{A}^{P}=\frac{BP}{B{{{P}_{0}}}}\frac{\alpha _{C}^{{P}_{0}}\sin A\sin {{\theta }_{B}}+\alpha _{A}^{{P}_{0}}\sin C\sin \left( B+{{\theta }_{B}} \right)}{\sin C\sin B}.\]
	
	3. $P\in {{\pi }_{ABC}}$, $\angle {{{P}_{0}}}CP={{\theta }_{C}}$, $C{{{P}_{0}}}\ne 0$, then the frame component of point $P$ is: 
	\[\alpha _{C}^{P}=1-\frac{CP}{C{{{P}_{0}}}}\frac{\left( \begin{aligned}
			& \left( \alpha _{A}^{{P}_{0}}{{\sin }^{2}}B-\alpha _{B}^{{P}_{0}}{{\sin }^{2}}A \right)\sin {{\theta }_{C}} \\ 
			& +\sin A\sin B\left( \alpha _{A}^{{P}_{0}}\sin \left( C-{{\theta }_{C}} \right)+\alpha _{B}^{{P}_{0}}\sin \left( C+{{\theta }_{C}} \right) \right) \\ 
		\end{aligned} \right)}{\sin A\sin B\sin C},\]
	\[\alpha _{A}^{P}=\frac{CP}{C{{{P}_{0}}}}\frac{\alpha _{A}^{{P}_{0}}\sin B\sin \left( C-{{\theta }_{C}} \right)-\alpha _{B}^{{P}_{0}}\sin A\sin {{\theta }_{C}}}{\sin B\sin C},\]
	\[\alpha _{B}^{P}=\frac{CP}{C{{{P}_{0}}}}\frac{\alpha _{A}^{{P}_{0}}\sin B\sin {{\theta }_{C}}+\alpha _{B}^{{P}_{0}}\sin A\sin \left( C+{{\theta }_{C}} \right)}{\sin C\sin A}.\]
\end{corollary}

\begin{proof}
	For $\triangle ABC$, using the sine theorem of triangles: $a=2R\sin A$, $b=2R\sin B$, $c=2R\sin C$, according to theorem \ref{thm:GenjuYizhidianDeBiaojiafenliangQiuWeizhidianDeBiaojiafenliangGongshi}, this corollary can be obtained.
\end{proof}
\hfill $\square$\par

The above theorem can also be written in the form of side length.

\begin{corollary}{Side-length form of theorem \ref{thm:GenjuYizhidianDeBiaojiafenliangQiuWeizhidianDeBiaojiafenliangGongshi}, Daiyuan Zhang}{GenjuYizhidianDeBiaojiafenliangQiuWeizhidianDeBiaojiafenliangGongshi_BianchangXingshi}\label{GenjuYizhidianDeBiaojiafenliangQiuWeizhidianDeBiaojiafenliangGongshi_BianchangXingshi}
	Given a $\triangle ABC$, the area of the triangle is $S$, the frame component of point ${{P}_{0}}$ is $\alpha _{A}^{{{P}_{0}}}$, $\alpha _{B}^{{{P}_{0}}}$, $\alpha _{C}^{{{P}_{0}}}$.
	
	1. $P\in {{\pi }_{ABC}}$, $\angle {{P}_{0}}AP={{\theta }_{A}}$, $A{{P}_{0}}\ne 0$, then the frame component of point $P$ is: 
	\[\alpha _{A}^{P}=1-\frac{AP}{A{{P}_{0}}}\frac{4S\left( \alpha _{B}^{{{P}_{0}}}+\alpha _{C}^{{{P}_{0}}} \right)\cos {{\theta }_{A}}+\left( \alpha _{B}^{{{P}_{0}}}\left( {{c}^{2}}+{{a}^{2}}-{{b}^{2}} \right)-\alpha _{C}^{{{P}_{0}}}\left( {{a}^{2}}+{{b}^{2}}-{{c}^{2}} \right) \right)\sin {{\theta }_{A}}}{4S},\]
	\[\alpha _{B}^{P}=\frac{AP}{A{{P}_{0}}}\frac{4\alpha _{B}^{{{P}_{0}}}S\cos {{\theta }_{A}}-\left( \alpha _{B}^{{{P}_{0}}}\left( {{b}^{2}}+{{c}^{2}}-{{a}^{2}} \right)+2\alpha _{C}^{{{P}_{0}}}{{b}^{2}} \right)\sin {{\theta }_{A}}}{4S},\]
	\[\alpha _{C}^{P}=\frac{AP}{A{{P}_{0}}}\frac{4S\alpha _{C}^{{{P}_{0}}}\cos {{\theta }_{A}}+\left( 2\alpha _{B}^{{{P}_{0}}}{{c}^{2}}+\alpha _{C}^{{{P}_{0}}}\left( {{b}^{2}}+{{c}^{2}}-{{a}^{2}} \right) \right)\sin {{\theta }_{A}}}{4S}.\]
	
	2. $P\in {{\pi }_{ABC}}$, $\angle {{P}_{0}}BP={{\theta }_{B}}$, $B{{P}_{0}}\ne 0$, then the frame component of point $P$ is: 
	\[\alpha _{B}^{P}=1-\frac{BP}{B{{P}_{0}}}\frac{4S\left( \alpha _{C}^{{{P}_{0}}}+\alpha _{A}^{{{P}_{0}}} \right)\cos {{\theta }_{B}}+\left( \alpha _{C}^{{{P}_{0}}}\left( {{a}^{2}}+{{b}^{2}}-{{c}^{2}} \right)-\alpha _{A}^{{{P}_{0}}}\left( {{b}^{2}}+{{c}^{2}}-{{a}^{2}} \right) \right)\sin {{\theta }_{B}}}{4S},\]
	\[\alpha _{C}^{P}=\frac{BP}{B{{P}_{0}}}\frac{4\alpha _{C}^{{{P}_{0}}}S\cos {{\theta }_{B}}-\left( \alpha _{C}^{{{P}_{0}}}\left( {{c}^{2}}+{{a}^{2}}-{{b}^{2}} \right)+2\alpha _{A}^{{{P}_{0}}}{{c}^{2}} \right)\sin {{\theta }_{B}}}{4S},\]
	\[\alpha _{A}^{P}=\frac{BP}{B{{P}_{0}}}\frac{4S\alpha _{A}^{{{P}_{0}}}\cos {{\theta }_{B}}+\left( 2\alpha _{C}^{{{P}_{0}}}{{a}^{2}}+\alpha _{A}^{{{P}_{0}}}\left( {{c}^{2}}+{{a}^{2}}-{{b}^{2}} \right) \right)\sin {{\theta }_{B}}}{4S}.\]
	
	3. $P\in {{\pi }_{ABC}}$, $\angle {{P}_{0}}CP={{\theta }_{C}}$, $C{{P}_{0}}\ne 0$, then the frame component of point $P$ is: 
	\[\alpha _{C}^{P}=1-\frac{CP}{C{{P}_{0}}}\frac{4S\left( \alpha _{A}^{{{P}_{0}}}+\alpha _{B}^{{{P}_{0}}} \right)\cos {{\theta }_{C}}+\left( \alpha _{A}^{{{P}_{0}}}\left( {{b}^{2}}+{{c}^{2}}-{{a}^{2}} \right)-\alpha _{B}^{{{P}_{0}}}\left( {{c}^{2}}+{{a}^{2}}-{{b}^{2}} \right) \right)\sin {{\theta }_{C}}}{4S},\]
	\[\alpha _{A}^{P}=\frac{CP}{C{{P}_{0}}}\frac{4\alpha _{A}^{{{P}_{0}}}S\cos {{\theta }_{C}}-\left( \alpha _{A}^{{{P}_{0}}}\left( {{a}^{2}}+{{b}^{2}}-{{c}^{2}} \right)+2\alpha _{B}^{{{P}_{0}}}{{a}^{2}} \right)\sin {{\theta }_{C}}}{4S},\]
	\[\alpha _{B}^{P}=\frac{CP}{C{{P}_{0}}}\frac{4S\alpha _{B}^{{{P}_{0}}}\cos {{\theta }_{C}}+\left( 2\alpha _{A}^{{{P}_{0}}}{{b}^{2}}+\alpha _{B}^{{{P}_{0}}}\left( {{a}^{2}}+{{b}^{2}}-{{c}^{2}} \right) \right)\sin {{\theta }_{C}}}{4S}.\]
\end{corollary}

\begin{proof}
	Only prove the first scenario, the rest are similar. Because
	\[\sin A=\frac{2S}{bc},\]
	\[\cos A=\frac{{{b}^{2}}+{{c}^{2}}-{{a}^{2}}}{2bc}.\]
	
	Therefore
	\begin{align*}
		\sin \left( A-{{\theta }_{A}} \right)& =\sin A\cos {{\theta }_{A}}-\cos A\sin {{\theta }_{A}} \\ 
		& =\frac{2S}{bc}\cos {{\theta }_{A}}-\frac{{{b}^{2}}+{{c}^{2}}-{{a}^{2}}}{2bc}\sin {{\theta }_{A}} \\ 
		& =\frac{1}{2bc}\left( 4S\cos {{\theta }_{A}}-\left( {{b}^{2}}+{{c}^{2}}-{{a}^{2}} \right)\sin {{\theta }_{A}} \right),  
	\end{align*}
	\begin{align*}
		\sin \left( A+{{\theta }_{A}} \right)& =\sin A\cos {{\theta }_{A}}+\cos A\sin {{\theta }_{A}} \\ 
		& =\frac{2S}{bc}\cos {{\theta }_{A}}+\frac{{{b}^{2}}+{{c}^{2}}-{{a}^{2}}}{2bc}\sin {{\theta }_{A}} \\ 
		& =\frac{1}{2bc}\left( 4S\cos {{\theta }_{A}}+\left( {{b}^{2}}+{{c}^{2}}-{{a}^{2}} \right)\sin {{\theta }_{A}} \right).  
	\end{align*}
	
	And
	\begin{align*}
		& bc\left( \alpha _{B}^{{{P}_{0}}}\sin \left( A-{{\theta }_{A}} \right)+\alpha _{C}^{{{P}_{0}}}\sin \left( A+{{\theta }_{A}} \right) \right) \\ 
		& =\frac{\alpha _{B}^{{{P}_{0}}}}{2}\left( 4S\cos {{\theta }_{A}}-\left( {{b}^{2}}+{{c}^{2}}-{{a}^{2}} \right)\sin {{\theta }_{A}} \right)+\frac{\alpha _{C}^{{{P}_{0}}}}{2}\left( 4S\cos {{\theta }_{A}}+\left( {{b}^{2}}+{{c}^{2}}-{{a}^{2}} \right)\sin {{\theta }_{A}} \right) \\ 
		& =2S\left( \alpha _{B}^{{{P}_{0}}}+\alpha _{C}^{{{P}_{0}}} \right)\cos {{\theta }_{A}}+\frac{1}{2}\left( \alpha _{C}^{{{P}_{0}}}\left( {{b}^{2}}+{{c}^{2}}-{{a}^{2}} \right)-\alpha _{B}^{{{P}_{0}}}\left( {{b}^{2}}+{{c}^{2}}-{{a}^{2}} \right) \right)\sin {{\theta }_{A}}. \\ 
	\end{align*}
	
	Therefore
	\begin{align*}
		& \left( \alpha _{B}^{{{P}_{0}}}{{c}^{2}}-\alpha _{C}^{{{P}_{0}}}{{b}^{2}} \right)\sin {{\theta }_{A}}+bc\left( \alpha _{B}^{{{P}_{0}}}\sin \left( A-{{\theta }_{A}} \right)+\alpha _{C}^{{{P}_{0}}}\sin \left( A+{{\theta }_{A}} \right) \right) \\ 
		& =2S\left( \alpha _{B}^{{{P}_{0}}}+\alpha _{C}^{{{P}_{0}}} \right)\cos {{\theta }_{A}}+\frac{1}{2}\left( \begin{aligned}
			& 2\left( \alpha _{B}^{{{P}_{0}}}{{c}^{2}}-\alpha _{C}^{P}{{b}^{2}} \right)+\alpha _{C}^{{{P}_{0}}}\left( {{b}^{2}}+{{c}^{2}}-{{a}^{2}} \right) \\ 
			& -\alpha _{B}^{{{P}_{0}}}\left( {{b}^{2}}+{{c}^{2}}-{{a}^{2}} \right) \\ 
		\end{aligned} \right)\sin {{\theta }_{A}} \\ 
		& =2S\left( \alpha _{B}^{{{P}_{0}}}+\alpha _{C}^{{{P}_{0}}} \right)\cos {{\theta }_{A}}+\frac{1}{2}\left( \alpha _{C}^{{{P}_{0}}}\left( {{c}^{2}}-{{a}^{2}}-{{b}^{2}} \right)-\alpha _{B}^{{{P}_{0}}}\left( {{b}^{2}}-{{c}^{2}}-{{a}^{2}} \right) \right)\sin {{\theta }_{A}} \\ 
		& =2S\left( \alpha _{B}^{{{P}_{0}}}+\alpha _{C}^{{{P}_{0}}} \right)\cos {{\theta }_{A}}+\frac{1}{2}\left( \alpha _{B}^{{{P}_{0}}}\left( {{c}^{2}}+{{a}^{2}}-{{b}^{2}} \right)-\alpha _{C}^{{{P}_{0}}}\left( {{a}^{2}}+{{b}^{2}}-{{c}^{2}} \right) \right)\sin {{\theta }_{A}}. \\ 
	\end{align*}
	
	According to theorem \ref{thm:GenjuYizhidianDeBiaojiafenliangQiuWeizhidianDeBiaojiafenliangGongshi}, it can be obtained that:
	\begin{align*}
		\alpha _{A}^{P}& =1-\frac{AP}{A{{P}_{0}}}\frac{\left( \alpha _{B}^{{{P}_{0}}}{{c}^{2}}-\alpha _{C}^{{{P}_{0}}}{{b}^{2}} \right)\sin {{\theta }_{A}}+bc\left( \alpha _{B}^{{{P}_{0}}}\sin \left( A-{{\theta }_{A}} \right)+\alpha _{C}^{{{P}_{0}}}\sin \left( A+{{\theta }_{A}} \right) \right)}{bc\sin A} \\ 
		& =1-\frac{AP}{A{{P}_{0}}}\frac{2S\left( \alpha _{B}^{{{P}_{0}}}+\alpha _{C}^{{{P}_{0}}} \right)\cos {{\theta }_{A}}+\frac{1}{2}\left( \alpha _{B}^{{{P}_{0}}}\left( {{c}^{2}}+{{a}^{2}}-{{b}^{2}} \right)-\alpha _{C}^{{{P}_{0}}}\left( {{a}^{2}}+{{b}^{2}}-{{c}^{2}} \right) \right)\sin {{\theta }_{A}}}{2S} \\ 
		& =1-\frac{AP}{A{{P}_{0}}}\frac{4S\left( \alpha _{B}^{{{P}_{0}}}+\alpha _{C}^{{{P}_{0}}} \right)\cos {{\theta }_{A}}+\left( \alpha _{B}^{{{P}_{0}}}\left( {{c}^{2}}+{{a}^{2}}-{{b}^{2}} \right)-\alpha _{C}^{{{P}_{0}}}\left( {{a}^{2}}+{{b}^{2}}-{{c}^{2}} \right) \right)\sin {{\theta }_{A}}}{4S}.  
	\end{align*}
	
	And
	\begin{align*}
		\alpha _{B}^{P}& =\frac{AP}{A{{P}_{0}}}\frac{\alpha _{B}^{{{P}_{0}}}c\sin \left( A-{{\theta }_{A}} \right)-\alpha _{C}^{{{P}_{0}}}b\sin {{\theta }_{A}}}{c\sin A} \\ 
		& =\frac{AP}{A{{P}_{0}}}\frac{\alpha _{B}^{{{P}_{0}}}c\frac{1}{2bc}\left( 4S\cos {{\theta }_{A}}-\left( {{b}^{2}}+{{c}^{2}}-{{a}^{2}} \right)\sin {{\theta }_{A}} \right)-\alpha _{C}^{{{P}_{0}}}b\sin {{\theta }_{A}}}{c\sin A} \\ 
		& =\frac{AP}{A{{P}_{0}}}\frac{\alpha _{B}^{{{P}_{0}}}\left( 4S\cos {{\theta }_{A}}-\left( {{b}^{2}}+{{c}^{2}}-{{a}^{2}} \right)\sin {{\theta }_{A}} \right)-2\alpha _{C}^{{{P}_{0}}}{{b}^{2}}\sin {{\theta }_{A}}}{2bc\sin A} \\ 
		& =\frac{AP}{A{{P}_{0}}}\frac{4\alpha _{B}^{{{P}_{0}}}S\cos {{\theta }_{A}}-\left( \alpha _{B}^{{{P}_{0}}}\left( {{b}^{2}}+{{c}^{2}}-{{a}^{2}} \right)+2\alpha _{C}^{{{P}_{0}}}{{b}^{2}} \right)\sin {{\theta }_{A}}}{4S},  
	\end{align*}
	\begin{align*}
		\alpha _{C}^{P}& =\frac{AP}{A{{P}_{0}}}\frac{\alpha _{B}^{{{P}_{0}}}c\sin {{\theta }_{A}}+\alpha _{C}^{{{P}_{0}}}b\sin \left( A+{{\theta }_{A}} \right)}{b\sin A} \\ 
		& =\frac{AP}{A{{P}_{0}}}\frac{\alpha _{B}^{{{P}_{0}}}c\sin {{\theta }_{A}}+\alpha _{C}^{{{P}_{0}}}b\frac{1}{2bc}\left( 4S\cos {{\theta }_{A}}+\left( {{b}^{2}}+{{c}^{2}}-{{a}^{2}} \right)\sin {{\theta }_{A}} \right)}{b\sin A} \\ 
		& =\frac{AP}{A{{P}_{0}}}\frac{2\alpha _{B}^{{{P}_{0}}}{{c}^{2}}\sin {{\theta }_{A}}+\alpha _{C}^{{{P}_{0}}}\left( 4S\cos {{\theta }_{A}}+\left( {{b}^{2}}+{{c}^{2}}-{{a}^{2}} \right)\sin {{\theta }_{A}} \right)}{4S} \\ 
		& =\frac{AP}{A{{P}_{0}}}\frac{4S\alpha _{C}^{{{P}_{0}}}\cos {{\theta }_{A}}+2\alpha _{B}^{{{P}_{0}}}{{c}^{2}}\sin {{\theta }_{A}}+\alpha _{C}^{{{P}_{0}}}\left( {{b}^{2}}+{{c}^{2}}-{{a}^{2}} \right)\sin {{\theta }_{A}}}{4S} \\ 
		& =\frac{AP}{A{{P}_{0}}}\frac{4S\alpha _{C}^{{{P}_{0}}}\cos {{\theta }_{A}}+\left( 2\alpha _{B}^{{{P}_{0}}}{{c}^{2}}+\alpha _{C}^{{{P}_{0}}}\left( {{b}^{2}}+{{c}^{2}}-{{a}^{2}} \right) \right)\sin {{\theta }_{A}}}{4S}.  
	\end{align*}
	
	Therefore
	\[\alpha _{A}^{P}=1-\frac{AP}{A{{P}_{0}}}\frac{4S\left( \alpha _{B}^{{{P}_{0}}}+\alpha _{C}^{{{P}_{0}}} \right)\cos {{\theta }_{A}}+\left( \alpha _{B}^{{{P}_{0}}}\left( {{c}^{2}}+{{a}^{2}}-{{b}^{2}} \right)-\alpha _{C}^{{{P}_{0}}}\left( {{a}^{2}}+{{b}^{2}}-{{c}^{2}} \right) \right)\sin {{\theta }_{A}}}{4S},\]
	\[\alpha _{B}^{P}=\frac{AP}{A{{P}_{0}}}\frac{4\alpha _{B}^{{{P}_{0}}}S\cos {{\theta }_{A}}-\left( \alpha _{B}^{{{P}_{0}}}\left( {{b}^{2}}+{{c}^{2}}-{{a}^{2}} \right)+2\alpha _{C}^{{{P}_{0}}}{{b}^{2}} \right)\sin {{\theta }_{A}}}{4S},\]
	\[\alpha _{C}^{P}=\frac{AP}{A{{P}_{0}}}\frac{4S\alpha _{C}^{{{P}_{0}}}\cos {{\theta }_{A}}+\left( 2\alpha _{B}^{{{P}_{0}}}{{c}^{2}}+\alpha _{C}^{{{P}_{0}}}\left( {{b}^{2}}+{{c}^{2}}-{{a}^{2}} \right) \right)\sin {{\theta }_{A}}}{4S}.\]
	
	Similarly, other formulas can be proven.
\end{proof}
\hfill $\square$\par

%
%
%

\begin{corollary}{Frame component of stretching transformation, Daiyuan Zhang}{ShensuobianhuanDeBiaojiafenliang}\label{ShensuobianhuanDeBiaojiafenliang}
	Given a $\triangle ABC$, the frame component of point ${{P}_{0}}$ is $\alpha _{A}^{{{P}_{0}}}$, $\alpha _{B}^{{{P}_{0}}}$, $\alpha _{C}^{{{P}_{0}}}$.
	
	1. $A{{P}_{0}}\ne 0$, $P\in \overleftrightarrow{A{{P}_{0}}}$,  then the frame component of point $P$ is: 
	\[\alpha _{A}^{P}=1-\frac{AP}{A{{P}_{0}}}\left( \alpha _{B}^{{{P}_{0}}}+\alpha _{C}^{{{P}_{0}}} \right),\]
	\[\alpha _{B}^{P}=\frac{AP}{A{{P}_{0}}}\alpha _{B}^{{{P}_{0}}},\]
	\[\alpha _{C}^{P}=\frac{AP}{A{{P}_{0}}}\alpha _{C}^{{{P}_{0}}}.\]
	
	2. $B{{P}_{0}}\ne 0$, $P\in \overleftrightarrow{B{{P}_{0}}}$,  then the frame component of point $P$ is: 
	\[\alpha _{B}^{P}=1-\frac{BP}{B{{P}_{0}}}\left( \alpha _{C}^{{{P}_{0}}}+\alpha _{A}^{{{P}_{0}}} \right),\]
	\[\alpha _{C}^{P}=\frac{BP}{B{{P}_{0}}}\alpha _{C}^{{{P}_{0}}},\]
	\[\alpha _{A}^{P}=\frac{BP}{B{{P}_{0}}}\alpha _{A}^{{{P}_{0}}}.\]
	
	3. $C{{P}_{0}}\ne 0$, $P\in \overleftrightarrow{C{{P}_{0}}}$,  then the frame component of point $P$ is: 
	\[\alpha _{C}^{P}=1-\frac{CP}{C{{P}_{0}}}\left( \alpha _{A}^{{{P}_{0}}}+\alpha _{B}^{{{P}_{0}}} \right),\]
	\[\alpha _{A}^{P}=\frac{CP}{C{{P}_{0}}}\alpha _{A}^{{{P}_{0}}},\]
	\[\alpha _{B}^{P}=\frac{CP}{C{{P}_{0}}}\alpha _{B}^{{{P}_{0}}}.\]
\end{corollary}


\begin{proof}
	The result can be obtained by making ${{\theta }_{A}}=0$, ${{\theta }_{B}}=0$ and ${{\theta }_{C}}=0$ in corollary \ref{cor:GenjuYizhidianDeBiaojiafenliangQiuWeizhidianDeBiaojiafenliangGongshi_SanjiaohanshuXingshi}, respectively.
\end{proof}
\hfill $\square$\par

\begin{corollary}{Frame component of rotation transformation-trigonometric form, Daiyuan Zhang}{XuanzhuanbianhuanDeBiaojiafenliang_Sanjiaohanshuxingshi}\label{XuanzhuanbianhuanDeBiaojiafenliang_Sanjiaohanshuxingshi}
	Given a $\triangle ABC$, the frame component of point ${{P}_{0}}$ is $\alpha _{A}^{{{P}_{0}}}$, $\alpha _{B}^{{{P}_{0}}}$, $\alpha _{C}^{{{P}_{0}}}$.
	
	1. $\angle {{P}_{0}}AP={{\theta }_{A}}$, $A{{P}_{0}}\ne 0$,  then the frame component of point $P$ is: 
	\[\alpha _{A}^{P}=1-\frac{\left( \begin{aligned}
			& \left( \alpha _{B}^{{{P}_{0}}}{{\sin }^{2}}C-\alpha _{C}^{{{P}_{0}}}{{\sin }^{2}}B \right)\sin {{\theta }_{A}} \\ 
			& +\sin B\sin C\left( \alpha _{B}^{{{P}_{0}}}\sin \left( A-{{\theta }_{A}} \right)+\alpha _{C}^{{{P}_{0}}}\sin \left( A+{{\theta }_{A}} \right) \right) \\ 
		\end{aligned} \right)}{\sin A\sin B\sin C},\]
	\[\alpha _{B}^{P}=\frac{\alpha _{B}^{{{P}_{0}}}\sin C\sin \left( A-{{\theta }_{A}} \right)-\alpha _{C}^{{{P}_{0}}}\sin B\sin {{\theta }_{A}}}{\sin C\sin A},\]
	\[\alpha _{C}^{P}=\frac{\alpha _{B}^{{{P}_{0}}}\sin C\sin {{\theta }_{A}}+\alpha _{C}^{{{P}_{0}}}\sin B\sin \left( A+{{\theta }_{A}} \right)}{\sin A\sin B}.\]
	
	2. $\angle {{P}_{0}}BP={{\theta }_{B}}$, $B{{P}_{0}}\ne 0$,  then the frame component of point $P$ is: 
	\[\alpha _{B}^{P}=1-\frac{\left( \begin{aligned}
			& \left( \alpha _{C}^{{{P}_{0}}}{{\sin }^{2}}A-\alpha _{A}^{{{P}_{0}}}{{\sin }^{2}}C \right)\sin {{\theta }_{B}} \\ 
			& +\sin C\sin A\left( \alpha _{C}^{{{P}_{0}}}\sin \left( B-{{\theta }_{B}} \right)+\alpha _{A}^{{{P}_{0}}}\sin \left( B+{{\theta }_{B}} \right) \right) \\ 
		\end{aligned} \right)}{\sin B\sin C\sin A},\]
	\[\alpha _{C}^{P}=\frac{\alpha _{C}^{{{P}_{0}}}\sin A\sin \left( B-{{\theta }_{B}} \right)-\alpha _{A}^{{{P}_{0}}}\sin C\sin {{\theta }_{B}}}{\sin A\sin B},\]
	\[\alpha _{A}^{P}=\frac{\alpha _{C}^{{{P}_{0}}}\sin A\sin {{\theta }_{B}}+\alpha _{A}^{{{P}_{0}}}\sin C\sin \left( B+{{\theta }_{B}} \right)}{\sin C\sin B}.\]
	
	3. $\angle {{P}_{0}}CP={{\theta }_{C}}$, $C{{P}_{0}}\ne 0$,  then the frame component of point $P$ is: 
	\[\alpha _{C}^{P}=1-\frac{\left( \begin{aligned}
			& \left( \alpha _{A}^{{{P}_{0}}}{{\sin }^{2}}B-\alpha _{B}^{{{P}_{0}}}{{\sin }^{2}}A \right)\sin {{\theta }_{C}} \\ 
			& +\sin A\sin B\left( \alpha _{A}^{{{P}_{0}}}\sin \left( C-{{\theta }_{C}} \right)+\alpha _{B}^{{{P}_{0}}}\sin \left( C+{{\theta }_{C}} \right) \right) \\ 
		\end{aligned} \right)}{\sin A\sin B\sin C},\]
	\[\alpha _{A}^{P}=\frac{\alpha _{A}^{{{P}_{0}}}\sin B\sin \left( C-{{\theta }_{C}} \right)-\alpha _{B}^{{{P}_{0}}}\sin A\sin {{\theta }_{C}}}{\sin B\sin C},\]
	\[\alpha _{B}^{P}=\frac{\alpha _{A}^{{{P}_{0}}}\sin B\sin {{\theta }_{C}}+\alpha _{B}^{{{P}_{0}}}\sin A\sin \left( C+{{\theta }_{C}} \right)}{\sin C\sin A}.\]
\end{corollary}


\begin{proof}
	The result can be obtained by making $AP=A{{P}_{0}}$, $BP=B{{P}_{0}}$ and $CP=C{{P}_{0}}$ in corollary \ref{cor:GenjuYizhidianDeBiaojiafenliangQiuWeizhidianDeBiaojiafenliangGongshi_SanjiaohanshuXingshi}, respectively. 
\end{proof}
\hfill $\square$\par

\begin{corollary}{Frame component of rotation transformation-side length form, Daiyuan Zhang}{XuanzhuanbianhuanDeBiaojiafenliang_Bianchangxingshi}\label{XuanzhuanbianhuanDeBiaojiafenliang_Bianchangxingshi}
	Given a $\triangle ABC$, the area of the triangle is $S$, the frame component of point ${{P}_{0}}$ is $\alpha _{A}^{{{P}_{0}}}$, $\alpha _{B}^{{{P}_{0}}}$, $\alpha _{C}^{{{P}_{0}}}$.
	
	1. $P\in {{\pi }_{ABC}}$, $\angle {{P}_{0}}AP={{\theta }_{A}}$, $A{{P}_{0}}\ne 0$,  then the frame component of point $P$ is: 
	\[\alpha _{A}^{P}=1-\frac{4S\left( \alpha _{B}^{{{P}_{0}}}+\alpha _{C}^{{{P}_{0}}} \right)\cos {{\theta }_{A}}+\left( \alpha _{B}^{{{P}_{0}}}\left( {{c}^{2}}+{{a}^{2}}-{{b}^{2}} \right)-\alpha _{C}^{{{P}_{0}}}\left( {{a}^{2}}+{{b}^{2}}-{{c}^{2}} \right) \right)\sin {{\theta }_{A}}}{4S},\]
	\[\alpha _{B}^{P}=\frac{4S\alpha _{B}^{{{P}_{0}}}\cos {{\theta }_{A}}-\left( \alpha _{B}^{{{P}_{0}}}\left( {{b}^{2}}+{{c}^{2}}-{{a}^{2}} \right)+2\alpha _{C}^{{{P}_{0}}}{{b}^{2}} \right)\sin {{\theta }_{A}}}{4S},\]
	\[\alpha _{C}^{P}=\frac{4S\alpha _{C}^{{{P}_{0}}}\cos {{\theta }_{A}}+\left( 2\alpha _{B}^{{{P}_{0}}}{{c}^{2}}+\alpha _{C}^{{{P}_{0}}}\left( {{b}^{2}}+{{c}^{2}}-{{a}^{2}} \right) \right)\sin {{\theta }_{A}}}{4S}.\]
	
	2. $P\in {{\pi }_{ABC}}$, $\angle {{P}_{0}}BP={{\theta }_{B}}$, $B{{P}_{0}}\ne 0$,  then the frame component of point $P$ is: 
	\[\alpha _{B}^{P}=1-\frac{4S\left( \alpha _{C}^{{{P}_{0}}}+\alpha _{A}^{{{P}_{0}}} \right)\cos {{\theta }_{B}}+\left( \alpha _{C}^{{{P}_{0}}}\left( {{a}^{2}}+{{b}^{2}}-{{c}^{2}} \right)-\alpha _{A}^{{{P}_{0}}}\left( {{b}^{2}}+{{c}^{2}}-{{a}^{2}} \right) \right)\sin {{\theta }_{B}}}{4S},\]
	\[\alpha _{C}^{P}=\frac{4S\alpha _{C}^{{{P}_{0}}}\cos {{\theta }_{B}}-\left( \alpha _{C}^{{{P}_{0}}}\left( {{c}^{2}}+{{a}^{2}}-{{b}^{2}} \right)+2\alpha _{A}^{{{P}_{0}}}{{c}^{2}} \right)\sin {{\theta }_{B}}}{4S},\]
	\[\alpha _{A}^{P}=\frac{4S\alpha _{A}^{{{P}_{0}}}\cos {{\theta }_{B}}+\left( 2\alpha _{C}^{{{P}_{0}}}{{a}^{2}}+\alpha _{A}^{{{P}_{0}}}\left( {{c}^{2}}+{{a}^{2}}-{{b}^{2}} \right) \right)\sin {{\theta }_{B}}}{4S}.\]
	
	3. $P\in {{\pi }_{ABC}}$, $\angle {{P}_{0}}CP={{\theta }_{C}}$, $C{{P}_{0}}\ne 0$,  then the frame component of point $P$ is: 
	\[\alpha _{C}^{P}=1-\frac{4S\left( \alpha _{A}^{{{P}_{0}}}+\alpha _{B}^{{{P}_{0}}} \right)\cos {{\theta }_{C}}+\left( \alpha _{A}^{{{P}_{0}}}\left( {{b}^{2}}+{{c}^{2}}-{{a}^{2}} \right)-\alpha _{B}^{{{P}_{0}}}\left( {{c}^{2}}+{{a}^{2}}-{{b}^{2}} \right) \right)\sin {{\theta }_{C}}}{4S},\]
	\[\alpha _{A}^{P}=\frac{4S\alpha _{A}^{{{P}_{0}}}\cos {{\theta }_{C}}-\left( \alpha _{A}^{{{P}_{0}}}\left( {{a}^{2}}+{{b}^{2}}-{{c}^{2}} \right)+2\alpha _{B}^{{{P}_{0}}}{{a}^{2}} \right)\sin {{\theta }_{C}}}{4S},\]
	\[\alpha _{B}^{P}=\frac{4S\alpha _{B}^{{{P}_{0}}}\cos {{\theta }_{C}}+\left( 2\alpha _{A}^{{{P}_{0}}}{{b}^{2}}+\alpha _{B}^{{{P}_{0}}}\left( {{a}^{2}}+{{b}^{2}}-{{c}^{2}} \right) \right)\sin {{\theta }_{C}}}{4S}.\]
\end{corollary}


\begin{proof}
	The result can be obtained by making $AP=A{{P}_{0}}$, $BP=B{{P}_{0}}$ and $CP=C{{P}_{0}}$ in corollary \ref{cor:GenjuYizhidianDeBiaojiafenliangQiuWeizhidianDeBiaojiafenliangGongshi_BianchangXingshi}, respectively. 
\end{proof}
\hfill $\square$\par  

\begin{example}{}\label{XuanzhuanYuShensuoJiansuanLizi}
	Given a $\triangle ABC$, let the area of $\triangle ABC$ be $S$, ${{P}_{0}}$ be the incenter $I$.
	
	1. Rotate $A{{P}_{0}}$ (i.e. $AI$) counterclockwise by $90{}^\circ $ with $A$ as the center, transform the point ${{P}_{0}}$ into a new point ${{P}_{1}}$, and calculate the frame component of ${{P}_{1}}$;
	
	2. Rotate $A{{P}_{0}}$ (i.e. $AI$) counterclockwise by $60{}^\circ $ with $A$ as the center, and then perform a stretching transformation to move the point ${{P}_{0}}$ into a new point ${{P}_{2}}$, where $A{{P}_{2}}=2A{{P}_{0}}=2AI$, find the frame component of ${{P}_{2}}$.
\end{example}

\begin{solution}
	1. Establish a Cartesian coordinate system (see Figure \ref{fig:Bianzhouzuobiaoxi_A}), and based on the previous corollary and the frame component formula of incenter, let ${{\theta }_{A}}=90{}^\circ $, $AP=A{{P}_{0}}$, it is obtained that:
	\begin{align*}
		\alpha _{A}^{{{P}_{1}}}&=1-\frac{AP}{A{{P}_{0}}}\frac{4S\left( \alpha _{B}^{{{P}_{0}}}+\alpha _{C}^{{{P}_{0}}} \right)\cos {{\theta }_{A}}+\left( \begin{aligned}
				& \alpha _{B}^{{{P}_{0}}}\left( {{c}^{2}}+{{a}^{2}}-{{b}^{2}} \right) 
				-\alpha _{C}^{{{P}_{0}}}\left( {{a}^{2}}+{{b}^{2}}-{{c}^{2}} \right) \\ 
			\end{aligned} \right)\sin {{\theta }_{A}}}{4S} \\ 
		& =1-\frac{4S\left( b+c \right)\cos {{\theta }_{A}}+\left( b\left( {{c}^{2}}+{{a}^{2}}-{{b}^{2}} \right)-c\left( {{a}^{2}}+{{b}^{2}}-{{c}^{2}} \right) \right)\sin {{\theta }_{A}}}{4S\left( a+b+c \right)},  
	\end{align*}
	
	And
	\begin{align*}
		& b\left( {{c}^{2}}+{{a}^{2}}-{{b}^{2}} \right)-c\left( {{a}^{2}}+{{b}^{2}}-{{c}^{2}} \right) \\ 
		& =b{{c}^{2}}-c{{a}^{2}}+b{{a}^{2}}-c{{b}^{2}}-{{b}^{3}}+{{c}^{3}} \\ 
		& =bc\left( c-b \right)-{{a}^{2}}\left( c-b \right)+\left( c-b \right)\left( {{c}^{2}}+cb+{{b}^{2}} \right) \\ 
		& =\left( c-b \right)\left( bc-{{a}^{2}}+\left( {{c}^{2}}+cb+{{b}^{2}} \right) \right) \\ 
		& =\left( c-b \right)\left( {{\left( c+b \right)}^{2}}-{{a}^{2}} \right) \\ 
		& =\left( c-b \right)\left( b+c+a \right)\left( b+c-a \right) \\ 
	\end{align*}
	
	Therefore
	\begin{align*}
		\alpha _{A}^{{{P}_{1}}}&=1-\frac{4S\left( b+c \right)\cos {{\theta }_{A}}+\left( c-b \right)\left( b+c+a \right)\left( b+c-a \right)\sin {{\theta }_{A}}}{4S\left( a+b+c \right)} \\ 
		& =1-\frac{\left( c-b \right)\left( b+c-a \right)}{4S},  
	\end{align*}
	
	And
	\begin{align*}
		\alpha _{B}^{{{P}_{1}}}& =\frac{AP}{A{{P}_{0}}}\frac{4S\alpha _{B}^{{{P}_{0}}}\cos {{\theta }_{A}}-\left( \alpha _{B}^{{{P}_{0}}}\left( {{b}^{2}}+{{c}^{2}}-{{a}^{2}} \right)+2\alpha _{C}^{{{P}_{0}}}{{b}^{2}} \right)\sin {{\theta }_{A}}}{4S} \\ 
		& =\frac{4Sb\cos {{\theta }_{A}}-\left( b\left( {{b}^{2}}+{{c}^{2}}-{{a}^{2}} \right)+2c{{b}^{2}} \right)\sin {{\theta }_{A}}}{4S\left( a+b+c \right)} \\ 
		& =\frac{4Sb\cos {{\theta }_{A}}-b\left( \left( {{b}^{2}}+{{c}^{2}}-{{a}^{2}} \right)+2cb \right)\sin {{\theta }_{A}}}{4S\left( a+b+c \right)} \\ 
		& =\frac{4Sb\cos {{\theta }_{A}}-b\left( b+c+a \right)\left( b+c-a \right)\sin {{\theta }_{A}}}{4S\left( a+b+c \right)} \\ 
		& =\frac{-b\left( b+c-a \right)}{4S},  
	\end{align*}
	
	And
	\begin{align*}
		\alpha _{C}^{{{P}_{1}}}& =\frac{AP}{A{{P}_{0}}}\frac{4S\alpha _{C}^{{{P}_{0}}}\cos {{\theta }_{A}}+\left( 2\alpha _{B}^{{{P}_{0}}}{{c}^{2}}+\alpha _{C}^{{{P}_{0}}}\left( {{b}^{2}}+{{c}^{2}}-{{a}^{2}} \right) \right)\sin {{\theta }_{A}}}{4S} \\ 
		& =\frac{4Sc\cos {{\theta }_{A}}+\left( 2b{{c}^{2}}+c\left( {{b}^{2}}+{{c}^{2}}-{{a}^{2}} \right) \right)\sin {{\theta }_{A}}}{4S\left( a+b+c \right)} \\ 
		& =\frac{4Sc\cos {{\theta }_{A}}+c\left( \left( {{b}^{2}}+{{c}^{2}}-{{a}^{2}} \right)+2cb \right)\sin {{\theta }_{A}}}{4S\left( a+b+c \right)} \\ 
		& =\frac{4Sc\cos {{\theta }_{A}}+c\left( b+c+a \right)\left( b+c-a \right)\sin {{\theta }_{A}}}{4S\left( a+b+c \right)} \\ 
		& =\frac{c\left( b+c-a \right)}{4S}.  
	\end{align*}
	
	2. Establish a Cartesian coordinate system (see Figure \ref{fig:Bianzhouzuobiaoxi_A}), and based on the previous corollary and the frame component formula of incenter, let ${{\theta }_{A}}=60{}^\circ $, ${AP}/{A{{P}_{0}}}\;=2$, it is obtained that:
	\begin{align*}
		\alpha _{A}^{{{P}_{1}}}& =\alpha _{A}^{P}=1-\frac{AP}{A{{P}_{0}}}\frac{4S\left( b+c \right)\cos {{\theta }_{A}}+\left( c-b \right)\left( b+c+a \right)\left( b+c-a \right)\sin {{\theta }_{A}}}{4S\left( a+b+c \right)} \\ 
		& =1-2\frac{4S\left( b+c \right)\cos 60{}^\circ +\left( c-b \right)\left( b+c+a \right)\left( b+c-a \right)\sin 60{}^\circ }{4S\left( a+b+c \right)} \\ 
		& =1-2\frac{4S\left( b+c \right)\frac{1}{2}+\left( c-b \right)\left( b+c+a \right)\left( b+c-a \right)\frac{\sqrt{3}}{2}}{4S\left( a+b+c \right)} \\ 
		& =1-\frac{4S\left( b+c \right)+\sqrt{3}\left( c-b \right)\left( b+c+a \right)\left( b+c-a \right)}{4S\left( a+b+c \right)},  
	\end{align*}
	
	And
	\begin{align*}
		\alpha _{B}^{{{P}_{1}}}& =\alpha _{B}^{P}=\frac{AP}{A{{P}_{0}}}\frac{4Sb\cos {{\theta }_{A}}-b\left( b+c+a \right)\left( b+c-a \right)\sin {{\theta }_{A}}}{4S\left( a+b+c \right)} \\ 
		& =2\frac{4Sb\cos 60{}^\circ -b\left( b+c+a \right)\left( b+c-a \right)\sin 60{}^\circ }{4S\left( a+b+c \right)} \\ 
		& =\frac{4Sb-\sqrt{3}b\left( b+c+a \right)\left( b+c-a \right)}{4S\left( a+b+c \right)},  
	\end{align*}
	
	And
	\begin{align*}
		\alpha _{C}^{{{P}_{1}}}& =\alpha _{C}^{P}=\frac{AP}{A{{P}_{0}}}\frac{4Sc\cos {{\theta }_{A}}+c\left( b+c+a \right)\left( b+c-a \right)\sin {{\theta }_{A}}}{4S\left( a+b+c \right)} \\ 
		& =2\frac{4Sc\cos 60{}^\circ +c\left( b+c+a \right)\left( b+c-a \right)\sin 60{}^\circ }{4S\left( a+b+c \right)} \\ 
		& =\frac{4Sc+\sqrt{3}c\left( b+c+a \right)\left( b+c-a \right)}{4S\left( a+b+c \right)}.  
	\end{align*}
\end{solution}
\hfill $\diamond$\par

The above example shows that as long as the position relationship between a point and a known point on the plane is known, the frame component of any point on the plane can be obtained by using the frame component of the known point, and the frame component is represented by the lengths of the three sides of the triangle.

Intercenter Geometry usually does not use Cartesian coordinate systems, which are introduced here only to determine the positional relationship between the point to be solved and the known point.

%
%
%
%
%

\section{Calculate frame components based on coordinates}\label{GenjuZuobiaoJisuanBiaojiafenliang}	
The theorems and inferences studied earlier indicate that if the frame component of a point is known, the frame component of another point can be obtained, as long as the position relationship between this point and the known point is determined. Determining the positional relationship between two points requires a reference frame. For planar problems, polar and Cartesian coordinate systems are commonly used reference frames.

In the Edge-Axis coordinate system, using the notations commonly used in polar coordinates, let ${{\rho }_{AP}}:=AP$.  Number pair $\left( {{\rho }_{AP}},{{\theta }_{A}} \right)$ is the polar coordinate system $\mathcal{P}\left( AB \right)$ of $\triangle ABC$. ${{\rho }_{AP}}$ and ${{\theta }_{A}}$ are called the radial coordinate and polar angle (angular coordinate), respectively.

Similarly, $\left( {{\rho }_{BP}},{{\theta }_{B}} \right)$ and $\left( {{\rho }_{CP}},{{\theta }_{C}} \right)$ can be defined as the polar coordinates of $\mathcal{P}\left( BC \right)$ and $\mathcal{P}\left( CA \right)$, respectively

In the Edge-Axis coordinate system of $\triangle ABC$, polar coordinates are closely related to Cartesian coordinates.

If point $P$ coincides with a vertex of a triangle, an important special case will be obtained.

Based on the results obtained earlier, the following corollary can be derived.

%
%
%

\begin{corollary}{Calculate frame component by polar coordinates-Form 1, Daiyuan Zhang}{GenjuJizuobiaoJisuanBiaojiafenliang_Xingshi1}\label{GenjuJizuobiaoJisuanBiaojiafenliang_Xingshi1}
	Given a $\triangle ABC$ and a point $P$ on the $\triangle ABC$ plane.
	
	1. If the polar coordinates of point $P$ in the polar coordinate system $\mathcal{P}\left( AB \right)$ are $\left( {{\rho }_{AP}},{{\theta }_{A}} \right)$, then the frame component of the center $P$ is:
	\[\alpha _{A}^{P}=1-\frac{{{\rho }_{AP}}}{c}\frac{\left( \sin C\sin {{\theta }_{A}}+\sin B\sin \left( A-{{\theta }_{A}} \right) \right)}{\sin A\sin B},\]
	\[\alpha _{B}^{P}=\frac{{{\rho }_{AP}}}{c}\frac{\sin \left( A-{{\theta }_{A}} \right)}{\sin A},\]
	\[\alpha _{C}^{P}=\frac{{{\rho }_{AP}}}{c}\frac{\sin C\sin {{\theta }_{A}}}{\sin A\sin B}.\]
	
	2. If the polar coordinates of point $P$ in the polar coordinate system $\mathcal{P}\left( BC \right)$ are $\left( {{\rho }_{BP}},{{\theta }_{B}} \right)$, then the frame component of the center $P$ is:
	\[\alpha _{B}^{P}=1-\frac{{{\rho }_{BP}}}{a}\frac{\left( \sin A\sin {{\theta }_{B}}+\sin C\sin \left( B-{{\theta }_{B}} \right) \right)}{\sin B\sin C},\]
	\[\alpha _{C}^{P}=\frac{{{\rho }_{BP}}}{a}\frac{\sin \left( B-{{\theta }_{B}} \right)}{\sin B},\]
	\[\alpha _{A}^{P}=\frac{{{\rho }_{BP}}}{a}\frac{\sin A\sin {{\theta }_{B}}}{\sin C\sin B}.\]
	
	3. If the polar coordinates of point $P$ in the polar coordinate system $\mathcal{P}\left( CA \right)$ are $\left( {{\rho }_{CP}},{{\theta }_{C}} \right)$, then the frame component of the center $P$ is:
	\[\alpha _{C}^{P}=1-\frac{{{\rho }_{CP}}}{b}\frac{\left( \sin B\sin {{\theta }_{C}}+\sin A\sin \left( C-{{\theta }_{C}} \right) \right)}{\sin C\sin A},\]
	\[\alpha _{A}^{P}=\frac{{{\rho }_{CP}}}{b}\frac{\sin \left( C-{{\theta }_{C}} \right)}{\sin C},\]
	\[\alpha _{B}^{P}=\frac{{{\rho }_{CP}}}{b}\frac{\sin B\sin {{\theta }_{C}}}{\sin C\sin A}.\]
\end{corollary}

\begin{proof}
	Let's consider the first scenario. Select the Edge-Axis coordinate system ${{\mathcal{T}}_{A}}$ of vertex $A$, and use the polar coordinate system $\mathcal{P}\left( AB \right)$ in the Edge-Axis coordinate system ${{\mathcal{T}}_{A}}$, let ${{P}_{0}}=B$, then $A{{P}_{0}}=AB=c$, according to theorem \ref{thm:SanjiaoxingDingdianDeBiaojiafenliang} and corollary \ref{cor:GenjuYizhidianDeBiaojiafenliangQiuWeizhidianDeBiaojiafenliangGongshi_SanjiaohanshuXingshi}, it is obtained that:
	\begin{align*}
		\alpha _{A}^{P}& =1-\frac{AP}{A{{P}_{0}}}\frac{\left( \begin{aligned}
				& \left( \alpha _{B}^{{{P}_{0}}}{{\sin }^{2}}C-\alpha _{C}^{{{P}_{0}}}{{\sin }^{2}}B \right)\sin {{\theta }_{A}} \\ 
				& +\sin B\sin C\left( \alpha _{B}^{{{P}_{0}}}\sin \left( A-{{\theta }_{A}} \right)+\alpha _{C}^{{{P}_{0}}}\sin \left( A+{{\theta }_{A}} \right) \right) \\ 
			\end{aligned} \right)}{\sin A\sin B\sin C} \\ 
		& =1-\frac{{{\rho }_{AP}}}{c}\frac{\left( {{\sin }^{2}}C\sin {{\theta }_{A}}+\sin B\sin C\sin \left( A-{{\theta }_{A}} \right) \right)}{\sin A\sin B\sin C} \\ 
		& =1-\frac{{{\rho }_{AP}}}{c}\frac{\left( \sin C\sin {{\theta }_{A}}+\sin B\sin \left( A-{{\theta }_{A}} \right) \right)}{\sin A\sin B},  
	\end{align*}
	\begin{align*}
		\alpha _{B}^{P}& =\frac{AP}{A{{P}_{0}}}\frac{\alpha _{B}^{{{P}_{0}}}\sin C\sin \left( A-{{\theta }_{A}} \right)-\alpha _{C}^{{{P}_{0}}}\sin B\sin {{\theta }_{A}}}{\sin C\sin A} \\ 
		& =\frac{{{\rho }_{AP}}}{c}\frac{\sin C\sin \left( A-{{\theta }_{A}} \right)}{\sin C\sin A}=\frac{{{\rho }_{AP}}}{c}\frac{\sin \left( A-{{\theta }_{A}} \right)}{\sin A},  
	\end{align*}
	\[\alpha _{C}^{P}=\frac{AP}{A{{P}_{0}}}\frac{\sin C\sin {{\theta }_{A}}}{\sin A\sin B}=\frac{{{\rho }_{AP}}}{c}\frac{\sin C\sin {{\theta }_{A}}}{\sin A\sin B}.\]
	
	Let's consider the second scenario. Select the Edge-Axis coordinate system ${{\mathcal{T}}_{B}}$ of vertex $B$, and use the polar coordinate system $\mathcal{P}\left( BC \right)$ in the Edge-Axis coordinate system ${{\mathcal{T}}_{B}}$, let ${{P}_{0}}=C$, then $B{{P}_{0}}=BC=a$, according to theorem \ref{thm:SanjiaoxingDingdianDeBiaojiafenliang} and corollary \ref{cor:GenjuYizhidianDeBiaojiafenliangQiuWeizhidianDeBiaojiafenliangGongshi_SanjiaohanshuXingshi}, it is obtained that:
	\begin{align*}
		\alpha _{B}^{P}& =1-\frac{BP}{B{{P}_{0}}}\frac{\left( \begin{aligned}
				& \left( \alpha _{C}^{{{P}_{0}}}{{\sin }^{2}}A-\alpha _{A}^{{{P}_{0}}}{{\sin }^{2}}C \right)\sin {{\theta }_{B}} \\ 
				& +\sin C\sin A\left( \alpha _{C}^{{{P}_{0}}}\sin \left( B-{{\theta }_{B}} \right)+\alpha _{A}^{{{P}_{0}}}\sin \left( B+{{\theta }_{B}} \right) \right) \\ 
			\end{aligned} \right)}{\sin B\sin C\sin A} \\ 
		& =1-\frac{{{\rho }_{BP}}}{a}\frac{\left( {{\sin }^{2}}A\sin {{\theta }_{B}}+\sin C\sin A\sin \left( B-{{\theta }_{B}} \right) \right)}{\sin B\sin C\sin A} \\ 
		& =1-\frac{{{\rho }_{BP}}}{a}\frac{\left( \sin A\sin {{\theta }_{B}}+\sin C\sin \left( B-{{\theta }_{B}} \right) \right)}{\sin B\sin C},  
	\end{align*}
	\begin{align*}
		\alpha _{C}^{P}& =\frac{BP}{B{{P}_{0}}}\frac{\alpha _{C}^{{{P}_{0}}}\sin A\sin \left( B-{{\theta }_{B}} \right)-\alpha _{A}^{{{P}_{0}}}\sin C\sin {{\theta }_{B}}}{\sin A\sin B} \\ 
		& =\frac{{{\rho }_{BP}}}{a}\frac{\sin A\sin \left( B-{{\theta }_{B}} \right)}{\sin A\sin B}=\frac{{{\rho }_{BP}}}{a}\frac{\sin \left( B-{{\theta }_{B}} \right)}{\sin B},  
	\end{align*}
	\begin{align*}
		\alpha _{A}^{P}& =\frac{BP}{B{{P}_{0}}}\frac{\alpha _{C}^{{{P}_{0}}}\sin A\sin {{\theta }_{B}}+\alpha _{A}^{{{P}_{0}}}\sin C\sin \left( B+{{\theta }_{B}} \right)}{\sin C\sin B} \\ 
		& =\frac{{{\rho }_{BP}}}{a}\frac{\sin A\sin {{\theta }_{B}}}{\sin C\sin B}.  
	\end{align*}
	
	Finally, consider the third scenario. Select the Edge-Axis coordinate system $\mathcal{P}\left( CA \right)$ of vertex $C$, and use the polar coordinate system $\mathcal{P}\left( CA \right)$ in the Edge-Axis coordinate system ${{\mathcal{T}}_{C}}$, let ${{P}_{0}}=A$, then $C{{P}_{0}}=CA=b$, according to theorem \ref{thm:SanjiaoxingDingdianDeBiaojiafenliang} and corollary \ref{cor:GenjuYizhidianDeBiaojiafenliangQiuWeizhidianDeBiaojiafenliangGongshi_SanjiaohanshuXingshi}, it is obtained that:
	\begin{align*}
		\alpha _{C}^{P}& =1-\frac{CP}{C{{P}_{0}}}\frac{\left( \begin{aligned}
				& \left( \alpha _{A}^{{{P}_{0}}}{{\sin }^{2}}B-\alpha _{B}^{{{P}_{0}}}{{\sin }^{2}}A \right)\sin {{\theta }_{C}} \\ 
				& +\sin A\sin B\left( \alpha _{A}^{{{P}_{0}}}\sin \left( C-{{\theta }_{C}} \right)+\alpha _{B}^{{{P}_{0}}}\sin \left( C+{{\theta }_{C}} \right) \right) \\ 
			\end{aligned} \right)}{\sin A\sin B\sin C} \\ 
		& =1-\frac{{{\rho }_{CP}}}{b}\frac{\left( \sin B\sin {{\theta }_{C}}+\sin A\sin \left( C-{{\theta }_{C}} \right) \right)}{\sin C\sin A},  
	\end{align*}
	\begin{align*}
		\alpha _{A}^{P}& =\frac{CP}{C{{P}_{0}}}\frac{\alpha _{A}^{{{P}_{0}}}\sin B\sin \left( C-{{\theta }_{C}} \right)-\alpha _{B}^{{{P}_{0}}}\sin A\sin {{\theta }_{C}}}{\sin B\sin C} \\ 
		& =\frac{{{\rho }_{CP}}}{b}\frac{\sin \left( C-{{\theta }_{C}} \right)}{\sin C},  
	\end{align*}
	\begin{align*}
		& \alpha _{B}^{P}=\frac{CP}{C{{P}_{0}}}\frac{\alpha _{A}^{{{P}_{0}}}\sin B\sin {{\theta }_{C}}+\alpha _{B}^{{{P}_{0}}}\sin A\sin \left( C+{{\theta }_{C}} \right)}{\sin C\sin A} \\ 
		& =\frac{{{\rho }_{CP}}}{b}\frac{\sin B\sin {{\theta }_{C}}}{\sin C\sin A}.  
	\end{align*}
\end{proof}
\hfill $\square$\par

As an application, let's take an example.


\begin{example}{}\label{QiuzhengYixiaLiangshi}
	Prove the following two equations:
	\[\frac{\sin {{\theta }_{A}}}{\sin \left( A-{{\theta }_{A}} \right)}\frac{\sin {{\theta }_{B}}}{\sin \left( B-{{\theta }_{B}} \right)}\frac{\sin {{\theta }_{C}}}{\sin \left( C-{{\theta }_{C}} \right)}=1,\]
	\[\alpha _{A}^{P}\alpha _{B}^{P}\alpha _{C}^{P}=\frac{{{\rho }_{AP}}{{\rho }_{BP}}{{\rho }_{CP}}}{abc}\frac{\sin {{\theta }_{A}}\sin {{\theta }_{B}}\sin {{\theta }_{C}}}{\sin A\sin B\sin C}.\]
\end{example}

\begin{solution}
	Obviously, regardless of wherever the ${{P}_{0}}$ is located, it will not affect the frame component of point $P$. Choose the first set of frame components from the above corollary:
	\[\alpha _{B}^{P}=\frac{{{\rho }_{AP}}}{c}\frac{\sin \left( A-{{\theta }_{A}} \right)}{\sin A},\]
	\[\alpha _{C}^{P}=\frac{{{\rho }_{AP}}}{c}\frac{\sin C\sin {{\theta }_{A}}}{\sin A\sin B}.\]
	
	Therefore
	\[\frac{\alpha _{C}^{P}}{\alpha _{B}^{P}}=\frac{\frac{{{\rho }_{AP}}}{c}\frac{\sin C\sin {{\theta }_{A}}}{\sin A\sin B}}{\frac{{{\rho }_{AP}}}{c}\frac{\sin \left( A-{{\theta }_{A}} \right)}{\sin A}}=\frac{\sin {{\theta }_{A}}}{\sin \left( A-{{\theta }_{A}} \right)}\frac{\sin C}{\sin B}.\]
	
	Choose the second set of frame components from the above corollary:
	\[\alpha _{C}^{P}=\frac{{{\rho }_{BP}}}{a}\frac{\sin \left( B-{{\theta }_{B}} \right)}{\sin B},\]
	\[\alpha _{A}^{P}=\frac{{{\rho }_{BP}}}{a}\frac{\sin A\sin {{\theta }_{B}}}{\sin C\sin B}.\]
	
	Therefore
	\[\frac{\alpha _{A}^{P}}{\alpha _{C}^{P}}=\frac{\frac{{{\rho }_{BP}}}{a}\frac{\sin A\sin {{\theta }_{B}}}{\sin C\sin B}}{\frac{{{\rho }_{BP}}}{a}\frac{\sin \left( B-{{\theta }_{B}} \right)}{\sin B}}=\frac{\sin {{\theta }_{B}}}{\sin \left( B-{{\theta }_{B}} \right)}\frac{\sin A}{\sin C}.\]
	
	Choose the third set of frame components from the above corollary:
	\[\alpha _{A}^{P}=\frac{{{\rho }_{CP}}}{b}\frac{\sin \left( C-{{\theta }_{C}} \right)}{\sin C},\]
	\[\alpha _{B}^{P}=\frac{{{\rho }_{CP}}}{b}\frac{\sin B\sin {{\theta }_{C}}}{\sin C\sin A}.\]
	
	Therefore
	\[\frac{\alpha _{B}^{P}}{\alpha _{A}^{P}}=\frac{\frac{{{\rho }_{CP}}}{b}\frac{\sin B\sin {{\theta }_{C}}}{\sin C\sin A}}{\frac{{{\rho }_{CP}}}{b}\frac{\sin \left( C-{{\theta }_{C}} \right)}{\sin C}}=\frac{\sin {{\theta }_{C}}}{\sin \left( C-{{\theta }_{C}} \right)}\frac{\sin B}{\sin A}.\]
	
	Thus
	\begin{align*}
		\frac{\alpha _{C}^{P}}{\alpha _{B}^{P}}\frac{\alpha _{A}^{P}}{\alpha _{C}^{P}}\frac{\alpha _{B}^{P}}{\alpha _{A}^{P}}& =\frac{\sin {{\theta }_{A}}}{\sin \left( A-{{\theta }_{A}} \right)}\frac{\sin C}{\sin B}\frac{\sin {{\theta }_{B}}}{\sin \left( B-{{\theta }_{B}} \right)}\frac{\sin A}{\sin C}\frac{\sin {{\theta }_{C}}}{\sin \left( C-{{\theta }_{C}} \right)}\frac{\sin B}{\sin A} \\ 
		& =\frac{\sin {{\theta }_{A}}}{\sin \left( A-{{\theta }_{A}} \right)}\frac{\sin {{\theta }_{B}}}{\sin \left( B-{{\theta }_{B}} \right)}\frac{\sin {{\theta }_{C}}}{\sin \left( C-{{\theta }_{C}} \right)}\text{.}  
	\end{align*}
	
	The left side of the above equation is 1, thus obtaining:
	\[\frac{\sin {{\theta }_{A}}}{\sin \left( A-{{\theta }_{A}} \right)}\frac{\sin {{\theta }_{B}}}{\sin \left( B-{{\theta }_{B}} \right)}\frac{\sin {{\theta }_{C}}}{\sin \left( C-{{\theta }_{C}} \right)}=1.\]
	
	Choose another set of frame components:
	\[\alpha _{A}^{P}=\frac{{{\rho }_{CP}}}{b}\frac{\sin \left( C-{{\theta }_{C}} \right)}{\sin C},\]
	\[\alpha _{B}^{P}=\frac{{{\rho }_{AP}}}{c}\frac{\sin \left( A-{{\theta }_{A}} \right)}{\sin A},\]
	\[\alpha _{C}^{P}=\frac{{{\rho }_{BP}}}{a}\frac{\sin \left( B-{{\theta }_{B}} \right)}{\sin B}.\]
	
	Therefore
	\begin{align*}
		\alpha _{A}^{P}\alpha _{B}^{P}\alpha _{C}^{P}& =\frac{{{\rho }_{CP}}}{b}\frac{\sin \left( C-{{\theta }_{C}} \right)}{\sin C}\frac{{{\rho }_{AP}}}{c}\frac{\sin \left( A-{{\theta }_{A}} \right)}{\sin A}\frac{{{\rho }_{BP}}}{a}\frac{\sin \left( B-{{\theta }_{B}} \right)}{\sin B} \\ 
		& =\frac{{{\rho }_{AP}}{{\rho }_{BP}}{{\rho }_{CP}}}{abc}\frac{\sin \left( C-{{\theta }_{C}} \right)}{\sin C}\frac{\sin \left( A-{{\theta }_{A}} \right)}{\sin A}\frac{\sin \left( B-{{\theta }_{B}} \right)}{\sin B} \\ 
		& =\frac{{{\rho }_{AP}}{{\rho }_{BP}}{{\rho }_{CP}}}{abc}\frac{\sin {{\theta }_{A}}\sin {{\theta }_{B}}\sin {{\theta }_{C}}}{\sin A\sin B\sin C}.  
	\end{align*}
\end{solution}
\hfill $\diamond$\par

\begin{corollary}{Calculate frame component by polar coordinates-Form 2, Daiyuan Zhang}{GenjuJizuobiaoJisuanBiaojiafenliang_Xingshi2}\label{GenjuJizuobiaoJisuanBiaojiafenliang_Xingshi2}
	Given a $\triangle ABC$ and a point $P$ on the $\triangle ABC$ plane, the area of $\triangle ABC$ is $S$.
	
	1. Let $\angle BAP={{\theta }_{A}}$, $AP={{\rho }_{AP}}$, then the frame component of point $P$ is: 
	\[\alpha _{A}^{P}=1-\frac{1}{c}\left( {{\rho }_{AP}}\cos {{\theta }_{A}}+\frac{{{c}^{2}}+{{a}^{2}}-{{b}^{2}}}{4S}{{\rho }_{AP}}\sin {{\theta }_{A}} \right),\]
	\[\alpha _{B}^{P}=\frac{1}{c}\left( {{\rho }_{AP}}\cos {{\theta }_{A}}-\frac{{{b}^{2}}+{{c}^{2}}-{{a}^{2}}}{4S}{{\rho }_{AP}}\sin {{\theta }_{A}} \right),\]
	\[\alpha _{C}^{P}=\frac{c}{2S}{{\rho }_{AP}}\sin {{\theta }_{A}}.\]
	
	2. Let $\angle CBP={{\theta }_{B}}$, $BP={{\rho }_{BP}}$, then the frame component of point $P$ is: 
	\[\alpha _{B}^{P}=1-\frac{1}{a}\left( {{\rho }_{BP}}\cos {{\theta }_{B}}+\frac{{{a}^{2}}+{{b}^{2}}-{{c}^{2}}}{4S}{{\rho }_{BP}}\sin {{\theta }_{B}} \right),\]
	\[\alpha _{C}^{P}=\frac{1}{a}\left( {{\rho }_{BP}}\cos {{\theta }_{B}}-\frac{{{c}^{2}}+{{a}^{2}}-{{b}^{2}}}{4S}{{\rho }_{BP}}\sin {{\theta }_{B}} \right),\]
	\[\alpha _{A}^{P}=\frac{a}{2S}{{\rho }_{BP}}\sin {{\theta }_{B}}.\]
	
	3. Let $\angle ACP={{\theta }_{C}}$, $CP={{\rho }_{CP}}$, then the frame component of point $P$ is: 
	\[\alpha _{C}^{P}=1-\frac{1}{b}\left( {{\rho }_{CP}}\cos {{\theta }_{C}}+\frac{{{b}^{2}}+{{c}^{2}}-{{a}^{2}}}{4S}{{\rho }_{CP}}\sin {{\theta }_{C}} \right),\]
	\[\alpha _{A}^{P}=\frac{1}{b}\left( {{\rho }_{CP}}\cos {{\theta }_{C}}-\frac{{{a}^{2}}+{{b}^{2}}-{{c}^{2}}}{4S}{{\rho }_{CP}}\sin {{\theta }_{C}} \right),\]
	\[\alpha _{B}^{P}=\frac{b}{2S}{{\rho }_{CP}}\sin {{\theta }_{C}}.\]
\end{corollary}

\begin{proof}
	Let's consider the first scenario. Select the Edge-Axis coordinate system ${{\mathcal{T}}_{A}}$ of vertex $A$, and use the polar coordinate system $\mathcal{P}\left( AB \right)$ in the Edge-Axis coordinate system ${{\mathcal{T}}_{A}}$, let ${{P}_{0}}=B$, then $A{{P}_{0}}=AB=c$, according to theorem \ref{thm:SanjiaoxingDingdianDeBiaojiafenliang} and corollary \ref{cor:GenjuYizhidianDeBiaojiafenliangQiuWeizhidianDeBiaojiafenliangGongshi_BianchangXingshi}, it is obtained that:
	\begin{align*}
		\alpha _{A}^{P}& =1-\frac{AP}{A{{P}_{0}}}\frac{4S\left( \alpha _{B}^{{{P}_{0}}}+\alpha _{C}^{{{P}_{0}}} \right)\cos {{\theta }_{A}}+\left( \alpha _{B}^{{{P}_{0}}}\left( {{c}^{2}}+{{a}^{2}}-{{b}^{2}} \right)-\alpha _{C}^{{{P}_{0}}}\left( {{a}^{2}}+{{b}^{2}}-{{c}^{2}} \right) \right)\sin {{\theta }_{A}}}{4S} \\ 
		& =1-\frac{{{\rho }_{AP}}}{c}\frac{4S\cos {{\theta }_{A}}+\left( {{c}^{2}}+{{a}^{2}}-{{b}^{2}} \right)\sin {{\theta }_{A}}}{4S} \\ 
		& =1-\frac{1}{c}\left( {{\rho }_{AP}}\cos {{\theta }_{A}}+\frac{{{c}^{2}}+{{a}^{2}}-{{b}^{2}}}{4S}{{\rho }_{AP}}\sin {{\theta }_{A}} \right),  
	\end{align*}
	\begin{align*}
		\alpha _{B}^{P}& =\frac{AP}{A{{P}_{0}}}\frac{4\alpha _{B}^{{{P}_{0}}}S\cos {{\theta }_{A}}-\left( \alpha _{B}^{{{P}_{0}}}\left( {{b}^{2}}+{{c}^{2}}-{{a}^{2}} \right)+2\alpha _{C}^{{{P}_{0}}}{{b}^{2}} \right)\sin {{\theta }_{A}}}{4S} \\ 
		& =\frac{{{\rho }_{AP}}}{c}\frac{4S\cos {{\theta }_{A}}-\left( {{b}^{2}}+{{c}^{2}}-{{a}^{2}} \right)\sin {{\theta }_{A}}}{4S} \\ 
		& =\frac{1}{c}\left( {{\rho }_{AP}}\cos {{\theta }_{A}}-\frac{{{b}^{2}}+{{c}^{2}}-{{a}^{2}}}{4S}{{\rho }_{AP}}\sin {{\theta }_{A}} \right),  
	\end{align*}
	\begin{align*}
		\alpha _{C}^{P}& =\frac{AP}{A{{P}_{0}}}\frac{4S\alpha _{C}^{{{P}_{0}}}\cos {{\theta }_{A}}+\left( 2\alpha _{B}^{{{P}_{0}}}{{c}^{2}}+\alpha _{C}^{{{P}_{0}}}\left( {{b}^{2}}+{{c}^{2}}-{{a}^{2}} \right) \right)\sin {{\theta }_{A}}}{4S} \\ 
		& =\frac{{{\rho }_{AP}}}{c}\frac{2{{c}^{2}}\sin {{\theta }_{A}}}{4S}=\frac{c}{2S}{{\rho }_{AP}}\sin {{\theta }_{A}}.  
	\end{align*}
	
	Let's consider the second scenario. Select the Edge-Axis coordinate system ${{\mathcal{T}}_{B}}$ of vertex $B$, and use the polar coordinate system $\mathcal{P}\left( BC \right)$ in the Edge-Axis coordinate system ${{\mathcal{T}}_{B}}$, let ${{P}_{0}}=C$, then $B{{P}_{0}}=BC=a$, according to theorem \ref{thm:SanjiaoxingDingdianDeBiaojiafenliang} and corollary \ref{cor:GenjuYizhidianDeBiaojiafenliangQiuWeizhidianDeBiaojiafenliangGongshi_BianchangXingshi}, it is obtained that:
	\begin{align*}
		\alpha _{B}^{P}& =1-\frac{BP}{B{{P}_{0}}}\frac{4S\left( \alpha _{C}^{{{P}_{0}}}+\alpha _{A}^{{{P}_{0}}} \right)\cos {{\theta }_{B}}+\left( \alpha _{C}^{{{P}_{0}}}\left( {{a}^{2}}+{{b}^{2}}-{{c}^{2}} \right)-\alpha _{A}^{{{P}_{0}}}\left( {{b}^{2}}+{{c}^{2}}-{{a}^{2}} \right) \right)\sin {{\theta }_{B}}}{4S} \\ 
		& =1-\frac{{{\rho }_{BP}}}{a}\frac{4S\cos {{\theta }_{B}}+\left( {{a}^{2}}+{{b}^{2}}-{{c}^{2}} \right)\sin {{\theta }_{B}}}{4S} \\ 
		& =1-\frac{1}{a}\left( {{\rho }_{BP}}\cos {{\theta }_{B}}+\frac{{{a}^{2}}+{{b}^{2}}-{{c}^{2}}}{4S}{{\rho }_{BP}}\sin {{\theta }_{B}} \right),  
	\end{align*}
	\begin{align*}
		& \alpha _{C}^{P}=\frac{BP}{B{{P}_{0}}}\frac{4\alpha _{C}^{{{P}_{0}}}S\cos {{\theta }_{B}}-\left( \alpha _{C}^{{{P}_{0}}}\left( {{c}^{2}}+{{a}^{2}}-{{b}^{2}} \right)+2\alpha _{A}^{{{P}_{0}}}{{c}^{2}} \right)\sin {{\theta }_{B}}}{4S} \\ 
		& =\frac{{{\rho }_{BP}}}{a}\frac{4S\cos {{\theta }_{B}}-\left( {{c}^{2}}+{{a}^{2}}-{{b}^{2}} \right)\sin {{\theta }_{B}}}{4S} \\ 
		& =\frac{1}{a}\left( {{\rho }_{BP}}\cos {{\theta }_{B}}-\frac{{{c}^{2}}+{{a}^{2}}-{{b}^{2}}}{4S}{{\rho }_{BP}}\sin {{\theta }_{B}} \right),  
	\end{align*}
	\begin{align*}
		\alpha _{A}^{P}& =\frac{BP}{B{{P}_{0}}}\frac{4S\alpha _{A}^{{{P}_{0}}}\cos {{\theta }_{B}}+\left( 2\alpha _{C}^{{{P}_{0}}}{{a}^{2}}+\alpha _{A}^{{{P}_{0}}}\left( {{c}^{2}}+{{a}^{2}}-{{b}^{2}} \right) \right)\sin {{\theta }_{B}}}{4S} \\ 
		& =\frac{{{\rho }_{BP}}}{a}\frac{2{{a}^{2}}\sin {{\theta }_{B}}}{4S}=\frac{a}{2S}{{\rho }_{BP}}\sin {{\theta }_{B}}.  
	\end{align*}
	
	Finally, consider the third scenario. Select the Edge-Axis coordinate system $\mathcal{P}\left( CA \right)$ of vertex $C$, and use the polar coordinate system $\mathcal{P}\left( CA \right)$ in the Edge-Axis coordinate system ${{\mathcal{T}}_{C}}$, let ${{P}_{0}}=A$, then $C{{P}_{0}}=CA=b$, according to theorem \ref{thm:SanjiaoxingDingdianDeBiaojiafenliang} and corollary \ref{cor:GenjuYizhidianDeBiaojiafenliangQiuWeizhidianDeBiaojiafenliangGongshi_BianchangXingshi}, it is obtained that:
	\begin{align*}
		\alpha _{C}^{P}& =1-\frac{CP}{C{{P}_{0}}}\frac{4S\left( \alpha _{A}^{{{P}_{0}}}+\alpha _{B}^{{{P}_{0}}} \right)\cos {{\theta }_{C}}+\left( \alpha _{A}^{{{P}_{0}}}\left( {{b}^{2}}+{{c}^{2}}-{{a}^{2}} \right)-\alpha _{B}^{{{P}_{0}}}\left( {{c}^{2}}+{{a}^{2}}-{{b}^{2}} \right) \right)\sin {{\theta }_{C}}}{4S} \\ 
		& =1-\frac{{{\rho }_{CP}}}{b}\frac{4S\cos {{\theta }_{C}}+\left( {{b}^{2}}+{{c}^{2}}-{{a}^{2}} \right)\sin {{\theta }_{C}}}{4S} \\ 
		& =1-\frac{1}{b}\left( {{\rho }_{CP}}\cos {{\theta }_{C}}+\frac{{{b}^{2}}+{{c}^{2}}-{{a}^{2}}}{4S}{{\rho }_{CP}}\sin {{\theta }_{C}} \right),  
	\end{align*}
	\begin{align*}
		\alpha _{A}^{P}& =\frac{CP}{C{{P}_{0}}}\frac{4\alpha _{A}^{{{P}_{0}}}S\cos {{\theta }_{C}}-\left( \alpha _{A}^{{{P}_{0}}}\left( {{a}^{2}}+{{b}^{2}}-{{c}^{2}} \right)+2\alpha _{B}^{{{P}_{0}}}{{a}^{2}} \right)\sin {{\theta }_{C}}}{4S} \\ 
		& =\frac{{{\rho }_{CP}}}{b}\frac{4S\cos {{\theta }_{C}}-\left( {{a}^{2}}+{{b}^{2}}-{{c}^{2}} \right)\sin {{\theta }_{C}}}{4S} \\ 
		& =\frac{1}{b}\left( {{\rho }_{CP}}\cos {{\theta }_{C}}-\frac{{{a}^{2}}+{{b}^{2}}-{{c}^{2}}}{4S}{{\rho }_{CP}}\sin {{\theta }_{C}} \right),  
	\end{align*}
	\begin{align*}
		\alpha _{B}^{P}& =\frac{CP}{C{{P}_{0}}}\frac{4S\alpha _{B}^{{{P}_{0}}}\cos {{\theta }_{C}}+\left( 2\alpha _{A}^{{{P}_{0}}}{{b}^{2}}+\alpha _{B}^{{{P}_{0}}}\left( {{a}^{2}}+{{b}^{2}}-{{c}^{2}} \right) \right)\sin {{\theta }_{C}}}{4S} \\ 
		& =\frac{{{\rho }_{CP}}}{b}\frac{2{{b}^{2}}\sin {{\theta }_{C}}}{4S}=\frac{b}{2S}{{\rho }_{CP}}\sin {{\theta }_{C}}.  
	\end{align*}
\end{proof}
\hfill $\square$\par

\begin{corollary}{Find frame component by Cartesian coordinate-side length form, Daiyuan Zhang}{GenjuZhijiaozuobiaoJisuanBiaojiafenliang_BianchangXingshi}\label{GenjuZhijiaozuobiaoJisuanBiaojiafenliang_BianchangXingshi}
	Given a $\triangle ABC$ and a point $P$ on the $\triangle ABC$ plane, the area of $\triangle ABC$ is $S$.
	
	1. In the $\mathcal{R}\left( AB \right)$ coordinate system, the Cartesian coordinates of point (IC) $P$ are $\left( x_{\mathcal{R}\left( AB \right)}^{P},y_{\mathcal{R}\left( AB \right)}^{P} \right)$, and the frame component of point $P$ is: 
	\[\alpha _{A}^{P}=1-\frac{1}{c}\left( x_{\mathcal{R}\left( AB \right)}^{P}+\frac{{{c}^{2}}+{{a}^{2}}-{{b}^{2}}}{4S}y_{\mathcal{R}\left( AB \right)}^{P} \right),\]
	\[\alpha _{B}^{P}=\frac{1}{c}\left( x_{\mathcal{R}\left( AB \right)}^{P}-\frac{{{b}^{2}}+{{c}^{2}}-{{a}^{2}}}{4S}y_{\mathcal{R}\left( AB \right)}^{P} \right),\]
	\[\alpha _{C}^{P}=\frac{c}{2S}y_{\mathcal{R}\left( AB \right)}^{P}.\]
	
	2. In the $\mathcal{R}\left( BC \right)$ coordinate system, the Cartesian coordinates of point (IC) $P$ are $\left( x_{\mathcal{R}\left( BC \right)}^{P},y_{\mathcal{R}\left( BC \right)}^{P} \right)$, and the frame component of point $P$ is: 
	\[\alpha _{B}^{P}=1-\frac{1}{a}\left( x_{\mathcal{R}\left( BC \right)}^{P}+\frac{{{a}^{2}}+{{b}^{2}}-{{c}^{2}}}{4S}y_{\mathcal{R}\left( BC \right)}^{P} \right),\]
	\[\alpha _{C}^{P}=\frac{1}{a}\left( x_{\mathcal{R}\left( BC \right)}^{P}-\frac{{{c}^{2}}+{{a}^{2}}-{{b}^{2}}}{4S}y_{\mathcal{R}\left( BC \right)}^{P} \right),\]
	\[\alpha _{A}^{P}=\frac{a}{2S}y_{\mathcal{R}\left( BC \right)}^{P}.\]
	
	3. In the $\mathcal{R}\left( CA \right)$ coordinate system, the Cartesian coordinates of point (IC) $P$ are $\left( x_{\mathcal{R}\left( CA \right)}^{P},y_{\mathcal{R}\left( CA \right)}^{P} \right)$, and the frame component of point $P$ is:
	\[\alpha _{C}^{P}=1-\frac{1}{b}\left( x_{\mathcal{R}\left( CA \right)}^{P}+\frac{{{b}^{2}}+{{c}^{2}}-{{a}^{2}}}{4S}y_{\mathcal{R}\left( CA \right)}^{P} \right),\]
	\[\alpha _{A}^{P}=\frac{1}{b}\left( x_{\mathcal{R}\left( CA \right)}^{P}-\frac{{{a}^{2}}+{{b}^{2}}-{{c}^{2}}}{4S}y_{\mathcal{R}\left( CA \right)}^{P} \right),\]
	\[\alpha _{B}^{P}=\frac{b}{2S}y_{\mathcal{R}\left( CA \right)}^{P}.\]
\end{corollary}


\begin{proof}
	In the Edge-Axis coordinate system, the Cartesian coordinate system and the polar coordinate system are related, and this corollary is directly obtained from corollary \ref{cor:GenjuJizuobiaoJisuanBiaojiafenliang_Xingshi1}.   
\end{proof}
\hfill $\square$\par

\begin{theorem}{Calculate frame component by the polar angle, Daiyuan Zhang}{GenjuJijiaoJisuanBiaojiafenliang}\label{GenjuJijiaoJisuanBiaojiafenliang}
	Given a $\triangle ABC$, the point $P$ is on the plane where $\triangle ABC$ is located, and the polar angles of point $P$ are ${{\theta }_{A}}$, ${{\theta }_{B}}$, ${{\theta }_{C}}$, then the frame components of point $P$ are:
	\[\alpha _{A}^{P}=\frac{{{f}_{a}}}{{{f}_{a}}+{{f}_{b}}+{{f}_{c}}},\]
	\[\alpha _{B}^{P}=\frac{{{f}_{b}}}{{{f}_{a}}+{{f}_{b}}+{{f}_{c}}},\]
	\[\alpha _{C}^{P}=\frac{{{f}_{c}}}{{{f}_{a}}+{{f}_{b}}+{{f}_{c}}},\]
	
	where
	\[{{f}_{a}}=\sin A\left( \sin {{\theta }_{A}}\sin {{\theta }_{B}}+\sin {{\theta }_{B}}\sin \left( C-{{\theta }_{C}} \right)+\sin \left( C-{{\theta }_{C}} \right)\sin \left( A-{{\theta }_{A}} \right) \right),\]
	\[{{f}_{b}}=\sin B\left( \sin {{\theta }_{B}}\sin {{\theta }_{C}}+\sin {{\theta }_{C}}\sin \left( A-{{\theta }_{A}} \right)+\sin \left( A-{{\theta }_{A}} \right)\sin \left( B-{{\theta }_{B}} \right) \right),\]
	\[{{f}_{c}}=\sin C\left( \sin {{\theta }_{C}}\sin {{\theta }_{A}}+\sin {{\theta }_{A}}\sin \left( B-{{\theta }_{B}} \right)+\sin \left( B-{{\theta }_{B}} \right)\sin \left( C-{{\theta }_{C}} \right) \right);\]
	
	or
	\[{{f}_{a}}=a\left( \sin {{\theta }_{A}}\sin {{\theta }_{B}}+\sin {{\theta }_{B}}\sin \left( C-{{\theta }_{C}} \right)+\sin \left( C-{{\theta }_{C}} \right)\sin \left( A-{{\theta }_{A}} \right) \right),\]
	\[{{f}_{b}}=b\left( \sin {{\theta }_{B}}\sin {{\theta }_{C}}+\sin {{\theta }_{C}}\sin \left( A-{{\theta }_{A}} \right)+\sin \left( A-{{\theta }_{A}} \right)\sin \left( B-{{\theta }_{B}} \right) \right),\]
	\[{{f}_{c}}=c\left( \sin {{\theta }_{C}}\sin {{\theta }_{A}}+\sin {{\theta }_{A}}\sin \left( B-{{\theta }_{B}} \right)+\sin \left( B-{{\theta }_{B}} \right)\sin \left( C-{{\theta }_{C}} \right) \right).\]
\end{theorem}

\begin{proof}
	According to theorem \ref{thm:Thm6.4.1}, corollary \ref{cor:GenjuJizuobiaoJisuanBiaojiafenliang_Xingshi1} and sine theorem, select the first set of frame components in corollary \ref{cor:GenjuJizuobiaoJisuanBiaojiafenliang_Xingshi1}: 
	\[\alpha _{B}^{P}=\frac{{{\rho }_{AP}}}{c}\frac{\sin \left( A-{{\theta }_{A}} \right)}{\sin A},\]
	\[\alpha _{C}^{P}=\frac{{{\rho }_{AP}}}{c}\frac{\sin C\sin {{\theta }_{A}}}{\sin A\sin B}.\]
	
	Therefore
	\[\lambda _{BC}^{P}=\frac{\alpha _{C}^{P}}{\alpha _{B}^{P}}=\frac{\frac{{{\rho }_{AP}}}{c}\frac{\sin C\sin {{\theta }_{A}}}{\sin A\sin B}}{\frac{{{\rho }_{AP}}}{c}\frac{\sin \left( A-{{\theta }_{A}} \right)}{\sin A}}=\frac{\sin {{\theta }_{A}}}{\sin \left( A-{{\theta }_{A}} \right)}\frac{\sin C}{\sin B}=\frac{c\sin {{\theta }_{A}}}{b\sin \left( A-{{\theta }_{A}} \right)}.\]
	
	Select the second set of frame components in corollary \ref{cor:GenjuJizuobiaoJisuanBiaojiafenliang_Xingshi1}:
	\[\alpha _{C}^{P}=\frac{{{\rho }_{BP}}}{a}\frac{\sin \left( B-{{\theta }_{B}} \right)}{\sin B},\]
	\[\alpha _{A}^{P}=\frac{{{\rho }_{BP}}}{a}\frac{\sin A\sin {{\theta }_{B}}}{\sin C\sin B}.\]
	
	Therefore
	\[\lambda _{CA}^{P}=\frac{\alpha _{A}^{P}}{\alpha _{C}^{P}}=\frac{\frac{{{\rho }_{BP}}}{a}\frac{\sin A\sin {{\theta }_{B}}}{\sin C\sin B}}{\frac{{{\rho }_{BP}}}{a}\frac{\sin \left( B-{{\theta }_{B}} \right)}{\sin B}}=\frac{\sin {{\theta }_{B}}}{\sin \left( B-{{\theta }_{B}} \right)}\frac{\sin A}{\sin C}=\frac{a\sin {{\theta }_{B}}}{c\sin \left( B-{{\theta }_{B}} \right)}.\]
	
	Select the third set of frame components in corollary \ref{cor:GenjuJizuobiaoJisuanBiaojiafenliang_Xingshi1}:
	\[\alpha _{A}^{P}=\frac{{{\rho }_{CP}}}{b}\frac{\sin \left( C-{{\theta }_{C}} \right)}{\sin C},\]
	\[\alpha _{B}^{P}=\frac{{{\rho }_{CP}}}{b}\frac{\sin B\sin {{\theta }_{C}}}{\sin C\sin A}.\]
	
	Therefore
	\[\lambda _{AB}^{P}=\frac{\alpha _{B}^{P}}{\alpha _{A}^{P}}=\frac{\frac{{{\rho }_{CP}}}{b}\frac{\sin B\sin {{\theta }_{C}}}{\sin C\sin A}}{\frac{{{\rho }_{CP}}}{b}\frac{\sin \left( C-{{\theta }_{C}} \right)}{\sin C}}=\frac{\sin {{\theta }_{C}}}{\sin \left( C-{{\theta }_{C}} \right)}\frac{\sin B}{\sin A}=\frac{b\sin {{\theta }_{C}}}{a\sin \left( C-{{\theta }_{C}} \right)}.\]
	
	According to theorem \ref{thm:Thm6.1.1}, it is obtained that:
	\begin{align*}
		\alpha _{A}^{P}& =\frac{1}{1+\lambda _{AB}^{P}+\lambda _{AC}^{P}}=\frac{1}{1+\frac{b\sin {{\theta }_{C}}}{a\sin \left( C-{{\theta }_{C}} \right)}+\frac{c\sin \left( B-{{\theta }_{B}} \right)}{a\sin {{\theta }_{B}}}} \\ 
		& =\frac{a\sin {{\theta }_{B}}\sin \left( C-{{\theta }_{C}} \right)}{a\sin {{\theta }_{B}}\sin \left( C-{{\theta }_{C}} \right)+b\sin {{\theta }_{B}}\sin {{\theta }_{C}}+c\sin \left( B-{{\theta }_{B}} \right)\sin \left( C-{{\theta }_{C}} \right)},  
	\end{align*}
	\begin{align*}
		\alpha _{A}^{P}& =\frac{\lambda _{BA}^{P}}{1+\lambda _{BC}^{P}+\lambda _{BA}^{P}}=\frac{\frac{a\sin \left( C-{{\theta }_{C}} \right)}{b\sin {{\theta }_{C}}}}{1+\frac{c\sin {{\theta }_{A}}}{b\sin \left( A-{{\theta }_{A}} \right)}+\frac{a\sin \left( C-{{\theta }_{C}} \right)}{b\sin {{\theta }_{C}}}} \\ 
		& =\frac{a\sin \left( C-{{\theta }_{C}} \right)\sin \left( A-{{\theta }_{A}} \right)}{b\sin {{\theta }_{C}}\sin \left( A-{{\theta }_{A}} \right)+c\sin {{\theta }_{C}}\sin {{\theta }_{A}}+a\sin \left( C-{{\theta }_{C}} \right)\sin \left( A-{{\theta }_{A}} \right)},  
	\end{align*}
	\begin{align*}
		\alpha _{A}^{P}& =\frac{\lambda _{CA}^{P}}{1+\lambda _{CA}^{P}+\lambda _{CB}^{P}}=\frac{\frac{a\sin {{\theta }_{B}}}{c\sin \left( B-{{\theta }_{B}} \right)}}{1+\frac{a\sin {{\theta }_{B}}}{c\sin \left( B-{{\theta }_{B}} \right)}+\frac{b\sin \left( A-{{\theta }_{A}} \right)}{c\sin {{\theta }_{A}}}} \\ 
		& =\frac{a\sin {{\theta }_{A}}\sin {{\theta }_{B}}}{c\sin {{\theta }_{A}}\sin \left( B-{{\theta }_{B}} \right)+a\sin {{\theta }_{A}}\sin {{\theta }_{B}}+b\sin \left( A-{{\theta }_{A}} \right)\sin \left( B-{{\theta }_{B}} \right)}.  
	\end{align*}
	
	Sum up the numerator and denominator separately to obtain (symmetric form):
	\[\alpha _{A}^{P}=\frac{a\sin {{\theta }_{A}}\sin {{\theta }_{B}}+a\sin {{\theta }_{B}}\sin \left( C-{{\theta }_{C}} \right)+a\sin \left( C-{{\theta }_{C}} \right)\sin \left( A-{{\theta }_{A}} \right)}{\left( \begin{aligned}
			& c\sin {{\theta }_{A}}\sin \left( B-{{\theta }_{B}} \right)+a\sin {{\theta }_{A}}\sin {{\theta }_{B}}+b\sin \left( A-{{\theta }_{A}} \right)\sin \left( B-{{\theta }_{B}} \right) \\ 
			& +a\sin {{\theta }_{B}}\sin \left( C-{{\theta }_{C}} \right)+b\sin {{\theta }_{B}}\sin {{\theta }_{C}}+c\sin \left( B-{{\theta }_{B}} \right)\sin \left( C-{{\theta }_{C}} \right) \\ 
			& +b\sin {{\theta }_{C}}\sin \left( A-{{\theta }_{A}} \right)+c\sin {{\theta }_{C}}\sin {{\theta }_{A}}+a\sin \left( C-{{\theta }_{C}} \right)\sin \left( A-{{\theta }_{A}} \right) \\ 
		\end{aligned} \right)},\]
	i.e.
	\[\alpha _{A}^{P}=\frac{a\left( \sin {{\theta }_{A}}\sin {{\theta }_{B}}+\sin {{\theta }_{B}}\sin \left( C-{{\theta }_{C}} \right)+\sin \left( C-{{\theta }_{C}} \right)\sin \left( A-{{\theta }_{A}} \right) \right)}{\left( \begin{aligned}
			& a\left( \sin {{\theta }_{A}}\sin {{\theta }_{B}}+\sin {{\theta }_{B}}\sin \left( C-{{\theta }_{C}} \right)+\sin \left( C-{{\theta }_{C}} \right)\sin \left( A-{{\theta }_{A}} \right) \right) \\ 
			& +b\left( \sin {{\theta }_{B}}\sin {{\theta }_{C}}+\sin {{\theta }_{C}}\sin \left( A-{{\theta }_{A}} \right)+\sin \left( A-{{\theta }_{A}} \right)\sin \left( B-{{\theta }_{B}} \right) \right) \\ 
			& +c\left( \sin {{\theta }_{C}}\sin {{\theta }_{A}}+\sin {{\theta }_{A}}\sin \left( B-{{\theta }_{B}} \right)+\sin \left( B-{{\theta }_{B}} \right)\sin \left( C-{{\theta }_{C}} \right) \right) \\ 
		\end{aligned} \right)}.\]
	
	Using the sine theorem to obtain
	\[\alpha _{A}^{P}=\frac{\sin A\left( \sin {{\theta }_{A}}\sin {{\theta }_{B}}+\sin {{\theta }_{B}}\sin \left( C-{{\theta }_{C}} \right)+\sin \left( C-{{\theta }_{C}} \right)\sin \left( A-{{\theta }_{A}} \right) \right)}{\left( \begin{aligned}
			& \sin A\left( \sin {{\theta }_{A}}\sin {{\theta }_{B}}+\sin {{\theta }_{B}}\sin \left( C-{{\theta }_{C}} \right)+\sin \left( C-{{\theta }_{C}} \right)\sin \left( A-{{\theta }_{A}} \right) \right) \\ 
			& +\sin B\left( \sin {{\theta }_{B}}\sin {{\theta }_{C}}+\sin {{\theta }_{C}}\sin \left( A-{{\theta }_{A}} \right)+\sin \left( A-{{\theta }_{A}} \right)\sin \left( B-{{\theta }_{B}} \right) \right) \\ 
			& +\sin C\left( \sin {{\theta }_{C}}\sin {{\theta }_{A}}+\sin {{\theta }_{A}}\sin \left( B-{{\theta }_{B}} \right)+\sin \left( B-{{\theta }_{B}} \right)\sin \left( C-{{\theta }_{C}} \right) \right) \\ 
		\end{aligned} \right)}.\]
	
	Similarly, it can be obtained that:
	\[\alpha _{B}^{P}=\frac{\sin B\left( \sin {{\theta }_{B}}\sin {{\theta }_{C}}+\sin {{\theta }_{C}}\sin \left( A-{{\theta }_{A}} \right)+\sin \left( A-{{\theta }_{A}} \right)\sin \left( B-{{\theta }_{B}} \right) \right)}{\left( \begin{aligned}
			& \sin A\left( \sin {{\theta }_{A}}\sin {{\theta }_{B}}+\sin {{\theta }_{B}}\sin \left( C-{{\theta }_{C}} \right)+\sin \left( C-{{\theta }_{C}} \right)\sin \left( A-{{\theta }_{A}} \right) \right) \\ 
			& +\sin B\left( \sin {{\theta }_{B}}\sin {{\theta }_{C}}+\sin {{\theta }_{C}}\sin \left( A-{{\theta }_{A}} \right)+\sin \left( A-{{\theta }_{A}} \right)\sin \left( B-{{\theta }_{B}} \right) \right) \\ 
			& +\sin C\left( \sin {{\theta }_{C}}\sin {{\theta }_{A}}+\sin {{\theta }_{A}}\sin \left( B-{{\theta }_{B}} \right)+\sin \left( B-{{\theta }_{B}} \right)\sin \left( C-{{\theta }_{C}} \right) \right) \\ 
		\end{aligned} \right)},\]
	\[\alpha _{C}^{P}=\frac{\sin C\left( \sin {{\theta }_{C}}\sin {{\theta }_{A}}+\sin {{\theta }_{A}}\sin \left( B-{{\theta }_{B}} \right)+\sin \left( B-{{\theta }_{B}} \right)\sin \left( C-{{\theta }_{C}} \right) \right)}{\left( \begin{aligned}
			& \sin A\left( \sin {{\theta }_{A}}\sin {{\theta }_{B}}+\sin {{\theta }_{B}}\sin \left( C-{{\theta }_{C}} \right)+\sin \left( C-{{\theta }_{C}} \right)\sin \left( A-{{\theta }_{A}} \right) \right) \\ 
			& +\sin B\left( \sin {{\theta }_{B}}\sin {{\theta }_{C}}+\sin {{\theta }_{C}}\sin \left( A-{{\theta }_{A}} \right)+\sin \left( A-{{\theta }_{A}} \right)\sin \left( B-{{\theta }_{B}} \right) \right) \\ 
			& +\sin C\left( \sin {{\theta }_{C}}\sin {{\theta }_{A}}+\sin {{\theta }_{A}}\sin \left( B-{{\theta }_{B}} \right)+\sin \left( B-{{\theta }_{B}} \right)\sin \left( C-{{\theta }_{C}} \right) \right) \\ 
		\end{aligned} \right)}.\]
\end{proof}
\hfill $\square$\par

From the above formula, it can be seen that as long as the numerator is calculated, the denominator can be obtained by rotating and then summing up.

Here is an example.


\begin{example}{}\label{BuluokadianDeBiaojiafenliang}
	Given a $\triangle ABC$, find the frame components of the first and second Brocard points of $\triangle ABC$.
\end{example}

\begin{solution}
	Firstly, calculate the frame component of the first Brocard point.
	Let the three polar angles of $\triangle ABC$ be ${{\theta }_{A}}={{\theta }_{B}}={{\theta }_{C}}=\theta $. Then $\theta $ is called the first Brocard angle The first Brocard angle can be calculated according to the following formula:
	\[\cot \theta =\frac{{{a}^{2}}+{{b}^{2}}+{{c}^{2}}}{4S}.\]
	
	Let's denote the first Brocard point as ${{B}_{1}}$, According to theorem \ref{thm:GenjuJijiaoJisuanBiaojiafenliang}, let $P:={{B}_{1}}$, it is obtained that
	\[\alpha _{A}^{{{B}_{1}}}=\alpha _{A}^{P}=\frac{{{f}_{a}}}{{{f}_{a}}+{{f}_{b}}+{{f}_{c}}}.\]
	
	Where
	\[{{f}_{a}}=a\left( \sin {{\theta }_{A}}\sin {{\theta }_{B}}+\sin {{\theta }_{B}}\sin \left( C-{{\theta }_{C}} \right)+\sin \left( C-{{\theta }_{C}} \right)\sin \left( A-{{\theta }_{A}} \right) \right),\]
	\[{{f}_{b}}=b\left( \sin {{\theta }_{B}}\sin {{\theta }_{C}}+\sin {{\theta }_{C}}\sin \left( A-{{\theta }_{A}} \right)+\sin \left( A-{{\theta }_{A}} \right)\sin \left( B-{{\theta }_{B}} \right) \right),\]
	\[{{f}_{c}}=c\left( \sin {{\theta }_{C}}\sin {{\theta }_{A}}+\sin {{\theta }_{A}}\sin \left( B-{{\theta }_{B}} \right)+\sin \left( B-{{\theta }_{B}} \right)\sin \left( C-{{\theta }_{C}} \right) \right).\]
	
	Firstly, calculate the numerator of $\alpha _{A}^{P}$.
	\begin{align*}
		{{f}_{a}}& =a\left( \sin {{\theta }_{A}}\sin {{\theta }_{B}}+\sin {{\theta }_{B}}\sin \left( C-{{\theta }_{C}} \right)+\sin \left( C-{{\theta }_{C}} \right)\sin \left( A-{{\theta }_{A}} \right) \right) \\ 
		& =a\left( {{\sin }^{2}}\theta +\sin \theta \sin \left( C-\theta  \right)+\sin \left( C-\theta  \right)\sin \left( A-\theta  \right) \right) \\ 
		& =a\left( {{\sin }^{2}}\theta +\sin \left( C-\theta  \right)\left( \sin \theta +\sin \left( A-\theta  \right) \right) \right) \\ 
		& =a\left( {{\sin }^{2}}\theta +\sin \left( C-\theta  \right)\left( \sin \theta +\sin A\cos \theta -\cos A\sin \theta  \right) \right) \\ 
		& =a\left( {{\sin }^{2}}\theta +\left( \sin C\cos \theta -\cos C\sin \theta  \right)\left( \sin \theta +\sin A\cos \theta -\cos A\sin \theta  \right) \right) \\ 
		& =a\left( {{\sin }^{2}}\theta +{{\sin }^{2}}\theta \left( \sin C\cot \theta -\cos C \right)\left( 1+\sin A\cot \theta -\cos A \right) \right) \\ 
		& =a{{\sin }^{2}}\theta \left( 1+\left( \sin C\cot \theta -\cos C \right)\left( 1+\sin A\cot \theta -\cos A \right) \right),  
	\end{align*}
	i.e.
	\begin{align*}
		{{f}_{a}}& =a{{\sin }^{2}}\theta \left( \begin{aligned}
			& 1+\left( \frac{2S}{ab}\frac{{{a}^{2}}+{{b}^{2}}+{{c}^{2}}}{4S}-\frac{{{a}^{2}}+{{b}^{2}}-{{c}^{2}}}{2ab} \right) \\ 
			& \times \left( 1+\frac{2S}{bc}\frac{{{a}^{2}}+{{b}^{2}}+{{c}^{2}}}{4S}-\frac{{{b}^{2}}+{{c}^{2}}-{{a}^{2}}}{2bc} \right) \\ 
		\end{aligned} \right) \\ 
		& =a{{\sin }^{2}}\theta \left( 1+\frac{{{c}^{2}}}{ab}\left( 1+\frac{{{a}^{2}}}{bc} \right) \right)=a{{\sin }^{2}}\theta \left( 1+\frac{c\left( {{a}^{2}}+bc \right)}{a{{b}^{2}}} \right) \\ 
		& ={{\sin }^{2}}\theta \left( \frac{a{{b}^{2}}+b{{c}^{2}}+c{{a}^{2}}}{{{b}^{2}}} \right)=\frac{1}{1\text{+}{{\cot }^{2}}\theta }\left( \frac{a{{b}^{2}}+b{{c}^{2}}+c{{a}^{2}}}{{{b}^{2}}} \right) \\ 
		& =\frac{1}{1\text{+}{{\left( \frac{{{a}^{2}}+{{b}^{2}}+{{c}^{2}}}{4S} \right)}^{2}}}\left( \frac{a{{b}^{2}}+b{{c}^{2}}+c{{a}^{2}}}{{{b}^{2}}} \right) \\ 
		& =\frac{16{{S}^{2}}}{16{{S}^{2}}\text{+}{{\left( {{a}^{2}}+{{b}^{2}}+{{c}^{2}} \right)}^{2}}}\left( \frac{a{{b}^{2}}+b{{c}^{2}}+c{{a}^{2}}}{{{b}^{2}}} \right).  
	\end{align*}
	
	Similarly, it can be obtained that:
	\[{{f}_{b}}=\frac{16{{S}^{2}}}{16{{S}^{2}}\text{+}{{\left( {{a}^{2}}+{{b}^{2}}+{{c}^{2}} \right)}^{2}}}\left( \frac{b{{c}^{2}}+c{{a}^{2}}+a{{b}^{2}}}{{{c}^{2}}} \right),\]
	\[{{f}_{c}}=\frac{16{{S}^{2}}}{16{{S}^{2}}\text{+}{{\left( {{a}^{2}}+{{b}^{2}}+{{c}^{2}} \right)}^{2}}}\left( \frac{c{{a}^{2}}+a{{b}^{2}}+b{{c}^{2}}}{{{a}^{2}}} \right).\]
	
	So we obtain the frame component of the first Brocard point ${{B}_{1}}$:
	\[\alpha _{A}^{{{B}_{1}}}=\frac{{{f}_{a}}}{{{f}_{a}}+{{f}_{b}}+{{f}_{c}}}=\frac{\frac{1}{{{b}^{2}}}}{\frac{1}{{{b}^{2}}}+\frac{1}{{{c}^{2}}}+\frac{1}{{{a}^{2}}}}=\frac{{{c}^{2}}{{a}^{2}}}{{{a}^{2}}{{b}^{2}}+{{b}^{2}}{{c}^{2}}+{{c}^{2}}{{a}^{2}}},\]
	\[\alpha _{B}^{{{B}_{1}}}=\frac{{{f}_{b}}}{{{f}_{a}}+{{f}_{b}}+{{f}_{c}}}=\frac{\frac{1}{{{c}^{2}}}}{\frac{1}{{{b}^{2}}}+\frac{1}{{{c}^{2}}}+\frac{1}{{{a}^{2}}}}=\frac{{{a}^{2}}{{b}^{2}}}{{{a}^{2}}{{b}^{2}}+{{b}^{2}}{{c}^{2}}+{{c}^{2}}{{a}^{2}}},\]
	\[\alpha _{C}^{{{B}_{1}}}=\frac{{{f}_{c}}}{{{f}_{a}}+{{f}_{b}}+{{f}_{c}}}=\frac{\frac{1}{{{a}^{2}}}}{\frac{1}{{{b}^{2}}}+\frac{1}{{{c}^{2}}}+\frac{1}{{{a}^{2}}}}=\frac{{{b}^{2}}{{c}^{2}}}{{{a}^{2}}{{b}^{2}}+{{b}^{2}}{{c}^{2}}+{{c}^{2}}{{a}^{2}}}.\]
	
	Next, let's solve for the frame component of the second Brocard point.
	Due to the three polar angles of $\triangle ABC$ being ${{\theta }_{A}}={{\theta }_{B}}={{\theta }_{C}}=\theta $, let
	\[{{\varphi }_{A}}:=A-{{\theta }_{A}}=A-\theta ,\quad {{\varphi }_{B}}:=B-{{\theta }_{B}}=B-\theta ,\quad {{\varphi }_{C}}:=C-{{\theta }_{C}}=C-\theta .\]
	
	${{\varphi }_{A}}$, ${{\varphi }_{B}}$, ${{\varphi }_{C}}$ are called the second Brocard angle (where $\theta $ is the first Brocard angle).
	
	We denote the second Brocard point as ${{B}_{2}}$, According to theorem \ref{thm:GenjuJijiaoJisuanBiaojiafenliang}, we have
	\[\alpha _{A}^{P}=\alpha _{A}^{{{B}_{2}}}=\frac{{{f}_{a}}}{{{f}_{a}}+{{f}_{b}}+{{f}_{c}}}.\]
	
	Where
	\[{{f}_{a}}=a\left( \sin {{\theta }_{A}}\sin {{\theta }_{B}}+\sin {{\theta }_{B}}\sin \left( C-{{\theta }_{C}} \right)+\sin \left( C-{{\theta }_{C}} \right)\sin \left( A-{{\theta }_{A}} \right) \right),\]
	\[{{f}_{b}}=b\left( \sin {{\theta }_{B}}\sin {{\theta }_{C}}+\sin {{\theta }_{C}}\sin \left( A-{{\theta }_{A}} \right)+\sin \left( A-{{\theta }_{A}} \right)\sin \left( B-{{\theta }_{B}} \right) \right),\]
	\[{{f}_{c}}=c\left( \sin {{\theta }_{C}}\sin {{\theta }_{A}}+\sin {{\theta }_{A}}\sin \left( B-{{\theta }_{B}} \right)+\sin \left( B-{{\theta }_{B}} \right)\sin \left( C-{{\theta }_{C}} \right) \right).\]
	
	Firstly, calculate the numerator of $\alpha _{A}^{P}$.
	\begin{align*}
		{{f}_{a}}& =a\left( \sin {{\theta }_{A}}\sin {{\theta }_{B}}+\sin {{\theta }_{B}}\sin \left( C-{{\theta }_{C}} \right)+\sin \left( C-{{\theta }_{C}} \right)\sin \left( A-{{\theta }_{A}} \right) \right) \\ 
		& =a\left( \sin \left( A-\theta  \right)\sin \left( B-\theta  \right)+\sin \left( B-\theta  \right)\sin \theta +{{\sin }^{2}}\theta  \right) \\ 
		& =a\left( {{\sin }^{2}}\theta +\sin \left( B-\theta  \right)\left( \sin \theta +\sin \left( A-\theta  \right) \right) \right) \\ 
		& =a\left( {{\sin }^{2}}\theta +\sin \left( B-\theta  \right)\left( \sin \theta +\sin A\cos \theta -\cos A\sin \theta  \right) \right) \\ 
		& =a\left( {{\sin }^{2}}\theta +\left( \sin B\cos \theta -\cos B\sin \theta  \right)\left( \sin \theta +\sin A\cos \theta -\cos A\sin \theta  \right) \right) \\ 
		& =a\left( {{\sin }^{2}}\theta +{{\sin }^{2}}\theta \left( \sin B\cot \theta -\cos B \right)\left( 1+\sin A\cot \theta -\cos A \right) \right) \\ 
		& =a{{\sin }^{2}}\theta \left( 1+\left( \sin B\cot \theta -\cos B \right)\left( 1+\sin A\cot \theta -\cos A \right) \right),  
	\end{align*}
	i.e.
	\begin{align*}
		{{f}_{a}}& =a{{\sin }^{2}}\theta \left( \begin{aligned}
			& 1+\left( \frac{2S}{ca}\frac{{{a}^{2}}+{{b}^{2}}+{{c}^{2}}}{4S}-\frac{{{c}^{2}}+{{a}^{2}}-{{b}^{2}}}{2ca} \right) \\ 
			& \times \left( 1+\frac{2S}{bc}\frac{{{a}^{2}}+{{b}^{2}}+{{c}^{2}}}{4S}-\frac{{{b}^{2}}+{{c}^{2}}-{{a}^{2}}}{2bc} \right) \\ 
		\end{aligned} \right) \\ 
		& =a{{\sin }^{2}}\theta \left( 1+\frac{{{b}^{2}}}{ca}\left( 1+\frac{{{a}^{2}}}{bc} \right) \right)=a{{\sin }^{2}}\theta \left( 1+\frac{b\left( {{a}^{2}}+bc \right)}{{{c}^{2}}a} \right) \\ 
		& ={{\sin }^{2}}\theta \left( \frac{{{a}^{2}}b+{{b}^{2}}c+{{c}^{2}}a}{{{c}^{2}}} \right)=\frac{1}{1\text{+}{{\cot }^{2}}\theta }\left( \frac{{{a}^{2}}b+{{b}^{2}}c+{{c}^{2}}a}{{{c}^{2}}} \right) \\ 
		& =\frac{1}{1\text{+}{{\left( \frac{{{a}^{2}}+{{b}^{2}}+{{c}^{2}}}{4S} \right)}^{2}}}\left( \frac{{{a}^{2}}b+{{b}^{2}}c+{{c}^{2}}a}{{{c}^{2}}} \right) \\ 
		& =\frac{16{{S}^{2}}}{16{{S}^{2}}\text{+}{{\left( {{a}^{2}}+{{b}^{2}}+{{c}^{2}} \right)}^{2}}}\left( \frac{{{a}^{2}}b+{{b}^{2}}c+{{c}^{2}}a}{{{c}^{2}}} \right).  
	\end{align*}
	
	Similarly, it can be obtained that:
	\[{{f}_{b}}=\frac{16{{S}^{2}}}{16{{S}^{2}}\text{+}{{\left( {{a}^{2}}+{{b}^{2}}+{{c}^{2}} \right)}^{2}}}\left( \frac{{{b}^{2}}c+{{c}^{2}}a+{{a}^{2}}b}{{{a}^{2}}} \right),\]
	\[{{f}_{c}}=\frac{16{{S}^{2}}}{16{{S}^{2}}\text{+}{{\left( {{a}^{2}}+{{b}^{2}}+{{c}^{2}} \right)}^{2}}}\left( \frac{{{c}^{2}}a+{{a}^{2}}b+{{b}^{2}}c}{{{b}^{2}}} \right).\]
	
	So we obtain the frame component of the second Brocard point ${{B}_{2}}$:
	\[\alpha _{A}^{{{B}_{2}}}=\frac{{{f}_{a}}}{{{f}_{a}}+{{f}_{b}}+{{f}_{c}}}=\frac{\frac{1}{{{c}^{2}}}}{\frac{1}{{{c}^{2}}}+\frac{1}{{{a}^{2}}}+\frac{1}{{{b}^{2}}}}=\frac{{{a}^{2}}{{b}^{2}}}{{{a}^{2}}{{b}^{2}}+{{b}^{2}}{{c}^{2}}+{{c}^{2}}{{a}^{2}}},\]
	\[\alpha _{B}^{{{B}_{2}}}=\frac{{{f}_{b}}}{{{f}_{a}}+{{f}_{b}}+{{f}_{c}}}=\frac{\frac{1}{{{a}^{2}}}}{\frac{1}{{{c}^{2}}}+\frac{1}{{{a}^{2}}}+\frac{1}{{{b}^{2}}}}=\frac{{{b}^{2}}{{c}^{2}}}{{{a}^{2}}{{b}^{2}}+{{b}^{2}}{{c}^{2}}+{{c}^{2}}{{a}^{2}}},\]
	\[\alpha _{C}^{{{B}_{2}}}=\frac{{{f}_{c}}}{{{f}_{a}}+{{f}_{b}}+{{f}_{c}}}=\frac{\frac{1}{{{b}^{2}}}}{\frac{1}{{{c}^{2}}}+\frac{1}{{{a}^{2}}}+\frac{1}{{{b}^{2}}}}=\frac{{{c}^{2}}{{a}^{2}}}{{{a}^{2}}{{b}^{2}}+{{b}^{2}}{{c}^{2}}+{{c}^{2}}{{a}^{2}}}.\]
\end{solution}
\hfill $\diamond$\par

It can be seen that the introduction of Edge-Axis coordinate systems in Intercenter Geometry does not violate the original intention of Intercenter Geometry. Coordinates are parameters, which are known quantities in practical problems and can be regarded as constants (see previous examples). Therefore, the frame component can still be regarded as a function of the lengths of the three sides of a triangle.

On the other hand, the coordinates used in the above corollary are the coordinates of different vertex coordinate systems (i.e. Edge-Axis coordinate systems of different vertex), which brings great convenience to practical work. People can choose coordinate systems according to their needs, and there is no need to be limited to a single coordinate system like analytic geometry. One drawback of a single coordinate system in analytic geometry is that the coordinates obtained in different coordinate systems cannot be generalized, and coordinate transformations have to be performed. Transforming coordinates from one coordinate system to another coordinate system is undoubtedly cumbersome and requires additional computational costs. In Intercenter Geometry, the frame component is unique and independent of the coordinate system and the position of the origin (pole). The frame component obtained in any coordinate system is the final result. The biggest difference between Intercenter Geometry and analytic geometry is the use of frame components instead of coordinates. When studying the distance between two points later, the advantage of the frame component will be more prominent.


Furthermore, the frame component is independent of the position of the frame origin $O$. As long as the position relationship between point $P$ and $\triangle ABC$ remains unchanged, a beautiful formula can be obtained, which cannot be achieved by analytic geometry. 
\[\overrightarrow{OP}=\alpha _{A}^{P}\overrightarrow{OA}+\alpha _{B}^{P}\overrightarrow{OB}+\alpha _{C}^{P}\overrightarrow{OC},\]
\[\alpha _{A}^{P}+\alpha _{B}^{P}+\alpha _{C}^{P}=1.\]


\section{Calculate coordinates based on frame components}\label{GenjuBiaojiafenliangJisuanZuobiao}	
According to theorems and corollaries in section \ref{GenjuZuobiaoJisuanBiaojiafenliang}, it is known that as long as the coordinates in the Edge-Axis coordinate system are given, the frame components can be calculated. The opposite question is: How to calculate the coordinates in the Edge-Axis coordinate system based on the frame components? This section will solve this problem.

%
%
%

\begin{theorem}{Calculating Cartesian coordinates from frame components,  Daiyuan Zhang}{YouBiaojiafenliangJisuanZuobiao}\label{YouBiaojiafenliangJisuanZuobiao}
	Given a $\triangle ABC$ and a point $P$ on the plane of $\triangle ABC$ , the area of $\triangle ABC$ is $S$. If the frame components of $P$ are $\alpha _{A}^{P}$, $\alpha _{B}^{P}$ and $\alpha _{C}^{P}$, then
	
	1. The Cartesian coordinates of $P$ in $\mathcal{R}\left( AB \right)$ are:
	\[x_{\mathcal{R}\left( AB \right)}^{P}=\frac{1}{2c}\left( \left( {{b}^{2}}+{{c}^{2}}-{{a}^{2}} \right)\left( 1-\alpha _{A}^{P} \right)+\left( {{c}^{2}}+{{a}^{2}}-{{b}^{2}} \right)\alpha _{B}^{P} \right),\]
	\[y_{\mathcal{R}\left( AB \right)}^{P}=\frac{2S}{c}\alpha _{C}^{P}.\]
	
	2. The Cartesian coordinates of $P$ in $\mathcal{R}\left( BC \right)$ are:
	\[x_{\mathcal{R}\left( BC \right)}^{P}=\frac{1}{2a}\left( \left( {{c}^{2}}+{{a}^{2}}-{{b}^{2}} \right)\left( 1-\alpha _{B}^{P} \right)+\left( {{a}^{2}}+{{b}^{2}}-{{c}^{2}} \right)\alpha _{C}^{P} \right),\]
	\[y_{\mathcal{R}\left( BC \right)}^{P}=\frac{2S}{a}\alpha _{A}^{P}.\]
	
	3. The Cartesian coordinates of $P$ in $\mathcal{R}\left( CA \right)$ are:
	\[x_{\mathcal{R}\left( CA \right)}^{P}=\frac{1}{2b}\left( \left( {{a}^{2}}+{{b}^{2}}-{{c}^{2}} \right)\left( 1-\alpha _{C}^{P} \right)+\left( {{b}^{2}}+{{c}^{2}}-{{a}^{2}} \right)\alpha _{A}^{P} \right),\]
	\[y_{\mathcal{R}\left( CA \right)}^{P}=\frac{2S}{b}\alpha _{B}^{P}.\]
\end{theorem}

\begin{proof}
	Only prove the first scenario. 
	According to corollary \ref{cor:GenjuZhijiaozuobiaoJisuanBiaojiafenliang_BianchangXingshi}, it can be concluded that: 
	\[y_{\mathcal{R}\left( AB \right)}^{P}=\frac{2S}{c}\alpha _{C}^{P}.\]
	
	Therefore
	\begin{align*}
		\alpha _{B}^{P}& =\frac{1}{c}\left( x_{\mathcal{R}\left( AB \right)}^{P}-\frac{{{b}^{2}}+{{c}^{2}}-{{a}^{2}}}{4S}y_{\mathcal{R}\left( AB \right)}^{P} \right) \\ 
		& =\frac{1}{c}\left( x_{\mathcal{R}\left( AB \right)}^{P}-\frac{2bc\cos A}{4\left( \frac{1}{2} \right)bc\sin A}y_{\mathcal{R}\left( AB \right)}^{P} \right) \\ 
		& =\frac{1}{c}\left( x_{\mathcal{R}\left( AB \right)}^{P}-\cot A\cdot y_{\mathcal{R}\left( AB \right)}^{P} \right).  
	\end{align*}
	
	Therefore
	\begin{align*}
		x_{\mathcal{R}\left( AB \right)}^{P}& =c\alpha _{B}^{P}+\cot A\cdot y_{\mathcal{R}\left( AB \right)}^{P} \\ 
		& =c\alpha _{B}^{P}+\cot A\cdot \frac{2S}{c}\alpha _{C}^{P} \\ 
		& =\frac{1}{c}\left( {{c}^{2}}\alpha _{B}^{P}+2S\cot A\cdot \alpha _{C}^{P} \right) \\ 
		& =\frac{1}{c}\left( {{c}^{2}}\alpha _{B}^{P}+\frac{{{b}^{2}}+{{c}^{2}}-{{a}^{2}}}{2}\alpha _{C}^{P} \right) \\ 
		& =\frac{1}{c}\left( {{c}^{2}}\alpha _{B}^{P}+\frac{{{b}^{2}}+{{c}^{2}}-{{a}^{2}}}{2}\left( 1-\alpha _{A}^{P}-\alpha _{B}^{P} \right) \right) \\ 
		& =\frac{1}{c}\left( \frac{{{b}^{2}}+{{c}^{2}}-{{a}^{2}}}{2}\left( 1-\alpha _{A}^{P} \right)+{{c}^{2}}\alpha _{B}^{P}-\frac{{{b}^{2}}+{{c}^{2}}-{{a}^{2}}}{2}\alpha _{B}^{P} \right) \\ 
		& =\frac{1}{c}\left( \frac{{{b}^{2}}+{{c}^{2}}-{{a}^{2}}}{2}\left( 1-\alpha _{A}^{P} \right)+\left( {{c}^{2}}-\frac{{{b}^{2}}+{{c}^{2}}-{{a}^{2}}}{2} \right)\alpha _{B}^{P} \right) \\ 
		& =\frac{1}{c}\left( \frac{{{b}^{2}}+{{c}^{2}}-{{a}^{2}}}{2}\left( 1-\alpha _{A}^{P} \right)+\left( \frac{{{c}^{2}}+{{a}^{2}}-{{b}^{2}}}{2} \right)\alpha _{B}^{P} \right) \\ 
		& =\frac{1}{2c}\left( \left( {{b}^{2}}+{{c}^{2}}-{{a}^{2}} \right)\left( 1-\alpha _{A}^{P} \right)+\left( {{c}^{2}}+{{a}^{2}}-{{b}^{2}} \right)\alpha _{B}^{P} \right).  
	\end{align*}
	
	Similarly, other results can be proven.
\end{proof}
\hfill $\square$\par

Obviously, if point $P$ is located inside of $\triangle ABC$, then $y_{\mathcal{R}\left( AB \right)}^{P}$, $y_{\mathcal{R}\left( BC \right)}^{P}$ and $y_{\mathcal{R}\left( CA \right)}^{P}$ are all positive real numbers, and the geometric meaning is the distances from point $P$ to the edges of $AB$, $BC$ and $CA$, respectively. If point $P$ is located outside of $\triangle ABC$, then $\left| y_{\mathcal{R}\left( AB \right)}^{P} \right|$, $\left| y_{\mathcal{R}\left( BC \right)}^{P} \right|$ and $\left| y_{\mathcal{R}\left( CA \right)}^{P} \right|$ are the distances from point $P$ to the edges of $AB$, $BC$, and $CA$, respectively, or the distance between point $P$ and the vertical projection points of the three edges of $\triangle ABC$.

The formula for calculating Cartesian coordinates from the frame component is derived from the above theorem. By utilizing the relationship between polar coordinates and Cartesian coordinates, the formula for calculating polar coordinates from the frame component can also be obtained.

Below are the applications of the above theorems.


\begin{corollary}{Representing the area of a triangle using Edge-Axis coordinates}{YongBianzhouzuobiaoBiaoshiSanjiaoxingMianji}\label{YongBianzhouzuobiaoBiaoshiSanjiaoxingMianji}
	Given a $\triangle ABC$ and a point $P$ on the plane of $\triangle ABC$, the area of $\triangle ABC$ is $S$. The frame components of $P$ are $\alpha _{A}^{P}$, $\alpha _{B}^{P}$ and $\alpha _{C}^{P}$, respectively; In the Edge-Axis coordinate system $\mathcal{R}\left( AB \right)$, $\mathcal{R}\left( BC \right)$ and $\mathcal{R}\left( CA \right)$, the three vertical coordinates are $y_{\mathcal{R}\left( AB \right)}^{P}$, $y_{\mathcal{R}\left( BC \right)}^{P}$, $y_{\mathcal{R}\left( CA \right)}^{P}$, respectively, then
	\[S=\frac{1}{2}\left( cy_{\mathcal{R}\left( AB \right)}^{P}+ay_{\mathcal{R}\left( BC \right)}^{P}+by_{\mathcal{R}\left( CA \right)}^{P} \right).\]
\end{corollary}

%
%
%

\begin{proof}
	According to theorem \ref{thm:YouBiaojiafenliangJisuanZuobiao}, it is obtained that:
	\[y_{\mathcal{R}\left( AB \right)}^{P}=\frac{2S}{c}\alpha _{C}^{P},\]
	\[y_{\mathcal{R}\left( BC \right)}^{P}=\frac{2S}{a}\alpha _{A}^{P},\]
	\[y_{\mathcal{R}\left( CA \right)}^{P}=\frac{2S}{b}\alpha _{B}^{P}.\]
	
	Therefore
	\[\frac{cy_{\mathcal{R}\left( AB \right)}^{P}}{2S}=\alpha _{C}^{P},\]
	\[\frac{ay_{\mathcal{R}\left( BC \right)}^{P}}{2S}=\alpha _{A}^{P},\]
	\[\frac{by_{\mathcal{R}\left( CA \right)}^{P}}{2S}=\alpha _{B}^{P}.\]
	
	i.e.
	\[\frac{cy_{\mathcal{R}\left( AB \right)}^{P}}{2S}+\frac{ay_{\mathcal{R}\left( BC \right)}^{P}}{2S}+\frac{by_{\mathcal{R}\left( CA \right)}^{P}}{2S}=\alpha _{C}^{P}+\alpha _{A}^{P}+\alpha _{B}^{P}=1.\]
	
	i.e.
	\[S=\frac{1}{2}\left( cy_{\mathcal{R}\left( AB \right)}^{P}+ay_{\mathcal{R}\left( BC \right)}^{P}+by_{\mathcal{R}\left( CA \right)}^{P} \right).\]
\end{proof}
\hfill $\square$\par

The above area formula holds for any point $P$ on the plane of $\triangle ABC$.


\begin{corollary}{Summation formula for vertical coordinates of Edge-Axis,  Daiyuan Zhang}{Cor3.1.1}\label{Cor3.1.1}
	Given a $\triangle ABC$ and a point $P$ on the plane of $\triangle ABC$, the area of $\triangle ABC$ is $S$. The frame components of $P$ are $\alpha _{A}^{P}$, $\alpha _{B}^{P}$ and $\alpha _{C}^{P}$, respectively; In the Edge-Axis coordinate system $\mathcal{R}\left( AB \right)$, $\mathcal{R}\left( BC \right)$ and $\mathcal{R}\left( CA \right)$, the three vertical coordinates are $y_{\mathcal{R}\left( AB \right)}^{P}$, $y_{\mathcal{R}\left( BC \right)}^{P}$, $y_{\mathcal{R}\left( CA \right)}^{P}$, respectively, then
	\begin{align*}
		y_{\mathcal{R}\left( AB \right)}^{P}+y_{\mathcal{R}\left( BC \right)}^{P}+y_{\mathcal{R}\left( CA \right)}^{P}& =2S\left( \frac{\alpha _{A}^{P}}{a}+\frac{\alpha _{B}^{P}}{b}+\frac{\alpha _{C}^{P}}{c} \right) \\ 
		& =\frac{2S}{abc}\left( bc\alpha _{A}^{P}+ca\alpha _{B}^{P}+ab\alpha _{C}^{P} \right).  
	\end{align*}
\end{corollary}

%

\begin{proof}
	According to theorem \ref{thm:YouBiaojiafenliangJisuanZuobiao}, it is obtained that:
	\[y_{\mathcal{R}\left( AB \right)}^{P}=\frac{2S}{c}\alpha _{C}^{P},\]
	\[y_{\mathcal{R}\left( BC \right)}^{P}=\frac{2S}{a}\alpha _{A}^{P},\]
	\[y_{\mathcal{R}\left( CA \right)}^{P}=\frac{2S}{b}\alpha _{B}^{P}.\]
	
	Therefore
	\begin{align*}
		y_{\mathcal{R}\left( AB \right)}^{P}+y_{\mathcal{R}\left( BC \right)}^{P}+y_{\mathcal{R}\left( CA \right)}^{P}& =2S\left( \frac{\alpha _{A}^{P}}{a}+\frac{\alpha _{B}^{P}}{b}+\frac{\alpha _{C}^{P}}{c} \right) \\ 
		& =\frac{2S}{abc}\left( bc\alpha _{A}^{P}+ca\alpha _{B}^{P}+ab\alpha _{C}^{P} \right).  
	\end{align*}
\end{proof}
\hfill $\square$\par

The above formula is symmetrical and beautiful. In fact, if point $P$ is located inside a triangle, then all three vertical coordinates are positive real numbers, which means that the three vertical coordinates are the distances between point $P$ and the three sides of $\triangle ABC$, respectively.


\begin{example}{}\label{QiuZhongxin}
	Find the sum vertical coordinates of centroid $G$, circumcenter $Q$ and orthocenter $H$ of the Edge-Axis, respectively. 
\end{example}

%
%

\begin{solution}
	For the centroid $G$, it is obtained that:
	\begin{align*}
		y_{\mathcal{R}\left( BC \right)}^{G}+y_{\mathcal{R}\left( CA \right)}^{G}+y_{\mathcal{R}\left( AB \right)}^{G}& =\frac{2S}{abc}\left( bc\alpha _{A}^{G}+ca\alpha _{B}^{G}+ab\alpha _{C}^{G} \right) \\ 
		& =\frac{2S}{3abc}\left( bc+ca+ab \right).  
	\end{align*}
	
	For the circumcenter $Q$, it is obtained that:
	\begin{align*}
		y_{\mathcal{R}\left( BC \right)}^{Q}+y_{\mathcal{R}\left( CA \right)}^{Q}+y_{\mathcal{R}\left( AB \right)}^{Q}& =\frac{2S}{abc}\left( bc\alpha _{A}^{Q}+ca\alpha _{B}^{Q}+ab\alpha _{C}^{Q} \right) \\ 
		& =\frac{2S}{abc}\left( \frac{bc\sin 2A+ca\sin 2B+ab\sin 2C}{\sin 2A+\sin 2B+\sin 2C} \right).  
	\end{align*}
	
	For the orthocenter $H$, it is obtained that:
	\begin{align*}
		y_{\mathcal{R}\left( BC \right)}^{H}+y_{\mathcal{R}\left( CA \right)}^{H}+y_{\mathcal{R}\left( AB \right)}^{H}& =\frac{2S}{abc}\left( bc\alpha _{A}^{H}+ca\alpha _{B}^{H}+ab\alpha _{C}^{H} \right) \\ 
		& =\frac{2S}{abc}\left( \frac{bc\tan A+ca\tan B+ab\tan C}{\tan A+\tan B+\tan C} \right).  
	\end{align*}
\end{solution}
\hfill $\diamond$\par


\begin{theorem}{Denominator of frame component is not zero, Daiyuan Zhang}{BiaojiafenliangDeFenmuBuWeiLing}\label{BiaojiafenliangDeFenmuBuWeiLing}
	Given a $\triangle ABC$ and a point $P$ on the plane of $\triangle ABC$, the IR of $P$ is $\lambda _{AB}^{P}$, $\lambda _{BC}^{P}$ and  $\lambda _{CA}^{P}$ respectively, then the following inequality holds on the finite plane:
	\[1+\lambda _{AB}^{P}+\lambda _{AC}^{P}\ne 0,\]
	\[1+\lambda _{BC}^{P}+\lambda _{BA}^{P}\ne 0,\]
	\[1+\lambda _{CA}^{P}+\lambda _{CB}^{P}\ne 0.\]
\end{theorem}

%
%
%
%

\begin{proof}
	Only prove the first inequality, and the proof methods for the other two inequalities are similar. Let the frame components of $P$ be $\alpha _{A}^{P}$, $\alpha _{B}^{P}$, $\alpha _{C}^{P}$, respectively. According to theorem \ref{thm:Thm6.1.1}, it is obtained that: 
	\[\alpha _{A}^{P}=\frac{1}{1+\lambda _{AB}^{P}+\lambda _{AC}^{P}},\]
	\[\alpha _{B}^{P}=\frac{1}{1+\lambda _{BC}^{P}+\lambda _{BA}^{P}},\]	
	\[\alpha _{C}^{P}=\frac{1}{1+\lambda _{CA}^{P}+\lambda _{CB}^{P}}.\]	
	
	According to theorem \ref{thm:YouBiaojiafenliangJisuanZuobiao}, it is obtained that:
	\[y_{\mathcal{R}\left( AB \right)}^{P}=\frac{2S}{c}\alpha _{C}^{P},\]
	\[y_{\mathcal{R}\left( BC \right)}^{P}=\frac{2S}{a}\alpha _{A}^{P},\]
	\[y_{\mathcal{R}\left( CA \right)}^{P}=\frac{2S}{b}\alpha _{B}^{P}.\]
	
	Where $S$ is the area of $\triangle ABC$, hence
	\[1+\lambda _{AB}^{P}+\lambda _{AC}^{P}=\frac{1}{\alpha _{A}^{P}}=\frac{2S}{ay_{\mathcal{R}\left( BC \right)}^{P}}.\]
	
	Obviously, on a finite plane, the $y_{\mathcal{R}\left( BC \right)}^{P}$ of the Cartesian coordinate of the Edge-Axis coordinate system $\mathcal{R}\left( BC \right)$ is finite, that is, $\left| y_{\mathcal{R}\left( BC \right)}^{P} \right|<\infty $. And the edge length of $\triangle ABC$ is finite, that is, the edge length $a<\infty $, and the area of $\triangle ABC$ is $S\ne 0$, so
	\[1+\lambda _{AB}^{P}+\lambda _{AC}^{P}\ne 0.\]
	
	Similarly, the other two inequalities can be proven.
\end{proof}
\hfill $\square$\par

The above theorem answers the question raised at the end of section \ref{SanjiaoxingDingdianYuBianshangYidianDeBiaojiafenliang}.

To describe the geometric meaning of the frame components, alternative notation is used to represent $y_{\mathcal{R}\left( BC \right)}^{P}$, $y_{\mathcal{R}\left( CA \right)}^{P}$, $y_{\mathcal{R}\left( AB \right)}^{P}$, such that $h_{a}^{P}:=y_{\mathcal{R}\left( BC \right)}^{P}$, $h_{b}^{P}:=y_{\mathcal{R}\left( CA \right)}^{P}$, $h_{c}^{P}:=y_{\mathcal{R}\left( AB \right)}^{P}$.


\begin{corollary}{Geometric meaning of frame components, Daiyuan Zhang}{BiaojiafenliangDeJiheyiyi}\label{BiaojiafenliangDeJiheyiyi}
	Given a $\triangle ABC$, if ${{h}_{a}}$, ${{h}_{b}}$ and ${{h}_{c}}$ are the heights on the three edges $a$, $b$, and $c$ of $\triangle ABC$, respectively; the frame components of point $P$ are $\alpha _{A}^{P}$, $\alpha _{B}^{P}$ and $\alpha _{C}^{P}$, respectively; in the Edge-Axis coordinate system, the three vertical coordinates of point $P$ are $h_{a}^{P}$, $h_{b}^{P}$ and $h_{c}^{P}$, respectively; then
	\[\alpha _{A}^{P}=\frac{h_{a}^{P}}{{{h}_{a}}},\]
	\[\alpha _{B}^{P}=\frac{h_{b}^{P}}{{{h}_{b}}},\]
	\[\alpha _{C}^{P}=\frac{h_{c}^{P}}{{{h}_{c}}}.\]
\end{corollary}

%

\begin{proof}
	According to theorem \ref{thm:YouBiaojiafenliangJisuanZuobiao}, it is obtained that:
	\[y_{\mathcal{R}\left( BC \right)}^{P}=\frac{2S}{a}\alpha _{A}^{P},\]
	\[y_{\mathcal{R}\left( CA \right)}^{P}=\frac{2S}{b}\alpha _{B}^{P},\]
	\[y_{\mathcal{R}\left( AB \right)}^{P}=\frac{2S}{c}\alpha _{C}^{P}.\]
	
	Because  $S=\frac{1}{2}c{{h}_{c}}=\frac{1}{2}a{{h}_{a}}=\frac{1}{2}b{{h}_{b}}$,  where ${{h}_{a}}$, ${{h}_{b}}$, ${{h}_{c}}$ are the heights on the three edges $a$, $b$ and $c$ of $\triangle ABC$. The three real numbers $y_{\mathcal{R}\left( BC \right)}^{P}$, $y_{\mathcal{R}\left( CA \right)}^{P}$, $y_{\mathcal{R}\left( AB \right)}^{P}$ are the vertical Cartesian coordinates of three Edge-Axis coordinate systems of point $P$, and $S$ is the area of $\triangle ABC$. So it is obtained that:
	
	\[\alpha _{A}^{P}=\frac{ay_{\mathcal{R}\left( BC \right)}^{P}}{2S}=\frac{y_{\mathcal{R}\left( BC \right)}^{P}}{{{h}_{a}}}=\frac{h_{a}^{P}}{{{h}_{a}}},\]
	\[\alpha _{B}^{P}=\frac{by_{\mathcal{R}\left( CA \right)}^{P}}{2S}=\frac{y_{\mathcal{R}\left( CA \right)}^{P}}{{{h}_{b}}}=\frac{h_{b}^{P}}{{{h}_{b}}},\]
	\[\alpha _{C}^{P}=\frac{cy_{\mathcal{R}\left( AB \right)}^{P}}{2S}=\frac{y_{\mathcal{R}\left( AB \right)}^{P}}{{{h}_{c}}}=\frac{h_{c}^{P}}{{{h}_{c}}}.\]
\end{proof}
\hfill $\square$\par

From the above corollary, it can be seen that the geometric meaning of the frame components. And it can conveniently determine the positive or negative of the frame components. If both point $P$ and vertex $A$ are on the same side of $BC$ (i.e. $a$), then $\alpha _{A}^{P}$ is positive, negative on the opposite side, and 0 on $\overleftrightarrow{BC}$; If both point $P$ and vertex $B$ are on the same side of $CA$ (i.e. $b$), then $\alpha _{B}^{P}$ is positive, negative on the opposite side, and 0 on $\overleftrightarrow{CA}$; If both point $P$ and vertex $C$ are on the same side of $AB$ (i.e. $c$), then $\alpha _{C}^{P}$ is positive, negative on the opposite side, and 0 on $\overleftrightarrow{AB}$.

If ${{l}_{BC}}$ is a straight line on the $\triangle ABC$ plane, and ${{l}_{BC}}\parallel \overleftrightarrow{BC}$, then for any point $P\in {{l}_{BC}}$, the frame component $\alpha _{A}^{P}$ always remains unchanged; If ${{l}_{CA}}$ is a straight line on the $\triangle ABC$ plane, and ${{l}_{CA}}\parallel \overleftrightarrow{CA}$, then for any point $P\in {{l}_{CA}}$, the frame component $\alpha _{B}^{P}$ always remains unchanged; If ${{l}_{AB}}$ is a straight line on the $\triangle ABC$ plane, and ${{l}_{AB}}\parallel \overleftrightarrow{AB}$, then for any point $P\in {{l}_{AB}}$, the frame component $\alpha _{C}^{P}$ always remains unchanged.

%

\section{The polar angle of a triangle}
\begin{corollary}{Cotangent formula of polar angle-Form 1, Daiyuan Zhang}{JijiaoYuqieGongshi_Xingshi1}\label{JijiaoYuqieGongshi_Xingshi1}
	Given a $\triangle ABC$ and a point $P$ on the$\triangle ABC$ plane, the area of $\triangle ABC$ is $S$. The frame components of $P$ are $\alpha _{A}^{P}$, $\alpha _{A}^{P}$ and $\alpha _{A}^{P}$ respectively, then the cotangents of the polar angles of $\triangle ABC$ are:
	\[\cot {{\theta }_{A}}=\frac{\left( {{p}_{2}}-{{a}^{2}} \right)+{{c}^{2}}\frac{\alpha _{B}^{P}}{\alpha _{C}^{P}}}{2S},\]
	\[\cot {{\theta }_{B}}=\frac{\left( {{p}_{2}}-{{b}^{2}} \right)+{{a}^{2}}\frac{\alpha _{C}^{P}}{\alpha _{A}^{P}}}{2S},\]
	\[\cot {{\theta }_{C}}=\frac{\left( {{p}_{2}}-{{c}^{2}} \right)+{{b}^{2}}\frac{\alpha _{A}^{P}}{\alpha _{B}^{P}}}{2S}.\]
	
	Where
	\[{{p}_{2}}=\frac{1}{2}\left( {{a}^{2}}+{{b}^{2}}+{{c}^{2}} \right).\]
\end{corollary}

%

\begin{proof}
	According to corollary \ref{cor:GenjuJizuobiaoJisuanBiaojiafenliang_Xingshi2}, it is obtained that:
	\[\alpha _{A}^{P}=1-\frac{1}{c}\left( {{\rho }_{AP}}\cos {{\theta }_{A}}+\frac{{{c}^{2}}+{{a}^{2}}-{{b}^{2}}}{4S}{{\rho }_{AP}}\sin {{\theta }_{A}} \right),\]
	\[\alpha _{B}^{P}=\frac{1}{c}\left( {{\rho }_{AP}}\cos {{\theta }_{A}}-\frac{{{b}^{2}}+{{c}^{2}}-{{a}^{2}}}{4S}{{\rho }_{AP}}\sin {{\theta }_{A}} \right),\]
	\[\alpha _{C}^{P}=\frac{c}{2S}{{\rho }_{AP}}\sin {{\theta }_{A}}.\]
	\begin{align*}
		\frac{\alpha _{B}^{P}}{\alpha _{C}^{P}}& =\frac{\frac{1}{c}\left( {{\rho }_{AP}}\cos {{\theta }_{A}}-\frac{{{b}^{2}}+{{c}^{2}}-{{a}^{2}}}{4S}{{\rho }_{AP}}\sin {{\theta }_{A}} \right)}{\frac{c}{2S}{{\rho }_{AP}}\sin {{\theta }_{A}}} \\ 
		& =\frac{2S\left( \cot {{\theta }_{A}}-\frac{{{b}^{2}}+{{c}^{2}}-{{a}^{2}}}{4S} \right)}{{{c}^{2}}}  
	\end{align*}
	\[\cot {{\theta }_{A}}=\frac{{{c}^{2}}}{2S}\frac{\alpha _{B}^{P}}{\alpha _{C}^{P}}+\frac{{{b}^{2}}+{{c}^{2}}-{{a}^{2}}}{4S}=\frac{\left( {{p}_{2}}-{{a}^{2}} \right)+{{c}^{2}}\frac{\alpha _{B}^{P}}{\alpha _{C}^{P}}}{2S},\]
	
	Similarly, it can be proven that:
	\[\cot {{\theta }_{B}}=\frac{\left( {{p}_{2}}-{{b}^{2}} \right)+{{a}^{2}}\frac{\alpha _{C}^{P}}{\alpha _{A}^{P}}}{2S},\]
	\[\cot {{\theta }_{C}}=\frac{\left( {{p}_{2}}-{{c}^{2}} \right)+{{b}^{2}}\frac{\alpha _{A}^{P}}{\alpha _{B}^{P}}}{2S}.\]
\end{proof}
\hfill $\square$\par	


\begin{corollary}{Cotangent formula of polar angle-Form 2, Daiyuan Zhang}{JijiaoYuqieGongshi_Xingshi2}\label{JijiaoYuqieGongshi_Xingshi2}
	Given a $\triangle ABC$ and a point $P$ on the$\triangle ABC$ plane, the area of $\triangle ABC$ is $S$, The frame components of $P$ are $\alpha _{A}^{P}$, $\alpha _{A}^{P}$ and $\alpha _{A}^{P}$ respectively, and the polar angles are ${{\theta }_{A}}$, ${{\theta }_{B}}$ and ${{\theta }_{C}}$ respectively, the first Brocard angle is $\theta $, then
	\[\alpha _{A}^{P}\cot {{\theta }_{B}}+\alpha _{B}^{P}\cot {{\theta }_{C}}+\alpha _{C}^{P}\cot {{\theta }_{A}}=\cot \theta =\frac{{{a}^{2}}+{{b}^{2}}+{{c}^{2}}}{4S}.\]
\end{corollary}

%
%

\begin{proof}	
	Let ${{\theta }_{A}}={{\theta }_{B}}={{\theta }_{C}}=\theta $, $\theta $ is called the first Brocard angle, the corresponding point $P$ is called the first Brocard point, and the first Brocard point is denoted as ${{B}_{1}}$. According to corollary \ref{cor:GenjuJizuobiaoJisuanBiaojiafenliang_Xingshi2}, it is obtained that:
	\begin{align*}
		& \alpha _{A}^{P}\cot {{\theta }_{B}}+\alpha _{B}^{P}\cot {{\theta }_{C}}+\alpha _{C}^{P}\cot {{\theta }_{A}} \\ 
		& =\frac{{{c}^{2}}\alpha _{B}^{P}+\left( {{p}_{2}}-{{a}^{2}} \right)\alpha _{C}^{P}+{{a}^{2}}\alpha _{C}^{P}+\left( {{p}_{2}}-{{b}^{2}} \right)\alpha _{A}^{P}+{{b}^{2}}\alpha _{A}^{P}+\left( {{p}_{2}}-{{c}^{2}} \right)\alpha _{B}^{P}}{2S} \\ 
		& =\frac{{{p}_{2}}\left( \alpha _{A}^{P}+\alpha _{B}^{P}+\alpha _{C}^{P} \right)}{2S}=\frac{{{p}_{2}}}{2S}=\frac{{{a}^{2}}+{{b}^{2}}+{{c}^{2}}}{4S}. \\ 
	\end{align*}
	
	The above equation indicates that, for a given triangle, regardless of the position of the point (IC) $P$, $\alpha _{A}^{P}\cot {{\theta }_{B}}+\alpha _{B}^{P}\cot {{\theta }_{C}}+\alpha _{C}^{P}\cot {{\theta }_{A}}$ is always a constant. So, choose point (IC) $P$ in a special position ${{B}_{1}}$, i.e. when $P={{B}_{1}}$, ${{\theta }_{A}}={{\theta }_{B}}={{\theta }_{C}}=\theta $, then according to the above equation, it is obtained that
	\begin{align*}
		& \alpha _{A}^{P}\cot {{\theta }_{B}}+\alpha _{B}^{P}\cot {{\theta }_{C}}+\alpha _{C}^{P}\cot {{\theta }_{A}} \\ 
		& =\alpha _{A}^{{{B}_{1}}}\cot \theta +\alpha _{B}^{{{B}_{1}}}\cot \theta +\alpha _{C}^{{{B}_{1}}}\cot \theta  \\ 
		& =\left( \alpha _{A}^{{{B}_{1}}}+\alpha _{B}^{{{B}_{1}}}+\alpha _{C}^{{{B}_{1}}} \right)\cot \theta =\cot \theta . \\ 
	\end{align*}
	
	Therefore
	\[\alpha _{A}^{P}\cot {{\theta }_{B}}+\alpha _{B}^{P}\cot {{\theta }_{C}}+\alpha _{C}^{P}\cot {{\theta }_{A}}=\cot \theta =\frac{{{a}^{2}}+{{b}^{2}}+{{c}^{2}}}{4S}.\]
\end{proof}
\hfill $\square$\par

Choosing some special ICs can yield some interesting results. Below are some applications.


\begin{example}{}\label{QiuZhongxinJijiaoYuqieDePingjunzhi}
	Find the average value of the cotangents of polar angles of centroid.
\end{example}

%

\begin{solution}
	For the centroid, $\alpha _{A}^{G}=\alpha _{B}^{G}=\alpha _{C}^{G}={1}/{3}\;$, if the polar angles of the centroid are ${{\theta }_{A}}$, ${{\theta }_{B}}$ and ${{\theta }_{C}}$ respectively, then
	\[\alpha _{A}^{G}\cot {{\theta }_{B}}+\alpha _{B}^{G}\cot {{\theta }_{C}}+\alpha _{C}^{G}\cot {{\theta }_{A}}=\frac{1}{3}\left( \cot {{\theta }_{A}}+\cot {{\theta }_{B}}+\cot {{\theta }_{C}} \right)=\cot \theta =\frac{{{a}^{2}}+{{b}^{2}}+{{c}^{2}}}{4S}.\]	
	
	The arithmetic mean of the cotangents of the three polar angles of the centroid is exactly the cotangent of the first Brocard angle.
\end{solution}
\hfill $\diamond$\par


\begin{example}{}\label{QiuZheng}
	Proof that: 
	\[S={{p}^{2}}\tan \frac{A}{2}\tan \frac{B}{2}\tan \frac{C}{2}.\]
\end{example}

%
%
%

\begin{solution}
	For the incenter $I$, we have: ${{\theta }_{A}}={A}/{2}\;$, ${{\theta }_{B}}={B}/{2}\;$, ${{\theta }_{C}}={C}/{2}\;$, therefore
	\begin{align*}
		\cot \frac{A}{2}& =\cot {{\theta }_{A}}=\frac{\left( {{p}_{2}}-{{a}^{2}} \right)+{{c}^{2}}\frac{\alpha _{B}^{I}}{\alpha _{C}^{I}}}{2S}=\frac{\left( {{p}_{2}}-{{a}^{2}} \right)+{{c}^{2}}\frac{b}{c}}{2S} \\ 
		& =\frac{\left( {{p}_{2}}-{{a}^{2}} \right)+bc}{2S}=\frac{\left( {{b}^{2}}+{{c}^{2}}-{{a}^{2}} \right)+2bc}{4S}=\frac{{{\left( b+c \right)}^{2}}-{{a}^{2}}}{4S} \\ 
		& =\frac{\left( b+c+a \right)\left( b+c-a \right)}{4S}=\frac{p\left( p-a \right)}{S}.  
	\end{align*}
	
	Similarly, it is obtained that:
	\[\cot \frac{B}{2}=\cot {{\theta }_{B}}=\frac{p\left( p-b \right)}{S},\]
	\[\cot \frac{C}{2}=\cot {{\theta }_{C}}=\frac{p\left( p-c \right)}{S}.\]
	
	Therefore
	\[\tan \frac{A}{2}\tan \frac{B}{2}\tan \frac{C}{2}=\frac{S}{p\left( p-a \right)}\frac{S}{p\left( p-b \right)}\frac{S}{p\left( p-c \right)}=\frac{S}{{{p}^{2}}}.\]
	
	This yields an identity for the area of a triangle:
	\[S={{p}^{2}}\tan \frac{A}{2}\tan \frac{B}{2}\tan \frac{C}{2}.\]
\end{solution}
\hfill $\diamond$\par


\section{A sufficient and necessary condition for concurrent of three lines}
Firstly, I propose a theorem that can be used to solve a class of problems for concurrent of three lines. Then several examples will be given.

%

\begin{figure}[h]
	\centering
	\includegraphics[totalheight=8cm]{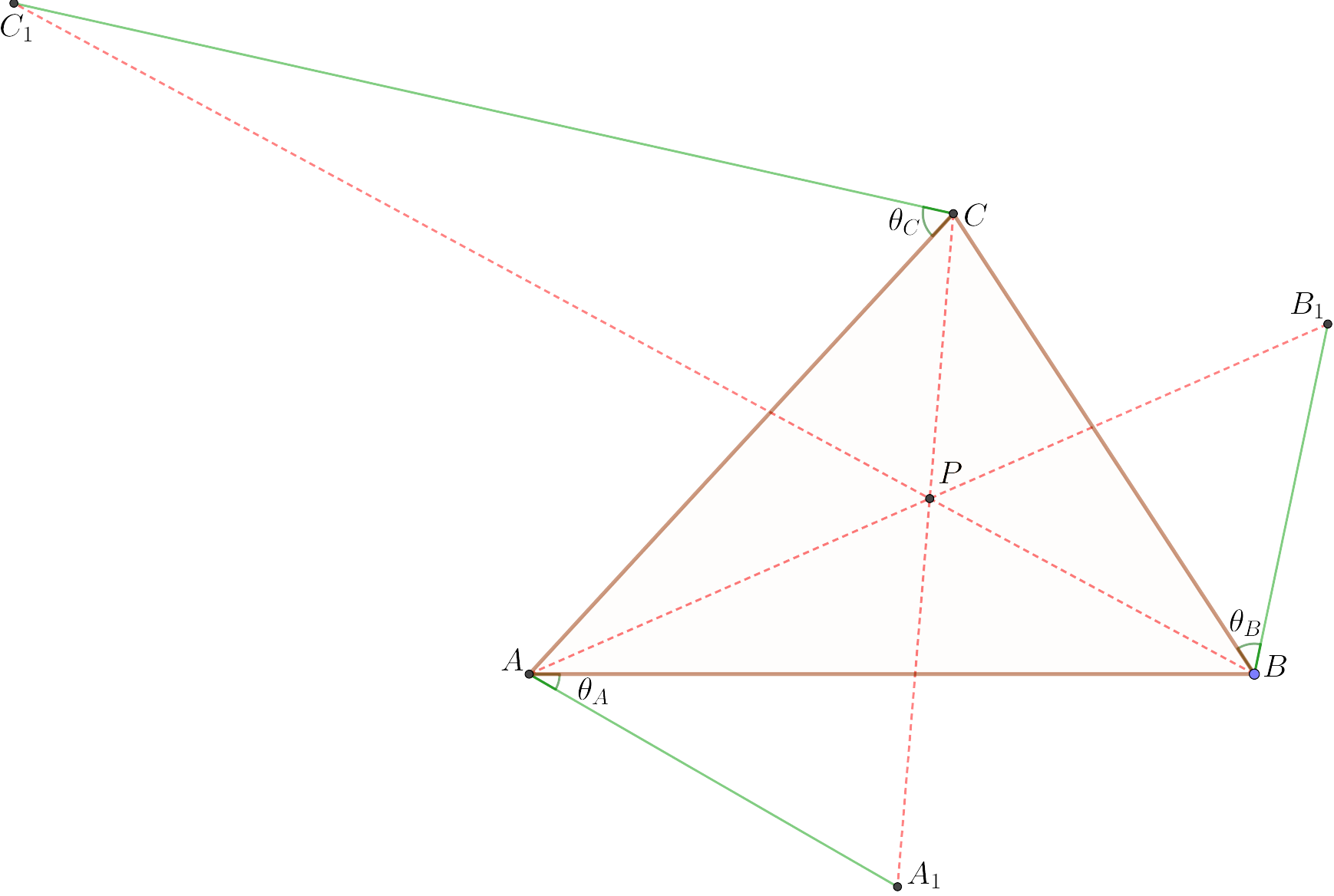}
	\caption{Figure of sufficient and necessary condition for concurrent of three lines} \label{fig:SanxiangongdianDeChongfenbiyaotiaojianTushi}
\end{figure}

\begin{theorem}{Sufficient and necessary condition for concurrent of three lines, Daiyuan Zhang}{YiGeSanxiangongdianDeChongfenbiyaotiaojian}\label{YiGeSanxiangongdianDeChongfenbiyaotiaojian}
	As shown in figure \ref{fig:SanxiangongdianDeChongfenbiyaotiaojianTushi}, given a $\triangle ABC$, and given the following quantities: ${{\theta }_{A}}$, ${{\theta }_{B}}$, ${{\theta }_{C}}$; $A{{A}_{1}}={{k }_{AB}}AB$, $B{{B}_{1}}={{k }_{BC}}BC$, $C{{C}_{1}}={{k }_{CA}}CA$; ${{k }_{AB}}\ne 0$, ${{k }_{BC}}\ne 0$, ${{k }_{CA}}\ne 0$.  Then the necessary and sufficient condition for the concurrent of three lines $\overleftrightarrow{A{{B}_{1}}}$, $\overleftrightarrow{B{{C}_{1}}}$ and $\overleftrightarrow{C{{A}_{1}}}$ is
	\begin{equation*}\label{SanxiangongdianDeChongfenbiyaotiaojianGongshi}
		\frac{{{k }_{AB}}\sin \left( A-{{\theta }_{A}} \right)}{\sin B-{{k }_{AB}}\sin \left( B+{{\theta }_{A}} \right)}\cdot \frac{{{k }_{BC}}\sin \left( B-{{\theta }_{B}} \right)}{\sin C-{{k }_{BC}}\sin \left( C+{{\theta }_{B}} \right)}\cdot \frac{{{k }_{CA}}\sin \left( C-{{\theta }_{C}} \right)}{\sin A-{{k }_{CA}}\sin \left( A+{{\theta }_{C}} \right)}=1.		
	\end{equation*}
	
	Or
	\[\frac{A{{A}_{1}}\sin \left( A-{{\theta }_{A}} \right)}{AB\sin B-A{{A}_{1}}\sin \left( B+{{\theta }_{A}} \right)}\cdot \frac{B{{B}_{1}}\sin \left( B-{{\theta }_{B}} \right)}{BC\sin C-B{{B}_{1}}\sin \left( C+{{\theta }_{B}} \right)}\cdot \frac{C{{C}_{1}}\sin \left( C-{{\theta }_{C}} \right)}{CA\sin A-C{{C}_{1}}\sin \left( A+{{\theta }_{C}} \right)}=1.\] 
\end{theorem}

\begin{proof}
	As shown in figure \ref{fig:SanxiangongdianDeChongfenbiyaotiaojianTushi}, assuming three straight lines $\overleftrightarrow{A{{B}_{1}}}$, $\overleftrightarrow{B{{C}_{1}}}$ and $\overleftrightarrow{C{{A}_{1}}}$ has a common intersection point, which is
	\[P=\overleftrightarrow{A{{B}_{1}}}\cap \overleftrightarrow{B{{C}_{1}}}\cap \overleftrightarrow{C{{A}_{1}}},\]
	
	For $\triangle ABC$, according to Ceva's theorem, the necessary and sufficient condition for the concurrent of three straight lines $\overleftrightarrow{A{{B}_{1}}}$, $\overleftrightarrow{B{{C}_{1}}}$ and $\overleftrightarrow{C{{A}_{1}}}$ is 
	\[\lambda _{AB}^{P}\cdot \lambda _{BC}^{P}\cdot \lambda _{CA}^{P}=1.\]
	
	And
	\[\lambda _{AB}^{P}=\lambda _{AB}^{{{A}_{1}}},\quad \lambda _{BC}^{P}=\lambda _{BC}^{{{B}_{1}}},\quad \lambda _{CA}^{P}=\lambda _{CA}^{{{C}_{1}}}.\]
	
	So the necessary and sufficient condition for the concurrent of three straight lines $\overleftrightarrow{A{{B}_{1}}}$, $\overleftrightarrow{B{{C}_{1}}}$ and $\overleftrightarrow{C{{A}_{1}}}$ is 
	\[\lambda _{AB}^{{{A}_{1}}}\cdot \lambda _{BC}^{{{B}_{1}}}\cdot \lambda _{CA}^{{{C}_{1}}}=1.\]
	
	Let the polar coordinates of point ${{A}_{1}}$ in the Edge-Axis coordinate system $\mathcal{P}\left( AB \right)$ of $\triangle ABC$ are $\left( {{\rho }_{A{{A}_{1}}}},{{\theta }_{A}} \right)$, according to theorem \ref{thm:Thm6.4.1} and corollary \ref{cor:GenjuJizuobiaoJisuanBiaojiafenliang_Xingshi1}, it is obtained that:
	\begin{align*}
		\alpha _{A}^{{{A}_{1}}}& =1-\frac{{{\rho }_{A{{A}_{1}}}}}{c}\frac{\left( \sin C\sin {{\theta }_{A}}+\sin B\sin \left( A-{{\theta }_{A}} \right) \right)}{\sin A\sin B} \\ 
		& =1-\frac{A{{A}_{1}}}{AB}\frac{\left( \sin C\sin {{\theta }_{A}}+\sin B\sin \left( A-{{\theta }_{A}} \right) \right)}{\sin A\sin B} \\ 
		& =1-{{k }_{AB}}\frac{\left( \sin C\sin {{\theta }_{A}}+\sin B\sin \left( A-{{\theta }_{A}} \right) \right)}{\sin A\sin B} \\ 
		& =\frac{\sin A\sin B-{{k }_{AB}}\left( \sin C\sin {{\theta }_{A}}+\sin B\sin \left( A-{{\theta }_{A}} \right) \right)}{\sin A\sin B},  
	\end{align*}
	\[\alpha _{B}^{{{A}_{1}}}=\frac{{{\rho }_{A{{A}_{1}}}}}{c}\frac{\sin \left( A-{{\theta }_{A}} \right)}{\sin A}=\frac{A{{A}_{1}}}{AB}\frac{\sin \left( A-{{\theta }_{A}} \right)}{\sin A}=\frac{{{k }_{AB}}\sin \left( A-{{\theta }_{A}} \right)}{\sin A}.\]
	
	Therefore
	\begin{align*}
		\lambda _{AB}^{{{A}_{1}}}& =\frac{\alpha _{B}^{{{A}_{1}}}}{\alpha _{A}^{{{A}_{1}}}}=\frac{{{k }_{AB}}\sin B\sin \left( A-{{\theta }_{A}} \right)}{\sin A\sin B-{{k }_{AB}}\left( \sin C\sin {{\theta }_{A}}+\sin B\sin \left( A-{{\theta }_{A}} \right) \right)} \\ 
		& =\frac{{{k }_{AB}}\sin B\sin \left( A-{{\theta }_{A}} \right)}{\sin A\sin B-{{k }_{AB}}\left( \sin \left( A+B \right)\sin {{\theta }_{A}}+\sin B\sin \left( A-{{\theta }_{A}} \right) \right)} \\ 
		& =\frac{{{k }_{AB}}\sin B\sin \left( A-{{\theta }_{A}} \right)}{\sin A\sin B-{{k }_{AB}}\left( \sin A\cos B\sin {{\theta }_{A}}+\sin B\sin A\cos {{\theta }_{A}} \right)} \\ 
		& =\frac{{{k }_{AB}}\sin B\sin \left( A-{{\theta }_{A}} \right)}{\sin A\left( \sin B-{{k }_{AB}}\left( \cos B\sin {{\theta }_{A}}+\sin B\cos {{\theta }_{A}} \right) \right)} \\ 
		& =\frac{{{k }_{AB}}\sin B\sin \left( A-{{\theta }_{A}} \right)}{\sin A\left( \sin B-{{k }_{AB}}\sin \left( B+{{\theta }_{A}} \right) \right)}.  
	\end{align*}
	
	Let the polar coordinates of point ${{B}_{1}}$ in the Edge-Axis coordinate system $\mathcal{P}\left( BC \right)$ of $\triangle ABC$ are $\left( {{\rho }_{B{{B}_{1}}}},{{\theta }_{B}} \right)$, according to theorem \ref{thm:Thm6.4.1} and corollary \ref{cor:GenjuJizuobiaoJisuanBiaojiafenliang_Xingshi1}, it is obtained that:
	\begin{align*}
		\alpha _{B}^{{{B}_{1}}}& =1-\frac{{{\rho }_{B{{B}_{1}}}}}{a}\frac{\left( \sin A\sin {{\theta }_{B}}+\sin C\sin \left( B-{{\theta }_{B}} \right) \right)}{\sin B\sin C} \\ 
		& =1-\frac{B{{B}_{1}}}{BC}\frac{\left( \sin A\sin {{\theta }_{B}}+\sin C\sin \left( B-{{\theta }_{B}} \right) \right)}{\sin B\sin C} \\ 
		& =1-{{k }_{BC}}\frac{\left( \sin A\sin {{\theta }_{B}}+\sin C\sin \left( B-{{\theta }_{B}} \right) \right)}{\sin B\sin C} \\ 
		& =\frac{\sin B\sin C-{{k }_{BC}}\left( \sin A\sin {{\theta }_{B}}+\sin C\sin \left( B-{{\theta }_{B}} \right) \right)}{\sin B\sin C},  
	\end{align*}
	\[\alpha _{C}^{{{B}_{1}}}=\frac{{{\rho }_{B{{B}_{1}}}}}{a}\frac{\sin \left( B-{{\theta }_{AB}} \right)}{\sin B}=\frac{B{{B}_{1}}}{BC}\frac{\sin \left( B-{{\theta }_{B}} \right)}{\sin B}=\frac{{{k }_{BC}}\sin \left( B-{{\theta }_{B}} \right)}{\sin B}.\]
	
	Therefore
	\begin{align*}
		\lambda _{BC}^{{{B}_{1}}}& =\frac{\alpha _{C}^{{{B}_{1}}}}{\alpha _{B}^{{{B}_{1}}}}=\frac{{{k }_{BC}}\sin C\sin \left( B-{{\theta }_{B}} \right)}{\sin B\sin C-{{k }_{BC}}\left( \sin A\sin {{\theta }_{B}}+\sin C\sin \left( B-{{\theta }_{B}} \right) \right)} \\ 
		& =\frac{{{k }_{BC}}\sin C\sin \left( B-{{\theta }_{B}} \right)}{\sin B\sin C-{{k }_{BC}}\left( \sin \left( B+C \right)\sin {{\theta }_{B}}+\sin C\sin \left( B-{{\theta }_{B}} \right) \right)} \\ 
		& =\frac{{{k }_{BC}}\sin C\sin \left( B-{{\theta }_{B}} \right)}{\sin B\sin C-{{k }_{BC}}\left( \sin B\cos C\sin {{\theta }_{B}}+\sin C\sin B\cos {{\theta }_{B}} \right)} \\ 
		& =\frac{{{k }_{BC}}\sin C\sin \left( B-{{\theta }_{B}} \right)}{\sin B\left( \sin C-{{k }_{BC}}\left( \cos C\sin {{\theta }_{B}}+\sin C\cos {{\theta }_{B}} \right) \right)} \\ 
		& =\frac{{{k }_{BC}}\sin C\sin \left( B-{{\theta }_{B}} \right)}{\sin B\left( \sin C-{{k }_{BC}}\sin \left( C+{{\theta }_{B}} \right) \right)}.  
	\end{align*}
	
	Let the polar coordinates of point ${{C}_{1}}$ in the Edge-Axis coordinate system $\mathcal{P}\left( CA \right)$ of $\triangle ABC$ are $\left( {{\rho }_{C{{C}_{1}}}},{{\theta }_{C}} \right)$, according to theorem \ref{thm:Thm6.4.1} and corollary \ref{cor:GenjuJizuobiaoJisuanBiaojiafenliang_Xingshi1}, it is obtained that:
	\begin{align*}
		\alpha _{C}^{{{C}_{1}}}& =1-\frac{{{\rho }_{C{{C}_{1}}}}}{b}\frac{\left( \sin B\sin {{\theta }_{C}}+\sin A\sin \left( C-{{\theta }_{C}} \right) \right)}{\sin C\sin A} \\ 
		& =1-\frac{C{{C}_{1}}}{CA}\frac{\left( \sin B\sin {{\theta }_{C}}+\sin A\sin \left( C-{{\theta }_{C}} \right) \right)}{\sin C\sin A} \\ 
		& =1-{{k }_{CA}}\frac{\left( \sin B\sin {{\theta }_{C}}+\sin A\sin \left( C-{{\theta }_{C}} \right) \right)}{\sin C\sin A} \\ 
		& =\frac{\sin C\sin A-{{k }_{CA}}\left( \sin B\sin {{\theta }_{C}}+\sin A\sin \left( C-{{\theta }_{C}} \right) \right)}{\sin C\sin A},  
	\end{align*}
	\[\alpha _{A}^{{{C}_{1}}}=\frac{{{\rho }_{C{{C}_{1}}}}}{b}\frac{\sin \left( C-{{\theta }_{C}} \right)}{\sin C}=\frac{C{{C}_{1}}}{CA}\frac{\sin \left( C-{{\theta }_{C}} \right)}{\sin C}=\frac{{{k }_{CA}}\sin \left( C-{{\theta }_{C}} \right)}{\sin C}.\]
	
	Therefore
	\begin{align*}
		\lambda _{CA}^{{{C}_{1}}}& =\frac{\alpha _{A}^{{{C}_{1}}}}{\alpha _{C}^{{{C}_{1}}}}=\frac{{{k }_{CA}}\sin A\sin \left( C-{{\theta }_{C}} \right)}{\sin C\sin A-{{k }_{CA}}\left( \sin B\sin {{\theta }_{C}}+\sin A\sin \left( C-{{\theta }_{C}} \right) \right)} \\ 
		& =\frac{{{k }_{CA}}\sin A\sin \left( C-{{\theta }_{C}} \right)}{\sin C\sin A-{{k }_{CA}}\left( \sin \left( C+A \right)\sin {{\theta }_{C}}+\sin A\sin \left( C-{{\theta }_{C}} \right) \right)} \\ 
		& =\frac{{{k }_{CA}}\sin A\sin \left( C-{{\theta }_{C}} \right)}{\sin C\sin A-{{k }_{CA}}\left( \sin C\cos A\sin {{\theta }_{C}}+\sin A\sin C\cos {{\theta }_{C}} \right)} \\ 
		& =\frac{{{k }_{CA}}\sin A\sin \left( C-{{\theta }_{C}} \right)}{\sin C\sin A-{{k }_{CA}}\left( \sin C\cos A\sin {{\theta }_{C}}+\sin A\sin C\cos {{\theta }_{C}} \right)} \\ 
		& =\frac{{{k }_{CA}}\sin A\sin \left( C-{{\theta }_{C}} \right)}{\sin C\left( \sin A-{{k }_{CA}}\left( \cos A\sin {{\theta }_{C}}+\sin A\cos {{\theta }_{C}} \right) \right)} \\ 
		& =\frac{{{k }_{CA}}\sin A\sin \left( C-{{\theta }_{C}} \right)}{\sin C\left( \sin A-{{k }_{CA}}\sin \left( A+{{\theta }_{C}} \right) \right)}.  
	\end{align*}
	
	So the necessary and sufficient condition for the three lines $\overleftrightarrow{A{{B}_{1}}}$, $\overleftrightarrow{B{{C}_{1}}}$ and $\overleftrightarrow{C{{A}_{1}}}$ to be concurrent is:
	\[\frac{{{k }_{AB}}\sin \left( A-{{\theta }_{A}} \right)}{\sin B-{{k }_{AB}}\sin \left( B+{{\theta }_{A}} \right)}\cdot \frac{{{k }_{BC}}\sin \left( B-{{\theta }_{B}} \right)}{\sin C-{{k }_{BC}}\sin \left( C+{{\theta }_{B}} \right)}\cdot \frac{{{k }_{CA}}\sin \left( C-{{\theta }_{C}} \right)}{\sin A-{{k }_{CA}}\sin \left( A+{{\theta }_{C}} \right)}=1.\]
	
	Or
	\[\frac{A{{A}_{1}}\sin \left( A-{{\theta }_{A}} \right)}{AB\sin B-A{{A}_{1}}\sin \left( B+{{\theta }_{A}} \right)}\cdot \frac{B{{B}_{1}}\sin \left( B-{{\theta }_{B}} \right)}{BC\sin C-B{{B}_{1}}\sin \left( C+{{\theta }_{B}} \right)}\cdot \frac{C{{C}_{1}}\sin \left( C-{{\theta }_{C}} \right)}{CA\sin A-C{{C}_{1}}\sin \left( A+{{\theta }_{C}} \right)}=1.\]
\end{proof}
\hfill $\square$\par

There are six independent variables in the necessary and sufficient conditions of the above theorem: ${{\theta }_{A}}$, ${{\theta }_{B}}$, ${{\theta }_{C}}$, $A{{A}_{1}}$, $B{{B}_{1}}$, $C{{C}_{1}}$. As long as you know five of them, you can find the other one.

Below are a few application examples.

\begin{figure}[h]
	\centering
	\includegraphics[totalheight=8cm]{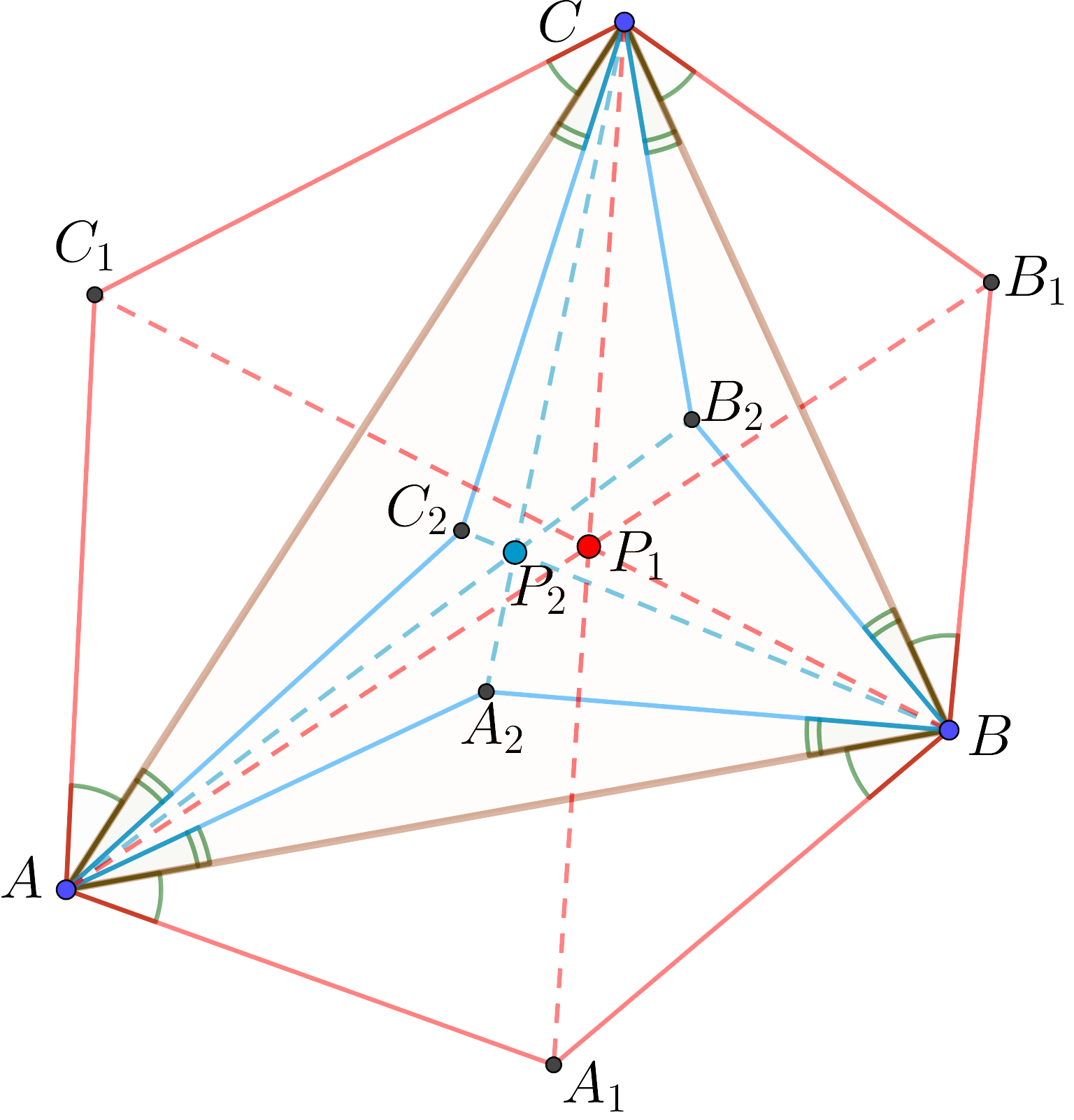}
	\caption{Concurrent of the outward and inward isosceles triangles} \label{fig:WaifansanjiaoxingHeNeifansanjiaoxingSanxiangongdian}
\end{figure}

\begin{example}{}\label{WaifansanjiaoxingHeNeifansanjiaoxingSanxiangongdian}
	As shown in figure \ref{fig:WaifansanjiaoxingHeNeifansanjiaoxingSanxiangongdian}, if the three outward flipped triangles of $\triangle ABC$ are isosceles triangles with equal base angles, then $\overleftrightarrow{A{{B}_{1}}}$, $\overleftrightarrow{B{{C}_{1}}}$ and $\overleftrightarrow{C{{A}_{1}}}$ are concurrent, that is ${{P}_{1}}=\overleftrightarrow{A{{B}_{1}}}\cap \overleftrightarrow{B{{C}_{1}}}\cap \overleftrightarrow{C{{A}_{1}}}$;  If the three inward flipped triangles of $\triangle ABC$ are isosceles triangles with equal base angles, then $\overleftrightarrow{A{{B}_{2}}}$, $\overleftrightarrow{B{{C}_{2}}}$ and $\overleftrightarrow{C{{A}_{2}}}$ are concurrent, that is ${{P}_{2}}=\overleftrightarrow{A{{B}_{2}}}\cap \overleftrightarrow{B{{C}_{2}}}\cap \overleftrightarrow{C{{A}_{2}}}$.
\end{example}

\begin{solution}
	Consider the situation of outward flipped triangles. According to the above theorem, let:
	\[\frac{{{k }_{AB}}\sin \left( A-{{\theta }_{A}} \right)}{\sin B-{{k }_{AB}}\sin \left( B+{{\theta }_{A}} \right)}=1,\]
	\[\frac{{{k }_{BC}}\sin \left( B-{{\theta }_{B}} \right)}{\sin C-{{k }_{BC}}\sin \left( C+{{\theta }_{B}} \right)}=1,\]
	\[\frac{{{k }_{CA}}\sin \left( C-{{\theta }_{C}} \right)}{\sin A-{{k }_{CA}}\sin \left( A+{{\theta }_{C}} \right)}=1.\]
	
	So the condition of concurrent of three lines $\overleftrightarrow{A{{B}_{1}}}$, $\overleftrightarrow{B{{C}_{1}}}$ and $\overleftrightarrow{C{{A}_{1}}}$ is satisfied (see theorem \ref{thm:YiGeSanxiangongdianDeChongfenbiyaotiaojian}). The above three formulas can be written as:
	\[{{k }_{AB}}\left( \sin \left( A-{{\theta }_{A}} \right)+\sin \left( B+{{\theta }_{A}} \right) \right)=\sin B,\]
	\[{{k }_{BC}}\left( \sin \left( B-{{\theta }_{B}} \right)+\sin \left( C+{{\theta }_{B}} \right) \right)=\sin C,\]
	\[{{k }_{CA}}\left( \sin \left( C-{{\theta }_{C}} \right)+\sin \left( A+{{\theta }_{C}} \right) \right)=\sin A.\]
	
	In the above equation, let ${{k }_{AB}}={{k }_{BC}}={{k }_{CA}}=k $, ${{\theta }_{A}}={{\theta }_{B}}={{\theta }_{C}}=\theta $, and obtain:
	\[k \left( \sin \left( A-\theta  \right)+\sin \left( B+\theta  \right) \right)=\sin B,\]
	\[k \left( \sin \left( B-\theta  \right)+\sin \left( C+\theta  \right) \right)=\sin C,\]
	\[k \left( \sin \left( C-\theta  \right)+\sin \left( A+\theta  \right) \right)=\sin A.\]
	
	That is
	\[k \left( \begin{aligned}
		& \sin \left( A-\theta  \right)+\sin \left( B+\theta  \right) \\ 
		& +\sin \left( B-\theta  \right)+\sin \left( C+\theta  \right) \\ 
		& +\sin \left( C-\theta  \right)+\sin \left( A+\theta  \right) \\ 
	\end{aligned} \right)=\sin A+\sin B+\sin C,\]
	i.e.
	\[2k \left( \sin A\cos \theta +\sin B\cos \theta +\sin C\cos \theta  \right)=\sin A+\sin B+\sin C,\]
	i.e.
	\begin{equation*}\label{YiGeChongfenBiyaotiaojianTeli1}
		k \cos \theta =\frac{1}{2}.	
	\end{equation*}
	
	Substitute the definitions of ${{k }_{AB}}$, ${{k }_{BC}}$, ${{k }_{CA}}$($=k $) into the above equation to obtain
	\begin{equation*}\label{YiGeChongfenBiyaotiaojianTeli2}
		\cos \theta =\frac{AB}{2A{{A}_{1}}},\quad \cos \theta =\frac{BC}{2B{{B}_{1}}},\quad \cos \theta =\frac{CA}{2C{{C}_{1}}}.	
	\end{equation*}
	
	The above equation indicates that if all three outward flipped triangles of $\triangle ABC$ are isosceles triangles, then the three straight lines $\overleftrightarrow{A{{B}_{1}}}$, $\overleftrightarrow{B{{C}_{1}}}$ and $\overleftrightarrow{C{{A}_{1}}}$ are concurrent, that is ${{P}_{1}}=\overleftrightarrow{A{{B}_{1}}}\cap \overleftrightarrow{B{{C}_{1}}}\cap \overleftrightarrow{C{{A}_{1}}}$.
	
	Similarly, it can be proven that if all three inward folded triangles of $\triangle ABC$ are isosceles triangles, then the three straight lines $\overleftrightarrow{A{{B}_{2}}}$, $\overleftrightarrow{B{{C}_{2}}}$ and $\overleftrightarrow{C{{A}_{2}}}$ are concurrent, that is ${{P}_{2}}=\overleftrightarrow{A{{B}_{2}}}\cap \overleftrightarrow{B{{C}_{2}}}\cap \overleftrightarrow{C{{A}_{2}}}$.
\end{solution}
\hfill $\diamond$\par

If the $A{{A}_{1}}B{{B}_{1}}C{{C}_{1}}$ (or $A{{A}_{2}}B{{B}_{2}}C{{C}_{2}}$) in figure \ref{fig:WaifansanjiaoxingHeNeifansanjiaoxingSanxiangongdian} is viewed as a hexagon, then the three diagonals $\overleftrightarrow{A{{B}_{1}}}$, $\overleftrightarrow{B{{C}_{1}}}$ and $\overleftrightarrow{C{{A}_{1}}}$ (or $\overleftrightarrow{A{{B}_{2}}}$, $\overleftrightarrow{B{{C}_{2}}}$ and $\overleftrightarrow{C{{A}_{2}}}$) of this hexagon are concurrent.

Through this example, it is known that theorem \ref{thm:GenjuYizhidianDeBiaojiafenliangQiuWeizhidianDeBiaojiafenliangGongshi} and its series of corollaries are important, as they can not only handle problems with triangles, but also with polygons.


\begin{example}{}\label{SanxiangongdianDeYingyong}
	Consider a special situation: ${{\theta }_{A}}=-B$, ${{\theta }_{B}}=-C$, ${{\theta }_{C}}=-A$.
\end{example}

%
%

\begin{figure}[h] 
	\centering
	\includegraphics[totalheight=6cm]{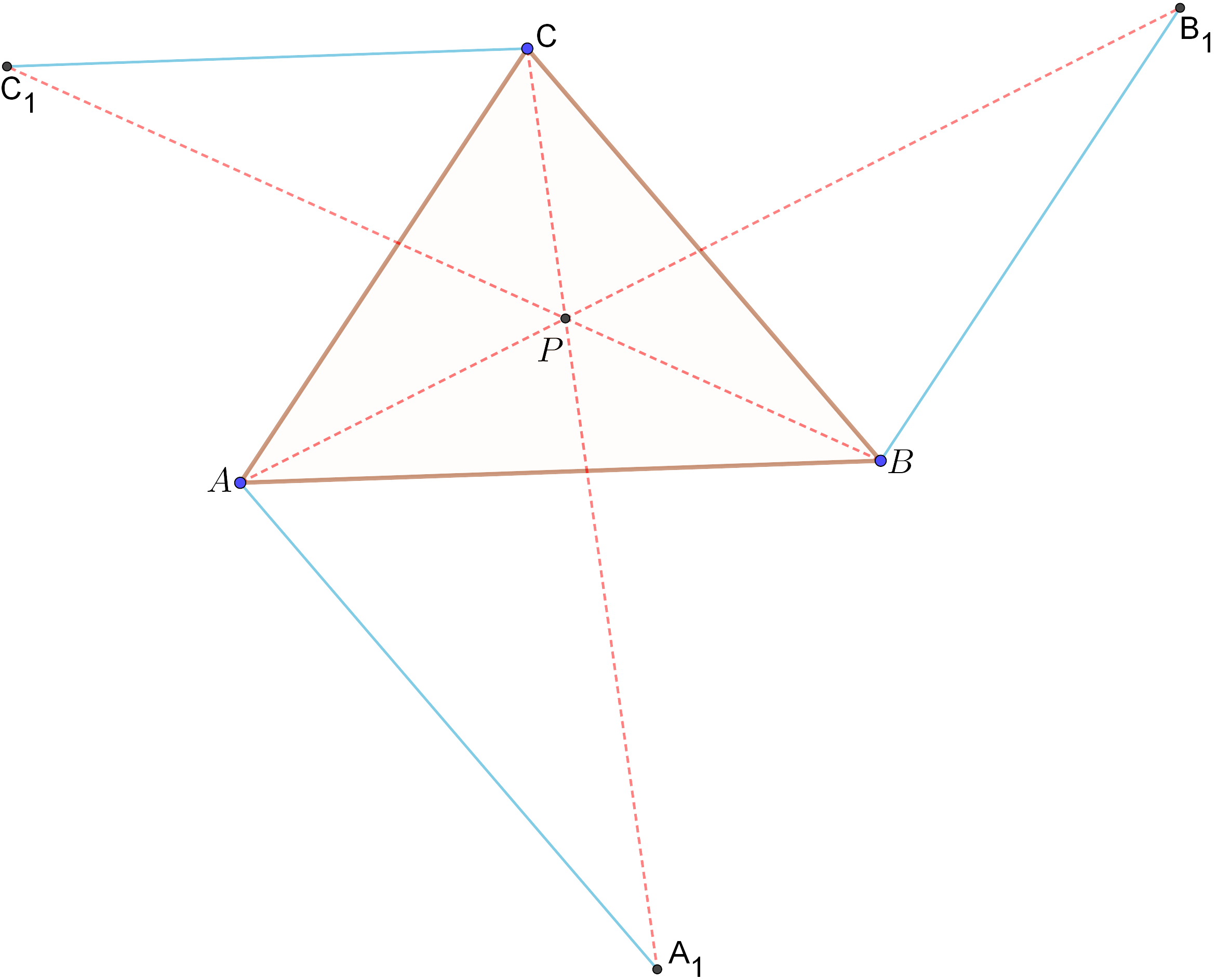}
	\caption{An application of concurrent condition of three straight lines 1} \label{fig:YizhongSanxiangongdianTiaojianDeYingyong1}
\end{figure}

\begin{figure}[h] 
	\centering
	\includegraphics[totalheight=6cm]{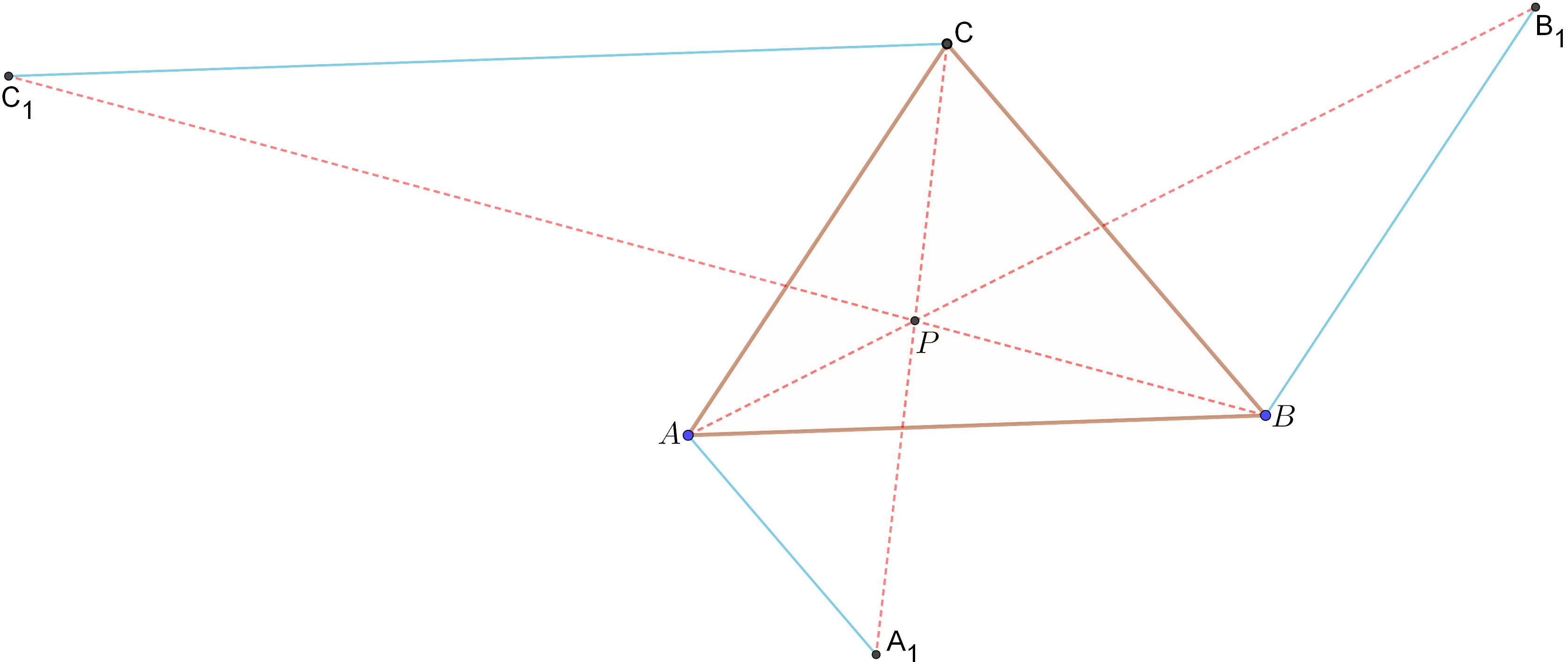}
	\caption{An application of concurrent condition of three straight lines 2} \label{fig:YizhongSanxiangongdianTiaojianDeYingyong2}
\end{figure}

\begin{solution}
	Let ${{\theta }_{A}}=-B$, ${{\theta }_{B}}=-C$, ${{\theta }_{C}}=-A$ in the left side of the formula in theorem \ref{thm:YiGeSanxiangongdianDeChongfenbiyaotiaojian},  it is obtained that
	\[\begin{aligned}
		& \frac{{{k }_{AB}}\sin \left( A-{{\theta }_{A}} \right)}{\sin B-{{k }_{AB}}\sin \left( B+{{\theta }_{A}} \right)}\cdot \frac{{{k }_{BC}}\sin \left( B-{{\theta }_{B}} \right)}{\sin C-{{k }_{BC}}\sin \left( C+{{\theta }_{B}} \right)}\cdot \frac{{{k }_{CA}}\sin \left( C-{{\theta }_{C}} \right)}{\sin A-{{k }_{CA}}\sin \left( A+{{\theta }_{C}} \right)} \\ 
		& =\frac{{{k }_{AB}}\sin \left( A+B \right)}{\sin B}\cdot \frac{{{k }_{BC}}\sin \left( B+C \right)}{\sin C}\cdot \frac{{{k }_{CA}}\sin \left( C+A \right)}{\sin A} \\ 
		& =\frac{{{k }_{AB}}\sin \left( \pi -C \right)}{\sin B}\cdot \frac{{{k }_{BC}}\sin \left( \pi -A \right)}{\sin C}\cdot \frac{{{k }_{CA}}\sin \left( \pi -B \right)}{\sin A} \\ 
		& =\frac{{{k }_{AB}}\sin C}{\sin B}\cdot \frac{{{k }_{BC}}\sin A}{\sin C}\cdot \frac{{{k }_{CA}}\sin B}{\sin A} \\ 
		& ={{k }_{AB}}{{k }_{BC}}{{k }_{CA}}. \\ 
	\end{aligned}\]
	
	So a condition for three lines to be concurrent is obtained.
	\[{{k }_{AB}}{{k }_{BC}}{{k }_{CA}}=1.\]
	
	In figure \ref{fig:YizhongSanxiangongdianTiaojianDeYingyong1}, ${{k }_{AB}}={{k }_{BC}}={{k }_{CA}}=1$. In figure \ref{fig:YizhongSanxiangongdianTiaojianDeYingyong2}, ${{k }_{AB}}={1}/{2}\;$, ${{k }_{BC}}=1$, ${{k }_{CA}}=2$.
\end{solution}
\hfill $\diamond$\par

\section{Application of frame component in trilinear coordinate}\label{BiaojiafenliangZaiSanxianzuobiaoZhongDeYingyong}
This section will discuss trilinear coordinates. Firstly, I introduce the concepts of true trilinear coordinates and trilinear coordinates. Then study the relationship between the trilinear coordinates and the frame components, and investigate the application of the frame components in the trilinear coordinates.

Based on the result of theorem \ref{thm:YouBiaojiafenliangJisuanZuobiao}, I propose the following concept of true trilinear coordinates.


\begin{definition}{True trilinear coordinate}{ZhenshiSanxianzuobiao}\label{ZhenshiSanxianzuobiao}
	Given a $\triangle ABC$ and a point $P$ on the $\triangle ABC$ plane, the area of $\triangle ABC$ is $S$. The frame components of the $P$ (an IC) are $y_{\mathcal{R}\left( BC \right)}^{P}$, $y_{\mathcal{R}\left( CA \right)}^{P}$, $y_{\mathcal{R}\left( AB \right)}^{P}$, respectively. Then, an ordered array $\left( y_{\mathcal{R}\left( BC \right)}^{P},y_{\mathcal{R}\left( CA \right)}^{P},y_{\mathcal{R}\left( AB \right)}^{P} \right)$ formed by three real numbers $y_{\mathcal{R}\left( BC \right)}^{P}$, $y_{\mathcal{R}\left( CA \right)}^{P}$, $y_{\mathcal{R}\left( AB \right)}^{P}$ is called the true trilinear coordinate of point $P$, abbreviated as the true trilinear coordinate. Among them, $y_{\mathcal{R}\left( BC \right)}^{P}$, $y_{\mathcal{R}\left( CA \right)}^{P}$, $y_{\mathcal{R}\left( AB \right)}^{P}$ are the vertical  Cartesian coordinates in the Edge-Axis coordinate systems $\mathcal{R}\left( BC \right)$, $\mathcal{R}\left( CA \right)$ and $\mathcal{R}\left( AB \right)$ respectively. The coordinates are given by the following equations:
	\[y_{\mathcal{R}\left( BC \right)}^{P}=\frac{2S}{a}\alpha _{A}^{P},\]
	\[y_{\mathcal{R}\left( CA \right)}^{P}=\frac{2S}{b}\alpha _{B}^{P},\]
	\[y_{\mathcal{R}\left( AB \right)}^{P}=\frac{2S}{c}\alpha _{C}^{P}.\]
\end{definition}

The true trilinear coordinates are uniquely determined by $\triangle ABC$ and point $P$, which means that the true trilinear coordinates are unique (see theorem \ref{thm:PingmianShangDeDianYuBiaojiafenliangYiyiduiying})

%

\begin{definition}{Trilinear coordinate}{Sanxianzuobiao}\label{Sanxianzuobiao}
	Given a $\triangle ABC$, the ordered array $\left( h_{\mathcal{R}\left( BC \right)}^{P},h_{\mathcal{R}\left( CA \right)}^{P},h_{\mathcal{R}\left( AB \right)}^{P} \right)$ obtained by reducing the common factor of three real numbers $y_{\mathcal{R}\left( BC \right)}^{P}$, $y_{\mathcal{R}\left( CA \right)}^{P}$, $y_{\mathcal{R}\left( AB \right)}^{P}$ in the true trilinear coordinates $\left( y_{\mathcal{R}\left( BC \right)}^{P},y_{\mathcal{R}\left( CA \right)}^{P},y_{\mathcal{R}\left( AB \right)}^{P} \right)$ is called the trilinear coordinates of point $P$ with respect to $\triangle ABC$, abbreviated as trilinear coordinates, where:
	\[kh_{\mathcal{R}\left( BC \right)}^{P}=y_{\mathcal{R}\left( BC \right)}^{P},\]
	\[kh_{\mathcal{R}\left( CA \right)}^{P}=y_{\mathcal{R}\left( CA \right)}^{P},\]
	\[kh_{\mathcal{R}\left( AB \right)}^{P}=y_{\mathcal{R}\left( AB \right)}^{P}.\]
	
	Where $k$ is called the common factor, $k\ne 0$.
\end{definition}

The essence of the above definition is to decompose a true trilinear coordinate into a product of a non-zero common factor and a trilinear coordinate. Due to the absence of a specific non-zero factor $k$ in the above definition, the trilinear coordinates are not unique. Obviously, if $k=1$, the true trilinear coordinates are also trilinear coordinates, which means that the true trilinear coordinates are also a type of trilinear coordinates. However, the opposite is not true, as trilinear coordinates may not necessarily be true trilinear coordinates.

The following theorem provides a new method for solving trilinear coordinates.


\begin{theorem}{Calculate trilinear coordinate by frame component, Daiyuan Zhang}{GenjuBiaojiafenliangQiuSanxianzuobiao}\label{GenjuBiaojiafenliangQiuSanxianzuobiao}
	Given a $\triangle ABC$ and a point $P$ on the plane of $\triangle ABC$, where the frame components of point $P$ are $\alpha _{A}^{P}$, $\alpha _{B}^{P}$, $\alpha _{C}^{P}$, then the trilinear coordinates of point $P$ are:
	\[h_{\mathcal{R}\left( BC \right)}^{P}:h_{\mathcal{R}\left( CA \right)}^{P}:h_{\mathcal{R}\left( AB \right)}^{P}=\frac{\alpha _{A}^{P}}{a}:\frac{\alpha _{B}^{P}}{b}:\frac{\alpha _{C}^{P}}{c}.\]
\end{theorem}

%

\begin{proof}
	According to the theorem \ref{thm:YouBiaojiafenliangJisuanZuobiao}, it is obtained that:
	\[y_{\mathcal{R}\left( BC \right)}^{P}=\frac{2S}{a}\alpha _{A}^{P},\]
	\[y_{\mathcal{R}\left( CA \right)}^{P}=\frac{2S}{b}\alpha _{B}^{P},\]
	\[y_{\mathcal{R}\left( AB \right)}^{P}=\frac{2S}{c}\alpha _{C}^{P}.\]
	
	By directly using the above definition, the common factor $2S$ can be reduced to obtain:
	\[h_{\mathcal{R}\left( BC \right)}^{P}:h_{\mathcal{R}\left( CA \right)}^{P}:h_{\mathcal{R}\left( AB \right)}^{P}=\frac{2S}{a}\alpha _{A}^{P}:\frac{2S}{b}\alpha _{B}^{P}:\frac{2S}{c}\alpha _{C}^{P}=\frac{\alpha _{A}^{P}}{a}:\frac{\alpha _{B}^{P}}{b}:\frac{\alpha _{C}^{P}}{c}.\]
\end{proof}
\hfill $\square$\par

The above theorem provides a new method for solving trilinear coordinates, which directly uses the frame component to solve trilinear coordinates, avoiding the difficulty of using Euclidean geometry methods to solve trilinear coordinates. It should be noted that the trilinear coordinates are not unique, so the trilinear coordinates obtained from the above theorem are only one form of them.


For example, if point $P$ is the centroid $G$ of $\triangle ABC$, then $\alpha _{A}^{G}=\alpha _{B}^{G}=\alpha _{C}^{G}={1}/{3}\;$, substituting them into the above theorem yields
\[h_{\mathcal{R}\left( BC \right)}^{G}:h_{\mathcal{R}\left( CA \right)}^{G}:h_{\mathcal{R}\left( AB \right)}^{G}=\frac{\alpha _{A}^{G}}{a}:\frac{\alpha _{B}^{G}}{b}:\frac{\alpha _{C}^{G}}{c}=\frac{1}{3a}:\frac{1}{3b}:\frac{1}{3c}=\frac{1}{a}:\frac{1}{b}:\frac{1}{c}.\]

Multiply $abc$ to the right of the above equation to obtain
\[h_{\mathcal{R}\left( BC \right)}^{G}:h_{\mathcal{R}\left( CA \right)}^{G}:h_{\mathcal{R}\left( AB \right)}^{G}=\frac{\alpha _{A}^{G}}{a}:\frac{\alpha _{B}^{G}}{b}:\frac{\alpha _{C}^{G}}{c}=\frac{1}{a}:\frac{1}{b}:\frac{1}{c}=bc:ca:ab.\]

Therefore, the trilinear coordinates of the centroid $G$ are not unique. They can be either $\left( {1}/{a}\;,{1}/{b}\;,{1}/{c}\; \right)$ or $\left( bc,ca,ab \right)$, and of course, they can have other forms.

But the true trilinear coordinates of the centroid $G$ are unique, determined solely by $\triangle ABC$ and the centroid $G$. The three quantities of the true trilinear coordinates of the centroid $G$ are:
\[y_{\mathcal{R}\left( BC \right)}^{P}=\frac{2S}{a}\alpha _{A}^{P}=\frac{2S}{3a},\]
\[y_{\mathcal{R}\left( CA \right)}^{P}=\frac{2S}{b}\alpha _{B}^{P}=\frac{2S}{3b},\]
\[y_{\mathcal{R}\left( AB \right)}^{P}=\frac{2S}{c}\alpha _{C}^{P}=\frac{2S}{3c}.\]


\begin{example}{}\label{QiuNeixinChuixinWaixinDeSanxianzuobiao}
	Find the trilinear coordinates of the incenter, orthocenter, circumcenter of $\triangle ABC$, only one form is required.
\end{example}

\begin{solution}
	For the incenter $I$, the frame components have already been calculated earlier:
	\[\alpha _{A}^{I}=\frac{a}{a+b+c},\]
	\[\alpha _{B}^{I}=\frac{b}{a+b+c},\]
	\[\alpha _{C}^{I}=\frac{c}{a+b+c}.\]
	
	So, according to the theorem \ref{thm:GenjuBiaojiafenliangQiuSanxianzuobiao}, the trilinear coordinates of incenter are:
	\begin{align*}
		& h_{\mathcal{R}\left( BC \right)}^{I}:h_{\mathcal{R}\left( CA \right)}^{I}:h_{\mathcal{R}\left( AB \right)}^{I} \\ 
		& =\frac{\alpha _{A}^{I}}{a}:\frac{\alpha _{B}^{I}}{b}:\frac{\alpha _{C}^{I}}{c} \\ 
		& =\frac{\frac{a}{a+b+c}}{a}:\frac{\frac{b}{a+b+c}}{b}:\frac{\frac{c}{a+b+c}}{c} \\ 
		& =\frac{1}{a+b+c}:\frac{1}{a+b+c}:\frac{1}{a+b+c} \\ 
		& =1:1:1.  
	\end{align*}
	
	For the orthocenter $H$, the frame components have already been calculated earlier:
	\[\alpha _{A}^{H}=\frac{\tan A}{\tan A+\tan B+\tan C},\]
	\[\alpha _{B}^{H}=\frac{\tan B}{\tan A+\tan B+\tan C},\]
	\[\alpha _{C}^{H}=\frac{\tan C}{\tan A+\tan B+\tan C}.\]
	
	So, according to the theorem \ref{thm:GenjuBiaojiafenliangQiuSanxianzuobiao} and the sine theorem, the trilinear coordinates of orthocenter are:
	\begin{align*}
		& h_{\mathcal{R}\left( BC \right)}^{H}:h_{\mathcal{R}\left( CA \right)}^{H}:h_{\mathcal{R}\left( AB \right)}^{H} \\ 
		& =\frac{\alpha _{A}^{H}}{a}:\frac{\alpha _{B}^{H}}{b}:\frac{\alpha _{C}^{H}}{c} \\ 
		& =\frac{\tan A}{a\left( \tan A+\tan B+\tan C \right)}:\frac{\tan B}{b\left( \tan A+\tan B+\tan C \right)}:\frac{\tan C}{c\left( \tan A+\tan B+\tan C \right)} \\ 
		& =\frac{\tan A}{a}:\frac{\tan B}{b}:\frac{\tan C}{c} \\ 
		& =\frac{\frac{\sin A}{\cos A}}{2R\sin A}:\frac{\frac{\sin B}{\cos B}}{2R\sin B}:\frac{\frac{\sin B}{\cos B}}{2R\sin B}\\ 
		& =\sec A:\sec B:\sec C. 
	\end{align*}
	
	For the circumcenter $Q$, the frame components have already been calculated earlier:
	\[\alpha _{A}^{Q}=\frac{\sin 2A}{\sin 2A+\sin 2B+\sin 2C},\]
	\[\alpha _{B}^{Q}=\frac{\sin 2B}{\sin 2A+\sin 2B+\sin 2C},\]
	\[\alpha _{C}^{Q}=\frac{\sin 2C}{\sin 2A+\sin 2B+\sin 2C},\]
	
	So, according to the theorem \ref{thm:GenjuBiaojiafenliangQiuSanxianzuobiao} and the sine theorem, the trilinear coordinates of circumcenter are:
	\begin{align*}
		& h_{\mathcal{R}\left( BC \right)}^{Q}:h_{\mathcal{R}\left( CA \right)}^{Q}:h_{\mathcal{R}\left( AB \right)}^{Q} \\ 
		& =\frac{\alpha _{A}^{Q}}{a}:\frac{\alpha _{B}^{Q}}{b}:\frac{\alpha _{C}^{Q}}{c} \\ 
		& =\frac{\sin 2A}{a\left( \sin 2A+\sin 2B+\sin 2C \right)}:\frac{\sin 2B}{b\left( \sin 2A+\sin 2B+\sin 2C \right)}:\frac{\sin 2C}{c\left( \sin 2A+\sin 2B+\sin 2C \right)} \\ 
		& =\frac{\sin 2A}{a}:\frac{\sin 2B}{b}:\frac{\sin 2C}{c} \\ 
		& =\frac{2\sin A\cos A}{2R\sin A}:\frac{2\sin B\cos B}{2R\sin B}:\frac{2\sin C\cos C}{2R\sin C} \\ 
		& =\cos A:\cos B:\cos C. 
	\end{align*}
\end{solution}
\hfill $\diamond$\par

Trilinear coordinates are a simplified form of true trilinear coordinates. When calculating the trilinear coordinates, one should try to obtain a simple form as much as possible. The advantage of trilinear coordinates is that the expression is concise, and in some applications, the relative quantity relationship between three coordinates can be quickly observed through proportional relationships. The disadvantage of trilinear coordinates is that important information may be lost, such as the dimensional information of the original geometric quantity, which means that trilinear coordinates may not have geometric meaning. Usually, obtaining trilinear coordinates through certain methods is not the true trilinear coordinates and cannot be directly used as a geometric quantity.

The true trilinear coordinates cannot be reduced by a common factor. The true trilinear coordinates reflect the true values of the Cartesian ordinate in the Edge-Axis coordinate system of $\triangle ABC$. If point $P$ is located inside $\triangle ABC$, then $y_{\mathcal{R}\left( AB \right)}^{P}$, $y_{\mathcal{R}\left( BC \right)}^{P}$, $y_{\mathcal{R}\left( CA \right)}^{P}$ are the distances between point $P$ and $AB$, $BC$, $CA$, respectively. they have accurate numerical values and clear geometric meanings.

In short, the trilinear coordinates are only the relative sizes of three quantities and usually do not have a clear geometric meaning; The true trilinear coordinates represent the true sizes of three quantities and have clear geometric meanings.

The current question is how to calculate the frame components and the true trilinear coordinates based on the given trilinear coordinates?


\begin{theorem}{Calculate frame component by trilinear coordinate, Daiyuan Zhang}{GenjuSanxianzuobiaoQiuBiaojiafenliang}\label{GenjuSanxianzuobiaoQiuBiaojiafenliang}
	Given a $\triangle ABC$ and a point $P$ on the plane of $\triangle ABC$, where the trilinear coordinates of point $P$ are $\left( h_{\mathcal{R}\left( BC \right)}^{P},h_{\mathcal{R}\left( CA \right)}^{P},h_{\mathcal{R}\left( AB \right)}^{P} \right)$, then the frame components of point $P$ are:
	\[\alpha _{A}^{P}=\frac{h_{\mathcal{R}\left( BC \right)}^{P}a}{h_{\mathcal{R}\left( BC \right)}^{P}a+h_{\mathcal{R}\left( CA \right)}^{P}b+h_{\mathcal{R}\left( AB \right)}^{P}c},\]
	\[\alpha _{B}^{P}=\frac{h_{\mathcal{R}\left( CA \right)}^{P}b}{h_{\mathcal{R}\left( BC \right)}^{P}a+h_{\mathcal{R}\left( CA \right)}^{P}b+h_{\mathcal{R}\left( AB \right)}^{P}c},\]
	\[\alpha _{C}^{P}=\frac{h_{\mathcal{R}\left( AB \right)}^{P}c}{h_{\mathcal{R}\left( BC \right)}^{P}a+h_{\mathcal{R}\left( CA \right)}^{P}b+h_{\mathcal{R}\left( AB \right)}^{P}c}.\]
\end{theorem}

\begin{proof}
	According to the definition of trilinear coordinates, assuming:
	\[kh_{\mathcal{R}\left( BC \right)}^{P}=y_{\mathcal{R}\left( BC \right)}^{P},\]
	\[kh_{\mathcal{R}\left( CA \right)}^{P}=y_{\mathcal{R}\left( CA \right)}^{P},\]
	\[kh_{\mathcal{R}\left( AB \right)}^{P}=y_{\mathcal{R}\left( AB \right)}^{P}.\]
	
	Where $k$ is the common factor, $k\ne 0$. The above equation can be understood as decomposing the true trilinear coordinates into the product of trilinear coordinates and a common non-zero factor. Next, we need to solve for the common factor. According to definition \ref{def:Sanxianzuobiao} and theorem \ref{thm:YouBiaojiafenliangJisuanZuobiao}, it can be concluded that:
	\[kh_{\mathcal{R}\left( BC \right)}^{P}=y_{\mathcal{R}\left( BC \right)}^{P}=\frac{2S}{a}\alpha _{A}^{P},\]
	\[kh_{\mathcal{R}\left( CA \right)}^{P}=y_{\mathcal{R}\left( CA \right)}^{P}=\frac{2S}{b}\alpha _{B}^{P},\]
	\[kh_{\mathcal{R}\left( AB \right)}^{P}=y_{\mathcal{R}\left( AB \right)}^{P}=\frac{2S}{c}\alpha _{C}^{P}.\]
	
	According to the above formulas, it can be obtained that:
	\[\frac{kh_{\mathcal{R}\left( BC \right)}^{P}a}{2S}=\alpha _{A}^{P},\]
	\[\frac{kh_{\mathcal{R}\left( CA \right)}^{P}b}{2S}=\alpha _{B}^{P},\]
	\[\frac{kh_{\mathcal{R}\left( AB \right)}^{P}c}{2S}=\alpha _{C}^{P}.\]
	
	Therefore
	\[\frac{kh_{\mathcal{R}\left( BC \right)}^{P}a+kh_{\mathcal{R}\left( CA \right)}^{P}b+kh_{\mathcal{R}\left( AB \right)}^{P}c}{2S}=\alpha _{A}^{P}+\alpha _{B}^{P}+\alpha _{C}^{P}=1,\]
	i.e.
	\[k=\frac{2S}{h_{\mathcal{R}\left( BC \right)}^{P}a+h_{\mathcal{R}\left( CA \right)}^{P}b+h_{\mathcal{R}\left( AB \right)}^{P}c}.\]
	
	So the frame components of point $P$ are:
	\[\alpha _{A}^{P}=\frac{kh_{\mathcal{R}\left( BC \right)}^{P}a}{2S}=\frac{h_{\mathcal{R}\left( BC \right)}^{P}a}{h_{\mathcal{R}\left( BC \right)}^{P}a+h_{\mathcal{R}\left( CA \right)}^{P}b+h_{\mathcal{R}\left( AB \right)}^{P}c},\]	
	\[\alpha _{B}^{P}=\frac{kh_{\mathcal{R}\left( CA \right)}^{P}b}{2S}=\frac{h_{\mathcal{R}\left( CA \right)}^{P}b}{h_{\mathcal{R}\left( BC \right)}^{P}a+h_{\mathcal{R}\left( CA \right)}^{P}b+h_{\mathcal{R}\left( AB \right)}^{P}c},\]
	\[\alpha _{C}^{P}=\frac{kh_{\mathcal{R}\left( AB \right)}^{P}c}{2S}=\frac{h_{\mathcal{R}\left( AB \right)}^{P}c}{h_{\mathcal{R}\left( BC \right)}^{P}a+h_{\mathcal{R}\left( CA \right)}^{P}b+h_{\mathcal{R}\left( AB \right)}^{P}c}.\]
\end{proof}
\hfill $\square$\par

So, if the trilinear coordinates are known, the frame components can be directly calculated based on the above formula. This provides another method for obtaining the frame components.


\begin{theorem}{Calculate true trilinear coordinates based on trilinear coordinates}{GenjuSanxianzuobiaoQiuZhenshisanxianzuobiao}\label{GenjuSanxianzuobiaoQiuZhenshisanxianzuobiao}
	Given a $\triangle ABC$ and a point $P$ on the plane of $\triangle ABC$, the trilinear coordinates of point $P$ are $\left( h_{\mathcal{R}\left( BC \right)}^{P},h_{\mathcal{R}\left( CA \right)}^{P},h_{\mathcal{R}\left( AB \right)}^{P} \right)$, then true trilinear coordinates of point $P$ are:
	\[y_{\mathcal{R}\left( BC \right)}^{P}=\frac{2Sh_{\mathcal{R}\left( BC \right)}^{P}}{h_{\mathcal{R}\left( BC \right)}^{P}a+h_{\mathcal{R}\left( CA \right)}^{P}b+h_{\mathcal{R}\left( AB \right)}^{P}c},\]
	\[y_{\mathcal{R}\left( CA \right)}^{P}=\frac{2Sh_{\mathcal{R}\left( CA \right)}^{P}}{h_{\mathcal{R}\left( BC \right)}^{P}a+h_{\mathcal{R}\left( CA \right)}^{P}b+h_{\mathcal{R}\left( AB \right)}^{P}c},\]
	\[y_{\mathcal{R}\left( AB \right)}^{P}=\frac{2Sh_{\mathcal{R}\left( AB \right)}^{P}}{h_{\mathcal{R}\left( BC \right)}^{P}a+h_{\mathcal{R}\left( CA \right)}^{P}b+h_{\mathcal{R}\left( AB \right)}^{P}c}.\]
\end{theorem}

\begin{proof}
	According to the definition of trilinear coordinates, assuming:
	\[kh_{\mathcal{R}\left( BC \right)}^{P}=y_{\mathcal{R}\left( BC \right)}^{P},\]
	\[kh_{\mathcal{R}\left( CA \right)}^{P}=y_{\mathcal{R}\left( CA \right)}^{P},\]
	\[kh_{\mathcal{R}\left( AB \right)}^{P}=y_{\mathcal{R}\left( AB \right)}^{P}.\]
	
	Where $k$ is the common factor, $k\ne 0$. The above equation can be understood as decomposing the true trilinear coordinates into the product of trilinear coordinates and a common non-zero factor. Next, we need to solve for the common factor. According to definition \ref{def:Sanxianzuobiao} and theorem \ref{thm:YouBiaojiafenliangJisuanZuobiao}, it can be concluded that:
	\[kh_{\mathcal{R}\left( BC \right)}^{P}=y_{\mathcal{R}\left( BC \right)}^{P}=\frac{2S}{a}\alpha _{A}^{P},\]
	\[kh_{\mathcal{R}\left( CA \right)}^{P}=y_{\mathcal{R}\left( CA \right)}^{P}=\frac{2S}{b}\alpha _{B}^{P},\]
	\[kh_{\mathcal{R}\left( AB \right)}^{P}=y_{\mathcal{R}\left( AB \right)}^{P}=\frac{2S}{c}\alpha _{C}^{P}.\]
	
	According to the above formulas, it can be obtained that:
	\[\frac{kh_{\mathcal{R}\left( BC \right)}^{P}a}{2S}=\alpha _{A}^{P},\]
	\[\frac{kh_{\mathcal{R}\left( CA \right)}^{P}b}{2S}=\alpha _{B}^{P},\]
	\[\frac{kh_{\mathcal{R}\left( AB \right)}^{P}c}{2S}=\alpha _{C}^{P}.\]
	
	Therefore
	\[\frac{kh_{\mathcal{R}\left( BC \right)}^{P}a+kh_{\mathcal{R}\left( CA \right)}^{P}b+kh_{\mathcal{R}\left( AB \right)}^{P}c}{2S}=\alpha _{A}^{P}+\alpha _{B}^{P}+\alpha _{C}^{P}=1,\]
	i.e.
	\[k=\frac{2S}{h_{\mathcal{R}\left( BC \right)}^{P}a+h_{\mathcal{R}\left( CA \right)}^{P}b+h_{\mathcal{R}\left( AB \right)}^{P}c}.\]
	
	So the true trilinear coordinates of point $P$ are:
	\[y_{\mathcal{R}\left( BC \right)}^{P}=kh_{\mathcal{R}\left( BC \right)}^{P}=\frac{2Sh_{\mathcal{R}\left( BC \right)}^{P}}{h_{\mathcal{R}\left( BC \right)}^{P}a+h_{\mathcal{R}\left( CA \right)}^{P}b+h_{\mathcal{R}\left( AB \right)}^{P}c},\]
	\[y_{\mathcal{R}\left( CA \right)}^{P}=kh_{\mathcal{R}\left( CA \right)}^{P}=\frac{2Sh_{\mathcal{R}\left( CA \right)}^{P}}{h_{\mathcal{R}\left( BC \right)}^{P}a+h_{\mathcal{R}\left( CA \right)}^{P}b+h_{\mathcal{R}\left( AB \right)}^{P}c},\]
	\[y_{\mathcal{R}\left( AB \right)}^{P}=kh_{\mathcal{R}\left( AB \right)}^{P}=\frac{2Sh_{\mathcal{R}\left( AB \right)}^{P}}{h_{\mathcal{R}\left( BC \right)}^{P}a+h_{\mathcal{R}\left( CA \right)}^{P}b+h_{\mathcal{R}\left( AB \right)}^{P}c}.\]
\end{proof}
\hfill $\square$\par

The trilinear coordinates can be seen as a proportion, and sometimes this proportion can have intuitive meaning.


\begin{example}{}\label{QiuNeixinDeBiqaojiafenliangHeZhenshisanxianzuobiao}
	Find the frame components and true trilinear coordinates of the incenter $I$ of $\triangle ABC$.
\end{example}

\begin{solution}
	Starting from basic geometric concepts, the distance between the incenter  $I$ of a triangle and its three sides is equal. Although it is difficult to see at a glance the actual value of the distance between a given triangle and its incenter $I$, these three equal distances can be immediately written as a continued proportion, which is $1:1:1 $. From this intuitive geometric concept, the trilinear coordinates of the incenter $I$ can be immediately obtained as $h_{\mathcal{R}\left( BC \right)}^{I}=h_{\mathcal{R}\left( CA \right)}^{I}=h_{\mathcal{R}\left( AB \right)}^{I}=1$. Therefore, according to theorem \ref{thm:GenjuSanxianzuobiaoQiuBiaojiafenliang}, the frame components of the incenter $I$ are:
	\[\alpha _{A}^{I}=\frac{h_{\mathcal{R}\left( BC \right)}^{I}a}{h_{\mathcal{R}\left( BC \right)}^{I}a+h_{\mathcal{R}\left( CA \right)}^{I}b+h_{\mathcal{R}\left( AB \right)}^{I}c}=\frac{a}{a+b+c}=\frac{a}{2p},\]
	\[\alpha _{B}^{I}=\frac{h_{\mathcal{R}\left( CA \right)}^{I}b}{h_{\mathcal{R}\left( BC \right)}^{I}a+h_{\mathcal{R}\left( CA \right)}^{I}b+h_{\mathcal{R}\left( AB \right)}^{I}c}=\frac{b}{a+b+c}=\frac{b}{2p},\]
	\[\alpha _{C}^{I}=\frac{h_{\mathcal{R}\left( AB \right)}^{I}c}{h_{\mathcal{R}\left( BC \right)}^{I}a+h_{\mathcal{R}\left( CA \right)}^{I}b+h_{\mathcal{R}\left( AB \right)}^{I}c}=\frac{c}{a+b+c}=\frac{c}{2p}.\]
	
	So, according to the above theorem, the true trilinear coordinates of incenter $I$ are:
	\[y_{\mathcal{R}\left( BC \right)}^{I}=\frac{2Sh_{\mathcal{R}\left( BC \right)}^{I}}{h_{\mathcal{R}\left( BC \right)}^{I}a+h_{\mathcal{R}\left( CA \right)}^{I}b+h_{\mathcal{R}\left( AB \right)}^{I}c}=\frac{2S}{a+b+c}=\frac{S}{p},\]
	\[y_{\mathcal{R}\left( CA \right)}^{I}=\frac{2Sh_{\mathcal{R}\left( CA \right)}^{I}}{h_{\mathcal{R}\left( BC \right)}^{I}a+h_{\mathcal{R}\left( CA \right)}^{I}b+h_{\mathcal{R}\left( AB \right)}^{I}c}=\frac{2S}{a+b+c}=\frac{S}{p},\]
	\[y_{\mathcal{R}\left( AB \right)}^{I}=\frac{2Sh_{\mathcal{R}\left( AB \right)}^{I}}{h_{\mathcal{R}\left( BC \right)}^{I}a+h_{\mathcal{R}\left( CA \right)}^{I}b+h_{\mathcal{R}\left( AB \right)}^{I}c}=\frac{2S}{a+b+c}=\frac{S}{p}.\]
\end{solution}
\hfill $\diamond$\par

For the frame components of incenter $I$, readers can compare the solution method here with the method given in Section \ref{Sec8.2}. Obviously, the method used here is simpler.

It should be noted that the trilinear coordinates of the incenter $I$ are $h_{\mathcal{R}\left( BC \right)}^{I}=h_{\mathcal{R}\left( CA \right)}^{I}=h_{\mathcal{R}\left( AB \right)}^{I}=1$. The fact that all three trilinear coordinates are 1 only indicates that the three trilinear coordinates are equal, and 1 is not the true trilinear coordinate, 1 is a pure number without dimensions, so the three trilinear coordinates do not represent the geometric quantities of a triangle. For Cartesian coordinate systems, the coordinates of points have a dimension of length. Therefore, the trilinear coordinates are not directly used for calculating geometric quantities (such as distance, area, etc.). The trilinear coordinates can also be seen as virtual coordinates. To be able to calculate geometric quantities, it is necessary to calculate the frame components or the true trilinear coordinates.


\begin{example}{}\label{QiuDuiyingDeBangxinDeBiaojiafenliangHeZhenshisanxianzuobiao}
	Find the frame components and the true trilinear coordinates of the excenter ${{E}_{A}}$ corresponding to $\angle A$ of $\triangle ABC$.
\end{example}

\begin{solution}
	According to the geometric concept of the excenter, for the excenter ${{E}_{A}}$ corresponding to $\angle A$, its three line coordinates are: $h_{\mathcal{R}\left( BC \right)}^{I}=-1$, $h_{\mathcal{R}\left( CA \right)}^{I}=1$, $h_{\mathcal{R}\left( AB \right)}^{I}=1$. Therefore, according to theorem \ref{thm:GenjuSanxianzuobiaoQiuBiaojiafenliang}, the frame components of the excenter ${{E}_{A}}$ are:
	\[\alpha _{A}^{{{E}_{A}}}=\frac{h_{\mathcal{R}\left( BC \right)}^{{{E}_{A}}}a}{h_{\mathcal{R}\left( BC \right)}^{{{E}_{A}}}a+h_{\mathcal{R}\left( CA \right)}^{{{E}_{A}}}b+h_{\mathcal{R}\left( AB \right)}^{{{E}_{A}}}c}=\frac{-a}{-a+b+c}=-\frac{a}{2\left( p-a \right)},\]
	\[\alpha _{B}^{{{E}_{A}}}=\frac{h_{\mathcal{R}\left( CA \right)}^{{{E}_{A}}}b}{h_{\mathcal{R}\left( BC \right)}^{{{E}_{A}}}a+h_{\mathcal{R}\left( CA \right)}^{{{E}_{A}}}b+h_{\mathcal{R}\left( AB \right)}^{{{E}_{A}}}c}=\frac{b}{-a+b+c}=\frac{b}{2\left( p-a \right)},\]
	\[\alpha _{C}^{{{E}_{A}}}=\frac{h_{\mathcal{R}\left( AB \right)}^{{{E}_{A}}}c}{h_{\mathcal{R}\left( BC \right)}^{{{E}_{A}}}a+h_{\mathcal{R}\left( CA \right)}^{{{E}_{A}}}b+h_{\mathcal{R}\left( AB \right)}^{{{E}_{A}}}c}=\frac{c}{-a+b+c}=\frac{c}{2\left( p-a \right)}.\]
	
	According to theorem \ref{thm:GenjuSanxianzuobiaoQiuZhenshisanxianzuobiao}, the true trilinear coordinates of the excenter ${{E}_{A}}$ are:
	\[y_{\mathcal{R}\left( BC \right)}^{{{E}_{A}}}=\frac{2Sh_{\mathcal{R}\left( BC \right)}^{{{E}_{A}}}}{h_{\mathcal{R}\left( BC \right)}^{{{E}_{A}}}a+h_{\mathcal{R}\left( CA \right)}^{{{E}_{A}}}b+h_{\mathcal{R}\left( AB \right)}^{{{E}_{A}}}c}=\frac{-2S}{-a+b+c}=-\frac{S}{p-a},\]
	\[y_{\mathcal{R}\left( CA \right)}^{{{E}_{A}}}=\frac{2Sh_{\mathcal{R}\left( CA \right)}^{{{E}_{A}}}}{h_{\mathcal{R}\left( BC \right)}^{{{E}_{A}}}a+h_{\mathcal{R}\left( CA \right)}^{{{E}_{A}}}b+h_{\mathcal{R}\left( AB \right)}^{{{E}_{A}}}c}=\frac{2S}{-a+b+c}=\frac{S}{p-a},\]
	\[y_{\mathcal{R}\left( AB \right)}^{{{E}_{A}}}=\frac{2Sh_{\mathcal{R}\left( AB \right)}^{{{E}_{A}}}}{h_{\mathcal{R}\left( BC \right)}^{{{E}_{A}}}a+h_{\mathcal{R}\left( CA \right)}^{{{E}_{A}}}b+h_{\mathcal{R}\left( AB \right)}^{{{E}_{A}}}c}=\frac{2S}{-a+b+c}=\frac{S}{p-a}.\]
\end{solution}
\hfill $\diamond$\par


Trilinear coordinates have been developed for many years, and people have accumulated rich achievements. The trilinear coordinates I introduced can only be given a new name, which is the true trilinear coordinates. I think theorems \ref{thm:GenjuSanxianzuobiaoQiuBiaojiafenliang} and \ref{thm:GenjuBiaojiafenliangQiuSanxianzuobiao} are both important. theorem \ref{thm:GenjuSanxianzuobiaoQiuBiaojiafenliang} establishes the relationship between trilinear coordinates and frame components. This means that Intercenter Geometry can fully utilize the rich achievements that people have accumulated. In fact, no matter what method is used to obtain the trilinear coordinates of two points, the frame components of these two points can be calculated using theorem \ref{thm:GenjuSanxianzuobiaoQiuBiaojiafenliang}, and then the distance between these two points can be calculated (see Chapter \ref{Ch12} below). theorem \ref{thm:GenjuBiaojiafenliangQiuSanxianzuobiao} provides a new method for solving trilinear coordinates. The theory and methods of Intercenter Geometry can be used to solve the trilinear coordinates. theorem \ref{thm:GenjuYizhidianDeBiaojiafenliangQiuWeizhidianDeBiaojiafenliangGongshi} and its related corollaries can theoretically determine the frame component of any point. Therefore, theorem \ref{thm:GenjuBiaojiafenliangQiuSanxianzuobiao} actually provides a method for solving the trilinear coordinates of any point, which promotes the research and application of trilinear coordinates.

\section{Application of frame component in analytical geometry}
	
In Plane Intercenter Geometry, the frame component is represented by the lengths of the three sides of a given triangle. Therefore, the distance between any two points can also be represented by the lengths of the three sides of the triangle (see Chapter \ref{Ch12} below). For impatient readers, you can first browse through the content of the following chapters to experience the charm of Intercenter Geometry. 

On the contrary, as long as the frame component is obtained, it can be easily converted into the representation of Cartesian coordinates.

Assuming a $\triangle ABC$ is given on a plane, a Cartesian coordinate system is established on the plane where $\triangle ABC$ is located. The Cartesian coordinates of the three vertices of $\triangle ABC$ are: $A\left( {{x}_{A}},{{y}_{A}} \right)$, $B\left( {{x}_{B}},{{y}_{B}} \right)$, $C\left( {{x}_{C}},{{y}_{C}} \right)$, then:
\[a:=BC=\sqrt{{{\left( {{x}_{B}}-{{x}_{C}} \right)}^{2}}+{{\left( {{y}_{B}}-{{y}_{C}} \right)}^{2}}},\]
\[b:=CA=\sqrt{{{\left( {{x}_{C}}-{{x}_{A}} \right)}^{2}}+{{\left( {{y}_{C}}-{{y}_{A}} \right)}^{2}}},\]
\[c:=AB=\sqrt{{{\left( {{x}_{A}}-{{x}_{B}} \right)}^{2}}+{{\left( {{y}_{A}}-{{y}_{B}} \right)}^{2}}}.\]


The above are the formulas for calculating the side lengths represented by the coordinates of the three vertices of the triangle. The above formulas are called the side lengths formulas of analytic geometry. Analytical geometry is also commonly referred to as coordinate geometry, and the above formulas are also known as the side lengths formulas of coordinate geometry.

The frame component obtained using Intercenter Geometry is a function of the lengths of the three sides, which can be written as:
\[\alpha _{A}^{P}=\alpha _{A}^{P}\left( a,b,c \right),\quad \alpha _{B}^{P}=\alpha _{B}^{P}\left( a,b,c \right),\quad \alpha _{C}^{P}=\alpha _{C}^{P}\left( a,b,c \right).\]	

\begin{align*}
	\alpha _{A}^{P}& =\alpha _{A}^{P}\left( a,b,c \right) \\ 
	& =\alpha _{A}^{P}\left( \sqrt{{{\left( {{x}_{B}}-{{x}_{C}} \right)}^{2}}+{{\left( {{y}_{B}}-{{y}_{C}} \right)}^{2}}},\sqrt{{{\left( {{x}_{C}}-{{x}_{A}} \right)}^{2}}+{{\left( {{y}_{C}}-{{y}_{A}} \right)}^{2}}},\sqrt{{{\left( {{x}_{A}}-{{x}_{B}} \right)}^{2}}+{{\left( {{y}_{A}}-{{y}_{B}} \right)}^{2}}} \right),  
\end{align*}
\begin{align*}
	\alpha _{B}^{P}& =\alpha _{B}^{P}\left( a,b,c \right) \\ 
	& =\alpha _{B}^{P}\left( \sqrt{{{\left( {{x}_{B}}-{{x}_{C}} \right)}^{2}}+{{\left( {{y}_{B}}-{{y}_{C}} \right)}^{2}}},\sqrt{{{\left( {{x}_{C}}-{{x}_{A}} \right)}^{2}}+{{\left( {{y}_{C}}-{{y}_{A}} \right)}^{2}}},\sqrt{{{\left( {{x}_{A}}-{{x}_{B}} \right)}^{2}}+{{\left( {{y}_{A}}-{{y}_{B}} \right)}^{2}}} \right),  
\end{align*}
\begin{align*}
	\alpha _{C}^{P}& =\alpha _{C}^{P}\left( a,b,c \right) \\ 
	& =\alpha _{C}^{P}\left( \sqrt{{{\left( {{x}_{B}}-{{x}_{C}} \right)}^{2}}+{{\left( {{y}_{B}}-{{y}_{C}} \right)}^{2}}},\sqrt{{{\left( {{x}_{C}}-{{x}_{A}} \right)}^{2}}+{{\left( {{y}_{C}}-{{y}_{A}} \right)}^{2}}},\sqrt{{{\left( {{x}_{A}}-{{x}_{B}} \right)}^{2}}+{{\left( {{y}_{A}}-{{y}_{B}} \right)}^{2}}} \right).  
\end{align*}


The above equation indicates that as long as the side length in the frame component is replaced by the side length formula of analytical geometry, the frame component can be represented by the coordinates of the triangle vertices.

That is to say, the results of Intercenter Geometry can be easily converted into the results of analytic geometry (in vertex coordinate form), but conversely, converting the results of analytic geometry (in vertex coordinate form) into the results of Intercenter Geometry is not an easy task. Interested readers can find a few examples to calculate for themselves.

In order to find the frame components of any point on the plane, I introduced an Edge-Axis coordinate system. In Intercenter Geometry, I recommend using the Edge-Axis coordinate system. But if a single Cartesian coordinate system must be used, Intercenter Geometry can also provide some solutions.


\begin{theorem}{Mixed formula of frame components and Cartesian coordinates, Daiyuan Zhang}{BiaojiafenliangYuZhijiaozuobiaoDeHunheGongshi}\label{BiaojiafenliangYuZhijiaozuobiaoDeHunheGongshi}
	Given a $\triangle ABC$, point $P$ is located on the plane of $\triangle ABC$, and the frame components of point $P$ are $\alpha _{A}^{P}$, $\alpha _{B}^{P}$, $\alpha _{C}^{P}$, respectively. Establish a Cartesian coordinate system on the plane where $\triangle ABC$ is located, with the Cartesian coordinates of the three vertices of $\triangle ABC$ being $A\left( {{x}_{A}},{{y}_{A}} \right)$, $B\left( {{x}_{B}},{{y}_{B}} \right)$, $C\left( {{x}_{C}},{{y}_{C}} \right)$, respectively. The Cartesian coordinates of point $P$ are $P\left( {{x}_{P}},{{y}_{P}} \right)$. Then:
	\[{{x}_{P}}=\alpha _{A}^{P}{{x}_{A}}+\alpha _{B}^{P}{{x}_{B}}+\alpha _{C}^{P}{{x}_{C}},\]
	\[{{y}_{P}}=\alpha _{A}^{P}{{y}_{A}}+\alpha _{B}^{P}{{y}_{B}}+\alpha _{C}^{P}{{y}_{C}}.\]
\end{theorem}

%
%
%

\begin{proof}
	Coinciding the origin of the frame with the origin of the Cartesian coordinate system yields:
	\[\overrightarrow{OP}=\alpha _{A}^{P}\overrightarrow{OA}+\alpha _{B}^{P}\overrightarrow{OB}+\alpha _{C}^{P}\overrightarrow{OC},\]
	\[\overrightarrow{OA}={{x}_{A}}{{\bm{e}}_{x}}+{{y}_{A}}{{\bm{e}}_{y}},\]
	\[\overrightarrow{OB}={{x}_{B}}{{\bm{e}}_{x}}+{{y}_{B}}{{\bm{e}}_{y}},\]
	\[\overrightarrow{OC}={{x}_{C}}{{\bm{e}}_{x}}+{{y}_{C}}{{\bm{e}}_{y}}.\]
	
	Among them, ${{\bm{e}}_{x}}$ and ${{\bm{e}}_{y}}$ are the unit vectors of the $OX$ axis and $OY$ axis in Cartesian coordinates, respectively. Therefore, the formula in Cartesian coordinates is obtained:
	\begin{align*}
		\overrightarrow{OP}& =\alpha _{A}^{P}\overrightarrow{OA}+\alpha _{B}^{P}\overrightarrow{OB}+\alpha _{C}^{P}\overrightarrow{OC} \\ 
		& =\alpha _{A}^{P}\left( {{x}_{A}}{{\bm{e}}_{x}}+{{y}_{A}}{{\bm{e}}_{y}} \right)+\alpha _{B}^{P}\left( {{x}_{B}}{{\bm{e}}_{x}}+{{y}_{B}}{{\bm{e}}_{y}} \right)+\alpha _{C}^{P}\left( {{x}_{C}}{{\bm{e}}_{x}}+{{y}_{C}}{{\bm{e}}_{y}} \right) \\ 
		& =\left( \alpha _{A}^{P}{{x}_{A}}+\alpha _{B}^{P}{{x}_{B}}+\alpha _{C}^{P}{{x}_{C}} \right){{\bm{e}}_{x}}+\left( \alpha _{A}^{P}{{y}_{A}}+\alpha _{B}^{P}{{y}_{B}}+\alpha _{C}^{P}{{y}_{C}} \right){{\bm{e}}_{y}}.  
	\end{align*}
	
	On the other hand,
	\[\overrightarrow{OP}={{x}_{P}}{{\bm{e}}_{x}}+{{y}_{P}}{{\bm{e}}_{y}}.\]
	
	Due to the linear independence between ${{\bm{e}}_{x}}$ and ${{\bm{e}}_{y}}$, it is obtained that:
	\[{{x}_{P}}=\alpha _{A}^{P}{{x}_{A}}+\alpha _{B}^{P}{{x}_{B}}+\alpha _{C}^{P}{{x}_{C}},\]
	\[{{y}_{P}}=\alpha _{A}^{P}{{y}_{A}}+\alpha _{B}^{P}{{y}_{B}}+\alpha _{C}^{P}{{y}_{C}}.\]
\end{proof}
\hfill $\square$\par

The mixed formula given by the above theorem is symmetric and elegant, with concise expression. I think it is also a good result.


\begin{example}{}\label{Geiding}
	Given a $\triangle ABC$, the Cartesian coordinates of the three vertices of $\triangle ABC$ are $A\left( {{x}_{A}},{{y}_{A}} \right)$, $B\left( {{x}_{B}},{{y}_{B}} \right)$, $C\left( {{x}_{C}},{{y}_{C}} \right)$, respectively, calculate the Cartesian coordinates of the centroid $G$, incenter $I$, orthocenter $H$, and circumcenter  $Q$ for $\triangle ABC$.
\end{example}

%
%
%
%

\begin{solution}
	Using the above theorem, for the centroid $G$, it is obtained that:
	\[{{x}_{G}}=\alpha _{A}^{G}{{x}_{A}}+\alpha _{B}^{G}{{x}_{B}}+\alpha _{C}^{G}{{x}_{C}}=\frac{1}{3}\left( {{x}_{A}}+{{x}_{B}}+{{x}_{C}} \right),\]
	\[{{y}_{G}}=\alpha _{A}^{G}{{y}_{A}}+\alpha _{B}^{G}{{y}_{B}}+\alpha _{C}^{G}{{y}_{C}}=\frac{1}{3}\left( {{y}_{A}}+{{y}_{B}}+{{y}_{C}} \right).\]
	
	For the incenter $I$, it is obtained that:
	\[{{x}_{I}}=\alpha _{A}^{I}{{x}_{A}}+\alpha _{B}^{I}{{x}_{B}}+\alpha _{C}^{I}{{x}_{C}}=\frac{a{{x}_{A}}+b{{x}_{B}}+c{{x}_{C}}}{a+b+c},\]
	\[{{y}_{I}}=\alpha _{A}^{I}{{y}_{A}}+\alpha _{B}^{I}{{y}_{B}}+\alpha _{C}^{I}{{y}_{C}}=\frac{a{{y}_{A}}+b{{y}_{B}}+c{{y}_{C}}}{a+b+c}.\]
	
	For the orthocenter $H$, it is obtained that:
	\[{{x}_{H}}=\alpha _{A}^{H}{{x}_{A}}+\alpha _{B}^{H}{{x}_{B}}+\alpha _{C}^{H}{{x}_{C}}=\frac{{{x}_{A}}\tan A+{{x}_{B}}\tan B+{{x}_{C}}\tan C}{\tan A+\tan B+\tan C},\]
	\[{{y}_{H}}=\alpha _{A}^{H}{{y}_{A}}+\alpha _{B}^{H}{{y}_{B}}+\alpha _{C}^{H}{{y}_{C}}=\frac{{{y}_{A}}\tan A+{{y}_{B}}\tan B+{{y}_{C}}\tan C}{\tan A+\tan B+\tan C}.\]
	
	For the circumcenter $Q$, it is obtained that:
	\[{{x}_{Q}}=\alpha _{A}^{Q}{{x}_{A}}+\alpha _{B}^{Q}{{x}_{B}}+\alpha _{C}^{Q}{{x}_{C}}=\frac{{{x}_{A}}\sin 2A+{{x}_{B}}\sin 2B+{{x}_{C}}\sin 2C}{\sin 2A+\sin 2B+\sin 2C}.\]
	\[{{y}_{Q}}=\alpha _{A}^{Q}{{y}_{A}}+\alpha _{B}^{Q}{{y}_{B}}+\alpha _{C}^{Q}{{y}_{C}}=\frac{{{y}_{A}}\sin 2A+{{y}_{B}}\sin 2B+{{y}_{C}}\sin 2C}{\sin 2A+\sin 2B+\sin 2C}.\]
	
	Trigonometric functions can be transformed into cosine functions using identity transformations, then expressed as functions of side lengths using cosine theorem, and finally expressed as functions of triangle vertex coordinates using analytical geometry's side length formula. Interested readers can try it out themselves.
\end{solution}
\hfill $\diamond$\par 

\section{Frame transformation}\label{Biaojiabianhuan}
This section studies the frame transformation of triangles.


\begin{theorem}{Frame transformation formula 1, Daiyuan Zhang}{BiaojiabianhuanGongshi1}\label{BiaojiabianhuanGongshi1}
	Let two triangles $\triangle {{A}_{1}}{{B}_{1}}{{C}_{1}}$, $\triangle {{A}_{2}}{{B}_{2}}{{C}_{2}}$ be on the same plane. 
	The frame components of the three vertices ${{A}_{2}}$, ${{B}_{2}}$, ${{C}_{2}}$ of $\triangle {{A}_{2}}{{B}_{2}}{{C}_{2}}$ on the frame $\left( O;{{A}_{1}},{{B}_{1}},{{C}_{1}} \right)$ are $\alpha _{{{A}_{1}}}^{{{A}_{2}}}$, $\alpha _{{{B}_{1}}}^{{{A}_{2}}}$, $\alpha _{{{C}_{1}}}^{{{A}_{2}}}$; $\alpha _{{{A}_{1}}}^{{{B}_{2}}}$, $\alpha _{{{B}_{1}}}^{{{B}_{2}}}$, $\alpha _{{{C}_{1}}}^{{{B}_{2}}}$;  $\alpha _{{{A}_{1}}}^{{{C}_{2}}}$, $\alpha _{{{B}_{1}}}^{{{C}_{2}}}$, $\alpha _{{{C}_{1}}}^{{{C}_{2}}}$, respectively. Then
	\[\left( \begin{matrix}
		\overrightarrow{O{{A}_{2}}}  \\
		\overrightarrow{O{{B}_{2}}}  \\
		\overrightarrow{O{{C}_{2}}}  \\
	\end{matrix} \right)=\left( \begin{matrix}
		\alpha _{{{A}_{1}}}^{{{A}_{2}}} & \alpha _{{{B}_{1}}}^{{{A}_{2}}} & \alpha _{{{C}_{1}}}^{{{A}_{2}}}  \\
		\alpha _{{{A}_{1}}}^{{{B}_{2}}} & \alpha _{{{B}_{1}}}^{{{B}_{2}}} & \alpha _{{{C}_{1}}}^{{{B}_{2}}}  \\
		\alpha _{{{A}_{1}}}^{{{C}_{2}}} & \alpha _{{{B}_{1}}}^{{{C}_{2}}} & \alpha _{{{C}_{1}}}^{{{C}_{2}}}  \\
	\end{matrix} \right)\left( \begin{matrix}
		\overrightarrow{O{{A}_{1}}}  \\
		\overrightarrow{O{{B}_{1}}}  \\
		\overrightarrow{O{{C}_{1}}}  \\
	\end{matrix} \right).\]
\end{theorem}

\begin{proof}
	Obviously,
	\[\overrightarrow{O{{A}_{2}}}=\alpha _{{{A}_{1}}}^{{{A}_{2}}}\overrightarrow{O{{A}_{1}}}+\alpha _{{{B}_{1}}}^{{{A}_{2}}}\overrightarrow{O{{B}_{1}}}+\alpha _{{{C}_{1}}}^{{{A}_{2}}}\overrightarrow{O{{C}_{1}}},\]
	\[\alpha _{{{A}_{1}}}^{{{A}_{2}}}+\alpha _{{{B}_{1}}}^{{{A}_{2}}}+\alpha _{{{C}_{1}}}^{{{A}_{2}}}=1;\]
	\[\overrightarrow{O{{B}_{2}}}=\alpha _{{{A}_{1}}}^{{{B}_{2}}}\overrightarrow{O{{A}_{1}}}+\alpha _{{{B}_{1}}}^{{{B}_{2}}}\overrightarrow{O{{B}_{1}}}+\alpha _{{{C}_{1}}}^{{{B}_{2}}}\overrightarrow{O{{C}_{1}}},\]
	\[\alpha _{{{A}_{1}}}^{{{B}_{2}}}+\alpha _{{{B}_{1}}}^{{{B}_{2}}}+\alpha _{{{C}_{1}}}^{{{B}_{2}}}=1;\]
	\[\overrightarrow{O{{C}_{2}}}=\alpha _{{{A}_{1}}}^{{{C}_{2}}}\overrightarrow{O{{A}_{1}}}+\alpha _{{{B}_{1}}}^{{{C}_{2}}}\overrightarrow{O{{B}_{1}}}+\alpha _{{{C}_{1}}}^{{{C}_{2}}}\overrightarrow{O{{C}_{1}}},\]
	\[\alpha _{{{A}_{1}}}^{{{C}_{2}}}+\alpha _{{{B}_{1}}}^{{{C}_{2}}}+\alpha _{{{C}_{1}}}^{{{C}_{2}}}=1.\]
	
	Writing it in matrix form is
	\[\left( \begin{matrix}
		\overrightarrow{O{{A}_{2}}}  \\
		\overrightarrow{O{{B}_{2}}}  \\
		\overrightarrow{O{{C}_{2}}}  \\
	\end{matrix} \right)=\left( \begin{matrix}
		\alpha _{{{A}_{1}}}^{{{A}_{2}}} & \alpha _{{{B}_{1}}}^{{{A}_{2}}} & \alpha _{{{C}_{1}}}^{{{A}_{2}}}  \\
		\alpha _{{{A}_{1}}}^{{{B}_{2}}} & \alpha _{{{B}_{1}}}^{{{B}_{2}}} & \alpha _{{{C}_{1}}}^{{{B}_{2}}}  \\
		\alpha _{{{A}_{1}}}^{{{C}_{2}}} & \alpha _{{{B}_{1}}}^{{{C}_{2}}} & \alpha _{{{C}_{1}}}^{{{C}_{2}}}  \\
	\end{matrix} \right)\left( \begin{matrix}
		\overrightarrow{O{{A}_{1}}}  \\
		\overrightarrow{O{{B}_{1}}}  \\
		\overrightarrow{O{{C}_{1}}}  \\
	\end{matrix} \right).\]
\end{proof}
\hfill $\square$\par

The above equation can be written as
\[{{\bm{F}}_{2}}={{\bm{A}}_{21}}{{\bm{F}}_{1}}.\]

Where:
\[{{\bm{F}}_{1}}=\left( \begin{matrix}
	\overrightarrow{O{{A}_{1}}}  \\
	\overrightarrow{O{{B}_{1}}}  \\
	\overrightarrow{O{{C}_{1}}}  \\
\end{matrix} \right)\text{,}\quad {{\bm{F}}_{2}}=\left( \begin{matrix}
	\overrightarrow{O{{A}_{2}}}  \\
	\overrightarrow{O{{B}_{2}}}  \\
	\overrightarrow{O{{C}_{2}}}  \\
\end{matrix} \right)\text{,}\]
\[{{\bm{A}}_{21}}=\left( \begin{matrix}
	\alpha _{{{A}_{1}}}^{{{A}_{2}}} & \alpha _{{{B}_{1}}}^{{{A}_{2}}} & \alpha _{{{C}_{1}}}^{{{A}_{2}}}  \\
	\alpha _{{{A}_{1}}}^{{{B}_{2}}} & \alpha _{{{B}_{1}}}^{{{B}_{2}}} & \alpha _{{{C}_{1}}}^{{{B}_{2}}}  \\
	\alpha _{{{A}_{1}}}^{{{C}_{2}}} & \alpha _{{{B}_{1}}}^{{{C}_{2}}} & \alpha _{{{C}_{1}}}^{{{C}_{2}}}  \\
\end{matrix} \right).\]

Among them, ${{\bm{F}}_{1}}$ and ${{\bm{F}}_{2}}$ are called \textbf{frame vectors}; and ${{\bm{A}}_{21}}$ is called the \textbf{transformation matrix} from ${{\bm{F}}_{1}}$ to ${{\bm{F}}_{2}}$, it is abbreviated as a transformation matrix.


\begin{theorem}{Frame transformation formula 2, Daiyuan Zhang}{BiaojiabianhuanGongshi2}\label{BiaojiabianhuanGongshi2}
	Let $n$ triangles $\triangle {{A}_{1}}{{B}_{1}}{{C}_{1}}$, $\triangle {{A}_{2}}{{B}_{2}}{{C}_{2}}$, …, $\triangle {{A}_{n}}{{B}_{n}}{{C}_{n}}$ be on the same plane. The transformation matrix from ${{\bm{F}}_{1}}$ to ${{\bm{F}}_{2}}$ is ${{\bm{A}}_{21}}$, the transformation matrix from ${{\bm{F}}_{2}}$ to ${{\bm{F}}_{3}}$ is ${{\bm{A}}_{32}}$,..., and the transformation matrix from ${{\bm{F}}_{n-1}}$ to ${{\bm{F}}_{n}}$ is ${{\bm{A}}_{n\left( n-1 \right)}}$. Then:
	\[{{\bm{F}}_{n}}={{\bm{A}}_{n\left( n-1 \right)}}\cdots {{\bm{A}}_{21}}{{\bm{F}}_{1}}.\]
\end{theorem}


\begin{proof}
	According to theorem \ref{thm:BiaojiabianhuanGongshi1}, it can be obtained that:
	\[{{\bm{F}}_{n}}={{\bm{A}}_{n\left( n-1 \right)}}{{\bm{F}}_{n-1}}={{\bm{A}}_{n\left( n-1 \right)}}{{\bm{A}}_{\left( n-1 \right)\left( n-2 \right)}}{{\bm{F}}_{n-2}}=\cdots ={{\bm{A}}_{n\left( n-1 \right)}}\cdots {{\bm{A}}_{21}}{{\bm{F}}_{1}}.\]
\end{proof}
\hfill $\square$\par

%

\begin{theorem}{Transformation of frame component, Daiyuan Zhang}{BiaojiafenliangbianhuanGongshi1}\label{BiaojiafenliangbianhuanGongshi1}
	Let two triangles $\triangle {{A}_{1}}{{B}_{1}}{{C}_{1}}$, $\triangle {{A}_{2}}{{B}_{2}}{{C}_{2}}$ be on the same plane. For the three vertices ${{A}_{2}}$, ${{B}_{2}}$, ${{C}_{2}}$ of $\triangle {{A}_{2}}{{B}_{2}}{{C}_{2}}$, their frame components on $\left( O;{{A}_{1}},{{B}_{1}},{{C}_{1}} \right)$ are $\alpha _{{{A}_{1}}}^{{{A}_{2}}}$, $\alpha _{{{B}_{1}}}^{{{A}_{2}}}$, $\alpha _{{{C}_{1}}}^{{{A}_{2}}}$; $\alpha _{{{A}_{1}}}^{{{B}_{2}}}$, $\alpha _{{{B}_{1}}}^{{{B}_{2}}}$, $\alpha _{{{C}_{1}}}^{{{B}_{2}}}$;  $\alpha _{{{A}_{1}}}^{{{C}_{2}}}$, $\alpha _{{{B}_{1}}}^{{{C}_{2}}}$, $\alpha _{{{C}_{1}}}^{{{C}_{2}}}$; respectively. The frame component of $P$ on $\left( O;{{A}_{1}},{{B}_{1}},{{C}_{1}} \right)$ is $\alpha _{{{A}_{1}}}^{P}$, $\alpha _{{{B}_{1}}}^{P}$, $\alpha _{{{C}_{1}}}^{P}$; the frame component of $P$ on $\left( O;{{A}_{2}},{{B}_{2}},{{C}_{2}} \right)$ is $\alpha _{{{A}_{2}}}^{P}$, $\alpha _{{{B}_{2}}}^{P}$, $\alpha _{{{C}_{2}}}^{P}$. Then
	\[\left( \begin{matrix}
		\alpha _{{{A}_{1}}}^{P}  \\
		\alpha _{{{B}_{1}}}^{P}  \\
		\alpha _{{{C}_{1}}}^{P}  \\
	\end{matrix} \right)=\left( \begin{matrix}
		\alpha _{{{A}_{1}}}^{A_2} & \alpha _{{{A}_{1}}}^{B_2} & \alpha _{{{A}_{1}}}^{C_2}  \\
		\alpha _{{{B}_{1}}}^{A_2} & \alpha _{{{B}_{1}}}^{B_2} & \alpha _{{{B}_{1}}}^{C_2}  \\
		\alpha _{{{C}_{1}}}^{A_2} & \alpha _{{{C}_{1}}}^{B_2} & \alpha _{{{C}_{1}}}^{C_2}  \\
	\end{matrix} \right)\left( \begin{matrix}
		\alpha _{{{A}_{2}}}^{P}  \\
		\alpha _{{{B}_{2}}}^{P}  \\
		\alpha _{{{C}_{2}}}^{P}  \\
	\end{matrix} \right).\]
	
	The sum of elements in each column of the matrix is 1.
\end{theorem}

\begin{proof}
	Given a point $P$ on the plane of $\triangle {{A}_{1}}{{B}_{1}}{{C}_{1}}$ and $\triangle {{A}_{2}}{{B}_{2}}{{C}_{2}}$, for the frame $\left( O;{{A}_{1}},{{B}_{1}},{{C}_{1}} \right)$, it can be obtained that:
	\[\overrightarrow{OP}=\alpha _{{{A}_{1}}}^{P}\overrightarrow{O{{A}_{1}}}+\alpha _{{{B}_{1}}}^{P}\overrightarrow{O{{B}_{1}}}+\alpha _{{{C}_{1}}}^{P}\overrightarrow{O{{C}_{1}}},\]
	\[\alpha _{{{A}_{1}}}^{P}+\alpha _{{{B}_{1}}}^{P}+\alpha _{{{C}_{1}}}^{P}=1.\]
	
	For the frame $\left( O;{{A}_{2}},{{B}_{2}},{{C}_{2}} \right)$, it can be obtained that:
	\[\overrightarrow{OP}=\alpha _{{{A}_{2}}}^{P}\overrightarrow{OA_2}+\alpha _{{{B}_{2}}}^{P}\overrightarrow{OB_2}+\alpha _{{{C}_{2}}}^{P}\overrightarrow{OC_2},\]
	\[\alpha _{{{A}_{2}}}^{P}+\alpha _{{{B}_{2}}}^{P}+\alpha _{{{C}_{2}}}^{P}=1.\]
	
	Therefore
	\begin{align*}
		\overrightarrow{OP}& =\alpha _{{{A}_{1}}}^{P}\overrightarrow{O{{A}_{1}}}+\alpha _{{{B}_{1}}}^{P}\overrightarrow{O{{B}_{1}}}+\alpha _{{{C}_{1}}}^{P}\overrightarrow{O{{C}_{1}}} \\ 
		& =\alpha _{{{A}_{2}}}^{P}\overrightarrow{OA_2}+\alpha _{{{B}_{2}}}^{P}\overrightarrow{OB_2}+\alpha _{{{C}_{2}}}^{P}\overrightarrow{OC_2},  
	\end{align*}
	i.e.
	\[\alpha _{{{A}_{1}}}^{P}\overrightarrow{O{{A}_{1}}}+\alpha _{{{B}_{1}}}^{P}\overrightarrow{O{{B}_{1}}}+\alpha _{{{C}_{1}}}^{P}\overrightarrow{O{{C}_{1}}}-\left( \alpha _{{{A}_{2}}}^{P}\overrightarrow{OA_2}+\alpha _{{{B}_{2}}}^{P}\overrightarrow{OB_2}+\alpha _{{{C}_{2}}}^{P}\overrightarrow{OC_2} \right)=\overrightarrow{0}.\]
	
	And
	\[\overrightarrow{OA_2}=\alpha _{{{A}_{1}}}^{A_2}\overrightarrow{O{{A}_{1}}}+\alpha _{{{B}_{1}}}^{A_2}\overrightarrow{O{{B}_{1}}}+\alpha _{{{C}_{1}}}^{A_2}\overrightarrow{O{{C}_{1}}},\]
	\[\overrightarrow{OB_2}=\alpha _{{{A}_{1}}}^{B_2}\overrightarrow{O{{A}_{1}}}+\alpha _{{{B}_{1}}}^{B_2}\overrightarrow{O{{B}_{1}}}+\alpha _{{{C}_{1}}}^{B_2}\overrightarrow{O{{C}_{1}}},\]
	\[\overrightarrow{OC_2}=\alpha _{{{A}_{1}}}^{C_2}\overrightarrow{O{{A}_{1}}}+\alpha _{{{B}_{1}}}^{C_2}\overrightarrow{O{{B}_{1}}}+\alpha _{{{C}_{1}}}^{C_2}\overrightarrow{O{{C}_{1}}}.\]
	
	Substituting into the above equation yields:
	\begin{align*}
		& \alpha _{{{A}_{1}}}^{P}\overrightarrow{O{{A}_{1}}}+\alpha _{{{B}_{1}}}^{P}\overrightarrow{O{{B}_{1}}}+\alpha _{{{C}_{1}}}^{P}\overrightarrow{O{{C}_{1}}} \\ 
		& -\left( \begin{aligned}
			& \alpha _{{{A}_{2}}}^{P}\left( \alpha _{{{A}_{1}}}^{A_2}\overrightarrow{O{{A}_{1}}}+\alpha _{{{B}_{1}}}^{A_2}\overrightarrow{O{{B}_{1}}}+\alpha _{{{C}_{1}}}^{A_2}\overrightarrow{O{{C}_{1}}} \right) \\ 
			& +\alpha _{{{B}_{2}}}^{P}\left( \alpha _{{{A}_{1}}}^{B_2}\overrightarrow{O{{A}_{1}}}+\alpha _{{{B}_{1}}}^{B_2}\overrightarrow{O{{B}_{1}}}+\alpha _{{{C}_{1}}}^{B_2}\overrightarrow{O{{C}_{1}}} \right) \\ 
			& +\alpha _{{{C}_{2}}}^{P}\left( \alpha _{{{A}_{1}}}^{C_2}\overrightarrow{O{{A}_{1}}}+\alpha _{{{B}_{1}}}^{C_2}\overrightarrow{O{{B}_{1}}}+\alpha _{{{C}_{1}}}^{C_2}\overrightarrow{O{{C}_{1}}} \right) \\ 
		\end{aligned} \right)=\overrightarrow{0}. \\ 
	\end{align*}

	After merging similar items in the above equation, it can be obtained that:
	\begin{align*}
		& \left( \alpha _{{{A}_{1}}}^{P}-\alpha _{{{A}_{2}}}^{P}\alpha _{{{A}_{1}}}^{A_2}-\alpha _{{{B}_{2}}}^{P}\alpha _{{{A}_{1}}}^{B_2}-\alpha _{{{C}_{2}}}^{P}\alpha _{{{A}_{1}}}^{C_2} \right)\overrightarrow{O{{A}_{1}}} \\ 
		& +\left( \alpha _{{{B}_{1}}}^{P}-\alpha _{{{A}_{2}}}^{P}\alpha _{{{B}_{1}}}^{A_2}-\alpha _{{{B}_{2}}}^{P}\alpha _{{{B}_{1}}}^{B_2}-\alpha _{{{C}_{2}}}^{P}\alpha _{{{B}_{1}}}^{C_2} \right)\overrightarrow{O{{B}_{1}}} \\ 
		& +\left( \alpha _{{{C}_{1}}}^{P}-\alpha _{{{A}_{2}}}^{P}\alpha _{{{C}_{1}}}^{A_2}-\alpha _{{{B}_{2}}}^{P}\alpha _{{{C}_{1}}}^{B_2}-\alpha _{{{C}_{2}}}^{P}\alpha _{{{C}_{1}}}^{C_2} \right)\overrightarrow{O{{C}_{1}}}=\overrightarrow{0}. \\ 
	\end{align*}
	
		
	It is easy to see that the sum of all coefficients in the above equation is 0.
	By theorem \ref{thm:Thm3.4.3}, it can be obtained that:  
	\[\left\{ \begin{aligned}
		& \alpha _{{{A}_{1}}}^{P}-\alpha _{{{A}_{2}}}^{P}\alpha _{{{A}_{1}}}^{A_2}-\alpha _{{{B}_{2}}}^{P}\alpha _{{{A}_{1}}}^{B_2}-\alpha _{{{C}_{2}}}^{P}\alpha _{{{A}_{1}}}^{C_2}=0 \\ 
		& \alpha _{{{B}_{1}}}^{P}-\alpha _{{{A}_{2}}}^{P}\alpha _{{{B}_{1}}}^{A_2}-\alpha _{{{B}_{2}}}^{P}\alpha _{{{B}_{1}}}^{B_2}-\alpha _{{{C}_{2}}}^{P}\alpha _{{{B}_{1}}}^{C_2}=0 \\ 
		& \alpha _{{{C}_{1}}}^{P}-\alpha _{{{A}_{2}}}^{P}\alpha _{{{C}_{1}}}^{A_2}-\alpha _{{{B}_{2}}}^{P}\alpha _{{{C}_{1}}}^{B_2}-\alpha _{{{C}_{2}}}^{P}\alpha _{{{C}_{1}}}^{C_2}=0. \\ 
	\end{aligned} \right.\]
	
	Or written in matrix form:
	\[\left( \begin{matrix}
		\alpha _{{{A}_{1}}}^{P}  \\
		\alpha _{{{B}_{1}}}^{P}  \\
		\alpha _{{{C}_{1}}}^{P}  \\
	\end{matrix} \right)=\left( \begin{matrix}
		\alpha _{{{A}_{1}}}^{A_2} & \alpha _{{{A}_{1}}}^{B_2} & \alpha _{{{A}_{1}}}^{C_2}  \\
		\alpha _{{{B}_{1}}}^{A_2} & \alpha _{{{B}_{1}}}^{B_2} & \alpha _{{{B}_{1}}}^{C_2}  \\
		\alpha _{{{C}_{1}}}^{A_2} & \alpha _{{{C}_{1}}}^{B_2} & \alpha _{{{C}_{1}}}^{C_2}  \\
	\end{matrix} \right)\left( \begin{matrix}
		\alpha _{{{A}_{2}}}^{P}  \\
		\alpha _{{{B}_{2}}}^{P}  \\
		\alpha _{{{C}_{2}}}^{P}  \\
	\end{matrix} \right).\]
\end{proof}
\hfill $\square$\par

The above theorem is very important.

Imagine $\triangle {{A}_{1}}{{B}_{1}}{{C}_{1}}$ as a reference triangle or an old triangle, and consider $\triangle {{A}_{2}}{{B}_{2}}{{C}_{2}}$ as a new triangle. The above equation indicates that if the frame components of a point $P$ on a new triangle are known, and the frame components of the three vertices of the new triangle on the old triangle are known, then the frame components of that point $P$ on the old triangle can be calculated.


\begin{example}{}\label{GeidingYigeSanjiaoxing}
	Given a triangle (referred to as the original triangle), the incenter, orthocenter and circumcenter of the triangle form a new triangle. Find the frame component of the centroid of the new triangle under the  frame of original triangle.
\end{example}

\begin{solution}
	Assuming the original triangle is $\triangle {{A}_{1}}{{B}_{1}}{{C}_{1}}$, its incenter is $I$, its orthocenter is $H$, and its circumcenter is $Q$. $\triangle IHQ$ is the new triangle, ${{G}_{2}}$ is the centroid of $\triangle IHQ$. So $\alpha _{I}^{{{G}_{2}}}=\alpha _{H}^{{{G}_{2}}}=\alpha _{Q}^{{{G}_{2}}}={1}/{3}\;$, According to the above equation, it can be obtained that
	\[\begin{aligned}
		\left( \begin{matrix}
			\alpha _{{{A}_{1}}}^{{{G}_{2}}}  \\
			\alpha _{{{B}_{1}}}^{{{G}_{2}}}  \\
			\alpha _{{{C}_{1}}}^{{{G}_{2}}}  \\
		\end{matrix} \right)& =\left( \begin{matrix}
			\alpha _{{{A}_{1}}}^{I} & \alpha _{{{A}_{1}}}^{H} & \alpha _{{{A}_{1}}}^{Q}  \\
			\alpha _{{{B}_{1}}}^{I} & \alpha _{{{B}_{1}}}^{H} & \alpha _{{{B}_{1}}}^{Q}  \\
			\alpha _{{{C}_{1}}}^{I} & \alpha _{{{C}_{1}}}^{H} & \alpha _{{{C}_{1}}}^{Q}  \\
		\end{matrix} \right)\left( \begin{matrix}
			\alpha _{I}^{{{G}_{2}}}  \\
			\alpha _{H}^{{{G}_{2}}}  \\
			\alpha _{Q}^{{{G}_{2}}}  \\
		\end{matrix} \right) \\ 
		& =\left( \begin{matrix}
			\frac{{{a}_{1}}}{{{a}_{1}}+{{b}_{1}}+{{c}_{1}}} & \frac{\tan {{A}_{1}}}{\tan {{A}_{1}}+\tan {{B}_{1}}+\tan {{C}_{1}}} & \frac{\sin 2{{A}_{1}}}{\sin 2{{A}_{1}}+\sin 2{{B}_{1}}+\sin 2{{C}_{1}}}  \\
			\frac{{{b}_{1}}}{{{a}_{1}}+{{b}_{1}}+{{c}_{1}}} & \frac{\tan {{B}_{1}}}{\tan {{A}_{1}}+\tan {{B}_{1}}+\tan {{C}_{1}}} & \frac{\sin 2{{B}_{1}}}{\sin 2{{A}_{1}}+\sin 2{{B}_{1}}+\sin 2{{C}_{1}}}  \\
			\frac{{{c}_{1}}}{{{a}_{1}}+{{b}_{1}}+{{c}_{1}}} & \frac{\tan {{C}_{1}}}{\tan {{A}_{1}}+\tan {{B}_{1}}+\tan {{C}_{1}}} & \frac{\sin 2{{C}_{1}}}{\sin 2{{A}_{1}}+\sin 2{{B}_{1}}+\sin 2{{C}_{1}}}  \\
		\end{matrix} \right)\left( \begin{matrix}
			\frac{1}{3}  \\
			\frac{1}{3}  \\
			\frac{1}{3}  \\
		\end{matrix} \right).  
	\end{aligned}\]
	
	Therefore
	\[\alpha _{{{A}_{1}}}^{{{G}_{2}}}=\frac{1}{3}\left( \frac{{{a}_{1}}}{{{a}_{1}}+{{b}_{1}}+{{c}_{1}}}+\frac{\tan {{A}_{1}}}{\tan {{A}_{1}}+\tan {{B}_{1}}+\tan {{C}_{1}}}+\frac{\sin 2{{A}_{1}}}{\sin 2{{A}_{1}}+\sin 2{{B}_{1}}+\sin 2{{C}_{1}}} \right),\]
	\[\alpha _{{{B}_{1}}}^{{{G}_{2}}}=\frac{1}{3}\left( \frac{{{b}_{1}}}{{{a}_{1}}+{{b}_{1}}+{{c}_{1}}}+\frac{\tan {{B}_{1}}}{\tan {{A}_{1}}+\tan {{B}_{1}}+\tan {{C}_{1}}}+\frac{\sin 2{{B}_{1}}}{\sin 2{{A}_{1}}+\sin 2{{B}_{1}}+\sin 2{{C}_{1}}} \right),\]
	\[\alpha _{{{C}_{1}}}^{{{G}_{2}}}=\frac{1}{3}\left( \frac{{{c}_{1}}}{{{a}_{1}}+{{b}_{1}}+{{c}_{1}}}+\frac{\tan {{C}_{1}}}{\tan {{A}_{1}}+\tan {{B}_{1}}+\tan {{C}_{1}}}+\frac{\sin 2{{C}_{1}}}{\sin 2{{A}_{1}}+\sin 2{{B}_{1}}+\sin 2{{C}_{1}}} \right).\]
\end{solution}
\hfill $\diamond$\par

\section{The centers of a triangle}\label{SanjiaoxingDeXin}

Many readers know some of the centers of a given triangle, such as the centroid, incenter, orthocenter, and so on. But are these centers just names for convenience or do they imply profound inner meanings? In other words, do these centers share any common characteristics? The answer is affirmative. In this section, I first introduced the concept of fixed point and used the method of Intercenter Geometry to classify the centers of triangle into several categories. These centers of the triangle are essentially fixed points under transformations. From the perspective of the fixed point, the centroid, incenter, orthocenter and circumcenter of the triangle are all of the same type of center, while the Brocard point is a different type of center. The intrinsic relationship between the centers of several types of triangle is given, and it is proved that there are infinitely many centers of each type of triangle. The center of the triangle cannot fill the entire plane, which means there are points that are not the center of the triangle. I refer to these points as the "non-center" of the triangle and prove that there are infinitely many non-centers of the triangle. I hope the conclusions derived in this section can have a positive promoting effect on the study of triangle geometry.

According to theorem \ref{thm:Thm6.1.1}, if $\triangle ABC$ is given, point $P$ is a point on the plane of $\triangle ABC$, and point $O$ is any given point in space, then:
\[\overrightarrow{OP}=\alpha _{A}^{P}\overrightarrow{OA}+\alpha _{B}^{P}\overrightarrow{OB}+\alpha _{C}^{P}\overrightarrow{OC},\]
\[\alpha _{A}^{P}+\alpha _{B}^{P}+\alpha _{C}^{P}=1.\]

According to theorem \ref{thm:PingmianShangDeDianYuBiaojiafenliangYiyiduiying}, point $P$ corresponds one-to-one with the frame component $\left( \alpha _{A}^{P},\alpha _{B}^{P},\alpha _{C}^{P} \right)$. That is to say, as long as the triangle remains unchanged and the frame components remain unchanged, point $P$ will not change, or in other words, the relative position of point $P$ and $\triangle ABC$ will remain unchanged.

The above equation also indicates that the vector $\overrightarrow{OP}$ can be viewed as a function of three vertices, i.e.
\[\overrightarrow{OP}=\overrightarrow{OP}\left( A,B,C \right).\] 

The frame component is a function of the side lengths, and can also be seen as a function of three vertices. Taking $\alpha _{A}^{P}$ as an example:
\[\alpha _{A}^{P}=\alpha _{A}^{P}\left( a,b,c \right)=\alpha _{A}^{P}\left( BC,CA,AB \right).\]

When $A$ and $B$ are exchanged, the right-hand side of the above equation becomes
\[\alpha _{A}^{P}\left( AC,CB,BA \right)=\alpha _{A}^{P}\left( CA,BC,AB \right)=\alpha _{A}^{P}\left( b,a,c \right).\]

So, for the frame component, the exchange of $A$ and $B$ is equivalent to the exchange of $a$ and $b$. Similarly, the exchange between $B$ and $C$ is equivalent to the exchange between $b$ and $c$, while the exchange between $C$ and $A$ is equivalent to the exchange between $c $ and $a $. On the contrary, exchanging $a$ and $b $ is equivalent to exchanging $A $ and $B $. Similarly, the exchange between $b $ and $c $ is equivalent to the exchange between $B $ and $C $, while the exchange between $c $ and $a $is equivalent to the exchange between $C $ and $A $.

Now let's take an example. Let:
\[\overrightarrow{OP}=\alpha _{A}^{P}\overrightarrow{OA}+\alpha _{B}^{P}\overrightarrow{OB}+\alpha _{C}^{P}\overrightarrow{OC},\]
\[\alpha _{A}^{P}+\alpha _{B}^{P}+\alpha _{C}^{P}=1.\]

Where:
\[\alpha _{A}^{P}=\frac{\tan A}{\tan A+\tan B+\tan C},\]
\[\alpha _{B}^{P}=\frac{\tan B}{\tan A+\tan B+\tan C},\]	
\[\,\alpha _{C}^{P}=\frac{\tan C}{\tan A+\tan B+\tan C}.\]

Therefore%
\[\overrightarrow{OP}\left( B,A,C \right)=\overrightarrow{OP}\left( A,B,C \right)=\overrightarrow{OP}.\] 

That is to say, after swapping $A $ and $B $, the vector $\overrightarrow{OP}$ does not change, or the position of point $P $ remains unchanged, that is, point $P $ is a fixed point under the transformation of $A $ and $B $ swapping. In fact, the above point $P $ is the orthocenter of $\triangle ABC$. Below is the definition of fixed points for a triangle.


\begin{definition}{Fixed point of triangle}{SanjiaoxingDeBudongdian}\label{SanjiaoxingDeBudongdian}
	If $\triangle ABC$ is given, point $P $ is a point on the plane of $\triangle ABC$, point $O $ is any given point in space, and:
	\[\overrightarrow{OP}=\alpha _{A}^{P}\overrightarrow{OA}+\alpha _{B}^{P}\overrightarrow{OB}+\alpha _{C}^{P}\overrightarrow{OC},\]
	\[\alpha _{A}^{P}+\alpha _{B}^{P}+\alpha _{C}^{P}=1.\]
	If there is a transformation $\bm{T}$ that:
	\[\bm{T}:\ \,\,\overrightarrow{OP}\to \overrightarrow{OP},\]
	The point $P$ is called the fixed point of $\triangle ABC$ under the transformation $\bm{T}$, abbreviated as the fixed point.
\end{definition}

The transformation in the above definition is not an identity transformation. Here are some examples to illustrate the meaning of the above definition.

In the first scenario, the frame changes and the frame components remain unchanged. For example, suppose the definition of the transformation $\bm{T}$ is as follows: changing the frame $\overrightarrow{OA}$ to $\overrightarrow{OB}$, changing the frame $\overrightarrow{OB}$ to $\overrightarrow{OA}$, keeping the frame $\overrightarrow{OC}$ unchanged, and keeping the frame components unchanged. Under the transformation $\bm{T}$, the above equation is written as:
\[\overrightarrow{O{{P}_{1}}}=\alpha _{A}^{P}\overrightarrow{OB}+\alpha _{B}^{P}\overrightarrow{OA}+\alpha _{C}^{P}\overrightarrow{OC},\]
\[\alpha _{A}^{P}+\alpha _{B}^{P}+\alpha _{C}^{P}=1.\]

Obviously, generally speaking, $\overrightarrow{O{{P}_{1}}}\ne \overrightarrow{OP}$. Or rather, at this point, the transformation $\bm{T}$ changes point $P$ to ${{P}_{1}}$. Or rather, point $P$ is not a fixed point for the transformation $\bm{T}$.

In the second scenario, the frame component changes while the frame remains unchanged. For example, if the definition of transformation $\bm{T}$ is as follows: the frame components $\alpha _{A}^{P}$ and $\alpha _{B}^{P}$ are exchanged, and the rest remain unchanged, then under transformation $\bm{T}$, the following results are obtained:
\[\overrightarrow{O{{P}_{2}}}=\alpha _{B}^{P}\overrightarrow{OA}+\alpha _{A}^{P}\overrightarrow{OB}+\alpha _{C}^{P}\overrightarrow{OC},\]
\[\alpha _{B}^{P}+\alpha _{A}^{P}+\alpha _{C}^{P}=1.\]

Obviously, generally speaking, $\overrightarrow{O{{P}_{2}}}\ne \overrightarrow{OP}$. This indicates that point $P $ is not a fixed point for transformation $\bm{T}$.

Changing the frame or changing the frame components can be seen as a transformation. So, is there any transformation that makes some of the transformed points fixed points? The answer is affirmative. There are two elements here, one is transformation and the other is point. Different transformations may turn different points into fixed points. I believe that the center of a triangle is a fixed point under a certain type of transformation by carefully comparing the frame component formulas of its incenter, orthocenter, circumcenter. I divide the centers of triangles into several categories: local symmetry centers, global symmetry centers, and rotation centers. Let's study this issue below.

%
%
%
%
%
%
%

\begin{definition}{The local symmetry center of a triangle}{SanjiaoxingDeJubuduichenxin}\label{SanjiaoxingDeJubuduichenxin}
	If $\triangle ABC$ is given, point $P $ is a point on the plane of $\triangle ABC$, point $O $ is any given point in space, and:
	\[\overrightarrow{OP}=\alpha _{A}^{P}\overrightarrow{OA}+\alpha _{B}^{P}\overrightarrow{OB}+\alpha _{C}^{P}\overrightarrow{OC},\]
	\[\alpha _{A}^{P}+\alpha _{B}^{P}+\alpha _{C}^{P}=1.\]
	
	1. If
	\[\underline{{{\overrightarrow{OP}}_{AB}}}=\overrightarrow{OP},\]
	then the point $P $ is called the center of $\triangle ABC$ under the transformation $\underline{{{\left( {\quad} \right)}_{AB}}}$, or the local symmetry center of $AB $. The meaning of $\underline{{{\left( {\quad} \right)}_{AB}}}$ is to swap $A $ and $B $ in $\left( {\quad} \right)$.
	
	2. If
	\[\underline{{{\overrightarrow{OP}}_{BC}}}=\overrightarrow{OP},\]
	then the  point $P $ is called the center of$\triangle ABC$ under the transformation $\underline{{{\left( {\quad} \right)}_{BC}}}$, or the local symmetry center of $BC $. The meaning of $\underline{{{\left( {\quad} \right)}_{BC}}}$ is to swap $B $ and $C $ in $\left( {\quad} \right)$.
	
	3. If
	\[\underline{{{\overrightarrow{OP}}_{CA}}}=\overrightarrow{OP},\]
	then the  point $P $ is called the center of $\triangle ABC$ under the transformation $\underline{{{\left( {\quad} \right)}_{CA}}}$, or the local symmetry center of $CA $. The meaning of $\underline{{{\left( {\quad} \right)}_{CA}}}$ is to swap $C $ and $A $ in $\left( {\quad} \right)$.

	The collective notation of the $AB $ local symmetry center is $\mathbf{L}{{\mathbf{S}}_{AB}}$, the collective notation of the $BC $ local symmetry center is $\mathbf{L}{{\mathbf{S}}_{BC}}$, and the collective notation of the $CA $ local symmetry center is $\mathbf{L}{{\mathbf{S}}_{CA}}$. The set composed of all three types of local symmetry centers mentioned above is denoted as $\mathbf{LS}$, and $\mathbf{LS}$ is called the complete set of local symmetry centers, abbreviated as local symmetry centers.
	\[\mathbf{LS}:=\mathbf{L}{{\mathbf{S}}_{AB}}\cup \mathbf{L}{{\mathbf{S}}_{BC}}\cup \mathbf{L}{{\mathbf{S}}_{CA}}.\]
\end{definition}

In the above definition, the details and correlation of transformations such as $\underline{{{\left( {\quad} \right)}_{AB}}}$ can be found in Appendix \ref{DuihuanLunhuan}. The local center of symmetry is a fixed point under the corresponding  transformation. 

%
%

\begin{definition}{The global center of symmetry of a triangle}{SanjiaoxingDeQuanjuduichenxin}\label{SanjiaoxingDeQuanjuduichenxin}
	If $\triangle ABC$ is given, point $P $ is a point on the plane of $\triangle ABC$, point $O$ is any given point in space, and:
	\[\overrightarrow{OP}=\alpha _{A}^{P}\overrightarrow{OA}+\alpha _{B}^{P}\overrightarrow{OB}+\alpha _{C}^{P}\overrightarrow{OC},\]
	\[\alpha _{A}^{P}+\alpha _{B}^{P}+\alpha _{C}^{P}=1.\]
	
	If point $P$ is both the local symmetry center of $AB$, the local symmetry center of $BC$, and the local symmetry center of $CA$, i.e
	\[\underline{{{\overrightarrow{OP}}_{AB}}}=\overrightarrow{OP},\quad \underline{{{\overrightarrow{OP}}_{BC}}}=\overrightarrow{OP},\quad \underline{{{\overrightarrow{OP}}_{CA}}}=\overrightarrow{OP},\]
	The point $P$ is called the global center of symmetry of $\triangle ABC$. The set composed of all global symmetry centers is denoted as $\mathbf{GS}$, i.e.
	\[\mathbf{GS}:=\mathbf{L}{{\mathbf{S}}_{AB}}\cap \mathbf{L}{{\mathbf{S}}_{BC}}\cap \mathbf{L}{{\mathbf{S}}_{CA}}.\]
\end{definition}

In the above definition, the details and correlation of transformations such as $\underline{\left( {\quad} \right)}$ can be found in Appendix \ref{DuihuanLunhuan}. The global center of symmetry is a fixed point under the $AB$ transformation, $BC$ transformation, and $CA$ transformation.


\begin{definition}{The symmetry center of a triangle}{SanjiaoxingDeduichenxin}\label{SanjiaoxingDeduichenxin}
	If point $P $ is the local or global center of symmetry of $\triangle ABC$, then point $P $ is called the center of symmetry of $\triangle ABC$, abbreviated as the center of symmetry. The set composed of all symmetry centers is denoted as $\mathbf{CS}$. 
	\[\mathbf{CS}:=\mathbf{LS}\cup \mathbf{GS}.\]
\end{definition}

The center of symmetry is for both local and global centers of symmetry.

%

\begin{definition}{The rotation center of a triangle}{SanjiaoxingDeLunhuanxin}\label{SanjiaoxingDeLunhuanxin}
	If $\triangle ABC$ is given, point $P$ is a point on the plane of $\triangle ABC$, point $O$ is any given point in space, and:
	\[\overrightarrow{OP}=\alpha _{A}^{P}\overrightarrow{OA}+\alpha _{B}^{P}\overrightarrow{OB}+\alpha _{C}^{P}\overrightarrow{OC},\]
	\[\alpha _{A}^{P}+\alpha _{B}^{P}+\alpha _{C}^{P}=1.\]
	
	If
	\[\underline{\overrightarrow{OP}}=\overrightarrow{OP},\]
	then the point $P$ is called the rotation center of the $\triangle ABC$. The set of all rotation centers is denoted as  $\mathbf{RC}$.
\end{definition}

Rotation center is a fixed point under rotation transformation. The concept and properties of rotation transformation can be found in Appendix \ref{DuihuanLunhuan}.


\begin{definition}{The centers of a triangle}{SanjiaoxingDeNaxieXin}\label{SanjiaoxingDeNaxieXin}
	If point $P$ is the center of symmetry or the center of rotation of a triangle, then point $P$ is called the center of the $\triangle ABC$, abbreviated as the center of a triangle. The set of all the centers of this triangle is denoted as $\mathbf{CT}$. 
	\[\mathbf{CT}:=\mathbf{CS}\cup \mathbf{RS}.\]
\end{definition}

Below are the necessary and sufficient conditions for determining the types of centers in a triangle.

%
%
%
%

\begin{theorem}{Sufficient and necessary conditions for local symmetry centers, Daiyuan Zhang}{PandingSanjiaoxingJubuduichenxinDeChongfenbiyaotiaojian}\label{PandingSanjiaoxingJubuduichenxinDeChongfenbiyaotiaojian}
	If $\triangle ABC$ is given, point $P$ is a point on the plane of $\triangle ABC$, point $O$ is any given point in space, and:
	\[\overrightarrow{OP}=\alpha _{A}^{P}\overrightarrow{OA}+\alpha _{B}^{P}\overrightarrow{OB}+\alpha _{C}^{P}\overrightarrow{OC},\]
	\[\alpha _{A}^{P}+\alpha _{B}^{P}+\alpha _{C}^{P}=1.\]
	
	Then:
	
	1. The necessary and sufficient condition for the point $P$ is the $AB$ local symmetry center of $\triangle ABC$ is:
	\[\underline{{{\left( \alpha _{A}^{P} \right)}_{AB}}}=\alpha _{B}^{P},\quad \underline{{{\left( \alpha _{B}^{P} \right)}_{AB}}}=\alpha _{A}^{P},\quad \underline{{{\left( \alpha _{C}^{P} \right)}_{AB}}}=\alpha _{C}^{P}.\]
	
	2. The necessary and sufficient condition for the point $P$ is the $BC$ local symmetry center of $\triangle ABC$ is:
	\[\underline{{{\left( \alpha _{A}^{P} \right)}_{BC}}}=\alpha _{A}^{P},\quad \underline{{{\left( \alpha _{B}^{P} \right)}_{BC}}}=\alpha _{C}^{P},\quad \underline{{{\left( \alpha _{C}^{P} \right)}_{BC}}}=\alpha _{B}^{P};\]
	
	3. The necessary and sufficient condition for the point $P$ is the $CA$ local symmetry center of $\triangle ABC$ is:
	\[\underline{{{\left( \alpha _{A}^{P} \right)}_{CA}}}=\alpha _{C}^{P},\quad \underline{{{\left( \alpha _{B}^{P} \right)}_{CA}}}=\alpha _{B}^{P},\quad \underline{{{\left( \alpha _{C}^{P} \right)}_{CA}}}=\alpha _{A}^{P}.\]
\end{theorem}

\begin{proof}
	Only prove the first scenario, the rest are similar.
	
	Necessity. According to theorem \ref{thm:DuihuanDeDaishuYunsuanXingzhi} in the appendix, it can be concluded that:
	\begin{align*}
		\underline{{{\overrightarrow{OP}}_{AB}}}& =\underline{{{\left( \alpha _{A}^{P}\overrightarrow{OA}+\alpha _{B}^{P}\overrightarrow{OB}+\alpha _{C}^{P}\overrightarrow{OC} \right)}_{AB}}} \\ 
		& =\underline{{{\left( \alpha _{A}^{P}\overrightarrow{OA} \right)}_{AB}}}+\underline{{{\left( \alpha _{B}^{P}\overrightarrow{OB} \right)}_{AB}}}+\underline{{{\left( \alpha _{C}^{P}\overrightarrow{OC} \right)}_{AB}}} \\ 
		& =\underline{{{\left( \alpha _{A}^{P} \right)}_{AB}}}\cdot \underline{{{\left( \overrightarrow{OA} \right)}_{AB}}}+\underline{{{\left( \alpha _{B}^{P} \right)}_{AB}}}\cdot \underline{{{\left( \overrightarrow{OB} \right)}_{AB}}}+\underline{{{\left( \alpha _{C}^{P} \right)}_{AB}}}\cdot \underline{{{\left( \overrightarrow{OC} \right)}_{AB}}} \\ 
		& =\underline{{{\left( \alpha _{A}^{P} \right)}_{AB}}}\overrightarrow{OB}+\underline{{{\left( \alpha _{B}^{P} \right)}_{AB}}}\overrightarrow{OA}+\underline{{{\left( \alpha _{C}^{P} \right)}_{AB}}}\overrightarrow{OC} \\ 
		& =\underline{{{\left( \alpha _{B}^{P} \right)}_{AB}}}\overrightarrow{OA}+\underline{{{\left( \alpha _{A}^{P} \right)}_{AB}}}\overrightarrow{OB}+\underline{{{\left( \alpha _{C}^{P} \right)}_{AB}}}\overrightarrow{OC},  
	\end{align*}
	\[\overrightarrow{OP}=\alpha _{A}^{P}\overrightarrow{OA}+\alpha _{B}^{P}\overrightarrow{OB}+\alpha _{C}^{P}\overrightarrow{OC},\]
	
	Assuming $\underline{{{\overrightarrow{OP}}_{AB}}}=\overrightarrow{OP}$, then $\overrightarrow{OP}-\underline{{{\overrightarrow{OP}}_{AB}}}=\overrightarrow{0}$, that is
	\[\left( \alpha _{A}^{P}-\underline{{{\left( \alpha _{B}^{P} \right)}_{AB}}} \right)\overrightarrow{OA}+\left( \alpha _{B}^{P}-\underline{{{\left( \alpha _{A}^{P} \right)}_{AB}}} \right)\overrightarrow{OB}+\left( \alpha _{C}^{P}-\underline{{{\left( \alpha _{C}^{P} \right)}_{AB}}} \right)\overrightarrow{OC}=\overrightarrow{0}.\]
	
	And
	\[\underline{{{\left( \alpha _{B}^{P} \right)}_{AB}}}+\underline{{{\left( \alpha _{A}^{P} \right)}_{AB}}}+\underline{{{\left( \alpha _{C}^{P} \right)}_{AB}}}=\underline{{{\left( \alpha _{B}^{P}+\alpha _{A}^{P}+\alpha _{C}^{P} \right)}_{AB}}}=\underline{{{\left( 1 \right)}_{AB}}}=1.\]
	
	Therefore
	\begin{align*}
		& \left( \alpha _{A}^{P}-\underline{{{\left( \alpha _{B}^{P} \right)}_{AB}}} \right)+\left( \alpha _{B}^{P}-\underline{{{\left( \alpha _{A}^{P} \right)}_{AB}}} \right)+\left( \alpha _{C}^{P}-\underline{{{\left( \alpha _{C}^{P} \right)}_{AB}}} \right) \\ 
		& =\left( \alpha _{A}^{P}+\alpha _{B}^{P}+\alpha _{C}^{P} \right)-\left( \underline{{{\left( \alpha _{B}^{P} \right)}_{AB}}}+\underline{{{\left( \alpha _{A}^{P} \right)}_{AB}}}+\underline{{{\left( \alpha _{C}^{P} \right)}_{AB}}} \right) \\ 
		& =1-1 \\ 
		& =0. 
	\end{align*}
	
	Using the triangle frame theorem (theorem \ref{thm:Thm3.4.3}), we obtain:
	\[\underline{{{\left( \alpha _{B}^{P} \right)}_{AB}}}=\alpha _{A}^{P},\quad \underline{{{\left( \alpha _{A}^{P} \right)}_{AB}}}=\alpha _{B}^{P},\quad \underline{{{\left( \alpha _{C}^{P} \right)}_{AB}}}=\alpha _{C}^{P}.\]
	
	Sufficiency. Assuming the following conditions hold:
	\[\underline{{{\left( \alpha _{A}^{P} \right)}_{AB}}}=\alpha _{B}^{P},\quad \underline{{{\left( \alpha _{B}^{P} \right)}_{AB}}}=\alpha _{A}^{P},\quad \underline{{{\left( \alpha _{C}^{P} \right)}_{AB}}}=\alpha _{C}^{P}.\]
	
	Based on the previous results, it can be concluded that:
	\begin{align*}
		& \underline{{{\overrightarrow{OP}}_{AB}}}=\underline{{{\left( \alpha _{B}^{P} \right)}_{AB}}}\overrightarrow{OA}+\underline{{{\left( \alpha _{A}^{P} \right)}_{AB}}}\overrightarrow{OB}+\underline{{{\left( \alpha _{C}^{P} \right)}_{AB}}}\overrightarrow{OC} \\ 
		& =\alpha _{A}^{P}\overrightarrow{OA}+\alpha _{B}^{P}\overrightarrow{OB}+\alpha _{C}^{P}\overrightarrow{OC} \\ 
		& =\overrightarrow{OP}.  
	\end{align*}
\end{proof}
\hfill $\square$\par

%

\begin{theorem}{Sufficient and necessary conditions for global symmetry centers, Daiyuan Zhang}{PanduanSanjiaoxingQuanjuduichenxinDeChongfenbiyaotiaojian}\label{PanduanSanjiaoxingQuanjuduichenxinDeChongfenbiyaotiaojian}
	If $\triangle ABC$ is given, point $P$ is a point on the plane of $\triangle ABC$, point $O$ is any given point in space, and:
	\[\overrightarrow{OP}=\alpha _{A}^{P}\overrightarrow{OA}+\alpha _{B}^{P}\overrightarrow{OB}+\alpha _{C}^{P}\overrightarrow{OC},\]
	\[\alpha _{A}^{P}+\alpha _{B}^{P}+\alpha _{C}^{P}=1.\]
	
	Then the necessary and sufficient condition for point $P$ to be the global center of symmetry of $\triangle ABC$ is:
	\[\underline{{{\left( \alpha _{A}^{P} \right)}_{AB}}}=\alpha _{B}^{P},\quad \underline{{{\left( \alpha _{B}^{P} \right)}_{AB}}}=\alpha _{A}^{P},\quad \underline{{{\left( \alpha _{C}^{P} \right)}_{AB}}}=\alpha _{C}^{P};\]
	\[\underline{{{\left( \alpha _{A}^{P} \right)}_{BC}}}=\alpha _{A}^{P},\quad \underline{{{\left( \alpha _{B}^{P} \right)}_{BC}}}=\alpha _{C}^{P},\quad \underline{{{\left( \alpha _{C}^{P} \right)}_{BC}}}=\alpha _{B}^{P};\]
	\[\underline{{{\left( \alpha _{A}^{P} \right)}_{CA}}}=\alpha _{C}^{P},\quad \underline{{{\left( \alpha _{B}^{P} \right)}_{CA}}}=\alpha _{B}^{P},\quad \underline{{{\left( \alpha _{C}^{P} \right)}_{CA}}}=\alpha _{A}^{P}.\]
\end{theorem}


\begin{proof}
	According to theorem\ref{PandingSanjiaoxingJubuduichenxinDeChongfenbiyaotiaojian}, this theorem is directly obtained.
\end{proof}
\hfill $\square$\par

%

\begin{theorem}{Sufficient and necessary conditions for centers of rotation, Daiyuan Zhang}{PandingSanjiaoxingLunhuanxinDeChongfenbiyaotiaojian}\label{PandingSanjiaoxingLunhuanxinDeChongfenbiyaotiaojian}
	If $\triangle ABC$ is given, point $P$ is a point on the plane of $\triangle ABC$, point $O$ is any given point in space, and:
	\[\overrightarrow{OP}=\alpha _{A}^{P}\overrightarrow{OA}+\alpha _{B}^{P}\overrightarrow{OB}+\alpha _{C}^{P}\overrightarrow{OC},\]
	\[\alpha _{A}^{P}+\alpha _{B}^{P}+\alpha _{C}^{P}=1.\]
	
	Then the necessary and sufficient condition for point $P$ to be the center of rotation of $\triangle ABC$ is:
	\[\underline{\alpha _{A}^{P}}=\alpha _{B}^{P},\quad \underline{\alpha _{B}^{P}}=\alpha _{C}^{P},\quad \underline{\alpha _{C}^{P}}=\alpha _{A}^{P}.\]
\end{theorem}

\begin{proof}
	Sufficiency. Let
	\[\underline{\alpha _{A}^{P}}=\alpha _{B}^{P},\quad \underline{\alpha _{B}^{P}}=\alpha _{C}^{P},\quad \underline{\alpha _{C}^{P}}=\alpha _{A}^{P}.\]
	
	According to the theorem \ref{thm:LunhuanDeDaishuYunsuanXingzhi}, it is obtained that
	\begin{align*}
		\underline{\overrightarrow{OP}}& =\underline{\alpha _{A}^{P}\overrightarrow{OA}+\alpha _{B}^{P}\overrightarrow{OB}+\alpha _{C}^{P}\overrightarrow{OC}} \\ 
		& =\underline{\alpha _{A}^{P}}\underline{\overrightarrow{OA}}+\underline{\alpha _{B}^{P}}\underline{\overrightarrow{OB}}+\underline{\alpha _{C}^{P}}\underline{\overrightarrow{OC}} \\ 
		& =\underline{\alpha _{A}^{P}}\overrightarrow{OB}+\underline{\alpha _{B}^{P}}\overrightarrow{OC}+\underline{\alpha _{C}^{P}}\overrightarrow{OA} \\ 
		& =\alpha _{B}^{P}\overrightarrow{OB}+\alpha _{C}^{P}\overrightarrow{OC}+\alpha _{A}^{P}\overrightarrow{OA} \\ 
		& =\overrightarrow{OP}.  
	\end{align*}
	
	Necessity. Let
	\[\underline{\overrightarrow{OP}}=\overrightarrow{OP}.\]
	
	Because
	\begin{align*}
		\underline{\overrightarrow{OP}}& =\underline{\alpha _{A}^{P}\overrightarrow{OA}+\alpha _{B}^{P}\overrightarrow{OB}+\alpha _{C}^{P}\overrightarrow{OC}} \\ 
		& =\underline{\alpha _{A}^{P}}\underline{\overrightarrow{OA}}+\underline{\alpha _{B}^{P}}\underline{\overrightarrow{OB}}+\underline{\alpha _{C}^{P}}\underline{\overrightarrow{OC}} \\ 
		& =\underline{\alpha _{A}^{P}}\overrightarrow{OB}+\underline{\alpha _{B}^{P}}\overrightarrow{OC}+\underline{\alpha _{C}^{P}}\overrightarrow{OA} \\ 
		& =\underline{\alpha _{C}^{P}}\overrightarrow{OA}+\underline{\alpha _{A}^{P}}\overrightarrow{OB}+\underline{\alpha _{B}^{P}}\overrightarrow{OC},  
	\end{align*}
	so
	\begin{align*}
		\overrightarrow{OP}-\underline{\overrightarrow{OP}}& =\left( \alpha _{A}^{P}\overrightarrow{OA}+\alpha _{B}^{P}\overrightarrow{OB}+\alpha _{C}^{P}\overrightarrow{OC} \right)-\left( \underline{\alpha _{C}^{P}}\overrightarrow{OA}+\underline{\alpha _{A}^{P}}\overrightarrow{OB}+\underline{\alpha _{B}^{P}}\overrightarrow{OC} \right) \\ 
		& =\left( \alpha _{A}^{P}-\underline{\alpha _{C}^{P}} \right)\overrightarrow{OA}+\left( \alpha _{B}^{P}-\underline{\alpha _{A}^{P}} \right)\overrightarrow{OB}+\left( \alpha _{C}^{P}-\underline{\alpha _{B}^{P}} \right)\overrightarrow{OC} \\ 
		& =\overrightarrow{0}.  
	\end{align*}
	
	And
	\[\underline{\alpha _{A}^{P}}+\underline{\alpha _{B}^{P}}+\underline{\alpha _{C}^{P}}=\underline{\alpha _{A}^{P}+\alpha _{B}^{P}+\alpha _{C}^{P}}=\underline{1}=1,\]
	therefore
	\begin{align*}
		& \left( \alpha _{A}^{P}-\underline{\alpha _{C}^{P}} \right)+\left( \alpha _{B}^{P}-\underline{\alpha _{A}^{P}} \right)+\left( \alpha _{C}^{P}-\underline{\alpha _{B}^{P}} \right) \\ 
		& =\left( \alpha _{A}^{P}+\alpha _{B}^{P}+\alpha _{C}^{P} \right)-\left( \underline{\alpha _{A}^{P}}+\underline{\alpha _{B}^{P}}+\underline{\alpha _{C}^{P}} \right) \\ 
		& =1-1 \\ 
		& =0. \\ 
	\end{align*}
	
	According to theorem \ref{thm:Thm3.4.3}, it can be obtained that:
	\[\alpha _{A}^{P}-\underline{\alpha _{C}^{P}}=0,\quad \alpha _{B}^{P}-\underline{\alpha _{A}^{P}}=0,\quad \alpha _{C}^{P}-\underline{\alpha _{B}^{P}}=0.\]
	
	i.e.
	\[\underline{\alpha _{A}^{P}}=\alpha _{B}^{P},\quad \underline{\alpha _{B}^{P}}=\alpha _{C}^{P},\quad \underline{\alpha _{C}^{P}}=\alpha _{A}^{P}.\]
\end{proof}
\hfill $\square$\par

The above equation indicates that if the point $P$ is known to be a rotation center of $\triangle ABC$, then as long as one frame component is known, the other two frame components can be directly obtained through rotation.

%

\begin{example}{}\label{Qiuzheng}
	Assuming that one of the IC of $\triangle ABC$ is $P$, and its frame components are:
	\[\alpha _{A}^{P}=\frac{BC}{BC+CA+AB}=\frac{a}{a+b+c},\]
	\[\alpha _{B}^{P}=\frac{CA}{CA+AB+BC}=\frac{b}{b+c+a},\]
	\[\alpha _{C}^{P}=\frac{AB}{AB+BC+CA}=\frac{c}{c+a+b}.\]
	
	Prove that the point $P$ is  the center of rotation of $\triangle ABC$.
\end{example}

%

\begin{solution}
	Since
	\[\underline{\alpha _{A}^{P}}=\frac{CA}{CA+AB+BC}=\alpha _{B}^{P},\]
	\[\underline{\alpha _{B}^{P}}=\frac{AB}{AB+BC+CA}=\alpha _{C}^{P},\]
	\[\underline{\alpha _{C}^{P}}=\frac{BC}{BC+CA+AB}=\alpha _{A}^{P}.\]
	
	According to theorem \ref{thm:PandingSanjiaoxingLunhuanxinDeChongfenbiyaotiaojian}, the IC $P$ is a rotation center of $\triangle ABC$, and in fact, this center is the incenter of $\triangle ABC$.
\end{solution}
\hfill $\diamond$\par

According to the above theorems, it is easy to prove that the centroid, incenter, orthocenter, circumcenter, etc. of a triangle are all local symmetry centers, global symmetry centers, and rotation centers of the triangle.

The following is a study of the relationship between the local symmetry center, global symmetry center, and rotation center of a triangle.


Obviously, according to the definition, the global symmetry center must be a local symmetry center, that is, the global symmetry center is a subset of the local symmetry center. Now consider the following question: Is there a local center of symmetry that is not the global center of symmetry? The answer is affirmative, and the following example provides an explanation.

According to example \ref{YuShangyigeLiziJinxingBijiao}, the frame component of the midpoint between the vertex $C $ and the orthocenter $H $ of $\triangle ABC$ is:
\[\alpha _{A}^{{{P}_{m\left( CH \right)}}}=\frac{\tan A}{2\left( \tan A+\tan B+\tan C \right)},\]
\[\alpha _{B}^{{{P}_{m\left( CH \right)}}}=\frac{\tan B}{2\left( \tan A+\tan B+\tan C \right)},\]
\[\alpha _{C}^{{{P}_{m\left( CH \right)}}}=\frac{\tan A+\tan B+2\tan C}{2\left( \tan A+\tan B+\tan C \right)}.\]

Therefore
\begin{align*}
	\underline{{{\left( \alpha _{B}^{{{P}_{m\left( CH \right)}}} \right)}_{AB}}}& =\underline{{{\left( \frac{\tan B}{2\left( \tan A+\tan B+\tan C \right)} \right)}_{AB}}} \\ 
	& =\underline{\frac{\tan A}{2\left( \tan B+\tan A+\tan C \right)}} \\ 
	& =\alpha _{A}^{{{P}_{m\left( CH \right)}}},  
\end{align*}
\begin{align*}
	\underline{{{\left( \alpha _{A}^{{{P}_{m\left( CH \right)}}} \right)}_{AB}}}& =\underline{{{\left( \frac{\tan A}{2\left( \tan A+\tan B+\tan C \right)} \right)}_{AB}}} \\ 
	& =\underline{\frac{\tan B}{2\left( \tan B+\tan A+\tan C \right)}} \\ 
	& =\alpha _{B}^{{{P}_{m\left( CH \right)}}},  
\end{align*}
\begin{align*}
	\underline{{{\left( \alpha _{C}^{{{P}_{m\left( CH \right)}}} \right)}_{AB}}}& =\underline{{{\left( \frac{\tan A+\tan B+2\tan C}{2\left( \tan A+\tan B+\tan C \right)} \right)}_{AB}}} \\ 
	& =\frac{\tan B+\tan A+2\tan C}{2\left( \tan B+\tan A+\tan C \right)} \\ 
	& =\alpha _{C}^{{{P}_{m\left( CH \right)}}}.  
\end{align*}

According to theorem \ref{thm:PandingSanjiaoxingJubuduichenxinDeChongfenbiyaotiaojian}, point ${{P}_{m\left( CH \right)}}$ is the $AB$ local symmetry center of $\triangle ABC$, so point ${{P}_{m\left( CH \right)}}$ is the local symmetry center of $\triangle ABC$.

And
\begin{align*}
	\underline{{{\left( \alpha _{C}^{{{P}_{m\left( CH \right)}}} \right)}_{CA}}}& =\underline{{{\left( \frac{\tan A+\tan B+2\tan C}{2\left( \tan A+\tan B+\tan C \right)} \right)}_{CA}}} \\ 
	& =\underline{{{\left( \frac{\tan C+\tan B+2\tan A}{2\left( \tan C+\tan B+\tan A \right)} \right)}_{CA}}} \\ 
	& \ne \frac{\tan A}{2\left( \tan A+\tan B+\tan C \right)}=\alpha _{A}^{{{P}_{m\left( CH \right)}}}  
\end{align*}

So, according to theorem\ref{PanduanSanjiaoxingQuanjuduichenxinDeChongfenbiyaotiaojian}, point ${{P}_{m\left( CH \right)}}$ is not the global center of symmetry of $\triangle ABC$.

Summarize the above discussion into the following theorem.


\begin{theorem}{Relationship between global and local symmetry center, Daiyuan Zhang}{SanjiaoxingDeQuanjuduichenxinYuJubuduichenxinDeGuanxi}\label{SanjiaoxingDeQuanjuduichenxinYuJubuduichenxinDeGuanxi}
 	The global center of symmetry is a true subset of the local center of symmetry of a triangle.
\end{theorem}

The following is a study of the relationship between the global center of symmetry and the rotation center of a triangle.


\begin{theorem}{Relationship between global center of symmetry and center of rotation, Daiyuan Zhang}{SanjiaoxingDeQuanjuduichenxinYuLunhuanxinDeGuanxi}\label{SanjiaoxingDeQuanjuduichenxinYuLunhuanxinDeGuanxi}
	The global center of symmetry of a triangle must be the center of rotation.
\end{theorem}

%
%
%

\begin{proof}
	According to definition \ref{def:SanjiaoxingDeQuanjuduichenxin}, if point $P$ is the global center of symmetry of $\triangle ABC$, then there must be
	\[\underline{{{\overrightarrow{OP}}_{AB}}}=\underline{{{\overrightarrow{OP}}_{BC}}}=\underline{{{\overrightarrow{OP}}_{CA}}}=\overrightarrow{OP},\]
	where
	\[\overrightarrow{OP}=\alpha _{A}^{P}\overrightarrow{OA}+\alpha _{B}^{P}\overrightarrow{OB}+\alpha _{C}^{P}\overrightarrow{OC},\]
	\[\alpha _{A}^{P}+\alpha _{B}^{P}+\alpha _{C}^{P}=1.\]
	
	And
	\[\underline{{{\left( \underline{{{\overrightarrow{OP}}_{AB}}} \right)}_{CA}}}=\underline{{{\left( \underline{{{\overrightarrow{OP}}_{CA}}} \right)}_{CA}}}=\underline{{{\overrightarrow{OP}}_{CA}}}=\overrightarrow{OP}.\]
	
	Calculate the left side of the above equation. denote $\overrightarrow{OP}:=\overrightarrow{OP}\left( A,B,C \right)$, so
	\[\underline{{{\overrightarrow{OP}}_{AB}}}=\overrightarrow{OP}\left( B,A,C \right),\]
	\[\underline{{{\left( \underline{{{\overrightarrow{OP}}_{AB}}} \right)}_{CA}}}=\underline{{{\left( \overrightarrow{OP}\left( B,A,C \right) \right)}_{CA}}}=\overrightarrow{OP}\left( B,C,A \right)=\underline{\overrightarrow{OP}},\]
	therefore
	\[\underline{\overrightarrow{OP}}=\overrightarrow{OP}.\]
	
	According to theorem \ref{thm:PandingSanjiaoxingLunhuanxinDeChongfenbiyaotiaojian}, point $P$ is the rotation center of $\triangle ABC$.
\end{proof}
\hfill $\square$\par

The above theorem indicates that the global center of symmetry is a subset of the rotational center, i.e. $\mathbf{GS}\subseteq \mathbf{RS}$.


Many centers of the triangle are both the global symmetry center and the rotation center, such as the centroid, incenter, circumcenter, etc. So, is there a center that is a rotational center but not a global symmetry center? The answer is affirmative. For example, the Brocard points are rotation centers rather than global symmetry centers.

According to example \ref{BuluokadianDeBiaojiafenliang}, it is known that the frame component of the first Brocard point ${{B}_{1}}$ is:
\[\alpha _{A}^{{{B}_{1}}}=\frac{\frac{1}{{{b}^{2}}}}{\frac{1}{{{b}^{2}}}+\frac{1}{{{c}^{2}}}+\frac{1}{{{a}^{2}}}},\]
\[\alpha _{B}^{{{B}_{1}}}=\frac{\frac{1}{{{c}^{2}}}}{\frac{1}{{{c}^{2}}}+\frac{1}{{{a}^{2}}}+\frac{1}{{{b}^{2}}}},\]
\[\alpha _{C}^{{{B}_{1}}}=\frac{\frac{1}{{{a}^{2}}}}{\frac{1}{{{a}^{2}}}+\frac{1}{{{b}^{2}}}+\frac{1}{{{c}^{2}}}}.\]

And
\[\overrightarrow{O{{B}_{1}}}=\alpha _{A}^{{{B}_{1}}}\overrightarrow{OA}+\alpha _{B}^{{{B}_{1}}}\overrightarrow{OB}+\alpha _{C}^{{{B}_{1}}}\overrightarrow{OC},\]
\[\alpha _{A}^{{{B}_{1}}}+\alpha _{B}^{{{B}_{1}}}+\alpha _{C}^{{{B}_{1}}}=1.\]

Therefore
\[\underline{\alpha _{A}^{{{B}_{1}}}}=\underline{\left( \frac{\frac{1}{{{b}^{2}}}}{\frac{1}{{{b}^{2}}}+\frac{1}{{{c}^{2}}}+\frac{1}{{{a}^{2}}}} \right)}=\frac{\frac{1}{{{c}^{2}}}}{\frac{1}{{{c}^{2}}}+\frac{1}{{{a}^{2}}}+\frac{1}{{{b}^{2}}}}=\alpha _{B}^{{{B}_{1}}},\]
\[\underline{\alpha _{B}^{{{B}_{1}}}}=\underline{\left( \frac{\frac{1}{{{c}^{2}}}}{\frac{1}{{{c}^{2}}}+\frac{1}{{{a}^{2}}}+\frac{1}{{{b}^{2}}}} \right)}=\frac{\frac{1}{{{a}^{2}}}}{\frac{1}{{{a}^{2}}}+\frac{1}{{{b}^{2}}}+\frac{1}{{{c}^{2}}}}=\alpha _{C}^{{{B}_{1}}},\]
\[\underline{\alpha _{C}^{{{B}_{1}}}}=\underline{\left( \frac{\frac{1}{{{a}^{2}}}}{\frac{1}{{{a}^{2}}}+\frac{1}{{{b}^{2}}}+\frac{1}{{{c}^{2}}}} \right)}=\frac{\frac{1}{{{b}^{2}}}}{\frac{1}{{{b}^{2}}}+\frac{1}{{{c}^{2}}}+\frac{1}{{{a}^{2}}}}=\alpha _{A}^{{{B}_{1}}}.\]

According to theorem \ref{thm:PandingSanjiaoxingLunhuanxinDeChongfenbiyaotiaojian}, the first Brocard point ${{B}_{1}}$ is a rotation center. Similarly, it can be proven that the second Brocard point ${{B}_{2}}$ is also a rotation center.

Here is proof of the first Brocard point ${{B}_{1}}$ is not the center of symmetry.
\[\underline{{{\left( \alpha _{A}^{{{B}_{1}}} \right)}_{AB}}}=\underline{{{\left( \frac{\frac{1}{{{b}^{2}}}}{\frac{1}{{{b}^{2}}}+\frac{1}{{{c}^{2}}}+\frac{1}{{{a}^{2}}}} \right)}_{AB}}}=\frac{\frac{1}{{{a}^{2}}}}{\frac{1}{{{a}^{2}}}+\frac{1}{{{c}^{2}}}+\frac{1}{{{b}^{2}}}}\ne \alpha _{B}^{{{B}_{1}}},\]

According to theorem \ref{thm:PandingSanjiaoxingJubuduichenxinDeChongfenbiyaotiaojian}, the first Brocard point ${{B}_{1}}$ is not a local symmetry center of $AB$. Similarly, it can be proven that the first Brocard point ${{B}_{1}}$ is not a local symmetry center of $BC$, nor is it a local symmetry center of $CA$.

It can also be proven that the first Brocard point ${{B}_{1}}$ is not a global center of symmetry.

Therefore, it can be concluded that the first Brocard point ${{B}_{1}}$ is not the center of symmetry.

Similarly, the second Brocard point ${{B}_{2}}$ is the rotation center rather than the symmetry center.


So we obtained the following theorem.


\begin{theorem}{Relationship between global center of symmetry and center of rotation, Daiyuan Zhang}{SanjiaoxingDeQuanjuduichenxinHeLunhuanduichenxinDeGuanxi}\label{SanjiaoxingDeQuanjuduichenxinHeLunhuanduichenxinDeGuanxi}
	The global center of symmetry of a triangle is a true subset of the rotation center. Namely,
	\[\mathbf{GS}\subset \mathbf{RS}.\]
\end{theorem}

According to the previous discussion, the Brocard point is the center of rotation rather than the center of symmetry, the ${{P}_{m\left( CH \right)}}$ is the center of symmetry rather than the rotation center. So we obtained the following theorem.


\begin{theorem}{Relationship between rotation center and symmetry center, Daiyuan Zhang}{SanjiaoxingDeLunhuanxinYuDuichenxinDeGuanxi}\label{SanjiaoxingDeLunhuanxinYuDuichenxinDeGuanxi}
	The rotation centers and symmetry centers of $\triangle ABC$ do not contain each other, that is
	$\mathbf{CS}\not\subset \mathbf{RS},\quad \mathbf{RS}\not\subset \mathbf{CS}.$
\end{theorem}

The above theorems clearly outline the global center of symmetry, the relationship between the center of symmetry and the rotation center of $\triangle ABC$. The current question is, how many are the local symmetry centers, global symmetry centers, and rotation centers of $\triangle ABC$? To answer these questions, first prove the following theorem.


\begin{theorem}{Midpoint of two centers is also a center, Daiyuan Zhang}{SanjiaoxingZhongdianDeXin}\label{SanjiaoxingZhongdianDeXin}
	Given a $\triangle ABC$, points ${{P}_{1}}$ and ${{P}_{2}}$ are located on the plane where $\triangle ABC$ is located, assuming point $P$ is the midpoint of the line segment ${{P}_{1}}{{P}_{2}}$, if the points ${{P}_{1}}$ and ${{P}_{2}}$ are both the local center of symmetry (global center of symmetry or rotation center) of $\triangle ABC$, then point $P$ is also the local center of symmetry (global center of symmetry or rotation center) of $\triangle ABC$.
\end{theorem}

%
%
%

\begin{proof}
	Take the example of rotation center to prove.
	
	Assuming point ${{P}_{1}}$ and point ${{P}_{2}}$ are all rotation centers of $\triangle ABC$, assuming the frame component of point ${{P}_{1}}$ is $\alpha _{A}^{{{P}_{1}}}$, $\alpha _{B}^{{{P}_{1}}}$, $\alpha _{C}^{{{P}_{1}}}$;  The frame component of point ${{P}_{2}}$ is $\alpha _{A}^{{{P}_{2}}}$, $\alpha _{B}^{{{P}_{2}}}$, $\alpha _{C}^{{{P}_{2}}}$. According to corollary \ref{cor:ZhongdianDeBiaojiafenliang}, the frame component of midpoint $P$ is:
	\[\alpha _{A}^{P}=\frac{1}{2}\left( \alpha _{A}^{{{P}_{1}}}+\alpha _{A}^{{{P}_{2}}} \right),\]
	\[\alpha _{B}^{P}=\frac{1}{2}\left( \alpha _{B}^{{{P}_{1}}}+\alpha _{B}^{{{P}_{2}}} \right),\]
	\[\alpha _{C}^{P}=\frac{1}{2}\left( \alpha _{C}^{{{P}_{1}}}+\alpha _{C}^{{{P}_{2}}} \right).\]
	
	So, according to theorem \ref{thm:LunhuanDeDaishuYunsuanXingzhi} in the appendix, the following conclusion can be obtained:
	\[\underline{\alpha _{A}^{P}}=\frac{1}{2}\left( \underline{\alpha _{A}^{{{P}_{1}}}}+\underline{\alpha _{A}^{{{P}_{2}}}} \right)=\frac{1}{2}\left( \alpha _{B}^{{{P}_{1}}}+\alpha _{B}^{{{P}_{2}}} \right)=\alpha _{B}^{P},\]
	\[\underline{\alpha _{B}^{P}}=\frac{1}{2}\left( \underline{\alpha _{B}^{{{P}_{1}}}}+\underline{\alpha _{B}^{{{P}_{2}}}} \right)=\frac{1}{2}\left( \alpha _{C}^{{{P}_{1}}}+\alpha _{C}^{{{P}_{2}}} \right)=\alpha _{C}^{P},\]
	\[\underline{\alpha _{C}^{P}}=\frac{1}{2}\left( \underline{\alpha _{C}^{{{P}_{1}}}}+\underline{\alpha _{C}^{{{P}_{2}}}} \right)=\frac{1}{2}\left( \alpha _{A}^{{{P}_{1}}}+\alpha _{A}^{{{P}_{2}}} \right)=\alpha _{A}^{P}.\]
	
	Therefore,  the midpoint $P$ of line segments ${{P}_{1}}{{P}_{2}}$ is the rotation center of $\triangle ABC$.

	For local symmetry centers and global symmetry centers, the proof method is similar.
\end{proof}
\hfill $\square$\par


If there are two different rotation centers ${{P}_{1}}$ and ${{P}_{2}}$ on the plane where $\triangle ABC$ is located (assume ${{P}_{1}}$ and ${{P}_{2}}$   do not coincide with each other) The midpoint of ${{P}_{1}}$ and ${{P}_{2}}$ is ${{P}_{3}}$. The midpoint of ${{P}_{1}}$ and ${{P}_{3}}$ is ${{P}_{4}}$, and so on. This results in an infinite point sequence ${{P}_{1}}$, ${{P}_{2}}$, ${{P}_{3}}$, ${{P}_{4}}$, ……, and any two points in this infinite point do not coincide with each other. Therefore, there are infinitely many rotation centers for $\triangle ABC$. The current question is whether there are two different rotation centers ${{P}_{1}}$ and ${{P}_{2}}$ for $\triangle ABC$? The answer is affirmative. For example, the first and second Brocard points are both rotation center. For non-equilateral triangles, the first and second Brocard points do not coincide, so an infinite number of rotation centers can be constructed using the above method.

Similarly, the centroid and incenter of a triangle are both global centers of symmetry. For non-equilateral triangles, if the center of centroid and incenter do not coincide, an infinite number of global symmetry centers can also be constructed using the above method.

It is also possible to construct an infinite number of local symmetry centers.

An equilateral triangle is an extremely special type of triangle that lacks generality and is not specifically discussed. I am studying general triangles. Unless otherwise specified, the triangles mentioned in this book are all general triangles.

Summarize the results of the above discussion into the following theorem.


\begin{theorem}{A triangle has infinite centers, Daiyuan Zhang}{SanjiaoxingYouWuqiongduogeXin}\label{SanjiaoxingYouWuqiongduogeXin}
	A triangle has an infinite number of global symmetry centers (local symmetry centers, rotation centers).
\end{theorem}

Since a triangle has infinite centers, can the centers of the triangle cover the entire plane in which the triangle is located? The answer is negative.


\begin{definition}{The non-center of a triangle}{SanjiaoxingDeFeixin}\label{SanjiaoxingDeFeixin}
	If point $P $ is not the center of $\triangle ABC$, i.e. $P\notin \mathbf{CT}$, then point $P $ is called the non-center of $\triangle ABC$, abbreviated as the non-center of the triangle. The set of all non-centers of the triangle is denoted as $\mathbf{NC}$.
\end{definition}

Consider the following examples. Assuming:
\[\overrightarrow{OP}=\alpha _{A}^{P}\overrightarrow{OA}+\alpha _{B}^{P}\overrightarrow{OB}+\alpha _{C}^{P}\overrightarrow{OC},\]
\[\alpha _{A}^{P}+\alpha _{B}^{P}+\alpha _{C}^{P}=1.\]

Where
\[\alpha _{A}^{P}=\frac{a}{a+b+c},\]
\[\alpha _{B}^{P}=\frac{b+ka}{a+b+c},\]
\[\alpha _{C}^{P}=\frac{c-ka}{a+b+c}.\]

Where $k $ is a positive integer, hence
\begin{align*}
	\underline{\overrightarrow{OP}}& =\underline{\alpha _{A}^{P}}\overrightarrow{OB}+\underline{\alpha _{B}^{P}}\overrightarrow{OC}+\underline{\alpha _{C}^{P}}\overrightarrow{OA} \\ 
	& =\frac{b}{a+b+c}\overrightarrow{OB}+\frac{c+kb}{a+b+c}\overrightarrow{OC}+\frac{a-kb}{a+b+c}\overrightarrow{OA}.  
\end{align*}

And
\begin{align*}
	\overrightarrow{OP}& =\alpha _{A}^{P}\overrightarrow{OA}+\alpha _{B}^{P}\overrightarrow{OB}+\alpha _{C}^{P}\overrightarrow{OC} \\ 
	& =\frac{a}{a+b+c}\overrightarrow{OA}+\frac{b+ka}{a+b+c}\overrightarrow{OB}+\frac{c-ka}{a+b+c}\overrightarrow{OC}.  
\end{align*}

Therefore, $\underline{\overrightarrow{OP}}\ne \overrightarrow{OP}$, so  $P$ is not the rotation center of $\triangle ABC$. So, $P$ is not the global center of symmetry of $\triangle ABC$.

And
\begin{align*}
	\underline{{{\overrightarrow{OP}}_{AB}}}& =\underline{{{\left( \alpha _{A}^{P}\overrightarrow{OA}+\alpha _{B}^{P}\overrightarrow{OB}+\alpha _{C}^{P}\overrightarrow{OC} \right)}_{AB}}} \\ 
	& =\underline{{{\left( \alpha _{A}^{P} \right)}_{AB}}}\overrightarrow{OB}+\underline{{{\left( \alpha _{B}^{P} \right)}_{AB}}}\overrightarrow{OA}+\underline{{{\left( \alpha _{C}^{P} \right)}_{AB}}}\overrightarrow{OC} \\ 
	& =\frac{b}{a+b+c}\overrightarrow{OB}+\frac{a+kb}{a+b+c}\overrightarrow{OA}+\frac{c-kb}{a+b+c}\overrightarrow{OC} \\ 
	& \ne \overrightarrow{OP}.  
\end{align*}

According to definition \ref{def:SanjiaoxingDeJubuduichenxin}, point $P$ is not the local symmetry center of $AB$  of $\triangle ABC$.

Similarly,
\begin{align*}
	\underline{{{\overrightarrow{OP}}_{BC}}}& =\underline{{{\left( \alpha _{A}^{P}\overrightarrow{OA}+\alpha _{B}^{P}\overrightarrow{OB}+\alpha _{C}^{P}\overrightarrow{OC} \right)}_{BC}}} \\ 
	& =\underline{{{\left( \alpha _{A}^{P} \right)}_{BC}}}\overrightarrow{OA}+\underline{{{\left( \alpha _{B}^{P} \right)}_{BC}}}\overrightarrow{OC}+\underline{{{\left( \alpha _{C}^{P} \right)}_{BC}}}\overrightarrow{OB} \\ 
	& =\frac{a}{a+b+c}\overrightarrow{OA}+\frac{c+ka}{a+b+c}\overrightarrow{OC}+\frac{b-ka}{a+b+c}\overrightarrow{OB} \\ 
	& \ne \overrightarrow{OP}.  
\end{align*}
\begin{align*}
	\underline{{{\overrightarrow{OP}}_{CA}}}& =\underline{{{\left( \alpha _{A}^{P}\overrightarrow{OA}+\alpha _{B}^{P}\overrightarrow{OB}+\alpha _{C}^{P}\overrightarrow{OC} \right)}_{CA}}} \\ 
	& =\underline{{{\left( \alpha _{A}^{P} \right)}_{CA}}}\overrightarrow{OC}+\underline{{{\left( \alpha _{B}^{P} \right)}_{CA}}}\overrightarrow{OB}+\underline{{{\left( \alpha _{C}^{P} \right)}_{CA}}}\overrightarrow{OA} \\ 
	& =\frac{c}{a+b+c}\overrightarrow{OC}+\frac{b+kc}{a+b+c}\overrightarrow{OB}+\frac{a-kc}{a+b+c}\overrightarrow{OA} \\ 
	& \ne \overrightarrow{OP}.  
\end{align*}

So according to definition \ref{def:SanjiaoxingDeJubuduichenxin}, point $P $ is not the local symmetry center of $BC$ of $\triangle ABC$, nor is it the local symmetry center of $CA$ of $\triangle ABC$.

The above results indicate that point $P$ is not the center of the triangle, or that point $P$ is the non-center of the triangle.

Considering
\begin{align*}
	\overrightarrow{O{{P}_{k}}}& =\alpha _{A}^{P}\overrightarrow{OA}+\alpha _{B}^{P}\overrightarrow{OB}+\alpha _{C}^{P}\overrightarrow{OC} \\ 
	& =\frac{a}{a+b+c}\overrightarrow{OA}+\frac{b+ka}{a+b+c}\overrightarrow{OB}+\frac{c-ka}{a+b+c}\overrightarrow{OC}.  
\end{align*}

Where $k $ is a positive integer, hence
\begin{align*}
	\overrightarrow{O{{P}_{m}}}-\overrightarrow{O{{P}_{n}}}	&=\frac{a}{a+b+c}\overrightarrow{OA}+\frac{b+ma}{a+b+c}\overrightarrow{OB}+\frac{c-ma}{a+b+c}\overrightarrow{OC} \\ 
	& -\left( \frac{a}{a+b+c}\overrightarrow{OA}+\frac{b+na}{a+b+c}\overrightarrow{OB}+\frac{c-na}{a+b+c}\overrightarrow{OC} \right) \\ 
	& =0\cdot \overrightarrow{OA}+\frac{\left( m-n \right)a}{a+b+c}\overrightarrow{OB}+\frac{\left( n-m \right)a}{a+b+c}\overrightarrow{OC}.  
\end{align*}

Where $m $ and $n $ are positive integers and $m\ne n$. And
\[0+\frac{\left( m-n \right)a}{a+b+c}+\frac{\left( n-m \right)a}{a+b+c}=0,\]
If
\[\overrightarrow{O{{P}_{m}}}-\overrightarrow{O{{P}_{n}}}=\overrightarrow{0},\]

According to the triangle frame theorem (theorem \ref{thm:Thm3.4.3}), $\left( m-n \right)a=\left( n-m \right)a=0$, which can only be $m=n$, which is contradictory. Therefore, $\overrightarrow{O{{P}_{m}}}-\overrightarrow{O{{P}_{n}}}\ne \overrightarrow{0}$. This means that when $m\ne n$, the points ${{P}_{m}}$ and ${{P}_{n}}$ are two different points. So we get an infinite point sequence ${{P}_{1}}$, ${{P}_{2}}$, ……, where each point is the non-center of the triangle, and any two points are different from each other. Therefore, the following theorem is obtained.


\begin{theorem}{A triangle has infinite non-centroids, Daiyuan Zhang}{SanjiaoxingYouWuqiongduogeFwixin}\label{SanjiaoxingYouWuqiongduogeFwixin}
	A given triangle ($\triangle ABC$) has an infinite number of non-centers.
\end{theorem}

I studied the center of a triangle using the method of Intercenter Geometry in this section. The center of a triangle is a type of fixed point. The global centroid is not only the true subset of the local center, but also the true subset of the rotation center, and can be regarded as the most core part of the center of a triangle. This most core subset has an infinite number of members, and the familiar center of centroid, incenter, orthocenter, circumcenter are all the most core members. But Brocard point, excenter, ${{P}_{m\left( CH \right)}}$ and so on are not the most core members.

I divided the triangle into center and non-center. The center of the triangle is further divided into a center of symmetry and a center of rotation. Symmetry center is further divided into local symmetry center and global symmetry center. The local symmetry center can be further divided into local symmetry center of $AB$, local symmetry center of $BC$, and local symmetry center of $CA$. It can be summarized into the following expressions:

Assuming that the set of all points on the $\triangle ABC$ plane is denoted as $\mathbf{P}$, then
\[\mathbf{P}=\mathbf{NC}\cup \mathbf{CT},\]
\[\mathbf{CT}:=\mathbf{CS}\cup \mathbf{RS},\]
\[\mathbf{CS}:=\mathbf{LS}\cup \mathbf{GS},\]
\[\mathbf{LS}:=\mathbf{L}{{\mathbf{S}}_{AB}}\cup \mathbf{L}{{\mathbf{S}}_{BC}}\cup \mathbf{L}{{\mathbf{S}}_{CA}},\]
\[\mathbf{GS}:=\mathbf{L}{{\mathbf{S}}_{AB}}\cap \mathbf{L}{{\mathbf{S}}_{BC}}\cap \mathbf{L}{{\mathbf{S}}_{CA}}.\]


\chapter{Vector of special intersecting centers}\label{Ch9}
\thispagestyle{empty}

This chapter is one of the applications of theorem \ref{thm:Thm6.1.1} and theorem \ref{thm:Thm6.2.1}. All the conclusions in this chapter are the corollaries of theorem \ref{thm:Thm6.1.1} and theorem \ref{thm:Thm6.2.1}.

The vector between two points is the basis of geometric quantities (such as distance, etc.). In this chapter, we use theorem \ref{thm:Thm6.1.1} and theorem \ref{thm:Thm6.2.1} to study the vector from origin to intersecting center (abbreviated as VOIC), and the vector of two intersecting centers (abbreviated as VTICs), such as the vector between circumcenter and centroid, the vector between circumcenter and incenter, the vector between circumcenter and orthocenter of a triangle; In this chapter, we all study the vectors of two intersecting centers (abbreviated as VTICs) among the centroid, incenter and orthocenter. The properties of VOIC and VTICs are compared.

\section{Vector from origin to intersecting center on triangular frame}\label{Sec9.1}
This section is one of the applications of theorem \ref{thm:Thm6.1.1}. All the conclusions in this section are corollaries of theorem \ref{thm:Thm6.1.1}. This section studies some special vectors of intersecting centers, specifically the centroid vector, incenter vector, orthocenter vector and circumcenter vector on the triangular frame, and gives the formulas of those vectors of intersecting centers on the triangular frame. Through the comparison with the vector in analytic geometry, the readers can further understand the characteristics of Intercenter Geometry.

\subsection{Centroid vector of a triangle}\label{Subsec9.1.1}
Let a $\triangle ABC$ and a triangular frame $\left( O;A,B,C \right)$ be given, and point $O$ is an arbitrary point in space. According to theorem \ref{thm:Thm6.1.1} and section \ref{Sec8.1}, it is obtained that:
\[\begin{aligned}
	\overrightarrow{OG}=\alpha _{A}^{G}\overrightarrow{OA}+\alpha _{B}^{G}\overrightarrow{OB}+\alpha _{C}^{G}\overrightarrow{OC}=\frac{1}{3}\left( \overrightarrow{OA}+\overrightarrow{OB}+\overrightarrow{OC} \right).  
\end{aligned}\]	

The above formula is the vector from any point $O$ in space to the centroid $G$. Point $O$ is the origin of the frame $\left( O;A,B,C \right)$. The above formula indicates that the frame components of the centroid vector $\overrightarrow{OG}$ of the triangle are all ${1}/{3}\;$, and are independent of the position of the origin $O$ of the frame $\left( O;A,B,C \right)$.

\subsection{Incenter vector of a triangle}\label{Subsec9.1.2}
Let a $\triangle ABC$ and a triangular frame $\left( O;A,B,C \right)$ be given, and point $O$ is an arbitrary point in space. According to theorem \ref{thm:Thm6.1.1} and section \ref{Sec8.2}, it is obtained that:
\[\begin{aligned}
	\overrightarrow{OI}=\alpha _{A}^{I}\overrightarrow{OA}+\alpha _{B}^{I}\overrightarrow{OB}+\alpha _{C}^{I}\overrightarrow{OC}=\frac{a\overrightarrow{OA}+b\overrightarrow{OB}+c\overrightarrow{OC}}{a+b+c}.  
\end{aligned}\]	


The above formula is the vector from any point $O$ in space to the incenter $I$. Point $O$ is the origin of the frame $\left( O;A,B,C \right)$. As explained above, the frame components of the incenter vector $\overrightarrow{OI}$ of a triangle are ${a}/{\left(a+b+c \right)}\;$, ${b}/{\left (a+b+c \right)}\;$ and ${c}/{\left (a+b+c \right)}\;$, and are independent of the position of the origin $O$ of the frame $\left( O;A,B,C \right)$.


It is necessary to make a comparison with analytic geometry. Let's start with analytic geometry. First of all, two terms will be explained in the following: in analytic geometry, the coordinate system is equivalent to the frame in this book, and the coordinate (or coordinate component) is equivalent to the frame component in this book. Suppose that the coordinate of incenter $I$ in Cartesian rectangular coordinate system is $I\left( x,y \right)$, and the origin of the coordinate is $O$, then the vector $\overrightarrow{OI}$ is
\[\overrightarrow{OI}=x{{\mathbf{e}}_{x}}+y{{\mathbf{e}}_{y}}.\]	

Where ${{\mathbf{e}}_{x}}$ and ${{\mathbf{e}}_{y}}$ are the unit vectors in the positive direction of $\overrightarrow{OX}$ axis and $\overrightarrow{OY}$ axis, respectively. Obviously, with the change of $\triangle ABC$ (translation or rotation), the frame components $x$ and $y$ of the incenter $I$ usually change. Moreover, the frame components $x$ and $y$ of the incenter $I$ are related to the measuring unit of $\overrightarrow{OX}$ axis and  $\overrightarrow{OY}$ axis. On the other hand, if a coordinate system changes, such as a polar coordinate system instead of a Cartesian coordinate system, the coordinates (the frame components) also change. The coordinates (frame components) in analytic geometry are closely related to the coordinate system (frame) and effected a lot by human factors.

Look again at Intercenter Geometry. As we have just discussed, the frame components of the incenter vector $\overrightarrow{OI}$ of a triangle are ${a}/{\left(a+b+c \right)}\;$, ${b}/{\left (a+b+c \right)}\;$, ${c}/{\left (a+b+c \right)}\;$. This is a set of “natural” parameters that relate only to the lengths of triangular sides and do not contain any frame information. Moreover, no matter how $\triangle ABC$ is translated, rotated or even flipped, the three frame components of the incenter vector $\overrightarrow{OI}$ are always ${a}/{\left(a+b+c \right)}\;$, ${b}/{\left (a+b+c \right)}\;$, ${c}/{\left (a+b+c \right)}\;$. On the other hand, if the position of the origin of frame $O$ changes, the frame $\left( O;A,B,C \right)$ usually changes, but the components of the frame are always constant. This is very convenient in application, because the position of the origin of the frame $O$ can be selected according to the needs, without affecting its frame components.

To sum up, Intercenter Geometry is different from analytic geometry, and I believe that readers have had a preliminary understanding.

\subsection{Orthocenter vector of a triangle}\label{Subsec9.1.3}
Let a $\triangle ABC$ and a triangular frame $\left( O;A,B,C \right)$ be given, and point $O$ is an arbitrary point in space. According to theorem \ref{thm:Thm6.1.1} and section \ref{Sec8.3}, it is obtained that:
\[\begin{aligned}
	\overrightarrow{OH}&=\alpha _{A}^{H}\overrightarrow{OA}+\alpha _{B}^{H}\overrightarrow{OB}+\alpha _{C}^{H}\overrightarrow{OC} \\ 
	& =\frac{\tan A\cdot \overrightarrow{OA}+\tan B\cdot \overrightarrow{OB}+\tan C\cdot \overrightarrow{OC}}{\tan A+\tan B+\tan C}.  
\end{aligned}\]	


The above formula is the vector from any point $O$ in space to the orthocenter $H$. Point $O$ is the origin of the frame $\left( O;A,B,C \right)$. The above formula indicates that the frame components of the orthocenter vector $\overrightarrow{OH}$ of the triangle are independent of the position of the origin $O$ of the frame $\left( O;A,B,C \right)$.

\subsection{Circumcenter vector of a triangle}\label{Subsec9.1.4}
Let a $\triangle ABC$ and a triangular frame $\left( O;A,B,C \right)$ be given, and point $O$ is an arbitrary point in space. According to theorem \ref{thm:Thm6.1.1} and section \ref{Sec8.4}, it is obtained that:
\[\begin{aligned}
	\overrightarrow{OQ}&=\alpha _{A}^{Q}\overrightarrow{OA}+\alpha _{B}^{Q}\overrightarrow{OB}+\alpha _{C}^{Q}\overrightarrow{OC} \\ 
	& =\frac{\sin 2A\cdot \overrightarrow{OA}+\sin 2B\cdot \overrightarrow{OB}+\sin 2C\cdot \overrightarrow{OC}}{\sin 2A+\sin 2B+\sin 2C}.  
\end{aligned}\]	

%
The above formula is the vector from any point $O$ in space to the circumcenter $Q$. Point $O$ is the origin of the frame $\left( O;A,B,C \right)$. The above formula indicates that the frame components of the circumcenter vector $\overrightarrow{OQ}$ of the triangle are independent of the position of the origin $O$ of the frame $\left( O;A,B,C \right)$.

An excenter vector of a triangle can be derived similarly and left to the reader as an exercise.

\section{Vector of two intersecting centers on triangular frame}\label{Sec9.2}

Now we study the expression of the vector of two intersecting centers (abbreviated as VTICs). In theorem \ref{thm:Thm6.2.1}, we need to know the frame components of the two intersecting centers.

\subsection{Vector of two intersecting centers between circumcenter and centroid}\label{Subsec9.2.1}
From Chapter \ref{Ch8}, we have:
\[\alpha _{A}^{QG}=\alpha _{A}^{G}-\alpha _{A}^{Q}=\frac{1}{3}-\frac{\sin 2A}{W}=\frac{W-3\sin 2A}{3W},\]	
\[\alpha _{B}^{QG}=\alpha _{B}^{G}-\alpha _{B}^{Q}=\frac{1}{3}-\frac{\sin 2B}{W}=\frac{W-3\sin 2B}{3W},\]	
\[\alpha _{C}^{QG}=\alpha _{C}^{G}-\alpha _{C}^{Q}=\frac{1}{3}-\frac{\sin 2C}{W}=\frac{W-3\sin 2C}{3W}.\]	

If the starting IC is the circumcenter $Q$ and the ending point is the centroid $G$, then from theorem \ref{thm:Thm6.2.1}, we have
\begin{equation}\label{Eq9.2.1}
	\begin{aligned}
		\overrightarrow{QG}&=\alpha _{A}^{QG}\overrightarrow{OA}+\alpha _{B}^{QG}\overrightarrow{OB}+\alpha _{C}^{QG}\overrightarrow{OC}\\
		&=\frac{1}{3W}\left( \begin{aligned}
			\left( W-3\sin 2A \right)\overrightarrow{OA}+\left( W-3\sin 2B \right)\overrightarrow{OB}+\left( W-3\sin 2C \right)\overrightarrow{OC}
		\end{aligned} \right).  
	\end{aligned}		
\end{equation}

Comparing the VOIC and the VTICs between circumcenter and centroid, we can find that the expression of the VOIC is simple, but the origin of the frame is limited to the circumcenter $Q$. However, the expression of the VTICs is more complex, but the origin of the frame $O$ can be placed at any position. For example, if the origin of the frame $O$ is placed at the vertex of the triangle $A$, $B$, $C$, it can be obtained respectively:

\[\overrightarrow{QG}=\frac{1}{3W}\left( \left( W-3\sin 2B \right)\overrightarrow{AB}+\left( W-3\sin 2C \right)\overrightarrow{AC} \right),\]	\[\overrightarrow{QG}=\frac{1}{3W}\left( \left( W-3\sin 2A \right)\overrightarrow{BA}+\left( W-3\sin 2C \right)\overrightarrow{BC} \right),\]	\[\overrightarrow{QG}=\frac{1}{3W}\left( \left( W-3\sin 2A \right)\overrightarrow{CA}+\left( W-3\sin 2B \right)\overrightarrow{CB} \right).\]

\subsection{Vector of two intersecting centers between circumcenter and incenter}\label{Subsec9.2.2}
If the starting point of IC is the circumcenter $Q$ and the ending point is the incenter $I$, then from theorem \ref{thm:Thm6.2.1}, we have

\[\alpha _{A}^{QI}=\alpha _{A}^{I}-\alpha _{A}^{Q}=\frac{a}{2p}-\frac{\sin 2A}{W}=\frac{aW-2p\sin 2A}{2pW},\]
\[\alpha _{B}^{QI}=\alpha _{B}^{I}-\alpha _{B}^{Q}=\frac{b}{2p}-\frac{\sin 2B}{W}=\frac{bW-2p\sin 2B}{2pW},\]
\[\alpha _{C}^{QI}=\alpha _{C}^{I}-\alpha _{C}^{Q}=\frac{c}{2p}-\frac{\sin 2C}{W}=\frac{cW-2p\sin 2C}{2pW},\]
\[\begin{aligned}
	\overrightarrow{QI}&=\alpha _{A}^{QI}\overrightarrow{OA}+\alpha _{B}^{QI}\overrightarrow{OB}+\alpha _{C}^{QI}\overrightarrow{OC} \\ 
	&=\frac{1}{2pW}\left( \begin{aligned}
		\left( aW-2p\sin 2A \right)\overrightarrow{OA}+\left( bW-2p\sin 2B \right)\overrightarrow{OB}+\left( cW-2p\sin 2C \right)\overrightarrow{OC} \\ 
	\end{aligned} \right).  
\end{aligned}\]

\subsection{Vector of two intersecting centers between circumcenter and orthocenter}\label{Subsec9.2.3}
If the starting IC is the circumcenter $Q$ and the ending point is the orthocenter $H$, then from theorem \ref{thm:Thm6.2.1}, we have

\[\alpha _{A}^{QH}=\alpha _{A}^{H}-\alpha _{A}^{Q}=\frac{\tan A}{T}-\frac{\sin 2A}{W}=\frac{W\tan A-T\sin 2A}{TW},\]	
\[\alpha _{B}^{QH}=\alpha _{B}^{H}-\alpha _{B}^{Q}=\frac{\tan B}{T}-\frac{\sin 2B}{W}=\frac{W\tan B-T\sin 2B}{TW},\]	
\[\alpha _{C}^{QH}=\alpha _{C}^{H}-\alpha _{C}^{Q}=\frac{\tan C}{T}-\frac{\sin 2C}{W}=\frac{W\tan C-T\sin 2C}{TW},\]	
\[\begin{aligned}
	\overrightarrow{QH}&=\alpha _{A}^{QH}\overrightarrow{OA}+\alpha _{B}^{QH}\overrightarrow{OB}+\alpha _{C}^{QH}\overrightarrow{OC} \\ 
	& =\frac{1}{TW}\left( \begin{aligned}
		\left( W\tan A-T\sin 2A \right)\overrightarrow{OA}+\left( W\tan B-T\sin 2B \right)\overrightarrow{OB}+\left( W\tan C-T\sin 2C \right)\overrightarrow{OC} 
	\end{aligned} \right).  
\end{aligned}\]

\subsection{Vector of two intersecting centers between centroid and incenter}\label{Subsec9.2.4}
If the starting IC is the centroid $G$ and the ending point is the incenter $I$, then from Chapter \ref{Ch8} and theorem \ref{thm:Thm6.2.1}, we have

\[\alpha _{A}^{GI}=\alpha _{A}^{I}-\alpha _{A}^{G}=\frac{a}{2p}-\frac{1}{3}=\frac{3a-2p}{6p},\]
\[\alpha _{B}^{GI}=\alpha _{B}^{I}-\alpha _{B}^{G}=\frac{b}{2p}-\frac{1}{3}=\frac{3b-2p}{6p},\]
\[\alpha _{C}^{GI}=\alpha _{C}^{I}-\alpha _{C}^{G}=\frac{c}{2p}-\frac{1}{3}=\frac{3c-2p}{6p},\]
\[\begin{aligned}
	\overrightarrow{GI}&=\alpha _{A}^{GI}\overrightarrow{OA}+\alpha _{B}^{GI}\overrightarrow{OB}+\alpha _{C}^{GI}\overrightarrow{OC} \\ 
	& =\frac{1}{6p}\left( \left( 3a-2p \right)\overrightarrow{OA}+\left( 3b-2p \right)\overrightarrow{OB}+\left( 3c-2p \right)\overrightarrow{OC} \right).  
\end{aligned}\]

\subsection{Vector of two intersecting centers between centroid and orthocenter}\label{Subsec9.2.5}
If the starting IC is the centroid $G$ and the ending point is the orthocenter $H$, then from Chapter \ref{Ch8} and theorem \ref{thm:Thm6.2.1}, we have

\[\alpha _{A}^{GH}=\alpha _{A}^{H}-\alpha _{A}^{G}=\frac{\tan A}{T}-\frac{1}{3}=\frac{3\tan A-T}{3T},\]
\[\alpha _{B}^{GH}=\alpha _{B}^{H}-\alpha _{B}^{G}=\frac{\tan B}{T}-\frac{1}{3}=\frac{3\tan B-T}{3T},\]
\[\alpha _{C}^{GH}=\alpha _{C}^{H}-\alpha _{C}^{G}=\frac{\tan C}{T}-\frac{1}{3}=\frac{3\tan C-T}{3T},\]
\begin{equation}\label{Eq9.2.2}
	\begin{aligned}
		\overrightarrow{GH}&=\alpha _{A}^{GH}\overrightarrow{OA}+\alpha _{B}^{GH}\overrightarrow{OB}+\alpha _{C}^{GH}\overrightarrow{OC} \\ 
		& =\frac{1}{3T}\left( \begin{aligned}
			\left( 3\tan A-T \right)\overrightarrow{OA}+\left( 3\tan B-T \right)\overrightarrow{OB}+\left( 3\tan C-T \right)\overrightarrow{OC} 
		\end{aligned} \right).  
	\end{aligned}
\end{equation}

\subsection{Vector of two intersecting centers between incenter and orthocenter}\label{Subsec9.2.6}
If the starting IC is the incenter $I$ and the ending point is the orthocenter $H$, then from Chapter \ref{Ch8} and theorem \ref{thm:Thm6.2.1}, we have

\[\alpha _{A}^{IH}=\alpha _{A}^{H}-\alpha _{A}^{I}=\frac{\tan A}{T}-\frac{a}{2p}=\frac{2p\tan A-aT}{2pT},\]
\[\alpha _{B}^{IH}=\alpha _{B}^{H}-\alpha _{B}^{I}=\frac{\tan B}{T}-\frac{b}{2p}=\frac{2p\tan B-bT}{2pT},\]
\[\alpha _{C}^{IH}=\alpha _{C}^{H}-\alpha _{C}^{I}=\frac{\tan C}{T}-\frac{c}{2p}=\frac{2p\tan C-cT}{2pT},\]
\[\begin{aligned}
	\overrightarrow{IH}&=\alpha _{A}^{IH}\overrightarrow{OA}+\alpha _{B}^{IH}\overrightarrow{OB}+\alpha _{C}^{IH}\overrightarrow{OC} \\ 
	& =\frac{1}{2pT}\left( \begin{aligned}
		\left( 2p\tan A-aT \right)\overrightarrow{OA}+\left( 2p\tan B-bT \right)\overrightarrow{OB}+\left( 2p\tan C-cT \right)\overrightarrow{OC} 
	\end{aligned} \right).  
\end{aligned}\]	

\section{Vector from origin to special intersecting centers on triangular frame of circumcenter}\label{Sec9.3}

It is an important special case that the origin of the frame $O$ coincides with the circumcenter $Q$ because the distances from the circumcenter $Q$ to the three vertices of the triangle are equal, which leads to some simple expressions when discussing distances later.

If the starting point $O$ coincides with the circumcenter $Q$, and the ending point of IC is $P$, then according to theorem \ref{thm:Thm6.1.1}, we have

\[\overrightarrow{QP}=\alpha _{A}^{P}\overrightarrow{QA}+\alpha _{B}^{P}\overrightarrow{QB}+\alpha _{C}^{P}\overrightarrow{QC}.\]	

By using the frame components of each IC of the triangle (see Chapter \ref{Ch8}), the vector between the circumcenter of the triangle $Q$ and the other ICs can be obtained.

\subsection{Vector between circumcenter and centroid (from origin to intersecting center)}\label{Subsec9.3.1}
If the starting point is the circumcenter $Q$, and the ending point is the centroid $G$, then from Chapter \ref{Ch8} we have

\[\begin{aligned}
	\overrightarrow{QG}=\alpha _{A}^{G}\overrightarrow{QA}+\alpha _{B}^{G}\overrightarrow{QB}+\alpha _{C}^{G}\overrightarrow{QC}=\frac{1}{3}\left( \overrightarrow{QA}+\overrightarrow{QB}+\overrightarrow{QC} \right).  
\end{aligned}\]

\subsection{Vector between circumcenter and incenter (from origin to intersecting center)}\label{Subsec9.3.2}
If the starting point is the circumcenter $Q$, and the ending point is the incenter $I$, then from Chapter \ref{Ch8} we have

\[\overrightarrow{QI}=\alpha _{A}^{I}\overrightarrow{QA}+\alpha _{B}^{I}\overrightarrow{QB}+\alpha _{C}^{I}\overrightarrow{QC}=\frac{a\overrightarrow{QA}+b\overrightarrow{QB}+c\overrightarrow{QC}}{a+b+c}.\]	

\subsection{Vector between circumcenter and orthocenter (from origin to intersecting center)}\label{Subsec9.3.3}
If the starting point is the circumcenter $Q$, and the ending point is the orthocenter $H$, then from Chapter \ref{Ch8} we have

\[\begin{aligned}
	\overrightarrow{QH}=\alpha _{A}^{H}\overrightarrow{QA}+\alpha _{B}^{H}\overrightarrow{QB}+\alpha _{C}^{H}\overrightarrow{QC}=\frac{\tan A\cdot \overrightarrow{QA}+\tan B\cdot \overrightarrow{QB}+\tan C\cdot \overrightarrow{QC}}{\tan A+\tan B+\tan C}.  
\end{aligned}\]


\chapter{Proof of some famous theorems using Intercenter Geometry}\label{Ch10}
\thispagestyle{empty}
To demonstrate the application of Intercenter Geometry, this chapter uses the method of Intercenter Geometry to prove some well-known theorems.
\section{Proof of Ceva's theorem}\label{Sec10.1}
Ceva's theorem can be proved by the traditional Euclidean geometry method. Here, the method of Intercenter Geometry is used to prove this theorem.
\begin{theorem}{Ceva's theorem, Giovanni Ceva}{Thm10.1.1}\label{Thm10.1.1} 
	Given a $\triangle ABC$, $P\in {{\pi }_{ABC}}$, then
	\[\lambda _{AB}^{P}\lambda _{BC}^{P}\lambda _{CA}^{P}=1.\]	
\end{theorem}

\begin{proof}
	
	Let the point $O$ coincide with the intersecting foot $L$ of vertex $A$ (see Figure \ref{fig:tu6.1.1}), that is, $O=L\in \overleftrightarrow{PA}\bigcap \overleftrightarrow{BC}$, that is, the starting point $L$ is both on $\overleftrightarrow{PA}$ and $\overleftrightarrow{BC}$, and let the ending intersecting center be at the point $P$.
	
	Since the starting point $O$ is the intersecting foot $L$ of the point $A$, according to the concept of fractional ratio and integral ratio (see Section \ref{Sec2.1}), it is obtained that:	
	
	\[\overrightarrow{LP}={{\kappa }_{LA}}\overrightarrow{LA},\]
	\[\overrightarrow{LB}={{\kappa }_{LC}}\overrightarrow{LC},\]
	\[\overrightarrow{BL}=\lambda _{BC}^{P}\overrightarrow{LC},\]
	\[\overrightarrow{LB}=-\lambda _{BC}^{P}\overrightarrow{LC}.\]
	
	Therefore, formula (\ref{Eq7.1.4}) can be written as follows:
	\[\overrightarrow{LP}-\frac{\overrightarrow{LA}}{1+\lambda _{AB}^{P}+\lambda _{AC}^{P}}=\frac{\lambda _{AB}^{P}\overrightarrow{LB}+\lambda _{AC}^{P}\overrightarrow{LC}}{1+\lambda _{AB}^{P}+\lambda _{AC}^{P}},\]	
	i.e.
	\[\left( {{\kappa }_{LA}}-\frac{1}{1+\lambda _{AB}^{P}+\lambda _{AC}^{P}} \right)\overrightarrow{LA}=\frac{-\lambda _{BC}^{P}\lambda _{AB}^{P}+\lambda _{AC}^{P}}{1+\lambda _{AB}^{P}+\lambda _{AC}^{P}}\overrightarrow{LC}.\]
	
	Since the point $L$ is the intersecting foot of the vertex $A$ of the triangle, therefore, $\overrightarrow{LA}$ is linearly independent of $\overrightarrow{LC}$, from the above formula, the following must be holds:
	
	\[{{\kappa }_{LA}}-\frac{1}{1+\lambda _{AB}^{P}+\lambda _{AC}^{P}}=0,\]
	\[\frac{-\lambda _{BC}^{P}\lambda _{AB}^{P}+\lambda _{AC}^{P}}{1+\lambda _{AB}^{P}+\lambda _{AC}^{P}}=0,\]
	i.e.
	\[-\lambda _{BC}^{P}\lambda _{AB}^{P}+\lambda _{AC}^{P}=0.\]	
	
	Multiply both sides of the above formula by $\lambda _{CA}^{P}$($\lambda _{CA}^{P}\ne 0$), and according to the formula (\ref{Eq2.4.1}), we know $\lambda _{AC}^{P}\lambda _{CA}^{P}=1$, then we obtain the following theorem:
	
	\[\lambda _{AB}^{P}\lambda _{BC}^{P}\lambda _{CA}^{P}=1.\]	
\end{proof}
\hfill $\square$\par

Using the formula (\ref{Eq2.4.5}), we obtain Ceva's theorem in the form of integral ratio:
\[\frac{{{\kappa }_{AB}}}{{{\kappa }_{BA}}}\frac{{{\kappa }_{BC}}}{{{\kappa }_{CB}}}\frac{{{\kappa }_{CA}}}{{{\kappa }_{AC}}}=1,\]
or
\[{{\kappa }_{AB}}{{\kappa }_{BC}}{{\kappa }_{CA}}={{\kappa }_{AC}}{{\kappa }_{BA}}{{\kappa }_{CB}}.\]	

\section{Proof of Menelaus's theorem}\label{Sec10.2}
Menelaus's theorem can be proved by the traditional Euclidean geometry method. Here, the method of Intercenter Geometry is used to prove this theorem.

\begin{theorem}{Menelaus' theorem, Menelaus}{Thm10.2.1}\label{Thm10.2.1} 
	As shown in Figure \ref{fig:tu10.2.1}, take a point $D$ on the extension line of $BC$ side of $\triangle ABC$, make a straight line through $D$, intersect the $AB$ side of $\triangle ABC$ at the point $N$, and $AC$ side of $\triangle ABC$ at the point $P$, then there must be
		
	\[{{\lambda }_{AB}}{{\lambda }_{BC}}{{\lambda }_{CA}}=-1.\]
	Where
	\begin{flalign*}
		{{\lambda }_{AB}}=\frac{\overrightarrow{AN}}{\overrightarrow{NB}}, {{\lambda }_{BC}}=\frac{\overrightarrow{BD}}{\overrightarrow{DC}}, {{\lambda }_{CA}}=\frac{\overrightarrow{CP}}{\overrightarrow{PA}}.
	\end{flalign*}
\end{theorem}

\begin{proof}
	\begin{figure}[h]
		\centering
		\includegraphics[totalheight=4cm]{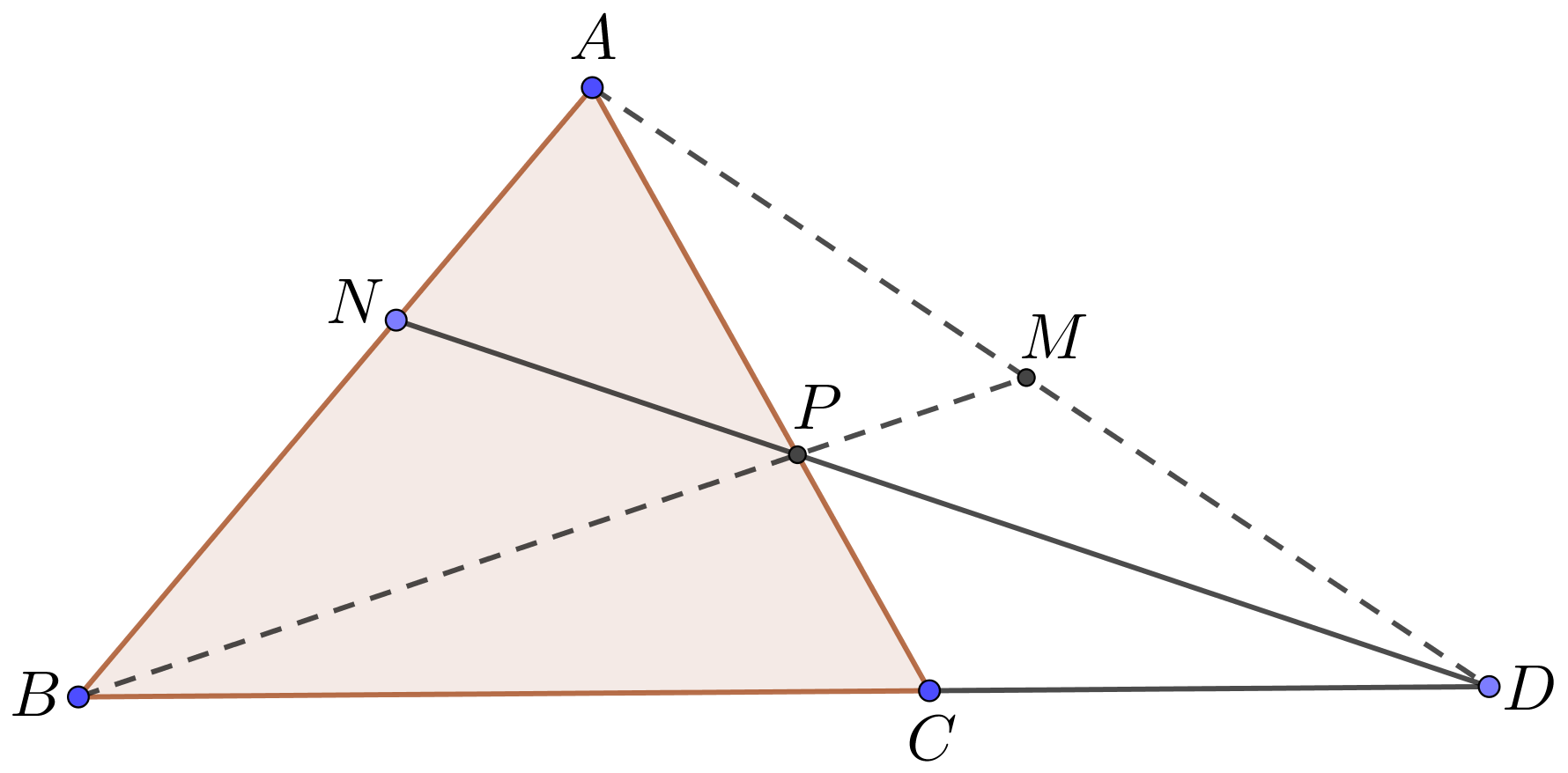}
		\caption{Geometric diagram of Menelaus' theorem} \label{fig:tu10.2.1}
	\end{figure}
	
	Here with Intercenter Geometry, a vector approach is adopted to prove this theorem. As shown in Figure \ref{fig:tu10.2.1}, considering the large triangle $\triangle ABD$, the point $P$ is the IC of $\triangle ABD$, and let the point $O$ coincide with the point $C$ of $\triangle ABC$, that is, $O=C\in \overleftrightarrow{PA}\bigcap \overleftrightarrow{BC}$, with the ending IC located at point $P$. For $\triangle ABD$, using theorem \ref{thm:Thm7.1.1} and theorem \ref{thm:Thm2.4.1}, we have
		
	\[{{\lambda }_{CA}}=\frac{1}{{{\lambda }_{AC}}}=\frac{1}{\lambda _{AB}^{P}+\lambda _{AD}^{P}}.\]	
	Where
	\begin{flalign*}
		\lambda _{AB}^{P}={{\lambda }_{AB}}=\frac{\overrightarrow{AN}}{\overrightarrow{NB}}, \lambda _{AD}^{P}=\frac{\overrightarrow{AM}}{\overrightarrow{MD}}.
	\end{flalign*}
	Therefore
	\[{{\lambda }_{CA}}\left( \lambda _{AB}^{P}+\lambda _{AD}^{P} \right)=1,\]
	\[{{\lambda }_{CA}}\left( \lambda _{AB}^{P}+\frac{1}{\lambda _{DA}^{P}} \right)=1,\]
	\[{{\lambda }_{CA}}\left( \lambda _{AB}^{P}\lambda _{DA}^{P}+1 \right)=\lambda _{DA}^{P},\]
	\[{{\lambda }_{CA}}\left( \frac{\lambda _{AB}^{P}\lambda _{BD}^{P}\lambda _{DA}^{P}}{\lambda _{BD}^{P}}+1 \right)=\lambda _{DA}^{P}.\]
	Where
	\[\lambda _{BD}^{P}=\frac{\overrightarrow{BC}}{\overrightarrow{CD}}.\]
	
	
	Using Ceva's theorem $\lambda _{AB}^{P}\lambda _{BD}^{P}\lambda _{DA}^{P}=1$ in $\triangle ABD$, the following results are obtained:
		
	\[{{\lambda }_{CA}}\left( \frac{1}{\lambda _{BD}^{P}}+1 \right)=\lambda _{DA}^{P},\]
	\[{{\lambda }_{CA}}\left( \frac{1}{\lambda _{BD}^{P}}+1 \right)=\frac{1}{\lambda _{AB}^{P}\lambda _{BD}^{P}},\]
	\[\lambda _{AB}^{P}{{\lambda }_{CA}}\left( 1+\lambda _{BD}^{P} \right)=1,\]
	\[\lambda _{AB}^{P}\left( 1+\frac{\overrightarrow{BC}}{\overrightarrow{CD}} \right){{\lambda }_{CA}}=1,\]
	\[\lambda _{AB}^{P}\left( \frac{\overrightarrow{CD}+\overrightarrow{BC}}{\overrightarrow{CD}} \right){{\lambda }_{CA}}=1,\]
	\[\lambda _{AB}^{P}\left( \frac{\overrightarrow{BD}}{\overrightarrow{CD}} \right){{\lambda }_{CA}}=1,\]
	\[\lambda _{AB}^{P}\left( \frac{\overrightarrow{BD}}{\overrightarrow{DC}} \right){{\lambda }_{CA}}=-1,\]
	i.e.
	\[{{\lambda }_{AB}}{{\lambda }_{BC}}{{\lambda }_{CA}}=-1.\]	
\end{proof}
\hfill $\square$\par

	I obtain Ceva's theorem and Menelaus' theorem in a unified approach by Intercenter Geometry.

\section{Proof of Euler's theorem}\label{Sec10.3}
	
\begin{theorem}{Euler's theorem, Leonhard Euler}{Thm10.3.1}\label{Thm10.3.1} 
	The centroid, the orthocenter and the circumcenter of a triangle are collinear, and the distance from the centroid to the circumcenter is half of the distance from the centroid to the orthocenter.
\end{theorem}

Although there are other methods to prove Euler's theorem, this book uses theorem \ref{thm:Thm6.2.1} to prove Euler's theorem, so as to illustrate the wide application of theorem \ref{thm:Thm6.2.1}. Moreover, one of the advantages of proving Euler's theorem by using theorem \ref{thm:Thm6.2.1} is that it is easy to see that the centroid lies between the orthocenter and the circumcenter (see the following formula (\ref{Eq10.3.2})). The specific proof method is given below.

\begin{proof}
If the ratio from the frame component (coefficient) of $\overrightarrow{OA}$ in formula (\ref{Eq9.2.2}) to the frame component (coefficient) of $\overrightarrow{OA}$ in formula (\ref{Eq9.2.1}) is denoted as $\gamma$, then:

\[\gamma =\frac{\alpha _{A}^{H}-\alpha _{A}^{G}}{\alpha _{A}^{Q}-\alpha _{A}^{G}}=\frac{\frac{\tan A}{\tan A+\tan B+\tan C}-\frac{1}{3}}{\frac{\sin 2A}{\sin 2A+\sin 2B+\sin 2C}-\frac{1}{3}},\]
i.e.
\[\gamma =J\cdot \frac{2\tan A-\tan B-\tan C}{2\sin 2A-\sin 2B-\sin 2C}=J\cdot \frac{\tan A-\tan B+\tan A-\tan C}{\sin 2A-\sin 2B+\sin 2A-\sin 2C},\]
i.e.
\begin{equation}\label{Eq10.3.1}
	\gamma =\frac{J}{2}\cdot \frac{\frac{a}{\cos A}-\frac{b}{\cos B}+\frac{a}{\cos A}-\frac{c}{\cos C}}{a\cos A-b\cos B+a\cos A-c\cos C},		
\end{equation}
where
\[J=\frac{\sin 2A+\sin 2B+\sin 2C}{\tan A+\tan B+\tan C}.\]	

Without considering $J/2$ for the time being, the first part of the numerator in formula (\ref{Eq10.3.1}) is
	\[\begin{aligned}
	\frac{a}{\cos A}-\frac{b}{\cos B}&=\frac{2abc}{{{b}^{2}}+{{c}^{2}}-{{a}^{2}}}-\frac{2abc}{{{c}^{2}}+{{a}^{2}}-{{b}^{2}}} \\ 
	& =2abc\cdot \frac{{{c}^{2}}+{{a}^{2}}-{{b}^{2}}-\left( {{b}^{2}}+{{c}^{2}}-{{a}^{2}} \right)}{\left( {{b}^{2}}+{{c}^{2}}-{{a}^{2}} \right)\left( {{c}^{2}}+{{a}^{2}}-{{b}^{2}} \right)} \\ 
	& =4abc\cdot \frac{{{a}^{2}}-{{b}^{2}}}{\left( {{b}^{2}}+{{c}^{2}}-{{a}^{2}} \right)\left( {{c}^{2}}+{{a}^{2}}-{{b}^{2}} \right)}. \\ 
\end{aligned}\]

The second part of the numerator in formula (\ref{Eq10.3.1}) is
\[\begin{aligned}
	\frac{a}{\cos A}-\frac{c}{\cos C}&=\frac{2abc}{{{b}^{2}}+{{c}^{2}}-{{a}^{2}}}-\frac{2abc}{{{a}^{2}}+{{b}^{2}}-{{c}^{2}}} \\ 
	& =2abc\cdot \frac{{{a}^{2}}+{{b}^{2}}-{{c}^{2}}-\left( {{b}^{2}}+{{c}^{2}}-{{a}^{2}} \right)}{\left( {{b}^{2}}+{{c}^{2}}-{{a}^{2}} \right)\left( {{a}^{2}}+{{b}^{2}}-{{c}^{2}} \right)} \\ 
	& =4abc\cdot \frac{{{a}^{2}}-{{c}^{2}}}{\left( {{b}^{2}}+{{c}^{2}}-{{a}^{2}} \right)\left( {{a}^{2}}+{{b}^{2}}-{{c}^{2}} \right)}. \\ 
\end{aligned}\]

So the numerator in formula (\ref{Eq10.3.1}) is
\[\begin{aligned}
	&\frac{a}{\cos A}-\frac{b}{\cos B}+\frac{a}{\cos A}-\frac{c}{\cos C} \\ 
	& =\frac{4abc}{{{b}^{2}}+{{c}^{2}}-{{a}^{2}}}\left( \frac{{{a}^{2}}-{{b}^{2}}}{{{c}^{2}}+{{a}^{2}}-{{b}^{2}}}+\frac{{{a}^{2}}-{{c}^{2}}}{{{a}^{2}}+{{b}^{2}}-{{c}^{2}}} \right) \\ 
	& =K\left( \left( {{a}^{2}}-{{b}^{2}} \right)\left( {{a}^{2}}+{{b}^{2}}-{{c}^{2}} \right)+\left( {{a}^{2}}-{{c}^{2}} \right)\left( {{c}^{2}}+{{a}^{2}}-{{b}^{2}} \right) \right), \\ 
\end{aligned}\]
i.e.
\[\begin{aligned}
	& \frac{a}{\cos A}-\frac{b}{\cos B}+\frac{a}{\cos A}-\frac{c}{\cos C} \\ 
	& =K\left( {{a}^{4}}-{{b}^{4}}-{{c}^{2}}{{a}^{2}}+{{b}^{2}}{{c}^{2}}+{{a}^{4}}-{{c}^{4}}-{{a}^{2}}{{b}^{2}}+{{b}^{2}}{{c}^{2}} \right) \\ 
	& =K\left( 2{{a}^{4}}-{{b}^{4}}-{{c}^{4}}-{{a}^{2}}{{b}^{2}}+2{{b}^{2}}{{c}^{2}}-{{c}^{2}}{{a}^{2}} \right), \\ 
\end{aligned}\]
where
\[K=\frac{4abc}{\left( {{b}^{2}}+{{c}^{2}}-{{a}^{2}} \right)\left( {{c}^{2}}+{{a}^{2}}-{{b}^{2}} \right)\left( {{a}^{2}}+{{b}^{2}}-{{c}^{2}} \right)}.\]

The first part of the denominator in formula (\ref{Eq10.3.1}) is:
\[\begin{aligned}
	a\cos A-b\cos B & =\frac{a}{2bc}\left( {{b}^{2}}+{{c}^{2}}-{{a}^{2}} \right)-\frac{b}{2ca}\left( {{c}^{2}}+{{a}^{2}}-{{b}^{2}} \right) \\ 
	& =\frac{1}{2abc}\left( {{a}^{2}}\left( {{b}^{2}}+{{c}^{2}}-{{a}^{2}} \right)-{{b}^{2}}\left( {{c}^{2}}+{{a}^{2}}-{{b}^{2}} \right) \right) \\ 
	& =\frac{1}{2abc}\left( {{b}^{4}}-{{a}^{4}}+{{a}^{2}}{{c}^{2}}-{{b}^{2}}{{c}^{2}} \right). \\ 
\end{aligned}\]

The second part of the denominator in formula (\ref{Eq10.3.1}) is:
\[\begin{aligned}
	a\cos A-c\cos C & =\frac{a}{2bc}\left( {{b}^{2}}+{{c}^{2}}-{{a}^{2}} \right)-\frac{c}{2ab}\left( {{a}^{2}}+{{b}^{2}}-{{c}^{2}} \right) \\ 
	& =\frac{1}{2abc}\left( {{a}^{2}}\left( {{b}^{2}}+{{c}^{2}}-{{a}^{2}} \right)-{{c}^{2}}\left( {{a}^{2}}+{{b}^{2}}-{{c}^{2}} \right) \right) \\ 
	& =\frac{1}{2abc}\left( {{c}^{4}}-{{a}^{4}}+{{a}^{2}}{{b}^{2}}-{{b}^{2}}{{c}^{2}} \right). \\ 
\end{aligned}\]

The denominator in formula (\ref{Eq10.3.1}) is
\[\begin{aligned}
	& a\cos A-b\cos B+a\cos A-c\cos C \\ 
	& =\frac{1}{2abc}\left( {{b}^{4}}-{{a}^{4}}+{{a}^{2}}{{c}^{2}}-{{b}^{2}}{{c}^{2}}+{{c}^{4}}-{{a}^{4}}+{{a}^{2}}{{b}^{2}}-{{b}^{2}}{{c}^{2}} \right) \\ 
	& =\frac{1}{2abc}\left( {{b}^{4}}+{{c}^{4}}-2{{a}^{4}}+{{a}^{2}}{{b}^{2}}-2{{b}^{2}}{{c}^{2}}+{{c}^{2}}{{a}^{2}} \right). \\ 
\end{aligned}\]

Therefore
\[\frac{\alpha _{A}^{H}-\alpha _{A}^{G}}{\alpha _{A}^{Q}-\alpha _{A}^{G}}=\frac{J}{2}\cdot \frac{K\left( 2{{a}^{4}}-{{b}^{4}}-{{c}^{4}}-{{a}^{2}}{{b}^{2}}+2{{b}^{2}}{{c}^{2}}-{{c}^{2}}{{a}^{2}} \right)}{\frac{1}{2abc}\left( {{b}^{4}}+{{c}^{4}}-2{{a}^{4}}+{{a}^{2}}{{b}^{2}}-2{{b}^{2}}{{c}^{2}}+{{c}^{2}}{{a}^{2}} \right)}=-abcJK.\]		

The expression $-abcJK$ on the right side of the above formula has rotation invariance ($A$→$B$→$C$→$A$), indicating that the following formula holds:

\[\frac{\alpha _{A}^{H}-\alpha _{A}^{G}}{\alpha _{A}^{Q}-\alpha _{A}^{G}}=\frac{\alpha _{B}^{H}-\alpha _{B}^{G}}{\alpha _{B}^{Q}-\alpha _{B}^{G}}=\frac{\alpha _{C}^{H}-\alpha _{C}^{G}}{\alpha _{C}^{Q}-\alpha _{C}^{G}}=-abcJK,\]	
i.e.
\[\overrightarrow{GH}=-abcJK\overrightarrow{GQ}.\]	

This proves that the centroid is collinear with the orthocenter and the circumcenter.

In the following, we will calculate the relationship between $\overrightarrow{GH}$ and $\overrightarrow{GQ}$, namely, calculat the coefficient $-abcJK$. The ratio of the coefficients of $\overrightarrow{OA}$ in formula (\ref{Eq9.2.2}) and (\ref{Eq9.2.1}) is as follows:

\[\begin{aligned}
	\gamma &=\frac{\alpha _{A}^{H}-\alpha _{A}^{G}}{\alpha _{A}^{Q}-\alpha _{A}^{G}}=-abcJK \\ 
	& =-\frac{4{{a}^{2}}{{b}^{2}}{{c}^{2}}}{\left( {{b}^{2}}+{{c}^{2}}-{{a}^{2}} \right)\left( {{c}^{2}}+{{a}^{2}}-{{b}^{2}} \right)\left( {{a}^{2}}+{{b}^{2}}-{{c}^{2}} \right)}\frac{\sin 2A+\sin 2B+\sin 2C}{\tan A+\tan B+\tan C}. \\ 
\end{aligned}\]

Using cosine theorem, the following formula is obtained:
\[\begin{aligned}
	\gamma &=-\frac{1}{2}\frac{1}{\cos A\cos B\cos C}\frac{\sin 2A+\sin 2B+\sin 2C}{\tan A+\tan B+\tan C} \\ 
	& =-\frac{1}{2}\frac{\sin 2A+\sin 2B+\sin 2C}{\sin A\cos B\cos C+\sin B\cos C\cos A+\sin C\cos A\cos B} \\ 
	& =-\frac{1}{2}\frac{\sin 2A+\sin 2B+\sin 2C}{\left( \sin A\cos B+\sin B\cos A \right)\cos C+\sin C\cos A\cos B}. \\ 
\end{aligned}\]

By using trigonometric identity, we get the following result:
\[\begin{aligned}
	\gamma &=-\frac{1}{2}\frac{\sin 2A+\sin 2B+\sin 2C}{\sin \left( A+B \right)\cos C+\sin C\cos A\cos B}=-\frac{1}{2}\frac{\sin 2A+\sin 2B+\sin 2C}{\sin C\left( \cos C+\cos A\cos B \right)} \\ 
	& =-\frac{1}{2}\frac{\sin 2A+\sin 2B+\sin 2C}{\sin C\left( -\cos \left( A+B \right)+\cos A\cos B \right)}=-\frac{1}{2}\frac{\sin 2A+\sin 2B+\sin 2C}{\sin A\sin B\sin C}, \\ 
\end{aligned}\]
i.e.
\[\begin{aligned}
	\gamma &=-\frac{1}{2}\frac{2\sin \frac{2A+2B}{2}\cos \frac{2A-2B}{2}+2\sin C\cos C}{\sin A\sin B\sin C} \\ 
	& =-\frac{\sin \left( \pi -C \right)\cos \left( A-B \right)+\sin C\cos C}{\sin A\sin B\sin C} \\
	& =-\frac{\cos \left( A-B \right)+\cos C}{\sin A\sin B}, \\ 
\end{aligned}\]
i.e.
\[\begin{aligned}
	\gamma &=-\frac{2\cos \left( \frac{A-B+C}{2} \right)\cos \left( \frac{A-B-C}{2} \right)}{\sin A\sin B} \\ 
	& =-\frac{2\cos \left( \frac{\pi -2B}{2} \right)\cos \left( \frac{A-\left( \pi -A \right)}{2} \right)}{\sin A\sin B} \\ 
	& =-\frac{2\cos \left( \frac{\pi }{2}-B \right)\cos \left( A-\frac{\pi }{2} \right)}{\sin A\sin B}=-2. \\ 
\end{aligned}\]

Therefore
\begin{equation}\label{Eq10.3.2}
	\overrightarrow{GH}=-2\overrightarrow{GQ}
\end{equation}
\end{proof}
\hfill $\square$\par

The above formula not only shows that the centroid, the orthocenter and the circumcenter of a triangle are collinear, and the distance from the centroid to the orthocenter is twice that from the centroid to the circumcenter, but also shows that the centroid lies between the orthocenter and the circumcenter. The last point is not easy to make clear by the proof method of traditional comprehensive geometry.

\chapter{Frame equation of intersecting center}\label{Ch11}
\thispagestyle{empty}

In this chapter, we discuss the frame equations of some special intersecting centers of a triangle. According to the frame equation of single intersecting center (theorem \ref{Thm6.3.1}), the calculation formula of frame components and the intersecting ratios (see Chapter \ref{Ch5} and Chapter \ref{Ch8}), the frame equations of intersecting centers can be obtained.

\section{Frame equation of centroid of a triangle}\label{Sec11.1}
The FESIC of centroid $G$ is

\[\overrightarrow{GA}+\overrightarrow{GB}+\overrightarrow{GC}=\overrightarrow{0}.\]	

\section{Frame equation of incenter of a triangle}\label{Sec11.2}
The FESIC of incenter $I$ is

\[a\overrightarrow{IA}+b\overrightarrow{IB}+c\overrightarrow{IC}=\overrightarrow{0}.\]	

\section{Frame equation of orthocenter of a triangle}\label{Sec11.3}
The FESIC of orthocenter $H$ is

\[\frac{\tan A\cdot \overrightarrow{HA}+\tan B\cdot \overrightarrow{HB}+\tan C\cdot \overrightarrow{HC}}{\tan A+\tan B+\tan C}=\overrightarrow{0},\]	
or
\[\tan A\cdot \overrightarrow{HA}+\tan B\cdot \overrightarrow{HB}+\tan C\cdot \overrightarrow{HC}=\overrightarrow{0}.\]	

\section{Frame equation of circumcenter of a triangle}\label{Sec11.4}
The FESIC of circumcenter $Q$ is

\[\frac{\sin 2A\cdot \overrightarrow{QA}+\sin 2B\cdot \overrightarrow{QB}+\sin 2C\cdot \overrightarrow{QC}}{\sin 2A+\sin 2B+\sin 2C}=\overrightarrow{0},\]	
or
\[\sin 2A\cdot \overrightarrow{QA}+\sin 2B\cdot \overrightarrow{QB}+\sin 2C\cdot \overrightarrow{QC}=\overrightarrow{0}.\]

\chapter{Distance between two points on a plane}\label{Ch12}
\thispagestyle{empty}


The distance between two points is one of the basic geometric elements. The “natural defect” of Euclidean geometry makes it lack of quantitative method to deal with the distance between two points. In other words, What parameters must be known to determine the distance between two points? Euclidean geometry doesn't answer this question, which leads to the confusion that there is a lack of unified thinking and solution method, but usually one solution for one problem. In analytic geometry, the distance between two points can be uniquely determined by the coordinates of two points, which is determined by an algebraic formula representing the distance. However, the coordinate system of analytic geometry is set artificially, the point coordinate is related to the choice of coordinate system, and the component of point coordinate in the frame is related to the origin position of coordinate system. Moreover, because the basic element of analytic geometry is the coordinates of points, the geometric parameters of some graphics (such as distance, etc.) are measured by the coordinates of points, which makes the results of analytic geometry usually described by the coordinates of points, and it is difficult to establish a direct connection with the “natural” parameters of some graphics (for example, the lengths of the three sides of a triangle is the most natural parameter, which is the easiest to measure and has nothing to do with the position of the triangle). The reason why I put forward Intercenter Geometry is to overcome the “natural defect” of Euclidean geometry and the “unnatural” attribute of analytic geometry. This is well explained by the distance between origin and intersecting center, and the distance between two intersecting centers studied in this chapter. In the Intercenter Geometry proposed by the author, the component of the frame is independent of the position of the origin of the frame. In Plane Intercenter Geometry, the distance between two points is calculated by giving the lengths of three sides of a triangle and the frame components. Therefore, many results of Intercenter Geometry are expressed by the natural parameters (side length) of basic figures (triangles).

Using the method of vector algebra, I obtain some very important theorems about the distance between origin and intersecting center, and the distance between two intersecting centers, which lays the foundation for the following chapters.

\section{Theorem of distance between origin and intersecting center on triangular frame}\label{Sec12.1}
This section deals with an important theorem. Given a $\triangle ABC$ and an intersecting center $P$ located on the $\triangle ABC$ plane, select point $\, O$ as the origin of frame $\left(O;A, B, C \right)$, we can find the distance between two points $O$ and $P$, as shown in the following formula \ref{Eq12.1.1}. Since the location of point $O$ is arbitrary and does not necessarily lie on the plane of $\triangle ABC$, point $O$ is not necessarily an IC, and the distance between such two points is called \textbf{distance between origin and intersecting center} (abbreviated as DOIC). If the IC is selected at some special points in the triangle, some special conclusions will be obtained.

\begin{theorem}{Distance between origin and intersecting center on triangular frame, Daiyuan Zhang}{Thm12.1.1}\label{Thm12.1.1} 
	Let point $\,O$ be the the origin of frame $\left( O;A,B,C \right)$, and point $P$ be the IC on the plane of $\triangle ABC$, then	
	\begin{equation}\label{Eq12.1.1}
		\begin{aligned}
			O{{P}^{2}}=\alpha _{A}^{P}O{{A}^{2}}+\alpha _{B}^{P}O{{B}^{2}}+\alpha _{C}^{P}O{{C}^{2}}-\left( \alpha _{B}^{P}\alpha _{C}^{P}{{a}^{2}}+\alpha _{C}^{P}\alpha _{A}^{P}{{b}^{2}}+\alpha _{A}^{P}\alpha _{B}^{P}{{c}^{2}} \right). \\ 
		\end{aligned}		
	\end{equation}
	Where $\alpha _{A}^{P}$, $\alpha _{B}^{P}$, $\alpha _{C}^{P}$ are the frame components of $\overrightarrow{OA}$, $\overrightarrow{OB}$, $\overrightarrow{OC}$ at point $P$ on the frame $\left( O;A,B,C \right)$, respectively.	
\end{theorem}

\begin{proof}
	From theorem \ref{thm:Thm6.1.1}, the following formula is obtained:
	\[\overrightarrow{OP}=\alpha _{A}^{P}\overrightarrow{OA}+\alpha _{B}^{P}\overrightarrow{OB}+\alpha _{C}^{P}\overrightarrow{OC}.\]
		
	Thus, the distance from any point of $\, O$ to the intersecting center $P$ is:	
	\[O{{P}^{2}}=\left( \alpha _{A}^{P}\overrightarrow{OA}+\alpha _{B}^{P}\overrightarrow{OB}+\alpha _{C}^{P}\overrightarrow{OC} \right)\cdot \left( \alpha _{A}^{P}\overrightarrow{OA}+\alpha _{B}^{P}\overrightarrow{OB}+\alpha _{C}^{P}\overrightarrow{OC} \right),\]
	i.e.
	\[\begin{aligned}
		O{{P}^{2}}&={{\left( \alpha _{A}^{P} \right)}^{2}}O{{A}^{2}}+{{\left( \alpha _{B}^{P} \right)}^{2}}O{{B}^{2}}+{{\left( \alpha _{C}^{P} \right)}^{2}}O{{C}^{2}} \\ 
		& +2\alpha _{A}^{P}\alpha _{B}^{P}\overrightarrow{OA}\cdot \overrightarrow{OB}+2\alpha _{A}^{P}\alpha _{C}^{P}\overrightarrow{OA}\cdot \overrightarrow{OC}+2\alpha _{B}^{P}\alpha _{C}^{P}\overrightarrow{OB}\cdot \overrightarrow{OC}.  
	\end{aligned}\]
		
	Using the cosine theorem in vector form, the following formula is obtained:
	\[\begin{aligned}
		O{{P}^{2}}&={{\left( \alpha _{A}^{P} \right)}^{2}}O{{A}^{2}}+{{\left( \alpha _{B}^{P} \right)}^{2}}O{{B}^{2}}+{{\left( \alpha _{C}^{P} \right)}^{2}}O{{C}^{2}}+\alpha _{A}^{P}\alpha _{B}^{P}\left( O{{A}^{2}}+O{{B}^{2}}-A{{B}^{2}} \right) \\ 
		& +\alpha _{A}^{P}\alpha _{C}^{P}\left( O{{A}^{2}}+O{{C}^{2}}-A{{C}^{2}} \right)+\alpha _{B}^{P}\alpha _{C}^{P}\left( O{{B}^{2}}+O{{C}^{2}}-B{{C}^{2}} \right), \\ 
	\end{aligned}\]
	i.e.
	\[\begin{aligned}
		O{{P}^{2}}&=\alpha _{A}^{P}O{{A}^{2}}+\alpha _{B}^{P}O{{B}^{2}}+\alpha _{C}^{P}O{{C}^{2}}-\left( \alpha _{B}^{P}\alpha _{C}^{P}B{{C}^{2}}+\alpha _{C}^{P}\alpha _{A}^{P}C{{A}^{2}}+\alpha _{A}^{P}\alpha _{B}^{P}A{{B}^{2}} \right) \\ 
		& =\alpha _{A}^{P}O{{A}^{2}}+\alpha _{B}^{P}O{{B}^{2}}+\alpha _{C}^{P}O{{C}^{2}}-\left( \alpha _{B}^{P}\alpha _{C}^{P}{{a}^{2}}+\alpha _{C}^{P}\alpha _{A}^{P}{{b}^{2}}+\alpha _{A}^{P}\alpha _{B}^{P}{{c}^{2}} \right). \\ 
	\end{aligned}\]	
\end{proof}
\hfill $\square$\par

In formula \ref{Eq12.1.1}, the distance $OP$ is determined by the frame magnitudes (lengths) (i.e. $OA$, $OB$ and $OC$) and the frame components ($\alpha _{A}^{P}$, $\alpha _{B}^{P}$ and $\alpha _{C}^{P}$) and the lengths of the three sides of the triangle.

\begin{corollary}{Weighted sum of squares of distance from vertex to IC, Daiyuan Zhang}{Cor12.1.1}\label{Cor12.1.1} 
	If the point $P$ is the IC on the plane of $\triangle ABC$, then
	\[\alpha _{A}^{P}A{{P}^{2}}+\alpha _{B}^{P}B{{P}^{2}}+\alpha _{C}^{P}C{{P}^{2}}=\left( \alpha _{B}^{P}\alpha _{C}^{P}{{a}^{2}}+\alpha _{C}^{P}\alpha _{A}^{P}{{b}^{2}}+\alpha _{A}^{P}\alpha _{B}^{P}{{c}^{2}} \right).\]	
	Where $\alpha _{A}^{P}$, $\alpha _{B}^{P}$, $\alpha _{C}^{P}$ are the frame components of $\overrightarrow{OA}$, $\overrightarrow{OB}$, $\overrightarrow{OC}$ at point $P$ on the frame $\left( O;A,B,C \right)$, respectively.
\end{corollary}

\begin{proof}
	Let the point $\,O$ and $P$ coincide, and use the above theorem \ref{thm:Thm12.1.1} to get the result.
\end{proof}
\hfill $\square$\par

\subsection{Distance between vertex and intersecting center of a triangle}\label{Subsec12.1.1}
This section discusses the distance between vertex and intersecting center of a triangle (abbreviated as DVIC).

\begin{theorem}{Distance between vertex and intersecting center of a triangle, Daiyuan Zhang}{Thm12.1.2}\label{Thm12.1.2} 
	Given a $\triangle ABC$ and an intersecting center $P$ of the $\triangle ABC$, then
	\[A{{P}^{2}}=\left( 1-\alpha _{A}^{P} \right)\left( \alpha _{B}^{P}{{c}^{2}}+\alpha _{C}^{P}{{b}^{2}} \right)-\alpha _{B}^{P}\alpha _{C}^{P}{{a}^{2}},\]  
	\[B{{P}^{2}}=\left( 1-\alpha _{B}^{P} \right)\left( \alpha _{C}^{P}{{a}^{2}}+\alpha _{A}^{P}{{c}^{2}} \right)-\alpha _{C}^{P}\alpha _{A}^{P}{{b}^{2}},\]  
	\[C{{P}^{2}}=\left( 1-\alpha _{C}^{P} \right)\left( \alpha _{A}^{P}{{b}^{2}}+\alpha _{B}^{P}{{a}^{2}} \right)-\alpha _{A}^{P}\alpha _{B}^{P}{{c}^{2}}.\]	
	Where $\alpha _{A}^{P}$, $\alpha _{B}^{P}$, $\alpha _{C}^{P}$ are the frame components of $\overrightarrow{OA}$, $\overrightarrow{OB}$, $\overrightarrow{OC}$ at point $P$ on the frame $\left( O;A,B,C \right)$, respectively.
\end{theorem}

\begin{proof}
	Let the origin $O$ coincide with a vertex $A$ of $\triangle ABC$, and $P$ be the IC of $\triangle ABC$, then the square of the distance from the vertex $A$ of $\triangle ABC$ to the point $P$ is obtained by using formula \ref{Eq12.1.1}:
		
	\[\begin{aligned}
		A{{P}^{2}}&=\alpha _{B}^{P}A{{B}^{2}}+\alpha _{C}^{P}A{{C}^{2}}-\left( \alpha _{B}^{P}\alpha _{C}^{P}{{a}^{2}}+\alpha _{C}^{P}\alpha _{A}^{P}{{b}^{2}}+\alpha _{A}^{P}\alpha _{B}^{P}{{c}^{2}} \right) \\ 
		& =\alpha _{B}^{P}{{c}^{2}}+\alpha _{C}^{P}{{b}^{2}}-\left( \alpha _{B}^{P}\alpha _{C}^{P}{{a}^{2}}+\alpha _{C}^{P}\alpha _{A}^{P}{{b}^{2}}+\alpha _{A}^{P}\alpha _{B}^{P}{{c}^{2}} \right) \\ 
		& =\left( \alpha _{B}^{P}-\alpha _{A}^{P}\alpha _{B}^{P} \right){{c}^{2}}+\left( \alpha _{C}^{P}-\alpha _{A}^{P}\alpha _{C}^{P} \right){{b}^{2}}-\alpha _{B}^{P}\alpha _{C}^{P}{{a}^{2}} \\ 
		& =\left( 1-\alpha _{A}^{P} \right)\left( \alpha _{B}^{P}{{c}^{2}}+\alpha _{C}^{P}{{b}^{2}} \right)-\alpha _{B}^{P}\alpha _{C}^{P}{{a}^{2}}.  
	\end{aligned}\]
	
	Similar results are obtained for $B{{P}^{2}}$ and $C{{P}^{2}}$.
\end{proof}
\hfill $\square$\par

The intersecting center $P$ can be selected as some special points of the triangle, such as centroid, incenter, orthocenter, circumcenter, excenter, etc..


According to the above theory, we can directly calculate $A{{P}^{2}}+B{{P}^{2}}+C{{P}^{2}}$ and obtain the following corollary.


\begin{corollary}{Sum of squares for distances from vertexs to IC, Daiyuan Zhang}{Cor12.1.2}\label{Cor12.1.2} 
	Given a $\triangle ABC$ and an intersecting center $P$ of the $\triangle ABC$, then
	\begin{align*}
		A{{P}^{2}}+B{{P}^{2}}+C{{P}^{2}}& =\left( \alpha _{B}^{P}+\alpha _{C}^{P}-3\alpha _{B}^{P}\alpha _{C}^{P} \right){{a}^{2}} \\ 
		& +\left( \alpha _{C}^{P}+\alpha _{A}^{P}-3\alpha _{C}^{P}\alpha _{A}^{P} \right){{b}^{2}} \\ 
		& +\left( \alpha _{A}^{P}+\alpha _{B}^{P}-3\alpha _{A}^{P}\alpha _{B}^{P} \right){{c}^{2}}.  
	\end{align*}	
\end{corollary}

\begin{example}{(Example of calculation of the DVIC)}\label{Exm12.1.1} 
	Find the DVIC of the centroid.
\end{example}

\begin{solution}
	Each of IR of the centroid is 1, and the frame components are ${1}/{3}\;$. So the DVIC is	
	\[A{{G}^{2}}=\left( 1-\frac{1}{3} \right)\left( \frac{1}{3}{{c}^{2}}+\frac{1}{3}{{b}^{2}} \right)-\frac{1}{3}\times \frac{1}{3}{{a}^{2}}=\frac{1}{9}\left( 2{{c}^{2}}+2{{b}^{2}}-{{a}^{2}} \right),\]
	i.e.
	\[AG=\frac{1}{3}\sqrt{2{{b}^{2}}+2{{c}^{2}}-{{a}^{2}}}.\]
	
	Similarly, 
	\[BG=\frac{1}{3}\sqrt{2{{c}^{2}}+2{{a}^{2}}-{{b}^{2}}},\]
	\[CG=\frac{1}{3}\sqrt{2{{a}^{2}}+2{{b}^{2}}-{{c}^{2}}}.\]
\end{solution}
\hfill $\diamond$\par

For some other centers of the triangle (such as incenter, circumcenter, orthocenter, etc.) the DVIC can also be similarly calculated.

\subsection{Distance between vertex and intersecting foot}\label{Subsec12.1.2}
The distance from the vertex of a triangle to the intersecting foot is called the distance between vertex and intersecting foot (abbreviated as DVIF).

\begin{theorem}{Distance between vertex and intersecting foot, Daiyuan Zhang}{Thm12.1.3}\label{Thm12.1.3} 
	Given a $\triangle ABC$, the intersecting center $P$ and the intersecting foot $L$, $M$ and $N$, then each of the distance between vertex and intersecting foot is
		
	\[AL=\frac{AP}{\left| 1-\alpha _{A}^{P} \right|}=\frac{\sqrt{\left( 1-\alpha _{A}^{P} \right)\left( \alpha _{B}^{P}{{c}^{2}}+\alpha _{C}^{P}{{b}^{2}} \right)-\alpha _{B}^{P}\alpha _{C}^{P}{{a}^{2}}}}{\left| 1-\alpha _{A}^{P} \right|},\]
	\[BM=\frac{BP}{\left| 1-\alpha _{B}^{P} \right|}=\frac{\sqrt{\left( 1-\alpha _{B}^{P} \right)\left( \alpha _{C}^{P}{{a}^{2}}+\alpha _{A}^{P}{{c}^{2}} \right)-\alpha _{C}^{P}\alpha _{A}^{P}{{b}^{2}}}}{\left| 1-\alpha _{B}^{P} \right|},\]
	\[CN=\frac{CP}{\left| 1-\alpha _{C}^{P} \right|}=\frac{\sqrt{\left( 1-\alpha _{C}^{P} \right)\left( \alpha _{A}^{P}{{b}^{2}}+\alpha _{B}^{P}{{a}^{2}} \right)-\alpha _{A}^{P}\alpha _{B}^{P}{{c}^{2}}}}{\left| 1-\alpha _{C}^{P} \right|}.\]
	Where $\alpha _{A}^{P}$, $\alpha _{B}^{P}$, $\alpha _{C}^{P}$ are the frame components of $\overrightarrow{OA}$, $\overrightarrow{OB}$, $\overrightarrow{OC}$ at point $P$ respectively on the frame $\left( O;A,B,C \right)$.
\end{theorem}

\begin{proof}
	According to theorem \ref{thm:Thm6.1.4}, this theorem can be directly obtained by taking the magnitude on both sides of the vector expression and using theorem \ref{thm:Thm12.1.2}.
\end{proof}
\hfill $\square$\par

\begin{example}{(Example of calculation of the DVIF)}\label{Exm12.1.2} 
	Find the DVIF of centroid.
\end{example}

\begin{solution}
	For the centroid, $\lambda _{AB}^{G}=\lambda _{BC}^{G}=\lambda _{CA}^{G}=1$, and the frame components are ${1}/{3}\;$, From the previous example, the following result is obtained:
	\[AL=\frac{AG}{\left| 1-\alpha _{A}^{P} \right|}=\frac{AG}{\left| 1-\frac{1}{3} \right|}=\frac{3AG}{2}=\frac{1}{2}\sqrt{2{{b}^{2}}+2{{c}^{2}}-{{a}^{2}}}.\]
	
	Similarly,
	\[BM=\frac{3BG}{2}=\frac{1}{2}\sqrt{2{{c}^{2}}+2{{a}^{2}}-{{b}^{2}}},\]	
	\[CN=\frac{3CG}{2}=\frac{1}{2}\sqrt{2{{a}^{2}}+2{{b}^{2}}-{{c}^{2}}}.\]	
\end{solution}
\hfill $\diamond$\par

By the way, the following formula is obtained:
\[A{{L}^{2}}+B{{M}^{2}}+C{{N}^{2}}=\frac{3}{4}\left( {{a}^{2}}+{{b}^{2}}+{{c}^{2}} \right).\]	

That is, the sum of the squares of the three medians of a triangle is equal to three-quarters of the sum of the squares of the three sides.

\subsection{Distance between origin and intersecting center on frame of circumcenter}\label{Subsec12.1.3}

\begin{theorem}{Distance between origin and IC on frame of circumcenter, Daiyuan Zhang}{Thm12.1.4}\label{Thm12.1.4}
	Let the circumcenter $Q$ of $\triangle ABC$ be the origin of the frame $\left( Q;A,B,C \right)$, $P$ is the IC on the plane of $\triangle ABC$, $R$ is the radius of circumscribed circle of $\triangle ABC$, then
	
	\begin{equation}\label{Eq12.1.2}
		Q{{P}^{2}}={{R}^{2}}-\left( \alpha _{B}^{P}\alpha _{C}^{P}{{a}^{2}}+\alpha _{C}^{P}\alpha _{A}^{P}{{b}^{2}}+\alpha _{A}^{P}\alpha _{B}^{P}{{c}^{2}} \right).	
	\end{equation}
	Where $\alpha _{A}^{P}$, $\alpha _{B}^{P}$, $\alpha _{C}^{P}$ are the frame components of $\overrightarrow{OA}$, $\overrightarrow{OB}$, $\overrightarrow{OC}$ at point $P$ respectively on the frame $\left( O;A,B,C \right)$.
\end{theorem}
	
\begin{proof}
	Since $QA=QB=QC=R$, in theorem \ref{thm:Thm12.1.1}, if the origin of the frame $\,O$ coincides with the circumcenter $Q$ of $\triangle ABC$, the desired result can be obtained.
\end{proof}
\hfill $\square$\par

\section{Theorem of distance between two intersecting centers}\label{Sec12.2}
The \textbf{distance between two intersecting centers} (abbreviated as DTICs) of a triangle refers to the distance between two intersecting centers on the same triangular frame, and is also called the distance of two intersecting centers. The following theorem is one of the core theorems of this book.

\begin{theorem}{Distance between two intersecting centers, Daiyuan Zhang}{Thm12.2.1}\label{Thm12.2.1} 
	Assume a given $\triangle ABC$, point $O$ is any point of space, ${{P}_{1}}\in {{\pi }_{ABC}}$, ${{P}_{2}}\in {{\pi }_{ABC}}$, then
				
	\[{{P}_{1}}{{P}_{2}}=\sqrt{-\alpha _{B}^{{{P}_{1}}{{P}_{2}}}\alpha _{C}^{{{P}_{1}}{{P}_{2}}}{{a}^{2}}-\alpha _{C}^{{{P}_{1}}{{P}_{2}}}\alpha _{A}^{{{P}_{1}}{{P}_{2}}}{{b}^{2}}-\alpha _{A}^{{{P}_{1}}{{P}_{2}}}\alpha _{B}^{{{P}_{1}}{{P}_{2}}}{{c}^{2}}},\]	
	\[\alpha _{A}^{{{P}_{1}}{{P}_{2}}}+\alpha _{B}^{{{P}_{1}}{{P}_{2}}}+\alpha _{C}^{{{P}_{1}}{{P}_{2}}}=0.\]	
	Where
	\[\alpha _{A}^{{{P}_{1}}{{P}_{2}}}=\alpha _{A}^{{{P}_{2}}}-\alpha _{A}^{{{P}_{1}}},\]	
	\[\alpha _{B}^{{{P}_{1}}{{P}_{2}}}=\alpha _{B}^{{{P}_{2}}}-\alpha _{B}^{{{P}_{1}}},\]	
	\[\alpha _{C}^{{{P}_{1}}{{P}_{2}}}=\alpha _{C}^{{{P}_{2}}}-\alpha _{C}^{{{P}_{1}}}.\]	
	Where $\alpha _{A}^{{{P}_{1}}}$, $\alpha _{B}^{{{P}_{1}}}$, $\alpha _{C}^{{{P}_{1}}}$ are the frame components of the frame $\overrightarrow{OA}$, $\overrightarrow{OB}$, $\overrightarrow{OC}$ on the frame system $\left( O;A,B,C \right)$ at point ${P}_{1}$, respectively; $\alpha _{A}^{{{P}_{2}}}$, $\alpha _{B}^{{{P}_{2}}}$, $\alpha _{C}^{{{P}_{2}}}$ are the frame components of the frame $\overrightarrow{OA}$, $\overrightarrow{OB}$, $\overrightarrow{OC}$ on the frame system $\left( O;A,B,C \right)$ at point ${P}_{2}$, respectively.
\end{theorem}

\begin{proof}
	From the inner product operation of the vectors and from equations \ref{Eq6.2.1} and \ref{Eq6.2.5} in theorem \ref{thm:Thm6.2.2}, the following formula is obtained:	
	\[\begin{aligned}
		{{P}_{1}}{{P}_{2}}^{2}&=\left( \alpha _{B}^{{{P}_{1}}{{P}_{2}}}\overrightarrow{AB}+\alpha _{C}^{{{P}_{1}}{{P}_{2}}}\overrightarrow{AC} \right)\cdot \left( \alpha _{B}^{{{P}_{1}}{{P}_{2}}}\overrightarrow{AB}+\alpha _{C}^{{{P}_{1}}{{P}_{2}}}\overrightarrow{AC} \right) \\ 
		& ={{\left( \alpha _{B}^{{{P}_{1}}{{P}_{2}}} \right)}^{2}}A{{B}^{2}}+{{\left( \alpha _{C}^{{{P}_{1}}{{P}_{2}}} \right)}^{2}}A{{C}^{2}}+2\alpha _{B}^{{{P}_{1}}{{P}_{2}}}\alpha _{C}^{{{P}_{1}}{{P}_{2}}}\overrightarrow{AB}\cdot \overrightarrow{AC}.  
	\end{aligned}\]
	
	Using the cosine theorem in vector form, the following formula is obtained:
	\[\begin{aligned}
		{{P}_{1}}{{P}_{2}}^{2}&={{\left( \alpha _{B}^{{{P}_{1}}{{P}_{2}}} \right)}^{2}}{{c}^{2}}+{{\left( \alpha _{C}^{{{P}_{1}}{{P}_{2}}} \right)}^{2}}{{b}^{2}}+\alpha _{B}^{{{P}_{1}}{{P}_{2}}}\alpha _{C}^{{{P}_{1}}{{P}_{2}}}\left( {{c}^{2}}+{{b}^{2}}-{{a}^{2}} \right) \\ 
		& =\left( {{\left( \alpha _{C}^{{{P}_{1}}{{P}_{2}}} \right)}^{2}}+\alpha _{B}^{{{P}_{1}}{{P}_{2}}}\alpha _{C}^{{{P}_{1}}{{P}_{2}}} \right){{b}^{2}}+\left( {{\left( \alpha _{B}^{{{P}_{1}}{{P}_{2}}} \right)}^{2}}+\alpha _{B}^{{{P}_{1}}{{P}_{2}}}\alpha _{C}^{{{P}_{1}}{{P}_{2}}} \right){{c}^{2}}-\alpha _{B}^{{{P}_{1}}{{P}_{2}}}\alpha _{C}^{{{P}_{1}}{{P}_{2}}}{{a}^{2}},  
	\end{aligned}\]
	i.e.
	\[\begin{aligned}
		{{P}_{1}}{{P}_{2}}^{2}&=\alpha _{C}^{{{P}_{1}}{{P}_{2}}}\left( \alpha _{C}^{{{P}_{1}}{{P}_{2}}}+\alpha _{B}^{{{P}_{1}}{{P}_{2}}} \right){{b}^{2}}+\alpha _{B}^{{{P}_{1}}{{P}_{2}}}\left( \alpha _{B}^{{{P}_{1}}{{P}_{2}}}+\alpha _{C}^{{{P}_{1}}{{P}_{2}}} \right){{c}^{2}}-\alpha _{B}^{{{P}_{1}}{{P}_{2}}}\alpha _{C}^{{{P}_{1}}{{P}_{2}}}{{a}^{2}} \\ 
		& =\left( \alpha _{C}^{{{P}_{1}}{{P}_{2}}}+\alpha _{B}^{{{P}_{1}}{{P}_{2}}} \right)\left( \alpha _{C}^{{{P}_{1}}{{P}_{2}}}{{b}^{2}}+\alpha _{B}^{{{P}_{1}}{{P}_{2}}}{{c}^{2}} \right)-\alpha _{B}^{{{P}_{1}}{{P}_{2}}}\alpha _{C}^{{{P}_{1}}{{P}_{2}}}{{a}^{2}},  
	\end{aligned}\]
	i.e.
	\[\begin{aligned}
		{{P}_{1}}{{P}_{2}}^{2}&=-\alpha _{A}^{{{P}_{1}}{{P}_{2}}}\left( \alpha _{C}^{{{P}_{1}}{{P}_{2}}}{{b}^{2}}+\alpha _{B}^{{{P}_{1}}{{P}_{2}}}{{c}^{2}} \right)-\alpha _{B}^{{{P}_{1}}{{P}_{2}}}\alpha _{C}^{{{P}_{1}}{{P}_{2}}}{{a}^{2}} \\ 
		& =-\alpha _{B}^{{{P}_{1}}{{P}_{2}}}\alpha _{C}^{{{P}_{1}}{{P}_{2}}}{{a}^{2}}-\alpha _{C}^{{{P}_{1}}{{P}_{2}}}\alpha _{A}^{{{P}_{1}}{{P}_{2}}}{{b}^{2}}-\alpha _{A}^{{{P}_{1}}{{P}_{2}}}\alpha _{B}^{{{P}_{1}}{{P}_{2}}}{{c}^{2}},  
	\end{aligned}\]
i.e.
	\[{{P}_{1}}{{P}_{2}}=\sqrt{-\alpha _{B}^{{{P}_{1}}{{P}_{2}}}\alpha _{C}^{{{P}_{1}}{{P}_{2}}}{{a}^{2}}-\alpha _{C}^{{{P}_{1}}{{P}_{2}}}\alpha _{A}^{{{P}_{1}}{{P}_{2}}}{{b}^{2}}-\alpha _{A}^{{{P}_{1}}{{P}_{2}}}\alpha _{B}^{{{P}_{1}}{{P}_{2}}}{{c}^{2}}}.\]	
	And
	\[\begin{aligned}
		\alpha _{A}^{{{P}_{1}}{{P}_{2}}}+\alpha _{B}^{{{P}_{1}}{{P}_{2}}}+\alpha _{C}^{{{P}_{1}}{{P}_{2}}}&=\left( \alpha _{A}^{{{P}_{2}}}-\alpha _{A}^{{{P}_{1}}} \right)+\left( \alpha _{B}^{{{P}_{2}}}-\alpha _{B}^{{{P}_{1}}} \right)+\left( \alpha _{C}^{{{P}_{2}}}-\alpha _{C}^{{{P}_{1}}} \right) \\ 
		& =\left( \alpha _{A}^{{{P}_{2}}}+\alpha _{B}^{{{P}_{2}}}+\alpha _{C}^{{{P}_{2}}} \right)-\left( \alpha _{A}^{{{P}_{1}}}+\alpha _{B}^{{{P}_{1}}}+\alpha _{C}^{{{P}_{1}}} \right)=1-1=0. \\ 
	\end{aligned}\]	
\end{proof}
\hfill $\square$\par


The above theorem shows that if we know the lengths of three sides of a triangle and the frame components (or IRs) of two ICs: ${{P}_{1}}$ and ${{P}_{2}}$, then the length of the vector $\overrightarrow{{{P}_{1}}{{P}_{2}}}$ can be calculated, i.e., the distance between the points ${{P}_{1}}$ and ${{P}_{2}}$ can be calculated.

Let's compare the difference of the distance between two points in Euclidean geometry, analytic geometry and Intercenter Geometry.

In Euclidean geometry, there is no corresponding relationship between a point and a real number, so there is no algebraic formula to describe the distance between two points. Although the distance between two points in some problems can be obtained by means of auxiliary lines, rotation transformation and other methods, it is usually one problem with one strategy, and there is no unified method.


In analytic geometry, if the coordinates of two points in Cartesian coordinate system are ${{P}_{1}}\left( {{x}_{1}},{{y}_{1}} \right)$ and ${{P}_{2}}\left( {{x}_{2}},{{y}_{2}} \right)$, then the distance between the two points is
\[{{P}_{ 1}}{{P}_{ 2}}=\sqrt{{{\left( {{x}_{ 1}}-{{x}_{ 2}} \right)}^{2}}+{{\left( {{y}_{ 1}}-{{y}_{ 2}} \right)}^{2}}}.\] 	 	


The coordinates of two points ${{x}_{1}}$, ${{y}_{1}}$ and ${{x}_{2}}$, ${{y}_{2}}$ usually changes with the origin of the coordinate. In other words, if you translate or rotate the coordinate system, then ${{x}_{1}}$, ${{y}_{1}}$ and ${{x}_{2}}$, ${{y}_{2}}$ usually changes, that is, the parameters in the above formula of distance ${{P}_{1}}{{P}_{2}}$ will be various with the change of coordinate system or coordinate origin.


For example, if you want to calculate the distance between the centroid and the incenter of $\triangle ABC$, you need to establish an appropriate coordinate system, then calculate the coordinates of the centroid and the incenter, and then calculate the distance between the two points according to the above formula. At this time, the parameters in the distance formula are the coordinates of centroid and incenter, which are not only related to the choice of the coordinate system, but also related to the origin position of the coordinate system. In particular, it is difficult to establish a direct relationship between the coordinates of the distance formula and the most “natural” parameters of the triangle (side lengths).

In Intercenter Geometry, the distance between two points is given by the following formula (see the previous theorem):
\[{{P}_{1}}{{P}_{2}}=\sqrt{-\alpha _{B}^{{{P}_{1}}{{P}_{2}}}\alpha _{C}^{{{P}_{1}}{{P}_{2}}}{{a}^{2}}-\alpha _{C}^{{{P}_{1}}{{P}_{2}}}\alpha _{A}^{{{P}_{1}}{{P}_{2}}}{{b}^{2}}-\alpha _{A}^{{{P}_{1}}{{P}_{2}}}\alpha _{B}^{{{P}_{1}}{{P}_{2}}}{{c}^{2}}}.\]	


The above formula shows that the distance between two points is expressed by the lengths of the three sides of the triangle and the frame components. Obviously, the lengths of the three sides of the triangle has nothing to do with the position of the origin of the frame. For a given $\triangle ABC$ and two given points ${{P}_{1}}$ and ${{P}_{2}}$ on the $\triangle ABC$ plane, each of the frame component is only related to the IRs of point ${{P}_{1}}$ and ${{P}_{2}}$, and is not related to the position of the origin of the frame. In other words, if the given points ${{P}_{1}}$, ${{P}_{2}}$ are fixed with the position of $\triangle ABC$, then, no matter how to translate, rotate or even flip the $\triangle ABC$ in space, each quantity in the formula of the distance ${{P}_{1}}{{P}_{2}}$ will keep unchanged, which brings great convenience for the analysis and application of some problems, because the origin of the frame can be selected according to the needs without changing the parameters in the formula. We will gradually realize this superiority later.


For example, if you want to calculate the distance between the centroid and the incenter of $\triangle ABC$, you only need to calculate the IRs of the centroid and the incenter (see Chapter \ref{Ch5}), then use the IRs to calculate the frame components (see Section \ref{Sec8.1} and \ref{Sec8.2}), and then directly calculate the distance between the two points according to the above formula. At this time, the parameters in the distance formula are independent of the origin of the frame, and include the most “natural” parameters of the triangle: the lengths of three sides.


\chapter{Distance between special intersecting centers}\label{Ch13}
\thispagestyle{empty}

This chapter is an application of the method of Intercenter Geometry.


A triangle has some special ICs, which are the circumcenter, centroid, incenter, orthocenter, and three excenters. There are seven such special ICs, which are called the seven centers for a triangle. This chapter studies the distance between those special ICs. Intercenter Geometry uses a unified method to solve those distances. Some of those distances are directly expressed by the most natural parameter of the triangle: lengths of three sides. Some distance expressions contain trigonometric functions of the internal angle of the triangle. Those formulas are usually symmetrical and graceful, reflecting the intrinsic information of the triangle. As we all know, trigonometric function can be transformed into cosine function by trigonometric identity, and then expressed as the function of the lengths of three sides of the triangle according to the cosine theorem, so in essence, it is the function of the lengths of the three sides of triangle. This is difficult to achieve by using analytic geometry method. The result obtained by analytic geometry method is usually the function of point coordinates, not the function of lengths of three sides. The functions of those points often lose the important internal information of the triangle (such as lengths of three sides, angles, etc.). Although the function of point coordinates can be transformed into the function of distance by transformation, the calculation is usually lack of guidance, difficult and very complicated. The biggest disadvantage of calculating the distance between those special ICs by Euclidean geometry is that there is no unified solving strategy. Because of the natural defect of Euclidean geometry, Euclidean geometry often needs to draw different auxiliary lines for different problems, and often needs different solving strategies for different problems. There is no unified method, only one solution for one problem. What does that mean? I believe this is a scene that no mathematician would like to see.

In this chapter, the distance formulas between the seven special centers of a triangle will be obtained by using the method of Intercenter Geometry. 


In Intercenter Geometry, the distance between two points has special terms, that is, the distance between origin and intersecting center (abbreviated as DOIC) and the distance between two intersecting centers (abbreviated as DTICs) (see the previous chapter). There are some connections and differences between them. They all calculate the distance between two points. The DOIC (see theorem \ref{thm:Thm12.1.1}) refers to the distance between the origin of the frame and an IC, that is, one point must be the origin of the frame, and the other point must be an IC on the triangular plane. The DOIC emphasizes the freedom of the origin of the frame, that is, the origin of the frame can be any point in space. The DOIC is related to the lengths of the frame. When the origin of the frame is also an IC on the triangle plane, the distance between the two ICs is calculated. On the surface, it seems to be the distance between the two ICs, but in fact it is still the distance between origin and intersecting center, because the distance is calculated according to theorem \ref{thm:Thm12.1.1}, the distance is related to the frame length, and the formula of distance only contains the frame components of one IC. The distance between two ICs refers to the distance between two ICs on the triangular plane. The distance between the two ICs has nothing to do with the length of frame. In the formula of distance between two ICs (see theorem \ref{thm:Thm12.2.1}), the formula of distance contains the frame components of two ICs.


Because the distance between the circumcenter and the three vertices of the triangle is equal, the circumcenter is special. When discussing the distance between the circumcenter and other special ICs, the formulas for the DOIC and the DTICs are given. When discussing the distance between other special ICs, only the formulas for the DTICs are given. Interested readers can deduce the formulas for the DOIC by themselves.

In order to simplify the writing process, most results only give the square of the distance.


This book uses the following notations to represent the seven centers of a triangle: circumcenter: $Q$; centroid: $G$; incenter: $I$; orthocenter: $H$; excenter corresponding to $\angle A$: ${{E}_{A}}$, excenter corresponding to $\angle B$: ${{E}_{B}}$, excenter corresponding to $\angle C$: ${{E}_{C}}$.   

Some geometric quantities are represented by the following notations:
\[p=\frac{1}{2}\left( a+b+c \right),\]
\[W=\sin 2A+\sin 2B+\sin 2C,\]
\[T=\tan A+\tan B+\tan C.\]

\section{Distance between origin and intersecting center on frame of circumcenter}\label{Sec13.1}
This section solves the distance between the circumcenter and centroid, incenter, orthocenter and excenter of a triangle. The calculation formula is formula (\ref{Eq12.1.2}), which is a direct inference of theorem \ref{thm:Thm12.1.1} on the triangular frame. Formula (\ref{Eq12.1.2}) contains the circumscribed radius $R$ of $\triangle ABC$.

\subsection{DOIC between circumcenter and centroid}\label{Subsec13.1.1}
From Chapter \ref{Ch8} and formula (\ref{Eq12.1.2}), we have
\[\begin{aligned}
	Q{{G}^{2}}&=\overrightarrow{QG}\cdot \overrightarrow{QG} \\ 
	& ={{R}^{2}}-\left( \alpha _{B}^{G}\alpha _{C}^{G}{{a}^{2}}+\alpha _{C}^{G}\alpha _{A}^{G}{{b}^{2}}+\alpha _{A}^{G}\alpha _{B}^{I}{{c}^{2}} \right) \\ 
	& ={{R}^{2}}-\frac{1}{9}\left( {{a}^{2}}+{{b}^{2}}+{{c}^{2}} \right).  
\end{aligned}\]

\subsection{DOIC between circumcenter and incenter}\label{Subsec13.1.2}
From Chapter \ref{Ch8} and formula (\ref{Eq12.1.2}), we have
\[\alpha _{B}^{I}\alpha _{C}^{I}=\frac{b}{2p}\cdot \frac{c}{2p}=\frac{bc}{4{{p}^{2}}},\]
\[\alpha _{C}^{I}\alpha _{A}^{I}=\frac{c}{2p}\cdot \frac{a}{2p}=\frac{ca}{4{{p}^{2}}},\]
\[\alpha _{A}^{I}\alpha _{B}^{I}=\frac{a}{2p}\cdot \frac{b}{2p}=\frac{ab}{4{{p}^{2}}}.\]

Therefore,
\[Q{{I}^{2}}=\overrightarrow{QI}\cdot \overrightarrow{QI}  
	={{R}^{2}}-\frac{1}{4{{p}^{2}}}\left( ab{{c}^{2}}+bc{{a}^{2}}+ca{{b}^{2}} \right) 
	={{R}^{2}}-\frac{abc}{2p}.\]

\subsection{DOIC between circumcenter and orthocenter}\label{Subsec13.1.3}
From Chapter \ref{Ch8} and formula (\ref{Eq12.1.2}), we have
\begin{flalign*}
	\alpha _{A}^{H}=\frac{\tan A}{T}, \alpha _{B}^{H}=\frac{\tan B}{T}, \alpha _{C}^{H}=\frac{\tan C}{T}
\end{flalign*}
\[\begin{aligned}
	Q{{H}^{2}}&=\overrightarrow{QH}\cdot \overrightarrow{QH}={{R}^{2}}-\left( \alpha _{B}^{H}\alpha _{C}^{H}{{a}^{2}}+\alpha _{C}^{H}\alpha _{A}^{H}{{b}^{2}}+\alpha _{A}^{H}\alpha _{B}^{H}{{c}^{2}} \right) \\ 
	& ={{R}^{2}}-\frac{\tan B\cdot \tan C\cdot {{a}^{2}}+\tan C\cdot \tan A\cdot {{b}^{2}}+\tan A\cdot \tan B\cdot {{c}^{2}}}{{{\left( \tan A+\tan B+\tan C \right)}^{2}}}.  
\end{aligned}\]	

\subsection{DOIC between circumcenter and excenter}\label{Subsec13.1.4}
\subsubsection{DOIC between circumcenter and excenter ${{E}_{A}}$ corresponding to $\angle A$} 
From Chapter \ref{Ch8} and formula (\ref{Eq12.1.2}), we have
\begin{flalign*}
	\alpha _{A}^{{{E}_{A}}}=-\frac{a}{2\left( p-a \right)}, \alpha _{B}^{{{E}_{A}}}=\frac{b}{2\left( p-a \right)}, \alpha _{C}^{{{E}_{A}}}=\frac{c}{2\left( p-a \right)}.
\end{flalign*}
\[\begin{aligned}
	Q{{E}_{A}}^{2}&={{R}^{2}}-\left( \alpha _{B}^{{{E}_{A}}}\alpha _{C}^{{{E}_{A}}}{{a}^{2}}+\alpha _{C}^{{{E}_{A}}}\alpha _{A}^{{{E}_{A}}}{{b}^{2}}+\alpha _{A}^{{{E}_{A}}}\alpha _{B}^{{{E}_{A}}}{{c}^{2}} \right) \\ 
	& ={{R}^{2}}-\frac{bc{{a}^{2}}-ca{{b}^{2}}-ab{{c}^{2}}}{4{{\left( p-a \right)}^{2}}}={{R}^{2}}+\frac{abc}{2\left( p-a \right)}.  
\end{aligned}\]	
\subsubsection{DOIC between circumcenter and excenter ${{E}_{B}}$ corresponding to $\angle B$}
From Chapter \ref{Ch8} and formula (\ref{Eq12.1.2}), we have
\begin{flalign*}
	\alpha _{A}^{{{E}_{B}}}=\frac{a}{2\left( p-b \right)}, \alpha _{B}^{{{E}_{B}}}=-\frac{b}{2\left( p-b \right)}, \alpha _{C}^{{{E}_{B}}}=\frac{c}{2\left( p-b \right)}.
\end{flalign*}
\[\begin{aligned}
	Q{{E}_{B}}^{2}&=\overrightarrow{Q{{E}_{B}}}\cdot \overrightarrow{Q{{E}_{B}}}={{R}^{2}}-\left( \alpha _{B}^{{{E}_{B}}}\alpha _{C}^{{{E}_{B}}}{{a}^{2}}+\alpha _{C}^{{{E}_{B}}}\alpha _{A}^{{{E}_{B}}}{{b}^{2}}+\alpha _{A}^{{{E}_{B}}}\alpha _{B}^{{{E}_{B}}}{{c}^{2}} \right) \\ 
	& ={{R}^{2}}+\frac{bc{{a}^{2}}}{4{{\left( p-b \right)}^{2}}}-\frac{ca{{b}^{2}}}{4{{\left( p-b \right)}^{2}}}+\frac{ab{{c}^{2}}}{4{{\left( p-b \right)}^{2}}}={{R}^{2}}+\frac{abc}{2\left( p-b \right)}.  
\end{aligned}\]	
\subsubsection{DOIC between circumcenter and excenter ${{E}_{C}}$ corresponding to $\angle C$}
From Chapter \ref{Ch8} and formula (\ref{Eq12.1.2}), we have
\begin{flalign*}
	\alpha _{A}^{{{E}_{C}}}=\frac{a}{2\left( p-c \right)}, \alpha _{B}^{{{E}_{C}}}=\frac{b}{2\left( p-c \right)}, \alpha _{C}^{{{E}_{C}}}=-\frac{c}{2\left( p-c \right)}.	
\end{flalign*}
\[\begin{aligned}
	Q{{E}_{C}}^{2}&=\overrightarrow{Q{{E}_{C}}}\cdot \overrightarrow{Q{{E}_{C}}}={{R}^{2}}-\left( \alpha _{B}^{{{E}_{C}}}\alpha _{C}^{{{E}_{C}}}{{a}^{2}}+\alpha _{C}^{{{E}_{C}}}\alpha _{A}^{{{E}_{C}}}{{b}^{2}}+\alpha _{A}^{{{E}_{C}}}\alpha _{B}^{{{E}_{C}}}{{c}^{2}} \right) \\ 
	& ={{R}^{2}}+\frac{bc{{a}^{2}}}{4{{\left( p-c \right)}^{2}}}+\frac{ca{{b}^{2}}}{4{{\left( p-c \right)}^{2}}}-\frac{ab{{c}^{2}}}{4{{\left( p-c \right)}^{2}}}={{R}^{2}}+\frac{abc}{2\left( p-c \right)}.  
\end{aligned}\]

\section{Distance between two intersecting centers on triangular frame}\label{Sec13.2}
\subsection{DTICs between circumcenter and centroid}\label{Subsec13.2.1}
On the other hand, by using Chapter \ref{Ch8} and theorem \ref{thm:Thm12.2.1}, the following formulas are obtained directly:
\[\alpha _{A}^{QG}=\alpha _{A}^{G}-\alpha _{A}^{Q}=\frac{1}{3}-\frac{\sin 2A}{W}=\frac{W-3\sin 2A}{3W},\]	
\[\alpha _{B}^{QG}=\alpha _{B}^{G}-\alpha _{B}^{Q}=\frac{1}{3}-\frac{\sin 2B}{W}=\frac{W-3\sin 2B}{3W},\]	
\[\alpha _{C}^{QG}=\alpha _{C}^{G}-\alpha _{C}^{Q}=\frac{1}{3}-\frac{\sin 2C}{W}=\frac{W-3\sin 2C}{3W},\]	
\[Q{{G}^{2}}=-\alpha _{B}^{QG}\alpha _{C}^{QG}{{a}^{2}}-\alpha _{C}^{QG}\alpha _{A}^{QG}{{b}^{2}}-\alpha _{A}^{QG}\alpha _{B}^{QG}{{c}^{2}}.\]	

Therefore,
\[\begin{aligned}
	Q{{G}^{2}}&=-\alpha _{B}^{QG}\alpha _{C}^{QG}{{a}^{2}}-\alpha _{C}^{QG}\alpha _{A}^{QG}{{b}^{2}}-\alpha _{A}^{QG}\alpha _{B}^{QG}{{c}^{2}} \\ 
	& =-\frac{1}{9{{W}^{2}}}\left( \begin{aligned}
		& \left( W-3\sin 2B \right)\left( W-3\sin 2C \right){{a}^{2}} \\ 
		& +\left( W-3\sin 2C \right)\left( W-3\sin 2A \right){{b}^{2}} \\ 
		& +\left( W-3\sin 2A \right)\left( W-3\sin 2B \right){{c}^{2}} \\ 
	\end{aligned} \right).  
\end{aligned}\]	

The distance between circumcenter and centroid is the same whether the DOIC or DTICs is used. Make the two formulas equal, and an expression of the circumscribed circle radius is obtained:
\[\begin{aligned}
	{{R}^{2}}&=\frac{1}{9}\left( {{a}^{2}}+{{b}^{2}}+{{c}^{2}} \right)-\frac{1}{9{{W}^{2}}}\left( \begin{aligned}
		& \left( W-3\sin 2B \right)\left( W-3\sin 2C \right){{a}^{2}} \\ 
		& +\left( W-3\sin 2C \right)\left( W-3\sin 2A \right){{b}^{2}} \\ 
		& +\left( W-3\sin 2A \right)\left( W-3\sin 2B \right){{c}^{2}} \\ 
	\end{aligned} \right) \\ 
	& =\frac{1}{3{{W}^{2}}}\left( \begin{aligned}
		& \left( W\left( \sin 2B+\sin 2C \right)-3\sin 2B\cdot \sin 2C \right){{a}^{2}} \\ 
		& +\left( W\left( \sin 2C+\sin 2A \right)-3\sin 2C\cdot \sin 2A \right){{b}^{2}} \\ 
		& +\left( W\left( \sin 2A+\sin 2B \right)-3\sin 2A\cdot \sin 2B \right){{c}^{2}} \\ 
	\end{aligned} \right).  
\end{aligned}\]	

\subsection{DTICs between circumcenter and incenter}\label{Subsec13.2.2}
The following formulas are obtained by using Chapter \ref{Ch8}:
\[\alpha _{A}^{QI}=\alpha _{A}^{I}-\alpha _{A}^{Q}=\frac{a}{2p}-\frac{\sin 2A}{W}=\frac{aW-2p\sin 2A}{2pW},\]
\[\alpha _{B}^{QI}=\alpha _{B}^{I}-\alpha _{B}^{Q}=\frac{b}{2p}-\frac{\sin 2B}{W}=\frac{bW-2p\sin 2B}{2pW},\]
\[\alpha _{C}^{QI}=\alpha _{C}^{I}-\alpha _{C}^{Q}=\frac{c}{2p}-\frac{\sin 2C}{W}=\frac{cW-2p\sin 2C}{2pW},\]
\[\begin{aligned}
	\alpha _{A}^{QI}\alpha _{B}^{QI}&=\frac{aW-2p\sin 2A}{2pW}\cdot \frac{bW-2p\sin 2B}{2pW} \\ 
	& =\frac{\left( aW-2p\sin 2A \right)\left( bW-2p\sin 2B \right)}{4{{p}^{2}}{{W}^{2}}},  
\end{aligned}\]
\[\begin{aligned}
	\alpha _{B}^{QI}\alpha _{C}^{QI}&=\frac{bW-2p\sin 2B}{2pW}\cdot \frac{cW-2p\sin 2C}{2pW} \\ 
	& =\frac{\left( bW-2p\sin 2B \right)\left( cW-2p\sin 2C \right)}{4{{p}^{2}}{{W}^{2}}},  
\end{aligned}\]
\[\begin{aligned}
	\alpha _{C}^{QI}\alpha _{A}^{QI}&=\frac{cW-2p\sin 2C}{2pW}\cdot \frac{aW-2p\sin 2A}{2pW} \\ 
	& =\frac{\left( cW-2p\sin 2C \right)\left( aW-2p\sin 2A \right)}{4{{p}^{2}}{{W}^{2}}}.  
\end{aligned}\]

Therefore, according to theorem \ref{thm:Thm12.2.1}, the following result is obtained
\[\begin{aligned}
	Q{{I}^{2}}&=-\alpha _{B}^{QI}\alpha _{C}^{QI}{{a}^{2}}-\alpha _{C}^{QI}\alpha _{A}^{QI}{{b}^{2}}-\alpha _{A}^{QI}\alpha _{B}^{QI}{{c}^{2}} \\ 
	& =-4{{p}^{2}}{{W}^{2}}\left( \begin{aligned}
		& \left( bW-2p\sin 2B \right)\left( cW-2p\sin 2C \right){{a}^{2}} \\ 
		& +\left( cW-2p\sin 2C \right)\left( aW-2p\sin 2A \right){{b}^{2}} \\ 
		& +\left( aW-2p\sin 2A \right)\left( bW-2p\sin 2B \right){{c}^{2}}  
	\end{aligned} \right).  
\end{aligned}\]

Although the above formula is more complex, it can lead to some trigonometric identities.

The results of the two expressions of $Q{{I}^{2}}$ from the DOIC and the DTICs should be the same, so we get
\[\begin{aligned}
	{{R}^{2}}-\frac{abc}{2p}&=\frac{2pW\left( b\sin 2C+c\sin 2B \right)-bc{{W}^{2}}-4{{p}^{2}}\sin 2B\cdot \sin 2C}{4{{p}^{2}}{{W}^{2}}}{{a}^{2}} \\ 
	& +\frac{2pW\left( c\sin 2A+a\sin 2C \right)-ca{{W}^{2}}-4{{p}^{2}}\sin 2C\cdot \sin 2A}{4{{p}^{2}}{{W}^{2}}}{{b}^{2}} \\ 
	& +\frac{2pW\left( a\sin 2B+b\sin 2A \right)-ab{{W}^{2}}-4{{p}^{2}}\sin 2A\cdot \sin 2B}{4{{p}^{2}}{{W}^{2}}}{{c}^{2}}.  
\end{aligned}\]

This gives a formula for the radius of the circumcircle of a triangle:
\[{{R}^{2}}=\frac{1}{2p{{W}^{2}}}\left( \begin{aligned}
	& \left( W\left( b\sin 2C+c\sin 2B \right)-2p\sin 2B\cdot \sin 2C \right){{a}^{2}} \\ 
	& +\left( W\left( c\sin 2A+a\sin 2C \right)-2p\sin 2C\cdot \sin 2A \right){{b}^{2}} \\ 
	& +\left( W\left( a\sin 2B+b\sin 2A \right)-2p\sin 2A\cdot \sin 2B \right){{c}^{2}}  
\end{aligned} \right).\]	

By applying the sine theorem to the left side of the equation above, some identities can be obtained.

\subsection{DTICs between circumcenter and orthocenter}\label{Subsec13.2.3}
By using Chapter \ref{Ch8} and theorem \ref{thm:Thm12.2.1}, the following formulas are obtained:
\[\alpha _{A}^{QH}=\alpha _{A}^{H}-\alpha _{A}^{Q}=\frac{\tan A}{T}-\frac{\sin 2A}{W}=\frac{W\tan A-T\sin 2A}{TW},\]	
\[\alpha _{B}^{QH}=\alpha _{B}^{H}-\alpha _{B}^{Q}=\frac{\tan B}{T}-\frac{\sin 2B}{W}=\frac{W\tan B-T\sin 2B}{TW},\]	
\[\alpha _{C}^{QH}=\alpha _{C}^{H}-\alpha _{C}^{Q}=\frac{\tan C}{T}-\frac{\sin 2C}{W}=\frac{W\tan C-T\sin 2C}{TW},\]	
\[\begin{aligned}
	\alpha _{A}^{QH}\alpha _{B}^{QH}&=\frac{W\tan A-T\sin 2A}{TW}\cdot \frac{W\tan B-T\sin 2B}{TW} \\ 
	& =\frac{\left( W\tan A-T\sin 2A \right)\left( W\tan B-T\sin 2B \right)}{{{T}^{2}}{{W}^{2}}},  
\end{aligned}\]	
\[\begin{aligned}
	\alpha _{B}^{QH}\alpha _{C}^{QH}&=\frac{W\tan B-T\sin 2B}{TW}\cdot \frac{W\tan C-T\sin 2C}{TW} \\ 
	& =\frac{\left( W\tan B-T\sin 2B \right)\left( W\tan C-T\sin 2C \right)}{{{T}^{2}}{{W}^{2}}},  
\end{aligned}\]	
\[\begin{aligned}
	\alpha _{C}^{QH}\alpha _{A}^{QH}&=\frac{W\tan C-T\sin 2C}{TW}\cdot \frac{W\tan A-T\sin 2A}{TW} \\ 
	& =\frac{\left( W\tan C-T\sin 2C \right)\left( W\tan A-T\sin 2A \right)}{{{T}^{2}}{{W}^{2}}}.  
\end{aligned}\]	

Therefore,
\[\begin{aligned}
	Q{{H}^{2}}&=-\alpha _{B}^{QH}\alpha _{C}^{QH}{{a}^{2}}-\alpha _{C}^{QH}\alpha _{A}^{QH}{{b}^{2}}-\alpha _{A}^{QH}\alpha _{B}^{QH}{{c}^{2}} \\ 
	& =-\frac{1}{{{T}^{2}}{{W}^{2}}}\left( \begin{aligned}
		& \left( W\tan B-T\sin 2B \right)\left( W\tan C-T\sin 2C \right){{a}^{2}} \\ 
		& +\left( W\tan C-T\sin 2C \right)\left( W\tan A-T\sin 2A \right){{b}^{2}} \\ 
		& +\left( W\tan A-T\sin 2A \right)\left( W\tan B-T\sin 2B \right){{c}^{2}} \\ 
	\end{aligned} \right)  
\end{aligned}.\]	

Similarly, the square of the radius of the circumscribed circle is
\[\begin{aligned}
	{{R}^{2}}&=\frac{\tan B\cdot \tan C\cdot {{a}^{2}}+\tan C\cdot \tan A\cdot {{b}^{2}}+\tan A\cdot \tan B\cdot {{c}^{2}}}{{{T}^{2}}} \\ 
	& -\frac{1}{{{T}^{2}}{{W}^{2}}}\left( \begin{aligned}
		& \left( W\tan B-T\sin 2B \right)\left( W\tan C-T\sin 2C \right){{a}^{2}} \\ 
		& +\left( W\tan C-T\sin 2C \right)\left( W\tan A-T\sin 2A \right){{b}^{2}} \\ 
		& +\left( W\tan A-T\sin 2A \right)\left( W\tan B-T\sin 2B \right){{c}^{2}} \\ 
	\end{aligned} \right) \\ 
	& =\frac{1}{T{{W}^{2}}}\left( \begin{aligned}
		& \left( W\left( \sin 2B\cdot \tan C+\tan B\cdot \sin 2C \right)-T\sin 2B\cdot \sin 2C \right){{a}^{2}} \\ 
		& +\left( W\left( \sin 2C\cdot \tan A+\tan C\cdot \sin 2A \right)-T\sin 2C\cdot \sin 2A \right){{b}^{2}} \\ 
		& +\left( W\left( \sin 2A\cdot \tan B+\tan A\cdot \sin 2B \right)-T\sin 2A\cdot \sin 2B \right){{c}^{2}} \\ 
	\end{aligned} \right).  
\end{aligned}\]

\subsection{DTICs between circumcenter and excenter}\label{Subsec13.2.4}
Similar to the previous method, leave it to the reader as an exercise.

The distance between the centroid and the circumcenter has been discussed before. The distance between the centroid and the other centers of the triangle (except the circumcenter) is discussed below.
The formula (\ref{Eq12.1.2}) is useless when the two ICs are neither circumcenter. In the following subsections, the formula of DTICs is used for calculation.

\subsection{DTICs between centroid and incenter}\label{Subsec13.2.5}
By using Chapter \ref{Ch8} and theorem \ref{thm:Thm12.2.1}, the following formulas are obtained:
\[\alpha _{A}^{GI}=\alpha _{A}^{I}-\alpha _{A}^{G}=\frac{a}{2p}-\frac{1}{3}=\frac{3a-2p}{6p},\]
\[\alpha _{B}^{GI}=\alpha _{B}^{I}-\alpha _{B}^{G}=\frac{b}{2p}-\frac{1}{3}=\frac{3b-2p}{6p},\]
\[\alpha _{C}^{GI}=\alpha _{C}^{I}-\alpha _{C}^{G}=\frac{c}{2p}-\frac{1}{3}=\frac{3c-2p}{6p}.\]
\[\begin{aligned}
	G{{I}^{2}}&=-\alpha _{B}^{GI}\alpha _{C}^{GI}{{a}^{2}}-\alpha _{C}^{GI}\alpha _{A}^{GI}{{b}^{2}}-\alpha _{A}^{GI}\alpha _{B}^{GI}{{c}^{2}} \\ 
	& =-\frac{1}{36{{p}^{2}}}\left( \begin{aligned}
		& \left( 3b-2p \right)\left( 3c-2p \right){{a}^{2}} \\ 
		& +\left( 3c-2p \right)\left( 3a-2p \right){{b}^{2}} \\ 
		& +\left( 3a-2p \right)\left( 3b-2p \right){{c}^{2}} \\ 
	\end{aligned} \right).  
\end{aligned}\]

\subsection{DTICs between centroid and orthocenter}\label{Subsec13.2.6}
By using Chapter \ref{Ch8} and theorem \ref{thm:Thm12.2.1}, the following formulas are obtained:
\[\alpha _{A}^{GH}=\alpha _{A}^{H}-\alpha _{A}^{G}=\frac{\tan A}{T}-\frac{1}{3}=\frac{3\tan A-T}{3T},\]
\[\alpha _{B}^{GH}=\alpha _{B}^{H}-\alpha _{B}^{G}=\frac{\tan B}{T}-\frac{1}{3}=\frac{3\tan B-T}{3T},\]
\[\alpha _{C}^{GH}=\alpha _{C}^{H}-\alpha _{C}^{G}=\frac{\tan C}{T}-\frac{1}{3}=\frac{3\tan C-T}{3T},\]
\[\begin{aligned}
	G{{H}^{2}}&=-\alpha _{B}^{GH}\alpha _{C}^{GH}{{a}^{2}}-\alpha _{C}^{GH}\alpha _{A}^{GH}{{b}^{2}}-\alpha _{A}^{GH}\alpha _{B}^{GH}{{c}^{2}} \\ 
	& =-\frac{1}{9{{T}^{2}}}\left( \begin{aligned}
		& \left( 3\tan B-T \right)\left( 3\tan C-T \right){{a}^{2}} \\ 
		& +\left( 3\tan C-T \right)\left( 3\tan A-T \right){{b}^{2}} \\ 
		& +\left( 3\tan A-T \right)\left( 3\tan B-T \right){{c}^{2}} \\ 
	\end{aligned} \right).  
\end{aligned}\]

\subsection{DTICs between centroid and excenter}\label{Subsec13.2.7}
\subsubsection{DTICs between centroid and excenter ${{E}_{A}}$ corresponding to $\angle A$}
By using Chapter \ref{Ch8} and theorem \ref{thm:Thm12.2.1}, the following formulas are obtained:
\[\begin{aligned}
	\alpha _{A}^{G{{E}_{A}}}&=\alpha _{A}^{{{E}_{A}}}-\alpha _{A}^{G}=-\frac{a}{2\left( p-a \right)}-\frac{1}{3} \\ 
	& =-\frac{3a+2\left( p-a \right)}{6\left( p-a \right)}=-\frac{2p+a}{6\left( p-a \right)},  
\end{aligned}\]
\[\begin{aligned}
	\alpha _{B}^{G{{E}_{A}}}&=\alpha _{B}^{{{E}_{A}}}-\alpha _{B}^{G}=\frac{b}{2\left( p-a \right)}-\frac{1}{3} \\ 
	& =\frac{3b-2\left( p-a \right)}{6\left( p-a \right)}=\frac{2b-c+a}{6\left( p-a \right)},  
\end{aligned}\]
\[\begin{aligned}
	\alpha _{C}^{G{{E}_{A}}}&=\alpha _{C}^{{{E}_{A}}}-\alpha _{C}^{G}=\frac{c}{2\left( p-a \right)}-\frac{1}{3} \\ 
	& =\frac{3c-2\left( p-a \right)}{6\left( p-a \right)}=\frac{2c+a-b}{6\left( p-a \right)},  
\end{aligned}\]
\[\begin{aligned}
	G{{E}_{A}}^{2}&=-\alpha _{B}^{G{{E}_{A}}}\alpha _{C}^{G{{E}_{A}}}{{a}^{2}}-\alpha _{C}^{G{{E}_{A}}}\alpha _{A}^{G{{E}_{A}}}{{b}^{2}}-\alpha _{A}^{G{{E}_{A}}}\alpha _{B}^{G{{E}_{A}}}{{c}^{2}} \\ 
	& =\frac{1}{36{{\left( p-a \right)}^{2}}}\left( \begin{aligned}
		& -\left( 2b-c+a \right)\left( 2c+a-b \right){{a}^{2}} \\ 
		& +\left( 2c+a-b \right)\left( 2p+a \right){{b}^{2}} \\ 
		& +\left( 2p+a \right)\left( 2b-c+a \right){{c}^{2}} \\ 
	\end{aligned} \right).  
\end{aligned}\]
\subsubsection{DTICs between centroid and excenter ${{E}_{B}}$ corresponding to $\angle B$}
By using Chapter \ref{Ch8} and theorem \ref{thm:Thm12.2.1}, the following formulas are obtained:
\[\begin{aligned}
	\alpha _{A}^{G{{E}_{B}}}&=\alpha _{A}^{{{E}_{B}}}-\alpha _{A}^{G}=\frac{a}{2\left( p-b \right)}-\frac{1}{3} \\ 
	& =\frac{3a-2\left( p-b \right)}{6\left( p-b \right)}=\frac{3a+2b-2p}{6\left( p-b \right)}=\frac{2a+b-c}{6\left( p-b \right)},  
\end{aligned}\]
\[\begin{aligned}
	\alpha _{B}^{G{{E}_{B}}}&=\alpha _{B}^{{{E}_{B}}}-\alpha _{B}^{G}=-\frac{b}{2\left( p-b \right)}-\frac{1}{3} \\ 
	& =-\frac{3b+2\left( p-b \right)}{6\left( p-b \right)}=-\frac{3b-2b+2p}{6\left( p-b \right)}=-\frac{2p+b}{6\left( p-b \right)},  
\end{aligned}\]
\[\begin{aligned}
	\alpha _{C}^{G{{E}_{B}}}&=\alpha _{C}^{{{E}_{B}}}-\alpha _{C}^{G}=\frac{c}{2\left( p-b \right)}-\frac{1}{3} \\ 
	& =\frac{3c-2\left( p-b \right)}{6\left( p-b \right)}=\frac{3c+2b-2p}{6\left( p-b \right)}=\frac{2c-a+b}{6\left( p-b \right)},  
\end{aligned}\]
\[\begin{aligned}
	G{{E}_{B}}^{2}&=-\alpha _{B}^{G{{E}_{B}}}\alpha _{C}^{G{{E}_{B}}}{{a}^{2}}-\alpha _{C}^{G{{E}_{B}}}\alpha _{A}^{G{{E}_{B}}}{{b}^{2}}-\alpha _{A}^{G{{E}_{B}}}\alpha _{B}^{G{{E}_{B}}}{{c}^{2}} \\ 
	& =\frac{1}{36{{\left( p-b \right)}^{2}}}\left( \begin{aligned}
		& \left( 2p+b \right)\left( 2c-a+b \right){{a}^{2}} \\ 
		& -\left( 2c-a+b \right)\left( 2a+b-c \right){{b}^{2}} \\ 
		& +\left( 2a+b-c \right)\left( 2p+b \right){{c}^{2}} \\ 
	\end{aligned} \right).  
\end{aligned}\]	
\subsubsection{DTICs between centroid and excenter ${{E}_{C}}$ corresponding to $\angle C$}
By using Chapter \ref{Ch8} and theorem \ref{thm:Thm12.2.1}, the following formulas are obtained:
\[\begin{aligned}
	\alpha _{A}^{G{{E}_{C}}}&=\alpha _{A}^{{{E}_{C}}}-\alpha _{A}^{G}=\frac{a}{2\left( p-c \right)}-\frac{1}{3} \\ 
	& =\frac{3a-2\left( p-c \right)}{6\left( p-c \right)}=\frac{3a+2c-2p}{6\left( p-c \right)}=\frac{2a-b+c}{6\left( p-c \right)},  
\end{aligned}\]
\[\begin{aligned}
	\alpha _{B}^{G{{E}_{C}}}&=\alpha _{B}^{{{E}_{C}}}-\alpha _{B}^{G}=\frac{b}{2\left( p-c \right)}-\frac{1}{3} \\ 
	& =\frac{3b-2\left( p-c \right)}{6\left( p-c \right)}=\frac{3b+2c-2p}{6\left( p-c \right)}=\frac{2b+c-b}{6\left( p-c \right)},  
\end{aligned}\]
\[\begin{aligned}
	\alpha _{C}^{G{{E}_{C}}}&=\alpha _{C}^{{{E}_{C}}}-\alpha _{C}^{G}=-\frac{c}{2\left( p-c \right)}-\frac{1}{3} \\ 
	& =-\frac{3c+2\left( p-c \right)}{6\left( p-c \right)}=-\frac{2p+c}{6\left( p-c \right)},  
\end{aligned}\]
\[\begin{aligned}
	G{{E}_{C}}^{2}&=-\alpha _{B}^{G{{E}_{C}}}\alpha _{C}^{G{{E}_{C}}}{{a}^{2}}-\alpha _{C}^{G{{E}_{C}}}\alpha _{A}^{G{{E}_{C}}}{{b}^{2}}-\alpha _{A}^{G{{E}_{C}}}\alpha _{B}^{G{{E}_{C}}}{{c}^{2}} \\ 
	& =\frac{1}{36{{\left( p-c \right)}^{2}}}\left( \begin{aligned}
		& \left( 2b+c-b \right)\left( 2p+c \right){{a}^{2}} \\ 
		& +\left( 2p+c \right)\left( 2a-b+c \right){{b}^{2}} \\ 
		& -\left( 2a-b+c \right)\left( 2b+c-b \right){{c}^{2}} \\ 
	\end{aligned} \right).  
\end{aligned}\]	

\subsection{DTICs between incenter and orthocenter}\label{Subsec13.2.8}
By using Chapter \ref{Ch8} and theorem \ref{thm:Thm12.2.1}, the following formulas are obtained:
\[\alpha _{A}^{HI}=\alpha _{A}^{I}-\alpha _{A}^{H}=\frac{a}{2p}-\frac{\tan A}{T}=\frac{aT-2p\tan A}{2pT},\]
\[\alpha _{B}^{HI}=\alpha _{B}^{I}-\alpha _{B}^{H}=\frac{b}{2p}-\frac{\tan B}{T}=\frac{bT-2p\tan B}{2pT},\]
\[\alpha _{C}^{HI}=\alpha _{C}^{I}-\alpha _{C}^{H}=\frac{c}{2p}-\frac{\tan C}{T}=\frac{cT-2p\tan C}{2pT},\]
\[\begin{aligned}
	\alpha _{A}^{HI}\alpha _{B}^{HI}&=\frac{aT-2p\tan A}{2pT}\cdot \frac{bT-2p\tan B}{2pT} \\ 
	& =\frac{1}{4{{p}^{2}}{{T}^{2}}}\left( aT-2p\tan A \right)\left( bT-2p\tan B \right),  
\end{aligned}\]
\[\begin{aligned}
	\alpha _{B}^{HI}\alpha _{C}^{HI}&=\frac{bT-2p\tan B}{2pT}\cdot \frac{cT-2p\tan C}{2pT} \\ 
	& =\frac{1}{4{{p}^{2}}{{T}^{2}}}\left( bT-2p\tan B \right)\left( cT-2p\tan C \right),  
\end{aligned}\]
\[\begin{aligned}
	\alpha _{C}^{HI}\alpha _{A}^{HI}&=\frac{cT-2p\tan C}{2pT}\cdot \frac{aT-2p\tan A}{2pT} \\ 
	& =\frac{1}{4{{p}^{2}}{{T}^{2}}}\left( cT-2p\tan C \right)\left( aT-2p\tan A \right).  
\end{aligned}\]

Therefore,
\[\begin{aligned}
	I{{H}^{2}}&=-\alpha _{B}^{HI}\alpha _{C}^{HI}{{a}^{2}}-\alpha _{C}^{HI}\alpha _{A}^{HI}{{b}^{2}}-\alpha _{A}^{HI}\alpha _{B}^{HI}{{c}^{2}} \\ 
	& =-\frac{1}{4{{p}^{2}}{{T}^{2}}}\left( \begin{aligned}
		& \left( bT-2p\tan B \right)\left( cT-2p\tan C \right){{a}^{2}} \\ 
		& +\left( cT-2p\tan C \right)\left( aT-2p\tan A \right){{b}^{2}} \\ 
		& +\left( aT-2p\tan A \right)\left( bT-2p\tan B \right){{c}^{2}} \\ 
	\end{aligned} \right).  
\end{aligned}\]

\subsection{DTICs between incenter and excenter}\label{Subsec13.2.9}
\subsubsection{DTICs between incenter and excenter ${{E}_{A}}$ corresponding to $\angle A$}
By using Chapter \ref{Ch8} and theorem \ref{thm:Thm12.2.1}, the following formulas are obtained:
\[\begin{aligned}
	\alpha _{A}^{I{{E}_{A}}}&=\alpha _{A}^{{{E}_{A}}}-\alpha _{A}^{I}=-\frac{a}{2\left( p-a \right)}-\frac{a}{2p} \\ 
	& =-\frac{a\left( p+\left( p-a \right) \right)}{2p\left( p-a \right)}=-\frac{a\left( 2p-a \right)}{2p\left( p-a \right)},  
\end{aligned}\]
\[\begin{aligned}
	\alpha _{B}^{I{{E}_{A}}}&=\alpha _{B}^{{{E}_{A}}}-\alpha _{B}^{I}=\frac{b}{2\left( p-a \right)}-\frac{b}{2p} \\ 
	& =\frac{b\left( p-\left( p-a \right) \right)}{2p\left( p-a \right)}=\frac{ab}{2p\left( p-a \right)},  
\end{aligned}\]
\[\begin{aligned}
	\alpha _{C}^{I{{E}_{A}}}&=\alpha _{C}^{{{E}_{A}}}-\alpha _{C}^{I}=\frac{c}{2\left( p-a \right)}-\frac{c}{2p} \\ 
	& =\frac{c\left( p-\left( p-a \right) \right)}{2p\left( p-a \right)}=\frac{ca}{2p\left( p-a \right)},  
\end{aligned}\]
\[\alpha _{A}^{I{{E}_{A}}}\alpha _{B}^{I{{E}_{A}}}=-\frac{a\left( 2p-a \right)}{2p\left( p-a \right)}\cdot \frac{ab}{2p\left( p-a \right)}=-\frac{{{a}^{2}}b\left( 2p-a \right)}{4{{p}^{2}}{{\left( p-a \right)}^{2}}},\]
\[\alpha _{B}^{I{{E}_{A}}}\alpha _{C}^{I{{E}_{A}}}=\frac{ab}{2p\left( p-a \right)}\cdot \frac{ca}{2p\left( p-a \right)}=\frac{bc{{a}^{2}}}{4{{p}^{2}}{{\left( p-a \right)}^{2}}},\]
\[\alpha _{C}^{I{{E}_{A}}}\alpha _{A}^{I{{E}_{A}}}=-\frac{ca}{2p\left( p-a \right)}\cdot \frac{a\left( 2p-a \right)}{2p\left( p-a \right)}=-\frac{c{{a}^{2}}\left( 2p-a \right)}{4{{p}^{2}}{{\left( p-a \right)}^{2}}}.\]

Therefore,
\[\begin{aligned}
	I{{E}_{A}}^{2}&=-\alpha _{B}^{I{{E}_{A}}}\alpha _{C}^{I{{E}_{A}}}{{a}^{2}}-\alpha _{C}^{I{{E}_{A}}}\alpha _{A}^{I{{E}_{A}}}{{b}^{2}}-\alpha _{A}^{I{{E}_{A}}}\alpha _{B}^{I{{E}_{A}}}{{c}^{2}} \\ 
	& =\frac{-bc{{a}^{4}}+c{{a}^{2}}\left( 2p-a \right){{b}^{2}}+{{a}^{2}}b\left( 2p-a \right){{c}^{2}}}{4{{p}^{2}}{{\left( p-a \right)}^{2}}} \\ 
	& =\frac{{{a}^{2}}bc}{4{{p}^{2}}{{\left( p-a \right)}^{2}}}\left( -{{a}^{2}}+b\left( 2p-a \right)+c\left( 2p-a \right) \right) \\ 
	& =\frac{{{a}^{2}}bc}{4{{p}^{2}}{{\left( p-a \right)}^{2}}}\left( {{\left( b+c \right)}^{2}}-{{a}^{2}} \right)=\frac{{{a}^{2}}bc}{4{{p}^{2}}{{\left( p-a \right)}^{2}}}\left( 2p\cdot \left( 2p-2a \right) \right)=\frac{{{a}^{2}}bc}{p\left( p-a \right)},  
\end{aligned}\]
\[I{{E}_{A}}=a\sqrt{\frac{bc}{p\left( p-a \right)}}.\]

\subsubsection{DTICs between incenter and excenter ${{E}_{B}}$ corresponding to $\angle B$}
By using Chapter \ref{Ch8} and theorem \ref{thm:Thm12.2.1}, the following formulas are obtained:
\[\begin{aligned}
	\alpha _{A}^{I{{E}_{B}}}&=\alpha _{A}^{{{E}_{B}}}-\alpha _{A}^{I}=\frac{a}{2\left( p-b \right)}-\frac{a}{2p} \\ 
	& =-\frac{a\left( p-\left( p-b \right) \right)}{2p\left( p-b \right)}=\frac{ab}{2p\left( p-b \right)},  
\end{aligned}\]
\[\begin{aligned}
	\alpha _{B}^{I{{E}_{B}}}&=\alpha _{B}^{{{E}_{B}}}-\alpha _{B}^{I}=-\frac{b}{2\left( p-b \right)}-\frac{b}{2p} \\ 
	& =-\frac{b\left( p+\left( p-b \right) \right)}{2p\left( p-b \right)}=-\frac{b\left( 2p-b \right)}{2p\left( p-b \right)},  
\end{aligned}\]
\[\begin{aligned}
	\alpha _{C}^{I{{E}_{B}}}&=\alpha _{C}^{{{E}_{B}}}-\alpha _{C}^{I}=\frac{c}{2\left( p-b \right)}-\frac{c}{2p} \\ 
	& =\frac{c\left( p-\left( p-b \right) \right)}{2p\left( p-b \right)}=\frac{bc}{2p\left( p-b \right)},  
\end{aligned}\]
\[\alpha _{A}^{I{{E}_{B}}}\alpha _{B}^{I{{E}_{B}}}=-\frac{ab}{2p\left( p-b \right)}\cdot \frac{b\left( 2p-b \right)}{2p\left( p-b \right)}=-\frac{a{{b}^{2}}\left( 2p-b \right)}{4{{p}^{2}}{{\left( p-b \right)}^{2}}},\]
\[\alpha _{B}^{I{{E}_{B}}}\alpha _{C}^{I{{E}_{B}}}=-\frac{b\left( 2p-b \right)}{2p\left( p-b \right)}\cdot \frac{bc}{2p\left( p-b \right)}=-\frac{{{b}^{2}}c\left( 2p-b \right)}{4{{p}^{2}}{{\left( p-b \right)}^{2}}},\]
\[\alpha _{C}^{I{{E}_{C}}}\alpha _{A}^{I{{E}_{C}}}=\frac{bc}{2p\left( p-b \right)}\cdot \frac{ab}{2p\left( p-b \right)}=\frac{ca{{b}^{2}}}{4{{p}^{2}}{{\left( p-b \right)}^{2}}}.\]

Therefore,
\[\begin{aligned}
	I{{E}_{B}}^{2}&=-\alpha _{B}^{I{{E}_{B}}}\alpha _{C}^{I{{E}_{B}}}{{a}^{2}}-\alpha _{C}^{I{{E}_{B}}}\alpha _{A}^{I{{E}_{B}}}{{b}^{2}}-\alpha _{A}^{I{{E}_{B}}}\alpha _{B}^{I{{E}_{B}}}{{c}^{2}} \\ 
	& =\frac{{{b}^{2}}c\left( 2p-b \right){{a}^{2}}-ca{{b}^{2}}\cdot {{b}^{2}}+a{{b}^{2}}\left( 2p-b \right){{c}^{2}}}{4{{p}^{2}}{{\left( p-b \right)}^{2}}} \\ 
	& =\frac{{{b}^{2}}ca}{4{{p}^{2}}{{\left( p-b \right)}^{2}}}\left( a\left( 2p-b \right)-{{b}^{2}}+c\left( 2p-b \right) \right) \\ 
	& =\frac{{{b}^{2}}ca}{4{{p}^{2}}{{\left( p-b \right)}^{2}}}\left( {{\left( c+a \right)}^{2}}-{{b}^{2}} \right)=\frac{{{b}^{2}}ca}{4{{p}^{2}}{{\left( p-b \right)}^{2}}}\left( 2p\left( 2p-b \right) \right)=\frac{{{b}^{2}}ca}{p\left( p-b \right)},  
\end{aligned}\]	
\[I{{E}_{B}}=b\sqrt{\frac{ca}{p\left( p-b \right)}}.\]

\subsubsection{DTICs between incenter and excenter ${{E}_{C}}$ corresponding to $\angle C$}
By using Chapter \ref{Ch8} and theorem \ref{thm:Thm12.2.1}, the following formulas are obtained:
\[\begin{aligned}
	\alpha _{A}^{I{{E}_{C}}}&=\alpha _{A}^{{{E}_{C}}}-\alpha _{A}^{I}=\frac{a}{2\left( p-c \right)}-\frac{a}{2p} \\ 
	& =\frac{a\left( p-\left( p-c \right) \right)}{2p\left( p-c \right)}=\frac{ca}{2p\left( p-c \right)},  
\end{aligned}\]
\[\begin{aligned}
	\alpha _{B}^{I{{E}_{C}}}&=\alpha _{B}^{{{E}_{C}}}-\alpha _{B}^{I}=\frac{b}{2\left( p-c \right)}-\frac{b}{2p} \\ 
	& =\frac{b\left( p-\left( p-c \right) \right)}{2p\left( p-c \right)}=\frac{bc}{2p\left( p-c \right)},  
\end{aligned}\]
\[\begin{aligned}
	\alpha _{C}^{I{{E}_{C}}}&=\alpha _{C}^{{{E}_{C}}}-\alpha _{C}^{I}=-\frac{c}{2\left( p-c \right)}-\frac{c}{2p} \\ 
	& =-\frac{c\left( p+\left( p-c \right) \right)}{2p\left( p-c \right)}=-\frac{c\left( 2p-c \right)}{2p\left( p-c \right)},  
\end{aligned}\]
\[\alpha _{A}^{I{{E}_{C}}}\alpha _{B}^{I{{E}_{C}}}=\frac{ca}{2p\left( p-c \right)}\cdot \frac{bc}{2p\left( p-c \right)}=\frac{ab{{c}^{2}}\left( 2p-c \right)}{4{{p}^{2}}{{\left( p-c \right)}^{2}}},\]
\[\alpha _{B}^{I{{E}_{C}}}\alpha _{C}^{I{{E}_{C}}}=-\frac{bc}{2p\left( p-c \right)}\cdot \frac{c\left( 2p-c \right)}{2p\left( p-c \right)}=-\frac{b{{c}^{2}}\left( 2p-c \right)}{4{{p}^{2}}{{\left( p-c \right)}^{2}}},\]
\[\alpha _{C}^{I{{E}_{C}}}\alpha _{A}^{I{{E}_{C}}}=-\frac{c\left( 2p-c \right)}{2p\left( p-c \right)}\cdot \frac{ca}{2p\left( p-c \right)}=-\frac{{{c}^{2}}a\left( 2p-c \right)}{4{{p}^{2}}{{\left( p-c \right)}^{2}}}.\]

Therefore,
\[\begin{aligned}
	I{{E}_{C}}^{2}&=-\alpha _{B}^{I{{E}_{C}}}\alpha _{C}^{I{{E}_{C}}}{{a}^{2}}-\alpha _{C}^{I{{E}_{C}}}\alpha _{A}^{I{{E}_{C}}}{{b}^{2}}-\alpha _{A}^{I{{E}_{C}}}\alpha _{B}^{I{{E}_{C}}}{{c}^{2}} \\ 
	& =\frac{b{{c}^{2}}\left( 2p-c \right){{a}^{2}}+{{c}^{2}}a\left( 2p-c \right){{b}^{2}}-ab{{c}^{2}}\left( 2p-c \right){{c}^{2}}}{4{{p}^{2}}{{\left( p-c \right)}^{2}}} \\ 
	& =\frac{{{c}^{2}}ab}{4{{p}^{2}}{{\left( p-c \right)}^{2}}}\left( a\left( 2p-c \right)+b\left( 2p-c \right)-{{c}^{2}} \right) \\ 
	& =\frac{{{c}^{2}}ab}{4{{p}^{2}}{{\left( p-c \right)}^{2}}}\left( {{\left( a+b \right)}^{2}}-{{c}^{2}} \right)=\frac{{{c}^{2}}ab}{4{{p}^{2}}{{\left( p-c \right)}^{2}}}\left( 2p\left( 2p-2c \right) \right)=\frac{{{c}^{2}}ab}{p\left( p-c \right)},  
\end{aligned}\]
\[I{{E}_{C}}=c\sqrt{\frac{ab}{p\left( p-c \right)}}.\]	

\subsection{DTICs between orthocenter and excenter}\label{Subsec13.2.10}
\subsubsection{DTICs between orthocenter and excenter ${{E}_{A}}$ corresponding to $\angle A$}
By using Chapter \ref{Ch8} and theorem \ref{thm:Thm12.2.1}, the following formulas are obtained:
\[\begin{aligned}
	\alpha _{A}^{H{{E}_{A}}}&=\alpha _{A}^{{{E}_{A}}}-\alpha _{A}^{H}=-\frac{a}{2\left( p-a \right)}-\frac{\tan A}{T} \\ 
	& =-\frac{aT+2\left( p-a \right)\tan A}{2\left( p-a \right)T},  
\end{aligned}\]
\[\begin{aligned}
	\alpha _{B}^{H{{E}_{A}}}&=\alpha _{B}^{{{E}_{A}}}-\alpha _{B}^{H}=\frac{b}{2\left( p-a \right)}-\frac{\tan B}{T} \\ 
	& =\frac{bT-2\left( p-a \right)\tan B}{2\left( p-a \right)T},  
\end{aligned}\]
\[\begin{aligned}
	\alpha _{C}^{H{{E}_{A}}}&=\alpha _{C}^{{{E}_{A}}}-\alpha _{C}^{H}=\frac{c}{2\left( p-a \right)}-\frac{\tan C}{T} \\ 
	& =\frac{cT-2\left( p-a \right)\tan C}{2\left( p-a \right)T},  
\end{aligned}\]
\[\begin{aligned}
	\alpha _{A}^{H{{E}_{A}}}\alpha _{B}^{H{{E}_{A}}}&=-\frac{aT+2\left( p-a \right)\tan A}{2\left( p-a \right)T}\cdot \frac{bT-2\left( p-a \right)\tan B}{2\left( p-a \right)T} \\ 
	& =-\frac{\left( aT+2\left( p-a \right)\tan A \right)\left( bT-2\left( p-a \right)\tan B \right)}{4{{\left( p-a \right)}^{2}}{{T}^{2}}},  
\end{aligned}\]
\[\begin{aligned}
	\alpha _{B}^{H{{E}_{A}}}\alpha _{C}^{H{{E}_{A}}}&=\frac{bT-2\left( p-a \right)\tan B}{2\left( p-a \right)T}\cdot \frac{cT-2\left( p-a \right)\tan C}{2\left( p-a \right)T} \\ 
	& =\frac{\left( bT-2\left( p-a \right)\tan B \right)\left( cT-2\left( p-a \right)\tan C \right)}{4{{\left( p-a \right)}^{2}}{{T}^{2}}},  
\end{aligned}\]
\[\begin{aligned}
	\alpha _{C}^{H{{E}_{A}}}\alpha _{A}^{H{{E}_{A}}}&=-\frac{cT-2\left( p-a \right)\tan C}{2\left( p-a \right)T}\cdot \frac{aT+2\left( p-a \right)\tan A}{2\left( p-a \right)T} \\ 
	& =-\frac{\left( cT-2\left( p-a \right)\tan C \right)\left( aT+2\left( p-a \right)\tan A \right)}{4{{\left( p-a \right)}^{2}}{{T}^{2}}}.  
\end{aligned}\]

Therefore,
\[\begin{aligned}
	H{{E}_{A}}^{2}&=-\alpha _{B}^{H{{E}_{A}}}\alpha _{C}^{H{{E}_{A}}}{{a}^{2}}-\alpha _{C}^{H{{E}_{A}}}\alpha _{A}^{H{{E}_{A}}}{{b}^{2}}-\alpha _{A}^{H{{E}_{A}}}\alpha _{B}^{H{{E}_{A}}}{{c}^{2}} \\ 
	& =\frac{1}{4{{T}^{2}}{{\left( p-a \right)}^{2}}}\left( \begin{aligned}
		& -\left( bT-2\left( p-a \right)\tan B \right)\left( cT-2\left( p-a \right)\tan C \right){{a}^{2}} \\ 
		& +\left( cT-2\left( p-a \right)\tan C \right)\left( aT+2\left( p-a \right)\tan A \right){{b}^{2}} \\ 
		& +\left( aT+2\left( p-a \right)\tan A \right)\left( bT-2\left( p-a \right)\tan B \right){{c}^{2}} \\ 
	\end{aligned} \right).  
\end{aligned}\]

\subsubsection{DTICs between orthocenter and excenter ${{E}_{B}}$ corresponding to $\angle B$}
By using Chapter \ref{Ch8} and theorem \ref{thm:Thm12.2.1}, the following formulas are obtained:
\[\begin{aligned}
	\alpha _{A}^{H{{E}_{B}}}&=\alpha _{A}^{{{E}_{B}}}-\alpha _{A}^{H}=\frac{a}{2\left( p-b \right)}-\frac{\tan A}{T} \\ 
	& =\frac{aT-2\left( p-b \right)\tan A}{2\left( p-b \right)T},  
\end{aligned}\]
\[\begin{aligned}
	\alpha _{B}^{H{{E}_{B}}}&=\alpha _{B}^{{{E}_{B}}}-\alpha _{B}^{H}=-\frac{b}{2\left( p-b \right)}-\frac{\tan B}{T} \\ 
	& =-\frac{bT+2\left( p-b \right)\tan B}{2\left( p-b \right)T},  
\end{aligned}\]
\[\begin{aligned}
	\alpha _{C}^{H{{E}_{B}}}&=\alpha _{C}^{{{E}_{B}}}-\alpha _{C}^{H}=\frac{c}{2\left( p-b \right)}-\frac{\tan C}{T} \\ 
	& =\frac{cT-2\left( p-b \right)\tan C}{2\left( p-b \right)T},  
\end{aligned}\]
\[\begin{aligned}
	\alpha _{A}^{H{{E}_{B}}}\alpha _{B}^{H{{E}_{B}}}&=-\frac{aT-2\left( p-b \right)\tan A}{2\left( p-b \right)T}\cdot \frac{bT+2\left( p-b \right)\tan B}{2\left( p-b \right)T} \\ 
	& =-\frac{\left( aT-2\left( p-b \right)\tan A \right)\left( bT+2\left( p-b \right)\tan B \right)}{4{{\left( p-b \right)}^{2}}{{T}^{2}}},  
\end{aligned}\]
\[\begin{aligned}
	\alpha _{B}^{H{{E}_{B}}}\alpha _{C}^{H{{E}_{B}}}&=-\frac{bT+2\left( p-b \right)\tan B}{2\left( p-b \right)T}\cdot \frac{cT-2\left( p-b \right)\tan C}{2\left( p-b \right)T} \\ 
	& =-\frac{\left( bT+2\left( p-b \right)\tan B \right)\left( cT-2\left( p-b \right)\tan C \right)}{4{{\left( p-b \right)}^{2}}{{T}^{2}}},  
\end{aligned}\]
\[\begin{aligned}
	\alpha _{C}^{H{{E}_{B}}}\alpha _{A}^{H{{E}_{B}}}&=-\frac{cT-2\left( p-b \right)\tan C}{2\left( p-b \right)T}\cdot \frac{aT-2\left( p-b \right)\tan A}{2\left( p-b \right)T} \\ 
	& =\frac{\left( cT-2\left( p-b \right)\tan C \right)\left( aT-2\left( p-b \right)\tan A \right)}{4{{\left( p-b \right)}^{2}}{{T}^{2}}}.  
\end{aligned}\]

Therefore,
\[\begin{aligned}
	H{{E}_{B}}^{2}&=-\alpha _{B}^{H{{E}_{B}}}\alpha _{C}^{H{{E}_{B}}}{{a}^{2}}-\alpha _{C}^{H{{E}_{B}}}\alpha _{A}^{H{{E}_{B}}}{{b}^{2}}-\alpha _{A}^{H{{E}_{B}}}\alpha _{B}^{H{{E}_{B}}}{{c}^{2}} \\ 
	& =\frac{1}{4{{T}^{2}}{{\left( p-b \right)}^{2}}}\left( \begin{aligned}
		& \left( bT+2\left( p-b \right)\tan B \right)\left( cT-2\left( p-b \right)\tan C \right){{a}^{2}} \\ 
		& -\left( cT-2\left( p-b \right)\tan C \right)\left( aT-2\left( p-b \right)\tan A \right){{b}^{2}} \\ 
		& +\left( aT-2\left( p-b \right)\tan A \right)\left( bT+2\left( p-b \right)\tan B \right){{c}^{2}} \\ 
	\end{aligned} \right).  
\end{aligned}\]

\subsubsection{DTICs between orthocenter and excenter ${{E}_{C}}$ corresponding to $\angle C$}
By using Chapter \ref{Ch8} and theorem \ref{thm:Thm12.2.1}, the following formulas are obtained:
\[\begin{aligned}
	\alpha _{A}^{H{{E}_{C}}}&=\alpha _{A}^{{{E}_{C}}}-\alpha _{A}^{H}=\frac{a}{2\left( p-c \right)}-\frac{\tan A}{T} \\ 
	& =\frac{aT-2\left( p-c \right)\tan A}{2\left( p-c \right)T},  
\end{aligned}\]
\[\begin{aligned}
	\alpha _{B}^{H{{E}_{C}}}&=\alpha _{B}^{{{E}_{C}}}-\alpha _{B}^{H}=\frac{b}{2\left( p-c \right)}-\frac{\tan B}{T} \\ 
	& =\frac{bT-2\left( p-c \right)\tan B}{2\left( p-c \right)T},  
\end{aligned}\]
\[\begin{aligned}
	\alpha _{C}^{H{{E}_{C}}}&=\alpha _{C}^{{{E}_{C}}}-\alpha _{C}^{H}=-\frac{c}{2\left( p-c \right)}-\frac{\tan C}{T} \\ 
	& =-\frac{cT+2\left( p-c \right)\tan C}{2\left( p-c \right)T},  
\end{aligned}\]
\[\begin{aligned}
	\alpha _{A}^{H{{E}_{C}}}\alpha _{B}^{H{{E}_{C}}}&=\frac{aT-2\left( p-c \right)\tan A}{2\left( p-c \right)T}\cdot \frac{bT-2\left( p-c \right)\tan B}{2\left( p-c \right)T} \\ 
	& =\frac{\left( aT-2\left( p-c \right)\tan A \right)\left( bT-2\left( p-c \right)\tan B \right)}{4{{\left( p-c \right)}^{2}}{{T}^{2}}},  
\end{aligned}\]
\[\begin{aligned}
	\alpha _{B}^{H{{E}_{C}}}\alpha _{C}^{H{{E}_{C}}}&=-\frac{bT-2\left( p-c \right)\tan B}{2\left( p-c \right)T}\cdot \frac{cT+2\left( p-c \right)\tan C}{2\left( p-c \right)T} \\ 
	& =-\frac{\left( bT-2\left( p-c \right)\tan B \right)\left( cT+2\left( p-c \right)\tan C \right)}{4{{\left( p-c \right)}^{2}}{{T}^{2}}},  
\end{aligned}\]
\[\begin{aligned}
	\alpha _{C}^{H{{E}_{C}}}\alpha _{A}^{H{{E}_{C}}}&=-\frac{cT+2\left( p-c \right)\tan C}{2\left( p-c \right)T}\cdot \frac{aT-2\left( p-c \right)\tan A}{2\left( p-c \right)T} \\ 
	& =-\frac{\left( cT+2\left( p-c \right)\tan C \right)\left( aT-2\left( p-c \right)\tan A \right)}{4{{\left( p-c \right)}^{2}}{{T}^{2}}}.  
\end{aligned}\]

Therefore,
\[\begin{aligned}
	H{{E}_{C}}^{2}&=-\alpha _{B}^{H{{E}_{C}}}\alpha _{C}^{H{{E}_{C}}}{{a}^{2}}-\alpha _{C}^{H{{E}_{C}}}\alpha _{A}^{H{{E}_{C}}}{{b}^{2}}-\alpha _{A}^{H{{E}_{C}}}\alpha _{B}^{H{{E}_{C}}}{{c}^{2}} \\ 
	& =\frac{1}{4{{T}^{2}}{{\left( p-c \right)}^{2}}}\left( \begin{aligned}
		& \left( bT-2\left( p-c \right)\tan B \right)\left( cT+2\left( p-c \right)\tan C \right){{a}^{2}} \\ 
		& +\left( cT+2\left( p-c \right)\tan C \right)\left( aT-2\left( p-c \right)\tan A \right){{b}^{2}} \\ 
		& -\left( aT-2\left( p-c \right)\tan A \right)\left( bT-2\left( p-c \right)\tan B \right){{c}^{2}} \\ 
	\end{aligned} \right).  
\end{aligned}\]

\subsection{DTICs between excenters}\label{Subsec13.2.11}
\subsubsection{DTICs between excenter ${{E}_{A}}$ corresponding to $\angle A$ and excenter ${{E}_{B}}$ corresponding to $\angle B$}
By using Chapter \ref{Ch8} and theorem \ref{thm:Thm12.2.1}, the following formulas are obtained:
\[\begin{aligned}
	\alpha _{A}^{{{E}_{A}}{{E}_{B}}}&=\alpha _{A}^{{{E}_{B}}}-\alpha _{A}^{{{E}_{A}}} \\ 
	& =\frac{a}{2\left( p-b \right)}+\frac{a}{2\left( p-a \right)}=\frac{a\left( p-a \right)+a\left( p-b \right)}{2\left( p-a \right)\left( p-b \right)} \\ 
	& =\frac{a\left( 2p-a-b \right)}{2\left( p-a \right)\left( p-b \right)}=\frac{ca}{2\left( p-a \right)\left( p-b \right)},  
\end{aligned}\]
\[\begin{aligned}
	\alpha _{B}^{{{E}_{A}}{{E}_{B}}}&=\alpha _{B}^{{{E}_{B}}}-\alpha _{B}^{{{E}_{A}}} \\ 
	& =-\frac{b}{2\left( p-b \right)}-\frac{b}{2\left( p-a \right)}=-\frac{b\left( 2p-a-b \right)}{2\left( p-a \right)\left( p-b \right)} \\ 
	& =-\frac{bc}{2\left( p-a \right)\left( p-b \right)},  
\end{aligned}\]
\[\begin{aligned}
	\alpha _{C}^{{{E}_{A}}{{E}_{B}}}&=\alpha _{C}^{{{E}_{B}}}-\alpha _{C}^{{{E}_{A}}} \\ 
	& =\frac{c}{2\left( p-b \right)}-\frac{c}{2\left( p-a \right)}=\frac{c\left( b-a \right)}{2\left( p-a \right)\left( p-b \right)}.  
\end{aligned}\]

Therefore,
\[\begin{aligned}
	{{E}_{A}}{{E}_{B}}^{2}&=-\alpha _{B}^{{{E}_{A}}{{E}_{B}}}\alpha _{C}^{{{E}_{A}}{{E}_{B}}}{{a}^{2}}-\alpha _{C}^{{{E}_{A}}{{E}_{B}}}\alpha _{A}^{{{E}_{A}}{{E}_{B}}}{{b}^{2}}-\alpha _{A}^{{{E}_{A}}{{E}_{B}}}\alpha _{B}^{{{E}_{A}}{{E}_{B}}}{{c}^{2}} \\ 
	& =\frac{b{{c}^{2}}\left( b-a \right){{a}^{2}}-{{c}^{2}}\left( b-a \right)a{{b}^{2}}+ab{{c}^{4}}}{4{{\left( p-a \right)}^{2}}{{\left( p-b \right)}^{2}}},  
\end{aligned}\]
\[\begin{aligned}
	{{E}_{A}}{{E}_{B}}^{2}&=\frac{ab{{c}^{2}}\left( b-a \right)\left( a-b \right)+ab{{c}^{4}}}{4{{\left( p-a \right)}^{2}}{{\left( p-b \right)}^{2}}}=\frac{ab{{c}^{2}}\left( {{c}^{2}}-{{\left( a-b \right)}^{2}} \right)}{4{{\left( p-a \right)}^{2}}{{\left( p-b \right)}^{2}}} \\ 
	& =\frac{ab{{c}^{2}}\left( c+a-b \right)\left( c-a+b \right)}{4{{\left( p-a \right)}^{2}}{{\left( p-b \right)}^{2}}}=\frac{ab{{c}^{2}}}{\left( p-a \right)\left( p-b \right)}.  
\end{aligned}\]

Therefore,
\[{{E}_{A}}{{E}_{B}}=c\cdot \sqrt{\frac{ab}{\left( p-a \right)\left( p-b \right)}}.\]	

\subsubsection{DTICs between excenter ${{E}_{B}}$ corresponding to $\angle B$ and excenter ${{E}_{C}}$ corresponding to $\angle C$}
By using Chapter \ref{Ch8} and theorem \ref{thm:Thm12.2.1}, the following formulas are obtained:
\[\begin{aligned}
	\alpha _{A}^{{{E}_{B}}{{E}_{C}}}=\alpha _{A}^{{{E}_{C}}}-\alpha _{A}^{{{E}_{B}}}
	=\frac{a}{2\left( p-c \right)}-\frac{a}{2\left( p-b \right)}=\frac{a\left( c-b \right)}{2\left( p-b \right)\left( p-c \right)},  
\end{aligned}\]
\[\begin{aligned}
	\alpha _{B}^{{{E}_{B}}{{E}_{C}}}=\alpha _{B}^{{{E}_{C}}}-\alpha _{B}^{{{E}_{B}}}
	=\frac{b}{2\left( p-c \right)}+\frac{b}{2\left( p-b \right)}=\frac{ab}{2\left( p-b \right)\left( p-c \right)},  
\end{aligned}\]
\[\begin{aligned}
	\alpha _{C}^{{{E}_{B}}{{E}_{C}}}=\alpha _{C}^{{{E}_{C}}}-\alpha _{C}^{{{E}_{B}}} =-\frac{c}{2\left( p-c \right)}-\frac{c}{2\left( p-b \right)}
	=-\frac{ca}{2\left( p-b \right)\left( p-c \right)}.  
\end{aligned}\]

Therefore,
\[\begin{aligned}
	{{E}_{B}}{{E}_{C}}^{2}&=-\alpha _{B}^{{{E}_{B}}{{E}_{C}}}\alpha _{C}^{{{E}_{B}}{{E}_{C}}}{{a}^{2}}-\alpha _{C}^{{{E}_{B}}{{E}_{C}}}\alpha _{A}^{{{E}_{B}}{{E}_{C}}}{{b}^{2}}-\alpha _{A}^{{{E}_{B}}{{E}_{C}}}\alpha _{B}^{{{E}_{B}}{{E}_{C}}}{{c}^{2}} \\ 
	& =-\frac{-{{a}^{2}}bc{{a}^{2}}-c{{a}^{2}}\left( c-b \right){{b}^{2}}+{{a}^{2}}b\left( c-b \right){{c}^{2}}}{4{{\left( p-b \right)}^{2}}{{\left( p-c \right)}^{2}}},  
\end{aligned}\]

\[\begin{aligned}
	{{E}_{B}}{{E}_{C}}^{2}&=\frac{{{a}^{2}}bc\left( {{a}^{2}}+b\left( c-b \right)-c\left( c-b \right) \right)}{4{{\left( p-b \right)}^{2}}{{\left( p-c \right)}^{2}}}=\frac{{{a}^{2}}bc\left( {{a}^{2}}-{{\left( b-c \right)}^{2}} \right)}{4{{\left( p-b \right)}^{2}}{{\left( p-c \right)}^{2}}} \\ 
	& =\frac{{{a}^{2}}bc\left( a+b-c \right)\left( a-b+c \right)}{4{{\left( p-b \right)}^{2}}{{\left( p-c \right)}^{2}}}=\frac{{{a}^{2}}bc}{\left( p-b \right)\left( p-c \right)}.  
\end{aligned}\]

Therefore,
\[{{E}_{B}}{{E}_{C}}=a\cdot \sqrt{\frac{bc}{\left( p-b \right)\left( p-c \right)}}.\]	
\subsubsection{DTICs between excenter ${{E}_{C}}$ corresponding to $\angle C$ and excenter ${{E}_{A}}$ corresponding to $\angle A$}
By using Chapter \ref{Ch8} and theorem \ref{thm:Thm12.2.1}, the following formulas are obtained:
\[\begin{aligned}
	\alpha _{A}^{{{E}_{C}}{{E}_{A}}}=\alpha _{A}^{{{E}_{A}}}-\alpha _{A}^{{{E}_{C}}}  =-\frac{a}{2\left( p-a \right)}-\frac{a}{2\left( p-c \right)}=-\frac{ab}{2\left( p-c \right)\left( p-a \right)},  
\end{aligned}\]
\[\begin{aligned}
	\alpha _{B}^{{{E}_{C}}{{E}_{A}}}=\alpha _{B}^{{{E}_{A}}}-\alpha _{B}^{{{E}_{C}}} =\frac{b}{2\left( p-a \right)}-\frac{b}{2\left( p-c \right)}=\frac{b\left( a-c \right)}{2\left( p-c \right)\left( p-a \right)},  
\end{aligned}\]
\[\begin{aligned}
	\alpha _{C}^{{{E}_{C}}{{E}_{A}}}=\alpha _{C}^{{{E}_{A}}}-\alpha _{C}^{{{E}_{C}}}=\frac{c}{2\left( p-a \right)}+\frac{c}{2\left( p-c \right)}=\frac{bc}{2\left( p-c \right)\left( p-a \right)},  
\end{aligned}\]

Therefore,
\[\begin{aligned}
	{{E}_{C}}{{E}_{A}}^{2}& =-\alpha _{B}^{{{E}_{C}}{{E}_{A}}}\alpha _{C}^{{{E}_{C}}{{E}_{A}}}{{a}^{2}}-\alpha _{C}^{{{E}_{C}}{{E}_{A}}}\alpha _{A}^{{{E}_{C}}{{E}_{A}}}{{b}^{2}}-\alpha _{A}^{{{E}_{C}}{{E}_{A}}}\alpha _{B}^{{{E}_{C}}{{E}_{A}}}{{c}^{2}} \\ 
	& =-\frac{{{b}^{2}}c\left( a-c \right){{a}^{2}}-a{{b}^{2}}c{{b}^{2}}-a{{b}^{2}}\left( a-c \right){{c}^{2}}}{4{{\left( p-c \right)}^{2}}{{\left( p-a \right)}^{2}}},  
\end{aligned}\]
\[\begin{aligned}
	{{E}_{C}}{{E}_{A}}^{2}& =\frac{{{b}^{2}}ca\left( -a\left( a-c \right)+{{b}^{2}}+c\left( a-c \right) \right)}{4{{\left( p-c \right)}^{2}}{{\left( p-a \right)}^{2}}}=\frac{{{b}^{2}}ca\left( {{b}^{2}}-{{\left( c-a \right)}^{2}} \right)}{4{{\left( p-c \right)}^{2}}{{\left( p-a \right)}^{2}}} \\ 
	& =\frac{{{b}^{2}}ca\left( b+c-a \right)\left( b-c+a \right)}{4{{\left( p-c \right)}^{2}}{{\left( p-a \right)}^{2}}}=\frac{{{b}^{2}}ca}{\left( p-c \right)\left( p-a \right)}.  
\end{aligned}\]

Therefore,
\[{{E}_{C}}{{E}_{A}}=b\cdot \sqrt{\frac{ca}{\left( p-c \right)\left( p-a \right)}}.\]	

For equilateral triangles, there is the following formula: ${{E}_{A}}{{E}_{B}}={{E}_{B}}{{E}_{C}}={{E}_{C}}{{E}_{A}}=2a$.

So far, the 21 distances between the centers of a triangle have been obtained, which are the direct results of the theorems proposed by the author.


\chapter{Area of intersecting center triangle}\label{Ch14}
\thispagestyle{empty}

This chapter deals with the area of the intersecting center triangle. Triangle is one of the basic plane geometric figures. In addition to distance and angle, the area of triangle is also a very important measure. Using the method of vector product, this chapter first studies the area formula of the intersecting center triangle in general, and calculates the height of the intersecting center triangle. The general results are applied to the some ICs of a triangle and the corresponding results are obtained.

\section{Definition of areal vector of intersecting center triangle}\label{Sec14.1}
For a given $\triangle ABC$, its IC and two vertices in $\triangle ABC$ form a triangle, which is called \textbf{intersecting center triangle} (abbreviated as ICT). For example, in Figure \ref{fig:tu6.1.1}, $\triangle ABP$, $\triangle BCP$, and $\triangle CAP$ are all ICTs. $\triangle ABC$ is called the primitive triangle of the ICT. The $\triangle ABP$ is also called the ICT corresponding to the vertex $A$ of the original triangle, and so on.

For $\triangle ABC$, half of the vector product (cross product) of its two sides can represent the area vector. The magnitude of the area vector is the area (scalar) of $\triangle ABC$, and the direction of the area vector is determined by the right-handed system. It is necessary to make some rules for the vector product of the ICT. For example, the definition of the area vector of $\triangle ABC$ is as follows.

\begin{definition}{Definition of area vector of triangle}{Defi14.1.1}\label{Defi14.1.1} 
	Given a $\triangle ABC$, the area vector of the $\triangle ABC$ is defined as
	\[{{\mathbf{S}}_{ABC}}=\frac{1}{2}\overrightarrow{AB}\times \overrightarrow{AC}=\frac{1}{2}\overrightarrow{BC}\times \overrightarrow{BA}=\frac{1}{2}\overrightarrow{CA}\times \overrightarrow{CB}.\]	
\end{definition}

For the three ICTs: $\triangle ABP$, $\triangle BCP$, $\triangle CAP$ of $\triangle ABC$, the area vectors of the three ICTs are defined as follows.

\begin{definition}{Definition of area vector of intersecting center triangle}{Defi14.1.2}\label{Defi14.1.2} 
	\begin{flalign*}
		{{\mathbf{S}}_{ABP}}=\frac{1}{2}\overrightarrow{AB}\times \overrightarrow{AP}, {{\mathbf{S}}_{BCP}}=\frac{1}{2}\overrightarrow{BC}\times \overrightarrow{BP}, {{\mathbf{S}}_{CAP}}=\frac{1}{2}\overrightarrow{CA}\times \overrightarrow{CP}.
	\end{flalign*}
\end{definition}

The starting vectors of the above formulas are $\overrightarrow{AB}$, $\overrightarrow{BC}$, $\overrightarrow{CA}$ respectively, and the corresponding ending vectors are $\overrightarrow{AP}$, $\overrightarrow{BP}$, $\overrightarrow{CP}$ respectively. It should be noted that the three starting points $A$, $B$, $C$ of the starting vector $\overrightarrow{AB}$, $\overrightarrow{BC}$, $\overrightarrow{CA}$ should be arranged in anti-clockwise.

According to the above definition, the magnitude of area vector of the triangle $\triangle ABC$ is the area of triangle $\triangle ABC$. The unit vector of area vector is as follows:
\[{{\mathbf{n}}_{ABC}}=\frac{{{\mathbf{S}}_{ABC}}}{\left| {{\mathbf{S}}_{ABC}} \right|}=\frac{\overrightarrow{AB}\times \overrightarrow{AC}}{\left| \overrightarrow{AB}\times \overrightarrow{AC} \right|}=\frac{\overrightarrow{BC}\times \overrightarrow{BA}}{\left| \overrightarrow{BC}\times \overrightarrow{BA} \right|}=\frac{\overrightarrow{CA}\times \overrightarrow{CB}}{\left| \overrightarrow{CA}\times \overrightarrow{CB} \right|},\]
\[{{\mathbf{n}}_{ABP}}=\frac{{{\mathbf{S}}_{ABP}}}{\left| {{\mathbf{S}}_{ABP}} \right|}=\frac{\overrightarrow{AB}\times \overrightarrow{AP}}{\left| \overrightarrow{AB}\times \overrightarrow{AP} \right|},\]	
\[{{\mathbf{n}}_{BCP}}=\frac{{{\mathbf{S}}_{BCP}}}{\left| {{\mathbf{S}}_{BCP}} \right|}=\frac{\overrightarrow{BC}\times \overrightarrow{BP}}{\left| \overrightarrow{BC}\times \overrightarrow{BP} \right|},\]	
\[{{\mathbf{n}}_{CAP}}=\frac{{{\mathbf{S}}_{CAP}}}{\left| {{\mathbf{S}}_{CAP}} \right|}=\frac{\overrightarrow{CA}\times \overrightarrow{CP}}{\left| \overrightarrow{CA}\times \overrightarrow{CP} \right|}.\]	

\section{Theorem of areal vector of intersecting center triangle}\label{Sec14.2}
Now we use the method of vector product to calculate the area of the ICT.

\begin{theorem}{Theorem of area vector of intersecting center triangle, Daiyuan Zhang}{Thm14.2.1}\label{Thm14.2.1} 
	For a $\triangle ABC$, $P\in {{\pi }_{ABC}}$, the area vector of the intersecting center triangle is equal to the area vector of the original triangle multiplied by the frame component of the corresponding vertex. i.e.,
	\begin{flalign}\label{Eq14.2.1}
		{{\mathbf{S}}_{ABP}}=\alpha _{C}^{P}{{\mathbf{S}}_{ABC}}, {{\mathbf{S}}_{BCP}}=\alpha _{A}^{P}{{\mathbf{S}}_{ABC}}, {{\mathbf{S}}_{CAP}}=\alpha _{B}^{P}{{\mathbf{S}}_{ABC}},	
	\end{flalign}	
	and
	\[{{\mathbf{S}}_{ABC}}={{\mathbf{S}}_{ABP}}+{{\mathbf{S}}_{BCP}}+{{\mathbf{S}}_{CAP}}.\]		
\end{theorem}

\begin{proof}
	\[\begin{aligned}
		 {{\mathbf{S}}_{ABP}}&=\frac{1}{2}\overrightarrow{AB}\times \overrightarrow{AP}=\frac{1}{2}\overrightarrow{AB}\times \left( \alpha _{B}^{P}\overrightarrow{AB}+\alpha _{C}^{P}\overrightarrow{AC} \right) \\ 
		& =\frac{1}{2}\alpha _{C}^{P}\overrightarrow{AB}\times \overrightarrow{AC}=\alpha _{C}^{P}{{\mathbf{S}}_{ABC}}.  
	\end{aligned}\]	
	
	Similarly,
	\[\begin{aligned}
		{{\mathbf{S}}_{BCP}}&=\frac{1}{2}\overrightarrow{BC}\times \overrightarrow{BP}=\frac{1}{2}\overrightarrow{BC}\times \left( \alpha _{C}^{P}\overrightarrow{BC}+\alpha _{A}^{P}\overrightarrow{BA} \right) \\ 
		& =\frac{1}{2}\alpha _{A}^{P}\overrightarrow{BC}\times \overrightarrow{BA}=\alpha _{A}^{P}{{\mathbf{S}}_{ABC}},  
	\end{aligned}\]	
	\[\begin{aligned}
		{{\mathbf{S}}_{CAP}}&=\frac{1}{2}\overrightarrow{CA}\times \overrightarrow{CP}=\frac{1}{2}\overrightarrow{CA}\times \left( \alpha _{A}^{P}\overrightarrow{CA}+\alpha _{B}^{P}\overrightarrow{CB} \right) \\ 
		& =\frac{1}{2}\alpha _{B}^{P}\overrightarrow{CA}\times \overrightarrow{CB}=\alpha _{B}^{P}{{\mathbf{S}}_{ABC}},  
	\end{aligned}\]	
		\[{{\mathbf{S}}_{ABP}}+{{\mathbf{S}}_{BCP}}+{{\mathbf{S}}_{CAP}} =\alpha _{C}^{P}{{\mathbf{S}}_{ABC}}+\alpha _{A}^{P}{{\mathbf{S}}_{ABC}}+\alpha _{B}^{P}{{\mathbf{S}}_{ABC}}={{\mathbf{S}}_{ABC}}.\] 
\end{proof}
\hfill $\square$\par
By introducing the following notation, ${{S}_{ABC}}=\left| {{\mathbf{S}}_{ABC}} \right|$, ${{S}_{ABP}}=\left| {{\mathbf{S}}_{ABP}} \right|$, ${{S}_{BCP}}=\left| {{\mathbf{S}}_{BCP}} \right|$, ${{S}_{CAP}}=\left| {{\mathbf{S}}_{CAP}} \right|$, we get the following result:

\begin{theorem}{Magnitude of area vector of intersecting center triangle, Daiyuan Zhang}{Thm14.2.2}\label{Thm14.2.2} 
	For a $\triangle ABC$, $P\in {{\pi }_{ABC}}$, we have
	\[\begin{aligned}
		{{S}_{ABC}}&=\left| {{\mathbf{S}}_{ABC}} \right|=\left| {{\mathbf{S}}_{ABP}}+{{\mathbf{S}}_{BCP}}+{{\mathbf{S}}_{CAP}} \right| \\ 
		& =\frac{1}{2}\left| \alpha _{C}^{P}bc\sin A+\alpha _{A}^{P}ca\sin B+\alpha _{B}^{P}ab\sin C \right|  
	\end{aligned}.\]	

	And
	\begin{flalign*}
		{{S}_{ABP}}=\left| \alpha _{C}^{P} \right|{{S}_{ABC}}, {{S}_{BCP}}=\left| \alpha _{A}^{P} \right|{{S}_{ABC}}, {{S}_{CAP}}=\left| \alpha _{B}^{P} \right|{{S}_{ABC}}.
	\end{flalign*}
\end{theorem}

\begin{proof}
	From theorem \ref{thm:Thm14.2.1}, the following formula is obtained:
	\begin{align*}
		{{\mathbf{S}}_{ABC}}&={{\mathbf{S}}_{ABP}}+{{\mathbf{S}}_{BCP}}+{{\mathbf{S}}_{CAP}} \\ 
		& =\alpha _{C}^{P}{{\mathbf{S}}_{ABC}}+\alpha _{A}^{P}{{\mathbf{S}}_{ABC}}+\alpha _{B}^{P}{{\mathbf{S}}_{ABC}} \\ 
		& =\left( \alpha _{C}^{P}{{S}_{ABC}}+\alpha _{A}^{P}{{S}_{ABC}}+\alpha _{B}^{P}{{S}_{ABC}} \right){{\mathbf{n}}_{ABC}} \\ 
		& =\frac{1}{2}\left( \alpha _{C}^{P}bc\sin A+\alpha _{A}^{P}ca\sin B+\alpha _{B}^{P}ab\sin C \right){{\mathbf{n}}_{ABC}}.  
	\end{align*}
	
	Therefore
	\[\begin{aligned}
		{{S}_{ABC}}&=\left| {{\mathbf{S}}_{ABC}} \right|=\left| {{\mathbf{S}}_{ABP}}+{{\mathbf{S}}_{BCP}}+{{\mathbf{S}}_{CAP}} \right| \\ 
		& =\frac{1}{2}\left| \alpha _{C}^{P}bc\sin A+\alpha _{A}^{P}ca\sin B+\alpha _{B}^{P}ab\sin C \right|.  
	\end{aligned}\]	
	
	Take the magnitude on both sides of formula (\ref{Eq14.2.1}) to get the following results:
	\begin{flalign*}
		{{S}_{ABP}}=\left| \alpha _{C}^{P} \right|{{S}_{ABC}}, {{S}_{BCP}}=\left| \alpha _{A}^{P} \right|{{S}_{ABC}}, {{S}_{CAP}}=\left| \alpha _{B}^{P} \right|{{S}_{ABC}}.
	\end{flalign*}
\end{proof}
\hfill $\square$\par

For the inner IC, if each frame component is a positive number, there will be
\begin{flalign*}
	{{S}_{ABP}}=\alpha _{C}^{P}{{S}_{ABC}}, {{S}_{BCP}}=\alpha _{A}^{P}{{S}_{ABC}}, {{S}_{CAP}}=\alpha _{B}^{P}{{S}_{ABC}}.
\end{flalign*}
\[\begin{aligned}
	{{S}_{ABC}}&={{S}_{ABP}}+{{S}_{BCP}}+{{S}_{CAP}} \\ 
	& =\frac{1}{2}\left( \alpha _{C}^{P}bc\sin A+\alpha _{A}^{P}ca\sin B+\alpha _{B}^{P}ab\sin C \right).  
\end{aligned}\]	

\section{Application of area vector of intersecting center triangle}\label{Sec14.3}
From theorem \ref{thm:Thm14.2.1}, the area formulas of special ICTs can be obtained directly.

For the centroid,
\[{{\mathbf{S}}_{ABG}}=\alpha _{C}^{G}{{\mathbf{S}}_{ABC}}=\frac{1}{3}{{\mathbf{S}}_{ABC}},\]	
\[{{\mathbf{S}}_{BCG}}=\alpha _{A}^{G}{{\mathbf{S}}_{ABC}}=\frac{1}{3}{{\mathbf{S}}_{ABC}},\]	
\[{{\mathbf{S}}_{CAG}}=\alpha _{B}^{G}{{\mathbf{S}}_{ABC}}=\frac{1}{3}{{\mathbf{S}}_{ABC}}.\]	

For the incenter,
\[{{\mathbf{S}}_{ABI}}=\alpha _{C}^{I}{{\mathbf{S}}_{ABC}}=\frac{c}{a+b+c}{{\mathbf{S}}_{ABC}},\]	
\[{{\mathbf{S}}_{BCI}}=\alpha _{A}^{I}{{\mathbf{S}}_{ABC}}=\frac{a}{a+b+c}{{\mathbf{S}}_{ABC}},\]	
\[{{\mathbf{S}}_{CAI}}=\alpha _{B}^{I}{{\mathbf{S}}_{ABC}}=\frac{b}{a+b+c}{{\mathbf{S}}_{ABC}}.\]	

For the orthocenter,
\[{{\mathbf{S}}_{ABH}}=\alpha _{C}^{H}{{\mathbf{S}}_{ABC}}=\frac{\tan C}{\tan A+\tan B+\tan C}{{\mathbf{S}}_{ABC}},\]	
\[{{\mathbf{S}}_{BCH}}=\alpha _{A}^{H}{{\mathbf{S}}_{ABC}}=\frac{\tan A}{\tan A+\tan B+\tan C}{{\mathbf{S}}_{ABC}},\]	
\[{{\mathbf{S}}_{CAH}}=\alpha _{B}^{H}{{\mathbf{S}}_{ABC}}=\frac{\tan B}{\tan A+\tan B+\tan C}{{\mathbf{S}}_{ABC}}.\]	

For the circumcenter,
\[{{\mathbf{S}}_{ABQ}}=\alpha _{C}^{Q}{{\mathbf{S}}_{ABC}}=\frac{\sin 2C}{\sin 2A+\sin 2B+\sin 2C}{{\mathbf{S}}_{ABC}},\]	
\[{{\mathbf{S}}_{BCQ}}=\alpha _{A}^{Q}{{\mathbf{S}}_{ABC}}=\frac{\sin 2A}{\sin 2A+\sin 2B+\sin 2C}{{\mathbf{S}}_{ABC}},\]	
\[{{\mathbf{S}}_{CAQ}}=\alpha _{B}^{Q}{{\mathbf{S}}_{ABC}}=\frac{\sin 2B}{\sin 2A+\sin 2B+\sin 2C}{{\mathbf{S}}_{ABC}}.\]

\section{Altitude of intersecting center triangle}\label{Sec14.4}
The altitudes of the ICTs $\triangle ABP$, $\triangle BCP$, $\triangle CAP$ are denoted as $h_{AB}^{P}$, $h_{BC}^{P}$, $h_{CA}^{P}$ respectively, according to theorem \ref{thm:Thm14.2.2}, we have 
\[h_{AB}^{P}=\frac{2{{S}_{ABP}}}{AB}=\frac{2\left| \alpha _{C}^{P} \right|}{c}{{S}_{ABC}}.\]	

Similarly,
\[h_{BC}^{P}=\frac{2{{S}_{BCP}}}{BC}=\frac{2\left| \alpha _{A}^{P} \right|}{a}{{S}_{ABC}},\]	
\[h_{CA}^{P}=\frac{2{{S}_{CAP}}}{CA}=\frac{2\left| \alpha _{B}^{P} \right|}{b}{{S}_{ABC}}.\]	

For the seven centers of a triangle, the altitude of the corresponding ICT can be obtained directly according to the above three formulas.


\chapter{New inequalities}\label{Ch15}
\thispagestyle{empty}

%
%
In addition to the previous results, some new inequalities can be obtained in Intercenter Geometry. In this chapter, I publish some new geometric inequalities for the first time.

Intercenter Geometry is good at calculating the distance between two points. Many new geometric inequalities can be created based on the principle that the distance between two points is non-negative. It can be said that Intercenter Geometry is a “generator” that creates a large class of geometric inequalities. In Plane Intercenter Geometry, geometric inequality mainly involves geometric inequalities on triangle.

Since the distance between the ICs is non-negative, according to the theorem of distance between origin and intersecting center on triangular frame (theorem \ref{thm:Thm12.1.1}), the theorem of distance between origin and intersecting center on
frame of circumcenter (theorem \ref{thm:Thm12.1.4}), and the theorem of distance between two intersecting centers (theorem \ref{thm:Thm12.2.1}), the following new inequalities can be obtained directly. Some traditional well-known inequalities are the corollaries of these new inequalities.

This chapter is divided into two sections, the geometric inequality of DOIC and the geometric inequality of DTICs.

\section{Geometric inequality of distance between origin and intersecting center}\label{Sec15.1}
\begin{theorem}{Geometric inequality of distance between origin and intersecting center, Daiyuan Zhang}{Thm15.1.1}\label{Thm15.1.1} 
	Let point $\,O$ be the origin of frame $\left( O;A,B,C \right)$, point $P$ be the intersecting center on the plane of $\triangle ABC$, then
	\[\alpha _{A}^{P}O{{A}^{2}}+\alpha _{B}^{P}O{{B}^{2}}+\alpha _{C}^{P}O{{C}^{2}}\ge \left( \alpha _{B}^{P}\alpha _{C}^{P}{{a}^{2}}+\alpha _{C}^{P}\alpha _{A}^{P}{{b}^{2}}+\alpha _{A}^{P}\alpha _{B}^{P}{{c}^{2}} \right).\]	
	
	The above equality occurs if and only if the origin of the frame $O$ coincides with the intersecting center $P$, where $\alpha _{A}^{P}$, $\alpha _{B}^{P}$, $\alpha _{C}^{P}$ are the frame components of the intersecting center $P$.
\end{theorem}

\begin{proof}
	From theorem \ref{thm:Thm12.1.1}, then
	\[\begin{aligned}
		O{{P}^{2}}=\alpha _{A}^{P}O{{A}^{2}}+\alpha _{B}^{P}O{{B}^{2}}+\alpha _{C}^{P}O{{C}^{2}}-\left( \alpha _{B}^{P}\alpha _{C}^{P}{{a}^{2}}+\alpha _{C}^{P}\alpha _{A}^{P}{{b}^{2}}+\alpha _{A}^{P}\alpha _{B}^{P}{{c}^{2}} \right). 
	\end{aligned}\]
	
	Since $O{{P}^{2}}\ge 0$, it follows that	
	\[\alpha _{A}^{P}O{{A}^{2}}+\alpha _{B}^{P}O{{B}^{2}}+\alpha _{C}^{P}O{{C}^{2}}\ge \left( \alpha _{B}^{P}\alpha _{C}^{P}{{a}^{2}}+\alpha _{C}^{P}\alpha _{A}^{P}{{b}^{2}}+\alpha _{A}^{P}\alpha _{B}^{P}{{c}^{2}} \right).\]
	
	Obviously, the above equality holds if and only if $OP=0$, that is, if and only if the origin $O$ of the frame coincides with the intersecting center $P$.
\end{proof}
\hfill $\square$\par

\begin{theorem}{Inequality of circumcircle radius-intersecting center of triangle, Daiyuan Zhang}{Thm15.1.2}\label{Thm15.1.2} 
	Let the point $P$ be the intersecting center on the plane of $\triangle ABC$, and $R$ be the circumcircle radius of $\triangle ABC$, then
	\[{{R}^{2}}\ge \alpha _{B}^{P}\alpha _{C}^{P}{{a}^{2}}+\alpha _{C}^{P}\alpha _{A}^{P}{{b}^{2}}+\alpha _{A}^{P}\alpha _{B}^{P}{{c}^{2}}.\]	
	
	The above equality occurs if and only if the circumcenter $Q$ coincides with the intersecting center $P$, where $\alpha _{A}^{P}$, $\alpha _{B}^{P}$, $\alpha _{C}^{P}$ are the frame components of the intersecting center $P$.
\end{theorem}

\begin{proof}
	From the theorem of distance between origin and intersecting center (theorem \ref{thm:Thm12.1.4}), then
	\[Q{{P}^{2}}={{R}^{2}}-\left( \alpha _{B}^{P}\alpha _{C}^{P}{{a}^{2}}+\alpha _{C}^{P}\alpha _{A}^{P}{{b}^{2}}+\alpha _{A}^{P}\alpha _{B}^{P}{{c}^{2}} \right).\]	
	
	Since $O{{P}^{2}}\ge 0$, then
	
	\[{{R}^{2}}\ge \left( \alpha _{B}^{P}\alpha _{C}^{P}{{a}^{2}}+\alpha _{C}^{P}\alpha _{A}^{P}{{b}^{2}}+\alpha _{A}^{P}\alpha _{B}^{P}{{c}^{2}} \right).\]	
	
	Obviously, the above equality holds if and only if $QP=0$, that is, if and only if the origin $O$ of the frame coincides with the intersecting center $P$.
\end{proof}
\hfill $\square$\par

When the above theorem is applied to some special ICs of triangles, such as circumcenter, centroid, incenter and orthocenter, a series of inequalities will be obtained. According to Euclidean geometry, the circumcenter, centroid, incenter and orthocenter will coincide with each other if and only if the triangle is an equilateral triangle. So we get the following corollaries.

\begin{corollary}{Inequality of circumcircle radius-centroid of triangle}{Cor15.1.1}\label{Cor15.1.1} 
	Let $R$ be the circumcircle radius of $\triangle ABC$, then	
	\[{{R}^{2}}\ge \frac{1}{9}\left( {{a}^{2}}+{{b}^{2}}+{{c}^{2}} \right).\]	
	
	The above equality occurs if and only if $a=b=c$.
\end{corollary}

\begin{proof}
	It is directly obtained from theorem \ref{thm:Thm15.1.2} and subsection \ref{Subsec13.1.1}.
	
	From Euclidean geometry, if and only if $a=b=c$, the circumcenter and the centroid will coincide with each other, the above equality occurs.
\end{proof}
\hfill $\square$\par

The above inequality is a well-known result, which has been discussed in many materials. Here, it is only a corollary of theorem \ref{thm:Thm15.1.2}.

\begin{corollary}{Inequality of circumcircle radius-Incenter of triangle}{Cor15.1.2}\label{Cor15.1.2} 
	Let $R$ be the circumcircle radius of $\triangle ABC$, then	
	\[{{R}^{2}}\ge \frac{abc}{2p}=\frac{abc}{a+b+c}\]	
	
	The above equality occurs if and only if $a=b=c$.
\end{corollary}

\begin{proof}
	It is directly obtained from theorem \ref{thm:Thm15.1.2} and subsection \ref{Subsec13.1.2}.
	
	From Euclidean geometry, if and only if $a=b=c$, the circumcenter and the incenter will coincide with each other, the above equality occurs.
\end{proof}
\hfill $\square$\par

\begin{corollary}{Inequality of circumcircle radius-orthocenter of triangle}{Cor15.1.3}\label{Cor15.1.3} 
	Let $R$ be the circumcircle radius of $\triangle ABC$, then	
	\[{{R}^{2}}\ge \frac{\tan B\cdot \tan C\cdot {{a}^{2}}+\tan C\cdot \tan A\cdot {{b}^{2}}+\tan A\cdot \tan B\cdot {{c}^{2}}}{{{\left( \tan A+\tan B+\tan C \right)}^{2}}}\]
	
	The above equality occurs if and only if $a=b=c$.
\end{corollary}

\begin{proof}
	It is directly obtained from theorem \ref{thm:Thm15.1.2} and subsection \ref{Subsec13.1.3}.
	
	From Euclidean geometry, if and only if $a=b=c$, the circumcenter and the orthocenter will coincide with each other, the above equality occurs.
\end{proof}
\hfill $\square$\par

\section{Geometric inequality of distance between two intersecting centers}\label{Sec15.2}
\begin{theorem}{Geometric inequality of distance between two intersecting centers, Daiyuan Zhang}{Thm15.2.1}\label{Thm15.2.1} 
	Assume a given $\triangle ABC$, intersecting centers ${{P}_{1}}\in {{\pi }_{ABC}}$, ${{P}_{2}}\in {{\pi }_{ABC}}$, then	
	\[\alpha _{B}^{{{P}_{1}}{{P}_{2}}}\alpha _{C}^{{{P}_{1}}{{P}_{2}}}{{a}^{2}}+\alpha _{C}^{{{P}_{1}}{{P}_{2}}}\alpha _{A}^{{{P}_{1}}{{P}_{2}}}{{b}^{2}}+\alpha _{A}^{{{P}_{1}}{{P}_{2}}}\alpha _{B}^{{{P}_{1}}{{P}_{2}}}{{c}^{2}}\le 0.\]	
	
	The above equality occurs if and only if the two intersecting centers coincides with each other. Where
	\[\alpha _{A}^{{{P}_{1}}{{P}_{2}}}=\alpha _{A}^{{{P}_{2}}}-\alpha _{A}^{{{P}_{1}}},\]	
	\[\alpha _{B}^{{{P}_{1}}{{P}_{2}}}=\alpha _{B}^{{{P}_{2}}}-\alpha _{B}^{{{P}_{1}}},\]	
	\[\alpha _{C}^{{{P}_{1}}{{P}_{2}}}=\alpha _{C}^{{{P}_{2}}}-\alpha _{C}^{{{P}_{1}}},\]	
	and $\alpha _{A}^{{{P}_{1}}}$, $\alpha _{B}^{{{P}_{1}}}$, $\alpha _{C}^{{{P}_{1}}}$ are the frame components of the frame $\overrightarrow{OA}$, $\overrightarrow{OB}$, $\overrightarrow{OC}$ on the frame system $\left( O;A,B,C \right)$ at point ${P}_{1}$, respectively; $\alpha _{A}^{{{P}_{2}}}$, $\alpha _{B}^{{{P}_{2}}}$, $\alpha _{C}^{{{P}_{2}}}$ are the frame components of the frame $\overrightarrow{OA}$, $\overrightarrow{OB}$, $\overrightarrow{OC}$ on the frame system $\left( O;A,B,C \right)$ at point ${P}_{2}$ respectively.
\end{theorem}

\begin{proof}
	From theorem \ref{thm:Thm12.2.1}, we have
	\[{{P}_{1}}{{P}_{2}}^{2}=-\alpha _{B}^{{{P}_{1}}{{P}_{2}}}\alpha _{C}^{{{P}_{1}}{{P}_{2}}}{{a}^{2}}-\alpha _{C}^{{{P}_{1}}{{P}_{2}}}\alpha _{A}^{{{P}_{1}}{{P}_{2}}}{{b}^{2}}-\alpha _{A}^{{{P}_{1}}{{P}_{2}}}\alpha _{B}^{{{P}_{1}}{{P}_{2}}}{{c}^{2}}.\]
	
	Since ${{P}_{1}}{{P}_{2}}^{2}\ge 0$, then
	
	\[\alpha _{B}^{{{P}_{1}}{{P}_{2}}}\alpha _{C}^{{{P}_{1}}{{P}_{2}}}{{a}^{2}}+\alpha _{C}^{{{P}_{1}}{{P}_{2}}}\alpha _{A}^{{{P}_{1}}{{P}_{2}}}{{b}^{2}}+\alpha _{A}^{{{P}_{1}}{{P}_{2}}}\alpha _{B}^{{{P}_{1}}{{P}_{2}}}{{c}^{2}}\le 0\]	
	
	
	Obviously, the above equality holds if and only if ${{P}_{1}}{{P}_{2}}=0$, that is, if and only if the two intersecting centers ${{P}_{1}}$ and ${{P}_{2}}$ coincides with each other.	
\end{proof}
\hfill $\square$\par

\begin{corollary}{Circumcenter-Centroid inequality of triangle, Daiyuan Zhang}{Cor15.2.1}\label{Cor15.2.1} 
	Given a $\triangle ABC$, then
	\[\begin{aligned}
		& \left( W-3\sin 2B \right)\left( W-3\sin 2C \right){{a}^{2}} \\ 
		& +\left( W-3\sin 2C \right)\left( W-3\sin 2A \right){{b}^{2}} \\ 
		& +\left( W-3\sin 2A \right)\left( W-3\sin 2B \right){{c}^{2}}\le 0. \\ 
	\end{aligned}\]	
	
	The above equality occurs if and only if $a=b=c$.
\end{corollary}

\begin{proof}
	It is directly obtained from theorem \ref{thm:Thm12.2.1}, theorem \ref{thm:Thm15.2.1} and subsection \ref{Subsec13.2.1}.
	
	From Euclidean geometry, if and only if $a=b=c$, the circumcenter and the centroid will coincide with each other, the above equality occurs.
\end{proof}
\hfill $\square$\par

\begin{corollary}{Circumcenter-Incenter inequality of triangle, Daiyuan Zhang}{Cor15.2.2}\label{Cor15.2.2} 
	Given a $\triangle ABC$, then	
	\[\begin{aligned}
		& \left( bW-2p\sin 2B \right)\left( cW-2p\sin 2C \right){{a}^{2}} \\ 
		& +\left( cW-2p\sin 2C \right)\left( aW-2p\sin 2A \right){{b}^{2}} \\ 
		& +\left( aW-2p\sin 2A \right)\left( bW-2p\sin 2B \right){{c}^{2}}\le 0. \\ 
	\end{aligned}\]	
	
	The above equality occurs if and only if $a=b=c$.
\end{corollary}

\begin{proof}
	It is directly obtained from theorem \ref{thm:Thm12.2.1}, theorem \ref{thm:Thm15.2.1} and subsection \ref{Subsec13.2.2}.
	
	From Euclidean geometry, if and only if $a=b=c$, the circumcenter and the incenter will coincide with each other, the above equality occurs.
\end{proof}
\hfill $\square$\par

\begin{corollary}{Circumcenter-Orthocenter inequality of triangle, Daiyuan Zhang}{Cor15.2.3}\label{Cor15.2.3} 
	Given a $\triangle ABC$, then
	\[\begin{aligned}
		& \left( W\tan B-T\sin 2B \right)\left( W\tan C-T\sin 2C \right){{a}^{2}} \\ 
		& +\left( W\tan C-T\sin 2C \right)\left( W\tan A-T\sin 2A \right){{b}^{2}} \\ 
		& +\left( W\tan A-T\sin 2A \right)\left( W\tan B-T\sin 2B \right){{c}^{2}}\le 0. \\ 
	\end{aligned}\]
	
	The above equality occurs if and only if $a=b=c$.	
\end{corollary}

\begin{proof}
	It is directly obtained from theorem \ref{thm:Thm12.2.1}, theorem \ref{thm:Thm15.2.1} and subsection \ref{Subsec13.2.3}.
	
	From Euclidean geometry, if and only if $a=b=c$, the circumcenter and the orthocenter will coincide with each other, the above equality occurs.
\end{proof}
\hfill $\square$\par

\begin{corollary}{Centroid-Incenter inequality of triangle, Daiyuan Zhang}{Cor15.2.4}\label{Cor15.2.4} 
	Given a $\triangle ABC$, then	
	\[\left( 3b-2p \right)\left( 3c-2p \right){{a}^{2}}+\left( 3c-2p \right)\left( 3a-2p \right){{b}^{2}}+\left( 3a-2p \right)\left( 3b-2p \right){{c}^{2}}\le 0.\]
	
	The above equality occurs if and only if $a=b=c$. 
\end{corollary}

\begin{proof}
	It is directly obtained from theorem \ref{thm:Thm12.2.1}, theorem \ref{thm:Thm15.2.1} and subsection \ref{Subsec13.2.5}.
	
	From Euclidean geometry, if and only if $a=b=c$, the centroid and the incenter will coincide with each other, the above equality occurs.
\end{proof}
\hfill $\square$\par

\begin{corollary}{Centroid-Orthocenter inequality of triangle, Daiyuan Zhang}{Cor15.2.5}\label{Cor15.2.5} 
	Given a $\triangle ABC$, then
	\[\begin{aligned}
		& \left( 3\tan B-T \right)\left( 3\tan C-T \right){{a}^{2}} \\ 
		& +\left( 3\tan C-T \right)\left( 3\tan A-T \right){{b}^{2}} \\ 
		& +\left( 3\tan A-T \right)\left( 3\tan B-T \right){{c}^{2}}\le 0. \\ 
	\end{aligned}\]	
	
	The above equality occurs if and only if $a=b=c$.
\end{corollary}

\begin{proof}
	It is directly obtained from theorem \ref{thm:Thm12.2.1}, theorem \ref{thm:Thm15.2.1} and subsection \ref{Subsec13.2.6}.
	
	From Euclidean geometry, if and only if $a=b=c$, the centroid and the orthocenter will coincide with each other, the above equality occurs.
\end{proof}
\hfill $\square$\par

\begin{corollary}{Incenter-Orthocenter inequality of triangle, Daiyuan Zhang}{Cor15.2.6}\label{Cor15.2.6} 
	Given a $\triangle ABC$, then
	\[\begin{aligned}
		& \left( bT-2p\tan B \right)\left( cT-2p\tan C \right){{a}^{2}} \\ 
		& +\left( cT-2p\tan C \right)\left( aT-2p\tan A \right){{b}^{2}} \\ 
		& +\left( aT-2p\tan A \right)\left( bT-2p\tan B \right){{c}^{2}}\le 0. \\ 
	\end{aligned}\]
	
	The above equality occurs if and only if $a=b=c$.
\end{corollary}

\begin{proof}
	It is directly obtained from theorem \ref{thm:Thm12.2.1}, theorem \ref{thm:Thm15.2.1} and subsection \ref{Subsec13.2.8}.
	
	From Euclidean geometry, if and only if $a=b=c$, the incenter and the orthocenter will coincide with each other, the above equality occurs.
\end{proof}
\hfill $\square$\par

The above inequalities in the corollaries seem to be various and unrelated, but they are all special cases of the original geometric inequality of the distance between origin and intersecting center (theorem \ref{thm:Thm15.1.1}) and the geometric inequality of the distance between two intersecting centers (theorem \ref{thm:Thm15.2.1}) proposed by the author. That is to say, the original geometric inequality of the DOIC (theorem \ref{thm:Thm15.1.1}) and the geometric inequality of the DTICs (theorem \ref{thm:Thm15.2.1}) are the “father” of those inequalities in the corollaries.

%
%

Some of the inequalities in corollaries are well-known, but many of them are published by the author for the first time.

According to Intercenter Geometry, Some other inequalities can also be created, which will not be introduced here.

In a word, Intercenter Geometry has a unique charm in the construction of new geometric inequalities.


\chapter*{Part Two \\ Space Intercenter Geometry}
\chapter{Basic concepts of Space Intercenter Geometry}\label{Ch16}

From this chapter, we start to study Space Intercenter Geometry. Firstly, the basic concept of tetrahedron is briefly reviewed, and then the concept of IC of tetrahedron is given. Similar to Plane Intercenter Geometry, I put forward many new concepts in this chapter. It is hoped that readers will be familiar with these concepts as soon as possible, and the calculation methods of quantities related to these new concepts will be studied later.


In order to make readers understand why I put forward these new concepts, I need to briefly explain the original idea of Intercenter Geometry. The original idea of Plane Intercenter Geometry has been introduced in Chapter \ref{Ch4}. The following briefly describes the original idea of Space  Intercenter Geometry.


The idea of Space Intercenter Geometry is that the geometric quantities in space are represented by the lengths of the six edges of a given tetrahedron. 


Space Intercenter Geometry uses a unified vector method to solve the calculation problem of geometric quantities related to the distance between two points in the space.


The key problem of Space Intercenter Geometry is to establish a frame and calculate the frame components, so that each of the frame components can be expressed as a function of the lengths of the six edges of a given tetrahedron, and the frame components are independent of the origin of the frame.


Suppose a given tetrahedron $ABCD$ is embedded in space, a point $P$ is given in space, and $O$ is an arbitrary point in space, we hope to have the following expression:
\[\overrightarrow{OP}={\beta _{A}^{P}}\overrightarrow{OA}+{\beta _{B}^{P}}\overrightarrow{OB}+{\beta _{C}^{P}}\overrightarrow{OC}+{\beta _{D}^{P}}\overrightarrow{OD},\]	
and the frame components ${\beta _{A}^{P}}$, ${\beta _{B}^{P}}$, ${\beta _{C}^{P}}$ and ${\beta _{D}^{P}}$  are such quantities that they are only related to point $P$ and not to point $O$.


The question now is, can the frame components ${\beta _{A}^{P}}$, ${\beta _{B}^{P}}$, ${\beta _{C}^{P}}$ and ${\beta _{D}^{P}}$ be related to the fractional ratio of some line segments similar to the plane case? My answer is yes. The frame components ${\beta _{A}^{P}}$, ${\beta _{B}^{P}}$, ${\beta _{C}^{P}}$ and ${\beta _{D}^{P}}$ can be expressed by the frame components of the  triangles of the tetrahedron, and the frame components of the triangles can be expressed by the IRs of the triangles. Moreover, if the point $P$ is selected as some special points of the tetrahedron $ABCD$ (for example, centroid, incenter, circumcenter, and so on), then its frame components ${\beta _{A}^{P}}$, ${\beta _{B}^{P}}$, ${\beta _{C}^{P}}$ and ${\beta _{D}^{P}}$ can be expressed by the lengths of six edges of tetrahedron $ABCD$. This means that the frame components ${\beta _{A}^{P}}$, ${\beta _{B}^{P}}$, ${\beta _{C}^{P}}$ and ${\beta _{D}^{P}}$ are only related to the lengths of the six edges of the tetrahedron $ABCD$, so that the distance between the two points can be expressed as a function of the lengths of the six edges of the tetrahedron $ABCD$. Therefore, those geometric quantities that can be expressed by the distance between two points in space can also be expressed by the lengths of the six edges of the tetrahedron $ABCD$. This is the unique charm of Intercenter Geometry, which is unmatched by Euclidean geometry and analytical geometry.

\section{Basic concepts of tetrahedron}\label{Sec16.1}

A \textbf{tetrahedron} in three-dimensional space is denoted as tetrahedron $ABCD$. A tetrahedron has four \textbf{vertices} $A$, $B$, $C$, $D$; six \textbf{edges} $AB$, $AC$, $AD$, $BC$, $BD$, $CD$; and four \textbf{triangular faces}, which are the triangles $ABC$, $BCD$, $CDA$ and $DAB$ with their interiors. A triangular face is also called \textbf{face} for short.

Four vertices of a tetrahedron are not coplanar, each vertex of the tetrahedron has a unique face that does not pass through it, which is called the \textbf{opposite face} of the vertex, and the vertex is called the \textbf{opposite vertex} of the face. For example, the face $ABC$ is the opposite face of the vertex $D$, and the point $D$ is the opposite vertex of the face $ABC$. The edges $AB$ and $CD$ are called \textbf{opposite edges} of the tetrahedron. A tetrahedron has three pairs of opposite edges. Sometimes we say that the face $ABC$ is \textbf{opposite} the vertex $D$.

A tetrahedron is also called a \textbf{triangular pyramid}. If it is called a triangular pyramid, one of its four faces is called \textbf{base}, and the rest are called \textbf{lateral faces}. For example, when viewed as a triangular pyramid with $A$ as its vertex, it can be denoted as a triangular pyramid $A-BCD$. In this case, the opposite face of vertex $A$ is the base, and the remaining faces are the lateral faces. If the base of a triangular pyramid is an equilateral triangle and the projection of the vertex on the base is the center of the triangle, then such a triangular pyramid is called the \textbf{equilateral triangular pyramid}; A tetrahedron composed of four congruent regular triangles is called the \textbf{regular tetrahedron}.

\section{Concept of intersecting center of a tetrahedron}\label{Sec16.2}
Like triangles, tetrahedrons also have concepts such as IC and IR.

\textbf{Intersecting center of tetrahedron}: If the four lines passing through the four vertices of a given tetrahedron intersect at a point, the intersection point is called intersecting center of the tetrahedron (abbreviated as IC-T).

\begin{figure}[h]
	\centering
	\includegraphics[totalheight=6cm]{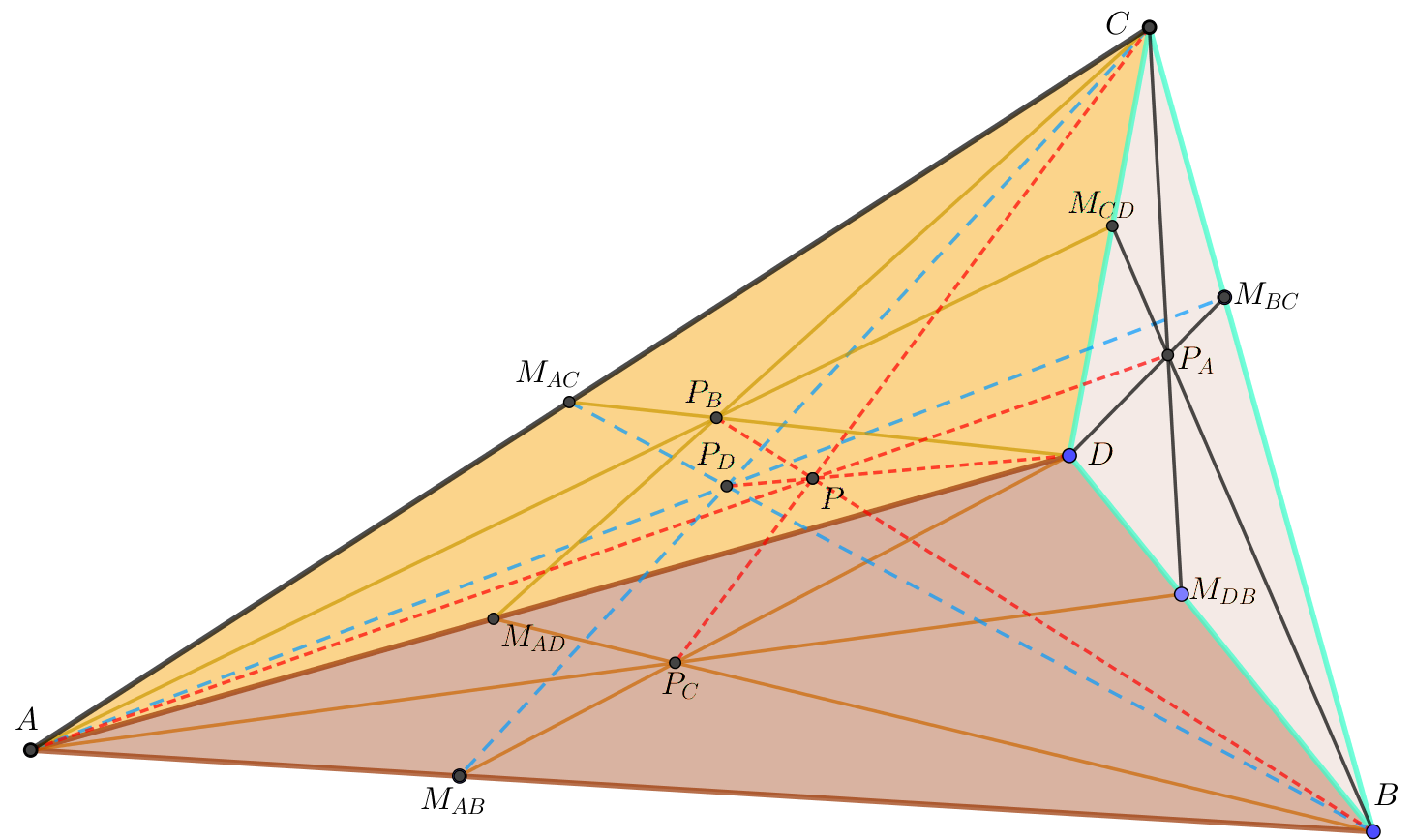}
	\caption{Intersecting center of tetrahedron} \label{fig:tu16.2.1}
\end{figure}

For example, the point $P$ in Figure \ref{fig:tu16.2.1} is the IC-T of the tetrahedron $ABCD$, and, 
	\[P=\overleftrightarrow{A{{P}_{A}}}\cap \overleftrightarrow{B{{P}_{B}}}\cap \overleftrightarrow{C{{P}_{C}}}\cap \overleftrightarrow{D{{P}_{D}}}.\]
	
When there is no confusion, IC-T is also called IC.

\textbf{Line through vertex and intersecting center} (abbreviated as LVIC): A straight line passing through the vertex and the IC-T of the tetrahedron.

\textbf{Intersecting center of face} (abbreviated as IC-F): The intersection of the LVIC and the base of the tetrahedron. If the intersection is the IC of the base triangle, it is called the intersecting center of face of the base.

In Figure \ref{fig:tu16.2.1}, four points ${{P}_{A}}$, ${{P}_{B}}$, ${{P}_{C}}$, ${{P}_{D}}$ are all IC-Fs of the tetrahedron, and,
\begin{flalign*}
	{{P}_{A}}=\overleftrightarrow{AP}\cap {P}_{\triangle BCD},
	{{P}_{B}}=\overleftrightarrow{BP}\cap {P}_{\triangle CDA},
	{{P}_{C}}=\overleftrightarrow{CP}\cap {P}_{\triangle DAB},
	{{P}_{D}}=\overleftrightarrow{DP}\cap {P}_{\triangle ABC}.
\end{flalign*}

In fact, ${{P}_{A}}$ is the IC of $\triangle BCD$, so the concept of IC discussed earlier can be applicable to the IC-F.

Of course, not every point has a IC-F. For example, points on the plane passing through the vertex $A$ and parallel to $\triangle BCD$ do not intersect the opposite face of the vertex $A$, so there is no IC-F for those points.

\textbf{Singular IC-T}: For a tetrahedron $ABCD$, if at least one IR of the IC-F of a point does not exist or is 0, then the point is called singular IC-T. The set of all singular ICs-T of the tetrahedron $ABCD$ is denoted as $\pi _{ABCD}^{\text{*}}$.


\textbf{Normal IC-T}: For a tetrahedron $ABCD$, if the IR of the four faces of a point exists and is not 0, then the point is called normal IC-T. The set of all normal ICs-T of the tetrahedron $ABCD$ is denoted as ${{\pi }_{ABCD}}$.

In the case of no confusion, the normal IC-T is also referred to as IC. Readers can distinguish triangle IC and tetrahedron IC according to context.

This book deals with the normal IC-T.

\textbf{Inner IC-T}: The ICs-T is within a tetrahedron.

For example, the point $P$ in Figure \ref{fig:tu16.2.1} is an inner IC-T.

\textbf{Outer IC-T}: The ICs-T is outside a tetrahedron.

\section{Concept of intersecting ratio of a tetrahedron}\label{Sec16.3}	
For the following concepts, see Figure \ref{fig:tu16.2.1}.


\textbf {Tensor of intersecting ratio of a tetrahedron (abbreviated as tensor of IR-T)}: Given a tetrahedron $ABCD$ and an IC-T $P$, the tensor of IR-T is defined by the following:
\[{\boldsymbol{\lambda}_{ABCD}^{P}}=\left( \begin{matrix}
	\boldsymbol{\lambda }_{\triangle BCD}^{{P}_A} & \boldsymbol{\lambda }_{\triangle CDA}^{{P}_B} & \boldsymbol{\lambda }_{\triangle DAB}^{{P}_C} & \boldsymbol{\lambda }_{\triangle ABC}^{{P}_D}  \\
\end{matrix} \right).\]


Where $\boldsymbol{\lambda }_{\triangle BCD}^{{P}_A}$, $\boldsymbol{\lambda }_{\triangle CDA}^{{P}_B}$, $\boldsymbol{\lambda }_{\triangle DAB}^{{P}_C}$ and $\boldsymbol{\lambda }_{\triangle ABC}^{{P}_D}$ are the vectors of IRs for the intersecting centers of faces (ICs-Fs) corresponding to $\triangle BCD $, $\triangle CDA $, $\triangle DAB $ and $\triangle ABC$ at point $P$, respectively, (see Chapter \ref{Ch4}), i.e. 

\[\boldsymbol{\lambda }_{\triangle BCD}^{{P}_A}=\left( \begin{matrix}
	\lambda _{BC}^{{P}_A} & \lambda _{CD}^{{P}_A} & \lambda _{DB}^{{P}_A}  \\
\end{matrix} \right),\]
\[\boldsymbol{\lambda }_{\triangle CDA}^{{P}_B}=\left( \begin{matrix}
	\lambda _{CD}^{{P}_B} & \lambda _{DA}^{{P}_B} & \lambda _{AC}^{{P}_B}  \\
\end{matrix} \right),\]
\[\boldsymbol{\lambda }_{\triangle DAB}^{{P}_C}=\left( \begin{matrix}
	\lambda _{DA}^{{P}_C} & \lambda _{AB}^{{P}_C} & \lambda _{BD}^{{P}_C}  \\
\end{matrix} \right),\]
\[\boldsymbol{\lambda }_{\triangle ABC}^{{P}_D}=\left( \begin{matrix}
	\lambda _{AB}^{{P}_D} & \lambda _{BC}^{{P}_D} & \lambda _{CA}^{{P}_D}  \\
\end{matrix} \right).\]


The tensor of IR-T is an ordered array of vectors of IRs for the $\boldsymbol{\lambda }_{\triangle BCD}^{{P}_A}$, $\boldsymbol{\lambda }_{\triangle CDA}^{{P}_B}$, $\boldsymbol{\lambda }_{\triangle DAB}^{{P}_C}$ and $\boldsymbol{\lambda }_{\triangle ABC}^{{P}_D}$. The tensor of IR-T (${{\boldsymbol{\lambda }}_{ABCD}^{P}}$) establishes a corresponding relationship with the point $P$ in space.

\section{Concept of tetrahedral frame}\label{Sec16.4}
\textbf{Tetrahedral frame}: Given a tetrahedron $ABCD$, point $O$ is any point in space, point $O$ and four vectors $\overrightarrow{OA}$, $\overrightarrow{OB}$, $\overrightarrow{OC}$  and $\overrightarrow{OD}$ together form a frame, and is called tetrahedral frame, denoted as $\left( O;A,B,C,D \right)$. Point $O$ is called the origin of the frame. Tetrahedral frame $\left( O;A,B,C,D \right)$ is also denoted as \textbf{tetrahedron \bm{$ABCD$} frame} or \textbf{\bm{$ABCD$} frame}. Vectors $\overrightarrow{OA}$, $\overrightarrow{OB}$, $\overrightarrow{OC}$ and $\overrightarrow{OD}$ are called \textbf{frame vectors}. Both the frame vectors and the tetrahedral frame are simply called frame when there is no need to distinguish them strictly.

Tetrahedral frame $\left( O;A,B,C,D \right)$ can be regarded as a coordinate system.

The position of origin $O$ of a tetrahedral frame $\left( O;A,B,C,D \right)$ can be arbitrary, which is also called free origin.

\textbf{Frame of IC-T}: The tetrahedral frame whose origin $O$ is an IC-T of the tetrahedron.

\textbf{Frame of circumcenter}: A frame of tetrahedron whose origin $O$ is the circumcenter of the tetrahedron.

The position of the origin $O$ for the frame of circumcenter cannot be arbitrary, but can only be the circumcenter of the tetrahedron $ABCD$. Similar concepts such as centroid frame can be obtained.

\textbf{Frame of edge} (abbreviated as FE): Given a tetrahedron, for a vertex $A$ and three vectors $\overrightarrow{AB}$, $\overrightarrow{AC}$, $\overrightarrow{AD}$, they form a frame together, this frame is denoted as $\left( A;B,C,D \right)$. Frame of edge is also called \textbf{edge frame}.

The origin of the edge frame $\left( A;B,C,D \right)$ is $A$, three vectors $\overrightarrow{AB}$, $\overrightarrow{AC}$, $\overrightarrow{AD}$ of the edge frame $\left( A;B,C,D \right)$ coincide with three edges $AB$, $AC$, $AD$ of the tetrahedron $ABCD$ respectively.

Similarly, we can define $\left( B;C,A,D \right)$, $\left( C;D,A,B \right)$, $\left( D;A,B,C \right)$.

\section{Concept on vector of intersecting center of a tetrahedron}\label{Sec16.5}
\textbf{Vector from origin to intersecting center of tetrahedron} (abbreviated as VOIC-T): Given a tetrahedron $ABCD$, point $O$ is the origin of the tetrahedral frame $\left( O;A,B,C,D \right)$ and point $P$ is an IC-T of $ABCD$, then vector $\overrightarrow{OP}$ is called the vector from origin to intersecting center of tetrahedron.

VOIC-T $\overrightarrow{OP}$ is the position vector where the origin of the frame $O$ points to $P$ (the IC-T). Since the origin of the frame $O$ can be arbitrary, the origin of the frame $O$ may also be an IC-T. The key point of the concept of a VOIC-T is to emphasize that the starting point of the vector $\overrightarrow{OP}$ must be the origin of tetrahedral frame $\left( O;A,B,C,D \right)$ and the terminal $P$ must be an intersecting center of the tetrahedron $ABCD$.

\textbf{Vector of two intersecting centers of tetrahedron} (abbreviated as VTICs-T): Given a $ABCD$, point $O$ is the origin of the tetrahedral frame $\left( O;A,B,C,D \right)$, point ${{P}_{1}}$, ${{P}_{2}}$ are the two ICs-T of $ABCD$, then the vector $\overrightarrow{{{P}_{1}}{{P}_{2}}}$ is called the vector of two intersecting centers of tetrahedron (abbreviated as VTICs-T).

The VTICs-T is a vector with two ICs-T. The position of the origin of the frame is arbitrary, that is, the free origin.

The main point of the concept VTICs-T is: for a vector $\overrightarrow{{{P}_{1}}{{P}_{2}}}$, its starting point ${{P}_{1}}$ and terminal point ${{P}_{2}}$ are both ICs-T of the tetrahedron $ABCD$, while the position of the origin $O$ is arbitrary and the origin $O$ is usually both different from point ${{P}_{1}}$ and point ${{P}_{2}}$.

Both the VOIC-T and the VTICs-T can be expressed linearly by different frames, which will be discussed later.

\textbf{Vector from vertex to intersecting center of tetrahedron} (abbreviated as VVIC-T): The position vector from the vertex to the IC-T of a given tetrahedron.


\chapter{Intersecting ratio of face on special intersecting center of a tetrahedron}\label{Ch17}
\thispagestyle{empty}
It has been pointed out that, in Intercenter Geometry, the IR is an important concept throughout. In this chapter, I give some methods for calculating the IR of some special intersecting centers of a tetrahedron, which include the centroid, incenter and excenter of the tetrahedron. The calculation formula of the IR of those special centers will be often used in the later research.

The intersecting ratio of a tetrahedron is related to the intersecting ratios of faces (abbreviated as IR-F), that is, the IR of the four triangular faces of the tetrahedron.

The intersecting center of face (abbreviated as IC-F) is the IC of four triangles on the surface of a tetrahedron, so it can be calculated by using the previous calculation formulas about the IR of triangles.

\section{IR-F of centroid of a tetrahedron}\label{Sec17.1}
For a given tetrahedron, the centroid of the tetrahedron has four ICs-F, which are the centroids of the four surface triangles of the tetrahedron. Each IR of the four centroids corresponding to the four triangles is 1, so each IR of the ICs-F of the corresponding four triangles of the tetrahedron is 1. So there is the following theorem.

\begin{theorem}{Centroid tensor of IR-T, Daiyuan Zhang}{Thm17.1.1}\label{Thm17.1.1} 
	Given the tetrahedron $ABCD$ and the centroid $G$ of the tetrahedron, then:
	\[\boldsymbol{\lambda }_{ABCD}^{G}=\left( \begin{matrix}
		\boldsymbol{\lambda }_{\triangle BCD}^{G_{A}} & \boldsymbol{\lambda }_{\triangle CDA}^{G_{B}} & \boldsymbol{\lambda }_{\triangle DAB}^{G_{C}} & \boldsymbol{\lambda }_{\triangle ABC}^{G_{D}}  \\
	\end{matrix} \right).\]
	Where
	\[\boldsymbol{\lambda }_{\triangle BCD}^{{G}_A}=\left( \begin{matrix}
		\lambda _{BC}^{{G}_A} & \lambda _{CD}^{{G}_A} & \lambda _{DB}^{{G}_A}  \\
	\end{matrix} \right)=\left( \begin{matrix}
		1 & 1 & 1  \\
	\end{matrix} \right),\]
	\[\boldsymbol{\lambda }_{\triangle CDA}^{{G}_B}=\left( \begin{matrix}
		\lambda _{CD}^{{G}_B} & \lambda _{DA}^{{G}_B} & \lambda _{AC}^{{G}_B}  \\
	\end{matrix} \right)=\left( \begin{matrix}
		1 & 1 & 1  \\
	\end{matrix} \right),\]
	\[\boldsymbol{\lambda }_{\triangle DAB}^{{G}_C}=\left( \begin{matrix}
		\lambda _{DA}^{{G}_C} & \lambda _{AB}^{{G}_C} & \lambda _{BD}^{{G}_C}  \\
	\end{matrix} \right)=\left( \begin{matrix}
		1 & 1 & 1  \\
	\end{matrix} \right),\]
	\[\boldsymbol{\lambda }_{\triangle ABC}^{{G}_D}=\left( \begin{matrix}
		\lambda _{AB}^{{G}_D} & \lambda _{BC}^{{G}_D} & \lambda _{CA}^{{G}_D}  \\
	\end{matrix} \right)=\left( \begin{matrix}
		1 & 1 & 1  \\
	\end{matrix} \right).\]
\end{theorem}

\section{IR-F of incenter of a tetrahedron}\label{Sec17.2}
\begin{figure}[h]
	\centering
	\includegraphics[totalheight=6cm]{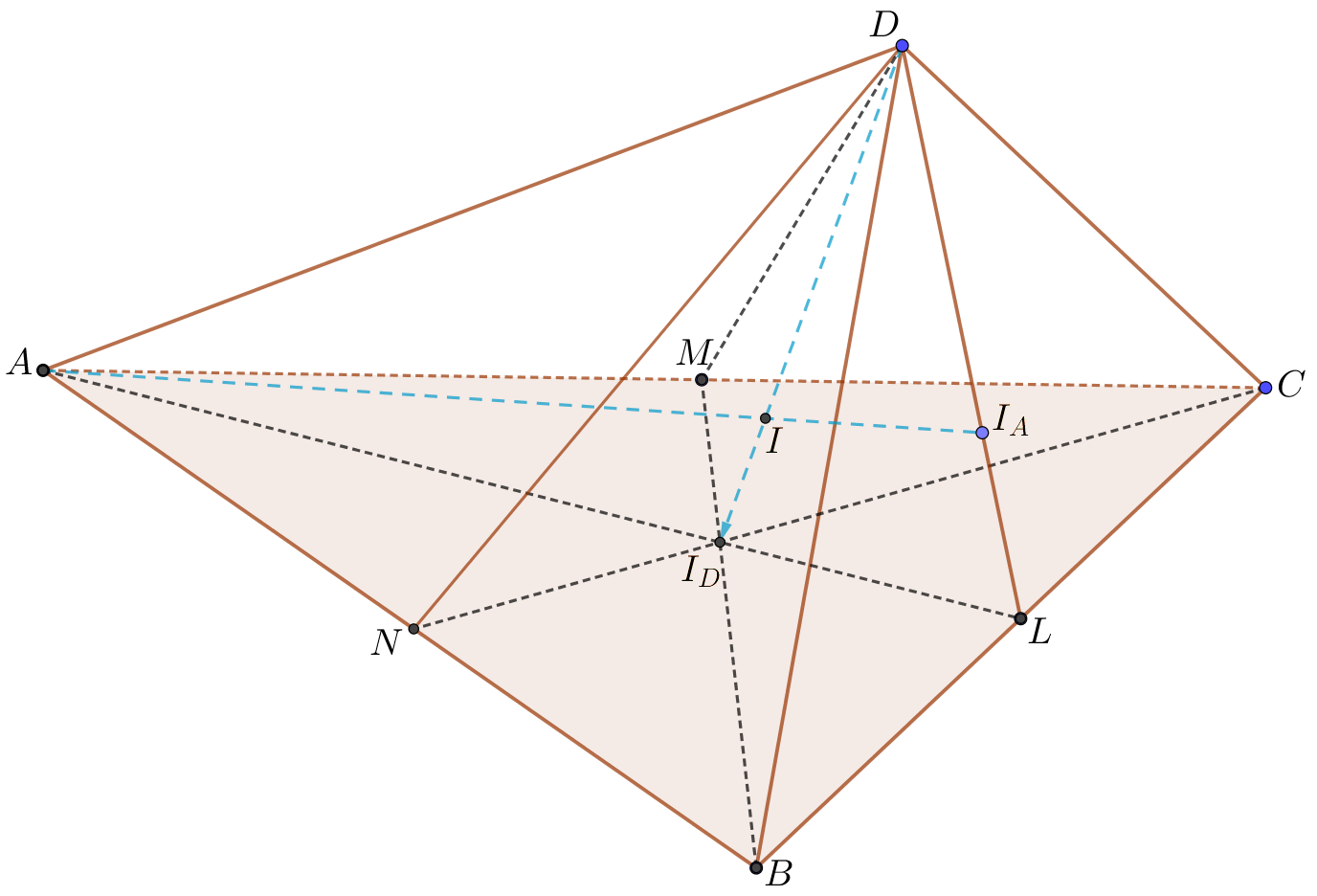}
	\caption{Diagram for calculating the IR-F of incenter of a tetrahedron} \label{fig:tu17.2.1}
\end{figure}

As shown in Figure \ref{fig:tu17.2.1}, the inscribed spherical center (or incenter) of the tetrahedron $ABCD$ is $I$, and the intersection of the straight line passing through the vertex $D$ and the incenter $I$ which intersect the opposite face is ${{I}_{D}}$, and ${{I}_{D}}$ is the IC-F of $\triangle ABC$. Similarly, we can get the ICs-F of the other three triangles.

Only two ICs-F ${{I}_{A}}$ and ${{I}_{D}}$ are shown in Figure \ref{fig:tu17.2.1}. The IRs-F will be calculated below.

\begin{itemize}
	\item IR-F of ${{I}_{D}}$
\end{itemize}

The bisector of dihedral angle $BDC-CD-ADC$ is $CDN$. The distances from the point $N$ on the angular bisector to the plane $BDC$ and $ADC$ are equal, and the result is:

\[\lambda _{AB}^{{{I}_{D}}}=\frac{\overrightarrow{AN}}{\overrightarrow{NB}}=\frac{AN}{NB}=\frac{{{S}_{ANC}}}{{{S}_{NBC}}}=\frac{{{V}_{ANCD}}}{{{V}_{NBCD}}}=\frac{{{S}_{ACD}}}{{{S}_{BDC}}}.\]	

Each of the four face areas of the tetrahedron can be calculated by the lengths of the six edges of the tetrahedron according to Helen's formula. To simplify the expression, some notations are introduced: ${{S}_{ACD}}={{S}^{B}}$, ${{S}_{BDC}}={{S}^{A}}$, and the above formula becomes

\[\lambda _{AB}^{{{I}_{D}}}=\frac{\overrightarrow{AN}}{\overrightarrow{NB}}=\frac{{{S}^{B}}}{{{S}^{A}}}.\]	

Similarly (the meaning of notations are the same as before, i.e. ${{S}_{ABD}}={{S}^{C}}$, ${{S}_{ABC}}={{S}^{D}}$),

\[\lambda _{BC}^{{{I}_{D}}}=\frac{\overrightarrow{BL}}{\overrightarrow{LC}}=\frac{{{S}^{C}}}{{{S}^{B}}},\]	
\[\lambda _{CA}^{{{I}_{D}}}=\frac{\overrightarrow{CM}}{\overrightarrow{MA}}=\frac{{{S}^{A}}}{{{S}^{C}}}.\]	

Similarly, the following results can be obtained.
\begin{itemize}
	\item IR-F of ${{I}_{A}}$
\end{itemize}

\begin{flalign*}
	\lambda _{BC}^{{{I}_{A}}}=\frac{{{S}^{C}}}{{{S}^{B}}}, \lambda _{CD}^{{{I}_{A}}}=\frac{{{S}^{D}}}{{{S}^{C}}}, \lambda _{DB}^{{{I}_{A}}}=\frac{{{S}^{B}}}{{{S}^{D}}}.	
\end{flalign*}

\begin{itemize}
	\item IR-F of ${{I}_{B}}$
\end{itemize}

\begin{flalign*}
	\lambda _{CD}^{{{I}_{B}}}=\frac{{{S}^{D}}}{{{S}^{C}}},\lambda _{DA}^{{{I}_{B}}}=\frac{{{S}^{A}}}{{{S}^{D}}},\lambda _{AC}^{{{I}_{B}}}=\frac{{{S}^{C}}}{{{S}^{A}}}.	
\end{flalign*}
	
\begin{itemize}
	\item IR-F of ${{I}_{C}}$
\end{itemize}

\begin{flalign*}
	\lambda _{DA}^{{{I}_{C}}}=\frac{{{S}^{A}}}{{{S}^{D}}},\lambda _{AB}^{{{I}_{C}}}=\frac{{{S}^{B}}}{{{S}^{A}}},\lambda _{BD}^{{{I}_{C}}}=\frac{{{S}^{D}}}{{{S}^{B}}}.	
\end{flalign*}

So there is the following theorem.
\begin{theorem}{Incenter tensor of IR-T, Daiyuan Zhang}{Thm17.2.1}\label{Thm17.2.1} 
	Given the tetrahedron $ABCD$ and the incenter $I$ of the tetrahedron, then:
\[\boldsymbol{\lambda }_{ABCD}^{I}=\left( \begin{matrix}
	\boldsymbol{\lambda }_{\triangle BCD}^{{I}_A} & \boldsymbol{\lambda }_{\triangle CDA}^{{I}_B} & \boldsymbol{\lambda }_{\triangle DAB}^{{I}_C} & \boldsymbol{\lambda }_{\triangle ABC}^{{I}_D}  \\
\end{matrix} \right).\]
Where
\[\boldsymbol{\lambda }_{\triangle BCD}^{{I}_A}=\left( \begin{matrix}
	\lambda _{BC}^{{I}_A} & \lambda _{CD}^{{I}_A} & \lambda _{DB}^{{I}_A}  \\
\end{matrix} \right)=\left( \begin{matrix}
	\frac{{{S}^{C}}}{{{S}^{B}}} & \frac{{{S}^{D}}}{{{S}^{C}}} & \frac{{{S}^{B}}}{{{S}^{D}}}  \\
\end{matrix} \right),\]
\[\boldsymbol{\lambda }_{\triangle CDA}^{{I}_B}=\left( \begin{matrix}
	\lambda _{CD}^{{I}_B} & \lambda _{DA}^{{I}_B} & \lambda _{AC}^{{I}_B}  \\
\end{matrix} \right)=\left( \begin{matrix}
	\frac{{{S}^{D}}}{{{S}^{C}}} & \frac{{{S}^{A}}}{{{S}^{D}}} & \frac{{{S}^{C}}}{{{S}^{A}}}  \\
\end{matrix} \right),\]
\[\boldsymbol{\lambda }_{\triangle DAB}^{{I}_C}=\left( \begin{matrix}
	\lambda _{DA}^{{I}_C} & \lambda _{AB}^{{I}_C} & \lambda _{BD}^{{I}_C}  \\
\end{matrix} \right)=\left( \begin{matrix}
	\frac{{{S}^{A}}}{{{S}^{D}}} & \frac{{{S}^{B}}}{{{S}^{A}}} & \frac{{{S}^{D}}}{{{S}^{B}}}  \\
\end{matrix} \right),\]
\[\boldsymbol{\lambda }_{\triangle ABC}^{{I}_D}=\left( \begin{matrix}
	\lambda _{AB}^{{I}_D} & \lambda _{BC}^{{I}_D} & \lambda _{CA}^{{I}_D}  \\
\end{matrix} \right)=\left( \begin{matrix}
	\frac{{{S}^{B}}}{{{S}^{A}}} & \frac{{{S}^{C}}}{{{S}^{B}}} & \frac{{{S}^{A}}}{{{S}^{C}}}  \\
\end{matrix} \right).\]
\end{theorem}

\section{Intersecting ratios of faces for excenters in the trihedral angles of a tetrahedron}\label{Sec17.3}
\begin{figure}[h]
	\centering
	\includegraphics[totalheight=6cm]{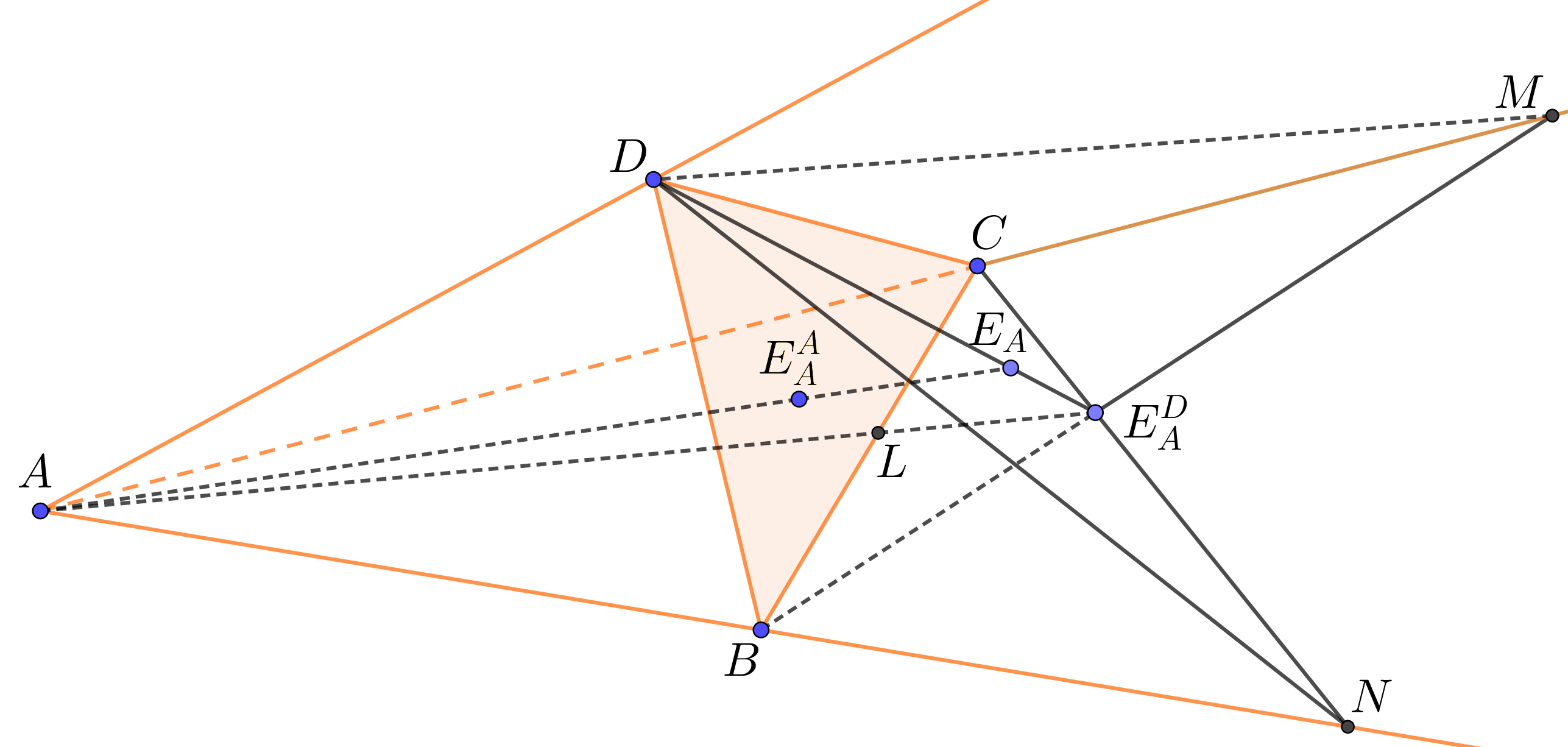}
	\caption{Diagram of calculating the IRs-Fs of the excenters in trihedral angles of a tetrahedron} \label{fig:tu17.3.1}
\end{figure}

In this section, the IRs-Fs of the excenters in the four trihedral angles of a tetrahedron are studied, see Figure \ref{fig:tu17.3.1}.

There are four trihedral angles in a tetrahedron, and each trihedral angle has an inscribed tangent sphere. Each center of inscribed tangent sphere can lead to a straight line with a vertex of the tetrahedron. There is a unique intersection point between the straight line and the opposite face of the vertex, and this intersection point is called the IC-F of the excenter.

According to solid geometry, the axis of a trihedral angle passes through the vertex, the incenter and the excenter of the trihedral angle.

Only two ICs-F: $E_{A}^{A}$ and $E_{A}^{D}$ are drawn in Figure \ref{fig:tu17.3.1}. 

%
The meanings of some notations are given below, with the notations $E_{A}^{D}$ and $E_{A}^{A}$ as examples to illustrate the meanings (see Figure \ref{fig:tu17.3.1}). 

First, the notation $E$ is used to denote the excenter. Since the tetrahedron $ABCD$ has four trihedral angles, they are the trihedral angles corresponding to the tetrahedral vertexs $A$, $B$, $C$, and $D$. Each trihedral angle has an inscribed tangent sphere center (excenter), so it needs to be denoted with different notations. 

Secondly, this book uses the method of adding subscript in the notation $E$ to distinguish the excenters of different trihedral angles, for example, ${{E}_{A}}$ is the center of the inscribed sphere in the trihedral angle corresponding to the vertex $A$ of the tetrahedron $ABCD$ (see Figure \ref{fig:tu17.3.1}).

Again, the four straight lines connecting the tangential spherical center ${{E}_{A}}$ with the four vertices $A$, $B$, $C$, $D$ of the tetrahedron $ABCD$ will intersect the opposite faces of the corresponding vertexs, and there are four such intersections. In order to distinguish those different intersections, on the basis of ${{E}_{A}}$, this book uses the superscript to distinguish them. For example, $E_{A}^{D}$ represents the intersection of the straight line $D{{E}_{A}}$ with the opposite face of the vertex $D$ (i.e., in the plane of $ABC$) (see Figure \ref{fig:tu17.3.1}), i.e.,
$E_{A}^{D}=\overleftrightarrow{D{{E}_{A}}}\cap {{P}_{\triangle ABC}}$.


To sum up, the notation $E_{A}^{D}$ is denoted by the subscript $A$ for the position of the excenter, i.e. the center of the inscribed tangent sphere in the trihedral angle of the vertex $A$ of the tetrahedron $ABCD$, and the superscript $D$ represents the opposite face of the vertex $D$ of the tetrahedron $ABCD$. The notation $E_{A}^{D}$ represents a point on the plane of $\triangle ABC$, so it may be an IC of $\triangle ABC$, i.e., an IC-F of the tetrahedron $ABCD$. In other words, the $E_{A}^{D}$ (IC-F) may only be on the plane of $\triangle ABC$, so the subscript of the IR for the $E_{A}^{D}$ (IC-F) of the $\triangle ABC$ cannot appear the symbol $D$, in fact, the vector of IR for the $E_{A}^{D}$ (IC-F) is:
\[\boldsymbol{\lambda }_{\triangle ABC}^{E_{A}^{D}}=\left( \begin{matrix}
	\lambda _{AB}^{E_{A}^{D}} & \lambda _{BC}^{E_{A}^{D}} & \lambda _{CA}^{E_{A}^{D}}  \\
\end{matrix} \right).\]


According to the meaning of this notation, we can infer the meaning of notation $E_{A}^{A}$. The subscript $A$ of $E_{A}^{A}$ denotes the excenter in the trihedral angle of the vertex $A$ of the tetrahedron $ABCD$. The superscript $A$ of $E_{A}^{A}$ represents the opposite face of the vertex $A$ of the tetrahedron $ABCD$. So the notation $E_{A}^{A}$ (IC-F) may only be the IC of $\triangle BCD$, and the subscript of the IR for the $E_{A}^{A}$ (IC-F) of the $\triangle BCD$ cannot appear the symbol $A$, in fact, the vector of IR for the $E_{A}^{A}$ (IC-F) is: 

\[\boldsymbol{\lambda }_{\triangle BCD}^{E_{A}^{A}}=\left( \begin{matrix}
	\lambda _{BC}^{E_{A}^{A}} & \lambda _{CD}^{E_{A}^{A}} & \lambda _{DB}^{E_{A}^{A}}  \\
\end{matrix} \right).\] 


The $E_{A}^{A}$ (IC-F) is located on the axis of the trihedral angle of the vertex $A$ of the tetrahedron $ABCD$, so it is located in the interior of $\triangle BCD$, and the three IRs are all positive numbers, that is to say, the IR must be positive as long as the subscripts $BC$, $CD$, $DB$ appear at the same time. In addition, according to Ceva's theorem, the product of the three IRs of a triangle is 1, so the three IRs are either all positive, or one is positive, and the other two are negative. According to this rule, we can write the vector of IRs for ${{E}_{B}}$ and ${{E}_{C}}$.


To sum up, each tensor of IR-T of the excenters is:

\[\boldsymbol{\lambda }_{ABCD}^{E_{A}^{A}}=\left( \begin{matrix}
	\boldsymbol{\lambda }_{\triangle BCD}^{E_{A}^{A}} & \boldsymbol{\lambda }_{\triangle CDA}^{E_{A}^{B}} & \boldsymbol{\lambda }_{\triangle DAB}^{E_{A}^{C}} & \boldsymbol{\lambda }_{\triangle ABC}^{E_{A}^{D}}  \\
\end{matrix} \right),\]
\[\boldsymbol{\lambda }_{\triangle BCD}^{E_{A}^{A}}=\left( \begin{matrix}
	\lambda _{BC}^{E_{A}^{A}} & \lambda _{CD}^{E_{A}^{A}} & \lambda _{DB}^{E_{A}^{A}}  \\
\end{matrix} \right),\]
\[\boldsymbol{\lambda }_{\triangle CDA}^{E_{A}^{B}}=\left( \begin{matrix}
	\lambda _{CD}^{E_{A}^{B}} & \lambda _{DA}^{E_{A}^{B}} & \lambda _{AC}^{E_{A}^{B}}  \\
\end{matrix} \right),\]
\[\boldsymbol{\lambda }_{\triangle DAB}^{E_{A}^{C}}=\left( \begin{matrix}
	\lambda _{DA}^{E_{A}^{C}} & \lambda _{AB}^{E_{A}^{C}} & \lambda _{BD}^{E_{A}^{C}}  \\
\end{matrix} \right),\]
\[\boldsymbol{\lambda }_{\triangle ABC}^{E_{A}^{D}}=\left( \begin{matrix}
	\lambda _{AB}^{E_{A}^{D}} & \lambda _{BC}^{E_{A}^{D}} & \lambda _{CA}^{E_{A}^{D}}  \\
\end{matrix} \right);\]

\[\boldsymbol{\lambda }_{ABCD}^{E_{B}^{A}}=\left( \begin{matrix}
	\boldsymbol{\lambda }_{\triangle BCD}^{E_{B}^{A}} & \boldsymbol{\lambda }_{\triangle CDA}^{E_{B}^{B}} & \boldsymbol{\lambda }_{\triangle DAB}^{E_{B}^{C}} & \boldsymbol{\lambda }_{\triangle ABC}^{E_{B}^{D}}  \\
\end{matrix} \right),\]
\[\boldsymbol{\lambda }_{\triangle BCD}^{E_{B}^{A}}=\left( \begin{matrix}
	\lambda _{BC}^{E_{B}^{A}} & \lambda _{CD}^{E_{B}^{A}} & \lambda _{DB}^{E_{B}^{A}}  \\
\end{matrix} \right),\]
\[\boldsymbol{\lambda }_{\triangle CDA}^{E_{B}^{B}}=\left( \begin{matrix}
	\lambda _{CD}^{E_{B}^{B}} & \lambda _{DA}^{E_{B}^{B}} & \lambda _{AC}^{E_{B}^{B}}  \\
\end{matrix} \right),\]
\[\boldsymbol{\lambda }_{\triangle DAB}^{E_{B}^{C}}=\left( \begin{matrix}
	\lambda _{DA}^{E_{B}^{C}} & \lambda _{AB}^{E_{B}^{C}} & \lambda _{BD}^{E_{B}^{C}}  \\
\end{matrix} \right),\]
\[\boldsymbol{\lambda }_{\triangle ABC}^{E_{B}^{D}}=\left( \begin{matrix}
	\lambda _{AB}^{E_{B}^{D}} & \lambda _{BC}^{E_{B}^{D}} & \lambda _{CA}^{E_{B}^{D}}  \\
\end{matrix} \right);\]

\[\boldsymbol{\lambda }_{ABCD}^{E_{C}^{A}}=\left( \begin{matrix}
	\boldsymbol{\lambda }_{\triangle BCD}^{E_{C}^{A}} & \boldsymbol{\lambda }_{\triangle CDA}^{E_{C}^{B}} & \boldsymbol{\lambda }_{\triangle DAB}^{E_{C}^{C}} & \boldsymbol{\lambda }_{\triangle ABC}^{E_{C}^{D}}  \\
\end{matrix} \right),\]
\[\boldsymbol{\lambda }_{\triangle BCD}^{E_{C}^{A}}=\left( \begin{matrix}
	\lambda _{BC}^{E_{C}^{A}} & \lambda _{CD}^{E_{C}^{A}} & \lambda _{DB}^{E_{C}^{A}}  \\
\end{matrix} \right),\]
\[\boldsymbol{\lambda }_{\triangle CDA}^{E_{C}^{B}}=\left( \begin{matrix}
	\lambda _{CD}^{E_{C}^{B}} & \lambda _{DA}^{E_{C}^{B}} & \lambda _{AC}^{E_{C}^{B}}  \\
\end{matrix} \right),\]
\[\boldsymbol{\lambda }_{\triangle DAB}^{E_{C}^{C}}=\left( \begin{matrix}
	\lambda _{DA}^{E_{C}^{C}} & \lambda _{AB}^{E_{C}^{C}} & \lambda _{BD}^{E_{C}^{C}}  \\
\end{matrix} \right),\]
\[\boldsymbol{\lambda }_{\triangle ABC}^{E_{C}^{D}}=\left( \begin{matrix}
	\lambda _{AB}^{E_{C}^{D}} & \lambda _{BC}^{E_{C}^{D}} & \lambda _{CA}^{E_{C}^{D}}  \\
\end{matrix} \right);\]

\[\boldsymbol{\lambda }_{ABCD}^{E_{D}^{A}}=\left( \begin{matrix}
	\boldsymbol{\lambda }_{\triangle BCD}^{E_{D}^{A}} & \boldsymbol{\lambda }_{\triangle CDA}^{E_{D}^{B}} & \boldsymbol{\lambda }_{\triangle DAB}^{E_{D}^{C}} & \boldsymbol{\lambda }_{\triangle ABC}^{E_{D}^{D}}  \\
\end{matrix} \right),\]
\[\boldsymbol{\lambda }_{\triangle BCD}^{E_{D}^{A}}=\left( \begin{matrix}
	\lambda _{BC}^{E_{D}^{A}} & \lambda _{CD}^{E_{D}^{A}} & \lambda _{DB}^{E_{D}^{A}}  \\
\end{matrix} \right),\]
\[\boldsymbol{\lambda }_{\triangle CDA}^{E_{D}^{B}}=\left( \begin{matrix}
	\lambda _{CD}^{E_{D}^{B}} & \lambda _{DA}^{E_{D}^{B}} & \lambda _{AC}^{E_{D}^{B}}  \\
\end{matrix} \right),\]
\[\boldsymbol{\lambda }_{\triangle DAB}^{E_{D}^{C}}=\left( \begin{matrix}
	\lambda _{DA}^{E_{D}^{C}} & \lambda _{AB}^{E_{D}^{C}} & \lambda _{BD}^{E_{D}^{C}}  \\
\end{matrix} \right),\]
\[\boldsymbol{\lambda }_{\triangle ABC}^{E_{D}^{D}}=\left( \begin{matrix}
	\lambda _{AB}^{E_{D}^{D}} & \lambda _{BC}^{E_{D}^{D}} & \lambda _{CA}^{E_{D}^{D}}  \\
\end{matrix} \right).\]

Now let's calculate the IRs-Fs.

\subsection{IR-F of excenter $E_A$}\label{Subsec17.3.1}

\subsubsection{IR of IC-F $E_{A}^{A}$ of excenter ${{E}_{A}}$}
From the discussion in the previous section it is known that: the IR of IC-F $E_{A}^{A}$ is
\begin{flalign*}
	\lambda _{BC}^{E_{A}^{A}}=\frac{{{S}^{C}}}{{{S}^{B}}},\lambda _{CD}^{E_{A}^{A}}=\frac{{{S}^{D}}}{{{S}^{C}}},\lambda _{DB}^{E_{A}^{A}}=\frac{{{S}^{B}}}{{{S}^{D}}}.	
\end{flalign*}

\subsubsection{IR of IC-F $E_{A}^{D}$ of excenter ${{E}_{A}}$}

The bisector of dihedral angle $BCD-CD-MCD$ is $NCD$. The distances from the point $N$ on the bisector of dihedral angle to the plane $BCD$ and the plane $MCD$ (that is, the plane $CDA$) are equal, then:
\[\begin{aligned}
	\lambda _{AB}^{E_{A}^{D}}&=\frac{\overrightarrow{AN}}{\overrightarrow{NB}}=-\frac{AN}{NB}=-\frac{{{S}_{DAN}}}{{{S}_{DNB}}}=-\frac{{{V}_{CDAN}}}{{{V}_{CDNB}}} \\ 
	& =-\frac{{{V}_{NCDA}}}{{{V}_{NBCD}}}=-\frac{{{S}_{CDA}}}{{{S}_{BCD}}}=-\frac{{{S}^{B}}}{{{S}^{A}}}.  
\end{aligned}\]	

According to the discussion in the previous section, we have
\[\lambda _{BC}^{E_{A}^{D}}=\frac{\overrightarrow{BL}}{\overrightarrow{LC}}=\frac{BL}{LC}=\frac{{{S}^{C}}}{{{S}^{B}}}.\]	


Similarly, for $\lambda _{CA}^{E_{A}^{D}}$, the bisector of dihedral angle $CBD-BD-NBD$ is plane $MBD$. The distances from the point $M$ on the bisector of dihedral angle to the plane $CBD$ and $NBD$ are equal, so that:
\[\begin{aligned}
	\lambda _{CA}^{E_{A}^{D}}&=\frac{\overrightarrow{CM}}{\overrightarrow{MA}}=-\frac{CM}{MA}=-\frac{{{S}_{DCM}}}{{{S}_{DMA}}}=-\frac{{{V}_{BDCM}}}{{{V}_{BDMA}}} \\ 
	& =-\frac{{{V}_{MBDC}}}{{{V}_{MABD}}}=-\frac{{{S}_{BDC}}}{{{S}_{ABD}}}=-\frac{{{S}^{A}}}{{{S}^{C}}}.  
\end{aligned}\]	

Similarly, the formula for the IRs of two other ICs-F can be obtained and written directly below.

\subsubsection{IR of IC-F $E_{A}^{B}$ of excenter ${{E}_{A}}$}
\begin{flalign*}
	\lambda _{CD}^{E_{A}^{B}}=\frac{{{S}^{D}}}{{{S}^{C}}},\lambda _{DA}^{E_{A}^{B}}=-\frac{{{S}^{A}}}{{{S}^{D}}},\lambda _{AC}^{E_{A}^{B}}=-\frac{{{S}^{C}}}{{{S}^{A}}}.	
\end{flalign*}
	
\subsubsection{IR of IC-F $E_{A}^{C}$ of excenter ${{E}_{A}}$}
\begin{flalign*}
	\lambda _{DA}^{E_{A}^{C}}=-\frac{{{S}^{A}}}{{{S}^{D}}},\lambda _{AB}^{E_{A}^{C}}=-\frac{{{S}^{B}}}{{{S}^{A}}},\lambda _{BD}^{E_{A}^{C}}=\frac{{{S}^{D}}}{{{S}^{B}}}.		
\end{flalign*}

Similar to the previous analysis, the IR formulas of the ICs-F corresponding to the centers (excenter) of inscribed tangent spheres of the three other trihedral angles of the tetrahedron $ABCD$ can be obtained, and they are written directly below.

\subsection{IR-F of excenter $E_B$}\label{Subsec17.3.2}
\subsubsection{IR of IC-F $E_{B}^{A}$ of excenter ${{E}_{B}}$}

\begin{flalign*}
	\lambda _{BC}^{E_{B}^{A}}=-\frac{{{S}^{C}}}{{{S}^{B}}},\lambda _{CD}^{E_{B}^{A}}=\frac{{{S}^{D}}}{{{S}^{C}}},\lambda _{DB}^{E_{B}^{A}}=-\frac{{{S}^{B}}}{{{S}^{D}}}.			
\end{flalign*}

\subsubsection{IR of IC-F $E_{B}^{B}$ of excenter ${{E}_{B}}$}
\begin{flalign*}
	\lambda _{CD}^{E_{B}^{B}}=\frac{{{S}^{D}}}{{{S}^{C}}},\lambda _{DA}^{E_{B}^{B}}=\frac{{{S}^{A}}}{{{S}^{D}}},\lambda _{AC}^{E_{B}^{B}}=\frac{{{S}^{C}}}{{{S}^{A}}}.			
\end{flalign*}

\subsubsection{IR of IC-F $E_{B}^{C}$ of excenter ${{E}_{B}}$}
\begin{flalign*}
	\lambda _{DA}^{E_{B}^{C}}=\frac{{{S}^{A}}}{{{S}^{D}}},\lambda _{AB}^{E_{B}^{C}}=-\frac{{{S}^{B}}}{{{S}^{A}}},\lambda _{BD}^{E_{B}^{C}}=-\frac{{{S}^{D}}}{{{S}^{B}}}.		
\end{flalign*}
	
\subsubsection{IR of IC-F $E_{B}^{D}$ of excenter ${{E}_{B}}$}
\begin{flalign*}
	\lambda _{AB}^{E_{B}^{D}}=-\frac{{{S}^{B}}}{{{S}^{A}}},\lambda _{BC}^{E_{B}^{D}}=-\frac{{{S}^{C}}}{{{S}^{B}}},\lambda _{CA}^{E_{B}^{D}}=\frac{{{S}^{A}}}{{{S}^{C}}}.		
\end{flalign*}

\subsection{IR-F of excenter $E_C$}\label{Subsec17.3.3}

\subsubsection{IR of IC-F $E_{C}^{A}$ of excenter ${{E}_{C}}$}
\begin{flalign*}
	\lambda _{BC}^{E_{C}^{A}}=-\frac{{{S}^{C}}}{{{S}^{B}}},\lambda _{CD}^{E_{C}^{A}}=-\frac{{{S}^{D}}}{{{S}^{A}}},\lambda _{DB}^{E_{C}^{A}}=\frac{{{S}^{B}}}{{{S}^{D}}}.
\end{flalign*}

\subsubsection{IR of IC-F $E_{C}^{B}$ of excenter ${{E}_{C}}$}
\begin{flalign*}
	\lambda _{CD}^{E_{C}^{B}}=-\frac{{{S}^{D}}}{{{S}^{C}}},\lambda _{DA}^{E_{C}^{B}}=\frac{{{S}^{A}}}{{{S}^{D}}},\lambda _{AC}^{E_{C}^{B}}=-\frac{{{S}^{C}}}{{{S}^{A}}}.		
\end{flalign*}
	
\subsubsection{IR of IC-F $E_{C}^{C}$ of excenter ${{E}_{C}}$}
\begin{flalign*}
	\lambda _{DA}^{E_{C}^{C}}=\frac{{{S}^{A}}}{{{S}^{D}}},\lambda _{AB}^{E_{C}^{C}}=\frac{{{S}^{B}}}{{{S}^{A}}},\lambda _{BD}^{E_{C}^{C}}=\frac{{{S}^{D}}}{{{S}^{B}}}.				
\end{flalign*}

\subsubsection{IR of IC-F $E_{C}^{D}$ of excenter ${{E}_{C}}$}
\begin{flalign*}
	\lambda _{AB}^{E_{C}^{D}}=\frac{{{S}^{B}}}{{{S}^{A}}},\lambda _{BC}^{E_{C}^{D}}=-\frac{{{S}^{C}}}{{{S}^{B}}},\lambda _{CA}^{E_{C}^{D}}=-\frac{{{S}^{A}}}{{{S}^{C}}}.		
\end{flalign*}

\subsection{IR-F of excenter $E_D$}\label{Subsec17.3.4}

\subsubsection{IR of IC-F $E_{C}^{A}$ of excenter ${{E}_{D}}$}
\begin{flalign*}
	\lambda _{BC}^{E_{D}^{A}}=\frac{{{S}^{C}}}{{{S}^{B}}},\lambda _{CD}^{E_{D}^{A}}=-\frac{{{S}^{D}}}{{{S}^{A}}},\lambda _{DB}^{E_{D}^{A}}=-\frac{{{S}^{B}}}{{{S}^{D}}}.		
\end{flalign*}
	
\subsubsection{IR of IC-F $E_{C}^{B}$ of excenter ${{E}_{D}}$}
\begin{flalign*}
	\lambda _{CD}^{E_{D}^{B}}=-\frac{{{S}^{D}}}{{{S}^{C}}},\lambda _{DA}^{E_{D}^{B}}=-\frac{{{S}^{A}}}{{{S}^{D}}},\lambda _{AC}^{E_{D}^{B}}=\frac{{{S}^{C}}}{{{S}^{A}}}.		
\end{flalign*}

\subsubsection{IR of IC-F $E_{C}^{C}$ of excenter ${{E}_{D}}$}
\begin{flalign*}
	\lambda _{DA}^{E_{D}^{C}}=-\frac{{{S}^{A}}}{{{S}^{D}}},\lambda _{AB}^{E_{D}^{C}}=\frac{{{S}^{B}}}{{{S}^{A}}},\lambda _{BD}^{E_{D}^{C}}=-\frac{{{S}^{D}}}{{{S}^{B}}}.		
\end{flalign*}

\subsubsection{IR of IC-F $E_{C}^{D}$ of excenter ${{E}_{D}}$}
\begin{flalign*}
	\lambda _{AB}^{E_{D}^{D}}=\frac{{{S}^{B}}}{{{S}^{A}}},\lambda _{BC}^{E_{D}^{D}}=\frac{{{S}^{C}}}{{{S}^{B}}},\lambda _{CA}^{E_{D}^{D}}=\frac{{{S}^{A}}}{{{S}^{C}}}.		
\end{flalign*}

\subsection{Tensors of IR-T of excenters}
Summarizing the previous results, the following theorem is obtained.


\begin{theorem}{Tensors of IR-T of excenters, Daiyuan Zhang}{Thm17.3.1}\label{Thm17.3.1} 
The tensors of IR-T of excenters for the tetrahedron $ABCD$ can be expressed by the following formulas.

The tensor of IR-T of excenter $E_{A}$:
\[\boldsymbol{\lambda }_{ABCD}^{E_{A}^{A}}=\left( \begin{matrix}
	\boldsymbol{\lambda }_{\triangle BCD}^{E_{A}^{A}} & \boldsymbol{\lambda }_{\triangle CDA}^{E_{A}^{B}} & \boldsymbol{\lambda }_{\triangle DAB}^{E_{A}^{C}} & \boldsymbol{\lambda }_{\triangle ABC}^{E_{A}^{D}}  \\
\end{matrix} \right),\]
\[\boldsymbol{\lambda }_{\triangle BCD}^{E_{A}^{A}}=\left( \begin{matrix}
	\lambda _{BC}^{E_{A}^{A}} & \lambda _{CD}^{E_{A}^{A}} & \lambda _{DB}^{E_{A}^{A}}  \\
\end{matrix} \right)=\left( \begin{matrix}
	\frac{{{S}^{C}}}{{{S}^{B}}} & \frac{{{S}^{D}}}{{{S}^{C}}} & \frac{{{S}^{B}}}{{{S}^{D}}}  \\
\end{matrix} \right),\]
\[\boldsymbol{\lambda }_{\triangle CDA}^{E_{A}^{B}}=\left( \begin{matrix}
	\lambda _{CD}^{E_{A}^{B}} & \lambda _{DA}^{E_{A}^{B}} & \lambda _{AC}^{E_{A}^{B}}  \\
\end{matrix} \right)=\left( \begin{matrix}
	\frac{{{S}^{D}}}{{{S}^{C}}} & -\frac{{{S}^{A}}}{{{S}^{D}}} & -\frac{{{S}^{C}}}{{{S}^{A}}}  \\
\end{matrix} \right),\]
\[\boldsymbol{\lambda }_{\triangle DAB}^{E_{A}^{C}}=\left( \begin{matrix}
	\lambda _{DA}^{E_{A}^{C}} & \lambda _{AB}^{E_{A}^{C}} & \lambda _{BD}^{E_{A}^{C}}  \\
\end{matrix} \right)=\left( \begin{matrix}
	-\frac{{{S}^{A}}}{{{S}^{D}}} & -\frac{{{S}^{B}}}{{{S}^{A}}} & \frac{{{S}^{D}}}{{{S}^{B}}}  \\
\end{matrix} \right),\]
\[\boldsymbol{\lambda }_{\triangle ABC}^{E_{A}^{D}}=\left( \begin{matrix}
	\lambda _{AB}^{E_{A}^{D}} & \lambda _{BC}^{E_{A}^{D}} & \lambda _{CA}^{E_{A}^{D}}  \\
\end{matrix} \right)=\left( \begin{matrix}
	-\frac{{{S}^{B}}}{{{S}^{A}}} & \frac{{{S}^{C}}}{{{S}^{B}}} & -\frac{{{S}^{A}}}{{{S}^{C}}}  \\
\end{matrix} \right);\]

The tensor of IR-T of excenter $E_{B}$:
\[\boldsymbol{\lambda }_{ABCD}^{E_{B}^{A}}=\left( \begin{matrix}
	\boldsymbol{\lambda }_{\triangle BCD}^{E_{B}^{A}} & \boldsymbol{\lambda }_{\triangle CDA}^{E_{B}^{B}} & \boldsymbol{\lambda }_{\triangle DAB}^{E_{B}^{C}} & \boldsymbol{\lambda }_{\triangle ABC}^{E_{B}^{D}}  \\
\end{matrix} \right),\]
\[\boldsymbol{\lambda }_{\triangle BCD}^{E_{B}^{A}}=\left( \begin{matrix}
	\lambda _{BC}^{E_{B}^{A}} & \lambda _{CD}^{E_{B}^{A}} & \lambda _{DB}^{E_{B}^{A}}  \\
\end{matrix} \right)=\left( \begin{matrix}
	-\frac{{{S}^{C}}}{{{S}^{B}}} & \frac{{{S}^{D}}}{{{S}^{C}}} & -\frac{{{S}^{B}}}{{{S}^{D}}}  \\
\end{matrix} \right),\]
\[\boldsymbol{\lambda }_{\triangle CDA}^{E_{B}^{B}}=\left( \begin{matrix}
	\lambda _{CD}^{E_{B}^{B}} & \lambda _{DA}^{E_{B}^{B}} & \lambda _{AC}^{E_{B}^{B}}  \\
\end{matrix} \right)=\left( \begin{matrix}
	\frac{{{S}^{D}}}{{{S}^{C}}} & \frac{{{S}^{A}}}{{{S}^{D}}} & \frac{{{S}^{C}}}{{{S}^{A}}}  \\
\end{matrix} \right),\]
\[\boldsymbol{\lambda }_{\triangle DAB}^{E_{B}^{C}}=\left( \begin{matrix}
	\lambda _{DA}^{E_{B}^{C}} & \lambda _{AB}^{E_{B}^{C}} & \lambda _{BD}^{E_{B}^{C}}  \\
\end{matrix} \right)=\left( \begin{matrix}
	\frac{{{S}^{A}}}{{{S}^{D}}} & -\frac{{{S}^{B}}}{{{S}^{A}}} & -\frac{{{S}^{D}}}{{{S}^{B}}}  \\
\end{matrix} \right),\]
\[\boldsymbol{\lambda }_{\triangle ABC}^{E_{B}^{D}}=\left( \begin{matrix}
	\lambda _{AB}^{E_{B}^{D}} & \lambda _{BC}^{E_{B}^{D}} & \lambda _{CA}^{E_{B}^{D}}  \\
\end{matrix} \right)=\left( \begin{matrix}
	-\frac{{{S}^{B}}}{{{S}^{A}}} & -\frac{{{S}^{C}}}{{{S}^{B}}} & \frac{{{S}^{A}}}{{{S}^{C}}}  \\
\end{matrix} \right);\]

The tensor of IR-T of excenter $E_{C}$:
\[\boldsymbol{\lambda }_{ABCD}^{E_{C}^{A}}=\left( \begin{matrix}
	\boldsymbol{\lambda }_{\triangle BCD}^{E_{C}^{A}} & \boldsymbol{\lambda }_{\triangle CDA}^{E_{C}^{B}} & \boldsymbol{\lambda }_{\triangle DAB}^{E_{C}^{C}} & \boldsymbol{\lambda }_{\triangle ABC}^{E_{C}^{D}}  \\
\end{matrix} \right),\]
\[\boldsymbol{\lambda }_{\triangle BCD}^{E_{C}^{A}}=\left( \begin{matrix}
	\lambda _{BC}^{E_{C}^{A}} & \lambda _{CD}^{E_{C}^{A}} & \lambda _{DB}^{E_{C}^{A}}  \\
\end{matrix} \right)=\left( \begin{matrix}
	-\frac{{{S}^{C}}}{{{S}^{B}}} & -\frac{{{S}^{D}}}{{{S}^{C}}} & \frac{{{S}^{B}}}{{{S}^{D}}}  \\
\end{matrix} \right),\]
\[\boldsymbol{\lambda }_{\triangle CDA}^{E_{C}^{B}}=\left( \begin{matrix}
	\lambda _{CD}^{E_{C}^{B}} & \lambda _{DA}^{E_{C}^{B}} & \lambda _{AC}^{E_{C}^{B}}  \\
\end{matrix} \right)=\left( \begin{matrix}
	-\frac{{{S}^{D}}}{{{S}^{C}}} & \frac{{{S}^{A}}}{{{S}^{D}}} & -\frac{{{S}^{C}}}{{{S}^{A}}}  \\
\end{matrix} \right),\]
\[\boldsymbol{\lambda }_{\triangle DAB}^{E_{C}^{C}}=\left( \begin{matrix}
	\lambda _{DA}^{E_{C}^{C}} & \lambda _{AB}^{E_{C}^{C}} & \lambda _{BD}^{E_{C}^{C}}  \\
\end{matrix} \right)=\left( \begin{matrix}
	\frac{{{S}^{A}}}{{{S}^{D}}} & \frac{{{S}^{B}}}{{{S}^{A}}} & \frac{{{S}^{D}}}{{{S}^{B}}}  \\
\end{matrix} \right),\]
\[\boldsymbol{\lambda }_{\triangle ABC}^{E_{C}^{D}}=\left( \begin{matrix}
	\lambda _{AB}^{E_{C}^{D}} & \lambda _{BC}^{E_{C}^{D}} & \lambda _{CA}^{E_{C}^{D}}  \\
\end{matrix} \right)=\left( \begin{matrix}
	\frac{{{S}^{B}}}{{{S}^{A}}} & -\frac{{{S}^{C}}}{{{S}^{B}}} & -\frac{{{S}^{A}}}{{{S}^{C}}}  \\
\end{matrix} \right);\]

The tensor of IR-T of excenter $E_{D}$:
\[\boldsymbol{\lambda }_{ABCD}^{E_{D}^{A}}=\left( \begin{matrix}
	\boldsymbol{\lambda }_{\triangle BCD}^{E_{D}^{A}} & \boldsymbol{\lambda }_{\triangle CDA}^{E_{D}^{B}} & \boldsymbol{\lambda }_{\triangle DAB}^{E_{D}^{C}} & \boldsymbol{\lambda }_{\triangle ABC}^{E_{D}^{D}}  \\
\end{matrix} \right),\]
\[\boldsymbol{\lambda }_{\triangle BCD}^{E_{D}^{A}}=\left( \begin{matrix}
	\lambda _{BC}^{E_{D}^{A}} & \lambda _{CD}^{E_{D}^{A}} & \lambda _{DB}^{E_{D}^{A}}  \\
\end{matrix} \right)=\left( \begin{matrix}
	\frac{{{S}^{C}}}{{{S}^{B}}} & -\frac{{{S}^{D}}}{{{S}^{C}}} & -\frac{{{S}^{B}}}{{{S}^{D}}}  \\
\end{matrix} \right),\]
\[\boldsymbol{\lambda }_{\triangle CDA}^{E_{D}^{B}}=\left( \begin{matrix}
	\lambda _{CD}^{E_{D}^{B}} & \lambda _{DA}^{E_{D}^{B}} & \lambda _{AC}^{E_{D}^{B}}  \\
\end{matrix} \right)=\left( \begin{matrix}
	-\frac{{{S}^{D}}}{{{S}^{C}}} & -\frac{{{S}^{A}}}{{{S}^{D}}} & \frac{{{S}^{C}}}{{{S}^{A}}}  \\
\end{matrix} \right),\]
\[\boldsymbol{\lambda }_{\triangle DAB}^{E_{D}^{C}}=\left( \begin{matrix}
	\lambda _{DA}^{E_{D}^{C}} & \lambda _{AB}^{E_{D}^{C}} & \lambda _{BD}^{E_{D}^{C}}  \\
\end{matrix} \right)=\left( \begin{matrix}
	-\frac{{{S}^{A}}}{{{S}^{D}}} & \frac{{{S}^{B}}}{{{S}^{A}}} & -\frac{{{S}^{D}}}{{{S}^{B}}}  \\
\end{matrix} \right),\]
\[\boldsymbol{\lambda }_{\triangle ABC}^{E_{D}^{D}}=\left( \begin{matrix}
	\lambda _{AB}^{E_{D}^{D}} & \lambda _{BC}^{E_{D}^{D}} & \lambda _{CA}^{E_{D}^{D}}  \\
\end{matrix} \right)=\left( \begin{matrix}
	\frac{{{S}^{B}}}{{{S}^{A}}} & \frac{{{S}^{C}}}{{{S}^{B}}} & \frac{{{S}^{A}}}{{{S}^{C}}}  \\
\end{matrix} \right).\]
\end{theorem}


\chapter{Theorem of tetrahedral frame and theorem of vector of intersecting center}\label{Ch18}
\thispagestyle{empty}


In Intercenter Geometry, vector of intersecting center of a tetrahedron plays an important role. The vector of intersecting center of a tetrahedron (abbreviated as VIC-T) is mainly divided into the vector from origin to intersecting center of the tetrahedron (abbreviated as VOIC-T) and the vector of two intersecting centers of the tetrahedron (abbreviated as VTICs-T). In this chapter, I will prove two important theorems: the theorem of VOIC-T and the theorem of VTICs-T. These two theorems are the core theorems of Space Intercenter Geometry.

Similar to Plane Intercenter Geometry, the VIC-T is composed of tetrahedral frame and frame component, which is the basis of other geometric quantities (such as distance, etc.). Many applications in the book are based on these two theorems in this chapter.

\section{Some basic theorems}\label{Sec18.1}
In order to deal with tetrahedral problems conveniently, a very useful theorem is given below. The author names the theorem as coplanar theorem of four vector terminal points, which is different from the unique representation theorem of coplanar vector (theorem \ref{thm:Thm3.3.1}).

\begin{theorem}{Theorem of coplanar vector of four terminal points, Daiyuan Zhang}{Thm18.1.1}\label{Thm18.1.1} 
	Suppose that four points $A$, $B$, $C$, $D$ do not coincide with each other, and any three of them are not collinear, $O\in {{\mathbb{R}}^{3}}$, then the necessary and sufficient conditions for four points $A$, $B$, $C$, $D$ to be coplanar are as follows:
		
	\[{{\beta}_{A}}\overrightarrow{OA}+{{\beta }_{B}}\overrightarrow{OB}+{{\beta }_{C}}\overrightarrow{OC}+{{\beta }_{D}}\overrightarrow{OD}=\vec{0},\]	
	\[{{\beta}_{A}}+{{\beta }_{B}}+{{\beta }_{C}}+{{\beta }_{D}}=0,\]	
	where the real numbers ${{\beta }_{A}}\ne 0$, ${{\beta }_{B}}\ne 0$, ${{\beta }_{C}}\ne 0$, ${{\beta }_{D}}\ne 0$, and if any one of four real numbers ${{\beta }_{A}}$, ${{\beta }_{B}}$, ${{\beta }_{C}}$ and ${{\beta }_{D}}$ is 1, the other three numbers are uniquely determined.
\end{theorem}

\begin{proof}
	\textbf{Sufficiency}. Based on the conditions of the theorem, we have 	
	\[\overrightarrow{OA}=\left( -\frac{{{\beta }_{B}}}{{{\beta }_{A}}} \right)\overrightarrow{OB}+\left( -\frac{{{\beta }_{C}}}{{{\beta }_{A}}} \right)\overrightarrow{OC}+\left( -\frac{{{\beta }_{D}}}{{{\beta }_{A}}} \right)\overrightarrow{OD},\]	
	\[\left( -\frac{{{\beta }_{B}}}{{{\beta }_{A}}} \right)+\left( -\frac{{{\beta }_{C}}}{{{\beta }_{A}}} \right)+\left( -\frac{{{\beta }_{D}}}{{{\beta }_{A}}} \right)=1.\]	
	
	According to theorem \ref{thm:Thm3.3.1}, the coplanarity of four terminal points $A$, $B$, $C$ and $D$ of vector $\overrightarrow{OA}$, $\overrightarrow{OB}$, $\overrightarrow{OC}$ and $\overrightarrow{OD}$ is obtained.
	
	
	\textbf{Necessity}. Suppose that given four endpoints $A$, $B$, $C$, $D$ are coplanar, then $\overrightarrow{DA}$ and $\overrightarrow{DB}$, $\overrightarrow{DC}$ are coplanar. And because the four endpoints $A$, $B$, $C$, $D$ do not coincide with each other, so $\overrightarrow{DA}\ne \vec{0}$, $\overrightarrow{DB}\ne \vec{0}$, $\overrightarrow{DC}\ne 0$. Because any three points are not collinear, so the three points $B$, $C$, $D$ are not collinear. Therefore, $\overrightarrow{DB}$ and $\overrightarrow{DC}$ are linearly independent. According to the fundamental theorem of plane vector, there is a unique set of real numbers ${{\gamma }_{1}}$ and ${{\gamma }_{2}}$ makes the following formula true:
		
	\[\overrightarrow{DA}={{\gamma }_{1}}\overrightarrow{DB}+{{\gamma }_{2}}\overrightarrow{DC}.\]	
	
	
	The following proves ${{\gamma }_{1}}\ne 0$ and ${{\gamma }_{2}}\ne 0$. In fact, if one of the two quantities of ${{\gamma }_{1}}$ and ${{\gamma }_{2}}$ is 0, without losing generality, we might as well assume that ${{\gamma }_{2}}=0$, the above expression becomes $\overrightarrow{DA}={{\gamma }_{1}}\overrightarrow{DB}$, which means $A$, $B$, $D$ are collinear and contradictory. Similarly, ${{\gamma}_{1}}$ cannot be 0. If ${{\gamma}_{1}}$ and ${{\gamma}_{2}}$ are all 0, then the above formula becomes $\overrightarrow{DA}=\overrightarrow{0}$, which means that point $A$ and $D$ coincide and contradict. So there's only ${{\gamma}_{1}}\ne 0$ and ${{\gamma}_{2}}\ne 0$.
	
	For the point $O\in {{\mathbb{R}}^{3}}$, the above formula can be written as
	\[\overrightarrow{OA}-\overrightarrow{OD}={{\gamma }_{1}}\left( \overrightarrow{OB}-\overrightarrow{OD} \right)+{{\gamma }_{2}}\left( \overrightarrow{OC}-\overrightarrow{OD} \right),\]
	i.e.
	\[\overrightarrow{OA}={{\gamma }_{1}}\overrightarrow{OB}+{{\gamma }_{2}}\overrightarrow{OC}+\left( 1-{{\gamma }_{1}}-{{\gamma }_{2}} \right)\overrightarrow{OD},\]
	i.e.
	\begin{equation}\label{Eq18.1.1}
		\overrightarrow{OA}-{{\gamma }_{1}}\overrightarrow{OB}-{{\gamma }_{2}}\overrightarrow{OC}-{{\gamma }_{3}}\overrightarrow{OD}=\vec{0},
	\end{equation}
	
	\begin{equation}\label{Eq18.1.2}
		1+\left( -{{\gamma }_{1}} \right)+\left( -{{\gamma }_{2}} \right)+\left( -{{\gamma }_{3}} \right)=0,
	\end{equation}
	where
	\[{{\gamma }_{3}}=1-{{\gamma }_{1}}-{{\gamma }_{2}}.\]	
	
	Because the nonzero real numbers ${{\gamma }_{1}}$ and ${{\gamma }_{2}}$ are unique, so ${{\gamma}_{3}}$ is also unique, so when the coefficient of $\overrightarrow{OA}$ in the formula (\ref{Eq18.1.1}) is 1, each of the remaining three coefficients ${{\gamma }_{1}}$, ${{\gamma }_{2}}$ and ${{\gamma }_{3}}$ is uniquely determined.
	
	Formulas (\ref{Eq18.1.1}) and (\ref{Eq18.1.2}) can also be written as:	
	\[{{\beta }_{A}}\overrightarrow{OA}+{{\beta }_{B}}\overrightarrow{OB}+{{\beta }_{C}}\overrightarrow{OC}+{{\beta }_{D}}\overrightarrow{OD}=\vec{0},\]	
	\[{{\beta }_{A}}+{{\beta }_{B}}+{{\beta }_{C}}+{{\beta }_{D}}=0,\]	
	where ${{\beta }_{A}}=1\ne 0$, ${{\beta }_{B}}=-{{\gamma }_{1}}\ne 0$, ${{\beta }_{C}}=-{{\gamma }_{2}}\ne 0$, ${{\beta }_{D}}=-{{\gamma }_{3}}$.
	
	Let's prove that ${{\beta }_{D}}\ne 0$. Suppose ${{\beta }_{D}}=0$, so, ${{\gamma }_{3}}=-{{\beta }_{D}}=0$, then	
	\[\overrightarrow{OA}=\left( -\frac{{{\beta }_{B}}}{{{\beta }_{A}}} \right)\overrightarrow{OB}+\left( -\frac{{{\beta }_{C}}}{{{\beta }_{A}}} \right)\overrightarrow{OC},\]	
	\[\left( -\frac{{{\beta }_{B}}}{{{\beta }_{A}}} \right)+\left( -\frac{{{\beta }_{C}}}{{{\beta }_{A}}} \right)=1.\]	
	
	According to theorem \ref{thm:Thm3.1.1}, we get that the three terminal points of $\overrightarrow{OA}$, $\overrightarrow{OB}$ and $\overrightarrow{OC}$ are collinear, which is contradictory. So there's only ${{\beta }_{D}}=0$.
	
	
	It is possible to multiply both sides of equations (\ref{Eq18.1.1}) and (\ref{Eq18.1.2}) by a real number $\mu $, where $\mu \ne 0$, then there will be: $\mu \overrightarrow{OA}+\mu \left( -{{\gamma }_{1}} \right)\overrightarrow{OB}+\mu \left( -{{\gamma }_{2}} \right)\overrightarrow{OC}=\vec{0}$, $\mu +\mu \left( -{{\gamma }_{1}} \right)+\mu \left( -{{\gamma }_{2}} \right)=0$,  In this case, the three coefficients are ${{\beta }_{A}}=\mu $, ${{\beta }_{B}}=\mu \left( -{{\gamma }_{1}} \right)$, ${{\beta }_{C}}=\mu \left( -{{\gamma }_{2}} \right)$. Obviously, due to the arbitrariness of the real number $\mu \ne 0$, ${{\beta }_{A}}$, ${{\beta }_{B}}$ and ${{\beta }_{C}}$ are not unique. However, according to the previous discussion, when the coefficient of $\overrightarrow{OA}$ is 1, the other coefficients are uniquely determined.
\end{proof}
\hfill $\square$\par
%
The coefficient of $\overrightarrow{OA}$ studied above is ${{\beta }_{A}}=1$, we can obtain a similar result when ${{\beta}_{B}}=1$, ${{\beta }_{C}}=1$, which is left to the readers as an exercise.

For a tetrahedron, its four vertices $A$, $B$, $C$, $D$ are not coincident with each other, and any three of them are not collinear. Therefore, the above theorem is applicable to tetrahedron. The above theorem has applications in dealing with tetrahedral problems, See the discussion later.

According to the above theorem, we can get the following theorem of vector with non-coplanar terminal points.

\begin{theorem}{Vector with non-coplanar four terminal points, Daiyuan Zhang}{Thm18.1.2}\label{Thm18.1.2} 
	Suppose that the given four points $A$, $B$, $C$, $D$ do not coincide with each other, and any three of them are not collinear, $O\in {{\mathbb{R}}^{3}}$, and the following two conditions are satisfied:	
	\begin{equation}\label{Eq18.1.3}
		{{\beta }_{A}}\overrightarrow{OA}+{{\beta }_{B}}\overrightarrow{OB}+{{\beta }_{C}}\overrightarrow{OC}+{{\beta }_{D}}\overrightarrow{OD}=\vec{0},
	\end{equation}	
	\begin{equation}\label{Eq18.1.4}
		{{\beta }_{A}}+{{\beta }_{B}}+{{\beta }_{C}}+{{\beta }_{D}}=0,
	\end{equation}
	then, if the four terminal points $A$, $B$, $C$ and $D$ of the four position vectors $\overrightarrow{OA}$, $\overrightarrow{OB}$, $\overrightarrow{OC}$ and $\overrightarrow{OD}$ with the same initial point are not coplanar, there must be ${{\beta }_{A}}={{\beta }_{B}}={{\beta }_{C}}={{\beta }_{D}}=0$.
\end{theorem}

\begin{proof}
	
	When the condition of this theorem is satisfied, from theorem \ref{thm:Thm18.1.1}, if four terminal points $A$, $B$, $C$ and $D$ are not coplanar, then at least one of the four coefficients ${{\beta}_{A}}$, ${{\beta}_{B}}$, ${{\beta}_{C}}$, ${{\beta}_{ D}}$ is zero. Without losing generality, suppose ${{\beta}_{A}}=0$. ${{\beta}_{A} }=0$ is substituted into formula (\ref{Eq18.1.4}) to obtain ${{\beta}_{D}}=-{{\beta }_{B}}-{{\beta }_{C} }$, and then substitute it into formula (\ref{Eq18.1.3}) to get ${{\beta }_{B}}\left( \overrightarrow{OB}-\overrightarrow{OD} \right)+{{\beta }_{C}}\left( \overrightarrow{OC}-\overrightarrow{OD} \right)=\vec{0}$, which is ${{\beta }_{B}}\overrightarrow{DB}+{{\beta }_{C}}\overrightarrow{DC}=\vec{0}$. Since points $C$, $B$ and $D$ do not coincide and are not colinear with each other, there can be only $\overrightarrow{DB}\ne \vec{0}$, $\overrightarrow{DC}\ne \vec{0}$, and $\overrightarrow{DB}$ and $\overrightarrow{DC}$ are linear independent, so there can be only ${{\beta}_{B}}={{\beta }_{C}}=0$. According to the formula (\ref{Eq18.1.3}), we obtain ${{\beta}_{D}}=0$.
	
	Obviously, the same conclusion is reached when more than one of the four coefficients ${{\beta }_{A}}$, ${{\beta }_{B}}$, ${{\beta }_{C}}$ and ${{\beta }_{D}}$ is 0.
\end{proof}
\hfill $\square$\par



The four vertices of a tetrahedron can not be coplanar, and any three of them will not be collinear. Therefore, the theorem of vector with non-coplanar four terminal points (theorem \ref{thm:Thm18.1.2}) has many applications in dealing with tetrahedron problems. This theorem will be cited many times later in this book.


The importance of theorem \ref{thm:Thm18.1.2} lies in: when the conditions of the theorem are satisfied, as long as it can be determined that the four terminal points $A$, $B$, $C$ and $D$ of the vectors $\overrightarrow{OA}$, $\overrightarrow{OB}$, $\overrightarrow{OC}$ and $\overrightarrow{OD}$ are not coplanar, and any three of them will not be collinear, Then there must be ${{\beta }_{A}}={{\beta }_{B}}={{\beta }_{C}}={{\beta }_{D}}=0$. For tetrahedral frames $\left( O;A,B,C,D \right)$, its four vertices $A$, $B$, $C$ and $D$ are certainly not coplanar, and any three of them are not collinear, so there can only be ${{\beta }_{A}}={{\beta }_{B}}={{\beta }_{C}}={{\beta }_{D}}=0$. According to the above theorem, we can get the following theorem:

\begin{theorem}{Theorem of tetrahedral frame, Daiyuan Zhang}{Thm18.1.3}\label{Thm18.1.3} 
	Given a tetrahedron $ABCD$ and the tetrahedral frame $\left( O;A,B,C,D \right)$, $O$ is any point in space and satisfies the following two conditions:
	\[{{\beta }_{A}}\overrightarrow{OA}+{{\beta }_{B}}\overrightarrow{OB}+{{\beta }_{C}}\overrightarrow{OC}+{{\beta }_{D}}\overrightarrow{OD}=\vec{0},\]			
	\[{{\beta }_{A}}+{{\beta }_{B}}+{{\beta }_{C}}+{{\beta }_{D}}=0,\]	
	then there must be ${{\beta }_{A}}={{\beta }_{B}}={{\beta }_{C}}={{\beta }_{D}}=0$.	
\end{theorem}
\hfill $\square$\par

The following theorem takes advantage of the above conclusion.
\begin{theorem}{Unique representation theorem of vector on tetrahedral frame, Daiyuan Zhang}{Thm18.1.4}\label{Thm18.1.4} 
	Given a tetrahedral frame $\left( O;A,B,C,D \right)$, point $P\in {{\mathbb{R}}^{3}}$, if there is a set of real numbers ${{\beta }_{A}^{P}}$, ${{\beta }_{B}^{P}}$, ${{\beta }_{C}^{P}}$ and ${{\beta }_{D}^{P}}$ make the following formulas true:	
	\[\overrightarrow{OP}={{\beta }_{A}^{P}}\overrightarrow{OA}+{{\beta }_{B}^{P}}\overrightarrow{OB}+{{\beta }_{C}^{P}}\overrightarrow{OC}+{{\beta }_{D}^{P}}\overrightarrow{OD},\]	
	\[{{\beta }_{A}^{P}}+{{\beta }_{B}^{P}}+{{\beta }_{C}^{P}}+{{\beta }_{D}^{P}}=1,\]
	then each of the real numbers ${{\beta }_{A}^{P}}$, ${{\beta }_{B}^{P}}$, ${{\beta }_{C}^{P}}$ and ${{\beta }_{D}^{P}}$ is unique.
\end{theorem}

\begin{proof}
	If two sets of real numbers ${{\beta }_{A}^{P}}$, ${{\beta }_{B}^{P}}$, ${{\beta }_{C}^{P}}$, ${{\beta }_{D}^{P}}$ and ${{\gamma }_{A}^{P}}$, ${{\gamma }_{B}^{P}}$, ${{\gamma }_{C}^{P}}$, ${{\gamma }_{D}^{P}}$ satisfy the above conditions at the same time, then
	\[\left( {{\beta }_{A}^{P}}-{{\gamma }_{A}^{P}} \right)\overrightarrow{OA}+\left( {{\beta }_{B}^{P}}-{{\gamma }_{B}^{P}} \right)\overrightarrow{OB}+\left( {{\beta }_{C}^{P}}-{{\gamma }_{C}^{P}} \right)\overrightarrow{OC}+\left( {{\beta }_{D}^{P}}-{{\gamma }_{D}^{P}} \right)\overrightarrow{OD}=\overrightarrow{0},\]	\[\left( {{\beta }_{A}^{P}}-{{\gamma }_{A}^{P}} \right)+\left( {{\beta }_{B}^{P}}-{{\gamma }_{B}^{P}} \right)+\left( {{\beta }_{C}^{P}}-{{\gamma }_{C}^{P}} \right)+\left( {{\beta }_{D}^{P}}-{{\gamma }_{D}^{P}} \right)=1-1=0.\]
	
	According to theorem \ref{thm:Thm18.1.3}, the following results are obtained: ${{\beta }_{A}^{P}}-{{\gamma }_{A}^{P}}=0$, ${{\beta }_{B}^{P}}-{{\gamma }_{B}^{P}}=0$, ${{\beta }_{C}^{P}}-{{\gamma }_{C}^{P}}=0$, ${{\beta }_{D}^{P}}-{{\gamma }_{D}^{P}}=0$.
\end{proof}
\hfill $\square$\par\

The above theorem shows that under the condition of the theorem, the vector $\overrightarrow{OP}$ can be represented by the tetrahedral frame $\left( O;A,B,C,D \right)$ uniquely and linearly, where ${{\beta }_{A}^{P}}$, ${{\beta }_{B}^{P}}$, ${{\beta }_{C}^{P}}$, ${{\beta }_{D}^{P}}$ are called the affine frame components of the frame $\overrightarrow{OA}$, $\overrightarrow{OB}$, $\overrightarrow{OC}$, $\overrightarrow{OD}$ at point $P$, which are called components or coefficients for short.

\section{Theorem of vector from origin to intersecting center of a tetrahedron}\label{Sec18.2}

The unique representation theorem of vector on tetrahedral frame (theorem \ref{thm:Thm18.1.4}) is the basis of studying tetrahedron, and is very important. Without confusion, the theorem is also called tetrahedral frame representation theorem.

Before introducing the theorem, we first introduce a concept. As shown in Figure \ref{fig:tu16.2.1}, let
\begin{flalign*}
	{{\kappa }_{A{{P}_{A}}}}=\frac{\overrightarrow{AP}}{\overrightarrow{A{{P}_{A}}}}, {{\kappa }_{B{{P}_{B}}}}=\frac{\overrightarrow{BP}}{\overrightarrow{B{{P}_{B}}}}, {{\kappa }_{C{{P}_{C}}}}=\frac{\overrightarrow{CP}}{\overrightarrow{C{{P}_{C}}}}, {{\kappa }_{D{{P}_{D}}}}=\frac{\overrightarrow{DP}}{\overrightarrow{D{{P}_{D}}}}.
\end{flalign*}

Its meaning is: the integral ratio of the line segment between the vertex of the tetrahedron and the intersection point of the opposite face (through the IC-T), which is called the \textbf{integral ratio of the line segment from vertex to intersecting foot} of the tetrahedron (abbreviated as IRLSVIF-T), which is called the \textbf{integral ratio from vertex to intersecting foot} of the tetrahedron for short (abbreviated as IRVIF-T).

The intersecting foot of a tetrahedron (abbreviated as IF-T) is on the opposite face of the vertex, while the intersecting foot of a triangle is on the opposite face of the vertex.

Before proving the following theorem, we should first give some notations (see Figure \ref{fig:tu16.2.1}).


In this book, the frame component of tetrahedron is denoted by $\beta $; The frame component of triangle is denoted by $\alpha $. A tetrahedron has one frame. A tetrahedron has four face triangles, each triangle has its own triangular frame, the frames of different triangles are different from each other.

The subscripts of $\beta $ and $\alpha $ correspond to the frame, and the superscript corresponds to the points.


For example, the point ${P}_{D} \in \triangle ABC$. The frame of $\triangle ABC$ is $\left( O;A,B,C \right)$. Those quantities $\alpha _{A}^{{{P}_{D}}}$, $\alpha _{B}^{{{P}_{D}}}$, $\alpha _{C}^{{{P}_{D}}}$ with superscript ${{P}_{D}}$ correspond to the frame components of $\overrightarrow{OA}$, $\overrightarrow{OB}$, $\overrightarrow{OC}$ at ${{P}_{D}}$ respectively. If each of the frame components has nothing to do with the position of the frame origin $O$, the $\alpha _{A}^{{{P}_{D}}}$, $\alpha _{B}^{{{P}_{D}}}$, $\alpha _{C}^{{{P}_{D}}}$ are also called the frame components of points $A$, $B$ and $C$ corresponding to point ${{P}_{D}}$ respectively. Obviously, in the tetrahedral frame $\left( O;A,B,C,D \right)$, there is a triangular frame $\left( O;A,B,C \right)$ corresponding to $\triangle ABC$, then the $\alpha _{A}^{{{P}_{D}}}$, $\alpha _{B}^{{{P}_{D}}}$, $\alpha _{C}^{{{P}_{D}}}$ are the frame components corresponding to the ${{P}_{D}}$ of $\triangle ABC$ on the triangular frame $\left( O;A,B,C \right)$. From theorem \ref{thm:Thm6.1.1}, the following formula is obtained:
\[\alpha _{A}^{{{P}_{D}}}+\alpha _{B}^{{{P}_{D}}}+\alpha _{C}^{{{P}_{D}}}=1.\]

The point ${{P}_{C}} \in \triangle DAB$. The frame of $\triangle DAB$ is $\left( O;D,A,B \right)$. Those quantities $\alpha _{D}^{{{P}_{C}}}$, $\alpha _{A}^{{{P}_{C}}}$, $\alpha _{B}^{{{P}_{C}}}$ with superscript ${{P}_{C}}$ correspond to the frame components of $\overrightarrow{OD}$, $\overrightarrow{OA}$, $\overrightarrow{OB}$ at point ${{P}_{C}}$ respectively. If each of the frame component has nothing to do with the position of the frame origin $O$, the $\alpha _{D}^{{{P}_{C}}}$, $\alpha _{A}^{{{P}_{C}}}$, $\alpha _{B}^{{{P}_{C}}}$ are also called the frame components of points $D$, $A$, $B$ corresponding to point ${{P}_{C}}$ respectively. Obviously, in the tetrahedral frame $\left( O;A,B,C,D \right)$, there is a triangular frame $\left( O;D,A,B \right)$ corresponding to the $\triangle DAB$, then the $\alpha _{D}^{{{P}_{C}}}$, $\alpha _{A}^{{{P}_{C}}}$, $\alpha _{B}^{{{P}_{C}}}$ are the frame components corresponding to the ${{P}_{C}}$ of $\triangle DAB$ on the triangular frame $\left( O;D,A,B \right)$. From theorem \ref{thm:Thm6.1.1}, the following formula is obtained:
\[\alpha _{D}^{{{P}_{C}}}+\alpha _{A}^{{{P}_{C}}}+\alpha _{B}^{{{P}_{C}}}=1.\]

The rest of the notations are similar in meaning.


%
%

\begin{theorem}{Vector from origin to intersecting center of a tetrahedron, Daiyuan Zhang}{Thm18.2.1}\label{Thm18.2.1} 
	Suppose that given a tetrahedron $ABCD$, the point $O$ is an arbitrary point, and the point $P\in {{\pi }_{ABCD}}$, then the vector $\overrightarrow{OP}$ can be expressed by the uniquely linear combination of the following form:
	\begin{equation}\label{Eq18.2.1}
		\overrightarrow{OP}=\beta _{A}^{P}\overrightarrow{OA}+\beta _{B}^{P}\overrightarrow{OB}+\beta _{C}^{P}\overrightarrow{OC}+\beta _{D}^{P}\overrightarrow{OD},
	\end{equation}
	\begin{equation}\label{Eq18.2.2}
		\beta _{A}^{P}+\beta _{B}^{P}+\beta _{C}^{P}+\beta _{D}^{P}=1.
	\end{equation}
	
	Where
	\begin{equation}\label{Eq18.2.3}
		\beta _{A}^{P}=\frac{1}{4}\left( 1-{{\kappa }_{A{{P}_{A}}}}+{{\kappa }_{B{{P}_{B}}}}\alpha _{A}^{{{P}_{B}}}+{{\kappa }_{C{{P}_{C}}}}\alpha _{A}^{{{P}_{C}}}+{{\kappa }_{D{{P}_{D}}}}\alpha _{A}^{{{P}_{D}}} \right),
	\end{equation}
	\begin{equation}\label{Eq18.2.4}
		\beta _{B}^{P}=\frac{1}{4}\left( 1-{{\kappa }_{B{{P}_{B}}}}+{{\kappa }_{C{{P}_{C}}}}\alpha _{B}^{{{P}_{C}}}+{{\kappa }_{D{{P}_{D}}}}\alpha _{B}^{{{P}_{D}}}+{{\kappa }_{A{{P}_{A}}}}\alpha _{B}^{{{P}_{A}}} \right),
	\end{equation}
	\begin{equation}\label{Eq18.2.5}
		\beta _{C}^{P}=\frac{1}{4}\left( 1-{{\kappa }_{C{{P}_{C}}}}+{{\kappa }_{D{{P}_{D}}}}\alpha _{C}^{{{P}_{D}}}+{{\kappa }_{A{{P}_{A}}}}\alpha _{C}^{{{P}_{A}}}+{{\kappa }_{B{{P}_{B}}}}\alpha _{C}^{{{P}_{B}}} \right),
	\end{equation}
	\begin{equation}\label{Eq18.2.6}
		\beta _{D}^{P}=\frac{1}{4}\left( 1-{{\kappa }_{D{{P}_{D}}}}+{{\kappa }_{A{{P}_{A}}}}\alpha _{D}^{{{P}_{A}}}+{{\kappa }_{B{{P}_{B}}}}\alpha _{D}^{{{P}_{B}}}+{{\kappa }_{C{{P}_{C}}}}\alpha _{D}^{{{P}_{C}}} \right).
	\end{equation}
	
	Or
	\begin{flalign}\label{Eq18.2.7}
		\beta _{A}^{P}=1-{{\kappa }_{A{{P}_{A}}}}, \beta _{B}^{P}={{\kappa }_{A{{P}_{A}}}}\alpha _{B}^{{{P}_{A}}}, \beta _{C}^{P}={{\kappa }_{A{{P}_{A}}}}\alpha _{C}^{{{P}_{A}}}, \beta _{D}^{P}={{\kappa }_{A{{P}_{A}}}}\alpha _{D}^{{{P}_{A}}},
	\end{flalign}
	\begin{flalign}\label{Eq18.2.8}
		\beta _{A}^{P}={{\kappa }_{B{{P}_{B}}}}\alpha _{A}^{{{P}_{B}}}, \beta _{B}^{P}=1-{{\kappa }_{B{{P}_{B}}}}, \beta _{C}^{P}={{\kappa }_{B{{P}_{B}}}}\alpha _{C}^{{{P}_{B}}}, \beta _{D}^{P}={{\kappa }_{B{{P}_{B}}}}\alpha _{D}^{{{P}_{B}}},
	\end{flalign}
	\begin{flalign}\label{Eq18.2.9}
		\beta _{A}^{P}={{\kappa }_{C{{P}_{C}}}}\alpha _{A}^{{{P}_{C}}}, \beta _{B}^{P}={{\kappa }_{C{{P}_{C}}}}\alpha _{B}^{{{P}_{C}}}, \beta _{C}^{P}=1-{{\kappa }_{C{{P}_{C}}}},\beta _{D}^{P}={{\kappa }_{C{{P}_{C}}}}\alpha _{D}^{{{P}_{C}}},
	\end{flalign}
	\begin{flalign}\label{Eq18.2.10}
		\beta _{A}^{P}={{\kappa }_{D{{P}_{D}}}}\alpha _{A}^{{{P}_{D}}},\beta _{B}^{P}={{\kappa }_{D{{P}_{D}}}}\alpha _{B}^{{{P}_{D}}},\beta _{C}^{P}={{\kappa }_{D{{P}_{D}}}}\alpha _{C}^{{{P}_{D}}},\beta _{D}^{P}=1-{{\kappa }_{D{{P}_{D}}}}.
	\end{flalign}
\end{theorem}

\begin{proof}
	As shown in Figure \ref{fig:tu16.2.1}, there are:	
	\[\begin{aligned}
		 \overrightarrow{OP}&=\overrightarrow{OA}+\overrightarrow{AP}=\overrightarrow{OA}+{{\kappa }_{A{{P}_{A}}}}\overrightarrow{A{{P}_{A}}} \\ 
		& =\overrightarrow{OA}+{{\kappa }_{A{{P}_{A}}}}\left( \alpha _{B}^{{{P}_{A}}}\overrightarrow{AB}+\alpha _{C}^{{{P}_{A}}}\overrightarrow{AC}+\alpha _{D}^{{{P}_{A}}}\overrightarrow{AD} \right),  
	\end{aligned}\]
	i.e.
	\[\begin{aligned}
		\overrightarrow{OP}&=\overrightarrow{OA}+{{\kappa }_{A{{P}_{A}}}}\left( \alpha _{B}^{{{P}_{A}}}\left( \overrightarrow{OB}-\overrightarrow{OA} \right)+\alpha _{C}^{{{P}_{A}}}\left( \overrightarrow{OC}-\overrightarrow{OA} \right)+\alpha _{D}^{{{P}_{A}}}\left( \overrightarrow{OD}-\overrightarrow{OA} \right) \right) \\ 
		& =\overrightarrow{OA}+{{\kappa }_{A{{P}_{A}}}}\left( \alpha _{B}^{{{P}_{A}}}\overrightarrow{OB}+\alpha _{C}^{{{P}_{A}}}\overrightarrow{OC}+\alpha _{D}^{{{P}_{A}}}\overrightarrow{OD}-\left( \alpha _{B}^{{{P}_{A}}}+\alpha _{C}^{{{P}_{A}}}+\alpha _{D}^{{{P}_{A}}} \right)\overrightarrow{OA} \right),  
	\end{aligned}\]
	i.e.
	\[\begin{aligned}
		\overrightarrow{OP}&=\overrightarrow{OA}+{{\kappa }_{A{{P}_{A}}}}\left( \alpha _{B}^{{{P}_{A}}}\overrightarrow{OB}+\alpha _{C}^{{{P}_{A}}}\overrightarrow{OC}+\alpha _{D}^{{{P}_{A}}}\overrightarrow{OD}-\overrightarrow{OA} \right) \\ 
		& =\left( 1-{{\kappa }_{A{{P}_{A}}}} \right)\overrightarrow{OA}+{{\kappa }_{A{{P}_{A}}}}\left( \alpha _{B}^{{{P}_{A}}}\overrightarrow{OB}+\alpha _{C}^{{{P}_{A}}}\overrightarrow{OC}+\alpha _{D}^{{{P}_{A}}}\overrightarrow{OD} \right).  
	\end{aligned}\]
	
	The sum of frame components of $\overrightarrow{OA}$, $\overrightarrow{OB}$, $\overrightarrow{OC}$, $\overrightarrow{OD}$ in the above formula is 1. In fact, according to theorem \ref{thm:Thm6.1.1}, we get		
	\[\left( 1-{{\kappa }_{A{{P}_{A}}}} \right)+{{\kappa }_{A{{P}_{A}}}}\left( \alpha _{B}^{{{P}_{A}}}+\alpha _{C}^{{{P}_{A}}}+\alpha _{D}^{{{P}_{A}}} \right)=1-{{\kappa }_{A{{P}_{A}}}}+{{\kappa }_{A{{P}_{A}}}}=1.\]
	
	Therefore, from theorem \ref{thm:Thm18.1.4}, the expression for vector $\overrightarrow{OP}$ is unique.
	
	Similarly,
	\[\begin{aligned}
		\overrightarrow{OP}&=\overrightarrow{OB}+\overrightarrow{BP}=\overrightarrow{OB}+{{\kappa }_{B{{P}_{B}}}}\overrightarrow{B{{P}_{B}}} \\ 
		& =\overrightarrow{OB}+{{\kappa }_{B{{P}_{B}}}}\left( \alpha _{C}^{{{P}_{B}}}\overrightarrow{BC}+\alpha _{D}^{{{P}_{B}}}\overrightarrow{BD}+\alpha _{A}^{{{P}_{B}}}\overrightarrow{BA} \right),  
	\end{aligned}\]
	i.e.
	\[\begin{aligned}
		\overrightarrow{OP}&=\overrightarrow{OB}+{{\kappa }_{B{{P}_{B}}}}\left( \alpha _{C}^{{{P}_{B}}}\left( \overrightarrow{OC}-\overrightarrow{OB} \right)+\alpha _{D}^{{{P}_{B}}}\left( \overrightarrow{OD}-\overrightarrow{OB} \right)+\alpha _{A}^{{{P}_{B}}}\left( \overrightarrow{OA}-\overrightarrow{OB} \right) \right) \\ 
		& =\overrightarrow{OB}+{{\kappa }_{B{{P}_{B}}}}\left( \alpha _{C}^{{{P}_{B}}}\overrightarrow{OC}+\alpha _{D}^{{{P}_{B}}}\overrightarrow{OD}+\alpha _{A}^{{{P}_{B}}}\overrightarrow{OA}-\left( \alpha _{C}^{{{P}_{B}}}+\alpha _{D}^{{{P}_{B}}}+\alpha _{A}^{{{P}_{B}}} \right)\overrightarrow{OB}, \right)  
	\end{aligned}\]
	i.e.
	\[\begin{aligned}
		\overrightarrow{OP}&=\overrightarrow{OB}+{{\kappa }_{B{{P}_{B}}}}\left( \alpha _{C}^{{{P}_{B}}}\overrightarrow{OC}+\alpha _{D}^{{{P}_{B}}}\overrightarrow{OD}+\alpha _{A}^{{{P}_{B}}}\overrightarrow{OA}-\overrightarrow{OB} \right) \\ 
		& =\left( 1-{{\kappa }_{B{{P}_{B}}}} \right)\overrightarrow{OB}+{{\kappa }_{B{{P}_{B}}}}\left( \alpha _{C}^{{{P}_{B}}}\overrightarrow{OC}+\alpha _{D}^{{{P}_{B}}}\overrightarrow{OD}+\alpha _{A}^{{{P}_{B}}}\overrightarrow{OA} \right).  
	\end{aligned}\]
	
	The above formulas can also be obtained directly by rotating the previous formulas. The sum of frame components of $\overrightarrow{OA}$, $\overrightarrow{OB}$, $\overrightarrow{OC}$, $\overrightarrow{OD}$ is 1, and the above formula is also unique.	
	\[\begin{aligned}
		\overrightarrow{OP}&=\overrightarrow{OC}+\overrightarrow{CP}=\overrightarrow{OC}+{{\kappa }_{C{{P}_{C}}}}\overrightarrow{C{{P}_{C}}} \\ 
		& =\overrightarrow{OC}+{{\kappa }_{C{{P}_{C}}}}\left( \alpha _{D}^{{{P}_{C}}}\overrightarrow{CD}+\alpha _{A}^{{{P}_{C}}}\overrightarrow{CA}+\alpha _{B}^{{{P}_{C}}}\overrightarrow{CB} \right),  
	\end{aligned}\]
	i.e.
	\[\begin{aligned}
		\overrightarrow{OP}&=\overrightarrow{OC}+{{\kappa }_{C{{P}_{C}}}}\left( \alpha _{D}^{{{P}_{C}}}\left( \overrightarrow{OD}-\overrightarrow{OC} \right)+\alpha _{A}^{{{P}_{C}}}\left( \overrightarrow{OA}-\overrightarrow{OC} \right)+\alpha _{B}^{{{P}_{C}}}\left( \overrightarrow{OB}-\overrightarrow{OC} \right) \right) \\ 
		& =\overrightarrow{OC}+{{\kappa }_{C{{P}_{C}}}}\left( \alpha _{D}^{{{P}_{C}}}\overrightarrow{OD}+\alpha _{A}^{{{P}_{C}}}\overrightarrow{OA}+\alpha _{B}^{{{P}_{C}}}\overrightarrow{OB}-\left( \alpha _{D}^{{{P}_{C}}}+\alpha _{A}^{{{P}_{C}}}+\alpha _{B}^{{{P}_{C}}} \right)\overrightarrow{OC} \right),  
	\end{aligned}\]
	i.e.
	\[\begin{aligned}
		\overrightarrow{OP}&=\overrightarrow{OC}+{{\kappa }_{C{{P}_{C}}}}\left( \alpha _{D}^{{{P}_{C}}}\overrightarrow{OD}+\alpha _{A}^{{{P}_{C}}}\overrightarrow{OA}+\alpha _{B}^{{{P}_{C}}}\overrightarrow{OB}-\overrightarrow{OC} \right) \\ 
		& =\left( 1-{{\kappa }_{C{{P}_{C}}}} \right)\overrightarrow{OC}+{{\kappa }_{C{{P}_{C}}}}\left( \alpha _{D}^{{{P}_{C}}}\overrightarrow{OD}+\alpha _{A}^{{{P}_{C}}}\overrightarrow{OA}+\alpha _{B}^{{{P}_{C}}}\overrightarrow{OB} \right).  
	\end{aligned}\]
	
	In above formula, the sum of frame components of $\overrightarrow{OA}$, $\overrightarrow{OB}$, $\overrightarrow{OC}$, $\overrightarrow{OD}$ is 1, and the above formula is also unique.	
	\[\begin{aligned}
		\overrightarrow{OP}&=\overrightarrow{OD}+\overrightarrow{DP}=\overrightarrow{OD}+{{\kappa }_{D{{P}_{D}}}}\overrightarrow{D{{P}_{D}}} \\ 
		& =\overrightarrow{OD}+{{\kappa }_{D{{P}_{D}}}}\left( \alpha _{A}^{{{P}_{D}}}\overrightarrow{DA}+\alpha _{B}^{{{P}_{D}}}\overrightarrow{DB}+\alpha _{C}^{{{P}_{D}}}\overrightarrow{DC} \right)  
	\end{aligned}\]
	i.e.
	\[\begin{aligned}
		\overrightarrow{OP}&=\overrightarrow{OD}+{{\kappa }_{D{{P}_{D}}}}\left( \alpha _{A}^{{{P}_{D}}}\left( \overrightarrow{OA}-\overrightarrow{OD} \right)+\alpha _{B}^{{{P}_{D}}}\left( \overrightarrow{OB}-\overrightarrow{OD} \right)+\alpha _{C}^{{{P}_{D}}}\left( \overrightarrow{OC}-\overrightarrow{OD} \right) \right) \\ 
		& =\overrightarrow{OD}+{{\kappa }_{D{{P}_{D}}}}\left( \alpha _{A}^{{{P}_{D}}}\overrightarrow{OA}+\alpha _{B}^{{{P}_{D}}}\overrightarrow{OB}+\alpha _{C}^{{{P}_{D}}}\overrightarrow{OC}-\left( \alpha _{A}^{{{P}_{D}}}+\alpha _{B}^{{{P}_{D}}}+\alpha _{C}^{{{P}_{D}}} \right)\overrightarrow{OD}, \right)  
	\end{aligned}\]
	i.e.
	\[\begin{aligned}
		\overrightarrow{OP}&=\overrightarrow{OD}+{{\kappa }_{D{{P}_{D}}}}\left( \alpha _{A}^{{{P}_{D}}}\overrightarrow{OA}+\alpha _{B}^{{{P}_{D}}}\overrightarrow{OB}+\alpha _{C}^{{{P}_{D}}}\overrightarrow{OC}-\overrightarrow{OD} \right) \\ 
		& =\left( 1-{{\kappa }_{D{{P}_{D}}}} \right)\overrightarrow{OD}+{{\kappa }_{D{{P}_{D}}}}\left( \alpha _{A}^{{{P}_{D}}}\overrightarrow{OA}+\alpha _{B}^{{{P}_{D}}}\overrightarrow{OB}+\alpha _{C}^{{{P}_{D}}}\overrightarrow{OC} \right).  
	\end{aligned}\]
	
	In above formula, the sum of frame components of $\overrightarrow{OA}$, $\overrightarrow{OB}$, $\overrightarrow{OC}$, $\overrightarrow{OD}$ is 1, and the above formula is also unique.
	
	So the above four formulas for calculating $\overrightarrow{OP}$ are added on the left and right sides to get the following result: 
	\[\overrightarrow{OP}=\frac{1}{4}\left( \begin{aligned}
		& \left( 1-{{\kappa }_{A{{P}_{A}}}} \right)\overrightarrow{OA}+{{\kappa }_{A{{P}_{A}}}}\left( \alpha _{B}^{{{P}_{A}}}\overrightarrow{OB}+\alpha _{C}^{{{P}_{A}}}\overrightarrow{OC}+\alpha _{D}^{{{P}_{A}}}\overrightarrow{OD} \right) \\ 
		&+ \left( 1-{{\kappa }_{B{{P}_{B}}}} \right)\overrightarrow{OB}+{{\kappa }_{B{{P}_{B}}}}\left( \alpha _{C}^{{{P}_{B}}}\overrightarrow{OC}+\alpha _{D}^{{{P}_{B}}}\overrightarrow{OD}+\alpha _{A}^{{{P}_{B}}}\overrightarrow{OA} \right) \\ 
		&+ \left( 1-{{\kappa }_{C{{P}_{C}}}} \right)\overrightarrow{OC}+{{\kappa }_{C{{P}_{C}}}}\left( \alpha _{D}^{{{P}_{C}}}\overrightarrow{OD}+\alpha _{A}^{{{P}_{C}}}\overrightarrow{OA}+\alpha _{B}^{{{P}_{C}}}\overrightarrow{OB} \right) \\ 
		&+ \left( 1-{{\kappa }_{D{{P}_{D}}}} \right)\overrightarrow{OD}+{{\kappa }_{D{{P}_{D}}}}\left( \alpha _{A}^{{{P}_{D}}}\overrightarrow{OA}+\alpha _{B}^{{{P}_{D}}}\overrightarrow{OB}+\alpha _{C}^{{{P}_{D}}}\overrightarrow{OC} \right) \\ 
	\end{aligned} \right),\]
	i.e.
	\[\overrightarrow{OP}=\frac{1}{4}\left( \begin{aligned}
		& \left( 1-{{\kappa }_{A{{P}_{A}}}}+{{\kappa }_{B{{P}_{B}}}}\alpha _{A}^{{{P}_{B}}}+{{\kappa }_{C{{P}_{C}}}}\alpha _{A}^{{{P}_{C}}}+{{\kappa }_{D{{P}_{D}}}}\alpha _{A}^{{{P}_{D}}} \right)\overrightarrow{OA} \\ 
		& +\left( 1-{{\kappa }_{B{{P}_{B}}}}+{{\kappa }_{C{{P}_{C}}}}\alpha _{B}^{{{P}_{C}}}+{{\kappa }_{D{{P}_{D}}}}\alpha _{B}^{{{P}_{D}}}+{{\kappa }_{A{{P}_{A}}}}\alpha _{B}^{{{P}_{A}}} \right)\overrightarrow{OB} \\ 
		& +\left( 1-{{\kappa }_{C{{P}_{C}}}}+{{\kappa }_{D{{P}_{D}}}}\alpha _{C}^{{{P}_{D}}}+{{\kappa }_{A{{P}_{A}}}}\alpha _{C}^{{{P}_{A}}}+{{\kappa }_{B{{P}_{B}}}}\alpha _{C}^{{{P}_{B}}} \right)\overrightarrow{OC} \\ 
		& +\left( 1-{{\kappa }_{D{{P}_{D}}}}+{{\kappa }_{A{{P}_{A}}}}\alpha _{D}^{{{P}_{A}}}+{{\kappa }_{B{{P}_{B}}}}\alpha _{D}^{{{P}_{B}}}+{{\kappa }_{C{{P}_{C}}}}\alpha _{D}^{{{P}_{C}}} \right)\overrightarrow{OD} \\ 
	\end{aligned} \right).\]
	
	Similarly, in above formula, the sum of frame components of $\overrightarrow{OA}$, $\overrightarrow{OB}$, $\overrightarrow{OC}$, $\overrightarrow{OD}$ is 1, and the above formula is also unique.
	
	Obviously, the expression of the vector is linear.	
\end{proof}
\hfill $\square$\par


In the above theorem, $\beta _{A}^{P}$, $\beta _{B}^{P}$, $\beta _{C}^{P}$, $\beta _{D}^{P}$ are the frame components of $A$, $B$, $C$, $D$ corresponding to point $P$ respectively, or the frame components of $\overrightarrow{OA}$, $\overrightarrow{OB}$, $\overrightarrow{OC}$, $\overrightarrow{OD}$ on $\left( O;A,B,C,D \right)$ corresponding to point $P$ respectively.

The above theorem points out that the vector of intersecting center of tetrahedron (abbreviated as VIC-T) can be expressed by the uniquely linear combination of the tetrahedral frame $\left( O;A,B,C,D \right)$, and each frame component is constant and independent of the origin $O$ of the frame, or independent of the tetrahedral frame $\left( O;A,B,C,D \right)$.


The above theorem is one of the core theorems in this book, and has many applications in dealing with tetrahedral problems. The first group of expressions of the frame component in the theorem are symmetrical and graceful, but the form is a little complicated. The expressions of the latter groups of frame components are simple and convenient for calculation, but the form is not symmetrical. According to theorem \ref{thm:Thm18.1.4}, no matter which group of formulas are used, the components of the same frame are equal.

\subsection{Representation theorem using frame of edge  for the vector from origin to intersecting center of a tetrahedron}\label{Subsec18.2.1}
Theorem \ref{thm:Thm18.2.1} states that the origin $O$ of the tetrahedral frame $\left( O;A,B,C,D \right)$ can be anywhere in space. Sometimes, for the convenience of analysis and calculation, the origin $O$ is selected in some special positions. In this section, the origin $O$ is placed at the vertex of the tetrahedron $ABCD$, so that the frame vectors coincide with the edges of the tetrahedron. Such a frame is called the frame of edge  of the tetrahedron, which is called the frame of edge  for short. For example, when the origin $O$ is placed at the vertex $A$ of the tetrahedron $ABCD$, the frame of edge  is the set of $\overrightarrow{AB}$, $\overrightarrow{AC}$ and $\overrightarrow{AD}$.

\begin{theorem}{Vector from origin to IC-T represented by edge frame, Daiyuan Zhang}{Thm18.2.2}\label{Thm18.2.2} 
	Suppose that given a tetrahedron $ABCD$, an intersecting center of the tetrahedron $P\in {{\pi }_{ABCD}}$, then the vector from the vertex of the tetrahedron to the intersecting center of the tetrahedron $P$ can be expressed by the uniquely linear combination of the frame of edge of the tetrahedron, i.e.
	\[\left\{ \begin{aligned}
		& \overrightarrow{AP}=\beta _{B}^{P}\overrightarrow{AB}+\beta _{C}^{P}\overrightarrow{AC}+\beta _{D}^{P}\overrightarrow{AD} \\ 
		& \beta _{A}^{P}+\beta _{B}^{P}+\beta _{C}^{P}+\beta _{D}^{P}=1, \\ 
	\end{aligned} \right.\]	
	\[\left\{ \begin{aligned}
		& \overrightarrow{BP}=\beta _{A}^{P}\overrightarrow{BA}+\beta _{C}^{P}\overrightarrow{BC}+\beta _{D}^{P}\overrightarrow{BD} \\ 
		& \beta _{A}^{P}+\beta _{B}^{P}+\beta _{C}^{P}+\beta _{D}^{P}=1, \\ 
	\end{aligned} \right.\]	
	\[\left\{ \begin{aligned}
		& \overrightarrow{CP}=\beta _{A}^{P}\overrightarrow{CA}+\beta _{B}^{P}\overrightarrow{CB}+\beta _{D}^{P}\overrightarrow{CD} \\ 
		& \beta _{A}^{P}+\beta _{B}^{P}+\beta _{C}^{P}+\beta _{D}^{P}=1, \\ 
	\end{aligned} \right.\]	
	\[\left\{ \begin{aligned}
		& \overrightarrow{DP}=\beta _{A}^{P}\overrightarrow{DA}+\beta _{B}^{P}\overrightarrow{DB}+\beta _{C}^{P}\overrightarrow{DC} \\ 
		& \beta _{A}^{P}+\beta _{B}^{P}+\beta _{C}^{P}+\beta _{D}^{P}=1. \\ 
	\end{aligned} \right.\]	
	
	Where $\beta _{A}^{P}$, $\beta _{B}^{P}$, $\beta _{C}^{P}$, $\beta _{D}^{P}$ are the frame components of $A$, $B$, $C$, $D$ corresponding to point $P$, respectively.	
\end{theorem}
\begin{proof}
	This theorem is obtained by replacing the point $O$ in theorem \ref{thm:Thm18.2.1} with point $A$, $B$, $C$, $D$, respectively.
\end{proof}
\hfill $\square$\par

Frame of edge can reduce the number of frames, thus reducing the amount of computation, so this is an important special case.

\subsection{Relationship between vector from vertex to intersecting center and vector from vertex to intersecting foot}\label{Subsec18.2.2}
The vector from the vertex of a tetrahedron to the intersecting center of the tetrahedron is called the vector from vertex to intersecting center (abbreviated as VVIC). In Figure \ref{fig:tu16.2.1}, each of  $\overrightarrow{AP}$, $\overrightarrow{BP}$, $\overrightarrow{CP}$ and $\overrightarrow{DP}$ is the VVIC.

The vector from the vertex of a tetrahedron to the intersecting foot of the tetrahedron is called the vector from vertex to intersecting foot (abbreviated as VVIF). In Figure \ref{fig:tu16.2.1}, each of  $\overrightarrow{A{{P}_{A}}}$, $\overrightarrow{B{{P}_{B}}}$, $\overrightarrow{C{{P}_{C}}}$ and $\overrightarrow{D{{P}_{D}}}$ is the VVIF.

The following theorem points out that there is a certain relationship between the VVIC and the VVIF.

\begin{theorem}{Relationship between VVIC and VVIF of a tetrahedron, Daiyuan Zhang}{Thm18.2.3}\label{Thm18.2.3} 
	For a tetrahedron $ABCD$, the intersecting center is $P\in {{\pi }_{ABCD}}$, then each vector from its vertex to corresponding intersecting foot is
	\begin{flalign*}
		\overrightarrow{A{{P}_{A}}}=\frac{\overrightarrow{AP}}{1-\beta _{A}^{P}},\overrightarrow{B{{P}_{B}}}=\frac{\overrightarrow{BP}}{1-\beta _{B}^{P}},\overrightarrow{C{{P}_{C}}}=\frac{\overrightarrow{CP}}{1-\beta _{C}^{P}},\overrightarrow{D{{P}_{D}}}=\frac{\overrightarrow{DP}}{1-\beta _{C}^{P}}.	
	\end{flalign*}	
		
	Where $\beta _{A}^{P}$, $\beta _{B}^{P}$, $\beta _{C}^{P}$, $\beta _{D}^{P}$ are the frame components of $A$, $B$, $C$, $D$ corresponding to point $P$, respectively.
\end{theorem}

\begin{proof}
	Obviously, the following formula can be obtained from (\ref{Eq18.2.7}):
	\[\overrightarrow{A{{P}_{A}}}=\frac{\overrightarrow{AP}}{{{\kappa }_{AP_A}}}=\frac{\overrightarrow{AP}}{1-\beta _{A}^{P}}.\]
	
	The other formulas can be proved similarly.
\end{proof}
\hfill $\square$\par

\begin{theorem}{Representation theorem using edge frame for VVIC-T, Daiyuan Zhang}{Thm18.2.4}\label{Thm18.2.4} 
	For a tetrahedron $ABCD$, the intersecting center is $P\in {{\pi }_{ABCD}}$, and its vector from vertex to intersecting foot can be expressed by frame of edge as:  	
	\[\overrightarrow{A{{P}_{A}}}=\frac{\beta _{B}^{P}\overrightarrow{AB}+\beta _{C}^{P}\overrightarrow{AC}+\beta _{D}^{P}\overrightarrow{AD}}{\beta _{B}^{P}+\beta _{C}^{P}+\beta _{D}^{P}},\]	
	\[\overrightarrow{B{{P}_{B}}}=\frac{\beta _{A}^{P}\overrightarrow{BA}+\beta _{C}^{P}\overrightarrow{BC}+\beta _{D}^{P}\overrightarrow{BD}}{\beta _{A}^{P}+\beta _{C}^{P}+\beta _{D}^{P}},\]	
	\[\overrightarrow{C{{P}_{C}}}=\frac{\beta _{A}^{P}\overrightarrow{CA}+\beta _{B}^{P}\overrightarrow{CB}+\beta _{D}^{P}\overrightarrow{CD}}{\beta _{A}^{P}+\beta _{B}^{P}+\beta _{D}^{P}},\]	
	\[\overrightarrow{D{{P}_{D}}}=\frac{\beta _{A}^{P}\overrightarrow{DA}+\beta _{B}^{P}\overrightarrow{DB}+\beta _{C}^{P}\overrightarrow{DC}}{\beta _{A}^{P}+\beta _{B}^{P}+\beta _{C}^{P}}.\]	
	
	Where $\beta _{A}^{P}$, $\beta _{B}^{P}$, $\beta _{C}^{P}$, $\beta _{D}^{P}$ are the frame components of $A$, $B$, $C$, $D$ corresponding to point $P$, respectively.	
\end{theorem}

\begin{proof}
	From theorem \ref{thm:Thm18.2.3} and theorem \ref{thm:Thm18.2.2}, the following results can be obtained:	
	\[\begin{aligned}
		\overrightarrow{A{{P}_{A}}}&=\frac{\overrightarrow{AP}}{1-\beta _{A}^{P}}=\frac{\beta _{B}^{P}\overrightarrow{AB}+\beta _{C}^{P}\overrightarrow{AC}+\beta _{D}^{P}\overrightarrow{AD}}{1-\beta _{A}^{P}} \\ 
		& =\frac{\beta _{B}^{P}\overrightarrow{AB}+\beta _{C}^{P}\overrightarrow{AC}+\beta _{D}^{P}\overrightarrow{AD}}{\beta _{B}^{P}+\beta _{C}^{P}+\beta _{D}^{P}}.  
	\end{aligned}\]	
	
	The other formulas can be obtained similarly.
\end{proof}
\hfill $\square$\par

\subsection{Representation of vector from origin to intersecting center of a tetrahedron using frame of circumcenter}\label{Subsec18.2.3}
In this section, the origin $O$ is placed at the circumcenter $Q$ of the tetrahedron $ABCD$ so that the magnitudes of the four frame vectors are equal, that is, the circumscribed sphere radius of the tetrahedron $ABCD$, which will bring convenience to the analysis of some problems. Such a frame is called the frame of circumcenter of the tetrahedron, which is called frame of circumcenter for short, and is denoted as $\left( Q;A,B,C,D \right)$.

\begin{theorem}{Representation theorem of VVIC-T using frame of circumcenter, Daiyuan Zhang}{Thm18.2.5}\label{Thm18.2.5} 
	Given the tetrahedron $ABCD$, the intersecting center of the tetrahedron is $P\in {{\pi }_{ABCD}}$, then the vector $\overrightarrow{QP}$ can be expressed by the uniquely linear combination of the frame of circumcenter $\left( Q;A,B,C,D \right)$, i.e.
	\[\overrightarrow{QP}=\beta _{A}^{P}\overrightarrow{QA}+\beta _{B}^{P}\overrightarrow{QB}+\beta _{C}^{P}\overrightarrow{QC}+\beta _{D}^{P}\overrightarrow{QD},\]	
	\[\beta _{A}^{P}+\beta _{B}^{P}+\beta _{C}^{P}+\beta _{D}^{P}=1.\]	
	
	Where $\beta _{A}^{P}$, $\beta _{B}^{P}$, $\beta _{C}^{P}$, $\beta _{D}^{P}$ are the frame components of $A$, $B$, $C$, $D$ corresponding to point $P$, respectively.
\end{theorem}

\begin{proof}
	This theorem can be obtained by replacing the point $O$ in theorem \ref{thm:Thm18.2.1} with the point $Q$.
\end{proof}
\hfill $\square$\par

\section{Theorem of vector of two intersecting centers of a tetrahedron}\label{Sec18.3}
As we know from the previous discussion, the initial point of the vector $\overrightarrow{OP}$ obtained in theorem \ref{thm:Thm18.2.1} is the origin $O$ of the tetrahedral frame $\left( O;A,B,C,D \right)$, the question now is how to use the frame vectors of $\overrightarrow{OA}$, $\overrightarrow{OB}$, $\overrightarrow{OC}$ and $\overrightarrow{OD}$ of the tetrahedral frame $\left( O;A,B,C,D \right)$ to represent the vector $\overrightarrow{{{P}_{1}}{{P}_{2}}}$ linearly (where points ${{P}_{1}}$ and ${{P}_{2}}$ does not necessarily coincide with the origin $O$)? This is the theorem to be studied below.

\begin{theorem}{Vector of two intersecting centers of a tetrahedron, Daiyuan Zhang}{Thm18.3.1}\label{Thm18.3.1} 
	Suppose that given a tetrahedron $ABCD$, the point $O$ is an arbitrary point, and ${{P}_{1}}\in {{\pi }_{ABCD}}$, ${{P}_{2}}\in {{\pi }_{ABCD}}$, then the vector $\overrightarrow{{{P}_{1}}{{P}_{2}}}$ can be expressed by the uniquely linear combination of the tetrahedral frame $\left( O;A,B,C,D \right)$ in the following:	
	\[\overrightarrow{{{P}_{1}}{{P}_{2}}}=\beta _{A}^{{{P}_{1}}{{P}_{2}}}\overrightarrow{OA}+\beta _{B}^{{{P}_{1}}{{P}_{2}}}\overrightarrow{OB}+\beta _{C}^{{{P}_{1}}{{P}_{2}}}\overrightarrow{OC}+\beta _{D}^{{{P}_{1}}{{P}_{2}}}\overrightarrow{OD},\]	
	\[\beta _{A}^{{{P}_{1}}{{P}_{2}}}+\beta _{B}^{{{P}_{1}}{{P}_{2}}}+\beta _{C}^{{{P}_{1}}{{P}_{2}}}+\beta _{D}^{{{P}_{1}}{{P}_{2}}}=0.\]  	
	Where
	\[\beta _{A}^{{{P}_{1}}{{P}_{2}}}=\beta _{A}^{{{P}_{2}}}-\beta _{A}^{{{P}_{1}}},\]	
	\[\beta _{B}^{{{P}_{1}}{{P}_{2}}}=\beta _{B}^{{{P}_{2}}}-\beta _{B}^{{{P}_{1}}},\]	
	\[\beta _{C}^{{{P}_{1}}{{P}_{2}}}=\beta _{C}^{{{P}_{2}}}-\beta _{C}^{{{P}_{1}}},\]	
	\[\beta _{D}^{{{P}_{1}}{{P}_{2}}}=\beta _{D}^{{{P}_{2}}}-\beta _{D}^{{{P}_{1}}}.\]	
	And the $\beta _{A}^{{{P}_{1}}}$, $\beta _{B}^{{{P}_{1}}}$, $\beta _{C}^{{{P}_{1}}}$, $\beta _{D}^{{{P}_{1}}}$ are the frame components of $A$, $B$, $C$, $D$ at the intersecting center ${{P}_{1}}$, respectively; the $\beta _{A}^{{{P}_{2}}}$, $\beta _{B}^{{{P}_{2}}}$, $\beta _{C}^{{{P}_{2}}}$, $\beta _{D}^{{{P}_{2}}}$ are the frame components of $A$, $B$, $C$, $D$ at the intersecting center ${{P}_{2}}$, respectively. 
\end{theorem}

\begin{proof}
	Obviously, for any point $O$, according to the basic operation of vector, the following result is obtained:
	\[\overrightarrow{{{P}_{1}}{{P}_{2}}}=\overrightarrow{O{{P}_{2}}}-\overrightarrow{O{{P}_{1}}}.\]
	
	According to theorem \ref{thm:Thm18.2.1}, there exists a unique set of $\beta _{A}^{{{P}_{1}}}$, $\beta _{B}^{{{P}_{1}}}$, $\beta _{C}^{{{P}_{1}}}$, $\beta _{D}^{{{P}_{1}}}$ such that		
	\[\overrightarrow{O{{P}_{1}}}=\beta _{A}^{{{P}_{1}}}\overrightarrow{OA}+\beta _{B}^{{{P}_{1}}}\overrightarrow{OB}+\beta _{C}^{{{P}_{1}}}\overrightarrow{OC}+\beta _{D}^{{{P}_{1}}}\overrightarrow{OD},\]	
	\[\beta _{A}^{{{P}_{1}}}+\beta _{B}^{{{P}_{1}}}+\beta _{C}^{{{P}_{1}}}+\beta _{D}^{{{P}_{1}}}=1.\]	
	
	There exists a unique set of $\beta _{A}^{{{P}_{2}}}$, $\beta _{B}^{{{P}_{2}}}$, $\beta _{C}^{{{P}_{2}}}$, $\beta _{D}^{{{P}_{2}}}$ such that	
	\[\overrightarrow{O{{P}_{2}}}=\beta _{A}^{{{P}_{2}}}\overrightarrow{OA}+\beta _{B}^{{{P}_{2}}}\overrightarrow{OB}+\beta _{C}^{{{P}_{2}}}\overrightarrow{OC}+\beta _{D}^{{{P}_{2}}}\overrightarrow{OD},\]	
	\[\beta _{A}^{{{P}_{2}}}+\beta _{B}^{{{P}_{2}}}+\beta _{C}^{{{P}_{2}}}+\beta _{D}^{{{P}_{2}}}=1.\]	
	
	Therefore,
	\[\begin{aligned}
		\overrightarrow{{{P}_{1}}{{P}_{2}}}&=\overrightarrow{O{{P}_{2}}}-\overrightarrow{O{{P}_{1}}}=\beta _{A}^{{{P}_{2}}}\overrightarrow{OA}+\beta _{B}^{{{P}_{2}}}\overrightarrow{OB}+\beta _{C}^{{{P}_{2}}}\overrightarrow{OC}+\beta _{D}^{{{P}_{2}}}\overrightarrow{OD} \\ 
		& -\left( \beta _{A}^{{{P}_{1}}}\overrightarrow{OA}+\beta _{B}^{{{P}_{1}}}\overrightarrow{OB}+\beta _{C}^{{{P}_{1}}}\overrightarrow{OC}+\beta _{D}^{{{P}_{1}}}\overrightarrow{OD} \right),  
	\end{aligned}\]	
	i.e.
	\[\begin{aligned}
		\overrightarrow{{{P}_{1}}{{P}_{2}}}&=\left( \beta _{A}^{{{P}_{2}}}-\beta _{A}^{{{P}_{1}}} \right)\overrightarrow{OA}+\left( \beta _{B}^{{{P}_{2}}}-\beta _{B}^{{{P}_{1}}} \right)\overrightarrow{OB} \\ 
		& +\left( \beta _{C}^{{{P}_{2}}}-\beta _{C}^{{{P}_{1}}} \right)\overrightarrow{OC}+\left( \beta _{D}^{{{P}_{2}}}-\beta _{D}^{{{P}_{1}}} \right)\overrightarrow{OD} \\ 
		& =\beta _{A}^{{{P}_{1}}{{P}_{2}}}\overrightarrow{OA}+\beta _{B}^{{{P}_{1}}{{P}_{2}}}\overrightarrow{OB}+\beta _{C}^{{{P}_{1}}{{P}_{2}}}\overrightarrow{OC}+\beta _{D}^{{{P}_{1}}{{P}_{2}}}\overrightarrow{OD}.  
	\end{aligned}\]
	and
	\[\begin{aligned}
		& \beta _{A}^{{{P}_{1}}{{P}_{2}}}+\beta _{B}^{{{P}_{1}}{{P}_{2}}}+\beta _{C}^{{{P}_{1}}{{P}_{2}}}+\beta _{D}^{{{P}_{1}}{{P}_{2}}} \\ 
		& =\left( \beta _{A}^{{{P}_{2}}}-\beta _{A}^{{{P}_{1}}} \right)+\left( \beta _{B}^{{{P}_{2}}}-\beta _{B}^{{{P}_{1}}} \right)+\left( \beta _{C}^{{{P}_{2}}}-\beta _{C}^{{{P}_{1}}} \right)+\left( \beta _{D}^{{{P}_{2}}}-\beta _{D}^{{{P}_{1}}} \right) \\ 
		& =\left( \beta _{A}^{{{P}_{2}}}+\beta _{B}^{{{P}_{2}}}+\beta _{C}^{{{P}_{2}}}+\beta _{D}^{{{P}_{2}}} \right)-\left( \beta _{A}^{{{P}_{1}}}+\beta _{B}^{{{P}_{1}}}+\beta _{C}^{{{P}_{1}}}+\beta _{D}^{{{P}_{1}}} \right)=0.  
	\end{aligned}\]	
	
	Since $\beta _{A}^{{{P}_{1}}}$, $\beta _{B}^{{{P}_{1}}}$, $\beta _{C}^{{{P}_{1}}}$, $\beta _{D}^{{{P}_{1}}}$ and $\beta _{A}^{{{P}_{2}}}$, $\beta _{B}^{{{P}_{2}}}$, $\beta _{C}^{{{P}_{2}}}$, $\beta _{D}^{{{P}_{2}}}$ are unique, $\beta _{A}^{{{P}_{1}}{{P}_{2}}}$, $\beta _{B}^{{{P}_{1}}{{P}_{2}}}$, $\beta _{C}^{{{P}_{1}}{{P}_{2}}}$, $\beta _{D}^{{{P}_{1}}{{P}_{2}}}$ are also unique.
\end{proof}
\hfill $\square$\par

The frame components of vector $\overrightarrow{{{P}_{1}}{{P}_{2}}}$ is related to the IRs of the two intersecting centers of the tetrahedron. Such a vector of intersecting center of the tetrahedron (abbreviated as VIC-T) is called the vector of two intersecting centers of the tetrahedron (abbreviated as VTICs-T).


The importance of theorem \ref{thm:Thm18.3.1} is that the vector $\overrightarrow{{{P}_{1}}{{P}_{2}}}$ can be expressed by the uniquely linear combination of the frame vectors $\overrightarrow{OA}$, $\overrightarrow{OB}$, $\overrightarrow{OC}$ and $\overrightarrow{OD}$ of the tetrahedral frame $\left( O;A,B,C,D \right)$, and its frame components (coefficients) are only related to the IRs of the two intersecting centers ${{P}_{1}}$ and ${{P}_{2}}$ of the tetrahedron, and has nothing to do with the position of the origin $O$ of the frame $\left( O;A,B,C \right)$.

\subsection{Representation theorem using frame of edge for the vector of two intersecting centers of a tetrahedron}\label{Subsec18.3.1}

\begin{theorem}{Representation theorem using edge frame for VTICs-T, Daiyuan Zhang}{Thm18.3.2}\label{Thm18.3.2} 
	Suppose that given a tetrahedron $ABCD$, and ${{P}_{1}}\in {{\pi }_{ABCD}}$, ${{P}_{2}}\in {{\pi }_{ABCD}}$, then the vector $\overrightarrow{{{P}_{1}}{{P}_{2}}}$ can be expressed by the uniquely linear combination of the frame of edge $\left( O;A,B,C,D \right)$ in the following:		
	\[\overrightarrow{{{P}_{1}}{{P}_{2}}}=\beta _{B}^{{{P}_{1}}{{P}_{2}}}\overrightarrow{AB}+\beta _{C}^{{{P}_{1}}{{P}_{2}}}\overrightarrow{AC}+\beta _{D}^{{{P}_{1}}{{P}_{2}}}\overrightarrow{AD},\]	
	\[\overrightarrow{{{P}_{1}}{{P}_{2}}}=\beta _{A}^{{{P}_{1}}{{P}_{2}}}\overrightarrow{BA}+\beta _{C}^{{{P}_{1}}{{P}_{2}}}\overrightarrow{BC}+\beta _{D}^{{{P}_{1}}{{P}_{2}}}\overrightarrow{BD},\]	
	\[\overrightarrow{{{P}_{1}}{{P}_{2}}}=\beta _{A}^{{{P}_{1}}{{P}_{2}}}\overrightarrow{CA}+\beta _{B}^{{{P}_{1}}{{P}_{2}}}\overrightarrow{CB}+\beta _{D}^{{{P}_{1}}{{P}_{2}}}\overrightarrow{CD},\]	
	\[\overrightarrow{{{P}_{1}}{{P}_{2}}}=\beta _{A}^{{{P}_{1}}{{P}_{2}}}\overrightarrow{DA}+\beta _{B}^{{{P}_{1}}{{P}_{2}}}\overrightarrow{DB}+\beta _{C}^{{{P}_{1}}{{P}_{2}}}\overrightarrow{DC}.\]	
	
	Or
	\[\overrightarrow{{{P}_{1}}{{P}_{2}}}=\frac{1}{4}\left( \begin{aligned}
		& \left( \beta _{B}^{{{P}_{1}}{{P}_{2}}}-\beta _{A}^{{{P}_{1}}{{P}_{2}}} \right)\overrightarrow{AB}+\left( \beta _{C}^{{{P}_{1}}{{P}_{2}}}-\beta _{A}^{{{P}_{1}}{{P}_{2}}} \right)\overrightarrow{AC}+\left( \beta _{D}^{{{P}_{1}}{{P}_{2}}}-\beta _{A}^{{{P}_{1}}{{P}_{2}}} \right)\overrightarrow{AD} \\ 
		& +\left( \beta _{C}^{{{P}_{1}}{{P}_{2}}}-\beta _{B}^{{{P}_{1}}{{P}_{2}}} \right)\overrightarrow{BC}+\left( \beta _{D}^{{{P}_{1}}{{P}_{2}}}-\beta _{B}^{{{P}_{1}}{{P}_{2}}} \right)\overrightarrow{BD}+\left( \beta _{D}^{{{P}_{1}}{{P}_{2}}}-\beta _{C}^{{{P}_{1}}{{P}_{2}}} \right)\overrightarrow{CD} \\ 
	\end{aligned} \right).\]  
	And		
	\[\beta _{A}^{{{P}_{1}}{{P}_{2}}}+\beta _{B}^{{{P}_{1}}{{P}_{2}}}+\beta _{C}^{{{P}_{1}}{{P}_{2}}}+\beta _{D}^{{{P}_{1}}{{P}_{2}}}=0,\]	
	where
	\[\beta _{A}^{{{P}_{1}}{{P}_{2}}}=\beta _{A}^{{{P}_{2}}}-\beta _{A}^{{{P}_{1}}},\]	
	\[\beta _{B}^{{{P}_{1}}{{P}_{2}}}=\beta _{B}^{{{P}_{2}}}-\beta _{B}^{{{P}_{1}}},\]	
	\[\beta _{C}^{{{P}_{1}}{{P}_{2}}}=\beta _{C}^{{{P}_{2}}}-\beta _{C}^{{{P}_{1}}},\]	
	\[\beta _{D}^{{{P}_{1}}{{P}_{2}}}=\beta _{D}^{{{P}_{2}}}-\beta _{D}^{{{P}_{1}}}.\]	
	And the $\beta _{A}^{{{P}_{1}}}$, $\beta _{B}^{{{P}_{1}}}$, $\beta _{C}^{{{P}_{1}}}$, $\beta _{D}^{{{P}_{1}}}$ are the frame components of $A$, $B$, $C$, $D$ at the intersecting center ${{P}_{1}}$, respectively; the $\beta _{A}^{{{P}_{2}}}$, $\beta _{B}^{{{P}_{2}}}$, $\beta _{C}^{{{P}_{2}}}$, $\beta _{D}^{{{P}_{2}}}$ are the frame components of $A$, $B$, $C$, $D$ at the intersecting center ${{P}_{2}}$, respectively. 	
\end{theorem}
\begin{proof}
	According to theorem \ref{thm:Thm18.3.1}, the results can be obtained by coincidence of point $O$ with the points $A$, $B$, $C$, $D$, respectively.
\end{proof}
\hfill $\square$\par

The The above theorem shows that: the vector $\overrightarrow{{{P}_{1}}{{P}_{2}}}$ can be expressed by the uniquely linear combination of the frame of edge, and its frame components (coefficients) are only related to the IRs of the two intersecting centers, and has nothing to do with the position of the origin of the frame $\left( O;A,B,C \right)$.

The above theorem can be conveniently used to calculate the magnitudes of the vector $\overrightarrow{{{P}_{1}}{{P}_{2}}}$ and $\overrightarrow{{{P}_{1}}{{P}_{2}}}$.

\subsection{Representation theorem using frame of circumcenter for the vector of two intersecting centers of a tetrahedron}\label{Subsec18.3.2}

\begin{theorem}{Representation using frame of circumcenter for VTICs-T, Daiyuan Zhang}{Thm18.3.3}\label{Thm18.3.3} 
	Suppose that given a tetrahedron $ABCD$, and ${{P}_{1}}\in {{\pi }_{ABCD}}$, ${{P}_{2}}\in {{\pi }_{ABCD}}$, then the vector $\overrightarrow{{{P}_{1}}{{P}_{2}}}$ can be expressed by the uniquely linear combination of the frame of circumcenter $\left(Q;A,B,C,D \right)$ in the following:
	\[\overrightarrow{{{P}_{1}}{{P}_{2}}}=\beta _{A}^{{{P}_{1}}{{P}_{2}}}\overrightarrow{QA}+\beta _{B}^{{{P}_{1}}{{P}_{2}}}\overrightarrow{QB}+\beta _{C}^{{{P}_{1}}{{P}_{2}}}\overrightarrow{QC}+\beta _{D}^{{{P}_{1}}{{P}_{2}}}\overrightarrow{QD}.\]	
	Where
	\[\beta _{A}^{{{P}_{1}}{{P}_{2}}}=\beta _{A}^{{{P}_{2}}}-\beta _{A}^{{{P}_{1}}},\]	
	\[\beta _{B}^{{{P}_{1}}{{P}_{2}}}=\beta _{B}^{{{P}_{2}}}-\beta _{B}^{{{P}_{1}}},\]	
	\[\beta _{C}^{{{P}_{1}}{{P}_{2}}}=\beta _{C}^{{{P}_{2}}}-\beta _{C}^{{{P}_{1}}},\]	
	\[\beta _{D}^{{{P}_{1}}{{P}_{2}}}=\beta _{D}^{{{P}_{2}}}-\beta _{D}^{{{P}_{1}}}.\]	
	And the $\beta _{A}^{{{P}_{1}}}$, $\beta _{B}^{{{P}_{1}}}$, $\beta _{C}^{{{P}_{1}}}$, $\beta _{D}^{{{P}_{1}}}$ are the frame components of $A$, $B$, $C$, $D$ at the intersecting center ${{P}_{1}}$, respectively; the $\beta _{A}^{{{P}_{2}}}$, $\beta _{B}^{{{P}_{2}}}$, $\beta _{C}^{{{P}_{2}}}$, $\beta _{D}^{{{P}_{2}}}$ are the frame components of $A$, $B$, $C$, $D$ at the intersecting center ${{P}_{2}}$, respectively. 	
\end{theorem}

\begin{proof}
	The origin of the tetrahedral frame is chosen at the circumcenter of the tetrahedron, and the result is obtained directly according to theorem \ref{thm:Thm18.3.1}.
\end{proof}
\hfill $\square$\par

\section{Frame equation of intersecting center of a tetrahedron}\label{Sec18.4}

The tetrahedral frame of intersecting center of a tetrahedron satisfies some equations.

\begin{theorem}{Frame equation of IC-T, Daiyuan Zhang}{Thm18.4.1}\label{Thm18.4.1} 
	Given the tetrahedron $ABCD$ and the point $P$ as the intersecting center of the tetrahedron, then
	\[\beta _{A}^{P}\overrightarrow{PA}+\beta _{B}^{P}\overrightarrow{PB}+\beta _{C}^{P}\overrightarrow{PC}+\beta _{D}^{P}\overrightarrow{PD}=\overrightarrow{0},\]	
	\[\beta _{A}^{P}+\beta _{B}^{P}+\beta _{C}^{P}+\beta _{D}^{P}=1.\]	
\end{theorem}

\begin{proof}
	For the tetrahedral frame $\left( O;A,B,C \right)$, according to theorem \ref{thm:Thm18.2.1}, let the point $O$ and $P$ coincide to get the result directly.
\end{proof}
\hfill $\square$\par

The frame equation for the vector of two intersecting centers is given below.
\begin{theorem}{Frame equation of VTICs-T, Daiyuan Zhang}{Thm18.4.2}\label{Thm18.4.2} 
	Given the tetrahedron $ABCD$ and the points ${{P}_{1}}$ and ${{P}_{2}}$ as the intersecting centers of the tetrahedron, then	
	\[\begin{aligned}
		& \beta _{A}^{{{P}_{1}}}\overrightarrow{{{P}_{2}}A}+\beta _{B}^{{{P}_{1}}}\overrightarrow{{{P}_{2}}B}+\beta _{C}^{{{P}_{1}}}\overrightarrow{{{P}_{2}}C}+\beta _{D}^{{{P}_{1}}}\overrightarrow{{{P}_{2}}D} \\ 
		& +\beta _{A}^{{{P}_{2}}}\overrightarrow{{{P}_{1}}A}+\beta _{B}^{{{P}_{2}}}\overrightarrow{{{P}_{1}}B}+\beta _{C}^{{{P}_{2}}}\overrightarrow{{{P}_{1}}C}+\beta _{D}^{{{P}_{2}}}\overrightarrow{{{P}_{1}}D}=\overrightarrow{0}, \\ 
	\end{aligned}\]	
	\[\beta _{A}^{{{P}_{1}}}+\beta _{B}^{{{P}_{1}}}+\beta _{C}^{{{P}_{1}}}+\beta _{D}^{{{P}_{1}}}=1,\]	
	\[\beta _{A}^{{{P}_{2}}}+\beta _{B}^{{{P}_{2}}}+\beta _{C}^{{{P}_{2}}}+\beta _{D}^{{{P}_{2}}}=1.\]	
\end{theorem}

\begin{proof}
	Because the ICs have equal status, it is obtained according to theorem \ref{thm:Thm18.2.1}:	
	\[\overrightarrow{{{P}_{1}}{{P}_{2}}}=\beta _{A}^{{{P}_{2}}}\overrightarrow{{{P}_{1}}A}+\beta _{B}^{{{P}_{2}}}\overrightarrow{{{P}_{1}}B}+\beta _{C}^{{{P}_{2}}}\overrightarrow{{{P}_{1}}C}+\beta _{D}^{{{P}_{2}}}\overrightarrow{{{P}_{1}}D},\]	
	\[\beta _{A}^{{{P}_{2}}}+\beta _{B}^{{{P}_{2}}}+\beta _{C}^{{{P}_{2}}}+\beta _{D}^{{{P}_{2}}}=1,\]	
	\[\overrightarrow{{{P}_{2}}{{P}_{1}}}=\beta _{A}^{{{P}_{1}}}\overrightarrow{{{P}_{2}}A}+\beta _{B}^{{{P}_{1}}}\overrightarrow{{{P}_{2}}B}+\beta _{C}^{{{P}_{1}}}\overrightarrow{{{P}_{2}}C}+\beta _{D}^{{{P}_{1}}}\overrightarrow{{{P}_{2}}D},\]	
	\[\beta _{A}^{{{P}_{1}}}+\beta _{B}^{{{P}_{1}}}+\beta _{C}^{{{P}_{1}}}+\beta _{D}^{{{P}_{1}}}=1.\]	
	
	Since $\overrightarrow{{{P}_{1}}{{P}_{2}}}+\overrightarrow{{{P}_{2}}{{P}_{1}}}=\overrightarrow{0}$, therefore
	\[\begin{aligned}
		& \beta _{A}^{{{P}_{1}}}\overrightarrow{{{P}_{2}}A}+\beta _{B}^{{{P}_{1}}}\overrightarrow{{{P}_{2}}B}+\beta _{C}^{{{P}_{1}}}\overrightarrow{{{P}_{2}}C}+\beta _{D}^{{{P}_{1}}}\overrightarrow{{{P}_{2}}D} \\ 
		& +\beta _{A}^{{{P}_{2}}}\overrightarrow{{{P}_{1}}A}+\beta _{B}^{{{P}_{2}}}\overrightarrow{{{P}_{1}}B}+\beta _{C}^{{{P}_{2}}}\overrightarrow{{{P}_{1}}C}+\beta _{D}^{{{P}_{2}}}\overrightarrow{{{P}_{1}}D}=\overrightarrow{0}. \\ 
	\end{aligned}\]	
\end{proof}
\hfill $\square$\par

%
The above equation gives the vector equation should be satisfied by the frame components of two ICs-T and the corresponding frames.

The frame equation of the vector from origin to intersecting center of the tetrahedron and the frame equation for the vector of two intersecting centers are also referred to as frame equation of intersecting center of the tetrahedron.


\chapter{New theorems of space geometry}\label{Ch19}
\thispagestyle{empty}

In this chapter, two new theorems of space geometry are proved according to theorem \ref{thm:Thm18.2.1}.
\section{Theorem of integral ratio of vertex to intersecting foot}\label{Sec19.1}

\begin{theorem}{Theorem of integral ratio of vertex to intersecting foot, Daiyuan Zhang}{Thm19.1.1}\label{Thm19.1.1} 
	For the tetrahedron $ABCD$, the following result is obtained on the tetrahedral frame $\left( O;A,B,C,D \right)$: 	
	\[{{\kappa }_{A{{P}_{A}}}}+{{\kappa }_{B{{P}_{B}}}}+{{\kappa }_{C{{P}_{C}}}}+{{\kappa }_{D{{P}_{D}}}}=3.\]	
\end{theorem}

\begin{proof}
	By using theorem \ref{thm:Thm18.2.1}, let the frame components (coefficients) be equal respectively, then
	\begin{equation}\label{Eq19.1.1}
		1-{{\kappa }_{A{{P}_{A}}}}={{\kappa }_{B{{P}_{B}}}}\alpha _{A}^{{{P}_{B}}}={{\kappa }_{C{{P}_{C}}}}\alpha _{A}^{{{P}_{C}}}={{\kappa }_{D{{P}_{D}}}}\alpha _{A}^{{{P}_{D}}},
	\end{equation}
	\begin{equation}\label{Eq19.1.2}
		{{\kappa }_{A{{P}_{A}}}}\alpha _{B}^{{{P}_{A}}}=1-{{\kappa }_{B{{P}_{B}}}}={{\kappa }_{C{{P}_{C}}}}\alpha _{B}^{{{P}_{C}}}={{\kappa }_{D{{P}_{D}}}}\alpha _{B}^{{{P}_{D}}},
	\end{equation}
	\begin{equation}\label{Eq19.1.3}
		{{\kappa }_{A{{P}_{A}}}}\alpha _{C}^{{{P}_{A}}}={{\kappa }_{B{{P}_{B}}}}\alpha _{C}^{{{P}_{B}}}=1-{{\kappa }_{C{{P}_{C}}}}={{\kappa }_{D{{P}_{D}}}}\alpha _{C}^{{{P}_{D}}},
	\end{equation}
	\begin{equation}\label{Eq19.1.4}
		{{\kappa }_{A{{P}_{A}}}}\alpha _{D}^{{{P}_{A}}}={{\kappa }_{B{{P}_{B}}}}\alpha _{D}^{{{P}_{B}}}={{\kappa }_{C{{P}_{C}}}}\alpha _{D}^{{{P}_{C}}}=1-{{\kappa }_{D{{P}_{D}}}}.
	\end{equation}
		
	According to formulas (\ref{Eq19.1.1})-(\ref{Eq19.1.3}), the following formula is obtained by using the properties of fraction and the formula of intersecting ratio of face $\alpha _{A}^{{{P}_{D}}}+\alpha _{B}^{{{P}_{D}}}+\alpha _{C}^{{{P}_{D}}}=1$:
	\[\begin{aligned}
		{{\kappa }_{D{{P}_{D}}}}&=\frac{1-{{\kappa }_{A{{P}_{A}}}}}{\alpha _{A}^{{{P}_{D}}}}=\frac{1-{{\kappa }_{B{{P}_{B}}}}}{\alpha _{B}^{{{P}_{D}}}}=\frac{1-{{\kappa }_{C{{P}_{C}}}}}{\alpha _{C}^{{{P}_{D}}}} =\frac{1-{{\kappa }_{A{{P}_{A}}}}+1-{{\kappa }_{B{{P}_{B}}}}+1-{{\kappa }_{C{{P}_{C}}}}}{\alpha _{A}^{{{P}_{D}}}+\alpha _{B}^{{{P}_{D}}}+\alpha _{C}^{{{P}_{D}}}}\\
		&=3-{{\kappa }_{A{{P}_{A}}}}-{{\kappa }_{B{{P}_{B}}}}-{{\kappa }_{C{{P}_{C}}}}.  
	\end{aligned}\]	
	
	Therefore,
	\[{{\kappa }_{A{{P}_{A}}}}+{{\kappa }_{B{{P}_{B}}}}+{{\kappa }_{C{{P}_{C}}}}+{{\kappa }_{D{{P}_{D}}}}=3.\]	
\end{proof}
\hfill $\square$\par

\section{Theorem of fractional ratio of vertex to intersecting foot of a tetrahedron}\label{Sec19.2}

\begin{theorem}{Fractional ratio of vertex to intersecting foot of a tetrahedron, Daiyuan Zhang}{Thm19.2.1}\label{Thm19.2.1} 
	For the tetrahedron $ABCD$, the following result is obtained on the tetrahedral frame $\left( O;A,B,C,D \right)$: 	
	\[\frac{1}{1+{{\lambda }_{A{{P}_{A}}}}}+\frac{1}{1+{{\lambda }_{B{{P}_{B}}}}}+\frac{1}{1+{{\lambda }_{C{{P}_{C}}}}}+\frac{1}{1+{{\lambda }_{D{{P}_{D}}}}}=1.\]	
\end{theorem}

\begin{proof}
	Using  formula (\ref{Eq2.4.3}), the results of theorem \ref{thm:Thm19.1.1} can also be written as follows:
	\[\frac{{{\lambda }_{A{{P}_{A}}}}}{1+{{\lambda }_{A{{P}_{A}}}}}+\frac{{{\lambda }_{B{{P}_{B}}}}}{1+{{\lambda }_{B{{P}_{B}}}}}+\frac{{{\lambda }_{C{{P}_{C}}}}}{1+{{\lambda }_{C{{P}_{C}}}}}+\frac{{{\lambda }_{D{{P}_{D}}}}}{1+{{\lambda }_{D{{P}_{D}}}}}=3,\]	
	i.e.
	\[1-\frac{1}{1+{{\lambda }_{A{{P}_{A}}}}}+1-\frac{1}{1+{{\lambda }_{B{{P}_{B}}}}}+1-\frac{1}{1+{{\lambda }_{C{{P}_{C}}}}}+1-\frac{1}{1+{{\lambda }_{D{{P}_{D}}}}}=3,\]	
	i.e.
	\[\frac{1}{1+{{\lambda }_{A{{P}_{A}}}}}+\frac{1}{1+{{\lambda }_{B{{P}_{B}}}}}+\frac{1}{1+{{\lambda }_{C{{P}_{C}}}}}+\frac{1}{1+{{\lambda }_{D{{P}_{D}}}}}=1.\]	
\end{proof}
\hfill $\square$\par


\chapter{Frame components on tetrahedral frame}\label{Ch20}
\thispagestyle{empty}


The frame component is very important in Intercenter Geometry. Compared with the triangular frame components, the tetrahedral frame components are more complex. This chapter studies the relationship between the frame components of tetrahedron and the frame components of triangle, the relationship between the frame components of triangle and the frame components of line segment, and the relationship between the frame components and some geometric quantities such as integral ratio and fractional ratio. These relations are of great value in the study of Intercenter Geometry.

\section{Relationship between integral ratio of vertex to intersecting foot and frame components of faces}\label{Sec20.1}

\subsection{Calculating integral ratio based on IR-F}\label{Subsec20.1.1}
\begin{theorem}{Calculating integral ratio (1) based on IR-F, Daiyuan Zhang}{Thm20.1.1}\label{Thm20.1.1} 
	Integral ratios of the vertex to its corresponding intersecting foot of tetrahedron $ABCD$ are
	\begin{flalign*}
		{{\kappa }_{A{{P}_{A}}}}=\frac{1-\alpha _{A}^{{{P}_{B}}}}{1-\alpha _{A}^{{{P}_{B}}}\alpha _{B}^{{{P}_{A}}}},\ {{\kappa }_{B{{P}_{B}}}}=\frac{1-\alpha _{B}^{{{P}_{C}}}}{1-\alpha _{B}^{{{P}_{C}}}\alpha _{C}^{{{P}_{B}}}},
	\end{flalign*}
	\begin{flalign*}
		{{\kappa }_{C{{P}_{C}}}}=\frac{1-\alpha _{C}^{{{P}_{D}}}}{1-\alpha _{C}^{{{P}_{D}}}\alpha _{D}^{{{P}_{C}}}},\ {{\kappa }_{D{{P}_{D}}}}=\frac{1-\alpha _{D}^{{{P}_{A}}}}{1-\alpha _{D}^{{{P}_{A}}}\alpha _{A}^{{{P}_{D}}}}.
	\end{flalign*}
\end{theorem}

\begin{proof}
	From the equation (\ref{Eq19.1.1}) - (\ref{Eq19.1.4}), a system of equations containing only ${{\kappa }_{A{{P}_{A}}}}$ and ${{\kappa }_{B{{P}_{B}}}}$ can be obtained:	
	\begin{equation}\label{Eq20.1.1}
		\left\{ \begin{aligned}
			& 1-{{\kappa }_{A{{P}_{A}}}}={{\kappa }_{B{{P}_{B}}}}\alpha _{A}^{{{P}_{B}}} \\ 
			& {{\kappa }_{A{{P}_{A}}}}\alpha _{B}^{{{P}_{A}}}=1-{{\kappa }_{B{{P}_{B}}}}.\\ 
		\end{aligned} \right.
	\end{equation}
		
	By solving the above system of equations, we get the following results:	
	\[1-{{\kappa }_{A{{P}_{A}}}}=\left( 1-{{\kappa }_{A{{P}_{A}}}}\alpha _{B}^{{{P}_{A}}} \right)\alpha _{A}^{{{P}_{B}}},\]	
	\[\left( 1-\alpha _{B}^{{{P}_{A}}}\alpha _{A}^{{{P}_{B}}} \right){{\kappa }_{A{{P}_{A}}}}=1-\alpha _{A}^{{{P}_{B}}}.\]	
	Therefore
	\begin{equation}\label{Eq20.1.2}
		{{\kappa }_{A{{P}_{A}}}}=\frac{1-\alpha _{A}^{{{P}_{B}}}}{1-\alpha _{A}^{{{P}_{B}}}\alpha _{B}^{{{P}_{A}}}}.
	\end{equation}
		
	By substituting the above results into equation (\ref{Eq20.1.1}), we can get ${{\kappa }_{B{{P}_{B}}}}$. But here we use the following rotation method to solve ${{\kappa }_{B{{P}_{B}}}}$ and so on.
	
	According to the equations (\ref{Eq19.1.1})-(\ref{Eq19.1.4}), we can get the corresponding equations which only contain ${{\kappa }_{B{{P}_{B}}}}$ and ${{\kappa }_{C{{P}_{C}}}}$, ${{\kappa }_{C{{P}_{C}}}}$ and ${{\kappa }_{D{{P}_{D}}}}$, ${{\kappa }_{D{{P}_{D}}}}$ and ${{\kappa }_{A{{P}_{A}}}}$, respectively:	
	\[\left\{ \begin{aligned}
		& 1-{{\kappa }_{B{{P}_{B}}}}={{\kappa }_{C{{P}_{C}}}}\alpha _{B}^{{{P}_{C}}} \\ 
		& {{\kappa }_{B{{P}_{B}}}}\alpha _{C}^{{{P}_{B}}}=1-{{\kappa }_{C{{P}_{C}}}}, \\ 
	\end{aligned} \right.\]	
	\[\left\{ \begin{aligned}
		& 1-{{\kappa }_{C{{P}_{C}}}}={{\kappa }_{D{{P}_{D}}}}\alpha _{C}^{{{P}_{D}}} \\ 
		& {{\kappa }_{C{{P}_{C}}}}\alpha _{D}^{{{P}_{C}}}=1-{{\kappa }_{D{{P}_{D}}}}, \\ 
	\end{aligned} \right.\]	
	\[\left\{ \begin{aligned}
		& 1-{{\kappa }_{D{{P}_{D}}}}={{\kappa }_{A{{P}_{A}}}}\alpha _{D}^{{{P}_{A}}} \\ 
		& {{\kappa }_{D{{P}_{D}}}}\alpha _{A}^{{{P}_{D}}}=1-{{\kappa }_{A{{P}_{A}}}}. \\ 
	\end{aligned} \right.\]	
	
	Obviously, the above equations are obtained by the rotation ($A$→$B$→$C$→$D$→$A$) of the equations (\ref{Eq20.1.1}), so their solutions can also be obtained by the rotation of the solutions (\ref{Eq20.1.2}) of the equations (\ref{Eq20.1.1}).
\end{proof}
\hfill $\square$\par

The above theorem shows that the IRVIF-T of an $P$ (IC-T) can be expressed by the frame components of the corresponding intersecting center of face (abbreviated as IC-F).

The expression of the above theorem is asymmetric, and the symmetric expression is given below.

\begin{theorem}{Symmetric form 1 for calculating integral ratio based on IR-F, Daiyuan Zhang}{Thm20.1.2}\label{Thm20.1.2} 
	Integral ratios of the vertex to its corresponding intersecting foot of tetrahedron $ABCD$ are	
	\[{{\kappa }_{A{{P}_{A}}}}=\frac{1-\alpha _{A}^{{{P}_{B}}}+\alpha _{A}^{{{P}_{B}}}\alpha _{B}^{{{P}_{C}}}-\alpha _{A}^{{{P}_{B}}}\alpha _{B}^{{{P}_{C}}}\alpha _{C}^{{{P}_{D}}}}{1-\alpha _{A}^{{{P}_{B}}}\alpha _{B}^{{{P}_{C}}}\alpha _{C}^{{{P}_{D}}}\alpha _{D}^{{{P}_{A}}}},\]	
	\[{{\kappa }_{B{{P}_{B}}}}=\frac{1-\alpha _{B}^{{{P}_{C}}}+\alpha _{B}^{{{P}_{C}}}\alpha _{C}^{{{P}_{D}}}-\alpha _{B}^{{{P}_{C}}}\alpha _{C}^{{{P}_{D}}}\alpha _{D}^{{{P}_{A}}}}{1-\alpha _{B}^{{{P}_{C}}}\alpha _{C}^{{{P}_{D}}}\alpha _{D}^{{{P}_{A}}}\alpha _{A}^{{{P}_{B}}}},\]	
	\[{{\kappa }_{C{{P}_{C}}}}=\frac{1-\alpha _{C}^{{{P}_{D}}}+\alpha _{C}^{{{P}_{D}}}\alpha _{D}^{{{P}_{A}}}-\alpha _{C}^{{{P}_{D}}}\alpha _{D}^{{{P}_{A}}}\alpha _{A}^{{{P}_{B}}}}{1-\alpha _{C}^{{{P}_{D}}}\alpha _{D}^{{{P}_{A}}}\alpha _{A}^{{{P}_{B}}}\alpha _{B}^{{{P}_{C}}}},\]	
	\[{{\kappa }_{D{{P}_{D}}}}=\frac{1-\alpha _{D}^{{{P}_{A}}}+\alpha _{D}^{{{P}_{A}}}\alpha _{A}^{{{P}_{B}}}-\alpha _{D}^{{{P}_{A}}}\alpha _{A}^{{{P}_{B}}}\alpha _{B}^{{{P}_{C}}}}{1-\alpha _{D}^{{{P}_{A}}}\alpha _{A}^{{{P}_{B}}}\alpha _{B}^{{{P}_{C}}}\alpha _{C}^{{{P}_{D}}}}.\]
\end{theorem}

\begin{proof}
	According to equations (\ref{Eq19.1.1})-(\ref{Eq19.1.4}), one of the linear equations can be listed:
	\[\left\{ \begin{aligned}
		& 1-{{\kappa }_{A{{P}_{A}}}}={{\kappa }_{B{{P}_{B}}}}\alpha _{A}^{{{P}_{B}}} \\ 
		& 1-{{\kappa }_{B{{P}_{B}}}}={{\kappa }_{C{{P}_{C}}}}\alpha _{B}^{{{P}_{C}}} \\ 
		& 1-{{\kappa }_{C{{P}_{C}}}}={{\kappa }_{D{{P}_{D}}}}\alpha _{C}^{{{P}_{D}}} \\ 
		& {{\kappa }_{A{{P}_{A}}}}\alpha _{D}^{{{P}_{A}}}=1-{{\kappa }_{D{{P}_{D}}}}. \\ 
	\end{aligned} \right.\]
	i.e.
	\[\left\{ \begin{aligned}
		& {{\kappa }_{A{{P}_{A}}}}+\alpha _{A}^{{{P}_{B}}}{{\kappa }_{B{{P}_{B}}}}=1 \\ 
		& {{\kappa }_{B{{P}_{B}}}}+\alpha _{B}^{{{P}_{C}}}{{\kappa }_{C{{P}_{C}}}}=1 \\ 
		& {{\kappa }_{C{{P}_{C}}}}+\alpha _{C}^{{{P}_{D}}}{{\kappa }_{D{{P}_{D}}}}=1 \\ 
		& \alpha _{D}^{{{P}_{A}}}{{\kappa }_{A{{P}_{A}}}}+{{\kappa }_{D{{P}_{D}}}}=1. \\ 
	\end{aligned} \right.\]
	
	The above equations can be written as matrix form and the following result is obtained:
	\[\left( \begin{matrix}
		1 & \alpha _{A}^{{{P}_{B}}} & {} & {}  \\
		{} & 1 & \alpha _{B}^{{{P}_{C}}} & {}  \\
		{} & {} & 1 & \alpha _{C}^{{{P}_{D}}}  \\
		\alpha _{D}^{{{P}_{A}}} & {} & {} & 1  \\
	\end{matrix} \right)\left( \begin{matrix}
		{{\kappa }_{A{{P}_{A}}}}  \\
		{{\kappa }_{B{{P}_{B}}}}  \\
		{{\kappa }_{C{{P}_{C}}}}  \\
		{{\kappa }_{D{{P}_{D}}}}  \\
	\end{matrix} \right)=\left( \begin{matrix}
		1  \\
		1  \\
		1  \\
		1  \\
	\end{matrix} \right).\]	
	
	The following equivalent equations are obtained by elementary row transformation:
	\[\left( \begin{matrix}
		1 & \alpha _{A}^{{{P}_{B}}} & {} & {}  \\
		{} & 1 & \alpha _{B}^{{{P}_{C}}} & {}  \\
		{} & {} & 1 & \alpha _{C}^{{{P}_{D}}}  \\
		{} & -\alpha _{D}^{{{P}_{A}}}\alpha _{A}^{{{P}_{B}}} & {} & 1  \\
	\end{matrix} \right)\left( \begin{matrix}
		{{\kappa }_{A{{P}_{A}}}}  \\
		{{\kappa }_{B{{P}_{B}}}}  \\
		{{\kappa }_{C{{P}_{C}}}}  \\
		{{\kappa }_{D{{P}_{D}}}}  \\
	\end{matrix} \right)=\left( \begin{matrix}
		1  \\
		1  \\
		1  \\
		1-\alpha _{D}^{{{P}_{A}}}  \\
	\end{matrix} \right),\]	
	\[\left( \begin{matrix}
		1 & \alpha _{A}^{{{P}_{B}}} & {} & {}  \\
		{} & 1 & \alpha _{B}^{{{P}_{C}}} & {}  \\
		{} & {} & 1 & \alpha _{C}^{{{P}_{D}}}  \\
		{} & {} & \alpha _{D}^{{{P}_{A}}}\alpha _{A}^{{{P}_{B}}}\alpha _{B}^{{{P}_{C}}} & 1  \\
	\end{matrix} \right)\left( \begin{matrix}
		{{\kappa }_{A{{P}_{A}}}}  \\
		{{\kappa }_{B{{P}_{B}}}}  \\
		{{\kappa }_{C{{P}_{C}}}}  \\
		{{\kappa }_{D{{P}_{D}}}}  \\
	\end{matrix} \right)=\left( \begin{matrix}
		1  \\
		1  \\
		1  \\
		1-\alpha _{D}^{{{P}_{A}}}+\alpha _{D}^{{{P}_{A}}}\alpha _{A}^{{{P}_{B}}} \\
	\end{matrix} \right),\]	
	\[\left( \begin{matrix}
		1 & \alpha _{A}^{{{P}_{B}}} & {} & {}  \\
		{} & 1 & \alpha _{B}^{{{P}_{C}}} & {}  \\
		{} & {} & 1 & \alpha _{C}^{{{P}_{D}}}  \\
		{} & {} & {} & 1-\alpha _{D}^{{{P}_{A}}}\alpha _{A}^{{{P}_{B}}}\alpha _{B}^{{{P}_{C}}}\alpha _{C}^{{{P}_{D}}}  \\
	\end{matrix} \right)\left( \begin{matrix}
		{{\kappa }_{A{{P}_{A}}}}  \\
		{{\kappa }_{B{{P}_{B}}}}  \\
		{{\kappa }_{C{{P}_{C}}}}  \\
		{{\kappa }_{D{{P}_{D}}}}  \\
	\end{matrix} \right)=\left( \begin{matrix}
		1  \\
		1  \\
		1  \\
		1-\alpha _{D}^{{{P}_{A}}}+\alpha _{D}^{{{P}_{A}}}\alpha _{A}^{{{P}_{B}}}-\alpha _{D}^{{{P}_{A}}}\alpha _{A}^{{{P}_{B}}}\alpha _{B}^{{{P}_{C}}} \\
	\end{matrix} \right).\]			
	
	The solution is in the following:
	\[{{\kappa }_{D{{P}_{D}}}}=\frac{1-\alpha _{D}^{{{P}_{A}}}+\alpha _{D}^{{{P}_{A}}}\alpha _{A}^{{{P}_{B}}}-\alpha _{D}^{{{P}_{A}}}\alpha _{A}^{{{P}_{B}}}\alpha _{B}^{{{P}_{C}}}}{1-\alpha _{D}^{{{P}_{A}}}\alpha _{A}^{{{P}_{B}}}\alpha _{B}^{{{P}_{C}}}\alpha _{C}^{{{P}_{D}}}}.\]	
	
	After rotation, the following results can be obtained:
	\[{{\kappa }_{A{{P}_{A}}}}=\frac{1-\alpha _{A}^{{{P}_{B}}}+\alpha _{A}^{{{P}_{B}}}\alpha _{B}^{{{P}_{C}}}-\alpha _{A}^{{{P}_{B}}}\alpha _{B}^{{{P}_{C}}}\alpha _{C}^{{{P}_{D}}}}{1-\alpha _{A}^{{{P}_{B}}}\alpha _{B}^{{{P}_{C}}}\alpha _{C}^{{{P}_{D}}}\alpha _{D}^{{{P}_{A}}}},\]	
	\[{{\kappa }_{B{{P}_{B}}}}=\frac{1-\alpha _{B}^{{{P}_{C}}}+\alpha _{B}^{{{P}_{C}}}\alpha _{C}^{{{P}_{D}}}-\alpha _{B}^{{{P}_{C}}}\alpha _{C}^{{{P}_{D}}}\alpha _{D}^{{{P}_{A}}}}{1-\alpha _{B}^{{{P}_{C}}}\alpha _{C}^{{{P}_{D}}}\alpha _{D}^{{{P}_{A}}}\alpha _{A}^{{{P}_{B}}}},\]	
	\[{{\kappa }_{C{{P}_{C}}}}=\frac{1-\alpha _{C}^{{{P}_{D}}}+\alpha _{C}^{{{P}_{D}}}\alpha _{D}^{{{P}_{A}}}-\alpha _{C}^{{{P}_{D}}}\alpha _{D}^{{{P}_{A}}}\alpha _{A}^{{{P}_{B}}}}{1-\alpha _{C}^{{{P}_{D}}}\alpha _{D}^{{{P}_{A}}}\alpha _{A}^{{{P}_{B}}}\alpha _{B}^{{{P}_{C}}}}.\]	
\end{proof}
\hfill $\square$\par

Several formulas, such as ${{\kappa }_{A{{P}_{A}}}}$ etc., given in theorem \ref{thm:Thm20.1.1}, are easy to calculate, but there are other expressions. See the following theorem. 


\begin{theorem}{Calculating integral ratio (2) based on IR-F, Daiyuan Zhang}{Thm20.1.3}\label{Thm20.1.3} 
	Integral ratios of vertex to intersecting foot of tetrahedron are:	
	\[{{\kappa }_{A{{P}_{A}}}}=\frac{1-\alpha _{A}^{{{P}_{B}}}}{1-\alpha _{A}^{{{P}_{B}}}\alpha _{B}^{{{P}_{A}}}}=\frac{1-\alpha _{A}^{{{P}_{C}}}}{1-\alpha _{A}^{{{P}_{C}}}\alpha _{C}^{{{P}_{A}}}}=\frac{1-\alpha _{A}^{{{P}_{D}}}}{1-\alpha _{A}^{{{P}_{D}}}\alpha _{D}^{{{P}_{A}}}},\]
	\[{{\kappa }_{B{{P}_{B}}}}=\frac{1-\alpha _{B}^{{{P}_{C}}}}{1-\alpha _{B}^{{{P}_{C}}}\alpha _{C}^{{{P}_{B}}}}=\frac{1-\alpha _{B}^{{{P}_{D}}}}{1-\alpha _{B}^{{{P}_{D}}}\alpha _{D}^{{{P}_{B}}}}=\frac{1-\alpha _{B}^{{{P}_{A}}}}{1-\alpha _{B}^{{{P}_{A}}}\alpha _{A}^{{{P}_{B}}}},\]
	\[{{\kappa }_{C{{P}_{C}}}}=\frac{1-\alpha _{C}^{{{P}_{D}}}}{1-\alpha _{C}^{{{P}_{D}}}\alpha _{D}^{{{P}_{C}}}}=\frac{1-\alpha _{C}^{{{P}_{A}}}}{1-\alpha _{C}^{{{P}_{A}}}\alpha _{A}^{{{P}_{C}}}}=\frac{1-\alpha _{C}^{{{P}_{B}}}}{1-\alpha _{C}^{{{P}_{B}}}\alpha _{B}^{{{P}_{C}}}},\]
	\[{{\kappa }_{D{{P}_{D}}}}=\frac{1-\alpha _{D}^{{{P}_{A}}}}{1-\alpha _{D}^{{{P}_{A}}}\alpha _{A}^{{{P}_{D}}}}=\frac{1-\alpha _{D}^{{{P}_{B}}}}{1-\alpha _{D}^{{{P}_{B}}}\alpha _{B}^{{{P}_{D}}}}=\frac{1-\alpha _{D}^{{{P}_{C}}}}{1-\alpha _{D}^{{{P}_{C}}}\alpha _{C}^{{{P}_{D}}}}.\]
\end{theorem}

\begin{proof}
	
	Only the proof of ${{\kappa }_{A{{P}_{A}}}}$ is given. The rest are similar. According to the equation (\ref{Eq19.1.1}) and (\ref{Eq19.1.2}), A system of functions containing only ${{\kappa }_{A{{P}_{A}}}}$ and ${{\kappa }_{B{{P}_{B}}}}$ can be obtained: 								
	\[\left\{ \begin{aligned}
		& 1-{{\kappa }_{A{{P}_{A}}}}={{\kappa }_{B{{P}_{B}}}}\alpha _{A}^{{{P}_{B}}} \\ 
		& {{\kappa }_{A{{P}_{A}}}}\alpha _{B}^{{{P}_{A}}}=1-{{\kappa }_{B{{P}_{B}}}}. \\ 
	\end{aligned} \right.\]	
	
	By solving the above equation, we have:
	\[{{\kappa }_{A{{P}_{A}}}}=\frac{1-\alpha _{A}^{{{P}_{B}}}}{1-\alpha _{A}^{{{P}_{B}}}\alpha _{B}^{{{P}_{A}}}}.\]	
	
	Similarly, according to equations (\ref{Eq19.1.1}) and (\ref{Eq19.1.3}), A system of functions containing only ${{\kappa }_{A{{P}_{A}}}}$ and ${{\kappa }_{C{{P}_{C}}}}$ can be obtained: 
	\[\left\{ \begin{aligned}
		& 1-{{\kappa }_{A{{P}_{A}}}}={{\kappa }_{C{{P}_{C}}}}\alpha _{A}^{{{P}_{C}}} \\ 
		& {{\kappa }_{A{{P}_{A}}}}\alpha _{C}^{{{P}_{A}}}=1-{{\kappa }_{C{{P}_{C}}}}. \\ 
	\end{aligned} \right.\]	
	
	By solving the above equation, we have:
	\[{{\kappa }_{A{{P}_{A}}}}=\frac{1-\alpha _{A}^{{{P}_{C}}}}{1-\alpha _{A}^{{{P}_{C}}}\alpha _{C}^{{{P}_{A}}}}.\]	
	
	According to the equation (\ref{Eq19.1.1}) and (\ref{Eq19.1.4}), A system of functions containing only ${{\kappa }_{A{{P}_{A}}}}$ and ${{\kappa }_{D{{P}_{D}}}}$ can be obtained: 	
	\[\left\{ \begin{aligned}
		& 1-{{\kappa }_{A{{P}_{A}}}}={{\kappa }_{D{{P}_{D}}}}\alpha _{A}^{{{P}_{D}}} \\ 
		& {{\kappa }_{A{{P}_{A}}}}\alpha _{D}^{{{P}_{A}}}=1-{{\kappa }_{D{{P}_{D}}}}. \\ 
	\end{aligned} \right.\]	
	
	By solving the above equation, we have:
	\[{{\kappa }_{A{{P}_{A}}}}=\frac{1-\alpha _{A}^{{{P}_{D}}}}{1-\alpha _{A}^{{{P}_{D}}}\alpha _{D}^{{{P}_{A}}}}.\]
	
	Therefore 
	\[{{\kappa }_{A{{P}_{A}}}}=\frac{1-\alpha _{A}^{{{P}_{B}}}}{1-\alpha _{A}^{{{P}_{B}}}\alpha _{B}^{{{P}_{A}}}}=\frac{1-\alpha _{A}^{{{P}_{C}}}}{1-\alpha _{A}^{{{P}_{C}}}\alpha _{C}^{{{P}_{A}}}}=\frac{1-\alpha _{A}^{{{P}_{D}}}}{1-\alpha _{A}^{{{P}_{D}}}\alpha _{D}^{{{P}_{A}}}}.\]
	
	Similarly, 
	\[{{\kappa }_{B{{P}_{B}}}}=\frac{1-\alpha _{B}^{{{P}_{C}}}}{1-\alpha _{B}^{{{P}_{C}}}\alpha _{C}^{{{P}_{B}}}}=\frac{1-\alpha _{B}^{{{P}_{D}}}}{1-\alpha _{B}^{{{P}_{D}}}\alpha _{D}^{{{P}_{B}}}}=\frac{1-\alpha _{B}^{{{P}_{A}}}}{1-\alpha _{B}^{{{P}_{A}}}\alpha _{A}^{{{P}_{B}}}},\]
	\[{{\kappa }_{C{{P}_{C}}}}=\frac{1-\alpha _{C}^{{{P}_{D}}}}{1-\alpha _{C}^{{{P}_{D}}}\alpha _{D}^{{{P}_{C}}}}=\frac{1-\alpha _{C}^{{{P}_{A}}}}{1-\alpha _{C}^{{{P}_{A}}}\alpha _{A}^{{{P}_{C}}}}=\frac{1-\alpha _{C}^{{{P}_{B}}}}{1-\alpha _{C}^{{{P}_{B}}}\alpha _{B}^{{{P}_{C}}}},\]
	\[{{\kappa }_{D{{P}_{D}}}}=\frac{1-\alpha _{D}^{{{P}_{A}}}}{1-\alpha _{D}^{{{P}_{A}}}\alpha _{A}^{{{P}_{D}}}}=\frac{1-\alpha _{D}^{{{P}_{B}}}}{1-\alpha _{D}^{{{P}_{B}}}\alpha _{B}^{{{P}_{D}}}}=\frac{1-\alpha _{D}^{{{P}_{C}}}}{1-\alpha _{D}^{{{P}_{C}}}\alpha _{C}^{{{P}_{D}}}}.\]
\end{proof}
\hfill $\square$\par
Another symmetrical form of calculating the integral ratio according to the IR-F is given below.

\begin{theorem}{Symmetric form 2 for calculating integral ratio based on IR-F, Daiyuan Zhang}{Thm20.1.4}\label{Thm20.1.4} 
	Integral ratios of vertex to intersecting foot of tetrahedron are:
	\[{{\kappa }_{A{{P}_{A}}}}=\frac{3-\left( \alpha _{A}^{{{P}_{B}}}+\alpha _{A}^{{{P}_{C}}}+\alpha _{A}^{{{P}_{D}}} \right)}{3-\left( \alpha _{A}^{{{P}_{B}}}\alpha _{B}^{{{P}_{A}}}+\alpha _{A}^{{{P}_{C}}}\alpha _{C}^{{{P}_{A}}}+\alpha _{A}^{{{P}_{D}}}\alpha _{D}^{{{P}_{A}}} \right)},\]
	\[{{\kappa }_{B{{P}_{B}}}}=\frac{3-\left( \alpha _{B}^{{{P}_{C}}}+\alpha _{B}^{{{P}_{D}}}+\alpha _{B}^{{{P}_{A}}} \right)}{3-\left( \alpha _{B}^{{{P}_{C}}}\alpha _{C}^{{{P}_{B}}}+\alpha _{B}^{{{P}_{D}}}\alpha _{D}^{{{P}_{B}}}+\alpha _{B}^{{{P}_{A}}}\alpha _{A}^{{{P}_{B}}} \right)},\]
	\[{{\kappa }_{C{{P}_{C}}}}=\frac{3-\left( \alpha _{C}^{{{P}_{D}}}+\alpha _{C}^{{{P}_{A}}}+\alpha _{C}^{{{P}_{B}}} \right)}{3-\left( \alpha _{C}^{{{P}_{D}}}\alpha _{D}^{{{P}_{C}}}+\alpha _{C}^{{{P}_{A}}}\alpha _{A}^{{{P}_{C}}}+\alpha _{C}^{{{P}_{B}}}\alpha _{B}^{{{P}_{C}}} \right)},\]
	\[{{\kappa }_{D{{P}_{D}}}}=\frac{3-\left( \alpha _{D}^{{{P}_{A}}}+\alpha _{D}^{{{P}_{B}}}+\alpha _{D}^{{{P}_{C}}} \right)}{3-\left( \alpha _{D}^{{{P}_{A}}}\alpha _{A}^{{{P}_{D}}}+\alpha _{D}^{{{P}_{B}}}\alpha _{B}^{{{P}_{D}}}+\alpha _{D}^{{{P}_{C}}}\alpha _{C}^{{{P}_{D}}} \right)}.\]
\end{theorem}

\begin{proof}
	From theorem \ref{thm:Thm20.1.2}, we have:
	\[\begin{aligned}
		& {{\kappa }_{A{{P}_{A}}}}=\frac{1-\alpha _{A}^{{{P}_{B}}}}{1-\alpha _{A}^{{{P}_{B}}}\alpha _{B}^{{{P}_{A}}}}=\frac{1-\alpha _{A}^{{{P}_{C}}}}{1-\alpha _{A}^{{{P}_{C}}}\alpha _{C}^{{{P}_{A}}}}=\frac{1-\alpha _{A}^{{{P}_{D}}}}{1-\alpha _{A}^{{{P}_{D}}}\alpha _{D}^{{{P}_{A}}}} \\ 
		& =\frac{1-\alpha _{A}^{{{P}_{B}}}+1-\alpha _{A}^{{{P}_{C}}}+1-\alpha _{A}^{{{P}_{D}}}}{1-\alpha _{A}^{{{P}_{B}}}\alpha _{B}^{{{P}_{A}}}+1-\alpha _{A}^{{{P}_{C}}}\alpha _{C}^{{{P}_{A}}}+1-\alpha _{A}^{{{P}_{D}}}\alpha _{D}^{{{P}_{A}}}} \\ 
		& =\frac{3-\left( \alpha _{A}^{{{P}_{B}}}+\alpha _{A}^{{{P}_{C}}}+\alpha _{A}^{{{P}_{D}}} \right)}{3-\left( \alpha _{A}^{{{P}_{B}}}\alpha _{B}^{{{P}_{A}}}+\alpha _{A}^{{{P}_{C}}}\alpha _{C}^{{{P}_{A}}}+\alpha _{A}^{{{P}_{D}}}\alpha _{D}^{{{P}_{A}}} \right)}.  
	\end{aligned}\]
	
	The other formulas can be proved similarly.
\end{proof}
\hfill $\square$\par
\subsection{Calculate IR-F according to integer ratio}\label{Subsec20.1.2}	
In the previous section, the method of calculating the integral ratio by using the IR-F is studied. In this section, the opposite problem is studied, that is, the IR-F is calculated by the integral ratio.

\begin{theorem}{Calculate IR-F by integer ratio, Daiyuan Zhang}{Thm20.1.5}\label{Thm20.1.5} 
	The IR-Fs of the tetrahedron are:
	\[\alpha _{B}^{{{P}_{A}}}=\frac{1-{{\kappa }_{B{{P}_{B}}}}}{{{\kappa }_{A{{P}_{A}}}}},\,\alpha _{C}^{{{P}_{A}}}=\frac{1-{{\kappa }_{C{{P}_{C}}}}}{{{\kappa }_{A{{P}_{A}}}}},\ \alpha _{D}^{{{P}_{A}}}=\frac{1-{{\kappa }_{D{{P}_{D}}}}}{{{\kappa }_{A{{P}_{A}}}}};\]
	\[\alpha _{C}^{{{P}_{B}}}=\frac{1-{{\kappa }_{C{{P}_{C}}}}}{{{\kappa }_{B{{P}_{B}}}}},\,\alpha _{D}^{{{P}_{B}}}=\frac{1-{{\kappa }_{D{{P}_{D}}}}}{{{\kappa }_{B{{P}_{B}}}}},\ \alpha _{A}^{{{P}_{B}}}=\frac{1-{{\kappa }_{A{{P}_{A}}}}}{{{\kappa }_{B{{P}_{B}}}}};\]
	\[\alpha _{D}^{{{P}_{C}}}=\frac{1-{{\kappa }_{D{{P}_{D}}}}}{{{\kappa }_{C{{P}_{C}}}}},\,\alpha _{A}^{{{P}_{C}}}=\frac{1-{{\kappa }_{A{{P}_{A}}}}}{{{\kappa }_{C{{P}_{C}}}}},\ \alpha _{B}^{{{P}_{C}}}=\frac{1-{{\kappa }_{B{{P}_{B}}}}}{{{\kappa }_{C{{P}_{C}}}}};\]
	\[\alpha _{A}^{{{P}_{D}}}=\frac{1-{{\kappa }_{A{{P}_{A}}}}}{{{\kappa }_{D{{P}_{D}}}}},\,\alpha _{B}^{{{P}_{D}}}=\frac{1-{{\kappa }_{B{{P}_{B}}}}}{{{\kappa }_{D{{P}_{D}}}}},\ \alpha _{C}^{{{P}_{D}}}=\frac{1-{{\kappa }_{C{{P}_{C}}}}}{{{\kappa }_{D{{P}_{D}}}}}.\]
\end{theorem}

\begin{proof}
	According to the formula (\ref{Eq19.1.2})-(\ref{Eq19.1.4}), the following result is obtained directly:
	\[\alpha _{B}^{{{P}_{A}}}=\frac{1-{{\kappa }_{B{{P}_{B}}}}}{{{\kappa }_{A{{P}_{A}}}}},\,\alpha _{C}^{{{P}_{A}}}=\frac{1-{{\kappa }_{C{{P}_{C}}}}}{{{\kappa }_{A{{P}_{A}}}}},\ \alpha _{D}^{{{P}_{A}}}=\frac{1-{{\kappa }_{D{{P}_{D}}}}}{{{\kappa }_{A{{P}_{A}}}}}.\]
	
	The other formulas can be proved similarly.
\end{proof}
\hfill $\square$\par
\begin{theorem}{Harmonic expression 1 on frame components of vertex at IC-F, Daiyuan Zhang}{Thm20.1.6}\label{Thm20.1.6} 
	\[\frac{1}{\alpha _{A}^{{{P}_{B}}}}+\frac{1}{\alpha _{A}^{{{P}_{C}}}}+\frac{1}{\alpha _{A}^{{{P}_{D}}}}=\frac{3-{{\kappa }_{A{{P}_{A}}}}}{1-{{\kappa }_{A{{P}_{A}}}}},\]
	\[\frac{1}{\alpha _{B}^{{{P}_{C}}}}+\frac{1}{\alpha _{B}^{{{P}_{D}}}}+\frac{1}{\alpha _{B}^{{{P}_{A}}}}=\frac{3-{{\kappa }_{B{{P}_{B}}}}}{1-{{\kappa }_{B{{P}_{B}}}}},\]
	\[\frac{1}{\alpha _{C}^{{{P}_{D}}}}+\frac{1}{\alpha _{C}^{{{P}_{A}}}}+\frac{1}{\alpha _{C}^{{{P}_{B}}}}=\frac{3-{{\kappa }_{C{{P}_{C}}}}}{1-{{\kappa }_{C{{P}_{C}}}}},\]
	\[\frac{1}{\alpha _{D}^{{{P}_{A}}}}+\frac{1}{\alpha _{D}^{{{P}_{B}}}}+\frac{1}{\alpha _{D}^{{{P}_{C}}}}=\frac{3-{{\kappa }_{D{{P}_{D}}}}}{1-{{\kappa }_{D{{P}_{D}}}}}.\]
\end{theorem}

\begin{proof}
	According to theorem \ref{thm:Thm20.1.4}, the following result is obtained:
	\[\frac{1}{\alpha _{A}^{{{P}_{B}}}}+\frac{1}{\alpha _{A}^{{{P}_{C}}}}+\frac{1}{\alpha _{A}^{{{P}_{D}}}}=\frac{{{\kappa }_{B{{P}_{B}}}}+{{\kappa }_{C{{P}_{C}}}}+{{\kappa }_{D{{P}_{D}}}}}{1-{{\kappa }_{A{{P}_{A}}}}}=\frac{3-{{\kappa }_{A{{P}_{A}}}}}{1-{{\kappa }_{A{{P}_{A}}}}}\]
	
	The other formulas can be proved similarly.
\end{proof}
\hfill $\square$\par

\section{The relationship of frame components between IC-Fs and IC-T}\label{Sec20.2}
\begin{theorem}{Calculate frame components (1) of IC-T by IC-Fs, Daiyuan Zhang}{Thm20.2.1}\label{Thm20.2.1} 
	The frame components for the given IC-T of a tetrahedron are	
	\begin{flalign*}
		\beta _{A}^{P}=\frac{\alpha _{A}^{{{P}_{B}}}\left( 1-\alpha _{B}^{{{P}_{A}}} \right)}{1-\alpha _{A}^{{{P}_{B}}}\alpha _{B}^{{{P}_{A}}}}, \beta _{B}^{P}=\frac{\alpha _{B}^{{{P}_{C}}}\left( 1-\alpha _{C}^{{{P}_{B}}} \right)}{1-\alpha _{B}^{{{P}_{C}}}\alpha _{C}^{{{P}_{B}}}},
	\end{flalign*}
	\begin{flalign*}
		\beta _{C}^{P}=\frac{\alpha _{C}^{{{P}_{D}}}\left( 1-\alpha _{D}^{{{P}_{C}}} \right)}{1-\alpha _{C}^{{{P}_{D}}}\alpha _{D}^{{{P}_{C}}}}, \beta _{D}^{P}=\frac{\alpha _{D}^{{{P}_{A}}}\left( 1-\alpha _{A}^{{{P}_{D}}} \right)}{1-\alpha _{D}^{{{P}_{A}}}\alpha _{A}^{{{P}_{D}}}}.
	\end{flalign*}
\end{theorem}

\begin{proof}
	From theorem \ref{thm:Thm18.2.1}, we have	
	\begin{flalign*}
		\beta _{A}^{P}=1-{{\kappa }_{A{{P}_{A}}}}, \beta _{B}^{P}=1-{{\kappa }_{B{{P}_{B}}}}, \beta _{C}^{P}=1-{{\kappa }_{C{{P}_{C}}}}, \beta _{D}^{P}={{\kappa }_{C{{P}_{C}}}}\alpha _{D}^{{{P}_{C}}}.
	\end{flalign*}
		
	Obviously, by the rotation of $A$→$B$→$C$→$D$→$A$, $\beta _{A}^{P}=1-{{\kappa }_{A{{P}_{A}}}}$ becomes $\beta _{B}^{P}=1-{{\kappa }_{B{{P}_{B}}}}$; $\beta _{B}^{P}=1-{{\kappa }_{B{{P}_{B}}}}$ becomes $\beta _{C}^{P}=1-{{\kappa }_{C{{P}_{C}}}}$; $\beta _{C}^{P}=1-{{\kappa }_{C{{P}_{C}}}}$ becomes $\beta _{D}^{P}=1-{{\kappa }_{D{{P}_{D}}}}$; $\beta _{D}^{P}=1-{{\kappa }_{D{{P}_{D}}}}$ becomes $\beta _{A}^{P}=1-{{\kappa }_{A{{P}_{A}}}}$. 
	
	This means that if we find the solution of one of the equations, we can get other solutions as long as we follow the same rotation ($A$→$B$→$C$→$D$→$A$) in the solution. 		
	
	From theorem \ref{thm:Thm20.1.1}, we have	
	\[\beta _{A}^{P}=1-{{\kappa }_{A{{P}_{A}}}}=1-\frac{1-\alpha _{A}^{{{P}_{B}}}}{1-\alpha _{A}^{{{P}_{B}}}\alpha _{B}^{{{P}_{A}}}}=\frac{\alpha _{A}^{{{P}_{B}}}\left( 1-\alpha _{B}^{{{P}_{A}}} \right)}{1-\alpha _{A}^{{{P}_{B}}}\alpha _{B}^{{{P}_{A}}}}.\]	
	
	By rotating the above results as $A$→$B$→$C$→$D$→$A$, the following results are obtained directly:
	\begin{flalign*}
		\beta _{B}^{P}=\frac{\alpha _{B}^{{{P}_{C}}}\left( 1-\alpha _{C}^{{{P}_{B}}} \right)}{1-\alpha _{B}^{{{P}_{C}}}\alpha _{C}^{{{P}_{B}}}},\beta _{C}^{P}=\frac{\alpha _{C}^{{{P}_{D}}}\left( 1-\alpha _{D}^{{{P}_{C}}} \right)}{1-\alpha _{C}^{{{P}_{D}}}\alpha _{D}^{{{P}_{C}}}},\beta _{D}^{P}=\frac{\alpha _{D}^{{{P}_{A}}}\left( 1-\alpha _{A}^{{{P}_{D}}} \right)}{1-\alpha _{D}^{{{P}_{A}}}\alpha _{A}^{{{P}_{D}}}}
	\end{flalign*}
\end{proof}
\hfill $\square$\par
%

The above theorem shows that the frame components of the intersecting center of a tetrahedron can be expressed by the corresponding frame components of the intersecting centers of faces.

Of course, the frame components of a tetrahedron can have other forms.

\begin{theorem}{Calculate frame components (2) of IC-T by IC-Fs, Daiyuan Zhang}{Thm20.2.2}\label{Thm20.2.2} 
	The frame components of the IC-Ts are:
	\[\beta _{A}^{P}=\frac{\alpha _{A}^{{{P}_{B}}}\left( 1-\alpha _{B}^{{{P}_{A}}} \right)}{1-\alpha _{A}^{{{P}_{B}}}\alpha _{B}^{{{P}_{A}}}}=\frac{\alpha _{A}^{{{P}_{C}}}\left( 1-\alpha _{C}^{{{P}_{A}}} \right)}{1-\alpha _{A}^{{{P}_{C}}}\alpha _{C}^{{{P}_{A}}}}=\frac{\alpha _{A}^{{{P}_{D}}}\left( 1-\alpha _{D}^{{{P}_{A}}} \right)}{1-\alpha _{A}^{{{P}_{D}}}\alpha _{D}^{{{P}_{A}}}}\text{,}\]
	\[\beta _{B}^{P}=\frac{\alpha _{B}^{{{P}_{C}}}\left( 1-\alpha _{C}^{{{P}_{B}}} \right)}{1-\alpha _{B}^{{{P}_{C}}}\alpha _{C}^{{{P}_{B}}}}=\frac{\alpha _{B}^{{{P}_{D}}}\left( 1-\alpha _{D}^{{{P}_{B}}} \right)}{1-\alpha _{B}^{{{P}_{D}}}\alpha _{D}^{{{P}_{B}}}}=\frac{\alpha _{B}^{{{P}_{A}}}\left( 1-\alpha _{A}^{{{P}_{B}}} \right)}{1-\alpha _{B}^{{{P}_{A}}}\alpha _{A}^{{{P}_{B}}}}\text{,}\]
	\[\beta _{C}^{P}=\frac{\alpha _{C}^{{{P}_{D}}}\left( 1-\alpha _{D}^{{{P}_{C}}} \right)}{1-\alpha _{C}^{{{P}_{D}}}\alpha _{D}^{{{P}_{C}}}}=\frac{\alpha _{C}^{{{P}_{A}}}\left( 1-\alpha _{A}^{{{P}_{C}}} \right)}{1-\alpha _{C}^{{{P}_{A}}}\alpha _{A}^{{{P}_{C}}}}=\frac{\alpha _{C}^{{{P}_{B}}}\left( 1-\alpha _{B}^{{{P}_{C}}} \right)}{1-\alpha _{C}^{{{P}_{B}}}\alpha _{B}^{{{P}_{C}}}}\text{,}\]
	\[\beta _{D}^{P}=\frac{\alpha _{D}^{{{P}_{A}}}\left( 1-\alpha _{A}^{{{P}_{D}}} \right)}{1-\alpha _{D}^{{{P}_{A}}}\alpha _{A}^{{{P}_{D}}}}=\frac{\alpha _{D}^{{{P}_{B}}}\left( 1-\alpha _{B}^{{{P}_{D}}} \right)}{1-\alpha _{D}^{{{P}_{B}}}\alpha _{B}^{{{P}_{D}}}}=\frac{\alpha _{D}^{{{P}_{C}}}\left( 1-\alpha _{C}^{{{P}_{D}}} \right)}{1-\alpha _{D}^{{{P}_{C}}}\alpha _{C}^{{{P}_{D}}}}\text{.}\]
\end{theorem}

\begin{proof}
	From theorem \ref{thm:Thm20.1.3} and \ref{thm:Thm18.2.1}, we have:
	\[{{\kappa }_{A{{P}_{A}}}}=\frac{1-\alpha _{A}^{{{P}_{B}}}}{1-\alpha _{A}^{{{P}_{B}}}\alpha _{B}^{{{P}_{A}}}}=\frac{1-\alpha _{A}^{{{P}_{C}}}}{1-\alpha _{A}^{{{P}_{C}}}\alpha _{C}^{{{P}_{A}}}}=\frac{1-\alpha _{A}^{{{P}_{D}}}}{1-\alpha _{A}^{{{P}_{D}}}\alpha _{D}^{{{P}_{A}}}},\]
	\[\beta _{A}^{P}=1-{{\kappa }_{A{{P}_{A}}}}=1-\frac{1-\alpha _{A}^{{{P}_{B}}}}{1-\alpha _{A}^{{{P}_{B}}}\alpha _{B}^{{{P}_{A}}}}=\frac{\alpha _{A}^{{{P}_{B}}}\left( 1-\alpha _{B}^{{{P}_{A}}} \right)}{1-\alpha _{A}^{{{P}_{B}}}\alpha _{B}^{{{P}_{A}}}}\text{,}\]
	\[\beta _{A}^{P}=1-{{\kappa }_{A{{P}_{A}}}}=1-\frac{1-\alpha _{A}^{{{P}_{C}}}}{1-\alpha _{A}^{{{P}_{C}}}\alpha _{C}^{{{P}_{A}}}}=\frac{\alpha _{A}^{{{P}_{C}}}\left( 1-\alpha _{C}^{{{P}_{A}}} \right)}{1-\alpha _{A}^{{{P}_{C}}}\alpha _{C}^{{{P}_{A}}}}\text{,}\]
	\[\beta _{A}^{P}=1-{{\kappa }_{A{{P}_{A}}}}=1-\frac{1-\alpha _{A}^{{{P}_{D}}}}{1-\alpha _{A}^{{{P}_{D}}}\alpha _{D}^{{{P}_{A}}}}=\frac{\alpha _{A}^{{{P}_{D}}}\left( 1-\alpha _{D}^{{{P}_{A}}} \right)}{1-\alpha _{A}^{{{P}_{D}}}\alpha _{D}^{{{P}_{A}}}}\text{.}\]
	
	Therefore 
	\[\beta _{A}^{P}=\frac{\alpha _{A}^{{{P}_{B}}}\left( 1-\alpha _{B}^{{{P}_{A}}} \right)}{1-\alpha _{A}^{{{P}_{B}}}\alpha _{B}^{{{P}_{A}}}}=\frac{\alpha _{A}^{{{P}_{C}}}\left( 1-\alpha _{C}^{{{P}_{A}}} \right)}{1-\alpha _{A}^{{{P}_{C}}}\alpha _{C}^{{{P}_{A}}}}=\frac{\alpha _{A}^{{{P}_{D}}}\left( 1-\alpha _{D}^{{{P}_{A}}} \right)}{1-\alpha _{A}^{{{P}_{D}}}\alpha _{D}^{{{P}_{A}}}}\text{.}\]
	
	The other formulas can be proved similarly.
\end{proof}
\hfill $\square$\par
The following theorem gives the symmetrical form of the frame components for the IC-T.
\begin{theorem}{Symmetrical form of frame components for IC-T, Daiyuan Zhang}{Thm20.2.3}\label{Thm20.2.3} 
	Symmetrical form of the frame components for the IC-T are as follows:
	\[\beta _{A}^{P}=\frac{\alpha _{A}^{{{P}_{B}}}+\alpha _{A}^{{{P}_{C}}}+\alpha _{A}^{{{P}_{D}}}-\left( \alpha _{A}^{{{P}_{B}}}\alpha _{B}^{{{P}_{A}}}+\alpha _{A}^{{{P}_{C}}}\alpha _{C}^{{{P}_{A}}}+\alpha _{A}^{{{P}_{D}}}\alpha _{D}^{{{P}_{A}}} \right)}{3-\left( \alpha _{A}^{{{P}_{B}}}\alpha _{B}^{{{P}_{A}}}+\alpha _{A}^{{{P}_{C}}}\alpha _{C}^{{{P}_{A}}}+\alpha _{A}^{{{P}_{D}}}\alpha _{D}^{{{P}_{A}}} \right)},\]
	\[\beta _{B}^{P}=\frac{\alpha _{B}^{{{P}_{C}}}+\alpha _{B}^{{{P}_{D}}}+\alpha _{B}^{{{P}_{A}}}-\left( \alpha _{B}^{{{P}_{C}}}\alpha _{C}^{{{P}_{B}}}+\alpha _{B}^{{{P}_{D}}}\alpha _{D}^{{{P}_{B}}}+\alpha _{B}^{{{P}_{A}}}\alpha _{A}^{{{P}_{B}}} \right)}{3-\left( \alpha _{B}^{{{P}_{C}}}\alpha _{C}^{{{P}_{B}}}+\alpha _{B}^{{{P}_{D}}}\alpha _{D}^{{{P}_{B}}}+\alpha _{B}^{{{P}_{A}}}\alpha _{A}^{{{P}_{B}}} \right)},\]
	\[\beta _{C}^{P}=\frac{\alpha _{C}^{{{P}_{D}}}+\alpha _{C}^{{{P}_{A}}}+\alpha _{C}^{{{P}_{B}}}-\left( \alpha _{C}^{{{P}_{D}}}\alpha _{D}^{{{P}_{C}}}+\alpha _{C}^{{{P}_{A}}}\alpha _{A}^{{{P}_{C}}}+\alpha _{C}^{{{P}_{B}}}\alpha _{B}^{{{P}_{C}}} \right)}{3-\left( \alpha _{C}^{{{P}_{D}}}\alpha _{D}^{{{P}_{C}}}+\alpha _{C}^{{{P}_{A}}}\alpha _{A}^{{{P}_{C}}}+\alpha _{C}^{{{P}_{B}}}\alpha _{B}^{{{P}_{C}}} \right)},\]
	\[\beta _{D}^{P}=\frac{\alpha _{D}^{{{P}_{A}}}+\alpha _{D}^{{{P}_{B}}}+\alpha _{D}^{{{P}_{C}}}-\left( \alpha _{D}^{{{P}_{A}}}\alpha _{A}^{{{P}_{D}}}+\alpha _{D}^{{{P}_{B}}}\alpha _{B}^{{{P}_{D}}}+\alpha _{D}^{{{P}_{C}}}\alpha _{C}^{{{P}_{D}}} \right)}{3-\left( \alpha _{D}^{{{P}_{A}}}\alpha _{A}^{{{P}_{D}}}+\alpha _{D}^{{{P}_{B}}}\alpha _{B}^{{{P}_{D}}}+\alpha _{D}^{{{P}_{C}}}\alpha _{C}^{{{P}_{D}}} \right)}.\]
\end{theorem}

\begin{proof}
	According to theorem \ref{thm:Thm20.2.2}, we get a symmetric form of $\beta _{A}^{P}$ in the following:
	\[\begin{aligned}
		\beta _{A}^{P}& =\frac{\alpha _{A}^{{{P}_{B}}}\left( 1-\alpha _{B}^{{{P}_{A}}} \right)+\alpha _{A}^{{{P}_{C}}}\left( 1-\alpha _{C}^{{{P}_{A}}} \right)+\alpha _{A}^{{{P}_{D}}}\left( 1-\alpha _{D}^{{{P}_{A}}} \right)}{1-\alpha _{A}^{{{P}_{B}}}\alpha _{B}^{{{P}_{A}}}+1-\alpha _{A}^{{{P}_{C}}}\alpha _{C}^{{{P}_{A}}}+1-\alpha _{A}^{{{P}_{D}}}\alpha _{D}^{{{P}_{A}}}} \\ 
		& =\frac{\alpha _{A}^{{{P}_{B}}}+\alpha _{A}^{{{P}_{C}}}+\alpha _{A}^{{{P}_{D}}}-\left( \alpha _{A}^{{{P}_{B}}}\alpha _{B}^{{{P}_{A}}}+\alpha _{A}^{{{P}_{C}}}\alpha _{C}^{{{P}_{A}}}+\alpha _{A}^{{{P}_{D}}}\alpha _{D}^{{{P}_{A}}} \right)}{3-\left( \alpha _{A}^{{{P}_{B}}}\alpha _{B}^{{{P}_{A}}}+\alpha _{A}^{{{P}_{C}}}\alpha _{C}^{{{P}_{A}}}+\alpha _{A}^{{{P}_{D}}}\alpha _{D}^{{{P}_{A}}} \right)}.  
	\end{aligned}\]
	
	Other formulas can be obtained similarly:
	\[\beta _{B}^{P}=\frac{\alpha _{B}^{{{P}_{C}}}+\alpha _{B}^{{{P}_{D}}}+\alpha _{B}^{{{P}_{A}}}-\left( \alpha _{B}^{{{P}_{C}}}\alpha _{C}^{{{P}_{B}}}+\alpha _{B}^{{{P}_{D}}}\alpha _{D}^{{{P}_{B}}}+\alpha _{B}^{{{P}_{A}}}\alpha _{A}^{{{P}_{B}}} \right)}{3-\left( \alpha _{B}^{{{P}_{C}}}\alpha _{C}^{{{P}_{B}}}+\alpha _{B}^{{{P}_{D}}}\alpha _{D}^{{{P}_{B}}}+\alpha _{B}^{{{P}_{A}}}\alpha _{A}^{{{P}_{B}}} \right)},\]
	\[\beta _{C}^{P}=\frac{\alpha _{C}^{{{P}_{D}}}+\alpha _{C}^{{{P}_{A}}}+\alpha _{C}^{{{P}_{B}}}-\left( \alpha _{C}^{{{P}_{D}}}\alpha _{D}^{{{P}_{C}}}+\alpha _{C}^{{{P}_{A}}}\alpha _{A}^{{{P}_{C}}}+\alpha _{C}^{{{P}_{B}}}\alpha _{B}^{{{P}_{C}}} \right)}{3-\left( \alpha _{C}^{{{P}_{D}}}\alpha _{D}^{{{P}_{C}}}+\alpha _{C}^{{{P}_{A}}}\alpha _{A}^{{{P}_{C}}}+\alpha _{C}^{{{P}_{B}}}\alpha _{B}^{{{P}_{C}}} \right)},\]
	\[\beta _{D}^{P}=\frac{\alpha _{D}^{{{P}_{A}}}+\alpha _{D}^{{{P}_{B}}}+\alpha _{D}^{{{P}_{C}}}-\left( \alpha _{D}^{{{P}_{A}}}\alpha _{A}^{{{P}_{D}}}+\alpha _{D}^{{{P}_{B}}}\alpha _{B}^{{{P}_{D}}}+\alpha _{D}^{{{P}_{C}}}\alpha _{C}^{{{P}_{D}}} \right)}{3-\left( \alpha _{D}^{{{P}_{A}}}\alpha _{A}^{{{P}_{D}}}+\alpha _{D}^{{{P}_{B}}}\alpha _{B}^{{{P}_{D}}}+\alpha _{D}^{{{P}_{C}}}\alpha _{C}^{{{P}_{D}}} \right)}.\]
\end{proof}
\hfill $\square$\par
\subsection{Calculate frame components of IC-Fs by frame components of IC-T}\label{Subsec20.2.2}
The above theorem is to calculate the frame components of the IC-T according to the frame components of the IC-Fs. Now, we study the opposite problem, that is, to calculate the frame components of the IC-Fs according to the frame components of the IC-T.

\begin{theorem}{Calculate frame components of IC-Fs by IC-T, Daiyuan Zhang}{Thm20.2.4}\label{Thm20.2.4} 
	The frame components of the IC-Fs for a given tetrahedron are:
	\[\alpha _{B}^{{{P}_{A}}}=\frac{\beta _{B}^{P}}{1-\beta _{A}^{P}},\ \alpha _{C}^{{{P}_{A}}}=\frac{\beta _{C}^{P}}{1-\beta _{A}^{P}},\ \alpha _{D}^{{{P}_{A}}}=\frac{\beta _{D}^{P}}{1-\beta _{A}^{P}};\]	
	\[\alpha _{C}^{{{P}_{B}}}=\frac{\beta _{C}^{P}}{1-\beta _{B}^{P}},\ \alpha _{D}^{{{P}_{B}}}=\frac{\beta _{D}^{P}}{1-\beta _{B}^{P}},\ \alpha _{A}^{{{P}_{B}}}=\frac{\beta _{A}^{P}}{1-\beta _{B}^{P}};\]	
	\[\alpha _{D}^{{{P}_{C}}}=\frac{\beta _{D}^{P}}{1-\beta _{C}^{P}},\ \alpha _{A}^{{{P}_{C}}}=\frac{\beta _{A}^{P}}{1-\beta _{C}^{P}},\ \alpha _{B}^{{{P}_{C}}}=\frac{\beta _{B}^{P}}{1-\beta _{C}^{P}};\]	
	\[\alpha _{A}^{{{P}_{D}}}=\frac{\beta _{A}^{P}}{1-\beta _{D}^{P}},\ \alpha _{B}^{{{P}_{D}}}=\frac{\beta _{B}^{P}}{1-\beta _{D}^{P}},\ \alpha _{C}^{{{P}_{D}}}=\frac{\beta _{C}^{P}}{1-\beta _{D}^{P}}.\]		
\end{theorem}

\begin{proof}
	Consider the formula in theorem \ref{thm:Thm20.2.1}:
	\[\beta _{A}^{P}=\frac{\alpha _{A}^{{{P}_{B}}}\left( 1-\alpha _{B}^{{{P}_{A}}} \right)}{1-\alpha _{A}^{{{P}_{B}}}\alpha _{B}^{{{P}_{A}}}}.\]
	
	Since the positions of the two vertices $A$ and $B$ are equivalent, the positions of the vertices $A$ and $B$ can be exchanged, and the following results are obtained:
	\[\beta _{B}^{P}=\frac{\alpha _{B}^{{{P}_{A}}}\left( 1-\alpha _{A}^{{{P}_{B}}} \right)}{1-\alpha _{B}^{{{P}_{A}}}\alpha _{A}^{{{P}_{B}}}}.\]
	
	In fact, according to formulas (\ref{Eq18.2.7}) and (\ref{Eq20.1.2}), the following results are obtained:
	\[\beta _{B}^{P}={{\kappa }_{A{{P}_{A}}}}\alpha _{B}^{{{P}_{A}}}=\frac{\alpha _{B}^{{{P}_{A}}}\left( 1-\alpha _{A}^{{{P}_{B}}} \right)}{1-\alpha _{A}^{{{P}_{B}}}\alpha _{B}^{{{P}_{A}}}},\]	
	
	This formula is the same as the above formula. Combine the above two formulas to obtain a system of equations:
	\begin{equation}\label{Eq20.2.1}
		\left\{ \begin{aligned}
			& \beta _{A}^{P}=\frac{\alpha _{A}^{{{P}_{B}}}\left( 1-\alpha _{B}^{{{P}_{A}}} \right)}{1-\alpha _{A}^{{{P}_{B}}}\alpha _{B}^{{{P}_{A}}}} \\ 
			& \beta _{B}^{P}=\frac{\alpha _{B}^{{{P}_{A}}}\left( 1-\alpha _{A}^{{{P}_{B}}} \right)}{1-\alpha _{B}^{{{P}_{A}}}\alpha _{A}^{{{P}_{B}}}}. \\ 
		\end{aligned} \right.
	\end{equation}
	
	From the first equation of the system of equations (\ref{Eq20.2.1}), we have
	\[\beta _{A}^{P}\left( 1-\alpha _{A}^{{{P}_{B}}}\alpha _{B}^{{{P}_{A}}} \right)=\alpha _{A}^{{{P}_{B}}}\left( 1-\alpha _{B}^{{{P}_{A}}} \right),\]
	i.e.
	\begin{equation}\label{Eq20.2.2}
		\left( 1-\beta _{A}^{P} \right)\alpha _{A}^{{{P}_{B}}}\alpha _{B}^{{{P}_{A}}}=\alpha _{A}^{{{P}_{B}}}-\beta _{A}^{P},
	\end{equation}
	\begin{equation}\label{Eq20.2.3}
		\alpha _{A}^{{{P}_{B}}}\alpha _{B}^{{{P}_{A}}}=\frac{\alpha _{A}^{{{P}_{B}}}-\beta _{A}^{P}}{1-\beta _{A}^{P}}.
	\end{equation}
	
	Similarly, from the second equation of the system of equations (\ref{Eq20.2.1}), we have
	\[\alpha _{B}^{{{P}_{A}}}\alpha _{A}^{{{P}_{B}}}=\frac{\alpha _{B}^{{{P}_{A}}}-\beta _{B}^{P}}{1-\beta _{B}^{P}},\]
	i.e.
	\[\frac{\alpha _{B}^{{{P}_{A}}}-\beta _{B}^{P}}{1-\beta _{B}^{P}}=\frac{\alpha _{A}^{{{P}_{B}}}-\beta _{A}^{P}}{1-\beta _{A}^{P}}.\]
	
	Theorefore
	\[\begin{aligned}
		\alpha _{B}^{{{P}_{A}}}& =\frac{1-\beta _{B}^{P}}{1-\beta _{A}^{P}}\left( \alpha _{A}^{{{P}_{B}}}-\beta _{A}^{P} \right)+\beta _{B}^{P}=\frac{\left( 1-\beta _{B}^{P} \right)\alpha _{A}^{{{P}_{B}}}}{1-\beta _{A}^{P}}-\frac{\left( 1-\beta _{B}^{P} \right)\beta _{A}^{P}}{1-\beta _{A}^{P}}+\beta _{B}^{P} \\ 
		& =\frac{\left( 1-\beta _{B}^{P} \right)\alpha _{A}^{{{P}_{B}}}}{1-\beta _{A}^{P}}-\frac{\left( 1-\beta _{B}^{P} \right)\beta _{A}^{P}-\left( 1-\beta _{A}^{P} \right)\beta _{B}^{P}}{1-\beta _{A}^{P}} \\ 
		& =\frac{\left( 1-\beta _{B}^{P} \right)\alpha _{A}^{{{P}_{B}}}}{1-\beta _{A}^{P}}-\frac{\beta _{A}^{P}-\beta _{B}^{P}}{1-\beta _{A}^{P}}=\frac{1}{1-\beta _{A}^{P}}\left( \left( 1-\beta _{B}^{P} \right)\alpha _{A}^{{{P}_{B}}}+\beta _{B}^{P}-\beta _{A}^{P} \right).  
	\end{aligned}\]
	
	Substitute the above formula into the formula (\ref{Eq20.2.2}) to obtain:
	\[\alpha _{A}^{{{P}_{B}}}\left( \left( 1-\beta _{B}^{P} \right)\alpha _{A}^{{{P}_{B}}}+\beta _{B}^{P}-\beta _{A}^{P} \right)=\alpha _{A}^{{{P}_{B}}}-\beta _{A}^{P},\]	
	i.e.
	\[\left( 1-\beta _{B}^{P} \right){{\left( \alpha _{A}^{{{P}_{B}}} \right)}^{2}}+\left( \beta _{B}^{P}-\beta _{A}^{P}-1 \right)\alpha _{A}^{{{P}_{B}}}+\beta _{A}^{P}=0.\]	
	
	Solving the above quadratic equation yields:
	\[\alpha _{A}^{{{P}_{B}}}=\frac{\beta _{A}^{P}-\beta _{B}^{P}+1\pm \sqrt{{{\left( \beta _{B}^{P}-\beta _{A}^{P}-1 \right)}^{2}}-4\left( 1-\beta _{B}^{P} \right)\beta _{A}^{P}}}{2\left( 1-\beta _{B}^{P} \right)}.\]	
	
	And
	\[\begin{aligned}
		& {{\left( \beta _{B}^{P}-\beta _{A}^{P}-1 \right)}^{2}}-4\left( 1-\beta _{B}^{P} \right)\beta _{A}^{P} \\ 
		& ={{\left( \beta _{B}^{P} \right)}^{2}}+{{\left( \beta _{A}^{P} \right)}^{2}}+1-2\beta _{B}^{P}\beta _{A}^{P}-2\beta _{B}^{P}+2\beta _{A}^{P}-4\beta _{A}^{P}+4\beta _{B}^{P}\beta _{A}^{P} \\ 
		& ={{\left( \beta _{A}^{P}+\beta _{B}^{P} \right)}^{2}}-2\left( \beta _{A}^{P}+\beta _{B}^{P} \right)+1={{\left( \beta _{A}^{P}+\beta _{B}^{P}-1 \right)}^{2}}. \\ 
	\end{aligned}\]
	
	Therefore
	\[\alpha _{A}^{{{P}_{B}}}=\frac{\beta _{A}^{P}-\beta _{B}^{P}+1\pm \left( \beta _{A}^{P}+\beta _{B}^{P}-1 \right)}{2\left( 1-\beta _{B}^{P} \right)}=\left\{ \begin{aligned}
		& \frac{\beta _{A}^{P}}{1-\beta _{B}^{P}} \\ 
		& 1 \\ 
	\end{aligned} \right..\]	
	
	According to the formula (\ref{Eq20.2.3}), if $\alpha _{A}^{{{P}_{B}}}={\beta _{A}^{P}}/{\left( 1-\beta _{B}^{P} \right)}\;$, we get:
	\[\alpha _{B}^{{{P}_{A}}}=\frac{\alpha _{A}^{{{P}_{B}}}-\beta _{A}^{P}}{\alpha _{A}^{{{P}_{B}}}\left( 1-\beta _{A}^{P} \right)}=\frac{\frac{\beta _{A}^{P}}{1-\beta _{B}^{P}}-\beta _{A}^{P}}{\frac{\beta _{A}^{P}}{1-\beta _{B}^{P}}\left( 1-\beta _{A}^{P} \right)}=\frac{\beta _{B}^{P}}{1-\beta _{A}^{P}};\]	
	
	If $\alpha _{A}^{{{P}_{B}}}=1$, we have
	\[\alpha _{B}^{{{P}_{A}}}=\frac{\alpha _{A}^{{{P}_{B}}}-\beta _{A}^{P}}{\alpha _{A}^{{{P}_{B}}}\left( 1-\beta _{A}^{P} \right)}=\frac{1-\beta _{A}^{P}}{\left( 1-\beta _{A}^{P} \right)}=1.\]	
	
	Let $\alpha _{A}^{{{P}_{B}}}=1$ and $\alpha _{B}^{{{P}_{A}}}=1$, and they are substituted into the original equations (\ref{Eq20.2.1}), then the denominators are 0, so this set of solutions $\alpha _{A}^{{{P}_{B}}}=1$ and $\alpha _{B}^{{{P}_{A}}}=1$ should be discarded, we get:
	\begin{equation}\label{Eq20.2.4}
		\alpha _{A}^{{{P}_{B}}}=\frac{\beta _{A}^{P}}{1-\beta _{B}^{P}},\ \alpha _{B}^{{{P}_{A}}}=\frac{\beta _{B}^{P}}{1-\beta _{A}^{P}}.
	\end{equation}
	
	In fact, replace $B$ on the right side of the expression $\beta _{A}^{P}$ with $C$, we get ($\beta _{A}^{P}$ corresponding edge $AC$):
	\[\beta _{A}^{P}=\frac{\alpha _{A}^{{{P}_{C}}}\left( 1-\alpha _{C}^{{{P}_{A}}} \right)}{1-\alpha _{A}^{{{P}_{C}}}\alpha _{C}^{{{P}_{A}}}}.\]
	
	The exchange of $A$ and $C$ in the above formula yields the following results:
	\[\beta _{C}^{P}=\frac{\alpha _{C}^{{{P}_{A}}}\left( 1-\alpha _{A}^{{{P}_{C}}} \right)}{1-\alpha _{C}^{{{P}_{A}}}\alpha _{A}^{{{P}_{C}}}}.\]
	
	So we get the system of equations:
	\begin{equation}\label{Eq20.2.5}
		\left\{ \begin{aligned}
			& \beta _{A}^{P}=\frac{\alpha _{A}^{{{P}_{C}}}\left( 1-\alpha _{C}^{{{P}_{A}}} \right)}{1-\alpha _{A}^{{{P}_{C}}}\alpha _{C}^{{{P}_{A}}}} \\ 
			& \beta _{C}^{P}=\frac{\alpha _{C}^{{{P}_{A}}}\left( 1-\alpha _{A}^{{{P}_{C}}} \right)}{1-\alpha _{C}^{{{P}_{A}}}\alpha _{A}^{{{P}_{C}}}}. \\ 
		\end{aligned} \right.
	\end{equation}
	
	The above system is obtained by replacing $B$ with $C$ of the system (\ref{Eq20.2.1}). Therefore, the solution of the equation system (\ref{Eq20.2.5}) is obtained after the solution $B$ of the equation system (\ref{Eq20.2.1}) is replaced by $C$, that is, according to the formula (\ref{Eq20.2.4}), $B$ is replaced by $C$ to obtain:
	\[\alpha _{A}^{{{P}_{C}}}=\frac{\beta _{A}^{P}}{1-\beta _{C}^{P}},\ \alpha _{C}^{{{P}_{A}}}=\frac{\beta _{C}^{P}}{1-\beta _{A}^{P}}.\]	
	
	Similarly:
	\[\alpha _{A}^{{{P}_{D}}}=\frac{\beta _{A}^{P}}{1-\beta _{D}^{P}},\ \alpha _{D}^{{{P}_{A}}}=\frac{\beta _{D}^{P}}{1-\beta _{A}^{P}};\]	
	\[\alpha _{B}^{{{P}_{C}}}=\frac{\beta _{B}^{P}}{1-\beta _{C}^{P}},\ \alpha _{C}^{{{P}_{B}}}=\frac{\beta _{C}^{P}}{1-\beta _{B}^{P}};\]	
	\[\alpha _{B}^{{{P}_{D}}}=\frac{\beta _{B}^{P}}{1-\beta _{D}^{P}},\ \alpha _{D}^{{{P}_{B}}}=\frac{\beta _{D}^{P}}{1-\beta _{B}^{P}};\]	
	\[\alpha _{C}^{{{P}_{D}}}=\frac{\beta _{C}^{P}}{1-\beta _{D}^{P}},\ \alpha _{D}^{{{P}_{C}}}=\frac{\beta _{D}^{P}}{1-\beta _{C}^{P}}.\]
\end{proof}
\hfill $\square$\par

\begin{theorem}{Harmonic expression 2 of frame components of vertex at IC-F, Daiyuan Zhang}{Thm20.2.5}\label{Thm20.2.5} 
	\[\frac{1}{\alpha _{A}^{{{P}_{B}}}}+\frac{1}{\alpha _{A}^{{{P}_{C}}}}+\frac{1}{\alpha _{A}^{{{P}_{D}}}}-\frac{2}{\beta _{A}^{P}}=1,\]
	\[\frac{1}{\alpha _{B}^{{{P}_{C}}}}+\frac{1}{\alpha _{B}^{{{P}_{D}}}}+\frac{1}{\alpha _{B}^{{{P}_{A}}}}-\frac{2}{\beta _{B}^{P}}=1,\]
	\[\frac{1}{\alpha _{C}^{{{P}_{D}}}}+\frac{1}{\alpha _{C}^{{{P}_{A}}}}+\frac{1}{\alpha _{C}^{{{P}_{B}}}}-\frac{2}{\beta _{C}^{P}}=1,\]
	\[\frac{1}{\alpha _{D}^{{{P}_{A}}}}+\frac{1}{\alpha _{D}^{{{P}_{B}}}}+\frac{1}{\alpha _{D}^{{{P}_{C}}}}-\frac{2}{\beta _{D}^{P}}=1.\]
\end{theorem}

\begin{proof}
	Only the first formula will be proved, and other formulas can be proved similarly. According to theorem \ref{thm:Thm20.2.4}, we can get:
	\[\begin{aligned}
		\frac{1}{\alpha _{A}^{{{P}_{B}}}}+\frac{1}{\alpha _{A}^{{{P}_{C}}}}+\frac{1}{\alpha _{A}^{{{P}_{D}}}}& =\frac{1-\beta _{B}^{P}}{\beta _{A}^{P}}+\frac{1-\beta _{C}^{P}}{\beta _{A}^{P}}+\frac{1-\beta _{D}^{P}}{\beta _{A}^{P}} \\ 
		& =\frac{3-\beta _{B}^{P}-\beta _{C}^{P}-\beta _{D}^{P}}{\beta _{A}^{P}}=\frac{\beta _{A}^{P}+2}{\beta _{A}^{P}}=1+\frac{2}{\beta _{A}^{P}}.  
	\end{aligned}\]
	
	Therefore 
	\[\frac{1}{\alpha _{A}^{{{P}_{B}}}}+\frac{1}{\alpha _{A}^{{{P}_{C}}}}+\frac{1}{\alpha _{A}^{{{P}_{D}}}}-\frac{2}{\beta _{A}^{P}}=1.\]
	
	Similarly, other formulas can be obtained.
\end{proof}
\hfill $\square$\par

\section{Relationship of frame components between triangle and line segment}\label{Sec20.3}

Firstly, we discuss the the frame components of line segment,  from theorem \ref{thm:Thm3.2.1}, we get:
\[\overrightarrow{OP}={{\alpha }_{A}}\overrightarrow{OA}+{{\alpha }_{B}}\overrightarrow{OB},\]	
\[{{\alpha }_{A}}+{{\alpha }_{B}}=1.\]	


The frame is composed of basis vectors $\overrightarrow{OA}$ and $\overrightarrow{OB}$, which is called \textbf{segmental frame} or \textbf{frame of line segment}, denoted as $\left( O;A,B \right)$. ${{\alpha }_{A}}$ and ${{\alpha }_{B}}$ are called the frame components of line segment on vectors $\overrightarrow{OA}$ and $\overrightarrow{OB}$ respectively, or referred to as \textbf{frame components of line segment} or \textbf{segmental frame components}.

We have discussed the concepts of intersecting center (abbreviated as IC), which is also known as intersecting center of triangle, and intersecting center of tetrahedron (abbreviated as IC-T). Now we introduce the concept of intersecting center of line segment (abbreviated as IC-L).

\textbf{Intersecting center of line segment}: For a given line segment $AB$ and a point $P$, if $P\in \overleftrightarrow{AB}$, the point $P$ is called the intersecting center of line segment $AB$, which is abbreviated as IC-L.

Obviously, each IC-T corresponds to four ICs-Fs, and each IC-F corresponds to three ICs-Ls.

\textbf{Intersecting ratio of line segment (abbreviated as IR-L)}: Assuming that the points $A$, $B$, $P$ are collinear, and $P$ is the IC-L of $AB$, then the following fractional ratio 
\[{{\lambda }_{AB}}=\frac{\overrightarrow{AP}}{\overrightarrow{PB}}\] 	
is called the intersecting ratio of line segment (abbreviated as IR-L).

\begin{theorem}{Relationship of frame components between triangle and line segment, Daiyuan Zhang}{Thm20.3.1}\label{Thm20.3.1} 
	If the IC $P\in \overleftrightarrow{AB}$, then the triangular frame $\left( O;A,B,C \right)$ will degenerate into the segmental frame $\left( O;A,B \right)$, and
	\[\alpha _{A}^{P}={{\alpha }_{A}},\ \alpha _{B}^{P}={{\alpha }_{B}},\ {{\alpha }_{C}^{P}}=0.\]
	
	If the IC $P\in \overleftrightarrow{BC}$, then the triangular frame $\left( O;A,B,C \right)$ will degenerate into the segmental frame $\left( O;B,C \right)$, and
	\[{{\alpha }_{A}^{P}}=0,\ \alpha _{B}^{P}={{\alpha }_{B}},\ \alpha _{C}^{P}={{\alpha }_{C}}.\]
	
	If the IC $P\in \overleftrightarrow{CA}$, then the triangular frame $\left( O;A,B,C \right)$ will degenerate into the segmental frame $\left( O;C,A \right)$, and
	\[\alpha _{A}^{P}={{\alpha }_{A}},\ {{\alpha }_{B}^{P}}=0,\ \alpha _{C}^{P}={{\alpha }_{C}}.\]
\end{theorem}

\begin{proof}
	
	Only the case of $P\in \overleftrightarrow{AB}$ will be proved, and other cases can be proved similarly.
	
	As shown in figure \ref{fig:tu6.1.1}, if point $P\in \overleftrightarrow{AB}$, for the triangular frame $\left( O;A,B,C \right)$, according to theorem \ref{thm:Thm6.1.1}, we have:	
	\[\overrightarrow{OP}=\alpha _{A}^{P}\overrightarrow{OA}+\alpha _{B}^{P}\overrightarrow{OB}+\alpha _{C}^{P}\overrightarrow{OC},\]
	\[\alpha _{A}^{P}+\alpha _{B}^{P}+\alpha _{C}^{P}=1.\]	
	
	For the the segmental frame $\left( O;A,B \right)$, according to theorem \ref{thm:Thm3.2.1}, we have:	
	\[\overrightarrow{OP}={{\alpha }_{A}}\overrightarrow{OA}+{{\alpha }_{B}}\overrightarrow{OB},\]	
	\[{{\alpha }_{A}}+{{\alpha }_{B}}=1.\]
	
	Therefore 
	\[\left( \alpha _{A}^{P}-{{\alpha }_{A}} \right)\overrightarrow{OA}+\left( \alpha _{B}^{P}-{{\alpha }_{B}} \right)\overrightarrow{OB}+\left( \alpha _{C}^{P}-0 \right)\overrightarrow{OC}=\overrightarrow{0},\]
	\[\left( \alpha _{A}^{P}-{{\alpha }_{A}} \right)+\left( \alpha _{B}^{P}-{{\alpha }_{B}} \right)+\alpha _{C}^{P}=\left( \alpha _{A}^{P}+\alpha _{B}^{P}+\alpha _{C}^{P} \right)-\left( {{\alpha }_{A}}+{{\alpha }_{B}} \right)=0.\]
	
	Apply theorem \ref{thm:Thm3.4.3} to obtain:
	\[\alpha _{A}^{P}-{{\alpha }_{A}}=0,\ \alpha _{B}^{P}-{{\alpha }_{B}}=0,\ {{\alpha }_{C}^{P}}=0,\]
	i.e. 
	\[\alpha _{A}^{P}={{\alpha }_{A}},\ \alpha _{B}^{P}={{\alpha }_{B}},\ {{\alpha }_{C}^{P}}=0.\]
\end{proof}
\hfill $\square$\par
The above theorem shows that if the IC $P$ of a triangle is located on the straight line of the opposite side of a vertex of the triangle, then the triangular frame component of the vertex is 0, and the triangular frame components of the other vertices of the triangle are equal to the corresponding segmental frame components of the segment (side) opposite the vertex.

\section{Relationship of frame components between tetrahedron and triangle}\label{Sec20.4}
A tetrahedron has four triangles, each of which can have its frame components. What is the relationship between the frame components of these triangles and the tetrahedron? This question will be discussed in this section.

\begin{theorem}{Relationship of frame components between tetrahedron and triangle, Daiyuan Zhang}{Thm20.4.1}\label{Thm20.4.1} 
	Assuming $O\in {{\mathbb{R}}^{3}}$, we have the following conclusions:
	
	1. If $P$ is an IC-T, $P\in {{\pi}_{ ABC}}$, then the tetrahedral frame $\left( O;A,B,C,D \right)$ will degenerate into the triangular frame $\left( O;A,B,C \right)$, and
	\[\beta _{A}^{P}=\alpha _{A}^{P},\ \beta _{B}^{P}=\alpha _{B}^{P},\ \beta _{C}^{P}=\alpha _{C}^{P},\ \beta _{D}^{P}=0.\ \]
	
	2. If $P$ is an IC-T, $P\in {{\pi}_{ BCD}}$, then the tetrahedral frame $\left( O;A,B,C,D \right)$ will degenerate into the triangular frame $\left( O;B,C,D \right)$, and	
	\[\beta _{A}^{P}=0,\ \beta _{B}^{P}=\alpha _{B}^{P},\ \beta _{C}^{P}=\alpha _{C}^{P},\ \beta _{D}^{P}=\alpha _{D}^{P}.\ \]
	
	3. If $P$ is an IC-T, $P\in {{\pi}_{ CDA}}$, then the tetrahedral frame $\left( O;A,B,C,D \right)$ will degenerate into the triangular frame $\left( O;C,D,A \right)$, and	
	\[\beta _{A}^{P}=\alpha _{A}^{P},\ \beta _{B}^{P}=0,\ \beta _{C}^{P}=\alpha _{C}^{P},\ \beta _{D}^{P}=\alpha _{D}^{P}.\ \]
	
	4. If $P$ is an IC-T, $P\in {{\pi}_{ DAB}}$, then the tetrahedral frame $\left( O;A,B,C,D \right)$ will degenerate into the triangular frame $\left( O;D,A,B \right)$, and	
	\[\beta _{A}^{P}=\alpha _{A}^{P},\ \beta _{B}^{P}=\alpha _{B}^{P},\ \beta _{C}^{P}=0,\ \beta _{D}^{P}=\alpha _{D}^{P}.\ \]
\end{theorem}

\begin{proof}
	In the following, only the case of $P\in {{\pi}_{ ABC}}$ will be proved, and other cases can be proved similarly.
	
	As shown in figure \ref{fig:tu6.1.1}, if point $P\in {{\pi}_{ ABC}}$, for the tetrahedral frame $\left( O;A,B,C,D \right)$, according to theorem \ref{thm:Thm18.1.4}, we have:	
	\[\overrightarrow{OP}=\beta _{A}^{P}\overrightarrow{OA}+\beta _{B}^{P}\overrightarrow{OB}+\beta _{C}^{P}\overrightarrow{OC}+\beta _{D}^{P}\overrightarrow{OD},\]
	\[\beta _{A}^{P}+\beta _{B}^{P}+\beta _{C}^{P}+\beta _{D}^{P}=1.\]	
	
	For the triangular frame $\left( O;A,B \right)$, according to theorem \ref{thm:Thm6.1.1}, we have:	
	\[\overrightarrow{OP}=\alpha _{A}^{P}\overrightarrow{OA}+\alpha _{B}^{P}\overrightarrow{OB}+\alpha _{C}^{P}\overrightarrow{OC},\]
	\[\alpha _{A}^{P}+\alpha _{B}^{P}+\alpha _{C}^{P}=1.\]
	
	Therefore 
	\[\left( \beta _{A}^{P}-\alpha _{A}^{P} \right)\overrightarrow{OA}+\left( \beta _{B}^{P}-\alpha _{B}^{P} \right)\overrightarrow{OB}+\left( \beta _{C}^{P}-\alpha _{C}^{P} \right)\overrightarrow{OC}+\left( \beta _{C}^{P}-\alpha _{C}^{P} \right)\overrightarrow{OD}=\overrightarrow{0},\]
	\[\begin{aligned}
		& \left( \beta _{A}^{P}-\alpha _{A}^{P} \right)+\left( \beta _{B}^{P}-\alpha _{B}^{P} \right)+\left( \beta _{C}^{P}-\alpha _{C}^{P} \right)+\beta _{D}^{P} \\ 
		& =\left( \beta _{A}^{P}+\beta _{B}^{P}+\beta _{C}^{P}+\beta _{D}^{P} \right)-\left( \alpha _{A}^{P}+\alpha _{B}^{P}+\alpha _{C}^{P} \right)=1-1=0. \\ 
	\end{aligned}\]
	
	Apply theorem \ref{thm:Thm18.1.2} to obtain:
	\[\beta _{A}^{P}-\alpha _{A}^{P}=0,\ \beta _{B}^{P}-\alpha _{B}^{P}=0,\ \beta _{C}^{P}-\alpha _{C}^{P}=0,\ \beta _{D}^{P}=0,\]
	i.e. 
	\[\beta _{A}^{P}=\alpha _{A}^{P},\ \beta _{B}^{P}=\alpha _{B}^{P},\ \beta _{C}^{P}=\alpha _{C}^{P},\ \beta _{D}^{P}=0.\]
\end{proof}
\hfill $\square$\par
The above theorem shows that if the IC-T $P$ is located on the plane of a triangle, then the tetrahedral frame component of the vertex opposite the triangle is 0, and the tetrahedral frame components of other vertices are equal to the frame components of the opposite triangle on the triangular frame.


\section{Application of the frame component of tetrahedra in analytical geometry}
In analytic geometry, the coordinates of other points are often calculated based on the coordinates of the four vertices of a given tetrahedron. I have provided a general method using Intercenter Geometry in this section.
\begin{theorem}{Mixed formula of frame component and Cartesian Coordinate, Daiyuan Zhang}\label{SimiantiDeBiaojiafenliangYuZhijiaozuobiaoDeHunhegongshi}
	Given a tetrahedron $ABCD$, $P$ is a given point in space, and the frame components of point $P$ are $\beta _{A}^{P}$, $\beta _{B}^{P}$, $\beta _{C}^{P}$, $\beta _{D}^{P}$, respectively. Establish a Cartesian coordinate system in space, with the Cartesian coordinates of the four vertices of the tetrahedron $ABCD$ are $A\left( {{x}_{A}},{{y}_{A}},{{z}_{A}} \right)$, $B\left( {{x}_{B}},{{y}_{B}},{{z}_{B}} \right)$, $C\left( {{x}_{C}},{{y}_{C}},{{z}_{C}} \right)$, $D\left( {{x}_{D}},{{y}_{D}},{{z}_{D}} \right)$,  respectively. The Cartesian coordinates of point $P$ are $P\left( {{x}_{P}},{{y}_{P}},{{z}_{P}} \right)$. Then:
	\[{{x}_{P}}=\beta _{A}^{P}{{x}_{A}}+\beta _{B}^{P}{{x}_{B}}+\beta _{C}^{P}{{x}_{C}}+\beta _{D}^{P}{{x}_{D}},\]
	\[{{y}_{P}}=\beta _{A}^{P}{{y}_{A}}+\beta _{B}^{P}{{y}_{B}}+\beta _{C}^{P}{{y}_{C}}+\beta _{D}^{P}{{y}_{D}},\]
	\[{{z}_{P}}=\beta _{A}^{P}{{z}_{A}}+\beta _{B}^{P}{{z}_{B}}+\beta _{C}^{P}{{z}_{C}}+\beta _{D}^{P}{{z}_{D}}.\]
\end{theorem}

\begin{proof}
	Coinciding the origin of the frame with the origin of the Cartesian coordinate system yields:
	\[\overrightarrow{OP}=\beta _{A}^{P}\overrightarrow{OA}+\beta _{B}^{P}\overrightarrow{OB}+\beta _{C}^{P}\overrightarrow{OC}+\beta _{D}^{P}\overrightarrow{OD},\]
	\[\overrightarrow{OA}={{x}_{A}}{{\bm{e}}_{x}}+{{y}_{A}}{{\bm{e}}_{y}}+{{z}_{A}}{{\bm{e}}_{z}},\]
	\[\overrightarrow{OB}={{x}_{B}}{{\bm{e}}_{x}}+{{y}_{B}}{{\bm{e}}_{y}}+{{z}_{B}}{{\bm{e}}_{z}},\]
	\[\overrightarrow{OC}={{x}_{C}}{{\bm{e}}_{x}}+{{y}_{C}}{{\bm{e}}_{y}}+{{z}_{C}}{{\bm{e}}_{z}},\]
	\[\overrightarrow{OD}={{x}_{D}}{{\bm{e}}_{x}}+{{y}_{D}}{{\bm{e}}_{y}}+{{z}_{D}}{{\bm{e}}_{z}}.\]
	
	Where ${{\bm{e}}_{x}}$, ${{\bm{e}}_{y}}$ and ${{\bm{e}}_{z}}$ are the unit vectors of the $OX$ axis, $OY$ axis and $OZ$ axis in Cartesian coordinates, respectively. Therefore, the formula in Cartesian coordinates is obtained:
	\begin{align*}
		\overrightarrow{OP}& =\beta _{A}^{P}\overrightarrow{OA}+\beta _{B}^{P}\overrightarrow{OB}+\beta _{C}^{P}\overrightarrow{OC}+\beta _{D}^{P}\overrightarrow{OD} \\ 
		& =\beta _{A}^{P}\left( {{x}_{A}}{{\bm{e}}_{x}}+{{y}_{A}}{{\bm{e}}_{y}}+{{z}_{A}}{{\bm{e}}_{z}} \right)+\beta _{B}^{P}\left( {{x}_{B}}{{\bm{e}}_{x}}+{{y}_{B}}{{\bm{e}}_{y}}+{{z}_{B}}{{\bm{e}}_{z}} \right) \\ 
		& +\beta _{C}^{P}\left( {{x}_{C}}{{\bm{e}}_{x}}+{{y}_{C}}{{\bm{e}}_{y}}+{{z}_{C}}{{\bm{e}}_{z}} \right)+\beta _{D}^{P}\left( {{x}_{D}}{{\bm{e}}_{x}}+{{y}_{D}}{{\bm{e}}_{y}}+{{z}_{D}}{{\bm{e}}_{z}} \right) \\ 
		& =\left( \beta _{A}^{P}{{x}_{A}}+\beta _{B}^{P}{{x}_{B}}+\beta _{C}^{P}{{x}_{C}}+\beta _{D}^{P}{{x}_{D}} \right){{\bm{e}}_{x}} \\ 
		& +\left( \beta _{A}^{P}{{y}_{A}}+\beta _{B}^{P}{{y}_{B}}+\beta _{C}^{P}{{y}_{C}}+\beta _{D}^{P}{{y}_{D}} \right){{\bm{e}}_{y}} \\ 
		& +\left( \beta _{A}^{P}{{z}_{A}}+\beta _{B}^{P}{{z}_{B}}+\beta _{C}^{P}{{z}_{C}}+\beta _{D}^{P}{{z}_{D}} \right){{\bm{e}}_{z}}.  
	\end{align*}
	
	On the other hand,
	\[\overrightarrow{OP}={{x}_{P}}{{\bm{e}}_{x}}+{{y}_{P}}{{\bm{e}}_{y}}+{{z}_{P}}{{\bm{e}}_{z}}.\]
	
	Due to the linear independence of ${{\bm{e}}_{x}}$, ${{\bm{e}}_{y}}$ and ${{\bm{e}}_{z}}$, we obtain:
	\[{{x}_{P}}=\beta _{A}^{P}{{x}_{A}}+\beta _{B}^{P}{{x}_{B}}+\beta _{C}^{P}{{x}_{C}}+\beta _{D}^{P}{{x}_{D}},\]
	\[{{y}_{P}}=\beta _{A}^{P}{{y}_{A}}+\beta _{B}^{P}{{y}_{B}}+\beta _{C}^{P}{{y}_{C}}+\beta _{D}^{P}{{y}_{D}},\]
	\[{{z}_{P}}=\beta _{A}^{P}{{z}_{A}}+\beta _{B}^{P}{{z}_{B}}+\beta _{C}^{P}{{z}_{C}}+\beta _{D}^{P}{{z}_{D}}.\]
\end{proof}
\hfill $\square$\par

The mixed formula given by the above theorem is symmetric and elegant, with concise expression. I think it is also a good result.

\begin{example}{}\label{YizhiSimiantiDeDingdianzuobiaoQiuZhongxinHeNeixinDeZhijiaozuobiao}
	Given the Cartesian coordinates $A\left( {{x}_{A}},{{y}_{A}},{{z}_{A}} \right)$, $B\left( {{x}_{B}},{{y}_{B}},{{z}_{B}} \right)$, $C\left( {{x}_{C}},{{y}_{C}},{{z}_{C}} \right)$, $D\left( {{x}_{D}},{{y}_{D}},{{z}_{D}} \right)$ of the four vertices of the tetrahedron $ABCD$, respectively. Find  the Cartesian coordinates of the centroid and incenter.
\end{example}

\begin{solution}
	Using the above theorem, for the centroid $G$, we obtain:
	\[{{x}_{G}}=\beta _{A}^{G}{{x}_{A}}+\beta _{B}^{G}{{x}_{B}}+\beta _{C}^{G}{{x}_{C}}+\beta _{D}^{G}{{x}_{D}}=\frac{1}{4}\left( {{x}_{A}}+{{x}_{B}}+{{x}_{C}}+{{x}_{D}} \right),\]
	\[{{y}_{G}}=\beta _{A}^{P}{{y}_{A}}+\beta _{B}^{P}{{y}_{B}}+\beta _{C}^{P}{{y}_{C}}+\beta _{D}^{P}{{y}_{D}}=\frac{1}{4}\left( {{y}_{A}}+{{y}_{B}}+{{y}_{C}}+{{y}_{D}} \right).\]
	\[{{z}_{G}}=\beta _{A}^{P}{{z}_{A}}+\beta _{B}^{P}{{z}_{B}}+\beta _{C}^{P}{{z}_{C}}+\beta _{D}^{P}{{z}_{D}}=\frac{1}{4}\left( {{z}_{A}}+{{z}_{B}}+{{z}_{C}}+{{z}_{D}} \right).\]
	
	For the incenter $I$, we obtain:
	\[{{x}_{I}}=\beta _{A}^{I}{{x}_{A}}+\beta _{B}^{I}{{x}_{B}}+\beta _{C}^{I}{{x}_{C}}+\beta _{D}^{I}{{x}_{D}}=\frac{{{S}^{A}}{{x}_{A}}+{{S}^{B}}{{x}_{B}}+{{S}^{C}}{{x}_{C}}+{{S}^{D}}{{x}_{D}}}{{{S}^{A}}+{{S}^{B}}+{{S}^{C}}+{{S}^{D}}},\]
	\[{{y}_{I}}=\beta _{A}^{I}{{y}_{A}}+\beta _{B}^{I}{{y}_{B}}+\beta _{C}^{I}{{y}_{C}}+\beta _{D}^{I}{{y}_{D}}=\frac{{{S}^{A}}{{y}_{A}}+{{S}^{B}}{{y}_{B}}+{{S}^{C}}{{y}_{C}}+{{S}^{D}}{{y}_{D}}}{{{S}^{A}}+{{S}^{B}}+{{S}^{C}}+{{S}^{D}}},\]
	\[{{z}_{I}}=\beta _{A}^{I}{{z}_{A}}+\beta _{B}^{I}{{z}_{B}}+\beta _{C}^{I}{{z}_{C}}+\beta _{D}^{I}{{z}_{D}}=\frac{{{S}^{A}}{{z}_{A}}+{{S}^{B}}{{z}_{B}}+{{S}^{C}}{{z}_{C}}+{{S}^{D}}{{z}_{D}}}{{{S}^{A}}+{{S}^{B}}+{{S}^{C}}+{{S}^{D}}}.\]
\end{solution}
\hfill $\diamond$\par

\section{Frame transformation}
This section studies the frame transformation of tetrahedra.

\begin{theorem}{Frame transformation formula 1, Daiyuan Zhang}{Biaojiabianhuangongshi1}\label{Biaojiabianhuangongshi1}
	Given two tetrahedra ${{A}_{1}}{{B}_{1}}{{C}_{1}}{{D}_{1}}$ and ${{A}_{2}}{{B}_{2}}{{C}_{2}}{{D}_{2}}$. The frame components on $\left( O;{{A}_{1}},{{B}_{1}},{{C}_{1}},{{D}_{1}} \right)$ of the four vertices ${{A}_{2}}$, ${{B}_{2}}$, ${{C}_{2}}$, ${{D}_{2}}$ of the tetrahedron ${{A}_{2}}{{B}_{2}}{{C}_{2}}{{D}_{2}}$ are $\beta _{{{A}_{1}}}^{{{A}_{2}}}$, $\beta _{{{B}_{1}}}^{{{A}_{2}}}$, $\beta _{{{C}_{1}}}^{{{A}_{2}}}$, $\beta _{{{D}_{1}}}^{{{A}_{2}}}$; $\beta _{{{A}_{1}}}^{{{B}_{2}}}$, $\beta _{{{B}_{1}}}^{{{B}_{2}}}$, $\beta _{{{C}_{1}}}^{{{B}_{2}}}$, $\beta _{{{D}_{1}}}^{{{B}_{2}}}$; $\beta _{{{A}_{1}}}^{{{C}_{2}}}$, $\beta _{{{B}_{1}}}^{{{C}_{2}}}$, $\beta _{{{C}_{1}}}^{{{C}_{2}}}$, $\beta _{{{D}_{1}}}^{{{C}_{2}}}$; $\beta _{{{A}_{1}}}^{{{D}_{2}}}$, $\beta _{{{B}_{1}}}^{{{D}_{2}}}$, $\beta _{{{C}_{1}}}^{{{D}_{2}}}$, $\beta _{{{D}_{1}}}^{{{D}_{2}}}$; respectively. Then
	\[\left( \begin{aligned}
		& \begin{matrix}
			\overrightarrow{O{{A}_{2}}}  \\
			\overrightarrow{O{{B}_{2}}}  \\
			\overrightarrow{O{{C}_{2}}}  \\
		\end{matrix} \\ 
		& \overrightarrow{O{{D}_{2}}} \\ 
	\end{aligned} \right)=\left( \begin{matrix}
		\beta _{{{A}_{1}}}^{{{A}_{2}}} & \beta _{{{B}_{1}}}^{{{A}_{2}}} & \beta _{{{C}_{1}}}^{{{A}_{2}}} & \beta _{{{D}_{1}}}^{{{A}_{2}}}  \\
		\beta _{{{A}_{1}}}^{{{B}_{2}}} & \beta _{{{B}_{1}}}^{{{B}_{2}}} & \beta _{{{C}_{1}}}^{{{B}_{2}}} & \beta _{{{D}_{1}}}^{{{B}_{2}}}  \\
		\beta _{{{A}_{1}}}^{{{C}_{2}}} & \beta _{{{B}_{1}}}^{{{C}_{2}}} & \beta _{{{C}_{1}}}^{{{C}_{2}}} & \beta _{{{D}_{1}}}^{{{C}_{2}}}  \\
		\beta _{{{A}_{1}}}^{{{D}_{2}}} & \beta _{{{B}_{1}}}^{{{D}_{2}}} & \beta _{{{C}_{1}}}^{{{D}_{2}}} & \beta _{{{D}_{1}}}^{{{D}_{2}}}  \\
	\end{matrix} \right)\left( \begin{aligned}
		& \begin{matrix}
			\overrightarrow{O{{A}_{1}}}  \\
			\overrightarrow{O{{B}_{1}}}  \\
			\overrightarrow{O{{C}_{1}}}  \\
		\end{matrix} \\ 
		& \overrightarrow{O{{D}_{1}}} \\ 
	\end{aligned} \right).\]
	
	The sum of elements in each row of the matrix is 1.
\end{theorem}

\begin{proof}
	Obviously,
	\[\overrightarrow{O{{A}_{2}}}=\beta _{{{A}_{1}}}^{{{A}_{2}}}\overrightarrow{O{{A}_{1}}}+\beta _{{{B}_{1}}}^{{{A}_{2}}}\overrightarrow{O{{B}_{1}}}+\beta _{{{C}_{1}}}^{{{A}_{2}}}\overrightarrow{O{{C}_{1}}}+\beta _{{{D}_{1}}}^{{{A}_{2}}}\overrightarrow{O{{D}_{1}}},\]
	\[\beta _{{{A}_{1}}}^{{{A}_{2}}}+\beta _{{{B}_{1}}}^{{{A}_{2}}}+\beta _{{{C}_{1}}}^{{{A}_{2}}}+\beta _{{{D}_{1}}}^{{{A}_{2}}}=1;\]
	\[\overrightarrow{O{{B}_{2}}}=\beta _{{{A}_{1}}}^{{{B}_{2}}}\overrightarrow{O{{A}_{1}}}+\beta _{{{B}_{1}}}^{{{B}_{2}}}\overrightarrow{O{{B}_{1}}}+\beta _{{{C}_{1}}}^{{{B}_{2}}}\overrightarrow{O{{C}_{1}}}+\beta _{{{D}_{1}}}^{{{B}_{2}}}\overrightarrow{O{{D}_{1}}},\]
	\[\beta _{{{A}_{1}}}^{{{B}_{2}}}+\beta _{{{B}_{1}}}^{{{B}_{2}}}+\beta _{{{C}_{1}}}^{{{B}_{2}}}+\beta _{{{D}_{1}}}^{{{B}_{2}}}=1;\]
	\[\overrightarrow{O{{C}_{2}}}=\beta _{{{A}_{1}}}^{{{C}_{2}}}\overrightarrow{O{{A}_{1}}}+\beta _{{{B}_{1}}}^{{{C}_{2}}}\overrightarrow{O{{B}_{1}}}+\beta _{{{C}_{1}}}^{{{C}_{2}}}\overrightarrow{O{{C}_{1}}}+\beta _{{{D}_{1}}}^{{{C}_{2}}}\overrightarrow{O{{D}_{1}}},\]
	\[\beta _{{{A}_{1}}}^{{{C}_{2}}}+\beta _{{{B}_{1}}}^{{{C}_{2}}}+\beta _{{{C}_{1}}}^{{{C}_{2}}}+\beta _{{{D}_{1}}}^{{{C}_{2}}}=1;\]
	\[\overrightarrow{O{{D}_{2}}}=\beta _{{{A}_{1}}}^{{{D}_{2}}}\overrightarrow{O{{A}_{1}}}+\beta _{{{B}_{1}}}^{{{D}_{2}}}\overrightarrow{O{{B}_{1}}}+\beta _{{{C}_{1}}}^{{{D}_{2}}}\overrightarrow{O{{C}_{1}}}+\beta _{{{D}_{1}}}^{{{D}_{2}}}\overrightarrow{O{{D}_{1}}},\]
	\[\beta _{{{A}_{1}}}^{{{D}_{2}}}+\beta _{{{B}_{1}}}^{{{D}_{2}}}+\beta _{{{C}_{1}}}^{{{D}_{2}}}+\beta _{{{D}_{1}}}^{{{D}_{2}}}=1.\]
	
	Written in matrix form, it is
	\[\left( \begin{aligned}
		& \begin{matrix}
			\overrightarrow{O{{A}_{2}}}  \\
			\overrightarrow{O{{B}_{2}}}  \\
			\overrightarrow{O{{C}_{2}}}  \\
		\end{matrix} \\ 
		& \overrightarrow{O{{D}_{2}}} \\ 
	\end{aligned} \right)=\left( \begin{matrix}
		\beta _{{{A}_{1}}}^{{{A}_{2}}} & \beta _{{{B}_{1}}}^{{{A}_{2}}} & \beta _{{{C}_{1}}}^{{{A}_{2}}} & \beta _{{{D}_{1}}}^{{{A}_{2}}}  \\
		\beta _{{{A}_{1}}}^{{{B}_{2}}} & \beta _{{{B}_{1}}}^{{{B}_{2}}} & \beta _{{{C}_{1}}}^{{{B}_{2}}} & \beta _{{{D}_{1}}}^{{{B}_{2}}}  \\
		\beta _{{{A}_{1}}}^{{{C}_{2}}} & \beta _{{{B}_{1}}}^{{{C}_{2}}} & \beta _{{{C}_{1}}}^{{{C}_{2}}} & \beta _{{{D}_{1}}}^{{{C}_{2}}}  \\
		\beta _{{{A}_{1}}}^{{{D}_{2}}} & \beta _{{{B}_{1}}}^{{{D}_{2}}} & \beta _{{{C}_{1}}}^{{{D}_{2}}} & \beta _{{{D}_{1}}}^{{{D}_{2}}}  \\
	\end{matrix} \right)\left( \begin{aligned}
		& \begin{matrix}
			\overrightarrow{O{{A}_{1}}}  \\
			\overrightarrow{O{{B}_{1}}}  \\
			\overrightarrow{O{{C}_{1}}}  \\
		\end{matrix} \\ 
		& \overrightarrow{O{{D}_{1}}} \\ 
	\end{aligned} \right).\]
\end{proof}
\hfill $\square$\par

The above equation can be written as
\[{{\bm{F}}_{2}}={{\bm{B}}_{21}}{{\bm{F}}_{1}}.\]

Where
\[{{\bm{F}}_{1}}=\left( \begin{aligned}
	& \begin{matrix}
		\overrightarrow{O{{A}_{1}}}  \\
		\overrightarrow{O{{B}_{1}}}  \\
		\overrightarrow{O{{C}_{1}}}  \\
	\end{matrix} \\ 
	& \overrightarrow{O{{D}_{1}}} \\ 
\end{aligned} \right)\text{,}\quad {{\bm{F}}_{2}}=\left( \begin{aligned}
	& \begin{matrix}
		\overrightarrow{O{{A}_{2}}}  \\
		\overrightarrow{O{{B}_{2}}}  \\
		\overrightarrow{O{{C}_{2}}}  \\
	\end{matrix} \\ 
	& \overrightarrow{O{{D}_{2}}} \\ 
\end{aligned} \right)\text{,}\]
\[{{\bm{B}}_{21}}=\left( \begin{matrix}
	\beta _{{{A}_{1}}}^{{{A}_{2}}} & \beta _{{{B}_{1}}}^{{{A}_{2}}} & \beta _{{{C}_{1}}}^{{{A}_{2}}} & \beta _{{{D}_{1}}}^{{{A}_{2}}}  \\
	\beta _{{{A}_{1}}}^{{{B}_{2}}} & \beta _{{{B}_{1}}}^{{{B}_{2}}} & \beta _{{{C}_{1}}}^{{{B}_{2}}} & \beta _{{{D}_{1}}}^{{{B}_{2}}}  \\
	\beta _{{{A}_{1}}}^{{{C}_{2}}} & \beta _{{{B}_{1}}}^{{{C}_{2}}} & \beta _{{{C}_{1}}}^{{{C}_{2}}} & \beta _{{{D}_{1}}}^{{{C}_{2}}}  \\
	\beta _{{{A}_{1}}}^{{{D}_{2}}} & \beta _{{{B}_{1}}}^{{{D}_{2}}} & \beta _{{{C}_{1}}}^{{{D}_{2}}} & \beta _{{{D}_{1}}}^{{{D}_{2}}}  \\
\end{matrix} \right).\]

Among them, ${{\bm{F}}_{1}}$ and ${{\bm{F}}_{2}}$ are called \textbf{frame vectors}, and ${{\bm{B}}_{21}}$ is called the \textbf{transformation matrix} from ${{\bm{F}}_{1}}$ to ${{\bm{F}}_{2}}$. It is abbreviated as a transformation matrix.


\begin{theorem}{Frame transformation formula 2, Daiyuan Zhang}{Biaojiabianhuangongshi2}\label{Biaojiabianhuangongshi2}
	Given $n$ tetrahedra ${{A}_{1}}{{B}_{1}}{{C}_{1}}{{D}_{1}}$, ${{A}_{2}}{{B}_{2}}{{C}_{2}}{{D}_{2}}$,…, ${{A}_{n}}{{B}_{n}}{{C}_{n}}{{D}_{n}}$. The transformation matrix from ${{\mathbf{F}}_{1}}$ to ${{\mathbf{F}}_{2}}$ is ${{\mathbf{B}}_{21}}$, the transformation matrix from ${{\mathbf{F}}_{2}}$ to ${{\mathbf{F}}_{3}}$ is ${{\mathbf{B}}_{32}}$,…, the transformation matrix from ${{\mathbf{F}}_{n-1}}$ to ${{\mathbf{F}}_{n}}$ is ${{\mathbf{B}}_{n\left( n-1 \right)}}$. Then:
	\[{{\mathbf{F}}_{n}}={{\mathbf{B}}_{n\left( n-1 \right)}}\cdots {{\mathbf{B}}_{21}}{{\mathbf{F}}_{1}}.\]
\end{theorem}


\begin{proof}
	According to theorem \ref{thm:Biaojiabianhuangongshi1}, it is obtained that
	\[{{\mathbf{F}}_{n}}={{\mathbf{B}}_{n\left( n-1 \right)}}{{\mathbf{F}}_{n-1}}={{\mathbf{B}}_{n\left( n-1 \right)}}{{\mathbf{B}}_{\left( n-1 \right)\left( n-2 \right)}}{{\mathbf{F}}_{n-2}}=\cdots ={{\mathbf{B}}_{n\left( n-1 \right)}}\cdots {{\mathbf{B}}_{21}}{{\mathbf{F}}_{1}}.\]
\end{proof}
\hfill $\square$\par

\begin{theorem}{Formula 1 for frame component transformation, Daiyuan Zhang}{BiaojiafenliangBianhuangongshi1}\label{BiaojiafenliangBianhuangongshi1}
	Given two tetrahedra ${{A}_{1}}{{B}_{1}}{{C}_{1}}{{D}_{1}}$ and ${{A}_{2}}{{B}_{2}}{{C}_{2}}{{D}_{2}}$. The frame components on $\left( O;{{A}_{1}},{{B}_{1}},{{C}_{1}},{{D}_{1}} \right)$ of the four vertices ${{A}_{2}}$, ${{B}_{2}}$, ${{C}_{2}}$, ${{D}_{2}}$ of the tetrahedron ${{A}_{2}}{{B}_{2}}{{C}_{2}}{{D}_{2}}$ are $\beta _{{{A}_{1}}}^{{{A}_{2}}}$, $\beta _{{{B}_{1}}}^{{{A}_{2}}}$, $\beta _{{{C}_{1}}}^{{{A}_{2}}}$, $\beta _{{{D}_{1}}}^{{{A}_{2}}}$; $\beta _{{{A}_{1}}}^{{{B}_{2}}}$, $\beta _{{{B}_{1}}}^{{{B}_{2}}}$, $\beta _{{{C}_{1}}}^{{{B}_{2}}}$, $\beta _{{{D}_{1}}}^{{{B}_{2}}}$; $\beta _{{{A}_{1}}}^{{{C}_{2}}}$, $\beta _{{{B}_{1}}}^{{{C}_{2}}}$, $\beta _{{{C}_{1}}}^{{{C}_{2}}}$, $\beta _{{{D}_{1}}}^{{{C}_{2}}}$; $\beta _{{{A}_{1}}}^{{{D}_{2}}}$, $\beta _{{{B}_{1}}}^{{{D}_{2}}}$, $\beta _{{{C}_{1}}}^{{{D}_{2}}}$, $\beta _{{{D}_{1}}}^{{{D}_{2}}}$; respectively. 
	The frame components of point $P$ on the frame $\left( O;{{A}_{1}},{{B}_{1}},{{C}_{1}},{{D}_{1}} \right)$ are $\beta _{{{A}_{1}}}^{P}$, $\beta _{{{B}_{1}}}^{P}$, $\beta _{{{C}_{1}}}^{P}$, $\beta _{{{D}_{1}}}^{P}$; the frame components of point $P$ on the frame $\left( O;{{A}_{2}},{{B}_{2}},{{C}_{2}},{{D}_{2}} \right)$ are $\beta _{{{A}_{2}}}^{P}$, $\beta _{{{B}_{2}}}^{P}$, $\beta _{{{C}_{2}}}^{P}$, $\beta _{{{D}_{2}}}^{P}$. Then
	\[\left( \begin{aligned}
		& \begin{matrix}
			\beta _{{{A}_{1}}}^{P}  \\
			\beta _{{{B}_{1}}}^{P}  \\
			\beta _{{{C}_{1}}}^{P}  \\
		\end{matrix} \\ 
		& \beta _{{{D}_{1}}}^{P} \\ 
	\end{aligned} \right)=\left( \begin{matrix}
		\beta _{{{A}_{1}}}^{{{A}_{2}}} & \beta _{{{A}_{1}}}^{{{B}_{2}}} & \beta _{{{A}_{1}}}^{{{C}_{2}}} & \beta _{{{A}_{1}}}^{{{D}_{2}}}  \\
		\beta _{{{B}_{1}}}^{{{A}_{2}}} & \beta _{{{B}_{1}}}^{{{B}_{2}}} & \beta _{{{B}_{1}}}^{{{C}_{2}}} & \beta _{{{B}_{1}}}^{{{D}_{2}}}  \\
		\beta _{{{C}_{1}}}^{{{A}_{2}}} & \beta _{{{C}_{1}}}^{{{B}_{2}}} & \beta _{{{C}_{1}}}^{{{C}_{2}}} & \beta _{{{C}_{1}}}^{{{D}_{2}}}  \\
		\beta _{{{D}_{1}}}^{{{A}_{2}}} & \beta _{{{D}_{1}}}^{{{B}_{2}}} & \beta _{{{D}_{1}}}^{{{C}_{2}}} & \beta _{{{D}_{1}}}^{{{D}_{2}}}  \\
	\end{matrix} \right)\left( \begin{aligned}
		& \begin{matrix}
			\beta _{{{A}_{2}}}^{P}  \\
			\beta _{{{B}_{2}}}^{P}  \\
			\beta _{{{C}_{2}}}^{P}  \\
		\end{matrix} \\ 
		& \beta _{{{D}_{2}}}^{P} \\ 
	\end{aligned} \right).\]
	
	The sum of elements in each column of the matrix is 1.
\end{theorem}

\begin{proof}
	Given a point $P$ in space, for the frame $\left( O;{{A}_{1}},{{B}_{1}},{{C}_{1}},{{D}_{1}} \right)$, it can be obtained that
	\[\overrightarrow{OP}=\beta _{{{A}_{1}}}^{P}\overrightarrow{O{{A}_{1}}}+\beta _{{{B}_{1}}}^{P}\overrightarrow{O{{B}_{1}}}+\beta _{{{C}_{1}}}^{P}\overrightarrow{O{{C}_{1}}}+\beta _{{{D}_{1}}}^{P}\overrightarrow{O{{D}_{1}}},\]
	\[\beta _{{{A}_{1}}}^{P}+\beta _{{{B}_{1}}}^{P}+\beta _{{{C}_{1}}}^{P}+\beta _{{{D}_{1}}}^{P}=1.\]
	
	For the frame $\left( O;{{A}_{2}},{{B}_{2}},{{C}_{2}},{{D}_{2}} \right)$, it can be obtained that
	\[\overrightarrow{OP}=\beta _{{{A}_{2}}}^{P}\overrightarrow{O{{A}_{2}}}+\beta _{{{B}_{2}}}^{P}\overrightarrow{O{{B}_{2}}}+\beta _{{{C}_{2}}}^{P}\overrightarrow{O{{C}_{2}}}+\beta _{{{D}_{2}}}^{P}\overrightarrow{O{{D}_{2}}},\]
	\[\beta _{{{A}_{2}}}^{P}+\beta _{{{B}_{2}}}^{P}+\beta _{{{C}_{2}}}^{P}+\beta _{{{D}_{2}}}^{P}=1.\]
	
	Therefore
	\begin{align*}
		\overrightarrow{OP}& =\beta _{{{A}_{1}}}^{P}\overrightarrow{O{{A}_{1}}}+\beta _{{{B}_{1}}}^{P}\overrightarrow{O{{B}_{1}}}+\beta _{{{C}_{1}}}^{P}\overrightarrow{O{{C}_{1}}}+\beta _{{{D}_{1}}}^{P}\overrightarrow{O{{D}_{1}}} \\ 
		& =\beta _{{{A}_{2}}}^{P}\overrightarrow{O{{A}_{2}}}+\beta _{{{B}_{2}}}^{P}\overrightarrow{O{{B}_{2}}}+\beta _{{{C}_{2}}}^{P}\overrightarrow{O{{C}_{2}}}+\beta _{{{D}_{2}}}^{P}\overrightarrow{O{{D}_{2}}},  
	\end{align*}
	i.e.
	\[\begin{aligned}
		& \left( \beta _{{{A}_{1}}}^{P}\overrightarrow{O{{A}_{1}}}+\beta _{{{B}_{1}}}^{P}\overrightarrow{O{{B}_{1}}}+\beta _{{{C}_{1}}}^{P}\overrightarrow{O{{C}_{1}}}+\beta _{{{D}_{1}}}^{P}\overrightarrow{O{{D}_{1}}} \right) \\ 
		& -\left( \beta _{{{A}_{2}}}^{P}\overrightarrow{O{{A}_{2}}}+\beta _{{{B}_{2}}}^{P}\overrightarrow{O{{B}_{2}}}+\beta _{{{C}_{2}}}^{P}\overrightarrow{O{{C}_{2}}}+\beta _{{{D}_{2}}}^{P}\overrightarrow{O{{D}_{2}}} \right)=\overrightarrow{0}. \\ 
	\end{aligned}\]
	
	And
	\[\overrightarrow{O{{A}_{2}}}=\beta _{{{A}_{1}}}^{{{A}_{2}}}\overrightarrow{O{{A}_{1}}}+\beta _{{{B}_{1}}}^{{{A}_{2}}}\overrightarrow{O{{B}_{1}}}+\beta _{{{C}_{1}}}^{{{A}_{2}}}\overrightarrow{O{{C}_{1}}}+\beta _{{{D}_{1}}}^{{{A}_{2}}}\overrightarrow{O{{D}_{1}}},\]
	\[\overrightarrow{O{{B}_{2}}}=\beta _{{{A}_{1}}}^{{{B}_{2}}}\overrightarrow{O{{A}_{1}}}+\beta _{{{B}_{1}}}^{{{B}_{2}}}\overrightarrow{O{{B}_{1}}}+\beta _{{{C}_{1}}}^{{{B}_{2}}}\overrightarrow{O{{C}_{1}}}+\beta _{{{D}_{1}}}^{{{B}_{2}}}\overrightarrow{O{{D}_{1}}},\]
	\[\overrightarrow{O{{C}_{2}}}=\beta _{{{A}_{1}}}^{{{C}_{2}}}\overrightarrow{O{{A}_{1}}}+\beta _{{{B}_{1}}}^{{{C}_{2}}}\overrightarrow{O{{B}_{1}}}+\beta _{{{C}_{1}}}^{{{C}_{2}}}\overrightarrow{O{{C}_{1}}}+\beta _{{{D}_{1}}}^{{{C}_{2}}}\overrightarrow{O{{D}_{1}}},\]
	\[\overrightarrow{O{{D}_{2}}}=\beta _{{{A}_{1}}}^{{{D}_{2}}}\overrightarrow{O{{A}_{1}}}+\beta _{{{B}_{1}}}^{{{D}_{2}}}\overrightarrow{O{{B}_{1}}}+\beta _{{{C}_{1}}}^{{{D}_{2}}}\overrightarrow{O{{C}_{1}}}+\beta _{{{D}_{1}}}^{{{D}_{2}}}\overrightarrow{O{{D}_{1}}}.\]
	
	Substituting into the above equation yields:
	\[\begin{aligned}
		& \left( \beta _{{{A}_{1}}}^{P}\overrightarrow{O{{A}_{1}}}+\beta _{{{B}_{1}}}^{P}\overrightarrow{O{{B}_{1}}}+\beta _{{{C}_{1}}}^{P}\overrightarrow{O{{C}_{1}}}+\beta _{{{D}_{1}}}^{P}\overrightarrow{O{{D}_{1}}} \right) \\ 
		& -\left( \begin{aligned}
			& \beta _{{{A}_{2}}}^{P}\left( \beta _{{{A}_{1}}}^{{{A}_{2}}}\overrightarrow{O{{A}_{1}}}+\beta _{{{B}_{1}}}^{{{A}_{2}}}\overrightarrow{O{{B}_{1}}}+\beta _{{{C}_{1}}}^{{{A}_{2}}}\overrightarrow{O{{C}_{1}}}+\beta _{{{D}_{1}}}^{{{A}_{2}}}\overrightarrow{O{{D}_{1}}} \right) \\ 
			& +\beta _{{{B}_{2}}}^{P}\left( \beta _{{{A}_{1}}}^{{{B}_{2}}}\overrightarrow{O{{A}_{1}}}+\beta _{{{B}_{1}}}^{{{B}_{2}}}\overrightarrow{O{{B}_{1}}}+\beta _{{{C}_{1}}}^{{{B}_{2}}}\overrightarrow{O{{C}_{1}}}+\beta _{{{D}_{1}}}^{{{B}_{2}}}\overrightarrow{O{{D}_{1}}} \right) \\ 
			& +\beta _{{{C}_{2}}}^{P}\left( \beta _{{{A}_{1}}}^{{{C}_{2}}}\overrightarrow{O{{A}_{1}}}+\beta _{{{B}_{1}}}^{{{C}_{2}}}\overrightarrow{O{{B}_{1}}}+\beta _{{{C}_{1}}}^{{{C}_{2}}}\overrightarrow{O{{C}_{1}}}+\beta _{{{D}_{1}}}^{{{C}_{2}}}\overrightarrow{O{{D}_{1}}} \right) \\ 
			& +\beta _{{{D}_{2}}}^{P}\left( \beta _{{{A}_{1}}}^{{{D}_{2}}}\overrightarrow{O{{A}_{1}}}+\beta _{{{B}_{1}}}^{{{D}_{2}}}\overrightarrow{O{{B}_{1}}}+\beta _{{{C}_{1}}}^{{{D}_{2}}}\overrightarrow{O{{C}_{1}}}+\beta _{{{D}_{1}}}^{{{D}_{2}}}\overrightarrow{O{{D}_{1}}} \right) \\ 
		\end{aligned} \right)=\overrightarrow{0}. \\ 
	\end{aligned}\]
	
	It is easy to see that the sum of all coefficients in the above equation is 0.
	
	After merging similar items in the above equation, it can be obtained that
	\[\begin{aligned}
		& \left( \beta _{{{A}_{1}}}^{P}-\beta _{{{A}_{2}}}^{P}\beta _{{{A}_{1}}}^{{{A}_{2}}}-\beta _{{{B}_{2}}}^{P}\beta _{{{A}_{1}}}^{{{B}_{2}}}-\beta _{{{C}_{2}}}^{P}\beta _{{{A}_{1}}}^{{{C}_{2}}}-\beta _{{{D}_{2}}}^{P}\beta _{{{A}_{1}}}^{{{D}_{2}}} \right)\overrightarrow{O{{A}_{1}}} \\ 
		& +\left( \beta _{{{B}_{1}}}^{P}-\beta _{{{A}_{2}}}^{P}\beta _{{{B}_{1}}}^{{{A}_{2}}}-\beta _{{{B}_{2}}}^{P}\beta _{{{B}_{1}}}^{{{B}_{2}}}-\beta _{{{C}_{2}}}^{P}\beta _{{{B}_{1}}}^{{{C}_{2}}}-\beta _{{{D}_{2}}}^{P}\beta _{{{B}_{1}}}^{{{D}_{2}}} \right)\overrightarrow{O{{B}_{1}}} \\ 
		& +\left( \beta _{{{C}_{1}}}^{P}-\beta _{{{A}_{2}}}^{P}\beta _{{{C}_{1}}}^{{{A}_{2}}}-\beta _{{{B}_{2}}}^{P}\beta _{{{C}_{1}}}^{{{B}_{2}}}-\beta _{{{C}_{2}}}^{P}\beta _{{{C}_{1}}}^{{{C}_{2}}}-\beta _{{{D}_{2}}}^{P}\beta _{{{C}_{1}}}^{{{D}_{2}}} \right)\overrightarrow{O{{C}_{1}}} \\ 
		& +\left( \beta _{{{D}_{1}}}^{P}-\beta _{{{A}_{2}}}^{P}\beta _{{{D}_{1}}}^{{{A}_{2}}}-\beta _{{{B}_{2}}}^{P}\beta _{{{D}_{1}}}^{{{B}_{2}}}-\beta _{{{C}_{2}}}^{P}\beta _{{{D}_{1}}}^{{{C}_{2}}}-\beta _{{{D}_{2}}}^{P}\beta _{{{D}_{1}}}^{{{D}_{2}}} \right)\overrightarrow{O{{D}_{1}}}=\overrightarrow{0}. \\ 
	\end{aligned}\]
	
	According to theorem \ref{thm:Thm18.1.3}, it is obtained that
	\[\left\{ \begin{aligned}
		& \beta _{{{A}_{1}}}^{P}-\beta _{{{A}_{2}}}^{P}\beta _{{{A}_{1}}}^{{{A}_{2}}}-\beta _{{{B}_{2}}}^{P}\beta _{{{A}_{1}}}^{{{B}_{2}}}-\beta _{{{C}_{2}}}^{P}\beta _{{{A}_{1}}}^{{{C}_{2}}}-\beta _{{{D}_{2}}}^{P}\beta _{{{A}_{1}}}^{{{D}_{2}}}=0 \\ 
		& \beta _{{{B}_{1}}}^{P}-\beta _{{{A}_{2}}}^{P}\beta _{{{B}_{1}}}^{{{A}_{2}}}-\beta _{{{B}_{2}}}^{P}\beta _{{{B}_{1}}}^{{{B}_{2}}}-\beta _{{{C}_{2}}}^{P}\beta _{{{B}_{1}}}^{{{C}_{2}}}-\beta _{{{D}_{2}}}^{P}\beta _{{{B}_{1}}}^{{{D}_{2}}}=0 \\ 
		& \beta _{{{C}_{1}}}^{P}-\beta _{{{A}_{2}}}^{P}\beta _{{{C}_{1}}}^{{{A}_{2}}}-\beta _{{{B}_{2}}}^{P}\beta _{{{C}_{1}}}^{{{B}_{2}}}-\beta _{{{C}_{2}}}^{P}\beta _{{{C}_{1}}}^{{{C}_{2}}}-\beta _{{{D}_{2}}}^{P}\beta _{{{C}_{1}}}^{{{D}_{2}}}=0 \\ 
		& \beta _{{{D}_{1}}}^{P}-\beta _{{{A}_{2}}}^{P}\beta _{{{D}_{1}}}^{{{A}_{2}}}-\beta _{{{B}_{2}}}^{P}\beta _{{{D}_{1}}}^{{{B}_{2}}}-\beta _{{{C}_{2}}}^{P}\beta _{{{D}_{1}}}^{{{C}_{2}}}-\beta _{{{D}_{2}}}^{P}\beta _{{{D}_{1}}}^{{{D}_{2}}}=0. \\ 
	\end{aligned} \right.\]
	
	Or written in matrix form:
	\[\left( \begin{aligned}
		& \begin{matrix}
			\beta _{{{A}_{1}}}^{P}  \\
			\beta _{{{B}_{1}}}^{P}  \\
			\beta _{{{C}_{1}}}^{P}  \\
		\end{matrix} \\ 
		& \beta _{{{D}_{1}}}^{P} \\ 
	\end{aligned} \right)=\left( \begin{matrix}
		\beta _{{{A}_{1}}}^{{{A}_{2}}} & \beta _{{{A}_{1}}}^{{{B}_{2}}} & \beta _{{{A}_{1}}}^{{{C}_{2}}} & \beta _{{{A}_{1}}}^{{{D}_{2}}}  \\
		\beta _{{{B}_{1}}}^{{{A}_{2}}} & \beta _{{{B}_{1}}}^{{{B}_{2}}} & \beta _{{{B}_{1}}}^{{{C}_{2}}} & \beta _{{{B}_{1}}}^{{{D}_{2}}}  \\
		\beta _{{{C}_{1}}}^{{{A}_{2}}} & \beta _{{{C}_{1}}}^{{{B}_{2}}} & \beta _{{{C}_{1}}}^{{{C}_{2}}} & \beta _{{{C}_{1}}}^{{{D}_{2}}}  \\
		\beta _{{{D}_{1}}}^{{{A}_{2}}} & \beta _{{{D}_{1}}}^{{{B}_{2}}} & \beta _{{{D}_{1}}}^{{{C}_{2}}} & \beta _{{{D}_{1}}}^{{{D}_{2}}}  \\
	\end{matrix} \right)\left( \begin{aligned}
		& \begin{matrix}
			\beta _{{{A}_{2}}}^{P}  \\
			\beta _{{{B}_{2}}}^{P}  \\
			\beta _{{{C}_{2}}}^{P}  \\
		\end{matrix} \\ 
		& \beta _{{{D}_{2}}}^{P} \\ 
	\end{aligned} \right).\]
\end{proof}
\hfill $\square$\par

Imagine the tetrahedron ${{A}_{1}}{{B}_{1}}{{C}_{1}}{{D}_{1}}$ as a reference tetrahedron or an old tetrahedron, and ${{A}_{2}}{{B}_{2}}{{C}_{2}}{{D}_{2}}$  as a new tetrahedron. The above equation indicates that if the frame components of a point $P$ on a new tetrahedron are known, and the frame components of the four vertices of the new tetrahedron on the old tetrahedron are known, then the frame components of that point $P$ on the old tetrahedron can be calculated.


\begin{example}{}\label{QiuZhongxinsimiantiDeBiaojiafenliang}
	Given a tetrahedron ${{A}_{1}}{{B}_{1}}{{C}_{1}}{{D}_{1}}$  (referred to as the original tetrahedron). Form a new tetrahedron using the centroid of each side triangle of the tetrahedron ${{A}_{2}}{{B}_{2}}{{C}_{2}}{{D}_{2}}$ (referred to as the centroid tetrahedron), where ${{A}_{2}}$ is the centroid of $\triangle {{B}_{1}}{{C}_{1}}{{D}_{1}}$, ${{B}_{2}}$ is the centroid of $\triangle {{C}_{1}}{{D}_{1}}{{A}_{1}}$, ${{C}_{2}}$ is the centroid of $\triangle {{D}_{1}}{{A}_{1}}{{B}_{1}}$, ${{D}_{2}}$ is the centroid of $\triangle {{A}_{1}}{{B}_{1}}{{C}_{1}}$. Find the tetrahedral frame components of the centroid of centroid tetrahedron ${{A}_{2}}{{B}_{2}}{{C}_{2}}{{D}_{2}}$ on the original tetrahedral frame $\left( O;{{A}_{1}},{{B}_{1}},{{C}_{1}},{{D}_{1}} \right)$.	
\end{example}

\begin{solution}
	Assuming the centroid of the centroid tetrahedron ${{A}_{2}}{{B}_{2}}{{C}_{2}}{{D}_{2}}$  is ${{G}_{2}}$. Because ${{A}_{2}}$ is the centroid of $\triangle {{B}_{1}}{{C}_{1}}{{D}_{1}}$, using theorem \ref{thm:Thm20.4.1}, the frame components of ${{A}_{2}}$ on tetrahedral frame $\left( O;{{A}_{1}},{{B}_{1}},{{C}_{1}},{{D}_{1}} \right)$ are: 
	\[\beta _{{{A}_{1}}}^{{{A}_{2}}}=0,\ \beta _{{{B}_{1}}}^{{{A}_{2}}}=\alpha _{{{B}_{1}}}^{{{A}_{2}}}={1}/{3}\;,\ \beta _{{{C}_{1}}}^{{{A}_{2}}}=\alpha _{{{C}_{1}}}^{{{A}_{2}}}={1}/{3}\;,\ \beta _{{{D}_{1}}}^{{{A}_{2}}}=\alpha _{{{D}_{1}}}^{{{A}_{2}}}={1}/{3}\;.\ \]
	
	The frame components of ${{B}_{2}}$ on tetrahedral frame $\left( O;{{A}_{1}},{{B}_{1}},{{C}_{1}},{{D}_{1}} \right)$ are:
	\[\beta _{{{A}_{1}}}^{{{B}_{2}}}=\alpha _{{{A}_{1}}}^{{{B}_{2}}}=1/3\ ,\ \beta _{{{B}_{1}}}^{{{B}_{2}}}=0,\ \beta _{{{C}_{1}}}^{{{B}_{2}}}=\alpha _{{{C}_{1}}}^{{{B}_{2}}}=1/3\ ,\ \beta _{{{D}_{1}}}^{{{B}_{2}}}=\alpha _{{{D}_{1}}}^{{{B}_{2}}}=1/3\ .\ \]
	
	The frame components of ${{C}_{2}}$ on tetrahedral frame $\left( O;{{A}_{1}},{{B}_{1}},{{C}_{1}},{{D}_{1}} \right)$ are:
	\[\beta _{{{A}_{1}}}^{{{C}_{2}}}=\alpha _{{{A}_{1}}}^{{{C}_{2}}}=1/3\ ,\ \beta _{{{B}_{1}}}^{{{C}_{2}}}=\alpha _{{{B}_{1}}}^{{{C}_{2}}}=1/3\ ,\ \beta _{{{C}_{1}}}^{{{C}_{2}}}=0\ ,\ \beta _{{{D}_{1}}}^{{{C}_{2}}}=\alpha _{{{D}_{1}}}^{{{C}_{2}}}=1/3\ .\ \]
	
	The frame components of ${{D}_{2}}$ on tetrahedral frame $\left( O;{{A}_{1}},{{B}_{1}},{{C}_{1}},{{D}_{1}} \right)$ are:
	\[\beta _{{{A}_{1}}}^{{{D}_{2}}}=\alpha _{{{A}_{1}}}^{{{D}_{2}}}=1/3\ ,\ \beta _{{{B}_{1}}}^{{{D}_{2}}}=\alpha _{{{B}_{1}}}^{{{D}_{2}}}=1/3\ ,\ \beta _{{{C}_{1}}}^{{{D}_{2}}}=\alpha _{{{C}_{1}}}^{{{D}_{2}}}=1/3\ ,\ \beta _{{{D}_{1}}}^{{{D}_{2}}}=0\ .\ \]
	
	The frame components of ${{G}_{2}}$ on tetrahedral frame $\left( O;{{A}_{2}},{{B}_{2}},{{C}_{2}},{{D}_{2}} \right)$ are:
	\[\beta _{{{A}_{2}}}^{{{G}_{2}}}={1}/{4}\;,\ \beta _{{{B}_{2}}}^{{{G}_{2}}}={1}/{4}\;,\ \beta _{{{C}_{2}}}^{{{G}_{2}}}={1}/{4}\;,\ \beta _{{{D}_{2}}}^{{{G}_{2}}}={1}/{4}\;.\ \]
	
	According to the above theorem, it can be concluded that
	\[\left( \begin{aligned}
		& \begin{matrix}
			\beta _{{{A}_{1}}}^{{{G}_{2}}}  \\
			\beta _{{{B}_{1}}}^{{{G}_{2}}}  \\
			\beta _{{{C}_{1}}}^{{{G}_{2}}}  \\
		\end{matrix} \\ 
		& \beta _{{{D}_{1}}}^{{{G}_{2}}} \\ 
	\end{aligned} \right)=\left( \begin{matrix}
		\beta _{{{A}_{1}}}^{{{A}_{2}}} & \beta _{{{A}_{1}}}^{{{B}_{2}}} & \beta _{{{A}_{1}}}^{{{C}_{2}}} & \beta _{{{A}_{1}}}^{{{D}_{2}}}  \\
		\beta _{{{B}_{1}}}^{{{A}_{2}}} & \beta _{{{B}_{1}}}^{{{B}_{2}}} & \beta _{{{B}_{1}}}^{{{C}_{2}}} & \beta _{{{B}_{1}}}^{{{D}_{2}}}  \\
		\beta _{{{C}_{1}}}^{{{A}_{2}}} & \beta _{{{C}_{1}}}^{{{B}_{2}}} & \beta _{{{C}_{1}}}^{{{C}_{2}}} & \beta _{{{C}_{1}}}^{{{D}_{2}}}  \\
		\beta _{{{D}_{1}}}^{{{A}_{2}}} & \beta _{{{D}_{1}}}^{{{B}_{2}}} & \beta _{{{D}_{1}}}^{{{C}_{2}}} & \beta _{{{D}_{1}}}^{{{D}_{2}}}  \\
	\end{matrix} \right)\left( \begin{aligned}
		& \begin{matrix}
			\beta _{{{A}_{2}}}^{{{G}_{2}}}  \\
			\beta _{{{B}_{2}}}^{{{G}_{2}}}  \\
			\beta _{{{C}_{2}}}^{{{G}_{2}}}  \\
		\end{matrix} \\ 
		& \beta _{{{D}_{2}}}^{{{G}_{2}}} \\ 
	\end{aligned} \right)=\left( \begin{matrix}
		0 & \frac{1}{3} & \frac{1}{3} & \frac{1}{3}  \\
		\frac{1}{3} & 0 & \frac{1}{3} & \frac{1}{3}  \\
		\frac{1}{3} & \frac{1}{3} & 0 & \frac{1}{3}  \\
		\frac{1}{3} & \frac{1}{3} & \frac{1}{3} & 0  \\
	\end{matrix} \right)\left( \begin{aligned}
		& \begin{matrix}
			\frac{1}{4}  \\
			\frac{1}{4}  \\
			\frac{1}{4}  \\
		\end{matrix} \\ 
		& \frac{1}{4} \\ 
	\end{aligned} \right).\]
	
	i.e.
	\[\beta _{{{A}_{1}}}^{{{G}_{2}}}=\beta _{{{B}_{1}}}^{{{G}_{2}}}=\beta _{{{C}_{1}}}^{{{G}_{2}}}=\beta _{{{D}_{1}}}^{{{G}_{2}}}=\frac{1}{4}.\ \]
	
	Assuming the centroid of tetrahedron ${{A}_{1}}{{B}_{1}}{{C}_{1}}{{D}_{1}}$ is ${{G}_{1}}$, the above equation explains ${{G}_{2}}={{G}_{1}}$, that is, the centroid of the centroid tetrahedron coincides with the centroid of the original tetrahedron.
\end{solution}
\hfill $\diamond$\par


%
%
%
%
%
%
%
%
%
%
%
%
%
%
%
%
%
%
%
%
%
%
%
%

\chapter{Frame components for special intersecting center of a tetrahedron}\label{Ch21}

In this chapter, the calculation on the frame components for the special intersecting centers of a tetrahedron is studied, which are the centroid and the incenter.

The results of the previous chapter show that the frame components of the intersecting center of a tetrahedron can be expressed by the corresponding frame components of the intersecting centers of faces. Therefore, first of all, it is needed to find the frame componentss of the corresponding intersecting center of faces (abbreviated as IC-F).

\section{Frame components of a tetrahedron corresponding to centroid}\label{Sec21.1}
\subsection{Frame components of intersecting center of face corresponding to centroid}\label{Subsec21.1.1}
According to the solid geometry, the centroids of faces corresponding to the centroid of a tetrahedron are the centroids of the corresponding triangles. From section \ref{Sec8.1}, the following results can be obtained:
\[\alpha _{B}^{{{G}_{A}}}=\alpha _{C}^{{{G}_{A}}}=\alpha _{D}^{{{G}_{A}}}=\frac{1}{3},\]
\[\alpha _{C}^{{{G}_{B}}}=\alpha _{D}^{{{G}_{B}}}=\alpha _{A}^{{{G}_{B}}}=\frac{1}{3},\]	
\[\alpha _{D}^{{{G}_{C}}}=\alpha _{A}^{{{G}_{C}}}=\alpha _{B}^{{{G}_{C}}}=\frac{1}{3},\]	
\[\alpha _{A}^{{{G}_{D}}}=\alpha _{B}^{{{G}_{D}}}=\alpha _{C}^{{{G}_{D}}}=\frac{1}{3}.\]

\subsection{Frame components of centroid of a tetrahedron}\label{Subsec21.1.2}
According to the previous section and theorem \ref{thm:Thm20.2.1}, the following results can be obtained:
\[\beta _{A}^{G}=\frac{\alpha _{A}^{{{G}_{B}}}\left( 1-\alpha _{B}^{{{G}_{A}}} \right)}{1-\alpha _{A}^{{{G}_{B}}}\alpha _{B}^{{{G}_{A}}}}=\frac{\frac{1}{3}\times \left( 1-\frac{1}{3} \right)}{1-\frac{1}{3}\times \frac{1}{3}}=\frac{1}{4},\]	
\[\beta _{B}^{G}=\frac{\alpha _{B}^{{{G}_{A}}}\left( 1-\alpha _{A}^{{{G}_{B}}} \right)}{1-\alpha _{A}^{{{G}_{B}}}\alpha _{B}^{{{G}_{A}}}}=\frac{\frac{1}{3}\times \left( 1-\frac{1}{3} \right)}{1-\frac{1}{3}\times \frac{1}{3}}=\frac{1}{4},\]	
\[\beta _{C}^{G}=\frac{\alpha _{C}^{{{G}_{A}}}\left( 1-\alpha _{A}^{{{G}_{B}}} \right)}{1-\alpha _{A}^{{{G}_{B}}}\alpha _{B}^{{{G}_{A}}}}=\frac{\frac{1}{3}\times \left( 1-\frac{1}{3} \right)}{1-\frac{1}{3}\times \frac{1}{3}}=\frac{1}{4},\]	
\[\beta _{D}^{G}=\frac{\alpha _{D}^{{{G}_{A}}}\left( 1-\alpha _{A}^{{{G}_{B}}} \right)}{1-\alpha _{A}^{{{G}_{B}}}\alpha _{B}^{{{G}_{A}}}}=\frac{\frac{1}{3}\times \left( 1-\frac{1}{3} \right)}{1-\frac{1}{3}\times \frac{1}{3}}=\frac{1}{4}.\]	

So the following theorem can be obtained:
\begin{theorem}{Frame components for centroid of a tetrahedron, Daiyuan Zhang}{Thm21.1.1}\label{Thm21.1.1} 
	The frame components for the centroid of a tetrahedron are all constant ${1}/{4}\;$.
\end{theorem}

\section{Frame components of a tetrahedron corresponding to  incenter}\label{Sec21.2}
The frame components of intersecting center of faces (intersecting center of faces is abbreviated as IC-Fs) are the triangular frame components, so it can be calculated by using the previous formulas about the frame components of the triangular frame.

According to theorem \ref{thm:Thm6.1.1} and the IR-F of the incenter of tetrahedron in subsection \ref{Sec17.2} (as shown in Figure \ref{fig:tu17.2.1}), the frame components of IC-Fs corresponding to the incenter can be obtained.

\subsection{Frame components of intersecting center of faces corresponding to the incenter of a tetrahedron}\label{Subsec21.2.1}
Let
\[S={{S}^{A}}+{{S}^{B}}+{{S}^{C}}+{{S}^{D}}.\]
\subsubsection{Frame components of intersecting center of face ${{I}_{A}}$}
\[\alpha _{B}^{{{I}_{A}}}=\frac{1}{1+\lambda _{BC}^{{{I}_{A}}}+\lambda _{BD}^{{{I}_{A}}}}=\frac{1}{1+\frac{{{S}^{C}}}{{{S}^{B}}}+\frac{{{S}^{D}}}{{{S}^{B}}}}=\frac{{{S}^{B}}}{S-{{S}^{A}}},\]	
\[\alpha _{C}^{{{I}_{A}}}=\frac{\lambda _{BC}^{{{I}_{A}}}}{1+\lambda _{BC}^{{{I}_{A}}}+\lambda _{BD}^{{{I}_{A}}}}=\frac{\frac{{{S}^{C}}}{{{S}^{B}}}}{1+\frac{{{S}^{C}}}{{{S}^{B}}}+\frac{{{S}^{D}}}{{{S}^{B}}}}=\frac{{{S}^{C}}}{S-{{S}^{A}}},\]
\[\alpha _{D}^{{{I}_{A}}}=\frac{\lambda _{BD}^{{{I}_{A}}}}{1+\lambda _{BC}^{{{I}_{A}}}+\lambda _{BD}^{{{I}_{A}}}}=\frac{\frac{{{S}^{D}}}{{{S}^{B}}}}{1+\frac{{{S}^{C}}}{{{S}^{B}}}+\frac{{{S}^{D}}}{{{S}^{B}}}}=\frac{{{S}^{D}}}{S-{{S}^{A}}}.\]

\subsubsection{Frame components of intersecting center of face ${{I}_{B}}$}
\[\alpha _{C}^{{{I}_{B}}}=\frac{1}{1+\lambda _{CD}^{{{I}_{B}}}+\lambda _{CA}^{{{I}_{B}}}}=\frac{1}{1+\frac{{{S}^{D}}}{{{S}^{C}}}+\frac{{{S}^{A}}}{{{S}^{C}}}}=\frac{{{S}^{C}}}{S-{{S}^{B}}},\]
\[\alpha _{D}^{{{I}_{B}}}=\frac{\lambda _{CD}^{{{I}_{B}}}}{1+\lambda _{CD}^{{{I}_{B}}}+\lambda _{CA}^{{{I}_{B}}}}=\frac{\frac{{{S}^{D}}}{{{S}^{C}}}}{1+\frac{{{S}^{D}}}{{{S}^{C}}}+\frac{{{S}^{A}}}{{{S}^{C}}}}=\frac{{{S}^{D}}}{S-{{S}^{B}}},\]
\[\alpha _{A}^{{{I}_{B}}}=\frac{\lambda _{CA}^{{{I}_{B}}}}{1+\lambda _{CD}^{{{I}_{B}}}+\lambda _{CA}^{{{I}_{B}}}}=\frac{\frac{{{S}^{A}}}{{{S}^{C}}}}{1+\frac{{{S}^{D}}}{{{S}^{C}}}+\frac{{{S}^{A}}}{{{S}^{C}}}}=\frac{{{S}^{A}}}{S-{{S}^{B}}}.\]	

\subsubsection{Frame components of intersecting center of face ${{I}_{C}}$}
\[\alpha _{D}^{{{I}_{C}}}=\frac{1}{1+\lambda _{DA}^{{{I}_{B}}}+\lambda _{DB}^{{{I}_{B}}}}=\frac{1}{1+\frac{{{S}^{A}}}{{{S}^{D}}}+\frac{{{S}^{B}}}{{{S}^{D}}}}=\frac{{{S}^{D}}}{S-{{S}^{C}}},\]
\[\alpha _{A}^{{{I}_{C}}}=\frac{\lambda _{DA}^{{{I}_{B}}}}{1+\lambda _{DA}^{{{I}_{B}}}+\lambda _{DB}^{{{I}_{B}}}}=\frac{\frac{{{S}^{A}}}{{{S}^{D}}}}{1+\frac{{{S}^{A}}}{{{S}^{D}}}+\frac{{{S}^{B}}}{{{S}^{D}}}}=\frac{{{S}^{A}}}{S-{{S}^{C}}},\]
\[\alpha _{B}^{{{I}_{C}}}=\frac{\lambda _{DB}^{{{I}_{B}}}}{1+\lambda _{DA}^{{{I}_{B}}}+\lambda _{DB}^{{{I}_{B}}}}=\frac{\frac{{{S}^{B}}}{{{S}^{D}}}}{1+\frac{{{S}^{A}}}{{{S}^{D}}}+\frac{{{S}^{B}}}{{{S}^{D}}}}=\frac{{{S}^{B}}}{S-{{S}^{C}}}.\]

\subsubsection{Frame components of intersecting center of face ${{I}_{D}}$}
\[\alpha _{A}^{{{I}_{D}}}=\frac{1}{1+\lambda _{AB}^{{{I}_{B}}}+\lambda _{AC}^{{{I}_{B}}}}=\frac{1}{1+\frac{{{S}^{B}}}{{{S}^{A}}}+\frac{{{S}^{C}}}{{{S}^{A}}}}=\frac{{{S}^{A}}}{S-{{S}^{D}}},\]
\[\alpha _{B}^{{{I}_{D}}}=\frac{\lambda _{AB}^{{{I}_{B}}}}{1+\lambda _{AB}^{{{I}_{B}}}+\lambda _{AC}^{{{I}_{B}}}}=\frac{\frac{{{S}^{B}}}{{{S}^{A}}}}{1+\frac{{{S}^{B}}}{{{S}^{A}}}+\frac{{{S}^{C}}}{{{S}^{A}}}}=\frac{{{S}^{B}}}{S-{{S}^{D}}},\]
\[\alpha _{C}^{{{I}_{D}}}=\frac{\lambda _{AC}^{{{I}_{B}}}}{1+\lambda _{AB}^{{{I}_{B}}}+\lambda _{AC}^{{{I}_{B}}}}=\frac{\frac{{{S}^{C}}}{{{S}^{A}}}}{1+\frac{{{S}^{B}}}{{{S}^{A}}}+\frac{{{S}^{C}}}{{{S}^{A}}}}=\frac{{{S}^{C}}}{S-{{S}^{D}}}.\]

Suppose the point $D$ is infinitely close to the point ${{I}_{D}}$, then another expression of component of triangular frame can be obtained (compared with the contents in the previous chapters).

\subsection{Frame components of a tetrahedron corresponding to incenter}\label{Subsec21.2.2}
Suppose the surface area of tetrahedron $ABCD$ is $S$; ${{S}^{A}}$ is the area of the triangle ($\triangle BCD$) opposite to the vertex of the tetrahedron $A$, and other notations are analogical. The following results can be obtained from the previous section and theorem \ref{thm:Thm20.2.1}:
\[\begin{aligned}
	\beta _{A}^{I}&=\frac{\alpha _{A}^{{{I}_{B}}}\left( 1-\alpha _{B}^{{{I}_{A}}} \right)}{1-\alpha _{A}^{{{I}_{B}}}\alpha _{B}^{{{I}_{A}}}}=\frac{\frac{{{S}^{A}}}{S-{{S}^{B}}}\left( 1-\frac{{{S}^{B}}}{S-{{S}^{A}}} \right)}{1-\frac{{{S}^{A}}}{S-{{S}^{B}}}\frac{{{S}^{B}}}{S-{{S}^{A}}}} \\ 
	& =\frac{{{S}^{A}}\left( S-{{S}^{A}}-{{S}^{B}} \right)}{\left( S-{{S}^{B}} \right)\left( S-{{S}^{A}} \right)-{{S}^{A}}{{S}^{B}}}=\frac{{{S}^{A}}\left( S-\left( {{S}^{A}}+{{S}^{B}} \right) \right)}{{{S}^{2}}-\left( {{S}^{A}}+{{S}^{B}} \right)S}=\frac{{{S}^{A}}}{S}.  
\end{aligned}\]	

After rotation ($A$→$B$→$C$→$D$→$A$), the other formulas are obtained as follows:
\[\beta _{B}^{I}=\frac{\alpha _{B}^{{{P}_{A}}}\left( 1-\alpha _{A}^{{{P}_{B}}} \right)}{1-\alpha _{A}^{{{P}_{B}}}\alpha _{B}^{{{P}_{A}}}}=\frac{{{S}^{B}}}{S},\]	
\[\beta _{C}^{I}=\frac{\alpha _{C}^{{{P}_{A}}}\left( 1-\alpha _{A}^{{{P}_{B}}} \right)}{1-\alpha _{A}^{{{P}_{B}}}\alpha _{B}^{{{P}_{A}}}}=\frac{{{S}^{C}}}{S},\]	
\[\beta _{D}^{I}=\frac{\alpha _{D}^{{{P}_{A}}}\left( 1-\alpha _{A}^{{{P}_{B}}} \right)}{1-\alpha _{A}^{{{P}_{B}}}\alpha _{B}^{{{P}_{A}}}}=\frac{{{S}^{D}}}{S}.\]	

Then the following theorem is obtained:
\begin{theorem}{Frame components of a tetrahedron corresponding to incenter, Daiyuan Zhang}{Thm21.2.1}\label{Thm21.2.1} 
	Frame components of incenter of tetrahedron $ABCD$ are	
	\begin{flalign*}
		\beta _{A}^{I}=\frac{{{S}^{A}}}{S}, \beta _{B}^{I}=\frac{{{S}^{B}}}{S}, \beta _{C}^{I}=\frac{{{S}^{C}}}{S}, \beta _{D}^{I}=\frac{{{S}^{D}}}{S}.
	\end{flalign*}
\end{theorem}

\section{Frame components of a tetrahedron corresponding to excenter}\label{Sec21.3}
\subsection{Frame components of intersecting center of faces corresponding to excenter of a tetrahedron}\label{Subsec21.3.1}
\subsubsection{Frame components of intersecting center of face $E_{A}^{A}$}
\[\alpha _{B}^{E_{A}^{A}}=\frac{1}{1+\lambda _{BC}^{E_{A}^{A}}+\lambda _{BD}^{E_{A}^{A}}}=\frac{1}{1+\frac{{{S}^{C}}}{{{S}^{B}}}+\frac{{{S}^{D}}}{{{S}^{B}}}}=\frac{{{S}^{B}}}{S-{{S}^{A}}},\]	
\[\alpha _{C}^{E_{A}^{A}}=\frac{\lambda _{BC}^{E_{A}^{A}}}{1+\lambda _{BC}^{E_{A}^{A}}+\lambda _{BD}^{E_{A}^{A}}}=\frac{\frac{{{S}^{C}}}{{{S}^{B}}}}{1+\frac{{{S}^{C}}}{{{S}^{B}}}+\frac{{{S}^{D}}}{{{S}^{B}}}}=\frac{{{S}^{C}}}{S-{{S}^{A}}},\]	
\[\alpha _{D}^{E_{A}^{A}}=\frac{\lambda _{BD}^{E_{A}^{A}}}{1+\lambda _{BC}^{E_{A}^{A}}+\lambda _{BD}^{E_{A}^{A}}}=\frac{\frac{{{S}^{D}}}{{{S}^{B}}}}{1+\frac{{{S}^{C}}}{{{S}^{B}}}+\frac{{{S}^{D}}}{{{S}^{B}}}}=\frac{{{S}^{D}}}{S-{{S}^{A}}}.\]	

\subsubsection{Frame components of intersecting center of face $E_{A}^{B}$}
\[\alpha _{C}^{E_{A}^{B}}=\frac{1}{1+\lambda _{CD}^{E_{A}^{B}}+\lambda _{CA}^{E_{A}^{B}}}=\frac{1}{1+\frac{{{S}^{D}}}{{{S}^{C}}}-\frac{{{S}^{A}}}{{{S}^{C}}}}=\frac{{{S}^{C}}}{{{S}^{C}}+{{S}^{D}}-{{S}^{A}}},\]
\[\alpha _{D}^{E_{A}^{B}}=\frac{\lambda _{CD}^{E_{A}^{B}}}{1+\lambda _{CD}^{E_{A}^{B}}+\lambda _{CA}^{E_{A}^{B}}}=\frac{\frac{{{S}^{D}}}{{{S}^{C}}}}{1+\frac{{{S}^{D}}}{{{S}^{C}}}-\frac{{{S}^{A}}}{{{S}^{C}}}}=\frac{{{S}^{D}}}{{{S}^{C}}+{{S}^{D}}-{{S}^{A}}},\]
\[\alpha _{A}^{E_{A}^{B}}=\frac{\lambda _{CA}^{E_{A}^{B}}}{1+\lambda _{CD}^{E_{A}^{B}}+\lambda _{CA}^{E_{A}^{B}}}=\frac{-\frac{{{S}^{A}}}{{{S}^{C}}}}{1+\frac{{{S}^{D}}}{{{S}^{C}}}-\frac{{{S}^{A}}}{{{S}^{C}}}}=-\frac{{{S}^{A}}}{{{S}^{C}}+{{S}^{D}}-{{S}^{A}}}.\]

\subsubsection{Frame components of intersecting center of face $E_{A}^{C}$}
\[\alpha _{D}^{E_{A}^{C}}=\frac{1}{1+\lambda _{DA}^{E_{A}^{C}}+\lambda _{DB}^{E_{A}^{C}}}=\frac{1}{1-\frac{{{S}^{A}}}{{{S}^{D}}}+\frac{{{S}^{B}}}{{{S}^{D}}}}=\frac{{{S}^{D}}}{{{S}^{D}}+{{S}^{B}}-{{S}^{A}}},\]
\[\alpha _{A}^{E_{A}^{C}}=\frac{\lambda _{DA}^{E_{A}^{C}}}{1+\lambda _{DA}^{E_{A}^{C}}+\lambda _{DB}^{E_{A}^{C}}}=\frac{-\frac{{{S}^{A}}}{{{S}^{D}}}}{1-\frac{{{S}^{A}}}{{{S}^{D}}}+\frac{{{S}^{B}}}{{{S}^{D}}}}=-\frac{{{S}^{A}}}{{{S}^{D}}+{{S}^{B}}-{{S}^{A}}},\]
\[\alpha _{B}^{E_{A}^{C}}=\frac{\lambda _{DB}^{E_{A}^{C}}}{1+\lambda _{DA}^{E_{A}^{C}}+\lambda _{DB}^{E_{A}^{C}}}=\frac{\frac{{{S}^{B}}}{{{S}^{D}}}}{1-\frac{{{S}^{A}}}{{{S}^{D}}}+\frac{{{S}^{B}}}{{{S}^{D}}}}=\frac{{{S}^{B}}}{{{S}^{D}}+{{S}^{B}}-{{S}^{A}}}.\]

\subsubsection{Frame components of intersecting center of face $E_{A}^{D}$}
\[\alpha _{A}^{E_{A}^{D}}=\frac{1}{1+\lambda _{AB}^{E_{A}^{D}}+\lambda _{AC}^{E_{A}^{D}}}=\frac{1}{1-\frac{{{S}^{B}}}{{{S}^{A}}}-\frac{{{S}^{C}}}{{{S}^{A}}}}=-\frac{{{S}^{A}}}{{{S}^{B}}+{{S}^{C}}-{{S}^{A}}},\]
\[\alpha _{B}^{E_{A}^{D}}=\frac{\lambda _{AB}^{E_{A}^{D}}}{1+\lambda _{AB}^{E_{A}^{D}}+\lambda _{AC}^{E_{A}^{D}}}=\frac{-\frac{{{S}^{B}}}{{{S}^{A}}}}{1-\frac{{{S}^{B}}}{{{S}^{A}}}-\frac{{{S}^{C}}}{{{S}^{A}}}}=\frac{{{S}^{B}}}{{{S}^{B}}+{{S}^{C}}-{{S}^{A}}},\]
\[\alpha _{C}^{E_{A}^{D}}=\frac{\lambda _{AC}^{E_{A}^{D}}}{1+\lambda _{AB}^{E_{A}^{D}}+\lambda _{AC}^{E_{A}^{D}}}=\frac{-\frac{{{S}^{C}}}{{{S}^{A}}}}{1-\frac{{{S}^{B}}}{{{S}^{A}}}-\frac{{{S}^{C}}}{{{S}^{A}}}}=\frac{{{S}^{C}}}{{{S}^{B}}+{{S}^{C}}-{{S}^{A}}}.\]

\subsection{Frame components for the center of escribed sphere in trihedral angle of a tetrahedron}\label{Subsec21.3.2}

According to theorem \ref{thm:Thm20.2.1}, the following results are obtained:
\[\beta _{A}^{{{E}_{A}}}=\frac{\alpha _{A}^{E_{A}^{B}}\left( 1-\alpha _{B}^{E_{A}^{A}} \right)}{1-\alpha _{A}^{E_{A}^{B}}\alpha _{B}^{E_{A}^{A}}}=-\frac{\frac{{{S}^{A}}}{{{S}^{C}}+{{S}^{D}}-{{S}^{A}}}\left( 1-\frac{{{S}^{B}}}{S-{{S}^{A}}} \right)}{1+\frac{{{S}^{A}}}{{{S}^{C}}+{{S}^{D}}-{{S}^{A}}}\frac{{{S}^{B}}}{S-{{S}^{A}}}},\]	
i.e.
\[\begin{aligned}
	\beta _{A}^{{{E}_{A}}}&=-\frac{{{S}^{A}}\left( S-{{S}^{A}}-{{S}^{B}} \right)}{\left( {{S}^{C}}+{{S}^{D}}-{{S}^{A}} \right)\left( S-{{S}^{A}} \right)+{{S}^{A}}{{S}^{B}}} \\ 
	& =-\frac{{{S}^{A}}\left( {{S}^{C}}+{{S}^{D}} \right)}{\left( {{S}^{C}}+{{S}^{D}}-{{S}^{A}} \right)\left( {{S}^{C}}+{{S}^{D}}+{{S}^{B}} \right)+{{S}^{A}}{{S}^{B}}},  
\end{aligned}\]
i.e.
\[\begin{aligned}
	\beta _{A}^{{{E}_{A}}}&=-\frac{{{S}^{A}}\left( {{S}^{C}}+{{S}^{D}} \right)}{{{\left( {{S}^{C}}+{{S}^{D}} \right)}^{2}}+\left( {{S}^{B}}-{{S}^{A}} \right)\left( {{S}^{C}}+{{S}^{D}} \right)} \\ 
	& =-\frac{{{S}^{A}}}{\left( {{S}^{C}}+{{S}^{D}} \right)+\left( {{S}^{B}}-{{S}^{A}} \right)}=-\frac{{{S}^{A}}}{S-2{{S}^{A}}}.  
\end{aligned}\]

And
\[\begin{aligned}
	\beta _{B}^{{{E}_{A}}}&=\frac{\alpha _{B}^{E_{A}^{C}}\left( 1-\alpha _{C}^{E_{A}^{B}} \right)}{1-\alpha _{B}^{E_{A}^{C}}\alpha _{C}^{E_{A}^{B}}}=\frac{\frac{{{S}^{B}}}{{{S}^{D}}+{{S}^{B}}-{{S}^{A}}}\left( 1-\frac{{{S}^{C}}}{{{S}^{C}}+{{S}^{D}}-{{S}^{A}}} \right)}{1-\frac{{{S}^{B}}}{{{S}^{D}}+{{S}^{B}}-{{S}^{A}}}\frac{{{S}^{C}}}{{{S}^{C}}+{{S}^{D}}-{{S}^{A}}}} \\ 
	& =\frac{{{S}^{B}}\left( {{S}^{D}}-{{S}^{A}} \right)}{\left( {{S}^{D}}+{{S}^{B}}-{{S}^{A}} \right)\left( {{S}^{C}}+{{S}^{D}}-{{S}^{A}} \right)-{{S}^{B}}{{S}^{C}}},  
\end{aligned}\]	
i.e.
\[\begin{aligned}
	\beta _{B}^{{{E}_{A}}}&=\frac{{{S}^{B}}\left( {{S}^{D}}-{{S}^{A}} \right)}{{{\left( {{S}^{D}}-{{S}^{A}} \right)}^{2}}+\left( {{S}^{B}}+{{S}^{C}} \right)\left( {{S}^{D}}-{{S}^{A}} \right)} \\ 
	& =\frac{{{S}^{B}}}{\left( {{S}^{D}}-{{S}^{A}} \right)+\left( {{S}^{B}}+{{S}^{C}} \right)}=\frac{{{S}^{B}}}{S-2{{S}^{A}}}.  
\end{aligned}\]

And
\[\begin{aligned}
	\beta _{C}^{{{E}_{A}}}&=\frac{\alpha _{C}^{E_{A}^{D}}\left( 1-\alpha _{D}^{E_{A}^{C}} \right)}{1-\alpha _{C}^{E_{A}^{D}}\alpha _{D}^{E_{A}^{C}}}=\frac{\frac{{{S}^{C}}}{{{S}^{B}}+{{S}^{C}}-{{S}^{A}}}\left( 1-\frac{{{S}^{D}}}{{{S}^{D}}+{{S}^{B}}-{{S}^{A}}} \right)}{1-\frac{{{S}^{C}}}{{{S}^{B}}+{{S}^{C}}-{{S}^{A}}}\frac{{{S}^{D}}}{{{S}^{D}}+{{S}^{B}}-{{S}^{A}}}} \\ 
	& =\frac{{{S}^{C}}\left( {{S}^{B}}-{{S}^{A}} \right)}{\left( {{S}^{B}}+{{S}^{C}}-{{S}^{A}} \right)\left( {{S}^{D}}+{{S}^{B}}-{{S}^{A}} \right)-{{S}^{C}}{{S}^{D}}},  
\end{aligned}\]
i.e.
\[\begin{aligned}
	\beta _{C}^{{{E}_{A}}}&=\frac{{{S}^{C}}\left( {{S}^{B}}-{{S}^{A}} \right)}{{{\left( {{S}^{B}}-{{S}^{A}} \right)}^{2}}+\left( {{S}^{D}}+{{S}^{C}} \right)\left( {{S}^{B}}-{{S}^{A}} \right)} \\ 
	& =\frac{{{S}^{C}}}{\left( {{S}^{B}}-{{S}^{A}} \right)+\left( {{S}^{D}}+{{S}^{C}} \right)}=\frac{{{S}^{C}}}{S-2{{S}^{A}}}.  
\end{aligned}\]

And
\[\begin{aligned}
	\beta _{D}^{{{E}_{A}}}&=\frac{\alpha _{D}^{E_{A}^{A}}\left( 1-\alpha _{A}^{E_{A}^{D}} \right)}{1-\alpha _{D}^{E_{A}^{A}}\alpha _{A}^{E_{A}^{D}}}=\frac{\frac{{{S}^{D}}}{S-{{S}^{A}}}\left( 1+\frac{{{S}^{A}}}{{{S}^{B}}+{{S}^{C}}-{{S}^{A}}} \right)}{1+\frac{{{S}^{D}}}{S-{{S}^{A}}}\frac{{{S}^{A}}}{{{S}^{B}}+{{S}^{C}}-{{S}^{A}}}} \\ 
	& =\frac{{{S}^{D}}\left( {{S}^{B}}+{{S}^{C}} \right)}{\left( {{S}^{B}}+{{S}^{C}}+{{S}^{D}} \right)\left( {{S}^{B}}+{{S}^{C}}-{{S}^{A}} \right)+{{S}^{D}}{{S}^{A}}},  
\end{aligned}\]
i.e.
\[\begin{aligned}
	\beta _{D}^{{{E}_{A}}}&=\frac{{{S}^{D}}\left( {{S}^{B}}+{{S}^{C}} \right)}{{{\left( {{S}^{B}}+{{S}^{C}} \right)}^{2}}+\left( {{S}^{B}}+{{S}^{C}} \right)\left( {{S}^{D}}-{{S}^{A}} \right)} \\ 
	& =\frac{{{S}^{D}}}{\left( {{S}^{B}}+{{S}^{C}} \right)+\left( {{S}^{D}}-{{S}^{A}} \right)}=\frac{{{S}^{D}}}{S-2{{S}^{A}}}.  
\end{aligned}\]

The frame components of the center of the escribed sphere in the trihedral angle of vertex $A$ of the tetrahedron have been studied above. Similarly, the frame components of the center of the escribed sphere in the trihedral angle of other vertexs can also be obtained. The following theorem is then obtained:

\begin{theorem}{Frame components for centers of escribed sphere, Daiyuan Zhang}{Thm21.3.1}\label{Thm21.3.1} 
	Given a tetrahedron $ABCD$, then its frame components of the centers of escribed spheres in the trihedral angles are	
\begin{flalign*}
	\beta _{A}^{{{E}_{A}}}=-\frac{{{S}^{A}}}{S-2{{S}^{A}}},\beta _{B}^{{{E}_{A}}}=\frac{{{S}^{B}}}{S-2{{S}^{A}}},\beta _{C}^{{{E}_{A}}}=\frac{{{S}^{C}}}{S-2{{S}^{A}}},\beta _{D}^{{{E}_{A}}}=\frac{{{S}^{D}}}{S-2{{S}^{A}}}.
\end{flalign*}
\begin{flalign*}
	\beta _{A}^{{{E}_{B}}}=\frac{{{S}^{A}}}{S-2{{S}^{B}}},\beta _{B}^{{{E}_{B}}}=-\frac{{{S}^{B}}}{S-2{{S}^{B}}},\beta _{C}^{{{E}_{B}}}=\frac{{{S}^{C}}}{S-2{{S}^{B}}},\beta _{D}^{{{E}_{B}}}=\frac{{{S}^{D}}}{S-2{{S}^{B}}}.	
\end{flalign*}
\begin{flalign*}
	\beta _{A}^{{{E}_{C}}}=\frac{{{S}^{A}}}{S-2{{S}^{C}}},\beta _{B}^{{{E}_{C}}}=\frac{{{S}^{B}}}{S-2{{S}^{C}}},\beta _{C}^{{{E}_{C}}}=-\frac{{{S}^{C}}}{S-2{{S}^{C}}},\beta _{D}^{{{E}_{C}}}=\frac{{{S}^{D}}}{S-2{{S}^{C}}}.	
\end{flalign*}
\begin{flalign*}
	\beta _{A}^{{{E}_{D}}}=\frac{{{S}^{A}}}{S-2{{S}^{D}}},\beta _{B}^{{{E}_{D}}}=\frac{{{S}^{B}}}{S-2{{S}^{D}}},\beta _{C}^{{{E}_{D}}}=\frac{{{S}^{C}}}{S-2{{S}^{D}}},\beta _{D}^{{{E}_{D}}}=-\frac{{{S}^{D}}}{S-2{{S}^{D}}}.	
\end{flalign*}	
\end{theorem}


\chapter{Special vectors of IC-T for a given tetrahedron}\label{Ch22}
\thispagestyle{empty}

This chapter is an important application of the results of previous chapters.

\section{Vector from origin to special intersecting center of a tetrahedron on tetrahedral frame}\label{Sec22.1}

\subsection{Centroid vector on tetrahedral frame}\label{Subsec22.1.1}
\begin{theorem}{Theorem of centroid vector on tetrahedral frame, Daiyuan Zhang}{Thm22.1.1}\label{Thm22.1.1} 
	If the point $O$ is any point and $G$ is the centroid of the tetrahedron $ABCD$, then the centroid vector $\overrightarrow{OG}$ can be uniquely expressed as follows:	
	\[\overrightarrow{OG}=\frac{1}{4}\left( \overrightarrow{OA}+\overrightarrow{OB}+\overrightarrow{OC}+\overrightarrow{OD} \right).\]
\end{theorem}

\begin{proof}
	This theorem is obtained by substituting the intersecting center $P$ of the tetrahedron in theorem \ref{thm:Thm18.2.1} with the centroid $G$ and substituting the frame components in theorem \ref{thm:Thm21.1.1} directly.
\end{proof}	
\hfill $\square$\par

\subsection{Incenter vector on tetrahedral frame}\label{Subsec22.1.2}
\begin{theorem}{Theorem of incenter vector on tetrahedral frame, Daiyuan Zhang}{Thm22.1.2}\label{Thm22.1.2} 
	If the point $O$ is any point and $I$ is the incenter of the tetrahedron $ABCD$, then the incenter vector $\overrightarrow{OI}$ can be uniquely expressed as follows:	
	\[\overrightarrow{OI}=\frac{1}{S}\left( {{S}^{A}}\overrightarrow{OA}+{{S}^{B}}\overrightarrow{OB}+{{S}^{C}}\overrightarrow{OC}+{{S}^{D}}\overrightarrow{OD} \right).\]	
\end{theorem}

\begin{proof}
	This theorem is obtained by substituting the intersecting center $P$ of the tetrahedron in theorem \ref{thm:Thm18.2.1} with the incenter $I$ and substituting the frame components in theorem \ref{thm:Thm21.2.1} directly.
\end{proof}
\hfill $\square$\par

\subsection{Excenter vector on tetrahedral frame}\label{Subsec22.1.3}
\begin{theorem}{Theorem of excenter vector in trihedral angle , Daiyuan Zhang}{Thm22.1.3}\label{Thm22.1.3} 
	If the point $O$ is any point, ${{E}_{A}}$, ${{E}_{B}}$, ${{E}_{C}}$, ${{E}_{D}}$ are the excenters in trihedral angles of vertexs $A$, $B$, $C$, $D$ of the tetrahedron, respectively, then the excenter vectors can be uniquely expressed in the following form:
	\[\overrightarrow{O{{E}_{A}}}=\frac{1}{S-2{{S}^{A}}}\left( -{{S}^{A}}\overrightarrow{OA}+{{S}^{B}}\overrightarrow{OB}+{{S}^{C}}\overrightarrow{OC}+{{S}^{D}}\overrightarrow{OD} \right),\]
	\[\overrightarrow{O{{E}_{B}}}=\frac{1}{S-2{{S}^{B}}}\left( {{S}^{A}}\overrightarrow{OA}-{{S}^{B}}\overrightarrow{OB}+{{S}^{C}}\overrightarrow{OC}+{{S}^{D}}\overrightarrow{OD} \right),\]	
	\[\overrightarrow{O{{E}_{C}}}=\frac{1}{S-2{{S}^{C}}}\left( {{S}^{A}}\overrightarrow{OA}+{{S}^{B}}\overrightarrow{OB}-{{S}^{C}}\overrightarrow{OC}+{{S}^{D}}\overrightarrow{OD} \right),\]	
	\[\overrightarrow{O{{E}_{D}}}=\frac{1}{S-2{{S}^{D}}}\left( {{S}^{A}}\overrightarrow{OA}+{{S}^{B}}\overrightarrow{OB}+{{S}^{C}}\overrightarrow{OC}-{{S}^{D}}\overrightarrow{OD} \right).\]		
\end{theorem}

\begin{proof}
	This theorem is obtained by substituting the intersecting center $P$ in theorem \ref{thm:Thm18.2.1} with the excenters ${{E}_{A}}$, ${{E}_{B}}$, ${{E}_{C}}$, ${{E}_{D}}$ in trihedral angles of vertexs $A$, $B$, $C$, $D$ of the tetrahedron $ABCD$, respectively, and substituting the frame components in theorem \ref{thm:Thm21.3.1} directly.
\end{proof}
\hfill $\square$\par

\section{Vector of two special intersecting centers of a tetrahedron on tetrahedral frame}\label{Sec22.2}
\subsection{Vector of two intersecting centers between centroid and incenter}\label{Subsec22.2.1}
\begin{theorem}{Theorem of vector from centroid to incenter of a tetrahedron, Daiyuan Zhang}{Thm22.2.1}\label{Thm22.2.1} 
	Given the tetrahedron $ABCD$ and let the point $O$ be an arbitrary point, then the vector $\overrightarrow{GI}$ from the centroid $G$ to the incenter $I$ can be expressed by the uniquely linear combination of the tetrahedral frame $\left( O;A,B,C,D \right)$, i.e.	
	\[\overrightarrow{GI}=\frac{1}{4S}\left( \begin{aligned}
		\left( 4{{S}^{A}}-S \right)\overrightarrow{OA}+\left( 4{{S}^{B}}-S \right)\overrightarrow{OB}+\left( 4{{S}^{C}}-S \right)\overrightarrow{OC}+\left( 4{{S}^{D}}-S \right)\overrightarrow{OD} \\ 
	\end{aligned} \right).\]	
\end{theorem}
\begin{proof}
	According to theorem \ref{thm:Thm18.3.1}, theorem \ref{thm:Thm21.1.1} and theorem \ref{thm:Thm21.2.1}, the following results are obtained:
	\[\overrightarrow{GI}=\beta _{A}^{GI}\overrightarrow{OA}+\beta _{B}^{GI}\overrightarrow{OB}+\beta _{C}^{GI}\overrightarrow{OC}+\beta _{D}^{GI}\overrightarrow{OD},\]	
	where
	\[\beta _{A}^{GI}=\beta _{A}^{I}-\beta _{A}^{G}=\frac{{{S}^{A}}}{S}-\frac{1}{4}=\frac{4{{S}^{A}}-S}{4S},\]	
	\[\beta _{B}^{GI}=\beta _{B}^{I}-\beta _{B}^{G}=\frac{{{S}^{B}}}{S}-\frac{1}{4}=\frac{4{{S}^{B}}-S}{4S},\]	
	\[\beta _{C}^{GI}=\beta _{C}^{I}-\beta _{C}^{G}=\frac{{{S}^{C}}}{S}-\frac{1}{4}=\frac{4{{S}^{C}}-S}{4S},\]	
	\[\beta _{D}^{GI}=\beta _{D}^{I}-\beta _{D}^{G}=\frac{{{S}^{D}}}{S}-\frac{1}{4}=\frac{4{{S}^{D}}-S}{4S}.\]	
	
	Substitute the above formula to get the result.
\end{proof}
\hfill $\square$\par

\subsection{Vector of two intersecting centers between centroid and excenter}\label{Subsec22.2.2}

\begin{theorem}{Theorem of vector from centroid to excenter of a tetrahedron, Daiyuan Zhang}{Thm22.2.2}\label{Thm22.2.2} 
	Suppose that given the tetrahedron $ABCD$ and the point $O$ is an arbitrary point, ${{E}_{A}}$, ${{E}_{B}}$, ${{E}_{C}}$, ${{E}_{D}}$ are the excenters in trihedral angles of vertexs $A$, $B$, $C$, $D$ of the tetrahedron respectively, then the vector from the centroid $G$ to each of the excenter ${{E}_{A}}$, ${{E}_{B}}$, ${{E}_{C}}$, ${{E}_{D}}$ can be expressed by the uniquely linear combination of the tetrahedral frame $\left( O;A,B,C,D \right)$, i.e.			
	\begin{flalign*}
		\overrightarrow{G{{E}_{A}}}=\frac{{{\mathbf{J}}_{A}}}{4\left( S-2{{S}^{A}} \right)},\overrightarrow{G{{E}_{B}}}=\frac{{{\mathbf{J}}_{B}}}{4\left( S-2{{S}^{B}} \right)},
	\end{flalign*}
	\begin{flalign*}
		\overrightarrow{G{{E}_{C}}}=\frac{{{\mathbf{J}}_{C}}}{4\left( S-2{{S}^{C}} \right)},\overrightarrow{G{{E}_{D}}}=\frac{{{\mathbf{J}}_{D}}}{4\left( S-2{{S}^{D}} \right)}.
	\end{flalign*}
	Where
	\[\begin{aligned}
		{{\mathbf{J}}_{A}}&=-\left( S+2{{S}^{A}} \right)\overrightarrow{OA}+\left( 4{{S}^{B}}+2{{S}^{A}}-S \right)\overrightarrow{OB} \\ 
		& +\left( 4{{S}^{C}}+2{{S}^{A}}-S \right)\overrightarrow{OC}+\left( 4{{S}^{D}}+2{{S}^{A}}-S \right)\overrightarrow{OD},  
	\end{aligned}\]	
	\[\begin{aligned}
		{{\mathbf{J}}_{B}}&=-\left( S+2{{S}^{B}} \right)\overrightarrow{OB}+\left( 4{{S}^{C}}+2{{S}^{B}}-S \right)\overrightarrow{OC} \\ 
		& +\left( 4{{S}^{D}}+2{{S}^{B}}-S \right)\overrightarrow{OD}+\left( 4{{S}^{A}}+2{{S}^{B}}-S \right)\overrightarrow{OA},  
	\end{aligned}\]	
	\[\begin{aligned}
		{{\mathbf{J}}_{C}}&=-\left( S+2{{S}^{C}} \right)\overrightarrow{OC}+\left( 4{{S}^{D}}+2{{S}^{C}}-S \right)\overrightarrow{OD} \\ 
		& +\left( 4{{S}^{A}}+2{{S}^{C}}-S \right)\overrightarrow{OA}+\left( 4{{S}^{B}}+2{{S}^{C}}-S \right)\overrightarrow{OB},  
	\end{aligned}\]	
	\[\begin{aligned}
		{{\mathbf{J}}_{D}}&=-\left( S+2{{S}^{D}} \right)\overrightarrow{OD}+\left( 4{{S}^{A}}+2{{S}^{D}}-S \right)\overrightarrow{OA} \\ 
		& +\left( 4{{S}^{B}}+2{{S}^{D}}-S \right)\overrightarrow{OB}+\left( 4{{S}^{C}}+2{{S}^{D}}-S \right)\overrightarrow{OC}.  
	\end{aligned}\]	
\end{theorem}

\begin{proof}
	First consider the vector from the centroid $G$ to the excenter ${{E}_{A}}$ in trihedral angle of vertex $A$ (excenter in trihedral angle of vertex $A$). From theorem \ref{thm:Thm18.3.1}, theorem \ref{thm:Thm21.1.1} and theorem \ref{thm:Thm21.3.1}, the vector $\overrightarrow{G{{E}_{A}}}$ from the centroid $G$ to the excenter ${{E}_{A}}$ can be expressed by the uniquely linear combination in the following, i.e.			
	\[\overrightarrow{G{{E}_{A}}}=\beta _{A}^{G{{E}_{A}}}\overrightarrow{OA}+\beta _{B}^{G{{E}_{A}}}\overrightarrow{OB}+\beta _{C}^{G{{E}_{A}}}\overrightarrow{OC}+\beta _{D}^{G{{E}_{A}}}\overrightarrow{OD}.\]	
	Where
	\[\beta _{A}^{G{{E}_{A}}}=\beta _{A}^{{{E}_{A}}}-\beta _{A}^{G}=-\frac{{{S}^{A}}}{S-2{{S}^{A}}}-\frac{1}{4}=-\frac{S+2{{S}^{A}}}{4\left( S-2{{S}^{A}} \right)},\]	
	\[\beta _{B}^{G{{E}_{A}}}=\beta _{B}^{{{E}_{A}}}-\beta _{B}^{G}=\frac{{{S}^{B}}}{S-2{{S}^{A}}}-\frac{1}{4}=\frac{4{{S}^{B}}+2{{S}^{A}}-S}{4\left( S-2{{S}^{A}} \right)},\]	
	\[\beta _{C}^{G{{E}_{A}}}=\beta _{C}^{{{E}_{A}}}-\beta _{C}^{G}=\frac{{{S}^{C}}}{S-2{{S}^{A}}}-\frac{1}{4}=\frac{4{{S}^{C}}+2{{S}^{A}}-S}{4\left( S-2{{S}^{A}} \right)},\]	
	\[\beta _{D}^{G{{E}_{A}}}=\beta _{D}^{{{E}_{A}}}-\beta _{D}^{G}=\frac{{{S}^{D}}}{S-2{{S}^{A}}}-\frac{1}{4}=\frac{4{{S}^{D}}+2{{S}^{A}}-S}{4\left( S-2{{S}^{A}} \right)}.\]	
	
	Substitute the previous formula to get the desired result. Similarly, the formulas of excenters in other trihedral angles of vertexs can be obtained.	
\end{proof}
\hfill $\square$\par

\subsection{Vector of two intersecting centers between incenter and excenter}\label{Subsec22.2.3}
\begin{theorem}{Theorem of vector from centroid to incenter of a tetrahedron, Daiyuan Zhang}{Thm22.2.3}\label{Thm22.2.3} 
	Suppose that given the tetrahedron $ABCD$ and the point $O$ is an arbitrary point, ${{E}_{A}}$, ${{E}_{B}}$, ${{E}_{C}}$, ${{E}_{D}}$ are the excenters in trihedral angles of vertexs $A$, $B$, $C$, $D$ of the tetrahedron, respectively, then the vector from the incenter $I$ to each of the excenter ${{E}_{A}}$, ${{E}_{B}}$, ${{E}_{C}}$, ${{E}_{D}}$ can be expressed by the uniquely linear combination of the tetrahedral frame $\left( O;A,B,C,D \right)$, i.e.
	\[\overrightarrow{I{{E}_{A}}}=\frac{2{{S}^{A}}}{S\left( S-2{{S}^{A}} \right)}\left( \left( {{S}^{A}}-S \right)\overrightarrow{OA}+{{S}^{B}}\overrightarrow{OB}+{{S}^{C}}\overrightarrow{OC}+{{S}^{D}}\overrightarrow{OD} \right),\]	
	\[\overrightarrow{I{{E}_{B}}}=\frac{2{{S}^{B}}}{S\left( S-2{{S}^{B}} \right)}\left( \left( {{S}^{B}}-S \right)\overrightarrow{OB}+{{S}^{C}}\overrightarrow{OC}+{{S}^{D}}\overrightarrow{OD}+{{S}^{A}}\overrightarrow{OA} \right),\]	
	\[\overrightarrow{I{{E}_{C}}}=\frac{2{{S}^{C}}}{S\left( S-2{{S}^{C}} \right)}\left( \left( {{S}^{C}}-S \right)\overrightarrow{OC}+{{S}^{D}}\overrightarrow{OD}+{{S}^{A}}\overrightarrow{OA}+{{S}^{B}}\overrightarrow{OB} \right),\]	
	\[\overrightarrow{I{{E}_{D}}}=\frac{2{{S}^{D}}}{S\left( S-2{{S}^{D}} \right)}\left( \left( {{S}^{D}}-S \right)\overrightarrow{OD}+{{S}^{A}}\overrightarrow{OA}+{{S}^{B}}\overrightarrow{OB}+{{S}^{C}}\overrightarrow{OC} \right).\]	
\end{theorem}

\begin{proof}
	First consider the vector from the incenter $I$ to the excenter ${{E}_{A}}$ in trihedral angle of vertex $A$ (excenter in trihedral angle of vertex $A$). From theorem \ref{thm:Thm18.3.1}, theorem \ref{thm:Thm21.2.1} and theorem \ref{thm:Thm21.3.1}, the vector $\overrightarrow{I{{E}_{A}}}$ from the incenter $I$ to the excenter ${{E}_{A}}$ can be expressed by the uniquely linear combination in the following, i.e.	
	\[\overrightarrow{I{{E}_{A}}}=\beta _{A}^{I{{E}_{A}}}\overrightarrow{OA}+\beta _{B}^{I{{E}_{A}}}\overrightarrow{OB}+\beta _{C}^{I{{E}_{A}}}\overrightarrow{OC}+\beta _{D}^{I{{E}_{A}}}\overrightarrow{OD}.\]	
	Where
	\[\beta _{A}^{I{{E}_{A}}}=\beta _{A}^{{{E}_{A}}}-\beta _{A}^{I}=-\frac{{{S}^{A}}}{S-2{{S}^{A}}}-\frac{{{S}^{A}}}{S}=\frac{2{{S}^{A}}\left( {{S}^{A}}-S \right)}{S\left( S-2{{S}^{A}} \right)},\]	
	\[\beta _{B}^{I{{E}_{A}}}=\beta _{B}^{{{E}_{A}}}-\beta _{B}^{I}=\frac{{{S}^{B}}}{S-2{{S}^{A}}}-\frac{{{S}^{B}}}{S}=\frac{2{{S}^{A}}{{S}^{B}}}{S\left( S-2{{S}^{A}} \right)},\]	
	\[\beta _{C}^{I{{E}_{A}}}=\beta _{C}^{{{E}_{A}}}-\beta _{C}^{I}=\frac{{{S}^{C}}}{S-2{{S}^{A}}}-\frac{{{S}^{C}}}{S}=\frac{2{{S}^{A}}{{S}^{C}}}{S\left( S-2{{S}^{A}} \right)},\]	
	\[\beta _{D}^{I{{E}_{A}}}=\beta _{D}^{{{E}_{A}}}-\beta _{D}^{I}=\frac{{{S}^{D}}}{S-2{{S}^{A}}}-\frac{{{S}^{D}}}{S}=\frac{2{{S}^{A}}{{S}^{D}}}{S\left( S-2{{S}^{A}} \right)}.\]	
	
	Substitute the previous formula to get the following result:
	\[\overrightarrow{I{{E}_{A}}}=\frac{2{{S}^{A}}}{S\left( S-2{{S}^{A}} \right)}\left( \left( {{S}^{A}}-S \right)\overrightarrow{OA}+{{S}^{B}}\overrightarrow{OB}+{{S}^{C}}\overrightarrow{OC}+{{S}^{D}}\overrightarrow{OD} \right).\]	
	Similarly, the formulas of excenters in trihedral angles of other vertexs can be obtained.
\end{proof}
\hfill $\square$\par

It can be seen that the frame components of the centroid, incenter and excenter of the tetrahedron are all constants, which are determined by the triangle areas of each face of the tetrahedron, and has nothing to do with the frame $\left( O;A,B,C,D \right)$.

\begin{example}\label{Exam22.2.1} 
	Suppose there is a regular pyramid $A-BCD$, the length of each base edge is $a=BC=CD=DB=2$, and the length of each face edge is $l=AB=AC=AD=3$, find the vectors $\overrightarrow{OG}$ and $\overrightarrow{OI}$, i.e., the vectors from any point $O$ of space to the centroid and to the center of inscribed sphere of the regular pyramid; and find the vector $\overrightarrow{GI}$ (vector from the centroid to the incenter)  on the tetrahedral frame and on the centroid frame, respectively.
\end{example}

\begin{solution}
	According to Helen's formula, the following results are obtained for the base triangle of:
	\[p=\frac{1}{2}\left( a+a+a \right)=\frac{1}{2}\left( 2+2+2 \right)=3,\]
	\[\begin{aligned}
		{{S}^{A}}&=\sqrt{p\left( p-a \right)\left( p-a \right)\left( p-a \right)} \\ 
		& =\sqrt{3\times \left( 3-2 \right)\times \left( 3-2 \right)\times \left( 3-2 \right)}=\sqrt{3}.  
	\end{aligned}\]
	
	For the base triangles, the following results are found:
	\[p=\frac{1}{2}\left( l+l+a \right)=\frac{1}{2}\left( 3+3+2 \right)=4,\]
	\[\begin{aligned}
		{{S}^{B}}={{S}^{C}}&={{S}^{D}}=\sqrt{p\left( p-l \right)\left( p-l \right)\left( p-a \right)} \\ 
		& =\sqrt{4\times \left( 4-3 \right)\times \left( 4-3 \right)\times \left( 4-2 \right)}=2\sqrt{2}.  
	\end{aligned}\]
	
	According to theorem \ref{thm:Thm22.1.1}, the following results can be obtained:
	\[\overrightarrow{OG}=\frac{1}{4}\left( \overrightarrow{OA}+\overrightarrow{OB}+\overrightarrow{OC}+\overrightarrow{OD} \right).\]
	
	According to theorem \ref{thm:Thm22.1.2}, the following results can be obtained:
	\begin{align*}
		\overrightarrow{OI}&=\frac{1}{S}\left( {{S}^{A}}\overrightarrow{OA}+{{S}^{B}}\overrightarrow{OB}+{{S}^{C}}\overrightarrow{OC}+{{S}^{D}}\overrightarrow{OD} \right) \\ 
		& =\frac{\sqrt{3}}{\sqrt{3}+6\sqrt{2}}\overrightarrow{OA}+\frac{2\sqrt{2}}{\sqrt{3}+6\sqrt{2}}\left( \overrightarrow{OB}+\overrightarrow{OC}+\overrightarrow{OD} \right).  
	\end{align*}	
	
	The surface area of the tetrahedron is: $S=\sqrt{3}+3\times 2\sqrt{2}=\sqrt{3}+6\sqrt{2}$. From theorem \ref{thm:Thm22.2.1}, the following results can be obtained on the tetrahedral frame:	
	\begin{align*}
		\overrightarrow{GI}&=\frac{1}{4S}\left( \begin{aligned}
			& \left( 4{{S}^{A}}-S \right)\overrightarrow{OA}+\left( 4{{S}^{B}}-S \right)\overrightarrow{OB} \\ 
			& +\left( 4{{S}^{C}}-S \right)\overrightarrow{OC}+\left( 4{{S}^{D}}-S \right)\overrightarrow{OD} \\ 
		\end{aligned} \right) \\ 
		& =\frac{1}{4\left( \sqrt{3}+6\sqrt{2} \right)}\left( \begin{aligned}
			& \left( 4\sqrt{3}-\left( \sqrt{3}+6\sqrt{2} \right) \right)\overrightarrow{OA}+\left( 16\sqrt{2}-\left( \sqrt{3}+6\sqrt{2} \right) \right)\overrightarrow{OB} \\ 
			& +\left( 16\sqrt{2}-\left( \sqrt{3}+6\sqrt{2} \right) \right)\overrightarrow{OC}+\left( 16\sqrt{2}-\left( \sqrt{3}+6\sqrt{2} \right) \right)\overrightarrow{OD} \\ 
		\end{aligned} \right) \\ 
		& =\frac{1}{4\left( \sqrt{3}+6\sqrt{2} \right)}\left( 3\left( \sqrt{3}-2\sqrt{2} \right)\overrightarrow{OA}+\left( 10\sqrt{2}-\sqrt{3} \right)\left( \overrightarrow{OB}+\overrightarrow{OC}+\overrightarrow{OD} \right) \right).  
	\end{align*}	
		
	In the case of centroid frame, using the results on the tetrahedral frame, replace $O$ with $G$ to get the following result:
		\[\begin{aligned}
			\overrightarrow{GI}&=\frac{1}{S}\left( {{S}^{A}}\overrightarrow{GA}+{{S}^{B}}\overrightarrow{GB}+{{S}^{C}}\overrightarrow{GC}+{{S}^{D}}\overrightarrow{GD} \right) \\ 
			& =\frac{\sqrt{3}}{\sqrt{3}+6\sqrt{2}}\overrightarrow{GA}+\frac{2\sqrt{2}}{\sqrt{3}+6\sqrt{2}}\left( \overrightarrow{GB}+\overrightarrow{GC}+\overrightarrow{GD} \right).  
		\end{aligned}\]
	
	In the next chapter, the frame equation of centroid will be given as follows:
	\[\overrightarrow{GA}+\overrightarrow{GB}+\overrightarrow{GC}+\overrightarrow{GD}=\overrightarrow{0}.\]
	
	The frame equation of centroid is substituted into the formula above to get:	
	\[\begin{aligned}
		\overrightarrow{GI}&=\frac{\sqrt{3}}{\sqrt{3}+6\sqrt{2}}\overrightarrow{GA}+\frac{2\sqrt{2}}{\sqrt{3}+6\sqrt{2}}\left( \overrightarrow{GB}+\overrightarrow{GC}+\overrightarrow{GD} \right) \\ 
		& =\frac{\sqrt{3}}{\sqrt{3}+6\sqrt{2}}\overrightarrow{GA}-\frac{2\sqrt{2}}{\sqrt{3}+6\sqrt{2}}\overrightarrow{GA}=\frac{\sqrt{3}-2\sqrt{2}}{\sqrt{3}+6\sqrt{2}}\overrightarrow{GA}.  
	\end{aligned}\]
	
	The above formula shows that the centroid $G$ is located between the incenter $I$ and the vertex $A$ of the regular pyramid ($\sqrt{3}-2\sqrt{2}<0$), and that the centroid $G$, the incenter $I$ and the vertex $A$ are in the same straight line.
\end{solution}
\hfill $\diamond$\par

Generally speaking, in application, if the origin of the frame is also an IC-T, the selection principle for the starting IC-T and the ending IC-T is: the position of the starting IC-T should make the vector from the starting IC-T to the vertexs of the tetrahedron easy to express and calculate; The position of the ending IC-T should make the frame components of the ending IC-T easy to calculate.

\section{Vector from origin to intersecting center of a tetrahedron on circumcenter frame}\label{Sec22.3}
\subsection{Centroid vector on circumcenter frame}\label{Subsec22.3.1}

\begin{theorem}{Representation of centroid vector on circumcenter frame, Daiyuan Zhang}{Thm22.3.1}\label{Thm22.3.1} 
	Let $G$ and $Q$ be the centroid and circumcenter of the tetrahedron $ABCD$ respectively, then the centroid vector $\overrightarrow{QG}$ on circumcenter frame is 
	\[\overrightarrow{QG}=\frac{1}{4}\left( \overrightarrow{QA}+\overrightarrow{QB}+\overrightarrow{QC}+\overrightarrow{QD} \right).\]	
\end{theorem}

\begin{proof}
	This theorem is obtained by substituting the intersecting center $P$ of the tetrahedron in theorem \ref{thm:Thm18.2.1} with the centroid $G$ and circumcenter $Q$ respectively, and substituting the frame components in theorem \ref{thm:Thm21.1.1} directly.
\end{proof}
\hfill $\square$\par

\subsection{Incenter vector on circumcenter frame}\label{Subsec22.3.2}

\begin{theorem}{Representation of incenter vector on circumcenter frame, Daiyuan Zhang}{Thm22.3.2}\label{Thm22.3.2} 
	Let $I$ and $Q$ be the incenter and circumcenter of the tetrahedron $ABCD$ respectively, then the incenter vector $\overrightarrow{QI}$ on circumcenter frame is: 	
	\[\overrightarrow{QI}=\frac{1}{S}\left( {{S}^{A}}\overrightarrow{QA}+{{S}^{B}}\overrightarrow{QB}+{{S}^{C}}\overrightarrow{QC}+{{S}^{D}}\overrightarrow{QD} \right).\]	
\end{theorem}

\begin{proof}
	This theorem is obtained by substituting the intersecting center $P$ of the tetrahedron in theorem \ref{thm:Thm18.2.1} with the incenter $I$ and circumcenter $Q$ respectively, and substituting the frame components in theorem \ref{thm:Thm21.2.1} directly.
\end{proof}
\hfill $\square$\par

\subsection{Excenter vector on circumcenter frame}\label{Subsec22.3.3}

\begin{theorem}{Representationof excenter vector on circumcenter frame, Daiyuan Zhang}{Thm22.3.3}\label{Thm22.3.3} 
	Let point $Q$ be the circumcenter of tetrahedron $ABCD$, ${{E}_{A}}$, ${{E}_{B}}$, ${{E}_{C}}$, ${{E}_{D}}$ are the excenters in trihedral angles of vertexs $A$, $B$, $C$, $D$ of the tetrahedron, respectively, then the excenter vectors are		
	\[\overrightarrow{Q{{E}_{A}}}=\frac{1}{S-2{{S}^{A}}}\left( -{{S}^{A}}\overrightarrow{QA}+{{S}^{B}}\overrightarrow{QB}+{{S}^{C}}\overrightarrow{QC}+{{S}^{D}}\overrightarrow{QD} \right),\]
	\[\overrightarrow{Q{{E}_{B}}}=\frac{1}{S-2{{S}^{B}}}\left( {{S}^{A}}\overrightarrow{QA}-{{S}^{B}}\overrightarrow{QB}+{{S}^{C}}\overrightarrow{QC}+{{S}^{D}}\overrightarrow{QD} \right),\]	
	\[\overrightarrow{Q{{E}_{C}}}=\frac{1}{S-2{{S}^{C}}}\left( {{S}^{A}}\overrightarrow{QA}+{{S}^{B}}\overrightarrow{QB}-{{S}^{C}}\overrightarrow{QC}+{{S}^{D}}\overrightarrow{QD} \right),\]	
	\[\overrightarrow{Q{{E}_{D}}}=\frac{1}{S-2{{S}^{D}}}\left( {{S}^{A}}\overrightarrow{QA}+{{S}^{B}}\overrightarrow{QB}+{{S}^{C}}\overrightarrow{QC}-{{S}^{D}}\overrightarrow{QD} \right).\]
\end{theorem}

\begin{proof}
	This theorem is obtained by substituting the intersecting center $P$ of the tetrahedron in theorem \ref{thm:Thm18.2.1} with the excenters ${{E}_{A}}$, ${{E}_{B}}$, ${{E}_{C}}$, ${{E}_{D}}$ and circumcenter $Q$ respectively, and substituting the frame components in theorem \ref{thm:Thm21.3.1} directly.
\end{proof}
\hfill $\square$\par


\chapter{Frame equation of special intersecting centers of a tetrahedron}\label{Ch23}
\thispagestyle{empty}

In some applications, it is necessary to consider the frame equations of some special intersecting centers. In this chapter, the frame equations of the centroid and the incenter of a tetrahedron are given.

\section{Frame equation of centroid of a tetrahedron}\label{Sec23.1}
\begin{theorem}{Frame equation of centroid of a tetrahedron, Daiyuan Zhang}{Thm23.1.1}\label{Thm23.1.1} 
	Suppose that $G$ is the centroid of tetrahedron $ABCD$, then:	
	\[\overrightarrow{GA}+\overrightarrow{GB}+\overrightarrow{GC}+\overrightarrow{GD}=\overrightarrow{0}.\]	
\end{theorem}

\begin{proof}
	In theorem \ref{thm:Thm18.2.1}, let $P$ and $O$ coincide with the centroid of tetrahedron $ABCD$ and then use the frame components of the centroid of tetrahedron (theorem \ref{thm:Thm21.1.1}), then	
	\[\frac{1}{4}\overrightarrow{GA}+\frac{1}{4}\overrightarrow{GB}+\frac{1}{4}\overrightarrow{GC}+\frac{1}{4}\overrightarrow{GD}=\overrightarrow{0},\]	
	i.e.
	\[\overrightarrow{GA}+\overrightarrow{GB}+\overrightarrow{GC}+\overrightarrow{GD}=\overrightarrow{0}.\]	
\end{proof}
\hfill $\square$\par

\section{Frame equation of incenter of a tetrahedron}\label{Sec23.2}
\begin{theorem}{Frame equation of incenter of a tetrahedron, Daiyuan Zhang}{Thm23.2.1}\label{Thm23.2.1} 
	Suppose that $I$ is the incenter of tetrahedron $ABCD$, then	
	\[{{S}^{A}}\overrightarrow{IA}+{{S}^{B}}\overrightarrow{IB}+{{S}^{C}}\overrightarrow{IC}+{{S}^{D}}\overrightarrow{ID}=\overrightarrow{0}.\]	
\end{theorem}

\begin{proof}
	In theorem \ref{thm:Thm18.2.1}, let $P$ and $O$ coincide with the incenter of tetrahedron $ABCD$ and then use the frame components of the incenter of tetrahedron (theorem \ref{thm:Thm21.2.1}), then	
	\[\frac{{{S}^{A}}}{S}\overrightarrow{IA}+\frac{{{S}^{B}}}{S}\overrightarrow{IB}+\frac{{{S}^{C}}}{S}\overrightarrow{IC}+\frac{{{S}^{D}}}{S}\overrightarrow{ID}=\overrightarrow{0},\]	
	i.e.
	\[{{S}^{A}}\overrightarrow{IA}+{{S}^{B}}\overrightarrow{IB}+{{S}^{C}}\overrightarrow{IC}+{{S}^{D}}\overrightarrow{ID}=\overrightarrow{0}.\]	
\end{proof}
\hfill $\square$\par

\section{Frame equation of excenter of a tetrahedron}\label{Sec23.3}
\begin{theorem}{Frame equation of excenter of a tetrahedron, Daiyuan Zhang}{Thm23.3.1}\label{Thm23.3.1} 
Let ${{E}_{A}}$, ${{E}_{B}}$, ${{E}_{C}}$, ${{E}_{D}}$ be the excenters in trihedral angles of vertexs $A$, $B$, $C$, $D$ of tetrahedron $ABCD$ respectively, then:
\[-{{S}^{A}}\overrightarrow{{{E}_{A}}A}+{{S}^{B}}\overrightarrow{{{E}_{A}}B}+{{S}^{C}}\overrightarrow{{{E}_{A}}C}+{{S}^{D}}\overrightarrow{{{E}_{A}}D}=\overrightarrow{0},\]
\[{{S}^{A}}\overrightarrow{{{E}_{B}}A}-{{S}^{B}}\overrightarrow{{{E}_{B}}B}+{{S}^{C}}\overrightarrow{{{E}_{B}}C}+{{S}^{D}}\overrightarrow{{{E}_{B}}D}=\overrightarrow{0},\]	
\[{{S}^{A}}\overrightarrow{{{E}_{C}}A}+{{S}^{B}}\overrightarrow{{{E}_{C}}B}-{{S}^{C}}\overrightarrow{{{E}_{C}}C}+{{S}^{D}}\overrightarrow{{{E}_{C}}D}=\overrightarrow{0},\]	
\[{{S}^{A}}\overrightarrow{{{E}_{D}}A}+{{S}^{B}}\overrightarrow{{{E}_{D}}B}+{{S}^{C}}\overrightarrow{{{E}_{D}}C}-{{S}^{D}}\overrightarrow{{{E}_{D}}D}=\overrightarrow{0}.\]	
\end{theorem}

\begin{proof}
	By theorem \ref{thm:Thm22.3.3}, this theorem is obtained by making the circumcenter $Q$ coincide with ${{E}_{A}}$, ${{E}_{B}}$, ${{E}_{C}}$, ${{E}_{D}}$ respectively.
\end{proof}
\hfill $\square$\par


\chapter{Distance between two points in space}\label{Ch24}
\thispagestyle{empty}

\section{Distance between origin and intersecting center of a tetrahedron on tetrahedral frame}\label{Sec24.1}

This section studies an important theorem. Given the tetrahedron $ABCD$ and a point $P$ (IC-T), select the point $\,O$ as the origin of the frame $\left( O;A,B,C,D \right)$, then the distance between the two points $O$ and $P$ can be obtained. Because the position of the point $O$ is arbitrary, the point $O$ is not necessarily the IC-T. The distance between the two points is called the distance between origin and intersecting center of the tetrahedron on the tetrahedral frame, which is called the distance between origin and intersecting center for short. If we choose the IC-T at some special points of the tetrahedron, we will get some special conclusions. In the previous discussion of the triangular frame, there is also the concept of the distance between origin and intersecting center. Please distinguish them according to the context.

\begin{theorem}{Distance between origin and IC-T on tetrahedral frame, Daiyuan Zhang}{Thm24.1.1}\label{Thm24.1.1} 
	Given the tetrahedron $ABCD$, let point $\,O$ be the origin of frame $\left( O;A,B,C,D \right)$, and point $P$ be the intersecting center of the tetrahedron, then	
	\[\begin{aligned}
		O{{P}^{2}}&=\beta _{A}^{P}O{{A}^{2}}+\beta _{B}^{P}O{{B}^{2}}+\beta _{C}^{P}O{{C}^{2}}+\beta _{D}^{P}O{{D}^{2}} \\ 
		& -\beta _{A}^{P}\beta _{B}^{P}A{{B}^{2}}-\beta _{A}^{P}\beta _{C}^{P}A{{C}^{2}}-\beta _{A}^{P}\beta _{D}^{P}A{{D}^{2}} \\ 
		& -\beta _{B}^{P}\beta _{C}^{P}B{{C}^{2}}-\beta _{B}^{P}\beta _{D}^{P}B{{D}^{2}}-\beta _{C}^{P}\beta _{D}^{P}C{{D}^{2}}.  
	\end{aligned}\]	
	Where $\beta _{A}^{P}$, $\beta _{B}^{P}$, $\beta _{C}^{P}$, $\beta _{D}^{P}$ are the frame components of $\overrightarrow{OA}$, $\overrightarrow{OB}$, $\overrightarrow{OC}$, $\overrightarrow{OD}$ at point $P$ respectively on the frame $\left( O;A,B,C \right)$.
\end{theorem}

\begin{proof}
	According to theorem \ref{thm:Thm18.2.1}, it is known that:
	\[\overrightarrow{OP}=\beta _{A}^{P}\overrightarrow{OA}+\beta _{B}^{P}\overrightarrow{OB}+\beta _{C}^{P}\overrightarrow{OC}+\beta _{D}^{P}\overrightarrow{OD}.\]	
	
	Therefore, the distance from any point $\,O$ to the intersecting center $P$ is
	\[\begin{aligned}
		O{{P}^{2}}&=\left( \beta _{A}^{P}\overrightarrow{OA}+\beta _{B}^{P}\overrightarrow{OB}+\beta _{C}^{P}\overrightarrow{OC}+\beta _{D}^{P}\overrightarrow{OD} \right) \cdot \left( \beta _{A}^{P}\overrightarrow{OA}+\beta _{B}^{P}\overrightarrow{OB}+\beta _{C}^{P}\overrightarrow{OC}+\beta _{D}^{P}\overrightarrow{OD} \right),  
	\end{aligned}\]
	i.e.
	\[\begin{aligned}
		O{{P}^{2}}&={{\left( \beta _{A}^{P} \right)}^{2}}O{{A}^{2}}+{{\left( \beta _{B}^{P} \right)}^{2}}O{{B}^{2}}+{{\left( \beta _{C}^{P} \right)}^{2}}O{{C}^{2}}+{{\left( \beta _{D}^{P} \right)}^{2}}O{{D}^{2}} \\ 
		& +2\beta _{A}^{P}\beta _{B}^{P}\overrightarrow{OA}\cdot \overrightarrow{OB}+2\beta _{A}^{P}\beta _{C}^{P}\overrightarrow{OA}\cdot \overrightarrow{OC}+2\beta _{A}^{P}\beta _{D}^{P}\overrightarrow{OA}\cdot \overrightarrow{OD} \\ 
		& +2\beta _{B}^{P}\beta _{C}^{P}\overrightarrow{OB}\cdot \overrightarrow{OC}+2\beta _{B}^{P}\beta _{D}^{P}\overrightarrow{OB}\cdot \overrightarrow{OD}+2\beta _{C}^{P}\beta _{D}^{P}\overrightarrow{OC}\cdot \overrightarrow{OD}.  
	\end{aligned}\]
	
	Using the cosine theorem in vector form, the following results can be obtained:
	\[2\overrightarrow{OA}\cdot \overrightarrow{OB}=O{{A}^{2}}+O{{B}^{2}}-A{{B}^{2}},\]
	\[2\overrightarrow{OA}\cdot \overrightarrow{OC}=O{{A}^{2}}+O{{C}^{2}}-A{{C}^{2}},\]
	\[2\overrightarrow{OA}\cdot \overrightarrow{OD}=O{{A}^{2}}+O{{D}^{2}}-A{{D}^{2}},\]
	\[2\overrightarrow{OB}\cdot \overrightarrow{OC}=O{{B}^{2}}+O{{C}^{2}}-B{{C}^{2}},\]
	\[2\overrightarrow{OB}\cdot \overrightarrow{OD}=O{{B}^{2}}+O{{D}^{2}}-B{{D}^{2}},\]
	\[2\overrightarrow{OC}\cdot \overrightarrow{OD}=O{{C}^{2}}+O{{D}^{2}}-C{{D}^{2}}.\]
	
	Therefore,
	\[O{{P}^{2}}={{J}_{1}}+{{J}_{2}}+{{J}_{3}}+{{J}_{4}}+{{K}_{1}}+{{K}_{2}},\]
	where
	\[{{J}_{1}}=\left( {{\left( \beta _{A}^{P} \right)}^{2}}+\beta _{A}^{P}\beta _{B}^{P}+\beta _{A}^{P}\beta _{C}^{P}+\beta _{A}^{P}\beta _{D}^{P} \right)O{{A}^{2}},\]
	\[{{J}_{2}}=\left( {{\left( \beta _{B}^{P} \right)}^{2}}+\beta _{A}^{P}\beta _{B}^{P}+\beta _{B}^{P}\beta _{C}^{P}+\beta _{B}^{P}\beta _{D}^{P} \right)O{{B}^{2}},\]
	\[{{J}_{3}}=\left( {{\left( \beta _{C}^{P} \right)}^{2}}+\beta _{A}^{P}\beta _{C}^{P}+\beta _{B}^{P}\beta _{C}^{P}+\beta _{C}^{P}\beta _{D}^{P} \right)O{{C}^{2}},\]
	\[{{J}_{4}}=\left( {{\left( \beta _{D}^{P} \right)}^{2}}+\beta _{A}^{P}\beta _{D}^{P}+\beta _{B}^{P}\beta _{D}^{P}+\beta _{C}^{P}\beta _{D}^{P} \right)O{{D}^{2}},\]
	\[{{K}_{1}}=-\beta _{A}^{P}\beta _{B}^{P}A{{B}^{2}}-\beta _{A}^{P}\beta _{C}^{P}A{{C}^{2}}-\beta _{A}^{P}\beta _{D}^{P}A{{D}^{2}},\]
	\[{{K}_{2}}=-\beta _{B}^{P}\beta _{C}^{P}B{{C}^{2}}-\beta _{B}^{P}\beta _{D}^{P}B{{D}^{2}}-\beta _{C}^{P}\beta _{D}^{P}C{{D}^{2}},\]
	i.e.
	\[\begin{aligned}
		O{{P}^{2}}&=\beta _{A}^{P}O{{A}^{2}}+\beta _{B}^{P}O{{B}^{2}}+\beta _{C}^{P}O{{C}^{2}}+\beta _{D}^{P}O{{D}^{2}} \\ 
		& -\beta _{A}^{P}\beta _{B}^{P}A{{B}^{2}}-\beta _{A}^{P}\beta _{C}^{P}A{{C}^{2}}-\beta _{A}^{P}\beta _{D}^{P}A{{D}^{2}} \\ 
		& -\beta _{B}^{P}\beta _{C}^{P}B{{C}^{2}}-\beta _{B}^{P}\beta _{D}^{P}B{{D}^{2}}-\beta _{C}^{P}\beta _{D}^{P}C{{D}^{2}}.  
	\end{aligned}\]
\end{proof}
\hfill $\square$\par

The distance $OP$ in the formula is expressed by the frame magnitudes (i.e. $OA$, $OB$, $OC$, $OD$), the frame components ($\beta _{A}^{P}$, $\beta _{B}^{P}$, $\beta _{C}^{P}$, $\beta _{D}^{P}$) of point $P$ and the length of each edge of six edges of the tetrahedron. 

The quantities $\beta _{A}^{P}$, $\beta _{B}^{P}$, $\beta _{C}^{P}$, $\beta _{D}^{P}$ are also called the frame components of point $P$.

\begin{corollary}{Relationship of distance between vetex and IC-T, frame components and lengths of edges, Daiyuan Zhang}{Cor24.1.1}\label{Cor24.1.1} 
	Given the tetrahedron $ABCD$, let point $P$ be the intersecting center, then 
	\[\begin{aligned}
		& \beta _{A}^{P}P{{A}^{2}}+\beta _{B}^{P}P{{B}^{2}}+\beta _{C}^{P}P{{C}^{2}}+\beta _{D}^{P}P{{D}^{2}} \\ 
		& =\beta _{A}^{P}\beta _{B}^{P}A{{B}^{2}}+\beta _{A}^{P}\beta _{C}^{P}A{{C}^{2}}+\beta _{A}^{P}\beta _{D}^{P}A{{D}^{2}} \\ 
		& +\beta _{B}^{P}\beta _{C}^{P}B{{C}^{2}}+\beta _{B}^{P}\beta _{D}^{P}B{{D}^{2}}+\beta _{C}^{P}\beta _{D}^{P}C{{D}^{2}},  
	\end{aligned}\]	
	where $\beta _{A}^{P}$, $\beta _{B}^{P}$, $\beta _{C}^{P}$, $\beta _{D}^{P}$ are the frame components of point $P$.
\end{corollary}

\begin{proof}
	Let the points $\,O$ and $P$ coincide, and use theorem \ref{thm:Thm24.1.1} to get the result.
\end{proof}
\hfill $\square$\par

\subsection{Distance between vertex and intersecting center of a tetrahedron}\label{Subsec24.1.1}
This subsection discusses the distance between vertex and intersecting center of a tetrahedron (abbreviated as DVIC-T).

\begin{theorem}{Distance between vertex and IC-T, Daiyuan Zhang}{Thm24.1.2}\label{Thm24.1.2} 
	Given the tetrahedron $ABCD$ and an intersecting center $P$ of the tetrahedron, then	
	\[\begin{aligned}
		A{{P}^{2}}&=\left( 1-\beta _{A}^{P} \right)\left( \beta _{B}^{P}A{{B}^{2}}+\beta _{C}^{P}A{{C}^{2}}+\beta _{D}^{P}A{{D}^{2}} \right) \\ 
		& -\beta _{B}^{P}\beta _{C}^{P}B{{C}^{2}}-\beta _{B}^{P}\beta _{D}^{P}B{{D}^{2}}-\beta _{C}^{P}\beta _{D}^{P}C{{D}^{2}},  
	\end{aligned}\]	
	\[\begin{aligned}
		B{{P}^{2}}&=\left( 1-\beta _{B}^{P} \right)\left( \beta _{C}^{P}B{{C}^{2}}+\beta _{D}^{P}B{{D}^{2}}+\beta _{A}^{P}B{{A}^{2}} \right) \\ 
		& -\beta _{C}^{P}\beta _{D}^{P}C{{D}^{2}}-\beta _{C}^{P}\beta _{A}^{P}C{{A}^{2}}-\beta _{D}^{P}\beta _{A}^{P}D{{A}^{2}},  
	\end{aligned}\]	
	\[\begin{aligned}
		C{{P}^{2}}&=\left( 1-\beta _{C}^{P} \right)\left( \beta _{D}^{P}C{{D}^{2}}+\beta _{A}^{P}C{{A}^{2}}+\beta _{B}^{P}C{{B}^{2}} \right) \\ 
		& -\beta _{D}^{P}\beta _{A}^{P}D{{A}^{2}}-\beta _{D}^{P}\beta _{B}^{P}D{{B}^{2}}-\beta _{A}^{P}\beta _{B}^{P}A{{B}^{2}},  
	\end{aligned}\]	
	\[\begin{aligned}
		D{{P}^{2}}&=\left( 1-\beta _{D}^{P} \right)\left( \beta _{A}^{P}D{{A}^{2}}+\beta _{B}^{P}D{{B}^{2}}+\beta _{C}^{P}D{{C}^{2}} \right) \\ 
		& -\beta _{A}^{P}\beta _{B}^{P}A{{B}^{2}}-\beta _{A}^{P}\beta _{C}^{P}A{{C}^{2}}-\beta _{B}^{P}\beta _{C}^{P}B{{C}^{2}}.  
	\end{aligned}\]	
	Where $\beta _{A}^{P}$, $\beta _{B}^{P}$, $\beta _{C}^{P}$, $\beta _{D}^{P}$ are the frame components of point $P$.
\end{theorem}

\begin{proof}
	Let the frame origin $O$ coincide with a vertex $A$ of tetrahedron $ABCD$, and $P$ be the intersecting center of the tetrahedron, then from theorem \ref{thm:Thm24.1.1}, the following result is obtained directly:	
	\[\begin{aligned}
		A{{P}^{2}}&=\beta _{B}^{P}A{{B}^{2}}+\beta _{C}^{P}A{{C}^{2}}+\beta _{D}^{P}A{{D}^{2}} \\ 
		& -\beta _{A}^{P}\beta _{B}^{P}A{{B}^{2}}-\beta _{A}^{P}\beta _{C}^{P}A{{C}^{2}}-\beta _{A}^{P}\beta _{D}^{P}A{{D}^{2}} \\ 
		& -\beta _{B}^{P}\beta _{C}^{P}B{{C}^{2}}-\beta _{B}^{P}\beta _{D}^{P}B{{D}^{2}}-\beta _{C}^{P}\beta _{D}^{P}C{{D}^{2}},  
	\end{aligned}\]	
	i.e.	
	\[\begin{aligned}
		A{{P}^{2}}&=\left( \beta _{B}^{P}-\beta _{A}^{P}\beta _{B}^{P} \right)A{{B}^{2}}+\left( \beta _{C}^{P}-\beta _{A}^{P}\beta _{C}^{P} \right)A{{C}^{2}}+\left( \beta _{D}^{P}-\beta _{A}^{P}\beta _{D}^{P} \right)A{{D}^{2}} \\ 
		& -\beta _{B}^{P}\beta _{C}^{P}B{{C}^{2}}-\beta _{B}^{P}\beta _{D}^{P}B{{D}^{2}}-\beta _{C}^{P}\beta _{D}^{P}C{{D}^{2}},  
	\end{aligned}\]
	i.e.	
	\[\begin{aligned}
		A{{P}^{2}}&=\beta _{B}^{P}\left( 1-\beta _{A}^{P} \right)A{{B}^{2}}+\beta _{C}^{P}\left( 1-\beta _{A}^{P} \right)A{{C}^{2}}+\beta _{D}^{P}\left( 1-\beta _{A}^{P} \right)A{{D}^{2}} \\ 
		& -\beta _{B}^{P}\beta _{C}^{P}B{{C}^{2}}-\beta _{B}^{P}\beta _{D}^{P}B{{D}^{2}}-\beta _{C}^{P}\beta _{D}^{P}C{{D}^{2}},  
	\end{aligned}\]
	i.e.
	\[\begin{aligned}
		A{{P}^{2}}&=\left( 1-\beta _{A}^{P} \right)\left( \beta _{B}^{P}A{{B}^{2}}+\beta _{C}^{P}A{{C}^{2}}+\beta _{D}^{P}A{{D}^{2}} \right) \\ 
		& -\beta _{B}^{P}\beta _{C}^{P}B{{C}^{2}}-\beta _{B}^{P}\beta _{D}^{P}B{{D}^{2}}-\beta _{C}^{P}\beta _{D}^{P}C{{D}^{2}}.  
	\end{aligned}\]	
	
	Similar results can be obtained for $B{{P}^{2}}$, $C{{P}^{2}}$ and $D{{P}^{2}}$.
\end{proof}
\hfill $\square$\par

The intersecting center $P$ of the tetrahedron can be selected as some special points of the tetrahedron, such as centroid, incenter, orthocenter, circumcenter, excenter, etc..


According to the above theorem, we directly calculate $A{{P}^{2}}+B{{P}^{2}}+C{{P}^{2}}+D{{P}^{2}}$, and obtain the following corollary.


\begin{corollary}{Sum of squares for distances from vertexs to intersecting center of a tetrahedron, Daiyuan Zhang}{Cor24.1.2}\label{Cor24.1.2} 
	Let the lengths of the six edges of the tetrahedron $ABCD$ be $AB$, $AC$, $AD$, $BC$, $CD$, $DB$, and the intersecting center of the tetrahedron $ABCD$ be $P$, then
	\[\begin{aligned}
		& A{{P}^{2}}+B{{P}^{2}}+C{{P}^{2}}+D{{P}^{2}} \\ 
		& =\left( \beta _{A}^{Q}+\beta _{B}^{Q}-4\beta _{A}^{Q}\beta _{B}^{Q} \right)A{{B}^{2}}+\left( \beta _{A}^{Q}+\beta _{C}^{Q}-4\beta _{A}^{Q}\beta _{C}^{Q} \right)A{{C}^{2}} \\ 
		& +\left( \beta _{A}^{Q}+\beta _{D}^{Q}-4\beta _{A}^{Q}\beta _{D}^{Q} \right)A{{D}^{2}}+\left( \beta _{B}^{Q}+\beta _{C}^{Q}-4\beta _{B}^{Q}\beta _{C}^{Q} \right)B{{C}^{2}} \\ 
		& +\left( \beta _{C}^{Q}+\beta _{D}^{Q}-4\beta _{C}^{Q}\beta _{D}^{Q} \right)C{{D}^{2}}+\left( \beta _{D}^{Q}+\beta _{B}^{Q}-4\beta _{D}^{Q}\beta _{B}^{Q} \right)D{{B}^{2}}.  
	\end{aligned}\]
\end{corollary}
\hfill $\square$\par

\begin{example}\label{Exam24.1.1} 
	Find the distance between vertex and centroid of a tetrahedron.
\end{example}

\begin{solution}
	Each of the frame components is ${1}/{4}\;$. So the distance between vertex and centroid of a tetrahedron is		
	\[\begin{aligned}
		A{{G}^{2}}&=\left( 1-\beta _{A}^{G} \right)\left( \beta _{B}^{G}A{{B}^{2}}+\beta _{C}^{G}A{{C}^{2}}+\beta _{D}^{G}A{{D}^{2}} \right) \\ 
		& -\beta _{B}^{G}\beta _{C}^{G}B{{C}^{2}}-\beta _{B}^{G}\beta _{D}^{G}B{{D}^{2}}-\beta _{C}^{G}\beta _{D}^{G}C{{D}^{2}},  
	\end{aligned}\]
	i.e.	
	\[\begin{aligned}
		A{{G}^{2}}&=\left( 1-\frac{1}{4} \right)\left( \frac{1}{4}A{{B}^{2}}+\frac{1}{4}A{{C}^{2}}+\frac{1}{4}A{{D}^{2}} \right) \\ 
		& -\frac{1}{4}\cdot \frac{1}{4}B{{C}^{2}}-\frac{1}{4}\cdot \frac{1}{4}B{{D}^{2}}-\frac{1}{4}\cdot \frac{1}{4}C{{D}^{2}} \\ 
		& =\frac{3}{16}\left( A{{B}^{2}}+A{{C}^{2}}+A{{D}^{2}} \right)-\frac{1}{16}\left( B{{C}^{2}}+B{{D}^{2}}+C{{D}^{2}} \right),  
	\end{aligned}\]
	i.e.	
	\[AG=\frac{1}{4}\sqrt{3\left( A{{B}^{2}}+A{{C}^{2}}+A{{D}^{2}} \right)-\left( B{{C}^{2}}+B{{D}^{2}}+C{{D}^{2}} \right)}.\]
	
	Similarly, 
	\[BG=\frac{1}{4}\sqrt{3\left( B{{C}^{2}}+B{{D}^{2}}+B{{A}^{2}} \right)-\left( C{{D}^{2}}+C{{A}^{2}}+D{{A}^{2}} \right)},\]
	\[CG=\frac{1}{4}\sqrt{3\left( C{{D}^{2}}+C{{A}^{2}}+C{{B}^{2}} \right)-\left( D{{A}^{2}}+D{{B}^{2}}+A{{B}^{2}} \right)},\]
	\[DG=\frac{1}{4}\sqrt{3\left( D{{A}^{2}}+D{{B}^{2}}+D{{C}^{2}} \right)-\left( A{{B}^{2}}+A{{C}^{2}}+B{{C}^{2}} \right)}.\]
\end{solution}
\hfill $\diamond$\par

From the results of the above example, the distance between vertex and centroid of a regular tetrahedron is obtained:
\[AG=BG=CG=DG=\frac{\sqrt{6}}{4}a,\]	
where $a$ is the edge length of a regular tetrahedron.

For some other centers of the tetrahedron (such as incenter, circumcenter, etc.) the DVIC-T can also be similarly calculated.

\subsection{Distance between vertex and intersecting foot of a tetrahedron}\label{Subsec24.1.2}
The distance from the vertex of a tetrahedron to the intersecting foot is called the distance between vertex and intersecting foot of the tetrahedron (abbreviated as DVIF-T).

\begin{theorem}{Distance between vertex and  IF-T, Daiyuan Zhang}{Thm24.1.3}\label{Thm24.1.3} 
	Given the tetrahedron $ABCD$, the intersecting center of the tetrahedron is $P$, then	
	\[A{{P}_{A}}=\frac{AP}{\left| 1-\beta _{A}^{P} \right|}=\frac{\sqrt{\begin{aligned}
				& \left( 1-\beta _{A}^{P} \right)\left( \beta _{B}^{P}A{{B}^{2}}+\beta _{C}^{P}A{{C}^{2}}+\beta _{D}^{P}A{{D}^{2}} \right) \\ 
				& -\beta _{B}^{P}\beta _{C}^{P}B{{C}^{2}}-\beta _{B}^{P}\beta _{D}^{P}B{{D}^{2}}-\beta _{C}^{P}\beta _{D}^{P}C{{D}^{2}}  
	\end{aligned}}}{\left| 1-\beta _{A}^{P} \right|},\]
	\[B{{P}_{B}}=\frac{BP}{\left| 1-\beta _{B}^{P} \right|}=\frac{\sqrt{\begin{aligned}
				& \left( 1-\beta _{B}^{P} \right)\left( \beta _{C}^{P}B{{C}^{2}}+\beta _{D}^{P}B{{D}^{2}}+\beta _{A}^{P}B{{A}^{2}} \right) \\ 
				& -\beta _{C}^{P}\beta _{D}^{P}C{{D}^{2}}-\beta _{C}^{P}\beta _{A}^{P}C{{A}^{2}}-\beta _{D}^{P}\beta _{A}^{P}D{{A}^{2}}  
	\end{aligned}}}{\left| 1-\beta _{B}^{P} \right|},\]	
	\[C{{P}_{C}}=\frac{CP}{\left| 1-\beta _{C}^{P} \right|}=\frac{\sqrt{\begin{aligned}
				& \left( 1-\beta _{C}^{P} \right)\left( \beta _{D}^{P}C{{D}^{2}}+\beta _{A}^{P}C{{A}^{2}}+\beta _{B}^{P}C{{B}^{2}} \right) \\ 
				& -\beta _{D}^{P}\beta _{A}^{P}D{{A}^{2}}-\beta _{D}^{P}\beta _{B}^{P}D{{B}^{2}}-\beta _{A}^{P}\beta _{B}^{P}A{{B}^{2}}  
	\end{aligned}}}{\left| 1-\beta _{C}^{P} \right|},\]	
	\[D{{P}_{D}}=\frac{DP}{\left| 1-\beta _{C}^{P} \right|}=\frac{\sqrt{\begin{aligned}
				& \left( 1-\beta _{D}^{P} \right)\left( \beta _{A}^{P}D{{A}^{2}}+\beta _{B}^{P}D{{B}^{2}}+\beta _{C}^{P}D{{C}^{2}} \right) \\ 
				& -\beta _{A}^{P}\beta _{B}^{P}A{{B}^{2}}-\beta _{A}^{P}\beta _{C}^{P}A{{C}^{2}}-\beta _{B}^{P}\beta _{C}^{P}B{{C}^{2}}  
	\end{aligned}}}{\left| 1-\beta _{C}^{P} \right|}.\]			
	Where $\beta _{A}^{P}$, $\beta _{B}^{P}$, $\beta _{C}^{P}$, $\beta _{D}^{P}$ are the frame components of point $P$.
\end{theorem}

\begin{proof}
	According to theorem \ref{thm:Thm24.1.2} and \ref{thm:Thm18.2.3}, this theorem can be directly obtained.
\end{proof}
\hfill $\square$\par

\begin{example}\label{Exam24.1.2} 
	Find the DVIF-T of centroid.
\end{example}

\begin{solution}
	According to theorem \ref{thm:Thm24.1.3} and theorem \ref{thm:Thm21.1.1}, the following result can be obtained:
	\begin{align*}
		A{{G}_{A}}&=\frac{AG}{\left| 1-\beta _{A}^{G} \right|}=\frac{\sqrt{\begin{aligned}
					& \left( 1-\frac{1}{4} \right)\left( \frac{1}{4}A{{B}^{2}}+\frac{1}{4}A{{C}^{2}}+\frac{1}{4}A{{D}^{2}} \right) \\ 
					& -\frac{1}{4}\times \frac{1}{4}B{{C}^{2}}-\frac{1}{4}\times \frac{1}{4}B{{D}^{2}}-\frac{1}{4}\times \frac{1}{4}C{{D}^{2}}  
		\end{aligned}}}{\left| 1-\frac{1}{4} \right|} \\ 
		& =\frac{1}{3}\sqrt{3\left( A{{B}^{2}}+A{{C}^{2}}+A{{D}^{2}} \right)-\left( B{{C}^{2}}+B{{D}^{2}}+C{{D}^{2}} \right)}. 
	\end{align*}
	
	Similarly,
	\[B{{G}_{B}}=\frac{BG}{\left| 1-\beta _{B}^{G} \right|}=\frac{1}{3}\sqrt{3\left( B{{C}^{2}}+B{{D}^{2}}+B{{A}^{2}} \right)-\left( C{{D}^{2}}+C{{A}^{2}}+D{{A}^{2}} \right)},\]
	\[C{{G}_{C}}=\frac{CG}{\left| 1-\beta _{C}^{G} \right|}=\frac{1}{3}\sqrt{3\left( C{{D}^{2}}+C{{A}^{2}}+C{{B}^{2}} \right)-\left( D{{A}^{2}}+D{{B}^{2}}+A{{B}^{2}} \right)},\]	
	\[D{{G}_{D}}=\frac{DG}{\left| 1-\beta _{C}^{G} \right|}=\frac{1}{3}\sqrt{3\left( D{{A}^{2}}+D{{B}^{2}}+D{{C}^{2}} \right)-\left( A{{B}^{2}}+A{{C}^{2}}+B{{C}^{2}} \right)}.\]
\end{solution}
\hfill $\diamond$\par

Sum the squares on both sides of the above four formulas to obtain the following conclusion.

\begin{corollary}{Sum of squares of DVIF-T for centroid of tetrahedron, Daiyuan Zhang}{Cor24.1.4}\label{Cor24.1.4} 
	Given the tetrahedron $ABCD$, the centroid of the tetrahedron is $G$, then
	\[\begin{aligned}
		A{{G}_{A}}^{2}+B{{G}_{B}}^{2}+C{{G}_{C}}^{2}+D{{G}_{D}}^{2} =\frac{4}{9}\left( A{{B}^{2}}+A{{C}^{2}}+A{{D}^{2}}+B{{C}^{2}}+B{{D}^{2}}+C{{D}^{2}} \right). \\ 
	\end{aligned}\]
\end{corollary}

\subsection{Distance between origin and intersecting center of a tetrahedron on frame of circumcenter}\label{Subsec24.1.3}

\begin{theorem}{Distance between origin and IC-T on frame of circumcenter, Daiyuan Zhang}{Thm24.1.4}\label{Thm24.1.4} 
	Given the tetrahedron $ABCD$, let the circumcenter $Q$ of tetrahedron $ABCD$ be the origin of the frame $\left( Q;A,B,C,D \right)$, point $P$ be the intersecting center of the tetrahedron, $R$ be the radius of circumscribed sphere of tetrahedron $ABCD$, then	
	\[Q{{P}^{2}}={{R}^{2}}-\left( \begin{aligned}
		& \beta _{A}^{P}\beta _{B}^{P}A{{B}^{2}}+\beta _{A}^{P}\beta _{C}^{P}A{{C}^{2}}+\beta _{A}^{P}\beta _{D}^{P}A{{D}^{2}} \\ 
		& +\beta _{B}^{P}\beta _{C}^{P}B{{C}^{2}}+\beta _{B}^{P}\beta _{D}^{P}B{{D}^{2}}+\beta _{C}^{P}\beta _{D}^{P}C{{D}^{2}}  
	\end{aligned} \right).\]
	Where $\beta _{A}^{P}$, $\beta _{B}^{P}$, $\beta _{C}^{P}$, $\beta _{D}^{P}$ are the frame components of point $P$.
\end{theorem}

\begin{proof}
	Since $QA=QB=QC=QD=R$, in theorem \ref{thm:Thm24.1.1}, if the origin of the frame $\,O$ coincides with the circumcenter $Q$ of tetrahedron $ABCD$, the desired result can be obtained.
\end{proof}
\hfill $\square$\par
%
\section{Distance between two intersecting centers on tetrahedral frame}\label{Sec24.2}
The distance between two intersecting centers of a tetrahedron (abbreviated as DTICs-T) refers to the distance between two intersecting centers of the same tetrahedron. The DTICs-T is also called the distance between intersecting centers of a tetrahedron (abbreviated as DICs-T). The following theorem is one of the core theorems of this book.

\begin{theorem}{Distance between two intersecting centers on tetrahedral frame, Daiyuan Zhang}{Thm24.2.1}\label{Thm24.2.1} 
	Suppose that given a tetrahedron $ABCD$, ${{P}_{1}}\in {{\pi }_{ABCD}}$, ${{P}_{2}}\in {{\pi }_{ABCD}}$, then 	
	\[\begin{aligned}
	{{P}_{1}}{{P}_{2}}^{2}= &-\beta _{A}^{{{P}_{1}}{{P}_{2}}}\beta _{B}^{{{P}_{1}}{{P}_{2}}}A{{B}^{2}}-\beta _{A}^{{{P}_{1}}{{P}_{2}}}\beta _{C}^{{{P}_{1}}{{P}_{2}}}A{{C}^{2}}-\beta _{A}^{{{P}_{1}}{{P}_{2}}}\beta _{D}^{{{P}_{1}}{{P}_{2}}}A{{D}^{2}} \\ 
		& -\beta _{B}^{{{P}_{1}}{{P}_{2}}}\beta _{C}^{{{P}_{1}}{{P}_{2}}}B{{C}^{2}}-\beta _{B}^{{{P}_{1}}{{P}_{2}}}\beta _{D}^{{{P}_{1}}{{P}_{2}}}B{{D}^{2}}-\beta _{C}^{{{P}_{1}}{{P}_{2}}}\beta _{D}^{{{P}_{1}}{{P}_{2}}}C{{D}^{2}},  
	\end{aligned}\]
	\[\beta _{A}^{{{P}_{1}}{{P}_{2}}}+\beta _{B}^{{{P}_{1}}{{P}_{2}}}+\beta _{C}^{{{P}_{1}}{{P}_{2}}}+\beta _{D}^{{{P}_{1}}{{P}_{2}}}=0.\]	
	
	Where
	\[\beta _{A}^{{{P}_{1}}{{P}_{2}}}=\beta _{A}^{{{P}_{2}}}-\beta _{A}^{{{P}_{1}}},\]	
	\[\beta _{B}^{{{P}_{1}}{{P}_{2}}}=\beta _{B}^{{{P}_{2}}}-\beta _{B}^{{{P}_{1}}},\]	
	\[\beta _{C}^{{{P}_{1}}{{P}_{2}}}=\beta _{C}^{{{P}_{2}}}-\beta _{C}^{{{P}_{1}}},\]	
	\[\beta _{D}^{{{P}_{1}}{{P}_{2}}}=\beta _{D}^{{{P}_{2}}}-\beta _{D}^{{{P}_{1}}}.\]	
	Where $\beta _{A}^{{{P}_{1}}}$, $\beta _{B}^{{{P}_{1}}}$, $\beta _{C}^{{{P}_{1}}}$, $\beta _{D}^{{{P}_{1}}}$ are the frame components of intersecting center ${{P}_{1}}$ on $\overrightarrow{OA}$, $\overrightarrow{OB}$, $\overrightarrow{OC}$, $\overrightarrow{OD}$ respectively; and $\beta _{A}^{{{P}_{2}}}$, $\beta _{B}^{{{P}_{2}}}$, $\beta _{C}^{{{P}_{2}}}$, $\beta _{D}^{{{P}_{2}}}$ are the frame components of intersecting center ${{P}_{2}}$  on $\overrightarrow{OA}$, $\overrightarrow{OB}$, $\overrightarrow{OC}$, $\overrightarrow{OD}$ respectively.
\end{theorem}

\begin{proof}
	From the inner product operation of the vectors and  theorem\ref{thm:Thm18.3.2}, the following formulas are obtained:	
	\[\begin{aligned}
		{{P}_{1}}{{P}_{2}}^{2}&=\overrightarrow{{{P}_{1}}{{P}_{2}}}\cdot \overrightarrow{{{P}_{1}}{{P}_{2}}} \\ 
		& =\left( \beta _{B}^{{{P}_{1}}{{P}_{2}}}\overrightarrow{AB}+\beta _{C}^{{{P}_{1}}{{P}_{2}}}\overrightarrow{AC}+\beta _{D}^{{{P}_{1}}{{P}_{2}}}\overrightarrow{AD} \right)\cdot \left( \beta _{B}^{{{P}_{1}}{{P}_{2}}}\overrightarrow{AB}+\beta _{C}^{{{P}_{1}}{{P}_{2}}}\overrightarrow{AC}+\beta _{D}^{{{P}_{1}}{{P}_{2}}}\overrightarrow{AD} \right) \\ 
		& ={{\left( \beta _{B}^{{{P}_{1}}{{P}_{2}}} \right)}^{2}}A{{B}^{2}}+{{\left( \beta _{C}^{{{P}_{1}}{{P}_{2}}} \right)}^{2}}A{{C}^{2}}+{{\left( \beta _{D}^{{{P}_{1}}{{P}_{2}}} \right)}^{2}}A{{D}^{2}} \\ 
		& +2\beta _{B}^{{{P}_{1}}{{P}_{2}}}\beta _{C}^{{{P}_{1}}{{P}_{2}}}\overrightarrow{AB}\cdot \overrightarrow{AC}+2\beta _{B}^{{{P}_{1}}{{P}_{2}}}\beta _{D}^{{{P}_{1}}{{P}_{2}}}\overrightarrow{AB}\cdot \overrightarrow{AD}+2\beta _{C}^{{{P}_{1}}{{P}_{2}}}\beta _{D}^{{{P}_{1}}{{P}_{2}}}\overrightarrow{AC}\cdot \overrightarrow{AD},  
	\end{aligned}\]
	i.e.
	\[\begin{aligned}
		{{P}_{1}}{{P}_{2}}^{2}&={{\left( \beta _{B}^{{{P}_{1}}{{P}_{2}}} \right)}^{2}}A{{B}^{2}}+{{\left( \beta _{C}^{{{P}_{1}}{{P}_{2}}} \right)}^{2}}A{{C}^{2}}+{{\left( \beta _{D}^{{{P}_{1}}{{P}_{2}}} \right)}^{2}}A{{D}^{2}} \\ 
		& +\beta _{B}^{{{P}_{1}}{{P}_{2}}}\beta _{C}^{{{P}_{1}}{{P}_{2}}}\left( A{{B}^{2}}+A{{C}^{2}}-B{{C}^{2}} \right)+\beta _{B}^{{{P}_{1}}{{P}_{2}}}\beta _{D}^{{{P}_{1}}{{P}_{2}}}\left( A{{B}^{2}}+A{{D}^{2}}-B{{D}^{2}} \right) \\ 
		& +\beta _{C}^{{{P}_{1}}{{P}_{2}}}\beta _{D}^{{{P}_{1}}{{P}_{2}}}\left( A{{C}^{2}}+A{{D}^{2}}-C{{D}^{2}} \right),  
	\end{aligned}\]
	i.e.
	\[\begin{aligned}
		{{P}_{1}}{{P}_{2}}^{2}&=\left( {{\left( \beta _{B}^{{{P}_{1}}{{P}_{2}}} \right)}^{2}}+\beta _{B}^{{{P}_{1}}{{P}_{2}}}\beta _{C}^{{{P}_{1}}{{P}_{2}}}+\beta _{B}^{{{P}_{1}}{{P}_{2}}}\beta _{D}^{{{P}_{1}}{{P}_{2}}} \right)A{{B}^{2}} \\ 
		& +\left( {{\left( \beta _{C}^{{{P}_{1}}{{P}_{2}}} \right)}^{2}}+\beta _{B}^{{{P}_{1}}{{P}_{2}}}\beta _{C}^{{{P}_{1}}{{P}_{2}}}+\beta _{C}^{{{P}_{1}}{{P}_{2}}}\beta _{D}^{{{P}_{1}}{{P}_{2}}} \right)A{{C}^{2}} \\ 
		& +\left( {{\left( \beta _{D}^{{{P}_{1}}{{P}_{2}}} \right)}^{2}}+\beta _{B}^{{{P}_{1}}{{P}_{2}}}\beta _{D}^{{{P}_{1}}{{P}_{2}}}+\beta _{C}^{{{P}_{1}}{{P}_{2}}}\beta _{D}^{{{P}_{1}}{{P}_{2}}} \right)A{{D}^{2}} \\ 
		& -\beta _{B}^{{{P}_{1}}{{P}_{2}}}\beta _{C}^{{{P}_{1}}{{P}_{2}}}B{{C}^{2}}-\beta _{B}^{{{P}_{1}}{{P}_{2}}}\beta _{D}^{{{P}_{1}}{{P}_{2}}}B{{D}^{2}}-\beta _{C}^{{{P}_{1}}{{P}_{2}}}\beta _{D}^{{{P}_{1}}{{P}_{2}}}C{{D}^{2}},  
	\end{aligned}\]
	i.e.
	\[\begin{aligned}
		{{P}_{1}}{{P}_{2}}^{2}&=\beta _{B}^{{{P}_{1}}{{P}_{2}}}\left( \beta _{B}^{{{P}_{1}}{{P}_{2}}}+\beta _{C}^{{{P}_{1}}{{P}_{2}}}+\beta _{D}^{{{P}_{1}}{{P}_{2}}} \right)A{{B}^{2}}+\beta _{C}^{{{P}_{1}}{{P}_{2}}}\left( \beta _{C}^{{{P}_{1}}{{P}_{2}}}+\beta _{B}^{{{P}_{1}}{{P}_{2}}}+\beta _{D}^{{{P}_{1}}{{P}_{2}}} \right)A{{C}^{2}} \\ 
		& +\beta _{D}^{{{P}_{1}}{{P}_{2}}}\left( \beta _{D}^{{{P}_{1}}{{P}_{2}}}+\beta _{B}^{{{P}_{1}}{{P}_{2}}}+\beta _{C}^{{{P}_{1}}{{P}_{2}}} \right)A{{D}^{2}}-\beta _{B}^{{{P}_{1}}{{P}_{2}}}\beta _{C}^{{{P}_{1}}{{P}_{2}}}B{{C}^{2}}-\beta _{B}^{{{P}_{1}}{{P}_{2}}}\beta _{D}^{{{P}_{1}}{{P}_{2}}}B{{D}^{2}}-\beta _{C}^{{{P}_{1}}{{P}_{2}}}\beta _{D}^{{{P}_{1}}{{P}_{2}}}C{{D}^{2}}  
	\end{aligned}\]
	i.e.
	\[\begin{aligned}
		{{P}_{1}}{{P}_{2}}^{2}&=-\beta _{A}^{{{P}_{1}}{{P}_{2}}}\beta _{B}^{{{P}_{1}}{{P}_{2}}}A{{B}^{2}}-\beta _{A}^{{{P}_{1}}{{P}_{2}}}\beta _{C}^{{{P}_{1}}{{P}_{2}}}A{{C}^{2}}-\beta _{A}^{{{P}_{1}}{{P}_{2}}}\beta _{D}^{{{P}_{1}}{{P}_{2}}}A{{D}^{2}} \\ 
		& -\beta _{B}^{{{P}_{1}}{{P}_{2}}}\beta _{C}^{{{P}_{1}}{{P}_{2}}}B{{C}^{2}}-\beta _{B}^{{{P}_{1}}{{P}_{2}}}\beta _{D}^{{{P}_{1}}{{P}_{2}}}B{{D}^{2}}-\beta _{C}^{{{P}_{1}}{{P}_{2}}}\beta _{D}^{{{P}_{1}}{{P}_{2}}}C{{D}^{2}}.  
	\end{aligned}\]
	
	And
	\[\begin{aligned}
		& \beta _{A}^{{{P}_{1}}{{P}_{2}}}+\beta _{B}^{{{P}_{1}}{{P}_{2}}}+\beta _{C}^{{{P}_{1}}{{P}_{2}}}+\beta _{D}^{{{P}_{1}}{{P}_{2}}} \\ 
		& =\left( \beta _{A}^{{{P}_{2}}}-\beta _{A}^{{{P}_{1}}} \right)+\left( \beta _{B}^{{{P}_{2}}}-\beta _{B}^{{{P}_{1}}} \right)+\left( \beta _{C}^{{{P}_{2}}}-\beta _{C}^{{{P}_{1}}} \right)+\left( \beta _{D}^{{{P}_{2}}}-\beta _{D}^{{{P}_{1}}} \right) \\ 
		& =\left( \beta _{A}^{{{P}_{2}}}+\beta _{B}^{{{P}_{2}}}+\beta _{C}^{{{P}_{2}}}+\beta _{D}^{{{P}_{2}}} \right)-\left( \beta _{A}^{{{P}_{1}}}+\beta _{B}^{{{P}_{1}}}+\beta _{C}^{{{P}_{1}}}+\beta _{D}^{{{P}_{1}}} \right)=0.  
	\end{aligned}\]
\end{proof}
\hfill $\square$\par

The above theorem shows that if we know each length of six edges of a tetrahedron and the frame components (or IRs) of two ICs-T: ${{P}_{1}}$ and ${{P}_{2}}$, then the magnitude of the vector $\overrightarrow{{{P}_{1}}{{P}_{2}}}$ can be calculated, i.e., the distance between the points ${{P}_{1}}$ and ${{P}_{2}}$ can be calculated.



Obviously, each length of the six edges of a tetrahedron has nothing to do with the location of the origin $O$ of the tetrahedron frame $\left( O;A,B,C,D \right)$. In other words, if the relative position between the given points ${{P}_{1}}$, ${{P}_{2}}$ and the tetrahedron $ABCD$ is fixed, that is, ${{P}_{1}}$, ${{P}_{2}}$ are the rigid points of tetrahedron $ABCD$, then any translation, rotation or even flipping of tetrahedron $ABCD$ in space, the parameters (frame components) representing the distance of ${{P}_{1}}{{P}_{2}}$ will keep unchanged in the formula, which is very convenient for the analysis and application of some problems, because the origin of the frame can be chosen as needed without changing the parameters (frame components) in the formula.


\chapter{Frame components of circumcenter of a tetrahedron}\label{Sec25}
\thispagestyle{empty}




%

%

%

%



%

%

%

%

%


Solving the frame components are the core problem of Intercenter Geometry. With the frame components, we can establish a vector (VOIC-T or VTICs-T), and then obtain the distance and other related quantities. The idea in the previous chapters is to first calculate the IR, and then calculate the frame components. However, in some cases, we can bypass the process of calculating the IR and directly solve the frame components, This method is used to find the frame components of the circumcenter of tetrahedron in this chapter. This chapter is divided into two sections. Section \ref{Sec25.1} first studies the frame components of the circumcenter of a tetrahedron; Section \ref{Sec25.2} is the application of section \ref{Sec25.1} to study the frame components of the circumcenter of regular triangular pyramid and regular tetrahedron.

\section{Frame components of the circumcenter of a tetrahedron}\label{Sec25.1}

The circumscribed sphere of tetrahedron has many applications. This section studies the calculation formula of the frame components of the circumcenter of general tetrahedron. The following theorems I give are very important for studying some measurement parameters of the circumcenter of tetrahedron.

\begin{theorem}{Frame component formulas for circumcenter of a tetrahedron-1, Daiyuan Zhang}{Thm25.1.1}\label{Thm25.1.1}
Let the lengths of the six edges of the tetrahedron $ABCD$ be $AB$, $AC$, $AD$, $BC$, $CD$, $DB$, and the point $Q$ be the circumcenter of the tetrahedron $ABCD$, then the frame components of the circumcenter $Q$ of the tetrahedron are
\begin{flalign}
	\beta _{A}^{Q}=\frac{{{U}_{A}}}{U},\beta _{B}^{Q}=\frac{{{U}_{B}}}{U},\beta _{C}^{Q}=\frac{{{U}_{C}}}{U},\beta _{D}^{Q}=\frac{{{U}_{D}}}{U}.	
\end{flalign}

Where
\[U={{U}_{A}}+{{U}_{B}}+{{U}_{C}}+{{U}_{D}},\]
\[{{U}_{A}}=\left| \begin{matrix}
	-A{{B}^{2}} & B{{C}^{2}}-A{{C}^{2}} & B{{D}^{2}}-A{{D}^{2}}  \\
	B{{C}^{2}} & -B{{C}^{2}} & C{{D}^{2}}-B{{D}^{2}}  \\
	B{{D}^{2}}-B{{C}^{2}} & C{{D}^{2}} & -C{{D}^{2}}  \\
\end{matrix} \right|,\]
\[{{U}_{B}}=\left| \begin{matrix}
	-A{{B}^{2}} & A{{C}^{2}}-B{{C}^{2}} & A{{D}^{2}}-B{{D}^{2}}  \\
	A{{C}^{2}}-A{{B}^{2}} & -B{{C}^{2}} & C{{D}^{2}}-B{{D}^{2}}  \\
	A{{D}^{2}}-A{{C}^{2}} & C{{D}^{2}} & -C{{D}^{2}}  \\
\end{matrix} \right|,\]
\[{{U}_{C}}=\left| \begin{matrix}
	-A{{B}^{2}} & A{{B}^{2}} & A{{D}^{2}}-B{{D}^{2}}  \\
	A{{B}^{2}}-A{{C}^{2}} & -B{{C}^{2}} & B{{D}^{2}}-C{{D}^{2}}  \\
	A{{D}^{2}}-A{{C}^{2}} & B{{D}^{2}}-B{{C}^{2}} & -C{{D}^{2}}  \\
\end{matrix} \right|,\]
\[{{U}_{D}}=\left| \begin{matrix}
	-A{{B}^{2}} & A{{B}^{2}} & A{{C}^{2}}-B{{C}^{2}}  \\
	A{{B}^{2}}-A{{C}^{2}} & -B{{C}^{2}} & B{{C}^{2}}  \\
	A{{C}^{2}}-A{{D}^{2}} & B{{C}^{2}}-B{{D}^{2}} & -C{{D}^{2}}  \\
\end{matrix} \right|.\]
\end{theorem}

\begin{proof}
	Refer to figure \ref{fig:tu16.2.1}, based on theorem \ref{thm:Thm24.1.2}, we get:
	\[\begin{aligned}
		A{{Q}^{2}}& =\left( 1-\beta _{A}^{Q} \right)\left( \beta _{B}^{Q}A{{B}^{2}}+\beta _{C}^{Q}A{{C}^{2}}+\beta _{D}^{Q}A{{D}^{2}} \right) \\ 
		& -\beta _{B}^{Q}\beta _{C}^{Q}B{{C}^{2}}-\beta _{C}^{Q}\beta _{D}^{Q}C{{D}^{2}}-\beta _{D}^{Q}\beta _{B}^{Q}D{{B}^{2}},  
	\end{aligned}\]	
	\[\begin{aligned}
		B{{Q}^{2}}& =\left( 1-\beta _{B}^{Q} \right)\left( \beta _{C}^{Q}B{{C}^{2}}+\beta _{D}^{Q}B{{D}^{2}}+\beta _{A}^{Q}B{{A}^{2}} \right) \\ 
		& -\beta _{C}^{Q}\beta _{D}^{Q}C{{D}^{2}}-\beta _{D}^{Q}\beta _{A}^{Q}D{{A}^{2}}-\beta _{A}^{Q}\beta _{C}^{Q}A{{C}^{2}},  
	\end{aligned}\]	
	\[\begin{aligned}
		C{{Q}^{2}}& =\left( 1-\beta _{C}^{Q} \right)\left( \beta _{D}^{Q}C{{D}^{2}}+\beta _{A}^{Q}C{{A}^{2}}+\beta _{B}^{Q}C{{B}^{2}} \right) \\ 
		& -\beta _{D}^{Q}\beta _{A}^{Q}D{{A}^{2}}-\beta _{A}^{Q}\beta _{B}^{Q}A{{B}^{2}}-\beta _{B}^{Q}\beta _{D}^{Q}B{{D}^{2}},  
	\end{aligned}\]	
	\[\begin{aligned}
		D{{Q}^{2}}& =\left( 1-\beta _{D}^{Q} \right)\left( \beta _{A}^{Q}D{{A}^{2}}+\beta _{B}^{Q}D{{B}^{2}}+\beta _{C}^{Q}D{{C}^{2}} \right) \\ 
		& -\beta _{A}^{Q}\beta _{B}^{Q}A{{B}^{2}}-\beta _{B}^{Q}\beta _{C}^{Q}B{{C}^{2}}-\beta _{C}^{Q}\beta _{A}^{Q}C{{A}^{2}}.  
	\end{aligned}\]	
	
	
	Where $\beta _{A}^{Q}$, $\beta _{B}^{Q}$, $\beta _{C}^{Q}$, $\beta _{D}^{Q}$ are the frame components of circumcenter $Q$. 
	For circumcenter $Q$, we have: $A{{Q}^{2}}=B{{Q}^{2}}=C{{Q}^{2}}=D{{Q}^{2}}={{R}^{2}}$. Using equation $A{{Q}^{2}}=B{{Q}^{2}}$, the following result is obtained:
	\[\begin{aligned}
		& \left( 1-\beta _{A}^{Q} \right)\left( \beta _{B}^{Q}A{{B}^{2}}+\beta _{C}^{Q}A{{C}^{2}}+\beta _{D}^{Q}A{{D}^{2}} \right) \\ 
		& -\beta _{B}^{Q}\beta _{C}^{Q}B{{C}^{2}}-\beta _{B}^{Q}\beta _{D}^{Q}B{{D}^{2}}-\beta _{C}^{Q}\beta _{D}^{Q}C{{D}^{2}} \\ 
		& =\left( 1-\beta _{B}^{Q} \right)\left( \beta _{C}^{Q}B{{C}^{2}}+\beta _{D}^{Q}B{{D}^{2}}+\beta _{A}^{Q}B{{A}^{2}} \right) \\ 
		& -\beta _{C}^{Q}\beta _{D}^{Q}C{{D}^{2}}-\beta _{C}^{Q}\beta _{A}^{Q}C{{A}^{2}}-\beta _{D}^{Q}\beta _{A}^{Q}D{{A}^{2}},  
	\end{aligned}\]
	i.e.
	\[\beta _{B}^{Q}A{{B}^{2}}+\beta _{C}^{Q}A{{C}^{2}}+\beta _{D}^{Q}A{{D}^{2}}=\beta _{C}^{Q}B{{C}^{2}}+\beta _{D}^{Q}B{{D}^{2}}+\beta _{A}^{Q}B{{A}^{2}},\]
	\[\beta _{A}^{Q}B{{A}^{2}}-\beta _{B}^{Q}A{{B}^{2}}+\beta _{C}^{Q}\left( B{{C}^{2}}-A{{C}^{2}} \right)+\beta _{D}^{Q}\left( B{{D}^{2}}-A{{D}^{2}} \right)=0.\]
	
	Using equation $B{{Q}^{2}}=C{{Q}^{2}}$ (or rotate the above formula), the following result is obtained:
	\[\beta _{B}^{Q}C{{B}^{2}}-\beta _{C}^{Q}B{{C}^{2}}+\beta _{D}^{Q}\left( C{{D}^{2}}-B{{D}^{2}} \right)+\beta _{A}^{Q}\left( C{{A}^{2}}-B{{A}^{2}} \right)=0.\]
	
	Using equation $C{{Q}^{2}}=D{{Q}^{2}}$ (or rotate the above formula), the following result is obtained:
	\[\beta _{C}^{Q}D{{C}^{2}}-\beta _{D}^{Q}C{{D}^{2}}+\beta _{A}^{Q}\left( D{{A}^{2}}-C{{A}^{2}} \right)+\beta _{B}^{Q}\left( D{{B}^{2}}-C{{B}^{2}} \right)=0.\]
	
	Therefore, using the condition $\beta _{A}^{Q}+\beta _{B}^{Q}+\beta _{C}^{Q}+\beta _{D}^{Q}=1$, the following linear equations are obtained:
	\[\left\{ \begin{aligned}
		& \beta _{A}^{Q}+\beta _{B}^{Q}+\beta _{C}^{Q}+\beta _{D}^{Q}=1 \\ 
		& B{{A}^{2}}\beta _{A}^{Q}-A{{B}^{2}}\beta _{B}^{Q}+\left( B{{C}^{2}}-A{{C}^{2}} \right)\beta _{C}^{Q}+\left( B{{D}^{2}}-A{{D}^{2}} \right)\beta _{D}^{Q}=0 \\ 
		& C{{B}^{2}}\beta _{B}^{Q}-B{{C}^{2}}\beta _{C}^{Q}+\left( C{{D}^{2}}-B{{D}^{2}} \right)\beta _{D}^{Q}+\left( C{{A}^{2}}-B{{A}^{2}} \right)\beta _{A}^{Q}=0 \\ 
		& D{{C}^{2}}\beta _{C}^{Q}-C{{D}^{2}}\beta _{D}^{Q}+\left( D{{A}^{2}}-C{{A}^{2}} \right)\beta _{A}^{Q}+\left( D{{B}^{2}}-C{{B}^{2}} \right)\beta _{B}^{Q}=0,\\ 
	\end{aligned} \right.\]
	i.e.
	\[\left\{ \begin{aligned}
		& \beta _{A}^{Q}+\beta _{B}^{Q}+\beta _{C}^{Q}+\beta _{D}^{Q}=1 \\ 
		& B{{A}^{2}}\beta _{A}^{Q}-A{{B}^{2}}\beta _{B}^{Q}+\left( B{{C}^{2}}-A{{C}^{2}} \right)\beta _{C}^{Q}+\left( B{{D}^{2}}-A{{D}^{2}} \right)\beta _{D}^{Q}=0 \\ 
		& \left( C{{A}^{2}}-B{{A}^{2}} \right)\beta _{A}^{Q}+C{{B}^{2}}\beta _{B}^{Q}-B{{C}^{2}}\beta _{C}^{Q}+\left( C{{D}^{2}}-B{{D}^{2}} \right)\beta _{D}^{Q}=0 \\ 
		& \left( D{{A}^{2}}-C{{A}^{2}} \right)\beta _{A}^{Q}+\left( D{{B}^{2}}-C{{B}^{2}} \right)\beta _{B}^{Q}+D{{C}^{2}}\beta _{C}^{Q}-C{{D}^{2}}\beta _{D}^{Q}=0. \\ 
	\end{aligned} \right.\]
	
	The following results are obtained by solving the above equations:
	\begin{flalign*}
		\beta _{A}^{Q}=\frac{{{U}_{A}}}{U},\beta _{B}^{Q}=\frac{{{U}_{B}}}{U},\beta _{C}^{Q}=\frac{{{U}_{C}}}{U},\beta _{D}^{Q}=\frac{{{U}_{D}}}{U}.	
	\end{flalign*}
	
	Where
	\[U=\left| \begin{matrix}
		1 & 1 & 1 & 1  \\
		A{{B}^{2}} & -A{{B}^{2}} & B{{C}^{2}}-A{{C}^{2}} & B{{D}^{2}}-A{{D}^{2}}  \\
		A{{C}^{2}}-A{{B}^{2}} & B{{C}^{2}} & -B{{C}^{2}} & C{{D}^{2}}-B{{D}^{2}}  \\
		A{{D}^{2}}-A{{C}^{2}} & B{{D}^{2}}-B{{C}^{2}} & C{{D}^{2}} & -C{{D}^{2}}  \\
	\end{matrix} \right|,\]
	\[\begin{aligned}
		{{U}_{A}}&=\left| \begin{matrix}
			1 & 1 & 1 & 1  \\
			0 & -A{{B}^{2}} & B{{C}^{2}}-A{{C}^{2}} & B{{D}^{2}}-A{{D}^{2}}  \\
			0 & B{{C}^{2}} & -B{{C}^{2}} & C{{D}^{2}}-B{{D}^{2}}  \\
			0 & B{{D}^{2}}-B{{C}^{2}} & C{{D}^{2}} & -C{{D}^{2}}  \\
		\end{matrix} \right|,
	\end{aligned}\]
	i.e.
	\[{{U}_{A}}=\left| \begin{matrix}
		-A{{B}^{2}} & B{{C}^{2}}-A{{C}^{2}} & B{{D}^{2}}-A{{D}^{2}}  \\
		B{{C}^{2}} & -B{{C}^{2}} & C{{D}^{2}}-B{{D}^{2}}  \\
		B{{D}^{2}}-B{{C}^{2}} & C{{D}^{2}} & -C{{D}^{2}}  \\
	\end{matrix} \right|,\]
	\[{{U}_{B}}=\left| \begin{matrix}
		1 & 1 & 1 & 1  \\
		A{{B}^{2}} & 0 & B{{C}^{2}}-A{{C}^{2}} & B{{D}^{2}}-A{{D}^{2}}  \\
		A{{C}^{2}}-A{{B}^{2}} & 0 & -B{{C}^{2}} & C{{D}^{2}}-B{{D}^{2}}  \\
		A{{D}^{2}}-A{{C}^{2}} & 0 & C{{D}^{2}} & -C{{D}^{2}}  \\
	\end{matrix} \right|,\]
	i.e.
	\[{{U}_{B}}=\left| \begin{matrix}
		-A{{B}^{2}} & A{{C}^{2}}-B{{C}^{2}} & A{{D}^{2}}-B{{D}^{2}}  \\
		A{{C}^{2}}-A{{B}^{2}} & -B{{C}^{2}} & C{{D}^{2}}-B{{D}^{2}}  \\
		A{{D}^{2}}-A{{C}^{2}} & C{{D}^{2}} & -C{{D}^{2}}  \\
	\end{matrix} \right|,\]
	\[{{U}_{C}}=\left| \begin{matrix}
		1 & 1 & 1 & 1  \\
		A{{B}^{2}} & -A{{B}^{2}} & 0 & B{{D}^{2}}-A{{D}^{2}}  \\
		A{{C}^{2}}-A{{B}^{2}} & B{{C}^{2}} & 0 & C{{D}^{2}}-B{{D}^{2}}  \\
		A{{D}^{2}}-A{{C}^{2}} & B{{D}^{2}}-B{{C}^{2}} & 0 & -C{{D}^{2}}  \\
	\end{matrix} \right|,\]
	i.e.
	\[{{U}_{C}}=\left| \begin{matrix}
		-A{{B}^{2}} & A{{B}^{2}} & A{{D}^{2}}-B{{D}^{2}}  \\
		A{{B}^{2}}-A{{C}^{2}} & -B{{C}^{2}} & B{{D}^{2}}-C{{D}^{2}}  \\
		A{{D}^{2}}-A{{C}^{2}} & B{{D}^{2}}-B{{C}^{2}} & -C{{D}^{2}}  \\
	\end{matrix} \right|,\]
	\[{{U}_{D}}=\left| \begin{matrix}
		1 & 1 & 1 & 1  \\
		A{{B}^{2}} & -A{{B}^{2}} & B{{C}^{2}}-A{{C}^{2}} & 0  \\
		A{{C}^{2}}-A{{B}^{2}} & B{{C}^{2}} & -B{{C}^{2}} & 0  \\
		A{{D}^{2}}-A{{C}^{2}} & B{{D}^{2}}-B{{C}^{2}} & C{{D}^{2}} & 0  \\
	\end{matrix} \right|,\]
	i.e.
	\[{{U}_{D}}=\left| \begin{matrix}
		-A{{B}^{2}} & A{{B}^{2}} & A{{C}^{2}}-B{{C}^{2}}  \\
		A{{B}^{2}}-A{{C}^{2}} & -B{{C}^{2}} & B{{C}^{2}}  \\
		A{{C}^{2}}-A{{D}^{2}} & B{{C}^{2}}-B{{D}^{2}} & -C{{D}^{2}}  \\
	\end{matrix} \right|.\]
	
	According to properties of determinant, expand $U$ by the first row to obtain the following results:
	\[U={{U}_{A}}+{{U}_{B}}+{{U}_{C}}+{{U}_{D}}.\]
	
	Therefore
	\[\beta _{A}^{Q}=\frac{{{U}_{A}}}{U}=\frac{{{U}_{A}}}{{{U}_{A}}+{{U}_{B}}+{{U}_{C}}+{{U}_{D}}},\]
	\[\beta _{B}^{Q}=\frac{{{U}_{B}}}{U}=\frac{{{U}_{B}}}{{{U}_{A}}+{{U}_{B}}+{{U}_{C}}+{{U}_{D}}},\]
	\[\beta _{C}^{Q}=\frac{{{U}_{C}}}{U}=\frac{{{U}_{C}}}{{{U}_{A}}+{{U}_{B}}+{{U}_{C}}+{{U}_{D}}},\]
	\[\beta _{D}^{Q}=\frac{{{U}_{D}}}{U}=\frac{{{U}_{D}}}{{{U}_{A}}+{{U}_{B}}+{{U}_{C}}+{{U}_{D}}}.\]	
\end{proof}
\hfill $\square$\par

The following theorem gives the symmetrical form of components of the circumcenter of tetrahedron. The results of the symmetrical form are more convenient for memory and application

In order to make the formula concise, some notations need to be introduced first.

The sum of the three edges of the opposite triangle ($\triangle BCD $) corresponding to vertex $A$ of the tetrahedron to the power of $n$ is denoted as:
\[\Delta _{n}^{A}=B{{C}^{n}}+C{{D}^{n}}+D{{B}^{n}}.\]

Similarly, we can get the notations of $\Delta _{n}^{B}$, $\Delta _{n}^{C}$ and $\Delta _{n}^{D}$. And $n = 2 $ or $n = 4 $ is a commonly used expression in this book. When $n = 1 $, omit $n$. 

The sum of the three edges connected to vertex $A$ of the tetrahedron to the power of $n$ is denoted as:
\[\Lambda _{n}^{A}=A{{B}^{n}}+A{{C}^{n}}+A{{D}^{n}}.\]

Similarly, we can get the notations of $\Lambda _{n}^{B}$, $\Lambda _{n}^{C}$, $\Lambda _{n}^{D}$. And $n = 2 $ or $n = 4 $ is a commonly used expression in this book. When $n = 1 $, omit $n$. 

\begin{theorem}{Frame component formulas for circumcenter of a tetrahedron-2, Daiyuan Zhang}{Thm25.1.2}\label{Thm25.1.2}
Let the lengths of the six edges of the tetrahedron $ABCD$ be $AB$, $AC$, $AD$, $BC$, $CD$, $DB$, and the point $Q$ be the circumcenter of the tetrahedron $ABCD$, then the frame components of the circumcenter of the tetrahedron $Q$ are
\begin{flalign*}
	\beta _{A}^{Q}=\frac{{{U}_{A}}}{U},\beta _{B}^{Q}=\frac{{{U}_{B}}}{U},\beta _{C}^{Q}=\frac{{{U}_{C}}}{U},\beta _{D}^{Q}=\frac{{{U}_{D}}}{U}.
\end{flalign*}

Where
\[U={{U}_{A}}+{{U}_{B}}+{{U}_{C}}+{{U}_{D}},\]

And
\[{{U}_{A}}=\frac{1}{3}\left( \begin{aligned}
	& A{{B}^{2}}\left( \Delta _{2}^{A}-2B{{C}^{2}} \right)\left( \Delta _{2}^{A}-2B{{D}^{2}} \right) \\ 
	& +A{{C}^{2}}\left( \Delta _{2}^{A}-2C{{D}^{2}} \right)\left( \Delta _{2}^{A}-2C{{B}^{2}} \right) \\ 
	& +A{{D}^{2}}\left( \Delta _{2}^{A}-2D{{B}^{2}} \right)\left( \Delta _{2}^{A}-2D{{C}^{2}} \right) \\ 
	& +B{{C}^{2}}\left( \Delta _{2}^{A}-2B{{C}^{2}} \right)\left( \Delta _{2}^{A}-\Lambda _{2}^{A}-\left( B{{C}^{2}}+A{{D}^{2}} \right) \right) \\ 
	& +C{{D}^{2}}\left( \Delta _{2}^{A}-2C{{D}^{2}} \right)\left( \Delta _{2}^{A}-\Lambda _{2}^{A}-\left( C{{D}^{2}}+A{{B}^{2}} \right) \right) \\ 
	& +D{{B}^{2}}\left( \Delta _{2}^{A}-2D{{B}^{2}} \right)\left( \Delta _{2}^{A}-\Lambda _{2}^{A}-\left( D{{B}^{2}}+A{{C}^{2}} \right) \right)  
\end{aligned} \right),\]
\[{{U}_{B}}=\frac{1}{3}\left( \begin{aligned}
	& B{{C}^{2}}\left( \Delta _{2}^{B}-2C{{D}^{2}} \right)\left( \Delta _{2}^{B}-2C{{A}^{2}} \right) \\ 
	& +B{{D}^{2}}\left( \Delta _{2}^{B}-2D{{A}^{2}} \right)\left( \Delta _{2}^{B}-2D{{C}^{2}} \right) \\ 
	& +B{{A}^{2}}\left( \Delta _{2}^{B}-2A{{C}^{2}} \right)\left( \Delta _{2}^{B}-2A{{D}^{2}} \right) \\ 
	& +C{{D}^{2}}\left( \Delta _{2}^{B}-2C{{D}^{2}} \right)\left( \Delta _{2}^{B}-\Lambda _{2}^{B}-\left( C{{D}^{2}}+B{{A}^{2}} \right) \right) \\ 
	& +D{{A}^{2}}\left( \Delta _{2}^{B}-2D{{A}^{2}} \right)\left( \Delta _{2}^{B}-\Lambda _{2}^{B}-\left( D{{A}^{2}}+B{{C}^{2}} \right) \right) \\ 
	& +A{{C}^{2}}\left( \Delta _{2}^{B}-2A{{C}^{2}} \right)\left( \Delta _{2}^{B}-\Lambda _{2}^{B}-\left( A{{C}^{2}}+B{{D}^{2}} \right) \right)  
\end{aligned} \right),\]
\[{{U}_{C}}=\frac{1}{3}\left( \begin{aligned}
	& C{{D}^{2}}\left( \Delta _{2}^{C}-2D{{A}^{2}} \right)\left( \Delta _{2}^{C}-2D{{B}^{2}} \right) \\ 
	& +C{{A}^{2}}\left( \Delta _{2}^{C}-2A{{B}^{2}} \right)\left( \Delta _{2}^{C}-2A{{D}^{2}} \right) \\ 
	& +C{{B}^{2}}\left( \Delta _{2}^{C}-2B{{D}^{2}} \right)\left( \Delta _{2}^{C}-2B{{A}^{2}} \right) \\ 
	& +D{{A}^{2}}\left( \Delta _{2}^{C}-2D{{A}^{2}} \right)\left( \Delta _{2}^{C}-\Lambda _{2}^{C}-\left( D{{A}^{2}}+C{{B}^{2}} \right) \right) \\ 
	& +A{{B}^{2}}\left( \Delta _{2}^{C}-2A{{B}^{2}} \right)\left( \Delta _{2}^{C}-\Lambda _{2}^{C}-\left( A{{B}^{2}}+C{{D}^{2}} \right) \right) \\ 
	& +B{{D}^{2}}\left( \Delta _{2}^{C}-2B{{D}^{2}} \right)\left( \Delta _{2}^{C}-\Lambda _{2}^{C}-\left( B{{D}^{2}}+C{{A}^{2}} \right) \right)  
\end{aligned} \right),\]
\[{{U}_{D}}=\frac{1}{3}\left( \begin{aligned}
	& D{{A}^{2}}\left( \Delta _{2}^{D}-2A{{B}^{2}} \right)\left( \Delta _{2}^{D}-2A{{C}^{2}} \right) \\ 
	& +D{{B}^{2}}\left( \Delta _{2}^{D}-2B{{C}^{2}} \right)\left( \Delta _{2}^{D}-2B{{A}^{2}} \right) \\ 
	& +D{{C}^{2}}\left( \Delta _{2}^{D}-2C{{A}^{2}} \right)\left( \Delta _{2}^{D}-2C{{B}^{2}} \right) \\ 
	& +A{{B}^{2}}\left( \Delta _{2}^{D}-2A{{B}^{2}} \right)\left( \Delta _{2}^{D}-\Lambda _{2}^{D}-\left( A{{B}^{2}}+D{{C}^{2}} \right) \right) \\ 
	& +B{{C}^{2}}\left( \Delta _{2}^{D}-2B{{C}^{2}} \right)\left( \Delta _{2}^{D}-\Lambda _{2}^{D}-\left( B{{C}^{2}}+D{{A}^{2}} \right) \right) \\ 
	& +C{{A}^{2}}\left( \Delta _{2}^{D}-2C{{A}^{2}} \right)\left( \Delta _{2}^{D}-\Lambda _{2}^{D}-\left( C{{A}^{2}}+D{{B}^{2}} \right) \right)  
\end{aligned} \right).\]
\end{theorem}	

\begin{proof}
	Firstly, let's calculate ${{U}_{A}}$. for the determinant ${{U}_{A}}$ in above theorem \ref{thm:Thm25.1.1}, using the properties of the determinant, three forms of the formula for calculating ${{U}_{A}}$ are obtained as follows.
	
	The first form of ${{U}_{A}}$:
	\[{{U}_{A}}=\left| \begin{matrix}
		-A{{B}^{2}} & B{{C}^{2}}-A{{C}^{2}} & B{{D}^{2}}-A{{D}^{2}}  \\
		B{{C}^{2}} & -B{{C}^{2}} & C{{D}^{2}}-B{{D}^{2}}  \\
		B{{D}^{2}}-B{{C}^{2}} & C{{D}^{2}} & -C{{D}^{2}}  \\
	\end{matrix} \right|,\]
	\[{{U}_{A}}=\left| \begin{matrix}
		-A{{B}^{2}} & B{{C}^{2}}-A{{C}^{2}}-A{{B}^{2}} & B{{D}^{2}}-A{{D}^{2}}-A{{B}^{2}}  \\
		B{{C}^{2}} & 0 & C{{D}^{2}}-B{{D}^{2}}+B{{C}^{2}}  \\
		B{{D}^{2}}-B{{C}^{2}} & C{{D}^{2}}+B{{D}^{2}}-B{{C}^{2}} & -C{{D}^{2}}+B{{D}^{2}}-B{{C}^{2}}  \\
	\end{matrix} \right|,\]
	\[{{U}_{A}}=\left| \begin{matrix}
		-A{{B}^{2}} & B{{C}^{2}}-A{{C}^{2}}-A{{B}^{2}} & B{{D}^{2}}-A{{D}^{2}}-A{{B}^{2}}  \\
		B{{C}^{2}} & 0 & C{{D}^{2}}-B{{D}^{2}}+B{{C}^{2}}  \\
		B{{D}^{2}} & C{{D}^{2}}+B{{D}^{2}}-B{{C}^{2}} & 0  \\
	\end{matrix} \right|,\]
	i.e.
	\[\begin{aligned}
		{{U}_{A}}& =A{{B}^{2}}\left( C{{D}^{2}}-B{{D}^{2}}+B{{C}^{2}} \right)\left( C{{D}^{2}}+B{{D}^{2}}-B{{C}^{2}} \right) \\ 
		& +B{{C}^{2}}\left( B{{D}^{2}}-A{{D}^{2}}-A{{B}^{2}} \right)\left( C{{D}^{2}}+B{{D}^{2}}-B{{C}^{2}} \right) \\ 
		& +B{{D}^{2}}\left( B{{C}^{2}}-A{{C}^{2}}-A{{B}^{2}} \right)\left( C{{D}^{2}}-B{{D}^{2}}+B{{C}^{2}} \right).  
	\end{aligned}\]
	
	The second form of ${{U}_{A}}$:
	\[{{U}_{A}}=\left| \begin{matrix}
		-A{{B}^{2}} & B{{C}^{2}}-A{{C}^{2}} & B{{D}^{2}}-A{{D}^{2}}  \\
		B{{C}^{2}} & -B{{C}^{2}} & C{{D}^{2}}-B{{D}^{2}}  \\
		B{{D}^{2}}-B{{C}^{2}} & C{{D}^{2}} & -C{{D}^{2}}  \\
	\end{matrix} \right|,\]
	\[{{U}_{A}}=\left| \begin{matrix}
		-A{{B}^{2}}+B{{C}^{2}}-A{{C}^{2}} & B{{C}^{2}}-A{{C}^{2}} & B{{D}^{2}}-A{{D}^{2}}  \\
		0 & -B{{C}^{2}} & C{{D}^{2}}-B{{D}^{2}}  \\
		B{{D}^{2}}-B{{C}^{2}}+C{{D}^{2}} & C{{D}^{2}} & -C{{D}^{2}}  \\
	\end{matrix} \right|,\]
	\[{{U}_{A}}=\left| \begin{matrix}
		-A{{B}^{2}}+B{{C}^{2}}-A{{C}^{2}} & B{{C}^{2}}-A{{C}^{2}}+B{{D}^{2}}-A{{D}^{2}} & B{{D}^{2}}-A{{D}^{2}}  \\
		0 & -B{{C}^{2}}+C{{D}^{2}}-B{{D}^{2}} & C{{D}^{2}}-B{{D}^{2}}  \\
		B{{D}^{2}}-B{{C}^{2}}+C{{D}^{2}} & 0 & -C{{D}^{2}}  \\
	\end{matrix} \right|,\]
	\[{{U}_{A}}=\left| \begin{matrix}
		-A{{B}^{2}}+B{{C}^{2}}-A{{C}^{2}} & -A{{C}^{2}}-A{{D}^{2}}+C{{D}^{2}} & C{{D}^{2}}-A{{D}^{2}}  \\
		0 & -B{{C}^{2}}+C{{D}^{2}}-B{{D}^{2}} & C{{D}^{2}}-B{{D}^{2}}  \\
		B{{D}^{2}}-B{{C}^{2}}+C{{D}^{2}} & 0 & -C{{D}^{2}}  \\
	\end{matrix} \right|,\]
	\[{{U}_{A}}=\left| \begin{matrix}
		A{{B}^{2}}-B{{C}^{2}}+A{{C}^{2}} & A{{C}^{2}}+A{{D}^{2}}-C{{D}^{2}} & A{{C}^{2}}  \\
		0 & B{{C}^{2}}-C{{D}^{2}}+B{{D}^{2}} & B{{C}^{2}}  \\
		-B{{D}^{2}}+B{{C}^{2}}-C{{D}^{2}} & 0 & -C{{D}^{2}}  \\
	\end{matrix} \right|,\]
	i.e.
	\[\begin{aligned}
		{{U}_{A}}& =-C{{D}^{2}}\left( A{{B}^{2}}-B{{C}^{2}}+A{{C}^{2}} \right)\left( B{{C}^{2}}-C{{D}^{2}}+B{{D}^{2}} \right) \\ 
		& +B{{C}^{2}}\left( A{{C}^{2}}+A{{D}^{2}}-C{{D}^{2}} \right)\left( -B{{D}^{2}}+B{{C}^{2}}-C{{D}^{2}} \right) \\ 
		& -A{{C}^{2}}\left( B{{C}^{2}}-C{{D}^{2}}+B{{D}^{2}} \right)\left( -B{{D}^{2}}+B{{C}^{2}}-C{{D}^{2}} \right).  
	\end{aligned}\]
	
	The third form of ${{U}_{A}}$: 
	\[{{U}_{A}}=\left| \begin{matrix}
		-A{{B}^{2}} & B{{C}^{2}}-A{{C}^{2}} & B{{D}^{2}}-A{{D}^{2}}  \\
		B{{C}^{2}} & -B{{C}^{2}} & C{{D}^{2}}-B{{D}^{2}}  \\
		B{{D}^{2}}-B{{C}^{2}} & C{{D}^{2}} & -C{{D}^{2}}  \\
	\end{matrix} \right|,\]
	\[{{U}_{A}}=\left| \begin{matrix}
		-A{{B}^{2}}+B{{D}^{2}}-A{{D}^{2}} & B{{C}^{2}}-A{{C}^{2}}+B{{D}^{2}}-A{{D}^{2}} & B{{D}^{2}}-A{{D}^{2}}  \\
		B{{C}^{2}}+C{{D}^{2}}-B{{D}^{2}} & -B{{C}^{2}}+C{{D}^{2}}-B{{D}^{2}} & C{{D}^{2}}-B{{D}^{2}}  \\
		B{{D}^{2}}-B{{C}^{2}}-C{{D}^{2}} & 0 & -C{{D}^{2}}  \\
	\end{matrix} \right|,\]
	\[{{U}_{A}}=\left| \begin{matrix}
		-A{{B}^{2}}+B{{D}^{2}}-A{{D}^{2}} & B{{C}^{2}}-A{{C}^{2}}+B{{D}^{2}}-A{{D}^{2}} & B{{D}^{2}}-A{{D}^{2}}  \\
		0 & -B{{C}^{2}}+C{{D}^{2}}-B{{D}^{2}} & -B{{D}^{2}}  \\
		B{{D}^{2}}-B{{C}^{2}}-C{{D}^{2}} & 0 & -C{{D}^{2}}  \\
	\end{matrix} \right|,\]
	\[{{U}_{A}}=\left| \begin{matrix}
		-A{{B}^{2}}+B{{D}^{2}}-A{{D}^{2}} & -A{{C}^{2}}-A{{D}^{2}}+C{{D}^{2}} & -A{{D}^{2}}  \\
		0 & -B{{C}^{2}}+C{{D}^{2}}-B{{D}^{2}} & -B{{D}^{2}}  \\
		B{{D}^{2}}-B{{C}^{2}}-C{{D}^{2}} & 0 & -C{{D}^{2}}  \\
	\end{matrix} \right|,\]
	that is
	\[\begin{aligned}
		{{U}_{A}}& =-C{{D}^{2}}\left( -A{{B}^{2}}+B{{D}^{2}}-A{{D}^{2}} \right)\left( -B{{C}^{2}}+C{{D}^{2}}-B{{D}^{2}} \right) \\ 
		& -B{{D}^{2}}\left( C{{D}^{2}}-A{{C}^{2}}-A{{D}^{2}} \right)\left( B{{D}^{2}}-B{{C}^{2}}-C{{D}^{2}} \right) \\ 
		& +A{{D}^{2}}\left( -B{{C}^{2}}+C{{D}^{2}}-B{{D}^{2}} \right)\left( B{{D}^{2}}-B{{C}^{2}}-C{{D}^{2}} \right).  
	\end{aligned}\]
	
	Thus, three formulas for calculating ${{U}_{A}}$ expressed by the lengths of the edges are obtained as follows:
	\[\begin{aligned}
		{{U}_{A}}&=A{{B}^{2}}\left( C{{D}^{2}}-B{{D}^{2}}+B{{C}^{2}} \right)\left( C{{D}^{2}}+B{{D}^{2}}-B{{C}^{2}} \right) \\ 
		& +B{{C}^{2}}\left( B{{D}^{2}}-A{{D}^{2}}-A{{B}^{2}} \right)\left( C{{D}^{2}}+B{{D}^{2}}-B{{C}^{2}} \right) \\ 
		& +B{{D}^{2}}\left( B{{C}^{2}}-A{{C}^{2}}-A{{B}^{2}} \right)\left( C{{D}^{2}}-B{{D}^{2}}+B{{C}^{2}} \right),  
	\end{aligned}\]
	\[\begin{aligned}
		{{U}_{A}}& =-C{{D}^{2}}\left( A{{B}^{2}}-B{{C}^{2}}+A{{C}^{2}} \right)\left( B{{C}^{2}}-C{{D}^{2}}+B{{D}^{2}} \right) \\ 
		& +B{{C}^{2}}\left( A{{C}^{2}}+A{{D}^{2}}-C{{D}^{2}} \right)\left( -B{{D}^{2}}+B{{C}^{2}}-C{{D}^{2}} \right) \\ 
		& -A{{C}^{2}}\left( B{{C}^{2}}-C{{D}^{2}}+B{{D}^{2}} \right)\left( -B{{D}^{2}}+B{{C}^{2}}-C{{D}^{2}} \right),  
	\end{aligned}\]
	\[\begin{aligned}
		{{U}_{A}}&=-C{{D}^{2}}\left( -A{{B}^{2}}+B{{D}^{2}}-A{{D}^{2}} \right)\left( -B{{C}^{2}}+C{{D}^{2}}-B{{D}^{2}} \right) \\ 
		& -B{{D}^{2}}\left( C{{D}^{2}}-A{{C}^{2}}-A{{D}^{2}} \right)\left( B{{D}^{2}}-B{{C}^{2}}-C{{D}^{2}} \right) \\ 
		& +A{{D}^{2}}\left( -B{{C}^{2}}+C{{D}^{2}}-B{{D}^{2}} \right)\left( B{{D}^{2}}-B{{C}^{2}}-C{{D}^{2}} \right).  
	\end{aligned}\]
	
	In order to obtain the symmetrical form of edges, sum the two sides of the above three formulas, and the items on the right can be grouped with like items. The specific calculation process is as follows.
	
	
	
	Grouping like terms for $B{{C}^{2}}$, we have:
	\[\begin{aligned}
		& B{{C}^{2}}\left( B{{D}^{2}}-A{{D}^{2}}-A{{B}^{2}} \right)\left( C{{D}^{2}}+B{{D}^{2}}-B{{C}^{2}} \right) \\ 
		& +B{{C}^{2}}\left( A{{C}^{2}}+A{{D}^{2}}-C{{D}^{2}} \right)\left( -B{{D}^{2}}+B{{C}^{2}}-C{{D}^{2}} \right), \\ 
	\end{aligned}\]
	i.e.
	\[B{{C}^{2}}\left( C{{D}^{2}}+B{{D}^{2}}-B{{C}^{2}} \right)\left( \left( B{{D}^{2}}-A{{D}^{2}}-A{{B}^{2}} \right)-\left( A{{C}^{2}}+A{{D}^{2}}-C{{D}^{2}} \right) \right),\]
	i.e.
	\[B{{C}^{2}}\left( C{{D}^{2}}+D{{B}^{2}}-B{{C}^{2}} \right)\left( B{{D}^{2}}+C{{D}^{2}}-A{{D}^{2}}-\left( A{{B}^{2}}+A{{C}^{2}}+A{{D}^{2}} \right) \right),\]
	i.e.
	\[B{{C}^{2}}\left( \left( B{{C}^{2}}+C{{D}^{2}}+D{{B}^{2}} \right)-2B{{C}^{2}} \right)\left( \begin{aligned}
		& \left( B{{C}^{2}}+C{{D}^{2}}+D{{B}^{2}} \right) \\ 
		& -\left( A{{B}^{2}}+A{{C}^{2}}+A{{D}^{2}} \right) \\ 
		& -\left( B{{C}^{2}}+A{{D}^{2}} \right) \\ 
	\end{aligned} \right),\]
	i.e.
	\[B{{C}^{2}}\left( {{\Delta }^{A}}-2B{{C}^{2}} \right)\left( {{\Delta }^{A}}-{{\Lambda }^{A}}-\left( B{{C}^{2}}+A{{D}^{2}} \right) \right).\]
	
	Grouping like terms for $C{{D}^{2}}$, we have:
	\[\begin{aligned}
		& -C{{D}^{2}}\left( A{{B}^{2}}-B{{C}^{2}}+A{{C}^{2}} \right)\left( B{{C}^{2}}-C{{D}^{2}}+B{{D}^{2}} \right) \\ 
		& -C{{D}^{2}}\left( -A{{B}^{2}}+B{{D}^{2}}-A{{D}^{2}} \right)\left( -B{{C}^{2}}+C{{D}^{2}}-B{{D}^{2}} \right), \\ 
	\end{aligned}\]
	i.e.
	\[C{{D}^{2}}\left( B{{C}^{2}}-C{{D}^{2}}+B{{D}^{2}} \right)\left( -\left( A{{B}^{2}}-B{{C}^{2}}+A{{C}^{2}} \right)+\left( -A{{B}^{2}}+B{{D}^{2}}-A{{D}^{2}} \right) \right),\]
	i.e.
	\[C{{D}^{2}}\left( B{{C}^{2}}-C{{D}^{2}}+B{{D}^{2}} \right)\left( B{{C}^{2}}+B{{D}^{2}}-A{{B}^{2}}-\left( A{{B}^{2}}+A{{C}^{2}}+A{{D}^{2}} \right) \right),\]
	i.e.
	\[C{{D}^{2}}\left( B{{C}^{2}}+C{{D}^{2}}+D{{B}^{2}}-2C{{D}^{2}} \right)\left( \begin{aligned}
		& B{{C}^{2}}+C{{D}^{2}}+D{{B}^{2}} \\ 
		& -\left( A{{B}^{2}}+A{{C}^{2}}+A{{D}^{2}} \right) \\ 
		& -\left( A{{B}^{2}}+C{{D}^{2}} \right) \\ 
	\end{aligned} \right),\]
	i.e.
	\[C{{D}^{2}}\left( {{\Delta }^{A}}-2C{{D}^{2}} \right)\left( {{\Delta }^{A}}-{{\Lambda }^{A}}-\left( C{{D}^{2}}+A{{B}^{2}} \right) \right).\]
	
	Grouping like terms for $D{{B}^{2}}$, we have:
	\[\begin{aligned}
		& B{{D}^{2}}\left( B{{C}^{2}}-A{{C}^{2}}-A{{B}^{2}} \right)\left( C{{D}^{2}}-B{{D}^{2}}+B{{C}^{2}} \right) \\ 
		& -B{{D}^{2}}\left( C{{D}^{2}}-A{{C}^{2}}-A{{D}^{2}} \right)\left( B{{D}^{2}}-B{{C}^{2}}-C{{D}^{2}} \right), \\ 
	\end{aligned}\]
	i.e.
	\[B{{D}^{2}}\left( C{{D}^{2}}-B{{D}^{2}}+B{{C}^{2}} \right)\left( \left( B{{C}^{2}}-A{{C}^{2}}-A{{B}^{2}} \right)+\left( C{{D}^{2}}-A{{C}^{2}}-A{{D}^{2}} \right) \right),\]
	i.e.
	\[B{{D}^{2}}\left( B{{C}^{2}}+C{{D}^{2}}+D{{B}^{2}}-2D{{B}^{2}} \right)\left( \begin{aligned}
		& \left( B{{C}^{2}}+C{{D}^{2}}+D{{B}^{2}} \right) \\ 
		& -\left( A{{B}^{2}}+A{{C}^{2}}+A{{D}^{2}} \right) \\ 
		& -\left( D{{B}^{2}}+A{{C}^{2}} \right) \\ 
	\end{aligned} \right).\]
	i.e.
	\[D{{B}^{2}}\left( {{\Delta }^{A}}-2D{{B}^{2}} \right)\left( {{\Delta }^{A}}-{{\Lambda }^{A}}-\left( D{{B}^{2}}+A{{C}^{2}} \right) \right).\]
	
	For items containing $A{{B}^{2}}$, $A{{C}^{2}}$, $A{{D}^{2}}$:
	\[\begin{aligned}
		& A{{B}^{2}}\left( C{{D}^{2}}-B{{D}^{2}}+B{{C}^{2}} \right)\left( C{{D}^{2}}+B{{D}^{2}}-B{{C}^{2}} \right) \\ 
		& -A{{C}^{2}}\left( B{{C}^{2}}-C{{D}^{2}}+B{{D}^{2}} \right)\left( -B{{D}^{2}}+B{{C}^{2}}-C{{D}^{2}} \right) \\ 
		& +A{{D}^{2}}\left( -B{{C}^{2}}+C{{D}^{2}}-B{{D}^{2}} \right)\left( B{{D}^{2}}-B{{C}^{2}}-C{{D}^{2}} \right), \\ 
	\end{aligned} \]
	i.e.
	\[\begin{aligned}
		& A{{B}^{2}}\left( B{{C}^{2}}+C{{D}^{2}}+D{{B}^{2}}-2B{{D}^{2}} \right)\left( B{{C}^{2}}+C{{D}^{2}}+D{{B}^{2}}-2B{{C}^{2}} \right) \\ 
		& +A{{C}^{2}}\left( B{{C}^{2}}+C{{D}^{2}}+D{{B}^{2}}-2C{{D}^{2}} \right)\left( B{{C}^{2}}+C{{D}^{2}}+D{{B}^{2}}-2B{{C}^{2}} \right) \\ 
		& +A{{D}^{2}}\left( B{{C}^{2}}+C{{D}^{2}}+D{{B}^{2}}-2C{{D}^{2}} \right)\left( B{{C}^{2}}+C{{D}^{2}}+D{{B}^{2}}-2B{{D}^{2}} \right), \\ 
	\end{aligned}\]
	i.e.
	\[\begin{aligned}
		& A{{B}^{2}}\left( {{\Delta }^{A}}-2B{{C}^{2}} \right)\left( {{\Delta }^{A}}-2B{{D}^{2}} \right) \\ 
		& +A{{C}^{2}}\left( {{\Delta }^{A}}-2C{{D}^{2}} \right)\left( {{\Delta }^{A}}-2C{{B}^{2}} \right) \\ 
		& +A{{D}^{2}}\left( {{\Delta }^{A}}-2D{{B}^{2}} \right)\left( {{\Delta }^{A}}-2D{{C}^{2}} \right). \\ 
	\end{aligned}\]
	
	Therefore
	\[{{U}_{A}}=\frac{1}{3}\left( \begin{aligned}
		& A{{B}^{2}}\left( {{\Delta }^{A}}-2B{{C}^{2}} \right)\left( {{\Delta }^{A}}-2B{{D}^{2}} \right) \\ 
		& +A{{C}^{2}}\left( {{\Delta }^{A}}-2C{{D}^{2}} \right)\left( {{\Delta }^{A}}-2C{{B}^{2}} \right) \\ 
		& +A{{D}^{2}}\left( {{\Delta }^{A}}-2D{{B}^{2}} \right)\left( {{\Delta }^{A}}-2D{{C}^{2}} \right) \\ 
		& +B{{C}^{2}}\left( {{\Delta }^{A}}-2B{{C}^{2}} \right)\left( {{\Delta }^{A}}-{{\Lambda }^{A}}-\left( B{{C}^{2}}+A{{D}^{2}} \right) \right) \\ 
		& +C{{D}^{2}}\left( {{\Delta }^{A}}-2C{{D}^{2}} \right)\left( {{\Delta }^{A}}-{{\Lambda }^{A}}-\left( C{{D}^{2}}+A{{B}^{2}} \right) \right) \\ 
		& +D{{B}^{2}}\left( {{\Delta }^{A}}-2D{{B}^{2}} \right)\left( {{\Delta }^{A}}-{{\Lambda }^{A}}-\left( D{{B}^{2}}+A{{C}^{2}} \right) \right)  
	\end{aligned} \right).\]
	
	Similar results can be obtained for ${{U}_{B}}$, ${{U}_{C}}$ and ${{U}_{D}}$. 
\end{proof}
\hfill $\square$\par

\section{Frame components of the circumcenter for some special tetrahedrons}\label{Sec25.2}
\subsection{Frame components of the circumcenter for regular triangular pyramid}\label{Subsec25.2.1}
As an application, let's find the frame components of the circumcenter of a regular triangular pyramid. For a regular triangular pyramid $A-BCD$, let $AB=AC=AD=l$, $BC=CD=DB=a$. Obviously, both the frame components of the circumcenter of a triangular pyramid can be calculated directly by using theorem \ref{thm:Thm25.1.1} and the formula in symmetric form (theorem \ref{thm:Thm25.1.2}). The following uses the formula of symmetrical form to calculate the frame components.

\[\Delta _{2}^{A}=B{{C}^{2}}+C{{D}^{2}}+D{{B}^{2}}=3{{a}^{2}},\]
\[\Lambda _{2}^{A}=A{{B}^{2}}+A{{C}^{2}}+A{{D}^{2}}=3{{l}^{2}},\]
\[{{U}_{A}}=\frac{1}{3}\left( \begin{aligned}
	& A{{B}^{2}}\left( \Delta _{2}^{A}-2B{{C}^{2}} \right)\left( \Delta _{2}^{A}-2B{{D}^{2}} \right) \\ 
	& +A{{C}^{2}}\left( \Delta _{2}^{A}-2C{{D}^{2}} \right)\left( \Delta _{2}^{A}-2C{{B}^{2}} \right) \\ 
	& +A{{D}^{2}}\left( \Delta _{2}^{A}-2D{{B}^{2}} \right)\left( \Delta _{2}^{A}-2D{{C}^{2}} \right) \\ 
	& +B{{C}^{2}}\left( \Delta _{2}^{A}-2B{{C}^{2}} \right)\left( \Delta _{2}^{A}-\Lambda _{2}^{A}-\left( B{{C}^{2}}+A{{D}^{2}} \right) \right) \\ 
	& +C{{D}^{2}}\left( \Delta _{2}^{A}-2C{{D}^{2}} \right)\left( \Delta _{2}^{A}-\Lambda _{2}^{A}-\left( C{{D}^{2}}+A{{B}^{2}} \right) \right) \\ 
	& +D{{B}^{2}}\left( \Delta _{2}^{A}-2D{{B}^{2}} \right)\left( \Delta _{2}^{A}-\Lambda _{2}^{A}-\left( D{{B}^{2}}+A{{C}^{2}} \right) \right)  
\end{aligned} \right).\]
i.e.
\[\begin{aligned}
	{{U}_{A}}& =\frac{1}{3}\left( \begin{aligned}
		& {{l}^{2}}\left( 3{{a}^{2}}-2{{a}^{2}} \right)\left( 3{{a}^{2}}-2{{a}^{2}} \right) \\ 
		& +{{l}^{2}}\left( 3{{a}^{2}}-2{{a}^{2}} \right)\left( 3{{a}^{2}}-2{{a}^{2}} \right) \\ 
		& +{{l}^{2}}\left( 3{{a}^{2}}-2{{a}^{2}} \right)\left( 3{{a}^{2}}-2{{a}^{2}} \right) \\ 
		& +{{a}^{2}}\left( 3{{a}^{2}}-2{{a}^{2}} \right)\left( 3{{a}^{2}}-3{{l}^{2}}-\left( {{a}^{2}}+{{l}^{2}} \right) \right) \\ 
		& +{{a}^{2}}\left( 3{{a}^{2}}-2{{a}^{2}} \right)\left( 3{{a}^{2}}-3{{l}^{2}}-\left( {{a}^{2}}+{{l}^{2}} \right) \right) \\ 
		& +{{a}^{2}}\left( 3{{a}^{2}}-2{{a}^{2}} \right)\left( 3{{a}^{2}}-3{{l}^{2}}-\left( {{a}^{2}}+{{l}^{2}} \right) \right)  
	\end{aligned} \right), \\ 
	& ={{a}^{4}}{{l}^{2}}+{{a}^{4}}\left( 2{{a}^{2}}-4{{l}^{2}} \right)={{a}^{4}}\left( 2{{a}^{2}}-3{{l}^{2}} \right)  
\end{aligned}\]

And 
\[\Delta _{2}^{B}=C{{D}^{2}}+D{{A}^{2}}+A{{C}^{2}}={{a}^{2}}+2{{l}^{2}},\]
\[\Lambda _{2}^{B}=B{{C}^{2}}+B{{D}^{2}}+B{{A}^{2}}=2{{a}^{2}}+{{l}^{2}},\]
\[{{U}_{B}}=\frac{1}{3}\left( \begin{aligned}
	& B{{C}^{2}}\left( \Delta _{2}^{B}-2C{{D}^{2}} \right)\left( \Delta _{2}^{B}-2C{{A}^{2}} \right) \\ 
	& +B{{D}^{2}}\left( \Delta _{2}^{B}-2D{{A}^{2}} \right)\left( \Delta _{2}^{B}-2D{{C}^{2}} \right) \\ 
	& +B{{A}^{2}}\left( \Delta _{2}^{B}-2A{{C}^{2}} \right)\left( \Delta _{2}^{B}-2A{{D}^{2}} \right) \\ 
	& +C{{D}^{2}}\left( \Delta _{2}^{B}-2C{{D}^{2}} \right)\left( \Delta _{2}^{B}-\Lambda _{2}^{B}-\left( C{{D}^{2}}+B{{A}^{2}} \right) \right) \\ 
	& +D{{A}^{2}}\left( \Delta _{2}^{B}-2D{{A}^{2}} \right)\left( \Delta _{2}^{B}-\Lambda _{2}^{B}-\left( D{{A}^{2}}+B{{C}^{2}} \right) \right) \\ 
	& +A{{C}^{2}}\left( \Delta _{2}^{B}-2A{{C}^{2}} \right)\left( \Delta _{2}^{B}-\Lambda _{2}^{B}-\left( A{{C}^{2}}+B{{D}^{2}} \right) \right)  
\end{aligned} \right),\]
i.e.
\[\begin{aligned}
	{{U}_{B}}& =\frac{1}{3}\left( \begin{aligned}
		& {{a}^{2}}\left( {{a}^{2}}+2{{l}^{2}}-2{{a}^{2}} \right)\left( {{a}^{2}}+2{{l}^{2}}-2{{l}^{2}} \right) \\ 
		& +{{a}^{2}}\left( {{a}^{2}}+2{{l}^{2}}-2{{l}^{2}} \right)\left( {{a}^{2}}+2{{l}^{2}}-2{{a}^{2}} \right) \\ 
		& +{{l}^{2}}\left( {{a}^{2}}+2{{l}^{2}}-2{{l}^{2}} \right)\left( {{a}^{2}}+2{{l}^{2}}-2{{l}^{2}} \right) \\ 
		& +{{a}^{2}}\left( {{a}^{2}}+2{{l}^{2}}-2{{a}^{2}} \right)\left( {{a}^{2}}+2{{l}^{2}}-\left( 2{{a}^{2}}+{{l}^{2}} \right)-\left( {{a}^{2}}+{{l}^{2}} \right) \right) \\ 
		& +{{l}^{2}}\left( {{a}^{2}}+2{{l}^{2}}-2{{l}^{2}} \right)\left( {{a}^{2}}+2{{l}^{2}}-\left( 2{{a}^{2}}+{{l}^{2}} \right)-\left( {{l}^{2}}+{{a}^{2}} \right) \right) \\ 
		& +{{l}^{2}}\left( {{a}^{2}}+2{{l}^{2}}-2{{l}^{2}} \right)\left( {{a}^{2}}+2{{l}^{2}}-\left( 2{{a}^{2}}+{{l}^{2}} \right)-\left( {{l}^{2}}+{{a}^{2}} \right) \right)  
	\end{aligned} \right) \\ 
	& =\frac{1}{3}\left( 2{{a}^{4}}\left( 2{{l}^{2}}-{{a}^{2}} \right)+{{l}^{2}}{{a}^{4}}-2{{a}^{4}}\left( 2{{l}^{2}}-{{a}^{2}} \right)-4{{a}^{4}}{{l}^{2}} \right)=-{{a}^{4}}{{l}^{2}}  
\end{aligned}.\]

Qbviously,
\[{{U}_{C}}={{U}_{D}}={{U}_{B}}=-{{a}^{4}}{{l}^{2}},\]

Therefore
\[\beta _{A}^{Q}=\frac{{{U}_{A}}}{U}=\frac{{{U}_{A}}}{{{U}_{A}}+{{U}_{B}}+{{U}_{C}}+{{U}_{D}}}=\frac{\left( 2{{a}^{2}}-3{{l}^{2}} \right){{a}^{4}}}{\left( 2{{a}^{2}}-3{{l}^{2}} \right){{a}^{4}}-3{{l}^{2}}{{a}^{4}}}=\frac{3{{l}^{2}}-2{{a}^{2}}}{6{{l}^{2}}-2{{a}^{2}}},\]
\[\begin{aligned}
	\beta _{B}^{Q}& =\beta _{C}^{Q}=\beta _{D}^{Q}=\frac{{{U}_{B}}}{U}=\frac{{{U}_{B}}}{{{U}_{A}}+{{U}_{B}}+{{U}_{C}}+{{U}_{D}}} \\ 
	& =\frac{-{{l}^{2}}{{a}^{4}}}{\left( 2{{a}^{2}}-3{{l}^{2}} \right){{a}^{4}}-3{{l}^{2}}{{a}^{4}}}=\frac{{{l}^{2}}}{6{{l}^{2}}-2{{a}^{2}}}.  
\end{aligned}\]

\subsection{Frame components of the circumcenter for regular tetrahedron}\label{Subsec25.2.2}
According to the previous results, let $l=a$, and the frame components of the circumcenter of the regular tetrahedron are obtained:
\[\beta _{A}^{Q}=\frac{3{{a}^{2}}-2{{a}^{2}}}{6{{a}^{2}}-2{{a}^{2}}}=\frac{1}{4},\]
\[\beta _{B}^{Q}=\beta _{C}^{Q}=\beta _{D}^{Q}=\frac{{{a}^{2}}}{6{{a}^{2}}-2{{a}^{2}}}=\frac{1}{4}.\]

This is precisely the frame components of the centroid of the tetrahedron, indicating that circumcenter of the regular tetrahedron coincides with the centroid.

\section{Frame components of the circumcenter of a tetrahedron(Continued)}\label{Sec25.3}
This section will give other calculation formulas of the frame components for the circumcenter of a given tetrahedron, which is more concise than the previous formulas.

\begin{theorem}{Frame component formulas for circumcenter of a tetrahedron-3, Daiyuan Zhang}{Thm25.3.1}\label{Thm25.3.1} 
	Let the lengths of the six edges of the tetrahedron $ABCD$ be $AB$, $AC$, $AD$, $BC$, $CD$, $DB$,  the point $Q$ be the circumcenter of the tetrahedron $ABCD$, and the frame components of point $Q$ on the tetrahedral frame $\left( O;A,B,C,D \right)$ be $\beta _{A}^{Q}$, $\beta _{B}^{Q}$, $\beta _{C}^{Q}$, $\beta _{D}^{Q}$, then:
	\begin{equation}\label{25.3.1}
		\beta _{A}^{Q}=\frac{{{U}_{A}}}{U},\ \beta _{B}^{Q}=\frac{{{U}_{B}}}{U},\ \beta _{C}^{Q}=\frac{{{U}_{C}}}{U},\ \beta _{D}^{Q}=\frac{{{U}_{D}}}{U}.
	\end{equation}
	Where 
	\[\begin{aligned}
		{{U}_{A}}& =\left( \Delta _{2}^{A}-B{{C}^{2}} \right)B{{C}^{2}}A{{D}^{2}}+\left( \Delta _{2}^{A}-C{{D}^{2}} \right)C{{D}^{2}}A{{B}^{2}} \\ 
		& +\left( \Delta _{2}^{A}-D{{B}^{2}} \right)D{{B}^{2}}A{{C}^{2}}-B{{C}^{2}}C{{D}^{2}}D{{B}^{2}},  
	\end{aligned}\]
	\[\Delta _{2}^{A}=\frac{1}{2}\left( B{{C}^{2}}+C{{D}^{2}}+D{{B}^{2}} \right);\]
	\[\begin{aligned}
		{{U}_{B}}& =\left( \Delta _{2}^{B}-C{{D}^{2}} \right)C{{D}^{2}}B{{A}^{2}}+\left( \Delta _{2}^{B}-D{{A}^{2}} \right)D{{A}^{2}}B{{C}^{2}} \\ 
		& +\left( \Delta _{2}^{B}-A{{C}^{2}} \right)A{{C}^{2}}B{{D}^{2}}-C{{D}^{2}}D{{A}^{2}}A{{C}^{2}},  
	\end{aligned}\]
	\[\Delta _{2}^{B}=\frac{1}{2}\left( C{{D}^{2}}+D{{A}^{2}}+A{{C}^{2}} \right);\]
	\[\begin{aligned}
		{{U}_{C}}& =\left( \Delta _{2}^{C}-D{{A}^{2}} \right)D{{A}^{2}}C{{B}^{2}}+\left( \Delta _{2}^{C}-A{{B}^{2}} \right)A{{B}^{2}}C{{D}^{2}} \\ 
		& +\left( \Delta _{2}^{C}-B{{D}^{2}} \right)B{{D}^{2}}C{{A}^{2}}-D{{A}^{2}}A{{B}^{2}}B{{D}^{2}},  
	\end{aligned}\]
	\[\Delta _{2}^{C}=\frac{1}{2}\left( D{{A}^{2}}+A{{B}^{2}}+B{{D}^{2}} \right);\]
	\[\begin{aligned}
		{{U}_{D}}& =\left( \Delta _{2}^{D}-A{{B}^{2}} \right)A{{B}^{2}}D{{C}^{2}}+\left( \Delta _{2}^{D}-B{{C}^{2}} \right)B{{C}^{2}}D{{A}^{2}} \\ 
		& +\left( \Delta _{2}^{D}-C{{A}^{2}} \right)C{{A}^{2}}D{{B}^{2}}-A{{B}^{2}}B{{C}^{2}}C{{A}^{2}},  
	\end{aligned}\]
	\[\Delta _{2}^{D}=\frac{1}{2}\left( A{{B}^{2}}+B{{C}^{2}}+C{{A}^{2}} \right);\]
	\[U={{U}_{A}}+{{U}_{B}}+{{U}_{C}}+{{U}_{D}}=4\left( {{t}_{1}}-{{t}_{2}}-{{t}_{3}} \right),\]
	\[{{t}_{1}}={{q}_{2}}{{\Delta }_{2}},\]
	\[{{q}_{2}}=\frac{1}{2}\left( A{{B}^{2}}C{{D}^{2}}+B{{C}^{2}}A{{D}^{2}}+C{{A}^{2}}B{{D}^{2}} \right),\]
	\[{{\Delta }_{2}}=\frac{1}{2}\left( A{{B}^{2}}+A{{C}^{2}}+A{{D}^{2}}+B{{C}^{2}}+C{{D}^{2}}+D{{B}^{2}} \right),\]
	\[{{t}_{2}}=\frac{1}{2}\left( A{{B}^{2}}C{{D}^{2}}\left( A{{B}^{2}}+C{{D}^{2}} \right)+B{{C}^{2}}A{{D}^{2}}\left( B{{C}^{2}}+A{{D}^{2}} \right)+C{{A}^{2}}B{{D}^{2}}\left( C{{A}^{2}}+B{{D}^{2}} \right) \right),\]
	\[{{t}_{3}}=\frac{1}{4}\left( A{{B}^{2}}B{{C}^{2}}C{{A}^{2}}+B{{C}^{2}}C{{D}^{2}}D{{B}^{2}}+C{{D}^{2}}D{{A}^{2}}A{{C}^{2}}+D{{A}^{2}}A{{B}^{2}}B{{D}^{2}} \right).\]

\end{theorem}

\begin{proof}
	Refer to figure \ref{fig:tu16.2.1}, based on theorem \ref{thm:Thm24.1.2}, we get:
	\[\begin{aligned}
		A{{Q}^{2}}& =\left( 1-\beta _{A}^{Q} \right)\left( \beta _{B}^{Q}A{{B}^{2}}+\beta _{C}^{Q}A{{C}^{2}}+\beta _{D}^{Q}A{{D}^{2}} \right) \\ 
		& -\beta _{B}^{Q}\beta _{C}^{Q}B{{C}^{2}}-\beta _{C}^{Q}\beta _{D}^{Q}C{{D}^{2}}-\beta _{D}^{Q}\beta _{B}^{Q}D{{B}^{2}},  
	\end{aligned}\]	
	\[\begin{aligned}
		B{{Q}^{2}}& =\left( 1-\beta _{B}^{Q} \right)\left( \beta _{C}^{Q}B{{C}^{2}}+\beta _{D}^{Q}B{{D}^{2}}+\beta _{A}^{Q}B{{A}^{2}} \right) \\ 
		& -\beta _{C}^{Q}\beta _{D}^{Q}C{{D}^{2}}-\beta _{D}^{Q}\beta _{A}^{Q}D{{A}^{2}}-\beta _{A}^{Q}\beta _{C}^{Q}A{{C}^{2}},  
	\end{aligned}\]	
	\[\begin{aligned}
		C{{Q}^{2}}& =\left( 1-\beta _{C}^{Q} \right)\left( \beta _{D}^{Q}C{{D}^{2}}+\beta _{A}^{Q}C{{A}^{2}}+\beta _{B}^{Q}C{{B}^{2}} \right) \\ 
		& -\beta _{D}^{Q}\beta _{A}^{Q}D{{A}^{2}}-\beta _{A}^{Q}\beta _{B}^{Q}A{{B}^{2}}-\beta _{B}^{Q}\beta _{D}^{Q}B{{D}^{2}},  
	\end{aligned}\]	
	\[\begin{aligned}
		D{{Q}^{2}}& =\left( 1-\beta _{D}^{Q} \right)\left( \beta _{A}^{Q}D{{A}^{2}}+\beta _{B}^{Q}D{{B}^{2}}+\beta _{C}^{Q}D{{C}^{2}} \right) \\ 
		& -\beta _{A}^{Q}\beta _{B}^{Q}A{{B}^{2}}-\beta _{B}^{Q}\beta _{C}^{Q}B{{C}^{2}}-\beta _{C}^{Q}\beta _{A}^{Q}C{{A}^{2}}.  
	\end{aligned}\]	
	
	For circumcenter $Q$, we have: $A{{Q}^{2}}=B{{Q}^{2}}=C{{Q}^{2}}=D{{Q}^{2}}={{R}^{2}}$. Using equation $A{{Q}^{2}}=B{{Q}^{2}}$, the following result is obtained:
	\[\begin{aligned}
		& \left( 1-\beta _{A}^{Q} \right)\left( \beta _{B}^{Q}A{{B}^{2}}+\beta _{C}^{Q}A{{C}^{2}}+\beta _{D}^{Q}A{{D}^{2}} \right) \\ 
		& -\beta _{B}^{Q}\beta _{C}^{Q}B{{C}^{2}}-\beta _{B}^{Q}\beta _{D}^{Q}B{{D}^{2}}-\beta _{C}^{Q}\beta _{D}^{Q}C{{D}^{2}} \\ 
		& =\left( 1-\beta _{B}^{Q} \right)\left( \beta _{C}^{Q}B{{C}^{2}}+\beta _{D}^{Q}B{{D}^{2}}+\beta _{A}^{Q}B{{A}^{2}} \right) \\ 
		& -\beta _{C}^{Q}\beta _{D}^{Q}C{{D}^{2}}-\beta _{C}^{Q}\beta _{A}^{Q}C{{A}^{2}}-\beta _{D}^{Q}\beta _{A}^{Q}D{{A}^{2}},  
	\end{aligned}\]
	i.e. 
	\[\beta _{B}^{Q}A{{B}^{2}}+\beta _{C}^{Q}A{{C}^{2}}+\beta _{D}^{Q}A{{D}^{2}}=\beta _{C}^{Q}B{{C}^{2}}+\beta _{D}^{Q}B{{D}^{2}}+\beta _{A}^{Q}B{{A}^{2}},\]
	\[\beta _{A}^{Q}B{{A}^{2}}-\beta _{B}^{Q}A{{B}^{2}}+\beta _{C}^{Q}\left( B{{C}^{2}}-A{{C}^{2}} \right)+\beta _{D}^{Q}\left( B{{D}^{2}}-A{{D}^{2}} \right)=0.\]
	
	Using equation $B{{Q}^{2}}=C{{Q}^{2}}$ (or rotate the above formula), the following result is obtained:
	\[\beta _{B}^{Q}C{{B}^{2}}-\beta _{C}^{Q}B{{C}^{2}}+\beta _{D}^{Q}\left( C{{D}^{2}}-B{{D}^{2}} \right)+\beta _{A}^{Q}\left( C{{A}^{2}}-B{{A}^{2}} \right)=0.\]
	
	Using equation $C{{Q}^{2}}=D{{Q}^{2}}$ (or rotate the above formula), the following result is obtained:
	\[\beta _{C}^{Q}D{{C}^{2}}-\beta _{D}^{Q}C{{D}^{2}}+\beta _{A}^{Q}\left( D{{A}^{2}}-C{{A}^{2}} \right)+\beta _{B}^{Q}\left( D{{B}^{2}}-C{{B}^{2}} \right)=0.\]
	
	Therefore, using the condition $\beta _{A}^{Q}+\beta _{B}^{Q}+\beta _{C}^{Q}+\beta _{D}^{Q}=1$, the following linear equations are obtained:	
	\[\left\{ \begin{aligned}
		& \beta _{A}^{Q}+\beta _{B}^{Q}+\beta _{C}^{Q}+\beta _{D}^{Q}=1 \\ 
		& B{{A}^{2}}\beta _{A}^{Q}-A{{B}^{2}}\beta _{B}^{Q}+\left( B{{C}^{2}}-A{{C}^{2}} \right)\beta _{C}^{Q}+\left( B{{D}^{2}}-A{{D}^{2}} \right)\beta _{D}^{Q}=0 \\ 
		& C{{B}^{2}}\beta _{B}^{Q}-B{{C}^{2}}\beta _{C}^{Q}+\left( C{{D}^{2}}-B{{D}^{2}} \right)\beta _{D}^{Q}+\left( C{{A}^{2}}-B{{A}^{2}} \right)\beta _{A}^{Q}=0 \\ 
		& D{{C}^{2}}\beta _{C}^{Q}-C{{D}^{2}}\beta _{D}^{Q}+\left( D{{A}^{2}}-C{{A}^{2}} \right)\beta _{A}^{Q}+\left( D{{B}^{2}}-C{{B}^{2}} \right)\beta _{B}^{Q}=0, \\ 
	\end{aligned} \right.\]
	i.e.
	\[\left\{ \begin{aligned}
		& \beta _{A}^{Q}+\beta _{B}^{Q}+\beta _{C}^{Q}+\beta _{D}^{Q}=1 \\ 
		& B{{A}^{2}}\beta _{A}^{Q}-A{{B}^{2}}\beta _{B}^{Q}+\left( B{{C}^{2}}-A{{C}^{2}} \right)\beta _{C}^{Q}+\left( B{{D}^{2}}-A{{D}^{2}} \right)\beta _{D}^{Q}=0 \\ 
		& \left( C{{A}^{2}}-B{{A}^{2}} \right)\beta _{A}^{Q}+C{{B}^{2}}\beta _{B}^{Q}-B{{C}^{2}}\beta _{C}^{Q}+\left( C{{D}^{2}}-B{{D}^{2}} \right)\beta _{D}^{Q}=0 \\ 
		& \left( D{{A}^{2}}-C{{A}^{2}} \right)\beta _{A}^{Q}+\left( D{{B}^{2}}-C{{B}^{2}} \right)\beta _{B}^{Q}+D{{C}^{2}}\beta _{C}^{Q}-C{{D}^{2}}\beta _{D}^{Q}=0. \\ 
	\end{aligned} \right.\]
	
	Solve the above equations to obtain:
	\[\beta _{A}^{Q}=\frac{{{U}_{A}}}{U}=\frac{{{U}_{A}}}{{{U}_{A}}+{{U}_{B}}+{{U}_{C}}+{{U}_{D}}},\]
	\[\beta _{B}^{Q}=\frac{{{U}_{B}}}{U}=\frac{{{U}_{B}}}{{{U}_{A}}+{{U}_{B}}+{{U}_{C}}+{{U}_{D}}},\]
	\[\beta _{C}^{Q}=\frac{{{U}_{C}}}{U}=\frac{{{U}_{C}}}{{{U}_{A}}+{{U}_{B}}+{{U}_{C}}+{{U}_{D}}},\]
	\[\beta _{D}^{Q}=\frac{{{U}_{D}}}{U}=\frac{{{U}_{D}}}{{{U}_{A}}+{{U}_{B}}+{{U}_{C}}+{{U}_{D}}}.\]
	Where 
	\[\begin{aligned}
		{{U}_{A}}& =\left( \Delta _{2}^{A}-B{{C}^{2}} \right)B{{C}^{2}}A{{D}^{2}}+\left( \Delta _{2}^{A}-C{{D}^{2}} \right)C{{D}^{2}}A{{B}^{2}} \\ 
		& +\left( \Delta _{2}^{A}-D{{B}^{2}} \right)D{{B}^{2}}A{{C}^{2}}-B{{C}^{2}}C{{D}^{2}}D{{B}^{2}},  
	\end{aligned}\]
	\[\Delta _{2}^{A}=\frac{1}{2}\left( B{{C}^{2}}+C{{D}^{2}}+D{{B}^{2}} \right);\]
	\[\begin{aligned}
		{{U}_{B}}& =\left( \Delta _{2}^{B}-C{{D}^{2}} \right)C{{D}^{2}}B{{A}^{2}}+\left( \Delta _{2}^{B}-D{{A}^{2}} \right)D{{A}^{2}}B{{C}^{2}} \\ 
		& +\left( \Delta _{2}^{B}-A{{C}^{2}} \right)A{{C}^{2}}B{{D}^{2}}-C{{D}^{2}}D{{A}^{2}}A{{C}^{2}},  
	\end{aligned}\]
	\[\Delta _{2}^{B}=\frac{1}{2}\left( C{{D}^{2}}+D{{A}^{2}}+A{{C}^{2}} \right);\]
	\[\begin{aligned}
		{{U}_{C}}& =\left( \Delta _{2}^{C}-D{{A}^{2}} \right)D{{A}^{2}}C{{B}^{2}}+\left( \Delta _{2}^{C}-A{{B}^{2}} \right)A{{B}^{2}}C{{D}^{2}} \\ 
		& +\left( \Delta _{2}^{C}-B{{D}^{2}} \right)B{{D}^{2}}C{{A}^{2}}-D{{A}^{2}}A{{B}^{2}}B{{D}^{2}},  
	\end{aligned}\]
	\[\Delta _{2}^{C}=\frac{1}{2}\left( D{{A}^{2}}+A{{B}^{2}}+B{{D}^{2}} \right);\]
	\[\begin{aligned}
		{{U}_{D}}& =\left( \Delta _{2}^{D}-A{{B}^{2}} \right)A{{B}^{2}}D{{C}^{2}}+\left( \Delta _{2}^{D}-B{{C}^{2}} \right)B{{C}^{2}}D{{A}^{2}} \\ 
		& +\left( \Delta _{2}^{D}-C{{A}^{2}} \right)C{{A}^{2}}D{{B}^{2}}-A{{B}^{2}}B{{C}^{2}}C{{A}^{2}},  
	\end{aligned}\]
	\[\Delta _{2}^{D}=\frac{1}{2}\left( A{{B}^{2}}+B{{C}^{2}}+C{{A}^{2}} \right).\]
	
	And 
	\[U={{U}_{A}}+{{U}_{B}}+{{U}_{C}}+{{U}_{D}}=4\left( {{t}_{1}}-{{t}_{2}}-{{t}_{3}} \right),\]
	where 
	\[{{t}_{1}}={{q}_{2}}{{\Delta }_{2}},\]
	\[{{q}_{2}}=\frac{1}{2}\left( A{{B}^{2}}C{{D}^{2}}+B{{C}^{2}}A{{D}^{2}}+C{{A}^{2}}B{{D}^{2}} \right),\]
	\[{{\Delta }_{2}}=\frac{1}{2}\left( A{{B}^{2}}+A{{C}^{2}}+A{{D}^{2}}+B{{C}^{2}}+C{{D}^{2}}+D{{B}^{2}} \right),\]
	\[{{t}_{2}}=\frac{1}{2}\left( A{{B}^{2}}C{{D}^{2}}\left( A{{B}^{2}}+C{{D}^{2}} \right)+B{{C}^{2}}A{{D}^{2}}\left( B{{C}^{2}}+A{{D}^{2}} \right)+C{{A}^{2}}B{{D}^{2}}\left( C{{A}^{2}}+B{{D}^{2}} \right) \right),\]
	\[{{t}_{3}}=\frac{1}{4}\left( A{{B}^{2}}B{{C}^{2}}C{{A}^{2}}+B{{C}^{2}}C{{D}^{2}}D{{B}^{2}}+C{{D}^{2}}D{{A}^{2}}A{{C}^{2}}+D{{A}^{2}}A{{B}^{2}}B{{D}^{2}} \right).\]

\end{proof}
\hfill $\square$\par

%
%
%
%
%
%
%
%
%
%
%
%
%
%
%
%
%
%
%
%
%
%
%
%

\chapter{Distance between special intersecting centers of a tetrahedron}\label{Ch26}
\thispagestyle{empty}

\section{Distance between origin and intersecting center of a tetrahedron on frame of circumcenter}\label{Sec26.1}

\subsection{Distance between circumcenter and centroid (DOIC-T)}\label{Subsec26.1.1}
\begin{theorem}{Distance between circumcenter and centroid, Daiyuan Zhang}{Thm26.1.1}\label{Thm26.1.1} 
	Given a tetrahedron $ABCD$, let points $Q$, $G$ be the circumcenter and centroid of the tetrahedron $ABCD$ respectively, and $R$ be the circumscribed radius of tetrahedron $ABCD$, then	
	\[QG=\frac{1}{4}\sqrt{16{{R}^{2}}-\left( A{{B}^{2}}+A{{C}^{2}}+A{{D}^{2}}+B{{C}^{2}}+B{{D}^{2}}+C{{D}^{2}} \right)}.\]			
	Where $AB$, $AC$, $AD$, $BC$, $BD$, $CD$ are the lengths of the six edges of the tetrahedron respectively.
\end{theorem}

\begin{proof}
	According to theorem \ref{thm:Thm24.1.4} and theorem \ref{thm:Thm21.1.1}, it is obtained that:	
	\[Q{{G}^{2}}={{R}^{2}}-\left( \begin{aligned}
		& \beta _{A}^{G}\beta _{B}^{G}A{{B}^{2}}+\beta _{A}^{G}\beta _{C}^{G}A{{C}^{2}}+\beta _{A}^{G}\beta _{D}^{G}A{{D}^{2}} \\ 
		& +\beta _{B}^{G}\beta _{C}^{G}B{{C}^{2}}+\beta _{B}^{G}\beta _{D}^{G}B{{D}^{2}}+\beta _{C}^{G}\beta _{D}^{G}C{{D}^{2}}  
	\end{aligned} \right),\]
	i.e.	
	\[\begin{aligned}
		Q{{G}^{2}}&={{R}^{2}}-\left( \begin{aligned}
			& \frac{1}{4}\times \frac{1}{4}A{{B}^{2}}+\frac{1}{4}\times \frac{1}{4}A{{C}^{2}}+\frac{1}{4}\times \frac{1}{4}A{{D}^{2}} \\ 
			& +\frac{1}{4}\times \frac{1}{4}B{{C}^{2}}+\frac{1}{4}\times \frac{1}{4}B{{D}^{2}}+\frac{1}{4}\times \frac{1}{4}C{{D}^{2}}  
		\end{aligned} \right) \\ 
		& ={{R}^{2}}-\frac{1}{16}\left( A{{B}^{2}}+A{{C}^{2}}+A{{D}^{2}}+B{{C}^{2}}+B{{D}^{2}}+C{{D}^{2}} \right).  
	\end{aligned}\]
	
	Therefore,	
	\[QG=\frac{1}{4}\sqrt{16{{R}^{2}}-\left( A{{B}^{2}}+A{{C}^{2}}+A{{D}^{2}}+B{{C}^{2}}+B{{D}^{2}}+C{{D}^{2}} \right)}.\]
\end{proof} \hfill $\square$\par

For a regular tetrahedron, $QG=0$, assume that each edge length of the tetrahedron is $a$, the following result is obtained from the above formula:
$16{{R}^{2}}-6{{a}^{2}}=0$

i.e.
\[R=\frac{\sqrt{6}}{4}a.\]	

The above formula is a well-known result.

\subsection{Distance between circumcenter and incenter (DOIC-T)}\label{Subsec26.1.2}
\begin{theorem}{Distance between circumcenter and incenter, Daiyuan Zhang}{Thm26.1.2}\label{Thm26.1.2} 
	Given a tetrahedron $ABCD$, let points $Q$, $I$ be the circumcenter and incenter of the tetrahedron $ABCD$ respectively, and $R$ be the circumscribed radius of tetrahedron $ABCD$, then	
	\[Q{{I}^{2}}={{R}^{2}}-\frac{1}{{{S}^{2}}}\left( \begin{aligned}
		& {{S}^{A}}{{S}^{B}}A{{B}^{2}}+{{S}^{A}}{{S}^{C}}A{{C}^{2}}+{{S}^{A}}{{S}^{C}}A{{D}^{2}} \\ 
		& +{{S}^{B}}{{S}^{C}}B{{C}^{2}}+{{S}^{B}}{{S}^{D}}B{{D}^{2}}+{{S}^{C}}{{S}^{D}}C{{D}^{2}}  
	\end{aligned} \right).\]
	
	Where ${{S}^{A}}$, ${{S}^{B}}$, ${{S}^{C}}$, ${{S}^{D}}$ are the areas of the triangle opposite to the four vertices $A$, $B$, $C$, $D$ of the tetrahedron $ABCD$ respectively; $S$ is the surface area of tetrahedron $ABCD$;	and $AB$, $AC$, $AD$, $BC$, $BD$, $CD$ are the lengths of the six edges of the tetrahedron respectively.
\end{theorem}

\begin{proof}
	According to theorem \ref{thm:Thm24.1.4} and theorem \ref{thm:Thm21.2.1}, it is obtained that:	
	\[Q{{I}^{2}}={{R}^{2}}-\left( \begin{aligned}
		& \beta _{A}^{I}\beta _{B}^{I}A{{B}^{2}}+\beta _{A}^{I}\beta _{C}^{I}A{{C}^{2}}+\beta _{A}^{I}\beta _{D}^{I}A{{D}^{2}} \\ 
		& +\beta _{B}^{I}\beta _{C}^{I}B{{C}^{2}}+\beta _{B}^{I}\beta _{D}^{I}B{{D}^{2}}+\beta _{C}^{I}\beta _{D}^{I}C{{D}^{2}}  
	\end{aligned} \right),\]
	\[Q{{I}^{2}}={{R}^{2}}-\left( \begin{aligned}
		& \frac{{{S}^{A}}}{S}\frac{{{S}^{B}}}{S}A{{B}^{2}}+\frac{{{S}^{A}}}{S}\frac{{{S}^{C}}}{S}A{{C}^{2}}+\frac{{{S}^{A}}}{S}\frac{{{S}^{D}}}{S}A{{D}^{2}} \\ 
		& +\frac{{{S}^{B}}}{S}\frac{{{S}^{C}}}{S}B{{C}^{2}}+\frac{{{S}^{B}}}{S}\frac{{{S}^{D}}}{S}B{{D}^{2}}+\frac{{{S}^{C}}}{S}\frac{{{S}^{D}}}{S}C{{D}^{2}}  
	\end{aligned} \right),\]
	\[Q{{I}^{2}}={{R}^{2}}-\frac{1}{{{S}^{2}}}\left( \begin{aligned}
		& {{S}^{A}}{{S}^{B}}A{{B}^{2}}+{{S}^{A}}{{S}^{C}}A{{C}^{2}}+{{S}^{A}}{{S}^{C}}A{{D}^{2}} \\ 
		& +{{S}^{B}}{{S}^{C}}B{{C}^{2}}+{{S}^{B}}{{S}^{D}}B{{D}^{2}}+{{S}^{C}}{{S}^{D}}C{{D}^{2}}  
	\end{aligned} \right).\]
\end{proof}
\hfill $\square$\par

\subsection{Distance between circumcenter and excenters (DOIC-T)}\label{Subsec26.1.2}

\begin{theorem}{Distance between circumcenter and excenters, Daiyuan Zhang}{Thm26.1.3}\label{Thm26.1.3}
Given a tetrahedron $ABCD$, let point $Q$ be the circumcenter, let ${{E}_{A}}$, ${{E}_{B}}$, ${{E}_{C}}$, ${{E}_{D}}$ be the excenters in trihedral angles of vertexs $A$, $B$, $C$, $D$ of tetrahedron $ABCD$ respectively, and $R$ be the circumscribed radius of tetrahedron $ABCD$, then

\[Q{{E}_{A}}^{2}={{R}^{2}}-\frac{1}{{{\left( S-2{{S}^{A}} \right)}^{2}}}\left( \begin{aligned}
	& -{{S}^{A}}{{S}^{B}}A{{B}^{2}}-{{S}^{A}}{{S}^{C}}A{{C}^{2}}-{{S}^{A}}{{S}^{D}}A{{D}^{2}} \\ 
	& +{{S}^{B}}{{S}^{C}}B{{C}^{2}}+{{S}^{B}}{{S}^{D}}B{{D}^{2}}+{{S}^{C}}{{S}^{D}}C{{D}^{2}}  
\end{aligned} \right),\]
\[Q{{E}_{B}}^{2}={{R}^{2}}-\frac{1}{{{\left( S-2{{S}^{B}} \right)}^{2}}}\left( \begin{aligned}
	& -{{S}^{A}}{{S}^{B}}A{{B}^{2}}+{{S}^{A}}{{S}^{C}}A{{C}^{2}}+{{S}^{A}}{{S}^{D}}A{{D}^{2}} \\ 
	& -{{S}^{B}}{{S}^{C}}B{{C}^{2}}-{{S}^{B}}{{S}^{D}}B{{D}^{2}}+{{S}^{C}}{{S}^{D}}C{{D}^{2}}  
\end{aligned} \right),\]
\[Q{{E}_{C}}^{2}={{R}^{2}}-\frac{1}{{{\left( S-2{{S}^{C}} \right)}^{2}}}\left( \begin{aligned}
	& {{S}^{A}}{{S}^{B}}A{{B}^{2}}-{{S}^{A}}{{S}^{C}}A{{C}^{2}}+{{S}^{A}}{{S}^{D}}A{{D}^{2}} \\ 
	& -{{S}^{B}}{{S}^{C}}B{{C}^{2}}+{{S}^{B}}{{S}^{D}}B{{D}^{2}}-{{S}^{C}}{{S}^{D}}C{{D}^{2}}  
\end{aligned} \right),\]
\[Q{{E}_{D}}^{2}={{R}^{2}}-\frac{1}{{{\left( S-2{{S}^{D}} \right)}^{2}}}\left( \begin{aligned}
	& {{S}^{A}}{{S}^{B}}A{{B}^{2}}+{{S}^{A}}{{S}^{C}}A{{C}^{2}}-{{S}^{A}}{{S}^{D}}A{{D}^{2}} \\ 
	& +{{S}^{B}}{{S}^{C}}B{{C}^{2}}-{{S}^{B}}{{S}^{D}}B{{D}^{2}}-{{S}^{C}}{{S}^{D}}C{{D}^{2}}  
\end{aligned} \right).\]
Where ${{S}^{A}}$, ${{S}^{B}}$, ${{S}^{C}}$, ${{S}^{D}}$ are the areas of the triangle opposite to the four vertices $A$, $B$, $C$, $D$ of the tetrahedron $ABCD$ respectively; $S$ is the surface area of tetrahedron $ABCD$;	and $AB$, $AC$, $AD$, $BC$, $BD$, $CD$ are the lengths of the six edges of the tetrahedron respectively.
\end{theorem}

\begin{proof}
	According to theorem \ref{thm:Thm24.1.4} and theorem \ref{thm:Thm21.3.1}, it is obtained that:	
	\[Q{{E}_{A}}^{2}={{R}^{2}}-\left( \begin{aligned}
		& \beta _{A}^{{{E}_{A}}}\beta _{B}^{{{E}_{A}}}A{{B}^{2}}+\beta _{A}^{{{E}_{A}}}\beta _{C}^{{{E}_{A}}}A{{C}^{2}}+\beta _{A}^{{{E}_{A}}}\beta _{D}^{{{E}_{A}}}A{{D}^{2}} \\ 
		& +\beta _{B}^{{{E}_{A}}}\beta _{C}^{{{E}_{A}}}B{{C}^{2}}+\beta _{B}^{{{E}_{A}}}\beta _{D}^{{{E}_{A}}}B{{D}^{2}}+\beta _{C}^{{{E}_{A}}}\beta _{D}^{{{E}_{A}}}C{{D}^{2}}  
	\end{aligned} \right),\]
	i.e.
	\[Q{{E}_{A}}^{2}={{R}^{2}}-\left( \begin{aligned}
		& \left( -\frac{{{S}^{A}}}{S-2{{S}^{A}}} \right)\left( \frac{{{S}^{B}}}{S-2{{S}^{A}}} \right)A{{B}^{2}}+\left( -\frac{{{S}^{A}}}{S-2{{S}^{A}}} \right)\left( \frac{{{S}^{C}}}{S-2{{S}^{A}}} \right)A{{C}^{2}} \\ 
		& +\left( -\frac{{{S}^{A}}}{S-2{{S}^{A}}} \right)\left( \frac{{{S}^{D}}}{S-2{{S}^{A}}} \right)A{{D}^{2}}+\left( \frac{{{S}^{B}}}{S-2{{S}^{A}}} \right)\left( \frac{{{S}^{C}}}{S-2{{S}^{A}}} \right)B{{C}^{2}} \\ 
		& +\left( \frac{{{S}^{B}}}{S-2{{S}^{A}}} \right)\left( \frac{{{S}^{D}}}{S-2{{S}^{A}}} \right)B{{D}^{2}}+\left( \frac{{{S}^{C}}}{S-2{{S}^{A}}} \right)\left( \frac{{{S}^{D}}}{S-2{{S}^{A}}} \right)C{{D}^{2}}  
	\end{aligned} \right),\]
	i.e.
	\[Q{{E}_{A}}^{2}={{R}^{2}}-\frac{1}{{{\left( S-2{{S}^{A}} \right)}^{2}}}\left( \begin{aligned}
		& -{{S}^{A}}{{S}^{B}}A{{B}^{2}}-{{S}^{A}}{{S}^{C}}A{{C}^{2}}-{{S}^{A}}{{S}^{D}}A{{D}^{2}} \\ 
		& +{{S}^{B}}{{S}^{C}}B{{C}^{2}}+{{S}^{B}}{{S}^{D}}B{{D}^{2}}+{{S}^{C}}{{S}^{D}}C{{D}^{2}}  
	\end{aligned} \right).\]
	
	Similarly:
	\[Q{{E}_{B}}^{2}={{R}^{2}}-\frac{1}{{{\left( S-2{{S}^{B}} \right)}^{2}}}\left( \begin{aligned}
		& -{{S}^{A}}{{S}^{B}}A{{B}^{2}}+{{S}^{A}}{{S}^{C}}A{{C}^{2}}+{{S}^{A}}{{S}^{D}}A{{D}^{2}} \\ 
		& -{{S}^{B}}{{S}^{C}}B{{C}^{2}}-{{S}^{B}}{{S}^{D}}B{{D}^{2}}+{{S}^{C}}{{S}^{D}}C{{D}^{2}}  
	\end{aligned} \right),\]
	\[Q{{E}_{C}}^{2}={{R}^{2}}-\frac{1}{{{\left( S-2{{S}^{C}} \right)}^{2}}}\left( \begin{aligned}
		& {{S}^{A}}{{S}^{B}}A{{B}^{2}}-{{S}^{A}}{{S}^{C}}A{{C}^{2}}+{{S}^{A}}{{S}^{D}}A{{D}^{2}} \\ 
		& -{{S}^{B}}{{S}^{C}}B{{C}^{2}}+{{S}^{B}}{{S}^{D}}B{{D}^{2}}-{{S}^{C}}{{S}^{D}}C{{D}^{2}}  
	\end{aligned} \right),\]
	\[Q{{E}_{D}}^{2}={{R}^{2}}-\frac{1}{{{\left( S-2{{S}^{D}} \right)}^{2}}}\left( \begin{aligned}
		& {{S}^{A}}{{S}^{B}}A{{B}^{2}}+{{S}^{A}}{{S}^{C}}A{{C}^{2}}-{{S}^{A}}{{S}^{D}}A{{D}^{2}} \\ 
		& +{{S}^{B}}{{S}^{C}}B{{C}^{2}}-{{S}^{B}}{{S}^{D}}B{{D}^{2}}-{{S}^{C}}{{S}^{D}}C{{D}^{2}}  
	\end{aligned} \right).\]	
\end{proof}
\hfill $\square$\par

\section{Distance between two intersecting centers of a tetrahedronon on tetrahedral frame}\label{Sec26.2}

\subsection{Distance between centroid and incenter of a tetrahedron (DTICs-T)}\label{Subsec26.2.1}
Let $S$ be the surface area of tetrahedron $ABCD$, i.e
\[S={{S}^{A}}+{{S}^{B}}+{{S}^{C}}+{{S}^{D}}.\]

Let $\bar{S}$ be the average area of four triangles on the surface of tetrahedron $ABCD$, i.e.
\[\bar{S}=\frac{1}{4}\left( {{S}^{A}}+{{S}^{B}}+{{S}^{C}}+{{S}^{D}} \right).\]

The following theorem can be directly derived from theorem \ref{thm:Thm24.2.1}. 

\begin{theorem}{Distance between centroid and incenter, Daiyuan Zhang}{Thm26.2.1}\label{Thm26.2.1} 
	Given a tetrahedron $ABCD$, let points $G$, $I$ be the centroid and incenter of the tetrahedron $ABCD$ respectively, then		
	\[G{{I}^{2}}=-\frac{1}{{{S}^{2}}}\left( \begin{aligned}
		& \left( {{S}^{A}}-\bar{S} \right)\left( {{S}^{B}}-\bar{S} \right)A{{B}^{2}}+\left( {{S}^{A}}-\bar{S} \right)\left( {{S}^{C}}-\bar{S} \right)A{{C}^{2}} \\ 
		& +\left( {{S}^{A}}-\bar{S} \right)\left( {{S}^{D}}-\bar{S} \right)A{{D}^{2}}+\left( {{S}^{B}}-\bar{S} \right)\left( {{S}^{C}}-\bar{S} \right)B{{C}^{2}} \\ 
		& +\left( {{S}^{B}}-\bar{S} \right)\left( {{S}^{D}}-\bar{S} \right)B{{D}^{2}}+\left( {{S}^{C}}-\bar{S} \right)\left( {{S}^{D}}-\bar{S} \right)C{{D}^{2}}  
	\end{aligned} \right).\]
		
	Where ${{S}^{A}}$, ${{S}^{B}}$, ${{S}^{C}}$, ${{S}^{D}}$ are the areas of the triangle opposite to the four vertices $A$, $B$, $C$, $D$ of the tetrahedron $ABCD$ respectively; $S$ is the surface area of tetrahedron $ABCD$;	$\bar{S}$ is the average area of four triangles on the surface of tetrahedron $ABCD$; and $AB$, $AC$, $AD$, $BC$, $BD$, $CD$ are the lengths of the six edges of the tetrahedron respectively.
\end{theorem}

\begin{proof}
	According to theorem \ref{thm:Thm21.2.1} and theorem \ref{thm:Thm21.1.1}, it is obtained that:	
	\[\beta _{A}^{GI}=\beta _{A}^{I}-\beta _{A}^{G}=\frac{{{S}^{A}}}{S}-\frac{1}{4}=\frac{4{{S}^{A}}-S}{4S},\]	
	\[\beta _{B}^{GI}=\beta _{B}^{I}-\beta _{B}^{G}=\frac{{{S}^{B}}}{S}-\frac{1}{4}=\frac{4{{S}^{B}}-S}{4S},\]	
	\[\beta _{C}^{GI}=\beta _{C}^{I}-\beta _{C}^{G}=\frac{{{S}^{C}}}{S}-\frac{1}{4}=\frac{4{{S}^{C}}-S}{4S},\]	
	\[\beta _{D}^{GI}=\beta _{D}^{I}-\beta _{D}^{G}=\frac{{{S}^{D}}}{S}-\frac{1}{4}=\frac{4{{S}^{D}}-S}{4S}.\]	
	
	Therefore, according to theorem \ref{thm:Thm24.2.1}, let ${{P}_{1}}$ be the centroid, and ${{P}_{2}}$ be the incenter, the following results can be obtained:
	\[\begin{aligned}
		G{{I}^{2}}=& -\beta _{A}^{GI}\beta _{B}^{GI}A{{B}^{2}}-\beta _{A}^{GI}\beta _{C}^{GI}A{{C}^{2}}-\beta _{A}^{GI}\beta _{D}^{GI}A{{D}^{2}} \\ 
		& -\beta _{B}^{GI}\beta _{C}^{GI}B{{C}^{2}}-\beta _{B}^{GI}\beta _{D}^{GI}B{{D}^{2}}-\beta _{C}^{GI}\beta _{D}^{GI}C{{D}^{2}},  
	\end{aligned}\]
	i.e.	
	\[G{{I}^{2}}=-\frac{1}{16{{S}^{2}}}\left( \begin{aligned}
		& \left( 4{{S}^{A}}-S \right)\left( 4{{S}^{B}}-S \right)A{{B}^{2}}+\left( 4{{S}^{A}}-S \right)\left( 4{{S}^{C}}-S \right)A{{C}^{2}} \\ 
		& +\left( 4{{S}^{A}}-S \right)\left( 4{{S}^{D}}-S \right)A{{D}^{2}}+\left( 4{{S}^{B}}-S \right)\left( 4{{S}^{C}}-S \right)B{{C}^{2}} \\ 
		& +\left( 4{{S}^{B}}-S \right)\left( 4{{S}^{D}}-S \right)B{{D}^{2}}+\left( 4{{S}^{C}}-S \right)\left( 4{{S}^{D}}-S \right)C{{D}^{2}}  
	\end{aligned} \right),\]
	i.e.
	\[G{{I}^{2}}=-\frac{1}{{{S}^{2}}}\left( \begin{aligned}
		& \left( {{S}^{A}}-\bar{S} \right)\left( {{S}^{B}}-\bar{S} \right)A{{B}^{2}}+\left( {{S}^{A}}-\bar{S} \right)\left( {{S}^{C}}-\bar{S} \right)A{{C}^{2}} \\ 
		& +\left( {{S}^{A}}-\bar{S} \right)\left( {{S}^{D}}-\bar{S} \right)A{{D}^{2}}+\left( {{S}^{B}}-\bar{S} \right)\left( {{S}^{C}}-\bar{S} \right)B{{C}^{2}} \\ 
		& +\left( {{S}^{B}}-\bar{S} \right)\left( {{S}^{D}}-\bar{S} \right)B{{D}^{2}}+\left( {{S}^{C}}-\bar{S} \right)\left( {{S}^{D}}-\bar{S} \right)C{{D}^{2}}  
	\end{aligned} \right).\]
\end{proof}
\hfill $\square$\par

The above theorem shows that the distance between the centroid and the incenter of a tetrahedron is related to the area of the four faces (triangles) and the lengths of six edges of the tetrahedron. If the lengths of six edges of a tetrahedron are given, the area of each face (triangle) can be calculated according to Helen's formula, and the distance between the centroid and the incenter can be obtained according to the above theorem. From the application point of view, the lengths of six edges of a tetrahedron are the most “natural” parameters, they are easy to be measured in practical applications. It can be seen that the above theorem not only have theoretical value, but also have application value.

It can be seen that the above formula is symmetrical and graceful. However, it is difficult to get the above formula by using analytic geometry or Euclidean geometry.


Suppose you ask the following question: solve the distance between the centroid and the incenter of a tetrahedron. If you don't know the formula I proposed here for calculating $G{{I}^{2}}$, what would you do to solve this problem?


If the Euclidean geometry method is used, the calculation will be difficult because Euclidean geometry is disjointed from algebraic methods and calculating the distance between two points is the natural soft rib of Euclidean geometry. Euclidean geometry has been developed for more than 2,000 years and no formula has been found for the calculation of $G{{I}^{2}}$. Does this imply that it is impossible to find a formula for the calculation of $G{{I}^{2}}$ using Euclidean geometry?


If you use analytic geometry (coordinate geometry), you must first set up a coordinate system, then give the coordinates of the four vertices of the tetrahedron, and take the coordinates of the four vertices as parameters to represent the centroid and the incenter. Finally, you calculate the distance between the two points using the formula in analytic geometry. Obviously, this will result in a complex system of nonlinear equations. How to find the analytical solution of this nonlinear system of equations? Analytic geometry has also been developed for centuries and no formula has been found for the calculation of $G{{I}^{2}}$. Does this also imply that no formula for the calculation of $G{{I}^{2}}$ can be found using analytic geometry? Even if an analytic solution can be found by analytic geometry, this result is expressed by the coordinates of the four vertices of the tetrahedron, so we can't see the “natural parameters” such as the length of edge or the area. You may say, I can find a way to make a transformation for seeing the “natural parameters”. Not to mention the complexity of this transformation, how would you transform it if you didn't know the formula I proposed in the above theorem? What is the transformation result do you want? Do you feel at a loss? It's like a boat floating on a vast sea, with no goal at all. But on the other hand, according to the formula I put forward to calculate $G{{I}^{2}}$, it is very convenient to get the result of the representation in analytic geometry.


In addition, the coordinates of points are the most microscopic (bottom) parameters, which are far away from people's visual senses and give people an illusory feeling. When you get a tetrahedron, you can see the “macro” quantities of the four triangular faces and six edges of the tetrahedron, rather than the coordinates of the four vertices. These “macro” quantities are the parameters that people's senses can directly touch.


Obviously, the coordinates of the points depend on the selection of the coordinate system and also depend on the selection of the coordinate origin. However, “natural parameters” (such as the length of edge, the area, etc.) are not related to the selection of coordinate system. This means that if the parameters of a geometric quantity are all “natural parameters”, then the geometric quantity will be able to more profoundly reveal the internal laws between the geometric quantity and “natural parameters”, and Intercenter Geometry has done this, but analytic geometry has not.

The above theorem I put forward once again explains the idea of Space Intercenter Geometry: the idea of Space Intercenter Geometry is that the geometric quantities in space
are expressed by the lengths of the six edges of a given tetrahedron.


The above formula I put forward solves the problem that cannot be solved by traditional methods. I am proud of it. The above theorem I put forward is only the tip of the iceberg in this book's fruitful and innovative achievements. If you can finish reading this book patiently, I believe you will have more surprises. I hope you like Intercenter Geometry.

\subsection{Distance between the centroid and the circumcenter of a tetrahedron (DTICs-T)}\label{Subsec26.2.2}

In subsection \ref{Subsec26.1.1}, the distance between the centroid and the circumcenter of the tetrahedron has been studied, but the circumcenter was selected as the origin of the frame when deriving the distance formula, and the given formula (see theorem \ref{thm:Thm26.1.1}) is the DOIC-T. The distance formula between the centroid and the circumcenter proposed in theorem \ref{thm:Thm26.1.1} has an important feature: the formula includes the circumscribed sphere radius. The distance formula between the centroid and the circumcenter given in this subsection is the DTICs-T, and the circumscribed sphere radius is not included in the formula. Obviously, the centroid and the circumcenter are both the ICs-T, and the frame components of the two ICs-T have been calculated, so the distance between the centroid and the circumcenter of the tetrahedron can be calculated directly.


\begin{theorem}{Distance between centroid and circumcenter, Daiyuan Zhang) }{Thm26.2.2}\label{Thm26.2.2}
Suppose that given a tetrahedron $ABCD$, point $Q$ is its circumcenter and point $G$ is its centroid, then
\[G{{Q}^{2}}=-\frac{1}{16{{U}^{2}}}\left( \begin{aligned}
	& \left( 4{{U}_{A}}-U \right)\left( 4{{U}_{B}}-U \right)A{{B}^{2}}+\left( 4{{U}_{A}}-U \right)\left( 4{{U}_{C}}-U \right)A{{C}^{2}} \\ 
	& +\left( 4{{U}_{A}}-U \right)\left( 4{{U}_{D}}-U \right)A{{D}^{2}}+\left( 4{{U}_{B}}-U \right)\left( 4{{U}_{C}}-U \right)B{{C}^{2}} \\ 
	& +\left( 4{{U}_{B}}-U \right)\left( 4{{U}_{D}}-U \right)B{{D}^{2}}+\left( 4{{U}_{C}}-U \right)\left( 4{{U}_{D}}-U \right)C{{D}^{2}}  
\end{aligned} \right).\]
Where $U$, ${{U}_{A}}$, ${{U}_{B}}$, ${{U}_{C}}$, ${{U}_{D}}$ are calculated by theorem \ref{thm:Thm25.3.1}, theorem \ref{thm:Thm25.1.2} or theorem \ref{thm:Thm25.1.1}, and $AB$, $AC$, $AD$, $BC$, $BD$, $CD$ are the lengths of the six edges of the tetrahedron respectively. 	
\end{theorem}

\begin{proof}
	From theorem \ref{thm:Thm21.1.1} and theorem \ref{thm:Thm25.3.1}, theorem \ref{thm:Thm25.1.2} or \ref{thm:Thm25.1.1}, we have
	\[\beta _{A}^{GQ}=\beta _{A}^{Q}-\beta _{A}^{G}=\frac{{{U}_{A}}}{U}-\frac{1}{4}=\frac{4{{U}_{A}}-U}{4U},\]	
	\[\beta _{B}^{GQ}=\beta _{B}^{Q}-\beta _{B}^{G}=\frac{{{U}_{B}}}{U}-\frac{1}{4}=\frac{4{{U}_{B}}-U}{4U},\]	
	\[\beta _{C}^{GQ}=\beta _{C}^{Q}-\beta _{C}^{G}=\frac{{{U}_{C}}}{U}-\frac{1}{4}=\frac{4{{U}_{C}}-U}{4U},\]	
	\[\beta _{D}^{GQ}=\beta _{D}^{Q}-\beta _{D}^{G}=\frac{{{U}_{D}}}{U}-\frac{1}{4}=\frac{4{{U}_{D}}-U}{4U}.\]	
	
	Therefore, from theorem \ref{thm:Thm24.2.1}, select ${{P}_{1}}$ as the centroid and ${{P}_{2}}$ as the circumcenter, then:
	\[\begin{aligned}
		G{{Q}^{2}} =&-\beta _{A}^{GQ}\beta _{B}^{GQ}A{{B}^{2}}-\beta _{A}^{GQ}\beta _{C}^{GQ}A{{C}^{2}}-\beta _{A}^{GQ}\beta _{D}^{GQ}A{{D}^{2}} \\ 
		& -\beta _{B}^{GQ}\beta _{C}^{GQ}B{{C}^{2}}-\beta _{B}^{GQ}\beta _{D}^{GQ}B{{D}^{2}}-\beta _{C}^{GQ}\beta _{D}^{GQ}C{{D}^{2}},  
	\end{aligned}\]
	i.e.
	\[G{{Q}^{2}}=-\frac{1}{16{{U}^{2}}}\left( \begin{aligned}
		& \left( 4{{U}_{A}}-U \right)\left( 4{{U}_{B}}-U \right)A{{B}^{2}}+\left( 4{{U}_{A}}-U \right)\left( 4{{U}_{C}}-U \right)A{{C}^{2}} \\ 
		& +\left( 4{{U}_{A}}-U \right)\left( 4{{U}_{D}}-U \right)A{{D}^{2}}+\left( 4{{U}_{B}}-U \right)\left( 4{{U}_{C}}-U \right)B{{C}^{2}} \\ 
		& +\left( 4{{U}_{B}}-U \right)\left( 4{{U}_{D}}-U \right)B{{D}^{2}}+\left( 4{{U}_{C}}-U \right)\left( 4{{U}_{D}}-U \right)C{{D}^{2}}  
	\end{aligned} \right).\]
\end{proof}
\hfill $\square$\par

\subsection{Distance between the incenter and the circumcenter of a tetrahedron}\label{Subsec26.2.3}
The incenter and the circumcenter of a tetrahedron are both the ICs-T, and the frame components of the two ICs-T have been calculated, so the distance between the incenter and the circumcenter of the tetrahedron can be calculated directly.

\begin{theorem}{Distance between incenter and circumcenter, Daiyuan Zhang) }{Thm26.2.3}\label{Thm26.2.3}
Suppose that given a tetrahedron $ABCD$, point $Q$ is its circumcenter and point $I$ is its incenter, then
\[I{{Q}^{2}}=-\frac{1}{{S}^{2}{U}^{2}}\left( \begin{aligned}
	& \left( S{{U}_{A}}-{{S}^{A}}U \right)\left( S{{U}_{B}}-{{S}^{B}}U \right)A{{B}^{2}} \\ 
	& +\left( S{{U}_{A}}-{{S}^{A}}U \right)\left( S{{U}_{C}}-{{S}^{C}}U \right)A{{C}^{2}} \\ 
	& +\left( S{{U}_{A}}-{{S}^{A}}U \right)\left( S{{U}_{D}}-{{S}^{D}}U \right)A{{D}^{2}} \\ 
	& +\left( S{{U}_{B}}-{{S}^{B}}U \right)\left( S{{U}_{C}}-{{S}^{C}}U \right)B{{C}^{2}} \\ 
	& +\left( S{{U}_{B}}-{{S}^{B}}U \right)\left( S{{U}_{D}}-{{S}^{D}}U \right)B{{D}^{2}} \\ 
	& +\left( S{{U}_{C}}-{{S}^{C}}U \right)\left( S{{U}_{D}}-{{S}^{D}}U \right)C{{D}^{2}}  
\end{aligned} \right).\]
where $U$, ${{U}_{A}}$, ${{U}_{B}}$, ${{U}_{C}}$, ${{U}_{D}}$ are calculated by theorem \ref{thm:Thm25.3.1}, theorem \ref{thm:Thm25.1.2} or theorem \ref{thm:Thm25.1.1}; ${{S}^{A}}$, ${{S}^{B}}$, ${{S}^{C}}$, ${{S}^{D}}$ are the areas of the triangle opposite to the four vertices $A$, $B$, $C$, $D$ of the tetrahedron $ABCD$ respectively; $S$ is the surface area of tetrahedron $ABCD$, i.e. $S={{S}^{A}}+{{S}^{B}}+{{S}^{C}}+{{S}^{D}}$; and $AB$, $AC$, $AD$, $BC$, $BD$, $CD$ are the lengths of the six edges of the tetrahedron respectively. 
\end{theorem}

\begin{proof}
	From theorem \ref{thm:Thm21.2.1} and theorem \ref{thm:Thm25.3.1}, theorem \ref{thm:Thm25.1.2} or \ref{thm:Thm25.1.1}, we have 
	\[\beta _{A}^{IQ}=\beta _{A}^{Q}-\beta _{A}^{I}=\frac{{{U}_{A}}}{U}-\frac{{{S}^{A}}}{S}=\frac{S{{U}_{A}}-{{S}^{A}}U}{SU},\]	
	\[\beta _{B}^{IQ}=\beta _{B}^{Q}-\beta _{B}^{I}=\frac{{{U}_{B}}}{U}-\frac{{{S}^{B}}}{S}=\frac{S{{U}_{B}}-{{S}^{B}}U}{SU},\]	
	\[\beta _{C}^{IQ}=\beta _{C}^{Q}-\beta _{C}^{I}=\frac{{{U}_{C}}}{U}-\frac{{{S}^{C}}}{S}=\frac{S{{U}_{C}}-{{S}^{C}}U}{SU},\]	
	\[\beta _{D}^{IQ}=\beta _{D}^{Q}-\beta _{D}^{I}=\frac{{{U}_{D}}}{U}-\frac{{{S}^{D}}}{S}=\frac{S{{U}_{D}}-{{S}^{D}}U}{SU}.\]	
	
	Therefore, from theorem \ref{thm:Thm24.2.1}, select ${{P}_{1}}$ as the incenter and ${{P}_{2}}$ as the circumcenter, then:
	\[\begin{aligned}
		I{{Q}^{2}}=& -\beta _{A}^{IQ}\beta _{B}^{IQ}A{{B}^{2}}-\beta _{A}^{IQ}\beta _{C}^{IQ}A{{C}^{2}}-\beta _{A}^{IQ}\beta _{D}^{IQ}A{{D}^{2}} \\ 
		& -\beta _{B}^{IQ}\beta _{C}^{IQ}B{{C}^{2}}-\beta _{B}^{IQ}\beta _{D}^{IQ}B{{D}^{2}}-\beta _{C}^{IQ}\beta _{D}^{IQ}C{{D}^{2}},  
	\end{aligned}\]
	i.e.
	\[I{{Q}^{2}}=-\frac{1}{{S}^{2}{U}^{2}}\left( \begin{aligned}
		& \left( S{{U}_{A}}-{{S}^{A}}U \right)\left( S{{U}_{B}}-{{S}^{B}}U \right)A{{B}^{2}} \\ 
		& +\left( S{{U}_{A}}-{{S}^{A}}U \right)\left( S{{U}_{C}}-{{S}^{C}}U \right)A{{C}^{2}} \\ 
		& +\left( S{{U}_{A}}-{{S}^{A}}U \right)\left( S{{U}_{D}}-{{S}^{D}}U \right)A{{D}^{2}} \\ 
		& +\left( S{{U}_{B}}-{{S}^{B}}U \right)\left( S{{U}_{C}}-{{S}^{C}}U \right)B{{C}^{2}} \\ 
		& +\left( S{{U}_{B}}-{{S}^{B}}U \right)\left( S{{U}_{D}}-{{S}^{D}}U \right)B{{D}^{2}} \\ 
		& +\left( S{{U}_{C}}-{{S}^{C}}U \right)\left( S{{U}_{D}}-{{S}^{D}}U \right)C{{D}^{2}}  
	\end{aligned} \right).\]
\end{proof}
\hfill $\square$\par

\subsection{(Distance between the centroid and the excenter of a tetrahedron (DTICs-T)}\label{Subsec26.2.4}
The following quantities are introduced:
\[S={{S}^{A}}+{{S}^{B}}+{{S}^{C}}+{{S}^{D}},\]
\[{{T}^{A}}=S-2{{S}^{A}}={{S}^{B}}+{{S}^{C}}+{{S}^{D}}-{{S}^{A}},\]
\[{{T}^{B}}=S-2{{S}^{B}}={{S}^{C}}+{{S}^{D}}+{{S}^{A}}-{{S}^{B}},\]
\[{{T}^{C}}=S-2{{S}^{C}}={{S}^{D}}+{{S}^{A}}+{{S}^{B}}-{{S}^{C}},\]
\[{{T}^{D}}=S-2{{S}^{D}}={{S}^{A}}+{{S}^{B}}+{{S}^{C}}-{{S}^{D}}.\]
The following theorem can be deduced directly by using theorem \ref{thm:Thm24.2.1} 
\begin{theorem}{Distance between centroid and excenter, Daiyuan Zhang}{Thm26.2.4}\label{Thm26.2.4} 
Suppose that given a tetrahedron $ABCD$, point $G$ is its centroid, and let ${{E}_{A}}$, ${{E}_{B}}$, ${{E}_{C}}$, ${{E}_{D}}$ be the centers in the trihedral of the vertex $A$, $B$, $C$, $D$ of tetrahedron $ABCD$ respectively, then
\[G{{E}_{A}}^{2}=\frac{1}{16{{\left( {{T}^{A}} \right)}^{2}}}\left( \begin{aligned}
	& \left( 4{{S}^{A}}+{{T}^{A}} \right)\left( 4{{S}^{B}}-{{T}^{A}} \right)A{{B}^{2}}+\left( 4{{S}^{A}}+{{T}^{A}} \right)\left( 4{{S}^{C}}-{{T}^{A}} \right)A{{C}^{2}} \\ 
	& +\left( 4{{S}^{A}}+{{T}^{A}} \right)\left( 4{{S}^{D}}-{{T}^{A}} \right)A{{D}^{2}}-\left( 4{{S}^{B}}-{{T}^{A}} \right)\left( 4{{S}^{C}}-{{T}^{A}} \right)B{{C}^{2}} \\ 
	& -\left( 4{{S}^{B}}-{{T}^{A}} \right)\left( 4{{S}^{D}}-{{T}^{A}} \right)B{{D}^{2}}-\left( 4{{S}^{C}}-{{T}^{A}} \right)\left( 4{{S}^{D}}-{{T}^{A}} \right)C{{D}^{2}}  
\end{aligned} \right),\]
\[G{{E}_{B}}^{2}=\frac{1}{16{{\left( {{T}^{B}} \right)}^{2}}}\left( \begin{aligned}
	& \left( 4{{S}^{B}}+{{T}^{B}} \right)\left( 4{{S}^{C}}-{{T}^{B}} \right)B{{C}^{2}}+\left( 4{{S}^{B}}+{{T}^{B}} \right)\left( 4{{S}^{D}}-{{T}^{B}} \right)B{{D}^{2}} \\ 
	& +\left( 4{{S}^{B}}+{{T}^{B}} \right)\left( 4{{S}^{A}}-{{T}^{B}} \right)B{{A}^{2}}-\left( 4{{S}^{C}}-{{T}^{B}} \right)\left( 4{{S}^{D}}-{{T}^{B}} \right)C{{D}^{2}} \\ 
	& -\left( 4{{S}^{C}}-{{T}^{B}} \right)\left( 4{{S}^{A}}-{{T}^{B}} \right)C{{A}^{2}}-\left( 4{{S}^{D}}-{{T}^{B}} \right)\left( 4{{S}^{A}}-{{T}^{B}} \right)D{{A}^{2}}  
\end{aligned} \right),\]
\[G{{E}_{C}}^{2}=\frac{1}{16{{\left( {{T}^{C}} \right)}^{2}}}\left( \begin{aligned}
	& \left( 4{{S}^{C}}+{{T}^{C}} \right)\left( 4{{S}^{D}}-{{T}^{C}} \right)C{{D}^{2}}+\left( 4{{S}^{C}}+{{T}^{C}} \right)\left( 4{{S}^{A}}-{{T}^{C}} \right)C{{A}^{2}} \\ 
	& +\left( 4{{S}^{C}}+{{T}^{C}} \right)\left( 4{{S}^{B}}-{{T}^{C}} \right)C{{B}^{2}}-\left( 4{{S}^{D}}-{{T}^{C}} \right)\left( 4{{S}^{A}}-{{T}^{C}} \right)D{{A}^{2}} \\ 
	& -\left( 4{{S}^{D}}-{{T}^{C}} \right)\left( 4{{S}^{B}}-{{T}^{C}} \right)D{{B}^{2}}-\left( 4{{S}^{A}}-{{T}^{C}} \right)\left( 4{{S}^{B}}-{{T}^{C}} \right)A{{B}^{2}}  
\end{aligned} \right),\]
\[G{{E}_{D}}^{2}=\frac{1}{16{{\left( {{T}^{D}} \right)}^{2}}}\left( \begin{aligned}
	& \left( 4{{S}^{D}}+{{T}^{D}} \right)\left( 4{{S}^{A}}-{{T}^{D}} \right)D{{A}^{2}}+\left( 4{{S}^{D}}+{{T}^{D}} \right)\left( 4{{S}^{B}}-{{T}^{D}} \right)D{{B}^{2}} \\ 
	& +\left( 4{{S}^{D}}+{{T}^{D}} \right)\left( 4{{S}^{C}}-{{T}^{D}} \right)D{{C}^{2}}-\left( 4{{S}^{A}}-{{T}^{D}} \right)\left( 4{{S}^{B}}-{{T}^{D}} \right)A{{B}^{2}} \\ 
	& -\left( 4{{S}^{A}}-{{T}^{D}} \right)\left( 4{{S}^{C}}-{{T}^{D}} \right)A{{C}^{2}}-\left( 4{{S}^{B}}-{{T}^{D}} \right)\left( 4{{S}^{C}}-{{T}^{D}} \right)B{{C}^{2}}  
\end{aligned} \right).\]

Where ${{S}^{A}}$, ${{S}^{B}}$, ${{S}^{C}}$, ${{S}^{D}}$ are the areas of the triangle opposite to the four vertices $A$, $B$, $C$, $D$ of the tetrahedron $ABCD$ respectively, and $AB$, $AC$, $AD$, $BC$, $BD$, $CD$ are the lengths of the six edges of the tetrahedron respectively. 
\end{theorem}

\begin{proof}
	From theorem \ref{thm:Thm21.3.1} and theorem \ref{thm:Thm21.1.1}, we have  
	\[\beta _{A}^{G{{E}_{A}}}=\beta _{A}^{{{E}_{A}}}-\beta _{A}^{G}=-\frac{{{S}^{A}}}{S-2{{S}^{A}}}-\frac{1}{4}=\frac{-4{{S}^{A}}-\left( S-2{{S}^{A}} \right)}{4\left( S-2{{S}^{A}} \right)}=-\frac{4{{S}^{A}}+{{T}^{A}}}{4{{T}^{A}}},\]
	\[\beta _{B}^{G{{E}_{A}}}=\beta _{B}^{{{E}_{A}}}-\beta _{B}^{G}=\frac{{{S}^{B}}}{S-2{{S}^{A}}}-\frac{1}{4}=\frac{4{{S}^{B}}-\left( S-2{{S}^{A}} \right)}{4\left( S-2{{S}^{A}} \right)}=\frac{4{{S}^{B}}-{{T}^{A}}}{4{{T}^{A}}},\]
	\[\beta _{C}^{G{{E}_{A}}}=\beta _{C}^{{{E}_{A}}}-\beta _{C}^{G}=\frac{{{S}^{C}}}{S-2{{S}^{A}}}-\frac{1}{4}=\frac{4{{S}^{C}}-\left( S-2{{S}^{A}} \right)}{4\left( S-2{{S}^{A}} \right)}=\frac{4{{S}^{C}}-{{T}^{A}}}{4{{T}^{A}}},\]
	\[\beta _{D}^{G{{E}_{A}}}=\beta _{D}^{{{E}_{A}}}-\beta _{D}^{G}=\frac{{{S}^{D}}}{S-2{{S}^{A}}}-\frac{1}{4}=\frac{4{{S}^{D}}-\left( S-2{{S}^{A}} \right)}{4\left( S-2{{S}^{A}} \right)}=\frac{4{{S}^{D}}-{{T}^{A}}}{4{{T}^{A}}}.\]
	
	Therefore, according to theorem \ref{thm:Thm24.2.1}, let ${{P}_{1}}$ be the centroid $G$, ${{P}_{2}}$ be the excenter ${{E}_{A}}$, the following results are obtained:
	\[\begin{aligned}
		G{{E}_{A}}^{2}=& -\beta _{A}^{G{{E}_{A}}}\beta _{B}^{G{{E}_{A}}}A{{B}^{2}}-\beta _{A}^{G{{E}_{A}}}\beta _{C}^{G{{E}_{A}}}A{{C}^{2}}-\beta _{A}^{G{{E}_{A}}}\beta _{D}^{G{{E}_{A}}}A{{D}^{2}} \\ 
		& -\beta _{B}^{G{{E}_{A}}}\beta _{C}^{G{{E}_{A}}}B{{C}^{2}}-\beta _{B}^{G{{E}_{A}}}\beta _{D}^{G{{E}_{A}}}B{{D}^{2}}-\beta _{C}^{G{{E}_{A}}}\beta _{D}^{G{{E}_{A}}}C{{D}^{2}},  
	\end{aligned}\]
	i.e.
	\[G{{E}_{A}}^{2}=\frac{1}{16{{\left( {{T}^{A}} \right)}^{2}}}\left( \begin{aligned}
		& \left( 4{{S}^{A}}+{{T}^{A}} \right)\left( 4{{S}^{B}}-{{T}^{A}} \right)A{{B}^{2}}+\left( 4{{S}^{A}}+{{T}^{A}} \right)\left( 4{{S}^{C}}-{{T}^{A}} \right)A{{C}^{2}} \\ 
		& +\left( 4{{S}^{A}}+{{T}^{A}} \right)\left( 4{{S}^{D}}-{{T}^{A}} \right)A{{D}^{2}}-\left( 4{{S}^{B}}-{{T}^{A}} \right)\left( 4{{S}^{C}}-{{T}^{A}} \right)B{{C}^{2}} \\ 
		& -\left( 4{{S}^{B}}-{{T}^{A}} \right)\left( 4{{S}^{D}}-{{T}^{A}} \right)B{{D}^{2}}-\left( 4{{S}^{C}}-{{T}^{A}} \right)\left( 4{{S}^{D}}-{{T}^{A}} \right)C{{D}^{2}}  
	\end{aligned} \right).\]
	
	Similarly:
	\[G{{E}_{B}}^{2}=\frac{1}{16{{\left( {{T}^{B}} \right)}^{2}}}\left( \begin{aligned}
		& \left( 4{{S}^{B}}+{{T}^{B}} \right)\left( 4{{S}^{C}}-{{T}^{B}} \right)B{{C}^{2}}+\left( 4{{S}^{B}}+{{T}^{B}} \right)\left( 4{{S}^{D}}-{{T}^{B}} \right)B{{D}^{2}} \\ 
		& +\left( 4{{S}^{B}}+{{T}^{B}} \right)\left( 4{{S}^{A}}-{{T}^{B}} \right)B{{A}^{2}}-\left( 4{{S}^{C}}-{{T}^{B}} \right)\left( 4{{S}^{D}}-{{T}^{B}} \right)C{{D}^{2}} \\ 
		& -\left( 4{{S}^{C}}-{{T}^{B}} \right)\left( 4{{S}^{A}}-{{T}^{B}} \right)C{{A}^{2}}-\left( 4{{S}^{D}}-{{T}^{B}} \right)\left( 4{{S}^{A}}-{{T}^{B}} \right)D{{A}^{2}}  
	\end{aligned} \right),\]
	\[G{{E}_{C}}^{2}=\frac{1}{16{{\left( {{T}^{C}} \right)}^{2}}}\left( \begin{aligned}
		& \left( 4{{S}^{C}}+{{T}^{C}} \right)\left( 4{{S}^{D}}-{{T}^{C}} \right)C{{D}^{2}}+\left( 4{{S}^{C}}+{{T}^{C}} \right)\left( 4{{S}^{A}}-{{T}^{C}} \right)C{{A}^{2}} \\ 
		& +\left( 4{{S}^{C}}+{{T}^{C}} \right)\left( 4{{S}^{B}}-{{T}^{C}} \right)C{{B}^{2}}-\left( 4{{S}^{D}}-{{T}^{C}} \right)\left( 4{{S}^{A}}-{{T}^{C}} \right)D{{A}^{2}} \\ 
		& -\left( 4{{S}^{D}}-{{T}^{C}} \right)\left( 4{{S}^{B}}-{{T}^{C}} \right)D{{B}^{2}}-\left( 4{{S}^{A}}-{{T}^{C}} \right)\left( 4{{S}^{B}}-{{T}^{C}} \right)A{{B}^{2}}  
	\end{aligned} \right),\]
	\[G{{E}_{D}}^{2}=\frac{1}{16{{\left( {{T}^{D}} \right)}^{2}}}\left( \begin{aligned}
		& \left( 4{{S}^{D}}+{{T}^{D}} \right)\left( 4{{S}^{A}}-{{T}^{D}} \right)D{{A}^{2}}+\left( 4{{S}^{D}}+{{T}^{D}} \right)\left( 4{{S}^{B}}-{{T}^{D}} \right)D{{B}^{2}} \\ 
		& +\left( 4{{S}^{D}}+{{T}^{D}} \right)\left( 4{{S}^{C}}-{{T}^{D}} \right)D{{C}^{2}}-\left( 4{{S}^{A}}-{{T}^{D}} \right)\left( 4{{S}^{B}}-{{T}^{D}} \right)A{{B}^{2}} \\ 
		& -\left( 4{{S}^{A}}-{{T}^{D}} \right)\left( 4{{S}^{C}}-{{T}^{D}} \right)A{{C}^{2}}-\left( 4{{S}^{B}}-{{T}^{D}} \right)\left( 4{{S}^{C}}-{{T}^{D}} \right)B{{C}^{2}}  
	\end{aligned} \right).\]	
\end{proof}
\hfill $\square$\par

\subsection{Distance between the incenter and the excenter of a tetrahedron (DTICs-T)}\label{Subsec26.2.5}
\begin{theorem}{Distance between incenter and excenter, Daiyuan Zhang}{Thm26.2.5}\label{Thm26.2.5} 
Suppose that given a tetrahedron $ABCD$, point $I$ is its incenter, and let ${{E}_{A}}$, ${{E}_{B}}$, ${{E}_{C}}$, ${{E}_{D}}$ be the centers in the trihedral of the vertex $A$, $B$, $C$, $D$ of tetrahedron $ABCD$ respectively, then
\[I{{E}_{A}}^{2}=\frac{1}{{{\left( S{{T}^{A}} \right)}^{2}}}\left( \begin{aligned}
	& \left( {{S}^{2}}-{{\left( {{T}^{A}} \right)}^{2}} \right)\left( {{S}^{A}}{{S}^{B}}A{{B}^{2}}+{{S}^{A}}{{S}^{C}}A{{C}^{2}}+{{S}^{A}}{{S}^{D}}A{{D}^{2}} \right) \\ 
	& -{{\left( S-{{T}^{A}} \right)}^{2}}\left( {{S}^{B}}{{S}^{C}}B{{C}^{2}}+{{S}^{B}}{{S}^{D}}B{{D}^{2}}+{{S}^{C}}{{S}^{D}}C{{D}^{2}} \right)  
\end{aligned} \right)\text{,}\]
\[I{{E}_{B}}^{2}=\frac{1}{{{\left( S{{T}^{B}} \right)}^{2}}}\left( \begin{aligned}
	& \left( {{S}^{2}}-{{\left( {{T}^{B}} \right)}^{2}} \right)\left( {{S}^{B}}{{S}^{C}}B{{C}^{2}}+{{S}^{B}}{{S}^{D}}B{{D}^{2}}+{{S}^{B}}{{S}^{A}}B{{A}^{2}} \right) \\ 
	& -{{\left( S-{{T}^{B}} \right)}^{2}}\left( {{S}^{C}}{{S}^{D}}C{{D}^{2}}+{{S}^{C}}{{S}^{A}}C{{A}^{2}}+{{S}^{D}}{{S}^{A}}D{{A}^{2}} \right)  
\end{aligned} \right)\text{,}\]
\[I{{E}_{C}}^{2}=\frac{1}{{{\left( S{{T}^{C}} \right)}^{2}}}\left( \begin{aligned}
	& \left( {{S}^{2}}-{{\left( {{T}^{C}} \right)}^{2}} \right)\left( {{S}^{C}}{{S}^{D}}C{{D}^{2}}+{{S}^{C}}{{S}^{A}}C{{A}^{2}}+{{S}^{C}}{{S}^{B}}C{{B}^{2}} \right) \\ 
	& -{{\left( S-{{T}^{C}} \right)}^{2}}\left( {{S}^{D}}{{S}^{A}}D{{A}^{2}}+{{S}^{D}}{{S}^{B}}D{{B}^{2}}+{{S}^{A}}{{S}^{B}}A{{B}^{2}} \right)  
\end{aligned} \right)\text{,}\]
\[I{{E}_{D}}^{2}=\frac{1}{{{\left( S{{T}^{D}} \right)}^{2}}}\left( \begin{aligned}
	& \left( {{S}^{2}}-{{\left( {{T}^{D}} \right)}^{2}} \right)\left( {{S}^{D}}{{S}^{A}}D{{A}^{2}}+{{S}^{D}}{{S}^{B}}D{{B}^{2}}+{{S}^{D}}{{S}^{C}}D{{C}^{2}} \right) \\ 
	& -{{\left( S-{{T}^{D}} \right)}^{2}}\left( {{S}^{A}}{{S}^{B}}A{{B}^{2}}+{{S}^{A}}{{S}^{C}}A{{C}^{2}}+{{S}^{B}}{{S}^{C}}B{{C}^{2}} \right)  
\end{aligned} \right).\]
\end{theorem}

\begin{proof}
	From theorem \ref{thm:Thm21.3.1} and theorem \ref{thm:Thm21.2.1}, we have
	\[\beta _{A}^{I{{E}_{A}}}=\beta _{A}^{{{E}_{A}}}-\beta _{A}^{I}=-\frac{{{S}^{A}}}{S-2{{S}^{A}}}-\frac{{{S}^{A}}}{S}=-\frac{{{S}^{A}}\left( S+\left( S-2{{S}^{A}} \right) \right)}{S\left( S-2{{S}^{A}} \right)}=-\frac{{{S}^{A}}\left( S+{{T}^{A}} \right)}{S{{T}^{A}}},\]
	\[\beta _{B}^{I{{E}_{A}}}=\beta _{B}^{{{E}_{A}}}-\beta _{B}^{I}=\frac{{{S}^{B}}}{S-2{{S}^{A}}}-\frac{{{S}^{B}}}{S}=\frac{{{S}^{B}}\left( S-\left( S-2{{S}^{A}} \right) \right)}{S\left( S-2{{S}^{A}} \right)}=\frac{{{S}^{B}}\left( S-{{T}^{A}} \right)}{S{{T}^{A}}},\]
	\[\beta _{C}^{I{{E}_{A}}}=\beta _{C}^{{{E}_{A}}}-\beta _{C}^{I}=\frac{{{S}^{C}}}{S-2{{S}^{A}}}-\frac{{{S}^{C}}}{S}=\frac{{{S}^{C}}\left( S-\left( S-2{{S}^{A}} \right) \right)}{S\left( S-2{{S}^{A}} \right)}=\frac{{{S}^{C}}\left( S-{{T}^{A}} \right)}{S{{T}^{A}}},\]
	\[\beta _{D}^{I{{E}_{A}}}=\beta _{D}^{{{E}_{A}}}-\beta _{D}^{I}=\frac{{{S}^{D}}}{S-2{{S}^{A}}}-\frac{{{S}^{D}}}{S}=\frac{{{S}^{D}}\left( S-\left( S-2{{S}^{A}} \right) \right)}{S\left( S-2{{S}^{A}} \right)}=\frac{{{S}^{D}}\left( S-{{T}^{A}} \right)}{S{{T}^{A}}}.\]
	
	Therefore, according to theorem \ref{thm:Thm24.2.1}, let ${{P}_{1}}$ be the incenter $I$, ${{P}_{2}}$ be the excenter ${{E}_{A}}$, the following results are obtained:
	\[\begin{aligned}
		I{{E}_{A}}^{2}=& -\beta _{A}^{I{{E}_{A}}}\beta _{B}^{I{{E}_{A}}}A{{B}^{2}}-\beta _{A}^{I{{E}_{A}}}\beta _{C}^{I{{E}_{A}}}A{{C}^{2}}-\beta _{A}^{I{{E}_{A}}}\beta _{D}^{I{{E}_{A}}}A{{D}^{2}} \\ 
		& -\beta _{B}^{I{{E}_{A}}}\beta _{C}^{I{{E}_{A}}}B{{C}^{2}}-\beta _{B}^{I{{E}_{A}}}\beta _{D}^{I{{E}_{A}}}B{{D}^{2}}-\beta _{C}^{I{{E}_{A}}}\beta _{D}^{I{{E}_{A}}}C{{D}^{2}},  
	\end{aligned}\]
	i.e.
	\[I{{E}_{A}}^{2}=\frac{1}{{{\left( S{{T}^{A}} \right)}^{2}}}\left( \begin{aligned}
		& {{S}^{A}}\left( S+{{T}^{A}} \right){{S}^{B}}\left( S-{{T}^{A}} \right)A{{B}^{2}}+{{S}^{A}}\left( S+{{T}^{A}} \right){{S}^{C}}\left( S-{{T}^{A}} \right)A{{C}^{2}} \\ 
		& +{{S}^{A}}\left( S+{{T}^{A}} \right){{S}^{D}}\left( S-{{T}^{A}} \right)A{{D}^{2}}-{{S}^{B}}\left( S-{{T}^{A}} \right){{S}^{C}}\left( S-{{T}^{A}} \right)B{{C}^{2}} \\ 
		& -{{S}^{B}}\left( S-{{T}^{A}} \right){{S}^{D}}\left( S-{{T}^{A}} \right)B{{D}^{2}}-{{S}^{C}}\left( S-{{T}^{A}} \right){{S}^{D}}\left( S-{{T}^{A}} \right)C{{D}^{2}}  
	\end{aligned} \right)\text{,}\]
	i.e.
	\[I{{E}_{A}}^{2}=\frac{1}{{{\left( S{{T}^{A}} \right)}^{2}}}\left( \begin{aligned}
		& \left( {{S}^{2}}-{{\left( {{T}^{A}} \right)}^{2}} \right)\left( {{S}^{A}}{{S}^{B}}A{{B}^{2}}+{{S}^{A}}{{S}^{C}}A{{C}^{2}}+{{S}^{A}}{{S}^{D}}A{{D}^{2}} \right) \\ 
		& -{{\left( S-{{T}^{A}} \right)}^{2}}\left( {{S}^{B}}{{S}^{C}}B{{C}^{2}}+{{S}^{B}}{{S}^{D}}B{{D}^{2}}+{{S}^{C}}{{S}^{D}}C{{D}^{2}} \right)  
	\end{aligned} \right).\]
	
	Similarly:
	\[I{{E}_{B}}^{2}=\frac{1}{{{\left( S{{T}^{B}} \right)}^{2}}}\left( \begin{aligned}
		& \left( {{S}^{2}}-{{\left( {{T}^{B}} \right)}^{2}} \right)\left( {{S}^{B}}{{S}^{C}}B{{C}^{2}}+{{S}^{B}}{{S}^{D}}B{{D}^{2}}+{{S}^{B}}{{S}^{A}}B{{A}^{2}} \right) \\ 
		& -{{\left( S-{{T}^{B}} \right)}^{2}}\left( {{S}^{C}}{{S}^{D}}C{{D}^{2}}+{{S}^{C}}{{S}^{A}}C{{A}^{2}}+{{S}^{D}}{{S}^{A}}D{{A}^{2}} \right)  
	\end{aligned} \right)\text{,}\]
	\[I{{E}_{C}}^{2}=\frac{1}{{{\left( S{{T}^{C}} \right)}^{2}}}\left( \begin{aligned}
		& \left( {{S}^{2}}-{{\left( {{T}^{C}} \right)}^{2}} \right)\left( {{S}^{C}}{{S}^{D}}C{{D}^{2}}+{{S}^{C}}{{S}^{A}}C{{A}^{2}}+{{S}^{C}}{{S}^{B}}C{{B}^{2}} \right) \\ 
		& -{{\left( S-{{T}^{C}} \right)}^{2}}\left( {{S}^{D}}{{S}^{A}}D{{A}^{2}}+{{S}^{D}}{{S}^{B}}D{{B}^{2}}+{{S}^{A}}{{S}^{B}}A{{B}^{2}} \right)  
	\end{aligned} \right)\text{,}\]
	\[I{{E}_{D}}^{2}=\frac{1}{{{\left( S{{T}^{D}} \right)}^{2}}}\left( \begin{aligned}
		& \left( {{S}^{2}}-{{\left( {{T}^{D}} \right)}^{2}} \right)\left( {{S}^{D}}{{S}^{A}}D{{A}^{2}}+{{S}^{D}}{{S}^{B}}D{{B}^{2}}+{{S}^{D}}{{S}^{C}}D{{C}^{2}} \right) \\ 
		& -{{\left( S-{{T}^{D}} \right)}^{2}}\left( {{S}^{A}}{{S}^{B}}A{{B}^{2}}+{{S}^{A}}{{S}^{C}}A{{C}^{2}}+{{S}^{B}}{{S}^{C}}B{{C}^{2}} \right)  
	\end{aligned} \right).\]	
\end{proof}
\hfill $\square$\par

\subsection{Distance between the circumcenter and the excenter of a tetrahedron (DTICs-T)}\label{Subsec26.2.6}
\begin{theorem}{Distance between circumcenter and excenter, Daiyuan Zhang}{Thm26.2.6}\label{Thm26.2.6} 
Suppose that given a tetrahedron $ABCD$, point $Q$ is its circumcenter, and let ${{E}_{A}}$, ${{E}_{B}}$, ${{E}_{C}}$, ${{E}_{D}}$ be the centers in the trihedral of the vertex $A$, $B$, $C$, $D$ of tetrahedron $ABCD$ respectively, then
\[Q{{E}_{A}}^{2}=\frac{1}{{{\left( {{T}^{A}}U \right)}^{2}}}\left( \begin{aligned}
	& \left( {{S}^{A}}U+{{T}^{A}}{{U}_{A}} \right)\left( {{S}^{B}}U-{{T}^{A}}{{U}_{B}} \right)A{{B}^{2}} \\ 
	& +\left( {{S}^{A}}U+{{T}^{A}}{{U}_{A}} \right)\left( {{S}^{C}}U-{{T}^{A}}{{U}_{C}} \right)A{{C}^{2}} \\ 
	& +\left( {{S}^{A}}U+{{T}^{A}}{{U}_{A}} \right)\left( {{S}^{D}}U-{{T}^{A}}{{U}_{D}} \right)A{{D}^{2}} \\ 
	& -\left( {{S}^{B}}U-{{T}^{A}}{{U}_{B}} \right)\left( {{S}^{C}}U-{{T}^{A}}{{U}_{C}} \right)B{{C}^{2}} \\ 
	& -\left( {{S}^{B}}U-{{T}^{A}}{{U}_{B}} \right)\left( {{S}^{D}}U-{{T}^{A}}{{U}_{D}} \right)B{{D}^{2}} \\ 
	& -\left( {{S}^{C}}U-{{T}^{A}}{{U}_{C}} \right)\left( {{S}^{D}}U-{{T}^{A}}{{U}_{D}} \right)C{{D}^{2}}  
\end{aligned} \right)\text{,}\]
\[Q{{E}_{B}}^{2}=\frac{1}{{{\left( {{T}^{B}}U \right)}^{2}}}\left( \begin{aligned}
	& \left( {{S}^{B}}U+{{T}^{B}}{{U}_{B}} \right)\left( {{S}^{C}}U-{{T}^{B}}{{U}_{C}} \right)B{{C}^{2}} \\ 
	& +\left( {{S}^{B}}U+{{T}^{B}}{{U}_{B}} \right)\left( {{S}^{D}}U-{{T}^{B}}{{U}_{D}} \right)B{{D}^{2}} \\ 
	& +\left( {{S}^{B}}U+{{T}^{B}}{{U}_{B}} \right)\left( {{S}^{A}}U-{{T}^{B}}{{U}_{A}} \right)B{{A}^{2}} \\ 
	& -\left( {{S}^{C}}U-{{T}^{B}}{{U}_{C}} \right)\left( {{S}^{D}}U-{{T}^{B}}{{U}_{D}} \right)C{{D}^{2}} \\ 
	& -\left( {{S}^{C}}U-{{T}^{B}}{{U}_{C}} \right)\left( {{S}^{A}}U-{{T}^{B}}{{U}_{A}} \right)C{{A}^{2}} \\ 
	& -\left( {{S}^{D}}U-{{T}^{B}}{{U}_{D}} \right)\left( {{S}^{A}}U-{{T}^{B}}{{U}_{A}} \right)D{{A}^{2}}  
\end{aligned} \right)\text{,}\]
\[Q{{E}_{C}}^{2}=\frac{1}{{{\left( {{T}^{C}}U \right)}^{2}}}\left( \begin{aligned}
	& \left( {{S}^{C}}U+{{T}^{C}}{{U}_{C}} \right)\left( {{S}^{D}}U-{{T}^{C}}{{U}_{D}} \right)C{{D}^{2}} \\ 
	& +\left( {{S}^{C}}U+{{T}^{C}}{{U}_{C}} \right)\left( {{S}^{A}}U-{{T}^{C}}{{U}_{A}} \right)C{{A}^{2}} \\ 
	& +\left( {{S}^{C}}U+{{T}^{C}}{{U}_{C}} \right)\left( {{S}^{B}}U-{{T}^{C}}{{U}_{B}} \right)C{{B}^{2}} \\ 
	& -\left( {{S}^{D}}U-{{T}^{C}}{{U}_{D}} \right)\left( {{S}^{A}}U-{{T}^{C}}{{U}_{A}} \right)D{{A}^{2}} \\ 
	& -\left( {{S}^{D}}U-{{T}^{C}}{{U}_{D}} \right)\left( {{S}^{B}}U-{{T}^{C}}{{U}_{B}} \right)D{{B}^{2}} \\ 
	& -\left( {{S}^{A}}U-{{T}^{C}}{{U}_{A}} \right)\left( {{S}^{B}}U-{{T}^{C}}{{U}_{B}} \right)A{{B}^{2}}  
\end{aligned} \right)\text{,}\]
\[Q{{E}_{D}}^{2}=\frac{1}{{{\left( {{T}^{D}}U \right)}^{2}}}\left( \begin{aligned}
	& \left( {{S}^{D}}U+{{T}^{D}}{{U}_{D}} \right)\left( {{S}^{A}}U-{{T}^{D}}{{U}_{A}} \right)D{{A}^{2}} \\ 
	& +\left( {{S}^{D}}U+{{T}^{D}}{{U}_{D}} \right)\left( {{S}^{B}}U-{{T}^{D}}{{U}_{B}} \right)D{{B}^{2}} \\ 
	& +\left( {{S}^{D}}U+{{T}^{D}}{{U}_{D}} \right)\left( {{S}^{C}}U-{{T}^{D}}{{U}_{C}} \right)D{{C}^{2}} \\ 
	& -\left( {{S}^{A}}U-{{T}^{D}}{{U}_{A}} \right)\left( {{S}^{B}}U-{{T}^{D}}{{U}_{B}} \right)A{{B}^{2}} \\ 
	& -\left( {{S}^{A}}U-{{T}^{D}}{{U}_{A}} \right)\left( {{S}^{C}}U-{{T}^{D}}{{U}_{C}} \right)A{{C}^{2}} \\ 
	& -\left( {{S}^{B}}U-{{T}^{D}}{{U}_{B}} \right)\left( {{S}^{C}}U-{{T}^{D}}{{U}_{C}} \right)B{{C}^{2}}  
\end{aligned} \right)\text{.}\]
\end{theorem}

\begin{proof}
	From theorem \ref{thm:Thm21.3.1} and theorem \ref{thm:Thm25.3.1} or theorem \ref{thm:Thm25.1.1}, we have
	\[\beta _{A}^{Q{{E}_{A}}}=\beta _{A}^{{{E}_{A}}}-\beta _{A}^{Q}=-\frac{{{S}^{A}}}{S-2{{S}^{A}}}-\frac{{{U}_{A}}}{U}=-\frac{{{S}^{A}}U+\left( S-2{{S}^{A}} \right){{U}_{A}}}{\left( S-2{{S}^{A}} \right)U}=-\frac{{{S}^{A}}U+{{T}^{A}}{{U}_{A}}}{{{T}^{A}}U},\]
	\[\beta _{B}^{Q{{E}_{A}}}=\beta _{B}^{{{E}_{A}}}-\beta _{B}^{Q}=\frac{{{S}^{B}}}{S-2{{S}^{A}}}-\frac{{{U}_{B}}}{U}=\frac{{{S}^{B}}U-\left( S-2{{S}^{A}} \right){{U}_{B}}}{\left( S-2{{S}^{A}} \right)U}=\frac{{{S}^{B}}U-{{T}^{A}}{{U}_{B}}}{{{T}^{A}}U},\]
	\[\beta _{C}^{Q{{E}_{A}}}=\beta _{C}^{{{E}_{A}}}-\beta _{C}^{Q}=\frac{{{S}^{C}}}{S-2{{S}^{A}}}-\frac{{{U}_{C}}}{U}=\frac{{{S}^{C}}U-\left( S-2{{S}^{A}} \right){{U}_{C}}}{\left( S-2{{S}^{A}} \right)U}=\frac{{{S}^{C}}U-{{T}^{A}}{{U}_{C}}}{{{T}^{A}}U},\]
	\[\beta _{D}^{Q{{E}_{A}}}=\beta _{D}^{{{E}_{A}}}-\beta _{D}^{Q}=\frac{{{S}^{D}}}{S-2{{S}^{A}}}-\frac{{{U}_{D}}}{U}=\frac{{{S}^{D}}U-\left( S-2{{S}^{A}} \right){{U}_{D}}}{\left( S-2{{S}^{A}} \right)U}=\frac{{{S}^{D}}U-{{T}^{A}}{{U}_{D}}}{{{T}^{A}}U}.\]
	
	Therefore, according to theorem \ref{thm:Thm24.2.1}, let ${{P}_{1}}$ be the circumcenter $Q$, ${{P}_{2}}$ be the excenter ${{E}_{A}}$, the following results are obtained:
	\[\begin{aligned}
		Q{{E}_{A}}^{2}=& -\beta _{A}^{Q{{E}_{A}}}\beta _{B}^{Q{{E}_{A}}}A{{B}^{2}}-\beta _{A}^{Q{{E}_{A}}}\beta _{C}^{Q{{E}_{A}}}A{{C}^{2}}-\beta _{A}^{Q{{E}_{A}}}\beta _{D}^{Q{{E}_{A}}}A{{D}^{2}} \\ 
		& -\beta _{B}^{Q{{E}_{A}}}\beta _{C}^{Q{{E}_{A}}}B{{C}^{2}}-\beta _{B}^{Q{{E}_{A}}}\beta _{D}^{Q{{E}_{A}}}B{{D}^{2}}-\beta _{C}^{Q{{E}_{A}}}\beta _{D}^{Q{{E}_{A}}}C{{D}^{2}}\text{,}  
	\end{aligned}\]
	i.e.
	\[Q{{E}_{A}}^{2}=\frac{1}{{{\left( {{T}^{A}}U \right)}^{2}}}\left( \begin{aligned}
		& \left( {{S}^{A}}U+{{T}^{A}}{{U}_{A}} \right)\left( {{S}^{B}}U-{{T}^{A}}{{U}_{B}} \right)A{{B}^{2}} \\ 
		& +\left( {{S}^{A}}U+{{T}^{A}}{{U}_{A}} \right)\left( {{S}^{C}}U-{{T}^{A}}{{U}_{C}} \right)A{{C}^{2}} \\ 
		& +\left( {{S}^{A}}U+{{T}^{A}}{{U}_{A}} \right)\left( {{S}^{D}}U-{{T}^{A}}{{U}_{D}} \right)A{{D}^{2}} \\ 
		& -\left( {{S}^{B}}U-{{T}^{A}}{{U}_{B}} \right)\left( {{S}^{C}}U-{{T}^{A}}{{U}_{C}} \right)B{{C}^{2}} \\ 
		& -\left( {{S}^{B}}U-{{T}^{A}}{{U}_{B}} \right)\left( {{S}^{D}}U-{{T}^{A}}{{U}_{D}} \right)B{{D}^{2}} \\ 
		& -\left( {{S}^{C}}U-{{T}^{A}}{{U}_{C}} \right)\left( {{S}^{D}}U-{{T}^{A}}{{U}_{D}} \right)C{{D}^{2}}  
	\end{aligned} \right)\text{,}\]
	
	Similarly:
	\[Q{{E}_{B}}^{2}=\frac{1}{{{\left( {{T}^{B}}U \right)}^{2}}}\left( \begin{aligned}
		& \left( {{S}^{B}}U+{{T}^{B}}{{U}_{B}} \right)\left( {{S}^{C}}U-{{T}^{B}}{{U}_{C}} \right)B{{C}^{2}} \\ 
		& +\left( {{S}^{B}}U+{{T}^{B}}{{U}_{B}} \right)\left( {{S}^{D}}U-{{T}^{B}}{{U}_{D}} \right)B{{D}^{2}} \\ 
		& +\left( {{S}^{B}}U+{{T}^{B}}{{U}_{B}} \right)\left( {{S}^{A}}U-{{T}^{B}}{{U}_{A}} \right)B{{A}^{2}} \\ 
		& -\left( {{S}^{C}}U-{{T}^{B}}{{U}_{C}} \right)\left( {{S}^{D}}U-{{T}^{B}}{{U}_{D}} \right)C{{D}^{2}} \\ 
		& -\left( {{S}^{C}}U-{{T}^{B}}{{U}_{C}} \right)\left( {{S}^{A}}U-{{T}^{B}}{{U}_{A}} \right)C{{A}^{2}} \\ 
		& -\left( {{S}^{D}}U-{{T}^{B}}{{U}_{D}} \right)\left( {{S}^{A}}U-{{T}^{B}}{{U}_{A}} \right)D{{A}^{2}}  
	\end{aligned} \right)\text{,}\]
	\[Q{{E}_{C}}^{2}=\frac{1}{{{\left( {{T}^{C}}U \right)}^{2}}}\left( \begin{aligned}
		& \left( {{S}^{C}}U+{{T}^{C}}{{U}_{C}} \right)\left( {{S}^{D}}U-{{T}^{C}}{{U}_{D}} \right)C{{D}^{2}} \\ 
		& +\left( {{S}^{C}}U+{{T}^{C}}{{U}_{C}} \right)\left( {{S}^{A}}U-{{T}^{C}}{{U}_{A}} \right)C{{A}^{2}} \\ 
		& +\left( {{S}^{C}}U+{{T}^{C}}{{U}_{C}} \right)\left( {{S}^{B}}U-{{T}^{C}}{{U}_{B}} \right)C{{B}^{2}} \\ 
		& -\left( {{S}^{D}}U-{{T}^{C}}{{U}_{D}} \right)\left( {{S}^{A}}U-{{T}^{C}}{{U}_{A}} \right)D{{A}^{2}} \\ 
		& -\left( {{S}^{D}}U-{{T}^{C}}{{U}_{D}} \right)\left( {{S}^{B}}U-{{T}^{C}}{{U}_{B}} \right)D{{B}^{2}} \\ 
		& -\left( {{S}^{A}}U-{{T}^{C}}{{U}_{A}} \right)\left( {{S}^{B}}U-{{T}^{C}}{{U}_{B}} \right)A{{B}^{2}}  
	\end{aligned} \right)\text{,}\]
	\[Q{{E}_{D}}^{2}=\frac{1}{{{\left( {{T}^{D}}U \right)}^{2}}}\left( \begin{aligned}
		& \left( {{S}^{D}}U+{{T}^{D}}{{U}_{D}} \right)\left( {{S}^{A}}U-{{T}^{D}}{{U}_{A}} \right)D{{A}^{2}} \\ 
		& +\left( {{S}^{D}}U+{{T}^{D}}{{U}_{D}} \right)\left( {{S}^{B}}U-{{T}^{D}}{{U}_{B}} \right)D{{B}^{2}} \\ 
		& +\left( {{S}^{D}}U+{{T}^{D}}{{U}_{D}} \right)\left( {{S}^{C}}U-{{T}^{D}}{{U}_{C}} \right)D{{C}^{2}} \\ 
		& -\left( {{S}^{A}}U-{{T}^{D}}{{U}_{A}} \right)\left( {{S}^{B}}U-{{T}^{D}}{{U}_{B}} \right)A{{B}^{2}} \\ 
		& -\left( {{S}^{A}}U-{{T}^{D}}{{U}_{A}} \right)\left( {{S}^{C}}U-{{T}^{D}}{{U}_{C}} \right)A{{C}^{2}} \\ 
		& -\left( {{S}^{B}}U-{{T}^{D}}{{U}_{B}} \right)\left( {{S}^{C}}U-{{T}^{D}}{{U}_{C}} \right)B{{C}^{2}}  
	\end{aligned} \right)\text{.}\]
\end{proof}
\hfill $\square$\par

\subsection{Distance between the excenters of a tetrahedron (DTICs-T)}\label{Subsec26.2.7}
\begin{theorem}{Distance between excenters, Daiyuan Zhang}{Thm26.2.7}\label{Thm26.2.7} 
Suppose that given a tetrahedron $ABCD$, let ${{E}_{A}}$, ${{E}_{B}}$, ${{E}_{C}}$, ${{E}_{D}}$ be the centers in the trihedral of the vertex $A$, $B$, $C$, $D$ of tetrahedron $ABCD$ respectively, then
\[{{E}_{A}}{{E}_{B}}^{2}=\frac{1}{{{\left( {{T}^{A}}{{T}^{B}} \right)}^{2}}}\left( \begin{aligned}
	& {{S}^{A}}{{S}^{B}}{{\left( {{T}^{A}}+{{T}^{B}} \right)}^{2}}A{{B}^{2}}-{{S}^{A}}{{S}^{C}}\left( {{\left( {{T}^{A}} \right)}^{2}}-{{\left( {{T}^{B}} \right)}^{2}} \right)A{{C}^{2}} \\ 
	& -{{S}^{A}}{{S}^{D}}\left( {{\left( {{T}^{A}} \right)}^{2}}-{{\left( {{T}^{B}} \right)}^{2}} \right)A{{D}^{2}}+{{S}^{B}}{{S}^{C}}\left( {{\left( {{T}^{A}} \right)}^{2}}-{{\left( {{T}^{B}} \right)}^{2}} \right)B{{C}^{2}} \\ 
	& +{{S}^{B}}{{S}^{D}}\left( {{\left( {{T}^{A}} \right)}^{2}}-{{\left( {{T}^{B}} \right)}^{2}} \right)B{{D}^{2}}-{{S}^{C}}{{S}^{D}}{{\left( {{T}^{A}}-{{T}^{B}} \right)}^{2}}C{{D}^{2}}  
\end{aligned} \right)\text{,}\]
\[{{E}_{A}}{{E}_{C}}^{2}=\frac{1}{{{\left( {{T}^{A}}{{T}^{C}} \right)}^{2}}}\left( \begin{aligned}
	& -{{S}^{A}}{{S}^{B}}\left( {{\left( {{T}^{A}} \right)}^{2}}-{{\left( {{T}^{C}} \right)}^{2}} \right)A{{B}^{2}}+{{S}^{A}}{{S}^{C}}{{\left( {{T}^{A}}+{{T}^{C}} \right)}^{2}}A{{C}^{2}} \\ 
	& -{{S}^{A}}{{S}^{D}}\left( {{\left( {{T}^{A}} \right)}^{2}}-{{\left( {{T}^{C}} \right)}^{2}} \right)A{{D}^{2}}+{{S}^{B}}{{S}^{C}}\left( {{\left( {{T}^{A}} \right)}^{2}}-{{\left( {{T}^{C}} \right)}^{2}} \right)B{{C}^{2}} \\ 
	& -{{S}^{B}}{{S}^{D}}{{\left( {{T}^{A}}-{{T}^{C}} \right)}^{2}}B{{D}^{2}}+{{S}^{C}}{{S}^{D}}\left( {{\left( {{T}^{A}} \right)}^{2}}-{{\left( {{T}^{C}} \right)}^{2}} \right)C{{D}^{2}}  
\end{aligned} \right)\text{,}\]
\[{{E}_{A}}{{E}_{D}}^{2}=\frac{1}{{{\left( {{T}^{A}}{{T}^{D}} \right)}^{2}}}\left( \begin{aligned}
	& -{{S}^{A}}{{S}^{B}}\left( {{\left( {{T}^{A}} \right)}^{2}}-{{\left( {{T}^{D}} \right)}^{2}} \right)A{{B}^{2}}-{{S}^{A}}{{S}^{C}}\left( {{\left( {{T}^{A}} \right)}^{2}}-{{\left( {{T}^{D}} \right)}^{2}} \right)A{{C}^{2}} \\ 
	& +{{S}^{A}}{{S}^{D}}{{\left( {{T}^{A}}+{{T}^{D}} \right)}^{2}}A{{D}^{2}}-{{S}^{B}}{{S}^{C}}{{\left( {{T}^{A}}-{{T}^{D}} \right)}^{2}}B{{C}^{2}} \\ 
	& +{{S}^{B}}{{S}^{D}}\left( {{\left( {{T}^{A}} \right)}^{2}}-{{\left( {{T}^{D}} \right)}^{2}} \right)B{{D}^{2}}+{{S}^{C}}{{S}^{D}}\left( {{\left( {{T}^{A}} \right)}^{2}}-{{\left( {{T}^{D}} \right)}^{2}} \right)C{{D}^{2}}  
\end{aligned} \right)\text{,}\]
\[{{E}_{B}}{{E}_{C}}^{2}=\frac{1}{{{\left( {{T}^{B}}{{T}^{C}} \right)}^{2}}}\left( \begin{aligned}
	& -{{S}^{A}}{{S}^{B}}\left( {{\left( {{T}^{B}} \right)}^{2}}-{{\left( {{T}^{C}} \right)}^{2}} \right)A{{B}^{2}}+{{S}^{A}}{{S}^{C}}\left( {{\left( {{T}^{B}} \right)}^{2}}-{{\left( {{T}^{C}} \right)}^{2}} \right)A{{C}^{2}} \\ 
	& -{{S}^{A}}{{S}^{D}}{{\left( {{T}^{B}}-{{T}^{C}} \right)}^{2}}A{{D}^{2}}+{{S}^{B}}{{S}^{C}}{{\left( {{T}^{B}}+{{T}^{C}} \right)}^{2}}B{{C}^{2}} \\ 
	& -{{S}^{B}}{{S}^{D}}\left( {{\left( {{T}^{B}} \right)}^{2}}-{{\left( {{T}^{C}} \right)}^{2}} \right)B{{D}^{2}}+{{S}^{C}}{{S}^{D}}\left( {{\left( {{T}^{B}} \right)}^{2}}-{{\left( {{T}^{C}} \right)}^{2}} \right)C{{D}^{2}}  
\end{aligned} \right)\text{,}\]
\[{{E}_{B}}{{E}_{D}}^{2}=\frac{1}{{{\left( {{T}^{B}}{{T}^{D}} \right)}^{2}}}\left( \begin{aligned}
	& -{{S}^{A}}{{S}^{B}}\left( {{\left( {{T}^{B}} \right)}^{2}}-{{\left( {{T}^{D}} \right)}^{2}} \right)A{{B}^{2}}-{{S}^{A}}{{S}^{C}}{{\left( {{T}^{B}}-{{T}^{D}} \right)}^{2}}A{{C}^{2}} \\ 
	& +{{S}^{A}}{{S}^{D}}\left( {{\left( {{T}^{B}} \right)}^{2}}-{{\left( {{T}^{D}} \right)}^{2}} \right)A{{D}^{2}}-{{S}^{B}}{{S}^{C}}\left( {{\left( {{T}^{B}} \right)}^{2}}-{{\left( {{T}^{D}} \right)}^{2}} \right)B{{C}^{2}} \\ 
	& +{{S}^{B}}{{S}^{D}}{{\left( {{T}^{B}}+{{T}^{D}} \right)}^{2}}B{{D}^{2}}+{{S}^{C}}{{S}^{D}}\left( {{\left( {{T}^{B}} \right)}^{2}}-{{\left( {{T}^{D}} \right)}^{2}} \right)C{{D}^{2}}  
\end{aligned} \right)\text{,}\]
\[{{E}_{B}}{{E}_{C}}^{2}=\frac{1}{{{\left( {{T}^{B}}{{T}^{C}} \right)}^{2}}}\left( \begin{aligned}
	& -{{S}^{A}}{{S}^{B}}\left( {{\left( {{T}^{B}} \right)}^{2}}-{{\left( {{T}^{C}} \right)}^{2}} \right)A{{B}^{2}}+{{S}^{A}}{{S}^{C}}\left( {{\left( {{T}^{B}} \right)}^{2}}-{{\left( {{T}^{C}} \right)}^{2}} \right)A{{C}^{2}} \\ 
	& -{{S}^{A}}{{S}^{D}}{{\left( {{T}^{B}}-{{T}^{C}} \right)}^{2}}A{{D}^{2}}+{{S}^{B}}{{S}^{C}}{{\left( {{T}^{B}}+{{T}^{C}} \right)}^{2}}B{{C}^{2}} \\ 
	& -{{S}^{B}}{{S}^{D}}\left( {{\left( {{T}^{B}} \right)}^{2}}-{{\left( {{T}^{C}} \right)}^{2}} \right)B{{D}^{2}}+{{S}^{C}}{{S}^{D}}\left( {{\left( {{T}^{B}} \right)}^{2}}-{{\left( {{T}^{C}} \right)}^{2}} \right)C{{D}^{2}}  
\end{aligned} \right).\]
\end{theorem}

\begin{proof}
	From theorem \ref{thm:Thm24.2.1}, let ${{P}_{1}}$ be the excenter ${{E}_{A}}$, ${{P}_{2}}$ be the excenter ${{E}_{B}}$, the following results are obtained:
	\[\begin{aligned}
		{{E}_{A}}{{E}_{B}}^{2}=& -\beta _{A}^{{{E}_{A}}{{E}_{B}}}\beta _{B}^{{{E}_{A}}{{E}_{B}}}A{{B}^{2}}-\beta _{A}^{{{E}_{A}}{{E}_{B}}}\beta _{C}^{{{E}_{A}}{{E}_{B}}}A{{C}^{2}}-\beta _{A}^{{{E}_{A}}{{E}_{B}}}\beta _{D}^{{{E}_{A}}{{E}_{B}}}A{{D}^{2}} \\ 
		& -\beta _{B}^{{{E}_{A}}{{E}_{B}}}\beta _{C}^{{{E}_{A}}{{E}_{B}}}B{{C}^{2}}-\beta _{B}^{{{E}_{A}}{{E}_{B}}}\beta _{D}^{{{E}_{A}}{{E}_{B}}}B{{D}^{2}}-\beta _{C}^{{{E}_{A}}{{E}_{B}}}\beta _{D}^{{{E}_{A}}{{E}_{B}}}C{{D}^{2}}.  
	\end{aligned}\]
	
	From theorem \ref{thm:Thm21.3.1}, we have
	\[\beta _{A}^{{{E}_{A}}{{E}_{B}}}=\beta _{A}^{{{E}_{B}}}-\beta _{A}^{{{E}_{A}}}=\frac{{{S}^{A}}}{S-2{{S}^{B}}}+\frac{{{S}^{A}}}{S-2{{S}^{A}}}=\frac{{{S}^{A}}\left( {{T}^{A}}+{{T}^{B}} \right)}{{{T}^{A}}{{T}^{B}}},\]	
	\[\beta _{B}^{{{E}_{A}}{{E}_{B}}}=\beta _{B}^{{{E}_{B}}}-\beta _{B}^{{{E}_{A}}}=-\frac{{{S}^{B}}}{S-2{{S}^{B}}}-\frac{{{S}^{B}}}{S-2{{S}^{A}}}=-\frac{{{S}^{B}}\left( {{T}^{A}}+{{T}^{B}} \right)}{{{T}^{A}}{{T}^{B}}},\]	\[\beta _{C}^{{{E}_{A}}{{E}_{B}}}=\beta _{C}^{{{E}_{B}}}-\beta _{C}^{{{E}_{A}}}=\frac{{{S}^{C}}}{S-2{{S}^{B}}}-\frac{{{S}^{C}}}{S-2{{S}^{A}}}=\frac{{{S}^{C}}\left( {{T}^{A}}-{{T}^{B}} \right)}{{{T}^{A}}{{T}^{B}}},\]	
	\[\beta _{D}^{{{E}_{A}}{{E}_{B}}}=\beta _{D}^{{{E}_{B}}}-\beta _{D}^{{{E}_{A}}}=\frac{{{S}^{D}}}{S-2{{S}^{B}}}-\frac{{{S}^{D}}}{S-2{{S}^{A}}}=\frac{{{S}^{D}}\left( {{T}^{A}}-{{T}^{B}} \right)}{{{T}^{A}}{{T}^{B}}}\text{.}\]
	
	Therefore
	\[{{E}_{A}}{{E}_{B}}^{2}=\frac{1}{{{\left( {{T}^{A}}{{T}^{B}} \right)}^{2}}}\left( \begin{aligned}
		& {{S}^{A}}{{S}^{B}}{{\left( {{T}^{A}}+{{T}^{B}} \right)}^{2}}A{{B}^{2}}-{{S}^{A}}{{S}^{C}}\left( {{\left( {{T}^{A}} \right)}^{2}}-{{\left( {{T}^{B}} \right)}^{2}} \right)A{{C}^{2}} \\ 
		& -{{S}^{A}}{{S}^{D}}\left( {{\left( {{T}^{A}} \right)}^{2}}-{{\left( {{T}^{B}} \right)}^{2}} \right)A{{D}^{2}}+{{S}^{B}}{{S}^{C}}\left( {{\left( {{T}^{A}} \right)}^{2}}-{{\left( {{T}^{B}} \right)}^{2}} \right)B{{C}^{2}} \\ 
		& +{{S}^{B}}{{S}^{D}}\left( {{\left( {{T}^{A}} \right)}^{2}}-{{\left( {{T}^{B}} \right)}^{2}} \right)B{{D}^{2}}-{{S}^{C}}{{S}^{D}}{{\left( {{T}^{A}}-{{T}^{B}} \right)}^{2}}C{{D}^{2}}  
	\end{aligned} \right)\text{.}\]
	
	Let ${{P}_{1}}$ be the excenter ${{E}_{A}}$, ${{P}_{2}}$ be the excenter ${{E}_{C}}$, the following results are obtained:
	\[\begin{aligned}
		{{E}_{A}}{{E}_{C}}^{2}=& -\beta _{A}^{{{E}_{A}}{{E}_{C}}}\beta _{B}^{{{E}_{A}}{{E}_{C}}}A{{B}^{2}}-\beta _{A}^{{{E}_{A}}{{E}_{C}}}\beta _{C}^{{{E}_{A}}{{E}_{C}}}A{{C}^{2}}-\beta _{A}^{{{E}_{A}}{{E}_{C}}}\beta _{D}^{{{E}_{A}}{{E}_{C}}}A{{D}^{2}} \\ 
		& -\beta _{B}^{{{E}_{A}}{{E}_{C}}}\beta _{C}^{{{E}_{A}}{{E}_{C}}}B{{C}^{2}}-\beta _{B}^{{{E}_{A}}{{E}_{C}}}\beta _{D}^{{{E}_{A}}{{E}_{C}}}B{{D}^{2}}-\beta _{C}^{{{E}_{A}}{{E}_{C}}}\beta _{D}^{{{E}_{A}}{{E}_{C}}}C{{D}^{2}}\text{.}  
	\end{aligned}\]
	
	From theorem \ref{thm:Thm21.3.1}, we have
	\[\beta _{A}^{{{E}_{A}}{{E}_{C}}}=\beta _{A}^{{{E}_{C}}}-\beta _{A}^{{{E}_{A}}}=\frac{{{S}^{A}}}{S-2{{S}^{C}}}+\frac{{{S}^{A}}}{S-2{{S}^{A}}}=\frac{{{S}^{A}}\left( {{T}^{A}}+{{T}^{C}} \right)}{{{T}^{A}}{{T}^{C}}},\]	
	\[\beta _{B}^{{{E}_{A}}{{E}_{C}}}=\beta _{B}^{{{E}_{C}}}-\beta _{B}^{{{E}_{A}}}=\frac{{{S}^{B}}}{S-2{{S}^{C}}}-\frac{{{S}^{B}}}{S-2{{S}^{A}}}=\frac{{{S}^{B}}\left( {{T}^{A}}-{{T}^{C}} \right)}{{{T}^{A}}{{T}^{C}}},\]	
	\[\beta _{C}^{{{E}_{A}}{{E}_{C}}}=\beta _{C}^{{{E}_{C}}}-\beta _{C}^{{{E}_{A}}}=-\frac{{{S}^{C}}}{S-2{{S}^{C}}}-\frac{{{S}^{C}}}{S-2{{S}^{A}}}=-\frac{{{S}^{C}}\left( {{T}^{A}}+{{T}^{C}} \right)}{{{T}^{A}}{{T}^{C}}},\]	\[\beta _{D}^{{{E}_{A}}{{E}_{C}}}=\beta _{D}^{{{E}_{C}}}-\beta _{D}^{{{E}_{A}}}=\frac{{{S}^{D}}}{S-2{{S}^{C}}}-\frac{{{S}^{D}}}{S-2{{S}^{A}}}=\frac{{{S}^{D}}\left( {{T}^{A}}-{{T}^{C}} \right)}{{{T}^{A}}{{T}^{C}}}\text{.}\]
	
	Therefore
	\[{{E}_{A}}{{E}_{C}}^{2}=\frac{1}{{{\left( {{T}^{A}}{{T}^{C}} \right)}^{2}}}\left( \begin{aligned}
		& -{{S}^{A}}{{S}^{B}}\left( {{\left( {{T}^{A}} \right)}^{2}}-{{\left( {{T}^{C}} \right)}^{2}} \right)A{{B}^{2}}+{{S}^{A}}{{S}^{C}}{{\left( {{T}^{A}}+{{T}^{C}} \right)}^{2}}A{{C}^{2}} \\ 
		& -{{S}^{A}}{{S}^{D}}\left( {{\left( {{T}^{A}} \right)}^{2}}-{{\left( {{T}^{C}} \right)}^{2}} \right)A{{D}^{2}}+{{S}^{B}}{{S}^{C}}\left( {{\left( {{T}^{A}} \right)}^{2}}-{{\left( {{T}^{C}} \right)}^{2}} \right)B{{C}^{2}} \\ 
		& -{{S}^{B}}{{S}^{D}}{{\left( {{T}^{A}}-{{T}^{C}} \right)}^{2}}B{{D}^{2}}+{{S}^{C}}{{S}^{D}}\left( {{\left( {{T}^{A}} \right)}^{2}}-{{\left( {{T}^{C}} \right)}^{2}} \right)C{{D}^{2}}  
	\end{aligned} \right)\text{.}\]
	
	Let ${{P}_{1}}$ be the excenter ${{E}_{A}}$, ${{P}_{2}}$ be the excenter ${{E}_{D}}$, the following results are obtained:
	\[\begin{aligned}
		{{E}_{A}}{{E}_{D}}^{2}=& -\beta _{A}^{{{E}_{A}}{{E}_{D}}}\beta _{B}^{{{E}_{A}}{{E}_{D}}}A{{B}^{2}}-\beta _{A}^{{{E}_{A}}{{E}_{D}}}\beta _{C}^{{{E}_{A}}{{E}_{D}}}A{{C}^{2}}-\beta _{A}^{{{E}_{A}}{{E}_{D}}}\beta _{D}^{{{E}_{A}}{{E}_{D}}}A{{D}^{2}} \\ 
		& -\beta _{B}^{{{E}_{A}}{{E}_{D}}}\beta _{C}^{{{E}_{A}}{{E}_{D}}}B{{C}^{2}}-\beta _{B}^{{{E}_{A}}{{E}_{D}}}\beta _{D}^{{{E}_{A}}{{E}_{D}}}B{{D}^{2}}-\beta _{C}^{{{E}_{A}}{{E}_{D}}}\beta _{D}^{{{E}_{A}}{{E}_{D}}}C{{D}^{2}}\text{.}  
	\end{aligned}\]
	
	From theorem \ref{thm:Thm21.3.1}, we have
	\[\beta _{A}^{{{E}_{A}}{{E}_{D}}}=\beta _{A}^{{{E}_{D}}}-\beta _{A}^{{{E}_{A}}}=\frac{{{S}^{A}}}{S-2{{S}^{D}}}+\frac{{{S}^{A}}}{S-2{{S}^{A}}}=\frac{{{S}^{A}}\left( {{T}^{A}}+{{T}^{D}} \right)}{{{T}^{A}}{{T}^{D}}},\]	
	\[\beta _{B}^{{{E}_{A}}{{E}_{D}}}=\beta _{B}^{{{E}_{D}}}-\beta _{B}^{{{E}_{A}}}=\frac{{{S}^{B}}}{S-2{{S}^{D}}}-\frac{{{S}^{B}}}{S-2{{S}^{A}}}=\frac{{{S}^{B}}\left( {{T}^{A}}-{{T}^{D}} \right)}{{{T}^{A}}{{T}^{D}}},\]
	\[\beta _{C}^{{{E}_{A}}{{E}_{D}}}=\beta _{C}^{{{E}_{D}}}-\beta _{C}^{{{E}_{A}}}=\frac{{{S}^{C}}}{S-2{{S}^{D}}}-\frac{{{S}^{C}}}{S-2{{S}^{A}}}=\frac{{{S}^{C}}\left( {{T}^{A}}-{{T}^{D}} \right)}{{{T}^{A}}{{T}^{D}}},\]	
	\[\beta _{D}^{{{E}_{A}}{{E}_{D}}}=\beta _{D}^{{{E}_{D}}}-\beta _{D}^{{{E}_{A}}}=-\frac{{{S}^{D}}}{S-2{{S}^{D}}}-\frac{{{S}^{D}}}{S-2{{S}^{A}}}=-\frac{{{S}^{D}}\left( {{T}^{A}}+{{T}^{D}} \right)}{{{T}^{A}}{{T}^{D}}}\text{.}\]
	
	Therefore
	\[{{E}_{A}}{{E}_{D}}^{2}=\frac{1}{{{\left( {{T}^{A}}{{T}^{D}} \right)}^{2}}}\left( \begin{aligned}
		& -{{S}^{A}}{{S}^{B}}\left( {{\left( {{T}^{A}} \right)}^{2}}-{{\left( {{T}^{D}} \right)}^{2}} \right)A{{B}^{2}}-{{S}^{A}}{{S}^{C}}\left( {{\left( {{T}^{A}} \right)}^{2}}-{{\left( {{T}^{D}} \right)}^{2}} \right)A{{C}^{2}} \\ 
		& +{{S}^{A}}{{S}^{D}}{{\left( {{T}^{A}}+{{T}^{D}} \right)}^{2}}A{{D}^{2}}-{{S}^{B}}{{S}^{C}}{{\left( {{T}^{A}}-{{T}^{D}} \right)}^{2}}B{{C}^{2}} \\ 
		& +{{S}^{B}}{{S}^{D}}\left( {{\left( {{T}^{A}} \right)}^{2}}-{{\left( {{T}^{D}} \right)}^{2}} \right)B{{D}^{2}}+{{S}^{C}}{{S}^{D}}\left( {{\left( {{T}^{A}} \right)}^{2}}-{{\left( {{T}^{D}} \right)}^{2}} \right)C{{D}^{2}}  
	\end{aligned} \right)\text{.}\]
	
	Let ${{P}_{1}}$ be the excenter ${{E}_{B}}$, ${{P}_{2}}$ be the excenter ${{E}_{C}}$, the following results are obtained:
	\[\begin{aligned}
		{{E}_{B}}{{E}_{C}}^{2}=& -\beta _{A}^{{{E}_{B}}{{E}_{C}}}\beta _{B}^{{{E}_{B}}{{E}_{C}}}A{{B}^{2}}-\beta _{A}^{{{E}_{B}}{{E}_{C}}}\beta _{C}^{{{E}_{B}}{{E}_{C}}}A{{C}^{2}}-\beta _{A}^{{{E}_{B}}{{E}_{C}}}\beta _{D}^{{{E}_{B}}{{E}_{C}}}A{{D}^{2}} \\ 
		& -\beta _{B}^{{{E}_{B}}{{E}_{C}}}\beta _{C}^{{{E}_{B}}{{E}_{C}}}B{{C}^{2}}-\beta _{B}^{{{E}_{B}}{{E}_{C}}}\beta _{D}^{{{E}_{B}}{{E}_{C}}}B{{D}^{2}}-\beta _{C}^{{{E}_{B}}{{E}_{C}}}\beta _{D}^{{{E}_{B}}{{E}_{C}}}C{{D}^{2}}\text{.}  
	\end{aligned}\]
	
	From theorem \ref{thm:Thm21.3.1}, we have
	\[\beta _{A}^{{{E}_{B}}{{E}_{C}}}=\beta _{A}^{{{E}_{C}}}-\beta _{A}^{{{E}_{B}}}=\frac{{{S}^{A}}}{S-2{{S}^{C}}}-\frac{{{S}^{A}}}{S-2{{S}^{B}}}=\frac{{{S}^{A}}\left( {{T}^{B}}-{{T}^{C}} \right)}{{{T}^{B}}{{T}^{C}}},\]	
	\[\beta _{B}^{{{E}_{B}}{{E}_{C}}}=\beta _{B}^{{{E}_{C}}}-\beta _{B}^{{{E}_{B}}}=\frac{{{S}^{B}}}{S-2{{S}^{C}}}+\frac{{{S}^{B}}}{S-2{{S}^{B}}}=\frac{{{S}^{B}}\left( {{T}^{B}}+{{T}^{C}} \right)}{{{T}^{B}}{{T}^{C}}},\]	
	\[\beta _{C}^{{{E}_{B}}{{E}_{C}}}=\beta _{C}^{{{E}_{C}}}-\beta _{C}^{{{E}_{B}}}=-\frac{{{S}^{C}}}{S-2{{S}^{C}}}-\frac{{{S}^{C}}}{S-2{{S}^{B}}}=-\frac{{{S}^{C}}\left( {{T}^{B}}+{{T}^{C}} \right)}{{{T}^{B}}{{T}^{C}}},\]	\[\beta _{D}^{{{E}_{B}}{{E}_{C}}}=\beta _{D}^{{{E}_{C}}}-\beta _{D}^{{{E}_{B}}}=\frac{{{S}^{D}}}{S-2{{S}^{C}}}-\frac{{{S}^{D}}}{S-2{{S}^{B}}}=\frac{{{S}^{D}}\left( {{T}^{B}}-{{T}^{C}} \right)}{{{T}^{B}}{{T}^{C}}}\text{.}\]
	
	Therefore
	\[{{E}_{B}}{{E}_{C}}^{2}=\frac{1}{{{\left( {{T}^{B}}{{T}^{C}} \right)}^{2}}}\left( \begin{aligned}
		& -{{S}^{A}}{{S}^{B}}\left( {{\left( {{T}^{B}} \right)}^{2}}-{{\left( {{T}^{C}} \right)}^{2}} \right)A{{B}^{2}}+{{S}^{A}}{{S}^{C}}\left( {{\left( {{T}^{B}} \right)}^{2}}-{{\left( {{T}^{C}} \right)}^{2}} \right)A{{C}^{2}} \\ 
		& -{{S}^{A}}{{S}^{D}}{{\left( {{T}^{B}}-{{T}^{C}} \right)}^{2}}A{{D}^{2}}+{{S}^{B}}{{S}^{C}}{{\left( {{T}^{B}}+{{T}^{C}} \right)}^{2}}B{{C}^{2}} \\ 
		& -{{S}^{B}}{{S}^{D}}\left( {{\left( {{T}^{B}} \right)}^{2}}-{{\left( {{T}^{C}} \right)}^{2}} \right)B{{D}^{2}}+{{S}^{C}}{{S}^{D}}\left( {{\left( {{T}^{B}} \right)}^{2}}-{{\left( {{T}^{C}} \right)}^{2}} \right)C{{D}^{2}}  
	\end{aligned} \right)\text{.}\]
	
	Let ${{P}_{1}}$ be the excenter ${{E}_{B}}$, ${{P}_{2}}$ be the excenter ${{E}_{D}}$, the following results are obtained:
	\[\begin{aligned}
		{{E}_{B}}{{E}_{D}}^{2}=& -\beta _{A}^{{{E}_{B}}{{E}_{D}}}\beta _{B}^{{{E}_{B}}{{E}_{D}}}A{{B}^{2}}-\beta _{A}^{{{E}_{B}}{{E}_{D}}}\beta _{C}^{{{E}_{B}}{{E}_{D}}}A{{C}^{2}}-\beta _{A}^{{{E}_{B}}{{E}_{D}}}\beta _{D}^{{{E}_{B}}{{E}_{D}}}A{{D}^{2}} \\ 
		& -\beta _{B}^{{{E}_{B}}{{E}_{D}}}\beta _{C}^{{{E}_{B}}{{E}_{D}}}B{{C}^{2}}-\beta _{B}^{{{E}_{B}}{{E}_{D}}}\beta _{D}^{{{E}_{B}}{{E}_{D}}}B{{D}^{2}}-\beta _{C}^{{{E}_{B}}{{E}_{D}}}\beta _{D}^{{{E}_{B}}{{E}_{D}}}C{{D}^{2}}\text{.}  
	\end{aligned}\]
	
	From theorem \ref{thm:Thm21.3.1}, we have
	\[\beta _{A}^{{{E}_{B}}{{E}_{D}}}=\beta _{A}^{{{E}_{D}}}-\beta _{A}^{{{E}_{B}}}=\frac{{{S}^{A}}}{S-2{{S}^{D}}}-\frac{{{S}^{A}}}{S-2{{S}^{B}}}=\frac{{{S}^{A}}\left( {{T}^{B}}-{{T}^{D}} \right)}{{{T}^{B}}{{T}^{D}}},\]	
	\[\beta _{B}^{{{E}_{B}}{{E}_{D}}}=\beta _{B}^{{{E}_{D}}}-\beta _{B}^{{{E}_{B}}}=\frac{{{S}^{B}}}{S-2{{S}^{D}}}+\frac{{{S}^{B}}}{S-2{{S}^{B}}}=\frac{{{S}^{B}}\left( {{T}^{B}}+{{T}^{D}} \right)}{{{T}^{B}}{{T}^{D}}},\]	
	\[\beta _{C}^{{{E}_{B}}{{E}_{D}}}=\beta _{C}^{{{E}_{D}}}-\beta _{C}^{{{E}_{B}}}=\frac{{{S}^{C}}}{S-2{{S}^{D}}}-\frac{{{S}^{C}}}{S-2{{S}^{B}}}=\frac{{{S}^{C}}\left( {{T}^{B}}-{{T}^{D}} \right)}{{{T}^{B}}{{T}^{D}}},\]	
	\[\beta _{D}^{{{E}_{B}}{{E}_{D}}}=\beta _{D}^{{{E}_{D}}}-\beta _{D}^{{{E}_{B}}}=-\frac{{{S}^{D}}}{S-2{{S}^{D}}}-\frac{{{S}^{D}}}{S-2{{S}^{B}}}=-\frac{{{S}^{D}}\left( {{T}^{B}}+{{T}^{D}} \right)}{{{T}^{B}}{{T}^{D}}}\text{.}\]
	
	Therefore
	\[{{E}_{B}}{{E}_{D}}^{2}=\frac{1}{{{\left( {{T}^{B}}{{T}^{D}} \right)}^{2}}}\left( \begin{aligned}
		& -{{S}^{A}}{{S}^{B}}\left( {{\left( {{T}^{B}} \right)}^{2}}-{{\left( {{T}^{D}} \right)}^{2}} \right)A{{B}^{2}}-{{S}^{A}}{{S}^{C}}{{\left( {{T}^{B}}-{{T}^{D}} \right)}^{2}}A{{C}^{2}} \\ 
		& +{{S}^{A}}{{S}^{D}}\left( {{\left( {{T}^{B}} \right)}^{2}}-{{\left( {{T}^{D}} \right)}^{2}} \right)A{{D}^{2}}-{{S}^{B}}{{S}^{C}}\left( {{\left( {{T}^{B}} \right)}^{2}}-{{\left( {{T}^{D}} \right)}^{2}} \right)B{{C}^{2}} \\ 
		& +{{S}^{B}}{{S}^{D}}{{\left( {{T}^{B}}+{{T}^{D}} \right)}^{2}}B{{D}^{2}}+{{S}^{C}}{{S}^{D}}\left( {{\left( {{T}^{B}} \right)}^{2}}-{{\left( {{T}^{D}} \right)}^{2}} \right)C{{D}^{2}}  
	\end{aligned} \right).\]
	
	Let ${{P}_{1}}$ be the excenter ${{E}_{C}}$, ${{P}_{2}}$ be the excenter ${{E}_{D}}$, the following results are obtained: 
	\[\begin{aligned}
		{{E}_{C}}{{E}_{D}}^{2}=& -\beta _{A}^{{{E}_{C}}{{E}_{D}}}\beta _{B}^{{{E}_{C}}{{E}_{D}}}A{{B}^{2}}-\beta _{A}^{{{E}_{C}}{{E}_{D}}}\beta _{C}^{{{E}_{C}}{{E}_{D}}}A{{C}^{2}}-\beta _{A}^{{{E}_{C}}{{E}_{D}}}\beta _{D}^{{{E}_{C}}{{E}_{D}}}A{{D}^{2}} \\ 
		& -\beta _{B}^{{{E}_{C}}{{E}_{D}}}\beta _{C}^{{{E}_{C}}{{E}_{D}}}B{{C}^{2}}-\beta _{B}^{{{E}_{C}}{{E}_{D}}}\beta _{D}^{{{E}_{C}}{{E}_{D}}}B{{D}^{2}}-\beta _{C}^{{{E}_{C}}{{E}_{D}}}\beta _{D}^{{{E}_{C}}{{E}_{D}}}C{{D}^{2}}\text{.}  
	\end{aligned}\]
	
	From theorem \ref{thm:Thm21.3.1}, we have
	\[\beta _{A}^{{{E}_{C}}{{E}_{D}}}=\beta _{A}^{{{E}_{D}}}-\beta _{A}^{{{E}_{C}}}=\frac{{{S}^{A}}}{S-2{{S}^{D}}}-\frac{{{S}^{A}}}{S-2{{S}^{C}}}=\frac{{{S}^{A}}\left( {{T}^{C}}-{{T}^{D}} \right)}{{{T}^{C}}{{T}^{D}}},\]	
	\[\beta _{B}^{{{E}_{C}}{{E}_{D}}}=\beta _{B}^{{{E}_{D}}}-\beta _{B}^{{{E}_{C}}}=\frac{{{S}^{B}}}{S-2{{S}^{D}}}-\frac{{{S}^{B}}}{S-2{{S}^{C}}}=\frac{{{S}^{B}}\left( {{T}^{C}}-{{T}^{D}} \right)}{{{T}^{C}}{{T}^{D}}},\]	
	\[\beta _{C}^{{{E}_{C}}{{E}_{D}}}=\beta _{C}^{{{E}_{D}}}-\beta _{C}^{{{E}_{C}}}=\frac{{{S}^{C}}}{S-2{{S}^{D}}}+\frac{{{S}^{C}}}{S-2{{S}^{C}}}=\frac{{{S}^{C}}\left( {{T}^{C}}+{{T}^{D}} \right)}{{{T}^{C}}{{T}^{D}}},\]	
	\[\beta _{D}^{{{E}_{C}}{{E}_{D}}}=\beta _{D}^{{{E}_{D}}}-\beta _{D}^{{{E}_{C}}}=-\frac{{{S}^{D}}}{S-2{{S}^{D}}}-\frac{{{S}^{D}}}{S-2{{S}^{C}}}=-\frac{{{S}^{D}}\left( {{T}^{C}}+{{T}^{D}} \right)}{{{T}^{C}}{{T}^{D}}}\text{.}\]
	
	Therefore
	\[{{E}_{C}}{{E}_{D}}^{2}=\frac{1}{{{\left( {{T}^{C}}{{T}^{D}} \right)}^{2}}}\left( \begin{aligned}
		& -{{S}^{A}}{{S}^{B}}{{\left( {{T}^{C}}-{{T}^{D}} \right)}^{2}}A{{B}^{2}}-{{S}^{A}}{{S}^{C}}{{\left( {{\left( {{T}^{C}} \right)}^{2}}-{{\left( {{T}^{D}} \right)}^{2}} \right)}^{2}}A{{C}^{2}} \\ 
		& +{{S}^{A}}{{S}^{D}}\left( {{\left( {{T}^{C}} \right)}^{2}}-{{\left( {{T}^{D}} \right)}^{2}} \right)A{{D}^{2}}-{{S}^{B}}{{S}^{C}}\left( {{\left( {{T}^{C}} \right)}^{2}}-{{\left( {{T}^{D}} \right)}^{2}} \right)B{{C}^{2}} \\ 
		& +{{S}^{B}}{{S}^{D}}\left( {{\left( {{T}^{C}} \right)}^{2}}-{{\left( {{T}^{D}} \right)}^{2}} \right)B{{D}^{2}}+{{S}^{C}}{{S}^{D}}{{\left( {{T}^{C}}+{{T}^{D}} \right)}^{2}}C{{D}^{2}}  
	\end{aligned} \right)\text{.}\]
\end{proof}
\hfill $\square$\par


\chapter{New inequalities of tetrahedron}\label{Ch27}
\thispagestyle{empty}

%
%
In addition to the previous achievements, some new inequalities can be obtained in Intercenter Geometry. In this chapter, I publish some new inequalities on tetrahedron for the first time.

Intercenter Geometry is good at calculating the distance between two points. Many new geometric inequalities can be created based on the principle that the distance between two points is non-negative. It can be said that Intercenter Geometry is a “generator” that creates a large class of geometric inequalities. In Space Intercenter Geometry, geometric inequality mainly involves geometric inequalities on tetrahedron.

Since the distance between the intersecting centers of a tetrahedron is non-negative, the following new geometric inequalities can be derived directly from the theorem of distance between two intersecting centers on tetrahedral frame (theorem \ref{thm:Thm24.2.1}), and the theorem of distance between origin and intersecting center on tetrahedral frame (theorem \ref{thm:Thm24.1.4}).

This chapter is divided into two sections, the geometric inequality of DOIC-T and the geometric inequality of DTICs-T.


\section{Geometric inequality of distance between origin and intersecting center of tetrahedron}\label{Sec27.1} 

\begin{theorem}{Geometric inequality of distance between origin and IC-T, Daiyuan Zhang}{Thm27.1.1}\label{Thm27.1.1} 
	Given a tetrahedron $ABCD$, let point $\,O$ be the origin of frame $\left( O;A,B,C,D \right)$, point $P$ be the intersecting center of the tetrahedron, then	
	\[\begin{aligned}
		& \beta _{A}^{P}O{{A}^{2}}+\beta _{B}^{P}O{{B}^{2}}+\beta _{C}^{P}O{{C}^{2}}+\beta _{D}^{P}O{{D}^{2}} \\ 
		& \ge \beta _{A}^{P}\beta _{B}^{P}A{{B}^{2}}+\beta _{A}^{P}\beta _{C}^{P}A{{C}^{2}}+\beta _{A}^{P}\beta _{D}^{P}A{{D}^{2}} \\ 
		& +\beta _{B}^{P}\beta _{C}^{P}B{{C}^{2}}+\beta _{B}^{P}\beta _{D}^{P}B{{D}^{2}}+\beta _{C}^{P}\beta _{D}^{P}C{{D}^{2}}.  
	\end{aligned}\]	

	The above equality occurs if and only if the origin $O$ of the frame coincides with the intersecting center $P$ of the tetrahedron, where $\beta _{A}^{P}$, $\beta _{B}^{P}$, $\beta _{C}^{P}$, $\beta _{D}^{P}$ are the frame components of $\overrightarrow{OA}$, $\overrightarrow{OB}$, $\overrightarrow{OC}$, $\overrightarrow{OD}$ at $P$ on the frme $\left( O;A,B,C,D \right)$.
\end{theorem}

\begin{proof}
	From theorem \ref{thm:Thm24.1.1}, 	
	\[\begin{aligned}
		O{{P}^{2}}&=\beta _{A}^{P}O{{A}^{2}}+\beta _{B}^{P}O{{B}^{2}}+\beta _{C}^{P}O{{C}^{2}}+\beta _{D}^{P}O{{D}^{2}} \\ 
		& -\beta _{A}^{P}\beta _{B}^{P}A{{B}^{2}}-\beta _{A}^{P}\beta _{C}^{P}A{{C}^{2}}-\beta _{A}^{P}\beta _{D}^{P}A{{D}^{2}} \\ 
		& -\beta _{B}^{P}\beta _{C}^{P}B{{C}^{2}}-\beta _{B}^{P}\beta _{D}^{P}B{{D}^{2}}-\beta _{C}^{P}\beta _{D}^{P}C{{D}^{2}}\ge 0. \\ 
	\end{aligned}\]	
	
	Since $O{{P}^{2}}\ge 0$, it follows that
	\[\begin{aligned}
		& \beta _{A}^{P}O{{A}^{2}}+\beta _{B}^{P}O{{B}^{2}}+\beta _{C}^{P}O{{C}^{2}}+\beta _{D}^{P}O{{D}^{2}} \\ 
		& \ge \beta _{A}^{P}\beta _{B}^{P}A{{B}^{2}}+\beta _{A}^{P}\beta _{C}^{P}A{{C}^{2}}+\beta _{A}^{P}\beta _{D}^{P}A{{D}^{2}} \\ 
		& +\beta _{B}^{P}\beta _{C}^{P}B{{C}^{2}}+\beta _{B}^{P}\beta _{D}^{P}B{{D}^{2}}+\beta _{C}^{P}\beta _{D}^{P}C{{D}^{2}}.  
	\end{aligned}\]
	
	Obviously, the above equality holds if and only if $OP=0$, that is, if and only if the origin $O$ of the frame coincides with the intersecting center $P$ of the tetrahedron.
\end{proof}
\hfill $\square$\par

\begin{theorem}{Inequality of radius of circumscribed sphere-IC-T, Daiyuan Zhang}{Thm27.1.2}\label{Thm27.1.2} 
	Given a tetrahedron $ABCD$, let point $\,O$ be the origin of frame $\left( O;A,B,C,D \right)$, point $P$ be the intersecting center of the tetrahedron, $R$ be the radius of the circumscribed sphere of tetrahedron $ABCD$, then 	
	\[{{R}^{2}}\ge \left( \begin{aligned}
		& \beta _{A}^{P}\beta _{B}^{P}A{{B}^{2}}+\beta _{A}^{P}\beta _{C}^{P}A{{C}^{2}}+\beta _{A}^{P}\beta _{D}^{P}A{{D}^{2}} \\ 
		& +\beta _{B}^{P}\beta _{C}^{P}B{{C}^{2}}+\beta _{B}^{P}\beta _{D}^{P}B{{D}^{2}}+\beta _{C}^{P}\beta _{D}^{P}C{{D}^{2}}  
	\end{aligned} \right).\]	
	
	The above equality occurs if and only if the circumcenter (origin of the frame) $Q$ coincides with the intersecting center $P$ of the tetrahedron, where $\beta _{A}^{P}$, $\beta _{B}^{P}$, $\beta _{C}^{P}$, $\beta _{D}^{P}$ are the frame components of $\overrightarrow{OA}$, $\overrightarrow{OB}$, $\overrightarrow{OC}$, $\overrightarrow{OD}$ at $P$ on the frme $\left( O;A,B,C,D \right)$.
\end{theorem}

\begin{proof}
	From the theorem of distance between origin and intersecting center of a tetrahedron (theorem \ref{thm:Thm24.1.4}), then	
	\[Q{{P}^{2}}={{R}^{2}}-\left( \begin{aligned}
		& \beta _{A}^{P}\beta _{B}^{P}A{{B}^{2}}+\beta _{A}^{P}\beta _{C}^{P}A{{C}^{2}}+\beta _{A}^{P}\beta _{D}^{P}A{{D}^{2}} \\ 
		& +\beta _{B}^{P}\beta _{C}^{P}B{{C}^{2}}+\beta _{B}^{P}\beta _{D}^{P}B{{D}^{2}}+\beta _{C}^{P}\beta _{D}^{P}C{{D}^{2}}  
	\end{aligned} \right).\]
	
	Since $Q{{P}^{2}}\ge 0$, then	
	\[{{R}^{2}}\ge \left( \begin{aligned}
		& \beta _{A}^{P}\beta _{B}^{P}A{{B}^{2}}+\beta _{A}^{P}\beta _{C}^{P}A{{C}^{2}}+\beta _{A}^{P}\beta _{D}^{P}A{{D}^{2}} \\ 
		& +\beta _{B}^{P}\beta _{C}^{P}B{{C}^{2}}+\beta _{B}^{P}\beta _{D}^{P}B{{D}^{2}}+\beta _{C}^{P}\beta _{D}^{P}C{{D}^{2}}  
	\end{aligned} \right).\]
	
	Obviously, the above equality holds if and only if $QP=0$, that is, if and only if the circumcenter $Q$ of the frame coincides with the intersecting center $P$ of the tetrahedron.
\end{proof}
\hfill $\square$\par

When the above theorem is applied to some special intersecting centers of a tetrahedron, such as circumcenter, centroid, and incenter, a series of inequalities will be obtained. According to Euclidean geometry, the circumcenter, centroid, incenter will coincide with each other If and only if the tetrahedron is a regular tetrahedron. So we get the following theorems.

\begin{theorem}{Inequality of radius-centroid of circumscribed sphere, Daiyuan Zhang}{Thm27.1.3}\label{Thm27.1.3} 
	Suppose that given a tetrahedron $ABCD$, $R$ is the radius of circumscribed sphere of tetrahedron $ABCD$, then
	\begin{equation}\label{Eq27.1.1}
		{{R}^{2}}\ge \frac{1}{16}\left( A{{B}^{2}}+A{{C}^{2}}+A{{D}^{2}}+B{{C}^{2}}+B{{D}^{2}}+C{{D}^{2}} \right).
	\end{equation}
		
	The above equality occurs if and only if the circumcenter coincides with the centroid, where $AB$, $AC$, $AD$, $BC$, $BD$, $CD$ are the lengths of the six edges of the tetrahedron, respectively.
\end{theorem}

\begin{proof}
	It is directly obtained from theorem \ref{thm:Thm26.1.1}:	
	\[QG=\frac{1}{4}\sqrt{16{{R}^{2}}-\left( A{{B}^{2}}+A{{C}^{2}}+A{{D}^{2}}+B{{C}^{2}}+B{{D}^{2}}+C{{D}^{2}} \right)}.\]
	
	Since $QG\ge 0$, then	
	\[{{R}^{2}}\ge \frac{1}{16}\left( A{{B}^{2}}+A{{C}^{2}}+A{{D}^{2}}+B{{C}^{2}}+B{{D}^{2}}+C{{D}^{2}} \right).\]	
	
	The above equality holds if and only if $QG=0$, i.e., if and only if the circumcenter coincides with the centroid.
\end{proof}
\hfill $\square$\par

It is difficult to find that it is related to the distance between the radius of the circumscribed sphere and the centroid of the tetrahedron from the inequality (\ref{Eq27.1.1}). I call this inequality “Inequality of radius of circumscribed sphere-centroid of a tetrahedron” to remind readers not to forget their origin. In fact, the inequality (\ref{Eq27.1.1}) and some inequalities to be given later in this section are derived from theorem \ref{thm:Thm27.1.2}.

Here's an example. For a triangular pyramid $A-BCD$, suppose that each length of lateral edges is $AB=AC=AD=l$, and each length of base edges is $BC=BD=CD=a$, then the inequality of circumscribed sphere radius is obtained according to the above formula:
\[R\ge \frac{\sqrt{3}}{4}\sqrt{{{l}^{2}}+{{a}^{2}}}.\]	

\begin{theorem}{Inequality of radius-incenter of circumscribed sphere of tetrahedron, Daiyuan Zhang}{Thm27.1.4}\label{Thm27.1.4} 
	Suppose that given a tetrahedron $ABCD$, $R$ is the radius of circumscribed sphere of tetrahedron $ABCD$, then
	\[{{R}^{2}}\ge \frac{1}{{{S}^{2}}}\left( \begin{aligned}
		& {{S}^{A}}{{S}^{B}}A{{B}^{2}}+{{S}^{A}}{{S}^{C}}A{{C}^{2}}+{{S}^{A}}{{S}^{C}}A{{D}^{2}} \\ 
		& +{{S}^{B}}{{S}^{C}}B{{C}^{2}}+{{S}^{B}}{{S}^{D}}B{{D}^{2}}+{{S}^{C}}{{S}^{D}}C{{D}^{2}}  
	\end{aligned} \right).\]	
	
	The above equality occurs if and only if the circumcenter coincides with the incenter. Where ${{S}^{A}}$, ${{S}^{B}}$, ${{S}^{C}}$, ${{S}^{D}}$ are the areas of the face triangles opposite to the four vertices $A$, $B$, $C$, $D$ of the tetrahedron $ABCD$ respectively; $S$ is the surface area of tetrahedron $ABCD$;	and $AB$, $AC$, $AD$, $BC$, $BD$, $CD$ are the lengths of the six edges of the tetrahedron respectively.
\end{theorem}

\begin{proof}
	It is directly obtained from theorem \ref{thm:Thm26.1.2}:	
	\[Q{{I}^{2}}={{R}^{2}}-\frac{1}{{{S}^{2}}}\left( \begin{aligned}
		& {{S}^{A}}{{S}^{B}}A{{B}^{2}}+{{S}^{A}}{{S}^{C}}A{{C}^{2}}+{{S}^{A}}{{S}^{C}}A{{D}^{2}} \\ 
		& +{{S}^{B}}{{S}^{C}}B{{C}^{2}}+{{S}^{B}}{{S}^{D}}B{{D}^{2}}+{{S}^{C}}{{S}^{D}}C{{D}^{2}}  
	\end{aligned} \right).\]
	
	Since $Q{{I}^{2}}\ge 0$, then	
	\[{{R}^{2}}\ge \frac{1}{{{S}^{2}}}\left( \begin{aligned}
		& {{S}^{A}}{{S}^{B}}A{{B}^{2}}+{{S}^{A}}{{S}^{C}}A{{C}^{2}}+{{S}^{A}}{{S}^{C}}A{{D}^{2}} \\ 
		& +{{S}^{B}}{{S}^{C}}B{{C}^{2}}+{{S}^{B}}{{S}^{D}}B{{D}^{2}}+{{S}^{C}}{{S}^{D}}C{{D}^{2}}  
	\end{aligned} \right).\]	
	
	The above equality occurs if and only if $Q{{I}^{2}}=0$, i.e., the equality occurs if and only if the circumcenter coincides with the incenter of the tetrahedron.
\end{proof}
\hfill $\square$\par

\section{Geometric inequality of distance between two intersecting centers of tetrahedron}\label{Sec27.2}
The following inequality is the origin from which some geometric inequalities can be constructed, that is, the “Father” inequalities of those geometric inequalities.

\begin{theorem}{Geometric inequality of distance between two ICs-T, Daiyuan Zhang}{Thm27.2.1}\label{Thm27.2.1} 
	Suppose that given a tetrahedron $ABCD$, ${{P}_{1}}\in {{\pi }_{ABCD}}$, ${{P}_{2}}\in {{\pi }_{ABCD}}$, then
	\[\begin{aligned}
		& \beta _{A}^{{{P}_{1}}{{P}_{2}}}\beta _{B}^{{{P}_{1}}{{P}_{2}}}A{{B}^{2}}+\beta _{A}^{{{P}_{1}}{{P}_{2}}}\beta _{C}^{{{P}_{1}}{{P}_{2}}}A{{C}^{2}}+\beta _{A}^{{{P}_{1}}{{P}_{2}}}\beta _{D}^{{{P}_{1}}{{P}_{2}}}A{{D}^{2}} \\ 
		& +\beta _{B}^{{{P}_{1}}{{P}_{2}}}\beta _{C}^{{{P}_{1}}{{P}_{2}}}B{{C}^{2}}+\beta _{B}^{{{P}_{1}}{{P}_{2}}}\beta _{D}^{{{P}_{1}}{{P}_{2}}}B{{D}^{2}}+\beta _{C}^{{{P}_{1}}{{P}_{2}}}\beta _{D}^{{{P}_{1}}{{P}_{2}}}C{{D}^{2}}\le 0. \\ 
	\end{aligned}\]	
	
	The above equality occurs if and only if the two intersecting centers ${{P}_{1}}$ and ${{P}_{2}}$ of the tetrahedron coincides with each other. Where	
	\[\beta _{A}^{{{P}_{1}}{{P}_{2}}}=\beta _{A}^{{{P}_{2}}}-\beta _{A}^{{{P}_{1}}},\]	
	\[\beta _{B}^{{{P}_{1}}{{P}_{2}}}=\beta _{B}^{{{P}_{2}}}-\beta _{B}^{{{P}_{1}}},\]	
	\[\beta _{C}^{{{P}_{1}}{{P}_{2}}}=\beta _{C}^{{{P}_{2}}}-\beta _{C}^{{{P}_{1}}},\]	
	\[\beta _{D}^{{{P}_{1}}{{P}_{2}}}=\beta _{D}^{{{P}_{2}}}-\beta _{D}^{{{P}_{1}}}.\]	
	Where $\beta _{A}^{{{P}_{1}}}$, $\beta _{B}^{{{P}_{1}}}$, $\beta _{C}^{{{P}_{1}}}$, $\beta _{D}^{{{P}_{1}}}$ are the frame components of the frame $\overrightarrow{OA}$, $\overrightarrow{OB}$, $\overrightarrow{OC}$, $\overrightarrow{OD}$ on the frame system $\left( O;A,B,C,D \right)$ at point ${P}_{1}$ respectively; $\beta _{A}^{{{P}_{2}}}$, $\beta _{B}^{{{P}_{2}}}$, $\beta _{C}^{{{P}_{2}}}$, $\beta _{D}^{{{P}_{2}}}$ are the frame components of the frame $\overrightarrow{OA}$, $\overrightarrow{OB}$, $\overrightarrow{OC}$, $\overrightarrow{OD}$ on the frame system $\left( O;A,B,C,D \right)$ at point ${P}_{2}$ respectively.
\end{theorem}

\begin{proof}
	According to theorem \ref{thm:Thm24.2.1}, the following result is obtained:	
	\[\begin{aligned}
		{{P}_{1}}{{P}_{2}}^{2}&=-\beta _{A}^{{{P}_{1}}{{P}_{2}}}\beta _{B}^{{{P}_{1}}{{P}_{2}}}A{{B}^{2}}-\beta _{A}^{{{P}_{1}}{{P}_{2}}}\beta _{C}^{{{P}_{1}}{{P}_{2}}}A{{C}^{2}}-\beta _{A}^{{{P}_{1}}{{P}_{2}}}\beta _{D}^{{{P}_{1}}{{P}_{2}}}A{{D}^{2}} \\ 
		& -\beta _{B}^{{{P}_{1}}{{P}_{2}}}\beta _{C}^{{{P}_{1}}{{P}_{2}}}B{{C}^{2}}-\beta _{B}^{{{P}_{1}}{{P}_{2}}}\beta _{D}^{{{P}_{1}}{{P}_{2}}}B{{D}^{2}}-\beta _{C}^{{{P}_{1}}{{P}_{2}}}\beta _{D}^{{{P}_{1}}{{P}_{2}}}C{{D}^{2}},  
	\end{aligned}\]
	
	Since ${{P}_{1}}{{P}_{2}}^{2}\ge 0$, then	
	\[\begin{aligned}
		& \beta _{A}^{{{P}_{1}}{{P}_{2}}}\beta _{B}^{{{P}_{1}}{{P}_{2}}}A{{B}^{2}}+\beta _{A}^{{{P}_{1}}{{P}_{2}}}\beta _{C}^{{{P}_{1}}{{P}_{2}}}A{{C}^{2}}+\beta _{A}^{{{P}_{1}}{{P}_{2}}}\beta _{D}^{{{P}_{1}}{{P}_{2}}}A{{D}^{2}} \\ 
		& +\beta _{B}^{{{P}_{1}}{{P}_{2}}}\beta _{C}^{{{P}_{1}}{{P}_{2}}}B{{C}^{2}}+\beta _{B}^{{{P}_{1}}{{P}_{2}}}\beta _{D}^{{{P}_{1}}{{P}_{2}}}B{{D}^{2}}+\beta _{C}^{{{P}_{1}}{{P}_{2}}}\beta _{D}^{{{P}_{1}}{{P}_{2}}}C{{D}^{2}}\le 0. \\ 
	\end{aligned}\]	
	
	Obviously, the above equality holds if and only if ${{P}_{1}}{{P}_{2}}=0$, that is, if and only if the two intersecting centers ${{P}_{1}}$ and ${{P}_{2}}$ coincides with each other.
\end{proof}
\hfill $\square$\par

As an application of the above theorem, several special cases are given below.
\begin{theorem}{Centroid-Incenter inequality of tetrahedron, Daiyuan Zhang}{Thm27.2.2}\label{Thm27.2.2} 
	Suppose that given a tetrahedron $ABCD$, then	
	\[\begin{aligned}
		& \left( 4{{S}^{A}}-S \right)\left( 4{{S}^{B}}-S \right)A{{B}^{2}}+\left( 4{{S}^{A}}-S \right)\left( 4{{S}^{C}}-S \right)A{{C}^{2}} \\ 
		& +\left( 4{{S}^{A}}-S \right)\left( 4{{S}^{D}}-S \right)A{{D}^{2}}+\left( 4{{S}^{B}}-S \right)\left( 4{{S}^{C}}-S \right)B{{C}^{2}} \\ 
		& +\left( 4{{S}^{B}}-S \right)\left( 4{{S}^{D}}-S \right)B{{D}^{2}}+\left( 4{{S}^{C}}-S \right)\left( 4{{S}^{D}}-S \right)C{{D}^{2}}\le 0. \\ 
	\end{aligned}\]			
	
	The above equality occurs if and only if the centroid coincides with the incenter. Where ${{S}^{A}}$, ${{S}^{B}}$, ${{S}^{C}}$, ${{S}^{D}}$ are the areas of the face triangles opposite to the four vertices $A$, $B$, $C$, $D$ of the tetrahedron $ABCD$ respectively; and $AB$, $AC$, $AD$, $BC$, $BD$, $CD$ are the lengths of the six edges of the tetrahedron respectively.
\end{theorem}

\begin{proof}
	According to theorem \ref{thm:Thm26.2.1}, the following result is obtained:	
	\[G{{I}^{2}}=-\frac{1}{16{{S}^{2}}}\left( \begin{aligned}
		& \left( 4{{S}^{A}}-S \right)\left( 4{{S}^{B}}-S \right)A{{B}^{2}}+\left( 4{{S}^{A}}-S \right)\left( 4{{S}^{C}}-S \right)A{{C}^{2}} \\ 
		& +\left( 4{{S}^{A}}-S \right)\left( 4{{S}^{D}}-S \right)A{{D}^{2}}+\left( 4{{S}^{B}}-S \right)\left( 4{{S}^{C}}-S \right)B{{C}^{2}} \\ 
		& +\left( 4{{S}^{B}}-S \right)\left( 4{{S}^{D}}-S \right)B{{D}^{2}}+\left( 4{{S}^{C}}-S \right)\left( 4{{S}^{D}}-S \right)C{{D}^{2}}  
	\end{aligned} \right).\]
	
	Since $G{{I}^{2}}\ge 0$, then	
	\[\begin{aligned}
		& \left( 4{{S}^{A}}-S \right)\left( 4{{S}^{B}}-S \right)A{{B}^{2}}+\left( 4{{S}^{A}}-S \right)\left( 4{{S}^{C}}-S \right)A{{C}^{2}} \\ 
		& +\left( 4{{S}^{A}}-S \right)\left( 4{{S}^{D}}-S \right)A{{D}^{2}}+\left( 4{{S}^{B}}-S \right)\left( 4{{S}^{C}}-S \right)B{{C}^{2}} \\ 
		& +\left( 4{{S}^{B}}-S \right)\left( 4{{S}^{D}}-S \right)B{{D}^{2}}+\left( 4{{S}^{C}}-S \right)\left( 4{{S}^{D}}-S \right)C{{D}^{2}}\le 0. \\ 
	\end{aligned}\]	
	
	Obviously, the above equality holds if and only if $G{{I}^{2}}=0$, that is, if and only if the centroid and incenter of the tetrahedron $ABCD$ coincides with each other.
\end{proof} \hfill $\square$\par	
	
	The above theorem can also be written in the following form:
	\[\begin{aligned}
		& \left( {{S}^{A}}-\bar{S} \right)\left( {{S}^{B}}-\bar{S} \right)A{{B}^{2}}+\left( {{S}^{A}}-\bar{S} \right)\left( {{S}^{C}}-\bar{S} \right)A{{C}^{2}} \\ 
		& +\left( {{S}^{A}}-\bar{S} \right)\left( {{S}^{D}}-\bar{S} \right)A{{D}^{2}}+\left( {{S}^{B}}-\bar{S} \right)\left( {{S}^{C}}-\bar{S} \right)B{{C}^{2}} \\ 
		& +\left( {{S}^{B}}-\bar{S} \right)\left( {{S}^{D}}-\bar{S} \right)B{{D}^{2}}+\left( {{S}^{C}}-\bar{S} \right)\left( {{S}^{D}}-\bar{S} \right)C{{D}^{2}}\le 0. \\ 
	\end{aligned}\]			
	
	Where $\bar{S}$ is the average value of tetrahedral surface area, i.e
	\[\bar{S}=\frac{1}{4}\left( {{S}^{A}}+{{S}^{B}}+{{S}^{C}}+{{S}^{D}} \right).\]	

\begin{theorem}{Centroid-Circumcenter inequality of tetrahedron, Daiyuan Zhang}{Thm27.2.3}\label{Thm27.2.3} 
Suppose that given a tetrahedron $ABCD$, then	
\[\begin{aligned}
	& \left( 4{{U}_{A}}-U \right)\left( 4{{U}_{B}}-U \right)A{{B}^{2}}+\left( 4{{U}_{A}}-U \right)\left( 4{{U}_{C}}-U \right)A{{C}^{2}} \\ 
	& +\left( 4{{U}_{A}}-U \right)\left( 4{{U}_{D}}-U \right)A{{D}^{2}}+\left( 4{{U}_{B}}-U \right)\left( 4{{U}_{C}}-U \right)B{{C}^{2}} \\ 
	& +\left( 4{{U}_{B}}-U \right)\left( 4{{U}_{D}}-U \right)B{{D}^{2}}+\left( 4{{U}_{C}}-U \right)\left( 4{{U}_{D}}-U \right)C{{D}^{2}}\le 0. \\ 
\end{aligned}\]	

The above equality occurs if and only if the centroid coincides with the circumcenter. Where $U$, ${{U}_{A}}$, ${{U}_{B}}$, ${{U}_{C}}$, ${{U}_{D}}$ are calculated by theorem \ref{thm:Thm25.3.1}, theorem \ref{thm:Thm25.1.2} or theorem \ref{thm:Thm25.1.1}; and $AB$, $AC$, $AD$, $BC$, $BD$, $CD$ are the lengths of the six edges of the tetrahedron respectively.	
\end{theorem}

\begin{proof}
	From theorem \ref{thm:Thm26.2.2} we have
	\[G{{Q}^{2}}=-\frac{1}{16{{U}^{2}}}\left( \begin{aligned}
		& \left( 4{{U}_{A}}-U \right)\left( 4{{U}_{B}}-U \right)A{{B}^{2}}+\left( 4{{U}_{A}}-U \right)\left( 4{{U}_{C}}-U \right)A{{C}^{2}} \\ 
		& +\left( 4{{U}_{A}}-U \right)\left( 4{{U}_{D}}-U \right)A{{D}^{2}}+\left( 4{{U}_{B}}-U \right)\left( 4{{U}_{C}}-U \right)B{{C}^{2}} \\ 
		& +\left( 4{{U}_{B}}-U \right)\left( 4{{U}_{D}}-U \right)B{{D}^{2}}+\left( 4{{U}_{C}}-U \right)\left( 4{{U}_{D}}-U \right)C{{D}^{2}}.  
	\end{aligned} \right)\]
	
	Because $G{{Q}^{2}}\ge 0$, then:
	\[\begin{aligned}
		& \left( 4{{U}_{A}}-U \right)\left( 4{{U}_{B}}-U \right)A{{B}^{2}}+\left( 4{{U}_{A}}-U \right)\left( 4{{U}_{C}}-U \right)A{{C}^{2}} \\ 
		& +\left( 4{{U}_{A}}-U \right)\left( 4{{U}_{D}}-U \right)A{{D}^{2}}+\left( 4{{U}_{B}}-U \right)\left( 4{{U}_{C}}-U \right)B{{C}^{2}} \\ 
		& +\left( 4{{U}_{B}}-U \right)\left( 4{{U}_{D}}-U \right)B{{D}^{2}}+\left( 4{{U}_{C}}-U \right)\left( 4{{U}_{D}}-U \right)C{{D}^{2}}\le 0. \\ 
	\end{aligned}\]
	
	Obviously, the above equality holds if and only if $G{{Q}^{2}}=0$, that is, if and only if the centroid and circumcenter of the tetrahedron $ABCD$ coincides with each other. 
\end{proof} \hfill $\square$\par

\begin{theorem}{Incenter-Circumcenter inequality of tetrahedron, Daiyuan Zhang}{Thm27.2.4}\label{Thm27.2.4} 
Suppose that given a tetrahedron $ABCD$, then	
\[\begin{aligned}
	& \left( S{{U}_{A}}-{{S}^{A}}U \right)\left( S{{U}_{B}}-{{S}^{B}}U \right)A{{B}^{2}}+\left( S{{U}_{A}}-{{S}^{A}}U \right)\left( S{{U}_{C}}-{{S}^{C}}U \right)A{{C}^{2}} \\ 
	& +\left( S{{U}_{A}}-{{S}^{A}}U \right)\left( S{{U}_{D}}-{{S}^{D}}U \right)A{{D}^{2}}+\left( S{{U}_{B}}-{{S}^{B}}U \right)\left( S{{U}_{C}}-{{S}^{C}}U \right)B{{C}^{2}} \\ 
	& +\left( S{{U}_{B}}-{{S}^{B}}U \right)\left( S{{U}_{D}}-{{S}^{D}}U \right)B{{D}^{2}}+\left( S{{U}_{C}}-{{S}^{C}}U \right)\left( S{{U}_{D}}-{{S}^{D}}U \right)C{{D}^{2}}\le 0. \\ 
\end{aligned}\]

The above equality occurs if and only if the incenter coincides with the circumcenter. Where $U$, ${{U}_{A}}$, ${{U}_{B}}$, ${{U}_{C}}$, ${{U}_{D}}$ are calculated by theorem \ref{thm:Thm25.3.1}, theorem \ref{thm:Thm25.1.2} or theorem \ref{thm:Thm25.1.1}; ${{S}^{A}}$, ${{S}^{B}}$, ${{S}^{C}}$, ${{S}^{D}}$ are the areas of the face triangles opposite to the four vertices $A$, $B$, $C$, $D$ of the tetrahedron $ABCD$ respectively; $S$ is the surface area of tetrahedron $ABCD$, i.e. $S={{S}^{A}}+{{S}^{B}}+{{S}^{C}}+{{S}^{D}}$; and $AB$, $AC$, $AD$, $BC$, $BD$, $CD$ are the lengths of the six edges of the tetrahedron respectively.
\end{theorem}

\begin{proof}
	From theorem \ref{thm:Thm26.2.3} we have
	\[I{{Q}^{2}}=-\frac{1}{{S}^{2}{U}^{2}}\left( \begin{aligned}
		& \left( S{{U}_{A}}-{{S}^{A}}U \right)\left( S{{U}_{B}}-{{S}^{B}}U \right)A{{B}^{2}} \\ 
		& +\left( S{{U}_{A}}-{{S}^{A}}U \right)\left( S{{U}_{C}}-{{S}^{C}}U \right)A{{C}^{2}} \\ 
		& +\left( S{{U}_{A}}-{{S}^{A}}U \right)\left( S{{U}_{D}}-{{S}^{D}}U \right)A{{D}^{2}} \\ 
		& +\left( S{{U}_{B}}-{{S}^{B}}U \right)\left( S{{U}_{C}}-{{S}^{C}}U \right)B{{C}^{2}} \\ 
		& +\left( S{{U}_{B}}-{{S}^{B}}U \right)\left( S{{U}_{D}}-{{S}^{D}}U \right)B{{D}^{2}} \\ 
		& +\left( S{{U}_{C}}-{{S}^{C}}U \right)\left( S{{U}_{D}}-{{S}^{D}}U \right)C{{D}^{2}}  
	\end{aligned} \right).\]
	
	Because $I{{Q}^{2}}\ge 0$, then:
	\[\begin{aligned}
		& \left( S{{U}_{A}}-{{S}^{A}}U \right)\left( S{{U}_{B}}-{{S}^{B}}U \right)A{{B}^{2}}+\left( S{{U}_{A}}-{{S}^{A}}U \right)\left( S{{U}_{C}}-{{S}^{C}}U \right)A{{C}^{2}} \\ 
		& +\left( S{{U}_{A}}-{{S}^{A}}U \right)\left( S{{U}_{D}}-{{S}^{D}}U \right)A{{D}^{2}}+\left( S{{U}_{B}}-{{S}^{B}}U \right)\left( S{{U}_{C}}-{{S}^{C}}U \right)B{{C}^{2}} \\ 
		& +\left( S{{U}_{B}}-{{S}^{B}}U \right)\left( S{{U}_{D}}-{{S}^{D}}U \right)B{{D}^{2}}+\left( S{{U}_{C}}-{{S}^{C}}U \right)\left( S{{U}_{D}}-{{S}^{D}}U \right)C{{D}^{2}}\le 0. \\ 
	\end{aligned}\]
	
	Obviously, the above equality holds if and only if $I{{Q}^{2}}=0$, that is, if and only if the incenter and circumcenter of the tetrahedron $ABCD$ coincides with each other. 
\end{proof} \hfill $\square$\par


\chapter{Frame components of projection point}\label{Ch28}
\thispagestyle{empty}




%

%

%

%



%

%

%

%

%


It has been pointed out that the frame component plays an important role in Intercenter Geometry. This chapter studies the frame components of projection point at any point in space, the frame components of projection point of tetrahedral vertex, the frame components of projection point of an IC-T, the relationship between the frame components of projection point and the frame components of the IC-T, and the application of the frame components of projection point.

\section{Basic concept of frame component of projection point}\label{Sec28.1}

In this chapter, we study the frame components of tetrahedral projection point.

For a given tetrahedron, any point $P$ in space can lead a vertical line to each triangular face (plane) of the tetrahedron. The intersection between the vertical line and the triangular face (plane) is called the \textbf{projection point} of point $P$ to the plane, and the plane where the projection point is located is called the \textbf{projection plane}. Point $P$ has four projection points and four projection distances for a tetrahedron. 

\begin{figure}[h]
	\centering
	\includegraphics[totalheight=6cm]{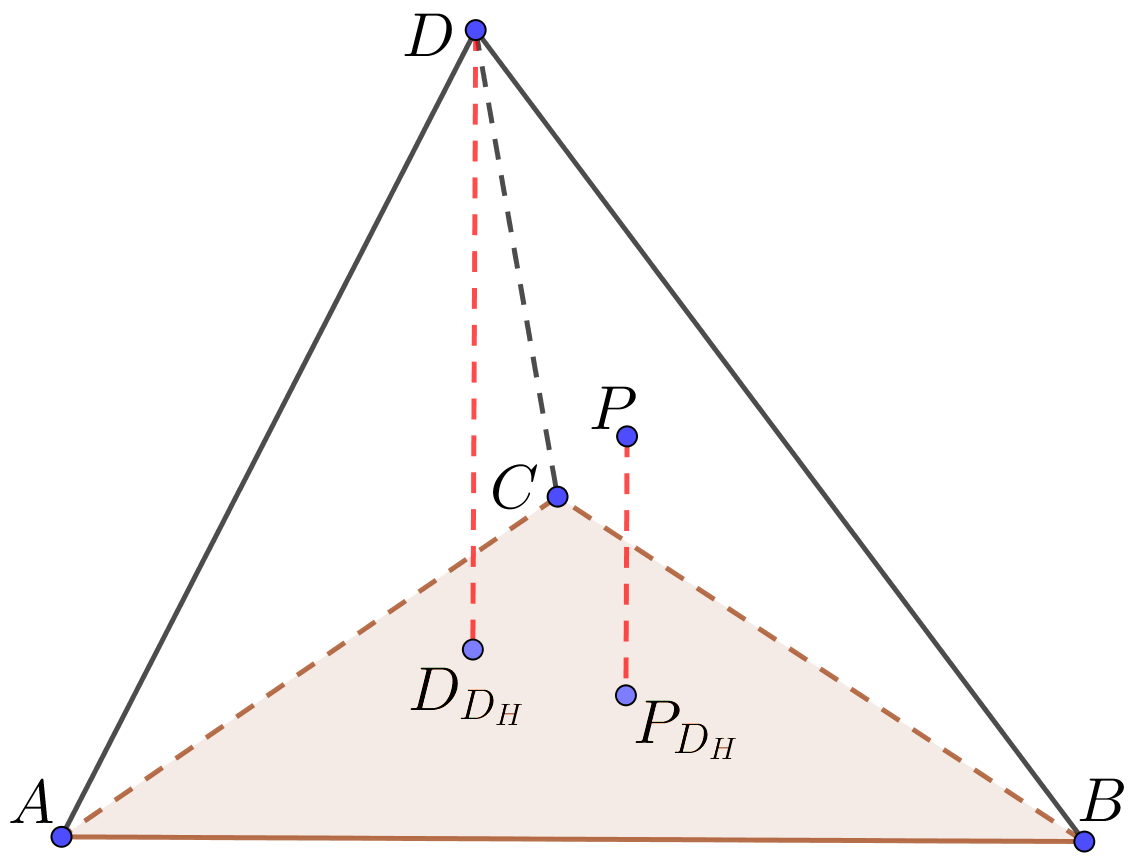}
	\caption{Projection points of a tetrahedron} \label{fig:tu28.1.1}
\end{figure}

Figure \ref{fig:tu28.1.1} shows the projection point ${{P}_{{{D}_{H}}}}$ of point $P$ on the plane of $\triangle ABC$ (projection plane). The concepts are given in the following. 

\textbf{Projection line}: the line between the point $P$ and the projection point ${{P}_{{{D}_{H}}}}$.

\textbf{Projection vector}: the position vector between the point $P$ and the projection point ${{P}_{{{D}_{H}}}}$, initial point is $P$.

\textbf{Projection distance}: the distance between the point $P$ and the projection point ${{P}_{{{D}_{H}}}}$.


\textbf{Frame components of projection point}: The frame components of the projection point ${{P}_{{{D}_{H}}}}$ (of point $P$) on the triangular frame $\left( O;A,B,C \right)$, where point $O$ is any point in space.

\section{Frame components of projection points of a point in space}\label{Sec28.2}
This section studies the frame components of the projection points of a point $P$ in space. The following theorem I give answers this question.
\begin{theorem}{Formulas of frame components for projection of a point in space-1, Daiyuan Zhang}{Thm28.2.1}\label{Thm28.2.1} 
	Given a tetrahedron $ABCD$, let point $O$ and point $P$ be a point in space respectively, then:
	
	
	1.Assuming that the projection point of the point $P$ to the $\triangle BCD$ plane is $ {{P}_{{{D}_{H}}}}$, ${{P}_{{{D}_{H}}}}\in {{\pi }_{ABC}}$, $\triangle BCD$ has an area of ${{S}^{A}}$, then the frame components of ${{P}_{{{D}_{H}}}}$ on the frame system $\left(O;A,B,C \right)$ are in the following:
	
	\[\alpha _{B}^{{{P}_{{{A}_{H}}}}}=\frac{1}{8{{\left( {{S}^{A}} \right)}^{2}}}\left( \begin{aligned}
		& \left( \Delta _{2}^{A}-C{{D}^{2}} \right)C{{D}^{2}}+\left( \Delta _{2}^{A}-D{{B}^{2}} \right)\left( D{{P}^{2}}-B{{P}^{2}} \right) \\ 
		& +\left( \Delta _{2}^{A}-B{{C}^{2}} \right)\left( C{{P}^{2}}-B{{P}^{2}} \right) \\ 
	\end{aligned} \right),\]
	\[\alpha _{C}^{{{P}_{{{A}_{H}}}}}=\frac{1}{8{{\left( {{S}^{A}} \right)}^{2}}}\left( \begin{aligned}
		& \left( \Delta _{2}^{A}-D{{B}^{2}} \right)D{{B}^{2}}+\left( \Delta _{2}^{A}-B{{C}^{2}} \right)\left( B{{P}^{2}}-C{{P}^{2}} \right) \\ 
		& +\left( \Delta _{2}^{A}-C{{D}^{2}} \right)\left( D{{P}^{2}}-C{{P}^{2}} \right) \\ 
	\end{aligned} \right),\]
	\[\alpha _{D}^{{{P}_{{{A}_{H}}}}}=\frac{1}{8{{\left( {{S}^{A}} \right)}^{2}}}\left( \begin{aligned}
		& \left( \Delta _{2}^{A}-B{{C}^{2}} \right)B{{C}^{2}}+\left( \Delta _{2}^{A}-C{{D}^{2}} \right)\left( C{{P}^{2}}-D{{P}^{2}} \right) \\ 
		& +\left( \Delta _{2}^{A}-D{{B}^{2}} \right)\left( B{{P}^{2}}-D{{P}^{2}} \right) \\ 
	\end{aligned} \right).\]
	Where 
	\[\Delta _{2}^{A}=\frac{1}{2}\left( B{{C}^{2}}+C{{D}^{2}}+D{{B}^{2}} \right).\]
	
	
	2.Assuming that the projection point of the point $P$ to the $\triangle CDA$ plane is ${{P}_{{{B}_{H}}}}$, ${{P}_{{{B}_{H}}}}\in {{\pi }_{CDA}}$, $\triangle CDA$ has an area of ${{S}^{B}}$, then the frame components of ${{P}_{{{B}_{H}}}}$ on the frame system $\left( O;C,D,A \right)$ are in the following:
	
	\[\alpha _{C}^{{{P}_{{{B}_{H}}}}}=\frac{1}{8{{\left( {{S}^{B}} \right)}^{2}}}\left( \begin{aligned}
		& \left( \Delta _{2}^{B}-D{{A}^{2}} \right)D{{A}^{2}}+\left( \Delta _{2}^{B}-A{{C}^{2}} \right)\left( A{{P}^{2}}-C{{P}^{2}} \right) \\ 
		& +\left( \Delta _{2}^{B}-C{{D}^{2}} \right)\left( D{{P}^{2}}-C{{P}^{2}} \right) \\ 
	\end{aligned} \right),\]
	\[\alpha _{D}^{{{P}_{{{B}_{H}}}}}=\frac{1}{8{{\left( {{S}^{B}} \right)}^{2}}}\left( \begin{aligned}
		& \left( \Delta _{2}^{B}-A{{C}^{2}} \right)A{{C}^{2}}+\left( \Delta _{2}^{B}-C{{D}^{2}} \right)\left( C{{P}^{2}}-D{{P}^{2}} \right) \\ 
		& +\left( \Delta _{2}^{B}-D{{A}^{2}} \right)\left( A{{P}^{2}}-D{{P}^{2}} \right) \\ 
	\end{aligned} \right),\]
	\[\alpha _{A}^{{{P}_{{{B}_{H}}}}}=\frac{1}{8{{\left( {{S}^{B}} \right)}^{2}}}\left( \begin{aligned}
		& \left( \Delta _{2}^{B}-C{{D}^{2}} \right)C{{D}^{2}}+\left( \Delta _{2}^{B}-D{{A}^{2}} \right)\left( D{{P}^{2}}-A{{P}^{2}} \right) \\ 
		& +\left( \Delta _{2}^{B}-A{{C}^{2}} \right)\left( C{{P}^{2}}-A{{P}^{2}} \right) \\ 
	\end{aligned} \right).\]
	Where 
	\[\Delta _{2}^{B}=\frac{1}{2}\left( C{{D}^{2}}+D{{A}^{2}}+A{{C}^{2}} \right).\]
	
	
	3.Assuming that the projection point of the point $P$ to the $\triangle DAB$ plane is ${{P}_{{{C}_{H}}}}$, ${{P}_{{{C}_{H}}}}\in {{\pi }_{DAB}}$, $\triangle DAB$ has an area of ${{S}^{C}}$, then the frame components of ${{P}_{{{C}_{H}}}}$ on the frame system $\left( O;D,A,B \right)$ are in the following:
	
	\[\alpha _{D}^{{{P}_{{{C}_{H}}}}}=\frac{1}{8{{\left( {{S}^{C}} \right)}^{2}}}\left( \begin{aligned}
		& \left( \Delta _{2}^{C}-A{{B}^{2}} \right)A{{B}^{2}}+\left( \Delta _{2}^{C}-B{{D}^{2}} \right)\left( B{{P}^{2}}-D{{P}^{2}} \right) \\ 
		& +\left( \Delta _{2}^{C}-D{{A}^{2}} \right)\left( A{{P}^{2}}-D{{P}^{2}} \right) \\ 
	\end{aligned} \right),\]
	\[\alpha _{A}^{{{P}_{{{C}_{H}}}}}=\frac{1}{8{{\left( {{S}^{C}} \right)}^{2}}}\left( \begin{aligned}
		& \left( \Delta _{2}^{C}-B{{D}^{2}} \right)B{{D}^{2}}+\left( \Delta _{2}^{C}-D{{A}^{2}} \right)\left( D{{P}^{2}}-A{{P}^{2}} \right) \\ 
		& +\left( \Delta _{2}^{C}-A{{B}^{2}} \right)\left( B{{P}^{2}}-A{{P}^{2}} \right) \\ 
	\end{aligned} \right),\]
	\[\alpha _{B}^{{{P}_{{{C}_{H}}}}}=\frac{1}{8{{\left( {{S}^{C}} \right)}^{2}}}\left( \begin{aligned}
		& \left( \Delta _{2}^{C}-D{{A}^{2}} \right)D{{A}^{2}}+\left( \Delta _{2}^{C}-A{{B}^{2}} \right)\left( A{{P}^{2}}-B{{P}^{2}} \right) \\ 
		& +\left( \Delta _{2}^{C}-B{{D}^{2}} \right)\left( D{{P}^{2}}-B{{P}^{2}} \right) \\ 
	\end{aligned} \right).\]
	Where 
	\[\Delta _{2}^{C}=\frac{1}{2}\left( D{{A}^{2}}+A{{B}^{2}}+B{{D}^{2}} \right).\]
	
	
	4.Assuming that the projection point of the point $P$ to the $\triangle ABC$ plane is ${{P}_{{{D}_{H}}}}$, ${{P}_{{{D}_{H}}}}\in {{\pi }_{ABC}}$, $\triangle ABC$ has an area of ${{S}^{D}}$, then the frame components of ${{P}_{{{D}_{H}}}}$ on the frame system $\left( O;A,B,C \right)$ are in the following:
	
	\[\alpha _{A}^{{{P}_{{{D}_{H}}}}}=\frac{1}{8{{\left( {{S}^{D}} \right)}^{2}}}\left( \begin{aligned}
		& \left( \Delta _{2}^{D}-B{{C}^{2}} \right)B{{C}^{2}}+\left( \Delta _{2}^{D}-C{{A}^{2}} \right)\left( C{{P}^{2}}-A{{P}^{2}} \right) \\ 
		& +\left( \Delta _{2}^{D}-A{{B}^{2}} \right)\left( B{{P}^{2}}-A{{P}^{2}} \right) \\ 
	\end{aligned} \right),\]
	\[\alpha _{B}^{{{P}_{{{D}_{H}}}}}=\frac{1}{8{{\left( {{S}^{D}} \right)}^{2}}}\left( \begin{aligned}
		& \left( \Delta _{2}^{D}-C{{A}^{2}} \right)C{{A}^{2}}+\left( \Delta _{2}^{D}-A{{B}^{2}} \right)\left( A{{P}^{2}}-B{{P}^{2}} \right) \\ 
		& +\left( \Delta _{2}^{D}-B{{C}^{2}} \right)\left( C{{P}^{2}}-B{{P}^{2}} \right) \\ 
	\end{aligned} \right),\]
	\[\alpha _{C}^{{{P}_{{{D}_{H}}}}}=\frac{1}{8{{\left( {{S}^{D}} \right)}^{2}}}\left( \begin{aligned}
		& \left( \Delta _{2}^{D}-A{{B}^{2}} \right)A{{B}^{2}}+\left( \Delta _{2}^{D}-B{{C}^{2}} \right)\left( B{{P}^{2}}-C{{P}^{2}} \right) \\ 
		& +\left( \Delta _{2}^{D}-C{{A}^{2}} \right)\left( A{{P}^{2}}-C{{P}^{2}} \right) \\ 
	\end{aligned} \right).\]
	Where 
	\[\Delta _{2}^{D}=\frac{1}{2}\left( A{{B}^{2}}+B{{C}^{2}}+C{{A}^{2}} \right).\]
\end{theorem}

\begin{proof}
	
	Take the projection point ${{P}_{{{D}_{H}}}}$ as an example to prove the theorem, which is the the projection point of $P$ on the plane of $\triangle ABC$. 
	
	According to theorem \ref{thm:Thm12.1.2}, it is obtained that the distance between the vertex and intersecting center of the projection point ${{P}_{{{D}_{H}}}}$ in the $\triangle ABC$ plane is
		
	\[AP_{{{D}_{H}}}^{2}=\left( 1-\alpha _{A}^{{{P}_{{{D}_{H}}}}} \right)\left( \alpha _{B}^{{{P}_{{{D}_{H}}}}}A{{B}^{2}}+\alpha _{C}^{{{P}_{{{D}_{H}}}}}C{{A}^{2}} \right)-\alpha _{B}^{{{P}_{{{D}_{H}}}}}\alpha _{C}^{{{P}_{{{D}_{H}}}}}B{{C}^{2}},\]
	\[BP_{{{D}_{H}}}^{2}=\left( 1-\alpha _{B}^{{{P}_{{{D}_{H}}}}} \right)\left( \alpha _{C}^{{{P}_{{{D}_{H}}}}}B{{C}^{2}}+\alpha _{A}^{{{P}_{{{D}_{H}}}}}A{{B}^{2}} \right)-\alpha _{C}^{{{P}_{{{D}_{H}}}}}\alpha _{A}^{{{P}_{{{D}_{H}}}}}C{{A}^{2}},\]
	\[CP_{{{D}_{H}}}^{2}=\left( 1-\alpha _{C}^{{{P}_{{{D}_{H}}}}} \right)\left( \alpha _{A}^{{{P}_{{{D}_{H}}}}}C{{A}^{2}}+\alpha _{B}^{{{P}_{{{D}_{H}}}}}B{{C}^{2}} \right)-\alpha _{A}^{{{P}_{{{D}_{H}}}}}\alpha _{B}^{{{P}_{{{D}_{H}}}}}A{{B}^{2}}.\]
	
	
	It is necessary to calculate the frame components of $\alpha _{A}^{{{P}_{{{D}_{H}}}}}$, $\alpha _{B}^{{{P}_{{{D}_{H}}}}}$ and $\alpha _{C}^{{{P}_{{{D}_{H}}}}}$ of point ${{P}_{{{D}_{H}}}}$ on the triangular frame $\left( O;A,B,C \right)$.
	
	
	Obviously, $\triangle AP{{P}_{{{D}_{H}}}}$, $\triangle BP{{P}_{{{D}_{H}}}}$ and $\triangle CP{{P}_{{{D}_{H}}}}$ are right triangles (see figure \ref{fig:tu28.1.1}), using Pythagorean theorem, we have:
		
	\begin{flalign*}
		A{{P}^{2}}=AP_{{{D}_{H}}}^{2}+PP_{{{D}_{H}}}^{2},B{{P}^{2}}=BP_{{{D}_{H}}}^{2}+PP_{{{D}_{H}}}^{2},C{{P}^{2}}=CP_{{{D}_{H}}}^{2}+PP_{{{D}_{H}}}^{2}.	
	\end{flalign*}
	i.e. 
	\[\left\{ \begin{aligned}
		& CP_{{{D}_{H}}}^{2}-AP_{{{D}_{H}}}^{2}=C{{P}^{2}}-A{{P}^{2}} \\ 
		& AP_{{{D}_{H}}}^{2}-BP_{{{D}_{H}}}^{2}=A{{P}^{2}}-B{{P}^{2}}. \\ 
	\end{aligned} \right.\]
	
	
	Add the condition $\alpha _{A}^{{{P}_{{{D}_{H}}}}}+\alpha _{B}^{{{P}_{{{D}_{H}}}}}+\alpha _{C}^{{{P}_{{{D}_{H}}}}}=1$, the following linear equations are obtained:
	
	\[\left\{ \begin{aligned}
		& \alpha _{A}^{{{P}_{{{D}_{H}}}}}+\alpha _{B}^{{{P}_{{{D}_{H}}}}}+\alpha _{C}^{{{P}_{{{D}_{H}}}}}=1 \\ 
		& C{{A}^{2}}\alpha _{A}^{{{P}_{{{D}_{H}}}}}+\left( B{{C}^{2}}-A{{B}^{2}} \right)\alpha _{B}^{{{P}_{{{D}_{H}}}}}-C{{A}^{2}}\alpha _{C}^{{{P}_{{{D}_{H}}}}}=C{{P}^{2}}-A{{P}^{2}} \\ 
		& -A{{B}^{2}}\alpha _{A}^{{{P}_{{{D}_{H}}}}}+A{{B}^{2}}\alpha _{B}^{{{P}_{{{D}_{H}}}}}+\left( C{{A}^{2}}-B{{C}^{2}} \right)\alpha _{C}^{{{P}_{{{D}_{H}}}}}=A{{P}^{2}}-B{{P}^{2}}. \\ 
	\end{aligned} \right.\]
	
	Written in matrix form is
	\[\left( \begin{matrix}
		1 & 1 & 1  \\
		C{{A}^{2}} & B{{C}^{2}}-A{{B}^{2}} & -C{{A}^{2}}  \\
		-A{{B}^{2}} & A{{B}^{2}} & C{{A}^{2}}-B{{C}^{2}}  \\
	\end{matrix} \right)\left( \begin{matrix}
		\alpha _{A}^{{{P}_{{{D}_{H}}}}}  \\
		\alpha _{B}^{{{P}_{{{D}_{H}}}}}  \\
		\alpha _{C}^{{{P}_{{{D}_{H}}}}}  \\
	\end{matrix} \right)=\left( \begin{matrix}
		1  \\
		C{{P}^{2}}-A{{P}^{2}}  \\
		A{{P}^{2}}-B{{P}^{2}}  \\
	\end{matrix} \right).\]
	i.e. 
	\[\left( \begin{matrix}
		\alpha _{A}^{{{P}_{{{D}_{H}}}}}  \\
		\alpha _{B}^{{{P}_{{{D}_{H}}}}}  \\
		\alpha _{C}^{{{P}_{{{D}_{H}}}}}  \\
	\end{matrix} \right)={{\left( \begin{matrix}
				1 & 1 & 1  \\
				C{{A}^{2}} & B{{C}^{2}}-A{{B}^{2}} & -C{{A}^{2}}  \\
				-A{{B}^{2}} & A{{B}^{2}} & C{{A}^{2}}-B{{C}^{2}}  \\
			\end{matrix} \right)}^{-1}}\left( \begin{matrix}
		1  \\
		C{{P}^{2}}-A{{P}^{2}}  \\
		A{{P}^{2}}-B{{P}^{2}}  \\
	\end{matrix} \right).\]
	
	By solving the inverse matrix, we have:
	\[\left( \begin{matrix}
		\alpha _{A}^{{{P}_{{{D}_{H}}}}}  \\
		\alpha _{B}^{{{P}_{{{D}_{H}}}}}  \\
		\alpha _{C}^{{{P}_{{{D}_{H}}}}}  \\
	\end{matrix} \right)=\frac{1}{8{{\left( {{S}^{D}} \right)}^{2}}}\left( \begin{matrix}
		\left( \Delta _{2}^{D}-B{{C}^{2}} \right)B{{C}^{2}} & \Delta _{2}^{D}-C{{A}^{2}} & A{{B}^{2}}-\Delta _{2}^{D}  \\
		\left( \Delta _{2}^{D}-C{{A}^{2}} \right)C{{A}^{2}} & \Delta _{2}^{D}-B{{C}^{2}} & 2C{{A}^{2}}  \\
		\left( \Delta _{2}^{D}-A{{B}^{2}} \right)A{{B}^{2}} & -2A{{B}^{2}} & B{{C}^{2}}-\Delta _{2}^{D}  \\
	\end{matrix} \right)\left( \begin{matrix}
		1  \\
		C{{P}^{2}}-A{{P}^{2}}  \\
		A{{P}^{2}}-B{{P}^{2}}  \\
	\end{matrix} \right).\]
	Where 
	\[\Delta _{2}^{D}=\frac{1}{2}\left( A{{B}^{2}}+B{{C}^{2}}+C{{A}^{2}} \right).\]
	
	So one of the solutions is
	\[\alpha _{A}^{{{P}_{{{D}_{H}}}}}=\frac{1}{8{{\left( {{S}^{D}} \right)}^{2}}}\left( \begin{aligned}
		& \left( \Delta _{2}^{D}-B{{C}^{2}} \right)B{{C}^{2}}+\left( \Delta _{2}^{D}-C{{A}^{2}} \right)\left( C{{P}^{2}}-A{{P}^{2}} \right) \\ 
		& +\left( \Delta _{2}^{D}-A{{B}^{2}} \right)\left( B{{P}^{2}}-A{{P}^{2}} \right) \\ 
	\end{aligned} \right).\]
	
	
	Similarly, $\triangle AP{{P}_{{{D}_{H}}}}$, $\triangle BP{{P}_{{{D}_{H}}}}$ and $\triangle CP{{P}_{{{D}_{H}}}}$ are right triangles (see figure \ref{fig:tu28.1.1}), using theorem \ref{thm:Thm12.1.2} and Pythagorean theorem, we have:
	
	\[\left\{ \begin{aligned}
		& AP_{{{D}_{H}}}^{2}-BP_{{{D}_{H}}}^{2}=A{{P}^{2}}-B{{P}^{2}} \\ 
		& BP_{{{D}_{H}}}^{2}-CP_{{{D}_{H}}}^{2}=B{{P}^{2}}-C{{P}^{2}}. \\ 
	\end{aligned} \right.\]
	
		
	Add the condition $\alpha _{B}^{{{P}_{{{D}_{H}}}}}+\alpha _{C}^{{{P}_{{{D}_{H}}}}}+\alpha _{A}^{{{P}_{{{D}_{H}}}}}=1$, the following linear equations are obtained:
	
	\[\left\{ \begin{aligned}
		& \alpha _{B}^{{{P}_{{{D}_{H}}}}}+\alpha _{C}^{{{P}_{{{D}_{H}}}}}+\alpha _{A}^{{{P}_{{{D}_{H}}}}}=1 \\ 
		& A{{B}^{2}}\alpha _{B}^{{{P}_{{{D}_{H}}}}}+\left( C{{A}^{2}}-B{{C}^{2}} \right)\alpha _{C}^{{{P}_{{{D}_{H}}}}}-A{{B}^{2}}\alpha _{A}^{{{P}_{{{D}_{H}}}}}=A{{P}^{2}}-B{{P}^{2}} \\ 
		& -B{{C}^{2}}\alpha _{B}^{{{P}_{{{D}_{H}}}}}+B{{C}^{2}}\alpha _{C}^{{{P}_{{{D}_{H}}}}}+\left( A{{B}^{2}}-C{{A}^{2}} \right)\alpha _{A}^{{{P}_{{{D}_{H}}}}}=B{{P}^{2}}-C{{P}^{2}}. \\ 
	\end{aligned} \right.\]
	
	Written in matrix form is
	\[\left( \begin{matrix}
		1 & 1 & 1  \\
		A{{B}^{2}} & C{{A}^{2}}-B{{C}^{2}} & -A{{B}^{2}}  \\
		-B{{C}^{2}} & B{{C}^{2}} & A{{B}^{2}}-C{{A}^{2}}  \\
	\end{matrix} \right)\left( \begin{matrix}
		\alpha _{B}^{{{P}_{{{D}_{H}}}}}  \\
		\alpha _{C}^{{{P}_{{{D}_{H}}}}}  \\
		\alpha _{A}^{{{P}_{{{D}_{H}}}}}  \\
	\end{matrix} \right)=\left( \begin{matrix}
		1  \\
		A{{P}^{2}}-B{{P}^{2}}  \\
		B{{P}^{2}}-C{{P}^{2}}  \\
	\end{matrix} \right).\]
	i.e. 
	\[\left( \begin{matrix}
		\alpha _{B}^{{{P}_{{{D}_{H}}}}}  \\
		\alpha _{C}^{{{P}_{{{D}_{H}}}}}  \\
		\alpha _{A}^{{{P}_{{{D}_{H}}}}}  \\
	\end{matrix} \right)={{\left( \begin{matrix}
				1 & 1 & 1  \\
				A{{B}^{2}} & C{{A}^{2}}-B{{C}^{2}} & -A{{B}^{2}}  \\
				-B{{C}^{2}} & B{{C}^{2}} & A{{B}^{2}}-C{{A}^{2}}  \\
			\end{matrix} \right)}^{-1}}\left( \begin{matrix}
		1  \\
		A{{P}^{2}}-B{{P}^{2}}  \\
		B{{P}^{2}}-C{{P}^{2}}  \\
	\end{matrix} \right).\]
	
	By solving the inverse matrix, we have:
	\[\left( \begin{matrix}
		\alpha _{B}^{{{P}_{{{D}_{H}}}}}  \\
		\alpha _{C}^{{{P}_{{{D}_{H}}}}}  \\
		\alpha _{A}^{{{P}_{{{D}_{H}}}}}  \\
	\end{matrix} \right)=\frac{1}{8{{\left( {{S}^{D}} \right)}^{2}}}\left( \begin{matrix}
		\left( \Delta _{2}^{D}-C{{A}^{2}} \right)C{{A}^{2}} & \Delta _{2}^{D}-A{{B}^{2}} & B{{C}^{2}}-\Delta _{2}^{D}  \\
		\left( \Delta _{2}^{D}-A{{B}^{2}} \right)A{{B}^{2}} & \Delta _{2}^{D}-C{{A}^{2}} & 2A{{B}^{2}}  \\
		\left( \Delta _{2}^{D}-B{{C}^{2}} \right)B{{C}^{2}} & -2B{{C}^{2}} & C{{A}^{2}}-\Delta _{2}^{D}  \\
	\end{matrix} \right)\left( \begin{matrix}
		1  \\
		A{{P}^{2}}-B{{P}^{2}}  \\
		B{{P}^{2}}-C{{P}^{2}}  \\
	\end{matrix} \right).\]
	
	So one of the solutions is
	\[\alpha _{B}^{{{P}_{{{D}_{H}}}}}=\frac{1}{8{{\left( {{S}^{D}} \right)}^{2}}}\left( \begin{aligned}
		& \left( \Delta _{2}^{D}-C{{A}^{2}} \right)C{{A}^{2}}+\left( \Delta _{2}^{D}-A{{B}^{2}} \right)\left( A{{P}^{2}}-B{{P}^{2}} \right) \\ 
		& +\left( \Delta _{2}^{D}-B{{C}^{2}} \right)\left( C{{P}^{2}}-B{{P}^{2}} \right) \\ 
	\end{aligned} \right).\]
	
	
	Similarly, $\triangle AP{{P}_{{{D}_{H}}}}$, $\triangle BP{{P}_{{{D}_{H}}}}$ and $\triangle CP{{P}_{{{D}_{H}}}}$ are right triangles (see figure \ref{fig:tu28.1.1}), using theorem \ref{thm:Thm12.1.2} and Pythagorean theorem, we have:
	
	\[\left\{ \begin{aligned}
		& BP_{{{D}_{H}}}^{2}-CP_{{{D}_{H}}}^{2}=B{{P}^{2}}-C{{P}^{2}} \\ 
		& CP_{{{D}_{H}}}^{2}-AP_{{{D}_{H}}}^{2}=C{{P}^{2}}-A{{P}^{2}}. \\ 
	\end{aligned} \right.\]
	
	
	After the condition $\alpha _{C}^{{{P}_{{{D}_{H}}}}}+\alpha _{A}^{{{P}_{{{D}_{H}}}}}+\alpha _{B}^{{{P}_{{{D}_{H}}}}}=1$ is added, the following linear equations are obtained:
	
	\[\left\{ \begin{aligned}
		& \alpha _{C}^{{{P}_{{{D}_{H}}}}}+\alpha _{A}^{{{P}_{{{D}_{H}}}}}+\alpha _{B}^{{{P}_{{{D}_{H}}}}}=1 \\ 
		& B{{C}^{2}}\alpha _{C}^{{{P}_{{{D}_{H}}}}}+\left( A{{B}^{2}}-C{{A}^{2}} \right)\alpha _{A}^{{{P}_{{{D}_{H}}}}}-B{{C}^{2}}\alpha _{B}^{{{P}_{{{D}_{H}}}}}=B{{P}^{2}}-C{{P}^{2}} \\ 
		& -C{{A}^{2}}\alpha _{C}^{{{P}_{{{D}_{H}}}}}+C{{A}^{2}}\alpha _{A}^{{{P}_{{{D}_{H}}}}}+\left( B{{C}^{2}}-A{{B}^{2}} \right)\alpha _{B}^{{{P}_{{{D}_{H}}}}}=C{{P}^{2}}-A{{P}^{2}}. \\ 
	\end{aligned} \right.\]
	
	Written in matrix form is
	\[\left( \begin{matrix}
		1 & 1 & 1  \\
		B{{C}^{2}} & A{{B}^{2}}-C{{A}^{2}} & -B{{C}^{2}}  \\
		-C{{A}^{2}} & C{{A}^{2}} & B{{C}^{2}}-A{{B}^{2}}  \\
	\end{matrix} \right)\left( \begin{matrix}
		\alpha _{C}^{{{P}_{{{D}_{H}}}}}  \\
		\alpha _{A}^{{{P}_{{{D}_{H}}}}}  \\
		\alpha _{B}^{{{P}_{{{D}_{H}}}}}  \\
	\end{matrix} \right)=\left( \begin{matrix}
		1  \\
		B{{P}^{2}}-C{{P}^{2}}  \\
		C{{P}^{2}}-A{{P}^{2}}  \\
	\end{matrix} \right).\]
	i.e. 
	\[\left( \begin{matrix}
		\alpha _{C}^{{{P}_{{{D}_{H}}}}}  \\
		\alpha _{A}^{{{P}_{{{D}_{H}}}}}  \\
		\alpha _{B}^{{{P}_{{{D}_{H}}}}}  \\
	\end{matrix} \right)={{\left( \begin{matrix}
				1 & 1 & 1  \\
				B{{C}^{2}} & A{{B}^{2}}-C{{A}^{2}} & -B{{C}^{2}}  \\
				-C{{A}^{2}} & C{{A}^{2}} & B{{C}^{2}}-A{{B}^{2}}  \\
			\end{matrix} \right)}^{-1}}\left( \begin{matrix}
		1  \\
		B{{P}^{2}}-C{{P}^{2}}  \\
		C{{P}^{2}}-A{{P}^{2}}  \\
	\end{matrix} \right).\]
	
	By solving the inverse matrix, we have:
	\[\left( \begin{matrix}
		\alpha _{C}^{{{P}_{{{D}_{H}}}}}  \\
		\alpha _{A}^{{{P}_{{{D}_{H}}}}}  \\
		\alpha _{B}^{{{P}_{{{D}_{H}}}}}  \\
	\end{matrix} \right)=\frac{1}{8{{\left( {{S}^{D}} \right)}^{2}}}\left( \begin{matrix}
		\left( \Delta _{2}^{D}-A{{B}^{2}} \right)A{{B}^{2}} & \Delta _{2}^{D}-B{{C}^{2}} & C{{A}^{2}}-\Delta _{2}^{D}  \\
		\left( \Delta _{2}^{D}-B{{C}^{2}} \right)B{{C}^{2}} & \Delta _{2}^{D}-A{{B}^{2}} & 2B{{C}^{2}}  \\
		\left( \Delta _{2}^{D}-C{{A}^{2}} \right)C{{A}^{2}} & -2C{{A}^{2}} & A{{B}^{2}}-\Delta _{2}^{D}  \\
	\end{matrix} \right)\left( \begin{matrix}
		1  \\
		B{{P}^{2}}-C{{P}^{2}}  \\
		C{{P}^{2}}-A{{P}^{2}}  \\
	\end{matrix} \right).\]
	
	So one of the solutions is
	\[\alpha _{C}^{{{P}_{{{D}_{H}}}}}=\frac{1}{8{{\left( {{S}^{D}} \right)}^{2}}}\left( \begin{aligned}
		& \left( \Delta _{2}^{D}-A{{B}^{2}} \right)A{{B}^{2}}+\left( \Delta _{2}^{D}-B{{C}^{2}} \right)\left( B{{P}^{2}}-C{{P}^{2}} \right) \\ 
		& +\left( \Delta _{2}^{D}-C{{A}^{2}} \right)\left( A{{P}^{2}}-C{{P}^{2}} \right) \\ 
	\end{aligned} \right).\]
	
	Similar, we can prove other formulas.
\end{proof}
\hfill $\square$\par

\begin{theorem}{Formulas of frame components for projection of a point in space-2, Daiyuan Zhang}{Thm28.2.2}\label{Thm28.2.2} 
	Given a tetrahedron $ABCD$, point $O$ is any point in space, and point $P$ is a given point in space, suppose that the projection point of point $P$ to the $\triangle BCD$ plane of the tetrahedron is ${{P}_{{{A}_{H}}}}$, the projection point of point $P$ on the $\triangle CDA$ plane of tetrahedron is ${{P}_{{{B}_{H}}}}$, the projection point of point $P$ on the $\triangle DAB$ plane of tetrahedron is ${{P}_{{{C}_{H}}}}$, the projection point of point $P$ on the $\triangle ABC$ plane of tetrahedron is ${{P}_{{{D}_{H}}}}$. Then, on the tetrahedral frame $\left( O;A,B,C,D \right)$, there are the following conclusions:
	
	1. 	The frame components of ${{P}_{{{A}_{H}}}}$ are:
			
	\[\beta _{B}^{{{P}_{{{A}_{H}}}}}=\frac{1}{8{{\left( {{S}^{A}} \right)}^{2}}}\left( \begin{aligned}
		& \left( \Delta _{2}^{A}-C{{D}^{2}} \right)C{{D}^{2}}+\left( \Delta _{2}^{A}-D{{B}^{2}} \right)\left( D{{P}^{2}}-B{{P}^{2}} \right) \\ 
		& +\left( \Delta _{2}^{A}-B{{C}^{2}} \right)\left( C{{P}^{2}}-B{{P}^{2}} \right) \\ 
	\end{aligned} \right),\]
	\[\beta _{C}^{{{P}_{{{A}_{H}}}}}=\frac{1}{8{{\left( {{S}^{A}} \right)}^{2}}}\left( \begin{aligned}
		& \left( \Delta _{2}^{A}-D{{B}^{2}} \right)D{{B}^{2}}+\left( \Delta _{2}^{A}-B{{C}^{2}} \right)\left( B{{P}^{2}}-C{{P}^{2}} \right) \\ 
		& +\left( \Delta _{2}^{A}-C{{D}^{2}} \right)\left( D{{P}^{2}}-C{{P}^{2}} \right) \\ 
	\end{aligned} \right),\]
	\[\beta _{D}^{{{P}_{{{A}_{H}}}}}=\frac{1}{8{{\left( {{S}^{A}} \right)}^{2}}}\left( \begin{aligned}
		& \left( \Delta _{2}^{A}-B{{C}^{2}} \right)B{{C}^{2}}+\left( \Delta _{2}^{A}-C{{D}^{2}} \right)\left( C{{P}^{2}}-D{{P}^{2}} \right) \\ 
		& +\left( \Delta _{2}^{A}-D{{B}^{2}} \right)\left( B{{P}^{2}}-D{{P}^{2}} \right) \\ 
	\end{aligned} \right),\]
	\[\beta _{A}^{{{P}_{{{A}_{H}}}}}=0.\]
	Where 
	\[\Delta _{2}^{A}=\frac{1}{2}\left( B{{C}^{2}}+C{{D}^{2}}+D{{B}^{2}} \right).\]
	
	2. 	The frame components of ${{P}_{{{B}_{H}}}}$ are:
			
	\[\beta _{C}^{{{P}_{{{B}_{H}}}}}=\frac{1}{8{{\left( {{S}^{B}} \right)}^{2}}}\left( \begin{aligned}
		& \left( \Delta _{2}^{B}-D{{A}^{2}} \right)D{{A}^{2}}+\left( \Delta _{2}^{B}-A{{C}^{2}} \right)\left( A{{P}^{2}}-C{{P}^{2}} \right) \\ 
		& +\left( \Delta _{2}^{B}-C{{D}^{2}} \right)\left( D{{P}^{2}}-C{{P}^{2}} \right) \\ 
	\end{aligned} \right),\]
	\[\beta _{D}^{{{P}_{{{B}_{H}}}}}=\frac{1}{8{{\left( {{S}^{B}} \right)}^{2}}}\left( \begin{aligned}
		& \left( \Delta _{2}^{B}-A{{C}^{2}} \right)A{{C}^{2}}+\left( \Delta _{2}^{B}-C{{D}^{2}} \right)\left( C{{P}^{2}}-D{{P}^{2}} \right) \\ 
		& +\left( \Delta _{2}^{B}-D{{A}^{2}} \right)\left( A{{P}^{2}}-D{{P}^{2}} \right) \\ 
	\end{aligned} \right),\]
	\[\beta _{A}^{{{P}_{{{B}_{H}}}}}=\frac{1}{8{{\left( {{S}^{B}} \right)}^{2}}}\left( \begin{aligned}
		& \left( \Delta _{2}^{B}-C{{D}^{2}} \right)C{{D}^{2}}+\left( \Delta _{2}^{B}-D{{A}^{2}} \right)\left( D{{P}^{2}}-A{{P}^{2}} \right) \\ 
		& +\left( \Delta _{2}^{B}-A{{C}^{2}} \right)\left( C{{P}^{2}}-A{{P}^{2}} \right) \\ 
	\end{aligned} \right),\]
	\[\beta _{B}^{{{P}_{{{B}_{H}}}}}=0.\]
	Where 
	\[\Delta _{2}^{B}=\frac{1}{2}\left( C{{D}^{2}}+D{{A}^{2}}+A{{C}^{2}} \right).\]
	
	3. 	The frame components of ${{P}_{{{C}_{H}}}}$ are:
		
	\[\beta _{D}^{{{P}_{{{C}_{H}}}}}=\frac{1}{8{{\left( {{S}^{C}} \right)}^{2}}}\left( \begin{aligned}
		& \left( \Delta _{2}^{C}-A{{B}^{2}} \right)A{{B}^{2}}+\left( \Delta _{2}^{C}-B{{D}^{2}} \right)\left( B{{P}^{2}}-D{{P}^{2}} \right) \\ 
		& +\left( \Delta _{2}^{C}-D{{A}^{2}} \right)\left( A{{P}^{2}}-D{{P}^{2}} \right) \\ 
	\end{aligned} \right),\]
	\[\beta _{A}^{{{P}_{{{C}_{H}}}}}=\frac{1}{8{{\left( {{S}^{C}} \right)}^{2}}}\left( \begin{aligned}
		& \left( \Delta _{2}^{C}-B{{D}^{2}} \right)B{{D}^{2}}+\left( \Delta _{2}^{C}-D{{A}^{2}} \right)\left( D{{P}^{2}}-A{{P}^{2}} \right) \\ 
		& +\left( \Delta _{2}^{C}-A{{B}^{2}} \right)\left( B{{P}^{2}}-A{{P}^{2}} \right) \\ 
	\end{aligned} \right),\]
	\[\beta _{B}^{{{P}_{{{C}_{H}}}}}=\frac{1}{8{{\left( {{S}^{C}} \right)}^{2}}}\left( \begin{aligned}
		& \left( \Delta _{2}^{C}-D{{A}^{2}} \right)D{{A}^{2}}+\left( \Delta _{2}^{C}-A{{B}^{2}} \right)\left( A{{P}^{2}}-B{{P}^{2}} \right) \\ 
		& +\left( \Delta _{2}^{C}-B{{D}^{2}} \right)\left( D{{P}^{2}}-B{{P}^{2}} \right) \\ 
	\end{aligned} \right),\]
	\[\beta _{C}^{{{P}_{{{C}_{H}}}}}=0.\]
	Where 
	\[\Delta _{2}^{C}=\frac{1}{2}\left( D{{A}^{2}}+A{{B}^{2}}+B{{D}^{2}} \right).\]
	
	4. 	The frame components of ${{P}_{{{D}_{H}}}}$ are:
		
	\[\beta _{A}^{{{P}_{{{D}_{H}}}}}=\frac{1}{8{{\left( {{S}^{D}} \right)}^{2}}}\left( \begin{aligned}
		& \left( \Delta _{2}^{D}-B{{C}^{2}} \right)B{{C}^{2}}+\left( \Delta _{2}^{D}-C{{A}^{2}} \right)\left( C{{P}^{2}}-A{{P}^{2}} \right) \\ 
		& +\left( \Delta _{2}^{D}-A{{B}^{2}} \right)\left( B{{P}^{2}}-A{{P}^{2}} \right) \\ 
	\end{aligned} \right),\]
	\[\beta _{B}^{{{P}_{{{D}_{H}}}}}=\frac{1}{8{{\left( {{S}^{D}} \right)}^{2}}}\left( \begin{aligned}
		& \left( \Delta _{2}^{D}-C{{A}^{2}} \right)C{{A}^{2}}+\left( \Delta _{2}^{D}-A{{B}^{2}} \right)\left( A{{P}^{2}}-B{{P}^{2}} \right) \\ 
		& +\left( \Delta _{2}^{D}-B{{C}^{2}} \right)\left( C{{P}^{2}}-B{{P}^{2}} \right) \\ 
	\end{aligned} \right),\]
	\[\beta _{C}^{{{P}_{{{D}_{H}}}}}=\frac{1}{8{{\left( {{S}^{D}} \right)}^{2}}}\left( \begin{aligned}
		& \left( \Delta _{2}^{D}-A{{B}^{2}} \right)A{{B}^{2}}+\left( \Delta _{2}^{D}-B{{C}^{2}} \right)\left( B{{P}^{2}}-C{{P}^{2}} \right) \\ 
		& +\left( \Delta _{2}^{D}-C{{A}^{2}} \right)\left( A{{P}^{2}}-C{{P}^{2}} \right) \\ 
	\end{aligned} \right),\]
	\[\beta _{D}^{{{P}_{{{D}_{H}}}}}=0.\]
	Where 
	\[\Delta _{2}^{D}=\frac{1}{2}\left( A{{B}^{2}}+B{{C}^{2}}+C{{A}^{2}} \right).\]
\end{theorem}

\begin{proof}
	This theorem can be obtained directly according to theorem \ref{thm:Thm28.2.1} and theorem \ref{thm:Thm20.4.1}.
\end{proof}
\hfill $\square$\par

\section{Frame components of projection points of tetrahedral vertexs}\label{Sec28.3}


Given a tetrahedron $ABCD$, the four projection points of the four vertices of the tetrahedron to the plane opposite the corresponding vertex are ${{A}_{{{A}_{H}}}}$, ${{B}_{{{B}_{H}}}}$, ${{C}_{{{C}_{H}}}}$, ${{D}_{{{D}_{H}}}}$ respectively. How to calculate the frame components of these projection points? The following theorem I give answers this question.

\begin{theorem}{Formulas of frame components for projection of vertex-1, Daiyuan Zhang}{Thm28.3.1}\label{Thm28.3.1} 
	Given a tetrahedron $ABCD$, point $O$ is any point in space. Suppose that the projection point of point $A$ to the $\triangle BCD$ plane of the tetrahedron is ${{A}_{{{A}_{H}}}}$, the projection point of point $B$ on the $\triangle CDA$ plane of tetrahedron is ${{B}_{{{B}_{H}}}}$, the projection point of point $C$ to the $\triangle DAB$ plane of tetrahedron is ${{C}_{{{C}_{H}}}}$, the projection point of point $D$ on the $\triangle ABC$ plane of tetrahedron is${{D}_{{{D}_{H}}}}$, then:
	
	1. 	The frame components of ${{A}_{{{A}_{H}}}}$ on the triangular frame $\left( O;B,C,D \right)$ is:
	
	\[\alpha _{B}^{{{A}_{{{A}_{H}}}}}=\frac{1}{8{{\left( {{S}^{A}} \right)}^{2}}}\left( \begin{aligned}
		& \left( \Delta _{2}^{A}-C{{D}^{2}} \right)C{{D}^{2}}+\left( \Delta _{2}^{A}-D{{B}^{2}} \right)\left( D{{A}^{2}}-B{{A}^{2}} \right) \\ 
		& +\left( \Delta _{2}^{A}-B{{C}^{2}} \right)\left( C{{A}^{2}}-B{{A}^{2}} \right) \\ 
	\end{aligned} \right),\]
	\[\alpha _{C}^{{{A}_{{{A}_{H}}}}}=\frac{1}{8{{\left( {{S}^{A}} \right)}^{2}}}\left( \begin{aligned}
		& \left( \Delta _{2}^{A}-D{{B}^{2}} \right)D{{B}^{2}}+\left( \Delta _{2}^{A}-B{{C}^{2}} \right)\left( B{{A}^{2}}-C{{A}^{2}} \right) \\ 
		& +\left( \Delta _{2}^{A}-C{{D}^{2}} \right)\left( D{{A}^{2}}-C{{A}^{2}} \right) \\ 
	\end{aligned} \right),\]
	\[\alpha _{D}^{{{A}_{{{A}_{H}}}}}=\frac{1}{8{{\left( {{S}^{A}} \right)}^{2}}}\left( \begin{aligned}
		& \left( \Delta _{2}^{A}-B{{C}^{2}} \right)B{{C}^{2}}+\left( \Delta _{2}^{A}-C{{D}^{2}} \right)\left( C{{A}^{2}}-D{{A}^{2}} \right) \\ 
		& +\left( \Delta _{2}^{A}-D{{B}^{2}} \right)\left( B{{A}^{2}}-D{{A}^{2}} \right) \\ 
	\end{aligned} \right).\]
	Where 
	\[\Delta _{2}^{A}=\frac{1}{2}\left( B{{C}^{2}}+C{{D}^{2}}+D{{B}^{2}} \right).\]
	
	2. 	The frame components of ${{B}_{{{B}_{H}}}}$ on the triangular frame $\left( O;C,D,A \right)$ is:
	
	\[\alpha _{C}^{{{B}_{{{B}_{H}}}}}=\frac{1}{8{{\left( {{S}^{B}} \right)}^{2}}}\left( \begin{aligned}
		& \left( \Delta _{2}^{B}-D{{A}^{2}} \right)D{{A}^{2}}+\left( \Delta _{2}^{B}-A{{C}^{2}} \right)\left( A{{B}^{2}}-C{{B}^{2}} \right) \\ 
		& +\left( \Delta _{2}^{B}-C{{D}^{2}} \right)\left( D{{B}^{2}}-C{{B}^{2}} \right) \\ 
	\end{aligned} \right),\]
	\[\alpha _{D}^{{{B}_{{{B}_{H}}}}}=\frac{1}{8{{\left( {{S}^{B}} \right)}^{2}}}\left( \begin{aligned}
		& \left( \Delta _{2}^{B}-A{{C}^{2}} \right)A{{C}^{2}}+\left( \Delta _{2}^{B}-C{{D}^{2}} \right)\left( C{{B}^{2}}-D{{B}^{2}} \right) \\ 
		& +\left( \Delta _{2}^{B}-D{{A}^{2}} \right)\left( A{{B}^{2}}-D{{B}^{2}} \right) \\ 
	\end{aligned} \right),\]
	\[\alpha _{A}^{{{B}_{{{B}_{H}}}}}=\frac{1}{8{{\left( {{S}^{B}} \right)}^{2}}}\left( \begin{aligned}
		& \left( \Delta _{2}^{B}-C{{D}^{2}} \right)C{{D}^{2}}+\left( \Delta _{2}^{B}-D{{A}^{2}} \right)\left( D{{B}^{2}}-A{{B}^{2}} \right) \\ 
		& +\left( \Delta _{2}^{B}-A{{C}^{2}} \right)\left( C{{B}^{2}}-A{{B}^{2}} \right) \\ 
	\end{aligned} \right).\]
	Where 
	\[\Delta _{2}^{B}=\frac{1}{2}\left( C{{D}^{2}}+D{{A}^{2}}+A{{C}^{2}} \right).\]
	
	3. 	The frame components of ${{C}_{{{C}_{H}}}}$ on the triangular frame $\left( O;D,A,B \right)$ is:
	
	\[\alpha _{D}^{{{C}_{{{C}_{H}}}}}=\frac{1}{8{{\left( {{S}^{C}} \right)}^{2}}}\left( \begin{aligned}
		& \left( \Delta _{2}^{C}-A{{B}^{2}} \right)A{{B}^{2}}+\left( \Delta _{2}^{C}-B{{D}^{2}} \right)\left( B{{C}^{2}}-D{{C}^{2}} \right) \\ 
		& +\left( \Delta _{2}^{C}-D{{A}^{2}} \right)\left( A{{C}^{2}}-D{{C}^{2}} \right) \\ 
	\end{aligned} \right),\]
	\[\alpha _{A}^{{{C}_{{{C}_{H}}}}}=\frac{1}{8{{\left( {{S}^{C}} \right)}^{2}}}\left( \begin{aligned}
		& \left( \Delta _{2}^{C}-B{{D}^{2}} \right)B{{D}^{2}}+\left( \Delta _{2}^{C}-D{{A}^{2}} \right)\left( D{{C}^{2}}-A{{C}^{2}} \right) \\ 
		& +\left( \Delta _{2}^{C}-A{{B}^{2}} \right)\left( B{{C}^{2}}-A{{C}^{2}} \right) \\ 
	\end{aligned} \right),\]
	\[\alpha _{B}^{{{C}_{{{C}_{H}}}}}=\frac{1}{8{{\left( {{S}^{C}} \right)}^{2}}}\left( \begin{aligned}
		& \left( \Delta _{2}^{C}-D{{A}^{2}} \right)D{{A}^{2}}+\left( \Delta _{2}^{C}-A{{B}^{2}} \right)\left( A{{C}^{2}}-B{{C}^{2}} \right) \\ 
		& +\left( \Delta _{2}^{C}-B{{D}^{2}} \right)\left( D{{C}^{2}}-B{{C}^{2}} \right) \\ 
	\end{aligned} \right).\]
	Where 
	\[\Delta _{2}^{C}=\frac{1}{2}\left( D{{A}^{2}}+A{{B}^{2}}+B{{D}^{2}} \right).\]
	
	4. 	The frame components of ${{D}_{{{D}_{H}}}}$ on the triangular frame $\left( O;A,B,C \right)$ is:
	
	\[\alpha _{A}^{{{D}_{{{D}_{H}}}}}=\frac{1}{8{{\left( {{S}^{D}} \right)}^{2}}}\left( \begin{aligned}
		& \left( \Delta _{2}^{D}-B{{C}^{2}} \right)B{{C}^{2}}+\left( \Delta _{2}^{D}-C{{A}^{2}} \right)\left( C{{D}^{2}}-A{{D}^{2}} \right) \\ 
		& +\left( \Delta _{2}^{D}-A{{B}^{2}} \right)\left( B{{D}^{2}}-A{{D}^{2}} \right) \\ 
	\end{aligned} \right),\]
	\[\alpha _{B}^{{{D}_{{{D}_{H}}}}}=\frac{1}{8{{\left( {{S}^{D}} \right)}^{2}}}\left( \begin{aligned}
		& \left( \Delta _{2}^{D}-C{{A}^{2}} \right)C{{A}^{2}}+\left( \Delta _{2}^{D}-A{{B}^{2}} \right)\left( A{{D}^{2}}-B{{D}^{2}} \right) \\ 
		& +\left( \Delta _{2}^{D}-B{{C}^{2}} \right)\left( C{{D}^{2}}-B{{D}^{2}} \right) \\ 
	\end{aligned} \right),\]
	\[\alpha _{C}^{{{D}_{{{D}_{H}}}}}=\frac{1}{8{{\left( {{S}^{D}} \right)}^{2}}}\left( \begin{aligned}
		& \left( \Delta _{2}^{D}-A{{B}^{2}} \right)A{{B}^{2}}+\left( \Delta _{2}^{D}-B{{C}^{2}} \right)\left( B{{D}^{2}}-C{{D}^{2}} \right) \\ 
		& +\left( \Delta _{2}^{D}-C{{A}^{2}} \right)\left( A{{D}^{2}}-C{{D}^{2}} \right) \\ 
	\end{aligned} \right).\]
	Where 
	\[\Delta _{2}^{D}=\frac{1}{2}\left( A{{B}^{2}}+B{{C}^{2}}+C{{A}^{2}} \right).\]
\end{theorem}

\begin{proof}
	According to theorem \ref{thm:Thm28.2.1} and theorem \ref{thm:Thm20.4.1}, the conclusion of the theorem can be obtained by replacing the point $P$ with points $A$, $B$ ,$C$, $D$ respectively.
\end{proof}
\hfill $\square$\par
\begin{theorem}{Formulas of frame components for projection of vertex-2, Daiyuan Zhang}{Thm28.3.2}\label{Thm28.3.2} 
	Given a tetrahedron $ABCD$, point $O$ is any point in space. Suppose that the projection point of point $A$ to the $\triangle BCD$ plane of the tetrahedron is ${{A}_{{{A}_{H}}}}$, the projection point of point $B$ on the $\triangle CDA$ plane of tetrahedron is ${{B}_{{{B}_{H}}}}$, the projection point of point $C$ to the $\triangle DAB$ plane of tetrahedron is ${{C}_{{{C}_{H}}}}$, the projection point of point $D$ on the $\triangle ABC$ plane of tetrahedron is ${{D}_{{{D}_{H}}}}$, then on the tetrahedral frame $\left( O;A,B,C,D \right)$, there are the following conclusions:
	
	1. 	The frame components of ${{A}_{{{A}_{H}}}}$ are:		
	\[\beta _{B}^{{{A}_{{{A}_{H}}}}}=\frac{1}{8{{\left( {{S}^{A}} \right)}^{2}}}\left( \begin{aligned}
		& \left( \Delta _{2}^{A}-C{{D}^{2}} \right)C{{D}^{2}}+\left( \Delta _{2}^{A}-D{{B}^{2}} \right)\left( D{{A}^{2}}-B{{A}^{2}} \right) \\ 
		& +\left( \Delta _{2}^{A}-B{{C}^{2}} \right)\left( C{{A}^{2}}-B{{A}^{2}} \right) \\ 
	\end{aligned} \right),\]
	\[\beta _{C}^{{{A}_{{{A}_{H}}}}}=\frac{1}{8{{\left( {{S}^{A}} \right)}^{2}}}\left( \begin{aligned}
		& \left( \Delta _{2}^{A}-D{{B}^{2}} \right)D{{B}^{2}}+\left( \Delta _{2}^{A}-B{{C}^{2}} \right)\left( B{{A}^{2}}-C{{A}^{2}} \right) \\ 
		& +\left( \Delta _{2}^{A}-C{{D}^{2}} \right)\left( D{{A}^{2}}-C{{A}^{2}} \right) \\ 
	\end{aligned} \right),\]
	\[\beta _{D}^{{{A}_{{{A}_{H}}}}}=\frac{1}{8{{\left( {{S}^{A}} \right)}^{2}}}\left( \begin{aligned}
		& \left( \Delta _{2}^{A}-B{{C}^{2}} \right)B{{C}^{2}}+\left( \Delta _{2}^{A}-C{{D}^{2}} \right)\left( C{{A}^{2}}-D{{A}^{2}} \right) \\ 
		& +\left( \Delta _{2}^{A}-D{{B}^{2}} \right)\left( B{{A}^{2}}-D{{A}^{2}} \right) \\ 
	\end{aligned} \right),\]
	\[\beta _{A}^{{{A}_{{{A}_{H}}}}}=0.\]
	Where 
	\[\Delta _{2}^{A}=\frac{1}{2}\left( B{{C}^{2}}+C{{D}^{2}}+D{{B}^{2}} \right).\]
	
	2. 	The frame components of ${{B}_{{{B}_{H}}}}$ are:
			
	\[\beta _{C}^{{{B}_{{{B}_{H}}}}}=\frac{1}{8{{\left( {{S}^{B}} \right)}^{2}}}\left( \begin{aligned}
		& \left( \Delta _{2}^{B}-D{{A}^{2}} \right)D{{A}^{2}}+\left( \Delta _{2}^{B}-A{{C}^{2}} \right)\left( A{{B}^{2}}-C{{B}^{2}} \right) \\ 
		& +\left( \Delta _{2}^{B}-C{{D}^{2}} \right)\left( D{{B}^{2}}-C{{B}^{2}} \right) \\ 
	\end{aligned} \right),\]
	\[\beta _{D}^{{{B}_{{{B}_{H}}}}}=\frac{1}{8{{\left( {{S}^{B}} \right)}^{2}}}\left( \begin{aligned}
		& \left( \Delta _{2}^{B}-A{{C}^{2}} \right)A{{C}^{2}}+\left( \Delta _{2}^{B}-C{{D}^{2}} \right)\left( C{{B}^{2}}-D{{B}^{2}} \right) \\ 
		& +\left( \Delta _{2}^{B}-D{{A}^{2}} \right)\left( A{{B}^{2}}-D{{B}^{2}} \right) \\ 
	\end{aligned} \right),\]
	\[\beta _{A}^{{{B}_{{{B}_{H}}}}}=\frac{1}{8{{\left( {{S}^{B}} \right)}^{2}}}\left( \begin{aligned}
		& \left( \Delta _{2}^{B}-C{{D}^{2}} \right)C{{D}^{2}}+\left( \Delta _{2}^{B}-D{{A}^{2}} \right)\left( D{{B}^{2}}-A{{B}^{2}} \right) \\ 
		& +\left( \Delta _{2}^{B}-A{{C}^{2}} \right)\left( C{{B}^{2}}-A{{B}^{2}} \right) \\ 
	\end{aligned} \right),\]
	\[\beta _{B}^{{{B}_{{{B}_{H}}}}}=0.\]
	Where 
	\[\Delta _{2}^{B}=\frac{1}{2}\left( C{{D}^{2}}+D{{A}^{2}}+A{{C}^{2}} \right).\]
	
	3. The frame components of ${{C}_{{{C}_{H}}}}$ are:
			
	\[\beta _{D}^{{{C}_{{{C}_{H}}}}}=\frac{1}{8{{\left( {{S}^{C}} \right)}^{2}}}\left( \begin{aligned}
		& \left( \Delta _{2}^{C}-A{{B}^{2}} \right)A{{B}^{2}}+\left( \Delta _{2}^{C}-B{{D}^{2}} \right)\left( B{{C}^{2}}-D{{C}^{2}} \right) \\ 
		& +\left( \Delta _{2}^{C}-D{{A}^{2}} \right)\left( A{{C}^{2}}-D{{C}^{2}} \right) \\ 
	\end{aligned} \right),\]
	\[\beta _{A}^{{{C}_{{{C}_{H}}}}}=\frac{1}{8{{\left( {{S}^{C}} \right)}^{2}}}\left( \begin{aligned}
		& \left( \Delta _{2}^{C}-B{{D}^{2}} \right)B{{D}^{2}}+\left( \Delta _{2}^{C}-D{{A}^{2}} \right)\left( D{{C}^{2}}-A{{C}^{2}} \right) \\ 
		& +\left( \Delta _{2}^{C}-A{{B}^{2}} \right)\left( B{{C}^{2}}-A{{C}^{2}} \right) \\ 
	\end{aligned} \right),\]
	\[\beta _{B}^{{{C}_{{{C}_{H}}}}}=\frac{1}{8{{\left( {{S}^{C}} \right)}^{2}}}\left( \begin{aligned}
		& \left( \Delta _{2}^{C}-D{{A}^{2}} \right)D{{A}^{2}}+\left( \Delta _{2}^{C}-A{{B}^{2}} \right)\left( A{{C}^{2}}-B{{C}^{2}} \right) \\ 
		& +\left( \Delta _{2}^{C}-B{{D}^{2}} \right)\left( D{{C}^{2}}-B{{C}^{2}} \right) \\ 
	\end{aligned} \right),\]
	\[\beta _{C}^{{{C}_{{{C}_{H}}}}}=0.\]
	Where 
	\[\Delta _{2}^{C}=\frac{1}{2}\left( D{{A}^{2}}+A{{B}^{2}}+B{{D}^{2}} \right).\]
	
	4. The frame components of ${{D}_{{{D}_{H}}}}$ are:		
	\[\beta _{A}^{{{D}_{{{D}_{H}}}}}=\frac{1}{8{{\left( {{S}^{D}} \right)}^{2}}}\left( \begin{aligned}
		& \left( \Delta _{2}^{D}-B{{C}^{2}} \right)B{{C}^{2}}+\left( \Delta _{2}^{D}-C{{A}^{2}} \right)\left( C{{D}^{2}}-A{{D}^{2}} \right) \\ 
		& +\left( \Delta _{2}^{D}-A{{B}^{2}} \right)\left( B{{D}^{2}}-A{{D}^{2}} \right) \\ 
	\end{aligned} \right),\]
	\[\beta _{B}^{{{D}_{{{D}_{H}}}}}=\frac{1}{8{{\left( {{S}^{D}} \right)}^{2}}}\left( \begin{aligned}
		& \left( \Delta _{2}^{D}-C{{A}^{2}} \right)C{{A}^{2}}+\left( \Delta _{2}^{D}-A{{B}^{2}} \right)\left( A{{D}^{2}}-B{{D}^{2}} \right) \\ 
		& +\left( \Delta _{2}^{D}-B{{C}^{2}} \right)\left( C{{D}^{2}}-B{{D}^{2}} \right) \\ 
	\end{aligned} \right),\]
	\[\beta _{C}^{{{D}_{{{D}_{H}}}}}=\frac{1}{8{{\left( {{S}^{D}} \right)}^{2}}}\left( \begin{aligned}
		& \left( \Delta _{2}^{D}-A{{B}^{2}} \right)A{{B}^{2}}+\left( \Delta _{2}^{D}-B{{C}^{2}} \right)\left( B{{D}^{2}}-C{{D}^{2}} \right) \\ 
		& +\left( \Delta _{2}^{D}-C{{A}^{2}} \right)\left( A{{D}^{2}}-C{{D}^{2}} \right) \\ 
	\end{aligned} \right),\]
	\[\beta _{D}^{{{D}_{{{D}_{H}}}}}=0.\]
	Where 
	\[\Delta _{2}^{D}=\frac{1}{2}\left( A{{B}^{2}}+B{{C}^{2}}+C{{A}^{2}} \right).\]
\end{theorem}

\begin{proof}
	According to theorem \ref{thm:Thm28.2.2} and theorem \ref{thm:Thm20.4.1}, the conclusion of the theorem can be obtained by replacing the point $P$ with points $A$, $B$ ,$C$, $D$ respectively.	 
\end{proof}
\hfill $\square$\par

\section{The relationship between the frame componentS of the projection point and the IC-T}\label{Sec28.4}


Given a tetrahedron $ABCD$, let $P$ be an IC-T. The point $P$ (IC-T) is projected to the four triangular faces, and four projection points can be obtained. The purpose is to calculate the frame components of projection points. I give the following theorem.

\begin{theorem}{Relationship 1 of frame components between the projection and IC-T, Daiyuan Zhang}{Thm28.4.1}\label{Thm28.4.1} 
	Given a tetrahedron $ABCD$, point $O$ is any point in space. Suppose that the projection point of point $A$ to the $\triangle BCD$ plane of the tetrahedron is ${{A}_{{{A}_{H}}}}$, the projection point of point $B$ on the $\triangle CDA$ plane of tetrahedron is ${{B}_{{{B}_{H}}}}$, the projection point of point $C$ to the $\triangle DAB$ plane of tetrahedron is ${{C}_{{{C}_{H}}}}$, the projection point of point $D$ on the $\triangle ABC$ plane of tetrahedron is ${{D}_{{{D}_{H}}}}$. The projection points of IC-T $P$ to the planes $\triangle BCD$, $\triangle CDA$, $\triangle DAB$, $\triangle ABC$ of the tetrahedron are ${{P}_{{{A}_{H}}}}$, ${{P}_{{{B}_{H}}}}$, ${{P}_{{{C}_{H}}}}$, ${{P}_{{{D}_{H}}}}$ respectively, then:	
	
	
	1. If ${{P}_{{{D}_{H}}}}\in {{\pi }_{ABC}}$, ${{D}_{{{D}_{H}}}}\in {{\pi }_{ABC}}$, then the frame components of ${{P}_{{{D}_{H}}}}$ in the frame $\left( O;A,B,C \right)$ is
	
	\[\left( \begin{matrix}
		\alpha _{A}^{{{P}_{{{D}_{H}}}}}  \\
		\alpha _{B}^{{{P}_{{{D}_{H}}}}}  \\
		\alpha _{C}^{{{P}_{{{D}_{H}}}}}  \\
	\end{matrix} \right)=\left( \begin{matrix}
		1 & 0 & 0 & \alpha _{A}^{{{D}_{{{D}_{H}}}}}  \\
		0 & 1 & 0 & \alpha _{B}^{{{D}_{{{D}_{H}}}}}  \\
		0 & 0 & 1 & \alpha _{C}^{{{D}_{{{D}_{H}}}}}  \\
	\end{matrix} \right)\left( \begin{matrix}
		\beta _{A}^{P}  \\
		\beta _{B}^{P}  \\
		\beta _{C}^{P}  \\
		\beta _{D}^{P}  \\
	\end{matrix} \right).\]	
	Where 
	\[\left( \begin{matrix}
		\alpha _{A}^{{{D}_{{{D}_{H}}}}}  \\
		\alpha _{B}^{{{D}_{{{D}_{H}}}}}  \\
		\alpha _{C}^{{{D}_{{{D}_{H}}}}}  \\
	\end{matrix} \right)={{\left( \begin{matrix}
				1 & 1 & 1  \\
				C{{A}^{2}} & B{{C}^{2}}-A{{B}^{2}} & -C{{A}^{2}}  \\
				-A{{B}^{2}} & A{{B}^{2}} & C{{A}^{2}}-B{{C}^{2}}  \\
			\end{matrix} \right)}^{-1}}\left( \begin{matrix}
		1  \\
		C{{D}^{2}}-A{{D}^{2}}  \\
		A{{D}^{2}}-B{{D}^{2}}  \\
	\end{matrix} \right).\]
	
	
	2. If ${{P}_{{{A}_{H}}}}\in {{\pi }_{BCD}}$, ${{A}_{{{A}_{H}}}}\in {{\pi }_{BCD}}$, then the frame components of ${{P}_{{{A}_{H}}}}$ in the frame $\left( O;B,C,D \right)$ is
	
	\[\left( \begin{matrix}
		\alpha _{B}^{{{P}_{{{A}_{H}}}}}  \\
		\alpha _{C}^{{{P}_{{{A}_{H}}}}}  \\
		\alpha _{D}^{{{P}_{{{A}_{H}}}}}  \\
	\end{matrix} \right)=\left( \begin{matrix}
		1 & 0 & 0 & \alpha _{B}^{{{A}_{{{A}_{H}}}}}  \\
		0 & 1 & 0 & \alpha _{C}^{{{A}_{{{A}_{H}}}}}  \\
		0 & 0 & 1 & \alpha _{D}^{{{A}_{{{A}_{H}}}}}  \\
	\end{matrix} \right)\left( \begin{matrix}
		\beta _{B}^{P}  \\
		\beta _{C}^{P}  \\
		\beta _{D}^{P}  \\
		\beta _{A}^{P}  \\
	\end{matrix} \right).\]	
	Where 
	\[\left( \begin{matrix}
		\alpha _{B}^{{{A}_{{{A}_{H}}}}}  \\
		\alpha _{C}^{{{A}_{{{A}_{H}}}}}  \\
		\alpha _{D}^{{{A}_{{{A}_{H}}}}}  \\
	\end{matrix} \right)={{\left( \begin{matrix}
				1 & 1 & 1  \\
				D{{B}^{2}} & C{{D}^{2}}-B{{C}^{2}} & -D{{B}^{2}}  \\
				-B{{C}^{2}} & B{{C}^{2}} & D{{B}^{2}}-C{{D}^{2}}  \\
			\end{matrix} \right)}^{-1}}\left( \begin{matrix}
		1  \\
		D{{A}^{2}}-B{{A}^{2}}  \\
		B{{A}^{2}}-C{{A}^{2}}  \\
	\end{matrix} \right).\]
	
	
	3. If ${{P}_{{{B}_{H}}}}\in {{\pi }_{CDA}}$, ${{B}_{{{B}_{H}}}}\in {{\pi }_{CDA}}$, then the frame components of ${{P}_{{{B}_{H}}}}$ in the frame $\left( O;C,D,A \right)$ is
	
	\[\left( \begin{matrix}
		\alpha _{C}^{{{P}_{{{B}_{H}}}}}  \\
		\alpha _{D}^{{{P}_{{{B}_{H}}}}}  \\
		\alpha _{A}^{{{P}_{{{B}_{H}}}}}  \\
	\end{matrix} \right)=\left( \begin{matrix}
		1 & 0 & 0 & \alpha _{C}^{{{B}_{{{B}_{H}}}}}  \\
		0 & 1 & 0 & \alpha _{D}^{{{B}_{{{B}_{H}}}}}  \\
		0 & 0 & 1 & \alpha _{A}^{{{B}_{{{B}_{H}}}}}  \\
	\end{matrix} \right)\left( \begin{matrix}
		\beta _{C}^{P}  \\
		\beta _{D}^{P}  \\
		\beta _{A}^{P}  \\
		\beta _{B}^{P}  \\
	\end{matrix} \right).\]	
	Where 
	\[\left( \begin{matrix}
		\alpha _{C}^{{{B}_{{{B}_{H}}}}}  \\
		\alpha _{D}^{{{B}_{{{B}_{H}}}}}  \\
		\alpha _{A}^{{{B}_{{{B}_{H}}}}}  \\
	\end{matrix} \right)={{\left( \begin{matrix}
				1 & 1 & 1  \\
				A{{C}^{2}} & D{{A}^{2}}-C{{D}^{2}} & -A{{C}^{2}}  \\
				-C{{D}^{2}} & C{{D}^{2}} & A{{C}^{2}}-D{{A}^{2}}  \\
			\end{matrix} \right)}^{-1}}\left( \begin{matrix}
		1  \\
		A{{B}^{2}}-C{{B}^{2}}  \\
		C{{B}^{2}}-D{{B}^{2}}  \\
	\end{matrix} \right).\]
	
	
	4. If ${{P}_{{{C}_{H}}}}\in {{\pi }_{DAB}}$,  ${{C}_{{{C}_{H}}}}\in {{\pi }_{DAB}}$, then the frame components of ${{P}_{{{C}_{H}}}}$ in the frame $\left( O;D,A,B \right)$ is
	
	\[\left( \begin{matrix}
		\alpha _{D}^{{{P}_{{{C}_{H}}}}}  \\
		\alpha _{A}^{{{P}_{{{C}_{H}}}}}  \\
		\alpha _{B}^{{{P}_{{{C}_{H}}}}}  \\
	\end{matrix} \right)=\left( \begin{matrix}
		1 & 0 & 0 & \alpha _{D}^{{{C}_{{{C}_{H}}}}}  \\
		0 & 1 & 0 & \alpha _{A}^{{{C}_{{{C}_{H}}}}}  \\
		0 & 0 & 1 & \alpha _{B}^{{{C}_{{{C}_{H}}}}}  \\
	\end{matrix} \right)\left( \begin{matrix}
		\beta _{D}^{P}  \\
		\beta _{A}^{P}  \\
		\beta _{B}^{P}  \\
		\beta _{C}^{P}  \\
	\end{matrix} \right).\]	
	Where 
	\[\left( \begin{matrix}
		\alpha _{D}^{{{C}_{{{C}_{H}}}}}  \\
		\alpha _{A}^{{{C}_{{{C}_{H}}}}}  \\
		\alpha _{B}^{{{C}_{{{C}_{H}}}}}  \\
	\end{matrix} \right)={{\left( \begin{matrix}
				1 & 1 & 1  \\
				B{{D}^{2}} & A{{B}^{2}}-D{{A}^{2}} & -B{{D}^{2}}  \\
				-D{{A}^{2}} & D{{A}^{2}} & B{{D}^{2}}-A{{B}^{2}}  \\
			\end{matrix} \right)}^{-1}}\left( \begin{matrix}
		1  \\
		B{{C}^{2}}-D{{C}^{2}}  \\
		D{{C}^{2}}-A{{C}^{2}}  \\
	\end{matrix} \right).\]
\end{theorem}

\begin{proof}
	
	Without losing generality, here calculate the frame components of the projection point ${{P}_{{{D}_{H}}}}$ of point $P$ to the $\triangle ABC$ plane of the tetrahedron. According to theorem \ref{thm:Thm24.1.2}, we have:	
	\[\begin{aligned}
		A{{P}^{2}}-B{{P}^{2}}& =\left( \beta _{B}^{P}A{{B}^{2}}+\beta _{C}^{P}A{{C}^{2}}+\beta _{D}^{P}A{{D}^{2}} \right)-\left( \beta _{C}^{P}B{{C}^{2}}+\beta _{D}^{P}B{{D}^{2}}+\beta _{A}^{P}B{{A}^{2}} \right) \\ 
		& =-B{{A}^{2}}\beta _{A}^{P}+A{{B}^{2}}\beta _{B}^{P}+\left( A{{C}^{2}}-B{{C}^{2}} \right)\beta _{C}^{P}+\left( A{{D}^{2}}-B{{D}^{2}} \right)\beta _{D}^{P},  
	\end{aligned}\]\[\begin{aligned}
		B{{P}^{2}}-C{{P}^{2}}& =\left( \beta _{C}^{P}B{{C}^{2}}+\beta _{D}^{P}B{{D}^{2}}+\beta _{A}^{P}B{{A}^{2}} \right)-\left( \beta _{D}^{P}C{{D}^{2}}+\beta _{A}^{P}C{{A}^{2}}+\beta _{B}^{P}C{{B}^{2}} \right) \\ 
		& =\left( B{{A}^{2}}-C{{A}^{2}} \right)\beta _{A}^{P}-C{{B}^{2}}\beta _{B}^{P}+B{{C}^{2}}\beta _{C}^{P}+\left( B{{D}^{2}}-C{{D}^{2}} \right)\beta _{D}^{P},  
	\end{aligned}\]\[\begin{aligned}
		C{{P}^{2}}-D{{P}^{2}}& =\left( \beta _{D}^{P}C{{D}^{2}}+\beta _{A}^{P}C{{A}^{2}}+\beta _{B}^{P}C{{B}^{2}} \right)-\left( \beta _{A}^{P}D{{A}^{2}}+\beta _{B}^{P}D{{B}^{2}}+\beta _{C}^{P}D{{C}^{2}} \right) \\ 
		& =\left( C{{A}^{2}}-D{{A}^{2}} \right)\beta _{A}^{P}+\left( C{{B}^{2}}-D{{B}^{2}} \right)\beta _{B}^{P}-D{{C}^{2}}\beta _{C}^{P}+C{{D}^{2}}\beta _{D}^{P}.  
	\end{aligned}\]
	
	
	Thus, the condition $\beta _{A}^{P}+\beta _{B}^{P}+\beta _{C}^{P}+\beta _{D}^{P}=1$ yields the following result:
	\begin{equation}\label{Eq28.4.1}
		\left( \begin{matrix}
			1  \\
			A{{P}^{2}}-B{{P}^{2}}  \\
			B{{P}^{2}}-C{{P}^{2}}  \\
			C{{P}^{2}}-D{{P}^{2}}  \\
		\end{matrix} \right)=\left( \begin{matrix}
			1 & 1 & 1 & 1  \\
			-B{{A}^{2}} & A{{B}^{2}} & A{{C}^{2}}-B{{C}^{2}} & A{{D}^{2}}-B{{D}^{2}}  \\
			B{{A}^{2}}-C{{A}^{2}} & -C{{B}^{2}} & B{{C}^{2}} & B{{D}^{2}}-C{{D}^{2}}  \\
			C{{A}^{2}}-D{{A}^{2}} & C{{B}^{2}}-D{{B}^{2}} & -D{{C}^{2}} & C{{D}^{2}}  \\
		\end{matrix} \right)\left( \begin{matrix}
			\beta _{A}^{P}  \\
			\beta _{B}^{P}  \\
			\beta _{C}^{P}  \\
			\beta _{D}^{P}  \\
		\end{matrix} \right).	
	\end{equation}
	
	
	The following formula appears in the proof process of theorem \ref{thm:Thm28.2.1}:
	\[\left( \begin{matrix}
		\alpha _{A}^{{{P}_{{{D}_{H}}}}}  \\
		\alpha _{B}^{{{P}_{{{D}_{H}}}}}  \\
		\alpha _{C}^{{{P}_{{{D}_{H}}}}}  \\
	\end{matrix} \right)={{\left( \begin{matrix}
				1 & 1 & 1  \\
				C{{A}^{2}} & B{{C}^{2}}-A{{B}^{2}} & -C{{A}^{2}}  \\
				-A{{B}^{2}} & A{{B}^{2}} & C{{A}^{2}}-B{{C}^{2}}  \\
			\end{matrix} \right)}^{-1}}\left( \begin{matrix}
		1  \\
		C{{P}^{2}}-A{{P}^{2}}  \\
		A{{P}^{2}}-B{{P}^{2}}  \\
	\end{matrix} \right).\]
	
	
	According to the formula (\ref{Eq28.4.1}), the following result is obtained:
	\[\left( \begin{matrix}
		1  \\
		C{{P}^{2}}-A{{P}^{2}}  \\
		A{{P}^{2}}-B{{P}^{2}}  \\
	\end{matrix} \right)=\left( \begin{matrix}
		1 & 1 & 1 & 1  \\
		C{{A}^{2}} & B{{C}^{2}}-A{{B}^{2}} & -C{{A}^{2}} & C{{D}^{2}}-A{{D}^{2}}  \\
		-B{{A}^{2}} & A{{B}^{2}} & A{{C}^{2}}-B{{C}^{2}} & A{{D}^{2}}-B{{D}^{2}}  \\
	\end{matrix} \right)\left( \begin{matrix}
		\beta _{A}^{P}  \\
		\beta _{B}^{P}  \\
		\beta _{C}^{P}  \\
		\beta _{D}^{P}  \\
	\end{matrix} \right).\]
	
	
	Thus, we obtaine the relationship of the frame components  between the IC-Fs related to the projection point ${{P}_{{{D}_{H}}}}$ of the IC-T $P$ on the $\triangle ABC$ plane and the IC-T of $P$. 
	
	\begin{equation}\label{Eq28.4.2}
		\begin{aligned}
			& \left( \begin{matrix}
				\alpha _{A}^{{{P}_{{{D}_{H}}}}}  \\
				\alpha _{B}^{{{P}_{{{D}_{H}}}}}  \\
				\alpha _{C}^{{{P}_{{{D}_{H}}}}}  \\
			\end{matrix} \right)={{\left( \begin{matrix}
						1 & 1 & 1  \\
						C{{A}^{2}} & B{{C}^{2}}-A{{B}^{2}} & -C{{A}^{2}}  \\
						-A{{B}^{2}} & A{{B}^{2}} & C{{A}^{2}}-B{{C}^{2}}  \\
					\end{matrix} \right)}^{-1}}\times  \\ 
			& \left( \begin{matrix}
				1 & 1 & 1 & 1  \\
				C{{A}^{2}} & B{{C}^{2}}-A{{B}^{2}} & -C{{A}^{2}} & C{{D}^{2}}-A{{D}^{2}}  \\
				-B{{A}^{2}} & A{{B}^{2}} & A{{C}^{2}}-B{{C}^{2}} & A{{D}^{2}}-B{{D}^{2}}  \\
			\end{matrix} \right)\left( \begin{matrix}
				\beta _{A}^{P}  \\
				\beta _{B}^{P}  \\
				\beta _{C}^{P}  \\
				\beta _{D}^{P}  \\
			\end{matrix} \right). \\ 
		\end{aligned}
	\end{equation}
	
	The following result can be obtained by block multiplication of matrix:
	\[\left( \begin{matrix}
		\alpha _{A}^{{{P}_{{{D}_{H}}}}}  \\
		\alpha _{B}^{{{P}_{{{D}_{H}}}}}  \\
		\alpha _{C}^{{{P}_{{{D}_{H}}}}}  \\
	\end{matrix} \right)=\left( \begin{matrix}
		1 & 0 & 0 & \alpha _{A}^{{{D}_{{{D}_{H}}}}}  \\
		0 & 1 & 0 & \alpha _{B}^{{{D}_{{{D}_{H}}}}}  \\
		0 & 0 & 1 & \alpha _{C}^{{{D}_{{{D}_{H}}}}}  \\
	\end{matrix} \right)\left( \begin{matrix}
		\beta _{A}^{P}  \\
		\beta _{B}^{P}  \\
		\beta _{C}^{P}  \\
		\beta _{D}^{P}  \\
	\end{matrix} \right).\]	
	Where 
	\[\left( \begin{matrix}
		\alpha _{A}^{{{D}_{{{D}_{H}}}}}  \\
		\alpha _{B}^{{{D}_{{{D}_{H}}}}}  \\
		\alpha _{C}^{{{D}_{{{D}_{H}}}}}  \\
	\end{matrix} \right)={{\left( \begin{matrix}
				1 & 1 & 1  \\
				C{{A}^{2}} & B{{C}^{2}}-A{{B}^{2}} & -C{{A}^{2}}  \\
				-A{{B}^{2}} & A{{B}^{2}} & C{{A}^{2}}-B{{C}^{2}}  \\
			\end{matrix} \right)}^{-1}}\left( \begin{matrix}
		1  \\
		C{{D}^{2}}-A{{D}^{2}}  \\
		A{{D}^{2}}-B{{D}^{2}}  \\
	\end{matrix} \right).\]
	
	
	According to the proof process of theorem \ref{thm:Thm28.2.1}, the above formula is exactly the frame components of the IC-F, i.e., the frame components of the projection point ${{D}_{{{D}_{H}}}}$ of the vertex $D$ of the tetrahedron on the $\triangle ABC$ plane. 
	
	The other formulas can be proved similarly.
\end{proof}
\hfill $\square$\par

The above theorem gives the relationship between the frame components of the projection point of the IC-T and the frame components of the projection point of the vertex and the frame component of the IC-T.

If the same tetrahedral frame is used, the following theorem is obtained.


\begin{theorem}{Relationship 2 of frame components between projection and IC-T, Daiyuan Zhang}{Thm28.4.2}\label{Thm28.4.2} 
	Given a tetrahedron $ABCD$, point $O$ is any point in space. Suppose that the projection point of point $A$ to the $\triangle BCD$ plane of the tetrahedron is ${{A}_{{{A}_{H}}}}$, the projection point of point $B$ on the $\triangle CDA$ plane of tetrahedron is ${{B}_{{{B}_{H}}}}$, the projection point of point $C$ to the $\triangle DAB$ plane of tetrahedron is ${{C}_{{{C}_{H}}}}$, the projection point of point $D$ on the $\triangle ABC$ plane of tetrahedron is ${{D}_{{{D}_{H}}}}$. The projection points of IC-T $P$ to the planes $\triangle BCD$, $\triangle CDA$, $\triangle DAB$, $\triangle ABC$ of the tetrahedron are ${{P}_{{{A}_{H}}}}$, ${{P}_{{{B}_{H}}}}$, ${{P}_{{{C}_{H}}}}$, ${{P}_{{{D}_{H}}}}$ respectively, then:	
		
	
	1. If ${{P}_{{{D}_{H}}}}\in {{\pi }_{ABC}}$, ${{D}_{{{D}_{H}}}}\in {{\pi }_{ABC}}$, then the frame components of ${{P}_{{{D}_{H}}}}$ in the frame $\left( O;A,B,C \right)$ is
	
	\[\left( \begin{matrix}
		\beta _{A}^{{{P}_{{{D}_{H}}}}}  \\
		\beta _{B}^{{{P}_{{{D}_{H}}}}}  \\
		\beta _{C}^{{{P}_{{{D}_{H}}}}}  \\
		\beta _{D}^{{{P}_{{{D}_{H}}}}}  \\
	\end{matrix} \right)=\left( \begin{matrix}
		1 & 0 & 0 & \beta _{A}^{{{D}_{{{D}_{H}}}}}  \\
		0 & 1 & 0 & \beta _{B}^{{{D}_{{{D}_{H}}}}}  \\
		0 & 0 & 1 & \beta _{C}^{{{D}_{{{D}_{H}}}}}  \\
		0 & 0 & 0 & \beta _{D}^{{{D}_{{{D}_{H}}}}}  \\
	\end{matrix} \right)\left( \begin{matrix}
		\beta _{A}^{P}  \\
		\beta _{B}^{P}  \\
		\beta _{C}^{P}  \\
		\beta _{D}^{P}  \\
	\end{matrix} \right).\]
	Where 
	\[\beta _{A}^{{{P}_{{{D}_{H}}}}}=\alpha _{A}^{{{P}_{{{D}_{H}}}}},\ \beta _{B}^{{{P}_{{{D}_{H}}}}}=\alpha _{B}^{{{P}_{{{D}_{H}}}}},\ \beta _{C}^{{{P}_{{{D}_{H}}}}}=\alpha _{C}^{{{P}_{{{D}_{H}}}}},\ \beta _{D}^{{{P}_{{{D}_{H}}}}}=0,\]
	\[\beta _{A}^{{{D}_{{{D}_{H}}}}}=\alpha _{A}^{{{D}_{{{D}_{H}}}}},\ \beta _{B}^{{{D}_{{{D}_{H}}}}}=\alpha _{B}^{{{D}_{{{D}_{H}}}}},\ \beta _{C}^{{{D}_{{{D}_{H}}}}}=\alpha _{C}^{{{D}_{{{D}_{H}}}}},\ \beta _{D}^{{{D}_{{{D}_{H}}}}}=0.\ \]
	
	
	2. If ${{P}_{{{A}_{H}}}}\in {{\pi }_{BCD}}$, ${{A}_{{{A}_{H}}}}\in {{\pi }_{BCD}}$, then the frame components of ${{P}_{{{A}_{H}}}}$ in the frame $\left( O;B,C,D \right)$ is
	
	\[\left( \begin{matrix}
		\beta _{B}^{{{P}_{{{A}_{H}}}}}  \\
		\beta _{C}^{{{P}_{{{A}_{H}}}}}  \\
		\beta _{D}^{{{P}_{{{A}_{H}}}}}  \\
		\beta _{A}^{{{P}_{{{A}_{H}}}}}  \\
	\end{matrix} \right)=\left( \begin{matrix}
		1 & 0 & 0 & \beta _{B}^{{{A}_{{{A}_{H}}}}}  \\
		0 & 1 & 0 & \beta _{C}^{{{A}_{{{A}_{H}}}}}  \\
		0 & 0 & 1 & \beta _{D}^{{{A}_{{{A}_{H}}}}}  \\
		0 & 0 & 0 & \beta _{A}^{{{A}_{{{A}_{H}}}}}  \\
	\end{matrix} \right)\left( \begin{matrix}
		\beta _{B}^{P}  \\
		\beta _{C}^{P}  \\
		\beta _{D}^{P}  \\
		\beta _{A}^{P}  \\
	\end{matrix} \right).\]
	Where 
	\[\beta _{B}^{{{P}_{{{A}_{H}}}}}=\alpha _{B}^{{{P}_{{{A}_{H}}}}},\ \beta _{C}^{{{P}_{{{A}_{H}}}}}=\alpha _{C}^{{{P}_{{{A}_{H}}}}},\ \beta _{D}^{{{P}_{{{A}_{H}}}}}=\alpha _{D}^{{{P}_{{{A}_{H}}}}},\ \beta _{A}^{{{P}_{{{A}_{H}}}}}=0,\]
	\[\beta _{B}^{{{A}_{{{A}_{H}}}}}=\alpha _{B}^{{{A}_{{{A}_{H}}}}},\ \beta _{C}^{{{A}_{{{A}_{H}}}}}=\alpha _{C}^{{{A}_{{{A}_{H}}}}},\ \beta _{D}^{{{A}_{{{A}_{H}}}}}=\alpha _{D}^{{{A}_{{{A}_{H}}}}},\ \beta _{A}^{{{A}_{{{A}_{H}}}}}=0.\ \]
	
	
	3. If ${{P}_{{{B}_{H}}}}\in {{\pi }_{CDA}}$, ${{B}_{{{B}_{H}}}}\in {{\pi }_{CDA}}$, then the frame components of ${{P}_{{{B}_{H}}}}$ in the frame $\left( O;C,D,A \right)$ is
	
	\[\left( \begin{matrix}
		\beta _{C}^{{{P}_{{{B}_{H}}}}}  \\
		\beta _{D}^{{{P}_{{{B}_{H}}}}}  \\
		\beta _{A}^{{{P}_{{{B}_{H}}}}}  \\
		\beta _{B}^{{{P}_{{{B}_{H}}}}}  \\
	\end{matrix} \right)=\left( \begin{matrix}
		1 & 0 & 0 & \beta _{C}^{{{B}_{{{B}_{H}}}}}  \\
		0 & 1 & 0 & \beta _{D}^{{{B}_{{{B}_{H}}}}}  \\
		0 & 0 & 1 & \beta _{A}^{{{B}_{{{B}_{H}}}}}  \\
		0 & 0 & 0 & \beta _{B}^{{{B}_{{{B}_{H}}}}}  \\
	\end{matrix} \right)\left( \begin{matrix}
		\beta _{C}^{P}  \\
		\beta _{D}^{P}  \\
		\beta _{A}^{P}  \\
		\beta _{B}^{P}  \\
	\end{matrix} \right).\]
	Where 
	\[\beta _{C}^{{{P}_{{{B}_{H}}}}}=\alpha _{C}^{{{P}_{{{B}_{H}}}}},\ \beta _{D}^{{{P}_{{{B}_{H}}}}}=\alpha _{D}^{{{P}_{{{B}_{H}}}}},\ \beta _{A}^{{{P}_{{{B}_{H}}}}}=\alpha _{A}^{{{P}_{{{B}_{H}}}}},\ \beta _{B}^{{{P}_{{{B}_{H}}}}}=0,\]
	\[\beta _{C}^{{{B}_{{{B}_{H}}}}}=\alpha _{C}^{{{B}_{{{B}_{H}}}}},\ \beta _{D}^{{{B}_{{{B}_{H}}}}}=\alpha _{D}^{{{B}_{{{B}_{H}}}}},\ \beta _{A}^{{{B}_{{{B}_{H}}}}}=\alpha _{A}^{{{B}_{{{B}_{H}}}}},\ \beta _{B}^{{{B}_{{{B}_{H}}}}}=0.\ \]
	
	
	4. If ${{P}_{{{C}_{H}}}}\in {{\pi }_{DAB}}$,  ${{C}_{{{C}_{H}}}}\in {{\pi }_{DAB}}$, then the frame components of ${{P}_{{{C}_{H}}}}$ in the frame $\left( O;D,A,B \right)$ is
	
	\[\left( \begin{matrix}
		\beta _{D}^{{{P}_{{{C}_{H}}}}}  \\
		\beta _{A}^{{{P}_{{{C}_{H}}}}}  \\
		\beta _{B}^{{{P}_{{{C}_{H}}}}}  \\
		\beta _{C}^{{{P}_{{{C}_{H}}}}}  \\
	\end{matrix} \right)=\left( \begin{matrix}
		1 & 0 & 0 & \beta _{D}^{{{C}_{{{C}_{H}}}}}  \\
		0 & 1 & 0 & \beta _{A}^{{{C}_{{{C}_{H}}}}}  \\
		0 & 0 & 1 & \beta _{B}^{{{C}_{{{C}_{H}}}}}  \\
		0 & 0 & 0 & \beta _{C}^{{{C}_{{{C}_{H}}}}}  \\
	\end{matrix} \right)\left( \begin{matrix}
		\beta _{D}^{P}  \\
		\beta _{A}^{P}  \\
		\beta _{B}^{P}  \\
		\beta _{C}^{P}  \\
	\end{matrix} \right).\]
	Where 
	\[\beta _{D}^{{{P}_{{{C}_{H}}}}}=\alpha _{D}^{{{P}_{{{C}_{H}}}}},\ \beta _{A}^{{{P}_{{{C}_{H}}}}}=\alpha _{A}^{{{P}_{{{C}_{H}}}}},\ \beta _{B}^{{{P}_{{{C}_{H}}}}}=\alpha _{B}^{{{P}_{{{C}_{H}}}}},\ \beta _{C}^{{{P}_{{{C}_{H}}}}}=0,\]
	\[\beta _{D}^{{{C}_{{{C}_{H}}}}}=\alpha _{D}^{{{C}_{{{C}_{H}}}}},\ \beta _{A}^{{{C}_{{{C}_{H}}}}}=\alpha _{A}^{{{C}_{{{C}_{H}}}}},\ \beta _{B}^{{{C}_{{{C}_{H}}}}}=\alpha _{B}^{{{C}_{{{C}_{H}}}}},\ \beta _{C}^{{{C}_{{{C}_{H}}}}}=0.\ \]
\end{theorem}

\begin{proof}
	According to theorem \ref{thm:Thm28.4.1} and theorem \ref{thm:Thm20.4.1}, the conclusion of the theorem can be obtained.		
\end{proof}
\hfill $\square$\par

\section{Frame components of projection points of special IC-T}\label{Sec28.5}

\subsection{Frame components of the projection point of the circumcenter of tetrahedron}\label{Subsec28.5.1}

The calculation foe the frame components of the projection point of the circumcenter of tetrahedron is simple.


\begin{theorem}{Formula 1 of frame components for projection of circumcenter, Daiyuan Zhang}{Thm28.5.1}\label{Thm28.5.1} 
	Given the tetrahedron $ABCD$, point $O$ is any point in space, and point $Q$ is the circumcenter of the circumscribed sphere of tetrahedron, then:
	
	
	1. Assume that the projection point of the circumcenter $Q$ to the plane of $\triangle BCD$ is ${{Q}_{{{A}_{H}}}}$, $\triangle BCD$ has an area of ${{S}^{A}}$, then the ocircumcenter $Q$ has the following frame components in the frame of $\left( O;B,C,D \right)$:
	
	\[\alpha _{B}^{{{Q}_{{{A}_{H}}}}}=\frac{1}{8{{\left( {{S}^{A}} \right)}^{2}}}\left( \Delta _{2}^{A}-C{{D}^{2}} \right)C{{D}^{2}},\]
	\[\alpha _{C}^{{{Q}_{{{A}_{H}}}}}=\frac{1}{8{{\left( {{S}^{A}} \right)}^{2}}}\left( \Delta _{2}^{A}-D{{B}^{2}} \right)D{{B}^{2}},\]
	\[\alpha _{D}^{{{Q}_{{{A}_{H}}}}}=\frac{1}{8{{\left( {{S}^{A}} \right)}^{2}}}\left( \Delta _{2}^{A}-B{{C}^{2}} \right)B{{C}^{2}}.\]
	Where 
	\[\Delta _{2}^{A}=\frac{1}{2}\left( B{{C}^{2}}+C{{D}^{2}}+D{{B}^{2}} \right).\]
	
	
	2. Assume that the projection point of the circumcenter $Q$ to the plane of $\triangle CDA$ is ${{Q}_{{{B}_{H}}}}$, $\triangle CDA$ has an area of ${{S}^{B}}$, then the ocircumcenter $Q$ has the following frame components in the frame of $\left( O;C,D,A \right)$:
	
	\[\alpha _{C}^{{{Q}_{{{B}_{H}}}}}=\frac{1}{8{{\left( {{S}^{B}} \right)}^{2}}}\left( \Delta _{2}^{B}-D{{A}^{2}} \right)D{{A}^{2}},\]
	\[\alpha _{D}^{{{Q}_{{{B}_{H}}}}}=\frac{1}{8{{\left( {{S}^{B}} \right)}^{2}}}\left( \Delta _{2}^{B}-A{{C}^{2}} \right)A{{C}^{2}},\]
	\[\alpha _{A}^{{{Q}_{{{B}_{H}}}}}=\frac{1}{8{{\left( {{S}^{B}} \right)}^{2}}}\left( \Delta _{2}^{B}-C{{D}^{2}} \right)C{{D}^{2}}.\]
	Where 
	\[\Delta _{2}^{B}=\frac{1}{2}\left( C{{D}^{2}}+D{{A}^{2}}+A{{C}^{2}} \right).\]
	
	
	3. Assume that the projection point of the circumcenter $Q$ to the plane of $\triangle DAB$ is ${{Q}_{{{C}_{H}}}}$, $\triangle DAB$ has an area of ${{S}^{C}}$, then the ocircumcenter $Q$ has the following frame components in the frame of $\left( O;D,A,B \right)$:
	
	\[\alpha _{D}^{{{Q}_{{{C}_{H}}}}}=\frac{1}{8{{\left( {{S}^{C}} \right)}^{2}}}\left( \Delta _{2}^{C}-A{{B}^{2}} \right)A{{B}^{2}},\]
	\[\alpha _{A}^{{{Q}_{{{C}_{H}}}}}=\frac{1}{8{{\left( {{S}^{C}} \right)}^{2}}}\left( \Delta _{2}^{C}-B{{D}^{2}} \right)B{{D}^{2}},\]
	\[\alpha _{B}^{{{Q}_{{{C}_{H}}}}}=\frac{1}{8{{\left( {{S}^{C}} \right)}^{2}}}\left( \Delta _{2}^{C}-D{{A}^{2}} \right)D{{A}^{2}}.\]
	Where 
	\[\Delta _{2}^{C}=\frac{1}{2}\left( D{{A}^{2}}+A{{B}^{2}}+B{{D}^{2}} \right).\]
	
	
	4. Assume that the projection point of the circumcenter $Q$ to the plane of $\triangle ABC$ is ${{Q}_{{{D}_{H}}}}$, $\triangle ABC$ has an area of ${{S}^{D}}$, then the ocircumcenter $Q$ has the following frame components in the frame of $\left( O;A,B,C \right)$:
	
	\[\alpha _{A}^{{{Q}_{{{D}_{H}}}}}=\frac{1}{8{{\left( {{S}^{D}} \right)}^{2}}}\left( \Delta _{2}^{D}-B{{C}^{2}} \right)B{{C}^{2}},\]
	\[\alpha _{B}^{{{Q}_{{{D}_{H}}}}}=\frac{1}{8{{\left( {{S}^{D}} \right)}^{2}}}\left( \Delta _{2}^{D}-C{{A}^{2}} \right)C{{A}^{2}},\]
	\[\alpha _{C}^{{{Q}_{{{D}_{H}}}}}=\frac{1}{8{{\left( {{S}^{D}} \right)}^{2}}}\left( \Delta _{2}^{D}-A{{B}^{2}} \right)A{{B}^{2}}.\]
	Where 
	\[\Delta _{2}^{D}=\frac{1}{2}\left( A{{B}^{2}}+B{{C}^{2}}+C{{A}^{2}} \right).\]
\end{theorem}

\begin{proof}
	Because $AQ=BQ=CQ=DQ$, the result can be obtained by replacing the point $Q$ with the circumcenter $Q$ in theorem \ref{thm:Thm28.2.1}
\end{proof}
\hfill $\square$\par

\begin{theorem}{Formula 2 of frame components for projection of circumcenter, Daiyuan Zhang}{Thm28.5.2}\label{Thm28.5.2} 
	Given the tetrahedron $ABCD$, point $O$ is any point in space, and point $Q$ is the circumcenter of the circumscribed sphere of tetrahedron, then:
		
	
	1. Assume that the projection point of the circumcenter $Q$ to the plane of $\triangle BCD$ is ${{Q}_{{{A}_{H}}}}$, $\triangle BCD$ has an area of ${{S}^{A}}$, then the ocircumcenter $Q$ has the following frame components in the frame of $\left( O;B,C,D \right)$:
	
	\[\beta _{B}^{{{Q}_{{{A}_{H}}}}}=\frac{1}{8{{\left( {{S}^{A}} \right)}^{2}}}\left( \Delta _{2}^{A}-C{{D}^{2}} \right)C{{D}^{2}},\]
	\[\beta _{C}^{{{Q}_{{{A}_{H}}}}}=\frac{1}{8{{\left( {{S}^{A}} \right)}^{2}}}\left( \Delta _{2}^{A}-D{{B}^{2}} \right)D{{B}^{2}},\]
	\[\beta _{D}^{{{Q}_{{{A}_{H}}}}}=\frac{1}{8{{\left( {{S}^{A}} \right)}^{2}}}\left( \Delta _{2}^{A}-B{{C}^{2}} \right)B{{C}^{2}},\]
	\[\beta _{A}^{{{Q}_{{{A}_{H}}}}}=0.\]
	Where 
	\[\Delta _{2}^{A}=\frac{1}{2}\left( B{{C}^{2}}+C{{D}^{2}}+D{{B}^{2}} \right).\]
	
	
	2. Assume that the projection point of the circumcenter $Q$ to the plane of $\triangle CDA$ is ${{Q}_{{{B}_{H}}}}$, $\triangle CDA$ has an area of ${{S}^{B}}$, then the ocircumcenter $Q$ has the following frame components in the frame of $\left( O;C,D,A \right)$:
	
	\[\beta _{C}^{{{Q}_{{{B}_{H}}}}}=\frac{1}{8{{\left( {{S}^{B}} \right)}^{2}}}\left( \Delta _{2}^{B}-D{{A}^{2}} \right)D{{A}^{2}},\]
	\[\beta _{D}^{{{Q}_{{{B}_{H}}}}}=\frac{1}{8{{\left( {{S}^{B}} \right)}^{2}}}\left( \Delta _{2}^{B}-A{{C}^{2}} \right)A{{C}^{2}},\]
	\[\beta _{A}^{{{Q}_{{{B}_{H}}}}}=\frac{1}{8{{\left( {{S}^{B}} \right)}^{2}}}\left( \Delta _{2}^{B}-C{{D}^{2}} \right)C{{D}^{2}},\]
	\[\beta _{B}^{{{Q}_{{{B}_{H}}}}}=0.\]
	Where 
	\[\Delta _{2}^{B}=\frac{1}{2}\left( C{{D}^{2}}+D{{A}^{2}}+A{{C}^{2}} \right).\]
	
	
	3. Assume that the projection point of the circumcenter $Q$ to the plane of $\triangle DAB$ is ${{Q}_{{{C}_{H}}}}$, $\triangle DAB$ has an area of ${{S}^{C}}$, then the ocircumcenter $Q$ has the following frame components in the frame of $\left( O;D,A,B \right)$:
	
	\[\beta _{D}^{{{Q}_{{{C}_{H}}}}}=\frac{1}{8{{\left( {{S}^{C}} \right)}^{2}}}\left( \Delta _{2}^{C}-A{{B}^{2}} \right)A{{B}^{2}},\]
	\[\beta _{A}^{{{Q}_{{{C}_{H}}}}}=\frac{1}{8{{\left( {{S}^{C}} \right)}^{2}}}\left( \Delta _{2}^{C}-B{{D}^{2}} \right)B{{D}^{2}},\]
	\[\beta _{B}^{{{Q}_{{{C}_{H}}}}}=\frac{1}{8{{\left( {{S}^{C}} \right)}^{2}}}\left( \Delta _{2}^{C}-D{{A}^{2}} \right)D{{A}^{2}},\]
	\[\beta _{C}^{{{Q}_{{{C}_{H}}}}}=0.\]
	Where 
	\[\Delta _{2}^{C}=\frac{1}{2}\left( D{{A}^{2}}+A{{B}^{2}}+B{{D}^{2}} \right).\]
	
	
	4. Assume that the projection point of the circumcenter $Q$ to the plane of $\triangle ABC$ is ${{Q}_{{{D}_{H}}}}$, $\triangle ABC$ has an area of ${{S}^{D}}$, then the ocircumcenter $Q$ has the following frame components in the frame of $\left( O;A,B,C \right)$:
	
	\[\beta _{A}^{{{Q}_{{{D}_{H}}}}}=\frac{1}{8{{\left( {{S}^{D}} \right)}^{2}}}\left( \Delta _{2}^{D}-B{{C}^{2}} \right)B{{C}^{2}},\]
	\[\beta _{B}^{{{Q}_{{{D}_{H}}}}}=\frac{1}{8{{\left( {{S}^{D}} \right)}^{2}}}\left( \Delta _{2}^{D}-C{{A}^{2}} \right)C{{A}^{2}},\]
	\[\beta _{C}^{{{Q}_{{{D}_{H}}}}}=\frac{1}{8{{\left( {{S}^{D}} \right)}^{2}}}\left( \Delta _{2}^{D}-A{{B}^{2}} \right)A{{B}^{2}},\]
	\[\beta _{D}^{{{Q}_{{{D}_{H}}}}}=0.\]
	Wherer 
	\[\Delta _{2}^{D}=\frac{1}{2}\left( A{{B}^{2}}+B{{C}^{2}}+C{{A}^{2}} \right).\]
\end{theorem}

\begin{proof}
	This theorem can be proved by using theorem \ref{thm:Thm28.5.1} and theorem \ref{thm:Thm20.4.1}
\end{proof}
\hfill $\square$\par

\subsection{Frame components of the projection points of the centroid of tetrahedron}\label{Subsec28.5.2}


Applying the theorems given above (theorem \ref{thm:Thm28.4.1} and theorem \ref{thm:Thm28.4.2}) to the centroid of the tetrahedron, the frame components of the projection points of the centroid of the tetrahedron can be obtained.
%

\begin{theorem}{Formula 1 of frame components for projection of centroid, Daiyuan Zhang}{Thm28.5.3}\label{Thm28.5.3} 
	Given a tetrahedron $ABCD$, point $O$ is any point in space. Suppose that the projection point of point $A$ to the $\triangle BCD$ plane of the tetrahedron is ${{A}_{{{A}_{H}}}}$, the projection point of point $B$ on the $\triangle CDA$ plane of tetrahedron is ${{B}_{{{B}_{H}}}}$, the projection point of point $C$ to the $\triangle DAB$ plane of tetrahedron is ${{C}_{{{C}_{H}}}}$, the projection point of point $D$ on the $\triangle ABC$ plane of tetrahedron is ${{D}_{{{D}_{H}}}}$. The projection points of centroid $G$ to the planes $\triangle BCD$, $\triangle CDA$, $\triangle DAB$, $\triangle ABC$ of the tetrahedron are ${{G}_{{{A}_{H}}}}$, ${{G}_{{{B}_{H}}}}$, ${{G}_{{{C}_{H}}}}$, ${{G}_{{{D}_{H}}}}$ respectively, then:	
	
	
	1. 	The frame components of ${{G}_{{{A}_{H}}}}$ on the triangular frame $\left( O;B,C,D \right)$ is:
	
	\[\alpha _{B}^{{{G}_{{{A}_{H}}}}}=\frac{1}{32{{\left( {{S}^{A}} \right)}^{2}}}\left( \begin{aligned}
		& C{{D}^{2}}\left( \Delta _{2}^{A}-C{{D}^{2}} \right)+\left( \Delta _{2}^{A}-D{{B}^{2}} \right)\left( D{{A}^{2}}-B{{A}^{2}} \right) \\ 
		& +\left( \Delta _{2}^{A}-B{{C}^{2}} \right)\left( C{{A}^{2}}-B{{A}^{2}} \right)+8{{\left( {{S}^{A}} \right)}^{2}} \\ 
	\end{aligned} \right),\]
	\[\alpha _{C}^{{{G}_{{{A}_{H}}}}}=\frac{1}{32{{\left( {{S}^{A}} \right)}^{2}}}\left( \begin{aligned}
		& D{{B}^{2}}\left( \Delta _{2}^{A}-D{{B}^{2}} \right)+\left( \Delta _{2}^{A}-B{{C}^{2}} \right)\left( B{{A}^{2}}-C{{A}^{2}} \right) \\ 
		& +\left( \Delta _{2}^{A}-C{{D}^{2}} \right)\left( D{{A}^{2}}-C{{A}^{2}} \right)+8{{\left( {{S}^{A}} \right)}^{2}} \\ 
	\end{aligned} \right),\]
	\[\alpha _{D}^{{{G}_{{{A}_{H}}}}}=\frac{1}{32{{\left( {{S}^{A}} \right)}^{2}}}\left( \begin{aligned}
		& B{{C}^{2}}\left( \Delta _{2}^{A}-B{{C}^{2}} \right)+\left( \Delta _{2}^{A}-C{{D}^{2}} \right)\left( C{{A}^{2}}-D{{A}^{2}} \right) \\ 
		& +\left( \Delta _{2}^{A}-D{{B}^{2}} \right)\left( B{{A}^{2}}-D{{A}^{2}} \right)+8{{\left( {{S}^{A}} \right)}^{2}} \\ 
	\end{aligned} \right).\]
	Where 
	\[\Delta _{2}^{A}=\frac{1}{2}\left( B{{C}^{2}}+C{{D}^{2}}+D{{B}^{2}} \right).\]
	
	
	2. 	The frame components of ${{G}_{{{B}_{H}}}}$ on the triangular frame $\left( O;C,D,A \right)$ is:
		
	\[\alpha _{C}^{{{G}_{{{B}_{H}}}}}=\frac{1}{32{{\left( {{S}^{B}} \right)}^{2}}}\left( \begin{aligned}
		& D{{A}^{2}}\left( \Delta _{2}^{B}-D{{A}^{2}} \right)+\left( \Delta _{2}^{B}-A{{C}^{2}} \right)\left( A{{B}^{2}}-C{{B}^{2}} \right) \\ 
		& +\left( \Delta _{2}^{B}-C{{D}^{2}} \right)\left( D{{B}^{2}}-C{{B}^{2}} \right)+8{{\left( {{S}^{B}} \right)}^{2}} \\ 
	\end{aligned} \right),\]
	\[\alpha _{D}^{{{G}_{{{B}_{H}}}}}=\frac{1}{32{{\left( {{S}^{B}} \right)}^{2}}}\left( \begin{aligned}
		& A{{C}^{2}}\left( \Delta _{2}^{B}-A{{C}^{2}} \right)+\left( \Delta _{2}^{B}-C{{D}^{2}} \right)\left( C{{B}^{2}}-D{{B}^{2}} \right) \\ 
		& +\left( \Delta _{2}^{B}-D{{A}^{2}} \right)\left( A{{B}^{2}}-D{{B}^{2}} \right)+8{{\left( {{S}^{B}} \right)}^{2}} \\ 
	\end{aligned} \right),\]
	\[\alpha _{A}^{{{G}_{{{B}_{H}}}}}=\frac{1}{32{{\left( {{S}^{B}} \right)}^{2}}}\left( \begin{aligned}
		& C{{D}^{2}}\left( \Delta _{2}^{B}-C{{D}^{2}} \right)+\left( \Delta _{2}^{B}-D{{A}^{2}} \right)\left( D{{B}^{2}}-A{{B}^{2}} \right) \\ 
		& +\left( \Delta _{2}^{B}-A{{C}^{2}} \right)\left( C{{B}^{2}}-A{{B}^{2}} \right)+8{{\left( {{S}^{B}} \right)}^{2}} \\ 
	\end{aligned} \right).\]
	Where 
	\[\Delta _{2}^{B}=\frac{1}{2}\left( C{{D}^{2}}+D{{A}^{2}}+A{{C}^{2}} \right).\]
	
	
	3. 	The frame components of ${{G}_{{{C}_{H}}}}$ on the triangular frame $\left( O;D,A,B \right)$ is:
		
	\[\alpha _{D}^{{{G}_{{{C}_{H}}}}}=\frac{1}{32{{\left( {{S}^{C}} \right)}^{2}}}\left( \begin{aligned}
		& A{{B}^{2}}\left( \Delta _{2}^{C}-A{{B}^{2}} \right)+\left( \Delta _{2}^{C}-B{{D}^{2}} \right)\left( B{{C}^{2}}-D{{C}^{2}} \right) \\ 
		& +\left( \Delta _{2}^{C}-D{{A}^{2}} \right)\left( A{{C}^{2}}-D{{C}^{2}} \right)+8{{\left( {{S}^{C}} \right)}^{2}} \\ 
	\end{aligned} \right),\]
	\[\alpha _{A}^{{{G}_{{{C}_{H}}}}}=\frac{1}{32{{\left( {{S}^{C}} \right)}^{2}}}\left( \begin{aligned}
		& B{{D}^{2}}\left( \Delta _{2}^{C}-B{{D}^{2}} \right)+\left( \Delta _{2}^{C}-D{{A}^{2}} \right)\left( D{{C}^{2}}-A{{C}^{2}} \right) \\ 
		& +\left( \Delta _{2}^{C}-A{{B}^{2}} \right)\left( B{{C}^{2}}-A{{C}^{2}} \right)+8{{\left( {{S}^{C}} \right)}^{2}} \\ 
	\end{aligned} \right),\]
	\[\alpha _{B}^{{{G}_{{{C}_{H}}}}}=\frac{1}{32{{\left( {{S}^{C}} \right)}^{2}}}\left( \begin{aligned}
		& D{{A}^{2}}\left( \Delta _{2}^{C}-D{{A}^{2}} \right)+\left( \Delta _{2}^{C}-A{{B}^{2}} \right)\left( A{{C}^{2}}-B{{C}^{2}} \right) \\ 
		& +\left( \Delta _{2}^{C}-B{{D}^{2}} \right)\left( D{{C}^{2}}-B{{C}^{2}} \right)+8{{\left( {{S}^{C}} \right)}^{2}} \\ 
	\end{aligned} \right).\]
	Where 
	\[\Delta _{2}^{C}=\frac{1}{2}\left( D{{A}^{2}}+A{{B}^{2}}+B{{D}^{2}} \right).\]
	
	
	4. 	The frame components of ${{G}_{{{D}_{H}}}}$ on the triangular frame $\left( O;A,B,C \right)$ is:
	
	\[\alpha _{A}^{{{G}_{{{D}_{H}}}}}=\frac{1}{32{{\left( {{S}^{D}} \right)}^{2}}}\left( \begin{aligned}
		& B{{C}^{2}}\left( \Delta _{2}^{D}-B{{C}^{2}} \right)+\left( \Delta _{2}^{D}-C{{A}^{2}} \right)\left( C{{D}^{2}}-A{{D}^{2}} \right) \\ 
		& +\left( \Delta _{2}^{D}-A{{B}^{2}} \right)\left( B{{D}^{2}}-A{{D}^{2}} \right)+8{{\left( {{S}^{D}} \right)}^{2}} \\ 
	\end{aligned} \right),\]
	\[\alpha _{B}^{{{G}_{{{D}_{H}}}}}=\frac{1}{32{{\left( {{S}^{D}} \right)}^{2}}}\left( \begin{aligned}
		& C{{A}^{2}}\left( \Delta _{2}^{D}-C{{A}^{2}} \right)+\left( \Delta _{2}^{D}-A{{B}^{2}} \right)\left( A{{D}^{2}}-B{{D}^{2}} \right) \\ 
		& +\left( \Delta _{2}^{D}-B{{C}^{2}} \right)\left( C{{D}^{2}}-B{{D}^{2}} \right)+8{{\left( {{S}^{D}} \right)}^{2}} \\ 
	\end{aligned} \right),\]
	\[\alpha _{C}^{{{G}_{{{D}_{H}}}}}=\frac{1}{32{{\left( {{S}^{D}} \right)}^{2}}}\left( \begin{aligned}
		& A{{B}^{2}}\left( \Delta _{2}^{D}-A{{B}^{2}} \right)+\left( \Delta _{2}^{D}-B{{C}^{2}} \right)\left( B{{D}^{2}}-C{{D}^{2}} \right) \\ 
		& +\left( \Delta _{2}^{D}-C{{A}^{2}} \right)\left( A{{D}^{2}}-C{{D}^{2}} \right)+8{{\left( {{S}^{D}} \right)}^{2}} \\ 
	\end{aligned} \right).\]
	Where 
	\[\Delta _{2}^{D}=\frac{1}{2}\left( A{{B}^{2}}+B{{C}^{2}}+C{{A}^{2}} \right).\]
\end{theorem}

\begin{proof}
	
		
	Here's the proof for the frame components of the projection point ${{G}_{{{D}_{H}}}}$ on the $\triangle ABC$ plane.
				
	According to theorem \ref{thm:Thm28.4.1} and theorem \ref{thm:Thm21.1.1}, we get:
		
	\[\left( \begin{matrix}
		\alpha _{A}^{{{G}_{{{D}_{H}}}}}  \\
		\alpha _{B}^{{{G}_{{{D}_{H}}}}}  \\
		\alpha _{C}^{{{G}_{{{D}_{H}}}}}  \\
	\end{matrix} \right)=\left( \begin{matrix}
		1 & 0 & 0 & \alpha _{A}^{{{D}_{{{D}_{H}}}}}  \\
		0 & 1 & 0 & \alpha _{B}^{{{D}_{{{D}_{H}}}}}  \\
		0 & 0 & 1 & \alpha _{C}^{{{D}_{{{D}_{H}}}}}  \\
	\end{matrix} \right)\left( \begin{matrix}
		\beta _{A}^{G}  \\
		\beta _{B}^{G}  \\
		\beta _{C}^{G}  \\
		\beta _{D}^{G}  \\
	\end{matrix} \right)=\frac{1}{4}\left( \begin{matrix}
		1 & 0 & 0 & \alpha _{A}^{{{D}_{{{D}_{H}}}}}  \\
		0 & 1 & 0 & \alpha _{B}^{{{D}_{{{D}_{H}}}}}  \\
		0 & 0 & 1 & \alpha _{C}^{{{D}_{{{D}_{H}}}}}  \\
	\end{matrix} \right)\left( \begin{matrix}
		1  \\
		1  \\
		1  \\
		1  \\
	\end{matrix} \right).\]
	
	
	Replace theorem \ref{thm:Thm28.3.1} and theorem \ref{thm:Thm21.1.1} directly into the above formula to obtain:
	
	\[\alpha _{A}^{{{G}_{{{D}_{H}}}}}=\frac{1}{32{{\left( {{S}^{D}} \right)}^{2}}}\left( \begin{aligned}
		& B{{C}^{2}}\left( \Delta _{2}^{D}-B{{C}^{2}} \right)+\left( \Delta _{2}^{D}-C{{A}^{2}} \right)\left( C{{D}^{2}}-A{{D}^{2}} \right) \\ 
		& +\left( \Delta _{2}^{D}-A{{B}^{2}} \right)\left( B{{D}^{2}}-A{{D}^{2}} \right)+8{{\left( {{S}^{D}} \right)}^{2}} \\ 
	\end{aligned} \right),\]
	\[\alpha _{B}^{{{G}_{{{D}_{H}}}}}=\frac{1}{32{{\left( {{S}^{D}} \right)}^{2}}}\left( \begin{aligned}
		& C{{A}^{2}}\left( \Delta _{2}^{D}-C{{A}^{2}} \right)+\left( \Delta _{2}^{D}-A{{B}^{2}} \right)\left( A{{D}^{2}}-B{{D}^{2}} \right) \\ 
		& +\left( \Delta _{2}^{D}-B{{C}^{2}} \right)\left( C{{D}^{2}}-B{{D}^{2}} \right)+8{{\left( {{S}^{D}} \right)}^{2}} \\ 
	\end{aligned} \right),\]
	\[\alpha _{C}^{{{G}_{{{D}_{H}}}}}=\frac{1}{32{{\left( {{S}^{D}} \right)}^{2}}}\left( \begin{aligned}
		& A{{B}^{2}}\left( \Delta _{2}^{D}-A{{B}^{2}} \right)+\left( \Delta _{2}^{D}-B{{C}^{2}} \right)\left( B{{D}^{2}}-C{{D}^{2}} \right) \\ 
		& +\left( \Delta _{2}^{D}-C{{A}^{2}} \right)\left( A{{D}^{2}}-C{{D}^{2}} \right)+8{{\left( {{S}^{D}} \right)}^{2}} \\ 
	\end{aligned} \right).\]
	
	Other formulas can be proved similarly.
\end{proof}
\hfill $\square$\par

\begin{theorem}{Formula 2 of frame components for projection of centroid,, Daiyuan Zhang}{Thm28.5.4}\label{Thm28.5.4} 
	Given a tetrahedron $ABCD$, point $O$ is any point in space. Suppose that the projection point of point $A$ to the $\triangle BCD$ plane of the tetrahedron is ${{A}_{{{A}_{H}}}}$, the projection point of point $B$ on the $\triangle CDA$ plane of tetrahedron is ${{B}_{{{B}_{H}}}}$, the projection point of point $C$ to the $\triangle DAB$ plane of tetrahedron is ${{C}_{{{C}_{H}}}}$, the projection point of point $D$ on the $\triangle ABC$ plane of tetrahedron is ${{D}_{{{D}_{H}}}}$. The projection points of centroid $G$ to the planes $\triangle BCD$, $\triangle CDA$, $\triangle DAB$, $\triangle ABC$ of the tetrahedron are ${{G}_{{{A}_{H}}}}$, ${{G}_{{{B}_{H}}}}$, ${{G}_{{{C}_{H}}}}$, ${{G}_{{{D}_{H}}}}$ respectively, then:	
	
	
	1. The frame components of ${{G}_{{{A}_{H}}}}$ on the tetrahedral frame $\left( O;A,B,C,D \right)$ is
	
	\[\beta _{B}^{{{G}_{{{A}_{H}}}}}=\frac{1}{32{{\left( {{S}^{A}} \right)}^{2}}}\left( \begin{aligned}
		& C{{D}^{2}}\left( \Delta _{2}^{A}-C{{D}^{2}} \right)+\left( \Delta _{2}^{A}-D{{B}^{2}} \right)\left( D{{A}^{2}}-B{{A}^{2}} \right) \\ 
		& +\left( \Delta _{2}^{A}-B{{C}^{2}} \right)\left( C{{A}^{2}}-B{{A}^{2}} \right)+8{{\left( {{S}^{A}} \right)}^{2}} \\ 
	\end{aligned} \right),\]
	\[\beta _{C}^{{{G}_{{{A}_{H}}}}}=\frac{1}{32{{\left( {{S}^{A}} \right)}^{2}}}\left( \begin{aligned}
		& D{{B}^{2}}\left( \Delta _{2}^{A}-D{{B}^{2}} \right)+\left( \Delta _{2}^{A}-B{{C}^{2}} \right)\left( B{{A}^{2}}-C{{A}^{2}} \right) \\ 
		& +\left( \Delta _{2}^{A}-C{{D}^{2}} \right)\left( D{{A}^{2}}-C{{A}^{2}} \right)+8{{\left( {{S}^{A}} \right)}^{2}} \\ 
	\end{aligned} \right),\]
	\[\beta _{D}^{{{G}_{{{A}_{H}}}}}=\frac{1}{32{{\left( {{S}^{A}} \right)}^{2}}}\left( \begin{aligned}
		& B{{C}^{2}}\left( \Delta _{2}^{A}-B{{C}^{2}} \right)+\left( \Delta _{2}^{A}-C{{D}^{2}} \right)\left( C{{A}^{2}}-D{{A}^{2}} \right) \\ 
		& +\left( \Delta _{2}^{A}-D{{B}^{2}} \right)\left( B{{A}^{2}}-D{{A}^{2}} \right)+8{{\left( {{S}^{A}} \right)}^{2}} \\ 
	\end{aligned} \right),\]
	\[\beta _{A}^{{{G}_{{{A}_{H}}}}}=0.\]
	Where 
	\[\Delta _{2}^{A}=\frac{1}{2}\left( B{{C}^{2}}+C{{D}^{2}}+D{{B}^{2}} \right).\]
	
	
	2. The frame components of ${{G}_{{{B}_{H}}}}$ on the tetrahedral frame $\left( O;A,B,C,D \right)$ is
	
	\[\beta _{C}^{{{G}_{{{B}_{H}}}}}=\frac{1}{32{{\left( {{S}^{B}} \right)}^{2}}}\left( \begin{aligned}
		& D{{A}^{2}}\left( \Delta _{2}^{B}-D{{A}^{2}} \right)+\left( \Delta _{2}^{B}-A{{C}^{2}} \right)\left( A{{B}^{2}}-C{{B}^{2}} \right) \\ 
		& +\left( \Delta _{2}^{B}-C{{D}^{2}} \right)\left( D{{B}^{2}}-C{{B}^{2}} \right)+8{{\left( {{S}^{B}} \right)}^{2}} \\ 
	\end{aligned} \right),\]
	\[\beta _{D}^{{{G}_{{{B}_{H}}}}}=\frac{1}{32{{\left( {{S}^{B}} \right)}^{2}}}\left( \begin{aligned}
		& A{{C}^{2}}\left( \Delta _{2}^{B}-A{{C}^{2}} \right)+\left( \Delta _{2}^{B}-C{{D}^{2}} \right)\left( C{{B}^{2}}-D{{B}^{2}} \right) \\ 
		& +\left( \Delta _{2}^{B}-D{{A}^{2}} \right)\left( A{{B}^{2}}-D{{B}^{2}} \right)+8{{\left( {{S}^{B}} \right)}^{2}} \\ 
	\end{aligned} \right),\]
	\[\beta _{A}^{{{G}_{{{B}_{H}}}}}=\frac{1}{32{{\left( {{S}^{B}} \right)}^{2}}}\left( \begin{aligned}
		& C{{D}^{2}}\left( \Delta _{2}^{B}-C{{D}^{2}} \right)+\left( \Delta _{2}^{B}-D{{A}^{2}} \right)\left( D{{B}^{2}}-A{{B}^{2}} \right) \\ 
		& +\left( \Delta _{2}^{B}-A{{C}^{2}} \right)\left( C{{B}^{2}}-A{{B}^{2}} \right)+8{{\left( {{S}^{B}} \right)}^{2}} \\ 
	\end{aligned} \right),\]
	\[\ \beta _{B}^{{{P}_{{{B}_{H}}}}}=0.\]
	Where 
	\[\Delta _{2}^{B}=\frac{1}{2}\left( C{{D}^{2}}+D{{A}^{2}}+A{{C}^{2}} \right).\]
	
	
	3. The frame components of ${{G}_{{{C}_{H}}}}$ on the tetrahedral frame $\left( O;A,B,C,D \right)$ is
	
	\[\beta _{D}^{{{G}_{{{C}_{H}}}}}=\frac{1}{32{{\left( {{S}^{C}} \right)}^{2}}}\left( \begin{aligned}
		& A{{B}^{2}}\left( \Delta _{2}^{C}-A{{B}^{2}} \right)+\left( \Delta _{2}^{C}-B{{D}^{2}} \right)\left( B{{C}^{2}}-D{{C}^{2}} \right) \\ 
		& +\left( \Delta _{2}^{C}-D{{A}^{2}} \right)\left( A{{C}^{2}}-D{{C}^{2}} \right)+8{{\left( {{S}^{C}} \right)}^{2}} \\ 
	\end{aligned} \right),\]
	\[\beta _{A}^{{{G}_{{{C}_{H}}}}}=\frac{1}{32{{\left( {{S}^{C}} \right)}^{2}}}\left( \begin{aligned}
		& B{{D}^{2}}\left( \Delta _{2}^{C}-B{{D}^{2}} \right)+\left( \Delta _{2}^{C}-D{{A}^{2}} \right)\left( D{{C}^{2}}-A{{C}^{2}} \right) \\ 
		& +\left( \Delta _{2}^{C}-A{{B}^{2}} \right)\left( B{{C}^{2}}-A{{C}^{2}} \right)+8{{\left( {{S}^{C}} \right)}^{2}} \\ 
	\end{aligned} \right),\]
	\[\beta _{B}^{{{G}_{{{C}_{H}}}}}=\frac{1}{32{{\left( {{S}^{C}} \right)}^{2}}}\left( \begin{aligned}
		& D{{A}^{2}}\left( \Delta _{2}^{C}-D{{A}^{2}} \right)+\left( \Delta _{2}^{C}-A{{B}^{2}} \right)\left( A{{C}^{2}}-B{{C}^{2}} \right) \\ 
		& +\left( \Delta _{2}^{C}-B{{D}^{2}} \right)\left( D{{C}^{2}}-B{{C}^{2}} \right)+8{{\left( {{S}^{C}} \right)}^{2}} \\ 
	\end{aligned} \right),\]
	\[\beta _{C}^{{{G}_{{{C}_{H}}}}}=0.\]
	Where 
	\[\Delta _{2}^{C}=\frac{1}{2}\left( D{{A}^{2}}+A{{B}^{2}}+B{{D}^{2}} \right).\]
	
	
	4. The frame components of ${{G}_{{{D}_{H}}}}$ on the tetrahedral frame $\left( O;A,B,C,D \right)$ is
	
	\[\beta _{A}^{{{G}_{{{D}_{H}}}}}=\frac{1}{32{{\left( {{S}^{D}} \right)}^{2}}}\left( \begin{aligned}
		& B{{C}^{2}}\left( \Delta _{2}^{D}-B{{C}^{2}} \right)+\left( \Delta _{2}^{D}-C{{A}^{2}} \right)\left( C{{D}^{2}}-A{{D}^{2}} \right) \\ 
		& +\left( \Delta _{2}^{D}-A{{B}^{2}} \right)\left( B{{D}^{2}}-A{{D}^{2}} \right)+8{{\left( {{S}^{D}} \right)}^{2}} \\ 
	\end{aligned} \right),\]
	\[\beta _{B}^{{{G}_{{{D}_{H}}}}}=\frac{1}{32{{\left( {{S}^{D}} \right)}^{2}}}\left( \begin{aligned}
		& C{{A}^{2}}\left( \Delta _{2}^{D}-C{{A}^{2}} \right)+\left( \Delta _{2}^{D}-A{{B}^{2}} \right)\left( A{{D}^{2}}-B{{D}^{2}} \right) \\ 
		& +\left( \Delta _{2}^{D}-B{{C}^{2}} \right)\left( C{{D}^{2}}-B{{D}^{2}} \right)+8{{\left( {{S}^{D}} \right)}^{2}} \\ 
	\end{aligned} \right),\]
	\[\beta _{C}^{{{G}_{{{D}_{H}}}}}=\frac{1}{32{{\left( {{S}^{D}} \right)}^{2}}}\left( \begin{aligned}
		& A{{B}^{2}}\left( \Delta _{2}^{D}-A{{B}^{2}} \right)+\left( \Delta _{2}^{D}-B{{C}^{2}} \right)\left( B{{D}^{2}}-C{{D}^{2}} \right) \\ 
		& +\left( \Delta _{2}^{D}-C{{A}^{2}} \right)\left( A{{D}^{2}}-C{{D}^{2}} \right)+8{{\left( {{S}^{D}} \right)}^{2}} \\ 
	\end{aligned} \right).\]
	\[\beta _{D}^{{{G}_{{{D}_{H}}}}}=0.\]
	Where 
	\[\Delta _{2}^{D}=\frac{1}{2}\left( A{{B}^{2}}+B{{C}^{2}}+C{{A}^{2}} \right).\]
\end{theorem}

\begin{proof}
	According to theorem \ref{thm:Thm28.5.3} and theorem \ref{thm:Thm20.4.1}, the conclusion of the theorem can be obtained. 
\end{proof}
\hfill $\square$\par

\subsection{Frame components of the projection points of the incenter of tetrahedron}\label{Subsec28.5.3}


Applying the theorems given above (theorem \ref{thm:Thm28.4.1} and theorem \ref{thm:Thm28.4.2}) to the incenter of the tetrahedron, the frame components of the projection points of the incenter of the tetrahedron can be obtained.


\begin{theorem}{Formula 1 of frame components for projection of incenter, Daiyuan Zhang}{Thm28.5.5}\label{Thm28.5.5} 
	Given a tetrahedron $ABCD$, point $O$ is any point in space. Suppose that the projection point of point $A$ to the $\triangle BCD$ plane of the tetrahedron is ${{A}_{{{A}_{H}}}}$, the projection point of point $B$ on the $\triangle CDA$ plane of tetrahedron is ${{B}_{{{B}_{H}}}}$, the projection point of point $C$ to the $\triangle DAB$ plane of tetrahedron is ${{C}_{{{C}_{H}}}}$, the projection point of point $D$ on the $\triangle ABC$ plane of tetrahedron is ${{D}_{{{D}_{H}}}}$. The projection points of incenter $I$ to the planes $\triangle BCD$, $\triangle CDA$, $\triangle DAB$, $\triangle ABC$ of the tetrahedron are ${{I}_{{{A}_{H}}}}$, ${{I}_{{{B}_{H}}}}$, ${{I}_{{{C}_{H}}}}$, ${{I}_{{{D}_{H}}}}$ respectively, then:	
		
	
	1. 	The frame components of ${{I}_{{{A}_{H}}}}$ on the triangular frame $\left( O;B,C,D \right)$ is:
	
	\[\alpha _{B}^{{{I}_{{{A}_{H}}}}}=\frac{1}{8S{{S}^{A}}}\left( \begin{aligned}
		& C{{D}^{2}}\left( \Delta _{2}^{A}-C{{D}^{2}} \right)+\left( \Delta _{2}^{A}-D{{B}^{2}} \right)\left( D{{A}^{2}}-B{{A}^{2}} \right) \\ 
		& +\left( \Delta _{2}^{A}-B{{C}^{2}} \right)\left( C{{A}^{2}}-B{{A}^{2}} \right)+8{{S}^{B}}{{S}^{A}} \\ 
	\end{aligned} \right),\]
	\[\alpha _{C}^{{{I}_{{{A}_{H}}}}}=\frac{1}{8S{{S}^{A}}}\left( \begin{aligned}
		& D{{B}^{2}}\left( \Delta _{2}^{A}-D{{B}^{2}} \right)+\left( \Delta _{2}^{A}-B{{C}^{2}} \right)\left( B{{A}^{2}}-C{{A}^{2}} \right) \\ 
		& +\left( \Delta _{2}^{A}-C{{D}^{2}} \right)\left( D{{A}^{2}}-C{{A}^{2}} \right)+8{{S}^{C}}{{S}^{A}} \\ 
	\end{aligned} \right),\]
	\[\alpha _{D}^{{{I}_{{{A}_{H}}}}}=\frac{1}{8S{{S}^{A}}}\left( \begin{aligned}
		& B{{C}^{2}}\left( \Delta _{2}^{A}-B{{C}^{2}} \right)+\left( \Delta _{2}^{A}-C{{D}^{2}} \right)\left( C{{A}^{2}}-D{{A}^{2}} \right) \\ 
		& +\left( \Delta _{2}^{A}-D{{B}^{2}} \right)\left( B{{A}^{2}}-D{{A}^{2}} \right)+8{{S}^{D}}{{S}^{A}} \\ 
	\end{aligned} \right).\]
	Where 
	\[\Delta _{2}^{A}=\frac{1}{2}\left( B{{C}^{2}}+C{{D}^{2}}+D{{B}^{2}} \right).\]
	
	
	2. 	The frame components of ${{I}_{{{B}_{H}}}}$ on the triangular frame $\left( O;C,D,A \right)$ is:
		
	\[\alpha _{C}^{{{I}_{{{B}_{H}}}}}=\frac{1}{8S{{S}^{B}}}\left( \begin{aligned}
		& D{{A}^{2}}\left( \Delta _{2}^{B}-D{{A}^{2}} \right)+\left( \Delta _{2}^{B}-A{{C}^{2}} \right)\left( A{{B}^{2}}-C{{B}^{2}} \right) \\ 
		& +\left( \Delta _{2}^{B}-C{{D}^{2}} \right)\left( D{{B}^{2}}-C{{B}^{2}} \right)+8{{S}^{C}}{{S}^{B}} \\ 
	\end{aligned} \right),\]
	\[\alpha _{D}^{{{I}_{{{B}_{H}}}}}=\frac{1}{8S{{S}^{B}}}\left( \begin{aligned}
		& A{{C}^{2}}\left( \Delta _{2}^{B}-A{{C}^{2}} \right)+\left( \Delta _{2}^{B}-C{{D}^{2}} \right)\left( C{{B}^{2}}-D{{B}^{2}} \right) \\ 
		& +\left( \Delta _{2}^{B}-D{{A}^{2}} \right)\left( A{{B}^{2}}-D{{B}^{2}} \right)+8{{S}^{D}}{{S}^{B}} \\ 
	\end{aligned} \right),\]
	\[\alpha _{A}^{{{I}_{{{B}_{H}}}}}=\frac{1}{8S{{S}^{B}}}\left( \begin{aligned}
		& C{{D}^{2}}\left( \Delta _{2}^{B}-C{{D}^{2}} \right)+\left( \Delta _{2}^{B}-D{{A}^{2}} \right)\left( D{{B}^{2}}-A{{B}^{2}} \right) \\ 
		& +\left( \Delta _{2}^{B}-A{{C}^{2}} \right)\left( C{{B}^{2}}-A{{B}^{2}} \right)+8{{S}^{A}}{{S}^{B}} \\ 
	\end{aligned} \right).\]
	Where 
	\[\Delta _{2}^{B}=\frac{1}{2}\left( C{{D}^{2}}+D{{A}^{2}}+A{{C}^{2}} \right).\]
	
	
	3. 	The frame components of ${{I}_{{{C}_{H}}}}$ on the triangular frame $\left( O;D,A,B \right)$ is:
		
	\[\alpha _{D}^{{{I}_{{{C}_{H}}}}}=\frac{1}{8S{{S}^{C}}}\left( \begin{aligned}
		& A{{B}^{2}}\left( \Delta _{2}^{C}-A{{B}^{2}} \right)+\left( \Delta _{2}^{C}-B{{D}^{2}} \right)\left( B{{C}^{2}}-D{{C}^{2}} \right) \\ 
		& +\left( \Delta _{2}^{C}-D{{A}^{2}} \right)\left( A{{C}^{2}}-D{{C}^{2}} \right)+8{{S}^{D}}{{S}^{C}} \\ 
	\end{aligned} \right),\]
	\[\alpha _{A}^{{{I}_{{{C}_{H}}}}}=\frac{1}{8S{{S}^{C}}}\left( \begin{aligned}
		& B{{D}^{2}}\left( \Delta _{2}^{C}-B{{D}^{2}} \right)+\left( \Delta _{2}^{C}-D{{A}^{2}} \right)\left( D{{C}^{2}}-A{{C}^{2}} \right) \\ 
		& +\left( \Delta _{2}^{C}-A{{B}^{2}} \right)\left( B{{C}^{2}}-A{{C}^{2}} \right)+8{{S}^{A}}{{S}^{C}} \\ 
	\end{aligned} \right),\]
	\[\alpha _{B}^{{{I}_{{{C}_{H}}}}}=\frac{1}{8S{{S}^{C}}}\left( \begin{aligned}
		& D{{A}^{2}}\left( \Delta _{2}^{C}-D{{A}^{2}} \right)+\left( \Delta _{2}^{C}-A{{B}^{2}} \right)\left( A{{C}^{2}}-B{{C}^{2}} \right) \\ 
		& +\left( \Delta _{2}^{C}-B{{D}^{2}} \right)\left( D{{C}^{2}}-B{{C}^{2}} \right)+8{{S}^{B}}{{S}^{C}} \\ 
	\end{aligned} \right).\]
	Where 
	\[\Delta _{2}^{C}=\frac{1}{2}\left( D{{A}^{2}}+A{{B}^{2}}+B{{D}^{2}} \right).\]
	
	
	4. 	The frame components of ${{I}_{{{D}_{H}}}}$ on the triangular frame $\left( O;A,B,C \right)$ is:
	
	\[\alpha _{A}^{{{I}_{{{D}_{H}}}}}=\frac{1}{8S{{S}^{D}}}\left( \begin{aligned}
		& B{{C}^{2}}\left( \Delta _{2}^{D}-B{{C}^{2}} \right)+\left( \Delta _{2}^{D}-C{{A}^{2}} \right)\left( C{{D}^{2}}-A{{D}^{2}} \right) \\ 
		& +\left( \Delta _{2}^{D}-A{{B}^{2}} \right)\left( B{{D}^{2}}-A{{D}^{2}} \right)+8{{S}^{A}}{{S}^{D}} \\ 
	\end{aligned} \right),\]
	\[\alpha _{B}^{{{I}_{{{D}_{H}}}}}=\frac{1}{8S{{S}^{D}}}\left( \begin{aligned}
		& C{{A}^{2}}\left( \Delta _{2}^{D}-C{{A}^{2}} \right)+\left( \Delta _{2}^{D}-A{{B}^{2}} \right)\left( A{{D}^{2}}-B{{D}^{2}} \right) \\ 
		& +\left( \Delta _{2}^{D}-B{{C}^{2}} \right)\left( C{{D}^{2}}-B{{D}^{2}} \right)+8{{S}^{B}}{{S}^{D}} \\ 
	\end{aligned} \right),\]
	\[\alpha _{C}^{{{I}_{{{D}_{H}}}}}=\frac{1}{8S{{S}^{D}}}\left( \begin{aligned}
		& A{{B}^{2}}\left( \Delta _{2}^{D}-A{{B}^{2}} \right)+\left( \Delta _{2}^{D}-B{{C}^{2}} \right)\left( B{{D}^{2}}-C{{D}^{2}} \right) \\ 
		& +\left( \Delta _{2}^{D}-C{{A}^{2}} \right)\left( A{{D}^{2}}-C{{D}^{2}} \right)+8{{S}^{C}}{{S}^{D}} \\ 
	\end{aligned} \right).\]
	Where 
	\[\Delta _{2}^{D}=\frac{1}{2}\left( A{{B}^{2}}+B{{C}^{2}}+C{{A}^{2}} \right).\]
\end{theorem}

\begin{proof}
	
	
	Here's the proof for the frame components of the projection point ${{I}_{{{D}_{H}}}}$ on the $\triangle ABC$ plane.
	
	According to theorem \ref{thm:Thm28.4.1} we get:
	
	\[\left( \begin{matrix}
		\alpha _{A}^{{{I}_{{{D}_{H}}}}}  \\
		\alpha _{B}^{{{I}_{{{D}_{H}}}}}  \\
		\alpha _{C}^{{{I}_{{{D}_{H}}}}}  \\
	\end{matrix} \right)=\left( \begin{matrix}
		1 & 0 & 0 & \alpha _{A}^{{{D}_{{{D}_{H}}}}}  \\
		0 & 1 & 0 & \alpha _{B}^{{{D}_{{{D}_{H}}}}}  \\
		0 & 0 & 1 & \alpha _{C}^{{{D}_{{{D}_{H}}}}}  \\
	\end{matrix} \right)\left( \begin{matrix}
		\beta _{A}^{I}  \\
		\beta _{B}^{I}  \\
		\beta _{C}^{I}  \\
		\beta _{D}^{I}  \\
	\end{matrix} \right).\]
	
	
	Replace theorem \ref{thm:Thm28.3.1} and theorem \ref{thm:Thm21.2.1} directly into the above formula to obtain:
	
	\[\alpha _{A}^{{{I}_{{{D}_{H}}}}}=\frac{1}{8S{{S}^{D}}}\left( \begin{aligned}
		& B{{C}^{2}}\left( \Delta _{2}^{D}-B{{C}^{2}} \right)+\left( \Delta _{2}^{D}-C{{A}^{2}} \right)\left( C{{D}^{2}}-A{{D}^{2}} \right) \\ 
		& +\left( \Delta _{2}^{D}-A{{B}^{2}} \right)\left( B{{D}^{2}}-A{{D}^{2}} \right)+8{{S}^{A}}{{S}^{D}} \\ 
	\end{aligned} \right),\]
	\[\alpha _{B}^{{{I}_{{{D}_{H}}}}}=\frac{1}{8S{{S}^{D}}}\left( \begin{aligned}
		& C{{A}^{2}}\left( \Delta _{2}^{D}-C{{A}^{2}} \right)+\left( \Delta _{2}^{D}-A{{B}^{2}} \right)\left( A{{D}^{2}}-B{{D}^{2}} \right) \\ 
		& +\left( \Delta _{2}^{D}-B{{C}^{2}} \right)\left( C{{D}^{2}}-B{{D}^{2}} \right)+8{{S}^{B}}{{S}^{D}} \\ 
	\end{aligned} \right),\]
	\[\alpha _{C}^{{{I}_{{{D}_{H}}}}}=\frac{1}{8S{{S}^{D}}}\left( \begin{aligned}
		& A{{B}^{2}}\left( \Delta _{2}^{D}-A{{B}^{2}} \right)+\left( \Delta _{2}^{D}-B{{C}^{2}} \right)\left( B{{D}^{2}}-C{{D}^{2}} \right) \\ 
		& +\left( \Delta _{2}^{D}-C{{A}^{2}} \right)\left( A{{D}^{2}}-C{{D}^{2}} \right)+8{{S}^{C}}{{S}^{D}} \\ 
	\end{aligned} \right).\]
	
	Other formulas can be proved similarly.	
\end{proof}
\hfill $\square$\par

\begin{theorem}{Formula 2 of frame components for projection of incenter, Daiyuan Zhang}{Thm28.5.6}\label{Thm28.5.6} 
	Given a tetrahedron $ABCD$, point $O$ is any point in space. Suppose that the projection point of point $A$ to the $\triangle BCD$ plane of the tetrahedron is ${{A}_{{{A}_{H}}}}$, the projection point of point $B$ on the $\triangle CDA$ plane of tetrahedron is ${{B}_{{{B}_{H}}}}$, the projection point of point $C$ to the $\triangle DAB$ plane of tetrahedron is ${{C}_{{{C}_{H}}}}$, the projection point of point $D$ on the $\triangle ABC$ plane of tetrahedron is ${{D}_{{{D}_{H}}}}$. The projection points of incenter $I$ to the planes $\triangle BCD$, $\triangle CDA$, $\triangle DAB$, $\triangle ABC$ of the tetrahedron are ${{I}_{{{A}_{H}}}}$, ${{I}_{{{B}_{H}}}}$, ${{I}_{{{C}_{H}}}}$, ${{I}_{{{D}_{H}}}}$ respectively, then:	
	
	
	1. The frame components of ${{I}_{{{A}_{H}}}}$ on the tetrahedral frame $\left( O;A,B,C,D \right)$ is
	
	\[\beta _{B}^{{{I}_{{{A}_{H}}}}}=\frac{1}{8S{{S}^{A}}}\left( \begin{aligned}
		& C{{D}^{2}}\left( \Delta _{2}^{A}-C{{D}^{2}} \right)+\left( \Delta _{2}^{A}-D{{B}^{2}} \right)\left( D{{A}^{2}}-B{{A}^{2}} \right) \\ 
		& +\left( \Delta _{2}^{A}-B{{C}^{2}} \right)\left( C{{A}^{2}}-B{{A}^{2}} \right)+8{{S}^{B}}{{S}^{A}} \\ 
	\end{aligned} \right),\]
	\[\beta _{C}^{{{I}_{{{A}_{H}}}}}=\frac{1}{8S{{S}^{A}}}\left( \begin{aligned}
		& D{{B}^{2}}\left( \Delta _{2}^{A}-D{{B}^{2}} \right)+\left( \Delta _{2}^{A}-B{{C}^{2}} \right)\left( B{{A}^{2}}-C{{A}^{2}} \right) \\ 
		& +\left( \Delta _{2}^{A}-C{{D}^{2}} \right)\left( D{{A}^{2}}-C{{A}^{2}} \right)+8{{S}^{C}}{{S}^{A}} \\ 
	\end{aligned} \right),\]
	\[\beta _{D}^{{{I}_{{{A}_{H}}}}}=\frac{1}{8S{{S}^{A}}}\left( \begin{aligned}
		& B{{C}^{2}}\left( \Delta _{2}^{A}-B{{C}^{2}} \right)+\left( \Delta _{2}^{A}-C{{D}^{2}} \right)\left( C{{A}^{2}}-D{{A}^{2}} \right) \\ 
		& +\left( \Delta _{2}^{A}-D{{B}^{2}} \right)\left( B{{A}^{2}}-D{{A}^{2}} \right)+8{{S}^{D}}{{S}^{A}} \\ 
	\end{aligned} \right),\]
	\[\beta _{A}^{{{I}_{{{A}_{H}}}}}=0.\]
	Where 
	\[\Delta _{2}^{A}=\frac{1}{2}\left( B{{C}^{2}}+C{{D}^{2}}+D{{B}^{2}} \right).\]
	
	
	2. The frame components of ${{I}_{{{B}_{H}}}}$ on the tetrahedral frame $\left( O;A,B,C,D \right)$ is
	
	\[\beta _{C}^{{{I}_{{{B}_{H}}}}}=\frac{1}{8S{{S}^{B}}}\left( \begin{aligned}
		& D{{A}^{2}}\left( \Delta _{2}^{B}-D{{A}^{2}} \right)+\left( \Delta _{2}^{B}-A{{C}^{2}} \right)\left( A{{B}^{2}}-C{{B}^{2}} \right) \\ 
		& +\left( \Delta _{2}^{B}-C{{D}^{2}} \right)\left( D{{B}^{2}}-C{{B}^{2}} \right)+8{{S}^{C}}{{S}^{B}} \\ 
	\end{aligned} \right),\]
	\[\beta _{D}^{{{I}_{{{B}_{H}}}}}=\frac{1}{8S{{S}^{B}}}\left( \begin{aligned}
		& A{{C}^{2}}\left( \Delta _{2}^{B}-A{{C}^{2}} \right)+\left( \Delta _{2}^{B}-C{{D}^{2}} \right)\left( C{{B}^{2}}-D{{B}^{2}} \right) \\ 
		& +\left( \Delta _{2}^{B}-D{{A}^{2}} \right)\left( A{{B}^{2}}-D{{B}^{2}} \right)+8{{S}^{D}}{{S}^{B}} \\ 
	\end{aligned} \right),\]
	\[\beta _{A}^{{{I}_{{{B}_{H}}}}}=\frac{1}{8S{{S}^{B}}}\left( \begin{aligned}
		& C{{D}^{2}}\left( \Delta _{2}^{B}-C{{D}^{2}} \right)+\left( \Delta _{2}^{B}-D{{A}^{2}} \right)\left( D{{B}^{2}}-A{{B}^{2}} \right) \\ 
		& +\left( \Delta _{2}^{B}-A{{C}^{2}} \right)\left( C{{B}^{2}}-A{{B}^{2}} \right)+8{{S}^{A}}{{S}^{B}} \\ 
	\end{aligned} \right),\]
	\[\ \beta _{B}^{{{P}_{{{B}_{H}}}}}=0.\]
	Where 
	\[\Delta _{2}^{B}=\frac{1}{2}\left( C{{D}^{2}}+D{{A}^{2}}+A{{C}^{2}} \right).\]
	
	
	3. The frame components of ${{I}_{{{C}_{H}}}}$ on the tetrahedral frame $\left( O;A,B,C,D \right)$ is
	
	\[\beta _{D}^{{{I}_{{{C}_{H}}}}}=\frac{1}{8S{{S}^{C}}}\left( \begin{aligned}
		& A{{B}^{2}}\left( \Delta _{2}^{C}-A{{B}^{2}} \right)+\left( \Delta _{2}^{C}-B{{D}^{2}} \right)\left( B{{C}^{2}}-D{{C}^{2}} \right) \\ 
		& +\left( \Delta _{2}^{C}-D{{A}^{2}} \right)\left( A{{C}^{2}}-D{{C}^{2}} \right)+8{{S}^{D}}{{S}^{C}} \\ 
	\end{aligned} \right),\]
	\[\beta _{A}^{{{I}_{{{C}_{H}}}}}=\frac{1}{8S{{S}^{C}}}\left( \begin{aligned}
		& B{{D}^{2}}\left( \Delta _{2}^{C}-B{{D}^{2}} \right)+\left( \Delta _{2}^{C}-D{{A}^{2}} \right)\left( D{{C}^{2}}-A{{C}^{2}} \right) \\ 
		& +\left( \Delta _{2}^{C}-A{{B}^{2}} \right)\left( B{{C}^{2}}-A{{C}^{2}} \right)+8{{S}^{A}}{{S}^{C}} \\ 
	\end{aligned} \right),\]
	\[\beta _{B}^{{{I}_{{{C}_{H}}}}}=\frac{1}{8S{{S}^{C}}}\left( \begin{aligned}
		& D{{A}^{2}}\left( \Delta _{2}^{C}-D{{A}^{2}} \right)+\left( \Delta _{2}^{C}-A{{B}^{2}} \right)\left( A{{C}^{2}}-B{{C}^{2}} \right) \\ 
		& +\left( \Delta _{2}^{C}-B{{D}^{2}} \right)\left( D{{C}^{2}}-B{{C}^{2}} \right)+8{{S}^{B}}{{S}^{C}} \\ 
	\end{aligned} \right),\]
	\[\beta _{C}^{{{P}_{{{C}_{H}}}}}=0.\]
	Where 
	\[\Delta _{2}^{C}=\frac{1}{2}\left( D{{A}^{2}}+A{{B}^{2}}+B{{D}^{2}} \right).\]
	
	
	4. The frame components of ${{I}_{{{D}_{H}}}}$ on the tetrahedral frame $\left( O;A,B,C,D \right)$ is
	
	\[\beta _{A}^{{{I}_{{{D}_{H}}}}}=\frac{1}{8S{{S}^{D}}}\left( \begin{aligned}
		& B{{C}^{2}}\left( \Delta _{2}^{D}-B{{C}^{2}} \right)+\left( \Delta _{2}^{D}-C{{A}^{2}} \right)\left( C{{D}^{2}}-A{{D}^{2}} \right) \\ 
		& +\left( \Delta _{2}^{D}-A{{B}^{2}} \right)\left( B{{D}^{2}}-A{{D}^{2}} \right)+8{{S}^{A}}{{S}^{D}} \\ 
	\end{aligned} \right),\]
	\[\beta _{B}^{{{I}_{{{D}_{H}}}}}=\frac{1}{8S{{S}^{D}}}\left( \begin{aligned}
		& C{{A}^{2}}\left( \Delta _{2}^{D}-C{{A}^{2}} \right)+\left( \Delta _{2}^{D}-A{{B}^{2}} \right)\left( A{{D}^{2}}-B{{D}^{2}} \right) \\ 
		& +\left( \Delta _{2}^{D}-B{{C}^{2}} \right)\left( C{{D}^{2}}-B{{D}^{2}} \right)+8{{S}^{B}}{{S}^{D}} \\ 
	\end{aligned} \right),\]
	\[\beta _{C}^{{{I}_{{{D}_{H}}}}}=\frac{1}{8S{{S}^{D}}}\left( \begin{aligned}
		& A{{B}^{2}}\left( \Delta _{2}^{D}-A{{B}^{2}} \right)+\left( \Delta _{2}^{D}-B{{C}^{2}} \right)\left( B{{D}^{2}}-C{{D}^{2}} \right) \\ 
		& +\left( \Delta _{2}^{D}-C{{A}^{2}} \right)\left( A{{D}^{2}}-C{{D}^{2}} \right)+8{{S}^{C}}{{S}^{D}} \\ 
	\end{aligned} \right),\]
	\[\beta _{D}^{{{I}_{{{D}_{H}}}}}=0.\]
	Where 
	\[\Delta _{2}^{D}=\frac{1}{2}\left( A{{B}^{2}}+B{{C}^{2}}+C{{A}^{2}} \right).\]
\end{theorem}

\begin{proof}
	According to theorem \ref{thm:Thm28.5.5} and theorem \ref{thm:Thm20.4.1}, the conclusion of the theorem can be obtained. 	
\end{proof}
\hfill $\square$\par

'
\chapter{Inscribed sphere of tetrahedron}\label{Ch29}
\thispagestyle{empty}


Although for some special tetrahedron, such as regular triangular pyramid or regular tetrahedron, the radius of the inscribed sphere can be obtained by a simple method (Euclidean geometry method), for a general tetrahedron, given the lengths of the six edges of the tetrahedron, how to calculate the radius of inscribed sphere of the tetrahedron? 


I answered this question in this chapter. The results show that the radius of inscribed sphere of the tetrahedron can be expressed by a formula, and this formula can be obtained by the lengths of the six edges of the tetrahedron through finite operations of addition, subtraction, multiplication, division and radical.

\section{Radius of inscribed sphere of tetrahedron}

As long as we know the lengths of the six edges of a tetrahedron, we can get the radius of its inscribed sphere. I give the following theorem.


\begin{theorem}{Formula of radius for inscribed sphere of tetrahedron, Daiyuan Zhang}{Thm29.1.1}\label{Thm29.1.1} 
	Given a tetrahedron $ABCD$, let $r$ be the radius of inscribed sphere of tetrahedron $ABCD$, $S$ be the surface area of tetrahedron $ABCD$, and $AB$, $AC$, $AD$, $BC$, $CD$, $DB$ be the lengths of the six edges of tetrahedron $ABCD$, then
	\[{{r}^{2}}=\frac{1}{4{{S}^{2}}}\left( {{t}_{1}}-{{t}_{2}}-{{t}_{3}} \right).\]
	Where 
	\[{{t}_{1}}={{q}_{2}}{{\Delta }_{2}},\]
	\[{{q}_{2}}=\frac{1}{2}\left( A{{B}^{2}}C{{D}^{2}}+B{{C}^{2}}A{{D}^{2}}+C{{A}^{2}}B{{D}^{2}} \right),\]
	\[{{\Delta }_{2}}=\frac{1}{2}\left( A{{B}^{2}}+A{{C}^{2}}+A{{D}^{2}}+B{{C}^{2}}+C{{D}^{2}}+D{{B}^{2}} \right),\]
	\[{{t}_{2}}=\frac{1}{2}\left( A{{B}^{2}}C{{D}^{2}}\left( A{{B}^{2}}+C{{D}^{2}} \right)+B{{C}^{2}}A{{D}^{2}}\left( B{{C}^{2}}+A{{D}^{2}} \right)+C{{A}^{2}}B{{D}^{2}}\left( C{{A}^{2}}+B{{D}^{2}} \right) \right),\]
	\[{{t}_{3}}=\frac{1}{4}\left( A{{B}^{2}}B{{C}^{2}}C{{A}^{2}}+B{{C}^{2}}C{{D}^{2}}D{{B}^{2}}+C{{D}^{2}}D{{A}^{2}}A{{C}^{2}}+D{{A}^{2}}A{{B}^{2}}B{{D}^{2}} \right).\]
\end{theorem}

\begin{proof}
	Consider The frame components of the projection point ${{I}_{{{D}_{H}}}}$ on the $\triangle ABC$ plane, according to theorem \ref{thm:Thm24.2.1}, the radius of the inscribed sphere can be obtained by the following formula:
	\[\begin{aligned}
		{{r}^{2}}& =II_{{{D}_{H}}}^{2}=-\beta _{A}^{I{{I}_{{{D}_{H}}}}}\beta _{B}^{I{{I}_{{{D}_{H}}}}}A{{B}^{2}}-\beta _{A}^{I{{I}_{{{D}_{H}}}}}\beta _{C}^{I{{I}_{{{D}_{H}}}}}A{{C}^{2}}-\beta _{A}^{I{{I}_{{{D}_{H}}}}}\beta _{D}^{I{{I}_{{{D}_{H}}}}}A{{D}^{2}} \\ 
		& -\beta _{B}^{I{{I}_{{{D}_{H}}}}}\beta _{C}^{I{{I}_{{{D}_{H}}}}}B{{C}^{2}}-\beta _{B}^{I{{I}_{{{D}_{H}}}}}\beta _{D}^{I{{I}_{{{D}_{H}}}}}B{{D}^{2}}-\beta _{C}^{I{{I}_{{{D}_{H}}}}}\beta _{D}^{I{{I}_{{{D}_{H}}}}}C{{D}^{2}}.  
	\end{aligned}\]
	Where 
	\[\beta _{A}^{I{{I}_{{{D}_{H}}}}}=\beta _{A}^{{{I}_{{{D}_{H}}}}}-\beta _{A}^{I},\]	
	\[\beta _{B}^{I{{I}_{{{D}_{H}}}}}=\beta _{B}^{{{I}_{{{D}_{H}}}}}-\beta _{B}^{I},\]	
	\[\beta _{C}^{I{{I}_{{{D}_{H}}}}}=\beta _{C}^{{{I}_{{{D}_{H}}}}}-\beta _{C}^{I},\]	
	\[\beta _{D}^{I{{I}_{{{D}_{H}}}}}=\beta _{D}^{{{I}_{{{D}_{H}}}}}-\beta _{D}^{I}.\]	
	
	
	Where $\beta _{A}^{I}$, $\beta _{B}^{I}$, $\beta _{C}^{I}$, $\beta _{D}^{I}$ are the frame components of incenter $I$ in tetrahedral frame $\left( O;A,B,C,D \right)$ respectively; and $\beta _{A}^{{{I}_{{{D}_{H}}}}}$, $\beta _{B}^{{{I}_{{{D}_{H}}}}}$, $\beta _{C}^{{{I}_{{{D}_{H}}}}}$, $\beta _{D}^{{{I}_{{{D}_{H}}}}}$ are the frame components of the projection point ${{I}_{{{D}_{H}}}}$ of incenter $I$ in tetrahedral frame $\left( O;A,B,C,D \right)$ respectively.
		
	
	According to theorem \ref{thm:Thm28.4.2}, theorem \ref{thm:Thm24.2.1} and theorem \ref{thm:Thm21.2.1}, we have:
	\[\beta _{A}^{I{{I}_{{{D}_{H}}}}}=\beta _{A}^{{{I}_{{{D}_{H}}}}}-\beta _{A}^{I}=\beta _{A}^{I}+\beta _{A}^{{{I}_{{{D}_{H}}}}}\beta _{D}^{I}-\beta _{A}^{I}=\beta _{A}^{{{I}_{{{D}_{H}}}}}\beta _{D}^{I}=\frac{\beta _{A}^{{{I}_{{{D}_{H}}}}}{{S}^{D}}}{S},\]	
	\[\beta _{B}^{I{{I}_{{{D}_{H}}}}}=\beta _{B}^{{{I}_{{{D}_{H}}}}}-\beta _{B}^{I}=\beta _{B}^{I}+\beta _{B}^{{{I}_{{{D}_{H}}}}}\beta _{D}^{I}-\beta _{B}^{I}=\beta _{B}^{{{I}_{{{D}_{H}}}}}\beta _{D}^{I}=\frac{\beta _{B}^{{{I}_{{{D}_{H}}}}}{{S}^{D}}}{S},\]	
	\[\beta _{C}^{I{{I}_{{{D}_{H}}}}}=\beta _{C}^{{{I}_{{{D}_{H}}}}}-\beta _{C}^{I}=\beta _{C}^{I}+\beta _{C}^{{{I}_{{{D}_{H}}}}}\beta _{D}^{I}-\beta _{C}^{I}=\beta _{C}^{{{I}_{{{D}_{H}}}}}\beta _{D}^{I}=\frac{\beta _{C}^{{{I}_{{{D}_{H}}}}}{{S}^{D}}}{S},\]	
	\[\beta _{D}^{I{{I}_{{{D}_{H}}}}}=\beta _{D}^{{{I}_{{{D}_{H}}}}}-\beta _{D}^{I}=0-\beta _{D}^{I}=-\beta _{D}^{I}=-\frac{{{S}^{D}}}{S}.\]	
	
	From theorem \ref{thm:Thm28.5.6} we have:
	\[\beta _{A}^{{{I}_{{{D}_{H}}}}}=\frac{1}{8S{{S}^{D}}}\left( \begin{aligned}
		& \left( \Delta _{2}^{D}-B{{C}^{2}} \right)B{{C}^{2}}+\left( \Delta _{2}^{D}-C{{A}^{2}} \right)\left( C{{D}^{2}}-A{{D}^{2}} \right) \\ 
		& +\left( \Delta _{2}^{D}-A{{B}^{2}} \right)\left( B{{D}^{2}}-A{{D}^{2}} \right)+8{{S}^{A}}{{S}^{D}} \\ 
	\end{aligned} \right),\]
	\[\beta _{B}^{{{I}_{{{D}_{H}}}}}=\frac{1}{8S{{S}^{D}}}\left( \begin{aligned}
		& \left( \Delta _{2}^{D}-C{{A}^{2}} \right)C{{A}^{2}}+\left( \Delta _{2}^{D}-A{{B}^{2}} \right)\left( A{{D}^{2}}-B{{D}^{2}} \right) \\ 
		& +\left( \Delta _{2}^{D}-B{{C}^{2}} \right)\left( C{{D}^{2}}-B{{D}^{2}} \right)+8{{S}^{B}}{{S}^{D}} \\ 
	\end{aligned} \right),\]
	\[\beta _{C}^{{{I}_{{{D}_{H}}}}}=\frac{1}{8S{{S}^{D}}}\left( \begin{aligned}
		& \left( \Delta _{2}^{D}-A{{B}^{2}} \right)A{{B}^{2}}+\left( \Delta _{2}^{D}-B{{C}^{2}} \right)\left( B{{D}^{2}}-C{{D}^{2}} \right) \\ 
		& +\left( \Delta _{2}^{D}-C{{A}^{2}} \right)\left( A{{D}^{2}}-C{{D}^{2}} \right)+8{{S}^{C}}{{S}^{D}} \\ 
	\end{aligned} \right),\]
	\[\beta _{D}^{{{I}_{{{D}_{H}}}}}=0.\]
	Where 
	\[\Delta _{2}^{D}=\frac{1}{2}\left( A{{B}^{2}}+B{{C}^{2}}+C{{A}^{2}} \right).\]
	
	Therefore
	\[{{r}^{2}}=\frac{1}{16{{S}^{2}}}\left( \begin{aligned}
		& \text{2}A{{B}^{2}}C{{D}^{2}}\left( {{\Delta }_{2}}-A{{B}^{2}}-C{{D}^{2}} \right)+2B{{C}^{2}}A{{D}^{2}}\left( {{\Delta }_{2}}-B{{C}^{2}}-A{{D}^{2}} \right) \\ 
		& +2C{{A}^{2}}B{{D}^{2}}\left( {{\Delta }_{2}}-C{{A}^{2}}-B{{D}^{2}} \right)-A{{B}^{2}}B{{C}^{2}}C{{A}^{2}}-B{{C}^{2}}C{{D}^{2}}D{{B}^{2}} \\ 
		& -C{{D}^{2}}D{{A}^{2}}A{{C}^{2}}-D{{A}^{2}}A{{B}^{2}}B{{D}^{2}}  
	\end{aligned} \right).\]
	Where ${{\Delta }_{2}}$ is half of the square sum of the six edges, i.e.
	\[{{\Delta }_{2}}=\frac{1}{2}\left( A{{B}^{2}}+A{{C}^{2}}+A{{D}^{2}}+B{{C}^{2}}+C{{D}^{2}}+D{{B}^{2}} \right).\]
	
	Let ${{q}_{2}}$ be half of the square sum of the opposite edges of the tetrahedron, i.e.
	\[{{q}_{2}}=\frac{1}{2}\left( A{{B}^{2}}C{{D}^{2}}+B{{C}^{2}}A{{D}^{2}}+C{{A}^{2}}B{{D}^{2}} \right).\]
	Then: 
	\[{{r}^{2}}=\frac{1}{16{{S}^{2}}}\left( \begin{aligned}
		& 2\left( A{{B}^{2}}C{{D}^{2}}+B{{C}^{2}}A{{D}^{2}}+C{{A}^{2}}B{{D}^{2}} \right){{\Delta }_{2}} \\ 
		& -2\left( \begin{aligned}
			& A{{B}^{2}}C{{D}^{2}}\left( A{{B}^{2}}+C{{D}^{2}} \right)+B{{C}^{2}}A{{D}^{2}}\left( B{{C}^{2}}+A{{D}^{2}} \right) \\ 
			& +C{{A}^{2}}B{{D}^{2}}\left( C{{A}^{2}}+B{{D}^{2}} \right) \\ 
		\end{aligned} \right) \\ 
		& -\left( \begin{aligned}
			& A{{B}^{2}}B{{C}^{2}}C{{A}^{2}}+B{{C}^{2}}C{{D}^{2}}D{{B}^{2}} \\ 
			& +C{{D}^{2}}D{{A}^{2}}A{{C}^{2}}+D{{A}^{2}}A{{B}^{2}}B{{D}^{2}} \\ 
		\end{aligned} \right)  
	\end{aligned} \right),\]
	i.e. 
	\[{{r}^{2}}=\frac{1}{4{{S}^{2}}}\left( {{t}_{1}}-{{t}_{2}}-{{t}_{3}} \right),\]
	Where 
	\[{{t}_{1}}={{q}_{2}}{{\Delta }_{2}},\]
	\[{{q}_{2}}=\frac{1}{2}\left( A{{B}^{2}}C{{D}^{2}}+B{{C}^{2}}A{{D}^{2}}+C{{A}^{2}}B{{D}^{2}} \right),\]
	\[{{\Delta }_{2}}=\frac{1}{2}\left( A{{B}^{2}}+A{{C}^{2}}+A{{D}^{2}}+B{{C}^{2}}+C{{D}^{2}}+D{{B}^{2}} \right),\]
	\[{{t}_{2}}=\frac{1}{2}\left( A{{B}^{2}}C{{D}^{2}}\left( A{{B}^{2}}+C{{D}^{2}} \right)+B{{C}^{2}}A{{D}^{2}}\left( B{{C}^{2}}+A{{D}^{2}} \right)+C{{A}^{2}}B{{D}^{2}}\left( C{{A}^{2}}+B{{D}^{2}} \right) \right),\]
	\[{{t}_{3}}=\frac{1}{4}\left( A{{B}^{2}}B{{C}^{2}}C{{A}^{2}}+B{{C}^{2}}C{{D}^{2}}D{{B}^{2}}+C{{D}^{2}}D{{A}^{2}}A{{C}^{2}}+D{{A}^{2}}A{{B}^{2}}B{{D}^{2}} \right).\]
\end{proof}
\hfill $\square$\par

Each quantity in theorem \ref{thm:Thm29.1.1} is only related to the lengths of the six edges of the tetrahedron (many quantities are related to the opposite edges) and the areas of the four triangles, which are “natural parameters” that can be touched by the senses. Each face of a tetrahedron is a triangle. As long as the length of each edge is known, the area of the triangle can be calculated by Helen's formula. Therefore, for any tetrahedron, as long as the lengths of the six edges of the tetrahedron are given, the radius of its inscribed sphere can be obtained. I am sure that the theorem I put forward here has not only theoretical value, but also application value.

A simple application of theorem \ref{thm:Thm29.1.1} is given below. For a regular tetrahedron with a length of edge $a$, the radius of its inscribed sphere is
\[\begin{aligned}
	& r=\frac{1}{2S}\sqrt{{{t}_{1}}-{{t}_{2}}-{{t}_{3}}}=\frac{{{a}^{3}}}{2S}\sqrt{\frac{6}{2}\times \frac{3}{2}-3-1} \\ 
	& =\frac{{{a}^{3}}}{2\sqrt{2}S}=\frac{{{a}^{3}}}{2\sqrt{2}\left( 4\times \frac{1}{2}\times \frac{\sqrt{3}}{2}{{a}^{2}} \right)}=\frac{1}{2\sqrt{6}}a=\frac{\sqrt{6}}{12}a. \\ 
\end{aligned}\]

\section{Volume of tetrahedron}\label{Sec29.2}


Given the lengths of six edges of tetrahedron $ABCD$, how to calculate the volume of the tetrahedron? Although the relationship between the volume of the tetrahedron and the radius of the inscribed sphere is not difficult to obtain, if the radius of the inscribed sphere cannot be calculated, then the volume of the tetrahedron cannot be calculated. Now that I have found the radius of the inscribed sphere of the tetrahedron (see theorem \ref{thm:Thm29.1.1}), the volume of the tetrahedron can be calculated.


\begin{theorem}{Volume formula of tetrahedron, Daiyuan Zhang}{Thm29.2.1}\label{Thm29.2.1} 
	Let the lengths of the six edges of the tetrahedron $ABCD$ be $AB$, $AC$, $AD$, $BC$, $CD$, $DB$, and the volume be $V$, then
	\[V=\frac{1}{6}\sqrt{{{t}_{1}}-{{t}_{2}}-{{t}_{3}}}.\]
	Where 
	\[{{t}_{1}}={{q}_{2}}{{\Delta }_{2}},\]
	\[{{q}_{2}}=\frac{1}{2}\left( A{{B}^{2}}C{{D}^{2}}+B{{C}^{2}}A{{D}^{2}}+C{{A}^{2}}B{{D}^{2}} \right),\]
	\[{{\Delta }_{2}}=\frac{1}{2}\left( A{{B}^{2}}+A{{C}^{2}}+A{{D}^{2}}+B{{C}^{2}}+C{{D}^{2}}+D{{B}^{2}} \right),\]
	\[{{t}_{2}}=\frac{1}{2}\left( A{{B}^{2}}C{{D}^{2}}\left( A{{B}^{2}}+C{{D}^{2}} \right)+B{{C}^{2}}A{{D}^{2}}\left( B{{C}^{2}}+A{{D}^{2}} \right)+C{{A}^{2}}B{{D}^{2}}\left( C{{A}^{2}}+B{{D}^{2}} \right) \right),\]
	\[{{t}_{3}}=\frac{1}{4}\left( A{{B}^{2}}B{{C}^{2}}C{{A}^{2}}+B{{C}^{2}}C{{D}^{2}}D{{B}^{2}}+C{{D}^{2}}D{{A}^{2}}A{{C}^{2}}+D{{A}^{2}}A{{B}^{2}}B{{D}^{2}} \right).\]
\end{theorem}

\begin{proof}
	Assuming that the radius of the inscribed sphere of the tetrahedron is $r$, it is obtained according to theorem \ref{thm:Thm29.1.1}:
	\[V=\frac{1}{3}\left( {{S}^{A}}+{{S}^{B}}+{{S}^{C}}+{{S}^{D}} \right)r=\frac{1}{3}Sr=\frac{1}{3}S\cdot \frac{\sqrt{{{t}_{1}}-{{t}_{2}}-{{t}_{3}}}}{2S}=\frac{1}{6}\sqrt{{{t}_{1}}-{{t}_{2}}-{{t}_{3}}}.\]
\end{proof}
\hfill $\square$\par




%

%

%

%



%

%

%

%

%
\chapter{Circumscribed sphere radius of tetrahedron}\label{Ch30}
\thispagestyle{empty}


This chapter mainly studies the calculation method of circumscribed sphere radius of an arbitrary tetrahedron.


Although for some special tetrahedron, such as regular triangular pyramid or regular tetrahedron, the radius of its circumscribed sphere can be obtained by a simple method (Euclidean geometry method), for a general tetrahedron, assuming the lengths of the six edges of the tetrahedron are given, how to calculate the circumscribed sphere radius of the tetrahedron? This chapter will answer the question.


In history, A.L. Crelle once studied the radius of circumscribed sphere of an arbitrary tetrahedron and obtained a formula for the radius of the circumscribed sphere, which means that if the volume of the tetrahedron and the lengths of six edges are given, the radius of the circumscribed sphere of the tetrahedron can be obtained. Because Crelle's formula needs to know the volume of the tetrahedron, it may not be considered that Crelle's formula completely solves the problem of calculating the radius of the circumscribed sphere of the tetrahedron.


In this chapter, the method of Intercenter Geometry proposed by me is successfully used to solve this problem. In this chapter, I give some formulas of the radius of the circumscribed sphere  of tetrahedron. The formula of the radius of the circumscribed sphere of tetrahedron can have many different forms, which can be said to be colorful.

\section{Colorful formula of radius of circumscribed sphere of tetrahedron}\label{Subsec30.1}


This section studies the radius of the circumscribed sphere of a tetrahedron. In this chapter, I give several different forms of formulas for the radius of circumscribed sphere of tetrahedron. 


\begin{theorem}{Basic formula of radius of circumscribed sphere of tetrahedron, Daiyuan Zhang}{Thm30.1.1}\label{Thm30.1.1} 
	Let the lengths of the six edges of tetrahedron $ABCD$ be $AB$, $AC$, $AD$, $BC$, $CD$, $DB$, and the radius of the circumscribed ball of tetrahedron $ABCD$ be $R$, then:
	\begin{align*}
		{{R}^{2}}& =\left( 1-\beta _{A}^{Q} \right)\left( \beta _{B}^{Q}A{{B}^{2}}+\beta _{C}^{Q}A{{C}^{2}}+\beta _{D}^{Q}A{{D}^{2}} \right) \\ 
		& -\beta _{B}^{Q}\beta _{C}^{Q}B{{C}^{2}}-\beta _{C}^{Q}\beta _{D}^{Q}C{{D}^{2}}-\beta _{D}^{Q}\beta _{B}^{Q}D{{B}^{2}},  
	\end{align*}	
	or 
	\begin{align*}
		{{R}^{2}}& =\left( 1-\beta _{B}^{Q} \right)\left( \beta _{C}^{Q}B{{C}^{2}}+\beta _{D}^{Q}B{{D}^{2}}+\beta _{A}^{Q}B{{A}^{2}} \right) \\ 
		& -\beta _{C}^{Q}\beta _{D}^{Q}C{{D}^{2}}-\beta _{D}^{Q}\beta _{A}^{Q}D{{A}^{2}}-\beta _{A}^{Q}\beta _{C}^{Q}A{{C}^{2}},  
	\end{align*}
	or 
	\begin{align*}
		{{R}^{2}}& =\left( 1-\beta _{C}^{Q} \right)\left( \beta _{D}^{Q}C{{D}^{2}}+\beta _{A}^{Q}C{{A}^{2}}+\beta _{B}^{Q}C{{B}^{2}} \right) \\ 
		& -\beta _{D}^{Q}\beta _{A}^{Q}D{{A}^{2}}-\beta _{A}^{Q}\beta _{B}^{Q}A{{B}^{2}}-\beta _{B}^{Q}\beta _{D}^{Q}B{{D}^{2}},  
	\end{align*}
	or 
	\begin{align*}
		{{R}^{2}}& =\left( 1-\beta _{D}^{Q} \right)\left( \beta _{A}^{Q}D{{A}^{2}}+\beta _{B}^{Q}D{{B}^{2}}+\beta _{C}^{Q}D{{C}^{2}} \right) \\ 
		& -\beta _{A}^{Q}\beta _{B}^{Q}A{{B}^{2}}-\beta _{B}^{Q}\beta _{C}^{Q}B{{C}^{2}}-\beta _{C}^{Q}\beta _{A}^{Q}C{{A}^{2}}.  
	\end{align*}	
	
	
	And the frame components of the circumcenter is given by theorem \ref{thm:Thm25.3.1}.
\end{theorem}

\begin{proof}
	
	The radius of the circumscribed sphere of a tetrahedron is the distance between the circumcenter and a vertexs. Now, let's calculate the distance between the radius of the circumscribed sphere of tetrahedron $ABCD$ and vertex $A$. According to theorem \ref{thm:Thm24.2.1}, the radius of the circumscribed sphere can be obtained by the following formula:
	
	\[\begin{aligned}
		{{R}^{2}}& =Q{{A}^{2}}=-\beta _{A}^{QA}\beta _{B}^{QA}A{{B}^{2}}-\beta _{A}^{QA}\beta _{C}^{QA}A{{C}^{2}}-\beta _{A}^{QA}\beta _{D}^{QA}A{{D}^{2}} \\ 
		& -\beta _{B}^{QA}\beta _{C}^{QA}B{{C}^{2}}-\beta _{C}^{QA}\beta _{D}^{QA}C{{D}^{2}}-\beta _{D}^{QA}\beta _{B}^{QA}D{{B}^{2}},  
	\end{aligned}\]
	where 
	\[\beta _{A}^{QA}=\beta _{A}^{A}-\beta _{A}^{Q},\]	
	\[\beta _{B}^{QA}=\beta _{B}^{A}-\beta _{B}^{Q},\]	
	\[\beta _{C}^{QA}=\beta _{C}^{A}-\beta _{C}^{Q},\]	
	\[\beta _{D}^{QA}=\beta _{D}^{A}-\beta _{D}^{Q}.\]	
	
	
	Here $\beta _{A}^{Q}$, $\beta _{B}^{Q}$, $\beta _{C}^{Q}$, $\beta _{D}^{Q}$ are the frame components of the circumcenter $Q$ in the tetrahedral frame $\left( O;A,B,C,D \right)$ respectively; and $\beta _{A}^{A}$, $\beta _{B}^{A}$, $\beta _{C}^{A}$, $\beta _{D}^{A}$ are the frame components of vertex $A$ in tetrahedral frame $\left( O;A,B,C,D \right)$ respectively.
	
	
	According to theorem \ref{thm:Thm20.4.1}, we get:
	
	\[\beta _{A}^{A}=1,\ \beta _{B}^{A}=0,\ \beta _{C}^{A}=0,\ \beta _{D}^{A}=0.\]
	
	Therefore 
	\[\beta _{A}^{QA}=\beta _{A}^{A}-\beta _{A}^{Q}=1-\beta _{A}^{Q},\]	
	\[\beta _{B}^{QA}=\beta _{B}^{A}-\beta _{B}^{Q}=-\beta _{B}^{Q},\]	
	\[\beta _{C}^{QA}=\beta _{C}^{A}-\beta _{C}^{Q}=-\beta _{C}^{Q},\]	
	\[\beta _{D}^{QA}=\beta _{D}^{A}-\beta _{D}^{Q}=-\beta _{D}^{Q}.\]	
	
	Therefore 
	\[\begin{aligned}
		{{R}^{2}}& =\left( 1-\beta _{A}^{Q} \right)\left( \beta _{B}^{Q}A{{B}^{2}}+\beta _{C}^{Q}A{{C}^{2}}+\beta _{D}^{Q}A{{D}^{2}} \right) \\ 
		& -\beta _{B}^{Q}\beta _{C}^{Q}B{{C}^{2}}-\beta _{C}^{Q}\beta _{D}^{Q}C{{D}^{2}}-\beta _{D}^{Q}\beta _{B}^{Q}D{{B}^{2}}.  
	\end{aligned}\]	
	
	
	If we calculate the distance between the radius of the circumscribed sphere of tetrahedron $ABCD$ and the vertices $B$, $C$, $D$ respectively, we can get the other formulas.
	
	
	Of course, the above theorem can also be obtained directly by using theorem \ref{thm:Thm24.1.2} and $A{{Q}^{2}}=B{{Q}^{2}}=C{{Q}^{2}}=D{{Q}^{2}}={{R}^{2}}$.
	
\end{proof}
\hfill $\square$\par

From the above theorem, we can directly see that the square of circumscribed sphere radius is the rational formula of the lengths of edges.


The following theorem can produce many different forms of the radius of circumscribed sphere.


\begin{theorem}{General formula of radius of circumscribed sphere of tetrahedron, Daiyuan Zhang}{Thm30.1.2}\label{Thm30.1.2} 
	Given a tetrahedron $ABCD$, let the point $Q$ be the circumscribed sphere center (circumcenter) of the tetrahedron, the point $P$ be the IC-T, $R$ be the radius of the circumscribed sphere of tetrahedron $ABCD$, and the lengths of the six edges of tetrahedron $ABCD$ be $AB$, $AC$, $AD$, $BC$, $CD$, $DB$ respectively, then
	\[\begin{aligned}
		{{R}^{2}}& =\left( \beta _{A}^{P}\beta _{B}^{Q}+\beta _{A}^{Q}\beta _{B}^{P}-\beta _{A}^{Q}\beta _{B}^{Q} \right)A{{B}^{2}}+\left( \beta _{A}^{P}\beta _{C}^{Q}+\beta _{A}^{Q}\beta _{C}^{P}-\beta _{A}^{Q}\beta _{C}^{Q} \right)A{{C}^{2}} \\ 
		& +\left( \beta _{A}^{P}\beta _{D}^{Q}+\beta _{A}^{Q}\beta _{D}^{P}-\beta _{A}^{Q}\beta _{D}^{Q} \right)A{{D}^{2}}+\left( \beta _{B}^{P}\beta _{C}^{Q}+\beta _{B}^{Q}\beta _{C}^{P}-\beta _{B}^{Q}\beta _{C}^{Q} \right)B{{C}^{2}} \\ 
		& +\left( \beta _{C}^{P}\beta _{D}^{Q}+\beta _{C}^{Q}\beta _{D}^{P}-\beta _{C}^{Q}\beta _{D}^{Q} \right)C{{D}^{2}}+\left( \beta _{D}^{Q}\beta _{B}^{P}+\beta _{D}^{P}\beta _{B}^{Q}-\beta _{D}^{Q}\beta _{B}^{Q} \right)D{{B}^{2}}.  
	\end{aligned}\] 
	Or 
	\[{{R}^{2}}=\frac{1}{{{U}^{2}}}\left( \begin{aligned}
		& \left( \beta _{B}^{P}{{U}_{A}}+\beta _{A}^{P}{{U}_{B}}-{{U}_{A}}{{U}_{B}} \right)A{{B}^{2}}+\left( \beta _{C}^{P}{{U}_{A}}+\beta _{A}^{P}{{U}_{C}}-{{U}_{A}}{{U}_{C}} \right)A{{C}^{2}} \\ 
		& +\left( \beta _{D}^{P}{{U}_{A}}+\beta _{A}^{P}{{U}_{D}}-{{U}_{A}}{{U}_{D}} \right)A{{D}^{2}}+\left( \beta _{C}^{P}{{U}_{B}}+\beta _{B}^{P}{{U}_{C}}-{{U}_{B}}{{U}_{C}} \right)B{{C}^{2}} \\ 
		& +\left( \beta _{D}^{P}{{U}_{C}}+\beta _{C}^{P}{{U}_{D}}-{{U}_{C}}{{U}_{D}} \right)C{{D}^{2}}+\left( \beta _{D}^{P}{{U}_{B}}+\beta _{B}^{P}{{U}_{D}}-{{U}_{D}}{{U}_{B}} \right)D{{B}^{2}}  
	\end{aligned} \right).\] 
	Where, the frame components of the circumscribed sphere center (circumcenter) are given by theorem \ref{thm:Thm25.3.1}.
\end{theorem}

\begin{proof}
	From theorem \ref{thm:Thm24.1.4}, we have
	\[Q{{P}^{2}}={{R}^{2}}-\left( \begin{aligned}
		& \beta _{A}^{P}\beta _{B}^{P}A{{B}^{2}}+\beta _{A}^{P}\beta _{C}^{P}A{{C}^{2}}+\beta _{A}^{P}\beta _{D}^{P}A{{D}^{2}} \\ 
		& +\beta _{B}^{P}\beta _{C}^{P}B{{C}^{2}}+\beta _{C}^{P}\beta _{D}^{P}C{{D}^{2}}+\beta _{D}^{P}\beta _{B}^{P}D{{B}^{2}}  
	\end{aligned} \right).\]
	
	According to theorem \ref{thm:Thm24.2.1}, for a given tetrahedron $ABCD$, $Q\in {{\pi }_{ABCD}}$, $P\in {{\pi }_{ABCD}}$, so
	\[\begin{aligned}
		Q{{P}^{2}}& =-\beta _{A}^{QP}\beta _{B}^{QP}A{{B}^{2}}-\beta _{A}^{QP}\beta _{C}^{QP}A{{C}^{2}}-\beta _{A}^{QP}\beta _{D}^{QP}A{{D}^{2}} \\ 
		& -\beta _{B}^{QP}\beta _{C}^{QP}B{{C}^{2}}-\beta _{C}^{QP}\beta _{D}^{QP}C{{D}^{2}}-\beta _{D}^{QP}\beta _{B}^{QP}D{{B}^{2}},  
	\end{aligned}\]
	\[\beta _{A}^{QP}+\beta _{B}^{QP}+\beta _{C}^{QP}+\beta _{D}^{QP}=0.\]	
	
	Where 
	\[\beta _{A}^{QP}=\beta _{A}^{P}-\beta _{A}^{Q},\]	
	\[\beta _{B}^{QP}=\beta _{B}^{P}-\beta _{B}^{Q},\]	
	\[\beta _{C}^{QP}=\beta _{C}^{P}-\beta _{C}^{Q},\]	
	\[\beta _{D}^{QP}=\beta _{D}^{P}-\beta _{D}^{Q}.\]	
	Where $\beta _{A}^{Q}$, $\beta _{B}^{Q}$, $\beta _{C}^{Q}$, $\beta _{D}^{Q}$ are the frame components of IC-T $Q$ in tetrahedral frame $\left( O;A,B,C,D \right)$ respectively; $\beta _{A}^{P}$, $\beta _{B}^{P}$, $\beta _{C}^{P}$, $\beta _{D}^{P}$ are the frame components of IC-T $P$ in tetrahedral frame $\left( O;A,B,C,D \right)$ respectively.
	
	Therefore 
	\[\begin{aligned}
		{{R}^{2}}& =Q{{P}^{2}}+\left( \begin{aligned}
			& \beta _{A}^{P}\beta _{B}^{P}A{{B}^{2}}+\beta _{A}^{P}\beta _{C}^{P}A{{C}^{2}}+\beta _{A}^{P}\beta _{D}^{P}A{{D}^{2}} \\ 
			& +\beta _{B}^{P}\beta _{C}^{P}B{{C}^{2}}+\beta _{C}^{P}\beta _{D}^{P}C{{D}^{2}}+\beta _{D}^{P}\beta _{B}^{P}D{{B}^{2}}  
		\end{aligned} \right) \\ 
		& =-\beta _{A}^{QP}\beta _{B}^{QP}A{{B}^{2}}-\beta _{A}^{QP}\beta _{C}^{QP}A{{C}^{2}}-\beta _{A}^{QP}\beta _{D}^{QP}A{{D}^{2}} \\ 
		& -\beta _{B}^{QP}\beta _{C}^{QP}B{{C}^{2}}-\beta _{C}^{QP}\beta _{D}^{QP}C{{D}^{2}}-\beta _{D}^{QP}\beta _{B}^{QP}D{{B}^{2}} \\ 
		& +\beta _{A}^{P}\beta _{B}^{P}A{{B}^{2}}+\beta _{A}^{P}\beta _{C}^{P}A{{C}^{2}}+\beta _{A}^{P}\beta _{D}^{P}A{{D}^{2}} \\ 
		& +\beta _{B}^{P}\beta _{C}^{P}B{{C}^{2}}+\beta _{C}^{P}\beta _{D}^{P}C{{D}^{2}}+\beta _{D}^{P}\beta _{B}^{P}D{{B}^{2}},  
	\end{aligned}\]
	i.e. 
	\[\begin{aligned}
		{{R}^{2}}& =\left( \beta _{A}^{P}\beta _{B}^{P}-\beta _{A}^{QP}\beta _{B}^{QP} \right)A{{B}^{2}}+\left( \beta _{A}^{P}\beta _{C}^{P}-\beta _{A}^{QP}\beta _{C}^{QP} \right)A{{C}^{2}}+\left( \beta _{A}^{P}\beta _{D}^{P}-\beta _{A}^{QP}\beta _{D}^{QP} \right)A{{D}^{2}} \\ 
		& +\left( \beta _{B}^{P}\beta _{C}^{P}-\beta _{B}^{QP}\beta _{C}^{QP} \right)B{{C}^{2}}+\left( \beta _{C}^{P}\beta _{D}^{P}-\beta _{C}^{QP}\beta _{D}^{QP} \right)C{{D}^{2}}+\left( \beta _{D}^{P}\beta _{B}^{P}-\beta _{D}^{QP}\beta _{B}^{QP} \right)D{{B}^{2}}.  
	\end{aligned}\]
	And 
	\[\begin{aligned}
		\beta _{A}^{P}\beta _{B}^{P}-\beta _{A}^{QP}\beta _{B}^{QP}& =\beta _{A}^{P}\beta _{B}^{P}-\left( \beta _{A}^{P}-\beta _{A}^{Q} \right)\left( \beta _{B}^{P}-\beta _{B}^{Q} \right) \\ 
		& =\beta _{A}^{P}\beta _{B}^{P}-\left( \beta _{A}^{P}\beta _{B}^{P}-\beta _{A}^{P}\beta _{B}^{Q}-\beta _{A}^{Q}\beta _{B}^{P}+\beta _{A}^{Q}\beta _{B}^{Q} \right) \\ 
		& =\beta _{A}^{P}\beta _{B}^{Q}+\beta _{A}^{Q}\beta _{B}^{P}-\beta _{A}^{Q}\beta _{B}^{Q},  
	\end{aligned}\]
	
	Similarly: 
	\[\beta _{A}^{P}\beta _{C}^{P}-\beta _{A}^{QP}\beta _{C}^{QP}=\beta _{A}^{P}\beta _{C}^{Q}+\beta _{A}^{Q}\beta _{C}^{P}-\beta _{A}^{Q}\beta _{C}^{Q},\]
	\[\beta _{A}^{P}\beta _{D}^{P}-\beta _{A}^{QP}\beta _{D}^{QP}=\beta _{A}^{P}\beta _{D}^{Q}+\beta _{A}^{Q}\beta _{D}^{P}-\beta _{A}^{Q}\beta _{D}^{Q},\]
	\[\beta _{B}^{P}\beta _{C}^{P}-\beta _{B}^{QP}\beta _{C}^{QP}=\beta _{B}^{P}\beta _{C}^{Q}+\beta _{B}^{Q}\beta _{C}^{P}-\beta _{B}^{Q}\beta _{C}^{Q},\]
	\[\beta _{C}^{P}\beta _{D}^{P}-\beta _{C}^{QP}\beta _{D}^{QP}=\beta _{C}^{P}\beta _{D}^{Q}+\beta _{C}^{Q}\beta _{D}^{P}-\beta _{C}^{Q}\beta _{D}^{Q},\]
	\[\beta _{D}^{P}\beta _{B}^{P}-\beta _{D}^{QP}\beta _{B}^{QP}=\beta _{D}^{Q}\beta _{B}^{P}+\beta _{D}^{P}\beta _{B}^{Q}-\beta _{D}^{Q}\beta _{B}^{Q}.\]
	
	Therefore 
	\[\begin{aligned}
		{{R}^{2}}& =\left( \beta _{A}^{P}\beta _{B}^{Q}+\beta _{A}^{Q}\beta _{B}^{P}-\beta _{A}^{Q}\beta _{B}^{Q} \right)A{{B}^{2}}+\left( \beta _{A}^{P}\beta _{C}^{Q}+\beta _{A}^{Q}\beta _{C}^{P}-\beta _{A}^{Q}\beta _{C}^{Q} \right)A{{C}^{2}} \\ 
		& +\left( \beta _{A}^{P}\beta _{D}^{Q}+\beta _{A}^{Q}\beta _{D}^{P}-\beta _{A}^{Q}\beta _{D}^{Q} \right)A{{D}^{2}}+\left( \beta _{B}^{P}\beta _{C}^{Q}+\beta _{B}^{Q}\beta _{C}^{P}-\beta _{B}^{Q}\beta _{C}^{Q} \right)B{{C}^{2}} \\ 
		& +\left( \beta _{C}^{P}\beta _{D}^{Q}+\beta _{C}^{Q}\beta _{D}^{P}-\beta _{C}^{Q}\beta _{D}^{Q} \right)C{{D}^{2}}+\left( \beta _{D}^{Q}\beta _{B}^{P}+\beta _{D}^{P}\beta _{B}^{Q}-\beta _{D}^{Q}\beta _{B}^{Q} \right)D{{B}^{2}}.  
	\end{aligned}\] 
	
	From theorem \ref{thm:Thm25.3.1}, we have
	\[{{R}^{2}}=\frac{1}{{{U}^{2}}}\left( \begin{aligned}
		& \left( \beta _{B}^{P}{{U}_{A}}+\beta _{A}^{P}{{U}_{B}}-{{U}_{A}}{{U}_{B}} \right)A{{B}^{2}}+\left( \beta _{C}^{P}{{U}_{A}}+\beta _{A}^{P}{{U}_{C}}-{{U}_{A}}{{U}_{C}} \right)A{{C}^{2}} \\ 
		& +\left( \beta _{D}^{P}{{U}_{A}}+\beta _{A}^{P}{{U}_{D}}-{{U}_{A}}{{U}_{D}} \right)A{{D}^{2}}+\left( \beta _{C}^{P}{{U}_{B}}+\beta _{B}^{P}{{U}_{C}}-{{U}_{B}}{{U}_{C}} \right)B{{C}^{2}} \\ 
		& +\left( \beta _{D}^{P}{{U}_{C}}+\beta _{C}^{P}{{U}_{D}}-{{U}_{C}}{{U}_{D}} \right)C{{D}^{2}}+\left( \beta _{D}^{P}{{U}_{B}}+\beta _{B}^{P}{{U}_{D}}-{{U}_{D}}{{U}_{B}} \right)D{{B}^{2}}  
	\end{aligned} \right).\] 
\end{proof}
\hfill $\square$\par
The above formula is a general formula for finding the radius of the circumscribed sphere of the tetrahedron. Different forms of formulas can be obtained by choosing different $P$. The following theorems give different forms of circumscribed sphere radius of tetrahedron.

If we choose point $P$ as the circumcenter, we get the following theorem.


\begin{theorem}{Formula 1 of Circumscribed sphere radius, Daiyuan Zhang}{Thm30.1.3}\label{Thm30.1.3} 
	Given a tetrahedron $ABCD$, let the point $Q$ be the circumscribed sphere center (circumcenter) of the tetrahedron, the point $P$ be the IC-T, $R$ be the radius of the circumscribed sphere  of tetrahedron $ABCD$, and the lengths of the six edges of tetrahedron $ABCD$ be $AB$, $AC$, $AD$, $BC$, $CD$, $DB$ respectively, then
	\[\begin{aligned}
		{{R}^{2}}& =\beta _{A}^{Q}\beta _{B}^{Q}A{{B}^{2}}+\beta _{A}^{Q}\beta _{C}^{Q}A{{C}^{2}}+\beta _{A}^{Q}\beta _{D}^{Q}A{{D}^{2}} \\ 
		& +\beta _{B}^{Q}\beta _{C}^{Q}B{{C}^{2}}+\beta _{C}^{Q}\beta _{D}^{Q}C{{D}^{2}}+\beta _{D}^{Q}\beta _{B}^{Q}D{{B}^{2}}.  
	\end{aligned}\]
	Or 
	\[{{R}^{2}}=\frac{1}{{{U}^{2}}}\left( \begin{aligned}
		& {{U}_{A}}{{U}_{B}}A{{B}^{2}}+{{U}_{A}}{{U}_{C}}A{{C}^{2}}+{{U}_{A}}{{U}_{D}}A{{D}^{2}} \\ 
		& +{{U}_{B}}{{U}_{C}}B{{C}^{2}}+{{U}_{C}}{{U}_{D}}C{{D}^{2}}+{{U}_{D}}{{U}_{B}}D{{B}^{2}}  
	\end{aligned} \right).\]
	
	Where, the frame components of the circumscribed sphere center (circumcenter) are given by theorem \ref{thm:Thm25.3.1}.
\end{theorem}

\begin{proof}
	In theorem \ref{thm:Thm30.1.2}, let $P=Q$, select the point $P$ as the circumcenter, then	
	\[\begin{aligned}
		{{R}^{2}}& =\beta _{A}^{Q}\beta _{B}^{Q}A{{B}^{2}}+\beta _{A}^{Q}\beta _{C}^{Q}A{{C}^{2}}+\beta _{A}^{Q}\beta _{D}^{Q}A{{D}^{2}} \\ 
		& +\beta _{B}^{Q}\beta _{C}^{Q}B{{C}^{2}}+\beta _{C}^{Q}\beta _{D}^{Q}C{{D}^{2}}+\beta _{D}^{Q}\beta _{B}^{Q}D{{B}^{2}}.  
	\end{aligned}\]
	
	According to theorem \ref{thm:Thm25.3.1}, the following result is obtained directly.
	\[{{R}^{2}}=\frac{1}{{{U}^{2}}}\left( \begin{aligned}
		& {{U}_{A}}{{U}_{B}}A{{B}^{2}}+{{U}_{A}}{{U}_{C}}A{{C}^{2}}+{{U}_{A}}{{U}_{D}}A{{D}^{2}} \\ 
		& +{{U}_{B}}{{U}_{C}}B{{C}^{2}}+{{U}_{C}}{{U}_{D}}C{{D}^{2}}+{{U}_{D}}{{U}_{B}}D{{B}^{2}}  
	\end{aligned} \right).\]
\end{proof}
\hfill $\square$\par
If we choose point $P$ as the centroid, we can get the following theorem.


\begin{theorem}{Formula 2 of Circumscribed sphere radius, Daiyuan Zhang, Daiyuan Zhang}{Thm30.1.4}\label{Thm30.1.4} 
	Given a tetrahedron $ABCD$, let the point $Q$ be the circumscribed sphere center (circumcenter) of the tetrahedron, $R$ be the radius of the circumscribed sphere  of tetrahedron $ABCD$, and the lengths of the six edges of tetrahedron $ABCD$ be $AB$, $AC$, $AD$, $BC$, $CD$, $DB$ respectively, then
	\[{{R}^{2}}=\frac{1}{4}\left( \begin{aligned}
		& \left( \beta _{A}^{Q}+\beta _{B}^{Q}-4\beta _{A}^{Q}\beta _{B}^{Q} \right)A{{B}^{2}}+\left( \beta _{A}^{Q}+\beta _{C}^{Q}-4\beta _{A}^{Q}\beta _{C}^{Q} \right)A{{C}^{2}} \\ 
		& +\left( \beta _{A}^{Q}+\beta _{D}^{Q}-4\beta _{A}^{Q}\beta _{D}^{Q} \right)A{{D}^{2}}+\left( \beta _{B}^{Q}+\beta _{C}^{Q}-4\beta _{B}^{Q}\beta _{C}^{Q} \right)B{{C}^{2}} \\ 
		& +\left( \beta _{C}^{Q}+\beta _{D}^{Q}-4\beta _{C}^{Q}\beta _{D}^{Q} \right)C{{D}^{2}}+\left( \beta _{D}^{Q}+\beta _{B}^{Q}-4\beta _{D}^{Q}\beta _{B}^{Q} \right)D{{B}^{2}}  
	\end{aligned} \right).\]
	Or 
	\[{{R}^{2}}=\frac{1}{4{{U}^{2}}}\left( \begin{aligned}
		& \left( U\left( {{U}_{A}}+{{U}_{B}} \right)-4{{U}_{A}}{{U}_{B}} \right)A{{B}^{2}}+\left( U\left( {{U}_{A}}+{{U}_{C}} \right)-4{{U}_{A}}{{U}_{C}} \right)A{{C}^{2}} \\ 
		& +\left( U\left( {{U}_{A}}+{{U}_{D}} \right)-4{{U}_{A}}{{U}_{D}} \right)A{{D}^{2}}+\left( U\left( {{U}_{B}}+{{U}_{C}} \right)-4{{U}_{B}}{{U}_{C}} \right)B{{C}^{2}} \\ 
		& +\left( U\left( {{U}_{C}}+{{U}_{D}} \right)-4{{U}_{C}}{{U}_{D}} \right)C{{D}^{2}}+\left( U\left( {{U}_{D}}+{{U}_{B}} \right)-4{{U}_{D}}{{U}_{B}} \right)D{{B}^{2}}  
	\end{aligned} \right).\]
	
	Where, the frame components of the circumscribed sphere center (circumcenter) are given by theorem \ref{thm:Thm25.3.1}.
\end{theorem}

\begin{proof}
	If we choose point $P$ as the centroid, then $\beta _{A}^{P}=\beta _{B}^{P}=\beta _{C}^{P}=\beta _{D}^{P}={1}/{4}\;$, we can get the following theorem by theorem \ref{thm:Thm30.1.2} directly.
\end{proof}
\hfill $\square$\par
Another form of formula can also be obtained as follows.


\begin{theorem}{Formula 3 of Circumscribed sphere radius, Daiyuan Zhang}{Thm30.1.5}\label{Thm30.1.5} 
	Given a tetrahedron $ABCD$, let the point $Q$ be the circumscribed sphere center (circumcenter) of the tetrahedron, $R$ be the radius of the circumscribed sphere  of tetrahedron $ABCD$, and the lengths of the six edges of tetrahedron $ABCD$ be $AB$, $AC$, $AD$, $BC$, $CD$, $DB$ respectively, then
	\[{{R}^{2}}=\frac{1}{8}\left( \begin{aligned}
		& \left( \beta _{A}^{Q}+\beta _{B}^{Q} \right)A{{B}^{2}}+\left( \beta _{A}^{Q}+\beta _{C}^{Q} \right)A{{C}^{2}}+\left( \beta _{A}^{Q}+\beta _{D}^{Q} \right)A{{D}^{2}} \\ 
		& +\left( \beta _{B}^{Q}+\beta _{C}^{Q} \right)B{{C}^{2}}+\left( \beta _{C}^{Q}+\beta _{D}^{Q} \right)C{{D}^{2}}+\left( \beta _{D}^{Q}+\beta _{B}^{Q} \right)D{{B}^{2}}  
	\end{aligned} \right).\] 
	Or
	\[{{R}^{2}}=\frac{1}{8U}\left( \begin{aligned}
		& \left( {{U}_{A}}+{{U}_{B}} \right)A{{B}^{2}}+\left( {{U}_{A}}+{{U}_{C}} \right)A{{C}^{2}}+\left( {{U}_{A}}+{{U}_{D}} \right)A{{D}^{2}} \\ 
		& +\left( {{U}_{B}}+{{U}_{C}} \right)B{{C}^{2}}+\left( {{U}_{C}}+{{U}_{D}} \right)C{{D}^{2}}+\left( {{U}_{D}}+{{U}_{B}} \right)D{{B}^{2}}  
	\end{aligned} \right).\] 
	
	Where, the frame components of the circumscribed sphere center (circumcenter) are given by theorem \ref{thm:Thm25.3.1}.
\end{theorem}

\begin{proof}
	Add both sides of the formulas in theorem \ref{thm:Thm30.1.2} and theorem \ref{thm:Thm30.1.3} to obtain:
	\[{{R}^{2}}=\frac{1}{2}\left( \begin{aligned}
		& \left( \beta _{A}^{P}\beta _{B}^{Q}+\beta _{A}^{Q}\beta _{B}^{P} \right)A{{B}^{2}}+\left( \beta _{A}^{P}\beta _{C}^{Q}+\beta _{A}^{Q}\beta _{C}^{P} \right)A{{C}^{2}}+\left( \beta _{A}^{P}\beta _{D}^{Q}+\beta _{A}^{Q}\beta _{D}^{P} \right)A{{D}^{2}} \\ 
		& +\left( \beta _{B}^{P}\beta _{C}^{Q}+\beta _{B}^{Q}\beta _{C}^{P} \right)B{{C}^{2}}+\left( \beta _{C}^{P}\beta _{D}^{Q}+\beta _{C}^{Q}\beta _{D}^{P} \right)C{{D}^{2}}+\left( \beta _{D}^{Q}\beta _{B}^{P}+\beta _{D}^{P}\beta _{B}^{Q} \right)D{{B}^{2}}  
	\end{aligned} \right).\] 
	
	Let $P=G$, that is, select point $P$ as the centroid, then
	\[{{R}^{2}}=\frac{1}{8}\left( \begin{aligned}
		& \left( \beta _{A}^{Q}+\beta _{B}^{Q} \right)A{{B}^{2}}+\left( \beta _{A}^{Q}+\beta _{C}^{Q} \right)A{{C}^{2}}+\left( \beta _{A}^{Q}+\beta _{D}^{Q} \right)A{{D}^{2}} \\ 
		& +\left( \beta _{B}^{Q}+\beta _{C}^{Q} \right)B{{C}^{2}}+\left( \beta _{C}^{Q}+\beta _{D}^{Q} \right)C{{D}^{2}}+\left( \beta _{D}^{Q}+\beta _{B}^{Q} \right)D{{B}^{2}}  
	\end{aligned} \right).\] 
	
	According to theorem \ref{thm:Thm25.3.1}, the following result is obtained directly:
	\[{{R}^{2}}=\frac{1}{8U}\left( \begin{aligned}
		& \left( {{U}_{A}}+{{U}_{B}} \right)A{{B}^{2}}+\left( {{U}_{A}}+{{U}_{C}} \right)A{{C}^{2}}+\left( {{U}_{A}}+{{U}_{D}} \right)A{{D}^{2}} \\ 
		& +\left( {{U}_{B}}+{{U}_{C}} \right)B{{C}^{2}}+\left( {{U}_{C}}+{{U}_{D}} \right)C{{D}^{2}}+\left( {{U}_{D}}+{{U}_{B}} \right)D{{B}^{2}}  
	\end{aligned} \right).\] 
\end{proof}
\hfill $\square$\par

\begin{theorem}{Formula 4 of Circumscribed sphere radius, Daiyuan Zhang}{Thm30.1.6}\label{Thm30.1.6} 
	Given a tetrahedron $ABCD$, let $R$ be the radius of the circumscribed sphere  of tetrahedron $ABCD$, and the lengths of the six edges of tetrahedron $ABCD$ be $AB$, $AC$, $AD$, $BC$, $CD$, $DB$ respectively, then
	\[{{R}^{2}}=\frac{q\left( q-AB\cdot CD \right)\left( q-BC\cdot AD \right)\left( q-CA\cdot BD \right)}{{{t}_{1}}-{{t}_{2}}-{{t}_{3}}}.\]	
	Where 
	\[q=\frac{1}{2}\left( AB\cdot CD+BC\cdot AD+CA\cdot BD \right),\]
	\[{{t}_{1}}={{q}_{2}}{{\Delta }_{2}},\]
	\[{{q}_{2}}=\frac{1}{2}\left( A{{B}^{2}}C{{D}^{2}}+B{{C}^{2}}A{{D}^{2}}+C{{A}^{2}}B{{D}^{2}} \right),\]
	\[{{\Delta }_{2}}=\frac{1}{2}\left( A{{B}^{2}}+A{{C}^{2}}+A{{D}^{2}}+B{{C}^{2}}+C{{D}^{2}}+D{{B}^{2}} \right),\]
	\[{{t}_{2}}=\frac{1}{2}\left( A{{B}^{2}}C{{D}^{2}}\left( A{{B}^{2}}+C{{D}^{2}} \right)+B{{C}^{2}}A{{D}^{2}}\left( B{{C}^{2}}+A{{D}^{2}} \right)+C{{A}^{2}}B{{D}^{2}}\left( C{{A}^{2}}+B{{D}^{2}} \right) \right),\]
	\[{{t}_{3}}=\frac{1}{4}\left( A{{B}^{2}}B{{C}^{2}}C{{A}^{2}}+B{{C}^{2}}C{{D}^{2}}D{{B}^{2}}+C{{D}^{2}}D{{A}^{2}}A{{C}^{2}}+D{{A}^{2}}A{{B}^{2}}B{{D}^{2}} \right).\]
\end{theorem}

\begin{proof}
	Substitute the result of theorem \ref{thm:Thm25.3.1} into theorem \ref{thm:Thm30.1.5}, we have:
	\[{{R}^{2}}=-\frac{\left( \begin{aligned}
			& A{{B}^{4}}C{{D}^{4}}-2A{{B}^{2}}A{{D}^{2}}B{{C}^{2}}C{{D}^{2}}-2A{{B}^{2}}B{{D}^{2}}C{{A}^{2}}C{{D}^{2}} \\ 
			& +A{{D}^{4}}B{{C}^{4}}-2A{{D}^{2}}B{{C}^{2}}B{{D}^{2}}C{{A}^{2}}+B{{D}^{4}}C{{A}^{4}} \\ 
		\end{aligned} \right)}{16\left( {{t}_{1}}-{{t}_{2}}-{{t}_{3}} \right)},\]
	wehre 
	\[{{t}_{1}}={{q}_{2}}{{\Delta }_{2}},\]
	\[{{q}_{2}}=\frac{1}{2}\left( A{{B}^{2}}C{{D}^{2}}+B{{C}^{2}}A{{D}^{2}}+C{{A}^{2}}B{{D}^{2}} \right),\]
	\[{{\Delta }_{2}}=\frac{1}{2}\left( A{{B}^{2}}+A{{C}^{2}}+A{{D}^{2}}+B{{C}^{2}}+C{{D}^{2}}+D{{B}^{2}} \right),\]
	\[{{t}_{2}}=\frac{1}{2}\left( A{{B}^{2}}C{{D}^{2}}\left( A{{B}^{2}}+C{{D}^{2}} \right)+B{{C}^{2}}A{{D}^{2}}\left( B{{C}^{2}}+A{{D}^{2}} \right)+C{{A}^{2}}B{{D}^{2}}\left( C{{A}^{2}}+B{{D}^{2}} \right) \right),\]
	\[{{t}_{3}}=\frac{1}{4}\left( A{{B}^{2}}B{{C}^{2}}C{{A}^{2}}+B{{C}^{2}}C{{D}^{2}}D{{B}^{2}}+C{{D}^{2}}D{{A}^{2}}A{{C}^{2}}+D{{A}^{2}}A{{B}^{2}}B{{D}^{2}} \right).\]
	i.e. 
	\[{{R}^{2}}=-\frac{\left( \begin{aligned}
			& \left( AB\cdot CD-AD\cdot BC-BD\cdot CA \right)\left( AB\cdot CD-AD\cdot BC+BD\cdot CA \right) \\ 
			& \cdot \left( AD\cdot BC+AB\cdot CD+BD\cdot CA \right)\left( AD\cdot BC+AB\cdot CD-BD\cdot CA \right) \\ 
		\end{aligned} \right)}{16\left( {{t}_{1}}-{{t}_{2}}-{{t}_{3}} \right)}.\]
	
	So we get a beautiful formula for the radius of the circumscribed sphere:
	\[{{R}^{2}}=\frac{q\left( q-AB\cdot CD \right)\left( q-BC\cdot AD \right)\left( q-CA\cdot BD \right)}{{{t}_{1}}-{{t}_{2}}-{{t}_{3}}}.\]
	Where
	\[q=\frac{1}{2}\left( AB\cdot CD+BC\cdot AD+CA\cdot BD \right).\]	
\end{proof}
\hfill $\square$\par

\section{The relationship between the volume of tetrahedron and the radius of circumscribed sphere}\label{Sec30.2}


\begin{corollary}{Crelle's formula, A.L.Crelle}{Cor30.2.1}\label{Cor30.2.1} 
	Given a tetrahedron $ABCD$, let the radius of the circumscribed sphere of tetrahedron $ABCD$ be $R$, the volume of tetrahedron $ABCD$ be $V$, and the lengths of the six edges of tetrahedron $ABCD$ be $AB$, $AC$, $AD$, $BC$, $CD$, $DB$ respectively, then
	\[36{{V}^{2}}{{R}^{2}}=q\left( q-AB\cdot CD \right)\left( q-BC\cdot AD \right)\left( q-CA\cdot BD \right),\]
	where 
	\[q=\frac{1}{2}\left( AB\cdot CD+BC\cdot AD+CA\cdot BD \right).\]
\end{corollary}

\begin{proof}
	According to theorem \ref{thm:Thm30.1.6} and theorem \ref{thm:Thm29.2.1}, the following result is obtained:
	\[{{R}^{2}}=\frac{q\left( q-AB\cdot CD \right)\left( q-BC\cdot AD \right)\left( q-CA\cdot BD \right)}{{{t}_{1}}-{{t}_{2}}-{{t}_{3}}},\]	\[{{r}^{2}}=\frac{1}{4{{S}^{2}}}\left( {{t}_{1}}-{{t}_{2}}-{{t}_{3}} \right).\]
	
	Therefore 
	\[{{r}^{2}}{{R}^{2}}=\frac{q\left( q-AB\cdot CD \right)\left( q-BC\cdot AD \right)\left( q-CA\cdot BD \right)}{4{{S}^{2}}}.\]
	i.e.
	\[4{{S}^{2}}{{r}^{2}}{{R}^{2}}=q\left( q-AB\cdot CD \right)\left( q-BC\cdot AD \right)\left( q-CA\cdot BD \right).\]
	
	And 
	\[V=\frac{1}{3}\left( {{S}^{A}}+{{S}^{B}}+{{S}^{C}}+{{S}^{D}} \right)r=\frac{1}{3}Sr,\]
	i.e. 
	\[{{S}^{2}}{{r}^{2}}=9{{V}^{2}}.\]
	
	Therefore 
	\[36{{V}^{2}}{{R}^{2}}=q\left( q-AB\cdot CD \right)\left( q-BC\cdot AD \right)\left( q-CA\cdot BD \right).\]
\end{proof}
\hfill $\square$\par

Using Euclidean geometry, Crelle proved the above theorem by using auxiliary surfaces and auxiliary lines. Crelle's proof was very personalized and highly skilled, but lack of generality. Obviously, It was impossible for Crelle to know  theorem \ref{thm:Thm30.1.6} and theorem \ref{thm:Thm29.2.1} in his time, so he was forced to choose the method of highly skilled but lack of generality.

Using the Intercenter Geometry method that I proposed, theorem \ref{thm:Thm30.1.6} and theorem \ref{thm:Thm29.2.1} have been obtained. Using these two theorems, we can easily get Crelle's formula, so this book takes Crelle's formula as an corollary.

Theorem \ref{thm:Thm30.1.6} and theorem \ref{thm:Thm29.2.1} I proposed are the latest results in Intercenter Geometry. These two graceful and symmetric formulas are presented to the world for the first time, and I hope the readers will like them.

After the radius of the circumscribed sphere of the tetrahedron is obtained, the surface area and volume of the circumscribed sphere of the tetrahedron can be easily obtained.

\chapter{Necessary and sufficient conditions of IC-T}\label{Ch31}
\thispagestyle{empty}

In this chapter, I give the necessary and sufficient conditions of intersecting center of tetrahedron (abbreviated as IC-T), which is very important for some researches and applications.

\section{Calculation of the IRs by the frame components of the IC-T}\label{Sec31.1}
\begin{theorem}{The calculation of IRs by the frame components of the IC-T, Daiyuan Zhang}{Thm31.1.1}\label{Thm31.1.1} 
	Given a tetrahedron $ABCD$, point $P$ is the IC-T of tetrahedron $ABCD$, and the corresponding ICs-Fs of tetrahedron $ABCD$ are ${{P}_{A}}$, ${{P}_{B}}$, ${{P}_{C}}$, ${{P}_{D}}$ respectively; and ${{\beta }_{A}}$, ${{\beta }_{B}}$, ${{\beta }_{C}}$, ${{\beta }_{D}}$ are the frame components of point $P$ in tetrahedron $ABCD$ respectively, then the ICs-Fs have the following properties:
	
	\[\lambda _{CD}^{{{P}_{A}}}=\lambda _{CD}^{{{P}_{B}}}=\frac{\beta _{D}^{P}}{\beta _{C}^{P}},\ \lambda _{BD}^{{{P}_{A}}}=\lambda _{BD}^{{{P}_{C}}}=\frac{\beta _{D}^{P}}{\beta _{B}^{P}},\ \lambda _{BC}^{{{P}_{A}}}=\lambda _{BC}^{{{P}_{D}}}=\frac{\beta _{C}^{P}}{\beta _{B}^{P}},\]
	\[\lambda _{AD}^{{{P}_{B}}}=\lambda _{AD}^{{{P}_{C}}}=\frac{\beta _{D}^{P}}{\beta _{A}^{P}},\ \lambda _{AC}^{{{P}_{B}}}=\lambda _{AC}^{{{P}_{D}}}=\frac{\beta _{C}^{P}}{\beta _{A}^{P}},\ \lambda _{AB}^{{{P}_{C}}}=\lambda _{AB}^{{{P}_{D}}}=\frac{\beta _{B}^{P}}{\beta _{A}^{P}}.\]
\end{theorem}

\begin{proof}
	According to theorem \ref{thm:Thm6.4.1} and theorem \ref{thm:Thm20.2.4}, the following results are obtained:
	\[\lambda _{AB}^{{{P}_{C}}}=\frac{\alpha _{B}^{{{P}_{C}}}}{\alpha _{A}^{{{P}_{C}}}}=\frac{\beta _{B}^{P}}{\beta _{A}^{P}},\ \lambda _{AB}^{{{P}_{D}}}=\frac{\alpha _{B}^{{{P}_{D}}}}{\alpha _{A}^{{{P}_{D}}}}=\frac{\beta _{B}^{P}}{\beta _{A}^{P}},\]
	therefore 
	\[\lambda _{AB}^{{{P}_{C}}}=\lambda _{AB}^{{{P}_{D}}}=\frac{\beta _{B}^{P}}{\beta _{A}^{P}}.\]
	
	Similarly: 
	\[\lambda _{DA}^{{{P}_{B}}}=\lambda _{DA}^{{{P}_{C}}}=\frac{\beta _{D}^{P}}{\beta _{A}^{P}},\ \lambda _{AB}^{{{P}_{C}}}=\lambda _{AB}^{{{P}_{D}}}=\frac{\beta _{B}^{P}}{\beta _{A}^{P}},\ \lambda _{DB}^{{{P}_{A}}}=\lambda _{DB}^{{{P}_{C}}}=\frac{\beta _{B}^{P}}{\beta _{D}^{P}},\ \]
	\[\lambda _{BC}^{{{P}_{A}}}=\lambda _{BC}^{{{P}_{D}}}=\frac{\beta _{C}^{P}}{\beta _{B}^{P}},\ \lambda _{AC}^{{{P}_{B}}}=\lambda _{AC}^{{{P}_{D}}}=\frac{\beta _{C}^{P}}{\beta _{A}^{P}}.\]
\end{proof}
\hfill $\square$\par

\section{Necessary and sufficient conditions of IC-T}\label{Sec31.2}
If we select one point on each of the four faces of a given tetrahedron and connect the point to the vertex opposite the face in a straight line, respectively, are the four straight lines concurrent? The following theorem that I put forward answers this question.


\begin{theorem}{Necessary and sufficient conditions 1 of IC-T, Daiyuan Zhang}{Thm31.2.1}\label{Thm31.2.1} 
	
	Given a tetrahedron $ABCD$, take a point on the plane of its four face triangles $\triangle BCD$, $\triangle CDA$, $\triangle DAB$ and $\triangle ABC$, which are ${{P}_{A}}$, ${{P}_{B}}$, ${{P}_{C}}$, ${{P}_{D}}$ (see figure \ref{fig:tu31.2.1}), and $\overleftrightarrow{A{{P}_{A}}}$, $\overleftrightarrow{B{{P}_{B}}}$, $\overleftrightarrow{C{{P}_{C}}}$, $\overleftrightarrow{D{{P}_{D}}}$ are not parallel to each other, 
	then the necessary and sufficient conditions for the four straight lines $\overleftrightarrow{A{{P}_{A}}}$, $\overleftrightarrow{B{{P}_{B}}}$, $\overleftrightarrow{C{{P}_{C}}}$ and  $\overleftrightarrow{D{{P}_{D}}}$ concurrent at the point $P$ are as follows:
	\[\lambda _{BC}^{{{P}_{A}}}\lambda _{CD}^{{{P}_{A}}}\lambda _{DB}^{{{P}_{A}}}=1,\ \lambda _{CD}^{{{P}_{B}}}\lambda _{DA}^{{{P}_{B}}}\lambda _{AC}^{{{P}_{B}}}=1,\ \lambda _{DA}^{{{P}_{C}}}\lambda _{AB}^{{{P}_{C}}}\lambda _{BD}^{{{P}_{C}}}=1,\ \lambda _{AB}^{{{P}_{D}}}\lambda _{BC}^{{{P}_{D}}}\lambda _{CA}^{{{P}_{D}}}=1;\]
	\[\lambda _{CD}^{{{P}_{A}}}=\lambda _{CD}^{{{P}_{B}}},\ \lambda _{DB}^{{{P}_{A}}}=\lambda _{DB}^{{{P}_{C}}},\ \lambda _{BC}^{{{P}_{A}}}=\lambda _{BC}^{{{P}_{D}}},\ \lambda _{DA}^{{{P}_{B}}}=\lambda _{DA}^{{{P}_{C}}},\ \lambda _{AB}^{{{P}_{C}}}=\lambda _{AB}^{{{P}_{D}}},\ \lambda _{AC}^{{{P}_{D}}}=\lambda _{AC}^{{{P}_{B}}}.\]
	Where 
	\[{{P}_{A}}=\overleftrightarrow{B{{M}_{CD}}}\cap \overleftrightarrow{C{{M}_{BD}}}\cap \overleftrightarrow{D{{M}_{BC}}},\]
	\[{{P}_{B}}=\overleftrightarrow{C{{M}_{AD}}}\cap \overleftrightarrow{D{{M}_{CA}}}\cap \overleftrightarrow{A{{M}_{CD}}},\]
	\[{{P}_{C}}=\overleftrightarrow{D{{M}_{AB}}}\cap \overleftrightarrow{A{{M}_{BD}}}\cap \overleftrightarrow{B{{M}_{AD}}},\]
	\[{{P}_{D}}=\overleftrightarrow{A{{M}_{CD}}}\cap \overleftrightarrow{B{{M}_{CA}}}\cap \overleftrightarrow{C{{M}_{AB}}}.\]
	Where ${{M}_{AB}}\in \overleftrightarrow{AB}$, ${{M}_{AC}}\in \overleftrightarrow{AC}$, ${{M}_{AD}}\in \overleftrightarrow{AD}$, ${{M}_{BC}}\in \overleftrightarrow{BC}$,  ${{M}_{CD}}\in \overleftrightarrow{CD}$, ${{M}_{DB}}\in \overleftrightarrow{DB}$.
\end{theorem}

\begin{proof}
	\textbf{Sufficiency.} According to the inverse theorem of Ceva's theorem, for the face triangle $\triangle BCD$ (as shown in figure \ref{fig:tu31.2.1}), if $\lambda _{BC}^{{{P}_{A}}}\lambda _{CD}^{{{P}_{A}}}\lambda _{DB}^{{{P}_{A}}}=1$, then the three straight lines $\overleftrightarrow{B{{M}_{CD}}}$, $\overleftrightarrow{C{{M}_{DB}}}$ and $\overleftrightarrow{D{{M}_{BC}}}$ are concurrent, denoted as ${{P}_{A}}$, i.e. 	
	\[{{P}_{A}}=\overleftrightarrow{B{{M}_{CD}}}\cap \overleftrightarrow{C{{M}_{BD}}}\cap \overleftrightarrow{D{{M}_{BC}}}.\]
	\begin{figure}[h]
		\centering
		\includegraphics[totalheight=6cm]{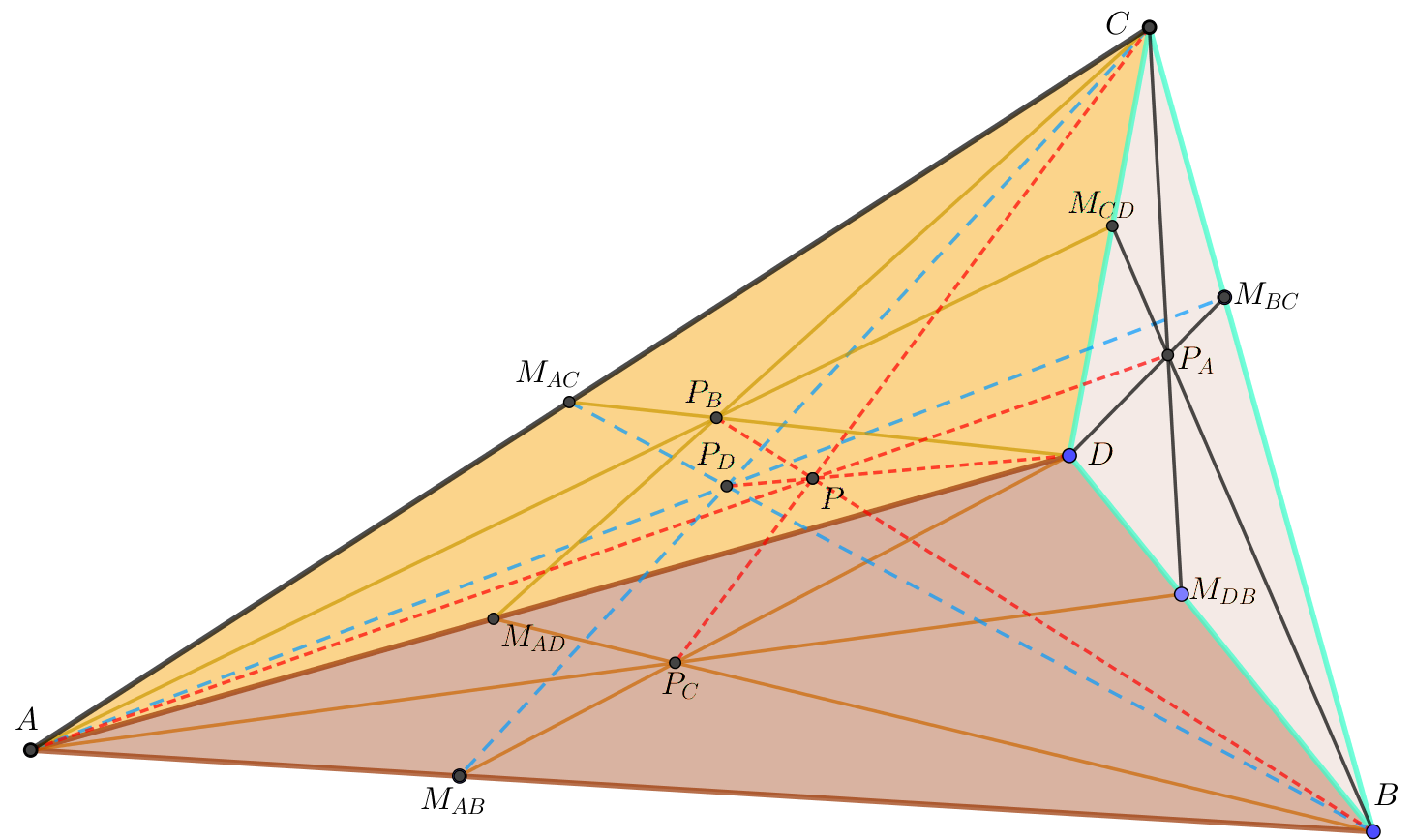}
		\caption{The geometric relationship between the IC-T and the ICs-Fs of a tetrahedron} \label{fig:tu31.2.1}
	\end{figure} 
	
	For the face triangle $\triangle CDA$ (as shown in figure \ref{fig:tu31.2.1}), if $\lambda _{CD}^{{{P}_{B}}}\lambda _{DA}^{{{P}_{B}}}\lambda _{AC}^{{{P}_{B}}}=1$, then three straight lines $\overleftrightarrow{C{{M}_{AD}}}$, $\overleftrightarrow{D{{M}_{AC}}}$ and $\overleftrightarrow{A{{M}_{CD}}}$ are concurrent, denoted as ${{P}_{B}}$, i.e. 
	\[{{P}_{B}}=\overleftrightarrow{C{{M}_{AD}}}\cap \overleftrightarrow{D{{M}_{AC}}}\cap \overleftrightarrow{A{{M}_{CD}}}.\]
	
	For the face triangle $\triangle DAB$ (as shown in figure \ref{fig:tu31.2.1}), if $\lambda _{DA}^{{{P}_{C}}}\lambda _{AB}^{{{P}_{C}}}\lambda _{BD}^{{{P}_{C}}}=1$, then three straight lines $\overleftrightarrow{D{{M}_{AB}}}$, $\overleftrightarrow{A{{M}_{DB}}}$ and $\overleftrightarrow{B{{M}_{AD}}}$ are concurrent, denoted as ${{P}_{C}}$, i.e.
	\[{{P}_{C}}=\overleftrightarrow{D{{M}_{AB}}}\cap \overleftrightarrow{A{{M}_{DB}}}\cap \overleftrightarrow{B{{M}_{AD}}}.\]
	
	For the face triangle $\triangle ABC$ (as shown in figure \ref{fig:tu31.2.1}), if $\lambda _{AB}^{{{P}_{D}}}\lambda _{BC}^{{{P}_{D}}}\lambda _{CA}^{{{P}_{D}}}=1$, then three straight lines $\overleftrightarrow{A{{M}_{BC}}}$, $\overleftrightarrow{B{{M}_{AC}}}$ and $\overleftrightarrow{C{{M}_{AB}}}$ are concurrent, denoted as ${{P}_{D}}$, i.e
	\[{{P}_{D}}=\overleftrightarrow{A{{M}_{BC}}}\cap \overleftrightarrow{B{{M}_{AC}}}\cap \overleftrightarrow{C{{M}_{AB}}}.\]
	
	The following results can be obtained from the definition of $\lambda_{CD}^{{{P}_{A}}}$ and $\lambda_{CD}^{{{P}_{B}}}$:
	\[\lambda _{CD}^{{{P}_{A}}}=\frac{\overrightarrow{CX}}{\overrightarrow{XD}},\quad X=\overleftrightarrow{B{{P}_{A}}}\cap \overleftrightarrow{CD};\]
	\[\lambda _{CD}^{{{P}_{B}}}=\frac{\overrightarrow{CY}}{\overrightarrow{YD}},\quad Y=\overleftrightarrow{A{{P}_{B}}}\cap \overleftrightarrow{CD}.\]
	
	According to section \ref{Sec2.3}, the IR is a single valued function, using the condition $\lambda _{CD}^{{{P}_{A}}}=\lambda _{CD}^{{{P}_{B}}}$ given in this theorem, we have $X=Y$, that is, point $X$ and point $Y$ coincide with each other. Let $X=Y={{M}_{CD}}$, obviously, $X=Y={{M}_{CD}}\in \overleftrightarrow{CD}$, so $\overleftrightarrow{A{{P}_{A}}}\in {{\pi }_{AB{{M}_{CD}}}}$, $\overleftrightarrow{B{{P}_{B}}}\in {{\pi }_{AB{{M}_{CD}}}}$.
	Similarly, using the condition $\lambda _{DB}^{{{P}_{A}}}=\lambda _{DB}^{{{P}_{C}}}$ given in this theorem, we have $\overleftrightarrow{A{{P}_{A}}}\in {{\pi }_{AC{{M}_{DB}}}}$, $\overleftrightarrow{C{{P}_{C}}}\in {{\pi }_{AC{{M}_{DB}}}}$; Using the condition $\lambda _{BC}^{{{P}_{A}}}=\lambda _{BC}^{{{P}_{D}}}$ given in this theorem, we have $\overleftrightarrow{A{{P}_{A}}}\in {{\pi }_{AD{{M}_{BC}}}}$, $\overleftrightarrow{D{{P}_{D}}}\in {{\pi }_{AD{{M}_{BC}}}}$. So we have	
	\begin{equation}\label{Eq31.2.1}
		\overleftrightarrow{A{{P}_{A}}}\in {{\pi }_{AB{{M}_{CD}}}}\cap {{\pi }_{AC{{M}_{DB}}}}\cap {{\pi }_{AD{{M}_{BC}}}}.
	\end{equation}
	
	Similarly, use the conditions $\lambda _{DA}^{{{P}_{B}}}=\lambda _{DA}^{{{P}_{C}}}$, $\lambda _{AC}^{{{P}_{D}}}=\lambda _{AC}^{{{P}_{B}}}$ and $\lambda _{CD}^{{{P}_{A}}}=\lambda _{CD}^{{{P}_{B}}}$ given in this theorem, we have
	\begin{equation}\label{Eq31.2.2}
		\overleftrightarrow{B{{P}_{B}}}\in {{\pi }_{BC{{M}_{AD}}}}\cap {{\pi }_{BD{{M}_{AC}}}}\cap {{\pi }_{BA{{M}_{CD}}}}.
	\end{equation}
	
	Using the conditions $\lambda _{DB}^{{{P}_{A}}}=\lambda _{DB}^{{{P}_{C}}}$, $\lambda _{DA}^{{{P}_{B}}}=\lambda _{DA}^{{{P}_{C}}}$ and $\lambda _{AB}^{{{P}_{C}}}=\lambda _{AB}^{{{P}_{D}}}$ given in this theorem, we have
	\begin{equation}\label{Eq31.2.3}
		\overleftrightarrow{C{{P}_{C}}}\in {{\pi }_{CD{{M}_{AB}}}}\cap {{\pi }_{CA{{M}_{DB}}}}\cap {{\pi }_{CB{{M}_{AD}}}}.	
	\end{equation}
	
	Using the conditions $\lambda _{BC}^{{{P}_{A}}}=\lambda _{BC}^{{{P}_{D}}}$, $\lambda _{AC}^{{{P}_{D}}}=\lambda _{AC}^{{{P}_{B}}}$ and $\lambda _{AB}^{{{P}_{C}}}=\lambda _{AB}^{{{P}_{D}}}$ given in this theorem, we have
	\begin{equation}\label{Eq31.2.4}
		\overleftrightarrow{D{{P}_{D}}}\in {{\pi }_{DA{{M}_{BC}}}}\cap {{\pi }_{DB{{M}_{AC}}}}\cap {{\pi }_{DC{{M}_{AB}}}}.
	\end{equation}
	
	According to formulas (\ref{Eq31.2.1}) and (\ref{Eq31.2.2}), $\overleftrightarrow{A{{P}_{A}}}\in {{\pi }_{AB{{M}_{CD}}}}$, $\overleftrightarrow{B{{P}_{B}}}\in {{\pi }_{BA{{M}_{CD}}}}$, so $\overleftrightarrow{A{{P}_{A}}}$ and $\overleftrightarrow{B{{P}_{B}}}$ are in the same plane ${{\pi }_{AB{{M}_{CD}}}}$. Since $\overleftrightarrow{A{{P}_{A}}}$ and $\overleftrightarrow{B{{P}_{B}}}$ are not parallel, so $\overleftrightarrow{A{{P}_{A}}}$ and $\overleftrightarrow{B{{P}_{B}}}$ must have a unique intersection. Let the two straight lines be concurrent at ${{P}_{AB}}$, i.e.
	\begin{equation}\label{Eq31.2.5}
		{{P}_{AB}}=\overleftrightarrow{A{{P}_{A}}}\cap \overleftrightarrow{B{{P}_{B}}}.
	\end{equation}
		
	According to formulas (\ref{Eq31.2.1}) and (\ref{Eq31.2.3}), $\overleftrightarrow{A{{P}_{A}}}\in {{\pi }_{AC{{M}_{DB}}}}$, $\overleftrightarrow{C{{P}_{C}}}\in {{\pi }_{CA{{M}_{DB}}}}$, so $\overleftrightarrow{A{{P}_{A}}}$ and $\overleftrightarrow{B{{P}_{B}}}$ are in the same plane ${{\pi }_{AC{{M}_{DB}}}}$. Since $\overleftrightarrow{A{{P}_{A}}}$ and $\overleftrightarrow{C{{P}_{C}}}$ are not parallel, so $\overleftrightarrow{A{{P}_{A}}}$ and $\overleftrightarrow{C{{P}_{C}}}$ must have a unique intersection. Let the two straight lines be concurrent at ${{P}_{AC}}$, i.e.
	\begin{equation}\label{Eq31.2.6}
		{{P}_{AC}}=\overleftrightarrow{A{{P}_{A}}}\cap \overleftrightarrow{C{{P}_{C}}}.	
	\end{equation}
	
	According to formulas (\ref{Eq31.2.2}) and (\ref{Eq31.2.3}), $\overleftrightarrow{B{{P}_{B}}}\in {{\pi }_{BC{{M}_{AD}}}}$,$\overleftrightarrow{C{{P}_{C}}}\in {{\pi }_{CB{{M}_{AD}}}}$, so $\overleftrightarrow{B{{P}_{B}}}$ and $\overleftrightarrow{C{{P}_{C}}}$ are in the same plane ${{\pi }_{BC{{M}_{AD}}}}$. Since $\overleftrightarrow{B{{P}_{B}}}$ and $\overleftrightarrow{C{{P}_{C}}}$ are not parallel, so $\overleftrightarrow{B{{P}_{B}}}$ and $\overleftrightarrow{C{{P}_{C}}}$ must have a unique intersection. Let the two straight lines be concurrent at ${{P}_{BC}}$, i.e.
	\[{{P}_{BC}}=\overleftrightarrow{B{{P}_{B}}}\cap \overleftrightarrow{C{{P}_{C}}}.\]
	
	Since $\overleftrightarrow{A{{P}_{A}}}$ and $\overleftrightarrow{B{{P}_{B}}}$ have a unique intersection ${{P}_{AB}}$, $\overleftrightarrow{A{{P}_{A}}}$ and $\overleftrightarrow{C{{P}_{C}}}$ also have a unique intersection ${{P}_{AC}}$, $\overleftrightarrow{B{{P}_{B}}}$ and $\overleftrightarrow{C{{P}_{C}}}$ are in the same plane ${{\pi }_{BC{{M}_{AD}}}}$, so the line $\overleftrightarrow{A{{P}_{A}}}$ and plane ${{\pi }_{BC{{M}_{AD}}}}$ cannot be parallel (Because if the line $\overleftrightarrow{A{{P}_{A}}}$ and the plane ${{\pi }_{BC{{M}_{AD}}}}$ are parallel, then $\overleftrightarrow{A{{P}_{A}}}$ and $\overleftrightarrow{B{{P}_{B}}}$ are impossible to be concurrent, which contradicts the formula (\ref{Eq31.2.5}); Similarly, $\overleftrightarrow{A{{P}_{A}}}$ and $\overleftrightarrow{C{{P}_{C}}}$ can not be concurrent, which contradicts the formula (\ref{Eq31.2.6})). It is that straight line $\overleftrightarrow{A{{P}_{A}}}$ and plane ${{\pi }_{BC{{M}_{AD}}}}$ can have one and only one intersection, so three intersections ${{P}_{AB}}$, ${{P}_{AC}}$ and ${{P}_{BC}}$ must be concurrent. Let this point be ${{P}_{1}}$, then
	\begin{equation}\label{Eq31.2.7}
		{{P}_{1}}={{P}_{AB}}={{P}_{AC}}={{P}_{BC}}.
	\end{equation}

	According to formulas (\ref{Eq31.2.2}) and (\ref{Eq31.2.4}), $\overleftrightarrow{B{{P}_{B}}}\in {{\pi }_{BD{{M}_{AC}}}}$, $\overleftrightarrow{D{{P}_{D}}}\in {{\pi }_{DB{{M}_{AC}}}}$, so $\overleftrightarrow{B{{P}_{B}}}$ and $\overleftrightarrow{C{{P}_{C}}}$ are in the same plane ${{\pi }_{BD{{M}_{AC}}}}$. Since $\overleftrightarrow{B{{P}_{B}}}$ and $\overleftrightarrow{D{{P}_{D}}}$ are not parallel, so $\overleftrightarrow{B{{P}_{B}}}$ and $\overleftrightarrow{D{{P}_{D}}}$ must have a unique intersection. Let the intersection be ${{P}_{BD}}$, i.e.
	\[{{P}_{BD}}=\overleftrightarrow{B{{P}_{B}}}\cap \overleftrightarrow{D{{P}_{D}}}.\]
		
	According to formulas (\ref{Eq31.2.3}) and (\ref{Eq31.2.4}), $\overleftrightarrow{C{{P}_{C}}}\in {{\pi }_{CD{{M}_{AB}}}}$,$\overleftrightarrow{D{{P}_{D}}}\in {{\pi }_{DC{{M}_{AB}}}}$, so $\overleftrightarrow{C{{P}_{C}}}$ and $\overleftrightarrow{D{{P}_{D}}}$ are in the same plane ${{\pi }_{CD{{M}_{AB}}}}$. Since $\overleftrightarrow{C{{P}_{C}}}$ and $\overleftrightarrow{D{{P}_{D}}}$ are not parallel, so $\overleftrightarrow{C{{P}_{C}}}$ and $\overleftrightarrow{D{{P}_{D}}}$ must have a unique intersection. Let the intersection be ${{P}_{CD}}$, i.e.
	\[{{P}_{CD}}=\overleftrightarrow{C{{P}_{C}}}\cap \overleftrightarrow{D{{P}_{D}}}.\]
	

	Since $\overleftrightarrow{B{{P}_{B}}}$ and $\overleftrightarrow{C{{P}_{C}}}$ have a unique intersection ${{P}_{BC}}$, $\overleftrightarrow{B{{P}_{B}}}$ and $\overleftrightarrow{D{{P}_{D}}}$ also have a unique intersection ${{P}_{BD}}$, $\overleftrightarrow{C{{P}_{C}}}$ and $\overleftrightarrow{D{{P}_{D}}}$ are in the same plane ${{\pi }_{CD{{M}_{AB}}}}$, so the straight line $\overleftrightarrow{B{{P}_{B}}}$ and plane ${{\pi }_{CD{{M}_{AB}}}}$ cannot be parallel (Because if the straight line $\overleftrightarrow{B{{P}_{B}}}$ and the plane ${{\pi }_{CD{{M}_{AB}}}}$ are parallel, then $\overleftrightarrow{B{{P}_{B}}}$ and $\overleftrightarrow{C{{P}_{C}}}$ are impossible to be concurrent, Similarly, $\overleftrightarrow{B{{P}_{B}}}$ and $\overleftrightarrow{D{{P}_{D}}}$ can not be concurrent). It is that the straight line $\overleftrightarrow{B{{P}_{B}}}$ and plane ${{\pi }_{CD{{M}_{AB}}}}$ can have one and only one intersection, so three intersections ${{P}_{BC}}$, ${{P}_{BD}}$ and ${{P}_{CD}}$ must be concurrent. Let this point be ${{P}_{2}}$, then
	\begin{equation}\label{Eq31.2.8}
		{{P}_{2}}={{P}_{BC}}={{P}_{BD}}={{P}_{CD}},
	\end{equation}
	
	According to formulas (\ref{Eq31.2.7}) and (\ref{Eq31.2.8}), we get ${{P}_{1}}={{P}_{2}}$. Let the point ${{P}_{1}}$ and ${{P}_{2}}$ be denoted as $P$, we have
	\[P={{P}_{1}}={{P}_{2}}={{P}_{AB}}={{P}_{AC}}={{P}_{BC}}={{P}_{BD}}={{P}_{CD}},\]
	i.e. 
	\[\begin{aligned}
		P=\overleftrightarrow{A{{P}_{A}}}\cap \overleftrightarrow{B{{P}_{B}}}=\overleftrightarrow{A{{P}_{A}}}\cap \overleftrightarrow{C{{P}_{C}}}=\overleftrightarrow{B{{P}_{B}}}\cap \overleftrightarrow{C{{P}_{C}}}=\overleftrightarrow{B{{P}_{B}}}\cap \overleftrightarrow{D{{P}_{D}}}=\overleftrightarrow{C{{P}_{C}}}\cap \overleftrightarrow{D{{P}_{D}}},  
	\end{aligned}\]
	i.e. 
	\[P\in \overleftrightarrow{A{{P}_{A}}}\cap \overleftrightarrow{B{{P}_{B}}}\cap \overleftrightarrow{C{{P}_{C}}}\cap \overleftrightarrow{D{{P}_{D}}}.\]
	
	The above formula indicates that the four straight lines  $\overleftrightarrow{A{{P}_{A}}}$, $\overleftrightarrow{B{{P}_{B}}}$, $\overleftrightarrow{C{{P}_{C}}}$ and $\overleftrightarrow{D{{P}_{D}}}$ are concurrent at the point $P$.
	
	\textbf {Necessity.} Let the IC-T be
	\[P=\overleftrightarrow{A{{P}_{A}}}\cap \overleftrightarrow{B{{P}_{B}}}\cap \overleftrightarrow{C{{P}_{C}}}\cap \overleftrightarrow{D{{P}_{D}}}.\]
	
	Select ${{P}_{A}}$, ${{P}_{B}}$, ${{P}_{C}}$, ${{P}_{D}}$ are the ICs-Fs of the corresponding triangles, respectively.  From Ceva's theorem we have:
	\[\lambda _{BC}^{{{P}_{A}}}\lambda _{CD}^{{{P}_{A}}}\lambda _{DB}^{{{P}_{A}}}=1,\ \lambda _{CD}^{{{P}_{B}}}\lambda _{DA}^{{{P}_{B}}}\lambda _{AC}^{{{P}_{B}}}=1,\ \lambda _{DA}^{{{P}_{C}}}\lambda _{AB}^{{{P}_{C}}}\lambda _{BD}^{{{P}_{C}}}=1,\ \lambda _{AB}^{{{P}_{D}}}\lambda _{BC}^{{{P}_{D}}}\lambda _{CA}^{{{P}_{D}}}=1.\]
	
	According to theorem \ref{thm:Thm31.1.1}, the following result is obtained
	\[\lambda _{CD}^{{{P}_{A}}}=\lambda _{CD}^{{{P}_{B}}}=\frac{\beta _{D}^{P}}{\beta _{C}^{P}}.\]
	
	Similarly, it can be proved that:
	\[\lambda _{DA}^{{{P}_{B}}}=\lambda _{DA}^{{{P}_{C}}},\ \lambda _{AB}^{{{P}_{C}}}=\lambda _{AB}^{{{P}_{D}}},\ \lambda _{DB}^{{{P}_{A}}}=\lambda _{DB}^{{{P}_{C}}},\ \lambda _{BC}^{{{P}_{A}}}=\lambda _{BC}^{{{P}_{D}}},\ \lambda _{AC}^{{{P}_{B}}}=\lambda _{AC}^{{{P}_{D}}}.\]
\end{proof}
\hfill $\square$\par

The above theorem shows that if the tetrahedron and the above conditions are given, there must be an IC-T:
\[P=\overleftrightarrow{A{{P}_{A}}}\cap \overleftrightarrow{B{{P}_{B}}}\cap \overleftrightarrow{C{{P}_{C}}}\cap \overleftrightarrow{D{{P}_{D}}}.\]


\begin{theorem}{Necessary and sufficient conditions 2 of IC-T, Daiyuan Zhang}{Thm31.2.2}\label{Thm31.2.2} 
	
	Given a tetrahedron $ABCD$, take a point on the plane of its four face triangles $\triangle BCD$, $\triangle CDA$, $\triangle DAB$ and $\triangle ABC$, which are ${{P}_{A}}$, ${{P}_{B}}$, ${{P}_{C}}$, ${{P}_{D}}$ (see figure \ref{fig:tu31.2.1}), and $\overleftrightarrow{A{{P}_{A}}}$, $\overleftrightarrow{B{{P}_{B}}}$, $\overleftrightarrow{C{{P}_{C}}}$, $\overleftrightarrow{D{{P}_{D}}}$ are not parallel to each other, 
	then the necessary and sufficient conditions for the four straight lines $\overleftrightarrow{A{{P}_{A}}}$, $\overleftrightarrow{B{{P}_{B}}}$, $\overleftrightarrow{C{{P}_{C}}}$ and  $\overleftrightarrow{D{{P}_{D}}}$ concurrent at the point $P$ are as follows:
	\[\frac{\alpha _{B}^{{{P}_{A}}}}{\beta _{B}^{P}}=\frac{\alpha _{C}^{{{P}_{A}}}}{\beta _{C}^{P}}=\frac{\alpha _{D}^{{{P}_{A}}}}{\beta _{D}^{P}},\,\frac{\alpha _{C}^{{{P}_{B}}}}{\beta _{C}^{P}}=\frac{\alpha _{D}^{{{P}_{B}}}}{\beta _{D}^{P}}=\frac{\alpha _{A}^{{{P}_{B}}}}{\beta _{A}^{P}},\]
	\[\frac{\alpha _{D}^{{{P}_{C}}}}{\beta _{D}^{P}}=\frac{\alpha _{A}^{{{P}_{C}}}}{\beta _{A}^{P}}=\frac{\alpha _{B}^{{{P}_{C}}}}{\beta _{B}^{P}},\,\frac{\alpha _{A}^{{{P}_{D}}}}{\beta _{A}^{P}}=\frac{\alpha _{B}^{{{P}_{D}}}}{\beta _{B}^{P}}=\frac{\alpha _{C}^{{{P}_{D}}}}{\beta _{C}^{P}}.\]
	Where $\alpha _{B}^{{{P}_{A}}}$, $\alpha _{C}^{{{P}_{A}}}$,$\alpha _{D}^{{{P}_{A}}}$ are the frame components of $\triangle BCD$; $\alpha _{C}^{{{P}_{B}}}$, $\alpha _{D}^{{{P}_{B}}}$,$\alpha _{A}^{{{P}_{B}}}$ are the frame components of $\triangle CDA$; $\alpha _{D}^{{{P}_{C}}}$, $\alpha _{A}^{{{P}_{C}}}$,$\alpha _{B}^{{{P}_{C}}}$ are the frame components of $\triangle DAB$; $\alpha _{A}^{{{P}_{D}}}$, $\alpha _{B}^{{{P}_{D}}}$,$\alpha _{C}^{{{P}_{D}}}$ are the frame components of $\triangle ABC$.
\end{theorem}

\begin{proof}
	\textbf{Sufficiency}. According to the conditions of this theorem, we can get:
	\[\frac{\alpha _{D}^{{{P}_{A}}}}{\alpha _{C}^{{{P}_{A}}}}=\frac{\alpha _{D}^{{{P}_{B}}}}{\alpha _{C}^{{{P}_{B}}}},\ \frac{\alpha _{B}^{{{P}_{A}}}}{\alpha _{D}^{{{P}_{A}}}}=\frac{\alpha _{B}^{{{P}_{C}}}}{\alpha _{D}^{{{P}_{C}}}},\ \frac{\alpha _{C}^{{{P}_{A}}}}{\alpha _{B}^{{{P}_{A}}}}=\frac{\alpha _{C}^{{{P}_{D}}}}{\alpha _{B}^{{{P}_{D}}}},\]
	\[\frac{\alpha _{A}^{{{P}_{B}}}}{\alpha _{D}^{{{P}_{B}}}}=\frac{\alpha _{A}^{{{P}_{C}}}}{\alpha _{D}^{{{P}_{C}}}},\ \frac{\alpha _{B}^{{{P}_{C}}}}{\alpha _{A}^{{{P}_{C}}}}=\frac{\alpha _{B}^{{{P}_{D}}}}{\alpha _{A}^{{{P}_{D}}}},\ \frac{\alpha _{C}^{{{P}_{D}}}}{\alpha _{A}^{{{P}_{D}}}}=\frac{\alpha _{C}^{{{P}_{B}}}}{\alpha _{A}^{{{P}_{B}}}}.\]
	
	According to theorem \ref{thm:Thm6.4.1}, we get:
	\[\lambda _{CD}^{{{P}_{A}}}=\lambda _{CD}^{{{P}_{B}}},\ \lambda _{DB}^{{{P}_{A}}}=\lambda _{DB}^{{{P}_{C}}},\ \lambda _{BC}^{{{P}_{A}}}=\lambda _{BC}^{{{P}_{D}}},\ \lambda _{DA}^{{{P}_{B}}}=\lambda _{DA}^{{{P}_{C}}},\ \lambda _{AB}^{{{P}_{C}}}=\lambda _{AB}^{{{P}_{D}}},\ \lambda _{AC}^{{{P}_{D}}}=\lambda _{AC}^{{{P}_{B}}}.\]
	
	Obviously,
	\[\lambda _{BC}^{{{P}_{A}}}\lambda _{CD}^{{{P}_{A}}}\lambda _{DB}^{{{P}_{A}}}=\frac{\alpha _{C}^{{{P}_{A}}}}{\alpha _{B}^{{{P}_{A}}}}\frac{\alpha _{D}^{{{P}_{A}}}}{\alpha _{C}^{{{P}_{A}}}}\frac{\alpha _{B}^{{{P}_{A}}}}{\alpha _{D}^{{{P}_{A}}}}=1,\ \]
	\[\lambda _{CD}^{{{P}_{B}}}\lambda _{DA}^{{{P}_{B}}}\lambda _{AC}^{{{P}_{B}}}=\frac{\alpha _{D}^{{{P}_{B}}}}{\alpha _{C}^{{{P}_{B}}}}\frac{\alpha _{A}^{{{P}_{B}}}}{\alpha _{D}^{{{P}_{B}}}}\frac{\alpha _{C}^{{{P}_{B}}}}{\alpha _{A}^{{{P}_{B}}}}=1,\ \]
	\[\lambda _{DA}^{{{P}_{C}}}\lambda _{AB}^{{{P}_{C}}}\lambda _{BD}^{{{P}_{C}}}=\frac{\alpha _{A}^{{{P}_{C}}}}{\alpha _{D}^{{{P}_{C}}}}\frac{\alpha _{B}^{{{P}_{C}}}}{\alpha _{A}^{{{P}_{C}}}}\frac{\alpha _{D}^{{{P}_{C}}}}{\alpha _{B}^{{{P}_{C}}}}=1,\]
	\[\lambda _{AB}^{{{P}_{D}}}\lambda _{BC}^{{{P}_{D}}}\lambda _{CA}^{{{P}_{D}}}=\frac{\alpha _{B}^{{{P}_{D}}}}{\alpha _{A}^{{{P}_{D}}}}\frac{\alpha _{C}^{{{P}_{D}}}}{\alpha _{B}^{{{P}_{D}}}}\frac{\alpha _{A}^{{{P}_{D}}}}{\alpha _{C}^{{{P}_{D}}}}=1.\]
		
	According to theorem \ref{thm:Thm31.2.1}, there is a point $P$ such that:
	\[P=\overleftrightarrow{A{{P}_{A}}}\cap \overleftrightarrow{B{{P}_{B}}}\cap \overleftrightarrow{C{{P}_{C}}}\cap \overleftrightarrow{D{{P}_{D}}}.\]
	Thus, the sufficiency is proved.
	
	\textbf{Necessity}. Suppose there is a point $P$ that makes
	\[P=\overleftrightarrow{A{{P}_{A}}}\cap \overleftrightarrow{B{{P}_{B}}}\cap \overleftrightarrow{C{{P}_{C}}}\cap \overleftrightarrow{D{{P}_{D}}}.\]
	
	According to theorem \ref{thm:Thm31.2.1}, we get:
	\[\lambda _{CD}^{{{P}_{A}}}=\lambda _{CD}^{{{P}_{B}}},\ \lambda _{DB}^{{{P}_{A}}}=\lambda _{DB}^{{{P}_{C}}},\ \lambda _{BC}^{{{P}_{A}}}=\lambda _{BC}^{{{P}_{D}}},\ \lambda _{DA}^{{{P}_{B}}}=\lambda _{DA}^{{{P}_{C}}},\ \lambda _{AB}^{{{P}_{C}}}=\lambda _{AB}^{{{P}_{D}}},\ \lambda _{AC}^{{{P}_{D}}}=\lambda _{AC}^{{{P}_{B}}}.\]
	
	According to theorem \ref{thm:Thm6.4.1}, we get:
	\[\frac{\alpha _{D}^{{{P}_{A}}}}{\alpha _{C}^{{{P}_{A}}}}=\frac{\alpha _{D}^{{{P}_{B}}}}{\alpha _{C}^{{{P}_{B}}}},\ \frac{\alpha _{B}^{{{P}_{A}}}}{\alpha _{D}^{{{P}_{A}}}}=\frac{\alpha _{B}^{{{P}_{C}}}}{\alpha _{D}^{{{P}_{C}}}},\ \frac{\alpha _{C}^{{{P}_{A}}}}{\alpha _{B}^{{{P}_{A}}}}=\frac{\alpha _{C}^{{{P}_{D}}}}{\alpha _{B}^{{{P}_{D}}}},\]
	\[\frac{\alpha _{A}^{{{P}_{B}}}}{\alpha _{D}^{{{P}_{B}}}}=\frac{\alpha _{A}^{{{P}_{C}}}}{\alpha _{D}^{{{P}_{C}}}},\ \frac{\alpha _{B}^{{{P}_{C}}}}{\alpha _{A}^{{{P}_{C}}}}=\frac{\alpha _{B}^{{{P}_{D}}}}{\alpha _{A}^{{{P}_{D}}}},\ \frac{\alpha _{C}^{{{P}_{D}}}}{\alpha _{A}^{{{P}_{D}}}}=\frac{\alpha _{C}^{{{P}_{B}}}}{\alpha _{A}^{{{P}_{B}}}}.\]
	So the necessity is proved.
\end{proof}
\hfill $\square$\par

\begin{theorem}{Necessary and sufficient conditions 3 of IC-T, Daiyuan Zhang}{Thm31.2.3}\label{Thm31.2.3} 
	
	Given a tetrahedron $ABCD$, take a point on the plane of its four face triangles $\triangle BCD$, $\triangle CDA$, $\triangle DAB$ and $\triangle ABC$, which are ${{P}_{A}}$, ${{P}_{B}}$, ${{P}_{C}}$, ${{P}_{D}}$ (see figure \ref{fig:tu31.2.1}), and $\overleftrightarrow{A{{P}_{A}}}$, $\overleftrightarrow{B{{P}_{B}}}$, $\overleftrightarrow{C{{P}_{C}}}$, $\overleftrightarrow{D{{P}_{D}}}$ are not parallel to each other, 
	then the necessary and sufficient conditions for the four straight lines $\overleftrightarrow{A{{P}_{A}}}$, $\overleftrightarrow{B{{P}_{B}}}$, $\overleftrightarrow{C{{P}_{C}}}$ and  $\overleftrightarrow{D{{P}_{D}}}$ concurrent at the point $P$ are as follows:
	\[\frac{\beta _{B}^{{{P}_{A}}}}{\beta _{B}^{P}}=\frac{\beta _{C}^{{{P}_{A}}}}{\beta _{C}^{P}}=\frac{\beta _{D}^{{{P}_{A}}}}{\beta _{D}^{P}},\,\frac{\beta _{C}^{{{P}_{B}}}}{\beta _{C}^{P}}=\frac{\beta _{D}^{{{P}_{B}}}}{\beta _{D}^{P}}=\frac{\beta _{A}^{{{P}_{B}}}}{\beta _{A}^{P}},\]
	\[\frac{\beta _{D}^{{{P}_{C}}}}{\beta _{D}^{P}}=\frac{\beta _{A}^{{{P}_{C}}}}{\beta _{A}^{P}}=\frac{\beta _{B}^{{{P}_{C}}}}{\beta _{B}^{P}},\,\frac{\beta _{A}^{{{P}_{D}}}}{\beta _{A}^{P}}=\frac{\beta _{B}^{{{P}_{D}}}}{\beta _{B}^{P}}=\frac{\beta _{C}^{{{P}_{D}}}}{\beta _{C}^{P}}.\]
	Where $\beta _{B}^{{{P}_{A}}}$, $\beta _{C}^{{{P}_{A}}}$, $\beta _{D}^{{{P}_{A}}}$ ($\beta _{A}^{{{P}_{A}}}=0$) are the frame components in tetrahedral frame of tetrahedron $ABCD$ at point ${{P}^{A}}$; $\beta _{C}^{{{P}_{B}}}$, $\beta _{D}^{{{P}_{B}}}$, $\beta _{A}^{{{P}_{B}}}$ ($\beta _{B}^{{{P}_{B}}}=0$) are the frame components in tetrahedral frame of tetrahedron $ABCD$ at point ${{P}^{B}}$; $\beta _{D}^{{{P}_{C}}}$, $\beta _{A}^{{{P}_{C}}}$,$\beta _{B}^{{{P}_{C}}}$ ($\beta _{C}^{{{P}_{C}}}=0$) are the frame components in tetrahedral frame of tetrahedron $ABCD$ at point ${{P}^{C}}$; $\beta _{A}^{{{P}_{D}}}$, $\beta _{B}^{{{P}_{D}}}$,$\beta _{C}^{{{P}_{D}}}$ ($\beta _{D}^{{{P}_{D}}}=0$) are the frame components in tetrahedral frame of tetrahedron $ABCD$ at point ${{P}^{D}}$.
\end{theorem}

\begin{proof}
	This theorem can be obtained directly according to theorem \ref{thm:Thm20.4.1} and theorem \ref{thm:Thm31.2.2}.
\end{proof}
\hfill $\square$\par

\begin{theorem}{Necessary and sufficient conditions 4 of IC-T, Daiyuan Zhang}{Thm31.2.4}\label{Thm31.2.4} 
	
	Given a tetrahedron $ABCD$, take one point on the plane of its four face triangles $\triangle BCD$, $\triangle CDA$, $\triangle DAB$ and $\triangle ABC$, which are ${{P}_{A}}$, ${{P}_{B}}$, ${{P}_{C}}$, ${{P}_{D}}$ (see figure \ref{fig:tu31.2.1}), and $\overleftrightarrow{A{{P}_{A}}}$, $\overleftrightarrow{B{{P}_{B}}}$, $\overleftrightarrow{C{{P}_{C}}}$, $\overleftrightarrow{D{{P}_{D}}}$ are not parallel to each other, 
	Suppose point $P$ is a given point in space, point $O$ is any point in space, and
	\[\overrightarrow{OP}=\beta _{A}^{P}\overrightarrow{OA}+\beta _{B}^{P}\overrightarrow{OB}+\beta _{C}^{P}\overrightarrow{OC}+\beta _{D}^{P}\overrightarrow{OD}\text{.}\]	
	where 
	\[\beta _{A}^{P}=\frac{{{f}^{A}}}{{{f}^{A}}+{{f}^{B}}+{{f}^{C}}+{{f}^{D}}},\]	
	\[\beta _{B}^{P}=\frac{{{f}^{B}}}{{{f}^{A}}+{{f}^{B}}+{{f}^{C}}+{{f}^{D}}},\]	
	\[\beta _{C}^{P}=\frac{{{f}^{C}}}{{{f}^{A}}+{{f}^{B}}+{{f}^{C}}+{{f}^{D}}},\]	
	\[\beta _{D}^{P}=\frac{{{f}^{D}}}{{{f}^{A}}+{{f}^{B}}+{{f}^{C}}+{{f}^{D}}}.\]	
	Then the necessary and sufficient conditions for the four straight lines $\overleftrightarrow{A{{P}_{A}}}$, $\overleftrightarrow{B{{P}_{B}}}$, $\overleftrightarrow{C{{P}_{C}}}$ and  $\overleftrightarrow{D{{P}_{D}}}$ concurrent at the point $P$ are as follows:
	\[\alpha _{B}^{{{P}_{A}}}=\frac{{{f}^{B}}}{{{f}^{B}}+{{f}^{C}}+{{f}^{D}}},\,\alpha _{C}^{{{P}_{A}}}=\frac{{{f}^{C}}}{{{f}^{B}}+{{f}^{C}}+{{f}^{D}}},\,\alpha _{D}^{{{P}_{A}}}=\frac{{{f}^{D}}}{{{f}^{B}}+{{f}^{C}}+{{f}^{D}}};\]\[\alpha _{C}^{{{P}_{B}}}=\frac{{{f}^{C}}}{{{f}^{C}}+{{f}^{D}}+{{f}^{A}}},\,\alpha _{D}^{{{P}_{B}}}=\frac{{{f}^{D}}}{{{f}^{C}}+{{f}^{D}}+{{f}^{A}}},\,\alpha _{A}^{{{P}_{B}}}=\frac{{{f}^{A}}}{{{f}^{C}}+{{f}^{D}}+{{f}^{A}}};\]\[\alpha _{D}^{{{P}_{C}}}=\frac{{{f}^{D}}}{{{f}^{D}}+{{f}^{A}}+{{f}^{B}}},\,\alpha _{A}^{{{P}_{C}}}=\frac{{{f}^{A}}}{{{f}^{D}}+{{f}^{A}}+{{f}^{B}}},\,\alpha _{B}^{{{P}_{C}}}=\frac{{{f}^{B}}}{{{f}^{D}}+{{f}^{A}}+{{f}^{B}}};\]\[\alpha _{A}^{{{P}_{D}}}=\frac{{{f}^{A}}}{{{f}^{A}}+{{f}^{B}}+{{f}^{C}}},\,\alpha _{B}^{{{P}_{D}}}=\frac{{{f}^{B}}}{{{f}^{A}}+{{f}^{B}}+{{f}^{C}}},\,\alpha _{C}^{{{P}_{D}}}=\frac{{{f}^{C}}}{{{f}^{A}}+{{f}^{B}}+{{f}^{C}}}.\]
	Where $\alpha _{B}^{{{P}_{A}}}$, $\alpha _{C}^{{{P}_{A}}}$,$\alpha _{D}^{{{P}_{A}}}$ are the frame components of $\triangle BCD$; $\alpha _{C}^{{{P}_{B}}}$, $\alpha _{D}^{{{P}_{B}}}$,$\alpha _{A}^{{{P}_{B}}}$ are the frame components of $\triangle CDA$; $\alpha _{D}^{{{P}_{C}}}$, $\alpha _{A}^{{{P}_{C}}}$,$\alpha _{B}^{{{P}_{C}}}$ are the frame components of $\triangle DAB$; $\alpha _{A}^{{{P}_{D}}}$, $\alpha _{B}^{{{P}_{D}}}$,$\alpha _{C}^{{{P}_{D}}}$ are the frame components of $\triangle ABC$. ${{f}^{A}}$, ${{f}^{A}}$,${{f}^{A}}$,${{f}^{A}}$ are a set of real numbers. The denominators in the above formulas are not zero.
\end{theorem}	
	
\begin{proof}
	\textbf{Sufficiency}. According to the conditions in this theorem, we have:
	\[\frac{\alpha _{B}^{{{P}_{A}}}}{\beta _{B}^{P}}=\frac{\alpha _{C}^{{{P}_{A}}}}{\beta _{C}^{P}}=\frac{\alpha _{D}^{{{P}_{A}}}}{\beta _{D}^{P}}=\frac{{{f}^{A}}+{{f}^{B}}+{{f}^{C}}+{{f}^{D}}}{{{f}^{B}}+{{f}^{C}}+{{f}^{D}}},\]
	\[\frac{\alpha _{C}^{{{P}_{B}}}}{\beta _{C}^{P}}=\frac{\alpha _{D}^{{{P}_{B}}}}{\beta _{D}^{P}}=\frac{\alpha _{A}^{{{P}_{B}}}}{\beta _{A}^{P}}=\frac{{{f}^{A}}+{{f}^{B}}+{{f}^{C}}+{{f}^{D}}}{{{f}^{C}}+{{f}^{D}}+{{f}^{A}}},\]
	\[\frac{\alpha _{D}^{{{P}_{C}}}}{\beta _{D}^{P}}=\frac{\alpha _{A}^{{{P}_{C}}}}{\beta _{A}^{P}}=\frac{\alpha _{B}^{{{P}_{C}}}}{\beta _{B}^{P}}=\frac{{{f}^{A}}+{{f}^{B}}+{{f}^{C}}+{{f}^{D}}}{{{f}^{D}}+{{f}^{A}}+{{f}^{B}}},\]
	\[\frac{\alpha _{A}^{{{P}_{D}}}}{\beta _{A}^{P}}=\frac{\alpha _{B}^{{{P}_{D}}}}{\beta _{B}^{P}}=\frac{\alpha _{C}^{{{P}_{D}}}}{\beta _{C}^{P}}=\frac{{{f}^{A}}+{{f}^{B}}+{{f}^{C}}+{{f}^{D}}}{{{f}^{A}}+{{f}^{B}}+{{f}^{C}}}.\]
		
	According to theorem \ref{thm:Thm31.2.2}, there is a point $P$ such that
	\[P=\overleftrightarrow{A{{P}_{A}}}\cap \overleftrightarrow{B{{P}_{B}}}\cap \overleftrightarrow{C{{P}_{C}}}\cap \overleftrightarrow{D{{P}_{D}}}.\]
	Thus, the sufficiency is proved.	
	
	\textbf{Necessity}. Suppose there is a point $P$ that makes
	\[P=\overleftrightarrow{A{{P}_{A}}}\cap \overleftrightarrow{B{{P}_{B}}}\cap \overleftrightarrow{C{{P}_{C}}}\cap \overleftrightarrow{D{{P}_{D}}}.\]
	Then according to theorem \ref{thm:Thm20.2.4}, the following results are obtained:
	\[\alpha _{B}^{{{P}_{A}}}=\frac{\beta _{B}^{P}}{1-\beta _{A}^{P}}=\frac{\frac{{{f}^{B}}}{{{f}^{A}}+{{f}^{B}}+{{f}^{C}}+{{f}^{D}}}}{1-\frac{{{f}^{A}}}{{{f}^{A}}+{{f}^{B}}+{{f}^{C}}+{{f}^{D}}}}=\frac{{{f}^{B}}}{{{f}^{B}}+{{f}^{C}}+{{f}^{D}}},\]
	\[\alpha _{C}^{{{P}_{A}}}=\frac{\beta _{C}^{P}}{1-\beta _{A}^{P}}=\frac{\frac{{{f}^{C}}}{{{f}^{A}}+{{f}^{B}}+{{f}^{C}}+{{f}^{D}}}}{1-\frac{{{f}^{A}}}{{{f}^{A}}+{{f}^{B}}+{{f}^{C}}+{{f}^{D}}}}=\frac{{{f}^{C}}}{{{f}^{B}}+{{f}^{C}}+{{f}^{D}}},\]
	\[\alpha _{D}^{{{P}_{A}}}=\frac{\beta _{D}^{P}}{1-\beta _{A}^{P}}=\frac{\frac{{{f}^{D}}}{{{f}^{A}}+{{f}^{B}}+{{f}^{C}}+{{f}^{D}}}}{1-\frac{{{f}^{A}}}{{{f}^{A}}+{{f}^{B}}+{{f}^{C}}+{{f}^{D}}}}=\frac{{{f}^{D}}}{{{f}^{B}}+{{f}^{C}}+{{f}^{D}}};\]
	\[\alpha _{C}^{{{P}_{B}}}=\frac{\beta _{C}^{P}}{1-\beta _{B}^{P}}=\frac{\frac{{{f}^{C}}}{{{f}^{A}}+{{f}^{B}}+{{f}^{C}}+{{f}^{D}}}}{1-\frac{{{f}^{B}}}{{{f}^{A}}+{{f}^{B}}+{{f}^{C}}+{{f}^{D}}}}=\frac{{{f}^{C}}}{{{f}^{C}}+{{f}^{D}}+{{f}^{A}}},\]
	\[\alpha _{D}^{{{P}_{B}}}=\frac{\beta _{D}^{P}}{1-\beta _{B}^{P}}=\frac{\frac{{{f}^{D}}}{{{f}^{A}}+{{f}^{B}}+{{f}^{C}}+{{f}^{D}}}}{1-\frac{{{f}^{B}}}{{{f}^{A}}+{{f}^{B}}+{{f}^{C}}+{{f}^{D}}}}=\frac{{{f}^{D}}}{{{f}^{C}}+{{f}^{D}}+{{f}^{A}}},\]\[\alpha _{A}^{{{P}_{B}}}=\frac{\beta _{A}^{P}}{1-\beta _{B}^{P}}=\frac{\frac{{{f}^{A}}}{{{f}^{A}}+{{f}^{B}}+{{f}^{C}}+{{f}^{D}}}}{1-\frac{{{f}^{B}}}{{{f}^{A}}+{{f}^{B}}+{{f}^{C}}+{{f}^{D}}}}=\frac{{{f}^{A}}}{{{f}^{C}}+{{f}^{D}}+{{f}^{A}}};\]
	\[\alpha _{D}^{{{P}_{C}}}=\frac{\beta _{D}^{P}}{1-\beta _{C}^{P}}=\frac{\frac{{{f}^{D}}}{{{f}^{A}}+{{f}^{B}}+{{f}^{C}}+{{f}^{D}}}}{1-\frac{{{f}^{C}}}{{{f}^{A}}+{{f}^{B}}+{{f}^{C}}+{{f}^{D}}}}=\frac{{{f}^{D}}}{{{f}^{D}}+{{f}^{A}}+{{f}^{B}}},\]
	\[\alpha _{A}^{{{P}_{C}}}=\frac{\beta _{A}^{P}}{1-\beta _{C}^{P}}=\frac{\frac{{{f}^{A}}}{{{f}^{A}}+{{f}^{B}}+{{f}^{C}}+{{f}^{D}}}}{1-\frac{{{f}^{C}}}{{{f}^{A}}+{{f}^{B}}+{{f}^{C}}+{{f}^{D}}}}=\frac{{{f}^{A}}}{{{f}^{D}}+{{f}^{A}}+{{f}^{B}}},\]
	\[\alpha _{B}^{{{P}_{C}}}=\frac{\beta _{B}^{P}}{1-\beta _{C}^{P}}=\frac{\frac{{{f}^{B}}}{{{f}^{A}}+{{f}^{B}}+{{f}^{C}}+{{f}^{D}}}}{1-\frac{{{f}^{C}}}{{{f}^{A}}+{{f}^{B}}+{{f}^{C}}+{{f}^{D}}}}=\frac{{{f}^{B}}}{{{f}^{D}}+{{f}^{A}}+{{f}^{B}}};\]
	\[\alpha _{A}^{{{P}_{D}}}=\frac{\beta _{A}^{P}}{1-\beta _{D}^{P}}=\frac{\frac{{{f}^{A}}}{{{f}^{A}}+{{f}^{B}}+{{f}^{C}}+{{f}^{D}}}}{1-\frac{{{f}^{D}}}{{{f}^{A}}+{{f}^{B}}+{{f}^{C}}+{{f}^{D}}}}=\frac{{{f}^{A}}}{{{f}^{A}}+{{f}^{B}}+{{f}^{C}}},\]
	\[\alpha _{B}^{{{P}_{D}}}=\frac{\beta _{B}^{P}}{1-\beta _{D}^{P}}=\frac{\frac{{{f}^{B}}}{{{f}^{A}}+{{f}^{B}}+{{f}^{C}}+{{f}^{D}}}}{1-\frac{{{f}^{D}}}{{{f}^{A}}+{{f}^{B}}+{{f}^{C}}+{{f}^{D}}}}=\frac{{{f}^{B}}}{{{f}^{A}}+{{f}^{B}}+{{f}^{C}}},\]
	\[\alpha _{C}^{{{P}_{D}}}=\frac{\beta _{C}^{P}}{1-\beta _{D}^{P}}=\frac{\frac{{{f}^{C}}}{{{f}^{A}}+{{f}^{B}}+{{f}^{C}}+{{f}^{D}}}}{1-\frac{{{f}^{D}}}{{{f}^{A}}+{{f}^{B}}+{{f}^{C}}+{{f}^{D}}}}=\frac{{{f}^{C}}}{{{f}^{A}}+{{f}^{B}}+{{f}^{C}}}.\]
	So the necessity is proved.
\end{proof}
\hfill $\square$\par
		

Obviously, in the above theorem, if ${{f}^{A}}$, ${{f}^{B}}$, ${{f}^{C}}$, ${{f}^{D}}$ are all functions of the six edges of a tetrahedron, then the frame components $\alpha _{B}^{{{P}_{A}}}$, $\alpha _{C}^{{{P}_{A}}}$,$\alpha _{D}^{{{P}_{A}}}$; $\alpha _{C}^{{{P}_{B}}}$, $\alpha _{D}^{{{P}_{B}}}$,$\alpha _{A}^{{{P}_{B}}}$; $\alpha _{D}^{{{P}_{C}}}$, $\alpha _{A}^{{{P}_{C}}}$,$\alpha _{B}^{{{P}_{C}}}$; $\alpha _{A}^{{{P}_{D}}}$, $\alpha _{B}^{{{P}_{D}}}$,$\alpha _{C}^{{{P}_{D}}}$ of ${{P}_{A}}$, ${{P}_{B}}$, ${{P}_{C}}$, ${{P}_{D}}$ are also functions of the lengths of the six edges of the tetrahedron.


According to theorem \ref{thm:Thm31.1.1}, the IRs of points ${{P}_{A}}$, ${{P}_{B}}$, ${{P}_{C}}$, ${{P}_{D}}$ in the corresponding triangles $\triangle BCD$, $\triangle CDA$, $\triangle DAB$, $\triangle ABC$ are, respectively:

\[\lambda _{CD}^{{{P}_{A}}}=\lambda _{CD}^{{{P}_{B}}}=\frac{{{f}^{D}}}{{{f}^{C}}},\ \lambda _{DB}^{{{P}_{A}}}=\lambda _{DB}^{{{P}_{C}}}=\frac{{{f}^{B}}}{{{f}^{D}}},\ \lambda _{BC}^{{{P}_{A}}}=\lambda _{BC}^{{{P}_{D}}}=\frac{{{f}^{C}}}{{{f}^{B}}},\]
\[\lambda _{DA}^{{{P}_{B}}}=\lambda _{DA}^{{{P}_{C}}}=\frac{{{f}^{A}}}{{{f}^{D}}},\ \lambda _{AB}^{{{P}_{C}}}=\lambda _{AB}^{{{P}_{D}}}=\frac{{{f}^{B}}}{{{f}^{A}}},\ \lambda _{AC}^{{{P}_{D}}}=\lambda _{AC}^{{{P}_{B}}}=\frac{{{f}^{C}}}{{{f}^{A}}}.\]
Therefore, each IR is also a function of the lengths of the six edges of the tetrahedron.

An application example of the above theorem is given below.


Given a tetrahedron $ABCD$, suppose point $P$ is a given point in space, point $O$ is any point in space, and
\[\overrightarrow{OP}=\beta _{A}^{P}\overrightarrow{OA}+\beta _{B}^{P}\overrightarrow{OB}+\beta _{C}^{P}\overrightarrow{OC}+\beta _{D}^{P}\overrightarrow{OD}.\]

If the frame components of tetrahedron $ABCD$ are selected in the following:
\[\beta _{A}^{P}=\frac{{{\left( {{S}^{A}} \right)}^{n}}}{{{\left( {{S}^{A}} \right)}^{n}}+{{\left( {{S}^{B}} \right)}^{n}}+{{\left( {{S}^{C}} \right)}^{n}}+{{\left( {{S}^{D}} \right)}^{n}}},\]	
\[\beta _{B}^{P}=\frac{{{\left( {{S}^{B}} \right)}^{n}}}{{{\left( {{S}^{A}} \right)}^{n}}+{{\left( {{S}^{B}} \right)}^{n}}+{{\left( {{S}^{C}} \right)}^{n}}+{{\left( {{S}^{D}} \right)}^{n}}},\]	
\[\beta _{C}^{P}=\frac{{{\left( {{S}^{C}} \right)}^{n}}}{{{\left( {{S}^{A}} \right)}^{n}}+{{\left( {{S}^{B}} \right)}^{n}}+{{\left( {{S}^{C}} \right)}^{n}}+{{\left( {{S}^{D}} \right)}^{n}}},\]	
\[\beta _{D}^{P}=\frac{{{\left( {{S}^{D}} \right)}^{n}}}{{{\left( {{S}^{A}} \right)}^{n}}+{{\left( {{S}^{B}} \right)}^{n}}+{{\left( {{S}^{C}} \right)}^{n}}+{{\left( {{S}^{D}} \right)}^{n}}}.\]

Therefore
\[\beta _{A}^{P}+\beta _{B}^{P}+\beta _{C}^{P}+\beta _{D}^{P}=1.\]	

Then, according to theorem\ref{thm:Thm18.1.4}, this set of frame components of tetrahedron $ABCD$ is unique.


Take one point ${{P}_{A}}$, ${{P}_{B}}$, ${{P}_{C}}$ and ${{P}_{D}}$ on each plane of the four face triangles $\triangle BCD$, $\triangle CDA$, $\triangle DAB$ and $\triangle ABC$ respectively (see figure \ref{fig:tu31.2.1}), let the corresponding frame components be:
\[\alpha _{B}^{{{P}_{A}}}=\frac{{{\left( {{S}^{B}} \right)}^{n}}}{{{\left( {{S}^{B}} \right)}^{n}}+{{\left( {{S}^{C}} \right)}^{n}}+{{\left( {{S}^{D}} \right)}^{n}}},\]
\[\alpha _{C}^{{{P}_{A}}}=\frac{{{\left( {{S}^{C}} \right)}^{n}}}{{{\left( {{S}^{B}} \right)}^{n}}+{{\left( {{S}^{C}} \right)}^{n}}+{{\left( {{S}^{D}} \right)}^{n}}},\]
\[\alpha _{D}^{{{P}_{A}}}=\frac{{{\left( {{S}^{D}} \right)}^{n}}}{{{\left( {{S}^{B}} \right)}^{n}}+{{\left( {{S}^{C}} \right)}^{n}}+{{\left( {{S}^{D}} \right)}^{n}}};\]
\[\alpha _{C}^{{{P}_{B}}}=\frac{{{\left( {{S}^{C}} \right)}^{n}}}{{{\left( {{S}^{C}} \right)}^{n}}+{{\left( {{S}^{D}} \right)}^{n}}+{{\left( {{S}^{A}} \right)}^{n}}},\]
\[\alpha _{D}^{{{P}_{B}}}=\frac{{{\left( {{S}^{D}} \right)}^{n}}}{{{\left( {{S}^{C}} \right)}^{n}}+{{\left( {{S}^{D}} \right)}^{n}}+{{\left( {{S}^{A}} \right)}^{n}}},\]
\[\alpha _{A}^{{{P}_{B}}}=\frac{{{\left( {{S}^{A}} \right)}^{n}}}{{{\left( {{S}^{C}} \right)}^{n}}+{{\left( {{S}^{D}} \right)}^{n}}+{{\left( {{S}^{A}} \right)}^{n}}};\]
\[\alpha _{D}^{{{P}_{C}}}=\frac{{{\left( {{S}^{D}} \right)}^{n}}}{{{\left( {{S}^{D}} \right)}^{n}}+{{\left( {{S}^{A}} \right)}^{n}}+{{\left( {{S}^{B}} \right)}^{n}}},\]
\[\alpha _{A}^{{{P}_{C}}}=\frac{{{\left( {{S}^{A}} \right)}^{n}}}{{{\left( {{S}^{D}} \right)}^{n}}+{{\left( {{S}^{A}} \right)}^{n}}+{{\left( {{S}^{B}} \right)}^{n}}},\]
\[\alpha _{B}^{{{P}_{C}}}=\frac{{{\left( {{S}^{B}} \right)}^{n}}}{{{\left( {{S}^{D}} \right)}^{n}}+{{\left( {{S}^{A}} \right)}^{n}}+{{\left( {{S}^{B}} \right)}^{n}}};\]
\[\alpha _{A}^{{{P}_{D}}}=\frac{{{\left( {{S}^{A}} \right)}^{n}}}{{{\left( {{S}^{A}} \right)}^{n}}+{{\left( {{S}^{B}} \right)}^{n}}+{{\left( {{S}^{C}} \right)}^{n}}},\]
\[\alpha _{B}^{{{P}_{D}}}=\frac{{{\left( {{S}^{B}} \right)}^{n}}}{{{\left( {{S}^{A}} \right)}^{n}}+{{\left( {{S}^{B}} \right)}^{n}}+{{\left( {{S}^{C}} \right)}^{n}}},\]
\[\alpha _{C}^{{{P}_{D}}}=\frac{{{\left( {{S}^{C}} \right)}^{n}}}{{{\left( {{S}^{A}} \right)}^{n}}+{{\left( {{S}^{B}} \right)}^{n}}+{{\left( {{S}^{C}} \right)}^{n}}}.\]
Then according to theorem \ref{thm:Thm31.2.4}, $\overleftrightarrow{A{{P}_{A}}}$, $\overleftrightarrow{B{{P}_{B}}}$, $\overleftrightarrow{C{{P}_{C}}}$ and $\overleftrightarrow{D{{P}_{D}}}$ must be concurrent at point $P$, i.e.
\[P\in \overleftrightarrow{A{{P}_{A}}}\cap \overleftrightarrow{B{{P}_{B}}}\cap \overleftrightarrow{C{{P}_{C}}}\cap \overleftrightarrow{D{{P}_{D}}}.\]

According to theorem \ref{thm:Thm31.1.1}, the IRs of points ${{P}_{A}}$, ${{P}_{B}}$, ${{P}_{C}}$, ${{P}_{D}}$ in the corresponding triangles $\triangle BCD$, $\triangle CDA$, $\triangle DAB$, $\triangle ABC$ are
\[\lambda _{BC}^{{{P}_{A}}}=\frac{{{\left( {{S}^{C}} \right)}^{n}}}{{{\left( {{S}^{B}} \right)}^{n}}},\ \lambda _{CD}^{{{P}_{A}}}=\frac{{{\left( {{S}^{D}} \right)}^{n}}}{{{\left( {{S}^{C}} \right)}^{n}}},\ \lambda _{DB}^{{{P}_{A}}}=\frac{{{\left( {{S}^{B}} \right)}^{n}}}{{{\left( {{S}^{D}} \right)}^{n}}};\]
\[\lambda _{CD}^{{{P}_{B}}}=\frac{{{\left( {{S}^{D}} \right)}^{n}}}{{{\left( {{S}^{C}} \right)}^{n}}},\ \lambda _{DA}^{{{P}_{B}}}=\frac{{{\left( {{S}^{A}} \right)}^{n}}}{{{\left( {{S}^{D}} \right)}^{n}}},\ \lambda _{AC}^{{{P}_{B}}}=\frac{{{\left( {{S}^{C}} \right)}^{n}}}{{{\left( {{S}^{A}} \right)}^{n}}};\]
\[\lambda _{DA}^{{{P}_{C}}}=\frac{{{\left( {{S}^{A}} \right)}^{n}}}{{{\left( {{S}^{D}} \right)}^{n}}},\ \lambda _{AB}^{{{P}_{C}}}=\frac{{{\left( {{S}^{B}} \right)}^{n}}}{{{\left( {{S}^{A}} \right)}^{n}}},\ \lambda _{BD}^{{{P}_{C}}}=\frac{{{\left( {{S}^{D}} \right)}^{n}}}{{{\left( {{S}^{B}} \right)}^{n}}};\]
\[\lambda _{AB}^{{{P}_{D}}}=\frac{{{\left( {{S}^{B}} \right)}^{n}}}{{{\left( {{S}^{A}} \right)}^{n}}},\ \lambda _{BC}^{{{P}_{D}}}=\frac{{{\left( {{S}^{C}} \right)}^{n}}}{{{\left( {{S}^{B}} \right)}^{n}}},\ \lambda _{CA}^{{{P}_{D}}}=\frac{{{\left( {{S}^{A}} \right)}^{n}}}{{{\left( {{S}^{C}} \right)}^{n}}}.\]


The above results show that the point $P$ is an IC-T, and the points ${{P}_{A}}$, ${{P}_{B}}$, ${{P}_{C}}$ and ${{P}_{D}}$ are all ICs-Fs corresponding to the IC-T $P$. I call this IC-T $P$ the \textbf{\bm{$n$}-th power of incenter}, which is also called the power of incenter for short.


Obviously, the number of the $n$-th power of incenters is infinite, in fact, all of the above formulas hold true as long as $n$ is a real number. When $n$ is even, the tetrahedral frame components of the $n$-th power of incenter are rational functions of the lengths of six edges. Obviously, when $n=1$, the power of incenter is the incenter, so it can be seen that the incenter is the 1-th power of incenter. When $n = 0$, the power of incenter corresponds to the centroid, or the centroid is the 0-th power of incenter. It can be seen that the concept of the power of incenter is the extension of the incenter, and the centroid is also a special power of incenter.


Since even numbers are infinite, there are infinite numbers of such points in a tetrahedron, where the square of the distance between any two points can be expressed as a rational function of the length of the six edges of the tetrahedron (see theorem \ref{thm:Thm24.2.1}).

\chapter*{Concluding remarks}

%
%

The dinner of Intercenter Geometry is coming to an end. I believe that the readers who have read this book should have realized the characteristics and originality of Intercenter Geometry.

How to evaluate the value of a new theory and method? I think it should be evaluated from the following aspects: 1. Whether there is innovation in thought, theory, method and viewpoint; 2. Whether we can get some new theoretical results; 3. Whether it can solve some problems which are difficult to be solved by traditional theories and methods; 4. Whether it has unique advantages in some specific fields; 5. Whether a kind of problems can be solved by a unified theoretical method; 6. Whether it has a good application prospect.

Intercenter Geometry has come to the field of mathematics. I hope this new born little life can flourish, spread its branches and leaves, and produce more fruitful results.

\begin{figure}[h]
	\centering
	\includegraphics[totalheight=8cm]{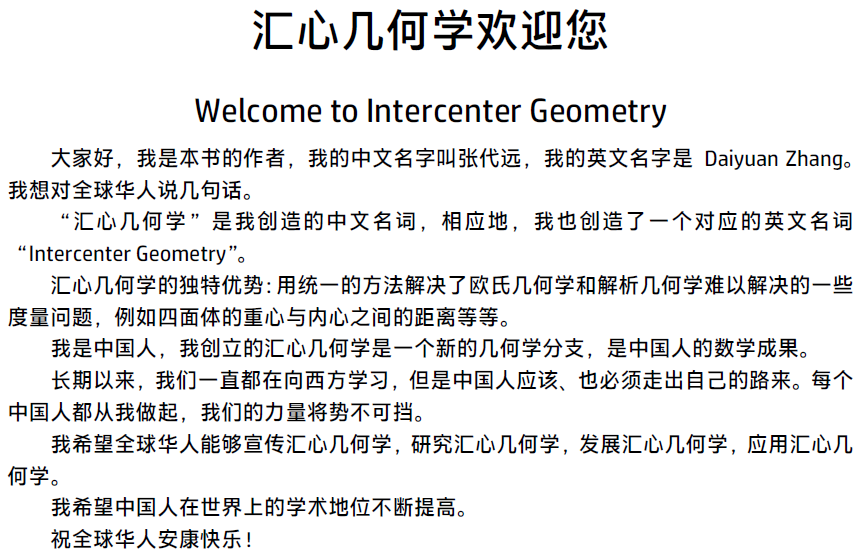}
\end{figure}



\appendix
\chapter{Preliminaries of Triangle}\label{SanjiaoxingXianguanZhishi} 
%
%
%
%


This chapter briefly introduces some of the basic knowledge needed in this book and explains some concepts which are usually confused. Readers who need a systematic understanding of these basics can consult the relevant materials.

\section{Fundamentals of vector}\label{XiangliangJibenzhishi}

The vector method is used in the research of Intercenter Geometry. This section lays the necessary foundation for the vector knowledge used in the research of Intercenter Geometry.

This section does not intend to introduce the theory and method of vector, and assumes that the reader is familiar with the basic concepts of vector, the algebraic operations of vector, the basic theorems of vector and so on.

We need to describe some confusing concepts of vector.


The following concepts need attention.

\textbf{Vector}: A quantity that has definite size and definite direction.

\textbf{Collinear vectors}: Vectors that are parallel to each other.

\textbf{Position vector}: A vector that can determine the size, direction and starting position.

The relations and differences between vector and directed line segment are as follows:

1. Relations. The directed line segment can represent the position vector determined by the \textbf{starting point} (or \textbf{initial point}) and a given direction. A vector can be expressed by a directed line segment.

2. Differences. The vector has only two elements: size and direction. The directed line segment not only has size and direction, but also has a certain starting position (starting point). That is to say, vector has only two elements: size and direction. The directed line segment has three elements: size, direction and starting point.

In this book, position vector and directed line segment represent the same concept. Vectors are denoted by bold characters, for example, $\mathbf{a}$, $\mathbf {b} $, and so on; Given two points $A$, $B$, the position vector (directed line segment) from point $A$ to point $B$ is denoted by $\overrightarrow{AB}$.

\section{Cosine theorem in the form of position vector}\label{WeizhixiangliangXingshiDeYuxiandingli}

In this book, the symbol $AB$ not only represents the side $AB$ of $\triangle ABC$ but also represents the length of the side $AB $ of $\triangle ABC$, and other symbols have similar meanings. Readers can distinguish them according to the context. The symbol $a $ denotes the length of side $BC $, i.e. $a = BC $, similarly, $b=CA$, $c=AB$.

\begin{theorem} {Cosine theorem in the form of position vector}{Thm1.2.1}\label{Thm1.2.1} 
	For a given $\triangle ABC$, we have
	\[2\overrightarrow{AB}\cdot \overrightarrow{AC}=A{{C}^{2}}+A{{B}^{2}}-B{{C}^{2}}={{b}^{2}}+{{c}^{2}}-{{a}^{2}}.\]
\end{theorem}

\begin{proof}
	In $\triangle ABC$,
	\[\overrightarrow{BC}=\overrightarrow{AC}-\overrightarrow{AB},\]
	then
	\[\begin{aligned}
		B{{C}^{2}}=\overrightarrow{BC}\cdot \overrightarrow{BC}=\left( \overrightarrow{AC}-\overrightarrow{AB} \right)\cdot \left( \overrightarrow{AC}-\overrightarrow{AB} \right) =A{{C}^{2}}+A{{B}^{2}}-2\overrightarrow{AB}\cdot \overrightarrow{AC}.  
	\end{aligned}\]
\end{proof}
\hfill $\square$\par

The above theorem states that twice the inner product of two adjacent edge vectors of a triangle equals the square of the two adjacent edges and subtracts the square of the opposite edge. The above theorem is also called cosine theorem in the form of position vector.

\section{Directed angle and vector product}\label{YouxiangjiaoYuXiangliangji}

Given a ray in a plane, the angle can be regarded as the rotation of the ray around its original point. The initial position of the ray is called the starting edge of the angle; The end position of the ray is called the ending edge of the angle; The original point of a ray is called the vertex of the angle.

Obviously, there are two directions of rotation from the beginning to the end, and a reference direction must be specified.

\textbf{Directed angle and right-handed system}. Given a plane $\alpha $, Suppose that the three points $O$, $A$ and $B$ on the plane $\alpha$ do not coincide with each other, make a unit position vector $\overrightarrow{OZ}\bot \alpha $, let $\mathbf{n}=\overrightarrow{OZ}$ ($\mathbf{n}$ is called normal vector of plane $\alpha $). The observer stands at point $Z $ and looks in the direction of $-\mathbf {n} $, if the starting edge of the position vector is $\overrightarrow {OA}$, and the ending edge of the position vector is $\overrightarrow {OB}$, then the amount of counterclockwise rotation from the starting edge $\overrightarrow {OA}$ to the ending edge $\overrightarrow {OB}$ is called the \textbf{directed angle}. Among them,  $\overrightarrow {OA}$, $\overrightarrow {OB}$ and $\mathbf {n}$ are called \textbf{right-handed system}. The right-handed system is denoted as $R\left( \overrightarrow{OA},\overrightarrow{OB},\mathbf{n} \right)$.

If a triangle is denoted as $\triangle ABC$, the order of its three vertices along the counter clockwise direction can be $A\to B\to C$; $B\to C\to A$ or $C\to A\to B$.

A directed angle is an angle that specifies the direction of rotation.

A directed angle can also be expressed by a directed line segment (position vector). In this book, the symbol $\angle AOB $ denotes a directed angle. Its starting edge is $\overrightarrow {OA} $, its ending edge is $\overrightarrow {OB} $, and its vertex is $O$.

\textbf{Vector product}. In the right-handed system of $R\left( \overrightarrow{OA},\overrightarrow{OB},\mathbf{n} \right)$, $\overrightarrow{OA}$ and $\overrightarrow{OB}$ is defined as

\begin{flalign*}
	\overrightarrow{OA}\times \overrightarrow{OB}=\left| \overrightarrow{OA} \right|\cdot \left| \overrightarrow{OB} \right|\sin \angle AOB\cdot \mathbf{n},0\le \angle AOB<2\pi.
\end{flalign*}

In the above definition, $\mathbf {n} $ is the positive reference direction of the vector product.

\section{Reflection on the formula of triangle area}\label{SanjiaoxingMianjiGongshiYinfaDeSikao}
Given the lengths of the three sides of a triangle, the area of the triangle can be derived from Helen's formula, this is well known. I have proposed a new formula to represent the area of a triangle by a 3 by 3 determinant.

\begin{theorem}{Determinant formula for triangle area, Daiyuan Zhang}{Thm1.5.1}\label{Thm1.5.1} 
	Assuming the lengths of the three sides of a given triangle are $a$, $b$, $c$, and the area is $S$, then:
	\[{{S}^{2}}=\frac{K}{16}=\frac{1}{16}\left| \begin{matrix}
		1 & 1 & 1  \\
		-{{a}^{2}} & {{b}^{2}}-{{c}^{2}} & {{a}^{2}}  \\
		{{b}^{2}} & -{{b}^{2}} & {{c}^{2}}-{{a}^{2}}  \\
	\end{matrix} \right|,\]
	or
	\[{{S}^{2}}=\frac{K}{16}=\frac{1}{16}\left| \begin{matrix}
		1 & 1 & 1  \\
		-{{b}^{2}} & {{c}^{2}}-{{a}^{2}} & {{b}^{2}}  \\
		{{c}^{2}} & -{{c}^{2}} & {{a}^{2}}-{{b}^{2}}  \\
	\end{matrix} \right|,\]
	or
	\[{{S}^{2}}=\frac{K}{16}=\frac{1}{16}\left| \begin{matrix}
		1 & 1 & 1  \\
		-{{c}^{2}} & {{a}^{2}}-{{b}^{2}} & {{c}^{2}}  \\
		{{a}^{2}} & -{{a}^{2}} & {{b}^{2}}-{{c}^{2}}  \\
	\end{matrix} \right|.\]
\end{theorem}

\begin{proof}
	According to theorem \ref{thm:YigeSanjiehanglieshiBubianliang}, $K$ is a cyclic invariant, so only the first formula needs to be proved, and the latter two formulas are also valid.  According to the formula (\ref{Eq1.4.1}):
	\[K={{\left( {{a}^{2}}+{{b}^{2}}+{{c}^{2}} \right)}^{2}}-2\left( {{a}^{4}}+{{b}^{4}}+{{c}^{4}} \right).\]
	
	According to Helen's formula, the area of the triangle is
	\[S=\sqrt{p\left( p-a \right)\left( p-b \right)\left( p-c \right)},\]
	where 
	\[p=\frac{1}{2}\left( a+b+c \right).\]
	
	Thererore 
	\[{{S}^{2}}=\frac{1}{16}H,\]
	where 
	\[H=\left( a+b+c \right)\left( b+c-a \right)\left( c+a-b \right)\left( a+b-c \right).\]
	
	And 
	\[\begin{aligned}
		H& =\left( a+b+c \right)\left( b+c-a \right)\left( c+a-b \right)\left( a+b-c \right) \\ 
		& =\left( {{\left( b+c \right)}^{2}}-{{a}^{2}} \right)\left( a-\left( b-c \right) \right)\left( a+b-c \right),  
	\end{aligned}\]
	i.e. 
	\[\begin{aligned}
		H& =\left( {{\left( b+c \right)}^{2}}-{{a}^{2}} \right)\left( {{a}^{2}}-{{\left( b-c \right)}^{2}} \right) \\ 
		& ={{\left( b+c \right)}^{2}}{{a}^{2}}-{{\left( b+c \right)}^{2}}{{\left( b-c \right)}^{2}}-{{a}^{4}}+{{a}^{2}}{{\left( b-c \right)}^{2}},  
	\end{aligned}\]
	i.e. 
	\[\begin{aligned}
		H& =\left( {{\left( b+c \right)}^{2}}+{{\left( b-c \right)}^{2}} \right){{a}^{2}}-{{\left( {{b}^{2}}-{{c}^{2}} \right)}^{2}}-{{a}^{4}} \\ 
		& =\left( 2{{b}^{2}}+2{{c}^{2}} \right){{a}^{2}}-\left( {{b}^{4}}-2{{b}^{2}}{{c}^{2}}+{{c}^{4}} \right)-{{a}^{4}} \\ 
		& ={{\left( {{a}^{2}}+{{b}^{2}}+{{c}^{2}} \right)}^{2}}-2\left( {{a}^{4}}+{{b}^{4}}+{{c}^{4}} \right).  
	\end{aligned}\]
	
	Therefore $K=H$.
\end{proof}
\hfill $\square$\par

%
Helen's formula is symmetrical, graceful and easy to remember. The determinant expression of triangle area I proposed is also symmetrical, graceful and easy to remember. 

The area formula that I proposed represents the independent variables in the formula are the lengths of the three sides of the given triangle without coordinate information. In analytic geometry, the area formula indicates that the independent variables in the formula are the coordinates of the three vertices of the triangle, as follows:
\[S=\frac{1}{2}\left| J \right|,\]
where 
\[J=\frac{1}{2}\left| \begin{matrix}
	1 & 1 & 1  \\
	{{x}_{A}} & {{x}_{B}} & {{x}_{C}}  \\
	{{y}_{A}} & {{y}_{B}} & {{y}_{C}}  \\
\end{matrix} \right|.\]

The above formula is completely different from theorem \ref{thm:Thm1.5.1}. $\left( {{x}_{A}},{{y}_{A}} \right)$, $\left( {{x}_{B}},{{y}_{B}} \right)$ and $\left( {{x}_{C}},{{y}_{C}} \right)$ in the above formula are the coordinates of the three vertices $A$, $B$ and $C$ of the triangle respectively. The coordinates of vertices lack clear geometric meaning. In particular, the formula in my theorem \ref{thm:Thm1.5.1} can deeply show the internal relationship between the area of a triangle and the lengths of three sides, which cannot be achieved by the area formula of analytical geometry (coordinate geometry). In addition, the area formula in the above analytic geometry is only applicable to the rectangular coordinate system, but not to other coordinate systems (such as polar coordinate system), that is, in the polar coordinate system, the area formula of triangle needs to be deduced again. Even in the rectangular coordinate system, due to the great human factors in the selection of coordinate origin, different people will have different selection habits, resulting in the inconsistency of coordinate values $\left( {{x}_{A}},{{y}_{A}} \right)$, $\left( {{x}_{B}},{{y}_{B}} \right)$ and $\left( {{x}_{C}},{{y}_{C}} \right)$. These are the shortcomings of analytical geometry (coordinate geometry). But in theorem \ref{thm:Thm1.5.1}, for a given triangle, the lengths of its three sides in the formula ${{S}^{2}}$ are constants, independent of coordinates, and there is no trouble in choosing the coordinate system. 

%

Readers may ask that all variables in Helen's formula are also the lengths of three sides of a triangle, which are also the geometric quantities with clear meaning, such a result can be obtained by directly using Euclidean geometry, why study Intercenter Geometry? In fact, few computational formulas of this kind (e.g., Helen formula, etc.) can be available using Euclidean geometry. This is because Euclidean geometry has no connection with algebra and does not give the calculation formula of the distance between two points, which makes it very difficult to find the calculation formulas of geometric quantities related to the distance between two points. However, the independent variables of all formulas in Plane Intercenter Geometry is the lengths of the three sides of the triangle.

The core idea of Intercenter Geometry is: the geometric quantities on the plane will be expressed by the lengths of the three sides of a given triangle; The geometric quantities in space will be expressed by the lengths of six edges of a given tetrahedron. Intercenter Geometry solves the problem of calculating the distance between two points in plane and space and the geometric quantities related to the distance between two points in a unified way by using vector methods, and many innovative achievements have been obtained.

\chapter{Exchange and Rotation}\label{DuihuanLunhuan}
\section{Exchange}\label{Duihuan}
Firstly, some concepts and symbols are introduced through the trivariate function $f=f\left( x,y,z \right)$. These concepts are also applicable to multivariate functions.


\begin{definition}{Exchange and exchange invariant}{Defi1.4.1}\label{Defi1.4.1} %
	Suppose a function is $f\left( x,y,z \right)$, the arguments are exchanged as follows: $x$ into $y$, $y$ into $x$, $z$ remains unchanged. This exchange is called the $xy$ exchange of the arguments, which is also called $xy$ \textbf {exchange}. The function after the $xy$ exchange of the arguments is denoted as $\underline{{{f}_{xy}}}\left( x,y,z \right)$, which means:
	\[\underline{{{f}_{xy}}}\left( x,y,z \right)={{f}_{xy}}\left( f\left( x,y,z \right) \right)=f\left( y,x,z \right),\]
	or 
	\[\underline{{{f}_{xy}}}\left( x,y,z \right)=f\left( y,x,z \right),\]
	then $\underline{{{f}_{xy}}}\left( x,y,z \right)$ is called the \textbf{exchange function} about $x $ and $y $ of function $f\left( x,y,z \right)$, which is also called \textbf {exchange}. If
	\[\underline{{{f}_{xy}}}\left( x,y,z \right)=f\left( y,x,z \right)=f\left( x,y,z \right),\] 
	i.e. $\underline{{{f}_{xy}}}=f$, 
	then the function $\underline{{{f}_{xy}}}\left( x,y,z \right)$ is called \textbf {function of exchange invariant} (or exchange invariant function) about $x $ and $y $ for function $f\left( x,y,z \right)$, which is also called \textbf{invariant function} for short.
\end{definition}

For example, if 
\[f\left( x,y,z \right)=z\left( x+y \right),\]
then 
\[\underline{{{f}_{xy}}}\left( x,y,z \right)=f\left( y,x,z \right)=z\left( y+x \right)=z\left( x+y \right)=f\left( x,y,z \right).\]

But 
\[\underline{{{f}_{yz}}}\left( x,y,z \right)=f\left( x,z,y \right)=y\left( z+x \right)\ne z\left( x+y \right).\]

Therefore, $\underline{{{f}_{xy}}}=f$, so the function $z\left(x+y \right)$ is a exchange invariant about $x$ and $y$; But $\underline{{{f}_{yz}}}\ne f$, that is, the function $z\left( x+y \right)$ is not a exchange invariant about $y$ and $z$.


\begin{theorem}{Algebraic properties of exchange, Daiyuan Zhang}{DuihuanDeDaishuYunsuanXingzhi}\label{DuihuanDeDaishuYunsuanXingzhi}
	Let $f=f\left( x,y,z \right)$, $g=g\left( x,y,z \right)$, ${{c}_{1}}$ and ${{c}_{2}}$ be a constant, then the algebraic operation of exchange has the following basic properties:
	\[\underline{{{\left( {{c}_{1}}f+{{c}_{2}}g \right)}_{xy}}}={{c}_{1}}\underline{{{f}_{xy}}}+{{c}_{2}}\underline{{{g}_{xy}}},\]
	\[\underline{{{\left( f\cdot g \right)}_{xy}}}=\underline{{{f}_{xy}}}\cdot \underline{{{g}_{xy}}},\]
	\[\underline{{{\left( \frac{f}{g} \right)}_{xy}}}=\frac{\underline{{{f}_{xy}}}}{\underline{{{g}_{xy}}}},\quad g\ne 0\cdot \]
\end{theorem}

%
%
%
%
%
%
%

\begin{proof}
	Using the addition formula as an example to prove. Let
	\[h\left( x,y,z \right):={{c}_{1}}f\left( x,y,z \right)+{{c}_{2}}g\left( x,y,z \right).\]
	
	In the above equation, let $x=u$ and $y=v$ to obtain
	\[h\left( u,v,z \right)={{c}_{1}}f\left( u,v,z \right)+{{c}_{2}}g\left( u,v,z \right).\]
	
	In the above equation, let $u=y$ and $v=x$ to obtain
	\[h\left( y,x,z \right)={{c}_{1}}f\left( y,x,z \right)+{{c}_{2}}g\left( y,x,z \right).\]
	
	And
	\[\underline{{{\left( h\left( x,y,z \right) \right)}_{xy}}}=h\left( y,x,z \right),\]
	\[\underline{{{\left( f\left( x,y,z \right) \right)}_{xy}}}=f\left( y,x,z \right),\]
	\[\underline{{{\left( g\left( x,y,z \right) \right)}_{xy}}}=g\left( y,x,z \right).\]
	
	i.e.
	\[\underline{{{\left( h\left( x,y,z \right) \right)}_{xy}}}={{c}_{1}}\underline{{{\left( f\left( x,y,z \right) \right)}_{xy}}}+{{c}_{2}}\underline{{{\left( g\left( x,y,z \right) \right)}_{xy}}}.\]
	
	i.e.
	\[\underline{{{\left( {{c}_{1}}f\left( x,y,z \right)+{{c}_{2}}g\left( x,y,z \right) \right)}_{xy}}}={{c}_{1}}\underline{{{\left( f\left( x,y,z \right) \right)}_{xy}}}+{{c}_{2}}\underline{{{\left( g\left( x,y,z \right) \right)}_{xy}}}.\]
	
	Using simplified notation means
	\[\underline{{{\left( {{c}_{1}}f+{{c}_{2}}g \right)}_{xy}}}={{c}_{1}}\underline{{{f}_{xy}}}+{{c}_{2}}\underline{{{g}_{xy}}}.\]
	
	Similarly, other formulas can be proven.
\end{proof}
\hfill $\square$\par

Obviously, $yz$ exchange and $zx$ exchange can have similar results.

\section{Rotation}\label{Lunhuan}
\begin{definition}{Rotation and rotation invariants}{LunhuanYuLunhuanbubianliang}\label{LunhuanYuLunhuanbubianliang}
	Assuming that the trivariate function with arguments of $x$, $y$, $z$ is $f\left( x,y,z \right)$, exchange the arguments of the function as follows: $x$ into $y$, $y$ into $z$, $z$ into $x$. This transformation is called the $xyz$ rotation transformation of the function $f\left( x,y,z \right)$, abbreviated as rotation, denoted as:
	\[\underline{\left( {\quad} \right)}:f\left( x,y,z \right)\to f\left( y,z,x \right).\]
	
	Or
	\[\underline{f\left( x,y,z \right)}=f\left( y,z,x \right).\]
	
	$\underline{f\left( x,y,z \right)}$ is also denoted as $\underline{f}\left( x,y,z \right)$ or $\underline{f}$.
	
	If $\underline{f}=f$, then function $f$ is called the invariant under rotation.
\end{definition}

%
%
%

For example, suppose a function 
\[f=f\left( x,y,z \right)=x+y+z,\]
then 
\[\underline{f}=f\left( y,z,x \right)=y+z+x=f\left( x,y,z \right)=f.\]

Thus, the function $x+y+z$ is a cyclic invariant.
The above invariant is very obvious. However, whether many functions are invariants or not can not be seen at a glance. In geometry, many invariants have definite geometric meaning. Finding invariants are an important research work in geometry.

For example, if 
\[K\left( a,b,c \right)=\left| \begin{matrix}
	1 & 1 & 1  \\
	-{{a}^{2}} & {{b}^{2}}-{{c}^{2}} & {{a}^{2}}  \\
	{{b}^{2}} & -{{b}^{2}} & {{c}^{2}}-{{a}^{2}}  \\
\end{matrix} \right|,\]	
then $K$ is the invariant under the $abc$ rotation of first-order, i.e. $\underline{K}=K$. In fact, $K$ is related to the area of a triangle, which will be proved later.

Of course, many functions are not rotation invariants.


\begin{theorem}{Algebraic operation properties of rotation, Daiyuan Zhang}{LunhuanDeDaishuYunsuanXingzhi}\label{LunhuanDeDaishuYunsuanXingzhi}
	Let $f=f\left( x,y,z \right)$, $g=g\left( x,y,z \right)$, ${{c}_{1}}$ and ${{c}_{2}}$ is a constant, and the algebraic operation of rotation has the following basic properties:
	\[\underline{\left( {{c}_{1}}f+{{c}_{2}}g \right)}={{c}_{1}}\underline{f}+{{c}_{2}}\underline{g},\]
	\[\underline{\left( f\cdot g \right)}=\underline{f}\cdot \underline{g},\]
	\[\underline{\left( \frac{f}{g} \right)}=\frac{\underline{f}}{\underline{g}},\quad g\ne 0\cdot \]
\end{theorem}

%
%
%
%
%
%
%

\begin{proof}
	Using the addition formula as an example to prove. Let
	\[h\left( x,y,z \right):={{c}_{1}}f\left( x,y,z \right)+{{c}_{2}}g\left( x,y,z \right).\]
	
	In the above equation, let $x=u$, $y=v$, $z=w$ yield
	\[h\left( u,v,w \right)={{c}_{1}}f\left( u,v,w \right)+{{c}_{2}}g\left( u,v,w \right).\]
	
	In the above equation, let $u=y$, $v=z$, $w=x$ yield
	\[h\left( y,z,x \right)={{c}_{1}}f\left( y,z,x \right)+{{c}_{2}}g\left( y,z,x \right).\]
	
	And
	\[\underline{h\left( x,y,z \right)}=h\left( y,z,x \right),\]
	\[\underline{f\left( x,y,z \right)}=f\left( y,z,x \right),\]
	\[\underline{g\left( x,y,z \right)}=g\left( y,z,x \right).\]
	
	i.e.
	\[\underline{h\left( x,y,z \right)}={{c}_{1}}\underline{f\left( x,y,z \right)}+{{c}_{2}}\underline{g\left( x,y,z \right)}.\]
	
	i.e.
	\[\underline{{{c}_{1}}f\left( x,y,z \right)+{{c}_{2}}g\left( x,y,z \right)}={{c}_{1}}\underline{f\left( x,y,z \right)}+{{c}_{2}}\underline{g\left( x,y,z \right)}.\]
	
	Using simplified notation means
	\[\underline{{{c}_{1}}f+{{c}_{2}}g}={{c}_{1}}\underline{f}+{{c}_{2}}\underline{g}.\]
	
	Similarly, other formulas can be proven.
\end{proof}
\hfill $\square$\par

\section{The relationship between rotation and exchange}\label{LunhuanYuDuihuanDeGuanxi}
\begin{theorem}{Relationship between exchange invariants and rotational invariants, Daiyuan Zhang}{DuihuanbubianliangYuLunhuanbubianliangDeGuanxi}\label{DuihuanbubianliangYuLunhuanbubianliangDeGuanxi}
	If any two different exchanges in a trivariate function are exchange invariants, the trivariate function is rotational (cyclic) invariant.
\end{theorem}

\begin{proof}
	Any two of the three variables, $x$, $y$, $z$, can be grouped together in three different ways, namely, $xy$, $yz$, $zx$. 
	
	First, suppose that the $xy$ exchange and $yz$ exchange of the function $f\left( x,y,z \right)$ are invariants, then
	\[\underline{{{f}_{xy}}}\left( x,y,z \right)=f\left( y,x,z \right)=f\left( x,y,z \right),\]
	\[\underline{{{f}_{yz}}}\left( x,y,z \right)=f\left( x,z,y \right)=f\left( x,y,z \right).\]
	
	Since $f\left( x,y,z \right)=f\left( x,z,y \right)$, then
	\[\underline{{{f}_{xy}}}\left( x,y,z \right)=\underline{{{f}_{xy}}}\left( x,z,y \right)=f\left( y,z,x \right)=\underline{f}\left( x,y,z \right).\]
	
	On the other hand, since $\underline{{{f}_{xy}}}\left( x,y,z \right)=f\left( x,y,z \right)$, then
	\[f\left( x,y,z \right)=\underline{f}\left( x,y,z \right).\]
	
	Secondly, assuming that the $yz$ exchange and $zx$ exchange of the function $f\left( x,y,z \right)$ are invariants, then
	\[\underline{{{f}_{yz}}}\left( x,y,z \right)=f\left( x,z,y \right)=f\left( x,y,z \right),\]
	\[\underline{{{f}_{zx}}}\left( x,y,z \right)=f\left( z,y,x \right)=f\left( x,y,z \right).\]
	
	Since $f\left( x,y,z \right)=f\left( z,y,x \right)$, then
	\[\underline{{{f}_{yz}}}\left( x,y,z \right)=\underline{{{f}_{yz}}}\left( z,y,x \right)=f\left( y,z,x \right)=\underline{f}\left( x,y,z \right).\]
	
	On the other hand, since $\underline{{{f}_{yz}}}\left( x,y,z \right)=f\left( x,y,z \right)$, then
	\[f\left( x,y,z \right)=\underline{f}\left( x,y,z \right).\]
	
	Finally, assuming that the $zx$ exchange and $xy$ exchange of the function $f\left( x,y,z \right)$ are invariants, then
	\[\underline{{{f}_{zx}}}\left( x,y,z \right)=f\left( z,y,x \right)=f\left( x,y,z \right),\]
	\[\underline{{{f}_{xy}}}\left( x,y,z \right)=f\left( y,x,z \right)=f\left( x,y,z \right).\]
	
	Since $f\left( x,y,z \right)=f\left( y,x,z \right)$, then
	\[\underline{{{f}_{zx}}}\left( x,y,z \right)=\underline{{{f}_{zx}}}\left( y,x,z \right)=f\left( y,z,x \right)=\underline{f}\left( x,y,z \right).\]
	
	On the other hand, since $\underline{{{f}_{zx}}}\left( x,y,z \right)=f\left( x,y,z \right)$, then
	\[f\left( x,y,z \right)=\underline{f}\left( x,y,z \right).\]
\end{proof}
\hfill $\square$\par


\begin{theorem}{Rotation decompose into pairs of exchange operations, Daiyuan Zhang}{LunhuanFenjieWeiDuihuan}\label{LunhuanFenjieWeiDuihuan}
	A rotation can be decomposed into pairs of exchange operations, i.e.:
	\[\underline{f}\left( x,y,z \right)=\underline{{{f}_{xy}}}\left( \underline{{{f}_{yz}}}\left( x,y,z \right) \right)=\underline{{{f}_{yz}}}\left( \underline{{{f}_{zx}}}\left( x,y,z \right) \right)=\underline{{{f}_{zx}}}\left( \underline{{{f}_{xy}}}\left( x,y,z \right) \right).\]
\end{theorem}

%
%
%

\begin{proof}
	For function $f\left (x, z, y \right)$, we have:
	\[{{f}_{yz}}\left( x,y,z \right)=f\left( x,z,y \right).\]
	
	Therefore 
	\[{{f}_{xy}}\left( {{f}_{yz}}\left( x,y,z \right) \right)={{f}_{xy}}\left( x,z,y \right)=f\left( y,z,x \right)=\underline{f}\left( x,y,z \right).\]
	
	The formula above states that rotation can be decomposed into the product of exchanges.
	
	Similarly, other results can be proven.
\end{proof}
\hfill $\square$\par


The above theorem can also be written as:
\[\underline{f}=\underline{{{f}_{xy}}}\circ \underline{{{f}_{yz}}}=\underline{{{f}_{yz}}}\circ \underline{{{f}_{zx}}}=\underline{{{f}_{zx}}}\circ \underline{{{f}_{xy}}}.\]

The $\circ$ in the above equation can be seen as an abstract multiplication operator.


The above theorems in this appendix are not the results obtained using the theory of Intercenter Geometry, and the content is also quite basic. It's just that I haven't found any relevant discussions, so writing my name after the theorem is a bit timid. If readers can confirm that it is not appropriate, please let me know and I will make the necessary modifications. But I think the following third-order determinant invariants are my innovation. Here I would like to reiterate that almost all the theorems in this book that involve Intercenter Geometry are my creations, and I can confidently write my name after the theorems.


\begin{theorem}{An invariant of 3 by 3 determinant, Daiyuan Zhang}{YigeSanjiehanglieshiBubianliang}\label{YigeSanjiehanglieshiBubianliang}
	Let $a$, $b$, $c$ be arbitrary real numbers, respectively, and the function $K$ be
	\[K\left( a,b,c \right)=\left| \begin{matrix}
		1 & 1 & 1  \\
		-{{a}^{2}} & {{b}^{2}}-{{c}^{2}} & {{a}^{2}}  \\
		{{b}^{2}} & -{{b}^{2}} & {{c}^{2}}-{{a}^{2}}  \\
	\end{matrix} \right|,\]	
	then 
	${{K}_{ab}}={{K}_{bc}}={{K}_{ca}}=\underline{K}=K.$
	That is, $K$ is the invariant under both exchange and rotation.
\end{theorem}

%
%

\begin{proof}
	Expand the determinant to:
	\[\begin{aligned}
		K& =\left( {{b}^{2}}-{{c}^{2}} \right)\left( {{c}^{2}}-{{a}^{2}} \right)+{{a}^{2}}{{b}^{2}}+{{a}^{2}}{{b}^{2}} \\ 
		& -{{b}^{2}}\left( {{b}^{2}}-{{c}^{2}} \right)+{{a}^{2}}\left( {{c}^{2}}-{{a}^{2}} \right)+{{a}^{2}}{{b}^{2}},  
	\end{aligned}\]
	i.e. 
	\[\begin{aligned}
		K&={{b}^{2}}{{c}^{2}}-{{b}^{2}}{{a}^{2}}-{{c}^{4}}+{{c}^{2}}{{a}^{2}}+{{a}^{2}}{{b}^{2}}+{{a}^{2}}{{b}^{2}} \\ 
		& -{{b}^{4}}+{{b}^{2}}{{c}^{2}}+{{a}^{2}}{{c}^{2}}-{{a}^{4}}+{{a}^{2}}{{b}^{2}}\text{,}  
	\end{aligned}\]
	i.e. 
	\begin{equation}\label{Eq1.4.1}
		\begin{aligned}
			K& =2\left( {{a}^{2}}{{b}^{2}}+{{b}^{2}}{{c}^{2}}+{{c}^{2}}{{a}^{2}} \right)-{{a}^{4}}-{{b}^{4}}-{{c}^{4}} \\ 
			& ={{\left( {{a}^{2}}+{{b}^{2}}+{{c}^{2}} \right)}^{2}}-2\left( {{a}^{4}}+{{b}^{4}}+{{c}^{4}} \right).  
		\end{aligned}
	\end{equation}
	
	Obviously, ${{a}^{2}}+{{b}^{2}}+{{c}^{2}}$ and ${{a}^{4}}+{{b}^{4}}+{{c}^{4}}$ are all the invariants under both exchange and rotation. According to theorem  \ref{DuihuanDeDaishuYunsuanXingzhi} and theorem \ref{thm:LunhuanDeDaishuYunsuanXingzhi}, $K$ is also invariant under both exchange and rotation.
\end{proof}
\hfill $\square$\par

\chapter{Symbols and Explanations}\label{FuhaoJiqiHanyiShuoming}

\section{Notations}\label{ZhuyaoGongshiFuhaoHanyiJiqiShuoming}
A brief explanation of the main notations in this book is provided. Please refer to the book for details. 

\begin{longtable}{rp{12cm}}
	\caption{\textbf{A brief table of the meanings of the main formula symbols}} \\
	\hline \hline \hline
	\textbf{Notations} & \textbf{Brief Explanations}\\ 
	\hline 
	\endfirsthead
	\multicolumn{2}{c}%
	{\bfseries \tablename\ \thetable{: Brief Explanations} -- Continued from previous page} \\ 
	\hline \hline
	\textbf{Symbols} & \textbf{Brief explanations}\\ 
	\hline 
	\endhead
	\hline \multicolumn{2}{r}{{Continue to the next page}} 
	\endfoot
	\hline \hline 
	\endlastfoot
	
	$AB$ & The line segment connecting point $A$ and point $B$, also represents the length of the line segment $AB$ \\ \hline
	
	
	
	
	$\overleftrightarrow{AB}$ & A straight line passing through points $A$ and $B$ \\ \hline
	
	$\overline{OP}$ & The ray passing through points $O$ and $P$, the starting point is $O$. \\ \hline
	
	
	$\overrightarrow{OP}$ & The position vector of starting point is $O$ and the ending point is $P$ \\ \hline
	
	$a$, $b$, $c$ &  The three sides or side lengths of a triangle ($\triangle{ABC}$) \\ \hline 
	
	$AB$, $BC$, $CA$ & The three sides or side lengths of a triangle ($\triangle{ABC}$), $a:=BC$, $b:=CA$, $c:=AB$. \\ \hline 

	$\alpha_{A}^{P}$, $\alpha_{B}^{P}$, $\alpha_{C}^{P}$ & Frame components under triangular frame (or line segment frame)\\ \hline 
	
	$\mathcal{R}\left( AB \right)$ & Edge-Axis Cartesian coordinate system with origin $A$ and horizontal axis $AB$ \\ \hline 
	
	$\mathcal{R}\left( BC \right)$ & Edge-Axis Cartesian coordinate system with origin $B$ and horizontal axis $BC$ \\ \hline 
	
	$\mathcal{R}\left( CA \right)$ & Edge-Axis Cartesian coordinate system with origin $C$ and horizontal axis $CA$ \\ \hline 
	
	$\mathcal{P}\left( AB \right)$ & Edge-axis polar coordinate system with origin $A$ and polar axis $AB$  \\ \hline
	
	$\mathcal{P}\left( BC \right)$ & Edge-axis polar coordinate system with origin $B$ and polar axis $BC$ \\ \hline
	
	$\mathcal{P}\left( CA \right)$ & Edge-axis polar coordinate system with origin $C$ and polar axis $CA$ \\ \hline
	
	$\left(x_{\mathcal{R}\left( AB \right)}^{P},y_{\mathcal{R}\left( AB \right)}^{P} \right)$ & The Cartesian coordinates of point $P$ under the coordinate system $\mathcal{R}\left( AB \right)$, also denoted as $P\left(x_{\mathcal{R}\left( AB \right)}^{P},y_{\mathcal{R}\left( AB \right)}^{P} \right)$ \\ \hline 
		
	$\left( x_{\mathcal{R}\left( BC \right)}^{P},y_{\mathcal{R}\left( BC \right)}^{P} \right)$ & The Cartesian coordinates of point $P$ under the coordinate system $\mathcal{R}\left( BC \right)$, also denoted as $P\left( x_{\mathcal{R}\left( BC \right)}^{P},y_{\mathcal{R}\left( BC \right)}^{P} \right)$ \\ \hline
	
	$\left( x_{\mathcal{R}\left( CA \right)}^{P},y_{\mathcal{R}\left( CA \right)}^{P} \right)$ & The Cartesian coordinates of point $P$ under the coordinate system $\mathcal{R}\left( CA \right)$, also denoted as $P\left( x_{\mathcal{R}\left( CA \right)}^{P},y_{\mathcal{R}\left( CA \right)}^{P} \right)$ \\ \hline 
	
	$\left( {{\rho }_{AP}},{{\theta }_{\mathcal{R}\left( AB \right)}} \right)$ & The polar coordinates of point $P$ associated with $\mathcal{R}\left( AB \right)$, also abbreviated as $\left( {{\rho }_{AP}},{{\theta }_{A}} \right)$ \\ \hline
	
	$\left( {{\rho }_{BP}},{{\theta }_{\mathcal{R}\left( BC \right)}} \right)$ & The polar coordinates of point $P$ associated with $\mathcal{R}\left( BC \right)$, also abbreviated as $\left( {{\rho }_{BP}},{{\theta }_{B}} \right)$ \\ \hline
	
	$\left( {{\rho }_{CP}},{{\theta }_{\mathcal{R}\left( CA \right)}} \right)$ & The polar coordinates of point $P$ associated with $\mathcal{R}\left( CA \right)$, also abbreviated as $\left( {{\rho }_{CP}},{{\theta }_{C}} \right)$ \\ \hline

	$\mathcal{T}_{A}$, $\mathcal{T}_{B}$, $\mathcal{T}_{C}$ & The Edge-axis coordinate systems for vertex $A$, $B$ and $C$ of $\triangle ABC$ \\ \hline 
	
	$\mathcal{T}_{ABC}$ & The set of Edge-axis coordinate system of $\triangle ABC$ , abbreviated as the set of Edge-axis coordinate system. $\mathcal{T}_{ABC}:={{\mathcal{T}_{A}} \cup{} {\mathcal{T}_{B}} \cup{} {\mathcal{T}_{C}}}$.\\ \hline 
	
	$\beta_{A}^{P}$, $\beta_{B}^{P}$, $\beta_{C}^{P}$, $\beta_{D}^{P}$ & The frame components of point $P$ under the tetrahedral frame \\ \hline

%
%
%
%
\end{longtable}

%
%
%
%
%
%
%
%
%
%
%
%
%
%
%
%
%
%
%
%
%
%
%
%
%
\end{document}